\documentclass[12pt,openany,twoside]{book}

\usepackage{amssymb}
\usepackage{amsmath}
\usepackage{amsthm}
\usepackage{amsfonts}
\usepackage{mathrsfs}
\usepackage{makeidx}
\usepackage{color}

\usepackage{makeidx}     
\usepackage{multicol}    

\allowdisplaybreaks

\makeindex

\numberwithin{equation}{section}

\newtheorem{theorem}{Theorem}[section]
\newtheorem{lemma}[theorem]{Lemma}
\newtheorem{corollary}[theorem]{Corollary}
\newtheorem{proposition}[theorem]{Proposition}
\newtheorem{example}[theorem]{Example}
\theoremstyle{definition}
\newtheorem{remark}[theorem]{Remark}
\newtheorem{definition}[theorem]{Definition}
\newtheorem{assumption}[theorem]{Assumption}

\def\red{\color{red}}

\def\fz{\infty}

\def\rp{\mathbb{R}_{+}}
\def\rr{\mathbb{R}}
\def\rrn{\mathbb{R}^{n}}
\def\rrnp{\mathbb{R}_{+}^{n+1}}

\def\zp{\mathbb{Z}_{+}}

\def\eps{\varepsilon}

\def\lims{\varlimsup}
\def\limi{\varliminf}
\def\lf{\left}
\def\r{\right}
\def\gfz{\genfrac{}{}{0pt}{}}

\def\m0{m_{0}}
\def\mi{m_{\infty}}
\def\M0{M_{0}}
\def\MI{M_{\infty}}
\def\Mw{\max\{\M0(\omega),\MI(\omega)\}}
\def\mw{\min\{\m0(\omega),\mi(\omega)\}}

\def\mwr{\min\left\{\m0(\omega^{r}),
\mi(\omega^{r})\right\}}

\def\Mwr{\max\left\{\M0(\omega^{r}),
\MI(\omega^{r})\right\}}

\def\mwms{\min\left\{\m0(1/\omega^{s}),
\mi(1/\omega^{s})\right\}}
\def\Mwms{\max\left\{\M0(1/\omega^{s}),
\MI(1/\omega^{s})\right\}}

\def\1bf{\mathbf{1}}
\def\0bf{\mathbf{0}}

\def\Kmc{\mathcal{K}}
\def\mc{\mathcal{M}}
\def\Kmp{\dot{\mathcal{K}}}
\def\Msc{\mathscr{M}}

\def\NHerzSea{\mathcal{K}_{1/\omega}^{p',q'}(\rrn)}

\def\HerzS{\Kmp_{\omega}^{p,q}(\rrn)}
\def\HerzScr{\Kmp_{\omega^{r}}^{p/r,q/r}(\rrn)}
\def\HerzScs{\Kmp_{\omega^{s}}^{p/s,q/s}(\rrn)}
\def\HerzScso{\Kmp_{\omega^{s_0}}^{p/s_0,q/s_0}(\rrn)}

\def\HerzSo{\Kmp_{\omega,\0bf}^{p,q}(\rrn)}
\def\HaSaHo{H\HerzSo}
\def\HaSaH{H\HerzS}

\def\HerzSoa{\Kmp_{1/\omega,\0bf}^{p',q'}(\rrn)}

\def\HerzSea{\Kmp_{1/\omega}^{p',q'}(\rrn)}
\def\HerzSocr{\Kmp_{\omega^{r},\0bf}^{p/r,q/r}(\rrn)}

\def\HerzSocs{\Kmp_{\omega^{s},\0bf}^{p/s,q/s}(\rrn)}
\def\HerzSocr{\Kmp_{\omega^{r},\0bf}^{p/r,q/r}(\rrn)}

\def\HaSaHol{h\HerzSo}
\def\HaSaHl{h\HerzS}
\def\WHerzS{W\HerzS}
\def\WHerzSo{W\HerzSo}
\def\WHaSaH{W\HaSaH}
\def\WHaSaHo{W\HaSaHo}

\def\NHerzS{\mathcal{K}_{\omega}^{p,q}(\rrn)}
\def\NHerzSo{\mathcal{K}_{\omega,\0bf}^{p,q}(\rrn)}
\def\NHaSaHo{H\NHerzSo}
\def\NHaSaH{H\NHerzS}
\def\NHaSaHl{h\NHerzS}
\def\NHaSaHol{h\NHerzSo}
\def\WNHaSaHo{WH\NHerzSo}
\def\WNHaSaH{WH\NHerzS}

\def\NHerzSocs{\Kmc_{\omega^{s},\0bf}^{p/s,q/s}(\rrn)}
\def\NHerzScs{\Kmc_{\omega^{s}}^{p/s,q/s}(\rrn)}

\def\aHaSaHo{H\Kmp_{\omega,\0bf}^{p,q,r,d,s}(\rrn)}
\def\aHaSaH{H\Kmp_{\omega}^{p,q,r,d,s}(\rrn)}
\def\iaHaSaHo{H\Kmp_{\omega,\0bf,\mathrm{fin}}^{p,q,r,d,s}(\rrn)}
\def\iaHaSaH{{H\Kmp_{\omega,\mathrm{fin}}^{p,q,r,d,s}(\rrn)}}
\def\HaSaHoc{\dot{\mathcal{L}}_{\omega,\0bf}^{p,q,r,d,s}(\rrn)}
\def\HaSaHod{\dot{\mathcal{L}}_{\omega,\0bf}^{p,q,r',d,s}(\rrn)}
\def\NaHaSaHo{H\Kmc_{\omega,\0bf}^{p,q,r,d,s}(\rrn)}
\def\NaHaSaH{H\Kmc_{\omega}^{p,q,r,d,s}(\rrn)}
\def\NiaHaSaHo{H\Kmc_{\omega,\0bf,\mathrm{fin}}^{p,q,r,d,s}(\rrn)}
\def\NiaHaSaH{{H\Kmc_{\omega,\mathrm{fin}}^{p,q,r,d,s}(\rrn)}}
\def\NHaSaHoc{\mathcal{L}_{\omega,\0bf}^{p,q,r,d,s}(\rrn)}
\def\NHaSaHod{\mathcal{L}_{\omega,\0bf}^{p,q,r',d,s}(\rrn)}

\def\HerzSx{\Kmp_{\omega,\xi}^{p,q}(\rrn)}
\def\HerzSxl{\Kmp_{\omega,\xi_l}^{p,q}(\rrn)}
\def\NHerzSx{\Kmc_{\omega,\xi}^{p,q}(\rrn)}

\def\HerzSxd{\Kmp_{1/\omega,\xi}^{p',q'}(\rrn)}

\def\NHerzSxd{\Kmc_{1/\omega,\xi}^{p',q'}(\rrn)}

\def\aHaSaHol{h\Kmp_{\omega,\0bf}^{p,q,r,d,s}(\rrn)}
\def\aHaSaHl{h\Kmp_{\omega}^{p,q,r,d,s}(\rrn)}
\def\NaHaSaHol{h\Kmc_{\omega,\0bf}^{p,q,r,d,s}(\rrn)}
\def\NaHaSaHl{h\Kmc_{\omega}^{p,q,r,d,s}(\rrn)}

\def\MorrSo{\textbf{\textsl{M}}_{\omega,\0bf}^{p,q}(\rrn)}
\def\MorrS{\textbf{\textsl{M}}_{\omega}^{p,q}(\rrn)}
\def\aMorrSo{\textbf{\textsl{M}}_{\omega,\0bf}^{p,q,r,d,s}(\rrn)}
\def\aMorrS{\textbf{\textsl{M}}_{\omega}^{p,q,r,d,s}(\rrn)}

\def\Bspace{\dot{\mathcal{B}}_{1/\omega}^{p',q'}(\rrn)}
\def\Bspaced{(\dot{\mathcal{B}}_{\omega}^{p,q}(\rrn))^{*}}
\def\Bspaceds{(\dot{\mathcal{B}}_{\omega}^{p,q}(\rrn))^{*}}
\def\bspace{\dot{\mathcal{B}}_{\omega}^{p,q}(\rrn)}

\def\NBspace{\mathcal{B}_{1/\omega}^{p',q'}(\rrn)}

\def\Nbspace{\mathcal{B}_{\omega}^{p,q}(\rrn)}

\def\BMO{{\rm BMO}(\rrn)}
\def\CMO{{\rm CMO}(\rrn)}

\def\mvr{\left\{\sum_{j\in\mathbb{N}}
\left[\mc(f_{j})\right]^{r}\right\}^{\frac{1}{r}}}

\def\vr{\left(\sum_{j\in\mathbb{N}}|f_{j}|^{r}\right)^{\frac{1}{r}}}

\def\vry{\left[\sum_{j\in\mathbb{N}}|f_{j}(y)|^{r}\right]^{\frac{1}{r}}}
\def\mvrge{\left\{\sum_{j\in\mathbb{N}}
\left[\mc(g\eta_{j})\right]^{r}\right\}^{\frac{1}{r}}}
\def\mvrgex{\left\{\sum_{j\in\mathbb{N}}
\left[\mc(g\eta_{j})(x)\right]^{r}\right\}^{\frac{1}{r}}}
\def\vrge{\left(\sum_{j\in\mathbb{N}}|
g\eta_{j}|^{r}\right)^{\frac{1}{r}}}
\def\vrey{\left[\sum_{j\in\mathbb{N}}|
\eta_{j}(y)|^{r}\right]^{\frac{1}{r}}}
\def\mvrys{\{\sum_{j\in\mathbb{N}}
[\mc(f_{j})(y)]^{r}\}^{\frac{1}{r}}}
\def\mvro{\left\{\sum_{j\in\mathbb{N}}
\left[\mc(f_{j,k,1})\right]^{r}\right\}^{\frac{1}{r}}}
\def\mvryo{\{\sum_{j\in\mathbb{N}}
[\mc(f_{j,k,1})(y)]^{r}\}^{\frac{1}{r}}}

\def\mvrytw{\{\sum_{j\in\mathbb{N}}
[\mc(f_{j,k,2})(y)]^{r}\}^{\frac{1}{r}}}
\def\mvrth{\left\{\sum_{j\in\mathbb{N}}
\left[\mc(f_{j,k,3})\right]^{r}\right\}^{\frac{1}{r}}}
\def\mvryth{\{\sum_{j\in\mathbb{N}}[\mc(f_{j,k,3})(y)]^{r}\}^{\frac{1}{r}}}
\def\vrys{[\sum_{j\in\mathbb{N}}|f_{j}(y)|^{r}]^{\frac{1}{r}}}

\def\vrtw{\left(\sum_{j\in\mathbb{N}}|f_{j,k,2}|^{r}\right)^{\frac{1}{r}}}

\def\xil{\xi_{l}}

\def\tk{2^{k}}
\def\tka{2^{k+1}}
\def\tkm{2^{k-1}}
\def\tkmd{2^{k-2}}

\def\zmn{\mathbb{Z}\setminus\mathbb{N}}

\newcommand{\normmm}[1]
{\left\|#1\right\|^{\star}}

\def\supp{\mathop\mathrm{\,supp\,}}

\def\HerzScsod{\dot{\mathcal{B}}
_{1/\omega^{s_0}}^{(p/s_0)',(q/s_0)'}(\rrn)}
\def\HerzScsd{\dot{\mathcal{B}}
_{1/\omega^{s}}^{(p/s)',(q/s)'}(\rrn)}

\begin{document}

\title{\vspace{-10cm}
\vspace{5cm}\bf\huge Real-Variable Theory of
Hardy Spaces Associated with
Generalized Herz Spaces of Rafeiro and Samko}
\author{Yinqin Li, Dachun Yang and Long Huang}
\date{\red March 26, 2022}
\maketitle

\noindent Yinqin Li

\medskip

\noindent Laboratory of Mathematics and Complex Systems (Ministry of
Education), School of Mathematical Sciences,
Beijing Normal University, Beijing 100875,
People's Republic of China

\smallskip

\noindent{\it E-mail address:} \texttt{yinqinli@mail.bnu.edu.cn}

\bigskip

\noindent Dachun Yang (Corresponding author)

\medskip

\noindent Laboratory of Mathematics and Complex Systems (Ministry of
Education), School of Mathematical Sciences,
Beijing Normal University, Beijing 100875,
People's Republic of China

\smallskip

\noindent{\it E-mail address:} \texttt{dcyang@bnu.edu.cn}

\bigskip

\noindent Long Huang

\medskip

\noindent School of Mathematics and Information Science,
Key Laboratory of Mathematics and Interdisciplinary Sciences of the Guangdong
Higher Education Institute, Guangzhou University,
Guangzhou, 510006,
People's Republic of China

\smallskip

\noindent{\it E-mail address:} \texttt{longhuang@gzhu.edu.cn}

\vspace{5cm}


\medskip


\frontmatter

\markboth{\scriptsize\rm\sc Preface}{\scriptsize\rm\sc Preface}

\chapter*{Preface\label{s0}}
\addcontentsline{toc}{chapter}{Preface}

It is well known that Herz spaces certainly play an
important role in harmonic analysis and
partial differential equations,
and have been systematically studied and
developed so far; see,
for instance, \cite{f87,LiY,ly98a,KY,YH}
for classical Herz spaces,
\cite{in20,S13,RS17,wl20,wf21,wx21,yl21}
for variable Herz spaces,
\cite{dh19,ds16,GLY,MRZ18,xy15,xy16}
for Herz-type Hardy spaces,
\cite{d13,d19,dh17,xy051,x09,x14}
for Herz-type Besov spaces,
and \cite{dd19,xy031,xy03,xy05,xy18} for
Herz-type Triebel--Lizorkin spaces.
Observe that the classical
Herz space was originally introduced by Herz
\cite{CSH} in 1968 to study the
Bernstein theorem on absolutely
convergent Fourier transforms, while
the research on Herz spaces
can be traced back to the work
of Beurling \cite{Beu}. Indeed,
in 1964, to study some convolution algebras,
Beurling {\cite{Beu}} first introduced
a special Herz space $A^{p}(\rrn)$,
with $p\in (1,\infty)$
[see Remark \ref{ighr}(ii) for its definition],
which is also called Beurling algebra.
After that,
great progress has been made on
Herz spaces and their applications.
For instance, in 1985, Baernstein
and Sawyer \cite{bs85} generalized
these Herz spaces and gave many applications;
Feichtinger \cite{f87}
introduced another norm of $A^{p}(\rrn)$,
which is equivalent to the
norm defined by Beurling \cite{Beu}.

Moreover, Herz spaces
play a crucial role in the convergence
and the summability problems of
Fourier transforms. Recall that the
study of summability means was
originally motivated by the famous
convergence problem of Dirichlet integral operators.
As one of the deepest results in harmonic analysis,
in the celebrated works of both
Carleson \cite{c66} and Hunt \cite{h67},
they showed that Dirichlet integral
operators converge almost everywhere
in one-dimensional case. In recent decades,
via replacing
Dirichlet integral operators by
some other summability means,
the summability of Fourier transforms was
systematically studied by
Butzer and Nessel \cite{bune},
Trigub and Belinsky \cite{Trigub2004}, and Feichtinger
and Weisz \cite{fw06jga,fw06-1,fw06-2} as well as
Weisz \cite{w08,w09,w11,w12,w14,w15,w16,w17}. Particularly,
let $\theta$ be an integrable function.
In the recent book \cite{wbook},
Weisz showed that the $\theta$-means of
some functions converge to
themselves at all their Lebesgue
points if and only if the Fourier transform
of $\theta$ belongs to a
suitable Herz space. This means that
those Herz spaces are the best choice
in the study of the summability of
Fourier transforms. Furthermore, Herz
spaces also prove important in the recent article
of Sawano et al. \cite{SHYY}.
Indeed, in 2017, Sawano et al. \cite{SHYY} introduced
the ball quasi-Banach
function space $X$ and the associated
Hardy space $H_X(\rrn)$
via the grand maximal functions.
As was pointed out in \cite{SHYY},
Sawano et al. used a certain inhomogeneous
Herz space to overcome the difficulty
appearing in the proof of the
convergence of the atomic decomposition of $H_X(\rrn)$.
Also, Herz-type spaces prove useful
in partial differential equations.
For instance, Scapellato \cite{s19} showed that
the variable Herz spaces are
the key tools in the study of the regularity of
solutions to elliptic equations.
In addition, the Fourier--Herz space is one of
the most suitable spaces to investigate
the global stability for fractional
Navier--Stokes equations;
see, for instance, \cite{lyhh21,cs18,nz18}.
For more progress on applications of Herz spaces,
we refer the reader to
\cite{c18,fp18,HoKP20E,HoKP21S,hjmaa,hwyy21,MRZ19,NRZ,ny18,RS15,yz18}.

On the other hand, as a good substitute of Herz spaces,
Herz-type Hardy spaces
are also useful in many mathematical fields
such as harmonic analysis,
partial differential equations,
and geometric analysis,
which have been systematically studied and
developed so far; see, for instance,
\cite{CL,GCJ,d18,lywl20,LuY08,hwy96}.
Recall that the classical real Hardy space
$H^p(\rrn)$ was originally
initiated by Stein and Weiss \cite{SW} and then
systematically developed
by Fefferman and Stein \cite{FS72}.
For more developments of
Hardy spaces, we refer the reader to
\cite{HoKP12,HoKP13,lu,Ky,s93,NSa12,NSa14,w92,w94}.
In addition, as a variant of the
classical real Hardy spaces, the Hardy
spaces associated with the
Beurling algebras on the real line
were first introduced by
Chen and Lau \cite{CL} in 1989,
in which they studied the
dual spaces and the maximal function
characterizations of these
Herz--Hardy spaces. Later,
Garc\'{i}a-Cuerva \cite{GCJ}
generalized the results of
Chen and Lau \cite{CL} to
higher-dimensional case, and
Garc\'{i}a-Cuerva and Herrero \cite{GH} further studied
the maximal function, atomic,
and Littlewood--Paley function
characterizations of these Herz--Hardy spaces.
In 1995, Lu and Yang \cite{LuY95}
introduced Herz--Hardy spaces
with general indices and
established their atomic and molecular
characterizations.

Meanwhile, Lu and Yang \cite{ly95}
proved that some oscillatory
singular integral operators are bounded
from Herz--Hardy spaces to Herz spaces, and
showed that their boundedness
fails on Herz--Hardy spaces
via a counterexample.
Recall that, as was pointed out by Pan \cite{p91},
the oscillatory singular integral
operators may not bounded
from classical Hardy space
$H^1(\rrn)$ to Lebesgue space $L^1(\rrn)$.
Thus, in some sense, the results
of \cite{ly95} showed that
the Herz--Hardy space is a
proper substitution of $H^1(\rrn)$
in the study on oscillatory singular
integral operators.
In addition, the commutators and the
multiplier theorems on the Herz--Hardy spaces
were investigated by Lu and Yang, respectively,
in \cite{ly97} and \cite{ly98};
both the interpolation of generalized Herz spaces
and its applications were given
in Hern\'{a}ndez and Yang \cite{hy98,HY}. Furthermore,
Hern\'{a}ndez et al. \cite{hwy96} established
the $\phi$-transform and the wavelet
characterizations for some Herz and
Herz-type Hardy spaces by means of a
local version of the discrete tent spaces at the origin.
For more progress on the Herz--Hardy spaces
and their applications, we refer the reader to
\cite{hly99,ly96,WF,FW,LYH,IS10,HoKP18,rs21}.

Nowadays, more and more new function
spaces continually
spring up to meet the increasing requirements
arising in harmonic analysis and
partial differential equations; see, for instance,
\cite{ad12,mb03,d22,nrz21,HSamko,RS,yyh13}.
Particularly, extending the classical Herz spaces to
some more general settings also has
attracted considerable attention recently.
For instance, Rafeiro and Samko
\cite{HSamko} creatively introduced
local and global
generalized Herz spaces recently,
which are the generalization of
classical homogeneous Herz spaces and
connect with generalized Morrey type spaces.
To be precise, Rafeiro and Samko \cite{HSamko}
showed the scale of these Herz spaces include
Morrey type spaces and complementary
Morrey type spaces and, as
applications, they also obtained
the boundedness of a class of
sublinear operators on these
generalized Herz spaces. Note that
Morrey type spaces have been studied in
\cite{as17-na,as17-jfa,noss21,s14,su20}.
Moreover, as a generalization of
classical variable Herz spaces
introdeced by Izuki \cite{IM},
Rafeiro and Samko \cite{RS}
further extended the generalized Herz spaces
to the variable exponent setting. Furthermore,
in \cite{RS}, Rafeiro and Samko showed that
the generalized variable exponent
Herz spaces coincide with
the generalized variable exponent
Morrey type spaces, and also established
some boundedness of sublinear operators on these spaces.

Observe that Herz-type Hardy spaces have been
found lots of applications in many branches
of mathematics, and
that both the real-variable theory of
function spaces and their applications
are always one of the central topic of
harmonic analysis.
Then it is a natural and meaningful
topic to introduce and develop
the real-variable theory of Hardy spaces
associated with local and global
generalized Herz spaces of Rafeiro and
Samko \cite{HSamko}, which is just
the subject of the present book.
To achieve this, we first investigate
some basic properties of these
generalized Herz spaces and realize that
the Hardy spaces associated
with generalized Herz spaces completely
fall into the framework of
Hardy spaces associated with
ball quasi-Banach function spaces
which were first studied by Sawano et al. \cite{SHYY}
as is mentioned above.
Let $X$ be a ball quasi-Banach function space on $\rrn$
(see \cite{SHYY} or Definition \ref{Df1} below).
Recall that, if the powered
Hardy--Littlewood maximal operators are bounded on
both $X$ and its associate space,
and if $X$ supports a Fefferman--Stein
vector-valued inequality, then
Sawano et al. \cite{SHYY} established
various real-variable characterizations
of the Hardy space
$$H_X(\rrn)$$
associated with $X$, respectively,
in terms of atoms, molecules, and the Lusin area function.
After the work of Sawano et al. \cite{SHYY},
the real-variable theory of
function spaces associated with
ball quasi-Banach function spaces
is well developed by Ho in \cite{ho16,h21} and also
by others in
\cite{CWYZ,YYY,HCY,HYYZ,WYY,SHYY,WYYZ,ZWYY}.

This book is devoted to
exploring further properties of
the local and the global generalized Herz spaces
and establishing a complete
real-variable theory of Hardy spaces associated with
local and global generalized Herz
spaces via a totally
fresh perspective which means that we view these
generalized Herz spaces as special
cases of ball quasi-Banach
function spaces.
In this perspective, the real-variable
theory of Hardy spaces
associated with local generalized Herz spaces can
be deduced directly from the general framework of the
real-variable theory of $H_X(\rrn)$ because
the local generalized Herz spaces satisfy
all the assumptions of the results
about Hardy-type spaces
associated with ball quasi-Banach function spaces.
However, due to the deficiency of the
associate space of the global
generalized Herz space, the known
real-variable characterizations about Hardy-type spaces
associated with ball quasi-Banach function spaces
are not applicable to Hardy spaces associated with
global generalized Herz spaces. Therefore,
the study of the real-variable
theory of Hardy spaces associated with
global generalized Herz spaces is more difficult.
To overcome this obstacle,
via replacing the assumptions on associate
spaces by some weaker assumptions about
the integral representations of quasi-norms
of ball quasi-Banach function
spaces [see Theorem \ref{Atogx}(ii) below],
we first develop some new real-variable
characterizations of Hardy-type spaces
associated with ball quasi-Banach function spaces,
which even improve the known results on Hardy-type spaces
associated with ball quasi-Banach function spaces
and, moreover, can have more additional
anticipating applications.
In particular, applying these
improved conclusions, we further
obtain a complete real-variable theory of Hardy spaces
associated with global generalized Herz spaces.

Precisely, in this book, we first study
some basic properties of these generalized Herz spaces
and obtain boundedness and
compactness characterizations of commutators on them.
Then, based on these local and global
generalized Herz spaces,
we introduce associated Herz--Hardy spaces,
localized Herz--Hardy spaces, and weak Herz--Hardy spaces,
and develop a complete real-variable theory of
these Herz--Hardy spaces,
including various maximal function, atomic,
finite atomic, molecular as well as
various Littlewood--Paley function characterizations.
As applications, we establish the boundedness of
some important operators arising
from harmonic analysis on these Herz--Hardy spaces.
Finally, the inhomogeneous Herz--Hardy spaces and their
complete real-variable theory are also investigated.
We should point out that, with the fresh perspective and the improved
conclusions on the real-variable theory of Hardy
spaces associated with ball quasi-Banach function spaces,
the exponents in all the obtained conclusions of this book
are \emph{sharp}. Moreover, all of these results in the book
are \emph{new} and have never been published before.

To be precise, this book is organized as follows.

In Chapter \ref{sec3}, we first recall the
concepts of the function class $M(\rp)$ and
the Matuszewska--Orlicz indices, and
the definitions of the local
generalized Herz space $\HerzSo$ as well as the
global generalized Herz space $\HerzS$
introduced by Rafeiro and Samko \cite{HSamko},
where $p,\ q\in(0,\fz)$, $\omega\in M(\rp)$,
and $\0bf$ denotes the \emph{origin} of $\rrn$.
In addition, under some reasonable and
sharp assumptions,
we show that these generalized Herz spaces are
special ball quasi-Banach function spaces.
Then we establish some basic properties
about these generalized Herz spaces,
which include their convexity, the absolute continuity
of the quasi-norm of $\HerzSo$,
the boundedness criterion of sublinear operators
on them which was essentially obtained
by Rafeiro and Samko \cite[Theorem 4.3]{HSamko},
and Fefferman--Stein vector-valued inequalities.
Furthermore, we find the dual space
and the associate space of
$\HerzSo$ when $p$, $q\in(1,\infty)$.
Finally, we establish the extrapolation
theorems of local and
global generalized Herz spaces.

In Chapter \ref{sec4}, we
first introduce the block space
$$\bspace$$
based on the concepts of $(\omega,\,p)$-blocks.
Then we investigate the properties
of block spaces in two aspects.
On the one hand, by establishing an equivalent
characterization of block spaces
via local generalized Herz spaces $\HerzSo$
and borrowing some ideas from the proof of
\cite[Theorem 6.1]{gm13},
we show that
the global generalized
Herz space $\HerzSea$ is
just the dual space of the block space
$\bspace$ when $p$, $q\in(1,\infty)$.
This dual theorem plays an essential role
in the study of
the real-variable theory
of Hardy spaces
associated with global
generalized Herz spaces
in the subsequent chapters.
On the other hand, we establish
the boundedness of some sublinear
operators on block spaces.
In particular, the boundedness
of powered Hardy--Littlewood
maximal operators on block spaces
is obtained, which also
plays an important role
in the subsequent chapters.

The main target of Chapter \ref{sec5}
is to study the boundedness and the
compactness characterizations of commutators
on generalized Herz spaces.
Recall that Tao et al. \cite{TYYZ} established
the boundedness and the
compactness characterizations of commutators
on ball Banach function spaces. Combining
these and the fact mentioned above
that the generalized Herz spaces are special
ball Banach function spaces,
we first show that local generalized Herz spaces
satisfy all the assumptions
of the results obtained in \cite{TYYZ},
and then obtain the boundedness and
the compactness characterizations
of commutators on the local generalized
Herz space $\HerzSo$.
However, conclusions obtained in \cite{TYYZ}
are not applicable to show
the boundedness and the compactness
characterizations of commutators on
$\HerzS$ due to the deficiency of the
associate space of the global
generalized Herz space $\HerzS$.
Notice that the most important
usage of associate spaces in the proof
of \cite{TYYZ} is that,
under the assumption of the boundedness
of the Hardy--Littlewood maximal operator
on both ball Banach function spaces
and their associate spaces,
Tao et al.\ obtained the extrapolation theorem
of ball Banach function spaces.
Via this extrapolation theorem and
some other technical lemmas independent
of associate spaces,
they then established the
boundedness and the compactness
characterizations of commutators
on ball Banach function spaces in
\cite{TYYZ}. Therefore,
to overcome the difficulty caused
by the deficiency of associate spaces
of $\HerzS$,
we establish new boundedness
and compactness characterizations
of commutators on ball Banach function
spaces under the assumption that
the extrapolation theorem
holds true for ball Banach
function spaces instead of the
assumption about associate spaces,
which improves the corresponding results of \cite{TYYZ}.
Finally, applying these improved boundedness
and compactness characterizations
of commutators on ball Banach function
spaces and the extrapolation theorem
of $\HerzS$ obtained in Chapter \ref{sec3}, we then obtain
the boundedness and the compactness
characterizations of commutators
on global generalized Herz spaces.

Chapter \ref{sec6} is devoted
to introducing the generalized Herz--Hardy
space and then
establishing its complete
real-variable theory.
To be precise, we first introduce
the generalized Herz--Hardy spaces,
$$\HaSaHo\quad\mathrm{and}\quad\HaSaH,$$
associated, respectively, with the local generalized
Herz space $\HerzSo$ and the global
generalized Herz space $\HerzS$.
Then, using the known
real-variable characterizations of
Hardy spaces associated with ball quasi-Banach
function spaces, we establish
various maximal function,
(finite) atomic, molecular, and
Littlewood--Paley function characterizations
of the Herz--Hardy space $\HaSaHo$.
Moreover, the duality and the Fourier
transform properties of $\HaSaHo$ are
also obtained based on the corresponding
results about Hardy spaces associated with
ball quasi-Banach function spaces.
However, the study of $\HaSaH$ is more
difficult than that of $\HaSaHo$ due to
the deficiency of associate spaces
of global generalized Herz spaces.
To overcome this obstacle,
via replacing the assumptions
of the boundedness
of powered Hardy--Littlewood maximal operators
on associate spaces
(see Assumption \ref{assas} below)
used in \cite[Theorems 3.6,\ 3,7,\ and 3.9]{SHYY}
by some weaker assumptions about
the integral representations of quasi-norms
of ball quasi-Banach function spaces
as well as some boundedness of
powered Hardy--Littlewood maximal operators
[see both (ii) and (iii)
of Theorem \ref{Atogx} below],
we establish new (finite) atomic and
molecular characterizations
of the Hardy space $H_X(\rrn)$
associated with the ball
quasi-Banach function space $X$,
which improve the corresponding
results obtained by Sawano et al.\
\cite{SHYY}.
Using these improved characterizations
and making full use of the obtained duality between
block spaces and global generalized
Herz spaces in Chapter \ref{sec4}
as well as the construction of
the quasi-norm $\|\cdot\|_{\HerzS}$,
we then obtain the maximal function,
(finite) atomic,
molecular, the various Littlewood--Paley
function characterizations of $\HaSaH$
and also give some properties about Fourier transforms
of distributions in $\HaSaH$.
Finally, as applications,
via first showing two boundedness criteria
of Calder\'{o}n--Zygmund operators
on Hardy spaces associated with ball
quasi-Banach function spaces,
we establish the boundedness
of Calder\'{o}n--Zygmund
operators on generalized Herz--Hardy spaces.

In Chapter \ref{sec7},
we first introduce the localized
generalized Herz--Hardy space
and then establish its complete
real-variable theory. To achieve
this, we begin with showing its various
maximal function characterizations
via the known maximal function characterizations
of the local Hardy space $h_X(\rrn)$
associated with the ball quasi-Banach function
space $X$. Then, via establishing new
atomic and molecular characterizations of $h_X(\rrn)$
as well as the boundedness of
pseudo-differential operators on $h_X(\rrn)$
without assumption about associate spaces (see Theorems
\ref{atomxxl}, \ref{molexl},
and \ref{pseudox} below), we
show the atomic and molecular characterizations
of localized generalized Herz--Hardy spaces
and the boundedness of pseudo-differential
operators on localized generalized Herz--Hardy spaces.
In addition, to clarify the relation
between localized generalized Herz--Hardy spaces
and generalized Herz--Hardy spaces, we first find the
relation between $h_X(\rrn)$ and
Hardy spaces $H_X(\rrn)$ associated with
the ball quasi-Banach function space $X$.
This extends the results obtained
by Goldberg \cite[Lemma 4]{Gold} for
classical Hardy spaces and also Nakai and Sawano
\cite[Lemma 9.1]{NSa12} for variable
Hardy spaces. Applying this and some auxiliary lemmas about
generalized Herz spaces, we then
obtain the relation between
localized generalized Herz--Hardy spaces
and generalized Herz--Hardy spaces.
As applications, we also establish various
Littlewood--Paley function
characterizations of $h_X(\rrn)$ and hence, together
with the construction of the
quasi-norm $\|\cdot\|_{\HerzS}$, the Littlewood--Paley function
characterizations of localized generalized
Herz--Hardy spaces are obtained.

The main target of Chapter \ref{sec8}
is to introduce weak
generalized Herz--Hardy spaces and
establish their complete
real-variable theory. For this purpose,
recall that Zhang et al.\ \cite{ZWYY}
and Wang et al.\ \cite{WYYZ}
investigated the real-variable
theory of the weak Hardy space $WH_X(\rrn)$
associated with the
ball quasi-Banach function space $X$.
Via removing the assumption about associate spaces,
we establish new atomic and molecular
characterizations of $WH_X(\rrn)$ as well
as new real interpolation between the Hardy space $H_X(\rrn)$
associated with the ball quasi-Banach
function space $X$ and the Lebesgue
space $L^{\infty}(\rrn)$, which improve the corresponding
results obtained in \cite{ZWYY} and \cite{WYYZ}.
Then, using these improved real-variable characterizations of $WH_X(\rrn)$,
we obtain various maximal function,
atomic, and molecular characterizations
of weak generalized Herz--Hardy spaces
and also show that the real interpolation spaces
between generalized Herz--Hardy spaces and
the Lebesgue space $L^{\infty}(\rrn)$
are just the new introduced weak generalized Herz--Hardy spaces.
In addition, by establishing a technique lemma about the
quasi-norm $\|\cdot\|_{W\HerzS}$ and the Littlewood--Paley
function characterizations of $WH_X(\rrn)$
obtained in \cite{WYYZ}, we show various Littlewood--Paley
function characterizations of weak generalized
Herz--Hardy spaces. Furthermore, we establish
two boundedness criteria of Calder\'{o}n--Zygmund operators
from the Hardy space $H_X(\rrn)$
to the weak Hardy space $WH_X(\rrn)$ and,
as a consequence, we finally deduce the boundedness
of Calder\'{o}n--Zygmund operators
from generalized Herz--Hardy spaces
to weak generalized Herz--Hardy spaces
even in the critical case.

In Chapter \ref{sec9}, we first introduce the inhomogeneous
generalized Herz spaces and then establish their corresponding conclusions
obtained in Chapters \ref{sec3} through \ref{sec5}.

Furthermore, in Chapter \ref{sec10}, based on the inhomogeneous
generalized Herz spaces studied in Chapter \ref{sec9},
we introduce the inhomogeneous generalized
Herz--Hardy spaces, the inhomogeneous
localized generalized Herz--Hardy spaces,
and the inhomogeneous weak generalized Herz--Hardy
spaces. Then we establish their various
real-variable characterizations and also
give some applications,
which are the corresponding inhomogeneous variants
obtained, respectively, in
Chapters \ref{sec6} through \ref{sec8}.

Throughout this book, we always let $\mathbb{N}:=
\{1,2,\ldots\}$\index{$\mathbb{N}$},
$\mathbb{Z}_{+}:=\mathbb{N}\cup\{0\}$\index{$\mathbb{Z}_+$},
$\rr_+:=(0,\fz)$\index{$\rp$},
and $$\rrnp:=\{(x,t):\ x\in\rrn,\
t\in(0,\infty)\}\index{$\rrnp$}.$$
We also use $\mathbf{0}:=(0,\ldots,0)$
to denote the \emph{origin} of $\rrn$.
For any $x:=(x_{1},\dots,x_{n})\in\rrn$ and
$\theta:=(\theta_{1},\dots,\theta_{n})\in(\mathbb{Z}_+)^n
=:\mathbb{Z}_{+}^{n}$,
let $|\theta|:=\theta_{1}+\cdots+\theta_{n}$,
$x^{\theta}:=x_{1}^{\theta_{1}}\cdots x_{n}^{\theta_{n}}$,
and $\partial^\gamma:=(\frac {\partial}{\partial x_1})^{\gamma_1}
\cdots(\frac {\partial}{\partial x_n})^{\gamma_n}$\index{$\partial^\gamma$}.
We always denote by $C$ a \emph{positive constant}
which is independent of the main parameters,
but it may vary from line to line.
We use $C_{(\alpha,\beta,\cdots)}$
to denote a positive constant depending
on the indicated parameters $\alpha$, $\beta$,
$\ldots$. The notation $f\lesssim g$
means $f\leq Cg$ and, if $f
\lesssim g\lesssim f$, then we write $
f\sim g$. If $f\leq Cg$ and $g=h$ or $g\leq h$,
we then write $f\lesssim g\sim h$
or $f\lesssim g\lesssim h$,
rather than $f\lesssim g=h$ or
$f\lesssim g\leq h$. For any $s\in\rr$,
the symbol
$\lceil s\rceil$\index{$\lceil s\rceil$}
denotes the smallest integer not less than $s$,
and the symbol
$\lfloor s\rfloor$\index{$\lfloor s\rfloor$} denotes
the largest integer not greater than $s$.
For any set $E\subset\rrn$, we denote
the set $\rrn\setminus E$ by
$E^\complement$\index{$E^\complement$},
its \emph{characteristic function} by $\1bf_{E}$,
and its \emph{$n$-dimensional Lebesgue measure} by $|E|$.
For any $q\in[1,\infty]$, we denote by $q'$ its \emph
{conjugate exponent}, namely, $1/q+1/q'=1$.
In addition, we use $\mathbb{S}^{n-1}:=\{x\in\rrn:\
|x|=1\}$\index{$\mathbb{S}^{n-1}$}
to denote the \emph{unit sphere} in $\rrn$
and $d\sigma$\index{$d\sigma$}
the \emph{area measure} on $\mathbb{S}^{n-1}$.
Furthermore, we always use the \emph{symbol}
$\Msc(\rrn)$\index{$\Msc(\rrn)$}
to denote the
set of all measurable functions on $\rrn$.
The symbol $\mathcal{Q}$\index{$\mathcal{Q}$}
denotes the set of all cubes
with edges parallel to the coordinate axes.
Finally, for any cube $Q\in\mathcal{Q}$,
$rQ$\index{$rQ$} means a cube with the same
center as $Q$ and $r$ times the edge length
of $Q$.\index{$|E|$}

This research of the authors is supported by the National
Key Research and Development Program of China
(Grant No. 2020YFA0712900)
and the National Natural Science Foundation of China
(Grant Nos. 11971058 and 12071197).
Long Huang is also supported by Guangdong Basic and
Applied Basic Research Foundation (Grant No. 2021A1515110905).
Yinqin Li would like to express
his deep thanks to Dr. Hongchao Jia,
Dr. Yangyang Zhang, and Dr. Yirui Zhao
for some helpful discussions on the subject of
this book.

\vspace{1cm}

\noindent Beijing, China
\hfill Yinqin Li

\noindent Beijing, China
\hfill Dachun Yang

\noindent Guangzhou, China
\hfill Long Huang

\smallskip

\noindent March 2022

\markboth{\scriptsize\rm\sc Contents}{\scriptsize\rm\sc Contents}
\tableofcontents

\mainmatter

\chapter{Generalized Herz Spaces of Rafeiro and Samko}\label{sec3}
\markboth{\scriptsize\rm\sc
Generalized Herz Spaces of Rafeiro and Samko}{\scriptsize\rm\sc
Generalized Herz Spaces of Rafeiro and Samko}

In this chapter, we first recall the concepts
of both the function class $M(\rp)$ and
the Matuszewska--Orlicz
indices, and the definitions of the local
generalized Herz space
$$\HerzSo$$
as well as the
global generalized Herz space
$$\HerzS$$
creatively
introduced by Rafeiro and Samko \cite{HSamko},
where $p,\ q\in(0,\fz)$ and $\omega\in M(\rp)$.
Next, we recall some basic concepts about
ball quasi-Banach function spaces
introduced by Sawano et al.\ \cite{SHYY}.
Moreover, under some sharp assumptions,
we show that these generalized Herz spaces are
special ball quasi-Banach function spaces.
Then we establish some basic properties about these
generalized Herz spaces, which include their convexity,
the absolutely continuity of the
quasi-norm of $\HerzSo$,
the boundedness criterion of
sublinear operators on them
which was essentially
obtained by Rafeiro and Samko
\cite[Theorem 4.3]{HSamko},
and Fefferman--Stein vector-valued inequalities.
We should point out that the
global generalized Herz space $\HerzS$
may not have the absolutely continuous quasi-norm.
Indeed, using the quasi-norm of
$\Kmp_{\omega}^{p,q}(\rr)$ and borrowing some ideas
from \cite[Example 5.1]{ST}, we
construct a special set $E$ and show that
its characteristic function $\mathbf{1}_E$
belongs to
certain global generalized Herz space.
But, $\mathbf{1}_E$ does not
have an absolutely continuous quasi-norm
in this global generalized Herz space.
Finally, by introducing the local
generalized Herz space $\HerzSx$
for any given $\xi\in\rrn$ and
establishing a dual result of it when $p$,
$q\in(1,\infty)$,
we find the associate space of $\HerzSo$
under some reasonable and sharp assumptions.
As an application, we establish the extrapolation
theorem of both $\HerzSo$ and $\HerzS$ at the
end of this chapter.

\section{Matuszewska--Orlicz Indices}

In this section, we first recall the concept
of a function class $M(\rp)$
given in \cite{HSamko} and the concept of
Matuszewska--Orlicz indices
originally introduced by Matuszewska and Orlicz in \cite{MO,mo65}.
Then we present some fundamental properties related
to $M(\rp)$ and Matuszewska--Orlicz indices, which are
widely used throughout this book.

To begin with, let $\omega$ be a nonnegative function on $\rp$.
Then the function $\omega$ is said to be
\emph{almost increasing}\index{almost increasing}
(resp.,\ \emph{almost decreasing}\index{almost decreasing})
on $\rp$ if there exists
a constant $C\in[1,\infty)$ such that,
for any $t,\ \tau\in(0,\infty)$ satisfying
$t\leq\tau$ (resp., $t\geq\tau$),
$$\omega(t)\leq C\omega(\tau)$$
(see, for instance, \cite[p.\,30]{KMRS}).
Now, we recall the concept of the function class $M(\rp)$
given in \cite[Definition 2.1]{HSamko} as follows.

\begin{definition}\label{mrp}
The \emph{function class} $M(\rp)$\index{$M(\rp)$}
is defined to be the set of all the positive functions
$\omega$ on $\rp$ such that, for any $0<\delta<N<\infty$,
$$0<\inf_{t\in(\delta,N)}\omega(t)\leq\sup_
{t\in(\delta,N)}\omega(t)<\infty$$ and
there exist four constants $\alpha_{0},\
  \beta_{0},\ \alpha_{\infty},\ \beta_{\infty}\in\rr$ such that
\begin{enumerate}
  \item[(i)] for any $t\in(0,1]$,
  $\omega(t)t^{-\alpha_{0}}$ is almost increasing and
  $\omega(t)t^{-\beta_{0}}$ is almost decreasing;
  \item[(ii)] for any $t\in[1,\infty)$,
  $\omega(t)t^{-\alpha_{\infty}}$ is
  almost increasing and  $\omega(t)t^{-\beta_{\infty}}$ is
  almost decreasing.
\end{enumerate}
\end{definition}

\begin{remark} Let $\omega\in M(\rp)$, with $\alpha_0$,
$\beta_0$, $\alpha_\infty$,
and $\beta_\infty$, be as in Definition \ref{mrp}.
\begin{enumerate}
\item[{\rm(i)}] Applying Definition \ref{mrp}(i),
we conclude that,
for any $t\in(0,1]$,
\begin{equation}\label{remark 1.1.3e1}
t^{\beta_0}\lesssim\omega(t)\lesssim t^{\alpha_0},
\end{equation}
which further implies that $\alpha_0\leq\beta_0$.
However, we should point out that, for any
positive function $\omega$ on $\rp$, the condition
\eqref{remark 1.1.3e1} may
not imply Definition \ref{mrp}(i),
and hence both \eqref{remark 1.1.3e1} and Definition
\ref{mrp}(i) may not be equivalent.
Indeed, for any given
$-\infty<\alpha_0<\beta_0<\infty$, let
\begin{align*}\omega(t):=&\,
2n(2n-1)\left[(2n-1)^{-\alpha_0}-(2n)^{-\beta_0}\right]
\left(t-\frac{1}{2n-1}\right)\\
&+(2n-1)^{-\alpha_0}
\end{align*}
for any $t\in(\frac{1}{2n},
\frac{1}{2n-1}]$ with $n\in\mathbb{N}$,
$$
\omega(t):=2n(2n+1)\left[(2n)^{-\beta_0}-
(2n+1)^{-\alpha_0}\right]
\left(t-\frac{1}{2n}\right)+(2n)^{-\beta_0}
$$
for any $t\in(\frac{1}{2n+1},\frac{1}{2n}]$
with $n\in\mathbb{N}$,
and $\omega(t):=1$ for any $t\in(1,\infty)$.
Then, obviously, for any $t\in(0,1]$, we have
$$
t^{\beta_0}\leq\omega(t)\leq t^{\alpha_0}.
$$
However, it is easy to show that
\begin{align*}
\frac{\omega(\frac{1}{2n+1})}{\omega(
\frac{1}{2n})}\left(
\frac{\frac{1}{2n+1}}{\frac{1}{2n}}
\right)^{-\alpha_0}=(2n)^
{\beta_0-\alpha_0}\to\infty
\end{align*}
as $n\to\infty$, which implies that
$\omega(t)t^{-\alpha_0}$ is not
almost increasing when $t\in(0,1]$.
Similarly, it holds true that
\begin{align*}
\frac{\omega(\frac{1}{2n-1})}
{\omega(\frac{1}{2n})}
\left(\frac{\frac{1}{2n-1}}{\frac{1}{2n}}\right)
^{-\beta_0}=(2n-1)^{\beta_0-\alpha_0}\to\infty
\end{align*}
as $n\to\infty$, which further implies that
$\omega(t)t^{-\beta_0}$ is not
almost decreasing when $t\in(0,1]$.
\item[{\rm(ii)}] From Definition \ref{mrp}(ii),
we deduce that,
for any $t\in[1,\infty)$,
\begin{equation}\label{remark 1.1.3e2}
t^{\alpha_\infty}\lesssim\omega(t)\lesssim t^{\beta_\infty}.
\end{equation}
By this, we further find that $\alpha_\infty\leq\beta_\infty$.
However, we should point out that, for any
positive function $\omega$ on $\rp$, the condition
\eqref{remark 1.1.3e2} may
not imply Definition \ref{mrp}(ii), and hence both
\eqref{remark 1.1.3e2} and Definition \ref{mrp}(ii)
may not be equivalent. Indeed,
for any given
$-\infty<\alpha_\infty<\beta_\infty<\infty$,
let $\omega(t):=1$ for any $t\in(0,1]$,
$$
\omega(t):=\left[(2n)^{\beta_\infty}-
(2n-1)^{\alpha_\infty}\right](t-2n)+(2n)^{\beta_\infty}
$$
for any $t\in(2n-1,2n]$ with $n\in\mathbb{N}$, and
$$
\omega(t):=\left[(2n+1)^{\alpha_\infty}-(2n)
^{\beta_\infty}\right](t-2n-1)+(2n+1)^{\alpha_\infty}
$$
for any $t\in(2n,2n+1]$ with $n\in\mathbb{N}$.
Then, in this case, it holds true that,
for any $t\in[1,\infty)$,
$$
t^{\alpha_\infty}\leq\omega(t)\leq t^{\beta_\infty}.
$$
However, we can easily find that
\begin{align*}
\frac{\omega(2n)}{\omega(2n+1)}\left(
\frac{2n}{2n+1}\right)^{-\alpha_\infty}=(2n)^
{\beta_\infty-\alpha_\infty}\to\infty
\end{align*}
as $n\to\infty$. Thus, $\omega(t)t^{-\alpha_\infty}$ is not
almost increasing when $t\in[1,\infty)$.
On another hand, we have
\begin{align*}
\frac{\omega(2n)}{\omega(2n-1)}\left(
\frac{2n}{2n-1}\right)^{-\beta_\infty}=(2n-1)^
{\beta_\infty-\alpha_\infty}\to\infty
\end{align*}
as $n\to\infty$, which further implies that
$\omega(t)t^{-\beta_\infty}$ is not
almost decreasing when $t\in[1,\infty)$.
\end{enumerate}
\end{remark}

The following property of the function class
$M(\rp)$ was stated in \cite[(7)]{HSamko}
without the proof. For the convenience
of the reader, we present it as follows
and give its detailed proof.

\begin{lemma}\label{mono}
Let $\omega\in M(\rp)$. Then, for any given
$\lambda\in(0,\infty)$, there exists
a positive constant $C_{(\omega,\lambda)}$,
depending only on $\omega$ and $\lambda$,
such that, for any $t\in(0,\infty)$,
$$\omega(\lambda t)\leq
C_{(\omega,\lambda)}\omega(t).$$
\end{lemma}

\begin{proof}
Let $\omega\in M(\rp)$, with $\alpha_0$, $\beta_0$,
$\alpha_\infty$, $\beta_\infty\in\rr$, be
as in Definition \ref{mrp}.
Then, there exists a
constant $C\in[1,\infty)$ such that,
for any $h$, $t\in(0,1]$,
\begin{equation}\label{monoe1}
\omega(ht)(ht)^{-\alpha_0}\leq
C\omega(h)h^{-\alpha_0}
\end{equation}
and
\begin{equation}\label{monoe2}
\omega(h)h^{-\beta_0}\leq
C\omega(ht)(ht)^{-\beta_0}
\end{equation}
and, for any $h$, $t\in[1,\infty)$,
\begin{equation}\label{monoe3}
\omega(h)h^{-\alpha_\infty}\leq
C\omega(ht)(ht)^{-\alpha_\infty}
\end{equation}
and
\begin{equation}\label{monoe4}
\omega(ht)(ht)^{-\beta_\infty}\leq
C\omega(h)h^{-\beta_\infty}.
\end{equation}
Now, we show the desired inequality
by considering the following
two cases on $\lambda$.

\emph{Case 1)} $\lambda\in(0,1]$. In this case,
using \eqref{monoe1}, we conclude that,
for any $t\in(0,1)$,
$$
\omega(\lambda t)(\lambda t)^{-\alpha_0}
\leq C\omega(t)t^{-\alpha_0},
$$
which further implies that
\begin{align}\label{monoe5}
\omega(\lambda t)\leq C\lambda^{\alpha_0}
\omega(t)\leq C^2\lambda
^{\min\{\alpha_0,\alpha_\infty\}}\omega(t).
\end{align}
In addition, combining \eqref{monoe1}
and \eqref{monoe3},
we find that, for any $t\in[1,\frac{1}{\lambda}]$,
\begin{align*}
&\omega(\lambda t)(\lambda t)^
{-\min\{\alpha_0,\alpha_\infty\}}\\
&\quad\leq\omega(\lambda t)
(\lambda t)^{-\alpha_0}\leq C\omega(1)
\leq C^2\omega(t)t^{-\alpha_\infty}\notag\\
&\quad\leq C^2\omega(t)t^{-\min\{\alpha_0,\alpha_\infty\}}.
\end{align*}
This implies that, for any $t\in[1,\frac{1}{\lambda}]$,
\begin{equation}\label{monoe6}
\omega(\lambda t)\leq
C^2\lambda^{\min\{\alpha_0,\alpha_\infty\}}\omega(t).
\end{equation}
On the other hand, applying \eqref{monoe3}, we know that,
for any $t\in(\frac{1}{\lambda},\infty)$,
\begin{align*}
\omega(\lambda t)(\lambda t)^{-\alpha_\infty}
\leq C\omega(t)t^{-\alpha_\infty},
\end{align*}
which further implies that
\begin{align*}
\omega(\lambda t)\leq C\lambda
^{\alpha_\infty}\omega(t)\leq C^2\lambda
^{\min\{\alpha_0,\alpha_\infty\}}\omega(t).
\end{align*}
This, together with \eqref{monoe5} and \eqref{monoe6},
further implies that, in this case,
Lemma \ref{mono} holds true
with $C_{(\omega,\lambda)}:=C^2
\lambda^{\min\{\alpha_0,\alpha_\infty\}}$.

\emph{Case 2)} $\lambda\in(1,\infty)$. In this case,
applying \eqref{monoe2}, \eqref{monoe4},
and an argument similar to that used
in the proof of Case 1),
we conclude that Lemma \ref{mono} holds true with
$C_{(\omega,\lambda)}:=C^2
\lambda^{\max\{\beta_0,\beta_\infty\}}$,
which completes the proof of Lemma \ref{mono}.
\end{proof}

To describe the properties of positive functions at
origin and infinity, we present the
following Matuszewska--Orlicz indices
from \cite{MO,mo65} (see also \cite{HSamko}).
Throughout this book, the
\emph{symbol} $h\to0^+$\index{$h\to0^+$}
means that $h\in(0,\infty)$ and $h\to0$.

\begin{definition}\label{moindex}
Let $\omega$ be a positive function on $\rp$. Then
the \emph{Matuszewska--Orlicz
indices}\index{Matuszewska--Orlicz index}
$\m0(\omega)$\index{$\m0(\omega)$},
$\M0(\omega)$\index{$\M0(\omega)$},
$\mi(\omega)$\index{$\mi(\omega)$},
and $\MI(\omega)$\index{$\MI(\omega)$} of
$\omega$ are defined, respectively, by setting, for
any $h\in(0,\infty)$,
$$\m0(\omega):=\sup_{t\in(0,1)}
\frac{\ln(\lims\limits_{h\to0^+}\frac{\omega(ht)}
  {\omega(h)})}{\ln t},$$
  $$\M0(\omega):=\inf_{t\in(0,1)}\frac{\ln(
\limi\limits_{h\to0^+}
\frac{\omega(ht)}{\omega(h)})}{\ln t},$$
$$
\mi(\omega):=\sup_{t\in(1,\infty)}
\frac{\ln(\limi\limits_{h\to\infty}
\frac{\omega(ht)}{\omega(h)})}{\ln t},$$
and
$$\MI(\omega)
:=\inf_{t\in(1,\infty)}\frac{\ln(\lims\limits_
{h\to\infty}\frac{\omega(ht)}{\omega(h)})}{\ln t}.
$$
\end{definition}

\begin{remark}\label{remark 1.1.3}
\begin{enumerate}
  \item[(i)] By \cite[p.\,12]{HSamko}, we find that,
  for any $\omega\in M(\rp)$ and $h\in(0,\infty)$,
  $$
m_{0}(\omega)=\lim\limits_{t\to0^+}
\frac{\ln(\lims\limits_{h\to0^+}
\frac{\omega(ht)}{\omega(h)})}{\ln t},$$
$$M_{0}(\omega)=\lim\limits_{t\to0^+}
\frac{\ln(\limi\limits_{h\to0^+}
\frac{\omega(ht)}{\omega(h)})}{\ln t},
$$
$$
m_{\infty}(\omega)=
\lim_{t\to\infty}\frac{\ln(\limi\limits_
{h\to\infty}\frac{\omega(ht)}{\omega(h)})}{\ln t},$$
and
$$M_{\infty}(\omega)=\lim_{t\to\infty}
\frac{\ln(\lims\limits_{h\to\infty}
\frac{\omega(ht)}{\omega(h)})}{\ln t},$$
where, for any positive function $\omega$ on
$\rp$, $\m0(\omega)$, $\M0(\omega)$,
$\mi(\omega)$, and $\MI(\omega)$
are as in Definition \ref{moindex}.
  \item[(ii)] We point out that,
  for any positive function $\omega$,
  its Matuszewska--Orlicz indices
  $\m0(\omega)$, $\M0(\omega)$,
$\mi(\omega)$, and $\MI(\omega)$
may not be finite.
Indeed, let $\omega(t):=e^{1/t}$
for any $t\in(0,\infty)$.
Then, in this case, we have $\m0(\omega)=-\infty$.
However, for any positive
function $\omega\in M(\rp)$, applying
\cite[Theorem 11.3]{LMO} and \cite[(6.14)]{SWGM},
we conclude that
$$
-\infty<\m0(\omega)\leq\M0(\omega)<\infty
$$
and
$$
-\infty<\mi(\omega)\leq\MI(\omega)<\infty.
$$
\item[(iii)] For any given positive function
$\omega$ on $\rp$ and for any $t\in(0,\infty)$,
  let $\omega_{*}\index{$\omega_*$}(t):=\omega(\frac{1}{t})$.
  Then, from \cite[(6.14)]{SWGM}, it follows that
  $$
  \m0(\omega_*)=-\MI(\omega)$$
  and
  $$
  \M0(\omega_*)=-\mi(\omega),
  $$
  where, for any positive function
  $\phi$, $\m0(\phi)$, $\M0(\phi)$, $\mi(\phi)$,
  and $\MI(\phi)$ denote its
  Matuszewska--Orlicz indices.
\end{enumerate}
\end{remark}

For both the function class $M(\rp)$
and the Matuszewska--Orlicz indices,
we have the following conclusions which are just
\cite[(6.4), (6.5), and (6.14)]{SWGM}.

\begin{lemma}\label{rela}
Let $\omega\in M(\rp)$. Then, for any given
$s\in(0,\infty)$, it holds true that
$1/\omega,\ \omega^s\in M(\rp)$ and
\begin{enumerate}
  \item[{\rm(i)}] $\m0\left(\frac{1}{\omega}\right)
  =-\M0(\omega)
  \text{ and }\M0\left(\frac{1}{\omega}\right)
  =-\m0(\omega);$
  \item[{\rm(ii)}] $\mi\left(\frac{1}{\omega}\right)
  =-\MI(\omega)\text{ and }
  \MI\left(\frac{1}{\omega}\right)=-\mi(\omega);$
  \item[{\rm(iii)}] $\m0(\omega^{s})=s\m0(\omega)
  \text{ and }\M0(\omega^{s})=s\M0(\omega);$
  \item[{\rm(iv)}] $\mi(\omega^{s})=s\mi(\omega)
  \text{ and }\MI(\omega^{s})=s\MI(\omega).$
\end{enumerate}
Here, for any positive function
$\phi$ on $\rp$, $\m0(\phi)$, $\M0(\phi)$, $\mi(\phi)$,
and $\MI(\phi)$ denote its Matuszewska--Orlicz indices.
\end{lemma}

Now, we give three typical examples of
functions in $M(\rp)$ and
their concrete Matuszewska--Orlicz indices.

\begin{example}\label{ex113}
For any given $\alpha\in\rr$
and for any $t\in(0,\infty)$, let $\omega(t):=t^\alpha$.
Then, using the fact that, for any $t\in(0,\infty)$,
$\omega(t)t^{-\alpha}=1$, we find that $\omega\in M(\rp)$.
In addition, for any $h$, $t\in(0,\infty)$,
we have
\begin{align*}
\frac{\omega(ht)}{\omega(h)}=
\frac{(ht)^\alpha}{h^\alpha}=t^\alpha.
\end{align*}
By this and Definition \ref{moindex}, we conclude that
\begin{align*}
\m0(\omega)=\sup_{t\in(0,1)}
\frac{\ln(\lims\limits_{h\to0^+}
\frac{\omega(ht)}{\omega(h)})}{\ln t}=\sup_{t\in(0,1)}
\frac{\ln(t^\alpha)}{\ln t}=\alpha.
\end{align*}
Similarly, we easily obtain
$$\mi(\omega)=\M0(\omega)=\MI(\omega)=\alpha.$$
\end{example}

\begin{example}
For any $t\in(0,\infty)$,
let $$\omega(t):=\frac{t}{\ln(e+t)}.$$
Then notice that, for any $t\in(0,\infty)$, we have
\begin{align*}
\omega'(t)&=\frac{\ln(e+t)
-\frac{t}{e+t}}{[\ln(e+t)]^2}
=\frac{(e+t)\ln(e+t)-t}{(e+t)[\ln(e+t)]^2}\\
&>\frac{t\ln e-t}{(e+t)[\ln(e+t)]^2}=0,
\end{align*}
which implies that $\omega(t)$ is
increasing when $t\in(0,\infty)$.
On the other hand, it is easy to show that
$\omega(t)t^{-1}$ is decreasing when
$t\in(0,\infty)$. Therefore, we conclude that
$\omega\in M(\rp)$. Moreover,
for any $t\in(0,\infty)$, we have
\begin{align*}
\lim\limits_{h\to0^+}\frac{\omega(ht)}{\omega(h)}
=\frac{\frac{ht}{\ln(e+ht)}}{\frac{h}{\ln(e+h)}}
=t\lim\limits_{h\to0^+}\frac{\ln(e+h)}{\ln(e+ht)}=t.
\end{align*}
From this and Definition \ref{moindex},
we deduce that
\begin{align*}
\m0(\omega)=\sup_{t\in(0,1)}
\frac{\ln(\lims\limits_{h\to0^+}
\frac{\omega(ht)}{\omega(h)})}{\ln t}
=\sup_{t\in(0,1)}\frac{\ln t}{\ln t}=1
\end{align*}
and
\begin{align*}
\M0(\omega)=\inf_{t\in(0,1)}
\frac{\ln(\limi\limits_{h\to0^+}
\frac{\omega(ht)}{\omega(h)})}{\ln t}
=\inf_{t\in(0,1)}\frac{\ln t}{\ln t}=1.
\end{align*}
In addition, for any $t\in(0,\infty)$, we have
\begin{align*}
\lim\limits_{h\to\infty}
\frac{\omega(ht)}{\omega(h)}
&=\lim\limits_{h\to\infty}
\frac{\frac{ht}{\ln(e+ht)}}{\frac{
h}{\ln(e+h)}}
=t\lim\limits_{h\to\infty}
\frac{\ln(e+h)}{\ln(e+ht)}\\
&=t\lim\limits_{h\to\infty}
\frac{1+\frac{\ln(\frac{e}{h}+1)}{\ln h}}
{1+\frac{\ln(\frac{e}{h}+t)}{\ln h}}=t.
\end{align*}
This, together with Definition \ref{moindex},
further implies that
\begin{align*}
\mi(\omega)=\sup_{t\in(1,\infty)}
\frac{\ln(\limi\limits_{h\to\infty}
\frac{\omega(ht)}{\omega(h)})}{\ln t}
=\sup_{t\in(1,\infty)}\frac{\ln t}{\ln t}=1
\end{align*}
and
\begin{align*}
\MI(\omega)=\inf_{t\in(1,\infty)}
\frac{\ln(\lims\limits_{h\to\infty}
\frac{\omega(ht)}{\omega(h)})}{\ln t}
=\inf_{t\in(1,\infty)}\frac{\ln t}{\ln t}=1.
\end{align*}
\end{example}

\begin{example}\label{ex117}
For any given $\alpha_{1}$,
$\alpha_{2}\in\rr$ and for
any $t\in(0,\infty)$, let
$$
\omega(t):=\left\{\begin{aligned}
&t^{\alpha_{1}}(1-\ln t)
\hspace{0.3cm}\text{when }t\in(0,1],\\
&t^{\alpha_{2}}(1+\ln t)
\hspace{0.3cm}\text{when }t\in(1,\infty).
\end{aligned}\right.
$$
Then, we find that, for any $t\in(0,1)$,
\begin{align*}
\left(\omega(t)t^{-\alpha_1+1}\right)'
=\left(t(1-\ln t)\right)'=-\ln t>0
\end{align*}
and, for any $t\in(1,\infty)$,
\begin{align*}
\left(\omega(t)t^{-\alpha_2-1}\right)'
=\left(\frac{1+\ln t}{t}\right)'
=-\frac{\ln t}{t^2}<0.
\end{align*}
This further implies that
$\omega(t)t^{-\alpha_1+1}$ is increasing
when $t\in(0,1]$, and $\omega(t)^{-\alpha_2-1}$
is decreasing when $t\in[1,\infty)$.
On the other hand, we can easily know that
$\omega(t)t^{-\alpha_1}$ is decreasing when
$t\in(0,1]$, and
$\omega(t)t^{-\alpha_2}$ is increasing
when $t\in[1,\infty)$.
Thus, the positive function $\omega$ belongs
to $M(\rp)$. In addition, it holds true that,
for any $h$, $t\in(0,1)$,
\begin{align*}
\lim\limits_{h\to0^+}\frac{\omega(ht)}{\omega(h)}
=\lim\limits_{h\to0^+}\frac{(ht)^{\alpha_1}[1-\ln(ht)]}
{h^{\alpha_1}(1-\ln h)}
=t^{\alpha_1}\lim\limits_{h\to0^+}\left(
1-\frac{\ln t}{1-\ln h}\right)=t^{\alpha_1}.
\end{align*}
Combining this and Definition \ref{moindex}, we
conclude that
\begin{align*}
\m0(\omega)=\sup_{t\in(0,1)}\frac{
\ln(\lims\limits_{h\to0^+}
\frac{\omega(ht)}{\omega(h)})}{
\ln t}=\sup_{t\in(0,1)}
\frac{\ln(t^{\alpha_1})}{\ln t}=\alpha_1
\end{align*}
and
\begin{align*}
\M0(\omega)=\inf_{t\in(0,1)}\frac{
\ln(\limi\limits_{h\to0^+}
\frac{\omega(ht)}{\omega(h)})}{
\ln t}=\inf_{t\in(0,1)}
\frac{\ln(t^{\alpha_1})}{\ln t}=\alpha_1.
\end{align*}
Similarly, from Definition \ref{moindex},
we deduce that
$$
\mi(\omega)=\MI(\omega)=\alpha_{2}.
$$
\end{example}

The following equivalent formulae of
Matuszewska--Orlicz indices
are just \cite[Theorem 11.13]{LMO}.

\begin{lemma}\label{re}
Let $\omega\in M(\rp)$, and $\m0(\omega)$,
$\M0(\omega)$, $\mi(\omega)$, and $\MI(\omega)$
denote its Matuszewska--Orlicz indices. Then
\begin{equation}\label{re1}
\m0(\omega)=\sup\left\{\alpha_0\in\rr:\
\omega(t)t^{-\alpha_0}
\text{ is almost increasing
for any $t\in(0,1]$}\right\},
\end{equation}
$$
\M0(\omega)=\inf\left\{\beta_0\in\rr:\
\omega(t)t^{-\beta_0}
\text{ is almost decreasing
for any $t\in(0,1]$}\right\},
$$
$$
\mi(\omega)=\sup\left\{\alpha_\infty\in\rr:\
\omega(t)t^{-\alpha_\infty}
\text{ is almost increasing
for any $t\in[1,\infty)$}\right\},
$$
and
$$
\MI(\omega)=\inf\left\{\beta_\infty\in\rr:\
\omega(t)t^{-\beta_\infty}
\text{ is almost decreasing
for any $t\in[1,\infty)$}\right\}.
$$
\end{lemma}

\begin{remark}
We should point out that the suprema and the infima
in Lemma \ref{re} may not be achieved.
Indeed, for any $t\in(0,\infty)$,
let
$$
\omega(t):=\left\{\begin{aligned}
&t(1-\ln t),\hspace{0.3cm}\text{when }t\in(0,1],\\
&t(1+\ln t),\hspace{0.3cm}\text{when }t\in(1,\infty).
\end{aligned}\right.
$$
Then, by Example \ref{ex117}, we conclude that
$\m0(\omega)=1$. However, we have
\begin{align*}
\frac{\omega(\frac{1}{n})(\frac{1}{n})^{-1}}
{\omega(\frac{1}{2})(\frac{1}{2})^{-1}}
=\frac{1+\ln n}{1+\ln2}
\to\infty
\end{align*}
as $n\to\infty$. This implies that
$\omega(t)t^{-\m0(\omega)}$ is not
almost increasing when
$t\in(0,1]$. Thus, in this case,
the supremum in \eqref{re1}
of the Matuszewska--Orlicz index $\m0$
can not be achieved.
\end{remark}

Applying Lemma \ref{re}, we immediately obtain the
following estimates of the
positive function $\omega\in M(\rp)$,
which were also stated in
\cite[(25) and (26)]{HSamko}; we
omit the details.

\begin{lemma}\label{Th1}
Let $\omega\in M(\rp)$, and $\m0(\omega)$, $\M0(\omega)$,
$\mi(\omega)$, and $\MI(\omega)$ denote
its Matuszewska--Orlicz indices.
Then, for any given $\eps\in(0,\infty)$,
there exists a constant $C_{(\eps)}\in[1,\infty)$,
depending on $\eps$, such that,
for any t $\in(0,1]$,
\begin{equation*}
C_{(\eps)}^{-1}t^{\M0(\omega)+\eps}
\leq\inf_{\tau\in(0,1]}
\frac{\omega(t\tau)}{\omega(\tau)}
\leq\sup_{\tau\in(0,1]}
\frac{\omega(t\tau)}{\omega(\tau)}
\leq
C_{(\eps)}t^{\m0(\omega)-\eps}
\end{equation*}
and, for any t $\in[1,\infty)$,
\begin{equation*}
C_{(\eps)}^{-1}t^{\mi(\omega)-\eps}
\leq\inf_{\tau\in[1,\infty)}
\frac{\omega(t\tau)}{\omega(\tau)}
\leq\sup_{\tau\in[1,\infty)}
\frac{\omega(t\tau)}{\omega(\tau)}
\leq C_{(\eps)}t^{\MI(\omega)+\eps}.
\end{equation*}
\end{lemma}

\begin{remark}
We should point out that, in Lemma \ref{Th1},
the constant $C_{(\eps)}$ may not be uniform about
$\eps$. Indeed, for any $t\in(0,\infty)$, let
$$
\omega(t):=\left\{\begin{aligned}
&t(1-\ln t)\hspace{0.3cm}\text{when }t\in(0,1],\\
&t(1+\ln t)\hspace{0.3cm}\text{when }t\in(1,\infty).
\end{aligned}\right.
$$
Then, from Example \ref{ex117},
it follows that
$\m0(\omega)=1$. By this, we conclude that
\begin{align*}
&\sup_{t,\,\tau\in(0,1]}\left\{\frac{\omega(t\tau)}
{\omega(\tau)t^{\m0(\omega)-\eps}}\right\}\\
&\quad=\sup_{t,\,\tau\in(0,1]}\left\{t^{\eps}
\frac{1-\ln\tau-\ln t}{1-\ln\tau}\right\}=
\sup_{t\in(0,1]}\left\{t^{\eps}(1-\ln t)\right\}\\
&\quad=\frac{e^{\eps-1}}{\eps}\to\infty
\end{align*}
as $\eps\to0^+$, which further
implies that $C_{(\eps)}\to\infty$ as
$\eps\to0^+$, where $C_{(\eps)}$ is as
in Lemma \ref{Th1}. Here and thereafter,
$\eps\to0^+$ means $\eps\in(0,\infty)$
and $\eps\to0$.
\end{remark}

\section{Generalized Herz Spaces}

The targets of this section are threefold.
The first one is to recall the concept
of generalized Herz spaces of
Rafeiro and Samko given
in \cite[Definition 2.2]{HSamko}.
The second one is to recall
some concepts related to
ball quasi-Banach function spaces.
The last one is to show that
generalized Herz spaces of Rafeiro and Samko
are ball Banach function spaces or
ball quasi-Banach function spaces under some
reasonable and sharp assumptions.

First, we present the concept of generalized Herz spaces
which were introduced by
Rafeiro and Samko in \cite[Definition 2.2]{HSamko}.
In what follows, for any $x\in\rrn$ and $r\in(0,\infty)$,
let\index{$B(x,r)$} $$B(x,r):=\left\{y\in\rrn:\
|x-y|<r\right\}$$ and
\begin{equation}\label{Bset}
\mathbb{B}\index{$\mathbb{B}$}:=\left\{B(x,r):\
x\in\rrn\text{ and }r\in(0,\infty)\right\}.
\end{equation}
Moreover, for any ball $B:=B(x,r)$ with
$x\in\rrn$ and $r\in(0,\infty)$,
$r$ is called the \emph{radius}\index{radius}
of $B$, which is denoted by $r(B)$.

\begin{definition}\label{gh}
Let $p,\ q\in(0,\infty)$ and $\omega\in M(\rp)$.
\begin{enumerate}
\item[{\rm(i)}] The \emph{local generalized
Herz space}\index{local generalized Herz space}
  $\HerzSo$\index{$\HerzSo$}
  is defined to be the set of all
  the measurable functions $f$
  on $\rrn$ such that
\begin{equation}\label{ghl}
\|f\|_{\HerzSo}:=\left\{\sum_{k\in\mathbb{Z}}
\left[\omega(2^{k})\right]^{q}
\left\|f\1bf_{B(\0bf,2^{k})\setminus B(\0bf,2^{k-1})}
\right\|_{L^{p}(\rrn)}^{q}\right\}^{\frac{1}{q}}
\end{equation}
is finite.
\item[{\rm(ii)}] The
\emph{global generalized
Herz space}\index{global generalized Herz space}
$\HerzS$\index{$\HerzS$} is defined
to be the set of all the
measurable functions $f$ on $\rrn$ such that
\begin{align}\label{ghg}
\|f\|_{\HerzS}:=\sup_{\xi\in\rrn}
\left\{\sum_{k\in\mathbb{Z}}\left[\omega(2^{k})
\right]^{q}\left\|f\1bf_{B(\xi,2^{k})
\setminus B(\xi,2^{k-1})}
\right\|_{L^{p}(\rrn)}^{q}\right\}
^{\frac{1}{q}}
\end{align}
is finite.
\end{enumerate}
\end{definition}

\begin{remark}\label{remhs}
\begin{enumerate}
  \item[{\rm(i)}] Recall that, in \cite{HSamko},
  Rafeiro and Samko also introduced the
  following \emph{continuous versions} of
  both \eqref{ghl} and \eqref{ghg},
  respectively, by setting,
  for any measurable function $f$ on $\rrn$,
  \begin{align}\label{ghcl}
\|f\|_{\textbf{\textsl{H}}_{\omega,\0bf}
^{p,q}(\rrn)}\index{$\textbf{\textsl{H}}
_{\omega,\0bf}^{p,q}(\rrn)$}
:=\left\{\int_{0}^{\infty}
\left[\omega(t)\right]^{q}
\left[\int_{t<|y|<2t}|f(y)|^{p}\,dy
\right]^{\frac{q}{p}}\,\frac{dt}{t}
\right\}^{\frac{1}{q}}
\end{align}
and
\begin{equation}\label{ghcg}
\|f\|_{\textbf{\textsl{H}}_{\omega}^{p,q}
(\rrn)}:=\sup_{\xi\in\rrn}\index{$\textbf{\textsl{H}}
_{\omega}^{p,q}(\rrn)$}
\left\{\int_{0}^{\infty}
\left[\omega(t)\right]^{q}
\left[\int_{t<|y-\xi|<2t}|f(y)|^{p}\,dy
\right]^{\frac{q}{p}}\,\frac{dt}{t}
\right\}^{\frac{1}{q}},
\end{equation}
where $p$, $q\in(0,\infty)$ and $\omega\in M(\rp)$.
We point out that, from \cite[Lemma 2.3]{HSamko},
it follows that, when $p$, $q\in[1,\infty)$
and $\omega\in M(\rp)$, the continuous
version \eqref{ghcl} is equivalent to the
discrete version \eqref{ghl} and, similarly,
\eqref{ghcg} is equivalent to \eqref{ghg}.
  \item[{\rm(ii)}] Obviously,
  by Definition \ref{gh}, we conclude that,
 for any measurable function $f$ on $\rrn$,
  $$\|f\|_{\HerzS}=\sup_{\xi\in\rrn}
  \left\|f(\cdot+\xi)\right\|_{\HerzSo}.$$
  \item[{\rm(iii)}] Observe that, in Definition \ref{gh}(i),
  for any given $\alpha\in\rr$ and for
  any $t\in(0,\infty)$, let $\omega(t):=t^\alpha$.
  Then, in this case, the local generalized Herz space
  $\HerzSo$ goes back to the classical
  \emph{homogeneous Herz
  space}\index{homogeneous Herz space}
  $\dot{K}_{p}^{\alpha,q}
  (\rrn)$\index{$\dot{K}_{p}^{\alpha,q}(\rrn)$},
   which was originally introduced
   in \cite[Definition 1.1(a)]{LuY95}
   (see also \cite[Chapter 1]{LYH}),
   with the same quasi-norms.
  \item[{\rm(iv)}] Let $p$, $q\in[1,\infty)$
  and $\omega\in M(\rp)$
  with $\Mw\in(-\infty,0)$. Then, applying
  \cite[Theorem 3.2 and Remark 3.3]{HSamko},
  we conclude that, in this case, the
  generalized Herz spaces $\HerzSo$ and $\HerzS$,
  in the sense of equivalent norms, coincide,
  respectively, with the \emph{local generalized Morrey
  space\index{local generalized Morrey space}
  $\MorrSo$\index{$\MorrSo$}}
  and the \emph{global generalized Morrey
  space\index{global generalized Morrey space}
  $\MorrS$\index{$\MorrS$}},
  originally introduced by
  Guliev and Mustafaev in \cite{gm98}
  (see also \cite{HSamko}), which are defined,
  respectively, to be the sets of
  all the measurable functions $f$ on $\rrn$
  such that
  $$
  \|f\|_{\MorrSo}:=\left\{\int_{0}^{\infty}
  \left[\omega(t)\right]^{q}
\left[\int_{|y|<t}|f(y)|^{p}\,dy
\right]^{\frac{q}{p}}\,\frac{dt}{t}\right\}
^{\frac{1}{q}}<\infty$$
and
$$
\|f\|_{\MorrS}:=\sup_{\xi\in\rrn}
\left\{\int_{0}^{\infty}\left[\omega(t)\right]^{q}
\left[\int_{|y-\xi|<t}|f(y)|^{p}\,dy
\right]^{\frac{q}{p}}\,\frac{dt}{t}\right\}
^{\frac{1}{q}}<\infty.
  $$
\end{enumerate}
\end{remark}

Now, we recall the definition of
ball quasi-Banach function
spaces and some other concepts which are used later.
In what follows, we always use the
\emph{symbol} $\Msc(\rrn)$
to denote the set of all measurable functions on $\rrn$.

The following concept of ball quasi-Banach function
spaces was given in \cite[Definition 2.2]{SHYY}.

\begin{definition}\label{Df1}
A quasi-normed linear space
$X\subset\Msc(\rrn)$,
equipped with a quasi-norm $\|\cdot\|_X$
which makes sense for all functions
in $\Msc(\rrn)$, is called a
\emph{ball quasi-Banach function space}\index{ball
quasi-Banach function space} (for short,
BQBF space\index{BQBF space}) if
it satisfies
\begin{enumerate}
\item[(i)] for any $f\in\Msc(\rrn)$,
$\|f\|_{X}=0$ implies that
$f=0$ almost everywhere;
\item[(ii)] for any $f$, $g\in\Msc(\rrn)$,
$|g|\leq|f|$ almost everywhere implies
that $\|g\|_{X}\leq\|f\|_{X}$;
\item[(iii)] for any $\{f\}_{m\in\mathbb{N}}
\subset\Msc(\rrn)$ and $f\in\Msc(\rrn)$,
$0\leq f_{m}\uparrow f$ as
$m\to\infty$ almost everywhere
 implies that $\|f_{m}\|_{X}\uparrow\|f\|_{X}$
 as $m\to\infty$;
\item[(iv)] $B\in\mathbb{B}$ implies that
$\1bf_{B}\in X$, where $\mathbb{B}$
is as in \eqref{Bset}.
\end{enumerate}
\end{definition}

Moreover, a ball quasi-Banach
function space $X$ is called
a \emph{ball Banach function
space}\index{ball Banach function space}
(for short, BBF space\index{BBF space})
if the quasi-norm $\|\cdot\|_X$
of $X$ satisfies the triangle
inequality\index{triangle inequality}:
for any $f,\ g\in X$,
$$
\|f+g\|_{X}\leq\|f\|_{X}+\|g\|_{X}
$$
and, for any $B\in\mathbb{B}$,
there exists a positive constant $C_{(B)}$,
depending on
$B$, such that, for any $f\in X$,
\begin{equation}\label{bal}
\int_{B}|f(x)|\,dx\leq C_{(B)}\|f\|_{X}.
\end{equation}

\begin{remark}
\begin{enumerate}
  \item[{\rm(i)}] We point out that,
  in Definition \ref{Df1}(iv),
  if we replace any ball $B$
  by any bounded measurable set
  $E$, we obtain an equivalent formulation
  of ball quasi-Banach function spaces.

  \item[{\rm(ii)}] Recall that a
  quasi-Banach space $X\subset\Msc(\rrn)$
  is called a \emph{quasi-Banach
  function space}\index{quasi-Banach function space}
  if it is a ball quasi-Banach function space
  and it satisfies Definition \ref{Df1}(iv)
  with ball replaced by any
  measurable set of \emph{finite measure}
(see, for instance,
\cite[Chapter 1, Definitions 1.1 and 1.3]{BSIO}).
It is easy to see that every quasi-Banach
function space is a ball quasi-Banach
function space, and the converse
is not necessary to be true. As was mentioned in
\cite[p.\ 9]{SHYY} and \cite[Section 5]{WYYZ},
the family of ball Banach function spaces
includes Morrey spaces, mixed-norm Lebesgue spaces,
variable Lebesgue spaces, weighted Lebesgue spaces,
and Orlicz--slice spaces,
which are not necessary to
be Banach function spaces.
\end{enumerate}
\end{remark}

The following basic properties of
ball quasi-Banach function spaces
can be deduced directly from
Definition \ref{Df1}; we omit the details.
\begin{lemma}
Let $X\subset\Msc(\rrn)$ be a
ball quasi-Banach function space.
Then
\begin{enumerate}
  \item[{\rm(i)}] for any $f\in X$,
   $\|f\|_{X}=0$ if and only if $f=0$
   almost everywhere in $\rrn$;
  \item[{\rm(ii)}] for any $f\in X$,
  $\|f\|_{X}=\|\,|f|\,\|_{X}$.
\end{enumerate}
\end{lemma}

Next, we recall the concept of the associate space of
a ball Banach function space
(see, for instance, \cite[p.\,9]{SHYY}).

\begin{definition}\label{assd}
For any ball Banach function space $X$, the
\emph{associate space\index{associate space}
\rm{(also called the
\emph{K\"{o}the dual}\index{K\"{o}the dual})}}
$X'$\index{$X'$} is defined to be the set of
all the $f\in\Msc(\rrn)$ such that
$$\|f\|_{X'}:=\sup\left\{\|fg\|_{L^1(\rrn)}:\
g\in X,\ \|g\|_{X}=1\right\}<\infty.$$
\end{definition}

\begin{remark}\label{remark127}
By \cite[Proposition 2.3]{SHYY}, we know that,
if $X$ is a ball Banach function space, then
its associate space $X'$
is also a ball Banach function space.
\end{remark}

The following concept on absolutely
continuous quasi-norms
of ball quasi-Banach function spaces
given in \cite[Definition
2.5]{HYYZ} is similar
to the absolute continuity of
Banach function spaces (see, for instance,
\cite[Chapter 1, Definition 3.1]{BSIO}).

\begin{definition}\label{absolute}
Let $X$ be a ball quasi-Banach function space.
A measurable function $f\in X$ is said to have an
\emph{absolutely continuous quasi-norm in
$X$\index{absolutely continuous quasi-norm}}
if, for any sequence $\{E_{i}\}_{i\in\mathbb{N}}$
of measurable sets satisfying $\1bf_{E_{i}}\to0$
almost everywhere as $i\to
\infty$, $\|f\1bf_{E_{i}}\|_{X}\to0$ as $i\to\infty$.
Moreover, $X$ is said to have an
\emph{absolutely continuous quasi-norm} if,
for any $f\in X$, $f$ has an absolutely
continuous quasi-norm in $X$.
\end{definition}

Now, we recall the concepts of both the
convexity and the concavity
of ball quasi-Banach function spaces
(see, for instance, \cite[Definition 2.6]{SHYY}).

\begin{definition}\label{convex}
Let $X$ be a ball quasi-Banach
function space and $p\in(0,\infty)$.
\begin{enumerate}
\item[(i)] The
\emph{$p$-convexification}\index{$p$-convexification}
$X^{p}$ of $X$ is defined by setting
$$X^{p}:=\left\{f\in\Msc(\rrn):\ |f|^{p}\in
X\right\}$$
equipped with the quasi-norm $\|f\|_{X^{p}}:=
\left\||f|^{p}\right\|^{1/p}_{X}.$
\item[(ii)] The space $X$ is said to be
\emph{$p$-convex} if there exists
a positive constant $C$
such that, for any $\{f_{j}\}_{j\in\mathbb{N}}
\subset X^{1/p}$,
$$\left\|\sum_{j\in\mathbb{N}}|f_{j}|
\right\|_{X^{1/p}}\leq C
\sum_{j\in\mathbb{N}}\|f_{j}\|_{X^{1/p}}.$$
In particular, when $C=1,\ X$ is said to be
\emph{strictly
$p$-convex}\index{$p$-convex}\index{strictly $p$-convex}.
\item[(iii)] The space $X$ is said to be
\emph{$p$-concave}\index{$p$-concave} if there exists
a positive constant $C$
such that, for any
$\{f_{j}\}_{j\in\mathbb{N}}\subset X^{1/p}$,
$$\sum_{j\in\mathbb{N}}\|f_{j}\|_{X^{1/p}}
\leq C\left\|\sum_{j\in\mathbb{N}}
|f_{j}|\right\|_{X^{1/p}}.$$
In particular, when $C=1,\
X$ is said to be
\emph{strictly
p-concave}\index{strictly p-concave}.
\end{enumerate}
\end{definition}

In what follows, let $L_{\text{loc}}^{1}
(\rrn)$\index{$L_{\text{loc}}^{1}(\rrn)$}
denote the set of all
locally integrable functions on $\rrn$.
Recall that the
\emph{Hardy--Littlewood maximal
operator}\index{Hardy--Littlewood maximal operator}
$\mc$\index{$\mc$} is defined by setting,
for any $f\in L_{\text{loc}}^{1}(\rrn)$
and $x\in\rrn$,
\begin{equation}\label{hlmax}
\mc (f)(x):=\sup_{B\ni x}\frac{1}{|B|}
\int_{B}|f(y)|\,dy,
\end{equation}
where the supremum is taken over all the
balls $B\in\mathbb{B}$ containing $x$.

For any $\theta\in(0,\infty)$,
the
\emph{powered Hardy--Littlewood maximal
operator}\index{powered Hardy--Littlewood
maximal operator} $\mc^{(\theta)}$\index{$\mc^{(\theta)}$}
is defined by setting, for any
$f\in L_{\text{loc}}^{1}(\rrn)$ and $x\in\rrn$,
\begin{equation}\label{hlmaxp}
\mc^{(\theta)}(f)(x):=
\left[\mc\left(|f|^{\theta}\right)(x)
\right]^{\frac{1}{\theta}}.
\end{equation}

We also need the following two assumptions about the
\emph{Fefferman--Stein vector-valued
inequality}\index{Fefferman--Stein vector-valued \\inequality}
on ball quasi-Banach function spaces
and the boundedness of powered Hardy--Littlewood
maximal operators on associate spaces.

\begin{assumption}\label{assfs}
Let $X$ be a ball
quasi-Banach function space and
$\theta,\ s\in(0,1]$.
Assume that there exists a positive constant
$C$ such that, for any
$\{f_{j}\}_{j\in\mathbb{N}}
\subset L^1_{{\rm loc}}(\rrn)$,
\begin{equation}\label{FSine}
\left\|\left\{\sum_{j\in\mathbb{N}}
\left[\mc^{(\theta)}
(f_{j})\right]^s\right\}^{1/s}
\right\|_{X}\leq C
\left\|\left(\sum_{j\in\mathbb{N}}
|f_{j}|^{s}\right)^{1/s}\right\|_{X}.
\end{equation}
\end{assumption}

\begin{remark}\label{power}
\begin{enumerate}
 \item[(i)] We claim that, for any given
 ball quasi-Banach function space
 $X$ and any given $r\in(0,\infty)$,
 the powered Hardy--Littlewood
 maximal operator $\mc^{(r)}$ is bounded
  on $X$ if and only if
  the centered Hardy--Littlewood maximal
  operator $\mc$ is bounded
  on $X^{\frac{1}{r}}$. Indeed,
  for any given $r\in(0,\infty)$, by the definition
 of the powered Hardy--Littlewood maximal operator
 $\mc^{(r)}$, we conclude that, for any
 $f\in L^1_{{\rm loc}}(\rrn)$,
\begin{align}\label{1111}
\left\|\mc^{(r)}(f)\right\|_{X}
=\left\|\left[\mc(|f|^{r})
\right]^{\frac{1}{r}}
\right\|_{X}
=\left\|\mc(|f|^{r})\right
\|_{X^{\frac{1}{r}}}^{\frac{1}{r}}.
\end{align}
In addition, from
Definition \ref{convex}(i), it follows that
$$
\|f\|_{X}=\||f|^{r}\|_{X^{\frac{1}{r}}}
^{\frac{1}{r}},
$$
which, combined with \eqref{1111},
implies that this claim holds true.
 \item[(ii)] Similarly
 to Remark \ref{power}(i) above,
 we know that the inequality \eqref{FSine}
 holds true if and only if, for any given
$\theta$, $s\in(0,1]$ and any $\{
f_{j}\}_{j\in\mathbb{N}}
\subset L^1_{{\rm loc}}(\rrn)$,
\begin{equation*}
\left\|\left\{
\sum_{j\in\mathbb{N}}
\left[\mc(f_{j})\right]^{\frac{s}{
\theta}}\right\}^{\frac{\theta}{s}}
\right\|_{X^{\frac{1}{\theta}}}
\lesssim\left\|\left(
\sum_{j\in\mathbb{N}}
|f_{j}|
^{\frac{s}{\theta}}\right)^
{\frac{\theta}{s}}\right\|
_{X^{\frac{1}{\theta}}},
\end{equation*}
where the implicit positive
constant is independent of
$\{f_{j}\}_{j\in\mathbb{N}}$.
\end{enumerate}
\end{remark}

\begin{assumption}\label{assas}
Let $X$ be a ball quasi-Banach function
space, $s\in(0,1]$, and
$r\in(1,\infty]$. Assume that $X^{1/s}$
is also a ball Banach function space
and there exists a positive
constant $C$ such that,
for any $f\in L^1_{\mathrm{loc}}(\rrn)$,
\begin{equation*}
\left\|\mc^{((r/s)')}(f)
\right\|_{(X^{1/s})'}\leq C
\left\|f\right\|_{(X^{1/s})'}.
\end{equation*}
\end{assumption}

We point out that,
for any $p$, $q\in(0,\infty)$
and $\omega\in M(\rp)$,
local and global generalized Herz spaces
are both quasi-normed linear spaces and,
in particular,
for any $p$, $q\in[1,\infty)$ and
$\omega\in M(\rp)$, they are
normed linear spaces.
Therefore, we next investigate the
completeness of these generalized Herz spaces.
For this purpose, we need the following lemma
which is a part of \cite[Theorem 2]{dfmn20}.

\begin{lemma}\label{lem}
Let $X\subset\Msc(\rrn)$ be a
quasi-normed linear space such that
Definition \ref{Df1}(ii) holds true. If,
for any positive increasing Cauchy sequence
$\{f_{n}\}_{n\in\mathbb{N}}$ of $X$,
the function $f:=\sup_{n\in\mathbb{N}}f_{n}$
belongs to $X$, then $X$ is complete.
\end{lemma}

Due to Lemma \ref{lem}, we obtain the following
conclusion about the completeness
of quasi-normed linear spaces.

\begin{proposition}\label{xcom}
Let $X\subset\Msc(\rrn)$ be a
quasi-normed linear space satisfying
both (ii) and (iii) of
Definition \ref{Df1}. Then $X$ is complete.
\end{proposition}

\begin{proof}
Let $X\subset\Msc(\rrn)$ be a quasi-normed
linear space satisfying both
(ii) and (iii) of Definition \ref{Df1}, and
$\{f_n\}_{n\in\mathbb{N}}$ be a positive
increasing Cauchy sequence of $X$.
Obviously, the sequence
$\{f_n\}_{n\in\mathbb{N}}$ is bounded on
$X$. This,
together with Definition \ref{Df1}(iii),
implies that
\begin{align*}
\left\|\sup_{n\in\mathbb{N}}f_n\right\|_{X}
=\left\|\lim\limits_{n\to\infty}f_n\right\|_{X}
=\lim\limits_{n\to\infty}\|f_n\|_{X}
\leq\sup_{n\in\mathbb{N}}\|f_n\|_{X}<\infty.
\end{align*}
Thus, the function $f:=\sup_{n\in\mathbb{N}}f_n$
belongs to $X$. Using this and Lemma \ref{lem},
we further conclude that $X$ is complete.
This finishes the proof of Proposition \ref{xcom}.
\end{proof}

The following completeness of
ball quasi-Banach function spaces
can be deduced immediately from
Proposition \ref{xcom}; we
omit the details.

\begin{corollary}
Let $X\subset\Msc(\rrn)$ be a
ball quasi-Banach function space.
Then $X$ is complete.
\end{corollary}

Via Proposition \ref{xcom}, we now show that
local and global generalized Herz spaces
are both complete and hence
are quasi-Banach spaces.

\begin{theorem}\label{hlcom}
Let $p$, $q\in(0,\infty)$ and $\omega\in M(\rp)$.
Then the local generalized Herz space
$\HerzSo$ is a quasi-Banach space.
\end{theorem}

\begin{proof}
Let $p$, $q\in(0,\infty)$ and $\omega\in M(\rp)$.
Obviously, $\HerzSo$ is a quasi-normed linear space.
Moreover, it is easy to show that the conditions
(i) and (ii) of Definition \ref{Df1} hold true and
Definition \ref{Df1}(iii) is a simple corollary of
the monotone convergence theorem.
This, combined with Proposition \ref{xcom},
further implies that
the local generalized Herz space
$\HerzSo$ is complete and hence finishes the proof of
Theorem \ref{hlcom}.
\end{proof}

\begin{theorem}\label{hgcom}
Let $p$, $q\in(0,\infty)$ and $\omega\in M(\rp)$.
Then the global generalized Herz space
$\HerzS$ is a quasi-Banach space.
\end{theorem}

\begin{proof}
Let $p$, $q\in(0,\infty)$ and $\omega\in M(\rp)$.
We now show that
the Herz space $\HerzS$ satisfies (i), (ii),
and (iii) of Definition \ref{Df1}.
Indeed, the global generalized Herz space
$\HerzS$ satisfies both (i) and (ii) of
Definition \ref{Df1} obviously.
Next, we prove that Definition \ref{Df1}(iii)
holds true for $\HerzS$.
To this end, assume that $f\in \Msc(\rrn)$ and
$\{f_{m}\}_{m\in\mathbb{N}}\subset
\Msc(\rrn)$ satisfies that
$0\leq f_{m}\uparrow f$ almost
everywhere as $m\to\infty$. Then,
for any given $\alpha\in(0,\|f\|_{\HerzS})$,
by the definition of $\|\cdot\|_{\HerzS}$,
we find that there exists a
$\xi\in\rrn$ such that
$$\|f(\cdot+\xi)\|_{\HerzSo}>\alpha.$$
Applying this and the monotone convergence theorem,
we find that there exists
an $N\in\mathbb{N}$ such that,
for any $m\in\mathbb{N}\cap(N,\infty)$,
$$\left\|f_{m}(\cdot+\xi)
\right\|_{\HerzSo}>\alpha,$$
which, together with the
definition of $\|\cdot\|_{\HerzS}$,
further
implies that, for any $m\in\mathbb{N}\cap(N,\infty)$,
\begin{align*}
\alpha<\sup_{\xi\in\rrn}\left\|
f_{m}(\cdot+\xi)\right\|_{\HerzSo}
=\|f_{m}\|_{\HerzS}.
\end{align*}
Therefore, we conclude that
$$\|f_{m}\|_{\HerzS}
\uparrow\|f\|_{\HerzS}$$
as $m\to\infty$
and hence Definition \ref{Df1}(iii) holds true
for $\HerzS$. From this and
Proposition \ref{xcom}, we deduce that
the global generalized Herz space
$\HerzS$ is complete.
This finishes the proof of Theorem \ref{hgcom}.
\end{proof}

However, the following two examples show that
local and global
generalized Herz spaces may not be
a ball quasi-Banach function space.
Namely, they may not satisfy
Definition \ref{Df1}(iv).

\begin{example}\label{Ex1}
Let $p,\ q\in(0,\infty),\
\alpha\in(-\infty,-\frac{n}{p}]$, and
$\omega(t):=t^{\alpha}$ for any $t\in(0,\infty)$.
Then $\1bf_{B(\0bf,1)}\notin\HerzSo$ and
$\1bf_{B(\0bf,1)}\notin\HerzS$,
which implies that the local generalized
Herz space $\HerzSo$ and the
global generalized Herz space $\HerzS$
are not ball quasi-Banach function spaces.
\end{example}

\begin{proof}
Let $p$, $q\in(0,\infty)$ and,
for any $t\in(0,\infty)$, $\omega(t):=t^\alpha$
with some fixed $\alpha\in(-\infty,-\frac{n}{p}]$.
Then, from the definition
of $\|\cdot\|_{\HerzSo}$, it follows that
\begin{align*}
\left\|\1bf_{B(\0bf,1)}\right\|_{\HerzSo}&=
\left[\sum_{k\in\mathbb{Z}}2^{k\alpha q}
\left\|\1bf_{B(\0bf,1)}\1bf_{B(\0bf,2^{k})
\setminus B(\0bf,2^{k-1})}\right\|
_{L^{p}(\rrn)}^{q}\right]^{\frac{1}{q}}\\
&=\left[\sum_{k\in\zmn}2^{k\alpha q}
\left\|\1bf_{B(\0bf,2^{k})\setminus
B(\0bf,2^{k-1})}\right\|
_{L^{p}(\rrn)}^{q}\right]^{\frac{1}{q}}\\
&\sim\left[\sum_{k\in\zmn}2^{kq
(\alpha+\frac{n}{p})}\right]^{\frac{1}{q}}.
\end{align*}
Applying this and the assumption
$\alpha\in(-\infty,-\frac{n}{p}]$,
we conclude that
\begin{equation}\label{ifz}
\left\|\1bf_{B(\0bf,1)}\right\|_{\HerzSo}=\infty,
\end{equation}
which implies that $\1bf_{B(\0bf,1)}
\notin\HerzSo$. Thus,
the local generalized Herz space
$\HerzSo$ is not a ball
quasi-Banach function space.

In addition, by \eqref{ifz} and
the definition of $\|\cdot\|_{\HerzS}$,
we find that
\begin{align*}
\left\|\1bf_{B(\0bf,1)}\right\|_{\HerzS}\geq
\left\|\1bf_{B(\0bf,1)}\right\|_{\HerzSo}=\infty.
\end{align*}
Therefore, we conclude that
$\1bf_{B}\notin\HerzS$ and hence
the global generalized
Herz space $\HerzS$ is also not
a ball quasi-Banach function space.
This finishes the proof of Example \ref{Ex1}.
\end{proof}

\begin{example}\label{Ex2}
Let $p$, $q\in(0,\infty)$, $\alpha_1\in\rr$,
$\alpha_2\in[0,\infty)$, and
$$
\omega(t):=\left\{\begin{aligned}
&t^{\alpha_{1}}(1-\ln t)
\hspace{0.3cm}\text{when }t\in(0,1],\\
&t^{\alpha_{2}}(1+\ln t)
\hspace{0.3cm}\text{when }t\in(1,\infty).
\end{aligned}\right.
$$
Then $\1bf_{B(\0bf,1)}\notin\HerzS$
and hence, in this case, the
global generalized Herz space
$\HerzS$ is not a
ball quasi-Banach function space.
\end{example}

\begin{proof}
Let all the symbols be as in the present example
and, for any $k\in\mathbb{N}$, let $\xi_k\in\rrn$
satisfy $|\xi_k|=2^{k}+1$. Then we claim that,
for any $k\in\mathbb{N}$,
$$
B(\0bf,1)\subset B(\xi_k,2^{k+1})\setminus
B(\xi_{k},\tk).
$$
Indeed, for any $y\in B(\0bf,1)$, we have
\begin{align*}
|y-\xi_k|\leq|y|+|\xi_k|<1+|\xi_k|\leq2^{k+1}
\end{align*}
and
\begin{align*}
|y-\xi_k|\geq|\xi_k|-|y|>|\xi_k|-1=\tk.
\end{align*}
These imply that $y\in B(\xi_k,2^{k+1})\setminus
B(\xi_{k},\tk)$ and hence
the above claim holds true.
From this and the definition of
$\|\cdot\|_{\HerzS}$, it follows that
\begin{align*}
\left\|\1bf_{B(\0bf,1)}\right\|_{\HerzS}
&\geq\left\|\1bf_{B(\0bf,1)}
(\cdot+\xi_k)\right\|_{\HerzSo}\\
&=\omega(2^{k+1})\left\|
\1bf_{B(\0bf,1)}\right\|_{L^p(\rrn)}\\
&\sim2^{(k+1)\alpha_2}
\left[1+(k+1)\ln2\right]\to\infty
\end{align*}
as $k\to\infty$. Thus,
$\1bf_{B(\0bf,1)}\notin\HerzS$.
This further implies that $\HerzS$
is not a ball quasi-Banach
function space and hence finishes
the proof of Example \ref{Ex2}.
\end{proof}

The following two theorems show that
local and global generalized Herz
spaces are ball
quasi-Banach function spaces under
some additional and sharp assumptions
of exponent $\omega\in M(\rp)$.

\begin{theorem}\label{Th3}
Let $p,\ q\in(0,\infty)$ and $\omega\in M(\rp)$ with
$\m0(\omega)\in(-\frac{n}{p},\infty)$.
Then the local generalized Herz space $\HerzSo$
is a ball quasi-Banach function space.
\end{theorem}

\begin{proof}
Let $p$, $q\in(0,\infty)$ and
$\omega\in M(\rp)$ with
$\m0(\omega)\in(-\frac{n}{p},
\infty)$.
Then, by the proof of Theorem \ref{hlcom},
we know that the local generalized Herz
space $\HerzSo$ is a quasi-Banach space satisfying
(i), (ii), and (iii) of Definition \ref{Df1}.
Therefore, to prove the present theorem,
it suffices to show that
Definition \ref{Df1}(iv)
holds true for $\HerzSo$.

For this purpose, let $B(x_{0},r)\in\mathbb{B}$
with $x_{0}\in\rrn$ and
$r\in(0,\infty)$. Then we have
\begin{align}\label{esb1}
&\left\|\1bf_{B(x_{0},r)}\right\|_
{\HerzSo}^{q}\notag\\&\quad=
\sum_{k\in\mathbb{Z}}
\left[\omega(2^{k})\right]^{q}
\left\|\1bf_{B(x_{0},r)}
\1bf_{B(\0bf,2^{k})\setminus
B(\0bf,2^{k-1})}\right
\|_{L^{p}(\rrn)}^{q}\notag\\&\quad=
\sum_{k\in\zmn}\left[
\omega(2^{k})\right]^{q}
\left\|\1bf_{B(x_{0},r)}\1bf_{B(\0bf,2^{k})
\setminus B(\0bf,2^{k-1})}
\right\|_{L^{p}(\rrn)}^{q}
+\sum_{k\in\mathbb{N}}\cdots\notag\\
&\quad=:\rm{I_{1}}+{\rm I_{2}}.
\end{align}

For ${\rm I_{1}}$,
from Lemma \ref{Th1}, it follows that,
for any $k\in\zmn$,
$
\omega(2^{k})\lesssim2^{k[\m0(\omega)-\eps]},
$
where $\eps\in(0,\m0(\omega)+\frac{n}{p})$
is a fixed positive constant.
By this and the assumption $\m0(\omega)+
\frac{n}{p}-\eps\in(0,\infty)$, we conclude that
\begin{align}\label{EsoI}
{\rm I_{1}}
&\lesssim\sum_{k\in\zmn}2^{kq[\m0(\omega)-\eps]}
\left\|\1bf_{B(x_{0},r)}\1bf_{B(\xi,2^{k})
\setminus B(\xi,2^{k-1})}\right
\|_{L^{p}(\rrn)}^{q}\notag\\
&\lesssim\sum_{k\in\zmn}2^{kq[\m0(\omega)
+\frac{n}{p}-\eps]}<\infty.
\end{align}

Now, we deal with ${\rm I_{2}}$. To this end,
we first claim that, for any
$k\in\mathbb{N}\cap(\ln(r+|x_{0}|)/
\ln2+1,\infty)$,
$$B(x_{0},r)\cap\left[B(\0bf,\tk)\setminus
B(\0bf,\tkm)\right]=\emptyset.$$
Indeed, from the assumption $k>\ln(r+|x_{0}|)/\ln2+1$,
we deduce that $2^{k-1}>r+|x_{0}|$.
By this, we find that,
for any $x\in B(x_{0},r)$,
$$|x|
\leq|x-x_{0}|+|x_{0}|<r+|x_{0}|<\tkm,$$
which implies that $x\in B(\0bf,\tkm)$
and hence $B(x_{0},r)\subset B(\0bf,\tkm)$.
Therefore, in this case,
$$B(x_{0},r)\cap\left[B(\0bf,\tk)\setminus
B(\0bf,\tkm)\right]=\emptyset,$$
which completes the proof of the above claim.
Applying this claim, we conclude that
\begin{align*}
{\rm I_{2}}
=\sum_{k\in\mathbb{N}\cap[1,\frac{\ln(r+|x_{0}|)}
{\ln2}+1]}\left[\omega(2^{k})
\right]^{q}\left\|\1bf_{B(x_{0},r)}
\1bf_{B(\0bf,2^{k})
\setminus B(\0bf,2^{k-1})}\right\|
_{L^{p}(\rrn)}^{q}<\infty,
\end{align*}
which, together with \eqref{esb1} and
\eqref{EsoI}, implies that
$\1bf_{B(x_{0},r)}\in\HerzSo$. Thus,
Definition \ref{Df1}(iv) holds true for
$\HerzSo$,
which completes the proof of Theorem \ref{Th3}.
\end{proof}

\begin{remark}
Combining Examples \ref{ex113}
and \ref{Ex1}, we find that the assumption
$\m0\in(-\frac{n}{p},\infty)$ in
Theorem \ref{Th3} is sharp.
\end{remark}

\begin{theorem}\label{Th2}
Let $p,\ q\in(0,\infty)$ and $\omega\in M(\rp)$
with $\m0(\omega)\in(-\frac{n}{p},
\infty)$ and $\MI(\omega)\in(-\infty,0)$.
Then the global generalized Herz space $\HerzS$
is a ball quasi-Banach function space.
\end{theorem}

\begin{proof}
Let $p$, $q\in(0,\infty)$ and $\omega\in M(\rp)$
with $\m0(\omega)\in(-\frac{n}{p},\infty)$
and $\MI(\omega)\in(-\infty,0)$.
Then, from the proof of Theorem \ref{hgcom},
it follows that the
global generalized Herz space $\HerzS$
is a quasi-Banach space satisfying
(i), (ii), and (iii) of Definition \ref{Df1}.

Next, we prove that $\HerzS$ satisfies
Definition \ref{Df1}(iv).
For this purpose, let
$B(x_{0},r)\in\mathbb{B}$ with
$x_{0}\in\rrn$ and $r\in(0,\infty)$.
In addition, by Lemma \ref{Th1},
we conclude that, for any
$k\in\zmn$,
$$\omega(\tk)\lesssim2^{k[\m0(\omega)-\eps]}$$
and, for any $k\in\mathbb{N}$,
$$\omega(2^{k})\lesssim2^{k[\MI(\omega)+\eps]},$$
where $\eps\in(0,\min\{\m0(\omega)+\frac{n}{p},
-\MI(\omega)\})$
is a fixed positive constant.
Then, for any $\xi\in\rrn$, we have
\begin{align}\label{eso3}
&\left\|\1bf_{B(x_{0},r)}(\cdot+\xi)
\right\|_{\HerzSo}^{q}\notag\\
&\quad=\sum_{k\in\mathbb{Z}}
\left[\omega(2^{k})\right]^{q}
\left\|\1bf_{B(x_{0},r)}\1bf_{B(\xi,2^{k})
\setminus B(\xi,2^{k-1})}
\right\|_{L^{p}(\rrn)}^{q}\notag\notag\\&\quad=
\sum_{k\in\zmn}\left[\omega(2^{k})\right]^{q}
\left\|\1bf_{B(x_{0},r)}\1bf_{B(\xi,2^{k})
\setminus B(\xi,2^{k-1})}
\right\|_{L^{p}(\rrn)}^{q}+
\sum_{k\in\mathbb{N}}\cdots\notag\\
&\quad=:{\rm I_{\xi,1}}+{\rm I_{\xi,2}}.
\end{align}
We first deal with ${\rm I_{\xi,1}}$. Indeed,
applying an argument similar to that
used in the estimation
of \eqref{EsoI}, we find that,
for any $\xi\in\rrn$,
\begin{align}\label{esoi}
{\rm I_{\xi,1}}\lesssim\sum_{k\in\zmn}
2^{kq[\m0(\omega)+\frac{n}{p}-\eps]}<\infty.
\end{align}
For ${\rm I_{\xi,2}}$, from the assumption
$\omega(2^{k})\lesssim2^{k[\MI(\omega)+\eps]}$
for any $k\in\mathbb{Z}$,
and the assumption $\MI(\omega)+\eps\in(-\infty,0)$,
it follows that, for any $\xi\in\rrn$,
\begin{align*}
{\rm I_{\xi,2}}&=\sum_{k\in\mathbb{N}}
\left[\omega(2^{k})\right]^{q}
\left\|\1bf_{B(x_{0},r)}
\1bf_{B(\xi,2^{k})\setminus
B(\xi,2^{k-1})}\right\|_{L^{p}(\rrn)}^{q}\notag\\
&\leq\sum_{k\in\mathbb{N}}
\left[\omega(2^{k})\right]^{q}
\left\|\1bf_{B(x_{0},r)}
\right\|_{L^{p}(\rrn)}^{q}\\
&\lesssim\sum_{k\in\mathbb{N}}
2^{kq[\MI(\omega)+\eps]}<\infty,
\end{align*}
which, combined with \eqref{eso3}
and \eqref{esoi}, further implies
that, for any $\xi\in\rrn$,
\begin{align*}
&\left\|\1bf_{B(x_{0},r)}(\cdot+\xi)
\right\|_{\HerzSo}\\
&\quad\lesssim
\left\{\sum_{k\in\zmn}2^{
kq[\m0(\omega)+\frac{n}{p}-\eps]}
+\sum_{k\in\mathbb{N}}2^{kq[
\MI(\omega)+\eps]}
\right\}^{\frac{1}{q}}<\infty.
\end{align*}
This implies that $\|\1bf_{B(x_{0},r)}
\|_{\HerzS}<\infty$ and hence
$\1bf_{B(x_{0},r)}\in\HerzS$. Thus,
Definition \ref{Df1}(iv) also holds true,
which completes the proof of Theorem \ref{Th2}.
\end{proof}

\begin{remark}\label{remark1221}
Applying Examples \ref{ex113}, \ref{Ex1},
\ref{ex117}, and \ref{Ex2},
we conclude
that the assumptions $\m0(\omega)
\in(-\frac{n}{p},\infty)$
and $\MI(\omega)\in(-\infty,0)$
in Theorem \ref{Th2} are sharp.
\end{remark}

Moreover, by the following conclusion,
we know that
the local Herz spaces are
ball Banach function spaces
when $p$, $q\in[1,\infty)$ and
$\omega\in M(\rp)$ satisfy
some reasonable and sharp assumptions.

\begin{theorem}\label{balll}
Let $p,\ q\in[1,\infty)$ and $\omega\in M(\rp)$
satisfy $$-\frac{n}{p}<\m0(\omega)
\leq\M0(\omega)<\frac{n}{p'},$$
where $\frac{1}{p}+\frac{1}{p'}=1$.
Then the local generalized Herz space
$\HerzSo$ is a ball Banach function space.
\end{theorem}

\begin{proof}
Let all the symbols be as in the present theorem.
Then, applying Theorem \ref{Th3},
we know that the local generalized
Herz space $\HerzSo$ is a
ball quasi-Banach function
space. In addition, note that
the quasi-norm
$\|\cdot\|_{\HerzSo}$ satisfies the
triangle inequality when $p$, $q\in[1,\infty)$.
Thus, to complete the whole proof,
it suffices to show that, for any
$B\in\mathbb{B}$,
\eqref{bal} holds true with $X$
replaced by $\HerzSo$.
To this end, we claim that,
for any $k\in\mathbb{Z}$, \eqref{bal} holds true
with $X$ replaced by $\HerzSo$ and $B:=B(\0bf,\tk)$.
Assume that this claim holds true for the moment.
Then, for any $B\in\mathbb{B}$,
there exists a $k\in\mathbb{Z}$ such that
$B\subset B(\0bf,\tk)$. This,
together with the above claim,
implies that, for any $f\in\HerzSo$,
\begin{align*}
\int_{B}|f(y)|\,dy\leq\int_{
B(\0bf,\tk)}|f(y)|\,dy\leq
C_{(B(\0bf,\tk))}\|f\|_{\HerzSo},
\end{align*}
which further implies that \eqref{bal}
holds true for any $B\in\mathbb{B}$
and hence completes
the proof of Theorem \ref{balll}.

Thus, to finish the whole proof,
it remains to show the
above claim.
For this purpose, let
$B(\0bf,2^{k_{0}})\in\mathbb{B}$
with some $k_{0}\in\mathbb{Z}$.
We now prove this claim
by considering the following four cases
on $p$, $q$.

\emph{Case 1)} $p,\ q\in(1,\infty)$.
In this case, by Lemma \ref{rela}, we find
that
$$\m0(1/\omega)=-\M0(\omega)>-\frac{n}{p'}.$$
This, together with Theorem \ref{Th3},
implies that the local generalized Herz
space $\HerzSoa$ is a
ball quasi-Banach function space.
Thus, applying the
H\"{o}lder inequality\index{H\"{o}lder inequality}
and the fact
that $\1bf_{B(\0bf,2^{k_{0}})}\in\HerzSoa$,
we find that, for any $f\in\HerzSo$,
\begin{align*}
&\int_{B(\0bf,2^{
k_{0}})}|f(y)|\,dy\\&
\quad=\int_{\rrn}\left|f(y)
\1bf_{B(\0bf,2^{k_{0}})}(y)\right|\,dy\\
&\quad=\sum_{k\in\mathbb{Z}}
\int_{B(\0bf,\tk)\setminus
B(\0bf,\tkm)}\left|f(y)\1bf_{
B(\0bf,2^{k_{0}})}(y)\right|\,dy\\
&\quad\leq\sum_{k\in\mathbb{Z}}
\left\|f\1bf_{B(\0bf,\tk)\setminus
B(\0bf,\tkm)}\right\|_{L^{p}(\rrn)}
\left\|\1bf_{B(\0bf,2^{k_{0}})}\1bf_{B(\0bf,
\tk)\setminus B(\0bf,\tkm)}
\right\|_{L^{p'}(\rrn)}\\
&\quad\leq\left\{\sum_{k\in\mathbb{Z}}
\left[\omega(2^{k})\right]^{q}
\left\|f\1bf_{B(\0bf,\tk)\setminus
B(\0bf,\tkm)}\right\|_{L^{p}(\rrn)}
^{q}\right\}^{\frac{1}{q}}\\
&\quad\quad\times\left\{\sum_{k\in\mathbb{Z}}
\left[\omega(2^{k})\right]^{-q'}
\left\|\1bf_{B(\0bf,2^{k_{0}
})}\1bf_{B(\0bf,\tk)\setminus
B(\0bf,\tkm)}\right\|_{L^{p'}(\rrn)}
^{q'}\right\}^{\frac{1}{q'}}\\
&\quad=\|f\|_{\HerzSo}\|
\1bf_{B(\0bf,2^{k_{0}})}\|_{\HerzSoa},
\end{align*}
which further
implies that, in this case, \eqref{bal}
holds true with the positive constant
$$C_{(B)}:=\left\|
\1bf_{B(\0bf,2^{k_{0}})}\right\|_{\HerzSoa}.$$

\emph{Case 2)} $p\in(1,\infty)$ and $q=1$.
In this case,
we first show that,
for any $k\in\mathbb{Z}\cap(-\infty,k_{0}]$,
\begin{equation}\label{cl1}
\left\|\1bf_{B(\0bf,\tk)\setminus
B(\0bf,\tkm)}\right\|_{L^{p'}(\rrn)}
\lesssim\omega(\tk).
\end{equation}
We prove \eqref{cl1} by considering,
respectively, $k_{0}\in\zmn$ and $k_{0}
\in\mathbb{N}$.
If $k_{0}\in\zmn$,
then, by Lemma \ref{Th1},
we conclude that, for any
$k\in\zmn$,
$$\omega(\tk)
\gtrsim2^{k[\M0(\omega)+\eps]},$$
where $\eps\in(0,\frac{n}{p'}-\M0(\omega))$
is a fixed positive constant.
From this and the assumption
$\frac{n}{p'}-\M0(\omega)-\eps\in(0,\infty)$,
we deduce that,
for any $k\in\mathbb{Z}\cap(-\infty,k_{0}]$,
\begin{align}\label{q12}
&\left\|\1bf_{B(\0bf,\tk)\setminus
B(\0bf,\tkm)}\right\|_{L^{p'}(\rrn)}\notag\\
&\quad\sim2^{\frac{nk}{p'}}
\sim2^{k[\M0(\omega)+\eps]}
2^{k[\frac{n}{p'}-\M0(\omega)-\eps]}\notag\\
&\quad\lesssim\omega(\tk)2^{k_{0}
[\frac{n}{p'}-\M0(\omega)-\eps]}
\sim\omega(\tk),
\end{align}
which implies that \eqref{cl1} holds
true for any $k\in\mathbb{Z}
\cap(-\infty,k_{0}]$ with $k_{0}\in\zmn$.
If $k_{0}\in\mathbb{N}$, let $$C_{(k_{0})}:=
\min\left\{\omega(2),\ldots,
\omega(2^{k_{0}})\right\}.$$ Then we find that,
for any $k\in\mathbb{N}\cap(-\infty,k_{0}]$,
\begin{align*}
\left\|\1bf_{B(\0bf,\tk)\setminus
B(\0bf,\tkm)}\right\|_{L^{p'})
(\rrn)}\sim2^{\frac{nk}{p'}}
\lesssim2^{\frac{nk_{0}}{p'}}
\lesssim2^{\frac{nk_{0}}{p'}}
C_{(k_{0})}^{-1}\omega(\tk).
\end{align*}
By this and an argument similar to that used
in the estimation of \eqref{q12},
we conclude that \eqref{cl1} holds
true for any $k\in\mathbb{Z}\cap(-\infty,k_{0}]$
with $k_{0}\in\mathbb{N}$,
which completes the proof of \eqref{cl1}.
Using this and the H\"{o}lder inequality,
we find that, for any $f\in\HerzSo$,
\begin{align*}
&\int_{B(\0bf,2^{k_{0}})}|f(y)|\,dy\\
&\quad=\sum_{k=-\infty}^{k_{0}}
\int_{B(\0bf,\tk)\setminus
B(\0bf,\tkm)}|f(y)|\,dy\\
&\quad\leq\sum_{k=-\infty}^{k_{0}}
\left\|f\1bf_{B(\0bf,\tk)\setminus
B(\0bf,\tkm)}\right\|_{L^p(\rrn)}
\left\|\1bf_{B(\0bf,\tk)\setminus
B(\0bf,\tkm)}\right\|_{L^{p'}(\rrn)}\\
&\quad\lesssim\sum_{k\in\mathbb{Z}}
\omega(\tk)\left\|f\1bf_{B(\0bf,\tk)\setminus
B(\0bf,\tkm)}\right\|_{L^{p}(\rrn)}
\sim\|f\|_{\Kmp_{\omega,\0bf}^{p,1}(\rrn)},
\end{align*}
which further implies that,
in this case, \eqref{bal} holds true.

\emph{Case 3)} $p=1$ and $q=1$. In this case,
we first claim that, for any
$k\in\mathbb{Z}\cap(-\infty,k_{0}]$,
$1\lesssim\omega(\tk)$. Indeed,
applying Lemma \ref{Th1}, we conclude that,
for any $k\in\zmn$,
$$\omega(\tk)\gtrsim2^{k[\M0(\omega)+\eps]},$$
where $\eps\in(0,-\M0(\omega))$ is
a fixed positive constant.
Now, we prove this claim by
considering, respectively,
$k_{0}\in\zmn$ and $k_{0}\in\mathbb{N}$.
If $k_{0}\in\zmn$, then, from both the assumption that,
for any $k\in\zmn$,
$\omega(\tk)\gtrsim2^{k[\M0(
\omega)+\eps]}$ and the
assumption $\M0(\omega)+\eps\in(-\infty,0)$,
it follows that, for any
$k\in\mathbb{Z}\cap(-\infty,k_{0}]$,
\begin{align}\label{eq33}
1\sim2^{k_{0}[\M0(\omega)+\eps]}\lesssim
2^{k[\M0(\omega)+\eps]}\lesssim\omega(\tk).
\end{align}
This implies that the above
claim holds true when $k_{0}\in\zmn$.
If $k_{0}\in\mathbb{N}$, by the fact
that $\min\{\omega(2),\ldots,
\omega(2^{k_{0}})\}\leq\omega(\tk)$
for any $k\in\mathbb{N}\cap(-\infty,k_{0}]$,
and an argument similar to that
used in the estimation of \eqref{eq33},
we find that this claim holds true
under the assumption $k_{0}\in\mathbb{N}$.
Thus, this claim is proved.
Therefore, for any $f\in\HerzSo$,
we have
\begin{align*}
\int_{B(\0bf,2^{k_{0}})}|f(y)|
\,dy=\sum_{k=-\infty}^{k_{0}}
\left\|f\1bf_{B(\0bf,\tk)\setminus
B(\0bf,\tkm)}\right\|_{L^{1}(\rrn)}\lesssim
\|f\|_{\Kmp_{\omega,\0bf}^{1,1}(\rrn)},
\end{align*}
which further implies that, in this case,
\eqref{bal} holds true.

\emph{Case 4)} $p=1$ and $q\in(1,\infty)$.
In this case, we first show that
\begin{equation}\label{eq44}
\sum_{k=-\infty}^{
k_{0}}\left[\omega(\tk)\right]^{-q'}<\infty.
\end{equation}
Indeed, by Lemma \ref{Th1}, we find that,
for any $k\in\zmn$,
$$\omega(\tk)
\gtrsim2^{k[\M0(\omega)+\eps]},$$
where $\eps\in(0,-\M0(\omega))$ is
a fixed positive constant.
We now prove \eqref{eq44} by considering, respectively,
$k_{0}\in\zmn$ and $k_{0}\in\mathbb{N}$.
If $k_{0}\in\zmn$. Then, applying both the assumption that,
for any $k\in\zmn$, $\omega(\tk)\gtrsim2^{k[\M0(
\omega)+\eps]}$ and the assumption
$\M0(\omega)+\eps\in(-\infty,0)$,
we conclude that
\begin{align*}
\sum_{k=-\infty}^{k_{0}}\left[\omega(\tk)\right]^{-q'}
\lesssim\sum_{k=-\infty}^{
k_{0}}2^{-kq'[\M0(\omega)+\eps]}<\infty,
\end{align*}
which completes \eqref{eq44}
under the assumption $k_{0}\in\zmn$.
If $k_{0}\in\mathbb{N}$, then,
using both the assumption that,
for any $k\in\zmn$, $\omega(\tk)\gtrsim2^{k[\M0(
\omega)+\eps]}$ and the assumption
$\M0(\omega)+\eps\in(-\infty,0)$,
we find that
\begin{align*}
\sum_{k=-\infty}^{k_{0}}
\left[\omega(\tk)\right]^{-q'}\lesssim
\sum_{k\in\zmn}2^{-kq'[\M0(\omega)+\eps]}+
\sum_{k=1}^{k_{0}}
\left[\omega(\tk)\right]^{-q'}<\infty,
\end{align*}
which completes the proof of \eqref{eq44}
under the assumption $k_{0}\in\mathbb{N}$,
and hence \eqref{eq44} holds true
for any $k_0\in\mathbb{Z}$.
Thus, from the H\"{o}lder inequality and \eqref{eq44},
we deduce that, for any $f\in\HerzSo$,
\begin{align*}
\int_{B(\0bf,2^{k_{0}})}|f(y)|\,dy
&=\sum_{k=-\infty}^{k_{0}}
\left\|f\1bf_{B(\0bf,\tk)\setminus
B(\0bf,\tkm)}\right\|_{L^{1}(\rrn)}\\
&\leq\left\{\sum_{k=-\infty}^{k_{0}}
\left[\omega(\tk)\right]^{q}
\left\|f\1bf_{B(\0bf,\tk)\setminus B(\0bf,\tkm)}
\right\|_{L^{1}(\rrn)}^{q}\right\}^{\frac{1}{q}}\\
&\quad\times\left\{
\sum_{k=-\infty}^{k_{0}}
\left[\omega(\tk)\right]^{-q'}\right\}
^{\frac{1}{q'}}\\
&\lesssim\|f\|_{\Kmp_{\omega,\0bf}^{1,q}(\rrn)},
\end{align*}
which further implies that,
in this case, \eqref{bal} holds true.
This finishes the proof of the above claim,
and hence of Theorem \ref{balll}.
\end{proof}

The following examples show that
the assumptions in
Theorem \ref{balll} are sharp.

\begin{example}\label{ex1223}
Let $p$, $q\in[1,\infty)$ and
$\alpha\in[-\frac{n}{p'},\infty)$,
where $\frac{1}{p}+\frac{1}{p'}=1$.
For any $t\in(0,\infty)$,
let $\omega(t):=t^\alpha$.
\begin{enumerate}
  \item[{\rm(i)}] If $\alpha\in(\frac{n}{p'},\infty)$,
  then the local generalized Herz space $\HerzSo$ is not
  a ball Banach function space.
  \item[{\rm(ii)}] If $p=q=1$ and $\alpha=\frac{n}{p'}=0$,
  then the local generalized
  Herz space $\HerzSo$ coincides with
  $L^1(\rrn)$ and hence is a
  ball Banach function space.
  \item[{\rm(iii)}] If $p=q>1$ and
  $\alpha=\frac{n}{p'}$, then the local
  generalized Herz space $\HerzSo$
  is not a ball Banach function space.
\end{enumerate}
\end{example}

\begin{proof}
Let $p$, $q\in[1,\infty)$ and $\omega(t):=t^\alpha$
for any $t\in(0,\infty)$ and any
given $\alpha\in[-\frac{n}{p'},\infty)$,
where $\frac{1}{p}+\frac{1}{p'}=1$.
We first show (i). Indeed, for any $x\in\rrn$, let
$$
f(x):=\frac{1}{|x|^n}\1bf_{B(\0bf,1)\setminus\{\0bf\}}.
$$
Then, for any $k\in\zmn$ and $x\in
B(\xi,\tk)\setminus B(\xi,\tkm)$,
we have $f(x)\sim2^{-nk}$.
Using this and the assumption
$\alpha\in(\frac{n}{p'},\infty)$, we conclude
that
\begin{align*}
\|f\|_{\HerzSo}&=\left[\sum_{k\in\mathbb{Z}}
2^{k\alpha q}\left\|f\1bf_{B(\0bf,\tk)\setminus
B(\0bf,\tkm)}\right\|_{L^p(\rrn)}^
q\right]^{\frac{1}{q}}\\
&\sim\left[\sum_{k\in\zmn}
2^{k\alpha q}2^{-nkq}
\left|B(\0bf,\tk)\setminus
B(\0bf,\tkm)\right|^{\frac{q}{p}}\right]
^{\frac{1}{q}}\\
&\sim\left[\sum_{k\in\zmn}2
^{kq(\alpha-\frac{n}{p'})}\right]^{\frac{1}{q}}
\sim1,
\end{align*}
which implies that $f\in\HerzSo$.
However, we have
\begin{align*}
\int_{B(\0bf,1)}|f(y)|\,dy
=\int_{B(\0bf,1)}\frac{1}{|y|^n}\,dy
=\int_{0}^{1}r^{-1}\,dr
\int_{\mathbb{S}^{n-1}}\,d\sigma(x)=\infty.
\end{align*}
Therefore, in this case,
the local generalized Herz space $\HerzSo$ does not
satisfy \eqref{bal}
and hence $\HerzSo$ is not a
ball Banach function space.
This implies that (i) holds true.

We next prove (ii). Indeed, from the assumption $p=1$,
it follows that $\alpha=\frac{n}{p'}=0$.
This further implies that,
for any $f\in\Msc(\rrn)$,
\begin{align*}
\|f\|_{\HerzSo}=\sum_{k\in\mathbb{Z}}
\left\|f\1bf_{B(\0bf,\tk)\setminus
B(\xi,\tkm)}\right\|_
{L^1(\rrn)}=\|f\|_{L^1(\rrn)}.
\end{align*}
Thus, in this case, it holds true
that $\HerzSo=L^1(\rrn)$
and hence $\HerzSo$ is a
ball Banach function space,
which completes the proof of (ii).

Now, we show (iii). To this end, for
any given $\beta\in(0,n)$ and for
any $x\in\rrn$,
let
$$
f_\beta(x):=\frac{1}{|x|^\beta}\1bf_{B(\0bf,1)
\setminus\{\0bf\}}.
$$
Then, for any given $\beta\in(0,n)$
and for any $k\in\zmn$
and $x\in B(\0bf,\tk)\setminus B(\0bf,\tkm)$,
we have
$|x|\sim2^{\frac{np}{p'}}$. Using this,
we find that, for any $\beta\in(0,n)$,
\begin{align}\label{ex1223e1}
\left\|f_\beta\right\|_{\HerzSo}
&=\left[\sum_{k\in\mathbb{Z}}2^{k\alpha q}
\left\|f_\beta
\1bf_{B(\0bf,\tk)\setminus B(\0bf,\tkm)}
\right\|_{L^p}^p\right]^{\frac{1}{p}}\notag\\
&\sim\left[\sum_{k\in\zmn}|x|^{\frac{np}{p'}}
\int_{B(\0bf,\tk)\setminus B(\0bf,\tkm)}|
f_\beta(y)|^p\,dy\right]^{\frac{1}{p}}\notag\\
&\sim\left[\int_{B(\0bf,1)}\left|x\right|
^{p(\frac{n}{p'}-\beta)}
\,dy\right]^{\frac{1}{p}}\notag\\
&\sim\left[\int_{0}^{1}r^{p(n-\beta)-1}
\,dr\int_{\mathbb{S}^{n-1}}
\,d\sigma(x)\right]^{\frac{1}{p}}\notag\\
&\sim\frac{1}{(n-\beta)^{\frac{1}{p}}}.
\end{align}
On another hand, for any $\beta\in(0,n)$,
we have
\begin{align*}
\int_{B(\0bf,1)}|f_\beta(y)|\,dy
&=\int_{B(\0bf,1)}\frac{1}{|x|^\beta}\,dy
=\int_{0}^{1}r^{n-\beta-1}
\,dr\int_{\mathbb{S}^{n-1}}\,d
\sigma(x)\\
&\sim\frac{1}{n-\beta}.
\end{align*}
Combining this and \eqref{ex1223e1},
we further conclude that
\begin{align*}
\frac{\int_{B(\0bf,1)}|f_\beta(y)|\,dy}
{\|f_\beta\|_{\HerzSo}}
\sim(n-\beta)^{\frac{1}{p}-1}
\to\infty
\end{align*}
as $\beta\to n$. This implies that,
in this case,
\eqref{bal} does not hold true
for the local generalized Herz space
$\HerzSo$, and hence $\HerzSo$
is not a ball Banach function space.
Thus, we complete the proof of (iii),
and hence of Example \ref{ex1223}.
\end{proof}

Finally, we show that, under some reasonable
and sharp assumptions,
the global generalized Herz space
$\HerzS$ is a ball Banach function space
as follows.

\begin{theorem}\label{ball}
Let $p,\ q\in[1,\infty)$ and
$\omega\in M(\rp)$ satisfy
$\m0(\omega)\in(-\frac{n}{p},\infty)$ and
$\MI(\omega)\in(-\infty,0)$.
Then the global generalized Herz space
$\HerzS$ is a ball Banach function space.
\end{theorem}

\begin{proof}
Let all the symbols be as in the
present theorem. Then,
by Theorem \ref{Th2}, we know that
the global generalized
Herz space $\HerzS$ is a
ball quasi-Banach function space.
Moreover, we can easily show that the quasi-norm
$\|\cdot\|_{\HerzS}$ satisfies the triangle inequality
because $p$, $q\in[1,\infty)$.
Therefore, to finish the proof of
the present theorem, we only need to
prove that, for any $B\in\mathbb{B}$,
\eqref{bal} holds true with $X$ replaced by $\HerzS$.

To this end, let
$B(x_0,r)\in\mathbb{B}$ with $x_0\in\rrn$,
$k\in\mathbb{Z}$, and $r\in[2^{k-1},2^k)$.
Assume $\xi\in\rrn$ satisfies that
$|x_0-\xi|=3\cdot\tk$. Then we claim that
$$
B(x_0,r)\subset B(\xi,2^{k+2})\setminus
B(\xi,2^{k+1}).
$$
Indeed, applying the assumption
$|x_0-\xi|=3\cdot\tk$, we find that, for any
$y\in B(x_0,r)$,
\begin{align*}
|y-\xi|\leq|y-x_0|+|x_0-\xi|
<r+|x_0-\xi|<2^{k+2}
\end{align*}
and
\begin{align*}
|y-\xi|\geq|x_0-\xi|-|y-x_0|>|x_0-\xi|-r>2^{k+1}.
\end{align*}
These imply that $
B(x_0,r)\subset B(\xi,2^{k+2})\setminus
B(\xi,2^{k+1})
$ and hence the above claim holds true.
By this, the H\"{o}lder inequality,
Definition \ref{gh}(i),
and Remark \ref{remhs}(ii), we find that,
for any $f\in\HerzS$,
\begin{align*}
\int_{B(x_0,r)}|f(y)|\,dy&\leq
\int_{B(\xi,2^{k+2})\setminus
B(\xi,2^{k+1})}|f(y)|\,dy\\
&\lesssim\left\|f\1bf_{B(\xi,2^{k+2})
\setminus B(\xi,2^{k+1})}\right\|_{L^p(\rrn)}\\
&\lesssim\left\|f(\cdot+\xi)\right\|_{\HerzSo}
\lesssim\|f\|_{\HerzS},
\end{align*}
where the implicit positive constants depend
only on $p$, $\omega$, and $k$.
Thus, for any $B\in\mathbb{B}$,
\eqref{bal} holds true with
$X$ replaced by $\HerzS$,
which completes the proof of Theorem \ref{ball}.
\end{proof}

\begin{remark}
From Remark \ref{remark1221}, it follows that
Theorem \ref{ball} is sharp.
\end{remark}

\section{Convexities}

In this section, we discuss the convexity
of local and global generalized Herz spaces.
To this end, we first
investigate the relation between
generalized Herz spaces and their
convexification as follows.

\begin{lemma}\label{convexll}
Let $p,\ q,\ s\in(0,\infty)$ and $\omega\in M(\rp)$.
Then
$$\left[\HerzSo\right]^{1/s}=\HerzSocs$$
with the same quasi-norms.
\end{lemma}

\begin{proof}
Let $p,\ q,\ s\in(0,\infty)$
and $\omega\in M(\rp)$.
By Definition \ref{convex}(i) with $X:=\HerzSo$,
we conclude that, for any $f\in\Msc(\rrn)$,
\begin{align}\label{con1}
\|f\|_{[\HerzSo]
^{1/s}}&=\left\||f|^{\frac{1}{s}}\right\|
_{\HerzSo}^{s}\notag\\&=
\left\{\sum_{k\in\mathbb{Z}}
\left[\omega(\tk)\right]^{q}
\left[\int_{\tkm<|y|<\tk}
|f(y)|^{\frac{p}{s}}\,dy\right]
^{\frac{q}{p}}\right\}
^{\frac{s}{q}}\notag\\
&=\left\{\sum_{k\in\mathbb{Z}}
\left[\omega^{s}(\tk)\right]^{\frac{q}{s}}
\left[\int_{\tkm<|y|<\tk}
|f(y)|^{\frac{p}{s}}\,dy\right]
^{\frac{q/s}{p/s}}
\right\}^{\frac{s}{q}}\notag\\
&=\|f\|_{\HerzSocs}.
\end{align}
This implies that $[\HerzSo]^{1/s}=
\HerzSocs$ with the same
quasi-norms, and hence finishes
the proof of Lemma \ref{convexll}.
\end{proof}

\begin{lemma}\label{convexl}
Let $p,\ q,\ s\in(0,\infty)$ and $\omega\in M(\rp)$.
Then
$$\left[\HerzS\right]^{1/s}=\HerzScs$$
with the same quasi-norms.
\end{lemma}

\begin{proof}
Let $p,\ q,\ s\in(0,\infty)$ and $\omega\in M(\rp)$.
Note that, for any $f\in\Msc(\rrn)$,
\begin{align*}
&\left(\sup_{\xi\in\rrn}
\left\{\sum_{k\in\mathbb{Z}}
\left[\omega(\tk)\right]^{q}
\left[\int_{\tkm<|y-\xi|<\tk}
|f(y)|^{\frac{p}{s}}\,dy\right]
^{\frac{q}{p}}\right\}
^{\frac{1}{q}}\right)^{s}
\\&\quad=\sup_{\xi\in\rrn}
\left\{\sum_{k\in\mathbb{Z}}
\left[\omega(\tk)\right]^{q}
\left[\int_{\tkm<|y-\xi|<\tk}
|f(y)|^{\frac{p}{s}}\,dy\right]
^{\frac{q}{p}}\right\}
^{\frac{s}{q}}.
\end{align*}
Applying this, Definition
\ref{convex}(i) with $X:=\HerzS$, and
some arguments similar
to those used in the estimation of
\eqref{con1}, we find that,
for any $f\in\Msc(\rrn)$,
$$\left\|f\right\|_{
[\HerzS]^{1/s}}=\left\|f
\right\|_{\HerzScs},$$
which completes the proof of
Lemma \ref{convexl}.
\end{proof}

The following two results show that
local and global generalized Herz
spaces are strictly convex under
some assumptions of exponents $p$ and $q$.

\begin{theorem}\label{strc}
Let $p$, $q\in(0,\infty)$,
$s\in(0,\min\{p,q\}]$, and $\omega\in M(\rp)$.
Then the local generalized Herz space $\HerzSo$ is
strictly $s$-convex.
\end{theorem}

\begin{proof}
Let $p$, $q$, $s\in(0,\infty)$
satisfy $s\in(0,\min\{p,q\}]$ and $\omega\in M(\rp)$.
Then, for any given sequence
$\{f_{j}\}_{j\in\mathbb{N}}$ of
measurable functions in $[\HerzSo]^{1/s}$,
from the assumption $s\in(0,\min\{p,q\}]$,
we deduce that $p/s,\ q/s\in[1,\infty)$.
This, together with
Lemma \ref{convexll} and the
Minkowski inequality\index{Minkowski inequality},
further implies that, for any $N\in\mathbb{N}$,
\begin{align*}
\left\|\sum_{j=1}^{N}
|f_{j}|\right\|_{[\HerzSo]^{1/s}}\leq\sum_{j=1}^{N}
\|f_{j}\|_{[\HerzSo]^{1/s}}\leq\sum_{j\in\mathbb{N}}
\|f_{j}\|_{[\HerzSo]^{1/s}}.
\end{align*}
Applying this and the
monotone convergence
theorem\index{monotone convergence theorem}, we
conclude that
$$
\left\|\sum_{j\in\mathbb{N}}|f_{j}|\right\|
_{[\HerzSo]^{1/s}}
\leq\sum_{j\in\mathbb{N}}\|f_{j}\|_{[\HerzSo]^{1/s}},
$$
which implies that the
local generalized Herz space $\HerzSo$
is strictly
$s$-convex, and hence completes the proof of
Theorem \ref{strc}.
\end{proof}

\begin{theorem}\label{Th3.1}
Let $p,\ q\in(0,\infty)$, $s\in(0,\min\{p,q\}]$,
and $\omega\in M(\rp)$.
Then the global generalized Herz space $\HerzS$ is
strictly $s$-convex.
\end{theorem}

\begin{proof}
Let $p,\ q,\ s\in(0,\infty)$ satisfy
$s\in(0,\min\{p,q\}]$ and $\omega\in M(\rp)$.
Then, similarly to the proof of
Theorem \ref{strc}, by Lemma \ref{convexl},
we find that, for any given sequence
$\{f_{j}\}_{j\in\mathbb{N}}$ of
measurable functions in $[\HerzS]^{1/s}$,
$$
\left\|\sum_{j\in\mathbb{N}}|f_{j}|
\right\|_{[\HerzS]^{1/s}}\leq
\sum_{j\in\mathbb{N}}\|f_{j}\|_{[\HerzS]^{1/s}}.
$$
This implies that the global generalized
Herz space $\HerzS$ is strictly
$s$-convex, and hence finishes
the proof of Theorem \ref{Th3.1}.
\end{proof}

\section{Absolutely Continuous Quasi-Norms}

In this section, we investigate the absolute
continuity of quasi-norms of
local and global generalized Herz spaces.
Indeed, we show that the
local generalized Herz space has
an absolutely continuous quasi-norm, while
the global generalized Herz space
does not have one.

\begin{theorem}\label{abso}
Let $p,\ q\in(0,\infty)$ and $\omega\in M(\rp)$.
Then the local generalized
Herz space $\HerzSo$ has an
absolutely continuous quasi-norm.
\end{theorem}

\begin{proof}
Let $p,\ q\in(0,\infty)$ and $\omega\in M(\rp)$.
Assume that $f$ is any given measurable function
in $\HerzSo$ and $\{E_{i}\}_{i\in\mathbb{N}}$
a sequence of measurable sets satisfying
$\1bf_{E_{i}}\to0$ almost everywhere as $i\to\infty$.
For any $k\in\mathbb{Z}$ and $i\in\mathbb{N}$, let
$$
a_{k,i}:=\omega(\tk)\|f\1bf_{E_{i}}\1bf_{B(\0bf,
\tk)\setminus B(\0bf,\tkm)}\|_{L^{p}(\rrn)}
$$
and
$$
a_{k}:=\omega(\tk)\|f\1bf_{B(\0bf,
\tk)\setminus B(\0bf,\tkm)}\|_{L^{p}(\rrn)}.
$$
From $f\in\HerzSo$ and the definition of $\HerzSo$,
it follows that
\begin{equation}\label{le}
\left(\sum_{k\in\mathbb{Z}}
|a_{k}|^{q}\right)^{\frac{1}{q}}<\infty.
\end{equation}
Thus, for any $k\in\mathbb{Z}$,
$$\omega(\tk)\left\|f
\1bf_{B(\0bf,\tk)\setminus
B(\0bf,\tkm)}\right\|_{L^{p}(\rrn)}<\infty.$$
This implies that, for any $k\in\mathbb{Z}$,
$$f\1bf_{B(\0bf,\tk)\setminus
B(\0bf,\tkm)}\in L^{p}(\rrn).$$
By this, the fact that, for any $i\in\mathbb{N}$,
$$|f\1bf_{E_{i}}\1bf_{B(\0bf,\tk)\setminus
B(\0bf,\tkm)}|\leq|f\1bf_{B(\0bf,\tk)\setminus
B(\0bf,\tkm)}|,$$
and the dominated convergence
theorem\index{dominated convergence theorem},
we find that, for any $k\in\mathbb{Z}$,
$a_{k,i}\to0$ as $i\to\infty$.
Applying this, the estimate that,
for any $k\in\mathbb{Z}$ and $i\in\mathbb{N}$,
$|a_{k,i}|\leq|a_{k}|$, \eqref{le}, and
the dominated convergence theorem again,
we conclude that
\begin{align*}
\lim\limits_{i\to\infty}\|f\1bf_{E_{i}}\|_{\HerzSo}
=\lim\limits_{i\to\infty}
\left(\sum_{k\in\mathbb{Z}}
|a_{k,i}|^{q}\right)^{\frac{1}{q}}
=\left(\sum_{k\in\mathbb{Z}}
\lim\limits_{i\to\infty}
|a_{k,i}|^{q}\right)^{\frac{1}{q}}=0,
\end{align*}
which implies that $f$ has an
absolutely continuous quasi-norm
in $\HerzSo$, and hence the
local generalized Herz space
$\HerzSo$ has an absolutely quasi-norm.
Therefore, the proof of
Theorem \ref{abso} is completed.
\end{proof}

Via borrowing some ideas from the proof
of \cite[Example 5.1]{ST} and using
the definition of
global generalized Herz spaces, we construct
a special set $E$ and further
show that the global generalized Herz space
does not have an absolutely continuous quasi-norm.

\begin{example}\label{conter}
Let $p,\ q\in(0,\infty)$
and $\alpha\in(-\frac{1}{p},0)$.
For any $t\in(0,\infty)$,
let $\omega(t):=t^{\alpha}$, and let
$$E:=\bigcup\limits_{k\in\mathbb{N}}\left(
k-1+k^{-\frac{2}{\alpha p}},
k+k^{-\frac{2}{\alpha p}}\right).$$
Then the characteristic function
$\1bf_{E}\in\Kmp_{\omega}^{p,q}(\rr)$,
but $\1bf_{E}$
does not have an absolutely continuous
quasi-norm in $\Kmp_{\omega}^{p,q}(\rr)$.
This implies that the
global generalized Herz space
$\Kmp_{\omega}^{p,q}(\rr)$ may not
have an absolutely continuous quasi-norm.
\end{example}

\begin{proof}
Let $p,\ q\in(0,\infty)$ and $\alpha\in(-\frac{1}{p},0)$.
For any $t\in(0,\infty)$, let $\omega(t):=t^{\alpha}$,
and let
$$E:=\bigcup\limits_{k\in\mathbb{N}}
\left(k-1+k^{-\frac{2}
{\alpha p}},k+k^{-\frac{2}{\alpha p}}\right)$$
and, for any $k\in\mathbb{N}$,
$$E_{k}:=\left(k-1+k^{-\frac{2}{\alpha p}},
k+k^{-\frac{2}{\alpha p}}\right).$$
We now show that $\1bf_{E}\in
\Kmp_{\omega}^{p,q}(\rr)$.
Indeed, from the definition of $\|\cdot\|_{
\Kmp_{\omega,\0bf}^{p,q}(\rr)}$, it follows that,
for any $\xi\in\rr$,
\begin{align}\label{xx11}
&\left\|\1bf_{E}(\cdot+\xi)\right\|_{\Kmp_
{\omega,\0bf}^{p,q}(\rr)}\notag\\
&\quad=\left[\sum_{k\in\mathbb{Z}}2^{k\alpha q}
\left\|\1bf_{E}\1bf_{B(\xi,\tk)\setminus
B(\xi,\tkm)}\right\|_{L^{p}(\rr)}
^{q}\right]^{\frac{1}{q}}\notag\\
&\quad\leq\left[\sum_{k\in\mathbb{Z}}
2^{k\alpha q}
\left\|\1bf_{E}\1bf_{B(\xi,\tk)}
\right\|_{L^{p}(\rr)}^{q}\right]^{\frac{1}{q}}
\notag\\
&\quad\lesssim\left[\sum_{k\in\mathbb{N}}
2^{k\alpha q}
\left\|\1bf_{E}\1bf_{(\xi,\xi+\tk)}
\right\|_{L^{p}(\rr)}^{q}\right]^{\frac{1}{q}}
+\left[\sum_{k\in\mathbb{N}}2^{k\alpha q}
\left\|\1bf_{E}\1bf_{(\xi-\tk,\xi)}
\right\|_{L^{p}(\rr)}^{q}\right]^{\frac{1}{q}}\notag\\
&\quad\quad+\left[\sum_{k\in\zmn}
2^{k\alpha q}
\left\|\1bf_{E}\1bf_{B(\xi,\tk)}
\right\|_{L^{p}(\rr)}^{q}\right]^{\frac{1}{q}}\notag\\
&\quad=:{\rm J_{\xi,1}}+{\rm J_{\xi,2}}+{\rm J_{\xi,3}}.
\end{align}

For ${\rm J}_{\xi,1}$, we have
\begin{align}\label{121}
{\rm J}_{\xi,1}\sim\left[\sum_{k\in\mathbb{N}}
2^{k\alpha q}\left|
\bigcup\limits_{\widetilde{k}\in\mathbb{N}}
E_{\widetilde{k}}\cap(\xi,\xi+\tk)
\right|^{\frac{q}{p}}\right]^{\frac{1}{q}}.
\end{align}
Notice that, for any $k$,
$\widetilde{k}\in\mathbb{N}$ and any $\xi\in\rr$,
$E_{\widetilde{k}}\cap
(\xi,\xi+\tk)\neq\emptyset$ if and only if
\begin{equation}\label{x11}
\xi\in\left(
\widetilde{k}-1+\widetilde{k}^{-\frac{2}{\alpha p}}-\tk,
\widetilde{k}+\widetilde{k}^{-\frac{2}{\alpha p}}\right).
\end{equation}
Now, we estimate ${\rm J}_{\xi,1}$
by considering the following
three cases on $\xi$.

\emph{Case 1)} $\xi\in(-\infty,1]$.
In this case, by \eqref{x11}, we know that
\begin{align*}
\widetilde{k}-1+\widetilde{k}^{-\frac{2}{\alpha p}}
<\xi+\tk\leq1+\tk.
\end{align*}
It means
that $\widetilde{k}^{-\frac{2}{\alpha p}}
<\widetilde{k}+\widetilde{k}^{-\frac{2}{
\alpha p}}<2+\tk\leq\tka$ and
hence $\widetilde{k}\in
\mathbb{N}\cap(0,2^{-\frac{(k+1)\alpha p}{2}}).$
This, together with \eqref{121}
and the fact that, for any $\widetilde{k}
\in\mathbb{N}$, $|E_{\widetilde{k}}|=1$,
further implies that
\begin{align}\label{122}
{\rm J}_{\xi,1}
&\lesssim\left\{\sum_{k\in\mathbb{N}}2^{k\alpha q}
\left[\sum_{\widetilde{k}\in\mathbb{N}
\cap(0,2^{-\frac{(k+1)\alpha p}{2}})}
\left|E_{\widetilde{k}}\right|\right]^{\frac{q}{p}}\right\}
^{\frac{1}{q}}\notag\\&\lesssim\left[
\sum_{k\in\mathbb{N}}2^{k\alpha q}2^{-\frac{(
k+1)\alpha q}{2}}\right]^{\frac{1}{q}}
\sim1.
\end{align}

\emph{Case 2)} $\xi\in(k_{0}-1+k_{0}
^{-\frac{2}{\alpha p}},
k_{0}+k_{0}^{-\frac{2}{\alpha p}}]$
for some $k_{0}\in\mathbb{N}$.
In this case, from \eqref{x11} and the assumption
$\xi>k_{0}-1+k_{0}^{-\frac{2}{\alpha p}}$,
it follows that
$$k_{0}-1+k_{0}
^{-\frac{2}{\alpha p}}<\widetilde{k}+
\widetilde{k}^{-\frac{2}{\alpha p}}$$
and hence $\widetilde{k}\in
\mathbb{N}\cap[k_{0},\infty)$. On the
other hand, by \eqref{x11} and
the assumption $\xi\leq
k_{0}+k_{0}^{-\frac{2}{\alpha p}}$,
we conclude that
$$\widetilde{k}-1+\widetilde{k}^{-\frac{2}{\alpha p}}
<k_{0}+k_{0}^{-\frac{2}{\alpha p}}+\tk.$$
This, combined with
$\widetilde{k}\in\mathbb{N}\cap[k_{0},\infty)$,
further implies that
\begin{align*}\widetilde{k}
^{-\frac{2}{\alpha p}}-k_{0}^{-\frac{2}{\alpha p}}
\leq\widetilde{k}-k_{0}+\widetilde{k}
^{-\frac{2}{\alpha p}}-k_{0}
^{-\frac{2}{\alpha p}}
<\tk+1<\tka.\end{align*} Then, from
$\widetilde{k}^{-\frac{2}{\alpha p}}-k_{0}
^{-\frac{2}{\alpha p}}<\tka$
and the assumption
$\alpha\in(-\frac{1}{p},0)$, we deduce that
$$\left(\widetilde{k}-k_{0}\right)^{-\frac{2}{\alpha p}}
\leq\widetilde{k}^{-\frac{2}{\alpha p}}-k_{0}^{
-\frac{2}{\alpha p}}<\tka$$
and hence
$\widetilde{k}\in\mathbb{N}
\cap[k_{0},k_{0}+2^{-\frac{(k+1)\alpha p}{2}})$.
Applying this, \eqref{121}, and the fact that,
for any $\widetilde{k}
\in\mathbb{N}$, $|E_{\widetilde{k}}|=1$, we find that
\begin{align}\label{123}
{\rm J}_{\xi,1}
&\lesssim\left\{\sum_{k\in\mathbb{N}}2^{k\alpha q}
\left[\sum_{\widetilde{k}\in\mathbb{N}
\cap[k_{0},k_{0}+2^{-\frac{(k+1)\alpha p}{2}})}
\left|E_{\widetilde{k}}\right|\right]^{\frac{q}{p}}
\right\}^{\frac{1}{q}}\notag\\
&\lesssim\left[
\sum_{k\in\mathbb{N}}2^{k\alpha q}2^{-\frac{(
k+1)\alpha q}{2}}\right]^{\frac{1}{q}}\sim1.
\end{align}

\emph{Case 3)} $\xi\in(k_{0}+k_{0}
^{-\frac{2}{\alpha p}},
k_{0}+(k_{0}+1)^{-\frac{2}{\alpha p}}]$
for some $k_{0}\in\mathbb{N}$.
In this case, by an argument similar to
that used in Case 2), we conclude that,
for any $k$, $\widetilde{k}\in\mathbb{N}$,
if $E_{\widetilde{k}}\cap(\xi,\xi+\tk)\neq\emptyset$,
then
$$\widetilde{k}\in\mathbb{N}
\cap\left[k_{0}+1,k_{0}+1+2^{-\frac{k\alpha p}{2}}\right).$$
From this, \eqref{121}, and the fact that,
for any $\widetilde{k}
\in\mathbb{N}$, $|E_{\widetilde{k}}|=1$,
we deduce that
\begin{align*}
{\rm J}_{\xi,1}
&\lesssim\left\{\sum_{k\in\mathbb{N}}2^{k\alpha q}
\left[\sum_{\widetilde{k}\in\mathbb{N}
\cap[k_{0}+1,k_{0}+1+2^{-\frac{k\alpha p}{2}})}
\left|E_{\widetilde{k}}\right|\right]^{\frac{q}{p}}
\right\}^{\frac{1}{q}}\\
&\lesssim\left[
\sum_{k\in\mathbb{N}}2^{k\alpha q}2^{-\frac{
k\alpha q}{2}}\right]^{\frac{1}{q}}\sim1,
\end{align*}
which, combined with \eqref{122}
and \eqref{123}, further implies that,
for any $\xi\in\rr$, $${\rm J}_{\xi,1}\lesssim1.$$
Thus, we finish the estimation
of ${\rm J}_{\xi,1}$.

Now, we deal with ${\rm J}_{\xi,2}$. Notice that
\begin{align}\label{124}
{\rm J}_{\xi,2}
\sim\left[\sum_{k\in\mathbb{N}}
2^{k\alpha q}\left|
\bigcup\limits_{\widetilde{k}\in\mathbb{N}}
E_{\widetilde{k}}\cap(\xi-\tk,\xi)\right|
^{\frac{q}{p}}\right]^{\frac{1}{q}}
\end{align}
and, for any $k$, $\widetilde{k}\in\mathbb{N}$
and any $\xi\in\rr$, $E_{\widetilde{k}}
\cap(\xi-\tk,\xi)\neq\emptyset$ if and
only if
\begin{equation}\label{x21}
\xi\in\left(
\widetilde{k}-1+\widetilde{k}^{-\frac{2}{\alpha p}},
\widetilde{k}+\widetilde{k}
^{-\frac{2}{\alpha p}}+\tk\right).
\end{equation}
We next estimate ${\rm J}_{\xi,2}$ by considering
the following three cases on $\xi$.

\emph{Case i)} $\xi\in(-\infty,1]$.
In this case, notice that,
for any $\widetilde{k}\in\mathbb{N}$,
$\xi\leq1\leq\widetilde{k}-1+\widetilde{k}^{-
\frac{2}{\alpha p}}$, which contradicts
\eqref{x21}. Thus,
for any $\widetilde{k}\in\mathbb{N}$ and
$\xi\in(-\infty,1]$, we have
$$E_{\widetilde{k}}
\cap(\xi-\tk,\xi)=\emptyset$$
and hence, in this case,
\begin{equation}\label{ee12}
{\rm J_{\xi,2}}=0.
\end{equation}

\emph{Case ii)} $\xi\in(k_{0}-1+
k_{0}^{-\frac{2}{\alpha p}},k_{0}+
k_{0}^{-\frac{2}{\alpha p}}]$
for some $k_{0}\in\mathbb{N}$.
In this case, applying \eqref{x21}
and the assumption
$\xi\leq k_{0}+
k_{0}^{-\frac{2}{\alpha p}}$,
we conclude that
$$\widetilde{k}-1+
\widetilde{k}^{-\frac{2}{\alpha p}}<k_{0}+k_{0}
^{-\frac{2}{\alpha p}}$$
and hence
$\widetilde{k}\in\mathbb{N}\cap(0,k_{0}]$.
Moreover, from \eqref{x21} and the assumption
$\xi>k_{0}-1+k_{0}^{-\frac{2}{
\alpha p}}$, it follows that
$$k_{0}-1+k_{0}^
{-\frac{2}{\alpha p}}-\tk
<\widetilde{k}+\widetilde{k}^{-\frac{2}{\alpha p}}.$$
This, together with the fact that
$\widetilde{k}\in\mathbb{N}\cap
(0,k_{0}]$, further implies that
$$k_{0}^{-\frac{2}{\alpha p}}
-\widetilde{k}^{-\frac{2}{\alpha p}}\leq
k_{0}-\widetilde{k}+k_{0}^{-\frac{2}{\alpha p}}
-\widetilde{k}^{-\frac{
2}{\alpha p}}<1+\tk<\tka.$$
By this and the assumption $\alpha\in
(-\frac{1}{p},0)$,
we find that
$$\left(k_{0}-\widetilde{k}\right)^{-\frac{2}{\alpha p}}
\leq k_{0}^{-\frac{2}{\alpha p}}
-\widetilde{k}^{-\frac{2}{\alpha p}}<\tka$$
and hence
$\widetilde{k}>k_{0}-2^{-\frac{(k+1)\alpha p}{2}}$. Thus, from
\eqref{124}, $\widetilde{k}\in\mathbb{N}
\cap(k_{0}-2^{-\frac{(k+1)\alpha p}{2}},k_{0}]$,
and the fact that, for any $\widetilde{k}
\in\mathbb{N}$, $|E_{\widetilde{k}}|=1$,
we infer that
\begin{align}\label{125}
{\rm J}_{\xi,2}
&\lesssim\left\{\sum_{k\in\mathbb{N}}2^{k\alpha q}
\left[\sum_{\widetilde{k}\in\mathbb{N}\cap
(k_{0}-2^{-\frac{(k+1)\alpha p}{2}},k_{0}]}
\left|E_{\widetilde{k}}\right|\right]
^{\frac{q}{p}}\right\}
^{\frac{1}{q}}\notag\\
&\lesssim\left[
\sum_{k\in\mathbb{N}}2^{k\alpha q}2^{-\frac{(
k+1)\alpha q}{2}}\right]^
{\frac{1}{q}}\sim1.
\end{align}

\emph{Case iii)} $\xi\in(k_{0}+
k_{0}^{-\frac{2}{\alpha p}},
k_{0}+(k_{0}+1)^{-\frac{2}{\alpha p}}]$
for some $k_{0}
\in\mathbb{N}$. In this case,
by an argument similar to
that used in the proof of Case ii),
we conclude that, for any $k$,
$\widetilde{k}\in\mathbb{N}$, if
$E_{\widetilde{k}}\cap(\xi-\tk,\xi)\neq\emptyset$,
then
$$\widetilde{k}\in\mathbb{N}
\cap\lf(k_{0}-2^{-\frac{k\alpha p}{2}},k_{0}\r].$$
This, combined with \eqref{124} and the fact that,
for any $\widetilde{k}\in\mathbb{N}$,
$|E_{\widetilde{k}}|=1$, further implies that
\begin{align*}
{\rm J}_{\xi,2}
&\lesssim\left\{\sum_{k\in\mathbb{N}}2^{k\alpha q}
\left[\sum_{\widetilde{k}\in\mathbb{N}
\cap(k_{0}-2^{-\frac{k\alpha p}{2}},k_{0}]}\left|
E_{\widetilde{k}}\right|\right]^{\frac{q}{p}}\right\}
^{\frac{1}{q}}\\&\lesssim\left[
\sum_{k\in\mathbb{N}}2^{k\alpha q}2^{-\frac{
k\alpha q}{2}}\right]^{\frac{1}{q}}\sim1.
\end{align*}
Combining this,
\eqref{ee12}, and \eqref{125}, we further
conclude that, for any
$\xi\in\rr$, $${\rm J_{\xi,2}}\lesssim1$$
and hence we finish the
estimation of ${\rm J_{\xi,2}}$.

Now, to deal with ${\rm J_{\xi,3}}$,
notice that, for any $\xi\in\rr$,
\begin{align*}
{\rm J_{\xi,3}}
&\lesssim\left[\sum_{k\in\zmn}2^{k\alpha q}
\left\|\1bf_{B(\xi,\tk)}\right\|_{L^{p}(\rr)}
^{q}\right]^{\frac{1}{q}}\\
&\sim\left[\sum_{k\in\zmn}2^{k\alpha q}
\left|(\xi-\tk,\xi+\tk)\right|^{\frac{q}{p}}
\right]^{\frac{1}{q}}\\
&\sim\left[\sum_{k\in\zmn}2^{k\alpha q}
2^{\frac{(k+1)q}{p}}\right]^{\frac{1}{q}}\sim1,
\end{align*}
which, combined with \eqref{xx11} and
the estimates of both
${\rm J_{\xi,1}}$ and ${\rm J_{\xi,2}}$,
further implies that, for any $\xi\in\rr$,
$$
\left\|\1bf_{E}(\cdot+\xi)\right\|
_{\Kmp_{\omega,\0bf}^{p,q}(\rr)}
\lesssim1.
$$
From this, the definition of $\|\cdot\|
_{\Kmp_{\omega}^{p,q}(\rr)}$,
and Remark \ref{remhs}(ii),
it follows that
$\1bf_{E}\in\Kmp_{\omega}^{p,q}(\rr)$.

Finally, we show that $\1bf_{E}$ does not have an
absolutely continuous quasi-norm in the global
generalized Herz space $\Kmp_{\omega}^{p,q}(\rr)$.
To achieve this, for any $\widetilde{k}\in\mathbb{N}$,
let $F_{\widetilde{k}}:=(\widetilde{k},\infty)$ and
$\xi_{\widetilde{k}}:=
\widetilde{k}+\widetilde{k}
^{-\frac{2}{\alpha p}}-\frac{1}{2}$.
Then, obviously, for any $x\in\rr$,
we have $\1bf_{F_{\widetilde{k}}}(x)\to0$
as $\widetilde{k}\to\infty$.
In addition, notice that, for any
$\widetilde{k}\in\mathbb{N}$,
$$\widetilde{k}-1+\widetilde{k}^{-\frac{2}{
\alpha p}}\geq\widetilde{k}.$$ This implies that,
for any $\widetilde{k}\in\mathbb{N}$,
$E_{\widetilde{k}}\subset F_{\widetilde{k}}$.
Thus, for any $\widetilde{k}\in\mathbb{N}$, we have
\begin{align*}
&\left\|\1bf_{E}\1bf_{F_{\widetilde{k}}}\right\|
_{\Kmp_{\omega}^{p,q}(\rr)}\\
&\quad\geq\left\|\1bf_{E_{\widetilde{k}}}\right\|
_{\Kmp_{\omega}^{p,q}(\rr)}\\
&\quad\geq\left\|\1bf_{E_{\widetilde{k}}}
(\cdot+\xi_{\widetilde{k}})\right\|_
{\Kmp_{\omega,\0bf}^{p,q}(\rr)}\\
&\quad=\left[\sum_{k\in
\mathbb{Z}}2^{k\alpha q}\left\|
\1bf_{E_{\widetilde{k}}}\1bf_
{B(\xi_{\widetilde{k}},\tk)\setminus
B(\xi_{\widetilde{k}},\tkm)}\right\|
_{L^{p}(\rr)}^{q}\right]^{\frac{1}{q}}\\
&\quad\geq2^{-\alpha}
\left\|\1bf_{E_{\widetilde{k}}}
\1bf_{(\xi_{\widetilde{k}}+\frac{1}{4},
\xi_{\widetilde{k}}+\frac{1}{2})}
\right\|_{L^{p}(\rr)}\\
&\quad=2^{-\alpha}\left|
\left(\widetilde{k}-1+\widetilde{k}
^{-\frac{2}{\alpha p}},\widetilde{k}+
\widetilde{k}^{-\frac{
2}{\alpha p}}\right)\cap
\left(\widetilde{k}-\frac{1}{4}+\widetilde{k}
^{-\frac{2}{\alpha p}},\widetilde{k}+
\widetilde{k}^{-\frac{2}{\alpha p}}\right)
\right|^{\frac{1}{p}}\\
&\quad=2^{-\alpha-\frac{2}{p}}.
\end{align*}
Therefore, $\1bf_{E}$ does not have an
absolutely continuous
quasi-norm in $\Kmp_{\omega}^{p,q}(\rr)$,
which further implies that
the global generalized Herz space
$\Kmp_{\omega}^{p,q}(\rr)$
does not have an absolutely continuous
quasi-norm.
This finishes the proof of Example \ref{conter}.
\end{proof}

\begin{remark}
Let $p$, $q\in(0,\infty)$. We should point out
that the assumptions
$\alpha\in(-\frac{1}{p},0)$ and,
for any $t\in(0,\infty)$, $\omega(t)
:=t^\alpha$, in Example \ref{conter},
are reasonable. Namely, in this case,
the global generalized Herz space
$\Kmp_{\omega}^{p,q}(\rr)$
is a ball quasi-Banach function space.
Indeed, by Example \ref{ex113},
we know that, in this case, $\omega\in M(\rp)$ and
its Matuszewska--Orlicz indices satisfy
$\m0(\omega)\in(-\frac{1}{p},\infty)$
and $\MI(\omega)\in(-\infty,0)$.
From this and Theorem \ref{Th2} with $n=1$,
it follows that the global
generalized Herz space $\Kmp_{\omega}^{p,q}(\rr)$
as in Example \ref{conter} is a
ball quasi-Banach function space.
\end{remark}

\section{Boundedness of Sublinear Operators}

In this section, we establish a criterion
about the boundedness
of sublinear operators on local and
global generalized Herz spaces,
which was essentially obtained by
Rafeiro and Samko in \cite[Theorem 4.3]{HSamko}.
As applications, we obtain the
boundedness of both the Hardy--Littlewood
maximal operator and the
Calder\'{o}n--Zygmund operator on these Herz spaces.

Recall that an operator $T$ defined on $\Msc(\rrn)$
is called a \emph{sublinear
operator}\index{sublinear operator}
if, for any $f$, $g\in\Msc(\rrn)$ and
$\lambda\in\mathbb{C}$,
$$
\left|T(f+g)\right|\leq|T(f)|+|T(g)|
$$
and
$$
|T(\lambda f)|=|\lambda||T(f)|.
$$
In addition, for any
normed linear space $X$
and any operator $T$ on $X$, the \emph{operator
norm}\index{operator norm}
$\|T\|_{X\to X}$ of $T$ is defined by setting
\begin{equation}\label{eon}
\|T\|_{X\to X}:=\sup
_{\{x\in X,\,\|x\|_{X}=1\}}\left\|Tx\right\|_{X}.
\end{equation}
Then we present the boundedness criterion
on generalized Herz spaces as follows.

\begin{theorem}\label{vmaxh}
Let $p\in(1,\infty),\
q\in(0,\infty)$, and $\omega\in M(\rp)$ satisfy
\begin{equation}\label{1}
-\frac{n}{p}<\m0(\omega)\leq\M0(\omega)<\frac{n}{p'}
\end{equation}
and
\begin{equation}\label{2}
-\frac{n}{p}<\mi(\omega)\leq\MI(\omega)<\frac{n}{p'},
\end{equation}
where $\frac{1}{p}+\frac{1}{p'}=1$.
Assume that $T$ is a
sublinear operator satisfying that $T$
is bounded on $L^{p}(\rrn)$
and that there exists a positive
constant $\widetilde{C}$ such that, for any
$f\in\HerzSo$ and $x\notin\supp(f)$,
\begin{equation}\label{size}
|T(f)(x)|\leq\widetilde{C}\int_{\rrn}
\frac{|f(y)|}{|x-y|^{n}}\,dy.
\end{equation}
Then there exists a positive constant $C$,
independent of $T$, such that
\begin{enumerate}
  \item[{\rm(i)}] for any $f\in\HerzSo$,
  $$\left\|T(f)\right\|_{\HerzSo}\leq
  C\left[\widetilde{C}+
  \left\|T\right\|_{L^p(\rrn)
  \to L^p(\rrn)}\right]\|f\|_{\HerzSo};$$
  \item[{\rm(ii)}] for any $f\in\HerzS$,
  $$\left\|T(f)\right\|_{\HerzS}\leq
  C\left[\widetilde{C}+
  \left\|T\right\|_{L^p(\rrn)
  \to L^p(\rrn)}\right]\|f\|_{\HerzS}.$$
\end{enumerate}
\end{theorem}

To prove Theorem \ref{vmaxh}, we need the
following auxiliary estimates of the
positive function
$\omega\in M(\rp)$, which
is just \cite[Lemma 3.1]{HSamko}.

\begin{lemma}\label{La3.5}
Let $\omega\in M(\rp)$ and $\m0(\omega)$,
$\M0(\omega)$,
$\mi(\omega)$, and $\MI(\omega)$ denote its
Matuszewska--Orlicz indices. Then,
for any given $\eps\in(0,\infty)$,
there exists a positive constant $C_{(\eps)}$,
depending on $\eps$, such that,
for any $0<t<\tau<\infty$,
\begin{equation*}
\frac{\omega(t)}{\omega(\tau)}\leq C_{(\eps)}
\left(\frac{t}{\tau}\right)
^{\min\{\m0(\omega),\mi(\omega)\}-\eps}
\end{equation*}
and, for any $0<\tau<t<\infty$,
\begin{equation*}
\frac{\omega(t)}{\omega(\tau)}\leq C_{(\eps)}
\left(\frac{t}{\tau}\right)
^{\max\{\M0(\omega),\MI(\omega)\}+\eps}.
\end{equation*}
\end{lemma}

We also require the following basic inequality
which is a part of \cite[Exercise 1.1.4]{LGCF}.

\begin{lemma}\label{l4320}
Let $r\in(0,1]$.
Then, for any $\{a_j\}_{
j\in\mathbb{N}}\subset\mathbb{C}$,
$$
\left(\sum_{j\in\mathbb{N}}
|a_j|\right)^r\leq\sum_{j\in
\mathbb{N}}|a_j|^r.
$$
\end{lemma}

\begin{proof}[Proof of Theorem \ref{vmaxh}]
Let all the symbols be as in the present theorem.
We first show (i). To this end, for any given
$f\in\HerzSo$ and for any $k\in\mathbb{Z}$, let
$$
f_{k,1}:=f\1bf_{B(\0bf,2^{k-2})},\
f_{k,2}:=f\1bf_{B(\0bf,2^{k+1})\setminus
B(\0bf,2^{k-2})},$$
and
$$
f_{k,3}:=f\1bf_{[B(\0bf,2^{k+1})]^{\complement}}.
$$
Then, obviously, for any $k\in\mathbb{Z}$, we have
$$f=f_{k,1}+f_{k,2}+f_{k,3}.$$
From this and the sublinearity of $T$,
we deduce that
\begin{align}\label{mi}
\|T(f)\|_{\HerzSo}&\lesssim
\left\{\sum_{k\in\mathbb{
Z}}\left[\omega(\tk)\right]^{q}\left\|T\left(f_{k,1}
\right)\1bf_{B(
\0bf,\tk)\setminus B(\0bf,\tkm)}
\right\|_{L^{p}(\rrn)}^{q}
\right\}^{\frac{1}{q}}\notag\\
&\quad+\left\{\sum_{k\in\mathbb{
Z}}\left[\omega(\tk)\right]^{q}
\left\|T\left(f_{k,2}\right)\1bf_{B(
\0bf,\tk)\setminus B(\0bf,\tkm)}
\right\|_{L^{p}(\rrn)}^{q}
\right\}^{\frac{1}{q}}\notag\\
&\quad+\left\{\sum_{k\in\mathbb{
Z}}\left[\omega(\tk)\right]^{q}
\left\|T\left(f_{k,3}\right)\1bf_{B(
\0bf,\tk)\setminus B(\0bf,\tkm)}
\right\|_{L^{p}(\rrn)}^{q}
\right\}^{\frac{1}{q}}\notag\\
&=:{\rm J_{1}}+{\rm J_{2}}+{\rm J_{3}},
\end{align}
where the implicit positive constant
is independent of both $T$ and $f$.
We now deal with ${\rm J_{1}}$,
${\rm J_{2}}$,
and ${\rm J_{3}}$, separately. For this purpose,
let
$$\eps\in\left(0,\min\left\{\mw+\frac{n}{p},
\,\frac{n}{p'}-\Mw\right\}\right)$$
be a fixed positive constant.
Then, from Lemma \ref{La3.5}, it follows that,
for any $0<\tau<t<\infty$,
\begin{equation}\label{td1}
\frac{\omega(t)}{\omega(\tau)}\lesssim
\left(\frac{t}{\tau}\right)
^{\max\{\M0(\omega),\MI(\omega)\}+\eps}
\end{equation}
and, for any $0<t<\tau<\infty$,
\begin{equation}\label{tt1}
\frac{\omega(t)}{\omega(\tau)}\lesssim
\left(\frac{t}{\tau}\right)
^{\min\{\m0(\omega),\mi(\omega)\}-\eps}.
\end{equation}
For the simplicity of the presentation, let
$$m:=\mw-\eps$$
and
$$
M:=\Mw+\eps.$$
Obviously, we have $m\in(-\frac{n}{p},\infty)$
and $M\in(-\infty,\frac{n}{p'})$.

We first estimate ${\rm J_{1}}$.
Notice that, for any $k,\ i\in\mathbb{Z}$
satisfying $i\in(-\infty,k-2]$,
and for any $x,\ y\in\rrn$ satisfying
$\tkm\leq|x|<\tk$ and $2^{i-1}\leq|y|<2^{i}$,
we have $|y|<2^{i}\leq2^{k-2}$.
This implies that $|x-y|\geq|x|-|y|\geq2^{k-2}$.
Applying this, \eqref{size},
and the H\"{o}lder inequality, we find that,
for any $k\in\mathbb{Z}$ and $x\in\rrn$
with $2^{k-1}\leq|x|<2^{k}$,
\begin{align}\label{mi1}
\left|T\left(f_{k,1}\right)(x)\right|
&\leq\widetilde{C}\int_{\rrn}
\frac{|f_{k,1}(y)|}{|x-y|^{n}}\,dy\notag\\
&\leq\widetilde{C}\sum_{i=-\infty}^{k-2}
2^{-nk}\int_{B(\0bf,2^{i})\setminus
B(\0bf,2^{i-1})}|f(y)|\,dy\notag\\
&\lesssim\widetilde{C}\sum_{i=-\infty}^{k-2}
2^{-nk+\frac{ni}{p'}}\left\|
f\1bf_{B(\0bf,2^{i})\setminus
B(\0bf,2^{i-1})}\right\|_{L^{p}(\rrn)},
\end{align}
where the implicit positive constant
is independent of both $T$ and $f$.
On the other hand,
from both Lemma \ref{mono} and \eqref{td1},
it follows that, for any
$k,\ i\in\mathbb{Z}$ satisfying $i\in(-\infty,k-2]$,
\begin{align*}
\frac{\omega(2^{k})}{\omega(2^{i})}
\lesssim\frac{\omega(\tkmd)}{
\omega(2^{i})}\lesssim
\left(\frac{2^{k-2}}{2^{i}}\right)^{\Mw+\eps}
\sim2^{(k-i)M},
\end{align*}
which, combined with \eqref{mi1},
further implies that, for any
$k\in\mathbb{Z}$,
\begin{align}\label{mi12}
&\omega(\tk)\left\|T(f_{k,1})
\1bf_{B(\0bf,\tk)\setminus
B(\0bf,\tkm)}\right\|_{L^{p}(\rrn)}\notag\\
&\quad\lesssim\widetilde{C}\omega(\tk)
\sum_{i=-\infty}^{k-2}
2^{-nk+\frac{ni}{p'}}
\left\|f\1bf_{B(\0bf,2^{i})\setminus
B(\0bf,2^{i-1})}\right\|_{L^{p}(\rrn)}\notag\\
&\qquad\times\left\|\1bf_{B(\0bf,2^{k})\setminus
B(\0bf,2^{k-1})}\right\|_{L^{p}(\rrn)}\notag\\
&\quad\sim\widetilde{C}\sum_
{i=-\infty}^{k-2}2^{(i-k)\frac{n}{p'}}\frac{
\omega(\tk)}{\omega(2^{i})}\omega(2^{i})
\left\|f\1bf_{B(\0bf,2^{i})
\setminus B(\0bf,2^{i-1})}
\right\|_{L^{p}(\rrn)}\notag\\
&\quad\lesssim\widetilde{C}\sum_{i=-\infty}^{k-2}
2^{(k-i)(M-\frac{n}{p'})}
\omega(2^{i})\left\|f\1bf_{B(\0bf,2^{i})
\setminus B(\0bf,2^{i-1})}
\right\|_{L^{p}(\rrn)},
\end{align}
where the implicit positive constants
are independent of both $T$ and $f$.
By this and Lemma \ref{l4320},
we conclude that, for any $q\in(0,1]$,
\begin{align}\label{mi13}
{\rm J_{1}}
&\lesssim\widetilde{C}\left\{\sum_{k\in
\mathbb{Z}}\left[\sum_{i=-\infty}^{k-2}
2^{(k-i)(M-\frac{n}{p'})}\omega(2^{i})
\left\|f\1bf_{B(\0bf,2^{i})\setminus
B(\0bf,2^{i-1})}\right\|_{L^{p}(\rrn)}
\right]^{q}\right\}^{\frac{1}{q}}\notag\\
&\lesssim\widetilde{C}
\left\{\sum_{k\in\mathbb{Z}}
\sum_{i=-\infty}^{k-2}
2^{(k-i)(M-\frac{n}{p'})q}
\left[\omega(2^{i})\right]^{q}
\left\|f\1bf_{B(\0bf,2^{i})\setminus
B(\0bf,2^{i-1})}
\right\|_{L^{p}(\rrn)}^{q}
\right\}^{\frac{1}{q}}\notag\\
&\sim\widetilde{C}
\left\{\sum_{i\in\mathbb{Z}}
\left[\omega(2^{i})\right]^{q}
\left\|f\1bf_{B(\0bf,2^{i})\setminus
B(\0bf,2^{i-1})}
\right\|_{L^{p}(\rrn)}^{q}\sum_{k=i+2}^{\infty}
2^{(k-i)(M-\frac{n}{p'})q}
\right\}^{\frac{1}{q}}\notag
\\&\sim\widetilde{C}\|f\|_{\HerzSo}
\lesssim\left[\widetilde{C}+\|T\|_{
L^p(\rrn)\to L^p(\rrn)}\right]
\|f\|_{\HerzSo},
\end{align}
where the implicit positive constants
are independent of both $T$ and $f$.
Moreover,
if $q\in(1,\infty)$, from both \eqref{mi12} and the H\"{o}lder
inequality, we deduce that
\begin{align*}
&\omega(\tk)\left\|T(f_{k,1})
\1bf_{B(\0bf,\tk)\setminus
B(\0bf,\tkm)}\right\|_{L^{p}(\rrn)}\\
&\quad\lesssim
\widetilde{C}\left[\sum_{i=-\infty}^{k-2}
2^{\frac{(k-i)}{2}(M-\frac{n}{p'})q'}\right]
^{\frac{1}{q'}}\\
&\quad\quad\times\left\{\sum_{i=-\infty}^{k-2}
2^{\frac{(k-i)}{2}(M-\frac{n}{p'})q}
\left[\omega(\tk)\right]^{q}
\left\|f\1bf_{B(\0bf,\tk)\setminus B(\0bf,\tkm)}
\right\|_{L^{p}(\rrn)}^{q}\right\}
^{\frac{1}{q}}\\
&\quad\sim\widetilde{C}
\left\{\sum_{i=-\infty}^{k-2}
2^{\frac{(k-i)}{2}(M-\frac{n}{p'})q}
\left[\omega(\tk)\right]^{q}
\left\|f\1bf_{B(\0bf,\tk)\setminus
B(\0bf,\tkm)}
\right\|_{L^{p}(\rrn)}^{q}\right\}
^{\frac{1}{q}},
\end{align*}
where the implicit positive constants
are independent of both $T$ and $f$.
Therefore,
\begin{align}\label{mi14}
{\rm J_{1}}&\lesssim\widetilde{C}\left\{
\sum_{k\in\mathbb{Z}}\sum_{i=-\infty}^{k-2}
2^{\frac{(k-i)}{
2}(M-\frac{n}{p'})q}\left[\omega(2^{i})
\right]^{q}
\left\|f\1bf_{B(\0bf,2^{i})\setminus
B(\0bf,2^{i-1})}\right\|_{L^{p}(\rrn)}^{q}
\right\}^{\frac{1}{q}}\notag\\
&\sim\widetilde{C}\left\{\sum_{i\in\mathbb{Z}}
\left[\omega(2^{i})\right]^{q}
\left\|f\1bf_{B(\0bf,2^{i})
\setminus B(\0bf,2^{i-1})}\right\|_{L^{p}(\rrn)}
^{q}\sum_{k=i+2}^{\infty}
2^{\frac{(k-i)}{2}(M-\frac{n}{p'})q}
\right\}^{\frac{1}{q}}\notag\\&\sim
\widetilde{C}\|f\|_{\HerzSo}
\lesssim\left[\widetilde{C}+
\|T\|_{L^p(\rrn)\to L^p(\rrn)}\right]
\|f\|_{\HerzSo},
\end{align}
where the implicit positive constants
are independent of both $T$ and $f$,
which completes the estimation of ${\rm J_{1}}$.

Next, for ${\rm J_{2}}$, using the
boundedness of $T$ on
$L^{p}(\rrn)$, we find that,
for any $k\in\mathbb{Z}$,
\begin{align}\label{mi2}
&\left\|T(f_{k,2})\1bf_{B(\0bf,\tk)
\setminus B(\0bf,\tkm)}
\right\|_{L^{p}(\rrn)}\notag\\
&\quad\leq\left\|T(f_{k,2})\right\|_{L^{p}
(\rrn)}\leq\|T\|_{L^p(\rrn)\to L^p(\rrn)}
\|f_{k,2}\|_{L^{p}(\rrn)}\notag\\
&\quad\leq\|T\|_{L^p(\rrn)\to L^p(\rrn)}
\left[\left\|f\1bf_{B(\0bf,\tka)\setminus
B(\0bf,\tk)}\right\|_{L^{p}(\rrn)}
+\left\|f\1bf_{B(\0bf,\tk)\setminus
B(\0bf,\tkm)}\right\|_{L^{p}(\rrn)}\right.\notag
\\&\quad\quad\left.+\left\|f\1bf_{B(\0bf,\tkm)\setminus
B(\0bf,\tkmd)}\right\|_{L^{p}(\rrn)}\right].
\end{align}
In addition, by Lemma \ref{mono},
we conclude that, for any $k\in\mathbb{Z}$,
$$\omega(\tkm)\sim\omega(\tk)\sim\omega(\tka).$$
This, together with
\eqref{mi2}, further implies that
\begin{align}\label{mi21}
{\rm J_{2}}&\lesssim
\|T\|_{L^p(\rrn)\to L^p(\rrn)}
\left(\left\{\sum_{k\in\mathbb{
Z}}\left[\omega(\tka)\right]^{q}\left\|f\1bf_{B(
\0bf,\tka)\setminus B(\0bf,\tk)}
\right\|_{L^{p}(\rrn)}^{q}
\right\}^{\frac{1}{q}}\right.\notag\\
&\quad+\left\{\sum_{k\in\mathbb{
Z}}\left[\omega(\tk)\right]^{q}\left\|f\1bf_{B(
\0bf,\tk)\setminus B(\0bf,\tkm)}
\right\|_{L^{p}(\rrn)}^{q}
\right\}^{\frac{1}{q}}\notag\\
&\quad\left.+\left\{\sum_{k\in\mathbb{
Z}}\left[\omega(\tkm)\right]^{q}
\left\|f\1bf_{B(
\0bf,\tkm)\setminus B(\0bf,\tkmd)}
\right\|_{L^{p}(\rrn)}^{q}
\right\}^{\frac{1}{q}}\right)\notag\\
&\sim\|T\|_{L^p(\rrn)\to L^p(\rrn)}
\|f\|_{\HerzSo}\notag\\
&\lesssim
\left[\widetilde{C}+
\|T\|_{L^p(\rrn)\to L^p(\rrn)}\right]
\|f\|_{\HerzSo},
\end{align}
where the implicit positive are independent
of both $T$ and $f$.
Therefore, we complete the estimation
of ${\rm J_{2}}$.

Finally, we deal with ${\rm J_{3}}$.
Notice that, for any $k,\ i\in\mathbb{Z}$
satisfying $i\in[k+2,\infty)$,
and for any $x,\ y\in\rrn$ satisfying
$\tkm\leq|x|<\tk$ and $2^{i-1}\leq|y|<2^{i}$,
we have $|x|<2^{k}\leq2^{i-2}$,
which implies that $|x-y|\geq|x|-|y|\geq2^{i-2}$.
Applying this, \eqref{size},
and the H\"{o}lder inequality,
we find that, for any $k\in\mathbb{Z}$ and
$x\in\rrn$ with $2^{k-1}\leq|x|<2^{k}$,
\begin{align}\label{mi3}
|T(f_{k,3})(x)|&\leq\widetilde{C}
\int_{\rrn}\frac{|f_{k,3}(y)|}{|x-y|^{n}}
\,dy\notag\\&\leq\widetilde{C}\sum_{i=k+2}
^{\infty}2^{-ni}\int_{B(\0bf,
2^{i})\setminus B(\0bf,2^{i-1})}|f(y)|\,dy\notag\\
&\lesssim\widetilde{C}\sum_{i=k+2}^{\infty}
2^{-\frac{ni}{p}}\left\|
f\1bf_{B(\0bf,2^{i})\setminus
B(\0bf,2^{i-1})}\right\|_{L^{p}(\rrn)},
\end{align}
where the implicit positive constant
is independent of both $T$ and $f$.
On the other hand, from Lemma \ref{mono}
and \eqref{tt1}, it follows that, for any
$k,\ i\in\mathbb{Z}$ satisfying
$i\in[k+2,\infty)$,
\begin{align*}
\frac{\omega(2^{k})}{\omega(2^{i})}
\lesssim\frac{\omega(2^{k+2})}{
\omega(2^{i})}\lesssim\left(
\frac{2^{k+2}}{2^{i}}\right)^{\mw-\eps}
\sim2^{(k-i)m},
\end{align*}
which, together with \eqref{mi3},
further implies that, for any $k\in\mathbb{Z}$,
\begin{align}
&\omega(\tk)\left\|T(f_{k,3})
\1bf_{B(\0bf,\tk)\setminus
B(\0bf,\tkm)}\right\|_{L^{p}(\rrn)}\notag\\
&\quad\lesssim\widetilde{C}
\omega(\tk)\sum_{i=k+2}^{\infty}
2^{-\frac{ni}{p}}\left\|f\1bf_{
B(\0bf,2^{i})
\setminus B(\0bf,2^{i-1})}\right\|_
{L^{p}(\rrn)}\left\|\1bf_{B(\0bf,2^{k})
\setminus B(\0bf,2^{k-1})}
\right\|_{L^{p}(\rrn)}\notag\\
&\quad\sim\widetilde{C}
\sum_{i=k+2}^{\infty}2^{(k-i)
\frac{n}{p}}\frac{
\omega(\tk)}{\omega(2^{i})}\omega(2^{i})
\left\|f\1bf_{B(\0bf,2^{i})\setminus
B(\0bf,2^{i-1})}\right\|_{L^{p}(\rrn)}\notag\\
&\quad\lesssim\widetilde{C}
\sum_{i=k+2}^{\infty}2^{(k-i)(m+\frac{n}{p})}
\omega(2^{i})
\left\|f\1bf_{B(\0bf,2^{i})\setminus
B(\0bf,2^{i-1})}\right\|_{L^{p}(\rrn)},\notag
\end{align}
where the implicit positive constants are
independent of both $T$ and $f$.
By this and an argument similar to
that used in the estimations of
both \eqref{mi13} and \eqref{mi14},
we find that
\begin{equation}
{\rm J_{3}}\lesssim
\widetilde{C}\|f\|_{\HerzSo}
\lesssim
\left[\widetilde{C}+\|T\|_{
L^p(\rrn)\to L^p(\rrn)}
\right]\|f\|_{\HerzSo},\notag
\end{equation}
where the implicit positive constants
are independent of both $T$ and $f$,
which completes the estimation
of ${\rm J_{3}}$. Combining this,
\eqref{mi}, \eqref{mi13}, \eqref{mi14}, and
\eqref{mi21},
we conclude that there
exists a positive constant $C$,
independent of $T$,
such that, for any $f\in\HerzSo$,
\begin{equation}\label{Th3.4e3}
\|T(f)\|_{\HerzSo}\leq C
\left[\widetilde{C}+
\|T\|_{L^p(\rrn)\to
L^p(\rrn)}\right]\|f\|_{\HerzSo},
\end{equation}
which completes the proof of (i).

Next, we turn to prove (ii). For this purpose,
let $f\in\HerzS$.
By the definition of $\|\cdot\|_{\HerzS}$,
we find that, for any $\xi\in\rrn$,
\begin{align*}
\|f(\cdot+\xi)\|_{\HerzSo}
\leq\|f\|_{\HerzS}<\infty.
\end{align*}
Thus, for any $\xi\in\rrn$,
$f(\cdot+\xi)\in\HerzSo$.
From this and \eqref{Th3.4e3},
we deduce that
\begin{align*}
\left\|T(f)(\cdot+\xi)
\right\|_{\HerzSo}&\leq C
\left[
\widetilde{C}+\|T\|_{L^p(\rrn)
\to L^p(\rrn)}\right]
\|f(\cdot+\xi)\|_{\HerzSo}\\
&\leq C\left[
\widetilde{C}+\|T\|_{L^p(\rrn)
\to L^p(\rrn)}\right]\|f\|_{\HerzS},
\end{align*}
which, combined with
the definition of $\|\cdot\|_{\HerzS}$
and the arbitrariness of $\xi\in\rrn$,
further implies that, for any $f\in\HerzS$,
$$
\|T(f)\|_{\HerzS}\leq
C\left[
\widetilde{C}+\|T\|_{L^p(\rrn)
\to L^p(\rrn)}\right]\|f\|_{\HerzS},
$$
where $C$ is a positive constant
independent of both $T$ and $f$.
This finishes the proof of (ii),
and hence of Theorem \ref{vmaxh}.
\end{proof}

\begin{remark}
\begin{enumerate}
\item[(i)] We point out that, in Theorem \ref{vmaxh},
both the assumptions \eqref{1} and
\eqref{2} on $\omega$ are reasonable
and sharp. Indeed, for any given
$\alpha\in\rr$ and for any $t\in(0,\infty)$,
let $\omega(t):=t^\alpha$. Then,
from Example \ref{ex113}, it follows
that both the assumptions
\eqref{1} and \eqref{2}
imply that $\alpha\in(-\frac{n}{p},
\frac{n}{p'})$.
In this case, Theorem \ref{vmaxh}(i)
coincides with
the conclusions on the classical Herz spaces
obtained by Li and Yang
\cite[Corollary 2.1]{LiY}. Moreover,
by \cite[Examples 5.1.1 and 5.1.2]{LYH},
we find that the assumption
$\alpha\in(-\frac{n}{p},\frac{n}{p'})$ on
\cite[Corollary 2.1]{LiY} is sharp.
Therefore, in this sense,
both the assumptions \eqref{1} and \eqref{2}
in Theorem \ref{vmaxh} are sharp.
  \item[(ii)] Recall that,
  when $q\in(1,\infty)$, Rafeiro and Samko
established the boundedness of
sublinear operators on
the global generalized Herz spaces
$\HerzS$ in \cite[Theorem 4.3]{HSamko},
which is just a special case
of Theorem \ref{vmaxh}(ii).
\end{enumerate}
\end{remark}

The following two conclusions
show the boundedness of
the Hardy--Littlewood maximal operator,
respectively, on
local and global generalized Herz spaces.
Indeed, they are simple
corollaries of Theorem \ref{vmaxh}
because the Hardy--Littlewood maximal operator
$\mc$ satisfies the size condition
\eqref{size} and is also bounded on the
Lebesgue space $L^p(\rrn)$ when
$p\in(1,\infty)$; we omit the details.

\begin{corollary}\label{bhlo}
Let $p$, $q$, and $\omega$ be as in
Theorem \ref{vmaxh} and $\mc$ the
Hardy--Littlewood maximal operator as
in \eqref{hlmax}.
Then there exists a
positive constant $C$ such that,
for any $f\in L^1_{\mathrm{loc}}(\rrn)$,
$$
\left\|\mc(f)\right\|_{\HerzSo}
\leq C\|f\|_{\HerzSo}.
$$
\end{corollary}

\begin{remark}
Let $\omega(t):=t^\alpha$ for
any $t\in(0,\infty)$ and for any
given $\alpha\in\rr$
in Corollary \ref{bhlo}.
Then, in this case,
Corollary \ref{bhlo} goes back to
the boundedness of
the Hardy--Littlewood maximal
operator on the classical
homogeneous Herz spaces obtained in
\cite[Theorem 4.1]{IM}.
\end{remark}

\begin{corollary}\label{Th4}
Let $p$, $q$, and $\omega$ be as in
Theorem \ref{vmaxh} and $\mc$ the
Hardy--Littlewood
maximal operator as in \eqref{hlmax}.
Then there exists a positive constant
$C$ such that,
for any $f\in L^1_{\mathrm{loc}}(\rrn)$,
$$
\left\|\mc(f)\right\|_{\HerzS}
\leq C\|f\|_{\HerzS}.
$$
\end{corollary}

Next, we establish the boundedness of
Calder\'{o}n--Zygmund operators on
local and global generalized
Herz spaces.
For this purpose, we first
recall some basic concepts
(see, for instance, \cite[Chapter
III]{s93}). In what follows,
for any $\gamma=(\gamma_1,\ldots,
\gamma_n)\in \zp^n$,
any $\gamma$-order differentiable
function $F(\cdot,\cdot)$
on $\rrn\times\rrn$, and any $(x,y)\in
\rrn\times\rrn$, let
$$\index{$\partial_{(1)}^{\gamma}$}
\partial_{(1)}^{\gamma}F(x,y):
=\frac{\partial^{|\gamma|}}
{\partial x_1^{\gamma_1}
\cdots\partial x_n^{\gamma_n}}F(x,y)
$$
and
$$\index{$\partial_{(2)}^{\gamma}$}
\partial_{(2)}^{\gamma}
F(x,y):=\frac{\partial^{|\gamma|}}
{\partial y_1^{\gamma_1}
\cdots\partial y_n^{\gamma_n}}F(x,y).
$$
Moreover, let $$\index{$\Delta$}\Delta:=\left\{(x,y)
\in\rrn\times\rrn:\ x=y\right\}$$
be the \emph{diagonal}\index{diagonal} of $\rrn\times\rrn$.
Then the definition of
$d$-order standard kernels is as follows.

\begin{definition}\label{def-s-k}
Let $d\in\zp$ and $\delta\in(0,1]$.
A measurable function $K$ on $(\rrn\times\rrn)\setminus\Delta$
is called a \textit{$d$-order standard kernel}
\index{$d$-order standard kernel}if
there exists a positive constant $C$ such that
the following two statements hold true:
\begin{enumerate}
\item[\rm (i)]
	for any $\gamma\in\zp^n$ with  $|\gamma|\leq d$,
and for any $x$, $y\in\rrn$ with $x\neq y$,
    \begin{align}\label{size-s'}
	\left|\partial_{(1)}^{\gamma}
K(x,y)\right|+\left|\partial_{(2)}^{\gamma}
K(x,y)\right|\leq\frac{C}{|x-y|^{n+|\gamma|}};
	\end{align}
	\item[\rm (ii)] for any $\gamma\in\zp^n$ with
$|\gamma|=d$, and for any $x$, $y$, $z\in\rrn$
with $x\neq y$ and $|x-y|>2|y-z|$,
	\begin{align}\label{regular2-s}
\left|\partial_{(2)}^{\gamma}K(x,y)-
\partial_{(2)}^{\gamma}
K(x,z)\right|\leq C\frac{|y-z|^{\delta}}
{|x-y|^{n+d+\delta}}
	\end{align}
and
\begin{align*}
\left|\partial_{(1)}^{\gamma}K(y,x)-
\partial_{(1)}^{\gamma}
K(z,x)\right|\leq C\frac{|y-z|^{\delta}}
{|x-y|^{n+d+\delta}}.
\end{align*}
\end{enumerate}
\end{definition}

The following definition of
Calder\'{o}n--Zygmund operators
was given in \cite[Definition 4.1.8]{LGMF} (see also
\cite[Definition 5.11]{d01}).

\begin{definition}\label{defin-C-Z-s}
Let $d\in\zp$ and $K$ be a $d$-order
standard kernel as in Definition \ref{def-s-k}.
A linear operator $T$ is called a
\emph{$d$-order Calder\'on--Zygmund operator}
\index{$d$-order Calder\'on--Zygmund operator}with kernel $K$ if
\begin{enumerate}
  \item[{\rm(i)}] $T$ is bounded on $L^2(\rrn)$;
  \item[{\rm(ii)}] for any given $f\in L^2(\rrn)$
  with compact support, and for any $x\notin\supp(f)$,
  \begin{equation}\label{czp}
  T(f)(x)=\int_{\rrn}K(x,y)f(y)\,dy.
  \end{equation}
\end{enumerate}
\end{definition}

Recall that the following
$L^p(\rrn)$ boundedness of
Calder\'{o}n--Zygmund operators is
just \cite[Theorem 5.10]{d01}.

\begin{lemma}\label{czop}
Let $p\in(1,\infty)$, $d\in\zp$,
$K$ be a $d$-order standard
kernel, and $T$ a $d$-order
Calder\'{o}n--Zygmund operator with kernel $K$.
Then $T$ is well defined on $L^p(\rrn)$
and there exists a positive constant
$C$ such that, for any
$f\in L^p(\rrn)$,
$$
\|T(f)\|_{L^p(\rrn)}
\leq C\|f\|_{L^p(\rrn)}.
$$
\end{lemma}

Now, we show the boundedness of the
Calder\'{o}n--Zygmund operator on
local generalized Herz spaces as follows.
Recall that $C_{\mathrm{c}}
(\rrn)$\index{$C_{\mathrm{c}}(\rrn)$}
is defined to be the set of all the continuous
functions on $\rrn$ with
compact support.

\begin{corollary}\label{czoh}
Let $p\in(1,\infty)$, $q\in(0,\infty)$,
and $\omega\in M(\rp)$
satisfy both \eqref{1} and \eqref{2}.
Assume that $d\in\zp$ and
$T$ is a $d$-order
Calder\'{o}n--Zygmund
operator with
a $d$-order standard kernel $K$.
Then $T$ has an unique extension
on $\HerzSo$ and
there exists a positive
constant $C$ such that,
for any $f\in\HerzSo$,
$$
\left\|T(f)\right\|_{\HerzSo}
\leq C\|f\|_{\HerzSo}.
$$
\end{corollary}

\begin{proof}
Let all the symbols be as in the
present corollary. We first
assume $f\in C_{\mathrm{c}}(\rrn)$.
Then, applying both \eqref{czp}
and \eqref{size-s'} with $\gamma:=
\0bf$,
we find that, for any $x\notin
\supp(f)$,
$$
\left|T(f)(x)\right|=
\left|\int_{\rrn}K(x,y)f(y)\,dy\right|
\lesssim
\int_{\rrn}\frac{|f(y)|}{|x-y|^n}\,dy.
$$
This, together with both Lemma
\ref{czop} and
Theorem \ref{vmaxh},
further implies that
\begin{equation*}
\|T(f)\|_{\HerzSo}\lesssim\|f\|_{\HerzSo},
\end{equation*}
where the implicit positive constant
is independent of $f$.

In addition, by the assumption
$\m0(\omega)\in(-\frac{n}{p},\infty)$
and Theorems \ref{Th3}
and \ref{abso},
we conclude that the local generalized
Herz space $\HerzSo$ is a ball quasi-Banach
function space having an absolutely continuous
quasi-norm.
From this and \cite[Proposition 3.8]{TYYZ}
with $X$ therein replaced by $\HerzSo$,
we deduce that $C_{\mathrm{c}}(\rrn)$
is dense in $\HerzSo$. This,
combined with a standard density argument
(see, for instance, \cite[p.\,45]{LDY}),
further implies that, for any
$f\in\HerzSo$, $T(f)$ is well defined and
$$
\|T(f)\|_{\HerzSo}\lesssim\|f\|_{\HerzSo}
$$
with the implicit positive constant
independent of $f$, which completes the
proof of Corollary \ref{czoh}.
\end{proof}

\begin{remark}
We should point out that,
in Corollary \ref{czoh}, when $d=0$
and $\omega(t):=t^\alpha$ for any
$t\in(0,\infty)$
and for any given $\alpha\in\rr$, then
this corollary goes back to
\cite[Remark 5.1.1]{LYH}.
\end{remark}

Applying Corollary \ref{czoh}
and an argument similar to that
used in the proof of Theorem \ref{Th3.4}(ii),
we conclude the following boundedness
of the Calder\'{o}n--Zygmund operator
on global generalized Herz spaces immediately;
we omit the details.

\begin{corollary}\label{czohg}
Let $p\in(1,\infty)$, $q\in(0,\infty)$,
and $\omega\in M(\rp)$
satisfy \eqref{1} and \eqref{2}.
Assume that $d\in\zp$ and $T$ is a
$d$-order Calder\'{o}n--Zygmund
operator with a $d$-order
standard kernel $K$.
Then $T$ has a unique extension on $\HerzS$
and there exists a positive constant $C$ such that,
for any $f\in\HerzS$,
$$
\left\|T(f)\right\|_{\HerzS}\leq C\|f\|_{\HerzS}.
$$
\end{corollary}

\section{Fefferman--Stein Vector-Valued Inequalities}

In this section, we show the Fefferman--Stein
vector-valued inequalities on
local and global generalized Herz spaces,
which are important in the study of
Hardy spaces associated with these
generalized Herz spaces. For instance,
in order to establish
the atomic characterizations of
generalized Herz--Hardy spaces,
we need to prove some auxiliary lemmas
about the Fefferman--Stein
vector-valued inequality in Chapter \ref{sec6}
(see Lemmas \ref{vmbhl} and \ref{vmbhg}
below).

We first have the following Fefferman--Stein vector-valued
inequality\index{Fefferman--Stein vector-valued \\inequality}
on the local generalized Herz space $\HerzSo$.

\begin{theorem}\label{Th3.4}
Let $p,\ r\in(1,\infty)$, $q\in(0,\infty)$, and
$\omega\in M(\rp)$ satisfy
\begin{equation*}
-\frac{n}{p}<\m0(\omega)
\leq\M0(\omega)<\frac{n}{p'}
\end{equation*}
and
\begin{equation*}
-\frac{n}{p}<\mi(\omega)
\leq\MI(\omega)<\frac{n}{p'},
\end{equation*}
where $\frac{1}{p}+\frac{1}{p'}=1$.
Then there exists a positive constant
$C$ such that, for any $\{
f_{j}\}_{j\in\mathbb{N}}\subset
L_{{\rm loc}}^{1}(\rrn)$,
$$
\left\|\left\{\sum_{j\in\mathbb{N}}\left[\mc
(f_{j})\right]^{r}\right\}
^{\frac{1}{r}}\right\|_{\HerzSo}
\leq C\left\|\left(
\sum_{j\in\mathbb{N}}|f_{j}|^{r}\right)
^{\frac{1}{r}}\right\|_{\HerzSo}.
$$
\end{theorem}

To prove Theorem \ref{Th3.4}, we need the
following Fefferman--Stein
vector-valued inequality on the Lebesgue space
$L^{p}(\rrn)$ with $p\in(1,\infty)$,
which is a part of \cite[Theorem 5.6.6]{LGCF}.

\begin{lemma}\label{Thvvmf}
Let $p,\ r\in(1,\infty)$.
Then there exists a positive constant
$C$ such that, for any $\{
f_{j}\}_{j\in\mathbb{N}}\subset
L_{{\rm loc}}^{1}(\rrn)$,
$$
\left\|\left\{\sum_{j\in\mathbb{N}}\left[\mc
(f_{j})\right]^{r}\right\}
^{\frac{1}{r}}\right\|_{L^{p}(\rrn)}
\leq C\left\|\left(
\sum_{j\in\mathbb{N}}|f_{j}|^{r}\right)^
{\frac{1}{r}}\right\|_{L^{p}(\rrn)}.
$$
\end{lemma}

We now show Theorem \ref{Th3.4}
via Lemma \ref{Thvvmf}.

\begin{proof}[Proof of Theorem \ref{Th3.4}]
Let all the symbols be as in
the present theorem.
For any given sequence
$\{f_{j}\}_{j\in\mathbb{N}}$ of
locally integrable functions on
$\rrn$, and for any $g\in\Msc(\rrn)$ and any
$x\in\rrn$, let
$$
A(g)(x):=\mvrgex,
$$
where, for any $j\in\mathbb{N}$ and $y\in\rrn$,
$$
\eta_{i}(y):=\left\{
\begin{aligned}
&\frac{f_{i}(y)}
{[\sum_{j\in\mathbb{N}}|f_{j}(y)|^{r}]^
{\frac{1}{r}}},\hspace{.5cm}
\text{when $\vry\neq0$},\\
&0,\hspace{3.2cm}\text{when $\vry=0$}.
\end{aligned}
\right.
$$
Then it is easy to show that, for any
$\lambda\in\mathbb{C}$ and $g\in
\Msc(\rrn)$,
\begin{equation}\label{hh}
A(\lambda g)=|\lambda|A(g).
\end{equation}
Moreover, by the Minkowski
inequality, we conclude that,
for any $g_{1},\ g_{2}\in\Msc(\rrn)$,
$$
A(g_{1}+g_{2})\leq A(g_{1})+A(g_{2}),
$$
which, combined with \eqref{hh},
further implies that $A$
is a sublinear operator.
Now, we show that $A$ satisfies \eqref{size}.
Indeed, observe that, for any $y\in\rrn$,
\begin{equation}\label{vv}
\vrey\leq1.
\end{equation}
This, together with both the fact
that $\mc$ satisfies \eqref{size}
(see \cite[Remark 4.4]{HSamko}), and the
Minkowski integral
inequality\index{Minkowski integral inequality},
implies that, for any $x\notin
\supp(g)$,
\begin{align}\label{Th3.4e1}
|A(g)(x)|&=\mvrgex
\lesssim\left\{\sum_{j\in\mathbb{N}}
\left[\int_{\rrn}\frac{|g(y)\eta_{j}(y)|}{
|x-y|^{n}}\,dy\right]^{r}\right\}
^{\frac{1}{r}}\notag\\
&\lesssim\int_{\rrn}\left\{\sum_{j\in\mathbb{N}}
\left[\frac{|g(y)\eta_{j}(y)|}{|
x-y|^{n}}\right]^{r}\right\}^{\frac{1}{r}}
\,dy\lesssim\int_{\rrn}\frac{|g(y)|}{|x-y|^{n}}\,dy,
\end{align}
where the implicit positive constants
are independent of
$\{f_j\}_{j\in\mathbb{N}}$,
which further implies
that the sublinear operator $A$
satisfies \eqref{size}.
Then we show that $A$ is bounded on $L^{p}(\rrn)$.
Indeed, from both
Lemma \ref{Thvvmf} and \eqref{vv}, we deduce that,
for any $g\in\Msc(\rrn)$,
\begin{align}\label{Th3.4e2}
\|A(g)\|_{L^{p}(\rrn)}&=
\left\|\mvrge\right\|_{L^{p}(\rrn)}
\notag\\&\lesssim
\left\|\vrge\right\|_{L^{p}(\rrn)}\lesssim
\|g\|_{L^{p}(\rrn)},
\end{align}
where the implicit positive constants
are independent of $\{f_j\}_{j\in\mathbb{N}}$.
Thus, the sublinear operator $A$
is bounded on $L^p(\rrn)$ and hence
$A$ satisfies all the assumptions of
Theorem \ref{vmaxh}.
This, together with
the fact that the implicit
positive constants in both \eqref{Th3.4e1}
and \eqref{Th3.4e2} are independent
of $\{f_j\}_{j\in\mathbb{N}}$,
further implies that,
for any $g\in\Msc(\rrn)$,
$$
\left\|A(g)\right\|_{\HerzSo}
\lesssim\|g\|_{\HerzSo},
$$
where the implicit positive constant
is independent of $\{f_j\}_{j\in\mathbb{N}}$.
Using this with $g:=\{\sum_{j\in\mathbb{N}}
|f_{j}|^{r}\}^{1/r}$, we find that
\begin{align*}
&\left\|\mvr\right\|_{\HerzSo}\\
&\quad=\|A(g)\|_{\HerzSo}\lesssim
\|g\|_{\HerzSo}\sim\left\|\vr\right\|_{\HerzSo},
\end{align*}
where the implicit positive constants
are independent of $\{f_j\}
_{j\in\mathbb{N}}$,
which completes the proof of Theorem \ref{Th3.4}.
\end{proof}

\begin{remark}
In Theorem \ref{Th3.4}, when $\omega(t):=t^\alpha$
for any $t\in(0,\infty)$ and
for any given $\alpha\in\rr$,
then, in this case,
the local generalized Herz space
$\HerzSo$ coincides with
the classical homogeneous Herz space
$\dot{K}_p^{\alpha,q}(\rrn)$.
In this case,
Theorem \ref{Th3.4} goes back to
\cite[Corollary 4.5]{i10} with
$p(\cdot):=p\in(1,\infty)$.
\end{remark}

For the global generalized Herz space $\HerzS$,
applying an argument similar to
that used in the
proof of Theorem \ref{Th3.4}
with $\HerzSo$ replaced
by $\HerzS$, we obtain the following
Fefferman--Stein vector-valued
inequality\index{Fefferman--Stein vector-valued \\inequality}
on $\HerzS$; we omit the details.

\begin{theorem}\label{Th3.3}
Let $p$, $q$, $r$, and $\omega$
be as in Theorem \ref{Th3.4}.
Then there exists a positive
constant $C$ such that, for any $\{
f_{j}\}_{j\in\mathbb{N}}\subset
L_{{\rm loc}}^{1}(\rrn)$,
$$
\left\|\left\{\sum_{j\in\mathbb{N}}\left[\mc
(f_{j})\right]^{r}\right\}
^{\frac{1}{r}}\right\|_{\HerzS}
\leq C\left\|\left(
\sum_{j\in\mathbb{N}}|f_{j}|^{r}\right)
^{\frac{1}{r}}\right\|_{\HerzS}.
$$
\end{theorem}

\section{Dual and
Associate Spaces of Local Generalized Herz Spaces}

In this section, we
first find the dual space of the local
generalized Herz space $\HerzSo$ with
$p$, $q\in(1,\infty)$.
Then, using this dual result and
the relation between the dual spaces
and the associate
spaces of ball Banach function spaces,
we obtain the associate space of $\HerzSo$
under some reasonable and sharp assumptions.

To begin with, we show the dual theorem
of the local generalized Herz space $\HerzSo$.
Indeed, we can establish a more general conclusion,
namely, the duality of the local
generalized Herz space $\HerzSx$ with
given $\xi\in\rrn$.
For this purpose, we first introduce the local
generalized Herz space $\HerzSx$ as follows.

\begin{definition}
Let $p$, $q\in(0,\infty)$, $\omega\in M(\rp)$,
and $\xi\in\rrn$.
Then the \emph{local generalized Herz
space}\index{local generalized Herz space} $\HerzSx$\index{$\HerzSx$}
is defined to be the set of all the measurable
functions $f$ on $\rrn$ such that
$$
\|f\|_{\HerzSx}:=
\left\{\sum_{k\in\mathbb{Z}}\left[
\omega(\tk)\right]^q\left\|f\1bf_{B(\xi,\tk)
\setminus B(\xi,\tkm)}\right\|_{L^p(\rrn)}
^q\right\}^{\frac{1}{q}}<\infty.
$$
\end{definition}

\begin{remark}\label{remark213}
Obviously, for any given $\xi\in\rrn$ and
for any
$f\in\Msc(\rrn)$, we have
$$
\|f\|_{\HerzSx}=\|f(\cdot+\xi)\|_{\HerzSo}
$$
and
$$
\|f\|_{\HerzS}=\sup_{\xi\in\rrn}\|f\|_{\HerzSx}.
$$
\end{remark}

Then we have the following dual theorem
on the Herz space $\HerzSx$.
To be precise, we show
that, when
$p$, $q\in(1,\infty)$,
the dual space of $\HerzSx$ is just
$\HerzSxd$.

\begin{theorem}\label{dual}
Let $p,\ q\in(1,\infty)$, $\omega\in M(\rp)$,
and $\xi\in\rrn$.
Then the dual space\index{dual space} of $\HerzSx$,
denoted by $(\HerzSx)^*$\index{$(\HerzSx)^*$},
is $\HerzSxd$ in the following sense:
\begin{enumerate}
 \item[{\rm(i)}]
 Let $g\in\HerzSxd$. Then the
 linear functional
 \begin{equation}\label{funld}
 \phi_{g}:\ f\mapsto\phi_{g}(f):=
 \int_{\rrn}f(y)g(y)\,dy,
 \end{equation}
 defined for any $f\in\HerzSx$,
 is bounded on $\HerzSx$.
 \item[{\rm(ii)}]
 Conversely, any continuous linear functional
 on $\HerzSx$ arises as in \eqref{funld}
 with a unique $g\in\HerzSxd$.
\end{enumerate}
Moreover, for any $g\in\HerzSxd$,
$$\left\|g
\right\|_{\HerzSxd}=\left\|\phi_{g}
\right\|_{(\HerzSx)^*}.$$
\end{theorem}

\begin{proof}
Let all the symbols be as in
the present theorem.
We first show (i). Indeed, by
the monotone convergence theorem
and the H\"{o}lder inequality,
we find that, for any
$f\in\HerzSx$ and $g\in\HerzSxd$,
\begin{align}\label{duale0}
&\left|\int_{\rrn}
f(y)g(y)\,dy\right|\notag\\
&\quad\leq\int_{\rrn}
|f(y)g(y)|\,dy\notag\\
&\quad=\sum_{k\in\mathbb{Z}}
\int_{B(\xi,2^k)\setminus
B(\xi,2^{k-1})}|f(y)g(y)|\,dy\notag\\
&\quad\leq\sum_{k\in\mathbb{Z}}\left\|f
\1bf_{B(\xi,2^k)\setminus
B(\xi,2^{k-1})}\right\|_{L^p(\rrn)}
\left\|g
\1bf_{B(\xi,2^k)\setminus
B(\xi,2^{k-1})}\right\|_{L^{p'}(\rrn)}\notag\\
&\quad\leq\left\{\sum_{k\in\mathbb{Z}}
\left[\omega(2^k)\right]^q\left\|f
\1bf_{B(\xi,2^k)\setminus B(\xi,2^{k-1})}
\right\|_{L^p(\rrn)}^q\right\}^{\frac{1}{q}}\notag\\
&\quad\quad\times\left\{\sum_{k\in\mathbb{Z}}
\left[\omega(2^k)\right]^{-q'}\left\|f
\1bf_{B(\xi,2^k)\setminus
B(\xi,2^{k-1})}\right\|_{L^{p'}(\rrn)}
^{q'}\right\}^{\frac{1}{q'}}\notag\\
&\quad=\|f\|_{\HerzSx}\|g\|_{\HerzSxd},
\end{align}
which implies that the linear functional $\phi_g$
as in \eqref{funld} is bounded on $\HerzSx$ and
\begin{equation}\label{duale1}
\left\|\phi_g\right\|_{(\HerzSx)^*}\leq
\left\|g\right\|_{\HerzSxd}.
\end{equation}
This finishes the proof of (i).

Conversely, we show (ii).
To this end, let $\phi\in(\HerzSx)^*$.
We first claim that,
for any $k\in\mathbb{Z}$, $$\phi
\in\left(L^p(B(\xi,2^k)\setminus
B(\xi,2^{k-1}))\right)^*,$$
here and thereafter, for any measurable set $E$
in $\rrn$, $L^p(E)$\index{$L^p(E)$}
is defined to be the set of all
the measurable functions $f$
on $\rrn$ such that
$\supp(f)\subset E$ and
$$\|f\|_{L^p(E)}:=
\|f\|_{L^p(\rrn)}<\infty.$$
Indeed, for any $f\in
L^p(B(\xi,2^k)\setminus B(\xi,2^{k-1}))$
with $k\in\mathbb{Z}$, we have
\begin{align*}
\|f\|_{\HerzSx}&=\omega(2^k)
\|f\1bf_{B(\xi,2^k)
\setminus B(\xi,2^{k-1})}\|_{L^p(\rrn)}\\
&=\omega(2^k)\|f\|_{L^p(B(\xi,2^k)
\setminus B(\xi,2^{k-1}))}<\infty,
\end{align*}
which implies that $f\in\HerzSx$ and hence
\begin{align*}
\left|\phi(f)\right|&\leq\|\phi\|_{(\HerzSx)^*}
\|f\|_{\HerzSx}\\
&=\omega(2^k)\|\phi\|_{(\HerzSx)^*}
\|f\|_{L^p(B(\xi,2^k)\setminus B(\xi,2^{k-1}))}.
\end{align*}
Therefore, for any
$k\in\mathbb{Z}$, $$\phi\in
\left(L^p(B(\xi,2^k)\setminus
B(\xi,2^{k-1}))\right)^*.$$
From this and \cite[Theorem 4.11]{b11},
it follows that,
for any $k\in\mathbb{Z}$,
there exists a $g_k\in L^{p'}(B(\xi,2^k)\setminus
B(\xi,2^{k-1}))$ such that, for any
$f\in L^p(B(\xi,2^k)\setminus B(\xi,2^{k-1}))$,
\begin{equation}\label{duale2}
\phi(f)=\int_{B(\xi,2^k)\setminus B(\xi,2^{k-1})}
f(y)g_k(y)\,dy.
\end{equation}
We now show that
\begin{equation}\label{duale11}
g:=\sum_{k\in\mathbb{Z}}g_k\in\HerzSxd
\end{equation}
by considering the following
two cases on the operator $\phi$.

\emph{Case 1)} $\phi=0$. In this case,
applying \eqref{duale2},
we conclude that, for any
$k\in\mathbb{Z}$,
$g_k=0$ almost everywhere in
$B(\xi,2^k)\setminus B(\xi,2^{k-1})$,
which further implies that $g_k=0$
almost everywhere in $\rrn$.
Using this, we find that $g=0$ almost everywhere
in $\rrn$ and hence $g\in\HerzSxd$.
Moreover, we have
\begin{align}\label{duale3}
\|g\|_{\HerzSxd}=0=\|\phi\|_{(\HerzSx)^*}.
\end{align}

\emph{Case 2)} $\phi\neq0$. In this case,
we first claim
that there exists a $k_0\in\mathbb{Z}$
such that
\begin{equation}\label{duale12}
\left|\left\{x\in
B(\xi,2^{k_0})\setminus
B(\xi,2^{k_0-1}):\
|g_{k_0}(x)|\neq0\right\}\right|
\in(0,\infty).
\end{equation}
Indeed, by the assumption
$\phi\neq0$, we find that there
exists an $f_0\in\HerzSx$ such that
$\phi(f_0)\neq0$. In addition,
for any $f\in\HerzSx$, we have
\begin{align}\label{duale4}
&\left\|
f-\sum_{k=-N}^{N}f\1bf_{B(\xi,2^k)
\setminus B(\xi,2^{k-1})}
\right\|_{\HerzSx}\notag\\
&\quad=\left\{
\sum_{j\in\mathbb{Z}}\left[\omega(2^j)\right]^q
\left\|\sum_{k\in\{k\in\mathbb{Z}:
\,|k|>N\}}f\1bf_{B(\xi,2^k)
\setminus B(\xi,2^{k-1})}\1bf_{B(\xi,2^j)
\setminus B(\xi,2^{j-1})}
\right\|_{L^p(\rrn)}^q\right\}
^{\frac{1}{q}}\notag\\
&\quad=\left\{\sum_{j\in
\{j\in\mathbb{Z}:\,|j|>N\}}
\left[\omega(2^j)\right]^q\left\|
f\1bf_{B(\xi,2^j)\setminus
B(\xi,2^{j-1})}\right\|_{L^p(\rrn)}
^q\right\}^{\frac{1}{q}}\to0
\end{align}
as $N\to\infty$. From this and
the continuity of $\phi$
on $\HerzSx$, it follows that,
for any $f\in\HerzSx$,
$$
\lim\limits_{N\to\infty}\phi
\left(\sum_{k=-N}^{N}f\1bf_{B(\xi,2^k)
\setminus B(\xi,2^{k-1})}\right)=\phi(f).
$$
This further implies that
there exists an $N_0\in\mathbb{N}$
such that
\begin{align}\label{duale5}
\left|\phi\left(\sum_{k=-N_0}^{N_0}
f_0\1bf_{B(\xi,2^k)\setminus
B(\xi,2^{k-1})}\right)\right|
>\frac{|\phi(f_0)|}{2}>0.
\end{align}
Thus, there exists a $k_0
\in\mathbb{Z}\cap
[-N_0,N_0]$ such that
$\phi(f_0\1bf_{B(\xi,2^{k_0})
\setminus B(\xi,2^{k_0-1})})\neq0$.
Otherwise, for any $k\in\mathbb{Z}\cap
[-N_0,N_0]$, it holds true that
$\phi(f_0\1bf_{B(\xi,2^{k})
\setminus B(\xi,2^{k-1})})=0$. This,
together with the linearity of $\phi$,
further implies
that
$$\phi\left(\sum_{k=-N_0}^{N_0}
f_0\1bf_{B(\xi,2^k)\setminus
B(\xi,2^{k-1})}\right)=0,$$
which contradicts \eqref{duale5}.
Therefore, for this $k_0\in\mathbb{Z}$
such that $$\phi\left(f_0\1bf_{B(\xi,2^{k_0})
\setminus B(\xi,2^{k_0-1})}\right)\neq0,$$
we have
\begin{align*}
\left\|
f_0\1bf_{B(\xi,2^{k_0})\setminus
B(\xi,2^{k_0-1})}\right\|_{L^p(\rrn)}
\leq\left[\omega(2^{k_0})\right]^{-1}
\left\|f_0\right\|_{\HerzSx}<\infty.
\end{align*}
Using this, the assumption
$\phi(f_0\1bf_{B(\xi,2^{k_0})
\setminus B(\xi,2^{k_0-1})})\neq0$,
and \eqref{duale2},
we further conclude that
$$
\left|\left\{x\in
B(\xi,2^{k_0})\setminus
B(\xi,2^{k_0-1}):\
|g_{k_{0}}(x)|\neq0\right\}
\right|\in(0,\infty),
$$
which completes the
proof of the above claim.

Now, for any given $k\in\mathbb{Z}$,
when $g_k=0$ almost everywhere
in $B(\xi,2^k)\setminus
B(\xi,2^{k-1})$, let $f_k:=0$ and,
otherwise, let $$f_k:=
\left[\omega(2^k)\right]^{-q'}
\left\|g_k\1bf_{B(\xi,2^k)
\setminus B(\xi,2^{k-1})}\right\|_{
L^{p'}(\rrn)}^{q'-p'}\left|g_k\right|^{p'-1}
\text{sgn}\overline{g_k},$$
here and thereafter, the function
\text{sgn} is defined by setting
$$
\text{sgn }z:=\left\{
\begin{aligned}
&\frac{z}{|z|},\hspace{.35cm}\text{when
$z\in\mathbb{C}\setminus\{0\}$},\\
&0,\hspace{.7cm}\text{when $z=0$}.
\end{aligned}
\right.
$$
Then we show that, for any
$k\in\mathbb{Z}$,
\begin{equation}\label{duale6}
\left\|f_k\1bf_{B(\xi,2^k)\setminus
B(\xi,2^{k-1})}\right\|_{L^p(\rrn)}=\left[
\omega(2^k)\right]^{-q'}\left\|g_k\1bf_{B(\xi,2^k)
\setminus B(\xi,2^{k-1})}
\right\|_{L^{p'}(\rrn)}^{q'-1}
\end{equation}
and
\begin{equation}\label{duale7}
\int_{B(\xi,2^{k})\setminus
B(\xi,2^{k-1})}f_k(y)g_k(y)\,dy=
\left[\omega(2^k)\right]^{-q'}
\left\|g_k\1bf_{B(\xi,2^k)
\setminus B(\xi,2^{k-1})}\right\|
_{L^{p'}(\rrn)}^{q'}.
\end{equation}
Indeed, for any $k\in\mathbb{Z}$
satisfying $g_k=0$ almost everywhere
in $B(\xi,2^k)\setminus
B(\xi,2^{k-1})$, it is obvious that both
\eqref{duale6} and \eqref{duale7} hold true.
On the other hand, for any
given $k\in\mathbb{Z}$
such that
$$
\left|\left\{x\in
B(\xi,2^{k})\setminus
B(\xi,2^{k-1}):\
|g_{k}(x)|\neq0\right\}\right|\in(0,\infty),
$$
we have $\|g_k\1bf_{B(\xi,2^k)\setminus
B(\xi,2^{k-1})}\|_{L^{p'}(\rrn)}\in(0,\infty)$.
Thus,
\begin{align*}
&\left\|f_k\1bf_{B(\xi,2^k)\setminus
B(\xi,2^{k-1})}\right\|_{L^p(\rrn)}\\
&\quad=\left[\omega(2^k)\right]^{-q'}
\left\|g_k\1bf_{B(\xi,2^k)
\setminus B(\xi,2^{k-1})}\right\|
_{L^{p'}(\rrn)}^{q'-p'}
\left[\int_{B(\xi,2^k)\setminus B(\xi,2^{k-1})}
|g_k(y)|^{p'}\,dy\right]^{\frac{1}{p}}\\
&\quad=\left[\omega(2^k)\right]^{-q'}
\left\|g_k\1bf_{B(\xi,2^k)
\setminus B(\xi,2^{k-1})}\right\|
_{L^{p'}(\rrn)}^{q'-p'}
\left\|g_k\1bf_{B(\xi,2^k)
\setminus B(\xi,2^{k-1})}\right\|
_{L^{p'}(\rrn)}^{\frac{p'}{p}}\\
&\quad=\left[\omega(2^k)\right]^{-q'}
\left\|g_k\1bf_{B(\xi,2^k)
\setminus B(\xi,2^{k-1})}\right\|
_{L^{p'}(\rrn)}^{q'-1}
\end{align*}
and
\begin{align*}
&\int_{B(\xi,2^{k})\setminus
B(\xi,2^{k-1})}f_k(y)g_k(y)\,dy\\&\quad=
\left[\omega(2^k)\right]^{-q'}
\left\|g_k\1bf_{B(\xi,2^k)
\setminus B(\xi,2^{k-1})}
\right\|_{L^{p'}(\rrn)}^{q'-p'}
\int_{B(\xi,2^k)\setminus B(\xi,2^{k-1})}
|g_k(y)|^{p'}\,dy\\
&\quad=\left[\omega(2^k)
\right]^{-q'}\left\|g_k\1bf_{B(\xi,2^k)
\setminus B(\xi,2^{k-1})}
\right\|_{L^{p'}(\rrn)}^{q'-p'}
\left\|g_k\1bf_{B(\xi,2^k)
\setminus B(\xi,2^{k-1})}
\right\|_{L^{p'}(\rrn)}^{p'}\\
&\quad=\left[\omega(2^k)
\right]^{-q'}\left\|g_k\1bf_{B(\xi,2^k)
\setminus B(\xi,2^{k-1})}
\right\|_{L^{p'}(\rrn)}^{q'}.
\end{align*}
These further imply that both
\eqref{duale6} and \eqref{duale7}
hold true.
Moreover, applying \eqref{duale6}
and the assumption that,
for any $j\in\mathbb{N}$, $g_j\in
L^{p'}(B(\xi,2^j)\setminus
B(\xi,2^{j-1}))$, we conclude that,
for any $N\in\mathbb{N}$,
\begin{align}\label{duale8}
&\left\|\sum_{k=-N}^{N}f_k\right\|_{\HerzSx}
\notag\\&\quad=\left\{
\sum_{j\in\mathbb{Z}}
\left[\omega(2^j)\right]^q
\left\|\sum_{k=-N}^{N}
f_k\1bf_{B(\xi,2^j)\setminus
B(\xi,2^{j-1})}\right\|_{L^p(\rrn)}
^q\right\}^{\frac{1}{q}}\notag\\
&\quad=\left\{\sum_{j=-N}^{N}
\left[\omega(2^j)\right]^{q}\left\|
f_j\1bf_{B(\xi,2^j)\setminus
B(\xi,2^{j-1})}\right\|_{L^p(\rrn)}
^q\right\}^{\frac{1}{q}}\notag\\
&\quad=\left\{\sum_{j=-N}^{N}
\left[\omega(2^j)\right]^{-q'}
\left\|g_j\1bf_{B(\xi,2^j)\setminus
B(\xi,2^{j-1})}\right\|_{L^{p'}(\rrn)}
^{q'}\right\}^{\frac{1}{q}}
<\infty,
\end{align}
which implies that
$$\sum_{k=-N}^{N}f_k\in\HerzSx.$$
By this, the linearity of $\phi$,
\eqref{duale6},
\eqref{duale2}, and \eqref{duale7},
we find that,
for any $N\in\mathbb{N}$,
\begin{align*}
\phi\left(\sum_{k=-N}^{N}f_k\right)
&=\sum_{k=-N}^{N}\phi(f_k)\\
&=\sum_{k=-N}^{N}\int_{B(\xi,2^k)
\setminus B(\xi,2^{k-1})}f_k(y)g_k(y)\,dy\\
&=\sum_{k=-N}^{N}\left[\omega(2^k)\right]^{-q'}
\left\|g_k\1bf_{B(\xi,2^k)\setminus
B(\xi,2^{k-1})}\right\|_{L^{p'}(\rrn)}^{q'},
\end{align*}
which, together with \eqref{duale8},
implies that
\begin{align*}
&\sum_{k=-N}^{N}\left[\omega(2^k)\right]^{-q'}
\left\|g_k\1bf_{B(\xi,2^k)\setminus
B(\xi,2^{k-1})}\right\|_{L^{p'}(\rrn)}^{q'}\\
&\quad\leq\|\phi\|_{(\HerzSx)^*}
\left\|\sum_{k=-N}^{N}f_k\right\|_{\HerzSx}\\
&\quad=\|\phi\|_{(\HerzSx)^*}
\left\{\sum_{k=-N}^{N}
\left[\omega(2^k)\right]^{-q'}
\left\|g_k\1bf_{B(\xi,2^k)\setminus
B(\xi,2^{k-1})}\right\|_{L^{p'}(\rrn)}
^{q'}\right\}^{\frac{1}{q}}.
\end{align*}
Then, by this, \eqref{duale11},
and \eqref{duale12},
we conclude that, for any
$N\in\mathbb{N}\cap[|k_0|,\infty)$,
\begin{align*}
&\left\{\sum_{k=-N}^{N}
\left[\omega(2^k)\right]^{-q'}
\left\|g\1bf_{B(\xi,2^k)\setminus
B(\xi,2^{k-1})}\right\|_{L^{p'}(\rrn)}^{q'}
\right\}^{\frac{1}{q'}}\\
&\quad=\left\{\sum_{k=-N}^{N}
\left[\omega(2^k)\right]^{-q'}
\left\|g_k\1bf_{B(\xi,2^k)\setminus
B(\xi,2^{k-1})}\right\|_{L^{p'}
(\rrn)}^{q'}\right\}^{\frac{1}{q'}}\\
&\quad\leq\|\phi\|_{(\HerzSx)^*}.
\end{align*}
Letting $N\to\infty$, we obtain
\begin{equation}\label{duale10}
\|g\|_{\HerzSxd}\leq\|\phi\|_{(\HerzSx)^*},
\end{equation}
which further implies that
$g\in\HerzSxd$ in this case,
and hence completes the proof of
the claim $g\in\HerzSxd$.
Thus, from (i), it follows that the
linear functional $\phi_g$ as in \eqref{funld}
is bounded on $\HerzSx$.

Next, we show that
$\phi=\phi_g$. Indeed,
applying \eqref{duale0}, we find that,
for any $f\in\HerzSx$,
$$|fg|\in L^1(\rrn).$$
Combining this, \eqref{duale4},
both the continuity and the linearity
of $\phi$ on $\HerzSx$,
\eqref{duale2},
\eqref{duale11}, and the dominated
convergence theorem,
we further conclude that,
for any $f\in\HerzSx$,
\begin{align*}
\phi(f)&=\lim\limits_{N\to\infty}
\sum_{k=-N}^{N}\phi\left(
f\1bf_{B(\xi,2^k)\setminus
B(\xi,2^{k-1})}\right)\\
&=\lim\limits_{N\to\infty}\sum_{k=-N}^{N}
\int_{B(\xi,2^k)\setminus
B(\xi,2^{k-1})}f(y)\1bf_{B(\xi,2^k)
\setminus B(\xi,2^{k-1})}(y)g_k(y)\,dy\\
&=\lim\limits_{N\to\infty}\int_{\rrn}
f(y)\1bf_{B(\xi,2^N)\setminus
B(\xi,2^{-N-1})}(y)g(y)\,dy\\
&=\int_{\rrn}f(y)g(y)\,dy=\phi_g(f).
\end{align*}
This finishes the proof that $\phi=\phi_g$.
Moreover, from this, \eqref{duale1}, \eqref{duale3},
and \eqref{duale10}, we infer that
\begin{align*}
\|\phi\|_{(\HerzSx)^*}=\|\phi_g\|_{(\HerzSx)^*}
\leq\|g\|_{\HerzSxd}\leq\|\phi\|_{(\HerzSx)^*}.
\end{align*}
Then we have
$$\left\|\phi_g\right\|_{
(\HerzSx)^*}=\left\|g
\right\|_{\HerzSx}.$$
Using this and the linearity of $\phi_g$ about $g$,
we further find that the $g$ is unique.
This finishes the proof of (ii),
and hence of Theorem \ref{dual}.
\end{proof}

\begin{remark}
We should point out that,
in Theorem \ref{dual}, when $\xi=\0bf$
and $\omega(t):=t^{\alpha}$ for any
$t\in(0,\infty)$ and for any given
$\alpha\in\rr$, then the conclusion obtained
in this theorem goes back to that of \cite
[Corollary 1.2.1]{LYH}.
\end{remark}

Now, we turn to find the associate
space of $\HerzSo$.
Recall that Bennett and Sharpley
\cite[Chapter 1, Corollary 4.3]{BSIO}
showed the relation between dual spaces
and associate spaces of \emph{Banach
function spaces}\index{Banach function
space} [namely, quasi-Banach
function spaces satisfy \eqref{bal}
with ball replaced
by any measurable set
of finite measure].
Moreover, we find that
the proof of
\cite[Chapter 1, Corollary 4.3]{BSIO}
still works on ball
Banach function spaces. Thus,
the following conclusion holds true,
which is an essential tool
to study the associate space of
$\HerzSo$.

\begin{lemma}\label{gassociate}
Let $X$ be a ball Banach function space,
and let $X'$ and $X^*$\index{$X^*$} denote, respectively,
its associate space and dual space.
Then $X'=X^*$ with the same norms
if and only if $X$ has an absolutely
continuous norm.
\end{lemma}

Via Lemma \ref{gassociate}, we now
show that the associate space of
$\HerzSo$ is $\HerzSoa$ under some reasonable
and sharp assumptions. Concretely,
we have the following conclusion.

\begin{theorem}\label{associate}
Let $p,\ q\in(1,\infty)$ and
$\omega\in M(\rp)$ satisfy
$$
-\frac{n}{p}<\m0(\omega)\leq
\M0(\omega)<\frac{n}{p'}.
$$
Then
$$
\left(\HerzSo\right)'=\HerzSoa
$$
with the same norms,
where $(\HerzSo)'$ denotes
the associate space
of the local generalized Herz space
$\HerzSo$.
\end{theorem}

\begin{proof}
Let all the symbols be as in the present
theorem. Then, by Theorem
\ref{balll} and the assumptions
of the present theorem, we find that
the local generalized Herz space
$\HerzSo$ is a ball Banach function space.
In addition, from Theorem \ref{abso},
it follows that $\HerzSo$ has
an absolutely continuous norm.
Combining this, Lemma \ref{gassociate},
and Theorem \ref{dual},
we further conclude that
$$
\left(\HerzSo\right)'=
\left(\HerzSo\right)^*=
\HerzSoa
$$ with the same norms, which completes
the proof of Theorem \ref{associate}.
\end{proof}

\section{Extrapolation Theorems}

In this section, we establish
the extrapolation theorems
of local and global
generalized Herz spaces via
Theorem \ref{associate}
obtained in last section.
To this end,
we first recall the following definition
of the Muckenhoupt
weights\index{Muckenhoupt weight}
(see, for instance,
\cite[Chapter 7]{LGCF}).

\begin{definition}\label{ap}
An $A_p(\rrn)$\emph{-weight}\index{$A_p(\rrn)$}
$\upsilon$, with $p\in[1,\infty)$,
is a locally integrable and nonnegative
function on $\rrn$ satisfying
that, when $p\in(1,\infty)$,
$$
\left[\upsilon\right]_{A_p(\rrn)}:=
\sup_{B\in\mathbb{B}}\left[\frac{1}{|B|}
\int_{B}\upsilon(x)\,dx\right]
\left\{\frac{1}{|B|}
\int_{B}\left[\upsilon(x)\right]
^{\frac{1}{1-p}}\,dx\right\}^{p-1}<\infty
$$
and
$$
\left[\upsilon\right]_{A_1(\rrn)}
:=\sup_{B\in\mathbb{B}}
\frac{1}{|B|}\int_{B}\upsilon(x)\,dx
\left[\left\|\upsilon^{-1}
\right\|_{L^\infty(B)}\right]
<\infty.
$$
Moreover, let $$\index{$A_{\infty}(\rrn)$}
A_{\infty}(\rrn):=\bigcup_{p\in[1,\infty)}A_p(\rrn).$$
\end{definition}

Then we have the following
extrapolation\index{extrapolation}
theorem about the local
generalized
Herz space $\HerzSo$.

\begin{theorem}\label{th175}
Let $p$, $q\in(1,\infty)$,
$r_0\in[1,\infty)$, and $\omega\in M(\rp)$
satisfy
$$
-\frac{n}{p}<\m0(\omega)
\leq\M0(\omega)<\frac{n}{p'}
$$
and
$$
-\frac{n}{p}<\mi(\omega)
\leq\MI(\omega)<\frac{n}{p'}.
$$
Assume $\mathcal{F}$ is a set of
all pairs of nonnegative
measurable functions $(F,\,G)$ such that,
for any given
$\upsilon\in A_{r_0}(\rrn)$,
$$
\int_{\rrn}\left[F(x)\right]^{r_0}\upsilon(x)\,dx
\leq C_{(r,[\upsilon]_{A_{r_0}(\rrn)})}
\int_{\rrn}\left[G(x)\right]^{r_0}\upsilon(x)\,dx,
$$
where the positive constant
$C_{(r_0,[\upsilon]_{A_{r_0}(\rrn)})}$
is independent of $(F,\,G)$, but depending on
both $r_0$ and $[\upsilon]_{A_{r_0}(\rrn)}$.
Then there exists a
positive constant $C$ such that, for any
$(F,\,G)\in\mathcal{F}$ with
$\|F\|_{\HerzSo}<\infty$,
$$
\left\|F\right\|_{\HerzSo}\leq
C\left\|G\right\|_{\HerzSo}.
$$
\end{theorem}

Recall that Tao et al.\ \cite[Lemma 2.13]{TYYZ}
obtained the following extrapolation
theorem of general ball Banach function spaces.

\begin{lemma}\label{th175l1}
Let $X$ be a ball Banach function space satisfying
that the Hardy--Littlewood maximal operator $\mc$
is bounded on both $X$ and its associate space $X'$,
and let $r_0\in[1,\infty)$.
Assume $\mathcal{F}$ is a set of
all pairs of nonnegative
measurable functions $(F,\,G)$ such that,
for any given
$\upsilon\in A_{r_0}(\rrn)$,
$$
\int_{\rrn}\left[F(x)\right]^{r_0}\upsilon(x)\,dx
\leq C_{(r,[\upsilon]_{A_{r_0}(\rrn)})}
\int_{\rrn}\left[G(x)\right]^{r_0}\upsilon(x)\,dx,
$$
where the positive constant
$C_{(r_0,[\upsilon]_{A_{r_0}(\rrn)})}$
is independent of $(F,\,G)$, but depending on
both $r_0$ and $[\upsilon]_{A_{r_0}(\rrn)}$.
Then there exists a
positive constant $C$ such that, for any
$(F,\,G)\in\mathcal{F}$ with
$\|F\|_{X}<\infty$,
$$
\left\|F\right\|_{X}\leq
C\left\|G\right\|_{X}.
$$
\end{lemma}

Therefore, in order to prove Theorem \ref{th175},
it suffices to show that, under the assumptions of
Theorem \ref{th175},
the local generalized Herz space $\HerzSo$
satisfies all the assumptions of Lemma \ref{th175l1}.
To achieve this, we need the following two
auxiliary lemmas about the boundedness of
the Hardy--Littlewood maximal operator on
local generalized Herz spaces and their associate
spaces.

\begin{lemma}\label{mbhl}
Let $p,\ q\in(0,\infty)$ and $\omega\in M(\rp)$ satisfy
$\m0(\omega)\in(-\frac{n}{p},\infty)$ and
$\mi(\omega)\in(-\frac{n}{p},\infty)$.
Then, for any given $r\in(0,\min\{p,\frac{n}{\Mw+n/p}\})$,
there exists a positive constant $C$ such that,
for any $f\in L^1_{\mathrm{loc}}(\rrn)$,
$$
\left\|\mc(f)\right\|_
{[\HerzSo]^{1/r}}\leq C\|f\|_{[\HerzSo]^{1/r}}.
$$
\end{lemma}

\begin{proof}
Let all the symbols be as in the present lemma.
Then, applying Lemma \ref{rela} and
the assumptions $r\in(0,\frac{n}{\Mw+n/p})$ and
$\mw\in(-\frac{n}{p},\infty)$, we conclude that
\begin{align}\label{mbhle1}
\Mwr&=r\Mw\notag\\
&<r\left(\frac{n}{r}-\frac{n}{p}\right)
=\frac{n}{(p/r)'}
\end{align}
and
\begin{align}\label{mbhle2}
\mwr=r\mw>-\frac{n}{p/r}.
\end{align}
Using this, the fact that $\frac{p}{r}\in(1,\infty)$,
and Corollary \ref{bhlo},
we find that the Hardy--Littlewood maximal operator
$\mc$ is bounded on the local generalized
Herz space $\HerzSocr$.
This, combined with Lemma \ref{convexll},
further implies that $\mc$ is bounded on $[\HerzSo]^{1/r}$,
which completes the proof of Lemma \ref{mbhl}.
\end{proof}

\begin{lemma}\label{mbhal}
Let $p,\ q\in(0,\infty)$ and
$\omega\in M(\rp)$ satisfy
$\m0(\omega)\in(-\frac{n}{p},\infty)$
and $\mi(\omega)\in(-\frac{n}{p},\infty)$.
Then, for any given $$s\in\left(0,
\min\left\{p,q,\frac{n}{\Mw+n/p}\right\}
\right)$$
and $$r\in\left(\max\left\{
p,\frac{n}{\mw+n/p}\right\},\infty\right],$$
the Herz space $[\HerzSo]^{1/s}$
is a ball Banach function space and
there exists a positive
constant $C$ such that, for any $f\in
L_{{\rm loc}}^{1}(\rrn)$,
\begin{equation}\label{mbhal1}
\left\|\mc^{((r/s)')}(f)\right\|_
{([\HerzSo]^{1/s})'}\leq C
\left\|f\right\|_{([\HerzSo]^{1/s})'}.
\end{equation}
\end{lemma}

\begin{proof}
Let all the symbols be as in the present lemma.
Then, by Lemma \ref{rela} and the assumptions that
\begin{align*}
s<\frac{n}{\Mw+n/p}\leq\frac{n}{\M0(\omega)+n/p}
\end{align*}
and $\m0(\omega)\in(-n/p,\infty)$,
we find that
\begin{align}\label{mbhal2}
\M0\left(\omega^{s}\right)=s\M0(\omega)
<s\left(\frac{n}{s}-\frac{n}{p}\right)=\frac{n}{(p/s)'}
\end{align}
and
\begin{align*}
\m0\left(\omega^{s}\right)=s\m0(\omega)>-\frac{n}{p/s}.
\end{align*}
These, together with the assumptions
$p/s,\ q/s\in(1,\infty)$,
and Theorem \ref{balll}, imply that
the local generalized Herz space $\HerzSocs$
is a ball Banach function space.
Combining this and Lemma \ref{convexll},
we further conclude that the
Herz space $[\HerzSo]^{1/s}$
is a ball Banach function space.

Next, we show \eqref{mbhal1}.
Indeed, from Lemma \ref{convexll},
\eqref{mbhal2}, and Theorem \ref{associate},
we deduce that
\begin{equation}\label{mbhal3}
\left([\HerzSo]^{1/s}\right)'=
\left(\Kmp_{\omega^{s},\0bf}
^{p/s,q/s}(\rrn)\right)'=\Kmp_{1/\omega^{s},
\0bf}^{(p/s)',(q/s)'}(\rrn).
\end{equation}
In addition,
by Lemma \ref{rela} and the assumption
$s\in(0,\frac{n}{\Mw+n/p})$,
we conclude that
\begin{align}\label{mbhal4}
&\mwms\notag\\&\quad=-s\Mw>-s\left(
\frac{n}{s}-\frac{n}{p}\right)=-\frac{n}{(p/s)'}.
\end{align}
On the other hand, from Lemma \ref{rela}
again and the assumption
$$r\in\left(\frac{n}{\mw+n/p},\infty\right),$$
it follows that
\begin{align*}
&\Mwms\\&\quad=-s\mw<-s
\left(\frac{n}{r}-\frac{n}{p}\right)=
\frac{n}{(r/s)'}-\frac{n}{(p/s)'},
\end{align*}
which, combined with $r\in(p,\infty)$,
further implies that
$$
\left(\frac{r}{s}\right)'<
\frac{n}{\Mwms+\frac{n}{(p/s)'}}
$$
and
$$
\left(\frac{r}{s}\right)'<
\left(\frac{p}{s}\right)'.
$$
Applying these, \eqref{mbhal4},
and Lemma \ref{mbhl},
we find that the operator $\mc$
is bounded on $[\Kmp_{1/\omega^{s},
\0bf}^{(p/s)',(q/s)'}(\rrn)]^{1/(r/s)'}$.
Thus, from both Remark \ref{power}
and \eqref{mbhal3},
it follows that,
for any $f\in
L_{{\rm loc}}^{1}(\rrn)$,
\begin{equation*}
\left\|\mc^{((r/s)')}(f)
\right\|_{([\HerzSo]^{1/s})'}\lesssim
\left\|f\right\|_{([\HerzSo]^{1/s})'}.
\end{equation*}
This finishes the proof of Lemma \ref{mbhal}.
\end{proof}

Via Lemmas \ref{th175l1} and \ref{mbhal},
we now show Theorem \ref{th175}.

\begin{proof}[Proof of Theorem \ref{th175}]
Let all the symbols be as in the present theorem.
Then, from the assumptions of the present
theorem and Theorem \ref{balll},
we deduce that the local generalized
Herz space $\HerzSo$ is a ball Banach function space.

Thus, to finish the proof of the present
theorem, we only need
to show that $\HerzSo$ satisfies all the assumptions
of Lemma \ref{th175l1}.
On the one hand, applying Corollary \ref{bhlo},
we find that the Hardy--Littlewood maximal
operator $\mc$ is bounded on $\HerzSo$, namely,
for any $f\in L^1_{\mathrm{loc}}(\rrn)$,
\begin{equation}\label{th175e1}
\left\|\mc(f)\right\|_{\HerzSo}
\lesssim\left\|f\right\|_{\HerzSo}.
\end{equation}
On the other hand, by the assumption
$$\Mw\in\left(-\frac{n}{p},\frac{n}{p'}\right),$$
we conclude that
\begin{equation*}
\frac{n}{\Mw+n/p}>\frac{n}{n(1/p'+1/p)}=1.
\end{equation*}
From this and the assumptions $p$, $q\in(1,\infty)$,
it follows that
$$
\min\left\{p,q,\frac{n}{\Mw+n/p}\right\}\in(1,\infty).
$$
This, combined with Lemma \ref{mbhal},
further implies that the Hardy--Littlewood maximal
operator $\mc$ is bounded on $(\HerzSo)'$, namely,
for any
$f\in L^1_{\mathrm{loc}}(\rrn)$,
\begin{equation}\label{th175e2}
\left\|\mc(f)\right\|_{(\HerzSo)'}
\lesssim\left\|f\right\|_{(\HerzSo)'}.
\end{equation}
Using this and \eqref{th175e1},
we conclude that $\HerzSo$ satisfies all
the assumptions of Lemma \ref{th175l1},
which completes the proof of Theorem \ref{th175}.
\end{proof}

Next, we show the following
extrapolation\index{extrapolation}
theorem of global generalized Herz
spaces.

\begin{theorem}\label{extrag}
Let $p$, $q\in(1,\infty)$,
$r_0\in[1,\infty)$, and $\omega\in M(\rp)$
satisfy
$$
-\frac{n}{p}<\m0(\omega)
\leq\M0(\omega)<\frac{n}{p'}
$$
and
$$
-\frac{n}{p}<\mi(\omega)
\leq\MI(\omega)<\frac{n}{p'}.
$$
Assume $\mathcal{F}$ is a set of
all pairs of nonnegative
measurable functions $(F,\,G)$ such that,
for any given
$\upsilon\in A_{r_0}(\rrn)$,
\begin{equation}\label{extrage0}
\int_{\rrn}\left[F(x)\right]^{r_0}\upsilon(x)\,dx
\leq C_{(r,[\upsilon]_{A_{r_0}(\rrn)})}
\int_{\rrn}\left[G(x)\right]^{r_0}\upsilon(x)\,dx,
\end{equation}
where the positive constant
$C_{(r_0,[\upsilon]_{A_{r_0}(\rrn)})}$
is independent of $(F,\,G)$, but depending on
both $r_0$ and $[\upsilon]_{A_{r_0}(\rrn)}$.
Then there exists a
positive constant $C$ such that, for any
$(F,\,G)\in\mathcal{F}$ with
$\|F\|_{\HerzS}<\infty$,
$$
\left\|F\right\|_{\HerzS}\leq
C\left\|G\right\|_{\HerzS}.
$$
\end{theorem}

To prove this theorem,
we require the following
property of Muckenhoupt
weights, which
can be immediately inferred from
Definition \ref{ap};
we omit the details.

\begin{lemma}\label{extragl1}
Let $p\in[1,\infty)$,
$\xi\in\rrn$, and
$\upsilon$ be a locally
integrable and nonnegative
function on $\rrn$.
Then $\upsilon\in A_p(\rrn)$
if and only if $\upsilon(\cdot+\xi)
\in A_p(\rrn)$. Moreover,
for any locally
integrable and nonnegative
function $\upsilon$,
$$
\left[\upsilon(\cdot+\xi)\right]_{
A_p(\rrn)}=\left[\upsilon\right]_{A_p(\rrn)}.
$$
\end{lemma}

\begin{proof}[Proof
of Theorem \ref{extrag}]
Let all the symbols be as
in the present theorem.
Then, for any given
$\upsilon\in A_{r_0}(\rrn)$
and any given $\xi\in\rrn$,
we claim that
$$
\int_{\rrn}\left[
F(x+\xi)
\right]^{r_0}\upsilon(x)\,dx\lesssim
\int_{\rrn}\left[G(x+\xi)\right]^{
r_0}\upsilon(x)\,dx,
$$
where the implicit
positive constant is
independent of both $(F,\,G)$ and $\xi$.
Indeed, from Lemma \ref{extragl1},
we deduce that $\upsilon(\cdot-\xi)
\in A_{r_0}(\rrn)$ and
\begin{equation}\label{extrage1}
\left[\upsilon(\cdot-\xi)
\right]_{A_{r_0}(\rrn)}
=\left[\upsilon\right]_{A_{r_0}(\rrn)}.
\end{equation}
By this and \eqref{extrage0},
we find that
$$
\int_{\rrn}\left[
F(x)
\right]^{r_0}\upsilon(x-\xi)\,dx\lesssim
\int_{\rrn}\left[G(x)\right]^{
r_0}\upsilon(x-\xi)\,dx,
$$
where the implicit
positive constant is
independent of $(F,\,G)$.
This, together with \eqref{extrage1}
again, further implies that
$$
\int_{\rrn}\left[
F(x+\xi)
\right]^{r_0}\upsilon(x)\,dx\lesssim
\int_{\rrn}\left[G(x+\xi)\right]^{
r_0}\upsilon(x)\,dx,
$$
where the implicit
positive constant is
independent of both $(F,\,G)$ and $\xi$,
which completes the proof
of the above claim.

On the other hand, applying the
assumption $\|F\|_{\HerzS}<\infty$
and Remark \ref{remhs}(ii),
we conclude that,
for any $\xi\in\rrn$,
$$\left\|F(\cdot+\xi)
\right\|_{\HerzSo}<\infty.$$
Combining this,
the above claim, Theorem \ref{th175},
and Remark \ref{remhs}(ii),
we further find that,
for any $\xi\in\rrn$,
$$
\left\|F(\cdot+\xi)\right\|_{\HerzSo}
\lesssim\left\|G(\cdot+\xi)
\right\|_{\HerzSo}\lesssim\|G\|_{\HerzS},
$$
where the implicit positive constants
are independent of $(F,\,G)$ and $\xi$.
From this and Remark \ref{remhs}(ii)
again,
it follows that
$$
\|F\|_{\HerzS}\lesssim
\|G\|_{\HerzS},
$$
where the implicit positive constant
is independent of $(F,\,G)$,
which completes the proof of Theorem \ref{extrag}.
\end{proof}

\chapter{Block Spaces and Their Applications}\label{sec4}
\markboth{\scriptsize\rm\sc Block Spaces
and Their Applications}
{\scriptsize\rm\sc Block Spaces and
Their Applications}

In this chapter, we first
introduce the concepts of
both the $(\omega,\,p)$-block
and the block space $\bspace$.
Then we establish
the duality between
block spaces and
global generalized Herz spaces,
which plays a key role
in the study of
the real-variable theory
of Hardy spaces
associated with global
generalized Herz spaces
in the subsequent chapters.
To this end, we
first characterize
local generalized Herz spaces
via $(\omega,\,p)$-blocks
and then establish
an equivalent characterization
of block spaces via
local generalized Herz spaces.
Using this equivalent
characterization and borrowing some
ideas from the proof of
\cite[Theorem 6.1]{gm13}, we show that
the global generalized
Herz space $\HerzSea$ is
just the dual space of the block space
$\bspace$ when $p$, $q\in(1,\infty)$.
As applications, we obtain
the boundedness of some sublinear
operators on block spaces.
To be precise, the boundedness
criteria of both sublinear operators
and Calder\'{o}m--Zygmund operators
are established.
In particular, the boundedness
of powered Hardy--Littlewood
maximal operators, which can
be deduced directly
from that of sublinear operators,
plays an important role
in the subsequent chapters.

\section{Block Spaces}

In this section, we introduce
block spaces. To begin with, we
introduce the following concept
of $(\omega,\,p)$-blocks.

\begin{definition}\label{blo}
Let $p\in(0,\infty)$ and $\omega\in M(\rp)$.
Then a measurable function $b$ on $\rrn$ is called
an \emph{$(\omega,
\,p)$-block}\index{$(\omega,\,p)$-block} if
there exists a cube $Q\in\mathcal{Q}$
such that $$\supp(b):=
\left\{x\in\rrn:\
b(x)\neq0\right\}\subset Q$$
and
$$
\left\|b\right\|_{L^p(\rrn)}\leq\left[
\omega\left(|Q|^{\frac{1}{n}}\right)\right]^{-1}.
$$
\end{definition}

Before presenting the definition
of block spaces, we need some notation.
Recall that, for
any $k\in\mathbb{Z}$,
$\mathcal{D}_k$\index{$\mathcal{D}_k$} is defined
to be the set
of all the standard dyadic cubes\index{dyadic cube} on
$\rrn$ of level $-k$,
namely,
\begin{equation}\label{dyadic-k}
\mathcal{D}_k:=\left\{2^k\left\{j+[0,1)^n\right\}:\
j\in\mathbb{Z}^n\right\}.
\end{equation}
Moreover, let $$\index{$\mathcal{D}$}
\mathcal{D}:=\bigcup\limits_{k\in\mathbb{Z}}
\mathcal{D}_k.$$
In what follows, for any cube $Q\in\mathcal{Q}$,
we always use $l(Q)$\index{$l(Q)$} to
denote its edge length. Furthermore,
for any $x\in\rrn$
and $l\in(0,\infty)$, the \emph{symbol}
$Q(x,l)$\index{$Q(x,l)$}
denotes the cube with center $x$ and edge length $l$.
Then, for any $k\in\mathbb{Z}$
and $\xi\in\rrn$, let\index{$\mathcal{Q}_k^{(\xi)}$}
$$\mathcal{Q}_k^{(\0bf)}
:=\mathcal{D}_{k-1}\cap\left[Q(\0bf,2^{k+1})\setminus
Q(\0bf,2^k)\right]$$
and
\begin{equation}\label{qmc}
\mathcal{Q}_k^{(\xi)}:=\left\{Q\in\mathcal{Q}:\
Q-\{\xi\}\in\mathcal{Q}_k^{(\0bf)}\right\},
\end{equation}
here and thereafter, $Q-\{\xi\}:=
\{x-\xi:\ x\in Q\}$.

\begin{remark}\label{remark215}
We can easily show that, for any $\xi\in\rrn$
and $k\in\mathbb{Z}$,
$\sharp\mathcal{Q}_k^{(\xi)}=2^{2n}-2^n$.
Here and thereafter, for any set $E$, $\sharp E$
denotes its \emph{cardinality}\index{cardinality}.
\end{remark}

Now, we introduce the
block space $\bspace$ as follows.

\begin{definition}\label{blos}
Let $p,\ q\in(0,\infty)$ and $\omega\in M(\rp)$.
Then the \emph{block
space}\index{block space} $\bspace$\index{$\bspace$} is
defined to be the set of all the
measurable functions $f$ on $\rrn$ such that
\begin{equation*}
f=\sum_{l\in\mathbb{N}}\sum_{k\in\mathbb{Z}}
\sum_{Q\in\mathcal{Q}_{k}^{(\xil)}}\lambda_{
\xil,k,Q}b_{\xil,k,Q}
\end{equation*}
almost everywhere in $\rrn$, and
$$
\left\{\sum_{l\in\mathbb{N}}
\left[\sum_{k\in\mathbb{Z}}
\sum_{Q\in\mathcal{Q}_k^{(\xil)}}
\lambda_{\xil,k,Q}^{q}
\right]^{\frac{1}{q}}\right\}<\infty,
$$
where $\{\xil\}_{l\in\mathbb{N}}
\subset\rrn$,
$\{\lambda_{\xil,k,Q}\}_{l\in\mathbb{N},
\,k\in\mathbb{Z},\,Q\in\mathcal{Q}
_{k}^{(\xil)}}\subset[0,\infty)$, and, for any
$l\in\mathbb{N}$, $k\in\mathbb{Z}$, and
$Q\in\mathcal{Q}_k^{(\xil)}$, $b_{\xil,k,Q}$
is an $(\omega,\,p)$-block supported in the cube
$Q$. Moreover, for any $f\in\Msc(\rrn)$,
\begin{displaymath}
\|f\|_{\bspace}:=\inf\left\{\sum_{l\in\mathbb{N}}
\left[\sum_{k\in\mathbb{Z}}
\sum_{Q\in\mathcal{Q}_k^{(\xil)}}
\lambda_{\xil,k,Q}^{q}
\right]^{\frac{1}{q}}\right\},
\end{displaymath}
where the infimum is taken over all
the decompositions of $f$ as above.
\end{definition}

It is easy to show that,
for any $p$, $q\in(0,\infty)$ and $\omega\in M(\rp)$,
the block space $\bspace$ is a linear space
equipped with the quasi-seminorm $\|\cdot\|_{\bspace}$
and, in particular,
for any $p$, $q\in[1,\infty)$ and $\omega\in M(\rp)$,
$\|\cdot\|_{\bspace}$ is a seminorm; we omit the details.

\section{Duality}

In this section, we show that the
dual space of the block space
$\bspace$ is the global generalized
Herz space $\HerzSea$ with $p$, $q\in(1,\infty)$.
Namely, the following conclusion holds true.

\begin{theorem}\label{pre}
Let $p,\ q\in(1,\infty)$
and $\omega\in M(\rp)$.
Then the dual space\index{dual space}
of $\bspace$, denoted by $\Bspaced$,
is $\HerzSea$ in the following sense:
\begin{enumerate}
 \item[{\rm(i)}] Let $g\in\HerzSea$.
 Then the linear functional
\begin{equation}\label{funl}
\phi_{g}:\ f\mapsto\phi_{g}(f):=
\int_{\rrn}f(y)g(y)\,dy,
\end{equation}
defined for any $f\in\bspace$,
is bounded on $\bspace$.
 \item[{\rm(ii)}]
 Conversely, any continuous linear functional
 on $\bspace$ arises as
 in \eqref{funl} with a unique $g\in\HerzSea$.
\end{enumerate}
Moreover, there exists a constant
$C\in[1,\infty)$ such that,
for any $g\in\HerzSea$, $$C^{-1}
\left\|g\right\|_{\HerzSea}\leq
\left\|\phi_{g}\right\|_{
\Bspaceds}\leq C\left\|g\right\|_{\HerzSea}.$$
\end{theorem}

To show this dual theorem, we first
establish the following equivalent
characterizations
of both the local generalized Herz spaces
$\HerzSx$ and the block space
$\bspace$, respectively, via $(\omega,\,p)$-blocks
and $\HerzSx$.

\begin{lemma}\label{prel1}
Let $p$, $q\in(0,\infty)$ and $\omega\in M(\rp)$.
Then
\begin{enumerate}
  \item[{\rm(i)}] for any given
  $\xi\in\rrn$, a measurable function
  $f$ belongs to $\HerzSx$ if and only if
  $$
  f=\sum_{k\in\mathbb{Z}}\sum_{Q\in\mathcal{Q}_k^{(\xi)}}
  \lambda_{k,Q}b_{k,Q}
  $$
  almost everywhere in $\rrn$, and
  $$\left[
  \sum_{k\in\mathbb{Z}}\sum_{Q\in\mathcal{Q}_k
  ^{(\xi)}}\lambda_{k,Q}^q\right]
  ^{\frac{1}{q}}<\infty,$$
where $\{\lambda_{k,Q}\}_{
  k\in\mathbb{Z},\,Q\in\mathcal{Q}_k^{(\xi)}}
  \subset[0,\infty)$ and,
  for any $k\in\mathbb{Z}$ and
  $Q\in\mathcal{Q}_k^{(\xi)}$, $b_{k,Q}$ is an
  $(\omega,\,p)$-block
  supported in the cube $Q$. Moreover,
  for any $f\in\HerzSx$,
  $$
  \left\|f\right\|_{\HerzSx}\sim\left[
  \sum_{k\in\mathbb{Z}}\sum_{Q\in\mathcal{Q}_k
  ^{(\xi)}}\lambda_{k,Q}^q\right]^{\frac{1}{q}},
  $$
  where the positive equivalence constants are independent
  of both $\xi$ and $f$;
  \item[{\rm(ii)}] a measurable function
  $f$ belongs to $\bspace$
  if and only if
  $$
  \normmm{f}_{\bspace}:=\inf\left\{
  \sum_{l\in\mathbb{N}}\left\|f_{\xil}
  \right\|_{\HerzSxl}\right\}<\infty,
  $$
  where the infimum is taken over all the sequences
  $\{\xil\}_{l\in\mathbb{N}}\subset\rrn$ and
  $\{f_{\xil}\}_{l\in\mathbb{N}}\subset\Msc(\rrn)$
  such that, for any $l\in\mathbb{N}$,
  $f_{\xil}\in\HerzSxl$ and
  $$
  f=\sum_{l\in\mathbb{N}}f_{\xil}
  $$
  almost everywhere in $\rrn$. Moreover, for any
  $f\in\bspace$,
  $$
  \|f\|_{\bspace}\sim\normmm{f}_{\bspace},
  $$
  where the positive equivalence constants
  are independent of $f$.
\end{enumerate}
\end{lemma}

\begin{proof}
Let $p$, $q\in(0,\infty)$ and
$\omega\in M(\rp)$. We first show (i).
Indeed, to prove the necessity of (i),
let $\xi\in\rrn$ be a fixed
point and $f\in\HerzSx$. Then, for any
$k\in\mathbb{Z}$ and $Q\in\mathcal{Q}_k^{(\xi)}$,
let $$\lambda_{k,Q}:=\omega\left(|Q|^{1/n}\right)
\left\|f\1bf_Q\right\|_{
L^p(\rrn)}.$$ Moreover, for any
$k\in\mathbb{Z}$ and $Q\in\mathcal{Q}_k^{(\xi)}$,
let $b_{k,Q}:=0$ when $f=0$ almost everywhere in $Q$;
otherwise, let
$$
b_{k,Q}:=\left[\omega\left(|Q|^{\frac{1}{n}}\right)
\right]^{-1}
\left[\left\|f\1bf_{Q}
\right\|_{L^p(\rrn)}\right]^{-1}f\1bf_Q.
$$
Therefore, from the definition
of $\mathcal{Q}_k^{(\xi)}$, it follows that
$$
f=\sum_{k\in\mathbb{Z}}f\1bf_{Q(\xi,\tka)\setminus
Q(\xi,\tk)}=\sum_{k\in\mathbb{Z}}
\sum_{Q\in\mathcal{Q}_k^{(\xi)}}f\1bf_Q
=\sum_{k\in\mathbb{Z}}\sum_{
Q\in\mathcal{Q}_{k}^{(\xi)}}\lambda_{k,Q}b_{k,Q}.
$$
On the other hand, applying Definition \ref{blo},
Remark \ref{remark215}, the fact that, for any
$k\in\mathbb{Z}$,
$$
\bigcup_{Q\in\mathcal{Q}_k^{(\xi)}}Q=Q(\xi,\tka)\setminus
Q(\xi,\tk)\subset B(\xi,\tka)\setminus
B(\xi,\tkm),
$$
and Lemma \ref{mono}, we conclude that, for any
$k\in\mathbb{Z}$ and $Q\in\mathcal{Q}_k^{(\xi)}$,
$b_{k,Q}$ is an $(\omega,\,p)$-block supported
in $Q$, and
\begin{align}\label{prel1e1}
\left[\sum_{k\in\mathbb{Z}}
\sum_{Q\in\mathcal{Q}_k^{(\xi)}}
\lambda_{k,Q}^q\right]^{\frac{1}{q}}
&=\left\{\sum_{k\in\mathbb{Z}}
\left[\omega(\tkm)\right]^q
\sum_{Q\in\mathcal{Q}_k^{(\xi)}}
\left\|f\1bf_{Q}\right\|_{L^p(\rrn)}
^q\right\}^{\frac{1}{q}}\notag\\
&\sim\left\{\sum_{k\in\mathbb{Z}}
\left[\omega(\tkm)\right]^q
\left\|f\1bf_{Q(
\xi,\tka)\setminus Q(\xi,\tk)}
\right\|_{L^p(\rrn)}^q\right\}^{\frac{1}{q}}\notag\\
&\lesssim\left\{\sum_{k\in\mathbb{Z}}
\left[\omega(\tka)\right]^q
\left\|f\1bf_{B(
\xi,\tka)\setminus B(\xi,\tk)}
\right\|_{L^p(\rrn)}^q\right\}^{\frac{1}{q}}\notag\\
&\quad+\left\{\sum_{k\in\mathbb{Z}}
\left[\omega(\tk)\right]^q
\left\|f\1bf_{B(
\xi,\tk)\setminus B(\xi,\tkm)}
\right\|_{L^p(\rrn)}^q\right\}^{\frac{1}{q}}\notag\\
&\sim\|f\|_{\HerzSx}<\infty,
\end{align}
where the implicit positive constants
are independent of both $\xi$ and $f$.
This finishes the proof of
the necessity of (i).

Conversely, we deal with
the sufficiency of (i). To this end, assume
$f\in\Msc(\rrn)$ satisfying that
$$
  f=\sum_{k\in\mathbb{Z}}\sum_{Q\in\mathcal{Q}_k^{(\xi)}}
  \lambda_{k,Q}b_{k,Q}
  $$
  almost everywhere in $\rrn$, and
  $$\left[
  \sum_{k\in\mathbb{Z}}\sum_{Q\in\mathcal{Q}_k
  ^{(\xi)}}\lambda_{k,Q}^q\right]
  ^{\frac{1}{q}}<\infty,$$
where $\{\lambda_{k,Q}\}_{
  k\in\mathbb{Z},\,Q\in\mathcal{Q}_k^{(\xi)}}
  \subset[0,\infty)$ and,
  for any $k\in\mathbb{Z}$ and
  $Q\in\mathcal{Q}_k^{(\xi)}$, $b_{k,Q}$ is an
  $(\omega,\,p)$-block supported in the cube $Q$.
  Then, from the fact that, for any
  $j\in\mathbb{Z}$,
  $$
  B(\xi,2^j)\setminus B(\xi,2^{j-1})\subset
  Q(\xi,2^{j+1})\setminus
  Q(\xi,2^{j-1}),
  $$
  Remark \ref{remark215}, Definition \ref{blo},
  and Lemma \ref{mono}, it follows that,
  for any $j\in\mathbb{Z}$,
  \begin{align}\label{prel1e6}
&\left\|f\1bf_{B(\xi,2^j)
\setminus B(\xi,2^{j-1})}\right\|_{
L^p(\rrn)}^q\notag\\
&\quad\leq\left\|\sum_{k=j}^{j+1}
\sum_{Q\in\mathcal{Q}_k^{(\xi)}}
\lambda_{k,Q}b_{k,Q}\right\|_{L^p(\rrn)}^q
\notag\\&\quad\lesssim\sum_{k=j}^{j+1}
\sum_{Q\in\mathcal{Q}_k^{(\xi)}}
\lambda_{k,Q}^q
\left\|b_{k,Q}\right\|_{L^p(\rrn)}^q\notag\\
&\quad\lesssim\left[\omega(2^{j-2})\right]^{-q}
\sum_{Q\in\mathcal{Q}_{j-1}^{(\xi)}}\lambda
_{j-1,Q}^q+
\left[\omega(2^{j-1})\right]^{-q}
\sum_{Q\in\mathcal{Q}_{j}^{(\xi)}}
\lambda_{j,Q}^q\notag\\
&\quad\sim\left[\omega(2^{j})\right]^{-q}
\sum_{k=j}^{j+1}
\sum_{Q\in\mathcal{Q}_{k}^{(\xi)}}
\lambda_{k,Q}^q.
  \end{align}
This further implies that
\begin{align}\label{prel1e2}
\|f\|_{\HerzSx}\lesssim
\left[\sum_{j\in\mathbb{Z}}\sum_{k=j}^{j+1}
\sum_{Q\in\mathcal{Q}_k^{(\xi)}}
\lambda_{k,Q}^q
\right]^{\frac{1}{q}}\lesssim
\left[\sum_{j\in\mathbb{Z}}
\sum_{Q\in\mathcal{Q}_j^{(\xi)}}
\lambda_{j,Q}^q
\right]^{\frac{1}{q}}<\infty,
\end{align}
where the implicit positive constants
are independent of both $\xi$ and $f$.
Thus, $f\in\HerzSx$ and hence the proof
of the sufficiency of (i) is completed.
Moreover, combining \eqref{prel1e1}
and \eqref{prel1e2}, we conclude that,
for any $f\in\HerzSx$,
$$
\left\|f\right\|_{\HerzSx}\sim\left[
\sum_{k\in\mathbb{Z}}
\sum_{Q\in\mathcal{Q}_k
^{(\xi)}}\lambda_{k,Q}^q
\right]^{\frac{1}{q}},
$$
where the positive equivalence constants
are independent
of both $\xi$ and $f$, which completes the proof of (i).

Next, we show (ii). We first consider the necessity
of (ii). For this purpose, let
$f\in\bspace$ satisfy
$$
f=\sum_{l\in\mathbb{N}}\sum_{k\in\mathbb{Z}}
\sum_{Q\in\mathcal{Q}_k^{(\xil)}}\lambda_
{\xil,k,Q}b_{\xil,k,Q}
$$
almost everywhere in $\rrn$, where
$\{\xil\}_{l\in\mathbb{N}}\subset\rrn$,
$\{\lambda_{\xil,k,Q}\}_
{l\in\mathbb{N},\,k\in\mathbb{Z},
\,Q\in\mathcal{Q}
_{k}^{(\xil)}}\subset[0,\infty)$, and, for any
$l\in\mathbb{N}$, $k\in\mathbb{Z}$,
and $Q\in\mathcal{Q}_k^{(\xil)}$, $b_{\xil,k,Q}$
is an $(\omega,\,p)$-block supported in the cube $Q$,
and
$$
\sum_{l\in\mathbb{N}}
\left[\sum_{k\in\mathbb{Z}}
\sum_{Q\in\mathcal{Q}_k^{(\xil)}}
\lambda_{\xil,k,Q}^q\right]
^{\frac{1}{q}}<\infty.
$$
This implies that, for any $l\in\mathbb{N}$,
$$
\left[\sum_{k\in\mathbb{Z}}
\sum_{Q\in\mathcal{Q}_k^{(\xil)}}
\lambda_{\xil,k,Q}^q\right]
^{\frac{1}{q}}<\infty.
$$
By this and (i), we find that,
for any $l\in\mathbb{N}$,
$$
f_{\xil}:=\sum_{k\in\mathbb{Z}}
\sum_{Q\in\mathcal{Q}_k^{(\xil)}}
\lambda_{\xil,k,Q}b_{\xil,k,Q}\in\HerzSxl
$$
and
$$
\|f_{\xil}\|_{\HerzSxl}
\sim\left[\sum_{k\in\mathbb{Z}}
\sum_{Q\in\mathcal{Q}_k^{(\xil)}}
\lambda_{\xil,k,Q}^q\right]^{\frac{1}{q}},
$$
where the positive equivalence constants are independent
of both $l$ and $f$. Therefore, we have
\begin{equation}\label{prel1e3}
\normmm{f}_{\bspace}\leq\sum_{l\in\mathbb{N}}
\|f_{\xil}\|_{\HerzSxl}\sim
\sum_{l\in\mathbb{N}}
\left[\sum_{k\in\mathbb{Z}}
\sum_{Q\in\mathcal{Q}_k^{(\xil)}}
\lambda_{\xil,k,Q}^q\right]^{\frac{1}{q}}
<\infty,
\end{equation}
where the positive equivalence constants
are independent of $f$.
This finishes the proof of the necessity of (ii).
Moreover, applying \eqref{prel1e3}
and the choice of the sequence
$\{\lambda_{\xil,k,Q}\}_{l\in\mathbb{N},
\,k\in\mathbb{Z},\,Q\in\mathcal{Q}
_k^{(\xil)}}$, we further
conclude that
\begin{equation}\label{prel1e4}
\normmm{f}_{\bspace}\lesssim\|f\|_{\bspace},
\end{equation}
where the implicit positive constant is independent
of $f$.

Conversely, we show the sufficiency of (ii).
Indeed, assume that $f\in\Msc(\rrn)$
satisfies that there exist $\{\xil\}
_{l\in\mathbb{N}}\subset\rrn$ and $\{f_{\xil}\}
_{l\in\mathbb{N}}\subset\Msc(\rrn)$ such that
$f=\sum_{l\in\mathbb{N}}
f_{\xil}$ almost everywhere in $\rrn$, and
$$
\sum_{l\in\mathbb{N}}\left\|
f_{\xil}\right\|_{\HerzSxl}<\infty.
$$
Then, for any
$l\in\mathbb{N}$,
$f_{\xil}\in\HerzSxl$. This, together with (i),
further implies that,
for any $l\in\mathbb{N}$,
$$
f_{\xil}=\sum_{k\in\mathbb{Z}}\sum_{Q\in\mathcal{Q}
_{k}^{(\xil)}}\lambda_{\xil,k,Q}b_{\xil,k,Q}
$$
almost everywhere in $\rrn$, and
$$
\|f_{\xil}\|_{\HerzSxl}\sim
\left[\sum_{k\in\mathbb{Z}}
\sum_{Q\in\mathcal{Q}_k^{(\xil)}}
\lambda_{\xil,k,Q}^q\right]^{\frac{1}{q}},
$$
where the positive equivalence constants are independent
of both $l$ and $f$, and for any
$k\in\mathbb{Z}$ and $Q\in\mathcal{Q}_{k}^{(\xil)}$,
$\lambda_{\xil,k,Q}\in[0,\infty)$
and $b_{\xil,k,Q}$ is an $(\omega,\,p)$-block
supported in Q.
From this and Definition \ref{blos},
we deduce that
$$
f=\sum_{l\in\mathbb{N}}
\sum_{k\in\mathbb{Z}}
\sum_{Q\in\mathcal{Q}_k^{(\xil)}}
\lambda_{\xil,k,Q}b_{\xil,k,Q}
$$
almost everywhere in $\rrn$, and
\begin{equation}\label{prel1e5}
\left\|f\right\|_{\bspace}
\leq\sum_{l\in\mathbb{N}}
\left[\sum_{k\in\mathbb{Z}}
\sum_{Q\in\mathcal{Q}_k^{(\xil)}}
\lambda_{\xil,k,Q}^q\right]^{\frac{1}{q}}
\sim\sum_{l\in\mathbb{N}}\left\|f_{\xil}
\right\|_{\HerzSxl}<\infty,
\end{equation}
where the positive equivalence constants
are independent of $f$. Thus, we conclude that $f\in
\bspace$, and hence complete
the proof of the sufficiency of (ii).
In addition, using \eqref{prel1e5}
and the choice of $\{f_{\xil}\}_
{l\in\mathbb{N}}$, we further find that
\begin{equation*}
\|f\|_{\bspace}\lesssim\normmm{f}_{\bspace},
\end{equation*}
where the implicit positive
constant is independent of $f$.
Combining this and \eqref{prel1e4}, we obtain
$$
\|f\|_{\bspace}\sim\normmm{f}_{\bspace},
$$
where the positive equivalence constants are independent
of $f$. This finishes the proof of (ii),
and hence of Lemma \ref{prel1}.
\end{proof}

By Lemma \ref{prel1}(ii),
we immediately obtain the following
relation between local generalized
Herz spaces and block spaces,
which is also important
in the proof of Theorem \ref{pre}; we omit the
details.

\begin{lemma}\label{ee}
Let $p$, $q\in(0,\infty)$
and $\omega\in M(\rp)$.
Then, for any given $\xi\in\rrn$,
$\HerzSx\subset\bspace$.
Moreover, there exists a
positive constant $C$, independent
of $\xi$, such that,
for any $f\in\HerzSx$,
$$
\|f\|_{\bspace}\leq
C\|f\|_{\HerzSx}.
$$
\end{lemma}

Via Lemmas \ref{prel1} and \ref{ee},
we now show Theorem \ref{pre}.

\begin{proof}[Proof of Theorem \ref{pre}]
Let $p$, $q$, and $\omega$
be as in the present theorem.
We first show (i). To this end,
let $g\in\HerzSea$.
Then, for any $f\in\bspace$,
from Lemma \ref{prel1}(ii),
it follows that
\begin{equation}\label{12}
\|f\|_{\bspace}^{\star}=
\inf\left\{\sum_{l\in\mathbb{N}}
\left\|f_{\xil}
\right\|_{\HerzSxl}\right\}<\infty,
\end{equation}
where the infimum is taken over
all the sequences $\{\xil\}_{l\in\mathbb{N}}
\subset\rrn$ and
$\{f_{\xil}\}_{l\in\mathbb{N}}
\subset\Msc(\rrn)$ such that,
for any $l\in\mathbb{N}$, $f_{\xil}\in
\HerzSxl$ and
$$
f=\sum_{l\in\mathbb{N}}f_{\xil}
$$
almost everywhere in $\rrn$.
Therefore, applying \eqref{duale0}
and Remark \ref{remark213},
we conclude that
\begin{align}\label{arg}
\left|\int_{\rrn}f(y)g(y)\,dy\right|
&\leq\sum_{l\in\mathbb{N}}
\int_{\rrn}|f_{\xil}(y)g(y)|\,dy\notag\\
&\leq\sum_{l\in\mathbb{N}}
\|f_{\xil}\|_{\HerzSxl}
\|g\|_{\Kmp_{1/\omega,
\xil}^{p',q'}(\rrn)}\notag\\
&\leq\|g\|_{\HerzSea}
\sum_{l\in\mathbb{N}}
\left\|f_{\xil}\right\|_{\HerzSxl},
\end{align}
which, together with \eqref{12}
and Lemma \ref{prel1}(ii),
further implies that
\begin{equation}\label{ope}
\left|\int_{\rrn}f(y)g(y)\,dy
\right|\leq\|f\|_{\bspace}^{\star}\|g\|_{\HerzSea}
\sim\|f\|_{\bspace}\|g\|_{\HerzSea},
\end{equation}
where the positive equivalence constants
are independent of both $f$ and $g$.
Thus, the linear functional $\phi_{g}$ defined
as in \eqref{funl} is bounded on $\bspace$ and
\begin{equation}\label{e0}
\|\phi_{g}\|_{\Bspaced}\lesssim\|g\|_{\HerzSea},
\end{equation}
where the implicit positive
constant is independent of $g$. This finishes the
proof of (i).

Next, we show (ii). For this purpose,
let $\phi\in\Bspaced$. Then,
using Lemma \ref{prel1}(ii)
and Lemma \ref{ee}, we find that,
for any given $\xi\in\rrn$
and for any $f\in\HerzSx$,
$f\in\bspace$ and
\begin{align}\label{pree0}
|\phi(f)|&\leq\|\phi\|_{\Bspaced}
\|f\|_{\bspace}\lesssim
\|\phi\|_{\Bspaced}\|f\|_{\HerzSx},
\end{align}
where the implicit positive constant
is independent of $f$.
This implies that, for any $\xi\in\rrn$,
$$\phi\in\left(\HerzSx\right)^{*}.$$
By this and Theorem \ref{dual},
we conclude that, for any
$\xi\in\rrn$, there exists a
$g_{\xi}\in\HerzSxd$ such that,
for any $f\in\HerzSx$,
\begin{equation}\label{pree1}
\phi(f)=\int_{\rrn}f(y)g_{\xi}(y)\,dy.
\end{equation}
Now, we show that, for any $\xi_{1},\ \xi_{2}\in
\rrn$, $g_{\xi_1}=g_{\xi_2}$ almost everywhere
in $\rrn$. Indeed, when $\xi_1=\xi_2$,
from Theorem \ref{dual}(ii), we deduce that
$g_{\xi_1}=g_{\xi_2}$ almost everywhere in $\rrn$.
On the other hand, when $\xi_1\neq\xi_2$,
applying \eqref{pree1}, we conclude that,
for any $f\in\Kmp_{\omega,\xi_{1}}^{p,q}(\rrn)
\cap\Kmp_{\omega,\xi_{2}}^{p,q}(\rrn)$,
$$
\int_{\rrn}f(y)g_{\xi_1}(y)\,dy=\phi(f)
=\int_{\rrn}f(y)g_{\xi_2}(y)\,dy,
$$
which further implies that
\begin{equation}\label{pree2}
\int_{\rrn}f(y)\left[g_{\xi_1}(y)-
g_{\xi_2}(y)\right]\,dy=0.
\end{equation}
For any $k\in\mathbb{N}$, let
$$
f_k:=\mathrm{sgn}\left(\overline{g_{\xi_1}-g_{\xi_2}}
\right)\1bf_
{B(\xi_{1},\tk)\setminus B(\xi_{1},2^{-k})}
\1bf_{[B(\xi_{2},\frac{|\xi_{1}-\xi_{2}|}{2})]^
\complement}.
$$
Then we have, for any $k\in\mathbb{N}$,
\begin{align*}
&\|f_k\|_{\Kmp_{\omega,\xi_{1}}^{p,q}(\rrn)}\\
&\quad\leq\left\{\sum_{j\in\mathbb{Z}}
\left[\omega(2^j)\right]^q
\left\|\1bf_
{B(\xi_{1},\tk)\setminus B(\xi_{1},2^{-k})}
\1bf_{B(\xi_{1},2^j)\setminus B(\xi_{1},2^{j-1})}
\right\|_{L^p(\rrn)}^q\right\}^{\frac{1}{q}}\\
&\quad=\left\{
\sum_{j=-k+1}^{k}\left[
\omega(2^j)\right]^q\left\|
\1bf_{B(\xi_{1},2^j)
\setminus B(\xi_{1},2^{j-1})}
\right\|_{L^p(\rrn)}^q\right\}^{\frac{1}{q}}<\infty,
\end{align*}
which implies that $f_k\in
\Kmp_{\omega,\xi_{1}}^{p,q}(\rrn)$.
Next, we prove that, for any $k\in\mathbb{N}$,
$f_k\in\Kmp_{\omega,\xi_{2}}^{p,q}(\rrn)$.
To this end, for any $k\in\mathbb{N}$,
let
$$
j_{k,1}:=\left\lfloor\frac{\ln(|\xi_{1}
-\xi_{2}|)}{\ln2}\right\rfloor
$$
and
$$
j_{k,2}:=\left\lceil\frac{\ln(|\xi_{1}
-\xi_{2}|+\tk)}{\ln2}\right\rceil.
$$
We claim that, for any given $k\in\mathbb{N}$
and for any $j\in\mathbb{Z}\cap(-\infty,
j_{k,1})$,
$$
B(\xi_{2},2^{j})\setminus
B(\xi_{2},2^{j-1})
\subset B\left(\xi_{2},\frac{|\xi_{1}-
\xi_{2}|}{2}\right)
$$
and, for any $j\in\mathbb{Z}\cap(j_{k,2},\infty)$,
$$
B(\xi_{2},2^{j})\setminus
B(\xi_{2},2^{j-1})\subset\left[B(\xi_{1},\tk)
\right]^\complement.
$$
Indeed, for any given
$k\in\mathbb{N}$ and for any $j\in\mathbb{Z}\cap(-\infty,
j_{k,1})$, we have
$$
j\leq j_{k,1}-1\leq\frac{\ln(|\xi_{1}
-\xi_{2}|)}{\ln2}-1,
$$
which implies that $2^{j}\leq\frac{|\xi_{1}
-\xi_{2}|}{2}$.
This further implies that, for any given
$k\in\mathbb{N}$ and for any
$j\in\mathbb{Z}\cap(-\infty,j_{k,1})$,
\begin{equation*}
B(\xi_{2},2^{j})\setminus
B(\xi_{2},2^{j-1})\subset B(\xi_{2},2^i)
\subset B\left(\xi_{2},\frac{|\xi_{1}-
\xi_{2}|}{2}\right).
\end{equation*}
In addition, for any given $k\in\mathbb{N}$
and for any $j\in\mathbb{Z}\cap
(j_{k,2},\infty)$, we have
$$
j\geq j_{k,2}+1\geq\frac{\ln(|\xi_{1}
-\xi_{2}|+\tk)}{\ln2}+1,
$$
which further implies that
$2^{j-1}\geq|\xi_{1}-\xi_{2}|+\tk$.
Using this, we conclude that, for any
given $k\in\mathbb{N}$ and for any $j\in
\mathbb{Z}\cap(j_{k,2},\infty)$ and $x\in
B(\xi_{2},2^j)\setminus
B(\xi_{2},2^{j-1})$,
$$
\left|x-\xi_{1}\right|\geq
\left|x-\xi_{2}\right|-\left|
\xi_{1}-\xi_{2}\right|
\geq2^{j-1}-\left|
\xi_{1}-\xi_{2}\right|\geq\tk.
$$
From this, we further deduce that,
for any
given $k\in\mathbb{N}$ and for any $j\in
\mathbb{Z}\cap(j_{k,2},\infty)$,
$$
B(\xi_{2},2^{j})\setminus
B(\xi_{2},2^{j-1})\subset
\left[B(\xi_{1},\tk)\right]^\complement
$$
and hence the above claim holds true.
Thus, we find that,
for any $k\in\mathbb{N}$,
\begin{align*}
\|f_k\|_{\Kmp_{\omega,\xi_{2}}^
{p,q}(\rrn)}\leq\left\{
\sum_{j=j_{k,1}}^{j_{k,2}}
\left[\omega(2^j)\right]^q
\left\|
\1bf_{B(\xi_{2},2^j)
\setminus B(\xi_{2},2^{j-1})}
\right\|_{L^p(\rrn)}^q\right\}^{\frac{1}{q}}<\infty,
\end{align*}
which implies that $f_k\in
\Kmp_{\omega,\xi_{2}}^{p,q}(\rrn)$.
Then, applying the fact that,
for any $k\in\mathbb{N}$,
$f_k\in\Kmp_{\omega,\xi_{1}}^
{p,q}(\rrn)\cap\Kmp_{\omega,\xi_{2}}^
{p,q}(\rrn)$ and \eqref{pree2} with $f$
therein replaced by $f_k$,
we conclude that, for any
$k\in\mathbb{N}$,
$$
\int_{[B(\xi_{1},\tk)\setminus
B(\xi_{1},2^{-k})]\cap[B(\xi_{2},
\frac{|\xi_{1}-\xi_{2}|}{2})]^\complement}
\left|g_{\xi_1}(y)-g_{\xi_2}(y)\right|\,dy=0.
$$
This, together with the arbitrariness of $k$,
further implies that
$g_{\xi_1}=g_{\xi_2}$ almost everywhere in
$\rrn\setminus B(\xi_{2},
\frac{|\xi_{1}-\xi_{2}|}{2})$.
Similarly, we can obtain $g_{\xi_1}=g_{\xi_2}$
almost everywhere in
$\rrn\setminus B(\xi_{1},
\frac{|\xi_{1}-\xi_{2}|}{2})$. Therefore,
we further find that $g_{\xi_1}=g_{\xi_2}$
almost everywhere in $\rrn$.

Next, we show that $$g:=g_{\0bf}\in\HerzSea.$$
Indeed, by  Theorem \ref{dual} and \eqref{pree0},
we conclude that, for any $\xi\in\rrn$,
\begin{align*}
\|g\|_{\HerzSxd}&=
\|g_{\xi}\|_{\HerzSxd}=
\|\phi\|_{(\HerzSx)^*}\\&=
\sup_{\|f\|_{\HerzSx}=1}\left|
\phi(f)\right|\lesssim
\|\phi\|_{\Bspaced},
\end{align*}
where the implicit positive constant
is independent of both $\xi$ and $g$.
This, together with the arbitrariness
of $\xi$,
further implies that $g\in\HerzSea$ and
\begin{equation}\label{e5}
\|g\|_{\HerzSea}\lesssim\|\phi\|_{\Bspaced},
\end{equation}
where the implicit positive constant
is independent of $g$.
Thus, the function $g$ can induce a
$\phi_{g}$ as in \eqref{funl}.

Next, we show that $\phi=\phi_g$.
Assume that $f\in\bspace$. Then, from
Lemma \ref{prel1}(ii), we infer that there exist
sequences $\{\xil\}_{l\in\mathbb{N}}\subset\rrn$
and $\{f_{\xil}\}_{l\in\mathbb{N}}
\subset\Msc(\rrn)$ such that, for any
$l\in\mathbb{N}$, $f_{\xil}\in\HerzSxl$,
$f=\sum_{l\in\mathbb{N}}
f_{\xil}$ almost everywhere in $\rrn$, and
\begin{equation}\label{pree3}
\sum_{l\in\mathbb{N}}\left\|
f_{\xil}\right\|_{\HerzSxl}\lesssim
\|f\|_{\bspace}+1.
\end{equation}
In addition, for any $N\in\mathbb{N}$,
let
$$f^{(N)}:=
\sum_{l=1}^{N}f_{\xil}.$$
Then, applying Lemma \ref{prel1}(ii) and
\eqref{pree3}, we find that, for any $N\in
\mathbb{N}$,
\begin{align*}
\left\|f^{(N)}\right\|_{\bspace}
&\sim\normmm{f^{(N)}}_{\bspace}
\lesssim\sum_{l=1}^N\left\|
f_{\xil}\right\|_{\HerzSxl}\\&\lesssim
\|f\|_{\bspace}+1<\infty,
\end{align*}
which implies that $f^{(N)}\in\bspace$.
Therefore, we have
\begin{align*}
\left\|
f-f^{(N)}\right\|_{\bspace}
&\sim\left\|
f-f^{(N)}\right\|_{\bspace}^{\star}\\
&\lesssim\sum_{l=N+1}
^{\infty}
\left\|f_{\xil}\right\|_{\HerzSxl}
\to0
\end{align*}
as $N\to\infty$. This, combined with the continuity
of $\phi$ on $\bspace$, \eqref{pree1}, \eqref{arg},
and the dominated convergence theorem, further implies
that
\begin{align*}
\phi(f)&=\lim\limits_{N\to\infty}\phi
\left(f^{(N)}\right)
=\lim\limits_{N\to\infty}
\int_{\rrn}\sum_{l=1}^N
f_{\xil}(y)g(y)\,dy\\&=\int_{\rrn}f(y)g(y)
\,dy=\phi_{g}(f),
\end{align*}
which completes the proof that $\phi=\phi_g$.
Moreover, from \eqref{e0} and \eqref{e5},
it follows that
$$
\|\phi_g\|_{\Bspaced}\sim\|g\|_{\HerzSea}
$$
with the
positive equivalence constants independent
of $g$.
By this and the linearity of $\phi_{g}$ about $g$,
we conclude that $g$ is unique.
This finishes the proof of (ii), and
hence of Theorem \ref{pre}.
\end{proof}

\section{Boundedness of Sublinear Operators}

The targets of this section
are twofold.
The first one is to establish a
boundedness criterion
of sublinear operators on
the block space $\bspace$
under some reasonable and
sharp assumptions.
To achieve this, we show the
lattice property of
block spaces (see Lemma \ref{equa}
below).
Using this lattice property and
the characterization of
block spaces via local generalized
Herz spaces obtained
in the last section,
we conclude the boundedness
of sublinear operators.
As an application, we
give the boundedness of
powered Hardy--Littlewood maximal
operators on block spaces,
which plays an important role
in the subsequent chapters.
The second target is
to investigate
the boundedness of
Calder\'{o}n--Zygmund operators
on block spaces.
Indeed, by embedding $\bspace$
into certain weighted Lebesgue space
(see Lemma \ref{czool1} below),
the boundedness of
Calder\'{o}n--Zygmund operators
can be concluded directly by
that of sublinear operators
just proved.

To begin with, we present
the boundedness criterion
of sublinear operators
on $\bspace$ as follows.

\begin{theorem}\label{sub}
Let $p\in(1,\infty)$, $q\in(0,\infty)$,
and $\omega\in M(\rp)$ satisfy
$$-\frac{n}{p}<\m0(\omega)
\leq\M0(\omega)<\frac{n}{p'}$$
and
$$-\frac{n}{p}<\mi(\omega)
\leq\MI(\omega)<\frac{n}{p'},$$
where $\frac{1}{p}+\frac{1}{p'}=1$.
Let $T$ be a bounded sublinear operator
on $L^p(\rrn)$ satisfying that
there exists a positive constant
$\widetilde{C}$ such that, for any
$f\in\HerzSo$ and $x\notin\supp(f)$,
\begin{equation*}
|T(f)(x)|\leq\widetilde{C}
\int_{\rrn}\frac{|f(y)|}{|x-y|^n}\,dy
\end{equation*}
and, for any $\{f_{j}\}_{j\in\mathbb{N}}
\subset\Msc(\rrn)$ and
almost every $x\in\rrn$,
\begin{equation}\label{conti}
\left|T\left(\sum_{j\in\mathbb{N}}
f_{j}\right)(x)\right|\leq\sum_{j
\in\mathbb{N}}|T(f_{j})(x)|.
\end{equation}
If $T$ is well defined on $\bspace$,
then there exists a positive constant
$C$ such that, for any
$f\in\bspace$,
$$
\left\|T(f)\right\|_{\bspace}
\leq C\|f\|_{\bspace}.
$$
\end{theorem}

In order to
show Theorem \ref{sub}, we
require the following lattice
property\index{lattice property}
of block spaces.

\begin{lemma}\label{equa}
Let $p,\ q\in(0,\infty)$ and
$\omega\in M(\rp)$.
Then a measurable function $f$
on $\rrn$ belongs to the block space
$\bspace$ if and only if
there exists a measurable function $g\in
\bspace$ such that $|f|\leq g$
almost everywhere in $\rrn$. Moreover,
for these $f$ and $g$,
$$\|f\|_{\bspace}\leq\|g\|_{\bspace}.$$
\end{lemma}

\begin{proof}
Let $p$, $q\in(0,\infty)$ and $\omega\in M(\rp)$.
We first show the necessity.
To this end, let $f\in\bspace$.
Then, by Definition \ref{blos},
we find that
$$f=\sum_{l\in\mathbb{N}}
\sum_{k\in\mathbb{Z}}\sum_{
Q\in\mathcal{Q}_k^{(\xil)}}\lambda_{\xil,
k,Q}b_{\xil,k,Q}$$
almost everywhere in $\rrn$, and
\begin{equation*}
\sum_{l\in\mathbb{N}}
\left[\sum_{k\in\mathbb{Z}}
\sum_{Q\in\mathcal{Q}_k^{(\xil)}}
\lambda_{\xil,k,Q}^{q}\right]
^{\frac{1}{q}}\in\left[\|f\|_{\bspace},
\|f\|_{\bspace}+1\right),
\end{equation*}
where $\{\xil\}_{l\in\mathbb{N}}
\subset\rrn$, $\{\lambda_{\xil,k,Q}\}_{
l\in\mathbb{N},\,k\in\mathbb{Z},\,Q\in
\mathcal{Q}_k^{(\xil)}}\subset[0,\infty)$,
and, for any $l\in\mathbb{N}$,
$k\in\mathbb{Z}$, and $Q\in\mathcal{Q}
_k^{(\xil)}$,
$b_{\xil,k,Q}$ is an $(\omega,\,p)$-block
supported in the cube $Q$.
Let
$$g:=\sum_{l\in\mathbb{N}}
\sum_{k\in\mathbb{Z}}
\sum_{Q\in\mathcal{Q}_k^{(\xil)}}
\lambda_{\xil,k,Q}|b_{\xil,k,Q}|.$$
Obviously, $|f|\leq g$ and,
for any $l\in\mathbb{N}$, $k\in\mathbb{Z}$,
and $Q\in\mathcal{Q}_k^{(\xil)}$,
the function $|b_{\xil,k,Q}|$ is
also an $(\omega,\,p)$-block supported
in $Q$.
This further implies that
\begin{align*}
\|g\|_{\bspace}\leq\sum_{l\in\mathbb{N}}
\left[\sum_{k\in
\mathbb{Z}}
\sum_{Q\in\mathcal{Q}_k^{(\xil)}}
\lambda_{\xil,k,Q}^{q}\right]^{\frac{1}{q}}
<\|f\|_{\bspace}+1<\infty.
\end{align*}
Thus, $g\in\bspace$ and hence
we finish the proof of the necessity.

Now, we show the sufficiency.
To achieve this, let $f\in\Msc(\rrn)$
and $g\in\bspace$ be such that
$|f|\leq g$ almost everywhere in $\rrn$.
Assume that $$g=\sum_{l\in\mathbb{N}}
\sum_{k\in\mathbb{Z}}\sum_{
Q\in\mathcal{Q}_k^{(\xil)}}\lambda_{\xil,
k,Q}b_{\xil,k,Q}$$ almost everywhere in $\rrn$,
where $\{\xil\}_{l\in\mathbb{N}}
\subset\rrn$, $\{\lambda_{\xil,k,Q}\}_{
l\in\mathbb{N},\,k\in\mathbb{Z},\,Q\in
\mathcal{Q}_k^{(\xil)}}\subset[0,\infty)$,
and, for any $l\in\mathbb{N}$,
$k\in\mathbb{Z}$, and $Q\in\mathcal{Q}
_k^{(\xil)}$,
$b_{\xil,k,Q}$ is an $(\omega,\,p)$-block
supported in the cube $Q$.
Then we have
$$\1bf_{\{y\in\rrn:\ g(y)\neq0\}}
=\sum_{l\in\mathbb{N}}
\sum_{k\in\mathbb{Z}}
\sum_{Q\in\mathcal{Q}_k^{(\xil)}}
\lambda_{\xil,k,Q}
\frac{1}{g}b_{\xil,k,Q}\1bf_{\{y\in\rrn:\
g(y)\neq0\}},$$
which, together with the assumption
$|f|\leq g$ almost
everywhere in $\rrn$, implies that
\begin{align*}f=f\1bf_{\{y\in\rrn:\
g(y)\neq0\}}
=\sum_{l\in\mathbb{N}}
\sum_{k\in\mathbb{Z}}
\sum_{Q\in\mathcal{Q}_k^{(\xil)}}
\lambda_{\xil,k,Q}
\frac{f}{g}b_{\xil,k,Q}\1bf_{\{y\in\rrn:\
g(y)\neq0\}}.
\end{align*}
Now, for any $l\in\mathbb{N}$,
$k\in\mathbb{Z}$, and $Q\in\mathcal{Q}_k^{(\xil)}$,
let $$b_{\xil,k,Q}^{\star}
:=\frac{f}{g}b_{\xil,k,Q}
\1bf_{\{y\in\rrn:\ g(y)\neq0\}}.$$
Obviously, for any $l\in\mathbb{N}$,
$k\in\mathbb{Z}$, and $Q\in\mathcal{Q}_k
^{(\xil)}$,
$\supp(b_{\xil,k,Q}^{\star})
\subset Q$. In addition, from the
assumption $|f|\leq g$
almost everywhere in $\rrn$, it follows that,
for any $l\in\mathbb{N}$, $k\in\mathbb{Z}$,
and $Q\in\mathcal{Q}_k^{(\xil)}$,
\begin{align*}
\left\|b_{\xil,k,Q}^{\star}
\right\|_{L^{p}(\rrn)}\leq\left\|b_
{\xil,k,Q}\right\|_{L^{p}(\rrn)}
\leq\left[\omega\left(
|Q|^{\frac{1}{n}}
\right)\right]^{-1},
\end{align*}
which further implies that,
for any $l\in\mathbb{N}$,
$k\in\mathbb{Z}$, and $Q\in\mathcal{Q}
_{k}^{(\xil)}$,
the function $b_{\xil,k,Q}^{\star}$
is an $(\omega,\,p)$-block supported
in $Q$, and hence
\begin{equation*}
\|f\|_{\bspace}\leq\sum_{l\in\mathbb{N}}
\left[\sum_{k\in\mathbb{Z}}
\sum_{Q\in\mathcal{Q}_k^{(\xil)}}
\lambda_{\xil,k,Q}^{q}\right]^{\frac{1}{q}}.
\end{equation*}
This, combined with the choice of
$\{\lambda_{\xil,k,Q}\}_{l\in\mathbb{N},
\,k\in\mathbb{N},\,Q\in\mathcal{Q}_k^{(\xil)}}$,
and the assumption $g\in\bspace$,
further implies that
\begin{align*}
\|f\|_{\bspace}\leq\|g\|_{\bspace}<\infty,
\end{align*}
which completes the proof of the
sufficiency, and hence of
Lemma \ref{equa}.
\end{proof}

\begin{remark}\label{remark233}
Let $p$, $q\in(0,\infty)$ and $\omega\in M(\rp)$.
We point out that a measurable function
$f$ on $\rrn$ belongs to
the block space $\bspace$ if and only if
the function $|f|$ belongs to $\bspace$.
Indeed, the sufficiency of this
conclusion is deduced
from Lemma \ref{equa} directly.
On the other hand,
let $f\in\bspace$ and $g$ be as
in the proof of the necessity of Lemma \ref{equa}.
Then, from the
proof of the necessity of Lemma \ref{equa},
it follows that
$g\in\bspace$. Applying this, the fact that
$|f|\leq g$, and
Lemma \ref{equa} again, we
conclude that $|f|\in\bspace$,
which completes the proof of the
necessity of the above claim.
Moreover, repeating the argument similar to that
used in the proof of Lemma \ref{equa},
we easily find that, for any $f\in\bspace$,
$$
\|f\|_{\bspace}=\left\|\,|f|
\,\right\|_{\bspace}.
$$
\end{remark}

Via the equivalent characterization
of block spaces
established in Lemma \ref{prel1}(ii),
we obtain the following boundedness
of sublinear operators on block spaces.

\begin{proposition}\label{subr}
Let $p$, $q\in(0,\infty)$, $\omega\in M(\rp)$,
and $T$ be a bounded sublinear
operator on $\HerzSo$ such that
\eqref{conti} holds true.
If $T$ is well defined
on $\bspace$, then
there exists a positive constant $C$ such that,
for any $f\in\bspace$,
$$
\left\|T(f)\right\|_{\bspace}\leq C\|f\|_{\bspace}.
$$
\end{proposition}

\begin{proof}
Let $p$, $q\in(0,\infty)$, $\omega\in M(\rp)$,
and $f\in\bspace$. Then,
from Lemma \ref{prel1}(ii) and
Remark \ref{remark213}, it follows that,
for any $k\in\mathbb{N}$,
there exist sequences $\{\xil^{(k)}\}
_{l\in\mathbb{N}}\subset\rrn$ and
$\{f_{\xil^{(k)}}\}_{l\in\mathbb{N}}
\subset\Msc(\rrn)$
such that $f=\sum_{l\in\mathbb{N}}
f_{\xil^{(k)}}$
almost everywhere in $\rrn$, and
\begin{equation}\label{2.344}
\sum_{l\in\mathbb{N}}\left\|f_{\xil^{
(k)}}\left(\cdot+\xil^{(k)}\right)
\right\|_{\HerzSo}\lesssim
\|f\|_{\bspace}+\frac{1}{k}.
\end{equation}
Moreover, by \eqref{conti}, we find that, for any
$k\in\mathbb{N}$ and almost every $x\in\rrn$,
\begin{equation}\label{2.354}
|T(f)(x)|\leq\sum_{l\in\mathbb{N}}
\left|T\left(f_{\xil^{(k)}}\right)(x)\right|.
\end{equation}
Applying the boundedness of the operator
$T$ on $\HerzSo$, we conclude that,
for any $k,\ l\in\mathbb{N}$,
$T(f_{\xil^{(k)}})(\cdot+\xil^{(k)})\in\HerzSo$ and
\begin{equation*}
\left\|T\left(f_{\xil^{(k)}}\right)
\left(\cdot+\xil^{(k)}\right)
\right\|_{\HerzSo}\lesssim
\left\|f_{\xil^{(k)}}\left(\cdot+\xil^{(k)}
\right)\right\|_{\HerzSo},
\end{equation*}
where the implicit positive constant
is independent of $k$, $l$, and $f$.
This, combined with
Lemma \ref{equa}, Lemma \ref{prel1}(ii),
\eqref{2.354}, and \eqref{2.344}, further
implies that, for any $k\in\mathbb{N}$,
\begin{align}\label{2.355}
\left\|T(f)\right\|_{\bspace}&\leq\left\|
\sum_{l\in\mathbb{N}}T\left(f_{\xil^{(k)}}
\right)\right\|_{\bspace}
\sim\normmm{\sum_{l\in\mathbb{N}}T
\left(f_{\xil^{(k)}}\right)}_{\bspace}\notag\\
&\lesssim\sum_{l\in\mathbb{N}}\left\|T\left(
f_{\xil^{(k)}}\right)\left(\cdot+
\xil^{(k)}\right)\right\|_{\HerzSo}\notag\\
&\lesssim\sum_{l\in\mathbb{N}}
\left\|f_{\xil^{(k)}}\left(\cdot+\xil^{(k)}
\right)\right\|_{\HerzSo}\notag\\
&\lesssim\|f\|_{\bspace}+\frac{1}{k},
\end{align}
where $\|\sum_{l\in\mathbb{N}}
T(f_{\xi_l^{(k)}})\|^{\star}_{\bspace}$
is defined as in Lemma \ref{prel1}(ii)
with $f$ therein replaced by
$\sum_{l\in\mathbb{N}}
T(f_{\xi_l^{(k)}})$, and
the implicit positive constant
is independent of both $f$ and $k$.
Letting $k\to\infty$ in \eqref{2.355}, we obtain
$$
\left\|T(f)\right\|_{\bspace}\lesssim\|f\|_{\bspace},
$$
where the implicit positive constant
is independent of $f$.
This finishes the proof of Proposition \ref{subr}.
\end{proof}

We now show Theorem \ref{sub}.

\begin{proof}[Proof of Theorem \ref{sub}]
Let all the symbols be as in the present theorem.
Then, combining Theorem \ref{vmaxh}
and Proposition \ref{subr},
we find that there exists
a positive constant $C$ such that,
for any $f\in\bspace$,
$$
\left\|T(f)\right\|_{\bspace}\leq C\|f\|_{\bspace}.
$$
This finishes the proof of Theorem \ref{sub}.
\end{proof}

As a consequence,
we now give the boundedness
of powered Hardy--Littlewood
maximal operators on block spaces
as follows.

\begin{corollary}\label{maxbl}
Let $p\in(1,\infty)$, $q\in(0,\infty)$,
$r\in[1,p)$, and
$\omega\in M(\rp)$ satisfy
$$
-\frac{n}{p}<\m0(\omega)\leq\M0(\omega)
<n\left(\frac{1}{r}-\frac{1}{p}\right)
$$
and
$$
-\frac{n}{p}<\mi(\omega)\leq\MI(\omega)
<n\left(\frac{1}{r}-\frac{1}{p}\right).
$$
Then there exists a positive constant
$C$ such that,
for any $f\in\bspace$,
$$
\left\|\mc^{(r)}(f)\right\|_{\bspace}
\leq C\|f\|_{\bspace},
$$
where $\mc^{(r)}$ is as in \eqref{hlmaxp}
with $\theta$ therein replaced
by $r$.
\end{corollary}

\begin{proof}
Let all the symbols be as in the present corollary.
We first show that the operator
$\mc^{(r)}$ satisfies \eqref{conti}.
Indeed, applying the
Minkowski inequality of $L^r(\rrn)$,
we know that, for any
$\{f_{j}\}_{j\in\mathbb{N}}
\subset\Msc(\rrn)$, any $x\in\rrn$,
and any ball $B\in\mathbb{B}$ containing $x$,
$$
\left[\int_{B}\left|\sum_{j\in\mathbb{N}}
f_{j}(y)\right|^{r}
\,dy\right]^{\frac{1}{r}}\leq
\sum_{j\in\mathbb{N}}\left[
\int_{B}\left|f_{j}(y)\right|^{r}
\,dy\right]^{\frac{1}{r}},
$$
which implies that
\begin{align*}
\left[\frac{1}{|B|}\int_{B}\left|
\sum_{j\in\mathbb{N}}
f_{j}(y)\right|^{r}\,dy\right]^{\frac{1}{r}}
&\leq\sum_{j\in\mathbb{N}}
\left[\frac{1}{|B|}
\int_{B}\left|f_{j}(y)\right|^{r}
\,dy\right]^{\frac{1}{r}}\\
&\leq\sum_{j\in\mathbb{N}}\mc^{(r)}(f_{j})(x).
\end{align*}
This further implies that, for any
$\{f_{j}\}_{j\in\mathbb{N}}
\subset\Msc(\rrn)$ and $x\in\rrn$,
$$
\mc^{(r)}\left(\sum_{j\in\mathbb{N}}f_{j}\right)(x)
\leq\sum_{j\in\mathbb{N}}\mc^{(r)}(f_{j})(x).
$$
Thus, the powered Hardy--Littlewood
maximal operator $\mc^{(r)}$ is a
sublinear operator which satisfies \eqref{conti}.

Next, we prove that $\mc^{(r)}$ is
bounded on the local
generalized Herz space $\HerzSo$.
Indeed, from Lemma \ref{convexll},
it follows that
$$
\left[\HerzSo\right]^{1/r}=\HerzSocr.
$$
In addition, applying Lemma \ref{rela} and
$p/r\in(1,\infty)$, we conclude that
\begin{align*}
\mwr=r\mw>-\frac{n}{p/r}
\end{align*}
and
\begin{align*}
\Mwr&=r\Mw\\&<r\left(\frac{n}{r}-\frac{n}{p}\right)
=\frac{n}{(p/r)'}.
\end{align*}
This, together with $p/r\in(1,\infty)$
and Corollary
\ref{bhlo}, further implies
that $\mc$ is bounded
on the local generalized Herz space $\HerzSocr$.
Therefore, the operator $\mc$ is
bounded on $[\HerzSo]^{1/r}$.
From this and Remark \ref{power}(i),
it follows that
the powered Hardy--Littlewood maximal
operator $\mc^{(r)}$
is bounded on $\HerzSo$. Using this,
the fact that $\mc^{(r)}$
satisfies \eqref{conti}, and
Proposition \ref{subr},
we then conclude that $\mc^{(r)}$
is bounded on $\bspace$,
which completes the proof of
Corollary \ref{maxbl}.
\end{proof}

Next, we devote to
establishing the boundedness
of Calder\'{o}n--Zygmund operators
on block spaces. To be precise,
we have the following conclusion.

\begin{theorem}\label{czoo}
Let $p,\ q\in(1,\infty)$
and $\omega\in M(\rp)$ with
$$
-\frac{n}{p}<\m0(\omega)
\leq\M0(\omega)<\frac{n}{p'}
$$
and
$$
-\frac{n}{p}<\mi(\omega)
\leq\MI(\omega)<\frac{n}{p'},
$$
where $\frac{1}{p}+\frac{1}{p'}=1.$
Assume $d\in\zp$
and $T$ is a $d$-order Calder\'{o}n--Zygmund
operator as in Definition \ref{defin-C-Z-s}.
Then $T$ is well defined on $\bspace$
and there exists a positive constant
$C$ such that, for any $f\in\bspace$,
$$
\left\|T(f)\right\|_{\bspace}\leq C\|f\|_{\bspace}.
$$
\end{theorem}

To prove Theorem \ref{czoo},
we need an embedding
lemma about block spaces.
First, we present the concepts
of weighted (weak) Lebesgue
spaces as follows.

\begin{definition}\label{df237}
Let $p\in(0,\infty)$ and
$\upsilon\in A_\infty(\rrn)$. Then
\begin{enumerate}
  \item[{\rm(i)}] the
\emph{weighted Lebesgue
space}\index{weighted Lebesgue space}
$L^p_\upsilon(\rrn)$\index{$L^p_\upsilon(\rrn)$}
is defined to be the set of
all the measurable functions
$f$ on $\rrn$ such that
$$
\left\|f\right\|_{L^p_\upsilon
(\rrn)}:=\left[\int_{\rrn}
\left|f(x)\right|^p\upsilon(x)
\,dx\right]^{\frac{1}{p}}<\infty;
$$
  \item[{\rm(ii)}] the
  \emph{weighted weak Lebesgue
space}\index{weighted weak Lebesgue space}
$WL^p_\upsilon(\rrn)$\index{$WL^p_\upsilon(\rrn)$}
is defined to be the set of
all the measurable functions
$f$ on $\rrn$ such that
$$
\left\|f\right\|_{WL^p_\upsilon
(\rrn)}:=\sup_{\lambda\in(0,\infty)}
\left\{\lambda
\left[\int_{\{x\in\rrn:\
|f(x)|>\lambda\}}\upsilon(y)\,dy
\right]^{\frac{1}{p}}\right\}<\infty.
$$
\end{enumerate}
\end{definition}

Note that Chang et al.
\cite[Lemma 4.7]{CWYZ} gave
an embedding lemma about
ball quasi-Banach function spaces.
However, due to
the deficiency of the
Fatou property\index{Fatou property}
[namely, whether
Definition \ref{Df1}(iii) holds true
for block spaces is unknown],
the block space $\bspace$ may
not be a ball quasi-Banach function
space. Fortunately,
we find that the proof
of \cite[Lemma 4.7]{CWYZ}
is still valid for Lemma \ref{czool1}
below; we omit the details.

\begin{lemma}\label{czool1}
Let $p$, $q$, and $\omega$ be
as in Theorem \ref{czoo}.
Then there exist an
$\eps\in(0,1)$ and a positive constant
$C$ such that, for any
$f\in\Msc(\rrn)$,
$$
\|f\|_{L^1_{\upsilon}(\rrn)}
\leq C\|f\|_{\bspace},
$$
where $\upsilon:=
[\mc(\1bf_{B(\0bf,1)})]^{\eps}$
with $\mc$ being the Hardy--Littlewood
maximal operator as in \eqref{hlmax}.
\end{lemma}

In order to prove Theorem \ref{czoo},
we also require the following
property of $A_1(\rrn)$-weights,
which is just
\cite[Theorem 7.7(1)]{d01}.

\begin{lemma}\label{czool2}
Let $f\in L^1_{\mathrm{loc}}(\rrn)$
satisfy $\mc(f)<\infty$ almost everywhere
in $\rrn$, and $\delta\in[0,1)$.
Then $\upsilon:=[\mc(f)]^{\delta}\in
A_1(\rrn)$.
\end{lemma}

The following conclusion shows
the weak boundedness of
Calder\'{o}n--Zygmund
operators on weighted Lebesgue spaces,
which is just \cite[Theorem 7.12]{d01}.

\begin{lemma}\label{czool3}
Let $\upsilon\in A_1(\rrn)$, $d\in\zp$,
and $T$ be a $d$-order
Calder\'{o}n--Zygmund
operator as in Definition
\ref{defin-C-Z-s}.
Then there exists a positive constant
$C$ such that, for any $f\in
L^1_{\upsilon}(\rrn)$
and $\lambda\in(0,\infty)$,
$$
\int_{\{x\in\rrn:\ |T(f)|>\lambda\}}
\upsilon(y)\,dy\leq\frac{C}{\lambda}
\|f\|_{L^1_{\upsilon}(\rrn)}
$$
and hence
$$
\left\|T(f)\right\|_{WL^1_{\upsilon}(\rrn)}
\leq C\|f\|_{L^1_{\upsilon}(\rrn)}.
$$
\end{lemma}

To prove Theorem \ref{czoo},
we also need the following
technique estimate about the
Hardy--Littlewood maximal
operator of characteristic functions
of balls (see, for instance, \cite[Example
2.1.8]{LGCF}).

\begin{lemma}\label{2.1.6}
Let $r\in(0,\infty)$ and
$\mc$ be the Hardy--Littlewood
maximal operator as in \eqref{hlmax}.
Then, for any $x\in\rrn$,
\begin{equation*}
\frac{r^n}{(r+|x|)^n}
\leq\mc\left(\1bf_{B(\0bf,r)}\right)(x)
\leq\frac{6^nr^n}{(r+|x|)^n}.
\end{equation*}
\end{lemma}

Via above preparations,
we now show Theorem \ref{czoo}.

\begin{proof}[Proof of Theorem \ref{czoo}]
Let all the symbols be
as in the present theorem
and $f\in\bspace$. Then, from
Lemma \ref{czool1}, we deduce that
$f\in L^1_{\upsilon}(\rrn)$, where
$\upsilon:=
[\mc(\1bf_{B(\0bf,1)})]^{\eps}$
with an $\eps\in(0,1)$.
Applying Lemma \ref{czool2},
we find that $\upsilon\in A_1(\rrn)$.
This, together with
the fact that $f\in L^1_{\upsilon}(\rrn)$
and Lemma \ref{czool3},
further implies that $T(f)\in WL^1_
{\upsilon}(\rrn)$ and hence
$T(f)$ is well defined.

In addition, using Lemma \ref{prel1}(ii),
we conclude that
there exist $\{\xil\}_{l\in\mathbb{N}}
\subset\rrn$ and $\{f_{\xil}\}_{l\in\mathbb{N}}
\subset\Msc(\rrn)$
such that, for any $l\in\mathbb{N}$,
$f_{\xil}\in\HerzSxl$,
$f=\sum_{l\in\mathbb{N}}f_{\xil}$ almost
everywhere in $\rrn$, and
\begin{equation}\label{czooe1}
\sum_{l\in\mathbb{N}}\left\|
f_{\xil}\right\|_{\HerzSxl}<\infty.
\end{equation}
We next prove that, for almost
everywhere $x\in\rrn$,
\begin{align}\label{czooe4}
\left|T(f)(x)\right|\leq
\sum_{l\in\mathbb{N}}\left|
T\left(f_{\xil}\right)(x)\right|.
\end{align}
To achieve this, for any $N\in\mathbb{N}$,
let $$f^{(N)}:=\sum_{l=1}^{N}f_{\xil}.$$
Then, from \eqref{czooe1},
it follows that, for any $N\in\mathbb{N}$,
$$
\sum_{l=1}^{N}\left\|
f_{\xil}\right\|_{\HerzSxl}
\leq\sum_{l\in\mathbb{N}}\left\|
f_{\xil}\right\|_{\HerzSxl}<\infty,
$$
which, combined with Lemma \ref{prel1}(ii)
again, implies that
$f^{(N)}\in\bspace$. By this, Lemma \ref{czool1},
the linearity of $T$, Lemmas \ref{czool3}
and \ref{prel1}(ii), and \eqref{czooe1},
we conclude that, for any
$N\in\mathbb{N}$, $$f-f^{(N)}\in
L^1_{\upsilon}(\rrn)$$ and, for any
$\lambda\in(0,\infty)$,
\begin{align}\label{czooe2}
&\int_{\{x\in\rrn:\
|T(f)(x)-T(f^{(N)})(x)|>\lambda\}}
\left[\mc\left(\1bf_{B(\0bf,1)}\right)(y)\right]
^{\eps}\,dy\notag\\
&\quad\lesssim\frac{1}{\lambda}
\left\|f-f^{(N)}\right\|_
{L^1_{\upsilon}(\rrn)}\lesssim
\frac{1}{\lambda}
\left\|f-f^{(N)}\right\|_{\bspace}\notag\\
&\quad\lesssim\frac{1}{\lambda}
\sum_{l=N+1}^{\infty}\left\|
f_{\xil}\right\|_{\HerzSxl}\to0
\end{align}
as $N\to\infty$.

Now, we claim that there exists
a subsequence $\{f^{(N_j)}\}_{
j\in\mathbb{N}}\subset\{f^{(N)}\}_{
N\in\mathbb{N}}$ such that
$T(f^{(N_j)})\to T(f)$ almost
everywhere in $\rrn$ as $j\to\infty$.
Indeed, applying Lemma \ref{2.1.6}
with $r:=1$,
we find that, for any $x\in\rrn$,
$$
\left[\mc\left(
\1bf_{B(\0bf,1)}\right)(x)\right]^{\eps}
\sim\left(1+|x|\right)^{-\eps n}.
$$
By this and \eqref{czooe2},
we conclude that, for any $\lambda\in(0,\infty)$,
\begin{align*}
&\left|\left\{x\in B(\0bf,1):\
\left|T(f)(x)-T\left(f^{(N)}\right)(x)
\right|>\lambda\right\}\right|\\
&\quad\lesssim\int_{\{x\in\rrn:\
|T(f)(x)-T(f^{(N)})(x)|>\lambda\}}
\left[\mc\left(\1bf_{B(\0bf,1)}\right)(y)\right]
^{\eps}\,dy\to0
\end{align*}
as $N\to\infty$ and, for any
$k\in\mathbb{N}$,
\begin{align*}
&\left|\left\{x\in B(\0bf,\tk)\setminus
B(\0bf,\tkm):\
\left|T(f)(x)-T\left(f^{(N)}\right)(x)
\right|>\lambda\right\}\right|\\
&\quad\lesssim2^{k\eps n}\int_{\{x\in\rrn:\
|T(f)(x)-T(f^{(N)})(x)|>\lambda\}}
\left[\mc\left(\1bf_{B(\0bf,1)}\right)(y)\right]
^{\eps}\,dy\to0
\end{align*}
as $N\to\infty$. These further imply
that the sequence
$\{T(f^{(N)})\}_{N\in\mathbb{N}}$
converges in measure to $T(f)$ both in
$B(\0bf,1)$ and $B(\0bf,\tk)\setminus
B(\0bf,\tkm)$ for any
$k\in\mathbb{N}$. From this,
the Riesz theorem\index{Riesz Theorem}, and
a diagonalization
argument\index{diagonalization argument},
we deduce that there exists
a subsequence $\{f^{(N_j)}\}_{
j\in\mathbb{N}}\subset\{f^{(N)}\}_{
N\in\mathbb{N}}$ such that
$T(f^{(N_j)})\to T(f)$ almost
everywhere in $\rrn$ as $j\to\infty$,
which completes the proof of
the above claim.

On the other hand,
using the linearity of $T$,
we conclude that, for any $j\in\mathbb{N}$,
\begin{align}\label{czooe3}
\left|T\left(f^{(N_j)}\right)
\right|=\left|\sum_{l=1}^{N_j}
T\left(f_{\xil}\right)(x)\right|
\leq\sum_{l=1}^{N_j}\left|
T\left(f_{\xil}\right)\right|
\leq\sum_{l\in\mathbb{N}}\left|
T\left(f_{\xil}\right)\right|.
\end{align}
By the above claim and
letting $j\to\infty$ in \eqref{czooe3},
we further find that, for almost
every $x\in\rrn$,
$$
\left|T(f)(x)\right|\leq
\sum_{l\in\mathbb{N}}\left|
T\left(f_{\xil}\right)(x)\right|.
$$
This finishes the proof of \eqref{czooe4}.
Combining this,
Corollary \ref{czoh}, and an argument
similar to that
used in the proof of \eqref{2.355},
we conclude that
$$
\left\|T(f)\right\|_{\bspace}
\lesssim\|f\|_{\bspace},
$$
where the implicit positive constant
is independent of $f$, which
completes the proof of
Theorem \ref{czoo}.
\end{proof}

\chapter{Boundedness and Compactness Characterizations
of Commutators on Generalized Herz Spaces}\label{sec5}
\markboth{\scriptsize\rm\sc Boundedness and Compactness
Characterizations of Commutators}
{\scriptsize\rm\sc Boundedness Characterization
and Compactness Characterization of Commutators}

In this chapter, we investigate
the boundedness and the compactness
characterizations of commutators on
local and global generalized
Herz spaces.
Recall that the commutator plays an
important role in various branch
of mathematics, such as
harmonic analysis (see, for instance,
\cite{AN18,AN20,CH,CD,GOS,NS,SS})
and partial differential equations
(see, for instance, \cite{CDD,CDH,TYYa}). In 1976,
Coifman et al. \cite{CRW} first obtained
a boundedness
characterization of commutators on the Lebesgue
space $L^p(\rrn)$ with $p\in(1,\infty)$. To be precise,
let $b\in L^1_{\mathrm{loc}}(\rrn)$
and $T_{\Omega}$ be a singular
integral operator with homogeneous kernel $\Omega$.
Coifman et al. \cite{CRW} proved that,
for any given $p\in(1,\infty)$
and any given $b\in\BMO$,
the commutator $[b,T_{\Omega}]$ is
bounded on $L^p(\rrn)$ and also that,
if, for any given $b\in
L_{\mathrm{loc}}^1(\rrn)$ and any
Riesz transform $R_{j}$ with $j\in\{1,\ldots,n\}$,
$[b,R_{j}]$ is bounded on $L^p(\rrn)$,
then $b\in\BMO$. Furthermore, in 1978,
Uchiyama \cite{Uch} showed that,
for any given $p\in(1,\infty)$,
the commutator $[b,T_{\Omega}]$
is bounded on $L^p(\rrn)$
if and only if $b\in\BMO$.
Later on, as extensions of the
results in Lebesgue spaces,
such boundedness characterizations
were established on
various function spaces (see \cite{DiR}
for Morrey spaces,
\cite{LDY} for weighted Lebesgue spaces,
and \cite{KL} for variable Lebesgue spaces).

On the other hand, the compactness
characterizations of commutators
were also studied. In 1978,
Uchiyama \cite{Uch} first showed that,
for any given $p\in(1,\infty)$,
the commutator $[b,T_{\Omega}]$ is compact
on $L^p(\rrn)$ if and only if
$b\in\CMO$. After that,
this compactness characterization
was extended to Morrey spaces
(see \cite{CDW}) and to weighted
Lebesgue spaces (see \cite{CC,GWY}).
Moreover, very recently,
Tao et al. \cite{TYYZ} studied
the boundedness and the
compactness characterizations of commutators
on ball Banach function spaces.

The main target of this chapter is to study the
boundedness and the
compactness characterizations of commutators
on generalized Herz spaces.
As was proved in Chapter \ref{sec3},
the generalized Herz spaces
are special ball Banach function spaces under
some reasonable and sharp assumptions
of exponents. Thus,
we can obtain the boundedness and the
compactness characterizations
of commutators on
generalized Herz spaces via
proving that generalized Herz spaces
satisfy all the assumptions
of the results obtained in \cite{TYYZ}.
This approach is feasible for
the boundedness and the compactness
characterizations of commutators on
local generalized Herz spaces.
However, since the
associate space of the global
generalized Herz space $\HerzS$ is still
unknown, we can not establish
the boundedness and the compactness
characterizations of commutators on
$\HerzS$ by the aforementioned method.
To overcome this difficulty,
we replace the
assumptions of conclusions in \cite{TYYZ}
about associate spaces
by other assumptions about ball
Banach function spaces
which are more convenient to check
for $\HerzS$.
In particular, we point
out that the extrapolation theorem
of $\HerzS$
established in Chapter \ref{sec3}
is a key tool
in the proof of the boundedness
and the compactness characterizations
of commutators on $\HerzS$.

We now recall some basic concepts.
Let $\Omega$ be a Lipschitz
function\index{Lipschitz function}
on the unit sphere of
$\rrn$, which is homogeneous of degree
zero\index{homogeneous of degree zero} and
has mean value zero\index{mean value zero}, namely,
for any $x$,
$y\in\mathbb{S}^{n-1}$ and $\mu\in(0,\infty)$,
\begin{equation}\label{lips}
|\Omega(x)-\Omega(y)|\leq|x-y|,
\end{equation}
\begin{equation}\label{degrz}
\Omega(\mu x)=\Omega(x),
\end{equation}
and
\begin{equation}\label{mz}
\int_{\mathbb{S}^{n-1}}\Omega(x)\,d\sigma(x)=0,
\end{equation}
here and thereafter, $\mathbb{S}^{n-1}:=\{x\in\rrn:\
|x|=1\}$ denotes
the unit sphere of $\rrn$ and $d\sigma$ the area
 measure on $\mathbb{S}^{n-1}$.
Furthermore, we state the
following $L^\infty$-Dini condition
(see, for instance, \cite[p.\,93]{d01}).

\begin{definition}\label{dini}
A function $\Omega\in L^{\infty}
(\mathbb{S}^{n-1})$ is said
to satisfy the
\emph{$L^{\infty}$-Dini
condition}\index{$L^{\infty}$-Dini condition} if
$$
\int_{0}^{1}\frac{\omega_{\infty}
(\tau)}{\tau}\,d\tau<\infty,
$$
where $\omega_{\infty}$\index{$\omega_{\infty}$}
is defined by setting, for any $\tau\in(0,1)$,
$$
\omega_{\infty}(\tau):=\sup_{\{x,\
y\in\mathbb{S}^{n-1}:\
|x-y|<\tau\}}|\Omega(x)-\Omega(y)|.
$$
\end{definition}

Assuming that $\Omega$
satisfies \eqref{degrz}, \eqref{mz},
and the $L^{\infty}$-Dini condition,
a linear operator $T_{\Omega}$ is called
a \emph{singular
integral operator with homogeneous
kernel $\Omega$}\index{singular
integral operator with homogeneous kernel}
(see, for instance,
\cite[p.\,53, Corollary 2.1.1]{LDY}) if,
for any $f\in L^p(\rrn)$ with
$p\in[1,\infty)$, and for any $x\in\rrn$,
\begin{align}\label{czcc}
T_{\Omega}(f)(x):&={\rm p.\,v.}\int_{\rrn}\frac{
\Omega(x-y)}{|x-y|^n}f(y)\,dy\notag\\
:&=\lim\limits_{\eps\to0^+}
\int_{\eps<|x-y|<1/\eps}\frac{\Omega(x-y)}
{|x-y|^{n}}f(y)\,dy.
\end{align}
Now, we recall the concept
of the commutator $[b,T_{\Omega}]$.
For any given $b\in L^1_{{\rm loc}}(\rrn)$,
the commutator\index{commutator}
$[b,T_{\Omega}]$\index{$[b,T_{\Omega}]$}
is defined by setting, for any bounded function
$f$ with compact support, and for
any $x\in\rrn$,
\begin{equation}\label{comu}
[b,T_{\Omega}](f)(x):=b(x)T_{\Omega}(f)(x)-
T_{\Omega}(bf)(x).
\end{equation}

Finally, recall that the \emph{bounded
mean oscillation function
space}\index{bounded mean
oscillation function space}
$\BMO$\index{$\BMO$},
introduced by John and Nirenberg \cite{jn61},
is defined to be the set of all the
$f\in L^{1}_{\rm{loc}}(\rrn)$ such that
\begin{equation}\label{fbmo}
\|f\|_{\BMO}:=\sup_{B\in\mathbb{B}}\left\{\frac{1}{|B|}
\int_{B}|f(x)-f_{B}|\,dx\right\}<\infty,
\end{equation}
where the supremum is taken over
all the balls $B\in\mathbb{B}$
and, for any given $B\in\mathbb{B}$,\index{$f_B$}
\begin{equation}\label{fb}
f_{B}:=\frac{1}{|B|}\int_{B}f(x)\,dx.
\end{equation}

\section{Boundedness Characterizations}

The target of this section
is to establish the boundedness
characterization of commutators
on local and global generalized Herz spaces.
We first consider commutators
on local generalized Herz spaces.
Namely, we show the following conclusion.

\begin{theorem}\label{iffbddl}
Let $p,\ q\in(1,\infty)$ and
$\omega\in M(\rp)$ satisfy
$$
-\frac{n}{p}<\m0(\omega)
\leq\M0(\omega)<\frac{n}{p'}
$$
and
$$
-\frac{n}{p}<\mi(\omega)
\leq\MI(\omega)<\frac{n}{p'},
$$
where $\frac{1}{p}+\frac{1}{p'}=1.$
Assume that $\Omega$ is a homogeneous Lipschitz
function of degree zero on $\mathbb{S}^{n-1}$
satisfying \eqref{mz},
$T_{\Omega}$ a singular integral
operator with homogeneous kernel $\Omega$,
and $b\in L^1_{{\rm loc}}(\rrn)$.
Then the commutator $[b,T_{\Omega}]$
is bounded on the local generalized Herz space
$\HerzSo$ if and only if $b\in\BMO$.
Moreover, there exists a constant $C\in[1,\infty)$,
independent of $b$, such that
$$
C^{-1}\|b\|_{\BMO}\leq\|[b,T_{\Omega}]\|_{\HerzSo
\to\HerzSo}\leq C\|b\|_{\BMO},
$$
where $\|[b,T_{\Omega}]\|_{\HerzSo
\to\HerzSo}$ is as in
\eqref{eon} with $T$ and $X$
therein replaced, respectively,
by $[b,T_{\Omega}]$ and $\HerzSo$.
\end{theorem}

First, we show the sufficiency
of Theorem \ref{iffbddl}. Indeed,
we have the following
conclusion for the sufficiency
of the boundedness of commutators
on local generalized Herz spaces,
where the assumption about $\Omega$
is weaker than that of Theorem
\ref{iffbddl}.

\begin{proposition}\label{suffbddl}
Let $p,\ q\in(1,\infty)$, $\omega\in M(\rp)$
be as in Theorem \ref{iffbddl},
and $r\in(1,\infty]$.
Assume that $b\in\BMO$, $\Omega\in
L^r(\mathbb{S}^{n-1})$ satisfies both
\eqref{degrz} and \eqref{mz},
and $T_{\Omega}$ is a singular integral
operator with homogeneous kernel $\Omega$.
Then there exists a positive constant $C$
such that, for any $f\in\HerzSo$,
$$
\left\|[b,T_{\Omega}](f)\right\|_{\HerzSo}
\leq C\|b\|_{\BMO}\|f\|_{\HerzSo}.
$$
\end{proposition}

To prove Proposition \ref{suffbddl},
we require the following
boundedness of commutators on
ball Banach function spaces,
which was obtained
by Tao et al.\ in
\cite[Theorem 2.17]{TYYZ}.

\begin{lemma}\label{suffbddll1}
Let $X$ be a
ball Banach function space
satisfying that the
Hardy--Littlewood maximal operator
$\mc$ is bounded on both $X$ and its associate
space $X'$, and $r\in(1,\infty]$.
Assume that $b\in\BMO$, $\Omega\in
L^r(\mathbb{S}^{n-1})$ satisfies
both \eqref{degrz} and \eqref{mz},
and $T_{\Omega}$ is a singular
integral operator
with homogeneous kernel $\Omega$.
Then there exists a positive constant
$C$ such that, for any $f\in X$,
$$
\left\|\left[b,T_{\Omega}\right](f)
\right\|_{X}\leq C\|b\|_{\BMO}\|f\|_{X}.
$$
\end{lemma}

Via this general conclusion for
ball Banach function spaces,
we now show Proposition \ref{suffbddl}.

\begin{proof}[Proof of
Proposition \ref{suffbddl}]
Let all the symbols be as in
the present proposition.
Then, from Theorem \ref{balll},
we deduce that, under the assumptions
of the present proposition,
the local generalized
Herz space $\HerzSo$ is a
BBF space.

Therefore, in order to complete
the proof of the present proposition,
it suffices to show that $\HerzSo$ satisfies
all the assumptions
of Lemma \ref{suffbddll1}.
Indeed, applying the assumptions
of the present proposition,
\eqref{th175e1}, and \eqref{th175e2},
we find that the Hardy--Littlewood
maximal operator $\mc$ is bounded on
$\HerzSo$ and $(\HerzSo)'$.
This further implies that
$\HerzSo$ satisfies all
the assumptions of Lemma \ref{suffbddll1},
which completes the proof
of Proposition \ref{suffbddl}.
\end{proof}

Next, we deal with
the necessity
of the boundedness
of commutators on local
generalized Herz spaces.
Notice that
the following proposition
gives a more general conclusion
than the necessity
of Theorem \ref{iffbddl}.

\begin{proposition}\label{necebddl}
Let $p\in(1,\infty)$, $q\in[1,\infty)$,
and $\omega\in M(\rp)$
be as in Theorem \ref{iffbddl}.
Assume that $b\in L^1_{{\rm loc}}(\rrn)$
and $\Omega\in
L^\infty(\mathbb{S}^{n-1})$ satisfies that
there exists an open set $\Lambda\subset
\mathbb{S}^{n-1}$ such that $\Omega$
never vanishes and
never changes sign on $\Lambda$.
If the commutator $[b,T_{\Omega}]$
is bounded on the
local generalized Herz space $\HerzSo$
and, for any bounded
measurable set $F\subset\rrn$ and
almost every $x\in\rrn\setminus\overline{F}$,
\begin{equation}\label{necebddle1}
\left[b,T_\Omega\right]\left(\1bf_{F}\right)(x)
=\int_{F}\left[b(x)-b(y)\right]\frac{
\Omega(x-y)}{|x-y|^n}\,dy,
\end{equation}
then $b\in
\BMO$ and there exists a positive constant $C$,
independent of $b$, such that
$$
\|b\|_{\BMO}\leq C\|[b,T_{\Omega}]
\|_{\HerzSo\to\HerzSo},
$$
where $\|[b,T_{\Omega}]
\|_{\HerzSo\to\HerzSo}$ is as in
\eqref{eon} with $T$ and $X$
therein replaced, respectively,
by $[b,T_{\Omega}]$ and $\HerzSo$.
\end{proposition}

Recall that the following
necessity of commutators on
ball Banach function spaces
was obtained in \cite[Theorem 2.22]{TYYZ}.

\begin{lemma}\label{necebddll1}
Let $X$ be a ball Banach function space
satisfying that the Hardy--Littlewood
maximal operator $\mc$ is bounded
on $X$, and let $b\in L^1_{\mathrm{loc}}(\rrn)$.
Assume $\Omega\in L^{\infty}
(\mathbb{S}^{n-1})$ satisfies that there exists
an open set $\Lambda\subset\mathbb{S}^{n-1}$
such that $\Omega$
never vanishes and never changes sign on
$\Lambda$. If $[b,T_{\Omega}]$ is bounded on
$X$ and satisfies \eqref{necebddle1},
then $b\in\BMO$ and there exists
a positive constant $C$, independent of $b$, such
that
$$
\|b\|_{\BMO}\leq C\left\|\left[
b,T_{\Omega}
\right]\right\|_{X\to X}.
$$
\end{lemma}

We now show Proposition \ref{necebddl}
via using Lemma \ref{necebddll1}.
To this end, we only need to
show that,
under the assumptions
of Proposition \ref{necebddl},
the local generalized Herz space
$\HerzSo$ satisfies all the assumptions
of Lemma \ref{necebddll1}.

\begin{proof}[Proof of Proposition
\ref{necebddl}]
Let all the symbols be as in the
present proposition.
Then, combining the assumptions
of the present proposition,
and Theorem \ref{balll},
we conclude that the
local generalized Herz space
$\HerzSo$ is a BBF space.
On the other hand,
from Corollary \ref{bhlo}, we deduce
that the Hardy--Littlewood maximal operator
$\mc$ is bounded on the
Herz space $\HerzSo$. This implies
that $\HerzSo$ satisfies
all the assumptions of
Lemma \ref{necebddll1} and
hence finishes the proof
of Proposition \ref{necebddl}.
\end{proof}

Via the sufficiency
and the necessity of the boundedness
of commutators on local generalized
Herz spaces established above, we
now prove Theorem \ref{iffbddl}.
We first need the following integral
representation
of $[b,T_\Omega](\1bf_{F})$ for any
bounded measurable set $F$ of $\rrn$,
which might be well known. However,
for the convenience of the reader,
we give its detailed proof here.

\begin{lemma}\label{well}
Let $\Omega$ be a homogeneous Lipschitz
function of degree zero on $\mathbb{S}^{n-1}$
satisfying \eqref{mz},
$T_{\Omega}$ a singular integral
operator with homogeneous kernel $\Omega$,
and $b\in L^1_{{\rm loc}}(\rrn)$. Then, for any
bounded measurable set $F$ and
almost every $x\in\rrn\setminus
\overline{F}$,
$$
\left[b,T_\Omega\right]\left(\1bf_{F}\right)(x)
=\int_{F}\left[b(x)-b(y)\right]\frac{
\Omega(x-y)}{|x-y|^n}\,dy.
$$
\end{lemma}

\begin{proof}
Let $\Omega$, $T_{\Omega}$,
$b$, and $F$ be as in the
present lemma.
Then, applying \eqref{lips},
we conclude that, for any
$\tau\in(0,1)$,
\begin{align*}
\omega_\infty(\tau)\leq
\sup_{\{x,\,y\in\mathbb{S}^{n-1},\,
|x-y|<\tau\}}\left|x-y\right|=\tau.
\end{align*}
Thus,
\begin{equation}\label{welle1}
\int_{0}^{1}\frac{\omega_\infty
(\tau)}{\tau}\,d\tau\leq1.
\end{equation}
In addition, for any $x\in\rrn\setminus\overline{F}$
and $\eps\in(0,d(x,F))$, we have
$[B(x,\eps)]^\complement\cap F=F$, where, for any
$x\in\rrn$ and $F\subset\rrn$,
$$d(x,F):=\inf\left\{|x-y|:\ y\in F\right\}.$$
From this, \eqref{welle1},
and \cite[Corollary 2.1.1]{LDY},
it follows that, for almost every
$x\in\rrn\setminus\overline{F}$,
\begin{align}\label{welle2}
T_\Omega(\1bf_{F})(x)&=\lim\limits_{\eps\to0^+}
\int_{|x-y|\geq\eps}
\frac{\Omega(x-y)}{|x-y|^n}\1bf_{F}(y)\,dy\notag\\
&=\lim\limits_{\eps\to0^+}\int_{[B(x,\eps)]^\complement
\cap F}\frac{\Omega(x-y)}{|x-y|^n}\,dy\notag\\
&=\int_{F}\frac{\Omega(x-y)}{|x-y|^n}\,dy.
\end{align}
Similarly, we find that, for almost every $x\in\rrn
\setminus\overline{F}$,
\begin{align*}
T_\Omega(b\1bf_{F})(x)=\int_{F}b(y)
\frac{\Omega(x-y)}{|x-y|^n}\,dy.
\end{align*}
Combining this and \eqref{welle2}, we conclude that,
for almost every $x\in\rrn\setminus
\overline{F}$,
$$
\left[b,T_\Omega\right]\left(\1bf_{F}\right)(x)
=\int_{F}\left[b(x)-b(y)\right]\frac{
\Omega(x-y)}{|x-y|^n}\,dy,
$$
which completes the proof of Lemma \ref{well}.
\end{proof}

Via both Propositions \ref{suffbddl} and \ref{necebddl},
and Lemma \ref{well}, we next show Theorem \ref{iffbddl}.

\begin{proof}[Proof of Theorem \ref{iffbddl}]
Let all the symbols be as in the present theorem.
Notice that the assumption \eqref{lips} on $\Omega$
implies that, for any given
$r\in(0,\infty]$, $\Omega\in L^{r}(
\mathbb{S}^{n-1})$ and the function $\Omega$ is
continuous on $\mathbb{S}^{n-1}$.
This further implies that
there exists an open set $\Lambda\subset
\mathbb{S}^{n-1}$ such that
$\Omega$ never vanishes and never
changes sign on $\Lambda$.
By this, Lemma \ref{well}, and both
Propositions \ref{suffbddl} and
\ref{necebddl},
we obtain Theorem \ref{iffbddl}.
\end{proof}

Next, we establish
the boundedness characterization
of commutators on global
generalized Herz spaces as follows.

\begin{theorem}\label{iffbddg}
Let $p,\ q\in(1,\infty)$ and
$\omega\in M(\rp)$ satisfy
$$
-\frac{n}{p}<\m0(\omega)
\leq\M0(\omega)<\frac{n}{p'}
$$
and
$$
-\frac{n}{p}<\mi(\omega)
\leq\MI(\omega)<0,
$$
where $\frac{1}{p}+\frac{1}{p'}=1.$
Assume that $\Omega$ is a homogeneous
function of degree zero on $\mathbb{S}^{n-1}$
satisfying both \eqref{lips}
and \eqref{mz}, and $T_{\Omega}$ a
singular integral
operator with homogeneous kernel
$\Omega$. Then, for any
$b\in L^1_{{\rm loc}}(\rrn)$, the commutator
$[b,T_{\Omega}]$ is bounded on
the global generalized Herz space
$\HerzS$ if and only if $b\in\BMO$. Moreover,
there exists a constant $C\in[1,\infty)$,
independent of $b$, such that
$$
C^{-1}\|b\|_{\BMO}\leq\left\|[b,T_{\Omega}]
\right\|_{\HerzS\to\HerzS}
\leq C\|b\|_{\BMO},
$$
where $\|[b,T_{\Omega}]\|_{\HerzS\to\HerzS}$
is as in \eqref{eon} with
$T$ and $X$ therein replaced, respectively,
by $[b,T_{\Omega}]$ and $\HerzS$.
\end{theorem}

First, we consider the sufficiency
of the above theorem.
To this end, we establish the following
boundedness of commutators
on global generalized Herz spaces,
where the assumption about $\Omega$
is weaker than that of Theorem \ref{iffbddg}.

\begin{proposition}\label{suffbddg}
Let $p$, $q$,
and $\omega$ be as in Theorem
\ref{iffbddg}
and $r\in(1,\infty]$.
Assume that $\Omega\in
L^r(\mathbb{S}^{n-1})$ satisfies
both \eqref{degrz} and \eqref{mz},
and $T_{\Omega}$ is a singular integral
operator with homogeneous kernel $\Omega$.
Then there exists a positive constant $C$
such that, for any $b\in\BMO$ and $f\in\HerzS$,
$$
\|[b,T_{\Omega}](f)\|_{\HerzS}\leq
C\|b\|_{\BMO}\|f\|_{\HerzS}.
$$
\end{proposition}

To show Proposition \ref{suffbddg},
recall that Tao et al.\
\cite[Theorem 2.17]{TYYZ} proved
the boundedness of commutators
on the ball Banach function space $X$.
However, due to the deficiency of the associate space
of the global generalized Herz space $\HerzS$, we can
not prove the above boundedness of commutators
on $\HerzS$ by using \cite[Theorem 2.17]{TYYZ} directly.
To overcome this difficulty, via replacing
the assumption about associate spaces
by an assumption about the
extrapolation\index{extrapolation} inequality [see Proposition
\ref{suffbddgl1}(ii) below], we
now establish the following
boundedness of commutators on ball Banach
function spaces, which plays a key role
in the proof of Proposition \ref{suffbddg}.

\begin{proposition}\label{suffbddgl1}
Let $b\in\BMO$, $r\in(1,\infty]$,
$\Omega\in
L^r(\mathbb{S}^{n-1})$ satisfy
both \eqref{degrz} and \eqref{mz},
and $T_{\Omega}$ be a singular integral
operator with homogeneous kernel $\Omega$.
Let $X$ be a ball Banach
function space satisfying that
\begin{enumerate}
  \item[{\rm(i)}] the commutator
  $[b,T_{\Omega}]$ is well defined
  on $X$;
  \item[{\rm(ii)}] there exist
  an $r_0\in[r',\infty)$ and a positive
  constant $\widetilde{C}$ such that,
  for any $(F,\,G)\in\mathcal{F}$ with
  $\|F\|_{X}<\infty$,
  $$
  \|F\|_{X}\leq\widetilde{C}\|G\|_{X},
  $$
  where $\mathcal{F}$ is the same
  as in Theorem \ref{extrag}.
\end{enumerate}
Then there exists a positive constant
$C$, independent of $b$, such that,
for any $f\in X$,
$$
\left\|\left[b,T_{\Omega}
\right](f)\right\|_{X}
\leq C\|b\|_{\BMO}\|f\|_{X}.
$$
\end{proposition}

\begin{proof}
Let all the symbols be as in
the present proposition.
Then, using (i) and repeating
the proof of \cite[Proposition 2.14]{TYYZ}
with $\mathcal{T}$ and \cite[Lemma 2.13]{TYYZ}
therein replaced, respectively, by
$[b,T_{\Omega}]$ and (ii),
we conclude that, for any $f\in X$,
$$
\left\|\left[b,T_{\Omega}(f)\right]\right\|_X
\lesssim\left\|b\right\|_{\BMO}\left\|f\right\|_X.
$$
This finishes the proof of Proposition
\ref{suffbddgl1}.
\end{proof}

We now prove Proposition \ref{suffbddg}.

\begin{proof}[Proof of Proposition
\ref{suffbddg}]
Let all the symbols be as in the present
proposition. Then, combining
the assumptions of both the present
proposition and Theorem \ref{ball},
we find that the global generalized
Herz space $\HerzS$ is
a BBF space.
Thus, to complete the proof of the present
proposition, it suffices for us to show that
$\HerzS$ satisfies both (i) and (ii)
of Proposition \ref{suffbddgl1}.

Indeed, from Definition \ref{gh},
it follows that,
for any $f\in\HerzS$,
\begin{equation}\label{suffbddge1}
\left\|f\right\|_{\HerzSo}
\leq\left\|f\right\|_{\HerzS}<\infty.
\end{equation}
By both this and Proposition \ref{suffbddl},
we conclude that,
for any $f\in\HerzS$,
$[b,T_{\Omega}](f)\in\HerzSo$. Therefore,
$[b,T_{\Omega}]$ is
well defined on $\HerzS$.
This further implies that
$\HerzS$ satisfies Proposition \ref{suffbddgl1}(i).

On the other hand, using
Theorem \ref{extrag},
we find that Proposition \ref{suffbddgl1}(ii)
holds true for $\HerzS$.
Then we conclude that
all the assumptions
of Proposition \ref{suffbddgl1} are satisfied
for the global
generalized Herz space
$\HerzS$, which completes the
proof of Proposition \ref{suffbddg}.
\end{proof}

Next, we focus on the
necessity of Theorem \ref{iffbddg}.
Indeed, via the necessity of the
boundedness of
commutators on ball Banach function
spaces established in
\cite[Theorem 2.17]{TYYZ} (see also
Lemma \ref{necebddll1} above),
we obtain the following conclusion which
is more general than the necessity
of Theorem \ref{iffbddg}.

\begin{proposition}\label{necebddg}
Let $p\in(1,\infty)$,
$q\in[1,\infty)$, and $\omega\in M(\rp)$
be as in Theorem \ref{iffbddg}.
Assume that $b\in L^1_{{\rm loc}}
(\rrn)$ and $\Omega\in
L^\infty(\mathbb{S}^{n-1})$ satisfying
that there exists an open set
$\Lambda\subset\mathbb{S}^{n-1}$
such that $\Omega$
never vanishes and never changes sign
on $\Lambda$.
If the commutator $[b,T_{\Omega}]$
is bounded on the
global generalized Herz space $\HerzS$
and, for any bounded
measurable set $F\subset\rrn$
and almost every
$x\in\rrn\setminus\overline{F}$,
$$
\left[b,T_\Omega\right]
\left(\1bf_{F}\right)(x)
=\int_{F}\left[b(x)-b(y)\right]\frac{
\Omega(x-y)}{|x-y|^n}\,dy,
$$
then $b\in
\BMO$ and there exists
a positive constant $C$,
independent of $b$, such that
$$
\|b\|_{\BMO}\leq C\|[b,T_{\Omega}]
\|_{\HerzS\to\HerzS},
$$
where $\|[b,T_{\Omega}]
\|_{\HerzS\to\HerzS}$ is as in
\eqref{eon} with $T$ and $X$
therein replaced, respectively,
by $[b,T_{\Omega}]$ and $\HerzS$.
\end{proposition}

\begin{proof}
Let all the symbols
be as in the present proposition.
Then, combining both
Theorem \ref{ball} and
Corollary \ref{Th4},
we conclude that
$\HerzS$ is a BBF space and the
Hardy--Littlewood maximal operator
$\mc$ is bounded on $\HerzS$.
Using this and Lemma \ref{necebddll1},
we find that $b\in\BMO$ and
$$
\|b\|_{\BMO}\lesssim
\left\|
\left[b,T_{\Omega}\right]
\right\|_{\HerzS\to\HerzS},
$$
which completes the proof of
Proposition \ref{necebddg}.
\end{proof}

We finally show Theorem \ref{iffbddg}.

\begin{proof}[Proof of
Theorem \ref{iffbddg}]
Let all the symbols be as in
the present theorem.
Then, by \eqref{lips}, we conclude that,
for any given $r\in(0,\infty]$, $\Omega\in L^{r}(
\mathbb{S}^{n-1})$ and the function $\Omega$
is continuous on $\mathbb{S}^{n-1}$.
This further implies that there exists
an open set $\Lambda\subset
\mathbb{S}^{n-1}$ such that
$\Omega$ never vanishes and never changes
sign on $\Lambda$.
Using this, Lemma \ref{well}, and
Propositions \ref{suffbddg} and \ref{necebddg},
we then complete the proof
of Theorem \ref{iffbddg}.
\end{proof}

\section{Compactness Characterizations}

In this section, we establish the
compactness characterization
of commutators on local and global
generalized Herz spaces. We first present the
compactness characterization on
local generalized Herz spaces as follows.
Recall that $\CMO$\index{$\CMO$} is
defined to be the closure of infinitely
differentiable functions
with compact support in $\BMO$.

\begin{theorem}\label{iffcoml}
Let $p,\ q\in(1,\infty)$ and
$\omega\in M(\rp)$ satisfy
$$
-\frac{n}{p}<\m0(\omega)
\leq\M0(\omega)<\frac{n}{p'}
$$
and
$$
-\frac{n}{p}<\mi(\omega)\leq\MI(\omega)<\frac{n}{p'},
$$
where $\frac{1}{p}+\frac{1}{p'}=1.$
Assume that $\Omega$ is a homogeneous
function satisfying \eqref{lips}, \eqref{degrz},
and \eqref{mz}, and $T_{\Omega}$ a singular integral
operator with homogeneous kernel $\Omega$. Then,
for any $b\in L^1_{{\rm loc}}(\rrn)$,
the commutator $[b,T_{\Omega}]$
is compact on the local generalized Herz space
$\HerzSo$ if and only if $b\in\CMO$.
\end{theorem}

We prove this theorem by considering, respectively,
the sufficiency and the necessity of the
compactness of commutators on
$\HerzSo$. First, for
the sufficiency, we have the following
more general conclusion.

\begin{proposition}\label{compsul}
Let $p$, $q$, and $\omega$ be as in Theorem \ref{iffcoml}.
Assume that $b\in L^1_{\mathrm{loc}}(\rrn)$,
$\Omega\in L^\infty(\mathbb{S}^{n-1})$ satisfies
\eqref{degrz}, \eqref{mz},
and the $L^{\infty}$-Dini condition,
and $T_{\Omega}$ is a singular
integral operator with homogeneous kernel
$\Omega$. If $b\in\CMO$,
then the commutator $[b,T_{\Omega}]$ is compact on the
local generalized Herz space $\HerzSo$.
\end{proposition}

Indeed, in \cite[Theorem 3.1]{TYYZ},
Tao et al.\ established the following
compactness of commutators
on ball Banach function spaces.

\begin{lemma}\label{compsull1}
Let $X$ be a ball Banach function
space satisfying
that the Hardy--Littlewood
maximal operator $\mc$ is bounded
on $X$ and its associate space
$X'$.
Assume that $b\in L^1_{\mathrm{loc}}(\rrn)$,
$\Omega\in L^\infty(\mathbb{S}^{n-1})$ satisfies
\eqref{degrz}, \eqref{mz},
and the $L^{\infty}$-Dini condition,
and $T_{\Omega}$ is a singular
integral operator with homogeneous kernel
$\Omega$. If $b\in\CMO$,
then the commutator $[b,T_{\Omega}]$ is
compact on $X$.
\end{lemma}

We now show Proposition \ref{compsul}
by using Lemma \ref{compsull1}. To this end,
we only need to check that all
the assumptions of Lemma \ref{compsull1}
are satisfied for $\HerzSo$ under
the assumptions of Proposition \ref{compsul}.

\begin{proof}[Proof of Proposition \ref{compsul}]
Let all the symbols be as in the present proposition.
Then, from the assumptions $p,\ q\in(1,\infty)$ and
$-\frac{n}{p}<\m0(\omega)\leq\M0(\omega)<\frac{n}{p'}$,
and Theorem \ref{balll}, we deduce
that the local generalized Herz
space $\HerzSo$ is a BBF space.
In addition, by the assumptions
of the present proposition, \eqref{th175e1},
and \eqref{th175e2},
we find that $\mc$ is bounded
on both $\HerzSo$ and its associate space
$(\HerzSo)'$.
This further implies that
$\HerzSo$ satisfies all the assumptions
of Lemma \ref{compsull1} and hence
finishes the proof of
Proposition \ref{compsul}.
\end{proof}

Next, we consider the
necessity of Theorem \ref{iffcoml}.
Indeed, we have the following
more general conclusion
for the necessity of compactness
of commutators on local
generalized Herz spaces.

\begin{proposition}\label{compnecel}
Let $p$, $q$, and $\omega$ be as in Theorem \ref{iffcoml}.
Assume that $b\in L^1_{{\rm loc}}(\rrn)$ and $\Omega\in
L^\infty(\mathbb{S}^{n-1})$ satisfies that
there exists an open set $\Lambda\subset
\mathbb{S}^{n-1}$ such that
$\Omega$ never vanishes and never
changes sign on $\Lambda$. If
the commutator $[b,T_{\Omega}]$ is compact
on the local
generalized Herz space $\HerzSo$ and
\eqref{necebddle1} holds true,
then $b\in\CMO$.
\end{proposition}

To prove this proposition,
we require the following
necessity of the compactness of
commutators on ball Banach function
spaces,
which is just \cite[Theorem 3.2]{TYYZ}.

\begin{lemma}\label{compnecell1}
Let $X$ be a ball Banach function space
satisfying that
the Hardy--Littlewood $\mc$ is
bounded on both $X$ and its associate
space $X'$.
Assume that
$b\in L^1_{{\rm loc}}(\rrn)$
and $\Omega\in
L^\infty(\mathbb{S}^{n-1})$
satisfies that
there exists an open
set $\Lambda\subset
\mathbb{S}^{n-1}$ such that
$\Omega$ never vanishes and never
changes sign on $\Lambda$. If
the commutator $[b,T_{\Omega}]$
is compact on $X$ and
\eqref{necebddle1} holds true,
then $b\in\CMO$.
\end{lemma}

Via Lemma \ref{compnecell1},
we now show Proposition \ref{compnecel}.

\begin{proof}[Proof of
Proposition \ref{compnecel}]
Let all the symbols be as in the present proposition.
Then, from the assumptions
of the present proposition,
and Theorem \ref{balll}, we deduce that
the local generalized
Herz space $\HerzSo$
is a BBF space.

Therefore, to complete the
proof of the present proposition,
it suffices to show that $\HerzSo$ satisfies all
the assumptions of Lemma \ref{compnecell1}.
To be precise, we only need
to prove that the Hardy--Littlewood
maximal operator $\mc$ is bounded
on both $\HerzSo$ and $(\HerzSo)'$.
Indeed, applying
\eqref{th175e1} and \eqref{th175e2},
we find that, for any
$f\in L^1_{\mathrm{loc}}(\rrn)$,
$$
\left\|\mc(f)\right\|_{\HerzSo}
\lesssim\|f\|_{\HerzSo}
$$
and
$$
\left\|\mc(f)\right\|_{(\HerzSo)'}
\lesssim\|f\|_{(\HerzSo)'},
$$
which imply that
$\mc$ is bounded on both
$\HerzSo$ and $(\HerzSo)'$, and
hence complete the proof of
Proposition \ref{compnecel}.
\end{proof}

We now give the proof of Theorem \ref{iffcoml}.

\begin{proof}[Proof of Theorem \ref{iffcoml}]
Let all the symbols be as in the present theorem.
Then, from \eqref{lips}, it
follows that $\Omega\in L^\infty(\mathbb{S}^{n-1})$
and $\Omega$ is continuous on $\mathbb{S}^{n-1}$.
This further implies that
there exists an open set $\Lambda\subset
\mathbb{S}^{n-1}$ such that $\Omega$
never vanishes and
never changes sign on $\Lambda$.
In addition, by \eqref{lips} again,
we conclude that, for any $\tau\in(0,1)$,
$\omega_{\infty}(\tau)\leq\tau$, where
$\omega_{\infty}(\tau)$ is as in
Definition \ref{dini}. This further implies that
$$\int_{0}^{1}\frac{\omega_{\infty}
(\tau)}{\tau}\,d\tau\in[0,1].$$
Thus, the function $\Omega$ satisfies
the $L^\infty$-Dini
condition. Then, applying Lemma
\ref{well} and both Propositions
\ref{compsul} and \ref{compnecel},
we then complete
the proof of Theorem \ref{iffcoml}.
\end{proof}

Next, we establish
the compactness characterization of
commutators on global
generalized Herz spaces.
Namely, the following conclusion
holds true.

\begin{theorem}\label{iffcomg}
Let $p,\ q\in(1,\infty)$ and
$\omega\in M(\rp)$ satisfy
$$
-\frac{n}{p}<\m0(\omega)\leq\M0(\omega)<\frac{n}{p'},
$$
and
$$
-\frac{n}{p}<\mi(\omega)\leq\MI(\omega)<0,
$$
where $\frac{1}{p}+\frac{1}{p'}=1.$
Assume that $\Omega$ is a homogeneous
function of degree zero on $\mathbb{S}^{n-1}$
satisfying both \eqref{lips} and \eqref{mz},
and $T_{\Omega}$ a singular integral
operator with homogeneous kernel $\Omega$. Then, for
any $b\in L^1_{{\rm loc}}(\rrn)$,
the commutator $[b,T_{\Omega}]$
is compact on the global generalized
Herz space $\HerzS$ if and only if $b\in\CMO$.
\end{theorem}

To obtain this theorem,
we first show the following
compactness of commutators
on global generalized Herz spaces, which
is more general than the sufficiency
of Theorem \ref{iffcomg}.

\begin{proposition}\label{compsug}
Let $p$, $q$, and $\omega$ be
as in Theorem \ref{iffcomg}.
Assume that
$b\in L^1_{\mathrm{loc}}(\rrn)$,
$\Omega\in L^\infty(\mathbb{S}^{n-1})$
satisfies \eqref{degrz}, \eqref{mz},
and the $L^{\infty}$-Dini condition,
and $T_{\Omega}$ is a singular
integral operator with homogeneous
kernel $\Omega$. If $b\in\CMO$,
then the commutator $[b,T_{\Omega}]$
is compact on $\HerzS$.
\end{proposition}

To show Proposition \ref{compsug},
let $\Omega\in L^{\infty}(\mathbb{S}^
{n-1})$ satisfy
\eqref{degrz}, \eqref{mz},
and the $L^{\infty}$-Dini condition, and
$T_{\Omega}$ be a singular
integral operator with homogeneous
kernel $\Omega$. Recall that
the \emph{maximal operator}
$T_{\Omega}^*$\index{maximal
operator}\index{$T_{\Omega}^*$}
is defined by setting,
for any $f\in\Msc(\rrn)$ and
$x\in\rrn$,
\begin{equation}\label{maxs}
T_{\Omega}^{*}(f)(x):=
\sup_{\eps\in(0,\infty)}
\left|\int_{\{y\in\rrn:\ |x-y|>\eps\}}
\frac{\Omega(x-y)}{|x-y|^n}f(y)
\,dy\right|.
\end{equation}
This operator plays a vital role
in the proof of Proposition \ref{compsug}.
Indeed, we have the following
boundedness of $T_{\Omega}^*$ on $X$,
which is just \cite[Proposition 3.10]{TYYZ}.

\begin{lemma}\label{compsugl3}
Let $X$ be a ball Banach function space
satisfying that the Hardy--Littlewood
maximal operator $\mc$ is bounded
on both $X$ and $X'$. Assume
$\Omega\in L^{\infty}
(\mathbb{S}^{n-1})$ satisfy
\eqref{degrz}, \eqref{mz},
and the $L^{\infty}$-Dini condition
and $T_{\Omega}^{*}$ is as in \eqref{maxs}.
Then there exists a positive constant
$C$ such that, for any
$f\in X$,
$$
\left\|T_{\Omega}^*(f)\right\|_{X}
\leq C\|f\|_{X}.
$$
\end{lemma}

In addition, to show Proposition \ref{compsug},
due to the deficiency of the associate space
of the global generalized Herz space $\HerzS$, we can
not obtain the above compactness of commutators
on $\HerzS$ by using \cite[Theorem 3.1]{TYYZ} directly.
To overcome this obstacle, via replacing
the assumption about associate spaces
by an assumption about the
extrapolation inequality [see Proposition
\ref{compsugl1}(ii) below], we
now establish the following
boundedness of commutators on ball Banach
function spaces, which is important in the proof
of the compactness characterization
of commutators on global generalized Herz spaces.

\begin{proposition}\label{compsugl1}
Let $b\in\CMO$, $\Omega\in L^{\infty}
(\mathbb{S}^{n-1})$ satisfy
\eqref{degrz}, \eqref{mz},
and the $L^{\infty}$-Dini condition,
$T_{\Omega}$ be a singular
integral operator with homogeneous
kernel $\Omega$, and $T_{\Omega}^{*}$
as in \eqref{maxs}. Assume that $X$
is a ball Banach function space satisfying
the following three conditions:
\begin{enumerate}
  \item[{\rm(i)}] $[b,T_{\Omega}]$ and
  $T_{\Omega}^*$ are well defined on $X$;
  \item[{\rm(ii)}] the Hardy--Littlewood
  maximal operator $\mc$ is bounded on $X$;
  \item[{\rm(iii)}] there exist
  an $r_0\in(1,\infty)$ and a positive
  constant $C$ such that,
  for any $(F,\,G)\in\mathcal{F}$ with
  $\|F\|_{X}<\infty$,
  $$
  \|F\|_{X}\leq C\|G\|_{X},
  $$
  where $\mathcal{F}$ is defined
  as in Theorem \ref{extrag}.
\end{enumerate}
Then the commutator $[b,T_{\Omega}]$ is
compact on $X$.
\end{proposition}

\begin{proof}
Let all the symbols be as in
the present proposition. Then,
using (i) and (ii), and repeating
an argument similar to that
used in the proof of \cite[Proposition
3.14]{TYYZ} with $\mathcal{T}$
and \cite[Lemma 2.13]{TYYZ} therein
replaced, respectively,
by $T_{\Omega}^*$ and (iii) here,
we conclude that,
for any $f\in X$,
$$
\left\|T_{\Omega}^*(f)\right\|_{X}
\lesssim\left\|f\right\|_X.
$$
From this and the proof of
\cite[Theorem 3.1]{TYYZ} via replacing
\cite[Theorem 2.17]{TYYZ} therein
by Lemma \ref{suffbddgl1} here, we infer
that the commutator $[b,T_{\Omega}]$
is compact on $X$.
This then finishes the proof of
Proposition \ref{compsugl1}.
\end{proof}

Now, we can show
Proposition \ref{compsug} via checking
that the global generalized
Herz space $\HerzS$ satisfies all the
assumptions of Proposition \ref{compsugl1}.
For this purpose, we first establish
an auxiliary lemma about
the maximal operator $T_{\Omega}^*$
as follows.

\begin{lemma}\label{compsugl2}
Let $p$, $q$, and $\omega$ be
as in Theorem \ref{iffcomg},
$\Omega\in L^{\infty}
(\mathbb{S}^{n-1})$ satisfy
\eqref{degrz}, \eqref{mz},
and the $L^{\infty}$-Dini condition,
and $T_{\Omega}^{*}$ be
as in \eqref{maxs}.
Then $T_{\Omega}^{*}$ is well defined
on $\HerzS$.
\end{lemma}

\begin{proof}
Let all the symbols be as in the
present lemma. Then, using
the assumptions of the present lemma
and Theorem \ref{balll}, we conclude that
the local generalized Herz space
$\HerzSo$ is a BBF space. In addition,
from \eqref{th175e1} and \eqref{th175e2},
we deduce that $\mc$ is bounded on both
$\HerzSo$ and $(\HerzSo)'$. Therefore,
by Lemma \ref{compsugl3}, we find
that, for any $f\in\HerzSo$,
\begin{equation}\label{compsugl2e1}
\left\|T_{\Omega}^*(f)\right\|_{\HerzSo}
\lesssim\left\|f\right\|_{\HerzSo}.
\end{equation}

On the other hand, applying
\eqref{suffbddge1}, we conclude that
$\HerzS\subset\HerzSo$. This, combined
with \eqref{compsugl2e1},
further implies that, for any
$f\in\HerzS$, $T_{\Omega}^*(f)\in
\HerzSo$, which shows that
$T_{\Omega}^*$ is well defined
on $\HerzS$, and hence completes the
proof of Lemma \ref{compsugl2}.
\end{proof}

Via the above two conclusions, we
now show Proposition \ref{compsug}.

\begin{proof}[Proof of Proposition
\ref{compsug}]
Let all the symbols be as in the
present proposition. Then, by
the assumptions of the present
propositions, and Theorem \ref{ball},
we find that the global
generalized Herz space $\HerzS$ is
a BBF space. Thus, to finish
the proof of the present proposition,
we only need to show that $\HerzS$ satisfies
(i), (ii), and (iii) of Proposition \ref{compsugl1}.

We first prove that Proposition \ref{compsugl1}(i)
holds true for $\HerzS$. Indeed,
from \eqref{suffbddge1},
it follows that $\HerzS\subset\HerzSo$.
This, together with Proposition \ref{suffbddl},
implies that, for any $f\in\HerzS$,
$[b,T_{\Omega}](f)$ is well defined.
On the other hand, applying
Lemma \ref{compsugl2}, we conclude that
$T_{\Omega}^*$ is also well defined
on $\HerzS$. Therefore, $\HerzS$
satisfies Proposition \ref{compsugl1}(i)
under the assumptions of the present proposition.

Next, we show that
$\HerzS$ satisfies Proposition \ref{compsugl1}(ii),
namely, $\mc$ is bounded on $\HerzS$. Indeed,
from the assumptions of the
present proposition, and Corollary \ref{Th4},
we deduce that, for any
$f\in L^1_{\mathrm{loc}}(\rrn)$,
$$
\left\|\mc(f)\right\|_{\HerzS}
\lesssim\|f\|_{\HerzS}.
$$
This implies that
$\mc$ is bounded
on $\HerzS$, and hence
Proposition \ref{compsugl1}(ii) holds true for
$\HerzS$.

Finally, by Theorem \ref{extrag},
we find that $\HerzS$ satisfies
Proposition \ref{compsugl1}(iii). Thus,
under the assumptions
of the present proposition,
the global generalized Herz space
$\HerzS$ satisfies all the assumptions
of Proposition \ref{compsugl1}.
This implies that
the commutator $[b,T_{\Omega}]$ is
compact on $\HerzS$, and hence finishes
the proof of Proposition \ref{compsug}.
\end{proof}

Then we focus on the
necessity of Theorem \ref{iffcomg}.
Indeed, we have the following
more general conclusion.

\begin{proposition}\label{compneceg}
Let $p$, $q$, and $\omega$ be as
in Theorem \ref{iffcomg}.
Assume that $b\in L^1_{{\rm loc}}(\rrn)$
and $\Omega\in
L^\infty(\mathbb{S}^{n-1})$ satisfies that
there exists an open set $\Lambda\subset
\mathbb{S}^{n-1}$ such that
$\Omega$ never vanishes and never
changes sign on $\Lambda$. If
the commutator $[b,T_{\Omega}]$ is compact
on the global
generalized Herz space $\HerzS$, and
satisfies \eqref{necebddle1},
then $b\in\CMO$.
\end{proposition}

To prove this proposition,
Lemma \ref{compnecell1}
is unavailable
due to the deficiency
of associate spaces of
global generalized
Herz spaces. To overcome this
difficulty, via replacing
the assumption about associate spaces
by an assumption about the
extrapolation inequality [see Proposition
\ref{compnecegl1}(ii) below],
we now show the
necessity of the compactness
of commutators on ball Banach function
spaces as follows.

\begin{proposition}\label{compnecegl1}
Let $b\in L^1_{{\rm loc}}(\rrn)$,
$\Omega\in
L^\infty(\mathbb{S}^{n-1})$ satisfy that
there exists an open set $\Lambda\subset
\mathbb{S}^{n-1}$ such that
$\Omega$ never vanishes and never
changes sign on $\Lambda$, and
$X$ be a ball
Banach function space
satisfy the following
two conditions:
\begin{enumerate}
  \item[{\rm(i)}] the Hardy--Littlewood
  maximal operator $\mc$ is bounded
  on $X$;
  \item[{\rm(ii)}] there exists
  a positive constant $C$ such that,
  for any $(F,\,G)\in\mathcal{F}$ with
  $\|F\|_{X}<\infty$,
  $$
  \|F\|_{X}\leq C\|G\|_{X},
  $$
  where $\mathcal{F}$ is defined
  as in Theorem \ref{extrag} with
  $r_0:=1$.
\end{enumerate}
If the commutator $[b,T_{\Omega}]$
is compact on $X$, and the assumption
\eqref{necebddle1} holds true,
then $b\in\CMO$.
\end{proposition}

To prove this proposition, we need
the following technique lemma
which is a part of
\cite[Lemma 2.2 and Remark 2.3]{is17}.

\begin{lemma}\label{compnecegl3}
Let $X$ be a ball Banach function
space such that $\mc$ is bounded on $X$.
Then there exists a positive
constant $C$ such that, for any ball $B\in
\mathbb{B}$,
$$
\frac{1}{|B|}\left\|\1bf_B\right\|_X\left\|
\1bf_B\right\|_{X'}\leq C.
$$
\end{lemma}

Via this technique lemma, we now show
the following equivalent characterization
of $\BMO$ via ball Banach function
spaces, which is an essential tool
in the proof of Proposition \ref{compnecegl1} above.

\begin{lemma}\label{compnecegl2}
Let $X$ be a ball Banach function space
satisfying the following
two statements:
\begin{enumerate}
  \item[{\rm(i)}] the Hardy--Littlewood
  maximal operator $\mc$ is bounded
  on $X$;
  \item[{\rm(ii)}] there exists
  a positive constant $C$ such that,
  for any $(F,\,G)\in\mathcal{F}$ with
  $\|F\|_{X}<\infty$,
  $$
  \|F\|_{X}\leq C\|G\|_{X},
  $$
  where $\mathcal{F}$ is defined
  as in Theorem \ref{extrag} with
  $r_0:=1$.
\end{enumerate}
Let $f\in L^1_{\mathrm{loc}}(\rrn)$.
Then $f\in\BMO$ if and only if
$$\index{$\mathrm{BMO}_X(\rrn)$}
\left\|f\right\|_{\mathrm{BMO}_X(\rrn)}
:=\sup_{B\in\mathbb{B}}\left\{\frac{1}{
\|\1bf_B\|_X}\left\|\left|f-f_B
\right|\1bf_B\right\|_X\right\}<\infty,
$$
where $f_B$ is defined as in \eqref{fb}.
Moreover, there exists
a positive constant $C$ such that, for any
$f\in\BMO$,
$$
C^{-1}\|f\|_{\BMO}\leq
\|f\|_{\mathrm{BMO}_X(\rrn)}\leq C\|f\|_{\BMO}.
$$
\end{lemma}

\begin{proof}
Let all the symbols be as in
the present lemma and
$f\in L^1_{\mathrm{loc}}(\rrn)$.
We first prove the sufficiency.
To this end, assume that
$$\left\|f\right\|_{
\mathrm{BMO}_X(\rrn)}<\infty.$$
By this, Definition \ref{assd}, the assumption
(i) of the present lemma,
and Lemma \ref{compnecegl3},
we conclude that, for any ball $B\in\mathbb{B}$,
\begin{align*}
&\frac{1}{|B|}\int_{B}\left|f(x)
-f_B\right|\,dx\\
&\quad\leq\frac{1}{|B|}\left\|
\left|f-f_B\right|\1bf_B\right\|_X
\left\|\1bf_B\right\|_{X'}
\lesssim\frac{1}{\|\1bf_B\|_{X}}
\left\|\left|f-f_B\right|\1bf_B\right\|_X\\
&\quad\lesssim
\left\|f\right\|_{\mathrm{BMO}_X(\rrn)}<\infty,
\end{align*}
which, combined with \eqref{fbmo},
further implies that
\begin{equation}\label{compnecegl2e1}
\left\|f\right\|_{\BMO}\lesssim
\left\|f\right\|_{\mathrm{BMO}_X(\rrn)}<\infty.
\end{equation}
This finishes the proof of the sufficiency.

Conversely, we next show the necessity.
To achieve this, assume $f\in\BMO$.
Then, for any given $\upsilon\in A_1(\rrn)$ and for any
ball $B\in\mathbb{B}$,
repeating the proof of \cite[Theorem
1.2]{ins19} with $Q$ therein
replaced by $B$, we find that
$$
\int_{B}\left|f(x)-f_B\right|\upsilon(x)\,dx
\lesssim\left[\upsilon\right]_{A_1(\rrn)}
\|f\|_{\BMO}\int_{B}\upsilon(x)\,dx,
$$
where the implicit positive constant
depends only on $n$. From this and
(ii), we infer that,
for any $B\in\mathbb{B}$,
$$
\left\|\left|f-f_B\right|\1bf_B\right\|_X
\lesssim\left\|\1bf_B\right\|_X\|f\|_{\BMO}.
$$
This further implies that
$$
\left\|f\right\|_{\mathrm{BMO}_X(\rrn)}
\lesssim\left\|f\right\|_{\BMO}<\infty,
$$
which completes the proof of the necessity.
Combining this and \eqref{compnecegl2e1},
we find that, for any $f\in\BMO$,
$$
\left\|f\right\|_{\BMO}
\sim\left\|f\right\|_{\mathrm{BMO}_X(\rrn)}
$$
with the positive equivalence constants independent
of $f$. Thus, the proof of Lemma \ref{compnecegl2}
is then completed.
\end{proof}

Applying Lemma \ref{compnecegl2}, we
now show Proposition \ref{compnecegl1}.

\begin{proof}[Proof of Proposition \ref{compnecegl1}]
Let all the symbols be as in the
present proposition and the commutator
$[b,T_{\Omega}]$ is compact on the
ball Banach function space $X$ under consideration.
Then, repeating the proof of \cite[Theorem 3.2]{TYYZ}
via replacing \cite[Lemma 3.15]{TYYZ}
therein by Lemma \ref{compnecegl2} here,
we conclude that
$b\in\CMO$. This then finishes the proof of
Proposition \ref{compnecegl1}.
\end{proof}

Via Proposition \ref{compnecegl1}, we next
prove Proposition \ref{compneceg}.

\begin{proof}[Proof of Proposition
\ref{compneceg}]
Let all the symbols
be as in the present proposition.
Then, from Theorem \ref{ball},
we deduce that the global
generalized Herz space $\HerzS$
is a BBF space under the assumptions
of the present proposition.
This implies that,
to finish the proof of the present
proposition,
we only need to show
that $\HerzS$ satisfies all the
assumptions of Proposition \ref{compnecegl1}.

Indeed, applying
Corollary \ref{Th4}, we find that
$\mc$ is bounded on $\HerzS$,
which implies that
Proposition \ref{compnecegl1}(i) holds
true for $\HerzS$. On the other hand,
using Theorem \ref{extrag} with $r_0:=1$,
we conclude that $\HerzS$ satisfies
Proposition \ref{compnecegl1}(ii).
This further implies that
all the assumptions of Proposition \ref{compnecegl1}
are satisfied for $\HerzS$.
Thus, $b\in\CMO$ and
we then complete the proof
of Proposition \ref{compneceg}.
\end{proof}

We then show Theorem \ref{iffcomg}
via both Propositions \ref{compsug}
and \ref{compneceg}.

\begin{proof}[Proof of Theorem \ref{iffcomg}]
Let all the symbols be as in the present theorem.
Then, by \eqref{lips}, we find
that $\Omega$ is continuous
and hence bounded on
$\mathbb{S}^{n-1}$.
This further implies that $\Omega
\in L^{\infty}(\mathbb{S}^{n-1})$ and
there exists an open set $\Lambda\subset
\mathbb{S}^{n-1}$ such that $\Omega$
never vanishes and
never changes sign on $\Lambda$.
On the other hand, combining
\eqref{lips} and Definition \ref{dini},
we conclude that, for any $\tau\in(0,1)$,
$\omega_{\infty}(\tau)\leq\tau$.
This further implies that
$$\int_{0}^{1}\frac{\omega_{\infty}
(\tau)}{\tau}\,d\tau\in[0,1].$$
Then the function $\Omega$ satisfies
the $L^\infty$-Dini
condition. Thus, applying Lemma
\ref{well} and Propositions
\ref{compsug} and \ref{compneceg},
we then complete
the proof of Theorem \ref{iffcomg}.
\end{proof}

\chapter{Generalized Herz--Hardy Spaces\label{sec6}}
\markboth{\scriptsize\rm\sc Generalized Herz--Hardy Spaces}
{\scriptsize\rm\sc Generalized Herz--Hardy Spaces}

Let $p$, $q\in(0,\infty)$ and $\omega\in M(\rp)$.
In this chapter, we first
introduce the generalized Herz--Hardy
spaces, both $\HaSaHo$ and $\HaSaH$, by
the grand maximal function. Then,
via a totally fresh perspective,
we establish their complete real-variable theory.
To be precise, via viewing the local generalized
Herz space $\HerzSo$ as a special
case of ball quasi-Banach function spaces
and applying the known
real-variable characterizations of
Hardy spaces $H_X(\rrn)$
associated with ball quasi-Banach
function spaces $X$, we establish
various maximal function,
(finite) atomic, molecular, and
Littlewood--Paley function characterizations
of the Herz--Hardy space $\HaSaHo$.
Moreover, the duality and the Fourier
transform properties of $\HaSaHo$ are
also obtained based on the corresponding
results about $H_X(\rrn)$.
However, the study of $\HaSaH$ is more
difficult than that of $\HaSaHo$ due to
the deficiency of associate spaces
of global generalized Herz spaces.
To overcome this obstacle,
via replacing the assumptions
of the boundedness
of powered Hardy--Littlewood maximal operators
on associate spaces
(see Assumption \ref{assas} above)
by some weaker assumptions about
the integral representations of quasi-norms
of ball quasi-Banach function spaces
as well as some boundedness of
powered Hardy--Littlewood maximal operators
[see both (ii) and (iii)
of Theorem \ref{Atogx} below],
we establish some improved (finite) atomic and
molecular characterizations
of the Hardy space $H_X(\rrn)$,
which include the corresponding
results obtained by Sawano et al.\
\cite{SHYY}.
Using these improved characterizations
and making full use of the obtained duality between
block spaces and global generalized
Herz spaces in Chapter \ref{sec4}
as well as the construction of
the quasi-norm $\|\cdot\|_{\HerzS}$,
we then obtain the maximal function,
(finite) atomic,
molecular, the various Littlewood--Paley
function characterizations of $\HaSaH$
and also give some properties about Fourier transforms
of distributions in $\HaSaH$.
As applications,
via first showing two boundedness criteria
of Calder\'{o}n--Zygmund operators
on $H_X(\rrn)$,
we establish the boundedness
of Calder\'{o}n--Zygmund
operators on generalized Herz--Hardy spaces.
In addition, we also introduce
the generalized Morrey--Hardy spaces.
Using the fact that,
under some reasonable assumptions on the exponents,
the generalized Morrey spaces
coincide with the generalized Herz spaces
in the sense of equivalent norms,
we obtain the complete real-variable
characterizations of generalized
Morrey--Hardy spaces immediately.

To begin with, we introduce the
generalized Herz--Hardy spaces.
For this purpose, we recall some basic concepts.
Throughout this book, we use $\mathcal{S}(\rrn)$
\index{$\mathcal{S}(\rrn)$}to
denote the \emph{space of all
Schwartz functions}\index{Schwartz function}
on $\rrn$ equipped with
the well-known topology determined by a
countable family of norms, and $\mathcal{S}'(\rrn)$
\index{$\mathcal{S}'(\rrn)$}its \emph{topological
dual space}\index{topological dual space} equipped with
the weak-$\ast$ topology\index{weak-$\ast$ topology}.
For any $N\in\mathbb{N}$ and
$\phi\in\mathcal{S}(\rrn)$, let
\begin{equation}\label{defp}\index{$p_N$}
p_N(\phi):=\sum_{\alpha\in\mathbb{Z}
^{n}_{+},|\alpha|\leq N}
\sup_{x\in\rrn}\left\{(1+|x|)^{N}|
\partial^{\alpha}\phi(x)|\right\}
\end{equation}
and
\begin{equation}\label{defn}
\mathcal{F}_N(\rrn)\index{$\mathcal{
F}_N(\rrn)$}:=
\left\{\phi\in\mathcal{S}
(\rrn)\index{$\mathcal{S}(\rrn)$}:\
p_N(\phi)
\in[0,1]\right\}.
\end{equation}
Moreover, the \emph{non-tangential
grand maximal function}\index{non-tangential
grand maximal function}
$\mathcal{M}_{N}(f)$\index{$\mathcal{M}_{N}$}
of $f\in\mathcal{S}'(\rrn)$
is defined by setting, for any $x\in\rrn$,
\begin{equation}\label{sec6e1}
\mathcal{M}_{N}(f)(x):=
\sup\left\{|f*\phi_{t}(y)|:\
\phi\in\mathcal{F}_{N}(\rrn),\
t\in(0,\infty),\ |x-y|<t\right\}.
\end{equation}
We now introduce the generalized
Herz--Hardy spaces\index{generalized
Herz--Hardy space}
via the grand maximal function as follows.

\begin{definition}\label{herzh}
Let $p$, $q\in(0,\infty)$, $\omega\in M(\rp)$, and
$N\in\mathbb{N}$.
Then
\begin{enumerate}
  \item[{\rm(i)}] the \emph{generalized Herz--Hardy space}
  $\HaSaHo$\index{$\HaSaHo$}, associated
  with the local generalized Herz space $\HerzSo$,
  is defined to be the set of all the $f\in
  \mathcal{S}'(\rrn)$ such that
  $$\left\|f\right\|_{\HaSaHo}:=
  \left\|\mathcal{M}_{N}(f)\right\|_{\HerzSo}<\infty;$$
  \item[{\rm(ii)}] the \emph{generalized Herz--Hardy space}
  $\HaSaH$\index{$\HaSaH$}, associated
  with the global generalized Herz space $\HerzS$,
  is defined to be the set of all the
  $f\in\mathcal{S}'(\rrn)$ such that $$
\|f\|_{\HaSaH}:=\left\|\mathcal{M}_{N}(f)
\right\|_{\HerzS}<\infty.$$
\end{enumerate}
\end{definition}

\begin{remark}\label{remark4.0.8}
If, in Definition \ref{herzh},
let $\omega(t):=t^\alpha$
for any $t\in(0,\infty)$ and for
any given $\alpha\in\rr$, then, in this
case, the generalized Herz--Hardy space
$\HaSaHo$ goes back
to the classical \emph{homogeneous Herz-type
Hardy space\index{homogeneous
Herz-type Hardy \\space}
$H\dot{K}^{\alpha,q}_{p}
(\rrn)$}\index{$H\dot{K}^{\alpha,q}_{p}(\rrn)$}
which was originally introduced in
\cite[Definition 2.1]{LuY95}
(see also \cite[Chapter 2]{LYH}).
However, we point out that, to the
best of our knowledge, even in this case,
the generalized
Herz--Hardy space $\HaSaH$ is also new.
\end{remark}

Recall that the local and the
global generalized Morrey spaces
are defined as in Remark \ref{remhs}(iv). Then we
introduce the generalized
Morrey--Hardy spaces\index{generalized
Morrey--Hardy space} as follows.

\begin{definition}\label{morrhz}
Let $p$, $q\in(0,\infty)$, $\omega\in M(\rp)$, and
$N\in\mathbb{N}$.
Then
\begin{enumerate}
  \item[{\rm(i)}] the \emph{generalized
  Morrey--Hardy space}
  $H\MorrSo$\index{$H\MorrSo$},
  associated with the local generalized
  Morrey space $\MorrSo$,
  is defined to be the set of all the
  $f\in\mathcal{S}'(\rrn)$ such that
  $$\left\|f\right\|_{H\MorrSo}:=
  \left\|\mathcal{M}_{N}(f)\right\|_
  {\MorrSo}<\infty;$$
  \item[{\rm(ii)}] the \emph{generalized
  Morrey--Hardy space}
  $H\MorrS$\index{$H\MorrS$}, associated
  with the global generalized Morrey space
  $\MorrS$,
  is defined to be the set of all the
  $f\in\mathcal{S}'(\rrn)$ such that $$
\|f\|_{H\MorrS}:=\left\|\mathcal{M}_{N}
(f)\right\|_{\MorrS}<\infty.$$
\end{enumerate}
\end{definition}

\begin{remark}\label{remark4.10}
\begin{enumerate}
  \item[{\rm(i)}] To the best of our
  knowledge, in Definition \ref{morrhz},
  even when $\omega(t):=t^\alpha$
  for any $t\in(0,\infty)$ and
  for any given $\alpha\in\rr$,
  the generalized Morrey--Hardy
  spaces, both $H\MorrSo$ and $H\MorrS$,
  are also new.
  \item[{\rm(ii)}] In Definition \ref{morrhz},
  let $p$, $q\in[1,\infty)$
  and $\omega\in M(\rp)$ satisfy
  $$\Mw\in(-\infty,0).$$
  Then, by Remark \ref{remhs}(iv),
  we conclude that,
  in this case, the generalized
  Morrey--Hardy spaces $H\MorrSo$
  and $H\MorrS$ coincide, respectively,
  with the generalized Herz--Hardy spaces
  $\HaSaHo$ and $\HaSaH$ with equivalent norms.
\end{enumerate}
\end{remark}

\section{Maximal Function Characterizations}

In this section, we establish various
maximal function characterizations
of generalized Herz--Hardy spaces.
We first present the concepts of various
maximal functions as follows (see
\cite[Definition 2.12]{SHYY}).

\begin{definition}\label{smax}
Let $N\in\mathbb{N},\ a,\ b\in(0,\infty)$,
$\phi\in\mathcal{S}(\rrn)$,
and $f\in\mathcal{S}'(\rrn)$.
\begin{enumerate}
\item[(i)] The \emph{radial maximal
function}\index{radial maximal function}
$M(f,\phi)$\index{$M(f,\phi)$} is defined
by setting, for any $x\in\rrn$,
$$M(f,\phi)(x):=\sup_{t\in(0,\infty)}|f*\phi_{t}(x)|,$$
here and thereafter, for any $\phi\in\mathcal{S}(\rrn)$ and
$x\in\rrn$,
\begin{equation}\label{smaxe1}
\phi_{t}(x):=\frac{1}{t^n}\phi\left(\frac{x}{t}\right).
\end{equation}
\item[(ii)] The \emph{non-tangential maximal
function}\index{non-tangential
maximal function}
$M^{*}_{a}(f,\phi)$\index{$M^{*}_{a}$}, with aperture
$a\in(0,\infty)$, is defined by setting,
for any $x\in\rrn$,
$$M_{a}^{*}(f,\phi)(x):=\sup_{(y,t)
\in\Gamma_{a}(x)}|f*\phi_{t}(y)|,$$
here and thereafter, for any $x\in\rrn$
and $a\in(0,\infty)$,
\begin{equation}\label{gammaa}
\Gamma_{a}(x)\index{$\Gamma_{a}(x)$}:=
\{(y,t)\in\mathbb{R}^{n+1}_{+}:\ |y-x|<at\}.
\end{equation}
\item[(iii)] The \emph{maximal function
$M_{b}^{**}(f,\phi)$
of Peetre type}\index{maximal function of Peetre type}
is defined by setting, for any $x\in\rrn$,
$$M_{b}^{**}(f,\phi)(x)\index{$M_{b}^{**}$}:
=\sup_{(y,t)\in\rrnp}\frac{|f\ast\phi_{t}
(x-y)|}{(1+t^{-1}|y|)^{b}}.$$
\item[(iv)] The \emph{grand maximal
function $\mathcal{M}^{**}_{b,N}(f)$
of Peetre
type}\index{grand maximal function of Peetre type}
is defined
by setting, for any $x\in\rrn$,
$$\mathcal{M}^{**}_{b,N}(f)(x)\index{$\mathcal{M}^{**}_{b,N}$}
:=\sup_{\phi\in\mathcal{F}_{N}(\rrn)}
\sup_{(y,t)\in\mathbb{R}^{n+1}_{
+}}\frac{|f*\phi_{t}(x-y)|}
{(1+t^{-1}|y|)^{b}},$$
where $\mathcal{F}_{N}(\rrn)$ is
as in \eqref{defn}.
\end{enumerate}
\end{definition}

Via these maximal functions, we characterize
the generalized Herz--Hardy space
$\HaSaHo$ as
follows.\index{maximal function characterization}

\begin{theorem}\label{Th5.4}
Let $p$, $q$, $a,\ b\in(0,\infty)$, $\omega\in M(\rp)$,
$N\in\mathbb{N}$, and $\phi\in\mathcal{S}(\rrn)$ satisfy
$\int_{\rrn}\phi(x)\,dx\neq0.$
\begin{enumerate}
\item[\rm{(i)}] Let $N\in\mathbb{N}
\cap[\lfloor b+1\rfloor,\infty)$
and $\omega$ satisfy $\m0(\omega)
\in(-\frac{n}{p},\infty)$.
Then, for any $f\in\mathcal{S}'(\rrn)$,
$$
\|M(f,\phi)\|_{\HerzSo}\lesssim
\|M^*_{a}(f,\phi)\|_{\HerzSo}
\lesssim\|M^{**}_{b}(f,\phi)\|_{\HerzSo},
$$
\begin{align*}
\|M(f,\phi)\|_{\HerzSo}&\lesssim
\|\mathcal{M}_{N}(f)\|_{\HerzSo}
\lesssim\|\mathcal{M}_{\lfloor b+1\rfloor}(f)
\|_{\HerzSo}\\&\lesssim
\|M^{**}_{b}(f,\phi)\|_{\HerzSo},
\end{align*}
and
$$
\|M^{**}_{b}(f,\phi)\|_{\HerzSo}\sim
\|\mathcal{M}^{**}_{b,N}(f)\|_{\HerzSo},
$$
where the implicit positive
constants are independent of $f$.
\item[\rm{(ii)}] Let
$\omega\in M(\rp)$ satisfy
$\m0(\omega)\in(-\frac{n}{p},\infty)$
and $\mi(\omega)\in(-\frac{n}{p},\infty)$.
Assume $b\in(\max\{\frac{n}{p},\Mw+\frac{n}{p}\},\infty)$.
Then, for any $f\in\mathcal{S}'(\rrn),$
$$\|M^{**}_{b}(f,\phi)\|_{\HerzSo}
\lesssim\|M(f,\phi)\|_{\HerzSo},$$
where the implicit positive
constant is independent of $f$.
In particular,
when $N\in\mathbb{N}\cap[\lfloor b+1\rfloor,\infty)$,
if one of the quantities
$$\|M(f,\phi)\|_{\HerzSo},\
\|M^*_{a}(f,\phi)\|_{\HerzSo},\ \|\mathcal{M}_{N}
(f)\|_{\HerzSo},$$
$$
\|M^{**}_{b}(f,\phi)\|_{\HerzSo},
\ \text{and}\ \|\mathcal{M}^{**}_{b,N}(f)\|_{\HerzSo}$$
is finite, then the others are also finite and mutually
equivalent with the
positive equivalence
constants independent of $f$.
\end{enumerate}
\end{theorem}

In order to prove this
theorem, we need the following
maximal function characterizations
of Hardy spaces associated with
ball quasi-Banach function spaces,
which is just \cite[Theorem 3.1]{SHYY}.

\begin{lemma}\label{Th5.4l1}
Let $a,\ b\in(0,\infty)$,
$X$ be a ball quasi-Banach
function space, and $\phi\in
\mathcal{S}(\rrn)$ satisfy
$\int_{\rrn}\phi(x)\,dx\neq0.$
\begin{enumerate}
\item[\rm{(i)}] Let $N\in\mathbb{N}
\cap[\lfloor b+1\rfloor,\infty)$.
Then, for any $f\in\mathcal{S}'(\rrn)$,
$$
\|M(f,\phi)\|_{X}\lesssim
\|M^*_{a}(f,\phi)\|_{X}
\lesssim\|M^{**}_{b}(f,\phi)\|_{X},
$$
\begin{align*}
\|M(f,\phi)\|_{X}\lesssim
\|\mathcal{M}_{N}(f)\|_{X}
\lesssim\|\mathcal{M}_{\lfloor b+1\rfloor}(f)
\|_{X}\lesssim
\|M^{**}_{b}(f,\phi)\|_{X},
\end{align*}
and
$$
\|M^{**}_{b}(f,\phi)\|_{X}\sim
\|\mathcal{M}^{**}_{b,N}(f)\|_{X},
$$
where the implicit positive
constants are independent of $f$.
\item[\rm{(ii)}] Let $r\in(0,\infty)$.
Assume $b\in(\frac{n}{r},\infty)$ and
$\mc$ is bounded on $X^{1/r}$.
Then, for any $f\in\mathcal{S}'(\rrn),$
$$\|M^{**}_{b}(f,\phi)\|_{X}
\lesssim\|M(f,\phi)\|_{X},$$
where the implicit positive
constant is independent of $f$.
In particular,
when $N\in\mathbb{N}
\cap[\lfloor b+1\rfloor,\infty)$,
if one of the quantities
$$\|M(f,\phi)\|_{X},\
\|M^*_{a}(f,\phi)\|_{X},\ \|\mathcal{M}_{N}
(f)\|_{X},$$
$$
\|M^{**}_{b}(f,\phi)\|_{X},
\ \text{and}\ \|\mathcal{M}^{**}_{b,N}(f)\|_{X}$$
is finite, then the others are also finite and mutually
equivalent with the positive
equivalence constants independent of $f$.
\end{enumerate}
\end{lemma}

As was pointed out in Chapter \ref{sec3},
under some reasonable and sharp assumptions,
the local generalized Herz space $\HerzSo$
is a ball quasi-Banach function space.
This implies that, to show Theorem \ref{Th5.4},
we only need to prove that
$\HerzSo$ satisfies all the assumptions
of Lemma \ref{Th5.4l1}.

\begin{proof}[Proof of Theorem \ref{Th5.4}]
Let all the symbols be as in the present theorem.
Then, from the assumption
$\m0(\omega)\in(-\frac{n}{p},\infty)$ and
Theorem \ref{Th3}, it follows that the
local generalized Herz space $\HerzSo$ is a
BQBF space.
Applying this and Lemma \ref{Th5.4l1}(i),
we conclude that (i) holds true.

Next, we show (ii). To achieve this,
we only need to prove that,
under the assumptions
of the present theorem, $\HerzSo$
satisfies all the assumptions of
Lemma \ref{Th5.4l1}(ii), namely, there exists
an $r\in(0,\infty)$ such that
$b\in(\frac{n}{r},\infty)$
and the Hardy--Littlewood
maximal operator $\mc$ is
bounded on $[\HerzSo]^{1/r}$.
Indeed, from the assumptions
$\mw\in(-\frac{n}{p},\infty)$ and
$$b\in\left(\max\left\{\frac{n}{p},
\Mw+\frac{n}{p}\right\},\infty\right),$$
we deduce that
$$\frac{n}{b}\in\left(0,\min
\left\{p,\frac{n}{\Mw+n/p}\right\}\right).$$
Thus, we can choose an
\begin{equation}\label{Th5.4e1}
r\in\left(\frac{n}{b},\min
\left\{p,\frac{n}{\Mw+n/p}\right\}\right).
\end{equation}
By this and Lemma \ref{mbhl},
we find that $b\in(\frac{n}{r},\infty)$
and $\mc$ is bounded on the Herz space
$[\HerzSo]^{1/r}$.
This implies that all the assumptions
of Lemma \ref{Th5.4l1}(ii) hold true
for $\HerzSo$, which completes
the proof of (ii), and hence
of Theorem \ref{Th5.4}.
\end{proof}

\begin{remark}\label{remth5.4}
In Theorem \ref{Th5.4},
let $\omega(t):=t^\alpha$ for
any $t\in(0,\infty)$ and for any
given $\alpha\in\rr$.
Then, in this case, when
$p\in(1,\infty)$, the
maximal function characterizations
established in Theorem \ref{Th5.4}
widen the range of $\alpha\in(0,\infty)$
in \cite[Theorem 2.1]{LYH}
into $\alpha\in(-\frac{n}{p},\infty)$;
when $p\in(0,1]$, to the best of out
knowledge, the result of Theorem \ref{Th5.4}
is totally new.
\end{remark}

As an application, we obtain the
following maximal function
characterizations\index{maximal function characterization}
of the
generalized Morrey--Hardy space $H\MorrSo$, which
is immediately deduced from Theorem \ref{Th5.4}
and Remark \ref{remhs}(iv);
we omit the details.

\begin{corollary}
Let $a,\ b\in(0,\infty)$, $p$, $q\in[1,\infty)$,
$\omega\in M(\rp)$,
$N\in\mathbb{N}$, and $\phi\in\mathcal{S}(\rrn)$ satisfy
$\int_{\rrn}\phi(x)\,dx\neq0.$
\begin{enumerate}
\item[\rm{(i)}] Let
$N\in\mathbb{N}\cap[\lfloor b+1\rfloor,\infty)$
and $\omega$ satisfy $\MI(\omega)\in(-\infty,0)$ and
$$
-\frac{n}{p}<\m0(\omega)\leq\M0(\omega)<0.
$$
Then, for any $f\in\mathcal{S}'(\rrn)$,
$$
\|M(f,\phi)\|_{\MorrSo}\lesssim
\|M^*_{a}(f,\phi)\|_{\MorrSo}
\lesssim\|M^{**}_{b}(f,\phi)\|_{\MorrSo},
$$
\begin{align*}
\|M(f,\phi)\|_{\MorrSo}&\lesssim
\|\mathcal{M}_{N}(f)\|_{\MorrSo}
\lesssim\|\mathcal{M}_{\lfloor
b+1\rfloor}(f)\|_{\MorrSo}\\&\lesssim
\|M^{**}_{b}(f,\phi)\|_{\MorrSo},
\end{align*}
and
$$
\|M^{**}_{b}(f,\phi)\|_{\MorrSo}\sim\|
\mathcal{M}^{**}_{b,N}(f)\|_{\MorrSo},
$$
where the implicit positive
constants are independent of $f$.
\item[\rm{(ii)}] Let $\omega\in M(\rp)$ satisfy
$$
-\frac{n}{p}<\m0(\omega)\leq\M0(\omega)<0
$$
and
$$
-\frac{n}{p}<\m0(\omega)\leq\M0(\omega)<0.
$$
Assume $b\in(\frac{n}{p},\infty)$.
Then, for any $f\in\mathcal{S}'(\rrn),$
$$\|M^{**}_{b}(f,\phi)\|_{\MorrSo}
\lesssim\|M(f,\phi)\|_{\MorrSo},$$
where the implicit positive
constant is independent of $f$.
In particular, when $N\in\mathbb{N}
\cap[\lfloor b+1\rfloor,\infty)$,
if one of the quantities
$$\|M(f,\phi)\|_{\MorrSo},\
\|M^*_{a}(f,\phi)\|_{\MorrSo},\ \|\mathcal{M}_{N}
(f)\|_{\MorrSo},$$
$$
\|M^{**}_{b}(f,\phi)\|_{\MorrSo},
\ \text{and}\ \|\mathcal{M}^{**}
_{b,N}(f)\|_{\MorrSo}$$
is finite, then the others are also
finite and mutually
equivalent with the
positive equivalence
constants independent of $f$.
\end{enumerate}
\end{corollary}

Next, we show the following maximal
function
characterizations\index{maximal function characterization}
of the Hardy space $\HaSaH$
associated with the global generalized
Herz space $\HerzS$.

\begin{theorem}\label{Th5.6}
Let $p$, $q$, $a,\ b\in(0,\infty)$, $\omega\in M(\rp)$,
$N\in\mathbb{N}$, and $\phi\in\mathcal{S}(\rrn)$ satisfy
$\int_{\rrn}\phi(x)\,dx\neq0.$
\begin{enumerate}
\item[\rm{(i)}] Let $N\in\mathbb{N}
\cap[\lfloor b+1\rfloor,\infty)$
and $\omega$ satisfy $\m0(\omega)
\in(-\frac{n}{p},\infty)$ and
$\MI(\omega)\in(-\infty,0)$.
Then, for any $f\in\mathcal{S}'(\rrn)$,
$$
\|M(f,\phi)\|_{\HerzS}\lesssim
\|M^*_{a}(f,\phi)\|_{\HerzS}
\lesssim\|M^{**}_{b}(f,\phi)\|_{\HerzS},
$$
\begin{align*}
\|M(f,\phi)\|_{\HerzS}&\lesssim
\|\mathcal{M}_{N}(f)\|_{\HerzS}
\lesssim\|\mathcal{M}_{\lfloor
b+1\rfloor}(f)\|_{\HerzS}\\&\lesssim
\|M^{**}_{b}(f,\phi)\|_{\HerzS},
\end{align*}
and
$$
\|M^{**}_{b}(f,\phi)\|_{\HerzS}\sim
\|\mathcal{M}^{**}_{b,N}(f)\|_{\HerzS},
$$
where the implicit positive
constants are independent of $f$.
\item[\rm{(ii)}] Let $\omega\in M(\rp)$ satisfy
$\m0(\omega)\in(-\frac{n}{p},\infty)$ and
$$-\frac{n}{p}<\mi(\omega)\leq\MI(\omega)<0.$$
Assume $b\in(\max\{\frac{n}{p},\Mw+\frac{n}{p}\},\infty)$.
Then, for any $f\in\mathcal{S}'(\rrn),$
$$\|M^{**}_{b}(f,\phi)\|_{\HerzS}
\lesssim\|M(f,\phi)\|_{\HerzS},$$
where the implicit positive
constant is independent of $f$.
In particular, when $N\in\mathbb{N}
\cap[\lfloor b+1\rfloor,\infty)$,
if one of the quantities
$$\|M(f,\phi)\|_{\HerzS},\
\|M^*_{a}(f,\phi)\|_{\HerzS},\ \|\mathcal{M}_{N}
(f)\|_{\HerzS},$$
$$
\|M^{**}_{b}(f,\phi)\|_{\HerzS},
\ \text{and}\ \|\mathcal{M}^{**}_{b,N}(f)\|_{\HerzS}$$
is finite, then the others are also finite and mutually
equivalent with the
positive equivalence constants independent of $f$.
\end{enumerate}
\end{theorem}

To prove this theorem, we first
present a lemma about the boundedness of the
Hardy--Littlewood maximal operator as follows.

\begin{lemma}\label{mbhg}
Let $p,\ q\in(0,\infty)$ and $\omega\in M(\rp)$ satisfy
$\m0(\omega)\in(-\frac{n}{p},
\infty)$ and $\mi(\omega)\in
(-\frac{n}{p},\infty)$.
Then, for any given
$r\in(0,\min\{p,\frac{n}{\Mw+n/p}\})$,
there exists a positive constant $C$ such that,
for any $f\in L^1_{\mathrm{loc}}(\rrn)$,
$$
\left\|\mc(f)\right\|_{[\HerzS]^{1/r}}
\leq C\|f\|_{[\HerzS]^{1/r}}.
$$
\end{lemma}

\begin{proof}
Let all the symbols be as in the present lemma.
From \eqref{mbhle1} and \eqref{mbhle2},
it follows that
\begin{align*}
\Mwr<\frac{n}{(p/r)'}
\end{align*}
and
\begin{align*}
\mwr>-\frac{n}{p/r}.
\end{align*}
This, together with the fact that
$\frac{p}{r}\in(1,\infty)$
and Corollary \ref{Th4},
further implies that, for any $f\in
L^1_{\mathrm{loc}}(\rrn)$,
$$
\left\|\mc(f)\right\|_{\HerzScr}
\lesssim\|f\|_{\HerzScr}.$$
Therefore, by Lemma \ref{convexl},
we further find that the Hardy--Littlewood
maximal operator $\mc$ is bounded on $[\HerzS]^{1/r}$.
This finishes the proof of Lemma \ref{mbhg}.
\end{proof}

Via Lemma \ref{mbhg} and the
maximal function characterizations
of Hardy spaces associated with
ball quasi-Banach function spaces
presented in Lemma \ref{Th5.4l1},
we now show Theorem \ref{Th5.6}.

\begin{proof}[Proof of Theorem \ref{Th5.6}]
Let all the symbols be as in the present theorem.
Then, using the assumptions
$\m0(\omega)\in(-\frac{n}{p},\infty)$ and
$\MI(\omega)\in(-\infty,0)$, and
Theorem \ref{Th2}, we conclude that the
global generalized Herz space $\HerzS$
is a BQBF space.
This, together with Lemma \ref{Th5.4l1}(i),
then finishes the proof of (i).

Next, we show (ii).
Indeed, let
$r$ be as in \eqref{Th5.4e1}.
Then, by Lemma \ref{mbhg},
we find that the
Hardy--Littlewood maximal operator
$\mc$ is bounded on $[\HerzS]^{1/r}$.
From this and Lemma \ref{Th5.4l1}(ii),
it follows that
(ii) holds true, which then completes
the proof of Theorem \ref{Th5.6}.
\end{proof}

Moreover, using the theorem above and Remark \ref{remhs}(iv),
we conclude the following maximal function
characterizations\index{maximal function characterization}
of the generalized Morrey--Hardy
space $H\MorrS$; we omit the
details.

\begin{corollary}
Let $a,\ b\in(0,\infty)$, $p$,
$q\in[1,\infty)$, $\omega\in M(\rp)$,
$N\in\mathbb{N}$, and
$\phi\in\mathcal{S}(\rrn)$ satisfy
$\int_{\rrn}\phi(x)\,dx\neq0.$
\begin{enumerate}
\item[\rm{(i)}] Let $N\in\mathbb{N}
\cap[\lfloor b+1\rfloor,\infty)$
and $\omega$ satisfy
$\MI(\omega)\in(-\infty,0)$ and
$$
-\frac{n}{p}<\m0(\omega)\leq\M0(\omega)<0.
$$
Then, for any $f\in\mathcal{S}'(\rrn)$,
$$
\|M(f,\phi)\|_{\MorrS}\lesssim
\|M^*_{a}(f,\phi)\|_{\MorrS}
\lesssim\|M^{**}_{b}(f,\phi)\|_{\MorrS},
$$
\begin{align*}
\|M(f,\phi)\|_{\MorrS}&\lesssim
\|\mathcal{M}_{N}(f)\|_{\MorrS}
\lesssim\|\mathcal{M}_{\lfloor
b+1\rfloor}(f)\|_{\MorrS}\\&\lesssim
\|M^{**}_{b}(f,\phi)\|_{\MorrS},
\end{align*}
and
$$
\|M^{**}_{b}(f,\phi)\|_{\MorrS}\sim\|
\mathcal{M}^{**}_{b,N}(f)\|_{\MorrS},
$$
where the implicit positive
constants are independent of $f$.
\item[\rm{(ii)}] Let $\omega\in M(\rp)$ satisfy
$$
-\frac{n}{p}<\m0(\omega)\leq\M0(\omega)<0
$$
and
$$
-\frac{n}{p}<\m0(\omega)\leq\M0(\omega)<0.
$$
Assume $b\in(\frac{n}{p},\infty)$.
Then, for any $f\in\mathcal{S}'(\rrn),$
$$\|M^{**}_{b}(f,\phi)\|_{\MorrS}
\lesssim\|M(f,\phi)\|_{\MorrS},$$
where the implicit positive
constant is independent of $f$.
In particular, when
$N\in\mathbb{N}\cap[\lfloor b+1\rfloor,\infty)$,
if one of the quantities
$$\|M(f,\phi)\|_{\MorrS},\
\|M^*_{a}(f,\phi)\|_{\MorrS},\ \|\mathcal{M}_{N}
(f)\|_{\MorrS},$$
$$
\|M^{**}_{b}(f,\phi)\|_{\MorrS},
\ \text{and}\
\|\mathcal{M}^{**}_{b,N}(f)\|_{\MorrS}$$
is finite, then the others are
also finite and mutually
equivalent with the
positive equivalence constants
independent of $f$.
\end{enumerate}
\end{corollary}

\section{Relations with Generalized Herz Spaces}

In this section,
we investigate the relations between
generalized Herz spaces and the
associated Hardy spaces.
The following conclusion shows that,
under some assumptions on the exponents,
the local generalized Herz space $\HerzSo$
coincides with the generalized
Herz--Hardy space $\HaSaHo$
in the sense of equivalent quasi-norms.

\begin{theorem}\label{Th5.7}
Let $p\in(1,\infty)$,
$q\in(0,\infty)$, and
$\omega\in M(\rp)$ satisfy
$$-\frac{n}{p}<\m0(\omega)\leq\M0(\omega)<
\frac{n}{p'}$$
and
$$-\frac{n}{p}<\mi(\omega)
\leq\MI(\omega)<\frac{n}{p'},$$
where $\frac{1}{p}+\frac{1}{p'}=1$.
Then
\begin{enumerate}
\item[\rm{(i)}] $\HerzSo
\hookrightarrow\mathcal{S}'(\rrn)$.
\item[\rm{(ii)}] If $f\in\HerzSo$,
then $f\in\HaSaHo$
and there exists a positive constant $C$,
independent
of $f$, such that
$$\|f\|_{\HaSaHo}
\leq C\|f\|_{\HerzSo}.$$
\item[\rm{(iii)}] If $f\in\HaSaHo$,
then there exists a
locally integrable function $g\in\HerzSo$
such that $g$ represents $f$,
which means that $f=g$ in
$\mathcal{S}'(\rrn)$,
$$\|f\|_{\HaSaHo}=\|g\|_{\HaSaHo},$$
and there exists a positive constant $C$,
independent of $f$,
such that
$$\|g\|_{\HerzSo}\leq C\|f\|_{\HaSaHo}.$$
\end{enumerate}
\end{theorem}

To prove this theorem,
recall that Sawano et al.\
\cite[Theorem 3.4]{SHYY} obtained
the relation between the ball
quasi-Banach function space
$X$ and the associated Hardy space
$H_{X}(\rrn)$
(see also Lemma \ref{Th5.7l1}
below). Thus, to prove
Theorem \ref{Th5.7},
it suffices to show that
all the assumptions
of \cite[Theorem 3.4]{SHYY}
hold true for local generalized
Herz spaces. For this purpose,
we first present the
definition of the Hardy space
$H_X(\rrn)$ as follows,
which was given in
\cite[Definition 2.22]{SHYY}.

\begin{definition}\label{hardyx}
Let $X$ be a ball quasi-Banach
function space and $N\in
\mathbb{N}$. Then the \emph{Hardy
space} $H_X(\rrn)$\index{$H_X(\rrn)$} is defined to be
the set of all the $f\in\mathcal{S}'(\rrn)$
such that
$$
\left\|f\right\|_{H_X(\rrn)}
:=\left\|\mc_N(f)\right\|_{X}<\infty.
$$
\end{definition}

The following lemma is just
\cite[Theorem 3.4]{SHYY},
which is vital in the proof
of Theorem \ref{Th5.7}.

\begin{lemma}\label{Th5.7l1}
Let $r\in(1,\infty)$
and $X$ be a ball quasi-Banach
function space satisfy that $\mc$
is bounded on $X^{1/r}$.
Then
\begin{enumerate}
\item[\rm{(i)}] $X
\hookrightarrow\mathcal{S}'(\rrn)$.
\item[\rm{(ii)}] If $f\in X$,
then $f\in H_{X}(\rrn)$
and there exists a positive constant $C$,
independent
of $f$, such that
$$\|f\|_{H_{X}(\rrn)}
\leq C\|f\|_{X}.$$
\item[\rm{(iii)}] If $f\in H_{X}(\rrn)$,
then there exists a
locally integrable function $g\in X$
such that $g$ represents $f$,
which means that $f=g$ in
$\mathcal{S}'(\rrn)$,
$$\|f\|_{H_{X}(\rrn)}=\|g\|_{H_{X}(\rrn)},$$
and there exists a positive constant $C$,
independent of $f$,
such that
$$\|g\|_{X}\leq C\|f\|_{H_{X}(\rrn)}.$$
\end{enumerate}
\end{lemma}

Now, we show Theorem \ref{Th5.7}.

\begin{proof}[Proof of Theorem \ref{Th5.7}]
Let all the symbols be as in the present theorem.
Then, combining the assumption
$\m0(\omega)\in(-\frac{n}{p},\infty)$
and Theorem \ref{Th3}, we
find that the local generalized Herz space
$\HerzSo$ is a BQBF space.
This implies that,
to complete the proof of the present theorem,
we only need to show that
$\HerzSo$ satisfies
all the assumptions of Lemma \ref{Th5.7l1}
under the assumptions of
the present theorem.
Indeed, from the assumptions $p\in(1,\infty)$ and
$\Mw\in(-\frac{n}{p},\frac{n}{p'})$, we deduce that
$$
\min\left\{p,\frac{n}{\Mw+n/p}\right\}\in(1,\infty).
$$
Therefore, we can choose an
$r\in(1,\min\{p,\frac{n}{\Mw+n/p}\})$.
For this $r$, by Lemma \ref{mbhl},
we conclude that, for any
$f\in L^1_{\mathrm{loc}}(\rrn)$,
$$
\left\|\mc(f)\right\|_{[\HerzSo]^{1/r}}
\lesssim\|f\|_{[\HerzSo]^{1/r}},
$$
which further implies that
$\mc$ is bounded on
$[\HerzSo]^{1/r}$, and hence
all the assumptions of Lemma \ref{Th5.7l1}
are satisfied for $\HerzSo$.
This finishes the proof of
Theorem \ref{Th5.7}.
\end{proof}

\begin{remark}
We should point out that, in Theorem \ref{Th5.7},
if $\omega(t):=t^\alpha$ for any $t\in(0,\infty)$
and for any given $\alpha\in\rr$,
then Theorem \ref{Th5.7}
goes back to \cite[Proposition 2.1.1(1)]{LYH}.
\end{remark}

Via Theorem \ref{Th5.7}
and Remark \ref{remark4.10}(ii),
we immediately obtain the following
corollary which shows that,
under some assumptions,
$H\MorrSo=\MorrSo$
with equivalent norms; we omit the details.

\begin{corollary}\label{Th5.7m}
Let $p\in(1,\infty)$, $q\in[1,\infty)$, and
$\omega\in M(\rp)$ satisfy
$$-\frac{n}{p}<\m0(\omega)\leq\M0(\omega)<0$$
and
$$-\frac{n}{p}<\mi(\omega)\leq\MI(\omega)<0.$$
Then
\begin{enumerate}
\item[\rm{(i)}] $\MorrSo\hookrightarrow\mathcal{S}'(\rrn)$.
\item[\rm{(ii)}] If $f\in\MorrSo$, then $f\in H\MorrSo$
and there exists a positive constant $C$, independent
of $f$, such that
$$\|f\|_{H\MorrSo}\leq C\|f\|_{\MorrSo}.$$
\item[\rm{(iii)}] If $f\in H\MorrSo$, then there exists a
locally integrable function $g\in\MorrSo$
such that $g$ represents $f$, which means that $f=g$ in
$\mathcal{S}'(\rrn)$,
$$\|f\|_{H\MorrSo}=\|g\|_{H\MorrSo},$$
and there exists a positive constant $C$, independent of $f$,
such that
$$\|g\|_{\MorrSo}\leq C\|f\|_{H\MorrSo}.$$
\end{enumerate}
\end{corollary}

Applying the known result
of ball quasi-Banach function spaces
mentioned in Lemma \ref{Th5.7l1},
we also have the following conclusion
which shows
that $\HaSaH=\HerzS$ with equivalent
quasi-norms under some reasonable
and sharp assumptions.

\begin{theorem}\label{Th5.8}
Let $p\in(1,\infty)$, $q\in(0,\infty)$, and
$\omega\in M(\rp)$ satisfy
$$-\frac{n}{p}<\m0(\omega)\leq
\M0(\omega)<\frac{n}{p'}$$
and
$$
-\frac{n}{p}<\mi(\omega)
\leq\MI(\omega)<0,
$$
where $\frac{1}{p}+\frac{1}{p'}=1$.
Then
\begin{enumerate}
\item[\rm{(i)}] $\HerzS\hookrightarrow\mathcal{S}'(\rrn)$.
\item[\rm{(ii)}] If $f\in\HerzS$, then $f\in\HaSaH$
and there exists a positive constant $C$, independent
of $f$, such that
$$\|f\|_{\HaSaH}\leq C\|f\|_{\HerzS}.$$
\item[\rm{(iii)}] If $f\in\HaSaH$, then there exists a
locally integrable function $g\in\HerzS$
such that $g$ represents $f$, which means that $f=g$ in
$\mathcal{S}'(\rrn)$,
$$\|f\|_{\HaSaH}=\|g\|_{\HaSaH},$$
and there exists a positive constant $C$, independent of $f$,
such that
$$\|g\|_{\HerzS}\leq C\|f\|_{\HaSaH}.$$
\end{enumerate}
\end{theorem}

\begin{proof}
Let all the symbols be as in the present theorem.
Then, from the assumptions
of the present theorem
and Theorem \ref{Th2}, we deduce
that the global generalized Herz space
$\HerzS$ is a BQBF space.
In addition, combining the
assumptions $p\in(1,\infty)$ and
$\Mw\in(-\frac{n}{p},\frac{n}{p'})$, we conclude that
$$
\min\left\{p,\frac{n}{\Mw+n/p}\right\}\in(1,\infty).
$$
Let
$r\in(1,\min\{p,\frac{n}{\Mw+n/p}\})$.
Then, by Lemma \ref{mbhg},
we find that
$\mc$ is bounded on $[\HerzS]^{1/r}$.
Using this and Lemma \ref{Th5.7l1},
we then complete the proof of
Theorem \ref{Th5.8}.
\end{proof}

Moreover, from both Theorem \ref{Th5.8} and
Remark \ref{remark4.10}(ii),
we immediately deduce the following relation
between the global generalized Morrey space
$\MorrS$ and the
associated Hardy space $H\MorrS$; we omit the details.

\begin{corollary}
Let $p$, $q$, and $\omega$ be as in Corollary \ref{Th5.7m}.
Then
\begin{enumerate}
\item[\rm{(i)}] $\MorrS\hookrightarrow\mathcal{S}'(\rrn)$.
\item[\rm{(ii)}] If $f\in\MorrS$, then $f\in H\MorrS$
and there exists a positive constant $C$, independent
of $f$, such that
$$\|f\|_{H\MorrS}\leq C\|f\|_{\MorrS}.$$
\item[\rm{(iii)}] If $f\in H\MorrS$, then there exists a
locally integrable function $g\in\MorrS$
such that $g$ represents $f$, which means that $f=g$ in
$\mathcal{S}'(\rrn)$,
$$\|f\|_{H\MorrS}=\|g\|_{H\MorrS},$$
and there exists a positive constant $C$, independent of $f$,
such that
$$\|g\|_{\MorrS}\leq C\|f\|_{H\MorrS}.$$
\end{enumerate}
\end{corollary}

\section{Atomic Characterizations}

The main target of this section
is to characterize
generalized Herz--Hardy spaces via
atoms. By the known
atomic characterization
of Hardy spaces associated with
ball quasi-Banach function spaces,
we can directly obtain
the atomic characterization
of the Hardy space $\HaSaHo$.
However,
the proof of the
atomic characterization
of $\HaSaH$ is more complex because of
the deficiency of the associate space
of $\HerzS$. Indeed,
in order to show the atomic
characterization of $\HaSaH$,
we establish an improved atomic characterization
of Hardy spaces associated with ball
quasi-Banach function spaces
under no assumption about associate spaces.
Via this conclusion and
some other auxiliary lemmas about
global generalized Herz spaces
and block spaces,
we obtain the atomic characterization
of $\HaSaH$ and get rid of
the dependence on associate spaces.

First, we establish the
atomic characterization of
the generalized Herz--Hardy space
$\HaSaHo$. To achieve
this, we introduce the definition
of $(\HerzSo,\,r,\,d)$-atoms
as follows.

\begin{definition}\label{atom}
Let $p$, $q\in(0,\infty)$, $\omega\in M(\rp)$
with $\m0(\omega)
\in(-\frac{n}{p},\infty)$,
$r\in[1,\infty]$, and $d\in\zp$.
Then a measurable function
$a$ on $\rrn$ is called a
\emph{$(\HerzSo,\,r,
\,d)$-atom}\index{$(\HerzSo,\,r,\,d)$-atom}
if there exists a ball $B\in\mathbb{B}$
such that
\begin{enumerate}
  \item[{\rm(i)}] $\supp(a)
  :=\{x\in\rrn:\ a(x)\neq0\}\subset B$;
  \item[{\rm(ii)}] $\|a\|_{L^r(\rrn)}
  \leq\frac{|B|^{1/r}}{\|\1bf_{B}\|_{\HerzSo}}$;
  \item[{\rm(iii)}] for any $\alpha\in\zp^n$
  with $|\alpha|\leq d$,
  $$\int_{\rrn}a(x)x^\alpha\,dx=0.$$
\end{enumerate}
\end{definition}

Via $(\HerzSo,\,r,\,d)$-atoms,
we now introduce the following
generalized atomic
Herz--Hardy spaces associated with
local generalized Herz spaces.

\begin{definition}\label{atsl}
Let $p$, $q\in(0,\infty)$, $\omega\in M(\rp)$ with
$\m0(\omega)\in(-\frac{n}{p},\infty)$ and
$\mi(\omega)\in(-\frac{n}{p},\infty)$,
$r\in(\max\{1,p,\frac{n}{\mw+n/p}\},\infty]$,
$$s\in\left(0,\min\left\{
1,p,q,\frac{n}{\Mw+n/p}\right\}\right),$$
and $d\geq\lfloor n(
1/s-1)\rfloor$ be a fixed integer.
Then the \emph{generalized atomic
Herz--Hardy space}\index{generalized
atomic Herz--Hardy \\space} $\aHaSaHo$\index{$\aHaSaHo$}
is defined to be the set of all the
$f\in\mathcal{S}'(\rrn)$ such that there exist
$\{\lambda_j\}_{j\in\mathbb{N}}
\subset[0,\infty)$ and
a sequence $\{a_{j}\}_{j\in\mathbb{N}}$ of
$(\HerzSo,\,r,\,d)$-atoms supported,
respectively, in the balls $\{B_{j}\}_{j\in\mathbb{N}}
\subset\mathbb{B}$ such that
$$
f=\sum_{j\in\mathbb{N}}\lambda_{j}a_{j}
$$
in $\mathcal{S}'(\rrn)$,
and
$$
\left\|\left\{
\sum_{j\in\mathbb{N}}
\left[\frac{\lambda_j}{\|
\1bf_{B_j}\|_{\HerzSo}}\right]^s
\1bf_{B_j}\right\}^
{\frac{1}{s}}\right\|_{\HerzSo}<\infty.
$$
Moreover, for any $f\in\aHaSaHo$, let
$$
\|f\|_{\aHaSaHo}:=\inf\left\{
\left\|\left\{\sum_{j\in\mathbb{N}}
\left[\frac{\lambda_{j}}
{\|\1bf_{B_{j}}\|_{\HerzSo}}\right]^{s}
\1bf_{B_{j}}\right\}^{\frac{1}{s}}
\right\|_{\HerzSo}\right\},
$$
where the infimum is taken over all
the decompositions of $f$ as above.
\end{definition}

Then we have the following
atomic characterization\index{atomic
characterization} of the
generalized Herz--Hardy space $\HaSaHo$.

\begin{theorem}\label{Atol}
Let $p$, $q$, $\omega$, $d$, $s$,
and $r$ be as in Definition \ref{atsl}.
Then $$\HaSaHo=\aHaSaHo$$
with equivalent quasi-norms.
\end{theorem}

To show this theorem,
we first recall the definition
of atoms of general ball quasi-Banach
function spaces as follows, which
is just \cite[Definition 3.5]{SHYY}.

\begin{definition}\label{atomx}
Let $X$ be a ball quasi-Banach
function space,
$r\in[1,\infty]$, and $d\in\zp$.
Then a measurable function
$a$ on $\rrn$ is called an
\emph{$(X,\,r,
\,d)$-atom}\index{$(X,\,r,\,d)$-atom}
if there exists a ball $B\in\mathbb{B}$
such that
\begin{enumerate}
  \item[{\rm(i)}] $\supp(a)
  :=\{x\in\rrn:\ a(x)\neq0\}\subset B$;
  \item[{\rm(ii)}] $\|a\|_{L^r(\rrn)}
  \leq\frac{|B|^{1/r}}{\|\1bf_{B}\|_{X}}$;
  \item[{\rm(iii)}] for any $\alpha\in\zp^n$
  with $|\alpha|\leq d$,
  $$\int_{\rrn}a(x)x^\alpha\,dx=0.$$
\end{enumerate}
\end{definition}

The following technical
lemma about atoms can be concluded by
both Definition \ref{atomx} and the
H\"{o}lder inequality immediately;
we omit the details.

\begin{lemma}\label{Atoll2}
Let $X$ be a ball quasi-Banach
function space, $r,\ t\in[1,\infty]$
with $r<t$,
and $d\in\zp$.
Assume $a$ is an $(X,\,t,\,d)$-atom
supported in the ball $B\in\mathbb{B}$.
Then $a$ is an $(X,\,r,\,d)$-atom
supported in $B$.
\end{lemma}

Note that Sawano et al.\ \cite{SHYY}
established the atomic characterization
of Hardy spaces associated with ball
quasi-Banach function spaces.
Indeed, the following conclusion
is just the atomic reconstruction
theorem from \cite[Theorem 3.6]{SHYY}.
This conclusion plays a key
role in the proof of the
atomic characterization
of $\HaSaHo$.

\begin{lemma}\label{Atoll1}
Let $X$ be a ball quasi-Banach
function space satisfying Assumption
\ref{assfs}
with some $\theta$, $s\in(0,1]$,
$d\geq\lfloor n
(1/\theta-1)\rfloor$ be
a fixed integer, and $r\in(1,\infty]$.
Assume that $X^{1/s}$ is a
ball Banach function space and
there exists a positive
constant $C$ such that,
for any $f\in L^1_{\mathrm{loc}}(\rrn)$,
\begin{equation}\label{Atoll1e1}
\left\|\mc^{((r/s)')}(f)
\right\|_{(X^{1/s})'}\leq C
\left\|f\right\|_{(X^{1/s})'}.
\end{equation}
Let $\{a_j\}_{j\in\mathbb{N}}$
be a sequence of $(X,\,r,\,d)$-atoms
supported, respectively,
in the balls $\{B_j\}_{j\in\mathbb{N}}
\subset\mathbb{B}$, and $\{\lambda_j\}_{
j\in\mathbb{N}}\subset[0,\infty)$
such that
$f:=\sum_{j\in\mathbb{N}}\lambda_j
a_j$ converges in $\mathcal{S'}(\rrn)$,
and
$$
\left\|\left[\sum_{j\in\mathbb
N}\left(\frac{\lambda_j}{\|\1bf_{B_j}\|_{X}}
\right)^s\1bf_{B_j}
\right]^{\frac{1}{s}}\right\|_{X}<\infty.
$$
Then $f\in H_X(\rrn)$ and
$$
\|f\|_{H_{X}(\rrn)}\lesssim
\left\|\left[\sum_{j\in\mathbb
N}\left(\frac{\lambda_j}{\|\1bf_{B_j}\|_{X}}
\right)^s\1bf_{B_j}
\right]^{\frac{1}{s}}\right\|_{X},
$$
where the implicit positive constant
is independent of $f$.
\end{lemma}

Moreover, the following
atomic decomposition theorem
of $H_X(\rrn)$
was obtained in
\cite[Theorem 3.7]{SHYY},
which is also a vital tool in the
proof of the atomic characterization
of $\HaSaHo$.

\begin{lemma}\label{Atogl4}
Let $X$ be a ball quasi-Banach function space
satisfying Assumption \ref{assfs}
with some $\theta$, $s\in(0,1]$, and
$d\geq\lfloor n
(1/\theta-1)\rfloor$ be a fixed integer.
Then, for any $f\in H_X(\rrn)$,
there exist $\{\lambda_j\}_
{j\in\mathbb{N}}\subset[0,\infty)$ and
a sequence $\{a_j\}_{j\in\mathbb{N}}$
of $(X,\,\infty,\,d)$-atoms supported,
respectively, in the balls
$\{B_j\}_{j\in\mathbb{N}}\subset
\mathbb{B}$ such that
$f=\sum_{j\in\mathbb{N}}\lambda_ja_j$
in $\mathcal{S}'(\rrn)$, and
$$
\left\|\left[\sum_{j\in\mathbb
N}\left(\frac{\lambda_j}{\|\1bf_{B_j}\|_{X}}
\right)^s\1bf_{B_j}
\right]^{\frac{1}{s}}\right\|_{X}\lesssim
\|f\|_{H_X(\rrn)},
$$
where the implicit positive constant
is independent of $f$.
\end{lemma}

In order to prove Theorem \ref{Atol},
we also require two auxiliary lemmas
about the Fefferman--Stein
vector-valued inequality on
local generalized Herz spaces as follows.

\begin{lemma}\label{vmbhl}
Let $p,\ q\in(0,\infty)$ and $\omega\in M(\rp)$ satisfy
$\m0(\omega)\in(-\frac{n}{p},\infty)$
and $\mi(\omega)\in(-\frac{n}{p},\infty)$.
Then, for any given $u\in(1,\infty)$ and
$$r\in\left(0,\min\left\{p,\frac{n}
{\Mw+n/p}\right\}\right),$$ there exists
a positive constant $C$ such that, for any $\{f
_{j}\}_{j\in\mathbb{N}}\subset L_{{\rm loc}}^{1}(\rrn)$,
\begin{equation*}
\left\|\left\{\sum_{j\in\mathbb{N}}\left[\mc
(f_{j})\right]^{u}\right\}^
{\frac{1}{u}}\right\|_{[\HerzSo]^{1/r}}
\leq C\left\|\left\{
\sum_{j\in\mathbb{N}}|f_{j}|^{u}\right\}
^{\frac{1}{u}}\right\|_{[\HerzSo
]^{1/r}}.
\end{equation*}
\end{lemma}

\begin{proof}
Let all the symbols be as
in the present lemma.
Then, from \eqref{mbhle1}
and \eqref{mbhle2},
it follows that
\begin{align*}
\Mwr<\frac{n}{(p/r)'}
\end{align*}
and
\begin{align*}
\mwr>-\frac{n}{p/r}.
\end{align*}
By this, the assumption
$\frac{p}{r}\in(1,\infty)$, Theorem \ref{Th3.4},
and Lemma \ref{convexll}, we conclude that,
for any $\{f_{j}\}_{j\in\mathbb{N}}\subset
L_{{\rm loc}}^{1}(\rrn)$,
\begin{equation*}
\left\|\left\{\sum_{j\in\mathbb{N}}\left[\mc
(f_{j})\right]^{u}\right\}^{\frac{1}{u}}
\right\|_{[\HerzSo]^{1/r}}
\lesssim\left\|\left\{
\sum_{j\in\mathbb{N}}|f_{j}|^{u}\right\}
^{\frac{1}{u}}\right\|_{[\HerzSo
]^{1/r}},
\end{equation*}
which completes the proof of Lemma \ref{vmbhl}.
\end{proof}

\begin{lemma}\label{Atoll3}
Let $p$, $q\in(0,\infty)$,
$\omega\in M(\rp)$ satisfy
$\m0(\omega)\in(-\frac{n}{p},\infty)$
and $\mi(\omega)\in(-\frac{n}{p},
\infty)$, $s\in(0,\infty)$,
and $$
\theta\in\left(0,\min\left\{
s,p,\frac{n}{\Mw+n/p}\right\}\right).
$$
Then there exists a positive
constant $C$ such that, for any
$\{f_j\}_{j\in\mathbb{N}}\subset
L^1_{\mathrm{loc}}(\rrn)$,
$$
\left\|
\left\{
\sum_{j\in\mathbb{N}}
\left[\mc^{(\theta)}(f_j)\right]^{s}
\right\}^{1/s}
\right\|_{\HerzSo}
\leq C\left\|
\left(\sum_{j\in\mathbb{N}}
|f_j|^{s}
\right)^{1/s}
\right\|_{\HerzSo}.
$$
\end{lemma}

\begin{proof}
Let all the symbols be
as in the present lemma.
By the assumption $\theta\in(0,s)$,
we find that
$\frac{s}{\theta}\in(1,\infty)$.
Then, using Lemma \ref{vmbhl}
with $u=\frac{s}{\theta}$
and $r=\frac{1}{\theta}$,
we conclude that,
for any $\{f_{j}\}_{
j\subset\mathbb{N}}
\subset L^1_{\mathrm{loc}}(\rrn)$,
$$
\left\|
\left\{
\sum_{j\in\mathbb{N}}
\left[\mc(f_j)\right]^{s/\theta}
\right\}^{\theta/s}
\right\|_{[\HerzSo]^{1/\theta}}
\lesssim\left\|
\left(\sum_{j\in\mathbb{N}}
|f_j|^{s/\theta}
\right)^{\theta/s}
\right\|_{
[\HerzSo]^{1/\theta}
},
$$
which, combined with Remark \ref{power}(ii)
with $X$ therein replaced by $\HerzSo$,
further implies that
$$
\left\|
\left\{
\sum_{j\in\mathbb{N}}
\left[\mc^{(\theta)}(f_j)\right]^{s}
\right\}^{1/s}
\right\|_{\HerzSo}
\lesssim\left\|
\left(\sum_{j\in\mathbb{N}}
|f_j|^{s}
\right)^{1/s}
\right\|_{\HerzSo}.
$$
This finishes the proof of
Lemma \ref{Atoll3}.
\end{proof}

Via above lemmas, we
now show the atomic
characterization of the
generalized Herz--Hardy space
$\HaSaHo$.

\begin{proof}[Proof of Theorem \ref{Atol}]
Let $p$, $q$, $\omega$,
$d$, $s$, and $r$ be as in the
present theorem.
By the definition of
$\lfloor n(1/s-1)\rfloor$, we find that
$$
\left\lfloor
n\left(\frac{1}{s}-1\right)
\right\rfloor\leq n\left(
\frac{1}{s}-1
\right)<\left\lfloor
n\left(\frac{1}{s}-1\right)
\right\rfloor+1.
$$
Then we can choose a
$\theta\in(0,s)$ such that
$$
\left\lfloor
n\left(\frac{1}{s}-1\right)
\right\rfloor\leq n\left(
\frac{1}{\theta}-1
\right)<\left\lfloor
n\left(\frac{1}{s}-1\right)
\right\rfloor+1,
$$
which further implies that
$d\geq \lfloor n(1/s-1)\rfloor
=\lfloor n(1/\theta-1)\rfloor$.
Now, we claim that
$\HerzSo$ satisfies Assumption
\ref{assfs} for this $\theta$
and $s$, namely, $\HerzSo$
is a BQBF space and, for any
$\{f_{j}\}_{j\in\mathbb{N}}
\subset L^1_{{\rm loc}}(\rrn)$,
\begin{equation}\label{Atole1}
\left\|\left\{\sum_{j\in\mathbb{N}}
\left[\mc^{(\theta)}
(f_{j})\right]^s\right\}^{1/s}
\right\|_{\HerzSo}\lesssim
\left\|\left(\sum_{j\in\mathbb{N}}
|f_{j}|^{s}\right)^{1/s}\right\|_{\HerzSo}.
\end{equation}
Indeed, from the assumption
$\m0(\omega)\in(-\frac{n}{p},\infty)$
and Theorem \ref{Th3}, it
follows that the local generalized
Herz space $\HerzSo$ is a BQBF space.
In addition, using Lemma \ref{Atoll3},
we conclude that \eqref{Atole1}
holds true. This finishes
the proof of the above claim.

Next, we show that
$\aHaSaHo\subset\HaSaHo$. To
this end, we first prove that,
under the assumptions
of the present theorem,
$\HerzSo$ satisfies all the assumptions
of Lemma \ref{Atoll1}.
Indeed, applying Lemma \ref{mbhal},
we find that $[\HerzSo]^{1/s}$ is
a BBF space and, for any $f\in
L_{{\rm loc}}^{1}(\rrn)$,
\begin{equation*}
\left\|\mc^{((r/s)')}(f)\right\|_
{([\HerzSo]^{1/s})'}\lesssim
\left\|f\right\|_{([\HerzSo]^{1/s})'}.
\end{equation*}
Combining this and the
above claim, we conclude that
all the assumptions of Lemma \ref{Atoll1}
hold true for $\HerzSo$.
Therefore,
for any sequence $\{a_j\}_{
j\in\mathbb{N}}$ of $(\HerzSo,\,r,\,d)$-atoms
supported, respectively, in the balls
$\{B_j\}_{j\in\mathbb{N}}\subset\mathbb{B}$,
and sequence $\{\lambda_j\}_{j\in\mathbb{N}}
\subset[0,\infty)$ such that $f:=
\sum_{j\in\mathbb{N}}\lambda_ja_j$
in $\mathcal{S'}(\rrn)$, and
$$
\left\|\left\{\sum_{j\in\mathbb
N}\left[\frac{\lambda_j}{\|\1bf_{B_j}\|_{
\HerzSo}}\right]^s\1bf_{B_j}
\right\}^{\frac{1}{s}}\right\|_{\HerzSo}<\infty,
$$
it holds true that $f\in\HaSaHo$ and
$$
\|f\|_{\HaSaHo}
\lesssim\left\|\left\{\sum_{j\in\mathbb
N}\left[\frac{\lambda_j}{\|\1bf_{B_j}\|_{
\HerzSo}}\right]^s\1bf_{B_j}
\right\}^{\frac{1}{s}}\right\|_{\HerzSo}
$$
with the implicit positive constant
independent of $f$.
From this and Definition \ref{atsl},
we deduce that $\aHaSaHo\subset\HaSaHo$ and,
for any $f\in\aHaSaHo$,
\begin{equation}\label{Atole2}
\|f\|_{\HaSaHo}\lesssim\|f\|_{\aHaSaHo},
\end{equation}
where the implicit
positive constant is independent of $f$.

We now prove that $\HaSaHo\subset\aHaSaHo$.
For this purpose,
let $f\in\HaSaHo$.
Notice that the local generalized
Herz space $\HerzSo$ satisfies
Assumptions \ref{assfs} for
the above $\theta$ and $s$. By this, the
assumption $d\geq\lfloor
n(1/s-1)\rfloor=\lfloor n(1/\theta-1)\rfloor$,
and Lemma \ref{Atogl4},
we find that there exist
$\{\lambda_j\}_{j\in\mathbb{N}}
\subset[0,\infty)$ and a sequence
$\{a_j\}_{j\in\mathbb{N}}$ of
$(\HerzSo,\,\infty,\,d)$-atoms supported,
respectively, in the balls
$\{B_j\}_{j\in\mathbb{N}}\subset\mathbb{B}$
such that
\begin{equation}\label{Atole3}
f=\sum_{j\in\mathbb{N}}\lambda_ja_j
\end{equation}
in $\mathcal{S}'(\rrn)$, and
\begin{equation}\label{Atole4}
\left\|\left\{
\sum_{j\in\mathbb{N}}\left[
\frac{\lambda_j}{\|\1bf_{B_j}\|_{\HerzSo}}
\right]^s\1bf_{B_j}
\right\}^{\frac{1}{s}}\right\|_{\HerzSo}
\lesssim\|f\|_{\HaSaHo},
\end{equation}
where the implicit positive constant
is independent of $f$.
In addition, using Lemma \ref{Atoll2}
with $X$, $t$, and $a$ therein replaced,
respectively, by $\HerzSo$, $\infty$,
and $a_j$,
we conclude that, for any $j\in\mathbb{N}$,
$a_j$ is a $(\HerzSo,\,r,\,d)$-atom
supported in the ball $B_j$.
From this, \eqref{Atole3},
\eqref{Atole4}, and Definition \ref{atsl},
it follows that $f\in\aHaSaHo$ and
\begin{align}\label{Atole5}
&\|f\|_{\aHaSaHo}\notag\\&\quad\leq
\left\|\left\{
\sum_{j\in\mathbb{N}}\left[
\frac{\lambda_j}{\|\1bf_{B_j}\|_{\HerzSo}}
\right]^s\1bf_{B_j}
\right\}^{\frac{1}{s}}\right\|_{\HerzSo}
\lesssim\|f\|_{\HaSaHo},
\end{align}
where the implicit positive constant
is independent of $f$.
Thus, we have
$\HaSaH\subset\aHaSaH$.
This further implies that
$$
\HaSaHo=\aHaSaHo.
$$
Moreover, combining \eqref{Atole2}
and \eqref{Atole5}, we conclude that,
for any $f\in\HaSaHo$,
$$
\|f\|_{\HaSaHo}\sim\|f\|_{
\aHaSaHo}
$$
with the positive equivalence
constants independent of $f$,
which completes the proof
of Theorem \ref{Atol}.
\end{proof}

As an application, we next establish the
atomic characterization of the generalized
Morrey--Hardy space $H\MorrSo$. We first
introduce the $(\MorrSo,\,r,\,d)$-atoms as follows.

\begin{definition}\label{atomm}
Let $p$, $q\in[1,\infty)$, $\omega\in M(\rp)$ with
$\MI(\omega)\in(-\infty,0)$ and
$$
-\frac{n}{p}<\m0(\omega)\leq\M0(\omega)<0,
$$
$r\in[1,\infty]$, and $d\in\zp$.
A measurable function $a$ is called an
\emph{$(\MorrSo,\,r,\,d)$-atom}\index{$(\MorrSo,\,r,\,d)$-atom}
if there exists a ball $B\in\mathbb{B}$ such that
\begin{enumerate}
  \item[{\rm(i)}] $\supp(a):=
  \{x\in\rrn:\ a(x)\neq0\}\subset B$;
  \item[{\rm(ii)}] $\|a\|_{L^r(\rrn)}
  \leq\frac{|B|^{1/r}}{\|\1bf_{B}\|_{\MorrSo}}$;
  \item[{\rm(iii)}] for any $\alpha\in\zp^n$
  with $|\alpha|\leq d$,
  $$\int_{\rrn}a(x)x^\alpha\,dx=0.$$
\end{enumerate}
\end{definition}

Then, using Theorem \ref{Atol} and
Remarks \ref{remhs}(iv) and
\ref{remark4.10}(ii),
we immediately obtain the following
atomic
characterization\index{atomic characterization}
of the generalized Morrey--Hardy space
$H\MorrSo$; we omit the details.

\begin{corollary}\label{cor438}
Let $p$, $q\in[1,\infty)$, $\omega\in M(\rp)$ with
$$
-\frac{n}{p}<\m0(\omega)\leq\M0(\omega)<0
$$
and
$$
-\frac{n}{p}<\mi(\omega)\leq\MI(\omega)<0,
$$
$s\in(0,1)$, $d\geq\lfloor n(1/s-1)\rfloor$
be a fixed integer, and
$r\in(\frac{n}{\mw+n/p},\infty]$.
Then the \emph{generalized atomic
Morrey--Hardy space}\index{generalized
atomic Morrey--Hardy space}
$H\aMorrSo$\index{$H\aMorrSo$},
associated with the local generalized
Morrey space $\MorrSo$, is
defined to be the set of all the
$f\in\mathcal{S}'(\rrn)$ such that there exist
a sequence $\{a_{j}\}_{j\in\mathbb{N}}$ of
$(\MorrSo,\,r,\,d)$-atoms supported,
respectively, in the balls
$\{B_{j}\}_{j\in\mathbb{N}}\subset\mathbb{B}$,
and a sequence $\{\lambda_{j}\}_{j\in\mathbb{N}}
\subset[0,\infty)$ such that
$$
f=\sum_{j\in\mathbb{N}}\lambda_{j}a_{j}
$$
in $\mathcal{S}'(\rrn)$, and
$$
\left\|\left\{\sum_{j\in\mathbb{N}}
\left[\frac{\lambda_{j}}
{\|\1bf_{B_{j}}\|_{\MorrSo}}\right]^
{s}\1bf_{B_{j}}\right\}^{\frac{1}{s}}
\right\|_{\MorrSo}<\infty.
$$
Moreover, for any $f\in H\aMorrSo$,
$$
\|f\|_{H\aMorrSo}:=\inf\left\{
\left\|\left\{\sum_{j\in\mathbb{N}}
\left[\frac{\lambda_{j}}
{\|\1bf_{B_{j}}\|_{\MorrSo}}\right]^
{s}\1bf_{B_{j}}\right\}^{\frac{1}{s}}
\right\|_{\MorrSo}\right\},
$$
where the infimum is taken over
all the decompositions of $f$ as above.
Then $$H\MorrSo=H\aMorrSo$$ with
equivalent quasi-norms.
\end{corollary}

The remainder of this section
is devoted to establishing the
atomic characterization of
the generalized Herz--Hardy space
$\HaSaH$. To this end,
we first introduce the following
definition of $(\HerzS,
\,r,\,d)$-atoms.

\begin{definition}\label{deatomg}
Let $p$, $q\in(0,\infty)$, $\omega\in
M(\rp)$ with $\m0(\omega)
\in(\frac{n}{p},\infty)$
and $\MI(\omega)\in(-\infty,0)$,
$r\in[1,\infty]$, and $d\in\zp$.
Then a measurable function $a$ on $\rrn$ is called a
\emph{$(\HerzS,
\,r,\,d)$-atom}\index{$(\HerzS,\,r,\,d)$-atom}
if there exists a ball $B\in\mathbb{B}$ such that
\begin{enumerate}
  \item[{\rm(i)}] $\supp(a):=
  \{x\in\rrn:\ a(x)\neq0\}\subset B$;
  \item[{\rm(ii)}] $\|a\|_{L^r(\rrn)}
  \leq\frac{|B|^{1/r}}{\|\1bf_{B}\|_{\HerzS}}$;
  \item[{\rm(iii)}] for any $\alpha\in\zp^n$
  with $|\alpha|\leq d$,
  $$\int_{\rrn}a(x)x^\alpha\,dx=0.$$
\end{enumerate}
\end{definition}

Then we introduce the following
concept of the generalized atomic
Herz--Hardy space $\aHaSaH$.

\begin{definition}\label{atsg}
Let $p$, $q\in(0,\infty)$,
$\omega\in M(\rp)$ with
$\m0(\omega)\in(-\frac{n}{p},\infty)$ and
$$-\frac{n}{p}<\mi(\omega)\leq\MI(\omega)<0,$$
$r\in(\max\{1,p,\frac{n}{\mw+n/p}\},\infty]$,
$s\in(0,\min\{1,p,q,\frac{n}{\Mw+n/p}\})$,
and $d\geq\lfloor n(
1/s-1)\rfloor$ be a fixed integer.
Then the \emph{generalized atomic
Herz--Hardy space}\index{generalized
atomic Herz--Hardy \\space} $\aHaSaH$\index{$\aHaSaH$}
is defined to be the set of all the
$f\in\mathcal{S}'(\rrn)$ such that there exist
$\{\lambda_j\}_{j\in\mathbb{N}}
\subset[0,\infty)$ and
a sequence $\{a_{j}\}_{j\in\mathbb{N}}$ of
$(\HerzS,\,r,\,d)$-atoms supported,
respectively, in the balls $\{B_{j}\}_{j\in\mathbb{N}}
\subset\mathbb{B}$ such that
$$
f=\sum_{j\in\mathbb{N}}\lambda_{j}a_{j}
$$
in $\mathcal{S}'(\rrn)$,
and
$$
\left\|\left\{
\sum_{j\in\mathbb{N}}
\left[\frac{\lambda_j}{\|
\1bf_{B_j}\|_{\HerzS}}\right]^s
\1bf_{B_j}\right\}^
{\frac{1}{s}}\right\|_{\HerzS}<\infty.
$$
Moreover, for any $f\in\aHaSaH$,
$$
\|f\|_{\aHaSaH}:=\inf\left\{
\left\|\left\{\sum_{j\in\mathbb{N}}
\left[\frac{\lambda_{j}}
{\|\1bf_{B_{j}}\|_{\HerzS}}\right]^{s}
\1bf_{B_{j}}\right\}^{\frac{1}{s}}
\right\|_{\HerzS}\right\},
$$
where the infimum is taken over all
the decompositions of $f$ as above.
\end{definition}

We now state the atomic
characterization\index{atomic characterization}
of the generalized Herz--Hardy space $\HaSaH$ as follows.

\begin{theorem}\label{Atog}
Let $p$, $q$, $\omega$, $r$, $s$, and $d$
be as in Definition \ref{atsg}.
Then $$\HaSaH=\aHaSaH$$ with
equivalent quasi-norms.
\end{theorem}

To prove this theorem,
observe that the associate spaces
of global generalized Herz spaces
are still unknown. Thus, we can not
show Theorem \ref{Atog} directly via using
the known atomic characterization
of Hardy spaces associated with
ball quasi-Banach function spaces
(see Lemmas \ref{Atoll1} and \ref{Atogl4}
above). In order to overcome
this difficulty,
we first establish an improved atomic
characterization of Hardy spaces
associated with ball quasi-Banach
function spaces under no assumption
about associate spaces.
To achieve this,
we first present
the following definition
of atomic Hardy spaces introduced
in \cite{SHYY} (see
also \cite[Definition 3.2]{HCY})

\begin{definition}\label{atsgx}
Let $X$ be a ball quasi-Banach
function space, $r\in(1,\infty]$,
$$0<\theta<s\leq1,$$ and $d\geq
\lfloor n(1/\theta-1)\rfloor$ be
a fixed integer. Then
the \emph{atomic Hardy space}
$H^{X,r,d,s}(\rrn)$\index{$H^{X,r,d,s}(\rrn)$},
associated with $X$, is defined
to be the set of all the $f\in\mathcal{S}'(\rrn)$
such that there exist
$\{\lambda_{j}\}_{j\in\mathbb{N}}\subset
[0,\infty)$ and a sequence
$\{a_j\}_{j\in\mathbb{N}}$ of
$(X,\,r,\,d)$-atoms supported,
respectively, in the balls
$\{B_j\}_{j\in\mathbb{N}}\subset
\mathbb{B}$ such that
$$f=\sum_{j\in\mathbb{N}}\lambda_j
a_j$$
in $\mathcal{S}'(\rrn)$, and
$$
\left\|
\left[\sum_{j\in\mathbb{N}}
\left(\frac{\lambda_j}{
\|\1bf_{B_j}\|_X}\right)^s\1bf_{B_j}
\right]^{\frac{1}{s}}
\right\|_{X}<\infty.
$$
Moreover, for any
$f\in H^{X,r,d,s}(\rrn)$,
\begin{equation}\label{atom hardyx}
\|f\|_{H^{X,r,d,s}(\rrn)}:=
\inf\left\{\left\|
\left[\sum_{j\in\mathbb{N}}
\left(\frac{\lambda_j}{
\|\1bf_{B_j}\|_X}\right)^s\1bf_{B_j}
\right]^{\frac{1}{s}}
\right\|_{X}\right\},
\end{equation}
where the infimum is taken over all
the decompositions of $f$ as above.
\end{definition}

Then we have the following atomic
characterization of the Hardy space $H_X(\rrn)$
associated with the ball quasi-Banach
function space $X$, which is valid
even when the associate space
of $X$ is not clear.

\begin{theorem}\label{Atogx}
Let $X$ be a ball quasi-Banach function space
satisfying:
\begin{enumerate}
  \item[{\rm(i)}] there exist $0<\theta<
  s\leq1$ such that
  Assumption \ref{assfs} holds true;
  \item[{\rm(ii)}] for $s$ in
  (i), there exists a linear space
  $Y\subset\Msc(\rrn)$
  equipped with a seminorm $\|\cdot\|_{Y}$
  such that, for any $f\in\Msc(\rrn)$,
  $$
  \|f\|_{X^{1/s}}
  \sim\sup\left\{\|fg\|_{L^1(\rrn)}:\
  \|g\|_{Y}=1\right\},
  $$
  where the positive equivalence
  constants are independent of $f$;
  \item[{\rm(iii)}] for $s$ in (i) and
  $Y$ in (ii), there exist an $r\in(1,\infty]$
  and a positive constant $C$ such that,
  for any $f\in L^1_{\mathrm{loc}}(\rrn)$,
  $$
  \left\|\mc^{((r/s)')}
  (f)\right\|_{Y}\leq C\|f\|_{Y}.
  $$
\end{enumerate}
Then
$$
H_X(\rrn)=H^{X,r,d,s}(\rrn)
$$
with equivalent quasi-norms.
\end{theorem}

\begin{remark}
We should point out that
Theorem \ref{Atogx} is an improved
version of the known atomic characterization
established by Sawano et al.\ in \cite[Theorems
3.6 and 3.7]{SHYY}. Indeed,
if $Y\equiv(X^{1/s})'$ in Theorem
\ref{Atogx}, then this theorem goes back
to \cite[Theorems 3.6 and 3.7]{SHYY}.
\end{remark}

To show this atomic characterization
of $H_X(\rrn)$, we
need some preliminary lemmas.
First, from the assumption (ii) of
Theorem \ref{Atogx},
the following H\"{o}lder inequality
about $X$ can be deduced directly,
which is useful in the proof of
Theorem \ref{Atogx};
we omit the details.

\begin{lemma}\label{Atogxl2}
Let $X$, $s$, and $Y$ be as in Theorem
\ref{Atogx}. Then there exists a positive
constant $C$ such that, for any
$f,\ g\in\Msc(\rrn)$,
$$
\left\|fg\right\|_{L^1(\rrn)}
\leq C\|f\|_{X^{1/s}}\|g\|_{Y}.
$$
\end{lemma}

To prove Theorem \ref{Atogx},
we also need the following
estimate about convolution,
which is just
\cite[Corollary 2.1.12]{LGCF}.

\begin{lemma}\label{Atogxl3}
Let $\phi\in L^1(\rrn)$ and
$\Phi$ be a nonnegative radial
decreasing function on $\rrn$ such that
$|\phi|\leq\Phi$.
If $\Phi\in L^1(\rrn)$,
then, for any $f\in L^1_{\mathrm{loc}}(\rrn)$
and $x\in\rrn$,
$$
\sup_{t\in(0,\infty)}
\left|f\ast\phi_t(x)\right|
\leq\left\|\Phi\right\|_{L^1(\rrn)}
\mc(f)(x),
$$
where $\mc$ is the Hardy--Littlewood
maximal operator defined as in
\eqref{hlmax}.
\end{lemma}

Via the above estimate, we
now establish a technical
estimate about $(X,r,d)$-atoms as follows,
which is an essential tool in
the proof of Theorem \ref{Atogx}.

\begin{lemma}\label{Atogxl1}
Let $X$ be a ball quasi-Banach function
space, $r\in(1,\infty]$, $\theta\in(0,1]$,
and $d\geq n(1/\theta-1)-1$ be a
nonnegative integer.
Assume that $\phi\in\mathcal{S}(\rrn)$
satisfies $\supp(\phi)\subset B(\0bf,1)$.
Then there exists a positive constant
$C$ such that, for any
$(X,\,r,\,d)$-atom $a$ supported
in the ball $B$,
$$
M(a,\phi)\leq C\left[\mc(a)\1bf_{2B}+
\frac{1}{\|\1bf_B\|_{X}}\mc^{(\theta)}
(\1bf_B)\1bf_{(2B)^{\complement}}\right],
$$
where $M$ denotes the radial
maximal function as in Definition \ref{smax}(i).
\end{lemma}

\begin{proof}
Let all the symbols be as in
the present lemma. Then,
using Lemma \ref{Atogxl3}
with $\Phi:=\|\phi\|_{L^{\infty}(\rrn)}
\1bf_{B(\0bf,1)}$ and $f:=a$,
we find that, for any $x\in\rrn$,
\begin{equation}\label{Atogl3e1}
M(a,\phi)(x)=\sup_{t\in(0,\infty)}
\left|a\ast\phi_t(x)\right|\lesssim
\mc(a)(x),
\end{equation}
where the implicit positive constant is
independent of both $a$ and $x$.

In addition, let $x\in (2B)^{\complement}$
and $t\in(0,\infty)$. Then, from the
assumption that $\supp(\phi)\subset
B(\0bf,1)$, it follows that
$\supp(\phi_t(x-\cdot))\subset
B(x,t)$. Therefore, we have
\begin{equation}\label{Atogl3e2}
\left|a\ast\phi_t(x)\right|
=\left|\int_{B\cap B(x,t)}a(y)
\phi_t(x-y)\,dy\right|=0
\end{equation}
when $B\cap B(x,t)=\emptyset$.
On the other hand,
when $B\cap B(x,t)\neq\emptyset$,
assume $B:=B(x_B,r_B)$ with
$x_B\in\rrn$ and $r_B\in(0,\infty)$.
Then we claim that $|x-x_B|<t+r_B$.
Otherwise, assume that $|x-x_B|\geq t+r_B$.
Applying this, we find that,
for any $y\in B(x,t)$,
$$
|y-x_B|\geq|x-x_B|-|y-x|>r_B,
$$
which implies that $B(x,t)\subset[B(x_B,r_B)]^
{\complement}$ and hence contradicts
to the assumption
$B(x_B,r_B)\cap B(x,t)\neq\emptyset$.
Thus, we obtain $|x-x_B|<t+r_B$
and finish the proof of the
above claim. Combining this and
the assumption $x\in(2B)^{\complement}$,
we further conclude that
\begin{equation*}
t>|x-x_B|-r_B\geq\frac{|x-x_B|}{2}.
\end{equation*}
From this, Definition \ref{atomx}(iii),
and the Taylor remainder
theorem\index{Taylor remainder theorem},
we deduce that, for any $y\in B=B(x_B,r_B)$,
there exists a $t_y\in(0,1)$ such that
\begin{align}\label{Atogl3e4}
&\left|a\ast\phi_t(x)\right|\notag\\
&\quad=\left|\int_{\rrn}a(y)\left[\phi_t(x-y)
-\sum_{\gfz{\gamma\in\zp^n}{|\gamma|\leq d}}
\frac{\partial^{\gamma}
\phi_t(x-x_B)}{\gamma!}(y-x_B)^{\gamma}\right]
\,dy\right|\notag\\
&\quad=\left|\int_{B(x_B,r_B)}a(y)\left[
\sum_{\gfz{\gamma\in\zp^n}
{|\gamma|=d+1}}\frac{\partial^{\gamma}
\phi_t(x-t_yy-(1-t_y)x_B)}
{\gamma!}(y-x_B)^{\gamma}\right]\,dy\right|
\notag\\
&\quad\lesssim t^{-n-d-1}\int_{B(x_B,r_B)}
|a(y)||y-x_B|^{d+1}\,dy\notag\\
&\quad\lesssim\frac{r_B^{d+1}}{|x-x_B|^{n+d+1}}
\|a\|_{L^1(\rrn)}.
\end{align}

Next, we claim that
$$
B\left(t_0x+(1-t_0)x_B,\frac{r_B}{2}\right)
\subset\left[B(x_B,r_B)\cap B(x,
|x-x_B|)\right],
$$
where $t_0:=\frac{r_B}{2|x-x_B|}$.
Indeed, for any $y\in
B(t_0x+(1-t_0)x_B,\frac{r_B}{2})$, we have
\begin{align*}
|y-x_B|\leq\left|y-t_0x-(1-t_0)x_B\right|
+t_0\left|x-x_B\right|<\frac{r_B}{2}
+\frac{r_B}{2}=r_B
\end{align*}
and
\begin{align*}
|y-x|&\leq\left|y-t_0x-(1-t_0)x_B\right|
+(1-t_0)|x-x_B|\\&<
\frac{r_B}{2}+|x-x_B|-\frac{r_B}{2}
=|x-x_B|.
\end{align*}
These imply $y\in B(x_B,r_B)\cap B(x,
|x-x_B|)$ and finish the proof of the above
claim.
By this, \eqref{Atogl3e4},
the H\"{o}lder inequality, Definition
\ref{atomx}(ii), and the
assumptions $x\in(2B)^{\complement}$ and
$d\geq n(1/\theta-1)-1$, we
conclude that
\begin{align}\label{atogxl1e1}
&\left|a\ast\phi_t(x)\right|\notag\\
&\quad\lesssim\frac{r_B^{d+1}}{|x-x_B|^{n+d+1}}
|B(x_B,r_B)|
^{1-\frac{1}{r}}\|a\|_{L^r(\rrn)}\notag\\
&\quad\lesssim\frac{r_B^{d+1}}{|x-x_B|^{n+d+1}}
\frac{|B(x_B,r_B)|}{\|\1bf_{B(x_B,
r_B)}\|_{X}}\notag\\
&\quad\sim\frac{1}{\|\1bf_{B(x_B,r_B)}\|_{X}}
\left(\frac{r_B}{|x-x_B|}\right)^{n+d+1}\notag\\
&\quad\lesssim\frac{1}{\|\1bf_{B
(x_B,r_B)}\|_{X}}
\left(\frac{r_B}{|x-x_B|}
\right)^{\frac{n}{\theta}}\notag\\
&\quad\sim\frac{1}{\|\1bf_{B(x_B,r_B)}\|_{X}}
\left[\frac{|B(t_0x+(1-t_0)x_B,
\frac{r_B}{2})|}{|B(x,|x-x_B|)|}\right]
^{\frac{1}{\theta}}\notag\\
&\quad\lesssim\frac{1}{\|\1bf_{B(x_B,r_B)}
\|_{X}}
\left[\frac{1}{|B(x,|x-x_B|)}
\int_{B(x,|x-x_B|)}\left|
\1bf_{B(x_B,r_B)}(y)\right|^{\theta}\,dy\right]
^{\frac{1}{\theta}}\notag\\
&\quad\lesssim\frac{1}{\|\1bf_B\|}_{X}\mc
^{(\theta)}
(\1bf_B)(x).
\end{align}
This, combined with \eqref{Atogl3e1}
and \eqref{Atogl3e2},
further implies that, for any
$x\in\rrn$,
$$
M(a,\phi)(x)\lesssim
\mc(a)(x)\1bf_{2B}(x)+
\frac{1}{\|\1bf_B\|_{X}}
\mc^{(\theta)}(x)\1bf_{(2B)^{
\complement}}(x)
$$
with the implicit positive constant
independent of both $a$ and $x$,
which completes the proof of
Lemma \ref{Atogxl1}.
\end{proof}

We now show Theorem \ref{Atogx}.

\begin{proof}[Proof of Theorem \ref{Atogx}]
Let all the symbols
be as in the present theorem.
We first show
$$
H^{X,r,d,s}(\rrn)\subset
H_X(\rrn).
$$
To this end, let
$f\in H^{X,r,d,s}(\rrn)$ satisfy
$f=\sum_{j\in\mathbb{N}}\lambda_{j}
a_{j}$ in $\mathcal{S}'(\rrn)$,
and $\phi\in\mathcal{S}(\rrn)$ satisfy
$\supp(\phi)\subset B(\0bf,1)$ and
$$\int_{\rrn}\phi(y)\,dy\neq0,$$
where
$\{\lambda_j\}_{j\in\mathbb{N}}\subset[0,\infty)$
and, for any $j\in\mathbb{N}$, $a_j$ is a
$(X,\,r,\,d)$-atom supported in the ball
$B_j\in\mathbb{B}$.
Then we find that, for any $t\in(0,\infty)$
and $x\in\rrn$,
\begin{align}\label{Atogxe10}
\left|f\ast\phi_t(x)\right|
&=\left|\left\langle
f,\phi_t(x-\cdot)\right\rangle\right|
=\left|\sum_{j\in\mathbb{N}}\lambda_j
\left\langle a_j,\phi_t(x-\cdot)\right\rangle
\right|\notag\\
&=\left|\sum_{j\in\mathbb{N}}\lambda_j
a_j\ast\phi_t(x)\right|
\leq\sum_{j\in\mathbb{N}}
\lambda_j\left|a_j\ast\phi_t(x)\right|.
\end{align}
This implies that
\begin{equation}\label{Atogxe1}
M(f,\phi)\leq\sum_{j\in\mathbb{N}}
\lambda_jM(a_j,\phi).
\end{equation}
In addition, by the definition of
$\lfloor n(1/\theta-1)\rfloor$, we conclude that
$$
d>n\left(\frac{1}{\theta}-1\right)-1.
$$
From this, \eqref{Atogxe1},
and Lemma \ref{Atogxl1} with
$a$ therein replaced by
$a_j$,
we deduce that
\begin{align}\label{Atogxe2}
\left\|M(f,\phi)\right\|_{X}
&\leq\left\|\sum_{j\in\mathbb{N}}
\lambda_jM(a_j,\phi)\right\|_{X}\notag\\
&\lesssim\left\|\sum_{j\in\mathbb{N}}
\lambda_j\mc(a_j)\1bf_{2B_j}\right\|_{X}
+\left\|\sum_{j\in\mathbb{N}}
\frac{\lambda_j}{\|\1bf_{B_j}\|_{X}}
\mc^{(\theta)}\left(\1bf_{B_j}\right)
\right\|_{X}\notag\\
&=:\mathrm{II}_1+\mathrm{II}_2.
\end{align}

Then we deal with $\mathrm{II}_1$
and $\mathrm{II}_2$, respectively.
For $\mathrm{II}_1$,
let $g\in\Msc(\rrn)$ satisfy $\|g\|_{Y}=1$.
Then, using the Tonelli theorem,
the H\"{o}lder inequality, the boundedness
of the Hardy--Littlewood maximal operator
$\mc$ on $L^r(\rrn)$, and Definition
\ref{deatomg}(ii), we find that
\begin{align}\label{Atogxe3}
&\int_{\rrn}\sum_{j\in\mathbb{N}}
\lambda_j^s\left[\mc(a_j)(y)\right]^s\1bf_
{2B_j}(y)g(y)\,dy\notag\\
&\quad\leq\sum_{j\in\mathbb{N}}
\lambda_j^s\left\|
\left[\mc(a_j)\right]^s
\right\|_{L^{r/s}(\rrn)}\left\|g
\1bf_{2B_j}\right\|
_{L^{(r/s)'}(\rrn)}\notag\\
&\quad=\sum_{j\in\mathbb{N}}
\lambda_j^s\left\|\mc(a_j)
\right\|_{L^{r}(\rrn)}^s\left\|g
\1bf_{2B_j}\right\|
_{L^{(r/s)'}(\rrn)}\notag\\
&\quad\lesssim\sum_{j\in\mathbb{N}}
\lambda_j^s\left\|a_j\right\|_{L^r(\rrn)}
^s\left\|g
\1bf_{2B_j}\right\|_{L^{(r/s)'}(\rrn)}\notag\\
&\quad\lesssim\sum_{j\in\mathbb{N}}
\lambda_j^s\frac{|B_j|^{\frac{s}{r}}}
{\|\1bf_{B_j}\|_{X}^s}\left\|
g\1bf_{2B_j}\right\|_{
L^{(r/s)'}(\rrn)}.
\end{align}
In addition,
for any $j\in\mathbb{N}$ and
$x\in B_j$, we have
\begin{align*}
\mc^{((r/s)')}(g)(x)
&\geq\left[\frac{1}{|2B_j|}
\int_{2B_j}\left|
g(y)\right|^{(r/s)'}\,dy\right]^{1/(r/s)'}\\
&\sim\left|B_j\right|^{-1/(r/s)'}
\left\|g\1bf_{
2B_j}\right\|_{L^{(r/s)'}(\rrn)}.
\end{align*}
By this, \eqref{Atogxe3},
the Tonelli theorem, Lemma \ref{Atogxl2},
and Definition \ref{atsgx}(iii) together
with the assumption $\|g\|_Y=1$,
we conclude that
\begin{align*}
&\int_{\rrn}\sum_{j\in\mathbb{N}}
\lambda_j^s\left[\mc(a_j)(y)\right]^s\1bf_
{2B_j}(y)g(y)\,dy\\
&\quad\lesssim\sum_{j\in\mathbb{N}}
\frac{\lambda_j^s}{\|\1bf_{B_j}\|_{X}^s}
\int_{B_j}\left|B_j\right|^{-1/(r/s)'}
\left\|g\1bf_{2B_j}\right\|
_{L^{(r/s)'}(\rrn)}\,dx\\
&\quad\lesssim\sum_{j\in\mathbb{N}}
\frac{\lambda_j^s}{\|\1bf_{B_j}\|_{X}
^{s}}\int_{B_j}\mc^{((r/s)')}(g)(x)\,dx\\
&\quad\sim\int_{\rrn}
\sum_{j\in\mathbb{N}}\left(\frac{\lambda_j}
{\|\1bf_{B_j}\|_{X}}\right)^s\1bf_{B_j}(x)
\mc^{((r/s)')}(g)(x)\,dx\\
&\quad\lesssim\left\|\sum_{j\in\mathbb{N}}
\left(\frac{\lambda_j}{\|\1bf_
{B_j}\|_{X}}\right)^s\1bf_{B_j}
\right\|_{X^{1/s}}\left\|
\mc^{((r/s)')}(g)\right\|_{Y}\\
&\quad\lesssim\left\|
\left[\sum_{j\in\mathbb{N}}
\left(\frac{\lambda_j}{\|\1bf_
{B_j}\|_{X}}\right)^s\1bf_{B_j}
\right]^{\frac{1}{s}}
\right\|_{X}^s.
\end{align*}
From this and Definition \ref{atsgx}(ii),
it follows that
\begin{align*}
\left\|\sum_{j\in\mathbb{N}}
\lambda_j^s\left[
\mc(a_j)\right]^s\1bf_{2B_j}
\right\|_{X^{1/s}}\lesssim
\left\|\left[\sum_{j\in\mathbb{N}}
\left(\frac{\lambda_j}{\|\1bf_
{B_j}\|_{X}}\right)^s\1bf_{B_j}
\right]^{\frac{1}{s}}
\right\|_{X}^s.
\end{align*}
This, together with Lemma \ref{l4320}
with $r$ and $\{a_j\}_{j\in\mathbb{N}}$
therein replaced, respectively, by
$s$ and $\{\lambda_j\mc(a_j)\1bf_{2B_j}\}
_{j\in\mathbb{N}}$,
further implies that
\begin{align}\label{Atogxe5}
\mathrm{II}_1&\lesssim
\left\|\left\{
\sum_{j\in\mathbb{N}}\lambda_j^s
\left[\mc(a_j)\right]^s\1bf_{2B_j}
\right\}^{\frac{1}{s}}
\right\|_{X}\notag\\
&\sim\left\|\sum_{j\in\mathbb{N}}
\lambda_j^s\left[
\mc(a_j)\right]^s\1bf_{2B_j}
\right\|_{X^{1/s}}^{\frac{1}{s}}\notag\\
&\lesssim\left\|\left[\sum_{j\in\mathbb{N}}
\left(\frac{\lambda_j}{\|\1bf_
{B_j}\|_{X}}\right)^s\1bf_{B_j}
\right]^{\frac{1}{s}}
\right\|_{X},
\end{align}
which completes the estimation
of $\mathrm{II}_1$.

Next, we deal with $\mathrm{II}_2$.
Indeed, applying Definition \ref{atsgx}(i),
we find that Assumption \ref{assfs} holds true
for $X$, namely, for any $\{f_j\}_{
j\in\mathbb{N}}\subset L^1_{\mathrm{loc}}(
\rrn)$,
\begin{equation}\label{Atogxe51}
\left\|\left\{
\sum_{j\in\mathbb{N}}
\left[
\mc^{(\theta)}(f_j)
\right]^{s}\right\}^{1/s}\right\|_{X}
\lesssim\left\|\left(\sum_{j\in\mathbb{N}}
|f_j|^s\right)^{1/s}\right\|_{X}.
\end{equation}
From this and Lemma \ref{l4320} again,
we deduce that
\begin{align}\label{Atogxe6}
\mathrm{II}_2&\lesssim
\left\|\left\{
\sum_{j\in\mathbb{N}}
\left[\frac{\lambda_j}{\|\1bf_
{B_j}\|_{X}}\mc^{(\theta)}(\1bf_
{B_j})\right]^s
\right\}^{\frac{1}{s}}
\right\|_{X}\notag\\
&\sim\left\|
\left\{
\sum_{j\in\mathbb{N}}\left[
\mc^{(\theta)}\left(
\frac{\lambda_j}{\|\1bf_{B_j}\|_{X}}
\1bf_{B_j}\right)
\right]^s
\right\}^{\frac{1}{s}}
\right\|_{X}\notag\\
&\lesssim\left\|\left[
\sum_{j\in\mathbb{N}}\left(\frac{\lambda_j}
{\|\1bf_{B_j}\|_{X}}\right)^s\1bf_
{B_j}
\right]^{\frac{1}{s}}\right\|_{X}.
\end{align}
This finishes the estimation of
$\mathrm{II}_2$.
Combining Lemma \ref{Th5.4l1},
\eqref{Atogxe2}, \eqref{Atogxe5}, \eqref{Atogxe6},
Definition \ref{atsgx}, and
the choice of $\{\lambda_j\}_{j\in\mathbb{N}}$,
we further conclude
that
\begin{equation}\label{Atogxe7}
\|f\|_{H_X(\rrn)}\sim\left\|
M(f,\phi)
\right\|_{X}\lesssim\|f\|_{H^{X,r,d,s}(\rrn)}
<\infty,
\end{equation}
where the implicit positive
constants are independent of $f$.
This implies that
$f\in H_X(\rrn)$ and hence
$H^{X,r,d,s}(\rrn)\subset H_X(\rrn)$.

Conversely, we prove that
$$
H_X(\rrn)\subset H^{X,r,d,s}(\rrn).
$$
To achieve this, let $f\in H_X(\rrn)$.
Notice that $X$ is a BQBF space.
From this,
\eqref{Atogxe51}, the
assumption $d\geq\lfloor n(1/\theta-1)\rfloor$,
and Lemma \ref{Atogl4},
it follows that there exist
$\{\lambda_j\}_{j\in\mathbb{N}}
\subset[0,\infty)$ and a sequence
$\{a_j\}_{j\in\mathbb{N}}$ of
$(X,\,\infty,\,d)$-atoms supported,
respectively, in the balls
$\{B_j\}_{j\in\mathbb{N}}\subset\mathbb{B}$
such that
\begin{equation}\label{Atogxe71}
f=\sum_{j\in\mathbb{N}}\lambda_ja_j
\end{equation}
in $\mathcal{S}'(\rrn)$, and
\begin{equation}\label{Atogxe8}
\left\|\left[
\sum_{j\in\mathbb{N}}\left(
\frac{\lambda_j}{\|\1bf_{B_j}\|_{X}}
\right)^s\1bf_{B_j}
\right]^{\frac{1}{s}}\right\|_{X}
\lesssim\|f\|_{H_X(\rrn)}
\end{equation}
with the implicit positive constant
independent of $f$.
In addition, applying Lemma \ref{Atoll2}
with $t=\infty$ and
$a=a_j$ for any $j\in\mathbb{N}$,
we conclude that, for
any $j\in\mathbb{N}$,
$a_j$ is a $(X,\,r,\,d)$-atom
supported in the ball $B_j$.
This, together with \eqref{Atogxe71},
\eqref{Atogxe8}, and Definition \ref{atsgx},
further implies that
$f\in H^{X,r,d,s}(\rrn)$ and
\begin{align}\label{Atogxe9}
\|f\|_{H^{X,r,d,s}(\rrn)}\leq
\left\|\left[
\sum_{j\in\mathbb{N}}\left(
\frac{\lambda_j}{\|\1bf_{B_j}\|_{X}}
\right)^s\1bf_{B_j}
\right]^{\frac{1}{s}}\right\|_{X}
\lesssim\|f\|_{H_X(\rrn)},
\end{align}
where the implicit positive constant
is independent of $f$.
Therefore, we have $$H_X(\rrn)\subset
H^{X,r,d,s}(\rrn).$$
This further implies that
$
H_X(\rrn)=H^{X,r,d,s}(\rrn).
$
Moreover, by \eqref{Atogxe7}
and \eqref{Atogxe9}, we find that,
for any $f\in H_X(\rrn)$,
$$
\|f\|_{H_X(\rrn)}\sim\|f\|_{
H^{X,r,d,s}(\rrn)}
$$
with the positive equivalence constants
independent of $f$,
which completes the proof of
Theorem \ref{Atogx}.
\end{proof}

Now, in order to
establish the atomic characterization
of the generalized Herz--Hardy
space $\HaSaH$, we only need
to prove that all the assumptions
of Theorem \ref{Atogx} are
satisfied for the Herz space $\HerzS$
under consideration.
To this end,
we first show several
auxiliary lemmas about
block spaces and global generalized
Herz spaces.
The following one shows
the boundedness of
powered Hardy--Littlewood
maximal operators on block spaces.

\begin{lemma}\label{mbhag}
Let $p,\ q\in(0,\infty)$ and
$\omega\in M(\rp)$ satisfy
$\m0(\omega)\in(-\frac{n}{p},\infty)$ and
$\mi(\omega)\in(-\frac{n}{p},\infty)$.
Then, for any given $s\in
(0,\min\{p,q,\frac{n}{\Mw+n/p}\})$
and $$r\in\left(\max\left\{p,\frac{n}{\mw+n/p}
\right\},\infty\right],$$
there exists a positive constant
$C$ such that,
for any $f\in L_{{\rm loc}}^{1}(\rrn)$,
\begin{equation*}
\left\|\mc^{((r/s)')}(f)\right\|_{\dot{\mathcal{B}}
_{1/\omega^{s}}^{(p/s)',(q/s)'}(\rrn)}\leq C
\|f\|_{\dot{\mathcal{B}}_
{1/\omega^{s}}^{(p/s)',(q/s)'}(\rrn)}.
\end{equation*}
\end{lemma}

\begin{proof}
Let all the symbols be as in the present lemma.
Then, from Lemma \ref{rela} and the
assumption $$s\in\left(0,\frac{n}{\Mw+n/p}\right),$$
it follows that
\begin{align}\label{mbhag4}
&\min\left\{
\m0\left(\frac{1}{\omega^s}\right),\
\mi\left(\frac{1}{\omega^s}\right)
\right\}\notag\\&\quad=-s\Mw
>-s\left(\frac{n}{s}-\frac{n}{p}\right)
=-\frac{n}{(p/s)'}.
\end{align}
In addition, using Lemma \ref{rela} and
the assumption
$$r\in\left(\frac{n}{\mw+n/p},\infty\right),$$
we conclude that
\begin{align}\label{mbhag7}
&\max\left\{
\M0\left(\frac{1}{\omega^s}\right),\
\MI\left(\frac{1}{\omega^s}\right)
\right\}\notag\\&
\quad=-s\mw<-s\left(\frac{n}{r}-\frac{n}{p}\right)
=\frac{n}{(r/s)'}-\frac{n}{(p/s)'}.
\end{align}
By the assumptions $r\in(p,\infty)$
and $s\in(0,p)$,
we find that
$$
\left(\frac{r}{s}\right)'\in
\left(1,\left(\frac{p}{s}\right)'\right).
$$
Applying this, \eqref{mbhag4}, \eqref{mbhag7},
and Corollary \ref{maxbl} with
$p$, $q$, $\omega$, and $r$ therein replaced,
respectively, by
$(p/s)'$, $(q/s)'$,
$1/\omega^s$, and $(r/s)'$,
we conclude that,
for any $f\in L^1_{\mathrm{loc}}(\rrn)$,
$$
\left\|\mc^{((r/s)')}(f)\right\|_{\dot{\mathcal{B}}
_{1/\omega^{s}}^{(p/s)',(q/s)'}(\rrn)}\lesssim
\|f\|_{\dot{\mathcal{B}}_
{1/\omega^{s}}^{(p/s)',(q/s)'}(\rrn)},
$$
which completes the proof of Lemma \ref{mbhag}.
\end{proof}

We also need the following two
technical lemmas about
the Fefferman--Stein vector-valued
inequality on global generalized
Herz spaces.

\begin{lemma}\label{vmbhg}
Let $p,\ q\in(0,\infty)$ and $\omega\in M(\rp)$ satisfy
$\m0(\omega)\in(-\frac{n}{p},\infty)$
and $\mi(\omega)\in(-\frac{n}{p},\infty)$.
Then, for any given $u\in(1,\infty)$ and
$$r\in\left(0,\min\left\{p,\frac{n}
{\Mw+n/p}\right\}\right),$$ there exists
a positive constant $C$ such that, for any $\{f
_{j}\}_{j\in\mathbb{N}}\subset L_{{\rm loc}}^{1}(\rrn)$,
\begin{equation*}
\left\|\left\{\sum_{j\in\mathbb{N}}\left[\mc
(f_{j})\right]^{u}\right\}^{\frac{1}{u}}
\right\|_{[\HerzS]^{1/r}}
\leq C\left\|\left(
\sum_{j\in\mathbb{N}}|f_{j}|^{u}\right)^
{\frac{1}{u}}\right\|_{[\HerzS]^{1/r}}.
\end{equation*}
\end{lemma}

\begin{proof}
Let all the symbols be as in the present lemma.
Then, by \eqref{mbhle1} and \eqref{mbhle2},
we find that
\begin{align*}
\Mwr<\frac{n}{(p/r)'}
\end{align*}
and
\begin{align*}
\mwr>-\frac{n}{p/r}.
\end{align*}
This, together with the assumption
$\frac{p}{r}\in(1,\infty)$, Theorem \ref{Th3.3},
and Lemma \ref{convexl},
further implies that,
for any $\{f_{j}\}_{j\in\mathbb{N}}
\subset L_{{\rm loc}}^{1}(\rrn)$,
\begin{equation*}
\left\|\left\{\sum_{j\in\mathbb{N}}\left[\mc
(f_{j})\right]^{u}\right\}^{\frac{1}{u}}
\right\|_{[\HerzS]^{1/r}}
\lesssim\left\|\left(
\sum_{j\in\mathbb{N}}|f_{j}|^{u}\right)
^{\frac{1}{u}}\right\|_{[\HerzS
]^{1/r}},
\end{equation*}
which completes the proof of Lemma \ref{vmbhg}.
\end{proof}

\begin{lemma}\label{Atogl5}
Let $p$, $q\in(0,\infty)$,
$\omega\in M(\rp)$ satisfy
$\m0(\omega)\in(-\frac{n}{p},\infty)$
and $\mi(\omega)\in(-\frac{n}{p},
\infty)$, $s\in(0,\infty)$,
and $$
\theta\in\left(0,\min\left\{
s,p,\frac{n}{\Mw+n/p}\right\}\right).
$$
Then there exists a positive
constant $C$ such that, for any
$\{f_j\}_{j\in\mathbb{N}}\subset
L^1_{\mathrm{loc}}(\rrn)$,
$$
\left\|
\left\{
\sum_{j\in\mathbb{N}}
\left[\mc^{(\theta)}(f_j)\right]^{s}
\right\}^{1/s}
\right\|_{\HerzS}
\leq C\left\|
\left(\sum_{j\in\mathbb{N}}
|f_j|^{s}
\right)^{1/s}
\right\|_{\HerzS}.
$$
\end{lemma}

\begin{proof}
Let all the symbols be
as in the present lemma.
Then, from the assumption $\theta\in(0,s)$,
we deduce that
$\frac{s}{\theta}\in(1,\infty)$.
Combining this and Lemma \ref{vmbhg}
with $u$ and $r$ therein replaced,
respectively, by $\frac{s}{\theta}$
and $\frac{1}{\theta}$,
we conclude that,
for any $\{f_{j}\}_{
j\subset\mathbb{N}}
\subset L^1_{\mathrm{loc}}(\rrn)$,
$$
\left\|
\left\{
\sum_{j\in\mathbb{N}}
\left[\mc(f_j)\right]^{s/\theta}
\right\}^{\theta/s}
\right\|_{[\HerzS]^{1/\theta}}
\lesssim\left\|
\left(\sum_{j\in\mathbb{N}}
|f_j|^{s/\theta}
\right)^{\theta/s}
\right\|_{
[\HerzS]^{1/\theta}
}.
$$
By this and Remark \ref{power}(ii)
with $X=\HerzSo$, we find that,
for any $\{f_j\}_{j\in\mathbb{N}}
\subset L^1_{\mathrm{loc}}(\rrn)$,
$$
\left\|
\left\{
\sum_{j\in\mathbb{N}}
\left[\mc^{(\theta)}(f_j)\right]^{s}
\right\}^{1/s}
\right\|_{\HerzS}
\lesssim\left\|
\left(\sum_{j\in\mathbb{N}}
|f_j|^{s}
\right)^{1/s}
\right\|_{\HerzS}.
$$
This finishes the proof of
Lemma \ref{Atogl5}.
\end{proof}

Moreover, we require the
following variant of
the H\"{o}lder inequality
of global generalized Herz spaces,
which can be deduced
directly from \eqref{ope};
we omit the details.

\begin{lemma}\label{Atogl1}
Let $p$, $q\in(1,\infty)$ and $\omega\in M(\rp)$.
Then there exists a positive constant $C$ such that,
for any $f$, $g\in\Msc(\rrn)$,
$$
\|fg\|_{L^1(\rrn)}\leq C\|f\|_{\HerzS}
\|g\|_{\Bspace}.
$$
\end{lemma}

We point out that
the following equivalent characterization
of the global generalized Herz space
$\HerzS$ is an essential tool
used to overcome the deficiency of
the associate space of $\HerzS$,
which plays a key role
throughout this book.

\begin{lemma}\label{Atogl2}
Let $p$, $q\in(1,\infty)$ and $\omega
\in M(\rp)$ satisfy $\m0(\omega)\in(-\frac{n}{p},
\infty)$ and $\MI(\omega)\in(-\infty,0)$. Then
a measurable function $f$ belongs to
the global generalized Herz space $\HerzS$ if and
only if
$$
\|f\|_{\HerzS}^{\star}:=\sup
\left\{\|fg\|_{L^1(\rrn)}:\
\|g\|_{\Bspace}=1\right\}<\infty.
$$
Moreover, there exists a constant
$C\in[1,\infty)$ such that, for any $f\in\HerzS$,
$$
C^{-1}\|f\|_{\HerzS}\leq
\|f\|_{\HerzS}^{\star}\leq C\|f\|_{\HerzS}.
$$
\end{lemma}

\begin{proof}
Let all the symbols be as in the present lemma.
We first show the necessity. Indeed,
let $f\in\HerzS$. Then, from Lemma \ref{Atogl1},
we deduce that, for any $g\in\Msc(\rrn)$
with $\|g\|_{\Bspace}=1$,
$$
\|fg\|_{L^1(\rrn)}\lesssim\|f\|_{\HerzS}\|g\|
_{\Bspace}\sim\|f\|_{\HerzS},
$$
which further implies that
\begin{equation}\label{Atogl2e1}
\|f\|_{\HerzS}^{\star}\lesssim\|f\|_{\HerzS}<\infty.
\end{equation}
Thus, we complete the proof of
the necessity.

Conversely, we show the sufficiency. To this end,
assume $f\in\Msc(\rrn)$ satisfying $\|f\|_
{\HerzS}^{\star}<\infty$, and, for any
$k\in\mathbb{N}$, let
$$
f_{k}:=\min\{|f|,k\}\1bf_{B(\0bf,\tk)}.
$$
Obviously,
for any $k\in\mathbb{N}$, $|f_k|\leq
k\1bf_{B(\0bf,\tk)}$ and $0\leq f_k\uparrow
|f|$ as $k\to\infty$.
In addition, applying Theorem \ref{Th2},
we conclude that the global
generalized Herz space $\HerzS$ is a BQBF space.
Therefore,
for any $k\in\mathbb{N}$,
we have $f_k\in\HerzS$.
By this and Theorem \ref{pre} with
$\bspace$ replaced by $\Bspace$, we find that,
for any $k\in\mathbb{N}$,
\begin{align}\label{Atogl2e2}
\|f_k\|_{\HerzS}&\sim
\sup\left\{\left|\int_{\rrn}f_k(y)g(y)\,dy\right|
:\ \|g\|_{\Bspace}=1\right\}\notag\\
&\lesssim\sup\left\{\|fg\|_{L^1(\rrn)}
:\ \|g\|_{\Bspace}=1\right\}\sim\|f\|_{\HerzS}^{\star}.
\end{align}
On the other hand, from Definition \ref{Df1}(iii)
and the fact that $0\leq f_k\uparrow
|f|$ as $k\to\infty$, it follows that
$$
\lim\limits_{k\to\infty}
\|f_k\|_{\HerzS}=\|f\|_{\HerzS}.
$$
Combining this and \eqref{Atogl2e2},
we further find that
\begin{equation}\label{Atogl2e3}
\|f\|_{\HerzS}\lesssim\|f\|_{\HerzS}
^{\star}<\infty.
\end{equation}
Then we conclude that $f\in\HerzS$ and
hence complete the proof of the sufficiency.
Moreover, by both \eqref{Atogl2e1} and
\eqref{Atogl2e3}, we find that
$$
\|f\|_{\HerzS}\sim
\|f\|_{\HerzS}^{\star},
$$
where the positive equivalence constants
are independent of $f$.
This finishes the proof of Lemma \ref{Atogl2}.
\end{proof}

Via above lemmas,
we now show Theorem \ref{Atog}.

\begin{proof}[Proof of Theorem \ref{Atog}]
Let $p$, $q$, $\omega$, $r$, $s$, and $d$
be as in the present theorem. Then, applying the
assumptions $\m0(\omega)\in(-\frac{n}{p},
\infty)$ and $\MI(\omega)\in(-\infty,0)$,
and Theorem \ref{Th2}, we find that
the global generalized Herz space $\HerzS$
under consideration is a BQBF space.
Thus, by Theorem \ref{Atogx}, we conclude that,
to complete the proof of the present theorem,
it suffices to show that
the assumptions (i), (ii), and (iii)
of Theorem \ref{Atogx} are satisfied for
$\HerzS$.

Indeed, from the definition of
$\lfloor n(1/s-1)\rfloor$, we deduce that
$$
\left\lfloor
n\left(\frac{1}{s}-1\right)
\right\rfloor\leq n\left(
\frac{1}{s}-1
\right)<\left\lfloor
n\left(\frac{1}{s}-1\right)
\right\rfloor+1.
$$
Then we can choose a $\theta\in(0,s)$ such that
$$
\left\lfloor
n\left(\frac{1}{s}-1\right)
\right\rfloor\leq n\left(
\frac{1}{\theta}-1
\right)<\left\lfloor
n\left(\frac{1}{s}-1\right)
\right\rfloor+1,
$$
which further implies that
$d\geq \lfloor n(1/s-1)\rfloor
=\lfloor n(1/\theta-1)\rfloor$.
Now, we show that
Theorem \ref{Atogx}(i) holds true for
the above $\theta$ and $s$.
Namely, for the above $\theta$ and
$s$, $\HerzS$ satisfies Assumption
\ref{assfs}. Indeed,
applying Lemma \ref{Atogl5},
we find that, for any $\{f_j\}_{
j\in\mathbb{N}}\subset L^1_{\mathrm{loc}}(
\rrn)$,
\begin{equation}\label{atoge1}
\left\|\left\{
\sum_{j\in\mathbb{N}}
\left[
\mc^{(\theta)}(f_j)
\right]^{s}\right\}^{1/s}\right\|_{\HerzS}
\lesssim\left\|\left(\sum_{j\in\mathbb{N}}
|f_j|^s\right)^{1/s}\right\|_{\HerzS}.
\end{equation}
This implies that Assumption
\ref{assfs} is satisfied for $\HerzS$ and hence
finishes the proof that Theorem
\ref{Atogx}(i) holds true.

Next, we prove that $\HerzS$ satisfies Theorem
\ref{Atogx}(ii). Indeed, from
the assumptions $\m0(\omega)\in(-\frac{n}{p},\infty)$
and $\MI(\omega)\in(-\infty,0)$, and
Lemma \ref{rela}, it follows that
$$
\m0\left(\omega^s\right)
=s\m0(\omega)>-\frac{n}{p/s}
$$
and
$$
\MI\left(\omega^s\right)=s
\MI(\omega)<0.
$$
These, combined with the assumptions
$\frac{p}{s},\ \frac{q}{s}\in(1,\infty)$
and Lemma \ref{Atogl2}
with $p$, $q$, and $\omega$
therein replaced, respectively,
by $\frac{p}{s}$, $\frac{q}{s}$,
and $\omega^s$, further imply that,
for any $f\in\Msc(\rrn)$,
\begin{equation}\label{atoge2}
\|f\|_{\HerzScs}\sim
\sup\left\{\|fg\|_{L^1(\rrn)}:\
\|g\|_{\dot{\mathcal{B}}_{1/\omega^s}
^{(p/s)',(q/s)'}(\rrn)}=1\right\}
\end{equation}
with the positive equivalence constants
independent of $f$.
On the other hand, by Theorem \ref{convexl},
we find that $\HerzScs=[\HerzS]^{1/s}$.
Combining this and \eqref{atoge2}, we
conclude that, for any $f\in\Msc(\rrn)$,
\begin{equation}\label{atoge3}
\|f\|_{[\HerzS]^{1/s}}\sim
\sup\left\{\|fg\|_{L^1(\rrn)}:\
\|g\|_{\dot{\mathcal{B}}_{1/\omega^s}
^{(p/s)',(q/s)'}(\rrn)}=1\right\}.
\end{equation}
Therefore, $\HerzS$ satisfies Theorem
\ref{Atogx}(ii) with $Y:=
\dot{\mathcal{B}}_{1/\omega^s}
^{(p/s)',(q/s)'}(\rrn)$.

Finally, we show that Theorem
\ref{Atogx}(iii) is satisfied for $\HerzS$.
Indeed, from Lemma \ref{mbhag},
it follows that, for any
$f\in L_{{\rm loc}}^{1}(\rrn)$,
\begin{equation}\label{atoge4}
\left\|\mc^{((r/s)')}(f)\right\|_{\dot{\mathcal{B}}
_{1/\omega^{s}}^{(p/s)',(q/s)'}(\rrn)}\lesssim
\|f\|_{\dot{\mathcal{B}}_
{1/\omega^{s}}^{(p/s)',(q/s)'}(\rrn)}.
\end{equation}
This implies that Theorem \ref{Atogx}(iii)
holds true for $\HerzS$ with
$$Y:=\dot{\mathcal{B}}_{1/\omega^s}
^{(p/s)',(q/s)'}(\rrn).$$
Moreover, combining \eqref{atoge1},
\eqref{atoge3}, and \eqref{atoge4},
we further find that the Herz space
$\HerzS$ under consideration satisfies
all the assumptions of Theorem \ref{Atogx}
and then complete the proof of Theorem
\ref{Atog}.
\end{proof}

As an application of this
atomic characterization,
we give the
atomic characterization of the
generalized Morrey--Hardy space
$H\MorrS$.
First, we introduce the
the following $(\MorrS,\,r,
\,d)$-atoms.

\begin{definition}\label{deatomgm}
Let $p$, $q$, and $\omega$ be as in
Definition \ref{atomm},
$r\in[1,\infty]$, and $d\in\zp$.
Then a measurable function $a$
on $\rrn$ is called an
\emph{$(\MorrS,\,r,
\,d)$-atom}\index{$(\MorrS,\,r,\,d)$-atom}
if there exists a ball $B\in\mathbb{B}$
such that
\begin{enumerate}
  \item[{\rm(i)}] $\supp(a):=
  \{x\in\rrn:\ a(x)\neq0\}\subset B$;
  \item[{\rm(ii)}] $\|a\|_{L^r(\rrn)}
  \leq\frac{|B|^{1/r}}{\|\1bf_{B}\|_{\MorrS}}$;
  \item[{\rm(iii)}] for any $\alpha\in\zp^n$
  with $|\alpha|\leq d$,
  $$\int_{\rrn}a(x)x^\alpha\,dx=0.$$
\end{enumerate}
\end{definition}

By Theorem \ref{Atog} and Remarks \ref{remhs}(iv)
and \ref{remark4.10}(ii), we immediately obtain the
atomic characterization\index{atomic characterization}
of the generalized Morrey--Hardy space
$H\MorrS$ as follows; we omit the details.

\begin{corollary}
Let $p$, $q$, $\omega$,
$r$, $d$, and $s$ be as in Corollary
\ref{cor438}. The \emph{generalized atomic
Morrey--Hardy space}\index{generalized atomic Morrey--Hardy space}
$H\aMorrS$\index{$H\aMorrS$},
associated with the global generalized
Morrey space $\MorrS$, is defined to be
the set of all the
$f\in\mathcal{S}'(\rrn)$ such that there exist
a sequence $\{a_{j}\}_{j\in\mathbb{N}}$ of
$(\MorrS,\,r,\,d)$-atoms supported,
respectively, in the balls
$\{B_{j}\}_{j\in\mathbb{N}}\subset\mathbb{B}$,
and a sequence
$\{\lambda_{j}\}_{j\in\mathbb{N}}
\subset[0,\infty)$ such that
$$
f=\sum_{j\in\mathbb{N}}\lambda_{j}a_{j}
$$
in $\mathcal{S}'(\rrn)$,
and
$$
\left\|\left\{\sum_{j\in\mathbb{N}}
\left[\frac{\lambda_{j}}
{\|\1bf_{B_{j}}\|_{\MorrS}}
\right]^{s}\1bf_{B_{j}}\right\}^{\frac{1}{s}}
\right\|_{\MorrS}<\infty.
$$
Moreover,
for any $f\in H\aMorrS$,
$$
\|f\|_{H\aMorrS}:=\inf\left\{
\left\|\left\{\sum_{j\in\mathbb{N}}\left[\frac{\lambda_{j}}
{\|\1bf_{B_{j}}\|_{\MorrS}}\right]^{s}\1bf_{B_{j}}
\right\}^{\frac{1}{s}}
\right\|_{\MorrS}\right\},
$$
where the infimum is taken over all
the decompositions of $f$ as above.
Then $$H\MorrS=H\aMorrS$$ with equivalent quasi-norms.
\end{corollary}

\section{Finite Atomic Characterizations}

In this section, we investigate the
finite atomic characterizations
of generalized Herz--Hardy spaces.
We first consider the finite
atomic characterization of
$\HaSaHo$ via introducing the
following finite atomic Hardy space
associated with $\HerzSo$.

\begin{definition}\label{finiat}
Let $p$, $q\in(0,\infty)$, $\omega\in M(\rp)$ with
$\m0(\omega)\in(-\frac{n}{p},\infty)$ and
$\mi(\omega)\in(-\frac{n}{p},\infty)$,
$r\in(\max\{1,p,\frac{n}{\mw+n/p}\},\infty]$,
$$s\in\left(0,\min\left\{
1,p,q,\frac{n}{\Mw+n/p}\right\}\right),$$
and $d\geq\lfloor n(
1/s-1)\rfloor$ be a fixed integer.
Then the \emph{generalized finite atomic
Herz--Hardy space\index{generalized
finite atomic Herz--Hardy \\ space}
$\iaHaSaHo$\index{$\iaHaSaHo$}},
associated with $\HerzSo$, is defined
to be the set of all
finite linear combinations of
$(\HerzSo,\,r,\,d)$-atoms.
Moreover,
for any $f\in\iaHaSaHo$,
$$
\|f\|_{\iaHaSaHo}:=\left\{
\inf\left\|\left\{\sum_{j=1}
^{N}\left[\frac{\lambda_{j}}{
\|\1bf_{B_{j}}\|_{\HerzSo}}
\right]^{s}\1bf_{
B_{j}}\right\}^{\frac{1}{s}}
\right\|_{\HerzSo}\right\},
$$
where the infimum is taken over all
finite linear combinations of $f$,
namely, $N\in\mathbb{N}$,
$$f=\sum_{j=1}^{N}\lambda_ja_j,$$
$\{\lambda_j\}_{j=1}^N\subset[0,\infty)$,
and $\{a_j\}_{j=1}^N$
being $(\HerzSo,\,r,\,\,d)$-atoms supported,
respectively, in the balls
$\{B_j\}_{j=1}^N\subset\mathbb{B}$.
\end{definition}

Based on the known finite atomic
characterization of the
Hardy space $H_{X}(\rrn)$
associated with the ball
quasi-Banach function space $X$
(see \cite[Theorem 1.10]{YYY} or Lemma
\ref{finiatomll1} below),
we have the following
finite atomic characterization\index{finite
atomic characterization} of $\HaSaHo$.
Throughout this book,
\index{$C(\rrn)$}$C(\rrn)$ is defined
to be the set of all continuous
functions on $\rrn$.

\begin{theorem}\label{finiatoml}
Let $p$, $q$, $\omega$, $d$, $s$,
and $r$ be as in Definition
\ref{finiat}.
\begin{enumerate}
  \item[{\rm(i)}] If $$r\in\left(\max
  \left\{1,p,\frac{n}{\mw+n/p}\right\},
  \infty\right),$$ then
  $\|\cdot\|_{\iaHaSaHo}$ and
  $\|\cdot\|_{\HaSaHo}$
  are equivalent quasi-norms
on the generalized
finite atomic Herz--Hardy space $\iaHaSaHo$.
  \item[{\rm(ii)}] If $r=\infty$, then
  $\|\cdot\|_{H\Kmp_{\omega,
  \0bf,\mathrm{fin}}
  ^{p,q,\infty,d,s}(\rrn)}$
  and $\|\cdot\|_{\HaSaHo}$
  are equivalent quasi-norms
  on $H\Kmp_{\omega,\0bf,\mathrm{fin}}
  ^{p,q,\infty,d,s}(\rrn)
  \cap C(\rrn)$.
\end{enumerate}
\end{theorem}

To prove this theorem,
recall that Yan et al.\
\cite[Theorem 1.10]{YYY}
established the following finite atomic
characterization of Hardy spaces associated
with ball quasi-Banach function spaces.

\begin{lemma}\label{finiatomll1}
Let $X$, $r$, $d$, and $s$
be as in Lemma \ref{Atoll1}.
Then the \emph{finite atomic
Hardy space} $H_{\mathrm{fin}}^{X,r,d,s}
(\rrn)$\index{$H_{\mathrm{fin}}^{X,r,d,s}(\rrn)$},
associated with $X$, is defined
to be the set of all finite
linear combinations of
$(X,\,r,\,d)$-atoms. Moreover,
for any $f\in H_{\mathrm{fin}}^{
X,r,d,s}(\rrn)$,
\begin{equation}\label{fi}
\|f\|_{H_{\mathrm{fin}}^{
X,r,d,s}(\rrn)}:=\left\{
\inf\left\|\left\{\sum_{j=1}
^{N}\left(\frac{\lambda_{j}}{
\|\1bf_{B_{j}}\|_{X}}
\right)^{s}\1bf_{
B_{j}}\right\}^{\frac{1}{s}}
\right\|_{X}\right\},
\end{equation}
where the infimum is taken over all
finite linear combinations of $f$,
namely, $N\in\mathbb{N}$,
$$f=\sum_{j=1}^{N}\lambda_ja_j,$$
$\{\lambda_j\}_{j=1}^N\subset[0,\infty)$,
and $\{a_j\}_{j=1}^N$
being $(X,\,r,\,d)$-atoms supported,
respectively, in the balls
$\{B_j\}_{j=1}^N\subset\mathbb{B}$.
Then
\begin{enumerate}
  \item[{\rm(i)}] if $r\in(1,\infty)$,
  $\|\cdot\|_{H_{\mathrm{fin}}^{
  X,r,d,s}(\rrn)}$ and $\|\cdot\|_{H_X(\rrn)}$
  are equivalent quasi-norms
on $H_{\mathrm{fin}}^{
  X,r,d,s}(\rrn)$;
  \item[{\rm(ii)}] if $r=\infty$,
  $\|\cdot\|_{H_{\mathrm{fin}}^{
  X,\infty,d,s}(\rrn)}$
  and $\|\cdot\|_{H_{X}(\rrn)}$
  are equivalent quasi-norms
  on $H_{\mathrm{fin}}^{X,\infty,d,s}(\rrn)
  \cap C(\rrn)$.
\end{enumerate}
\end{lemma}

Lemma \ref{finiatomll1} implies that,
to show Theorem \ref{finiatoml},
we only need to prove that,
under the assumptions of
Theorem \ref{finiatoml},
$\HerzSo$ satisfies all the
assumptions of Lemma \ref{finiatomll1}.
Applying this idea
and some technical lemmas
obtained in the last
section, we now show Theorem \ref{finiatoml}.

\begin{proof}[Proof of Theorem
\ref{finiatoml}]
Let $p$, $q$, $\omega$, $d$,
$s$, and $r$ be as in the present theorem.
Then, using the assumption
$\m0(\omega)\in(-\frac{n}{p},\infty)$
and Theorem \ref{Th3}, we conclude that
the local generalized Herz space
$\HerzSo$ is a BQBF space.
Thus, in order to finish
the proof of the present theorem, it suffices
to show that all the assumptions
of Lemma \ref{finiatomll1}
hold true for $\HerzSo$.

First, let $\theta\in(0,s)$ be such
that
$$
\left\lfloor
n\left(\frac{1}{s}-1\right)
\right\rfloor\leq n\left(
\frac{1}{\theta}-1
\right)<\left\lfloor
n\left(\frac{1}{s}-1\right)
\right\rfloor+1.
$$
This implies that
$d\geq\lfloor n(1/s-1)\rfloor=
\lfloor n(1/\theta-1)\rfloor$.
We now prove that $\HerzSo$
satisfies Assumption \ref{assfs}
for the above $\theta$ and $s$.
Indeed, from Lemma \ref{Atoll3},
it follows that, for any
$\{f_j\}_{j\in\mathbb{N}}\subset
L^1_{\mathrm{loc}}(\rrn)$,
$$
\left\|
\left\{
\sum_{j\in\mathbb{N}}
\left[\mc^{(\theta)}(f_j)\right]^{s}
\right\}^{1/s}
\right\|_{\HerzSo}
\lesssim\left\|
\left(\sum_{j\in\mathbb{N}}
|f_j|^{s}
\right)^{1/s}
\right\|_{\HerzSo},
$$
which implies that $\HerzSo$
satisfies Assumption \ref{assfs}
for the above $\theta$ and $s$.

On the other hand, by Lemma \ref{mbhal},
we find that $[\HerzSo]^{1/s}$ is a
BBF space and,
for any $f\in L^1_{\mathrm{loc}}(\rrn)$,
\begin{equation*}
\left\|\mc^{((r/s)')}(f)
\right\|_{([\HerzSo]^{1/s})'}\lesssim
\left\|f\right\|_{([\HerzSo]^{1/s})'}.
\end{equation*}
This further implies that,
under the assumptions of
the present theorem, $\HerzSo$
satisfies all the assumptions
of Lemma \ref{finiatomll1},
which completes the proof of
Theorem \ref{finiatoml}.
\end{proof}

Then, as an application
of Theorem \ref{finiatoml},
we have the following finite atomic
characterization\index{finite atomic
characterization} of $H\MorrSo$,
which is a simple corollary
of Theorem \ref{finiatoml} and
Remarks \ref{remhs}(iv)
and \ref{remark4.10}(ii);
we omit the details.

\begin{corollary}
Let $p$, $q$, $\omega$,
$r$, $d$, and $s$ be as in Corollary
\ref{cor438}.
Then the \emph{generalized finite atomic Morrey--Hardy
space\index{generalized finite atomic Morrey--Hardy space}
$H\textsl{\textbf{M}}
_{\omega,\0bf,\mathrm{fin}}^
{p,q,r,d,s}(\rrn)$\index{$H\textsl{\textbf{M}}
_{\omega,\0bf,\mathrm{fin}}^{p,q,r,d,s}(\rrn)$}},
associated with the local generalized Morrey space
$\MorrSo$, is defined to be the set of all
finite linear combinations of
$(\MorrSo,\,r,\,d)$-atoms.
Moreover, for any $f\in H\textsl{\textbf{M}}
_{\omega,\0bf,\mathrm{fin}}^{p,q,r,d,s}(\rrn)$,
$$
\|f\|_{H\textsl{\textbf{M}}_{\omega,\0bf,
\mathrm{fin}}^{p,q,r,d,s}(\rrn)}:=\left\{
\inf\left\|\left\{\sum_{j=1}^{N}
\left[\frac{\lambda_{j}}{
\|\1bf_{B_{j}}\|_{\MorrSo}}\right]^{s}\1bf_{
B_{j}}\right\}^{\frac{1}{s}}
\right\|_{\MorrSo}\right\},
$$
where the infimum is taken over
all finite linear combinations
of $f$, namely, $N\in\mathbb{N}$,
$$f=\sum_{j=1}^{N}\lambda_ja_j,$$
$\{\lambda_j\}_{j=1}^N\subset[0,\infty)$,
and $\{a_j\}_{j=1}^N$
being $(\MorrSo,\,r,\,\,d)$-atoms supported,
respectively, in the balls
$\{B_j\}_{j=1}^N\subset\mathbb{B}$.
Then
\begin{enumerate}
  \item[{\rm(i)}] if $$r\in\left(
  \frac{n}{\mw+n/p},\infty\right),$$
  $\|\cdot\|_{H\textsl{\textbf{M}}_
  {\omega,\0bf,\mathrm{fin}}^{p,q,r,d,s}(\rrn)}$
  and $\|\cdot\|_{H\MorrSo}$ are equivalent quasi-norms
on the generalized finite atomic Morrey--Hardy space
$H\textsl{\textbf{M}}_{\omega,\0bf,\mathrm{fin}}
^{p,q,r,d,s}(\rrn)$;
  \item[{\rm(ii)}] if $r=\infty$,
  $\|\cdot\|_{H\textsl{\textbf{M}}_
  {\omega,\0bf,\mathrm{fin}}^{p,q,\infty,d,s}(\rrn)}$
  and $\|\cdot\|_{H\MorrSo}$ are equivalent quasi-norms
  on $H\textsl{\textbf{M}}_{\omega,
  \0bf,\mathrm{fin}}^{p,q,\infty,d,s}(\rrn)
  \cap C(\rrn)$.
\end{enumerate}
\end{corollary}

Next, we are devoted to
establishing the finite atomic
characterization of the generalized
Herz--Hardy space $\HaSaH$.
For this purpose, we first
introduce the following generalized
finite atomic Herz--Hardy space.

\begin{definition}\label{finiatg}
Let $p$, $q\in(0,\infty)$, $\omega\in M(\rp)$ with
$\m0(\omega)\in(-\frac{n}{p},\infty)$ and
$$
-\frac{n}{p}<\mi(\omega)\leq
\MI(\omega)<0,
$$
$r\in(\max\{1,p,\frac{n}{\mw+n/p}\},\infty]$,
$s\in(0,\min\{
1,p,q,\frac{n}{\Mw+n/p}\}),$
and $d\geq\lfloor n(
1/s-1)\rfloor$ be a fixed integer.
Then the \emph{generalized finite atomic
Herz--Hardy space\index{generalized
finite atomic Herz--Hardy \\ space}
$\iaHaSaH$\index{$\iaHaSaH$}},
associated with $\HerzS$, is defined
to be the set of all
finite linear combinations of
$(\HerzS,\,r,\,d)$-atoms.
Moreover,
for any $f\in\iaHaSaH$,
$$
\|f\|_{\iaHaSaH}:=\left\{
\inf\left\|\left\{\sum_{j=1}
^{N}\left[\frac{\lambda_{j}}{
\|\1bf_{B_{j}}\|_{\HerzS}}
\right]^{s}\1bf_{
B_{j}}\right\}^{\frac{1}{s}}
\right\|_{\HerzS}\right\},
$$
where the infimum is taken over all
finite linear combinations of $f$,
namely, $N\in\mathbb{N}$,
$$f=\sum_{j=1}^{N}
\lambda_ja_j,$$
$\{\lambda_j\}_{j=1}^N\subset[0,\infty)$,
and $\{a_j\}_{j=1}^N$
being $(\HerzS,\,r,\,\,d)$-atoms supported,
respectively, in the balls
$\{B_j\}_{j=1}^N\subset\mathbb{B}$.
\end{definition}

Then we have the following
finite atomic characterization
of $\HaSaH$.

\begin{theorem}\label{finiatomg}
Let $p$, $q$, $\omega$, $d$, $s$,
and $r$ be as in Definition
\ref{finiatg}.
\begin{enumerate}
  \item[{\rm(i)}] If $$r\in\left(\max
  \left\{1,p,\frac{n}{\mw+n/p}\right\},\infty\right),$$ then
  $\|\cdot\|_{\iaHaSaH}$ and $\|\cdot\|_{\HaSaH}$
  are equivalent quasi-norms
on the generalized
finite atomic Herz--Hardy space $\iaHaSaH$.
  \item[{\rm(ii)}] If $r=\infty$, then
  $\|\cdot\|_{H\Kmp_{\omega,\mathrm{fin}}
  ^{p,q,\infty,d,s}(\rrn)}$
  and $\|\cdot\|_{\HaSaH}$ are equivalent quasi-norms
  on $H\Kmp_{\omega,\mathrm{fin}}
  ^{p,q,\infty,d,s}(\rrn)
  \cap C(\rrn)$.
\end{enumerate}
\end{theorem}

To show Theorem \ref{finiatomg},
note that the associate space
of the global generalized Herz space
$\HerzS$ is still unknown. This means
that we can not prove Theorem \ref{finiatomg}
via applying Lemma \ref{finiatomll1}
directly. Fortunately, by
the proof of \cite[Theorem 1.10]{YYY},
we find that, if the atomic characterization
of $H_X(\rrn)$ holds true,
Lemma \ref{finiatomll1}
still holds true even when there is
not any assumption about
the associate space of $X$.
Namely, we have the following
finite atomic characterization of
$H_X(\rrn)$.

\begin{lemma}\label{finiatomgl1}
Let $X$ be a ball quasi-Banach function space
satisfy Assumption \ref{assfs} with
some $\theta$, $s\in(0,1]$,
$r\in(1,\infty]$, and $d\geq\lfloor
n(1/\theta-1)\rfloor$ be a fixed
integer such that
$$
H_X(\rrn)=H^{X,r,d,s}(\rrn)
$$
with equivalent quasi-norms, where
the atomic Hardy space $H^{X,r,d,s}(\rrn)$
is defined as in Definition \ref{atsgx}.
Then
\begin{enumerate}
  \item[{\rm(i)}] if $r\in(1,\infty)$,
  $\|\cdot\|_{H_{\mathrm{fin}}^{
  X,r,d,s}(\rrn)}$ and $\|\cdot\|_{H_X(\rrn)}$
  are equivalent quasi-norms
on $H_{\mathrm{fin}}^{
  X,r,d,s}(\rrn)$;
  \item[{\rm(ii)}] if $r=\infty$,
  $\|\cdot\|_{H_{\mathrm{fin}}^{
  X,\infty,d,s}(\rrn)}$
  and $\|\cdot\|_{H_{X}(\rrn)}$
  are equivalent quasi-norms
  on $H_{\mathrm{fin}}^{X,\infty,d,s}(\rrn)
  \cap C(\rrn)$.
\end{enumerate}
\end{lemma}

\begin{proof}
Let all the symbols be as in the
present lemma and $f\in H_{\mathrm{fin}}^{
X,r,d,s}(\rrn)$. Then,
applying Lemma \ref{Th5.4l1}(ii),
\eqref{fi}, \eqref{atom hardyx}, and the assumption
that $$
H_X(\rrn)=H^{X,r,d,s}(\rrn)
$$ with equivalent quasi-norms,
we
find that
\begin{equation}\label{finiatomgl1e1}
\left\|f\right\|_{H_X(\rrn)}
\sim\left\|f\right\|_{H^{X,r,d,s}(\rrn)}
\lesssim\left\|f\right\|_{H_{\mathrm{fin}}^{
X,r,d,s}(\rrn)}.
\end{equation}
Conversely, using Assumption \ref{assfs}
and repeating the proof of \cite[Theorem
1.10]{YYY} with $q$ therein replaced by
$r$, we conclude that, if
$r\in(1,\infty)$, then
$$
\left\|f\right\|_{H_{\mathrm{fin}}^{X,r,d,s}(\rrn)}
\lesssim\left\|f\right\|_{H_X(\rrn)}
$$
and, if $r=\infty$ and $f\in C(\rrn)$, then
$$
\left\|f\right\|_{H_{\mathrm{fin}}^{
X,\infty,d,s}(\rrn)}
\lesssim\left\|f\right\|_{H_X(\rrn)}.
$$
These, combined with \eqref{finiatomgl1e1},
further imply that both
(i) and (ii) of the present lemma hold true.
Thus, the proof of Lemma \ref{finiatomgl1}
is completed.
\end{proof}

Via this finite atomic characterization
of $H_X(\rrn)$, we now show Theorem \ref{finiatomg}.

\begin{proof}[Proof of Theorem
\ref{finiatomg}]
Let $p$, $q$, $\omega$,
$d$, $s$, and $r$ be as in the
present theorem. Then, applying
the assumptions $\m0(\omega)\in
(-\frac{n}{p},\infty)$ and
$\MI(\omega)\in(-\infty,0)$,
and Theorem \ref{Th2}, we find
that the global generalized Herz
space $\HerzS$ is a BQBF space.
Therefore, to complete the
proof of the present theorem,
we only need to show that
all the assumptions of
Lemma \ref{finiatomgl1} are satisfied
for $\HerzS$.

First, we prove that
there exists a $\theta\in(0,1]$ such that
$d\geq\lfloor n(1/\theta-1)\rfloor$
and,
for this $\theta$ and the $s$ same as in
the present theorem, $\HerzS$
satisfies Assumption \ref{assfs}, namely,
for any $\{f_j\}_{j\in\mathbb{N}}\subset
L^1_{\mathrm{loc}}(\rrn)$,
\begin{equation}\label{finiatomge1}
\left\|
\left\{
\sum_{j\in\mathbb{N}}
\left[\mc^{(\theta)}(f_j)\right]^{s}
\right\}^{1/s}
\right\|_{\HerzSo}
\lesssim\left\|
\left(\sum_{j\in\mathbb{N}}
|f_j|^{s}
\right)^{1/s}
\right\|_{\HerzSo}.
\end{equation}
Indeed, from the definition of
$\lfloor n(1/s-1)\rfloor$, it follows that
$$
\left\lfloor
n\left(\frac{1}{s}-1\right)
\right\rfloor\leq n\left(
\frac{1}{s}-1
\right)<\left\lfloor
n\left(\frac{1}{s}-1\right)
\right\rfloor+1.
$$
Thus, we can choose a
$\theta\in(0,s)$ such that
$$
\left\lfloor
n\left(\frac{1}{s}-1\right)
\right\rfloor\leq n\left(
\frac{1}{\theta}-1
\right)<\left\lfloor
n\left(\frac{1}{s}-1\right)
\right\rfloor+1,
$$
which further implies that
$d\geq \lfloor n(1/s-1)\rfloor
=\lfloor n(1/\theta-1)\rfloor$.
On the other hand, by
Lemma \ref{Atogl5}, we conclude
that \eqref{finiatomge1}
holds true. This
implies that Assumption \ref{assfs}
holds true for this $\theta$ and
$s$.

Next, using Theorem \ref{Atog},
we find that
$$
\HaSaH=\aHaSaH
$$
with equivalent quasi-norms.
This further implies that
$\HerzS$ satisfies all the
assumptions of Lemma \ref{finiatomgl1},
which completes the
proof of Theorem \ref{finiatomg}.
\end{proof}

Combining Theorem \ref{finiatomg} and
Remarks \ref{remhs}(iv)
and \ref{remark4.10}(ii), we immediately
obtain the following
finite atomic characterization of
the generalized Morrey--Hardy space
$H\MorrS$; we omit the details.

\begin{corollary}
Let $p$, $q$, $\omega$,
$r$, $d$, and $s$ be as in Corollary
\ref{cor438}.
Then the \emph{generalized
finite atomic Morrey--Hardy
space\index{generalized finite atomic Morrey--Hardy space}
$H\textsl{\textbf{M}}
_{\omega,\mathrm{fin}}^
{p,q,r,d,s}(\rrn)$\index{$H\textsl{\textbf{M}}
_{\omega,\mathrm{fin}}^{p,q,r,d,s}(\rrn)$}},
associated with the global generalized Morrey space
$\MorrS$, is defined to be the set of all
finite linear combinations of
$(\MorrS,\,r,\,d)$-atoms.
Moreover, for any $f\in H\textsl{\textbf{M}}
_{\omega,\mathrm{fin}}^{p,q,r,d,s}(\rrn)$,
$$
\|f\|_{H\textsl{\textbf{M}}_{\omega,
\mathrm{fin}}^{p,q,r,d,s}(\rrn)}:=\left\{
\inf\left\|\left\{\sum_{j=1}^{N}
\left[\frac{\lambda_{j}}{
\|\1bf_{B_{j}}\|_{\MorrS}}\right]^{s}\1bf_{
B_{j}}\right\}^{\frac{1}{s}}
\right\|_{\MorrS}\right\},
$$
where the infimum is taken over
all finite linear combinations
of $f$, namely, $N\in\mathbb{N}$,
$$f=\sum_{j=1}^{N}\lambda_ja_j,$$
$\{\lambda_j\}_{j=1}^N\subset[0,\infty)$,
and $\{a_j\}_{j=1}^N$
being $(\MorrS,\,r,\,\,d)$-atoms supported,
respectively, in the balls
$\{B_j\}_{j=1}^N\subset\mathbb{B}$.
Then
\begin{enumerate}
  \item[{\rm(i)}] if $$r\in\left(
  \frac{n}{\mw+n/p},\infty\right),$$
  then $\|\cdot\|_{H\textsl{\textbf{M}}_
  {\omega,\mathrm{fin}}^{p,q,r,d,s}(\rrn)}$
  and $\|\cdot\|_{H\MorrS}$ are equivalent quasi-norms
on the generalized finite atomic Morrey--Hardy space
$H\textsl{\textbf{M}}_{\omega,\mathrm{fin}}
^{p,q,r,d,s}(\rrn)$;
  \item[{\rm(ii)}] if $r=\infty$,
  then $\|\cdot\|_{H\textsl{\textbf{M}}_
  {\omega,\mathrm{fin}}^{p,q,\infty,d,s}(\rrn)}$
  and $\|\cdot\|_{H\MorrS}$ are equivalent quasi-norms
  on $H\textsl{\textbf{M}}_{\omega,
  \mathrm{fin}}^{p,q,\infty,d,s}(\rrn)
  \cap C(\rrn)$.
\end{enumerate}
\end{corollary}

\section{Molecular Characterizations}

The target of this section is to establish the
molecular characterization of
generalized Herz--Hardy spaces.
Indeed, we first show the molecular
characterization
of the generalized
Herz--Hardy space $\HaSaHo$ via
the known molecular
characterization of the Hardy space
$H_{X}(\rrn)$ associated with the
ball quasi-Banach function space $X$
(see Lemma \ref{Molell1} below).
However, we should point out that,
due to the deficiency of the associate space
of $\HerzS$, the molecular characterization
of the generalized Herz--Hardy space $\HaSaH$
can not be obtained by applying Lemma \ref{Molell1}.
To overcome this obstacle,
we first establish an improved
molecular characterization
of $H_X(\rrn)$ (see Theorem \ref{Molegx} below)
via borrowing some ideas from the
proof of \cite[Theorem 3.9]{SHYY} and
get the rid of associate spaces.
Then, by this molecular
characterization, we obtain the
molecular characterization of $\HaSaH$.

We first use molecules to characterize
the Hardy space $\HaSaHo$ associated with
the local generalized Herz space $\HerzSo$.
To this end, we now introduce
the following concept of
$(\HerzSo,\,r,\,d,\,\tau)$-molecules.
In what follows, for any $j\in\mathbb{N}$
and $B\in\mathbb{B}$, let
$$\index{$S_{j}(B)$}S_{j}(B):=
\left(2^{j}B\right)\setminus
\left(2^{j-1}B\right)
\text{ and }S_{0}(B):=B.$$

\begin{definition}\label{mole}
Let $p$, $q\in(0,\infty)$, $\omega\in M(\rp)$
with $\m0(\omega)\in(\frac{n}{p},\infty)$,
$\tau\in(0,\infty)$,
$r\in[1,\infty]$, and $d\in\zp$.
Then a measurable function
$m$ on $\rrn$ is called a
\emph{$(\HerzSo,\,r,
\,d,\,\tau)$-molecule}\index{$(\HerzSo,\,r,\,d,\,\tau)$-molecule}
centered at a ball $B\in\mathbb{B}$ if
\begin{enumerate}
  \item[{\rm(i)}] for any $j\in\zp$,
  $$\left\|m\1bf_{S_{j}(B)}
  \right\|_{L^{r}(\rrn)}\leq2^{-\tau j}
\frac{|B|^{1/r}}{\|\1bf_{B}\|_{\HerzSo}};$$
 \item[{\rm(ii)}] for any
 $\alpha\in\zp^n$ with $|\alpha|\leq d$,
  $$\int_{\rrn}m(x)x^\alpha\,dx=0.$$
\end{enumerate}
\end{definition}

Then we establish the molecular characterization
of the generalized Herz--Hardy space $\HaSaHo$
as follows\index{molecular \\characterization}.

\begin{theorem}\label{Molel}
Let $p,\ q\in(0,\infty)$,
$\omega\in M(\rp)$ with
$\m0(\omega)\in(-\frac{n}{p},\infty)$
and $\mi(\omega)\in(-\frac{n}{p},\infty)$,
$$s\in\left(0,\min\left\{1,p,q,
\frac{n}{\Mw+n/p}\right\}\right),$$
$d\geq\lfloor n(1/s-1)\rfloor$ be a fixed integer,
$$r\in\left(\max\left\{1,p,\frac{n}
{\mw+n/p}\right\},\infty\right],$$ and
$\tau\in(0,\infty)$ with $\tau>n(1/s-1/r)$.
Then $f\in\HaSaHo$ if and only
if $f\in\mathcal{S}'(\rrn)$
and there exist a sequence
$\{m_{j}\}_{j\in\mathbb{N}}$
of $(\HerzSo,\,r,\,d,\,\tau)$-molecules
centered, respectively,
at the balls $\{B_{j}\}_{j\in\mathbb{N}}
\subset\mathbb{B}$ and
a sequence $\{\lambda_{j}\}
_{j\in\mathbb{N}}\subset
[0,\infty)$ such that
$
f=\sum_{j\in\mathbb{N}}\lambda_{j}a_{j}
$
in $\mathcal{S}'(\rrn)$, and
$$
\left\|\left\{\sum_{j\in\mathbb{N}}
\left[\frac{\lambda_{j}}
{\|\1bf_{B_{j}}\|_{\HerzSo}}
\right]^{s}\1bf_{B_{j}}
\right\}^{\frac{1}{s}}
\right\|_{\HerzSo}<\infty.
$$
Moreover, there exists
a constant $C\in[1,\infty)$
such that, for
any $f\in\HaSaHo$,
\begin{align*}
C^{-1}\|f\|_{\HaSaHo}&\leq\inf\left\{
\left\|\left\{\sum_{j\in\mathbb{N}}
\left[\frac{\lambda_{i}}
{\|\1bf_{B_{j}}\|_{\HerzSo}}\right]^s\1bf_
{B_{j}}\right\}^{\frac{1}{s}}\right\|_{\HerzSo}
\right\}\\&\leq C\|f\|_{\HaSaHo},
\end{align*}
where the infimum is taken
over all the decompositions of $f$ as above.
\end{theorem}

To show this theorem, we
need the molecular characterization
of Hardy spaces associated with
ball quasi-Banach
function spaces. First, we present
the definition of $(X,\,r,\,d,\,\tau)$-molecules
introduced in \cite[Definition 3.8]{SHYY}.

\begin{definition}\label{molex}
Let $X$ be a ball quasi-Banach function
space, $r\in[1,\infty]$, $d\in\zp$,
and $\tau\in(0,\infty)$.
Then a measurable function
$m$ on $\rrn$ is called an
\emph{$(X,\,r,
\,d,\,\tau)$-molecule}\index{$(X,\,r,\,d,\,\tau)$-molecule}
centered at a ball $B\in\mathbb{B}$ if
\begin{enumerate}
  \item[{\rm(i)}] for any $j\in\zp$,
  $$\left\|m\1bf_{S_{j}(B)}
  \right\|_{L^{r}(\rrn)}\leq2^{-\tau j}
\frac{|B|^{1/r}}{\|\1bf_{B}\|_{X}};$$
 \item[{\rm(ii)}] for any
 $\alpha\in\zp^n$ with $|\alpha|\leq d$,
  $$\int_{\rrn}m(x)x^\alpha\,dx=0.$$
\end{enumerate}
\end{definition}

\begin{remark}\label{molexr}
Let $X$, $r$, $d$, and $\tau$ be
as in Definition \ref{molex}.
Then it is easy to
see that, for any $(X,\,r,\,d)$-atom
$a$ supported in the ball $B\in
\mathbb{B}$, $a$ is also
an $(X,\,r,\,d,\,\tau)$-molecule
centered at $B$.
\end{remark}

Via these molecules,
Sawano et al.\ \cite[Theorem 3.9]{SHYY}
established the molecular
characterization of the Hardy space $H_{X}(\rrn)$,
which plays a key role in the proof of
the molecular characterization of $\HaSaHo$.
Namely, the following conclusion holds true.

\begin{lemma}\label{Molell1}
Let $X$, $r$, $d$, $\theta$,
and $s$ be as in Lemma \ref{Atoll1}
and $\tau\in(0,\infty)$ with
$\tau>n(1/\theta-1/r)$.
Then $f\in H_X(\rrn)$ if and only
if $f\in\mathcal{S}'(\rrn)$
and there exist a sequence
$\{m_{j}\}_{j\in\mathbb{N}}$
of $(X,\,r,\,d,\,\tau)$-molecules
centered, respectively,
at the balls $\{B_{j}\}_{j\in\mathbb{N}}
\subset\mathbb{B}$ and
a sequence $\{\lambda_{j}\}
_{j\in\mathbb{N}}\subset
[0,\infty)$ such that
$$
f=\sum_{j\in\mathbb{N}}\lambda_{j}a_{j}
$$
in $\mathcal{S}'(\rrn)$, and
$$
\left\|\left[\sum_{j\in\mathbb{N}}
\left(\frac{\lambda_{j}}
{\|\1bf_{B_{j}}\|_{X}}
\right)^{s}\1bf_{B_{j}}
\right]^{\frac{1}{s}}
\right\|_{X}<\infty.
$$
Moreover, for any $f\in H_X(\rrn)$,
\begin{align*}
\|f\|_{H_X(\rrn)}\sim
\inf\left\{
\left\|\left[\sum_{j\in\mathbb{N}}
\left(\frac{\lambda_{i}}
{\|\1bf_{B_{j}}\|_{X}}\right)^s\1bf_
{B_{j}}\right]^{\frac{1}{s}}\right\|_{X}
\right\}
\end{align*}
with the positive equivalence
constants independent of $f$,
where the infimum is taken
over all the decompositions of $f$ as above.
\end{lemma}

Applying Lemma \ref{Molell1}, we next
show Theorem \ref{Molel}.

\begin{proof}[Proof of Theorem \ref{Molel}]
Let $p$, $q$, $\omega$, $r$, $d$,
$s$, and $\tau$ be as in the
present theorem. Then, from the
assumption $\m0(\omega)\in(-
\frac{n}{p},\infty)$ and Theorem \ref{Th3},
it follows that the local generalized
Herz space $\HerzSo$ is a BQBF space.
Therefore, to
complete the proof of the present theorem,
we only need to show that all the assumptions
of Lemma \ref{Molell1} are satisfied
for $\HerzSo$.

First, let $\theta\in(0,s)$ satisfy
$$
\left\lfloor n\left(\frac{1}{s}-1\right)
\right\rfloor\leq n\left(\frac{1}{\theta}
-1\right)<\left\lfloor n\left(
\frac{1}{s}-1\right)
\right\rfloor+1
$$
and
$$
\tau>n\left(\frac{1}{\theta}-
\frac{1}{r}\right).
$$
This implies that $d\geq\lfloor
n(1/s-1)\rfloor=\lfloor n(1/\theta-1)\rfloor$.
We now show that $\HerzSo$
satisfies Assumption \ref{assfs}
for the above $\theta$ and $s$.
Indeed, applying Lemma \ref{Atoll3},
we find that, for any
$\{f_{j}\}_{j\in\mathbb{N}}
\subset L^1_{{\rm loc}}(\rrn)$,
\begin{equation}\label{Molele1}
\left\|\left\{\sum_{j\in\mathbb{N}}
\left[\mc^{(\theta)}
(f_{j})\right]^s\right\}^{1/s}
\right\|_{\HerzSo}\lesssim
\left\|\left(\sum_{j\in\mathbb{N}}
|f_{j}|^{s}\right)^{1/s}\right\|_{\HerzSo},
\end{equation}
which implies that
Assumption \ref{assfs} holds
true for $\HerzSo$ with the
above $\theta$ and $s$.

On the other hand, by Lemma \ref{mbhal},
we conclude that $[\HerzSo]^{1/s}$
is a BBF space and, for any $f\in
L_{{\rm loc}}^{1}(\rrn)$,
\begin{equation*}
\left\|\mc^{((r/s)')}(f)\right\|_
{([\HerzSo]^{1/s})'}\lesssim
\left\|f\right\|_{([\HerzSo]^{1/s})'}.
\end{equation*}
This, together with
\eqref{Molele1}, further implies
that, under the assumptions
of the present theorem,
$\HerzSo$ satisfies all the
assumptions of Lemma \ref{Molell1}
and hence finishes the
proof of Theorem \ref{Molel}.
\end{proof}

As an application of Theorem
\ref{Molel}, we now consider the
molecular characterization of the generalized
Morrey--Hardy space $H\MorrSo$.
For this purpose, we first introduce
$(\MorrSo,\,r,\,d,\,\tau)$-molecules as follows.

\begin{definition}\label{molem}
Let $p$, $q\in[1,\infty)$,
$\omega\in M(\rp)$
with $\MI(\omega)\in(-\infty,0)$
and
$$
-\frac{n}{p}<\m0(\omega)
\leq\M0(\omega)<0,
$$
$\tau\in(0,\infty)$,
$r\in[1,\infty]$, and $d\in\zp$.
Then a measurable function
$m$ on $\rrn$ is called an
\emph{$(\MorrSo,\,r,\,d,\,
\tau)$-molecule}\index{$(\MorrSo,\,r,\,d,\,\tau)$-molecule}
centered at a ball $B\in\mathbb{B}$ if
\begin{enumerate}
  \item[{\rm(i)}] for any $j\in\zp$,
  $$\left\|m\1bf_{S_{j}(B)}
  \right\|_{L^{r}(\rrn)}\leq2^{-\tau j}
\frac{|B|^{1/r}}{\|\1bf_{B}\|_{\MorrSo}};$$
 \item[{\rm(ii)}] for any $\alpha\in\zp^n$
 with $|\alpha|\leq d$,
  $\int_{\rrn}a(x)x^\alpha\,dx=0$.
\end{enumerate}
\end{definition}

Then we have the following molecular characterization
of $H\MorrSo$, which can be deduced directly from
Theorem \ref{Molel} and
Remarks \ref{remhs}(iv) and \ref{remark4.10}(ii);
we omit the details.

\begin{corollary}\label{Molelm}
Let $p$, $q$, $\omega$, $r$,
$d$, and $s$ be as in
Corollary \ref{cor438} and
$\tau\in(0,\infty)$ with $\tau>n(1/s-1/r)$.
Then $f\in H\MorrSo$ if and only
if $f\in\mathcal{S}'(\rrn)$
and there exist a sequence
$\{m_{j}\}_{j\in\mathbb{N}}$
of $(\MorrSo,\,r,\,d,\,\tau)$-molecules
centered, respectively,
at the balls $\{B_{j}\}_{j\in\mathbb{N}}
\subset\mathbb{B}$ and
a sequence $\{\lambda_{j}\}
_{j\in\mathbb{N}}\subset
[0,\infty)$ such that
$$
f=\sum_{j\in\mathbb{N}}\lambda_{j}a_{j}
$$
in $\mathcal{S}'(\rrn)$, and
$$
\left\|\left\{\sum_{j\in\mathbb{N}}
\left[\frac{\lambda_{j}}
{\|\1bf_{B_{j}}\|_{\MorrSo}}\right]
^{s}\1bf_{B_{j}}\right\}^{\frac{1}{s}}
\right\|_{\MorrSo}<\infty.
$$
Moreover, for any $f\in H\MorrSo$,
\begin{align*}
\|f\|_{H\MorrSo}\sim\inf\left\{
\left\|\left\{\sum_{j\in\mathbb{N}}
\left[\frac{\lambda_{i}}
{\|\1bf_{B_{j}}\|_{\MorrSo}}\right]^s\1bf_
{B_{j}}\right\}^{\frac{1}{s}}\right\|_{\MorrSo}
\right\}
\end{align*}
with the positive equivalence
constants independent of $f$,
where the infimum is taken
over all the decompositions of $f$ as above.
\end{corollary}

Next, we establish the molecular characterization
of the generalized Herz--Hardy space
$\HaSaH$. To achieve this,
we first introduce the following
concept of $(\HerzS,\,r,\,d,\,\tau)$-molecules.

\begin{definition}\label{moleg}
Let $p$, $q\in(0,\infty)$, $\omega\in M(\rp)$
with $\m0(\omega)\in(-\frac{n}{p},\infty)$
and $\MI(\omega)\in(-\infty,0)$,
$\tau\in(0,\infty)$, $r\in[1,\infty]$,
and $d\in\zp$.
Then a measurable function $m$ on
$\rrn$ is called a
\emph{$(\HerzS,\,r,\,d,\,\tau)
$-molecule}\index{$(\HerzS,\,r,\,d,\,\tau)$-molecule}
centered at a ball $B\in\mathbb{B}$ if
\begin{enumerate}
  \item[{\rm(i)}] for any
  $j\in\zp$,
  $$
  \left\|m\1bf_{S_j(B)}\right\|_{L^r(\rrn)}
  \leq2^{-j\tau}\frac{|B|^{1/r}}
  {\|\1bf_B\|_{\HerzS}};
  $$
 \item[{\rm(ii)}] for any $\alpha\in
 \zp^n$ with $|\alpha|\leq d$,
  $$\int_{\rrn}m(x)x^\alpha\,dx=0.$$
\end{enumerate}
\end{definition}

Via these molecules, we characterize the
generalized
Herz--Hardy space $\HaSaH$ associated with
the global generalized Herz space $\HerzS$
as follows.

\begin{theorem}\label{Moleg}
Let $p,\ q\in(0,\infty)$,
$\omega\in M(\rp)$ with
$\m0(\omega)\in(-\frac{n}{p},\infty)$
and $$-\frac{n}{p}<\mi
(\omega)\leq\MI(\omega)<0,$$
$r\in(\max\{1,p,\frac{n}{\mw+n/p}\},\infty]$,
$s\in(0,\min\{1,p,q,\frac{n}{\Mw+n/p}\})$,
$d\geq\lfloor n(\frac{1}{s}-1)\rfloor$ be a
fixed integer, and
$
\tau\in(n(
\frac{1}{s}-\frac{1}{r}),\infty).
$
Then $f\in\HaSaH$ if and only if
$f\in\mathcal{S}'(\rrn)$
and there exist $\{\lambda_{j}\}
_{j\in\mathbb{N}}\subset[0,\infty)$ and
a sequence $\{m_{j}\}_{j\in\mathbb{N}}$
of $(\HerzS,\,r,\,d,\,\tau)$-molecules
centered, respectively,
at the balls $\{B_{j}\}_{j\in\mathbb{N}}
\subset\mathbb{B}$
such that
\begin{equation*}
f=\sum_{j\in\mathbb{N}}\lambda_{j}m_{j}
\end{equation*}
in $\mathcal{S}'(\rrn)$, and
\begin{equation*}
\left\|\left\{\sum_{j\in\mathbb{N}}
\left[\frac{\lambda_{j}}
{\|\1bf_{B_{j}}\|_{\HerzS}}
\right]^{s}\1bf_{B_{j}}\right\}^{\frac{1}{s}}
\right\|_{\HerzS}<\infty.
\end{equation*}
Moreover, there exists a
constant $C\in[1,\infty)$ such that, for
any $f\in\HaSaH$,
\begin{align*}
C^{-1}\|f\|_{\HaSaH}&\leq\inf\left\{
\left\|\left\{\sum_{j\in\mathbb{N}}
\left[\frac{\lambda_{j}}
{\|\1bf_{B_{j}}\|_{\HerzS}}\right]^s\1bf_
{B_{j}}\right\}^{\frac{1}{s}}\right\|_{\HerzS}
\right\}\\
&\leq C\|f\|_{\HaSaH},
\end{align*}
where the infimum is taken
over all the decompositions of $f$ as above.
\end{theorem}

Due to the deficiency of associate spaces,
we can not show Theorem \ref{Moleg}
via using the known molecular
characterization of the
Hardy space $H_X(\rrn)$ associated with
the ball quasi-Banach function space
$X$ (see Lemma \ref{Molell1} above)
directly. To overcome this obstacle,
we first give the following
molecular characterization
of $H_X(\rrn)$
when there is no assumption about the
associate space $X'$.

\begin{theorem}\label{Molegx}
Let $X$, $r$, $d$, $s$, and $\theta$
be as in Theorem \ref{Atogx}
and $\tau\in(0,\infty)$
with $\tau>n(\frac{1}{\theta}-\frac{1}{r})$.
Then $f\in H_X(\rrn)$ if and only if
$f\in\mathcal{S}'(\rrn)$
and there exist $\{\lambda_{j}\}
_{j\in\mathbb{N}}\subset[0,\infty)$ and
a sequence $\{m_{j}\}_{j\in\mathbb{N}}$
of $(X,\,r,\,d,\,\tau)$-molecules
centered, respectively,
at the balls $\{B_{j}\}_{j\in\mathbb{N}}
\subset\mathbb{B}$
such that
\begin{equation}\label{Molegxe01}
f=\sum_{j\in\mathbb{N}}\lambda_{j}m_{j}
\end{equation}
in $\mathcal{S}'(\rrn)$, and
\begin{equation}\label{Molegxe02}
\left\|\left[\sum_{j\in\mathbb{N}}
\left(\frac{\lambda_{j}}
{\|\1bf_{B_{j}}\|_{X}}
\right)^{s}\1bf_{B_{j}}\right]
^{\frac{1}{s}}
\right\|_{X}<\infty.
\end{equation}
Moreover, there exists a
constant $C\in[1,\infty)$ such that, for
any $f\in H_X(\rrn)$,
\begin{align}\label{Molegxe03}
C^{-1}\|f\|_{H_X(\rrn)}\leq\inf\left\{
\left\|\left[\sum_{j\in\mathbb{N}}
\left(\frac{\lambda_{i}}
{\|\1bf_{B_{j}}\|_{X}}\right)^s\1bf_
{B_{j}}\right]^{\frac{1}{s}}
\right\|_{X}
\right\}\leq C\|f\|_{H_X(\rrn)},
\end{align}
where the infimum is taken
over all the decompositions of $f$ as above.
\end{theorem}

\begin{remark}
We should point out that
Theorem \ref{Molegx} is an improved
version of the known
molecular characterization obtained
in \cite[Theorem 3.9]{SHYY}. Indeed,
if $Y\equiv(X^{1/s})'$ in Theorem \ref{Molegx},
then this theorem goes back to \cite[Theorem 3.9]{SHYY}.
\end{remark}

In order to show the above theorem,
we require the following pointwise
estimate on the radial maximal function
of $(X,\,r,\,d,\,\tau)$-molecules.

\begin{lemma}\label{Molegl1}
Let $X$ be a ball quasi-Banach
function space, $r\in(1,\infty]$,
$\theta\in(0,1]$,
$d\geq \lfloor n(1/\theta-1)\rfloor$ be
a fixed integer, and $\tau\in
(n(\frac{1}{\theta}-\frac{1}{r}),\infty)$.
Assume that $\phi\in\mathcal{S}(\rrn)$
satisfies $\supp(\phi)\subset
B(\0bf,1)$. Then there exists a positive
constant $C$ such that, for any $(X,\,r,
\,d,\,\tau)$-molecule $m$ centered at the ball
$B\in\mathbb{B}$,
\begin{align*}
M(m,\phi)&\leq
C\left[\mc(m)\1bf_{4B}+\sum_{j=3}^{\infty}
\mc\left(m\1bf_{(2^{j+1}B)\setminus
(2^{j-2}B)}\right)\1bf_{(2^jB)\setminus
(2^{j-1}B)}\right.\\&\quad+\left.
\frac{1}{\|\1bf_B\|_{X}}
\mc^{(\theta)}\left(\1bf_B\right)
\right],
\end{align*}
where $M$ is the radial maximal
function defined as in Definition
\ref{smax}(i).
\end{lemma}

\begin{proof}
Let all the symbols be as in the
present lemma and $B:=
B(x_0,r_0)$ with $x_0\in\rrn$
and $r_0\in(0,\infty)$. Then, from
Lemma \ref{Atogxl3} with
$\Phi$ and $f$ therein replaced,
respectively, by $\|\phi\|_
{L^{\infty}(\rrn)}\1bf_{B(\0bf,1)}$
and $m$,
we deduce that, for any $x\in 4B$,
\begin{equation}\label{Molegl1e1}
M(m,\phi)(x)=
\sup_{t\in(0,\infty)}
\left|m\ast\phi_t(x)\right|\lesssim\mc
(m)(x),
\end{equation}
which is the desired estimate of
$M(m,\phi)(x)$ when $x\in 4B$.

Next, we estimate $M(m,\phi)(x)$
for any $x\in S_j(B)$ with
$j\in\mathbb{N}\cap[3,\infty)$.
To this end, let $j\in\mathbb{N}
\cap[3,\infty)$, $x\in
S_j(B)$, $t\in(0,\infty)$, and
$$d_{\theta}:=\left\lfloor n\left(
\frac{1}{\theta}-1\right)\right\rfloor.$$
Then, by Definition \ref{molex}(ii),
we find that
\begin{align}\label{Molegl1e2}
&\left|m\ast\phi_t(x)\right|\notag\\
&\quad=\left|\int_{\rrn}\phi_t(x-y)m(y)\,dy\right|
\notag\\
&\quad=\left|
\frac{1}{t^n}\int_{\rrn}\left[\phi\left(
\frac{x-y}{t}\right)-\sum_{\gfz{\gamma\in\zp^n}
{|\gamma|\leq d_{\theta}}}
\frac{\partial^{\gamma}(\phi(\frac{x-\cdot}
{t}))(x_0)}{\gamma!}(y-x_0)^{\gamma}\right]
m(y)\,dy\right|\notag\\
&\quad\leq\frac{1}{t^n}\left[\left|
\int_{2^{j-2}B}\left[\phi\left(
\frac{x-y}{t}\right)\r.\r.\r.\notag\\
&\qquad\lf.\lf.\lf.-\sum_{\gfz{\gamma\in\zp^n}
{|\gamma|\leq d_{\theta}}}\frac{
\partial^{\gamma}(\phi(\frac{x-\cdot}
{t}))(x_0)}{\gamma!}(y-x_0)^{\gamma}\right]
m(y)\,dy\right|\right.\notag\\
&\qquad\left.+\left|\int_{(2^{j+1}B)\setminus
(2^{j-2}B)}\cdots\right|+\left|\int_{
(2^{j-2}B)^{\complement}}\cdots\right|\right]\notag\\
&\quad\leq\int_{2^{j-2}B}\left|\phi\left(
\frac{x-y}{t}\right)-\sum_{\gfz{\gamma\in\zp^n}
{|\gamma|\leq d_{\theta}}}
\frac{\partial^{\gamma}(\phi(\frac{x-\cdot}
{t}))(x_0)}{\gamma!}(y-x_0)^{\gamma}\right|
\left|m(y)\right|\,dy\notag\\
&\qquad+\left|\int_{(2^{j+1}B)\setminus
(2^{j-2}B)}\phi_t(x-y)m(y)\,dy\right|\notag\\
&\qquad+\frac{1}{t^n}\int_{(2^{j+1}B)^{\complement}}
\left|\phi\left(\frac{x-y}{t}
\right)\right|\left|m(y)\right|\,dy\notag\\
&\qquad+\frac{1}{t^n}
\int_{(2^{j-2}B)^{\complement}}
\left|\sum_{\gfz{\gamma\in\zp^n}{|\gamma|
\leq d_{\theta}}}\frac{
\partial^{\gamma}(\phi(\frac{x-\cdot}
{t}))(x_0)}{\gamma!}(y-x_0)^{\gamma}\right|
\left|m(y)\right|\,dy\notag\\
&\quad=:\mathrm{A}_1+\mathrm{A}_2
+\mathrm{A}_3+\mathrm{A}_4.
\end{align}
We next estimate $\mathrm{A}_1$,
$\mathrm{A}_2$, $\mathrm{A}_3$,
and $\mathrm{A}_4$, respectively.

First, we deal with $\mathrm{A}_1$.
Indeed, from the Taylor remainder theorem
and the Tonelli theorem, it follows
that, for any $y\in 2^{j-2}B$, there exists
a $t_y\in(0,1)$ such that
\begin{align*}
\mathrm{A}_1&=
\sum_{k=0}^{j-2}\frac{1}{t^n}
\int_{S_k(B)}\left|\sum_{\gfz{\gamma\in\zp^n}
{|\gamma|=d_{\theta}+1}}\frac{\partial^{\gamma}
(\phi(\frac{x-\cdot}{t}))(t_y+(1-t_y)x_0)}{\gamma!}
\left(y-x_0\right)^{\gamma}
\right|\\&\quad\times\left|m(y)\right|\,dy\\
&\lesssim\sum_{k=0}^{j-2}\sum_{\gfz{\gamma\in\zp^n}
{|\gamma|=d_{\theta}+1}}\frac{1}{
t^{n+d_{\theta}+1}}\\
&\quad\times\int_{S_{k}(B)}\left|
\partial^{\gamma}\phi\left(
\frac{x-t_y-(1-t_y)x_0}{t}\right)
\right|\left|y-x_0\right|^{d_{\theta}+1}
\left|m(y)\right|\,dy\\
&\lesssim\sum_{k=0}^{j-2}
\int_{S_k(B)}\frac{|y-x_0|^{d_{\theta}+1}}{
|x-t_y-(1-t_y)x_0|^{n+d_{\theta}+1}}\left|
m(y)\right|\,dy.
\end{align*}
Observe that, for any
$y\in2^{j-2}B$, we have
$$|y-x_0|<2^{j-2}r_0\leq\frac{1}{2}|x-x_0|,$$
which, combined with the assumption $t_y\in(0,1)$,
further implies that
$$
\left|x-t_y-(1-t_y)x_0\right|
\geq|x-x_0|-|y-y_0|>\frac{1}{2}|x-x_0|.
$$
Using this and the assumption
$n(1/\theta-1)<d_{\theta}+1$,
we conclude that
\begin{align}\label{Molegl1e4}
\mathrm{A}_1&\lesssim\sum_{k=0}^{j-2}
\int_{S_k(B)}\frac{|y-x_0|^{d_{\theta}+1}}{
|x-x_0|^{n+d_{\theta}+1}}\left|m(y)
\right|\,dy\notag\\
&\lesssim\sum_{k=0}^{j-2}\int_{S_k(B)}
\frac{|y-x_0|^{\frac{n}{\theta}-n}}{
|x-x_0|^{\frac{n}{\theta}}}\,dy
\lesssim\sum_{k=0}^{j-2}\frac{(2^kr_0)^{
\frac{n}{\theta}-n}}{|x-x_0|^{\frac{n}{
\theta}}}\left\|m\1bf_{S_k(B)}\right\|_{L^1(\rrn)}.
\end{align}
Notice that, for any
$k\in\zp$, by the H\"{o}lder
inequality and Definition \ref{molex}(i),
we find that
\begin{align}\label{Molegl1e5}
\left\|m\1bf_{S_k(B)}\right\|_{L^1(\rrn)}
&\leq\left|S_k(B)\right|^{1-\frac{1}{r}}
\left\|m\1bf_{S_k(B)}\right\|_{L^r(\rrn)}\notag\\
&\leq\left|S_k(B)\right|^{1-\frac{1}{r}}
2^{-k\tau}\frac{|B|^{\frac{1}{r}}}{\|\1bf_B\|_X}
\sim2^{-k(\tau-n+\frac{n}{r})}\frac{r_0^{n}}{
\|\1bf_B\|_X}.
\end{align}
From this, \eqref{Molegl1e4}, and the
assumption $\tau\in(n(\frac{1}{\theta}-
\frac{1}{r}),\infty)$,
we deduce that
\begin{align}\label{Molegl1e6}
\mathrm{A}_1\lesssim\frac{1}{\|\1bf_B\|_X}
\frac{r_0^{\frac{n}{\theta}}}{|x-x_0|^{
\frac{n}{\theta}}}\sum_{k=0}^{j-2}2^{-k(
\tau-\frac{n}{\theta}+\frac{n}{r})}
\sim\frac{1}{\|\1bf_B\|_X}
\frac{r_0^{\frac{n}{\theta}}}{|x-x_0|^{
\frac{n}{\theta}}},
\end{align}
which is the desired estimate of
$\mathrm{A}_1$.

Next, we estimate $\mathrm{A}_2$.
Indeed, applying Lemma \ref{Atogxl3}
with $f:=m\1bf_{(2^{j+1}B)\setminus
(2^{j-2}B)}$ and
$\Phi:=\|\phi\|_{L^{\infty}(\rrn)}
\1bf_{B(\0bf,1)}$,
we conclude that
\begin{align}\label{Molegl1e7}
\mathrm{A}_2=\left|\left[m\1bf_{
(2^{j+1}B)\setminus(2^{j-2}B)}\right]
\ast\phi(x)\right|\lesssim
\mc\left(m\1bf_{
(2^{j+1}B)\setminus(2^{j-2}B)}\right)(x).
\end{align}
This finishes the estimate of $\mathrm{A}_2$.

We now deal with $\mathrm{A}_3$.
Notice that, for any $y\in
(2^{j+1}B)^{\complement}$, we have
$$
|y-x_0|\geq2^{j+1}r_0>2|x-x_0|,
$$
which further implies that
$$
|x-y|\geq|y-x_0|-|x-x_0|>|x-x_0|.
$$
From this, the Tonelli theorem,
and \eqref{Molegl1e5}, it then follows
that
\begin{align*}
\mathrm{A}_3&\lesssim
\sum_{k=j+2}^{\infty}
\int_{S_k(B)}\frac{1}{|x-y|^n}
\left|m(y)\right|\,dy
\lesssim\sum_{k=j+2}^{\infty}
\frac{1}{|x-x_0|^n}\left\|
m\1bf_{S_k(B)}\right\|_{L^1(\rrn)}\\
&\lesssim\sum_{k=j+2}^{\infty}
2^{-k(\tau+\frac{n}{r})}\frac{1}{
\|\1bf_B\|_X}\frac{(2^kr_0)^{n}}{|x-x_0|^n}.
\end{align*}
In addition, for any $k\in\mathbb{N}
\cap[j+2,\infty)$, we have
$$|x-x_0|<2^jr_0\leq2^{k-2}r_0,$$
which implies that $\frac{\tk r_0}{|x-x_0|}
\in[4,\infty)$.
This, together with
the assumptions $\theta\in(0,1]$
and $\tau\in(n(\frac{1}{\theta}-\frac{1}{r})
,\infty)$,
further implies that
\begin{align}\label{Molegl1e9}
\mathrm{A}_3&\lesssim\sum_{k=j+2}^{\infty}
2^{-k(\tau+\frac{n}{r})}\frac{1}{\|\1bf_B\|_X}
\left(\frac{\tk r_0}{|x-x_0|}\right)^{
\frac{n}{\theta}}\notag\\
&\sim\frac{1}{\|\1bf_B\|_X}\frac{r_0^{\frac{n}{\theta}}}
{|x-x_0|^{\frac{n}{\theta}}}\sum_{k=j+2}^{\infty}
2^{-k(\tau-\frac{n}{\theta}+\frac{n}{r})}
\sim\frac{1}{\|\1bf_B\|_X}\frac{r_0^{\frac{n}{\theta}}}
{|x-x_0|^{\frac{n}{\theta}}},
\end{align}
which is the desired estimate of
$\mathrm{A}_3$.

Finally, we turn to estimate $\mathrm{A}_4$.
Indeed, for any $y\in(2^{j-2}B)^{\complement}$,
it holds true that
$$
|y-x_0|\geq2^{j-2}r_0>\frac{1}{4}|x-x_0|,
$$
which further implies that
$\frac{4|y-x_0|}{|x-x_0|}\in(1,\infty)$.
Applying this, the assumption
$d_{\theta}\leq n(1/\theta-1)$,
\eqref{Molegl1e5}, and the
assumption $\tau\in(n(\frac{1}{\theta}-
\frac{1}{r}),\infty)$, we
conclude that
\begin{align*}
\mathrm{A}_4&\lesssim
\sum_{k=j-1}^{\infty}\sum_{\gfz{\gamma\in\zp^n}
{|\gamma|\leq d_{\theta}}}\frac{1}{t^{n+|\gamma|}}
\int_{S_k(B)}\left|\partial^{\gamma}
\phi\left(\frac{x-x_0}{t}\right)\right|
\left|y-x_0\right|^{|\gamma|}\left|
m(y)\right|\,dy\\
&\lesssim\sum_{k=j-1}^{\infty}
\sum_{\gfz{\gamma\in\zp^n}{|\gamma|\leq
d_{\theta}}}\int_{S_{k}(B)}\frac{|y-x_0|^{|
\gamma|}}{|x-x_0|^{n+|\gamma|}}\,dy\\
&\lesssim\sum_{k=j-1}^{\infty}
\int_{S_k(B)}\frac{|y-x_0|^{\frac{n}{\theta}-n}}
{|x-x_0|^{\frac{n}{\theta}}}\left|m(y)\right|
\,dy\\&\lesssim
\sum_{k=j-1}^{\infty}\frac{(2^kr_0)^{
\frac{n}{\theta}-n}}{|x-x_0|^{\frac{n}{\theta}}}
\left\|m\1bf_{S_k(B)}\right\|_{L^1(\rrn)}\\
&\lesssim\frac{1}{\|\1bf_B\|_X}\frac{r_0^{
\frac{n}{\theta}}}{|x-x_0|^{\frac{n}{\theta}}}
\sum_{k=j-1}^{\infty}2^{-k(\tau-\frac{n}{\theta}
+\frac{n}{r})}\sim\frac{1}{\|\1bf_B\|_X}\frac{r_0^{
\frac{n}{\theta}}}{|x-x_0|^{\frac{n}{\theta}}},
\end{align*}
which is the desired estimate of
$\mathrm{A}_4$. Combining this,
\eqref{Molegl1e2}, \eqref{Molegl1e6},
\eqref{Molegl1e7}, \eqref{Molegl1e9},
and an argument similar to
that used in the estimation of \eqref{atogxl1e1}
with $r_B$ and $x_B$ therein replaced,
respectively, by $r_0$ and $x_0$,
we find that, for any
$x\in S_j(B)$ with $j\in\mathbb{N}
\cap[3,\infty)$ and $t\in(0,\infty)$,
\begin{align*}
\left|m\ast\phi_t(x)\right|
&\lesssim\mc\left(m\1bf_{(2^{j+1}B)
\setminus(2^{j-2}B)}\right)(x)
+\frac{1}{\|\1bf_B\|_X}
\frac{r_0^{\frac{n}{\theta}}}{
|x-x_0|^{\frac{n}{\theta}}}\\
&\lesssim\mc\left(m\1bf_{(2^{j+1}B)
\setminus(2^{j-2}B)}\right)(x)+\frac{1}{
\|\1bf_B\|_X}\mc^{(\theta)}\left(
\1bf_B\right)(x),
\end{align*}
which, together with the arbitrariness
of $t$ and \eqref{Molegl1e1},
further implies that
\begin{align*}
M(m,\phi)
&\lesssim\mc(m)\1bf_{4B}+\sum_{j=3}^{\infty}
\mc\left(m\1bf_{(2^{j+1}B)
\setminus(2^{j-2}B)}\right)\1bf_{(2^{j}B)\setminus
(2^{j-1}B)}\\&\quad+\frac{1}{
\|\1bf_B\|_X}\mc^{(\theta)}\left(
\1bf_B\right).
\end{align*}
This then finishes the proof of Lemma
\ref{Molegl1}.
\end{proof}

We now show Theorem \ref{Molegx}.

\begin{proof}[Proof of
Theorem \ref{Molegx}]
Let all the symbols be as in the present
theorem. We first prove the necessity.
Indeed, let $f\in H_X(\rrn)$.
Then, applying Theorem \ref{Atogx}, we find that,
under the assumptions of the present
theorem, $H_X(\rrn)=H^{X,r,d,s}(\rrn)$
with equivalent
quasi-norms. This implies that
$f\in H^{X,r,d,s}(\rrn)$. Therefore,
there exist a sequence $\{\lambda_j\}_{j\in\mathbb{N}}
\subset[0,\infty)$ and
$\{a_j\}_{j\in\mathbb{N}}$ of
$(X,\,r,\,d)$-atoms supported,
respectively, in the balls $\{B_j\}
_{j\in\mathbb{N}}\subset\mathbb{B}$ such that
\begin{equation}\label{Molegxe1}
f=\sum_{j\in\mathbb{N}}\lambda_ja_j
\end{equation}
in $\mathcal{S}'(\rrn)$, and
\begin{equation}\label{Molegxe2}
\left\|\left[\sum_{j\in\mathbb{N}}
\left(\frac{\lambda_{j}}
{\|\1bf_{B_{j}}\|_{X}}
\right)^{s}\1bf_{B_{j}}\right]^{\frac{1}{s}}
\right\|_{X}<\infty.
\end{equation}
Moreover, by
Remark \ref{molexr}, we conclude that,
for any $j\in\mathbb{N}$, $a_j$
is a $(X,\,r,\,d,\,\tau)$-molecule
centered at the ball $B_j$.
This, combined with \eqref{Molegxe1} and
\eqref{Molegxe2}, then finishes the proof
of the necessity. In addition,
from the choice of $\{\lambda_j\}
_{j\in\mathbb{N}}$, Definition
\ref{atsgx}, and Theorem \ref{Atogx},
we deduce that
\begin{align}\label{Molegxe3}
\inf\left\{
\left\|\left[\sum_{j\in\mathbb{N}}
\left(\frac{\lambda_{j}}
{\|\1bf_{B_{j}}\|_{X}}
\right)^{s}\1bf_{B_{j}}
\right]^{\frac{1}{s}}
\right\|_{X}\right\}
\leq\left\|f\right\|
_{H^{X,r,d,s}(\rrn)}
\sim\|f\|_{H_X(\rrn)},
\end{align}
where the infimum is taken
over all the sequences $\{\lambda_j\}_{
j\in\mathbb{N}}\subset[0,\infty)$
and $\{m_j\}_{j\in\mathbb{N}}$
of $(X,\,r,\,d,\,\tau)$-molecules
centered, respectively, at
the balls $\{B_j\}_{j\in\mathbb{N}}
\subset\mathbb{B}$ such that
\eqref{Molegxe01} and \eqref{Molegxe02}
hold true.

Next, we show the sufficiency. To this end,
let $f\in\mathcal{S}'(\rrn)$ satisfy
$f=\sum_{j\in\mathbb{N}}
\lambda_jm_j$ in $\mathcal{S}'(\rrn)$,
where $\{\lambda_j\}_{j\in\mathbb{N}}
\subset[0,\infty)$ and $\{m_j\}
_{j\in\mathbb{N}}$ is a sequence of
$(X,\,r,\,d,\,\tau)$-molecules
centered, respectively, at the
balls $\{B_j\}_{j\in\mathbb{N}}
\subset\mathbb{B}$ such that
\begin{equation}\label{Molegxe11}
\left\|\left[\sum_{j\in\mathbb{N}}
\left(\frac{\lambda_{j}}
{\|\1bf_{B_{j}}\|_{X}}
\right)^{s}\1bf_{B_{j}}\right]^{\frac{1}{s}}
\right\|_{X}<\infty.
\end{equation}
To show the sufficiency,
we only need to prove that
$f\in H_X(\rrn)$. To achieve this,
we choose a $\phi\in\mathcal{S}(\rrn)$
satisfying that $\supp(\phi)\subset
B(\0bf,1)$ and $\int_{\rrn}\phi(x)
\,dx\neq0$.
Then, by Lemma \ref{Th5.4l1}, we find that,
to prove $f\in\ H_X(\rrn)$,
it suffices to show
$\|M(f,\phi)\|_{X}<\infty$.
Indeed, applying the assumption
that $f=\sum_{j\in\mathbb{N}}\lambda_j
m_j$ in $\mathcal{S}'(\rrn)$, and
an argument similar to
that used in the proof of \eqref{Atogxe1}
with $a_j$ therein replaced by
$m_j$ for any $j\in\mathbb{N}$, we conclude that
\begin{equation*}
M(f,\phi)\leq\sum_{j\in\mathbb{N}}
\lambda_jM(m_j,\phi).
\end{equation*}
From this, the assumption $d\geq
\lfloor n(1/\theta-1)\rfloor$,
and Lemma \ref{Molegl1}
with $m$ replaced by $m_j$ for any
$j\in\mathbb{N}$, we deduce
that
\begin{align}\label{Molegxe5}
\left\|M(f,\phi)\right\|_{X}
&\lesssim\left\|\sum_{j\in\mathbb{N}}
\lambda_j\mc(m_j)\1bf_{4B_j}\right\|_{X}
\notag\\&\quad+\left\|
\sum_{j\in\mathbb{N}}
\sum_{k=3}^{\infty}
\lambda_j\mc\left(m_j\1bf_{(2^{k+1}B_j)
\setminus(2^{k-2}B_j)}\right)
\1bf_{(2^kB_j)\setminus(2^{k-1}B_j)}
\right\|_{X}\notag\\
&\quad+\left\|\sum_{j\in\mathbb{N}}
\frac{\lambda_j}{\|
\1bf_{B_j}\|_{X}}
\mc^{(\theta)}(\1bf_{B_j})\right\|_{X}
\notag\\&=:\mathrm{III}_1
+\mathrm{III}_2+\mathrm{III}_3.
\end{align}

We first estimate $\mathrm{III}_1$.
Indeed, applying Definition \ref{molex}(i)
and an argument similar
to that used in the
estimation of $\mathrm{II}_1$ in
the proof of Theorem \ref{Atogx} with
$\mc(a_j)\1bf_{2B_j}$
therein replaced by
$\mc(m_j)\1bf_{4B_j}$ for any
$j\in\mathbb{N}$, we find that
\begin{equation}\label{Molegxe6}
\mathrm{III}_1\lesssim
\left\|\left[\sum_{j\in\mathbb{N}}
\left(\frac{\lambda_{j}}
{\|\1bf_{B_{j}}\|_{X}}
\right)^{s}\1bf_{B_{j}}
\right]^{\frac{1}{s}}
\right\|_{X}.
\end{equation}
This is the desired estimate
of $\mathrm{III}_1$.

Next, we deal with $\mathrm{III}_2$.
Indeed, from the $L^r(\rrn)$
boundedness of $\mc$ and Definition
\ref{molex}(i), it follows that,
for any $j\in\mathbb{N}$ and
$k\in\mathbb{N}\cap[3,\infty)$,
\begin{align*}
&\left\|\mc\left(
m_j\1bf_{(2^{k+1}B_j)\setminus
(2^{k-2}B_j)}\right)
\1bf_{(2^kB_j)\setminus
(2^{k-1}B_j)}\right\|_{L^r(\rrn)}\\
&\quad\lesssim
\left\|m_j\1bf_{(2^{k+1}B_j)\setminus
(2^{k-2}B_j)}\right\|_{L^r(\rrn)}
\lesssim\sum_{l=k-1}^{k+1}
\left\|m_j\1bf_{(2^{l}B_j)\setminus
(2^{l-1}B_j)}\right\|_{L^r(\rrn)}\\
&\quad\lesssim2^{-k\tau}\frac{|B_j|^{1/r}}
{\|\1bf_{B_j}\|_{X}}
\sim\left[2^{-k(\tau+\frac{n}{r})}
\frac{\|\1bf_{2^kB_j}\|_{X}}{
\|\1bf_{B_j}\|_{X}}\right]
\frac{|2^kB_j|^{1/r}}{\|\1bf_{2^kB_j}\|_{X}},
\end{align*}
where the implicit positive
constants are independent of both $j$ and $k$.
This implies that there exists
a positive constant $C$ such that,
for any $j\in\mathbb{N}$ and $k\in
\mathbb{N}\cap[3,\infty)$,
\begin{align}\label{Molegxe7}
&\left\|\mc\left(
m_j\1bf_{(2^{k+1}B_j)\setminus
(2^{k-2}B_j)}\right)
\1bf_{(2^kB_j)\setminus
(2^{k-1}B_j)}\right\|_{L^r(\rrn)}\notag\\
&\quad\leq C\left[2^{-k(\tau+\frac{n}{r})}
\frac{\|\1bf_{2^kB_j}\|_{X}}{
\|\1bf_{B_j}\|_{X}}\right]
\frac{|2^kB_j|^{1/r}}{\|\1bf_{2^kB_j}\|_{X}}.
\end{align}
For any $j\in\mathbb{N}$ and $k\in
\mathbb{N}\cap[3,\infty)$, let
$$
\mu_{j,k}:=C2^{-k(\tau+\frac{n}{r})}
\frac{\|\1bf_{2^kB_j}\|_{X}}{
\|\1bf_{B_j}\|_{X}}
$$
and
$$
a_{j,k}:=C^{-1}2^{k(\tau+\frac{n}{r})}
\frac{\|\1bf_{B_j}\|_{X}}{
\|\1bf_{2^kB_j}\|_{X}}\mc\left(
m_j\1bf_{(2^{k+1}B_j)\setminus
(2^{k-2}B_j)}\right)
\1bf_{(2^kB_j)\setminus
(2^{k-1}B_j)}.
$$
Then, by \eqref{Molegxe7}, we conclude that,
for any $j\in\mathbb{N}$ and
$k\in\mathbb{N}\cap[3,\infty)$,
$\supp(a_{j,k})\subset2^{k}B_j$ and
$$
\left\|a_{j,k}\right\|_{L^r(\rrn)}
\leq\frac{|2^kB_j|^{1/r}}{\|2^kB_j\|_{X}}.
$$
These, together with the
definitions of both $\mu_{j,k}$ and
$a_{j,k}$, and an
argument similar to that used in
the estimation of $\mathrm{II}_1$ in
the proof of Theorem \ref{Atogx} with
$\{\lambda_j\}_{j\in\mathbb{N}}$
and $\{\mc(a_j)\1bf_{2B_j}\}_{j\in\mathbb{N}}$
therein replaced, respectively, by
$$\left\{\lambda_j\mu_{j,k}
\right\}_{j\in\mathbb{N},k\in
\mathbb{N}\cap[3,\infty)}$$ and
$$\left\{a_{j,k}\1bf_{2^kB_j}\right\}
_{j\in\mathbb{N},k\in
\mathbb{N}\cap[3,\infty)},$$ further imply
that
\begin{align}\label{Molegxe8}
\mathrm{III}_2&\sim\left\|
\sum_{j\in\mathbb{N}}\sum_{k=3}^{\infty}
\lambda_j\mu_{j,k}a_{j,k}\right\|_{X}\notag
\\&\lesssim\left\|\left[
\sum_{j\in\mathbb{N}}\sum_{k=3}^{\infty}
\left(\frac{\lambda_j\mu_{j,k}}{\|\1bf_{2^kB_j}\|_{
X}}\right)^s\1bf_{2^kB_j}
\right]^{\frac{1}{s}}\right\|_{X}\notag\\
&\sim\left\|\left[\sum_{j\in\mathbb{N}}
\left(\frac{\lambda_j}{\|\1bf_{B_j}\|_{X}}
\right)^{s}\sum_{k=3}^{\infty}
2^{-ks(\tau+\frac{n}{r})}\1bf_{2^kB_j}
\right]^{\frac{1}{s}}\right\|_{X}.
\end{align}
In order to complete the estimation
of $\mathrm{III}_2$,
we now estimate the characteristic
function $\1bf_{2^kB_j}$ with $j\in\mathbb{N}$
and $k\in\mathbb{N}\cap[3,\infty)$.
Indeed, from \eqref{hlmaxp}, we deduce that,
for any $j\in\mathbb{N}$,
$k\in\mathbb{N}\cap[3,\infty)$, and $x\in2^kB_j$,
\begin{align*}
\left[\mc^{(\theta)}
\left(\1bf_{B_j}\right)(x)\right]^{s}
&\geq\left\{\frac{1}{
|2^kB_j|}\int_{2^kB_j}
\left[\1bf_{B_j}(y)\right]
^\theta\,dy\right\}^{\frac{s}{\theta}}\\
&\sim\left(\frac{|B_j|}{|2^kB_j|}
\right)^{\frac{s}{\theta}}\sim2^{-\frac{nks}
{\theta}},
\end{align*}
which implies that
$$
\1bf_{2^kB_j}\lesssim2^{\frac{nks}{\theta}}
\left[\mc^{(\theta)}
\left(\1bf_{B_j}\right)\right]^{s}.
$$
Combining this and the assumption
$\tau\in(n(\frac{1}{\theta}-
\frac{1}{r}),\infty)$, we further
conclude that
\begin{align*}
\sum_{k=3}^{\infty}2^{-ks(\tau+
\frac{n}{r})}\1bf_{2^kB_j}\lesssim
\sum_{k=3}^{\infty}2^{-ks(\tau-\frac{n}{
\theta}+\frac{n}{r})}\left[
\mc^{(\theta)}\left(\1bf_{
B_j}\right)\right]^{s}
\sim\left[\mc^{(\theta)}\left(\1bf_{
B_j}\right)\right]^{s}.
\end{align*}
Using this, \eqref{Molegxe8},
and Definition \ref{atsgx}(i) with
$f_j$ therein replaced by
$\frac{\lambda_j}{\|
\1bf_{B_j}\|_{X}}\1bf_{B_j}$
for any $j\in\mathbb{N}$, we find that
\begin{align}\label{Molegxe9}
\mathrm{III}_2&\lesssim
\left\|\left\{
\sum_{j\in\mathbb{N}}\left[
\mc^{(\theta)}\left(\frac{\lambda_j}{
\|\1bf_{B_j}\|_{X}}\1bf
_{B_j}\right)\right]^s
\right\}^{\frac{1}{s}}\right\|_{X}\notag\\
&\lesssim
\left\|\left[\sum_{j\in\mathbb{N}}
\left(\frac{\lambda_{j}}
{\|\1bf_{B_{j}}\|_{X}}
\right)^{s}\1bf_{B_{j}}
\right]^{\frac{1}{s}}
\right\|_{X},
\end{align}
which completes the estimation of
$\mathrm{III}_2$.

Finally, for the term $\mathrm{III}_3$,
from \eqref{Atogxe6}, we infer that
\begin{equation}\label{Molegxe10}
\mathrm{III}_3\lesssim
\left\|\left[\sum_{j\in\mathbb{N}}
\left(\frac{\lambda_{j}}
{\|\1bf_{B_{j}}\|_{X}}
\right)^{s}\1bf_{B_{j}}
\right]^{\frac{1}{s}}
\right\|_{X}.
\end{equation}
This is the desired estimate of
$\mathrm{III}_3$. Thus, combining \eqref{Molegxe5},
\eqref{Molegxe6}, \eqref{Molegxe9},
\eqref{Molegxe10}, and \eqref{Molegxe11},
we conclude that
\begin{equation}\label{Molegxe12}
\|M(f,\phi)\|_{X}
\lesssim\left\|\left[\sum_{j\in\mathbb{N}}
\left(\frac{\lambda_{j}}
{\|\1bf_{B_{j}}\|_{X}}
\right)^{s}\1bf_{B_{j}}
\right]^{\frac{1}{s}}
\right\|_{X}<\infty,
\end{equation}
which further implies that
$f\in H_X(\rrn)$ and hence completes
the proof of the sufficiency.
Moreover, from Lemma \ref{Th5.4l1},
\eqref{Molegxe12}, and the
choice of $\{\lambda_j\}_{
j\in\mathbb{N}}$,
it follows that
\begin{align*}
\|f\|_{H_X(\rrn)}&\sim
\left\|M(f,\phi)\right\|_{X}\\
&\lesssim\inf\left\{
\left\|\left[\sum_{j\in\mathbb{N}}
\left(\frac{\lambda_{j}}
{\|\1bf_{B_{j}}\|_{X}}
\right)^{s}\1bf_{B_{j}}
\right]^{\frac{1}{s}}
\right\|_{X}
\right\},
\end{align*}
where the infimum is taken
over all the sequences $\{\lambda_j\}_{
j\in\mathbb{N}}\subset[0,\infty)$
and $\{m_j\}_{j\in\mathbb{N}}$
of $(X,\,r,\,d,\,\tau)$-molecules
centered, respectively, at
the balls $\{B_j\}_{j\in\mathbb{N}}
\subset\mathbb{B}$ such that
\eqref{Molegxe01} and \eqref{Molegxe02}
hold true.
This, together with \eqref{Molegxe3},
implies that \eqref{Molegxe03}
holds true, which completes the
proof of Theorem \ref{Molegx}.
\end{proof}

Via the above molecular
characterization of $H_X(\rrn)$,
we next prove Theorem \ref{Moleg}.

\begin{proof}[Proof of Theorem
\ref{Moleg}]
Let $p$, $q$, $\omega$, $r$, $d$,
$s$, and $\tau$ be as in the present theorem.
Then, combining the assumptions,
$\m0(\omega)\in(-\frac{n}{p},\infty)$
and $\MI(\omega)\in(-\infty,0)$,
and Theorem \ref{Th2}, we find
that the global generalized Herz space
$\HerzS$ under consideration
is a BQBF space. Therefore,
in order to show the present theorem,
we only need to prove that
all the assumptions of Theorem \ref{Molegx}
are satisfied.

Indeed, let $\theta\in(0,s)$ be such that
$$
\left\lfloor n\left(\frac{1}{s}-1\right)
\right\rfloor\leq n\left(\frac{1}{\theta}
-1\right)<\left\lfloor n\left(
\frac{1}{s}-1\right)
\right\rfloor+1
$$
and
$$
\tau>n\left(\frac{1}{\theta}-
\frac{1}{r}\right).
$$
Thus, we have $d\geq\lfloor n(1/s-1)\rfloor
=\lfloor n(1/\theta-1)\rfloor$.
Then, from
\eqref{atoge1}, \eqref{atoge3}, and
\eqref{atoge4},
it follows that the following three
statements hold true:
\begin{enumerate}
  \item[{\rm(i)}] for any $\{f_j\}_{
j\in\mathbb{N}}\subset L^1_{\mathrm{loc}}(
\rrn)$,
\begin{equation*}
\left\|\left\{
\sum_{j\in\mathbb{N}}
\left[
\mc^{(\theta)}(f_j)
\right]^{s}\right\}^{1/s}\right\|_{\HerzS}
\lesssim\left\|\left(\sum_{j\in\mathbb{N}}
|f_j|^s\right)^{1/s}\right\|_{\HerzS};
\end{equation*}
  \item[{\rm(ii)}] for any $f\in\Msc(\rrn)$,
\begin{equation*}
\|f\|_{[\HerzS]^{1/s}}\sim
\sup\left\{\|fg\|_{L^1(\rrn)}:\
\|g\|_{\dot{\mathcal{B}}_{1/\omega^s}
^{(p/s)',(q/s)'}(\rrn)}=1\right\};
\end{equation*}
  \item[{\rm(iii)}] for any
$f\in L_{{\rm loc}}^{1}(\rrn)$,
\begin{equation*}
\left\|\mc^{((r/s)')}(f)\right\|_{\dot{\mathcal{B}}
_{1/\omega^{s}}^{(p/s)',(q/s)'}(\rrn)}\lesssim
\|f\|_{\dot{\mathcal{B}}_
{1/\omega^{s}}^{(p/s)',(q/s)'}(\rrn)}.
\end{equation*}
\end{enumerate}
These further imply that,
under the assumptions of the present theorem,
$\HerzS$ satisfies all the assumptions of
Theorem \ref{Molegx},
which completes the proof of Theorem
\ref{Moleg}.
\end{proof}

Similarly, we now characterize
the generalized Morrey--Hardy
space $H\MorrS$ via molecules.
To this end, we first introduce
the following
$(\MorrS,\,r,\,d,\,\tau)$-molecules.

\begin{definition}
Let $p$, $q$, $\omega$,
$r$, $d$, and $\tau$ be as
in Definition \ref{molem}.
Then a measurable function
$m$ on $\rrn$ is called an
\emph{$(\MorrS,\,r,\,d,\,
\tau)$-molecule}\index{$(\MorrS,\,r,\,d,\,\tau)$-molecule}
centered at a ball $B\in\mathbb{B}$ if
\begin{enumerate}
  \item[{\rm(i)}] for any $j\in\zp$,
  $$\left\|m\1bf_{S_{j}(B)}
  \right\|_{L^{r}(\rrn)}\leq2^{-\tau j}
\frac{|B|^{1/r}}{\|\1bf_{B}\|_{\MorrS}};$$
 \item[{\rm(ii)}] for any $\alpha\in\zp^n$
 with $|\alpha|\leq d$,
  $$\int_{\rrn}m(x)x^\alpha\,dx=0.$$
\end{enumerate}
\end{definition}

Then we immediately obtain
the following
molecular characterization of the generalized
Morrey--Hardy space $H\MorrS$ via
Theorem \ref{Moleg} and Remarks \ref{remhs}(iv)
and \ref{remark4.10}(ii);
we omit the details.

\begin{corollary}
Let $p$, $q$, $\omega$, $r$,
$d$, and $s$ be as in
Corollary \ref{cor438} and
$\tau\in(0,\infty)$ with
$\tau>n(1/s-1/r)$.
Then $f\in H\MorrS$ if and only
if $f\in\mathcal{S}'(\rrn)$
and there exist a sequence
$\{m_{j}\}_{j\in\mathbb{N}}$
of $(\MorrS,\,r,\,d,\,\tau)$-molecules
centered, respectively,
at the balls $\{B_{j}\}_{j\in\mathbb{N}}
\subset\mathbb{B}$ and
a sequence $\{\lambda_{j}\}
_{j\in\mathbb{N}}\subset
[0,\infty)$ such that
$$
f=\sum_{j\in\mathbb{N}}\lambda_{j}a_{j}
$$
in $\mathcal{S}'(\rrn)$, and
$$
\left\|\left\{\sum_{j\in\mathbb{N}}
\left[\frac{\lambda_{j}}
{\|\1bf_{B_{j}}\|_{\MorrS}}\right]
^{s}\1bf_{B_{j}}\right\}^{\frac{1}{s}}
\right\|_{\MorrS}<\infty.
$$
Moreover, for any $f\in H\MorrS$,
\begin{align*}
\|f\|_{H\MorrS}\sim\inf\left\{
\left\|\left\{\sum_{j\in\mathbb{N}}
\left[\frac{\lambda_{i}}
{\|\1bf_{B_{j}}\|_{\MorrS}}\right]^s\1bf_
{B_{j}}\right\}^{\frac{1}{s}}
\right\|_{\MorrS}
\right\}
\end{align*}
with the positive equivalence
constants independent of $f$,
where the infimum is taken
over all the decompositions of $f$ as above.
\end{corollary}

\section{Littlewood--Paley Function
Characterizations}\label{sec4.5}

In this section, we establish
various Littlewood--Paley
function characterizations of
generalized Herz--Hardy spaces.
To be precise, we characterize generalized
Herz--Hardy spaces
via the Lusin area function,
the Littlewood--Paley $g$-function,
and the Littlewood--Paley $g^*_{\lambda}$-function.

Throughout this book, for any
$\varphi\in\mathcal{S}(\rrn)$, $\mathcal{F}\varphi$
or $\widehat{\varphi}$ always denotes its
\emph{Fourier transform}\index{Fourier transform},
namely, for any $x\in\rrn$,
$$\index{$\widehat{\varphi}$}\index{$\mathcal{F}$}
\mathcal{F}\varphi(x):=\widehat{\varphi}(x):=\int_{
\rrn}\varphi(\xi)e^{-2\pi ix\cdot\xi}\,d\xi,
$$
where $i:=\sqrt{-1}$ and
$x\cdot\xi:=\sum_{j=1}^{n}x_{j}\xi_{j}$
for any $x:=(x_{1},\ldots,x_{n})$, $\xi:=(\xi_{1},\ldots,
\xi_{n})\in\rrn$. Moreover,
for any $f\in\mathcal{S}'(\rrn)$,
its \emph{Fourier transform},
also denoted by $\mathcal{F}f$ or
$\widehat{f}$, is defined by setting,
for any $\phi\in\mathcal{S}(\rrn)$,
$$
\langle\mathcal{F}f,\phi\rangle:=
\langle\widehat{f},\phi\rangle
:=\langle f,\widehat{\phi}\rangle.
$$
Then we state
the definitions of the
Lusin area function, the
Littlewood--Paley $g$-function,
and the Littlewood--Paley
$g^*_{\lambda}$-function as follows
(see, for instance,
\cite[Definitions 4.1 and 4.2]{CWYZ}).

\begin{definition}\label{df451}
Let $\varphi\in\mathcal{S}(\rrn)$ satisfy
$\widehat{\varphi}(\0bf)=0$ and,
for any $\xi\in\rrn\setminus\{\0bf\}$,
there exists a $t\in(0,\infty)$
such that $\widehat{\varphi}(t\xi)\neq0$.
Then, for any $f\in\mathcal{S}'(\rrn)$,
the \emph{Lusin area function}\index{Lusin area function}
$S(f)$\index{$S(f)$} and the
\emph{Littlewood--Paley
$g^{*}_{\lambda}$-function}\index{Littlewood--Paley $g^{*}_{\lambda}$-function}
$g^{*}_{\lambda}(f)$\index{$g^{*}_{\lambda}$} with
$\lambda\in(0,\infty)$ are defined,
respectively, by setting, for any $x\in\rrn$,
$$S(f)(x):=\left[\int_{\Gamma(x)}|
f\ast\varphi_t(y)|^{2}
\,\frac{dy\,dt}{t^{n+1}}\right]^{\frac{1}{2}}$$
and
$$
g_{\lambda}^{*}(f)(x):=\left[\int_{0}
^{\infty}\int_{\rrn}
\left(\frac{t}{t+|x-y|}\right)^{\lambda n}|
f\ast\varphi_{t}(y)|^{2}
\,\frac{dy\,dt}{t^{n+1}}\right]^{\frac{1}{2}},
$$
here and thereafter, for any $x\in\rrn$,
$$\index{$\Gamma(x)$}
\Gamma(x):=\{(y,t)\in\mathbb{R}^{n+1}_{+}:\ |y-x|<t\}.
$$
\end{definition}

\begin{definition}\label{df452}
Let $\varphi\in\mathcal{S}(\rrn)$
satisfy $\widehat{\varphi}(\0bf)=0$
and, for any $x\in\rrn\setminus\{\0bf\}$,
there exists a $j\in\mathbb{Z}$
such that $\widehat{\varphi}(2^{j}x)\neq0$.
Then, for any $f\in\mathcal{S}'(\rrn)$, the
\emph{Littlewood--Paley $g$-function}\index{Littlewood--Paley $g$-function}
$g(f)$\index{$g(f)$} is defined by setting,
for any $x\in\rrn$,
$$
g(f)(x):=\left[\int_{0}^{\infty}|f\ast\varphi_{t}(
x)|^{2}\,\frac{dt}{t}\right]^{\frac{1}{2}}.
$$
\end{definition}

Recall that
$f\in\mathcal{S}'(\rrn)$ is said to
\emph{vanish weakly at
infinity}\index{vanish weakly at infinity}
if, for any $\phi\in\mathcal{S}(\rrn)$,
$f\ast\phi_{t}\to0$ in
$\mathcal{S}'(\rrn)$ as $t\to\infty$
with $\phi_t(\cdot)=t^{-n}\phi(\cdot/t)$.
Now, we establish various
Littlewood--Paley function
characterizations
of the generalized Herz--Hardy spaces
$\HaSaHo$ and $\HaSaH$ as follows. To
begin with, we give the Littlewood--Paley
function characterizations\index{Littlewood--Paley function characterization}
of $\HaSaHo$.

\begin{theorem}\label{lusin}
Let $p,\ q\in(0,\infty)$, $\omega\in
M(\rp)$ satisfy
$\m0(\omega)\in(-\frac{n}{p},\infty)$ and
$\mi(\omega)\in(-\frac{n}{p},\infty)$,
$$
s_{0}:=\min\left\{1,p,q,\frac{n}{\Mw+n/p}\right\},
$$
and $\lambda\in(\max\{1,2/s_{0}\},\infty)$.
Then the following four statements are
mutually equivalent:
\begin{enumerate}
  \item[{\rm(i)}] $f\in\HaSaHo$;
  \item[{\rm(ii)}] $f\in\mathcal{S}'(\rrn)$,
  $f$ vanishes weakly
  at infinity, and $S(f)\in\HerzSo$;
  \item[{\rm(iii)}] $f\in\mathcal{S}'(\rrn)$,
  $f$ vanishes weakly
  at infinity, and $g(f)\in\HerzSo$;
  \item[{\rm(iv)}] $f\in\mathcal{S}'(\rrn)$,
  $f$ vanishes weakly
  at infinity, and $g_{\lambda}^*(f)\in\HerzSo$.
\end{enumerate}
Moreover, for any $f\in\HaSaHo$,
$$
\|f\|_{\HaSaHo}\sim\|S(f)\|_{\HerzSo}
\sim\|g(f)\|_{\HerzSo}\sim
\left\|g_{\lambda}^*(f)\right\|_{\HerzSo},
$$
where the positive equivalence
constants are independent of $f$.
\end{theorem}

To show Theorem \ref{lusin},
recall that Chang et al.\
\cite[Theorems 4.9, 4.11, and 4.13]{CWYZ}
investigated the Lusin area function,
the $g$-function, and the
$g_{\lambda}^*$-function characterizations
of the Hardy space $H_X(\rrn)$ as follows, which
is vital in the proof of the
Littlewood--Paley function characterizations
of $\HaSaHo$.

\begin{lemma}\label{lusinl1}
Let $s\in(0,1]$, $\theta\in(0,s)$, $\lambda
\in(\max\{1,\frac{2}{s}\},\infty)$,
and $X$ be a ball quasi-Banach
function space.
Assume that Assumption \ref{assfs} holds true
for both $X$ and $X^{s/2}$ with
$\theta$ and $s$ as above and
$X$ satisfies Assumption \ref{assas}
with $s$ as above.
Then the following four statements are mutually
equivalent:
\begin{enumerate}
  \item[{\rm(i)}] $f\in H_X(\rrn)$;
  \item[{\rm(ii)}] $f\in\mathcal{S}'(\rrn)$,
  $f$ vanishes weakly
  at infinity, and $S(f)\in X$;
  \item[{\rm(iii)}] $f\in\mathcal{S}'(\rrn)$,
  $f$ vanishes weakly
  at infinity, and $g(f)\in X$;
  \item[{\rm(iv)}] $f\in\mathcal{S}'(\rrn)$,
  $f$ vanishes weakly
  at infinity, and $g_{\lambda}^*(f)\in X$.
\end{enumerate}
Moreover, for any $f\in H_X(\rrn)$,
$$
\|f\|_{H_X(\rrn)}\sim\|S(f)\|_{X}
\sim\|g(f)\|_{X}\sim
\left\|g_{\lambda}^*(f)\right\|_{X}
$$
with the positive equivalence
constants independent of $f$.
\end{lemma}

Via this lemma, we now show
Theorem \ref{lusin}.

\begin{proof}[Proof of
Theorem \ref{lusin}]
Let all the symbols be as in
the present theorem.
Then, combining the
assumption $\m0(\omega)
\in(-\frac{n}{p},\infty)$ and
Theorem \ref{Th3}, we find that
the local generalized Herz
space $\HerzSo$ is a BQBF space.
This implies that, to complete
the proof of the present theorem, we only
need to show that,
under the assumptions of the present
theorem, $\HerzSo$ satisfies all the
assumptions of Lemma \ref{lusinl1}.

For this purpose, let
$s\in(\frac{2}{\lambda},s_0)$ and
$\theta\in(0,\min\{s,\frac{s^2}{2}\})$.
We now prove that both
$\HerzSo$ and $[\HerzSo]^{s/2}$
satisfy Assumption \ref{assfs} with these $\theta$
and $s$.
Indeed, applying Lemma \ref{Atoll3},
we conclude that,
for any $\{f_{j}\}_{j\in\mathbb{N}}
\subset L^{1}_{\text{loc}}(\rrn)$,
\begin{equation*}
\left\|\left\{\sum_{j\in\mathbb{N}}
\left[\mc^{(\theta)}
(f_{j})\right]^s\right\}^{1/s}
\right\|_{\HerzSo}\lesssim
\left\|\left(\sum_{j\in\mathbb{N}}
|f_{j}|^{s}\right)^{1/s}\right\|_{\HerzSo}.
\end{equation*}
This implies that,
for the above $\theta$ and $s$,
Assumption \ref{assfs} holds true
for $\HerzSo$.
On the other hand,
from the assumptions
$\theta\in(0,\min\{s,\frac{s^{2}}{2}\})$ and
$$s\in\left(0,\min\left\{p,\frac{n}
{\Mw+n/p}\right\}\right),$$
we deduce that $\frac{s}{\theta}\in(1,\infty)$ and
\begin{align*}
\min\left\{p,\frac{n}{\Mw+n/p}
\right\}>s>\frac{2\theta}{s}.
\end{align*}
This, together with
Lemma \ref{Atoll3} with
$\HerzSo$ therein
replaced by $[\HerzSo]^{s/2}$,
further implies that,
for any $\{f_{j}\}_{j\in\mathbb{N}}
\subset L^1_{\mathrm{loc}}(\rrn)$,
\begin{equation*}
\left\|\left\{\sum_{j\in\mathbb{N}}
\left[\mc^{(\theta)}
(f_{j})\right]^s\right\}^{1/s}
\right\|_{[\HerzSo]^{s/2}}\lesssim
\left\|\left(\sum_{j\in
\mathbb{N}}|f_{j}|^{s}\right)^{1/s}\right\|
_{[\HerzSo]^{s/2}},
\end{equation*}
which completes the proof
that Assumption \ref{assfs} holds true
for both $\HerzSo$ and $[\HerzSo]^{s/2}$.

In addition,
notice that $\lambda\in
(\max\{1,\frac{2}{s}\},\infty)$
and we can choose an
$$r\in\left(\max\left\{1,p,\frac{n}
{\mw+n/p}\right\},\infty\right].$$
Then, by Lemma \ref{mbhal}, we
find that $[\HerzSo]^{1/s}$ is a BBF space
and, for any $f\in L_{{\rm loc}}^{1}(\rrn)$,
\begin{equation*}
\left\|\mc^{((r/s)')}(f)
\right\|_{([\HerzSo]^{1/s})'}\lesssim
\left\|f\right\|_{([\HerzSo]^{1/s})'},
\end{equation*}
which further implies that
all the assumptions of Lemma \ref{lusinl1}
are satisfied for $\HerzSo$.
This finishes the proof of Theorem \ref{lusin}.
\end{proof}

\begin{remark}
We point out that,
in Theorem \ref{lusin}, when $p\in(1,\infty)$,
$q\in[1,\infty)$, and $\omega(t):=t^\alpha$
for any $t\in(0,\infty)$ and for any given
$\alpha\in(-\frac{n}{p},\frac{n}{p})$,
from Remark \ref{remark4.0.8} and Theorem \ref{Th5.7},
it follows that
the generalized Herz--Hardy
space $\HaSaHo$ coincides with the
classical homogeneous Herz space
$\dot{K}_p^{\alpha,q}(\rrn)$ and hence
the conclusion obtained in this theorem
goes back to \cite[Theorem 1.1.1]{LYH}.
\end{remark}

Using Theorem \ref{lusin} and Remarks \ref{remhs}(iv)
and \ref{remark4.10}(ii), we immediately obtain
the following
Littlewood--Paley function characterizations
of the generalized Morrey--Hardy space $H\MorrSo$; we
omit the details.

\begin{corollary}\label{lusinm}
Let $p,\ q\in[1,\infty)$, $\omega\in
M(\rp)$ satisfy
$$-\frac{n}{p}<
\m0(\omega)\leq\M0(\omega)<0$$ and
$$-\frac{n}{p}<
\mi(\omega)\leq\MI(\omega)<0,$$
and $\lambda\in(2,\infty)$.
Then the following four statements are
mutually equivalent:
\begin{enumerate}
  \item[{\rm(i)}] $f\in H\MorrSo$;
  \item[{\rm(ii)}] $f\in\mathcal{S}'(\rrn)$,
  $f$ vanishes weakly
  at infinity, and $S(f)\in\MorrSo$;
  \item[{\rm(iii)}] $f\in\mathcal{S}'(\rrn)$,
  $f$ vanishes weakly
  at infinity, and $g(f)\in\MorrSo$;
  \item[{\rm(iv)}] $f\in\mathcal{S}'(\rrn)$,
  $f$ vanishes weakly
  at infinity, and $g_{\lambda}^*(f)\in\MorrSo$.
\end{enumerate}
Moreover, for any $f\in H\MorrSo$,
$$
\|f\|_{H\MorrSo}\sim\|S(f)\|_{\MorrSo}
\sim\|g(f)\|_{\MorrSo}\sim
\left\|g_{\lambda}^*(f)\right\|_{\MorrSo},
$$
where the positive equivalence
constants are independent of $f$.
\end{corollary}

On the other hand, we now show the following
Littlewood--Paley
function characterizations of the
generalized Herz--Hardy space $\HaSaH$.

\begin{theorem}\label{lusing}
Let $p,\ q\in(0,\infty)$, $\omega\in
M(\rp)$ satisfy
$\m0(\omega)\in(-\frac{n}{p},\infty)$ and
$$
-\frac{n}{p}<\mi(\omega)\leq\MI(\omega)<0,
$$
$$
s_{0}:=\min\left\{1,p,q,\frac{n}{\Mw+n/p}\right\},
$$
and $\lambda\in(\max\{1,2/s_{0}\},\infty)$.
Then the following four statements are
mutually equivalent:
\begin{enumerate}
  \item[{\rm(i)}] $f\in\HaSaH$;
  \item[{\rm(ii)}] $f\in\mathcal{S}'(\rrn)$,
  $f$ vanishes weakly
  at infinity, and $S(f)\in\HerzS$;
  \item[{\rm(iii)}] $f\in\mathcal{S}'(\rrn)$,
  $f$ vanishes weakly
  at infinity, and $g(f)\in\HerzS$;
  \item[{\rm(iv)}] $f\in\mathcal{S}'(\rrn)$,
  $f$ vanishes weakly
  at infinity, and $g_{\lambda}^*(f)\in\HerzS$.
\end{enumerate}
Moreover, for any $f\in\HaSaH$,
$$
\|f\|_{\HaSaH}\sim\|S(f)\|_{\HerzS}
\sim\|g(f)\|_{\HerzS}\sim
\left\|g_{\lambda}^*(f)\right\|_{\HerzS},
$$
where the positive equivalence
constants are independent of $f$.
\end{theorem}

To prove Theorem \ref{lusing}, we first
give some symbols.
Recall that, for any given
$\xi\in\rrn$, the \emph{translation operator}
\index{translation operator}$\tau_{\xi}$
\index{$\tau_{\xi}$}is defined
by setting, for any $f\in\Msc(\rrn)$
and $x\in\rrn$,
\begin{equation}\label{translation}
\tau_{\xi}(f)(x):=f(x-\xi).
\end{equation}
Furthermore, the \emph{translation
operator} $\tau_{\xi}$ of distributions with
$\xi\in\rrn$ is defined by setting,
for any $f\in\mathcal{S}'(\rrn)$ and
$\phi\in\mathcal{S}(\rrn)$,
$$
\left\langle\tau_{\xi}(f),\phi
\right\rangle:=
\left\langle f,\tau_{-\xi}(\phi)
\right\rangle
=\left\langle f,\phi(\cdot+\xi)
\right\rangle.
$$
The following
technique lemma
establishes the relations
among translations, convolutions,
and various Littlewood--Paley
functions, which plays a key role
in the proof of Theorem \ref{lusing}.

\begin{lemma}\label{lusingl1}
Let $f\in\mathcal{S}'(\rrn)$,
$\phi\in\mathcal{S}(\rrn)$, and $\xi\in\rrn$.
Then
\begin{enumerate}
  \item[{\rm(i)}] $[\tau_{\xi}(f)]\ast\phi=
  \tau_{\xi}(f\ast\phi)$;
  \item[{\rm(ii)}] $
  M(\tau_{\xi}(f),\phi)=
  \tau_{\xi}(M(f,\phi)),
  $
  where the radial maximal function $M$
  is defined as
  in Definition \ref{smax}(i);
  \item[{\rm(iii)}] $A(
  \tau_{\xi}(f))=\tau_{\xi}(A(f))$,
  where $A\in\{S,g,g_{\lambda}^*\}$ with
  $\lambda\in(0,\infty)$;
  \item[{\rm(iv)}] if $f$ vanishes
  weakly at infinity, then $\tau_{\xi}
  (f)$ vanishes weakly at infinity.
\end{enumerate}
\end{lemma}

\begin{proof}
Let all the symbols be as in the
present lemma. We first prove (i).
Indeed, for any $x\in\rrn$, we have
\begin{align*}
\left[\tau_{\xi}(f)\right]\ast
\phi(x)&=\left\langle
\tau_{\xi}(f),\phi(x-\cdot)\right\rangle
=\left\langle f,\phi\left(x-(\cdot+\xi)
\right)\right\rangle\\
&=\left\langle
f,\phi(x-\xi-\cdot)\right\rangle
=\left(f\ast\phi\right)
(x-\xi)\\&=\tau_{\xi}(f\ast\phi)
(x),
\end{align*}
which implies that (i) holds true.

Next, we show (ii). Applying (i)
with $\phi$ therein replaced
by $\phi_t$ for any $t\in(0,\infty)$,
we find that
\begin{align*}
M\left(\tau_{\xi}(f),\phi\right)
&=\sup_{t\in(0,\infty)}
\left\{\left|
\left[\tau_{\xi}(f)\right]
\ast\phi_{t}\right|\right\}=\sup_{t\in(0,\infty)}
\left\{\left|\tau_{\xi}
(f\ast\phi_t)\right|\right\}\\
&=\tau_{\xi}\left(\sup_{t\in(0,\infty)}
\left\{|f\ast\phi_t|\right\}\right)
=\tau_{\xi}\left(M(f,\phi)\right).
\end{align*}
This finishes the proof of (ii).

We then prove (iii). First, let
$A:=S$. Then, from Definition \ref{df451}
and (i) with $\phi=\varphi_t$ for any
$t\in(0,\infty)$, it follows that,
for any $x\in\rrn$,
\begin{align*}
\left[S\left(\tau_{\xi}(f)\right)
(x)\right]^2&=
\int_{\Gamma(x)}\left|
\left[\tau_{\xi}(f)\right]\ast\varphi_{t}
(y)\right|^2\,\frac{dy\,dt}{t^{n+1}}\\
&=\int_{\Gamma(x)}\left|
(f\ast\varphi_{t})
(y-\xi)\right|^2\,\frac{dy\,dt}{t^{n+1}}\\
&=\int_{\Gamma(x-\xi)}\left|
f\ast\varphi_{t}(y)\right|^2\,\frac{dy\,dt}{t^{n+1}}
=\left[\tau_{\xi}(S(f))(x)\right]^2.
\end{align*}
Thus,
$S(\tau_{\xi}(f))=\tau_{\xi}(S(f))$
holds true. Similarly, we can
obtain $g(\tau_{\xi}(f))=\tau_{\xi}
(g(f))$ and $g_{
\lambda}^*(\tau_{\xi}(f))=\tau_{\xi}
(g_{\lambda}^*(f))$
and hence (iii) holds true.

Finally, we prove (iv). Indeed, for
any $\psi$, using (i) with $\phi$ therein
replaced by $\psi_t$ for any $t\in(0,\infty)$,
and the assumption that
$f$ vanishes weakly at infinity,
we find that, for any $\eta\in\mathcal{S}(\rrn)$,
\begin{align*}
\int_{\rrn}
\left[\tau_{\xi}(f)\right]
\ast\psi_t(y)\eta(y)\,dy
&=\int_{\rrn}\left(f\ast\psi_{t}\right)
(y-\xi)\eta(y)\,dy\\
&=\int_{\rrn}\left(f\ast\psi_{t}\right)
(y)\eta(y+\xi)\,dy\to0
\end{align*}
as $t\to\infty$. This
further implies that, for
any $\psi\in\mathcal{S}(\rrn)$,
$[\tau_{\xi}(f)]\ast\psi\to0$
in $\mathcal{S}'(\rrn)$ as $t\to\infty$,
which completes the proof of (iv),
and hence of Lemma \ref{lusingl1}.
\end{proof}

Now, we show Theorem \ref{lusing} via
Lemma \ref{lusingl1}.

\begin{proof}[Proof of Theorem
\ref{lusing}]
Let $p$, $q$, $\omega$, and
$\lambda$ be as in the
present theorem, and let
$\phi\in\mathcal{S}(\rrn)$ satisfy that
$\int_{\rrn}\phi(x)\,dx\neq0$.
We now first show that (i)
implies (ii). To this end, let
$f\in\HaSaH$. Then, applying this and
Theorem \ref{Th5.6},
we find that $M(f,\phi)\in\HerzS$.
This, together with
Lemma \ref{lusingl1}(ii),
further implies that,
for any $\xi\in\rrn$,
\begin{align*}
\left\|M\left(\tau_{\xi}
(f),\phi\right)\right\|_{\HerzSo}
&=\left\|\tau_{\xi}\left(M(f,\phi)
\right)\right\|_{\HerzSo}\\&\leq\left\|
M(f,\phi)\right\|_{\HerzS}
\sim\|f\|_{\HaSaH}<\infty.
\end{align*}
Therefore, for any $\xi\in\rrn$,
$M(\tau_{\xi}(f),\phi)\in\HerzSo$.
From this and Theorem \ref{Th5.4}, we
deduce that, for any $\xi\in\rrn$,
$\tau_{\xi}(f)\in\HaSaHo$.
Combining this and Theorem \ref{lusin},
we conclude that,
for any $\xi\in\rrn$,
$\tau_{\xi}(f)$ vanishes weakly at infinity,
$S(\tau_{\xi}(f))\in\HerzSo$, and
\begin{equation}\label{lusinge1}
\left\|S\left(\tau_{\xi}(f)\right)
\right\|_{\HerzSo}\sim\left\|
\tau_{\xi}(f)\right\|_{\HaSaHo}.
\end{equation}
In particular, letting $\xi:=\0bf$, we have
$f$ vanishes weakly at infinity. Then,
from Lemma \ref{lusingl1}(iii), \eqref{lusinge1},
Theorem \ref{Th5.4}, and Lemma \ref{lusingl1}(ii),
we further deduce that,
for any $\xi\in\rrn$,
\begin{align*}
\left\|\tau_{\xi}\left(
S(f)\right)\right\|_{\HerzSo}
&=\left\|S\left(\tau_{\xi}(f)\right)
\right\|_{\HerzSo}\\
&\sim\left\|\tau_{\xi}(f)\right\|
_{\HaSaHo}\sim
\left\|M\left(\tau_{\xi}(f),
\phi\right)\right\|_{\HerzSo}\\
&\sim\left\|\tau_{\xi}\left(
M(f,\phi)\right)\right\|_{\HerzSo}.
\end{align*}
By this, the definition of $\|
\cdot\|_{\HerzS}$, and Theorem \ref{Th5.6},
we find that
\begin{equation}\label{lusinge2}
\|S(f)\|_{\HerzS}\sim\|M(f,\phi)\|_{\HerzS}
\sim\|f\|_{\HaSaH}<\infty,
\end{equation}
which further implies that $S(f)\in\HerzS$,
and hence (i) implies (ii).

Conversely, we show that (ii) implies (i),
namely, assume $f\in\mathcal{S}'(\rrn)$
vanishes weakly at infinity
and $S(f)\in\HerzS$, we need to prove
that $f\in\HaSaH$.
Indeed, from Lemma \ref{lusingl1}(iii),
we deduce that, for any $\xi\in\rrn$,
\begin{align}\label{lusinge3}
\left\|S\left(\tau_{\xi}(f)\right)
\right\|_{\HerzSo}=
\left\|\tau_{\xi}
\left(S(f)\right)\right\|_{\HerzSo}
\leq\|S(f)\|_{\HerzS}<\infty,
\end{align}
which implies that $S(\tau_{\xi}(f))
\in\HerzSo$. On the other hand,
by Lemma \ref{lusingl1}(iv), we conclude that,
for any $\xi\in\rrn$,
$\tau_{\xi}(f)$ vanishes weakly at infinity.
Combining this, the fact that $S(\tau_{\xi}(f))
\in\HerzSo$ for any $\xi\in\rrn$, and
Theorem \ref{lusin}, we further conclude
that, for any $\xi\in\rrn$, $\tau_{\xi}(f)
\in\HaSaHo$. Then, applying Theorem \ref{Th5.6},
Lemma \ref{lusingl1}(ii), Theorems \ref{Th5.4}
and \ref{lusin}, and \eqref{lusinge3}, we find that
\begin{align*}
\|f\|_{\HaSaH}&\sim\sup_{\xi\in\rrn}
\left\|\tau_{\xi}\left(
M(f,\phi)\right)\right\|_{\HerzSo}\\
&\sim\sup_{\xi\in\rrn}\left\|
M\left(\tau_{\xi}(f),\phi\right)\right\|_{\HerzSo}
\sim\sup_{\xi\in\rrn}\left\|\tau_{\xi}
(f)\right\|_{\HaSaHo}\\
&\sim\sup_{\xi\in\rrn}\left\|
S\left(\tau_{\xi}(f)\right)\right\|_{\HerzSo}
\lesssim\|S(f)\|_{\HerzS}<\infty,
\end{align*}
which further implies that $f\in\HaSaH$, and
hence (ii) implies (i). Moreover, from
\eqref{lusinge2}, it follows that, for any
$f\in\HaSaH$,
$$
\|f\|_{\HaSaH}\sim\|S(f)\|_{\HerzS}.
$$

Similarly, we can obtain (i)
is equivalent to both (iii)
and (iv) and, for any $f\in\HaSaH$,
$$
\|f\|_{\HaSaH}\sim\|A(f)\|_{\HerzS},
$$
where $A\in\{g,g_{\lambda}^*\}$. This
then finishes the proof of Theorem \ref{lusing}.
\end{proof}

As an application, we then characterize
the generalized Morrey--Hardy space
$H\MorrS$ via the Lusin area function,
the $g$-function, and the $g_{
\lambda}^*$-function. Namely,
the following conclusion holds true,
which is a direct corollary of
Theorem \ref{lusing} and Remarks \ref{remhs}(iv)
and \ref{remark4.10}(ii); we
omit the details.

\begin{corollary}
Let $p$, $q$, $\omega$, and $\lambda$
be as in Corollary \ref{lusinm}.
Then the following four statements are
mutually equivalent:
\begin{enumerate}
  \item[{\rm(i)}] $f\in H\MorrS$;
  \item[{\rm(ii)}] $f\in\mathcal{S}'(\rrn)$,
  $f$ vanishes weakly
  at infinity, and $S(f)\in\MorrS$;
  \item[{\rm(iii)}] $f\in\mathcal{S}'(\rrn)$,
  $f$ vanishes weakly
  at infinity, and $g(f)\in\MorrS$;
  \item[{\rm(iv)}] $f\in\mathcal{S}'(\rrn)$,
  $f$ vanishes weakly
  at infinity, and $g_{\lambda}^*(f)\in\MorrS$.
\end{enumerate}
Moreover, for any $f\in H\MorrS$,
$$
\|f\|_{H\MorrS}\sim\|S(f)\|_{\MorrS}
\sim\|g(f)\|_{\MorrS}\sim
\left\|g_{\lambda}^*(f)\right\|_{\MorrS}
$$
with the positive equivalence
constants independent of $f$.
\end{corollary}

\section{Dual Space of $\HaSaHo$}

In this section, we investigate
the dual space of the
generalized Herz--Hardy space $\HaSaHo$.
Throughout this book, for any
$d\in\zp$, the \emph{symbol}
$\mathcal{P}_d(\rrn)$\index{$\mathcal{P}_d(\rrn)$} denotes
the set of all polynomials on $\rrn$
with degree not greater than $d$. Moreover,
for any ball $B\in\mathbb{B}$
and any $g\in L^{1}_{\rm loc}(\rrn)$,
$P_{B}^{d}g$\index{$P_{B}^{d}$} denotes
the \emph{minimizing
polynomial}\index{minimizing polynomial}
of $g$ with degree not greater than $d$,
which means that $P_{B}^{d}g$ is the unique
polynomial
$f\in\mathcal{P}_{d}(\rrn)$ such
that, for any $h\in\mathcal{P}_{d}(\rrn)$,
$$
\int_{B}[g(x)-f(x)]h(x)\,dx=0.
$$
Furthermore, for any
$r\in(0,\infty)$, in what follows,
we use $L^r_\mathrm{loc}
(\rrn)$\index{$L^r_\mathrm{loc}(\rrn)$}
to denote the set of all
$r$-order locally integrable
functions on $\rrn$.

Now, we introduce the following
Campanato-type function spaces
associated with local
generalized Herz spaces, which were
originally introduced in
\cite[Definition 3.2]{HYYZ} for
any given general
ball quasi-Banach function space $X$.

\begin{definition}\label{cams}
Let $p,\ q,\ s\in(0,\infty)$,
$r\in[1,\infty)$, $d\in\mathbb{Z}_{+}$,
and $\omega\in M(\rp)$ with $\m0(\omega)
\in(-\frac{n}{p},\infty)$. Then the
\emph{Campanato-type function
space}\index{Campanato-type function \\space}
$\HaSaHoc$\index{$\HaSaHoc$}, associated
with the local generalized Herz space $\HerzSo$,
is defined to be the set of all the
$f\in L^{r}_{\rm loc}(\rrn)$
such that
\begin{align*}
&\|f\|_{\HaSaHoc}\\&\quad:=\sup\left\|
\left\{\sum_{i=1}^{m}\left[\frac{\lambda_{i}}{
\|\1bf_{B_{i}}\|_{\HerzSo}}\right]
^{s}\1bf_{B_{i}}\right\}^{\frac{1}{s}}
\right\|_{\HerzSo}^{-1}\\
&\qquad\times\sum_{j=1}^{m}\left\{
\frac{\lambda_{j}|B_{j}|}{
\|\1bf_{B_{j}}\|_{\HerzSo}}\left[
\frac{1}{|B_{j}|}\int_{
B_{j}}\left|f(x)-P_{B_{j}}^{d}f(x)
\right|^{r}\,dx\right]
^{\frac{1}{r}}\right\}
\end{align*}
is finite,
where the supremum is taken over
all $m\in\mathbb{N}$, $\{B_{j}
\}_{j=1}^{m}\subset\mathbb{B}$,
and $\{\lambda_{j}\}_{j=1}^{m}
\subset[0,\infty)$
with $\sum_{j=1}^{m}\lambda_{
j}\neq0$.
\end{definition}

\begin{remark}
By Definition \ref{cams},
we can easily show that
$\mathcal{P}_d(\rrn)\subset
\HaSaHoc$ and, for any $f\in\HaSaHoc$,
$\|f\|_{\HaSaHoc}=0$ if
and only if $f\in\mathcal{P}_d(\rrn)$.
Therefore, in what follows,
we always identify $f\in\HaSaHoc$
with $\{f+P:\ P\in\mathcal{P}_d(\rrn)\}$.
\end{remark}

Applying \cite[Remark 3.3(iii) and
Proposition 3.4]{HYYZ} with
$X$ therein replaced by $\HerzSo$,
we immediately obtain the following
equivalent characterizations
of these Campanato-type
function spaces; we omit the details.

\begin{proposition}
Let $p$, $q$, $\omega$, $r$, $d$,
and $s$ be as in Definition \ref{cams}.
Then the following three
statements are equivalent:
\begin{enumerate}
  \item[{\rm(i)}] $f\in\HaSaHoc$;
  \item[{\rm(ii)}] $f\in L^r_\mathrm{loc}(\rrn)$
  and
  \begin{align*}
  &\normmm{f}_{\HaSaHoc}\\&\quad:=\sup\inf\left\|
\left\{\sum_{i=1}^{m}\left[\frac{\lambda_{i}}{
\|\1bf_{B_{i}}\|_{\HerzSo}}\right]
^{s}\1bf_{B_{i}}\right\}^{\frac{1}{s}}
\right\|_{\HerzSo}^{-1}\\
&\qquad\times\sum_{j=1}^{m}
\left\{\frac{\lambda_{j}|B_{j}|}{
\|\1bf_{B_{j}}\|_{\HerzSo}}
\left[\frac{1}{|B_{j}|}\int_{
B_{j}}\left|f(x)-P(x)\right|
^{r}\,dx\right]
^{\frac{1}{r}}\right\}
\end{align*}
is finite, where the supremum is
the same as in Definition \ref{cams}
and the infimum is taken over all
$P\in\mathcal{P}_d(\rrn)$;
  \item[{\rm(iii)}] $f\in L^r_\mathrm{loc}(\rrn)$
  and
  \begin{align*}
&\widetilde{\|f\|}_{\HaSaHoc}\\&\quad:=\sup\left\|
\left\{\sum_{i\in\mathbb{N}}
\left[\frac{\lambda_{i}}{
\|\1bf_{B_{i}}\|_{\HerzSo}}
\right]^{s}\1bf_{B_{i}}\right\}^{\frac{1}{s}}
\right\|_{\HerzSo}^{-1}\\
&\qquad\times\sum_{j\in\mathbb{N}}
\left\{\frac{\lambda_{j}|B_{j}|}{
\|\1bf_{B_{j}}\|_{\HerzSo}}\left[
\frac{1}{|B_{j}|}\int_{
B_{j}}\left|f(x)-P_{B_{j}}^{d}f(x)
\right|^{r}\,dx\right]
^{\frac{1}{r}}\right\}
\end{align*}
is finite, where the supremum is
taken over all $\{B_j\}_{
j\in\mathbb{N}}\subset\mathbb{B}$ and
$\{\lambda_j\}_{j\in\mathbb{N}}
\subset[0,\infty)$ satisfying
$$
\left\|
\left\{\sum_{i\in\mathbb{N}}\left[
\frac{\lambda_{i}}{
\|\1bf_{B_{i}}\|_{\HerzSo}}\right]^{s}
\1bf_{B_{i}}\right\}^{\frac{1}{s}}
\right\|_{\HerzSo}\in(0,\infty).
$$
\end{enumerate}
Moreover, there exist a constant
$C\in[1,\infty)$ such that,
for any $f\in L^r_\mathrm{loc}(\rrn)$,
$$
C^{-1}\|f\|_{\HaSaHoc}\leq
\normmm{f}_{\HaSaHoc}\leq
C\|f\|_{\HaSaHoc}
$$
and
$$
\widetilde{\|f\|}_{\HaSaHoc}=
\|f\|_{\HaSaHoc}.
$$
\end{proposition}

Via the known dual theorem of the
Hardy space $H_X(\rrn)$ associated with
the ball quasi-Banach function space $X$,
we now show that the dual space of
the generalized Herz--Hardy space $\HaSaHo$
is just the Campanato-type
function space $\HaSaHod$.
Namely, we have the following theorem.

\begin{theorem}\label{dualh}
Let $p$, $q\in(0,\infty)$,
$\omega\in M(\rp)$ with
$\m0(\omega)\in(-\frac{n}{p},\infty)$
and $\mi(\omega)\in(-\frac{n}{p},\infty)$,
$$p_{-}:=\min\left\{1,p,\frac{n}
{\Mw+n/p}\right\},$$
$d\geq\lfloor n(1/p_{-}-1)\rfloor$ be
a fixed integer, $s\in(0,\min\{p_{-},q\})$,
and $$r\in\left(\max\left\{1,p,\frac{n}
{\mw+n/p}\right\},\infty\right].$$
Then $\HaSaHod$ is the dual space of
$\HaSaHo$ in the following sense:
\begin{enumerate}
  \item[{\rm(i)}]
  Let $g\in\HaSaHod$. Then the linear
  functional
  \begin{equation}\label{dualh2}
  L_{g}:\ f\mapsto L_{g}(f):=
  \int_{\rrn}f(x)g(x)\,dx,
  \end{equation}
  initially defined for any
  $f\in\iaHaSaHo$, has a bounded extension to
  the generalized Herz--Hardy space $\HaSaHo$.
  \item[{\rm(ii)}] Conversely,
  any continuous linear functional
  on $\HaSaHo$ arises as in
  \eqref{dualh2} with a unique
  $g\in\HaSaHod$.
\end{enumerate}
Moreover, there exists a constant
$C\in[1,\infty)$ such that,
for any $g\in\HaSaHod$,
$$
C^{-1}\|g\|_{\HaSaHod}\leq\|L_{g}\|_{(\HaSaHo)^*}
\leq C\|g\|_{\HaSaHod},
$$
where $(\HaSaHo)^*$ denotes the
dual space of $\HaSaHo$.
\end{theorem}

To show this theorem, we first
recall the definition of the
ball Campanato-type function space
$\mathcal{L}_{X,
r,d,s}(\rrn)$\index{$\mathcal{L}_
{X,r,d,s}(\rrn)$}
associated with the general
ball quasi-Banach function space $X$
as follows, which is just
\cite[Definition 3.2]{HYYZ}.

\begin{definition}\label{camsx}
Let $X$ be a ball quasi-Banach function space,
$r\in[1,\infty)$, $s\in(0,\infty)$,
and $d\in\mathbb{Z}_{+}$. Then the
\emph{ball Campanato-type function
space}\index{ball Campanato-type function \\space}
$\mathcal{L}_{X,
r,d,s}(\rrn)$\index{$\mathcal{L}_{X,r,d,s}(\rrn)$},
associated with $X$,
is defined to be the set of all
the $f\in L^{r}_{\mathrm{loc}}(\rrn)$
such that
\begin{align*}
&\|f\|_{\mathcal{L}_{X,r,d,s}(\rrn)}\\&
\quad:=\sup\left\|
\left[\sum_{i=1}^{m}\left(\frac{\lambda_{i}}{
\|\1bf_{B_{i}}\|_{X}}\right)
^{s}\1bf_{B_{i}}\right]^{\frac{1}{s}}
\right\|_{X}^{-1}\\
&\qquad\times\sum_{j=1}^{m}\left\{
\frac{\lambda_{j}|B_{j}|}{
\|\1bf_{B_{j}}\|_{X}}\left[
\frac{1}{|B_{j}|}\int_{
B_{j}}\left|f(x)-P_{B_{j}}^{d}f(x)
\right|^{r}\,dx\right]
^{\frac{1}{r}}\right\}
\end{align*}
is finite,
where the supremum is taken over
all $m\in\mathbb{N}$, $\{B_{j}
\}_{j=1}^{m}\subset\mathbb{B}$,
and $\{\lambda_{j}\}_{j=1}^{m}
\subset[0,\infty)$
with $\sum_{j=1}^{m}\lambda_{
j}\neq0$.
\end{definition}

Notice that Zhang et al.\
\cite[Theorem 3.14]{HYYZ} established the
following duality between the
Hardy space $H_X(\rrn)$ and
the ball Campanato-type function space
$\mathcal{L}_{X,r',d,s}(\rrn)$, which
is an essential tool
for us to show Theorem \ref{dualh}.

\begin{lemma}\label{dualhl1}
Let $X$ be a ball quasi-Banach function
space satisfy:
\begin{enumerate}
  \item[{\rm(i)}] there exists a
  $p_-\in(0,\infty)$ such that, for
  any given $\theta\in(0,p_-)$ and
  $u\in(1,\infty)$, and for any
  $\{f_{j}\}_{j\in\mathbb{N}}\subset
  L^1_{\mathrm{loc}}(\rrn)$,
  $$
  \left\|\left\{
  \sum_{j\in\mathbb{N}}\left[
  \mc(f_j)\right]^{u}
  \right\}^{\frac{1}{u}}\right\|_{X^{1/\theta}}
  \lesssim\left\|
  \left\{
  \sum_{j\in\mathbb{N}}|f_j|^u
  \right\}^{\frac{1}{u}}\right\|_{X^{1/\theta}},
  $$
  where the implicit positive constant
  is independent of $\{f_j\}_{j\in\mathbb{N}}$;
  \item[{\rm(ii)}] for the above $p_-$, there exist
  an $s_0\in(0,\min\{1,p_-\})$ and
  an $r_0\in(s_0,\infty)$ such that
  $X^{1/s_0}$ is a ball Banach function space
  and, for any $f\in L^1_{\mathrm{loc}}(\rrn)$,
  $$
  \left\|\mc^{((r_0/s_0)')}
  (f)\right\|_{(X^{1/s_0})'}\lesssim
  \|f\|_{(X^{1/s_0})'},
  $$
  where the implicit positive constant
  is independent of $f$;
  \item[{\rm(iii)}] $X$ has an absolutely
  continuous quasi-norm.
\end{enumerate}
Assume $d\geq\lfloor n(1/\min\{1,p_-\}-1)\rfloor$,
$s\in(0,s_0]$, and $r\in(\max\{1,r_0\},\infty]$.
Then the dual space of $H_X(\rrn)$, denoted by
$(H_X(\rrn))^*$, is $\mathcal{L}_{X,r',d,s}(\rrn)$
in the following sense:
\begin{enumerate}
\item[{\rm (a)}] Let
$g\in\mathcal{L}_{X,r',d,s}(\rrn)$.
Then the linear functional
\begin{align}\label{2te1}
L_g:\ f\mapsto
L_g(f):=\int_{\rrn}f(x)g(x)\,dx,
\end{align}
initially defined for any
$f\in H_{\rm fin}^{X,r,d,s}(\rrn)$,
has a bounded extension to $H_X(\rrn)$.

\item[{\rm (b)}] Conversely, any continuous linear
functional on $H_X(\rrn)$ arises as in \eqref{2te1}
with a unique $g\in\mathcal{L}_{X,r',d,s}(\rrn)$.
\end{enumerate}
Moreover,
$\|g\|_{\mathcal{L}_
{X,r',d,s}(\rrn)}\sim\|L_g\|_{(H_X(\rrn))^*}$,
where the positive equivalence constants
are independent of $g$.
\end{lemma}

With the help of the above
lemma, we then prove Theorem \ref{dualh}.

\begin{proof}[Proof of Theorem \ref{dualh}]
Let all the symbols be as in the present theorem.
Then, from the assumption
$\m0(\omega)\in(-\frac{n}{p},\infty)$
and Theorem \ref{Th3}, we deduce
that the local generalized
Herz space $\HerzSo$ is a BQBF space.
This implies that, to finish the proof
of the present theorem,
we only need to show that $\HerzSo$ satisfies
(i), (ii), and (iii) of Lemma \ref{dualhl1}.

First, we prove that Lemma \ref{dualhl1}(i)
holds true for $\HerzSo$.
Indeed, for any given $\theta\in(0,p_-)$
and $u\in(1,\infty)$,
by Lemma \ref{vmbhl} with
$r:=\theta$, we conclude that, for
any $\{f_{j}\}_{j\in\mathbb{N}}
\subset L^1_{\mathrm{loc}}(\rrn)$,
\begin{equation*}
\left\|\left\{\sum_{j\in\mathbb{N}}\left[\mc
(f_{j})\right]^{u}\right\}
^{\frac{1}{u}}\right\|_{[\HerzSo]^{1/\theta}}
\lesssim\left\|\left(
\sum_{j\in\mathbb{N}}|f_{j}|
^{u}\right)^{\frac{1}{u}}\right\|_{[\HerzSo
]^{1/\theta}},
\end{equation*}
which implies that Lemma \ref{dualhl1}(i)
holds true for $\HerzSo$.

Next, we show that $\HerzSo$
satisfies Lemma \ref{dualhl1}(ii).
To this end, let $s_{0}\in(s,\min\{p_-,q\})$ and
$$
r_{0}\in\left(\max\left
\{1,p,\frac{n}{\mw+n/p}\right\},r\right).
$$
Then, from Lemma \ref{mbhal} with
$s$ and $r$ therein replaced,
respectively, by $s_0$ and $r_0$,
it follows that
$[\HerzSo]^{1/s_{0}}$ is a BBF space and,
for any $f\in L^1_{\mathrm{loc}}(\rrn)$,
\begin{equation*}
\left\|\mc^{((r_{0}/s_{0})')}(f)
\right\|_{([\HerzSo]^{1/s_{0}})'}\lesssim
\left\|f\right\|_{([\HerzSo]^{1/s_{0}})'}.
\end{equation*}
This implies that, under the assumptions
of the present theorem, $\HerzSo$ satisfies
Lemma \ref{dualhl1}(ii).

Finally, applying Theorem \ref{abso}, we find that
the local generalized Herz space $\HerzSo$
has an absolutely continuous quasi-norm.
Therefore,
all the assumptions of Lemma \ref{dualhl1}
hold true for $\HerzSo$. Thus,
the proof of Theorem \ref{dualh} is
then completed.
\end{proof}

\begin{remark}
We point out that the absolutely continuous quasi-norm
of the ball quasi-Banach function space $X$ plays a key
role in the dual theorem of the associated Hardy space $H_X(\rrn)$.
However, by Example \ref{conter}, we find
that the global generalized
Herz space does not have an absolutely continuous quasi-norm. Thus, the dual
space of the generalized Herz--Hardy space $\HaSaH$ is still unknown.
\end{remark}

From the above dual theorem, we immediately deduce
the following equivalence of the Campanato-type
function space $\HaSaHoc$; we omit the details.

\begin{corollary}
Let $p$, $q$, $\omega$, $p_-$, $d$, and $s$ be as in Theorem \ref{dualh},
$$
p_+:=\max\left\{1,p,\frac{n}{\mw+n/p}\right\},
$$
$r\in[1,p_+')$, $d_0:=\lfloor n(1/p_--1)\rfloor$, and $s_0
\in(0,\min\{p_-,q\})$. Then
$$
\HaSaHoc=\dot{\mathcal{L}}_{\omega,\0bf}^{p,q,1,d_0,s_0}(\rrn)
$$
with equivalent quasi-norms.
\end{corollary}

As an application of Theorem \ref{dualh}, we next
investigate the dual
space of the generalized Morrey--Hardy space $H\MorrSo$ via
introducing the following Campanato-type function space
$\dot{\mathscr{L}}_{\omega,\0bf}^{p,q,r,d,s}(\rrn)$
associated with the local generalized Morrey space.

\begin{definition}\label{camsm}
Let $p$, $q$, $r\in[1,\infty)$,
$s\in(0,\infty)$, $d\in\mathbb{Z}_{+}$,
and $\omega\in M(\rp)$ satisfy $\MI(\omega)
\in(-\infty,0)$ and
$$
-\frac{n}{p}<\m0(\omega)\leq\M0(\omega)<0.
$$
Then the
\emph{Campanato-type
function space}\index{Campanato-type function \\space}
$\dot{\mathscr{L}}_{\omega,\0bf}^{
p,q,r,d,s}(\rrn)$\index{$\dot{\mathscr{L}}_{\omega,\0bf}^{p,q,r,d,s}(\rrn)$},
associated with the local generalized Morrey space
$\MorrSo$, is defined to be the set of all
the $f\in L^{r}_{\rm loc}(\rrn)$
such that
\begin{align*}
&\|f\|_{\dot{\mathscr{L}}_{
\omega,\0bf}^{p,q,r,d,s}(\rrn)}\\&\quad:=\sup\left\|
\left\{\sum_{i=1}^{m}\left[\frac{\lambda_{i}}{
\|\1bf_{B_{i}}\|_{\MorrSo}}\right]^{s}
\1bf_{B_{i}}\right\}^{\frac{1}{s}}
\right\|_{\MorrSo}^{-1}\\
&\qquad\times\sum_{j=1}^{m}\left\{\frac{\lambda_{j}|B_{j}|}{
\|\1bf_{B_{j}}\|_{\MorrSo}}\left[\frac{1}{|B_{j}|}\int_{
B_{j}}\left|f(x)-P_{B_{j}}^{d}f(x)\right|^{r}\,dx\right]
^{\frac{1}{r}}\right\}
\end{align*}
is finite,
where the supremum is taken over
all $m\in\mathbb{N}$, $\{B_{j}\}_{j=1}^{m}\subset\mathbb{B}$,
and $\{\lambda_{j}\}_{j=1}^{m}
\subset[0,\infty)$ with $\sum_{j=1}^{m}\lambda_{
j}\neq0$.
\end{definition}

Using Theorem \ref{dualh} and
Remarks \ref{remhs}(iv) and \ref{remark4.10}(ii),
we immediately obtain the following conclusion,
which shows that the dual space of the
generalized Morrey--Hardy
space $H\MorrSo$ is just the Campanato-type
function space $\dot{\mathscr{L}
}_{\omega,\0bf}^{p,q,r',d,s}(\rrn)$;
we omit the details.

\begin{corollary}
Let $p$, $q\in[1,\infty)$, $\omega\in M(\rp)$ with
$$
-\frac{n}{p}<\m0(\omega)\leq\M0(\omega)<0
$$
and
$$
-\frac{n}{p}<\mi(\omega)\leq\MI(\omega)<0,
$$
$d\in\zp$, $s\in(0,1)$, and $$r
\in\left(\frac{n}{\mw+n/p},\infty\right].$$
Then $\dot{\mathscr{L}}_{\omega,\0bf}^{p,q,r',d,s}(\rrn)$
is the dual space of $H\MorrSo$ in the following sense:
\begin{enumerate}
  \item[{\rm(i)}]
  Let $g\in\dot{\mathscr{L}}_{
  \omega,\0bf}^{p,q,r',d,s}(\rrn)$. Then the linear
  functional
  \begin{equation}\label{dualhm2}
  L_{g}:\ f\mapsto L_{g}(f):=\int_{\rrn}f(x)g(x)\,dx,
  \end{equation}
  initially defined for any $f\in H
  \textsl{\textbf{M}}_{\omega,\0bf,\mathrm{fin}}
  ^{p,q,r,d,s}(\rrn)$, has a bounded extension to
  the generalized Morrey--Hardy space $H\MorrSo$.
  \item[{\rm(ii)}] Conversely, any continuous linear functional
  $L\in(H\MorrSo)^*$ arises as in \eqref{dualhm2} with a unique
  $g\in\dot{\mathscr{L}}_{\omega,\0bf}^{p,q,r',d,s}(\rrn)$.
\end{enumerate}
Moreover, there exists a constant $C\in[1,\infty)$ such that,
for any $g\in\dot{\mathscr{L}}_{
\omega,\0bf}^{p,q,r',d,s}(\rrn)$,
$$
C^{-1}\|g\|_{\dot{\mathscr{L}}_{
\omega,\0bf}^{p,q,r',d,s}(\rrn)}
\leq\|L_{g}\|_{(H\MorrSo)^*}
\leq C\|g\|_{\dot{\mathscr{L}
}_{\omega,\0bf}^{p,q,r',d,s}(\rrn)},
$$
where $(H\MorrSo)^*$ denotes the dual space of $H\MorrSo$.
\end{corollary}

\section{Boundedness of Calder\'{o}n--Zygmund Operators}

The main target of this section is
to investigate the boundedness of
Calder\'{o}n--Zygmund operators
on generalized
Herz--Hardy spaces. To this end, we first
establish two general boundedness
criteria of Calder\'{o}n--Zygmund
operators on Hardy spaces associated with
ball quasi-Banach function spaces
(see Proposition \ref{bddoxl4} and Theorem
\ref{bddox} below)
under some reasonable assumptions.
Via these results and
the facts that both local
and global generalized Herz spaces
are ball quasi-Banach function spaces,
we then obtain the boundedness of
Calder\'{o}n--Zygmund
operators on generalized Herz--Hardy spaces.

Let $d\in\zp$ and $T$ be a $d$-order
Calder\'on--Zygmund operator defined
as in Definition \ref{defin-C-Z-s}.
Recall the well-known assumption on
$T$ that,
for any $\gamma\in\zp^n$ with $|\gamma|\leq d$,
$T^*(x^{\gamma})=0$\index{$T^*(x^{\gamma})=0$},
namely, for any $a\in L^2(\rrn)$ having compact support
and satisfying that, for any $\gamma\in\zp^n$ with
$|\gamma|\leq d$,
$\int_{\rrn} a(x)x^\gamma\,dx=0$, it holds true that
\begin{align*}
\int_{\rrn}T(a)(x)x^\gamma\,dx=0
\end{align*}
(see, for instance, \cite[p.\,119]{lu}).

\begin{definition}\label{Def-T-s-v}
Let $d\in\zp$. A $d$-order
Calder\'on--Zygmund operator $T$
is said to have the \emph{vanishing moments up to
order $d$}\index{vanishing moments up to order $d$}
if, for any $\gamma\in\zp^n$ with
$|\gamma|\leq d$,
$T^*(x^{\gamma})=0$.
\end{definition}

We should point out that
the assumption that the
$d$-order Calder\'{o}n--Zygmund
operator $T$ has the vanishing
moments up to order $d$ is reasonable.
Indeed, this assumption
holds true automatically
when $T$ is a
Calder\'{o}n--Zygmund operator with
kernel $K(x,y):=K_1(x-y)$
for a locally integrable
function $K_1$ on $\rrn\setminus
\{\0bf\}$.
To be precise, we have
the following interesting proposition
about convolutional
type Calder\'{o}n--Zygmund
operators\index{convolutional type Calder\'{o}n--Zygmund
\\ operator},
which might be well known.
But, we do not find
its detailed proof in the literature.
Thus, for the convenience
of the reader,
we present its detailed proof as
follows.

\begin{proposition}\label{con-T-s-v}
Let $d\in\zp$,
$K\in\mathcal{S}'(\rrn)$, and the operator
$T$ be defined by setting, for any $f\in
\mathcal{S}(\rrn)$, $T(f):=K\ast f$.
Assume that the following three statements
hold true:
 \begin{enumerate}
  \item[{\rm(i)}] $\widehat{K}
  \in L^{\infty}(\rrn)$;
  \item[{\rm(ii)}] $K$ coincides with a
  function belonging to $C^{d}(\rrn\setminus\{\0bf\})
  $ in the sense that, for any given
  $a\in L^2(\rrn)$
  with compact support and for any $x\notin\supp(a)$,
  $$
  T(a)(x)=\int_{\rrn}K(x-y)a(y)\,dy;
  $$
  \item[{\rm(iii)}] there exist a positive
  constant $C$ and a $\delta\in(0,1]$ such that,
  for any $\gamma\in\zp^n$ with $|\gamma|\leq d$
  and for any $x\in\rrn\setminus\{\0bf\}$,
  \begin{equation}\label{con-T-size}
  \left|\partial^{\gamma}K(x)\right|
  \leq\frac{C}{|x|^{n+|\gamma|}}
  \end{equation}
  and, for any $\gamma\in\zp^n$ with $|\gamma|
  =d$ and $x$, $y\in\rrn$ with $|x|>2|y|$,
  \begin{equation}\label{con-T-regular}
  \left|\partial^{\gamma}K(x-y)-
  \partial^{\gamma}K(x)\right|\leq C\frac{|y|
  ^{\delta}}{|x|^{n+d+\delta}}.
  \end{equation}
\end{enumerate}
Then $T$ has a vanishing moments up
to order $d$.
\end{proposition}

To show Proposition \ref{con-T-s-v}, we
need some preliminary lemmas. Recall
that the following conclusion gives
the $L^p(\rrn)$ boundedness
of convolutional type
Calder\'{o}n--Zygmund operators,
which is just \cite[Theorem 5.1]{d01}.

\begin{lemma}\label{con-T-bdd}
Let $p\in(1,\infty)$
and $T$ be as in Proposition \ref{con-T-s-v}.
Then $T$ is well defined on $L^p(\rrn)$ and
there exists a positive constant $C$
such that, for any $f\in L^p(\rrn)$,
$$
\left\|T(f)\right\|_{L^p(\rrn)}
\leq C\|f\|_{L^p(\rrn)}.
$$
\end{lemma}

Via borrowing some ideas
from the proof of \cite[p.\,117,\ Lemma]{s93},
we next establish the following
technique estimates about
the kernel $K$, which
plays an important role in the proof
of Proposition \ref{con-T-s-v}.

\begin{lemma}\label{stein-size}
Let $d$, $K$, and $\delta$ be as in Proposition
\ref{con-T-s-v}.
For any $t\in(0,\infty)$, let $K^{(t)}:=
K\ast\phi_t$, where $\phi\in\mathcal{S}(\rrn)$
with $\supp(\phi)\subset B(\0bf,1)$.
Then there exists
a positive constant $C$ such that
\begin{enumerate}
  \item[{\rm(i)}] for any $t\in(0,\infty)$,
  $\gamma\in\zp^n$
  with $|\gamma|\leq d$,
  and $x\in\rrn\setminus\{\0bf\}$,
  \begin{equation*}
  \left|\partial^{\gamma}K^{(t)}(x)\right|
  \leq\frac{C}{|x|^{n+|\gamma|}};
  \end{equation*}
  \item[{\rm(ii)}] for any $t\in(0,\infty)$,
  $\gamma\in\zp^n$ with $|
\gamma|=d$, and $x$, $y\in\rrn$ with
$|x|>4|y|$,
$$
\left|\partial^{\gamma}K^{(t)}(x-y)-
\partial^{\gamma}K^{(t)}(x)
\right|\leq C\frac{|y|
^{\delta}}{|x|^{n+d+\delta}}.
$$
\end{enumerate}
\end{lemma}

\begin{proof}
Let all the symbols be as in
the present lemma and $t\in(0,\infty)$.
We first prove (i). Indeed,
notice that, for any $x\in\rrn$,
\begin{equation*}
K^{(t)}(x)=\int_{\rrn}
e^{2\pi ix\cdot\xi}\widehat{K}(\xi)
\widehat{\phi}(t\xi)\,d\xi.
\end{equation*}
From this, we further deduce that,
for any $\gamma\in\zp^n$ with
$|\gamma|\leq d$ and for any $x\in\rrn$,
\begin{equation}\label{stein-sizee2}
\partial^{\gamma}K^{(t)}(x)=
\int_{\rrn}\left(2\pi i\xi\right)^{\gamma}
e^{2\pi ix\cdot\xi}\widehat{K}(\xi)
\widehat{\phi}(t\xi)\,d\xi.
\end{equation}
which, together with the assumption
$\widehat{K}\in L^{\infty}(\rrn)$,
further implies that, for any
$x\in\rrn$ with $0<|x|<2t$,
\begin{align}\label{stein-sizee3}
\left|\partial^{\gamma}K^{(t)}(x)\right|
&\lesssim\int_{\rrn}|\xi|^{\gamma}
\left|\widehat{\phi}(t\xi)\right|\,d\xi
\sim\frac{1}{t^{n+|\gamma|}}\int_{\rrn}
\left|\widehat{\phi}(\xi)\right|\,d\xi\notag\\
&\sim\frac{1}{t^{n+|\gamma|}}\lesssim
\frac{1}{|x|^{n+|\gamma|}},
\end{align}
where the implicit positive
constants are independent of $t$.

On the other hand, applying the both assumptions
Proposition \ref{con-T-s-v}(ii) and
$$\supp(\phi_t)\subset B(\0bf,t),$$ we
find that, for any $x\in\rrn$ with $|x|\geq2t$,
$$
K^{(t)}(x)=\int_{\rrn}K(x-y)\phi_t(y)\,dy.
$$
This implies that, for any $\gamma\in\zp^n$
with $|\gamma|\leq d$ and for any $x\in\rrn$ with
$|x|\geq2t$,
\begin{align}\label{stein-sizee4}
\partial^{\gamma}K^{(t)}(x)&=\int_{\rrn}
\partial^{\gamma}K(x-y)\phi_t(y)\,dy\notag\\
&=\int_{B(\0bf,t)}
\partial^{\gamma}K(x-y)\phi_t(y)\,dy.
\end{align}
In addition, for any $x\in\rrn$
with $|x|\geq2t$ and for any $y\in B(\0bf,t)$,
we have
$$
|x-y|\geq|x|-|y|>\frac{1}{2}|x|.
$$
Combining this, \eqref{stein-sizee4},
and \eqref{con-T-size} with $x$
therein replaced by $x-y$,
we further conclude that, for any $\gamma\in\zp^n$
with $|\gamma|\leq d$ and for any $x\in\rrn$ with
$|x|\geq2t$,
\begin{equation}\label{stein-sizee5}
\left|\partial^{\gamma}K^{(t)}(x)\right|
\lesssim\frac{1}{|x-y|^{n+|\gamma|}}
\int_{\rrn}|\phi(y)|\,dy\lesssim
\frac{1}{|x|^{n+|\gamma|}},
\end{equation}
where the implicit positive
constants are independent of $t$. Moreover,
from \eqref{stein-sizee3} and \eqref{stein-sizee5},
it follows that (i) holds true.

Next, we show (ii). To this end, fix
$x$, $y\in\rrn$ with $|x|>4|y|$. We now
prove (ii) by considering the following
two cases on $x$.

\emph{Case 1)} $x\in\rrn$ with $|x|<4t$.
In this case, we have $|y|<\frac{1}{4}|x|<t$.
Then, using \eqref{stein-sizee2},
the Lagrange mean value theorem,
and the assumptions $\widehat{K}\in L^{\infty}
(\rrn)$ and $\delta\in(0,1]$, we
find that there exists a $y_1\in\rrn$
such that, for any $\gamma\in\zp^n$
with $|\gamma|=d$,
\begin{align}\label{stein-sizee6}
&\left|\partial^{\gamma}K^{(t)}(x-y)
-\partial^{\gamma}K^{(t)}(x)\right|\notag\\
&\quad=\left|\int_{\rrn}\left(2\pi i\xi\right)
^{\gamma}\left[e^{2\pi i(x-y)\cdot\xi}-e^{
2\pi ix\cdot\xi}\right]\widehat{K}(\xi)
\widehat{\phi}(t\xi)\,d\xi\right|\notag\\
&\quad=\left|\int_{\rrn}\left(2\pi i\xi\right)
^{\gamma}(2\pi i\xi)\cdot ye^{2\pi i(x-y_1)\cdot\xi}
\widehat{K}(\xi)\widehat{\phi}(t\xi)\,d\xi\right|\notag\\
&\quad\lesssim|y|\int_{\rrn}|\xi|^{d+1}\left|
\widehat{\phi}(t\xi)\right|\,d\xi
\sim\frac{|y|}{t^{n+d+1}}\int_{\rrn}
\left|\widehat{\phi}(\xi)\right|\,d\xi\notag\\
&\quad\lesssim\frac{1}{t^{n+d}}\frac{|y|^{\delta}}
{t^{\delta}}\lesssim\frac{|y|^{\delta}}{
|x|^{n+d+\delta}},
\end{align}
where the implicit positive constants
are independent of $t$. This finishes the
proof of (ii) under the assumption
$|x|<4t$.

\emph{Case 2)} $|x|\geq4t$. In this case,
we first claim that $|x-y|>t$.
Indeed, when $|y|<t$, we have
$$
|x-y|\geq|x|-|y|>3t>t.
$$
On the other hand, if $|y|\geq t$,
from the assumption $|x|>4|y|$, we deduce that
$$
|x-y|\geq|x|-|y|>3|y|>t.
$$
Therefore, the above claim holds true.
Using this claim and the assumptions
Proposition \ref{con-T-s-v}(ii),
$|x|\geq 4t$, and
$\supp(\phi_t)\subset B(\0bf,t)$,
we conclude that
\begin{equation*}
K^{(t)}(x)=\int_{|z|<t}K(x-z)\phi_t(z)\,dz
\end{equation*}
and
\begin{equation*}
K^{(t)}(x-y)=\int_{|z|<t}K(x-y-z)\phi_t(z)\,dz,
\end{equation*}
which further imply that, for any $\gamma\in\zp^n$
with $|\gamma|=d$,
\begin{equation}\label{stein-sizee7}
\partial^{\gamma}
K^{(t)}(x)=\int_{|z|<t}\partial^{\gamma}
K(x-z)\phi_t(z)\,dz
\end{equation}
and
\begin{equation}\label{stein-sizee8}
\partial^{\gamma}
K^{(t)}(x-y)=\int_{|z|<t}\partial^{
\gamma}K(x-y-z)\phi_t(z)\,dz.
\end{equation}
We now estimate $|x-z|$ for any $|z|<t$.
Indeed, for any $z\in\rrn$ with
$|z|<t$, we have $|x|\geq4t>4|z|$.
This implies that, for any
$z\in\rrn$ with $|z|<t$,
$$
|x-z|\geq|x|-|z|>\frac{3}{4}|x|>3|y|.
$$
Applying this, \eqref{stein-sizee7},
\eqref{stein-sizee8}, and the assumption
\eqref{con-T-regular} with
$x$ and $y$ therein replaced,
respectively, by
$x-z$ and $y$,
we conclude that, for any $\gamma\in\zp^n$
with $|\gamma|=d$,
\begin{align}\label{stein-sizee9}
&\left|\partial^{\gamma}K^{(t)}(x-y)
-\partial^{\gamma}K^{(t)}(x)\right|\notag\\
&\quad\leq\int_{|z|<1}\left|
\partial^{\gamma}K(x-y-z)-\partial^{\gamma}
K(x-z)\right|\left|\phi_t(z)\right|\,dz\notag\\
&\quad\lesssim\int_{|z|<1}\frac{|y|^{\delta}}{|x-z|
^{n+d+\delta}}\left|\phi_t(z)\right|\,dz
\lesssim\|\phi\|_{L^1(\rrn)}\frac{|y|^{\delta}}{
|x|^{n+d+\delta}},
\end{align}
where the implicit positive constants
are independent of $t$. This finishes
the proof of (ii) under the assumption $|x|\geq4t$.
Combining both \eqref{stein-sizee6}
and \eqref{stein-sizee9}, we then
complete the proof of (ii), and hence
of Lemma \ref{stein-size}.
\end{proof}

In addition, via borrowing some ideas from
the arguments used in \cite[p.\,118]{s93},
we obtain the following auxiliary
lemma about convolutional
type Calder\'{o}n--Zygmund operators
and radial maximal functions,
which plays a key role in the proof of
Proposition \ref{con-T-s-v}.

\begin{lemma}\label{stein-hardy}
Let $\phi\in\mathcal{S}(\rrn)$
with $\supp(\phi)\subset B(\0bf,1)$, and
$d$, $K$, $\delta$,
and $T$ be as in Proposition
\ref{con-T-s-v}.
Assume $a\in L^2(\rrn)$, with compact
support, satisfying that,
for any $\gamma\in\zp^n$ with
$|\gamma|\leq d$, $\int_{\rrn}a(x)x^{
\gamma}\,dx=0$. Then
$
M(T(a),\phi)|\cdot|^{|\gamma|}\in
L^1(\rrn)$ for any $\gamma\in\zp^n$
with $|\gamma|\leq d$,
where, for any $f\in\mathcal{S}'(\rrn)$
and $\psi\in\mathcal{S}(\rrn)$,
$M(f,\psi)$ is defined as in Definition
\ref{smax}(i).
\end{lemma}

\begin{proof}
Let all the symbols be as
in the present lemma and
$\supp(a)\subset B(x_0,r_0)$
with $x_0\in\rrn$ and $r_0\in(0,\infty)$.
We first estimate $M(T(a),\phi)$.
Indeed, from the H\"{o}lder
inequality, \eqref{Atogl3e1},
the $L^2(\rrn)$ boundedness
of the Hardy--Littlewood maximal
operator $\mc$,
and Lemma \ref{con-T-bdd} with $p:=2$,
it follows that
\begin{align}\label{stein-hardye1}
&\int_{B(x_0,4r_0)}M\left(T(a),\phi
\right)(x)\,dx\notag\\
&\quad\lesssim\left\|
M\left(T(a),\phi\right)\1bf_{B(x_0,4r_0)}
\right\|_{L^2(\rrn)}
\lesssim\left\|\mc\left(T(a)\right)\right\|_{L^2(\rrn)}
\notag\\&\quad\lesssim\|T(a)\|_{L^2(\rrn)}\lesssim
\|a\|_{L^2(\rrn)}<\infty.
\end{align}
This is the desired estimate
of $M(T(a),\phi)(x)$ with
$x\in B(x_0,4r_0)$.

On the other hand, we estimate
$M(T(a),\phi)(x)$ when
$x\in[B(x_0,4r_0)]^{\complement}$.
To achieve this, fix a $t\in(0,\infty)$.
Then, by the assumption
that, for any $\gamma\in\zp^n$
with $|\gamma|\leq d$, $
\int_{\rrn}a(x)x^{\gamma}\,dx=0$
and the Taylor remainder theorem,
we find that, for any $y\in
B(x_0,r_0)$, there exists a $t_y\in
[0,1]$ such that,
for any $x\in\rrn$,
\begin{align}\label{stein-hardye2}
&T(a)\ast\phi_t(x)\notag\\
&\quad=a\ast K^{(t)}(x)=\int_{B(x_0,r_0)}
K^{(t)}(x-y)a(y)\,dy\notag\\
&\quad=\int_{B(x_0,r_0)}
\left[K^{(t)}(x-y)-\sum_{\gfz{\gamma\in\zp^n}
{|\gamma|\leq d}}\frac{\partial^{\gamma}
K^{(t)}(x-x_0)}{\gamma!}(y-x_0)
^{\gamma}\right]a(y)\,dy\notag\\
&\quad=\int_{B(x_0,r_0)}\sum_{\gfz{\gamma\in\zp^n}
{|\gamma|=d}}\frac{\partial^
{\gamma}K^{(t)}(x-t_yy-(1-t_y)x_0)
-\partial^{\gamma}K^{(t)}(x-x_0)}{\gamma!}\notag\\
&\qquad\times(y-x_0)^{\gamma}
a(y)\,dy.
\end{align}
Notice that, for any
$x\in[B(x_0,4r_0)]^{\complement}$ and
$y\in B(x_0,r_0)$, we have
$$|x-x_0|\geq4r_0>4|t_y(y-x_0)|,$$ where
$t_y$ is as in \eqref{stein-hardye2}.
This, together with
\eqref{stein-hardye2},
Lemma \ref{stein-size}(ii)
with $x$ and $y$ replaced, respectively,
by $x-x_0$ and $t_y(y-x_0)$ for any
$y\in B(x_0,r_0)$, and the H\"{o}lder
inequality, further implies that,
for any $x\in[B(x_0,4r_0)]^{\complement}$,
\begin{align*}
\left|T(a)\ast\phi_t(x)\right|
\lesssim\int_{B(x_0,r_0)}\frac{|y-x_0|^{d+\delta}}{
|x-x_0|^{n+d+\delta}}|a(y)|\,dy
\lesssim\frac{\|a\|_{L^2(\rrn)}}
{|x-x_0|^{n+d+\delta}},
\end{align*}
where the implicit positive constants
are independent of $t$. Therefore,
for any $x\in[B(x_0,4r_0)]^{\complement}$,
it holds true that
\begin{equation}\label{stein-hardye3}
M\left(T(a),\phi\right)(x)\lesssim
\frac{\|a\|_{L^2(\rrn)}}
{|x-x_0|^{n+d+\delta}},
\end{equation}
which is the desired estimate of
$M(T(a),\phi)(x)$ with
$x\in[B(x_0,4r_0)]^{\complement}$.
Then, for any $\gamma\in\zp^n$ with
$|\gamma|\leq d$,
applying \eqref{stein-hardye1},
\eqref{stein-hardye3}, and the
facts $n+d+\delta>n$ and
$n+d-|\gamma|+\delta>n$,
we further conclude that
\begin{align*}
&\int_{\rrn}M\left(T(a),\phi
\right)(x)|x|^{|\gamma|}\,dx\\
&\quad=\int_{B(x_0,4r_0)}M\left(T(a),\phi
\right)(x)|x|^{|\gamma|}\,dx+\int_{[B(x_0,4r_0)]
^{\complement}}\cdots\\&\quad\lesssim
(x_0+4r_0)^{|\gamma|}\int_{B(x_0,
4r_0)}M\left(T(a),\phi\right)(x)\,dx\\
&\qquad+\|a\|_{L^2(\rrn)}\left[\int_{[B(x_0,4r_0)]^{
\complement}}\left(\frac{|x_0|^{|\gamma|}}{
|x-x_0|^{n+d+\delta}}+\frac{1}{
|x-x_0|^{n+d-|\gamma|+\delta}}\right)\,dx\right]\\
&\quad\lesssim\|a\|
_{L^2(\rrn)}\left[1+
\int_{[B(x_0,4r_0)]^{
\complement}}\left(\frac{1}{
|x-x_0|^{n+d+\delta}}+\frac{1}{
|x-x_0|^{n+d-|\gamma|+\delta}}\right)\,dx\right]\\
&\quad<\infty.
\end{align*}
This implies that,
for any $\gamma\in\zp^n$ with
$|\gamma|\leq d$, $M(T(a),\phi)|\cdot|^{|\gamma|}
\in L^1(\rrn)$ and hence finishes the proof
of Lemma \ref{stein-hardy}.
\end{proof}

Recall that the \emph{Hardy space}
$H^1(\rrn)$\index{$H^1(\rrn)$} is defined as
in Definition \ref{hardyx}
with $X:=L^1(\rrn)$.
The following atomic
decomposition of $H^1(\rrn)$ can
be deduced from the proof of
\cite[Chapter 2, Proposition 3.3]{lu}
immediately,
which is an essential tool
in the proof of Proposition \ref{con-T-s-v};
we omit the details.

\begin{lemma}\label{atom-hardy}
Let $N\in\mathbb{N}$, $d\in\zp$,
and $f\in H^1(\rrn)\cap L^2(\rrn)$.
Then there exist $\{\lambda_{i,j}\}
_{i\in\mathbb{Z},\,j\in\mathbb{N}}
\subset[0,\infty)$ and
a sequence $\{a_{i,j}\}_{i\in\mathbb{Z},\,
j\in\mathbb{N}}$
of $(L^1(\rrn),\,\infty,\,d)$-atoms supported,
respectively, in the balls
$\{B_{i,j}\}_{i\in\mathbb{Z},\,j\in\mathbb{N}}
\subset\mathbb{B}$ such that
$f=\sum_{i\in\mathbb{Z}}\sum_{j\in\mathbb{N}}
\lambda_{i,j}a_{i,j}$ almost everywhere
in $\rrn$ and there exists a positive
constant $C$ such that, for any $i\in\mathbb{Z}$
and $j\in\mathbb{N}$, the following
three statements hold true:
\begin{enumerate}
  \item[{\rm(i)}] $\lambda_{i,j}
  |a_{i,j}|\leq C2^i$;
  \item[{\rm(ii)}] $\bigcup_{j\in\mathbb{N}}
  B_{i,j}=\Omega_i:=\{x\in\rrn:\ \mc_{N}(f)(x)>2^i\}$
  with $\mc_N$ as in \eqref{sec6e1};
  \item[{\rm(iii)}] $\sum_{j\in\mathbb{N}}
  \1bf_{B_{i,j}}\leq C$.
\end{enumerate}
\end{lemma}

In order to show Proposition
\ref{con-T-s-v},
we also require some auxiliary
conclusions about weighted function spaces.
The following one about powered weights
is a part of \cite[Example 7.1.7]{LGCF}.

\begin{lemma}\label{power-weight}
Let $a\in\rr$ and $p\in(1,\infty)$.
Then $|\cdot|^{a}\in A_p(\rrn)$ if and
only if
$-n<a<n(p-1)$.
\end{lemma}

Moreover, the following conclusion
gives the strong type
inequality characterization of
$A_p(\rrn)$-weights with $p\in(1,\infty)$,
which can be found in \cite[Theorem 7.3]{d01}.

\begin{lemma}\label{hl-weight}
Let $p\in(1,\infty)$ and
$\upsilon\in A_{\infty}(\rrn)$.
Then the Hardy--Littlewood maximal operator
$\mc$ is bounded on $L^p_{\upsilon}(\rrn)$
if and only if $\upsilon\in A_p(\rrn)$.
\end{lemma}

Via the above two lemmas,
we show the following
technique conclusion about
grand maximal functions
and radial maximal functions,
which is vital in the
proof of Proposition \ref{con-T-s-v}.

\begin{lemma}\label{max-hardy}
Let $d\in\zp$, $N\in\mathbb{N}\cap
[n+d+1,\infty)$, $f\in\mathcal{S'}(\rrn)$,
and $\phi\in\mathcal{S}(\rrn)$ with
$\supp(\phi)\subset B(\0bf,1)$ and
$\int_{\rrn}\phi(x)\,dx\neq0$.
Assume that $M(f,\phi)|\cdot|^d\in L^1(\rrn)$,
then $\mc_N(f)|\cdot|^d\in L^1(\rrn)$ and
$$
\left\|\mc_N(f)|\cdot|^d
\right\|_{L^1(\rrn)}
\sim\left\|M(f,\phi)|\cdot|^d
\right\|_{L^1(\rrn)}
$$
with the positive equivalence
constants independent of $f$.
\end{lemma}

\begin{proof}
Let all the symbols be as in the present lemma
and $r\in(0,\frac{n}{n+d})$. Then
we have $d<n(\frac{1}{r}-1)$. This,
together with Lemma \ref{power-weight}
with $a:=d$ and $p:=\frac{1}{r}$,
implies that $|\cdot|^{d}\in A_{1/r}(\rrn)$.
Applying this and Lemma \ref{hl-weight}
with $p:=\frac{1}{r}$
and $\upsilon:=|\cdot|^d$, we conclude
that, for any $g\in L^1_{\mathrm{loc}}(\rrn)$,
\begin{equation}\label{max-hardye1}
\left\|\mc(g)\right\|_{
L^{1/r}_{|\cdot|^d}(\rrn)}
\lesssim\|g\|_{L^{1/r}_{|\cdot|^d}(\rrn)}.
\end{equation}
In addition, observe that
$L^{1/r}_{|\cdot|^{d}}(\rrn)=
[L^1_{|\cdot|^{d}}(\rrn)]^{1/r}$.
Thus, from this, \eqref{max-hardye1},
and Lemma \ref{Th5.4l1} with
$X:=L^1_{|\cdot|^{d}}(\rrn)$,
it follows that
$$
\left\|\mc_N(f)|\cdot|^d
\right\|_{L^1(\rrn)}
\sim\left\|M(f,\phi)|\cdot|^d
\right\|_{L^1(\rrn)},
$$
where the positive
equivalence constants are independent
of $f$. This finishes the proof of
Lemma \ref{max-hardy}.
\end{proof}

Based on above preparations,
we now prove Proposition \ref{con-T-s-v}.

\begin{proof}[Proof of Proposition
\ref{con-T-s-v}]
Let all the symbols be as in the present
proposition and $a\in L^2(\rrn)$,
with compact support, satisfying that,
for any $\gamma\in\zp^n$ with
$|\gamma|\leq d$, $$\int_{\rrn}a(x)
x^{\gamma}\,dx=0.$$
Then, from Lemma \ref{con-T-bdd}
with $p=2$, we deduce that
$T(a)\in L^2(\rrn)$. In addition,
let $\phi\in\mathcal{S}(\rrn)$
with $\supp(\phi)\subset B(\0bf,1)$ and
$$\int_{\rrn}\phi(x)\,dx\neq0.$$
By Lemma \ref{stein-hardy} with $\gamma=
\0bf$, we find that $M(T(a),\phi)\in L^1(\rrn)$
and hence $T(a)\in H^1(\rrn)$.
Thus, we have $T(a)\in H^1(\rrn)\cap
L^2(\rrn)$. Choose an $N\in\mathbb{
N}\cap[n+d+1,\infty)$. Then,
applying Lemma \ref{atom-hardy}
with $f$ therein replaced by $T(a)$,
we conclude that there exist $\{\lambda_{i,j}\}
_{i\in\mathbb{Z},\,j\in\mathbb{N}}\subset[0,\infty)$ and
a sequence $\{a_{i,j}\}_{i\in\mathbb{Z},\,
j\in\mathbb{N}}$
of $(L^1(\rrn),\,\infty,\,d)$-atoms such that
\begin{equation}\label{con-T-s-ve1}
T(a)=\sum_{i\in\mathbb{Z}}\sum_{j\in\mathbb{N}}
\lambda_{i,j}a_{i,j}
\end{equation}
almost everywhere in $\rrn$ and
(i), (ii), and (iii) of Lemma \ref{atom-hardy}
hold true. Assume $E\subset\rrn$
such that $|E|=0$ and, for
any $x\in\rrn\setminus E$,
$$T(a)(x)=\sum_{i\in\mathbb{Z}}\sum_{
j\in\mathbb{N}}
\lambda_{i,j}a_{i,j}(x).$$
Now, we show that
$$\sum_{i\in\mathbb{Z},\,j\in\mathbb{N}}
\lambda_{i,j}
\left|a_{i,j}\right|\lesssim\mc_N\left(T
\left(a\right)\right)$$
almost everywhere in $\rrn$ by
considering the following two cases
on $x$.

\emph{Case 1)} $x\in(
\bigcup_{i\in\mathbb{Z}}\Omega
_{i})^{\complement}\setminus E$.
Here and thereafter,
for any $i\in\mathbb{Z}$,
$\Omega_i$ is defined as in Lemma
\ref{atom-hardy}(ii)
with $f$ replaced by $T(a)$.
In this case, we have $\mc_N(T(a))(x)=0$.
In addition, for any $i\in\mathbb{Z}$
and $j\in\mathbb{N}$,
by the assumption
$\supp(a_{i,j})\subset B_{i,j}$ and
Lemma \ref{atom-hardy}(ii),
we conclude that $a_{i,j}(x)=0$.
This, together with $\mc_N(T(a))(x)=0$,
further implies that
\begin{equation}\label{con-T-s-ve2}
\sum_{i\in\mathbb{Z},\,j\in\mathbb{N}}
\lambda_{i,j}
\left|a_{i,j}(x)\right|
=0=\mc_N\left(T(a)\right)(x).
\end{equation}
Thus, $$\sum_{i\in\mathbb{Z},\,j\in\mathbb{N}}
\lambda_{i,j}
\left|a_{i,j}\right|\lesssim\mc_N\left(T
\left(a\right)\right)$$
almost everywhere in $(
\bigcup_{i\in\mathbb{Z}}\Omega
_{i})^{\complement}$.

\emph{Case 2)} $x\in(\Omega_{i_0}\setminus
\Omega_{i_0+1})\setminus E$ for
some $i_0\in\mathbb{Z}$. In this case,
for any $i\in\mathbb{N}\cap[i_0+1,\infty)$
and $j\in\mathbb{N}$,
we have $a_{i,j}(x)=0$.
Applying this, both (i) and
(iii) of Lemma \ref{atom-hardy},
the assumption $\supp(a_{i,j})
\subset B_{i,j}$ for any
$i\in\mathbb{Z}$ and $j\in\mathbb{N}$,
and the definition of $\Omega_i$,
we find that, for any
$x\in(\Omega_{i_0}\setminus
\Omega_{i_0+1})\setminus E$,
\begin{align}\label{con-T-s-ve3}
\sum_{i\in\mathbb{Z},\,j\in\mathbb{N}}
\lambda_{i,j}\left|a_{i,j}(x)\right|
&=\sum_{i=-\infty}^{i_0}\sum_{
j\in\mathbb{N}}\lambda_{i,j}
\left|a_{i,j}(x)\right|\notag\\
&\lesssim\sum_{i=-\infty}^{i_0}
2^{i}\sum_{j\in\mathbb{N}}\1bf_{B_{i,j}}
\lesssim\sum_{i=-\infty}^{i_0}
2^i\sim2^{i_0}\notag\\&\sim\mc_N\left(T(a)\right)(x),
\end{align}
which further implies that
$$\sum_{i\in\mathbb{Z},\,j\in\mathbb{N}}
\lambda_{i,j}
\left|a_{i,j}\right|\lesssim\mc_N\left(
T\left(a\right)\right)$$
almost everywhere in $\Omega_{i_0}\setminus
\Omega_{i_0+1}$. Combining
\eqref{con-T-s-ve2} and \eqref{con-T-s-ve3},
we then conclude that
$$\sum_{i\in\mathbb{Z},\,j\in\mathbb{N}}
\lambda_{i,j}\left|a_{i,j}\right|\lesssim\mc_N(T(a))$$
almost everywhere in $\rrn$.

Therefore, for any $\gamma\in\zp^n$
with $|\gamma|\leq d$, from both Lemmas
\ref{stein-hardy} and
\ref{max-hardy} with $d$ therein
replaced by $|\gamma|$,
we deduce that
$$
\sum_{i\in\mathbb{Z},\,j\in\mathbb{N}}
\lambda_{i,j}\left|a_{i,j}\right||
\cdot|^{|\gamma|}\lesssim
\mc_N(T(a))|\cdot|^{|\gamma|}\in
L^1(\rrn).
$$
Using this, \eqref{con-T-s-ve1},
the dominated convergence theorem, and
the assumption that
$\{a_{i,j}\}_{i\in\mathbb{Z},\,
j\in\mathbb{N}}$ is a
sequence of $(L^1(\rrn),\,\infty,\,d)$-atoms,
we further
find that, for any $\gamma\in\zp^n$
with $|\gamma|\leq d$,
$$
\int_{\rrn}T(a)(x)x^{\gamma}\,dx=
\sum_{i\in\mathbb{Z}}\sum_{j\in\mathbb{N}}
\int_{\rrn}a_{i,j}(x)x^{\gamma}\,dx=0,
$$
which completes the proof of Proposition
\ref{con-T-s-v}.
\end{proof}

Under the reasonable assumption
that the Calder\'{o}n--Zygmund operator
$T$ has vanishing moments, we
then have the boundedness of
Calder\'{o}n--Zygmund operators
on the generalized Herz--Hardy space
$\HaSaHo$ as follows.

\begin{theorem}\label{bddol}
Let $d\in\zp$, $\delta\in(0,1]$,
$p$, $q\in(\frac{n}{n+d+\delta},\infty)$,
$K$ be a $d$-order standard kernel defined
as in Definition \ref{def-s-k}, $T$ a
$d$-order Calder\'{o}n--Zygmund operator
with kernel $K$ having the vanishing moments
up to order $d$, and $\omega\in
M(\rp)$ with
$$
-\frac{n}{p}<\m0(\omega)\leq
\M0(\omega)<n-\frac{n}{p}+d+\delta
$$
and
$$
-\frac{n}{p}<\mi(\omega)\leq
\MI(\omega)<n-\frac{n}{p}+d+\delta.
$$
Then $T$ has a unique extension
on $\HaSaHo$ and there exists
a positive constant $C$ such that,
for any $f\in\HaSaHo$,
$$
\|T(f)\|_{\HaSaHo}\leq C\|f\|_{
\HaSaHo}.
$$
\end{theorem}

\begin{remark}
We should point out that, in Theorem
\ref{bddol}, when $d=0$, $p\in(1,\infty)$,
and $\omega(t):=t^{\alpha}$ for
any $t\in(0,\infty)$
and for any given $\alpha\in[n(1-\frac{1}{p}),
n(1-\frac{1}{p})+\delta)$,
then Theorem \ref{bddol}
goes back to \cite[Theorem 1]{ll97}.
\end{remark}

To show Theorem \ref{bddol},
we first establish the
boundedness of Calder\'{o}n--Zygmund
operators on Hardy spaces
associated with ball quasi-Banach
function spaces as follows.

\begin{proposition}\label{bddoxl4}
Let $X$ be a ball quasi-Banach
function space satisfying Assumption
\ref{assfs} for some $0<\theta<s\leq1$.
Assume that $X^{1/s}$ is a
ball Banach function space and
there exist an $r_0\in(1,\infty)$ and
a positive
constant $\widetilde{C}$ such that,
for any $f\in L^1_{\mathrm{loc}}(\rrn)$,
\begin{equation}\label{bddoxl4e0}
\left\|\mc^{((r_0/s)')}(f)
\right\|_{(X^{1/s})'}\leq\widetilde{C}
\left\|f\right\|_{(X^{1/s})'}.
\end{equation}
Let $d\in\zp$, $\delta\in(0,1]$, $K$ be a $d$-order
standard kernel defined as in Definition
\ref{def-s-k}, and $T$ a $d$-order
Calder\'{o}n--Zygmund operator
with kernel $K$ having the vanishing moments
up to order $d$.
If $\theta\in(\frac{n}{n+d+\delta},
\frac{n}{n+d}]$ and $X$ has
an absolutely continuous quasi-norm,
then $T$ has a unique extension on $H_X(\rrn)$
and there exists a positive constant $C$
such that, for any $f\in H_X(\rrn)$,
$$
\|T(f)\|_{H_X(\rrn)}\leq C\|f\|_{H_X(\rrn)}.
$$
\end{proposition}

To show Proposition \ref{bddoxl4},
we need the following auxiliary
lemma which established the
relations among
Calder\'{o}n--Zygmund operators,
atoms, and molecules.

\begin{lemma}\label{bddoxl3}
Let $X$ be a ball quasi-Banach
function space,
$r\in[2,\infty)$, and $d\in\zp$.
Assume that $K$ is a $d$-order standard
kernel defined as in Definition
\ref{def-s-k} with some $\delta\in(0,1]$, and $T$ a
$d$-order Calder\'{o}n--Zygmund operator
with kernel $K$ having
the vanishing moments up to order $d$.
Then, for any $(X,\,r,\,d)$-atom $a$
supported in the ball $B\in\mathbb{B}$,
$T(a)$ is a harmless constant multiple of
an $(X,\,r,\,d,\,d+\delta+\frac{n}{r'})$-molecule
centered at $B$.
\end{lemma}

\begin{proof}
Let all the symbols be as in the
present lemma and
$a$ an $(X,\,r,\,d)$-atom
supported in the ball $B:=B(x_0,r_0)$
with $x_0\in\rrn$ and $r_0\in(0,\infty)$.
Then, combining Definitions
\ref{atomx}(iii) and \ref{Def-T-s-v},
we find that, for any $\gamma\in\zp^n$
with $|\gamma|\leq d$,
\begin{equation*}
\int_{\rrn}T(a)(x)x^{\gamma}\,dx=0.
\end{equation*}
This implies that $T(a)$
satisfies Definition \ref{molex}(ii).

Next, we show that Definition \ref{molex}(i)
holds true for a harmless constant
multiple of $T(a)$. Indeed, from
Lemma \ref{czop} with $p:=r$ and Definition
\ref{atomx}(ii), it follows that
\begin{align}\label{bddoxl3e2}
\left\|T(a)\1bf_{B}\right\|_{L^r(\rrn)}
\leq\left\|T(a)\right\|_{L^r(\rrn)}
\lesssim\|a\|_{L^r(\rrn)}\lesssim
\frac{|B|^{1/r}}{\|\1bf_B\|_X}
\end{align}
and
\begin{align}\label{bddoxl3e6}
&\left\|T(a)\1bf_{(2B)\setminus
B}\right\|_{L^r(\rrn)}\notag\\
&\quad\leq\left\|T(a)\right\|_{L^r(\rrn)}
\lesssim\|a\|_{L^r(\rrn)}\lesssim
\frac{|B|^{1/r}}{\|\1bf_B\|_X}
\sim2^{-(d+\delta+\frac{n}{r'})}\frac{
|B|^{1/r}}{\|\1bf_B\|_X}.
\end{align}
These are the desired estimates of
$\|T(a)\1bf_{B}\|_{L^r(\rrn)}$ and
$\|T(a)\1bf_{(2B)\setminus B}\|_{L^r(\rrn)}$,
respectively.

In addition, recall that, for any
$j\in\mathbb{N}$, $S_j(B):=
(2^{j}B)\setminus(2^{j-1}B)$.
We next estimate $\|T(a)\1bf_
{S_{j+1}(B)}\|_{L^r(\rrn)}$. Indeed,
by Definition \ref{defin-C-Z-s}(ii), we
find that,
for any $j\in\mathbb{N}$
and $x\in S_{j+1}(B)$,
\begin{equation*}
T(a)(x)=\int_{\rrn}K(x,y)a(y).
\end{equation*}
This, together with
Definition \ref{atomx}(iii) and the Taylor
remainder theorem,
further implies that, for any
$j\in\mathbb{N}$ and $y\in B$, there
exists a $t_y\in[0,1]$ such that
\begin{align}\label{bddoxl3e5}
&\left\|T(a)\1bf_{S_{j+1}(B)}
\right\|_{L^r(\rrn)}\notag\\
&\quad=\left[
\int_{S_{j+1}(B)}\left|
\int_{B}K(x,y)a(y)\,dy
\right|^r\,dx
\right]^{\frac{1}{r}}\notag\\
&\quad=\left\{
\int_{S_{j+1}(B)}\left|
\int_{B}\left[K(x,y)-
\sum_{\gfz{\gamma\in\zp^n}{|\gamma|\leq d}}
\frac{\partial^{\gamma}_{(2)}K(x,x_0)}{\gamma!}
(y-x_0)^{\gamma}\right]a(y)\,dy
\right|^r\,dx
\right\}^{\frac{1}{r}}\notag\\
&\quad=\left\{
\int_{S_{j+1}(B)}\left|
\int_{B}\sum_{\gfz{\gamma\in\zp^n}
{|\gamma|=d}}\frac{\partial_{(2)}^{\gamma}
K(x,t_yy+(1-t_y)x_0)-\partial^{\gamma}_{(2)}
K(x,x_0)}{\gamma!}\right.
\right.\notag\\
&\qquad\times\lf.
(y-x_0)^{\gamma}a(y)\,dy
\Bigg|^r\,dx
\right\}^{\frac{1}{r}}\notag\\
&\quad\lesssim
\int_{B}|y-x_0|^da(y)\left[\int_{S_{j+1}(B)}\r.
\notag\\
&\qquad\times\lf.
\sum_{\gfz{\gamma\in\zp^n}{|\gamma|=d}}\left|
\partial_{(2)}^{\gamma}K(x,t_yy+(1-t_y)x_0)-
\partial_{(2)}^{\gamma}K(x,x_0)\right|^r
\,dx\right]^{\frac{1}{r}}\,dy.
\end{align}
On the other hand,
for any $j\in\mathbb{N}$,
$x\in S_{j+1}(B)$, and
$y\in B$, we have
$$
|x-x_0|\geq2^jr_0>2|y-x_0|\geq2|t_y(y-x_0)|.
$$
Using this, \eqref{bddoxl3e5}, \eqref{regular2-s}
with $y$ and $z$ therein replaced, respectively,
by $x_0$ and $t_yy+(1-t_y)x_0$ for any $y\in B$,
the H\"{o}lder inequality, and Definition
\ref{atomx}(ii), we conclude that, for any
$j\in\mathbb{N}$,
\begin{align}\label{bddoxl3e7}
&\left\|T(a)\1bf_{S_{j+1}(B)}\right\|_{L^r(\rrn)}
\notag\\&\quad\lesssim\int_{B}r_0^d|a(y)|\left[
\int_{S_{j+1}(B)}\frac{|y-x_0|^{\delta r}}
{|x-x_0|^{(n+d+\delta)r}}\,dx
\right]^{\frac{1}{r}}\,dy\notag\\
&\quad\lesssim r_0^{d+\delta}(2^jr_0)^{-n+d+\delta}
\left|S_{j+1}(B)\right|^{\frac{1}{r}}
\int_{B}|a(y)|\,dy\notag\\
&\quad\lesssim2^{-j(d+\delta+\frac{n}{r'})}
\left\|a\right\|_{L^r(\rrn)}
\lesssim2^{-j(d+\delta+\frac{n}{r'})}\frac{
|B|^{1/r}}{\|\1bf_B\|_X}\notag\\
&\quad\sim2^{-(j+1)(d+\delta+\frac{n}{r'})}
\frac{|B|^{1/r}}{\|\1bf_B\|_X},
\end{align}
which is the desired estimate of
$\|T(a)\1bf_{S_{j+1}(B)}\|_
{L^r(\rrn)}$ with $j\in\mathbb{N}$.
Applying \eqref{bddoxl3e2}, \eqref{bddoxl3e6},
and \eqref{bddoxl3e7}, we find that
there exists a positive constant $C$,
independent of $a$, such that, for any
$j\in\zp$,
$$
\|CT(a)\1bf_{S_j(B)}\|_{L^r(\rrn)}\leq
2^{-j(d+\delta+\frac{n}{r'})}
\frac{|B|^{1/r}}{\|\1bf_B\|_X}.
$$
This further implies that $CT(a)$ satisfies
Definition \ref{molex}(ii) and hence
$CT(a)$ is an $(X,\,r,\,d,\,d+\delta+\frac{n}
{r'})$-molecule centered at $B$, which
then completes
the proof of Lemma \ref{bddoxl3}.
\end{proof}

Via Lemma \ref{bddoxl3},
we next show Proposition \ref{bddoxl4}.

\begin{proof}[Proof of Proposition
\ref{bddoxl4}]
Let all the symbols be as in
the present proposition and
$r:=\max\{2,r_0\}$.
Then, from the
H\"{o}lder inequality, it follows that,
for any $f\in L^1_{\mathrm{loc}}(\rrn)$,
$x\in\rrn$, and $B\in\mathbb{B}$ satisfying
$x\in B$,
\begin{align*}
\left[\frac{1}{|B|}
\int_{B}|f(y)|^{(r/s)'}
\,dy\right]^{\frac{1}{(r/s)'}}
&\leq\left[\frac{1}{|B|}
\int_{B}|f(y)|^{(r_0/s)'}
\,dy\right]^{\frac{1}{(r_0/s)'}}\\
&\leq\mc^{((r_0/s)')}(f)(x),
\end{align*}
which further implies that
$$
\mc^{((r/s)')}(f)(x)\leq
\mc^{((r_0/s)')}(f)(x).
$$
Using this, \eqref{bddoxl4e0},
the assumption that $X^{1/s}$
is a BBF space, Remark \ref{remark127},
and Definition \ref{Df1}(ii),
we find that, for any
$f\in L^1_{\mathrm{loc}}(\rrn)$,
\begin{equation}\label{bddoxl4e1}
\left\|\mc^{((r/s)')}(f)
\right\|_{(X^{1/s})'}\leq
\left\|\mc^{((r_0/s)')}(f)
\right\|_{(X^{1/s})'}\lesssim
\left\|f\right\|_{(X^{1/s})'}.
\end{equation}
On the other hand, by the assumptions
$\theta\in(\frac{n}{n+d+\delta},
\frac{n}{n+d}]$ and $\delta\in(0,1]$,
we conclude that
$$
d>n\left(\frac{1}{\theta}-1\right)
-\delta\geq n\left(\frac{1}{\theta}-1
\right)-1
$$
and
$$
d\leq n\left(\frac{1}{\theta}-1\right).
$$
These imply that $d=\lfloor n(
1/\theta-1)\rfloor$.
From this, \eqref{bddoxl4e1}, the
assumption that $X$ satisfies Assumption
\ref{assfs} for the above $\theta$ and $s$,
and Lemma \ref{finiatomll1},
we deduce that, for any
$g\in H^{X,r,d,s}_{\mathrm{fin}}(\rrn)$,
\begin{equation}\label{bddoxl4e2}
\|g\|_{H^{X,r,d,s}_{\mathrm{fin}}(\rrn)}
\sim\|g\|_{H_X(\rrn)}
\end{equation}
with the positive equivalence constants
independent of $g$.

Now, let $f\in
H^{X,r,d,s}_{\mathrm{fin}}(\rrn)$,
$m\in\mathbb{N}$,
$\{\lambda_{j}\}_{j=1}^m\subset
[0,\infty)$, and $\{a_j\}_{j=1}^m$
of $(X,\,r,\,d)$-atoms supported,
respectively, in the balls $\{
B_j\}_{j=1}^m\subset\mathbb{B}$
such that
$$f=\sum_{j=1}^{m}\lambda_ja_j.$$
This, combined with the linearity
of $T$, implies that
\begin{equation}\label{bddoxl4e3}
T(f)=\sum_{j=1}^{m}T(a_j).
\end{equation}
On the other hand, from
Lemma \ref{bddoxl3}, it follows that,
for any $j\in\{1,\ldots,m\}$,
$T(a_j)$ is a harmless constant
multiple of
an $(X,\,r,\,d,\,d+\delta+
\frac{n}{r'})$-molecule centered at $B_j$.
Moreover, by the assumption
$\theta>\frac{n}{n+d+\delta}$, we find
that
$$d+\delta+\frac{n}{r'}>n\lf(\frac{1}{\theta}-
\frac{1}{r}\r).$$
Combining this, the assumption that
$X$ satisfies Assumption \ref{assfs}
for the above $\theta$ and $s$,
\eqref{bddoxl4e1}, \eqref{bddoxl4e3},
the fact that, for any
$j\in\{1,\ldots,m\}$, $T(a_j)$ is
a harmless constant multiple of
an $(X,\,r,\,d,\,d+\delta+
\frac{n}{r'})$-molecule centered at $B_j$, and
Lemma \ref{Molell1}, we further
conclude that $T(f)\in H_X(\rrn)$
and
$$
\|T(f)\|_{H_X(\rrn)}
\lesssim\left\|
\left[\sum_{j=1}^{m}\left(\frac{\lambda_j}
{\|\1bf_{B_j}\|_X}\right)^s\1bf_{B_j}
\right]^{\frac{1}{s}}\right\|_{X}.
$$
This, together with the choice of
$\{\lambda_j\}_{j=1}^m$, \eqref{fi},
and \eqref{bddoxl4e2},
implies that
\begin{equation*}
\|T(f)\|_{H_X(\rrn)}
\lesssim\|f\|_{H^{X,r,d,s}_{\mathrm{fin}}(\rrn)}
\sim\|f\|_{H_X(\rrn)}.
\end{equation*}
Therefore, $T$ is bounded on
the finite atomic Hardy space $
H^{X,r,d,s}_{\mathrm{fin}}(\rrn)$.

Finally, from the assumption that
$X$ has an absolutely continuous
quasi-norm and \cite[Remark 3.12]{SHYY},
we deduce that
the finite atomic Hardy space $H^{X,r,d,s}_
{\mathrm{fin}}(\rrn)$ is dense in the
Hardy space $H_X(\rrn)$. Thus, by a
standard density argument, we find that
$T$ has a unique extension on
$H_X(\rrn)$ and, for any $f\in H_X(\rrn)$,
$$
\|T(f)\|_{H_X(\rrn)}\lesssim\|f\|_{H_X(\rrn)},
$$
which completes the proof of Proposition
\ref{bddoxl4}.
\end{proof}

Via the above boundedness
of Calder\'{o}n--Zygmund operators
on the Hardy space $H_X(\rrn)$,
we next prove Theorem \ref{bddol}.

\begin{proof}[Proof of Theorem
\ref{bddol}]
Let all the symbols be as in the
present theorem.
Then, from the assumption
$\m0(\omega)\in(-\frac{n}{p},\infty)$
and Theorem \ref{Th3}, it follows that
the local generalized
Herz space $\HerzSo$ is a
BQBF space. We now show Theorem \ref{bddol}
via proving that all the assumptions
of Proposition \ref{bddoxl4}
hold true for $\HerzSo$.

First, we show that there exist
$\theta$, $s\in(0,1]$ such that Assumption
\ref{assfs} holds true for $\HerzSo$.
Indeed, by the
assumption
$\Mw\in(-\frac{n}{p},
n-\frac{n}{p}+d+\delta)$,
we conclude that
$$\frac{n}{\Mw+n/p}
\in\left(\frac{n}{n+d+\delta},
\infty\right).$$
This, combined with the
assumptions $p$, $q\in(n/(n+d+\delta),\infty)$,
further implies that
$$
\min\left\{1,p,q,\frac{n}
{\Mw+n/p}\right\}\in\left(
\frac{n}{n+d+\delta},\infty\right).
$$
Therefore, we can choose an
$$s\in\left(\frac{n}{n+d+\delta},\,
\min\left\{1,p,q,\frac{n}
{\Mw+n/p}\right\}\right)$$
and a $$\theta\in\left(\frac{n}{n+d+\delta},\,
\min\left\{s,\frac{n}{n+d}\right\}\right).$$
Then, applying Lemma \ref{Atoll3},
we find that,
for any $\{f_{j}\}_
{j\in\mathbb{N}}\subset L^{1}_
{\text{loc}}(\rrn)$,
\begin{equation}\label{czoe2}
\left\|\left\{\sum_{j\in\mathbb{N}}
\left[\mc^{(\theta)}
(f_{j})\right]^s\right\}^{1/s}
\right\|_{\HerzSo}\lesssim
\left\|\left(\sum_{j\in\mathbb{N}}
|f_{j}|^{s}\right)^{1/s}\right\|_{\HerzSo}.
\end{equation}
This implies that, for the above
$\theta$ and $s$, $\HerzSo$ satisfies
Assumption \ref{assfs}.

Next, from
Lemma \ref{mbhal}, we deduce that
the Herz space $[\HerzSo]^{1/s}$ is
a BBF space.
In addition, let
$$r_0\in\left(\max\left\{1,p,
\frac{n}{\mw+n/p}\right\},\infty\right).$$
We now show that, for this $r_0$ and the
above $s$, \eqref{Atoll1e1} holds true.
Indeed, using Lemma \ref{mbhal} with
$r:=r_0$, we conclude that,
for any $f\in L_{{\rm loc}}^{1}(\rrn)$,
\begin{equation}\label{czoe1}
\left\|\mc^{((r_0/s)')}(f)
\right\|_{([\HerzSo]^{1/s})'}\lesssim
\left\|f\right\|_{([\HerzSo]^{1/s})'}.
\end{equation}

Finally, from Theorem \ref{abso},
it follows that the local
generalized Herz space $\HerzSo$
has an absolutely continuous quasi-norm.
Combining this,
\eqref{czoe2}, the fact that
$[\HerzSo]^{1/s}$ is a BBF space, and
\eqref{czoe1}, we find that
all the assumptions of Lemma \ref{bddoxl4}
hold true for $\HerzSo$
under consideration and then complete
the proof of Theorem \ref{bddol}.
\end{proof}

As an application, we have the following
boundedness of Calder\'{o}n--Zygmund
operators on the generalized Morrey--Hardy space
$H\MorrSo$, which can be deduced
from Theorem \ref{bddol} and Remarks \ref{remhs}(iv)
and \ref{remark4.10}(ii) immediately; we omit the details.

\begin{corollary}\label{bddolm}
Let $d\in\zp$, $p$, $q\in[1,\infty)$,
$K$ be a $d$-order standard kernel defined
as in Definition \ref{def-s-k}
with some $\delta\in(0,1]$, $T$ a
$d$-order Calder\'{o}n--Zygmund operator
with kernel $K$ having the vanishing moments
up to order $d$,
and $\omega\in
M(\rp)$ with
$$
-\frac{n}{p}<\m0(\omega)\leq
\M0(\omega)<0
$$
and
$$
-\frac{n}{p}<\mi(\omega)\leq
\MI(\omega)<0.
$$
Then $T$ has a unique extension
on $H\MorrSo$ and there exists
a positive constant $C$ such that,
for any $f\in H\MorrSo$,
$$
\|T(f)\|_{H\MorrSo}\leq C\|f\|_{
H\MorrSo}.
$$
\end{corollary}

The remainder of this section
is devoted to showing the
boundedness of
Calder\'{o}n--Zygmund operators
on the generalized Herz--Hardy
space $\HaSaH$. Precisely,
we turn to prove the following theorem.

\begin{theorem}\label{bddog}
Let $d\in\zp$, $\delta\in(0,1]$,
$p$, $q\in(\frac{n}{n+d+\delta},\infty)$,
$K$ be a $d$-order standard kernel with
defined as in Definition \ref{def-s-k}, $T$ a
$d$-order Calder\'{o}n--Zygmund operator
with kernel $K$ having the vanishing moments
up to order $d$, and $\omega\in
M(\rp)$ with
$$
-\frac{n}{p}<\m0(\omega)\leq
\M0(\omega)<n-\frac{n}{p}+d+\delta
$$
and
$$
-\frac{n}{p}<\mi(\omega)\leq
\MI(\omega)<0.
$$
Then $T$ can be extended into a bounded
linear operator on $\HaSaH$ and
there exists
a positive constant $C$ such that,
for any $f\in\HaSaH$,
$$
\|T(f)\|_{\HaSaH}\leq C\|f\|_{
\HaSaH}.
$$
\end{theorem}

To show this
theorem, we first establish
a general result about the
boundedness of Calder\'{o}n--Zygmund
operators on Hardy spaces
associated with ball quasi-Banach function
spaces as follows.

\begin{theorem}\label{bddox}
Let $X$ be a ball quasi-Banach function
space, $Y$ a linear space equipped with a
quasi-seminorm $\|\cdot\|_Y$,
$Y_0$ a linear space equipped with
a quasi-seminorm $\|\cdot\|_{Y_0}$,
$\eta\in(1,\infty)$,
and $0<\theta<s<s_0\leq1$ such that
\begin{enumerate}
  \item[{\rm(i)}] for the above $\theta$ and
  $s$, Assumption \ref{assfs} holds true;
  \item[{\rm(ii)}] both $\|\cdot\|_{Y}$
  and $\|\cdot\|_{Y_0}$
  satisfy Definition \ref{Df1}(ii);
  \item[{\rm(iii)}] $\1bf_{B(\0bf,1)}
  \in Y_0$;
  \item[{\rm(iv)}] for any $f\in\Msc(\rrn)$,
  $$
  \|f\|_{X^{1/s}}
  \sim\sup\left\{\|fg\|_{L^1(\rrn)}:\
  \|g\|_{Y}=1\right\}
  $$
  and
  $$
  \|f\|_{X^{1/s_{0}}}
  \sim\sup\left\{\|fg\|_{L^1(\rrn)}:\
  \|g\|_{Y_0}=1\right\}
  $$
  with the positive equivalence constants
  independent of $f$;
  \item[{\rm(v)}] $\mc^{(\eta)}$ is bounded
  on both $Y$ and $Y_0$.
\end{enumerate}
Assume that $d\in\zp$, $\delta\in(0,1]$,
$K$ is a $d$-order standard kernel
defined as in Definition \ref{def-s-k},
and $T$ a $d$-order Calder\'{o}n--Zygmund
operator with kernel $K$
having the vanishing moments up to $d$.
If $\theta\in(\frac{n}{n+d+\delta},
\frac{n}{n+d}]$, then $T$ can be
extended into a bounded linear operator
on $H_X(\rrn)$, namely, there exists a positive
constant $C$ such that,
for any $f\in H_X(\rrn)$,
$$
\|T(f)\|_{H_X(\rrn)}\leq C\|f\|_{H_X(\rrn)}.
$$
\end{theorem}

To show Theorem
\ref{bddox},
the following two embedding theorems
about $X$ and $H_X(\rrn)$ are
the essential tools.

\begin{lemma}\label{bddol1}
Let $X$ be a ball quasi-Banach
function space, $Y\subset\Msc(\rrn)$ a
linear space equipped with a
quasi-seminorm $\|\cdot\|_Y$,
$\theta\in(1,\infty)$, and $s\in(0,\infty)$
satisfy the following four statements:
\begin{enumerate}
  \item[{\rm(i)}] $\|\cdot\|_Y$ satisfies
  Definition \ref{Df1}(ii);
  \item[{\rm(ii)}] $\1bf_{B(\0bf,1)}\in Y$;
  \item[{\rm(iii)}] for any $f\in\Msc(\rrn)$,
  $$
  \|f\|_{X^{1/s}}
  \sim\sup\left\{\|fg\|_{L^1(\rrn)}:\
  \|g\|_{Y}=1\right\},
  $$
  where the positive equivalence
  constants are independent of $f$;
  \item[{\rm(iv)}] $\mc^{(\theta)}$ is
  bounded on $Y$.
\end{enumerate}
Assume $\eps\in(\frac{1}{\theta},1)$ and
$\upsilon:=[\mc(\1bf_{B(\0bf,1)}
)]^{\eps}$.
Then there exists a positive
constant $C$
such that, for any $f\in X$,
$$
\|f\|_{L^s_{\upsilon}(\rrn)}
\leq C\|f\|_{X}.
$$
\end{lemma}

\begin{proof}
Let all the symbols be as
in the present lemma and
$f\in X$.
Then, by Lemma \ref{2.1.6}
with $r:=1$, we find
that, for any $x\in\rrn$,
\begin{equation*}
\upsilon(x):=\left[\mc
\left(\1bf_{B(\0bf,1)}
\right)\right]^{\eps}\sim
\left(1+|x|\right)^{-\eps n}.
\end{equation*}
This implies that, for any
$x\in B(\0bf,1)$,
$
\upsilon(x)\sim1
$
and, for any $k\in\mathbb{N}$
and $x\in B(\0bf,\tk)
\setminus B(\0bf,\tkm)$,
$
\upsilon(x)\sim2^{-k\eps n}.
$
From these, we further deduce that
\begin{align}\label{bddol1e1}
&\int_{\rrn}|f(y)|^{s}
\upsilon(y)\,dy\notag\\
&\quad=\int_{B(\0bf,1)}|f(y)|^{s}
\upsilon(y)\,dy+
\sum_{k\in\mathbb{N}}
\int_{B(\0bf,\tk)\setminus
B(\0bf,\tkm)}|f(y)|^{s}
\upsilon(y)\,dy\notag\\
&\quad\sim\int_{B(\0bf,1)}|f(y)|^{s}\,dy+
\sum_{k\in\mathbb{N}}2^{-k\eps n}
\int_{B(\0bf,\tk)\setminus
B(\0bf,\tkm)}|f(y)|^{s}\,dy\notag\\
&\quad\lesssim\int_{\rrn}|f(y)|^{s}
\1bf_{B(\0bf,1)}(y)\,dy
+\sum_{k\in\mathbb{N}}2^{-k\eps n}
\int_{\rrn}|f(y)|^{s}
\1bf_{B(\0bf,\tk)}(y)\,dy.
\end{align}

Next, applying \eqref{hlmaxp},
we conclude that,
for any $k\in\mathbb{N}$ and
$x\in B(\0bf,\tk)$,
\begin{align*}
\mc^{(\theta)}\left(
\1bf_{B(\0bf,1)}\right)(x)
&\geq\left\{\frac{1}{|B(\0bf,\tk)|}
\int_{B(\0bf,\tk)}
\left[\1bf_{B(\0bf,1)}(y)\right]
^{\theta}\,dy\right\}
^{\frac{1}{\theta}}\\
&\sim\left[
\frac{|B(\0bf,1)|}{
|B(\0bf,\tk)|}\right]^{\frac{1}{\theta}}
\sim2^{-\frac{nk}{\theta}},
\end{align*}
which implies that
\begin{equation*}
\1bf_{B(\0bf,\tk)}
\lesssim2^{\frac{nk}{\theta}}\mc^{(
\theta)}\left(\1bf_{B(\0bf,1)}\right),
\end{equation*}
where the implicit positive
constant is independent of $k$.
Applying this and both the assumptions
(i) and (iv) of the present lemma,
we find that, for any $k\in\mathbb{N}$,
\begin{align}\label{bddol1e7}
\left\|\1bf_{B(\0bf,\tk)}\right\|_{Y}
\lesssim2^{\frac{nk}{\theta}}
\left\|\mc^{(\theta)}\left(
\1bf_{B(\0bf,1)}\right)\right\|_{Y}
\lesssim2^{\frac{nk}{\theta}}
\left\|\1bf_{B(\0bf,1)}\right\|_{Y},
\end{align}
which is the desired estimate of $\1bf_{B(\0bf,\tk)}$.
Thus, using \eqref{bddol1e1},
the assumption (iii) of the present lemma,
\eqref{bddol1e7},
and the assumption
$\eps-\frac{1}{\theta}\in(0,\infty)$,
we find that
\begin{align}\label{bddol1e4}
&\int_{\rrn}|f(y)|^{s}\upsilon(y)\,dy\notag\\
&\quad\lesssim\left\|\,|f|^{s}
\,\right\|_{X^{1/s}}\left[
\left\|\1bf_{B(\0bf,1)}\right\|_{Y}
+\sum_{k\in\mathbb{N}}2^{-k\eps n}
\left\|\1bf_{B(\0bf,\tk)}\right\|_{Y}
\right]\notag\\
&\quad\lesssim\|f\|_{X}^{s}
\left\|\1bf_{B(\0bf,1)}
\right\|_{Y}\left[1
+\sum_{k\in\mathbb{N}}2^{-kn(\eps-\frac{1}{
\theta})}\right]\notag\\
&\quad\sim\|f\|_{X}^{s}
\left\|\1bf_{B(\0bf,1)}
\right\|_{Y}.
\end{align}
On the other hand, by
the assumption (ii) of the present lemma,
we conclude that $$\left\|\1bf_{B(\0bf,1)}
\right\|_{Y}<\infty.$$
Therefore, from \eqref{bddol1e4}, we
infer that
$$
\|f\|_{L_{\upsilon}^{s}(\rrn)}^{s}
=\int_{\rrn}|f(y)|^{s}\upsilon(y)\,dy
\lesssim\|f\|_{X}^{s},
$$
which completes the proof of
Lemma \ref{bddol1}.
\end{proof}

\begin{lemma}\label{bddol2}
Let $X$ be a ball quasi-Banach
function space and $N\in\mathbb{N}$.
Then the Hardy space
$H_X(\rrn)$ embeds continuously
into $\mathcal{S}'(\rrn)$. Namely,
there exists a positive constant
$C$ such that, for any $f\in H_X(\rrn)$
and $\phi\in\mathcal{S}(\rrn)$,
$$
|\langle f,\phi\rangle|\leq
Cp_N(\phi)\|f\|_{H_X(\rrn)},
$$
where $p_N$ is defined as in
\eqref{defp}.
\end{lemma}

\begin{proof}
Let all the symbols be as
in the present lemma, $f\in H_X(\rrn)$,
and $\phi\in\mathcal{S}(\rrn)$.
Then we prove the present lemma by considering
the following two cases on $\phi$.

\emph{Case 1)} $p_N(\phi)=0$. In this
case, we have
$\phi=0$. This further implies that
\begin{equation}\label{bddol2e1}
\left|\left\langle
f,\phi\right\rangle\right|
=0=p_N(\phi)\|f\|_{H_X(\rrn)},
\end{equation}
which completes the proof of
Lemma \ref{bddol2} in this case.

\emph{Case 2)} $p_N(\phi)\neq0$. In this case,
from \eqref{defp}, we deduce that
$\frac{\phi(-\cdot)}{p_N(\phi)}
\in\mathcal{S}(\rrn)$ and
$$p_N\left(\frac{\phi(-\cdot)}
{p_N(\phi)}\right)=\frac{
p_N(\phi)}{p_N(\phi)}=1.$$
Therefore, $\frac{\phi(-\cdot)}{p_N(\phi)}
\in\mathcal{F}_N(\rrn)$ with $\mathcal{F}_N(\rrn)$
as in \eqref{defn}.
By this and \eqref{sec6e1}, we conclude that,
for any $x\in B(\0bf,1)$,
\begin{align*}
\left|\left\langle f,\phi
\right\rangle\right|
&=\left|f\ast\phi(-\cdot)(\0bf)
\right|\\
&=p_N(\phi)\left|f\ast
\left[\frac{\phi(-\cdot)}{p_N(\phi)}\right](
\0bf)\right|\leq p_N(\phi)\mc_{N}(f)(x).
\end{align*}
This further implies that
\begin{equation}\label{bddol2e2}
\left|\left\langle f,\phi
\right\rangle\right|\1bf_{B(\0bf,1)}
\leq p_N(\phi)\mc_{N}(f).
\end{equation}
In addition, since $X$ is
a BQBF space, from both (i)
and (iv) of Definition \ref{Df1}, we deduce that
$\|\1bf_{B(\0bf,1)}\|_{X}\in(0,\infty)$.
This, together with \eqref{bddol2e2}
and Definitions \ref{Df1}(ii) and \ref{hardyx},
further implies that
\begin{align}\label{bddol2e3}
\left|\left\langle
f,\phi\right\rangle\right|
&\leq\left\|\1bf_{B(\0bf,1)}
\right\|_{X}^{-1}p_N(\phi)
\left\|\mc_N(f)\right\|_{X}\notag\\
&=\left\|\1bf_{B(\0bf,1)}
\right\|_{X}^{-1}p_N(\phi)
\left\|f\right\|_{H_X(\rrn)}.
\end{align}
Using this, we find that the present
lemma holds true in this case.
Combining this and \eqref{bddol2e1},
we then complete the proof of
Lemma \ref{bddol2}.
\end{proof}

Furthermore,
the following auxiliary lemma
about function spaces
having absolutely continuous
quasi-norms is given in
\cite[Corollary 3.11(ii)]{SHYY},
which is important in the proof of
Theorem \ref{bddox}.

\begin{lemma}\label{bddoxl6}
Let $X$, $r$, $d$, $s$,
$\{\lambda_{j}\}_{j\in\mathbb{N}}$,
and $\{a_j\}_{j\in\mathbb{N}}$
be as in Lemma \ref{Atoll1}.
If $X$ has an absolutely continuous
quasi-norm, then $f\in H_X(\rrn)$ and
$f:=\sum_{j\in\mathbb{N}}
\lambda_ja_j$ converges in $H_X(\rrn)$.
\end{lemma}

To prove Theorem \ref{bddox},
we also need some auxiliary conclusions
about weighted function spaces.
In what follows, the
\emph{weighted Hardy space}\index{weighted Hardy space}
$H^p_{\upsilon}(\rrn)$\index{$H^p_{\upsilon}(\rrn)$},
with $p\in(0,\infty)$ and
$\upsilon\in A_{\infty}(\rrn)$,
is defined as in Definition
\ref{hardyx} with $X:=
L^p_{\upsilon}(\rrn)$.
Then, via Lemma \ref{bddoxl6} above,
we establish the following
atomic reconstruction theorem
about weighted Hardy spaces, which
is one of the key tools
used in the proof of Theorem \ref{bddox}.

\begin{lemma}\label{bddoxl2}
Let $0<\theta<s<s_0\leq1$, $\upsilon
\in A_1(\rrn)$, and $d\geq\lfloor n(1/\theta-1)
\rfloor$ be a fixed integer.
Assume that $\{a_j\}_{j\in\mathbb{N}}$
is a sequence of $(L^{s_0}_{\upsilon}(\rrn),
\,\infty,\,d)$-atoms supported, respectively,
in the balls $\{B_j\}_{j\in\mathbb{N}}
\subset\mathbb{B}$ and $\{\lambda_j\}
_{j\in\mathbb{N}}\subset[0,\infty)$
such that $f:=\sum_{j\in\mathbb{N}}
\lambda_ja_j$ in $\mathcal{S}'(\rrn)$,
and $$
\left\|
\left\{
\sum_{j\in\mathbb{N}}
\left[\frac{\lambda_j}{\|\1bf_{B_j}\|
_{L^{s_0}_{\upsilon}(\rrn)}}
\right]^s\1bf_{B_j}
\right\}^{\frac{1}{s}}
\right\|_{L^{s_0}_{\upsilon}(\rrn)}<\infty.
$$
Then $f\in H^{s_0}_{\upsilon}(\rrn)$
and $f=\sum_{j\in\mathbb{N}}\lambda_j
a_j$ holds true in $H^{s_0}_{\upsilon}(\rrn)$.
\end{lemma}

\begin{proof}
Let all the symbols be as
in the present lemma. Then,
applying \cite[Remarks 2.4(b),
2.7(b), and 3.4(i)]{WYY}, we find that
the following four statements hold true:
\begin{enumerate}
  \item[{\rm(i)}] $L^{s_0}_{
  \upsilon}(\rrn)$ is a BQBF space;
  \item[{\rm(ii)}] for any
$\{f_{j}\}_{j\in\mathbb{N}}
\subset L^1_{{\rm loc}}(\rrn)$,
\begin{equation*}
\left\|\left\{\sum_{j\in\mathbb{N}}
\left[\mc^{(\theta)}
(f_{j})\right]^s\right\}^{1/s}
\right\|_{L^{s_0}_{
\upsilon}(\rrn)}\lesssim
\left\|\left(\sum_{j\in\mathbb{N}}
|f_{j}|^{s}\right)^{1/s}\right\|_{
L^{s_0}_{\upsilon}(\rrn)};
\end{equation*}
  \item[{\rm(iii)}] $[L^{s_0}_{\upsilon}(\rrn)]
  ^{1/s}$ is a BBF space and,
for any $f\in L^1_{\mathrm{loc}}(\rrn)$,
\begin{equation*}
\left\|\mc(f)
\right\|_{([L^{s_0}_{\upsilon}(\rrn)]
^{1/s})'}\lesssim
\left\|f\right\|_{([L^{s_0}_{\upsilon}(\rrn)]
^{1/s})'};
\end{equation*}
  \item[{\rm(iv)}] $L^{s_0}_{\upsilon}(\rrn)$
  has an absolutely continuous quasi-norm.
\end{enumerate}
Thus, the weighted Lebesgue space $L^{s_0}_
{\upsilon}(\rrn)$ under consideration
satisfies all the assumptions
of Lemma \ref{bddoxl6} with $r$
therein replaced by $\infty$.
This finishes the proof of Lemma
\ref{bddoxl2}.
\end{proof}

On the other hand,
applying Proposition
\ref{bddoxl4}, we can obtain
the following boundedness of
Calder\'{o}n--Zygmund operators
on weighted Hardy spaces,
which is important
in the proof of the
boundedness of
Calder\'{o}n--Zygmund operators
on the Hardy space $H_X(\rrn)$.

\begin{proposition}\label{bddoxl5}
Let $s_0\in(0,1]$, $\upsilon\in A_1(\rrn)$,
$d\in\zp$, $K$ be a $d$-order
standard kernel defined
as in Definition \ref{def-s-k}
with some $\delta\in(0,1]$, and $T$ a $d$-order
Calder\'{o}n--Zygmund operator with kernel $K$
having the vanishing moments up to order $d$.
If $s_0\in(\frac{n}{n+d+\delta},1]$,
then $T$ has a unique extension on
$H^{s_0}_{\upsilon}(\rrn)$
and there exists a positive constant $C$
such that, for any $f\in
H^{s_0}_{\upsilon}(\rrn)$,
$$
\|T(f)\|_{H^{s_0}_{\upsilon}(\rrn)}
\leq C\|f\|_{
H^{s_0}_{\upsilon}(\rrn)}.
$$
\end{proposition}

\begin{proof}
Let all the symbols
be as in the present proposition,
$s\in(\frac{n}{n+d+\delta},s_0)$,
and $$\theta\in
\left(\frac{n}{n+d+\delta},\,
\min\left\{s,\frac{n}{n+d}\right\}\right).
$$
Then, by \cite[Remarks 2.4(b),
2.7(b), and 3.4(i)]{WYY}, we conclude that
the following four statements hold true:
\begin{enumerate}
  \item[{\rm(i)}] $L^{s_0}_{
  \upsilon}(\rrn)$ is a BQBF space;
  \item[{\rm(ii)}] for any
$\{f_{j}\}_{j\in\mathbb{N}}
\subset L^1_{{\rm loc}}(\rrn)$,
\begin{equation*}
\left\|\left\{\sum_{j\in\mathbb{N}}
\left[\mc^{(\theta)}
(f_{j})\right]^s\right\}^{1/s}
\right\|_{L^{s_0}_{
\upsilon}(\rrn)}\lesssim
\left\|\left(\sum_{j\in\mathbb{N}}
|f_{j}|^{s}\right)^{1/s}\right\|_{
L^{s_0}_{\upsilon}(\rrn)};
\end{equation*}
  \item[{\rm(iii)}] $[L^{s_0}_{\upsilon}(\rrn)]
  ^{1/s}$ is a BBF space and
there exists an $r_0\in(1,\infty)$ such that,
for any $f\in L^1_{\mathrm{loc}}(\rrn)$,
\begin{equation*}
\left\|\mc^{((r_0/s)')}(f)
\right\|_{([L^{s_0}_{\upsilon}(\rrn)]
^{1/s})'}\lesssim
\left\|f\right\|_{([L^{s_0}_{\upsilon}(\rrn)]
^{1/s})'};
\end{equation*}
  \item[{\rm(iv)}] $L^{s_0}_{\upsilon}(\rrn)$
  has an absolutely continuous quasi-norm.
\end{enumerate}
These imply that
the weighted Lebesgue space $L^{s_0}_
{\upsilon}(\rrn)$ under consideration
satisfies all the assumptions
of Proposition \ref{bddoxl4},
which completes the proof of Proposition
\ref{bddoxl5}.
\end{proof}

We now show Theorem \ref{bddox}.

\begin{proof}[Proof of Theorem
\ref{bddox}]
Let all the symbols be as in
the present theorem and $f\in
H_X(\rrn)$. Then, by the assumptions
$\theta\in(\frac{n}{n+d+\delta},
\frac{n}{n+d}]$ and $\delta\in(0,1]$,
we conclude that
$$
d>n\left(\frac{1}{\theta}-1\right)
-\delta\geq n\left(\frac{1}{\theta}-1
\right)-1
$$
and
$$
d\leq n\left(\frac{1}{\theta}-1\right).
$$
Thus, $d=\lfloor n(
1/\theta-1)\rfloor$.
From this,
the assumption (i) of the present
theorem, and Lemma \ref{Atogl4},
it follows that there exist
$\{\lambda_j\}_{j\in\mathbb{N}}
\subset[0,\infty)$ and $\{a_j\}_{j\in
\mathbb{N}}$ of $(X,\,\infty,\,d)$-atoms
supported, respectively, in
the balls $\{B_j\}_{j\in\mathbb{N}}\subset
\mathbb{B}$ such that
\begin{equation}\label{bddoxe1}
f=\sum_{j\in\mathbb{N}}\lambda_ja_j
\end{equation}
in $\mathcal{S}'(\rrn)$, and
\begin{equation}\label{bddoxe2}
\left\|\left[\sum_{j\in\mathbb
N}\left(\frac{\lambda_j}{\|\1bf_{B_j}\|_{X}}
\right)^s\1bf_{B_j}
\right]^{\frac{1}{s}}\right\|_{X}\lesssim
\|f\|_{H_X(\rrn)}.
\end{equation}

In addition, by the assumptions
(ii) through (v) of the present theorem
and Lemma \ref{bddol1} with
$Y:=Y_0$, $\theta:=\eta$, and $s:=s_0$,
we find that there exists an $\eps\in(0,1)$
such that, for any $h\in\Msc(\rrn)$,
\begin{equation*}
\|h\|_{L^{s_0}_{\upsilon}}(\rrn)
\lesssim\|h\|_X,
\end{equation*}
where $\upsilon:=[\mc(\1bf_{B(\0bf,1)})]
^{\eps}$. Combining this, \eqref{bddoxe1},
and \eqref{bddoxe2}, we conclude that
\begin{equation}\label{bddoxe5}
f=\sum_{j\in\mathbb{N}}\lambda_ja_j
=\sum_{j\in\mathbb{N}}\left[\lambda_j
\frac{\|\1bf_{B_j}\|_{L^{s_0}_{\upsilon}(\rrn)}}
{\|\1bf_{B_j}\|_X}\right]\left[\frac
{\|\1bf_{B_j}\|_X}{\|\1bf_{B_j}\|_{L^{s_0}
_{\upsilon}(\rrn)}}a_j\right]
\end{equation}
in $\mathcal{S}'(\rrn)$, and
\begin{align}\label{bddoxe6}
&\left\|\left\{
\sum_{j\in\mathbb{N}}\left[
\frac{\lambda_j\frac{\|\1bf_{B_j}
\|_{L^{s_0}_{\upsilon}(\rrn)}}
{\|\1bf_{B_j}\|_X}}{\|\1bf_{B_j}\|
_{L^{s_0}_{\upsilon}(\rrn)}}\right]^s
\1bf_{B_j}\right\}^{\frac{1}{s}}
\right\|_{L^{s_0}_{\upsilon}(\rrn)}\notag\\
&\quad\lesssim\left\|\left[\sum_{j\in\mathbb
N}\left(\frac{\lambda_j}{\|\1bf_{B_j}\|_{X}}
\right)^s\1bf_{B_j}
\right]^{\frac{1}{s}}\right\|_{X}\lesssim
\|f\|_{H_X(\rrn)}<\infty.
\end{align}
Observe that, for any
$j\in\mathbb{N}$, $\frac
{\|\1bf_{B_j}\|_X}{\|\1bf_{B_j}\|_{L^{s_0}
_{\upsilon}(\rrn)}}a_j$ is an $(L^{s_0}_{
\upsilon}(\rrn),\,\infty,\,d)$-atom
supported in $B_j$.
This, together with \eqref{bddoxe5},
\eqref{bddoxe6}, and Lemma \ref{bddoxl2},
further implies that $f\in H^{s_0}_{\upsilon}
(\rrn)$ and $f=\sum_{j\in\mathbb{N}}
\lambda_ja_j$ in $H^{s_0}_{\upsilon}(\rrn)$.
Applying this and
Proposition \ref{bddoxl5}, we find that
$T(f)=\sum_{j\in\mathbb{N}}
\lambda_jT(a_j)$ holds true in $H_{\upsilon}
^{s_0}(\rrn)$.
Therefore, from Lemma \ref{bddol2}
with $X$ replaced by $L^{s_0}_{
\upsilon}(\rrn)$, we further infer that
\begin{equation}\label{bddoxe4}
T(f)=\sum_{j\in\mathbb{N}}
\lambda_jT(a_j)
\end{equation}
in $\mathcal{S}'(\rrn)$.

Next, we prove that $T(f)\in H_X(\rrn)$.
To this end, choose an $r\in
(\max\{2,s\eta'\},\infty)$ with
$\frac{1}{\eta}+\frac{1}{\eta'}=1$.
Then, by Lemma \ref{Atoll2}
with $t=\infty$, we conclude that,
for any $j\in\mathbb{N}$, $a_j$
is an $(X,\,r,\,d)$-atom supported in
the ball $B_j$.
This, together with Lemma \ref{bddoxl3},
implies that,
for any $j\in\mathbb{N}$,
$T(a_j)$ is a harmless multiple
of an $(X,\,r,\,d,\,d+\delta+
\frac{n}{r'})$-molecule centered at
$B_j$.
Moreover, since $\theta>
\frac{n}{n+d+\delta}$, it follows
that
$$d+\delta+\frac{n}{r'}>
n\lf(\frac{1}{\theta}-\frac{1}{r}\r).$$
On the other hand,
using the assumption (ii) of
the present theorem and
an argument similar to that
used in the proof of \eqref{bddoxl4e1}
with $(r_0/s)'$ and $(X^{1/s})'$
therein replaced, respectively, by
$\eta$ and $Y$, we find that,
for any $f\in L^1_{\mathrm{loc}}(\rrn)$,
$$
\left\|\mc^{((r/s)')}(f)\right\|_{Y}
\lesssim\|f\|_Y.
$$
From this, both the assumptions (i)
and (iv) of the present theorem,
Theorem \ref{Molegx}, \eqref{bddoxe4},
and \eqref{bddoxe2},
we deduce that $T(f)\in H_X(\rrn)$
and
$$
\|T(f)\|_{H_X(\rrn)}
\lesssim\left\|\left[
\sum_{j\in\mathbb{N}}\left(
\frac{\lambda_j}{\|\1bf_{B_j}\|_X}\right)^s
\1bf_{B_j}\right]^{\frac{1}{s}}\right\|_{X}
\lesssim\|f\|_{H_X(\rrn)}.
$$
This finishes the proof that $T(f)\in H_X(\rrn)$
and further implies that $T$ is bounded on
$H_X(\rrn)$.
Thus, the proof of Theorem \ref{bddox}
is completed.
\end{proof}

Now, we prove Theorem \ref{bddog}.

\begin{proof}[Proof of Theorem
\ref{bddog}]
Let all the symbols be
as in the present theorem.
Then, using the assumptions
$\m0(\omega)\in(-\frac{n}{p},
\infty)$ and $\MI(\omega)\in
(-\infty,0)$, and Theorem \ref{Th2},
we find that,
under the assumptions of the present
theorem, the global generalized
Herz space $\HerzS$ is a BQBF space.
Therefore, to prove the present theorem,
it suffices to show that
all the assumptions of Theorem \ref{bddox}
hold true for $\HerzS$.

First, by the assumptions
$p,\ q\in(\frac{n}{n+d+\delta},\infty)$
and $$\Mw\in\left(-\frac{n}{p},n-\frac{n}{p}+d+
\delta\right),$$ we conclude that
$$
\frac{n}{n+d+\delta}<\min\left\{
1,p,q,\frac{n}{\Mw+n/p}\right\}.
$$
Let
$$
s\in\left(\frac{n}{n+d+\delta},
\min\left\{
1,p,q,\frac{n}
{\Mw+n/p}\right\}\right)
$$
and $$\theta\in\left(
\frac{n}{n+d+\delta},\,\min\left\{s,
\frac{n}{n+d}\right\}\right).$$
We now prove that $\HerzS$
satisfies Theorem \ref{bddox}(i) for the
above $\theta$ and $s$.
Indeed, from Lemma \ref{Atogl5},
it follows that, for any
$\{f_j\}_{j\in\mathbb{N}}\subset
L^1_{\mathrm{loc}}(\rrn)$,
\begin{equation*}
\left\|
\left\{
\sum_{j\in\mathbb{N}}
\left[\mc^{(\theta)}(f_j)\right]^{s}
\right\}^{1/s}
\right\|_{\HerzS}
\lesssim\left\|
\left(\sum_{j\in\mathbb{N}}
|f_j|^{s}
\right)^{1/s}
\right\|_{\HerzS}.
\end{equation*}
This implies that, for the above
$\theta$ and $s$,
Assumption \ref{assfs} holds
true for $\HerzS$ and hence
Theorem \ref{bddox}(i) is satisfied
for $\HerzS$.

We next show that there
exist two linear spaces $Y,\ Y_0$
and $\eta$, $s\in(0,\infty)$ such that
the assumptions (ii) through
(v) of Theorem \ref{bddox} hold true.
To this end, let
$Y:=\HerzScsd$, $s_0\in(0,1]$ be such that
$$
s_0\in\left(s,\min\left\{1,p,q,\frac{n}{\Mw+n/p}
\right\}\right),
$$
$Y_0:=\HerzScsod$, and $\eta\in(1,\infty)$
be such that
\begin{equation}\label{bddoge5}
\eta<\min\left\{
\frac{n}{n(1-s/p)-s\mw},\left(
\frac{p}{s}\right)'\right\}.
\end{equation}
Then, applying Lemma
\ref{equa}, we find that both
$$\HerzScsd\text{ and }\HerzScsod$$
satisfy Definition \ref{Df1}(ii).
This further implies that Theorem
\ref{bddox}(ii) holds true.

Now, we show that the assumption
(iii) of Theorem \ref{bddox} holds true
for the above $Y_0$. Namely, $\1bf_{B(\0bf,1)}
\in\HerzScsod$. Indeed, combining
Lemma \ref{rela} and the assumptions $
s<\frac{n}{\Mw+n/p}$ and $\M0(\omega)>-\frac{n}{p}$,
we conclude that
\begin{align*}
\m0\left(\omega^{-s_0}\right)
&=-s_0\M0(\omega)>-s_0\Mw\\&>
-n\left(1-\frac{s}{p}\right)
=-\frac{n}{(p/s)'}.
\end{align*}
Applying this and Theorem \ref{Th3}
with $p$, $q$, and $\omega$ therein
replaced, respectively, by
$(p/s_0)'$, $(q/s_0)'$, and $1/
\omega^{s_0}$,
we find that the Herz space $\Kmp
^{(p/s_0)',(q/s_0)'}_{
1/\omega^{s_0},\0bf}(\rrn)$ is a BQBF space.
This, together with Lemma
\ref{ee} with $p$, $q$, $\omega$, and $\xi$
therein replaced, respectively, by
$(p/s_0)'$, $(q/s_0)'$, $1/
\omega^{s_0}$, and $\0bf$, and Definition
\ref{Df1}(iv), further implies that
$$
\left\|\1bf_{B(\0bf,1)}
\right\|_{\HerzScsod}\lesssim
\left\|\1bf_{B(\0bf,1)}\right\|_{\Kmp
^{(p/s_0)',(q/s_0)'}_{
1/\omega^{s_0},\0bf}(\rrn)}<\infty,
$$
which completes the proof that
$\1bf_{B(\0bf,1)}\in\HerzScsod$.
Therefore, Theorem \ref{bddox}(iii)
holds true.

Next, we prove that both
$\HerzScsd$ and
$\HerzScsod$ satisfy Theorem \ref{bddox}(iv).
Indeed, using Lemma \ref{rela} and
the assumption $\m0(\omega)>-\frac{n}{p}$,
we find that
\begin{equation}\label{bddoge2}
\m0\left(\omega^s\right)
=s\m0(\omega)>-\frac{n}{p/s}.
\end{equation}
On the other hand, from Lemma
\ref{rela} again and the assumption
$\MI(\omega)<0$, we deduce that
\begin{equation*}
\MI\left(\omega^s\right)=
s\MI(\omega)<0.
\end{equation*}
Combining this, \eqref{bddoge2},
and Lemma \ref{Atogl2} with
$p$, $q$, and $\omega$ therein replaced,
respectively, by $p/s$, $q/s$,
and $\omega^s$, we conclude that,
for any $f\in\Msc(\rrn)$,
\begin{equation}\label{bddoge3}
\|f\|_{\HerzScs}\sim
\sup\left\{\|fg\|_{L^1(\rrn)}:\
\|g\|_{\HerzScsd}=1\right\},
\end{equation}
where the positive equivalence
constants are independent of $f$.
Similarly, repeating an argument
used in the proof of \eqref{bddoge3}
with $s$ replaced by $s_0$, we
find that, for any $f\in\Msc(\rrn)$,
\begin{equation*}
\|f\|_{\HerzScso}\sim
\sup\left\{\|fg\|_{L^1(\rrn)}:\
\|g\|_{\HerzScsod}=1\right\},
\end{equation*}
which, together with \eqref{bddoge3},
further implies that Theorem
\ref{bddog}(iv) holds true
with $Y=\HerzScsd$ and $
Y_0=\HerzScsod$.

Finally, we prove that the powered
Hardy--Littlewood maximal operator
$\mc^{(\eta)}$ is bounded on
both $\HerzScsd$ and $\HerzScsod$.
Indeed, from the assumption
$s<\frac{n}{\Mw+n/p}$,
we deduce that
\begin{align*}
&n\left(1-\frac{s}{p}\right)
-s\mw\\&\quad=n-s\left[\mw+
\frac{n}{p}\right]\\
&\quad\geq n-s\left[\Mw+
\frac{n}{p}\right]>0,
\end{align*}
which, combined with Lemma
\ref{rela} and \eqref{bddoge5},
further implies that
\begin{align}\label{bddoge6}
&\max\left\{
\M0\left(\omega^{-s}\right),
\MI\left(\omega^{-s}\right)\right\}
\notag\\&\quad=
-s\mw<n\left[\frac{1}{\eta}-
\frac{1}{(p/s)'}\right].
\end{align}
In addition, applying Lemma
\ref{rela} again and the assumptions
$$s<\frac{n}{\Mw+n/p}$$ and
$\Mw>-\frac{n}{p},$
we conclude that
\begin{align*}
&\min\left\{
\m0\left(\omega^{-s}\right),
\mi\left(\omega^{-s}\right)\right\}\\
&\quad=-s\Mw>-n\left(1-\frac{s}{p}\right)
=-\frac{n}{(p/s)'}.
\end{align*}
From this, \eqref{bddoge6}, the assumption
$\eta<(p/s)'$, and Corollary \ref{maxbl}
with $p$, $q$, $\omega$, and $r$
therein replaced, respectively,
by $(p/s)'$, $(q/s)'$, $1/\omega^s$,
and $\eta$,
it follows that, for any
$f\in L^1_{\mathrm{loc}}(\rrn)$,
\begin{equation}\label{bddoge7}
\left\|\mc^{(\eta)}(f)\right\|_{\HerzScsd}
\lesssim\left\|f\right\|_{\HerzScsd}.
\end{equation}
This further implies that $\mc^{(\eta)}$
is bounded on $Y=\HerzScsd$.
Moreover, repeating
an argument used in
the proof of \eqref{bddoge7} with
$s$ replaced by $s_0$, and
using the assumption $s_0>s$, we
find that, for any
$f\in L^1_{\mathrm{loc}}(\rrn)$,
\begin{equation*}
\left\|\mc^{(\eta)}(f)\right\|_{
\HerzScsod}\lesssim\left\|f\right\|_
{\HerzScsod},
\end{equation*}
which implies that $\mc^{(\eta)}$
is bounded on $Y_0=\HerzScsod$
and hence Theorem \ref{bddox}(v)
holds true.
Therefore, all the assumptions
of Theorem \ref{bddox} hold true
for $\HerzS$. Therefore, $T$ can be extended
into a bounded linear operator
on $\HaSaH$ and, for any $f\in\HaSaH$,
$$
\|T(f)\|_{\HaSaH}\lesssim
\|f\|_{\HaSaH}.
$$
This finishes the proof of Theorem
\ref{bddog}.
\end{proof}

Via Theorem \ref{bddog} and Remarks
\ref{remhs}(iv) and
\ref{remark4.10}(ii), we immediately obtain
the following boundedness of Calder\'{o}n--Zygmund
operators on the
generalized Morrey--Hardy space $H\MorrS$;
we omit the details.

\begin{corollary}
Let $d\in\zp$, $K$ be a
$d$-order standard kernel
defined as in Definition
\ref{def-s-k} with
some $\delta\in(0,1]$, $T$ a
$d$-order Calder\'{o}n--Zygmund operator
with kernel $K$ having the vanishing moments
up to order $d$, and
$p$, $q$, and $\omega$ as in Corollary
\ref{bddolm}.
Then $T$ can be extended
into a bounded linear operator
on $H\MorrS$ and there exists
a positive constant $C$ such that,
for any $f\in H\MorrS$,
$$
\|T(f)\|_{H\MorrS}\leq C\|f\|_{
H\MorrS}.
$$
\end{corollary}

\section{Fourier Transform}

The target of this section is
to investigate the Fourier transform of a distribution
in the generalized Herz--Hardy space $\HaSaHo$
or $\HaSaH$.
Recall that, in 1974,
Coifman \cite{Coi} characterized all $\widehat{f}$\,
via entire functions of exponential type for $n=1$, where
$f\in H^p(\rr)$ with $p\in(0,1]$.
Later, a number of authors investigated the
characterization of $\widehat{f}$ with distribution $f$ belonging to Hardy
spaces in higher dimensions (see,
for instance, \cite{BW,Col,HCY-an,HCY,TW}). In particular,
Huang et al. \cite{HCY} studied the Fourier transform of the
distribution belonging to the Hardy space associated
with the ball quasi-Banach function space, which plays an important role
in this section.

We first consider the Fourier transform properties
of the Hardy spaces associated with
the local generalized Herz spaces. Namely, we have
the following theorem.

\begin{theorem}\label{Fouriertl}
Let $p,\ q\in(0,1]$, $\omega\in M(\rp)$ with
$\m0(\omega)\in(0,\infty)$ and
$\mi(\omega)\in(0,\infty)$,
and $p_{-}\in(0,\frac{n}{\Mw+n/p})$.
Then, for any $f\in\HaSaHo$,
there exists a continuous function
$g$ on $\rrn$ such that
$\widehat{f}=g$ in $\mathcal{S}'(\rrn)$
and
$$
\lim\limits_{|x|\to0^{+}}\frac{|g(x)|}
{|x|^{n(\frac{1}{p_-}-1)}}=0.
$$
Moreover, there exists a positive
constant $C$, independent of both $f$ and $g$,
such that, for any $x\in\rrn$,
$$
|g(x)|\leq C\|f\|_{\HaSaHo}\max
\left\{1,|x|^{n(\frac{1}{p_-}-1)}\right\}
$$
and
$$
\int_{\rrn}|g(x)|\min\left\{
|x|^{-\frac{n}{p_{-}}},
|x|^{-n}\right\}\,dx\leq C\|f\|_{\HaSaHo}.
$$
\end{theorem}

To show this theorem, we
first investigate the properties
of the Fourier transform of the
distribution belonging to the
Hardy space $H_X(\rrn)$ associated with
the general ball quasi-Banach function space
$X$. Indeed, we have
the following technique lemma, which is
essential obtained by Huang et al.\
\cite[Theorems 2.1, 2.2, and 2.3]{HCY}.

\begin{lemma}\label{fouriertll2}
Let $X$ be a ball quasi-Banach
function space and $p_-\in(0,1]$ such that,
for any given $\theta\in(0,p_-)$ and
$u\in(1,\infty)$, there exists a positive
constant $C$ such that, for any
$\{f_j\}_{j\in\mathbb{N}}\subset
L^1_{\mathrm{loc}}(\rrn)$,
\begin{equation}\label{fouriertle14}
\left\|\left\{
\sum_{j\in\mathbb{N}}
\left[\mc(f_j)\right]^u
\right\}^{\frac{1}{u}}\right\|_{X^{1/\theta}}
\leq C\left\|\left(
\sum_{j\in\mathbb{N}}|f_j|^u
\right)^{\frac{1}{u}}\right\|_{X^{1/\theta}}.
\end{equation}
Assume that there exists a $p_0\in[p_-,1]$
such that $X$ is $p_0$-concave and
there exists a positive constant
$C$ such that, for any
$B\in\mathbb{B}$,
\begin{equation}\label{fouriertle13}
\left\|\1bf_B\right\|_X
\geq C\min\left\{|B|^{\frac{1}{p_0}},
|B|^{\frac{1}{p_-}}\right\}.
\end{equation}
Then, for any $f\in H_X(\rrn)$,
there exists a continuous function
$g$ on $\rrn$ such that
$\widehat{f}=g$ in $\mathcal{S}'(\rrn)$
and
$$
\lim\limits_{|x|\to0^{+}}
\frac{|g(x)|}{|x|^{n(\frac{1}{p_-}-1)}}=0.
$$
Moreover, there exists a positive constant
$C$, independent of both $f$ and $g$,
such that, for any $x\in\rrn$,
$$
|g(x)|\leq C\|f\|_{H_X(\rrn)}
\max\left\{|x|^{n(\frac{1}{p_0}-1)},
|x|^{n(\frac{1}{p_-}-1)}\right\}
$$
and
$$
\left[\int_{\rrn}|g(x)|^{p_0}\min
\left\{|x|^{n(p_0-1-\frac{p_0}{p_-})},
|x|^{n(p_0-2)}
\right\}\,dx\right]^{1/p_0}\leq
C\|f\|_{H_X(\rrn)}.
$$
\end{lemma}

\begin{proof}
Let all the symbols be as in
the present lemma, $d\geq\lfloor
n(1/p_{-}-1)\rfloor$ a fixed integer,
and $f\in\mathcal{S}'(\rrn)$. Then, by
\eqref{fouriertle14} and
the known atomic decomposition of
$H_X(\rrn)$ (see \cite[Theorem 3.7]{SHYY}
or Lemma \ref{Atogl4}), we find that
there exist $\{\lambda_{j}\}_{
j\in\mathbb{N}}\subset[0,\infty)$ and
a sequence $\{a_j\}_{j\in\mathbb{N}}$
of $(X,\infty,d)$-atoms supported, respectively,
in the balls $\{B_{j}\}_{j\in\mathbb{N}}
\subset\mathbb{B}$ such that
$$f=\sum_{j\in\mathbb{N}}\lambda_ja_j$$
in $\mathcal{S}'(\rrn)$, and
$$
\left\|\left[\sum_{j\in\mathbb
N}\left(\frac{\lambda_j}{\|\1bf_{B_j}\|_{X}}
\right)^s\1bf_{B_j}
\right]^{\frac{1}{s}}\right\|_{X}\lesssim
\|f\|_{H_X(\rrn)}.
$$
Using this and repeating the proof
of \cite[Theorems 2.1, 2.2, and 2.3]{HCY},
we conclude that there exists
a continuous function $g$ on $\rrn$ such that
$\widehat{f}=g$ in $\mathcal{S}'(\rrn)$,
$$
\lim\limits_{|x|\to0^{+}}
\frac{|g(x)|}{|x|^{n(\frac{1}{p_-}-1)}}=0,
$$
and, for any $x\in\rrn$,
$$
|g(x)|\lesssim\|f\|_{H_X(\rrn)}
\max\left\{|x|^{n(\frac{1}{p_0}-1)},
|x|^{n(\frac{1}{p_-}-1)}\right\}
$$
and
$$
\left[\int_{\rrn}|g(x)|^{p_0}\min
\left\{|x|^{n(p_0-1-\frac{p_0}{p_-})},
|x|^{n(p_0-2)}
\right\}\,dx\right]^{1/p_0}\lesssim
\|f\|_{H_X(\rrn)}
$$
with the implicit positive
constants independent of both $f$
and $g$.
This then finishes the
proof of Lemma \ref{fouriertll2}.
\end{proof}

In order to prove Theorem \ref{Fouriertl},
we also require the following
estimate about the quasi-norm of
the characteristic functions
of balls on the local
generalized Herz space $\HerzSo$,
which plays an important role
in the proof of Theorem \ref{Fouriertl}
and is also of independent
interest.

\begin{lemma}\label{Fouriertll1}
Let $p$, $q$, $\omega$, and $p_{-}$
be as in Theorem \ref{Fouriertl}.
Then there exists a positive constant
$C$ such that, for any $B\in\mathbb{B}$,
\begin{equation}\label{Fouriertle1}
\|\1bf_{B}\|_{\HerzSo}\geq
C\min\left\{|B|,|B|^{\frac{1}{p_{-}}}\right\}.
\end{equation}
\end{lemma}

\begin{proof}
Let all the symbols be as in the present
lemma. We first show that,
for any $B(x_{0},2^{k_{0}})\in\mathbb{B}$
with $x_{0}\in\rrn$
and $k_{0}\in\mathbb{Z}$,
\begin{equation}\label{Fouriertle2}
2^{\frac{nk_{0}}{p}}\omega(2^{k_{0}})\lesssim
\left\|\1bf_{B(x_{0},
2^{k_{0}})}\right\|_{\HerzSo}.
\end{equation}
To achieve this, we consider
the following four cases on $x_0$.

\emph{Case 1)} $x_{0}=\0bf$. In this case,
from the fact
that $B(\0bf,2^{k_{0}})\setminus
B(\0bf,2^{k_{0}-1})\subset
B(\0bf,2^{k_{0}})$ and the
definition of $\|\cdot\|_{\HerzSo}$,
we deduce that
\begin{align}\label{Fouriertle3}
2^{\frac{nk_{0}}{p}}\omega(2^{k_{0}})
&\sim\omega(2^{k_{0}})
\left\|\1bf_{B(\0bf,2^{k_{0}})
\setminus B(\0bf,2^{k_{0}-1})}
\right\|_{L^p(\rrn)}\notag
\\&
\lesssim\left\|\1bf_{B(\0bf,2^{k_{0}})}
\right\|_{\HerzSo},
\end{align}
which completes the proof of
\eqref{Fouriertle2} in this case.

\emph{Case 2)} $|x_{0}|\in(0,
3\cdot2^{k_{0}-1})$. In this case, we
first claim that
\begin{equation}\label{fouriertle11}
B\left(\frac{3\cdot2^{k_{0}-2}}
{|x_{0}|}x_0,2^{k_0-2}\right)
\subset\left[B(x_0,2^{k_0})
\cap\left(B(\0bf,2^{k_0})\setminus
B(\0bf,2^{k_0-1})\right)\right].
\end{equation}
Indeed, for any $y\in
B((3\cdot2^{k_{0}-2}/|x_{0}|)
x_0,2^{k_0-2})$, we have
\begin{align*}
|y|\leq\left|y-\frac{3\cdot2^
{k_{0}-2}}{|x_{0}|}x_0\right|+
\left|\frac{3\cdot2^{k_{0}-2}}
{|x_{0}|}x_0\right|<2^{k_0-2}+3\cdot2^{k_0-2}
=2^{k_0}
\end{align*}
and
\begin{align*}
|y|\geq\left|\frac{3\cdot2^{k_{0}-2}}
{|x_{0}|}x_0\right|-
\left|y-\frac{3\cdot2^{k_{0}-2}}
{|x_{0}|}x_0\right|>3\cdot2^{k_0-2}-2^{k_0-2}
=2^{k_0-1}.
\end{align*}
This implies that
\begin{equation}\label{Fouriertle31}
B\left(\frac{3\cdot2^{k_{0}-2}}
{|x_{0}|}x_0,2^{k_0-2}\right)
\subset\left[B(\0bf,2^{k_0})
\setminus B(\0bf,2^{k_0-1})\right].
\end{equation}
On the other hand, for any $y\in
B((3\cdot2^{k_{0}-2}/|x_{0}|)
x_0,2^{k_0-2})$, we have
\begin{align*}
|y-x_0|&\leq\left|y-\frac{3\cdot
2^{k_{0}-2}}{|x_{0}|}x_0\right|+
\left|\frac{3\cdot2^{k_{0}-2}}
{|x_{0}|}x_0-x_{0}\right|\\&<2^{k_0-2}+
\left|3\cdot2^{k_0-2}-|x_{0}|\right|.
\end{align*}
This further implies that, for any $y\in
B((3\cdot2^{k_{0}-2}/|x_{0}|)x_0,2^{k_0-2})$
with $x_{0}\in[3\cdot2^{k_0-2},3\cdot2^{k_0-1})$,
\begin{align*}
|y-x_0|<|x_0|-2^{k_0-1}<2^{k_0}
\end{align*}
and, for any $y\in B((3\cdot2^
{k_{0}-2}/|x_{0}|)x_0,2^{k_0-2})$
with $|x_{0}|\in(0,3\cdot2^{k_0-2})$,
\begin{align*}
|y-x_0|<2^{k_0}-|x_{0}|<2^{k_0}.
\end{align*}
Therefore,
$B((3\cdot2^{k_{0}-2}/|x_{0}|)x_0,
2^{k_0-2})\subset B(x_0,2^{k_0})$.
By this and \eqref{Fouriertle31}, we find that
$$
B\left(\frac{3\cdot2^{k_{0}-2}}
{|x_{0}|}x_0,2^{k_0-2}\right)
\subset\left[B(x_0,2^{k_0})\cap
\left\{B(\0bf,2^{k_0})\setminus
B(\0bf,2^{k_0-1})\right\}\right].
$$
This finishes the proof of the above claim.
Thus, from the definition
of the quasi-norm $\|\cdot
\|_{\HerzSo}$, it follows that
\begin{align}\label{Fouriertle6}
2^{\frac{nk_{0}}{p}}\omega(2^{k_0})
&\sim\omega(2^{k_0})\left|
B\left(\frac{3\cdot2^{k_{0}-2}}
{|x_{0}|}x_0,2^{k_0-2}\right)\right|
^{\frac{1}{p}}\notag\\
&\lesssim\omega(2^{k_0})\left|
B(x_0,2^{k_0})\cap\left[B(\0bf,2^{k_0})\setminus
B(\0bf,2^{k_0-1})\right]\right|^
{\frac{1}{p}}\notag\\&\lesssim\left\|
\1bf_{B(x_0,2^{k_0})}\right\|_{\HerzSo}
\end{align}
and hence \eqref{Fouriertle2} holds
true in this case.

\emph{Case 3)} $|x_0|\in[3\cdot2^
{k_0-1},2^{k_0+1})$.
In this case,
for any $y\in B(x_0,2^{k_0})$, we have
\begin{align*}
|y|\leq|y-x_0|+|x_0|<2^{k_0}+
2^{k_0+1}<2^{k_0+2}
\end{align*}
and
\begin{align*}
|y|\geq|x_0|-|y-x_0|>3\cdot2^
{k_0-1}-2^{k_0}=2^{k_0-1}.
\end{align*}
This implies that $B(x_{0},2^{k_0})
\subset B(\0bf,2^{k_0+2})\setminus
B(\0bf,2^{k_0-1})$. Applying this
and Lemma \ref{mono}, we conclude that
\begin{align}\label{Fouriertle7}
2^{\frac{nk_0q}{p}}\left[\omega
(2^{k_0})\right]^q&\sim\left[\omega(
2^{k_0})\right]^q\left|B(x_0,2^{k_0})
\right|^{\frac{q}{p}}\notag\\
&\lesssim\left[\omega(2^{k_0+2})
\right]^q\left|
B(x_0,2^{k_0})\cap\left[B(\0bf,
2^{k_0+2})\setminus
B(\0bf,2^{k_0+1})\right]\right|
^{\frac{q}{p}}\notag\\
&\quad+\left[\omega(2^
{k_0+1})\right]^q\left|
B(x_0,2^{k_0})\cap\left[B
(\0bf,2^{k_0+1})\setminus
B(\0bf,2^{k_0})\right]\right|
^{\frac{q}{p}}\notag\\
&\quad+\left[\omega(2^{k_0})
\right]^q\left|
B(x_0,2^{k_0})\cap\left
[B(\0bf,2^{k_0})\setminus
B(\0bf,2^{k_0-1})\right]\right|
^{\frac{q}{p}}\notag\\
&\lesssim\left\|\1bf_{B(x_0,2^
{k_0})}\right\|_{\HerzSo}^q,
\end{align}
which completes the estimation of
\eqref{Fouriertle2} in this case.

\emph{Case 4)} $|x_{0}|\in[2^k,
2^{k+1})$ with $k\in\mathbb{Z}\cap
[k_0+1,\infty)$. In this case,
from Lemma \ref{La3.5},
it follows that,
for any $0<t<\tau<\infty$,
\begin{equation*}
\frac{\omega(t)}{\omega(\tau)}\lesssim
\left(\frac{t}{\tau}\right)
^{\min\{\m0(\omega),\mi(\omega)\}-\eps},
\end{equation*}
where $\eps\in(0,\mw)$ is a fixed positive constant.
Using this, Lemma \ref{mono}, and the assumptions
$k\in[k_{0}+1,\infty)$ and $$\eps\in(0,\mw),$$
we further find that
\begin{align}\label{Fouriertle71}
\frac{\omega(2^{k_0})}{\omega(2^k)}
\sim\frac{\omega(2^{k_0+1})}{\omega(2^k)}
\lesssim\left(\frac{2^{k_0+1}}
{2^k}\right)^{\mw-\eps}\lesssim1.
\end{align}
On the other hand, applying the assumption
$|x_{0}|\in[2^k,2^{k+1})$ with
$k\in\mathbb{Z}\cap
[k_0+1,\infty)$, we conclude that,
for any $y\in B(x_0,2^{k_0})$,
\begin{align*}
|y|\leq|y-x_0|+|x_0|<2^{k_0}+2^{k+1}
\leq2^{k-1}+2^{k+1}<2^{k+2}
\end{align*}
and
\begin{align*}
|y|\geq|x_0|-|y-x_0|>2^k-2^{k_0}
\geq2^k-2^{k-1}=2^{k-1}.
\end{align*}
This implies that
\begin{equation}\label{fouriertle15}
B(x_0,2^{k_0})\subset
B(\0bf,2^{k+2})\setminus B(\0bf,2^{k-1}).
\end{equation}
Thus, from \eqref{Fouriertle71} and
an argument similar to that used in
the estimation of \eqref{Fouriertle7},
we deduce that
\begin{align*}
2^{\frac{nk_0}{p}}\omega
(2^{k_0})\lesssim\omega(\tk)
\left|B(x_0,2^{k_0})\right|^
{\frac{n}{p}}\lesssim\left\|
\1bf_{B(x_0,2^{k_0})}\right\|_{\HerzSo},
\end{align*}
which, together with \eqref{Fouriertle3},
\eqref{Fouriertle6},
and \eqref{Fouriertle7},
implies that \eqref{Fouriertle2} holds true.

Next, we show that, for any
$B(x_{0},r)\in\mathbb{B}$ with
$x_0\in\rrn$ and $r\in(0,\infty)$,
\begin{equation}\label{Fouriertle8}
r^{\frac{n}{p}}\omega(r)\lesssim
\left\|\1bf_{B(x_0,r)}\right\|_{\HerzSo}.
\end{equation}
Indeed, for any $r\in(0,\infty)$,
there exists a $k\in\mathbb{Z}$
such that $r\in[\tk,\tka)$.
Thus, applying Lemma \ref{mono}
and \eqref{Fouriertle2}, we find that
\begin{align*}
r^{\frac{n}{p}}\omega(r)\sim
2^{\frac{nk}{p}}\omega(\tk)
\lesssim\left\|\1bf_{B(x_0,\tk)}
\right\|_{\HerzSo}\lesssim
\left\|\1bf_{B(x_0,r)}\right\|_{\HerzSo}.
\end{align*}
This finishes the proof of \eqref{Fouriertle8}.

Finally, from the definition of
$p_{-}$ and Lemma \ref{Th1},
we infer that,
for any $r\in(0,1]$,
\begin{equation*}
\omega(r)\gtrsim r^{\M0(\omega)+\eps}
\end{equation*}
and, for any $r\in(1,\infty)$,
\begin{equation}\label{Fouriertle9}
\omega(r)\gtrsim r^{\mi(\omega)-\eps},
\end{equation}
where
$$
\eps\in\left(0,\min\left\{
\frac{n}{p_{-}}-\M0(\omega)-\frac{n}{p},
-n+\mi(\omega)+\frac{n}{p}\right\}\right)
$$
is a fixed positive constant.
Therefore, using the assumption
$\eps\in(0,\frac{n}{p_{-}}
-\M0(\omega)-\frac{n}{p})$
and \eqref{Fouriertle8}, we
conclude that, for any $x_0\in\rrn$
and $r\in(0,1]$,
\begin{align}\label{Fouriertle10}
\left|B(x_{0},r)\right|^{\frac{1}
{p_{-}}}\sim r^{\frac{n}{p_{-}}}
\lesssim r^{\M0(\omega)+\frac{n}{p}
+\eps}\lesssim r^{\frac{n}{p}}
\omega(r)\lesssim\left\|\1bf_
{B(x_{0},r)}\right\|_{\HerzSo}.
\end{align}
Moreover, applying the fact that $\eps\in
(0,-n+\mi(\omega)+\frac{n}{p})$,
\eqref{Fouriertle9}, and
\eqref{Fouriertle8}, we find that,
for any
$x_0\in\rrn$ and $r\in(1,\infty)$,
\begin{align}\label{fouriertle16}
\left|B(x_{0},r)\right|\sim r^{n}
\lesssim r^{\mi(\omega)+
\frac{n}{p}-\eps}
\lesssim r^{\frac{n}{p}}\omega(r)
\lesssim\left\|\1bf_{B(x_{0},r)}
\right\|_{\HerzSo},
\end{align}
which, together with \eqref{Fouriertle10},
further implies that
\eqref{Fouriertle1} holds true
and hence completes the proof
of Lemma \ref{Fouriertll1}.
\end{proof}

Via both Lemmas
\ref{fouriertll2} and \ref{Fouriertll1},
we now show Theorem \ref{Fouriertl}.

\begin{proof}[Proof of Theorem \ref{Fouriertl}]
Let all the symbols be as in the
present theorem.
Then, from the assumption
$\m0(\omega)\in(0,\infty)$ and
Theorem \ref{Th3}, it follows that
the local generalized
Herz space $\HerzSo$ is
a BQBF space. This implies that,
in order to complete
the proof of the present theorem,
we only need to show that the Herz space $\HerzSo$
under consideration satisfies all
the assumptions of Lemma
\ref{fouriertll2}. Namely,
$\HerzSo$ satisfies \eqref{fouriertle14},
the $p_0$-concavity for some $p_0\in
[p_-,1]$, and \eqref{fouriertle13}.

We first prove that \eqref{fouriertle14}
holds true for $\HerzSo$.
Indeed, by the fact that
$\mw\in(0,\infty)$ and Remark
\ref{remark 1.1.3}(iii),
we find that $$\Mw\in(0,\infty)$$
and hence
\begin{equation*}
\frac{n}{\Mw+n/p}\in(0,p).
\end{equation*}
Applying this and Lemma \ref{vmbhl}, we conclude that,
for any given $\theta\in(0,p_{-})$ and
$u\in(1,\infty)$,
and for any $\{f_{j}\}_{j\in\mathbb{N}}\subset
L^{1}_{\text{loc}}(\rrn)$,
\begin{equation}\label{fouriertl1}
\left\|\left\{\sum_{j\in\mathbb{N}}
\left[\mc
(f_{j})\right]^{u}\right\}^{\frac{1}{u}}
\right\|_{[\HerzSo]^{1/\theta}}
\lesssim\left\|\left(
\sum_{j\in\mathbb{N}}|f_{j}|^{u}
\right)^{\frac{1}{u}}\right\|
_{[\HerzSo]^{1/\theta}}.
\end{equation}
This further implies that
$\HerzSo$ satisfies \eqref{fouriertle14}.

In addition, from
$p,\ q\in(0,1]$, it follows that,
for any $\{f_{j}\}_{j\in\mathbb{N}}
\subset\HerzSo$,
\begin{equation}\label{Fouriertle12}
\sum_{j\in\mathbb{N}}\|f_{j}\|_{\HerzSo}\leq
\left\|\sum_{j\in\mathbb{N}}
|f_{j}|\right\|_{\HerzSo},
\end{equation}
which implies that
$\HerzSo$ is strictly 1-concave and
hence $\HerzSo$ satisfies
the $p_0$-concavity with $p_0=1$.

Finally, applying Lemma \ref{Fouriertll1},
we conclude that,
for any $B\in\mathbb{B}$,
$$
\left\|\1bf_{B}\right\|_{\HerzSo}
\gtrsim\min\left\{|B|,|
B|^{\frac{1}{p_-}}\right\}.
$$
Thus, \eqref{fouriertle13} also holds true
for $\HerzSo$ with $p_0=1$.
This, combined with both \eqref{fouriertl1}
and \eqref{Fouriertle12}, further implies that
all the assumptions of Lemma
\ref{fouriertll2} hold true for $\HerzSo$,
and hence finishes the proof of
Theorem \ref{Fouriertl}.
\end{proof}

Next, we investigate the Fourier transform
properties of the generalized
Herz--Hardy space $\HaSaH$. Indeed,
we have the following conclusion.

\begin{theorem}\label{Fouriert}
Let $p\in(0,1),\ q\in(0,1]$,
$\omega\in M(\rp)$ with
$\m0(\omega)\in(-\frac{n}{p},\infty)$ and
$$
n\left(1-\frac{1}{p}\right)<
\mi(\omega)\leq\MI(\omega)<0,
$$
and $p_{-}\in(0,\min\{p,\frac{n}{\Mw+n/p})$.
Then, for any $f\in\HaSaH$,
there exists a continuous function
$g$ on $\rrn$ such that
$\widehat{f}=g$ in $\mathcal{S}'(\rrn)$
and
$$
\lim\limits_{|x|\to0^{+}}\frac{|g(x)|}
{|x|^{n(\frac{1}{p_-}-1)}}=0.
$$
Moreover, there exists a positive
constant $C$, independent of both $f$ and $g$,
such that, for any $x\in\rrn$,
$$
|g(x)|\leq C\|f\|_{\HaSaH}\max
\left\{1,|x|^{n(\frac{1}{p_-}-1)}\right\}
$$
and
$$
\int_{\rrn}|g(x)|\min\left\{|x|^{-\frac{n}{p_{-}}},
|x|^{-n}\right\}\,dx\leq C\|f\|_{\HaSaH}.
$$
\end{theorem}

To show Theorem \ref{Fouriert},
we need the following technical
estimate for the quasi-norm of
the characteristic function of balls
on the global generalized Herz space
$\HerzS$.

\begin{lemma}\label{estaf}
Let $p$, $q$, $\omega$, and $p_{-}$
be as in Theorem \ref{Fouriert}.
Then there exists a positive
constant $C$ such that,
for any $B\in\mathbb{B}$,
\begin{equation}\label{F0}
\|\1bf_{B}\|_{\HerzS}\geq C\min
\left\{|B|,|B|^{\frac{1}{p_{-}}}\right\}.
\end{equation}
\end{lemma}

\begin{proof}
Let all the symbols be as in the present lemma and
$B:=B(x_{0},r)\in\mathbb{B}$
with $x_{0}\in\rrn$ and $r\in(0,\infty)$.
Then, there exists a $k\in\mathbb{Z}$ such that
$r\in[2^{k},2^{k+1})$. This implies that
$B(x_0,\tk)\setminus B(x_0,\tkm)\subset B(x_0,r)$.
Therefore, by Lemma
\ref{mono} and the definition of the
quasi-norm $\|\cdot\|_{\HerzS}$, we find that
\begin{align}\label{F1}
r^{\frac{n}{p}}\omega(r)
&\sim2^{\frac{nk}{p}}\omega(\tk)
\sim\omega(\tk)\left\|
\1bf_{B(x_{0},\tk)\setminus
B(x_0,\tkm)}\right\|_{L^p(\rrn)}\notag\\
&\lesssim\left\|\1bf_{B(x_0,r)}
\right\|_{\HerzS}.
\end{align}
Moreover, from the definition
of $p_{-}$ and Lemma \ref{Th1},
it follows that,
for any $r\in(0,1]$,
\begin{equation}\label{W1}
\omega(r)\gtrsim r^{\M0(\omega)+\eps}
\end{equation}
and, for any $r\in(1,\infty)$,
\begin{equation}\label{W2}
\omega(r)\gtrsim r^{\mi(\omega)-\eps},
\end{equation}
where
$$
\eps\in\left(0,\min\left\{\frac{n}
{p_{-}}-\M0(\omega)-\frac{n}{p},
-n+\mi(\omega)+\frac{n}{p}\right\}\right)
$$
is a fixed positive constant.
We now show \eqref{F0}
by considering the following two
cases on $r$.

\emph{Case 1)} $r\in(0,1]$. In this case,
using the assumption $\eps\in(0,
\frac{n}{p_{-}}-\M0(\omega)-\frac{n}{p})$,
\eqref{W1}, and \eqref{F1}, we conclude that
\begin{align}\label{w3}
\left|B(x_{0},r)\right|^{\frac{1}{p_{-}}}
\sim r^{\frac{n}{p_{-}}}
\lesssim r^{\M0(\omega)+\frac{n}
{p}+\eps}\lesssim r^{\frac{n}{p}}
\omega(r)\lesssim\left\|\1bf_
{B(x_{0},r)}\right\|_{\HerzS}.
\end{align}
This finishes the estimation of
\eqref{F0} in this case.

\emph{Case 2)} $r\in(1,\infty)$. In this case,
applying the fact that $\eps\in
(0,-n+\mi(\omega)+\frac{n}{p})$,
\eqref{W2}, and \eqref{F1}, we find that
\begin{align*}
\left|B(x_{0},r)\right|\sim r^{n}\lesssim
r^{\mi(\omega)+\frac{n}{p}-\eps}
\lesssim r^{\frac{n}{p}}\omega(r)\lesssim
\left\|\1bf_{B(x_{0},r)}\right\|_{\HerzS},
\end{align*}
which completes the proof of \eqref{F0}
in this case.
Combining this and \eqref{w3}, we conclude
that Lemma \ref{estaf} holds true.
\end{proof}

Based on the above lemma
and the Fourier transform properties
of Hardy spaces associated with
ball quasi-Banach function spaces
established in Lemma \ref{fouriertll2},
we next prove Theorem \ref{Fouriert}

\begin{proof}[Proof of Theorem
\ref{Fouriert}]
Let all the symbols be as in
the present theorem. Then,
by the assumptions $\m0
(\omega)\in(-\frac{n}{p},
\infty)$ and $\MI(\omega)\in
(-\infty,0)$, and Theorem
\ref{Th2}, we find that the global
generalized Herz space $\HerzS$
under consideration is a BQBF space.
Thus, to finish the proof of the present
theorem, it suffices to prove that
the Herz space $\HerzS$ satisfies
all the assumptions of Lemma
\ref{fouriertll2}.

First, we show that
\eqref{fouriertle14} holds true.
Indeed, from Lemma \ref{vmbhg} with $
r=\theta$, it follows that,
for any given $\theta
\in(0,p_-)$ and $u\in(1,\infty)$,
and for any $\{f
_{j}\}_{j\in\mathbb{N}}\subset
L_{{\rm loc}}^{1}(\rrn)$,
\begin{equation*}
\left\|\left\{\sum_{j\in\mathbb{N}}\left[\mc
(f_{j})\right]^{u}\right\}^{\frac{1}{u}}
\right\|_{[\HerzS]^{1/\theta}}
\leq C\left\|\left(
\sum_{j\in\mathbb{N}}|f_{j}|^{u}\right)^
{\frac{1}{u}}\right\|_{[\HerzS]^{1/\theta}},
\end{equation*}
which implies that
$\HerzS$ satisfies \eqref{fouriertle14}.

Next, we prove that
$\HerzS$ is strictly $1$-concave.
Indeed, applying the assumptions
$p$, $q\in(0,1]$, we conclude that,
for any $\{f_{j}\}_{j\in\mathbb{N}}
\subset\HerzS$,
\begin{equation*}
\sum_{j\in\mathbb{N}}\|f_{j}\|_{\HerzS}\leq
\left\|\sum_{j\in\mathbb{N}}|f_{j}|\right\|_{\HerzS}.
\end{equation*}
This implies that $\HerzS$ is
strictly $p_0$-concave
with $p_0:=1$.

Finally, using Lemma \ref{estaf},
we find that, for any $B\in\mathbb{B}$,
$$
\|\1bf_{B}\|_{\HerzS}\gtrsim\min
\left\{|B|,|B|^{\frac{1}{p_{-}}}\right\},
$$
which implies that
\eqref{fouriertle13} holds true
for $\HerzS$ with $p_0:=1$. Thus,
the Herz space $\HerzS$
under consideration satisfies \eqref{fouriertle14},
$1$-concavity, and \eqref{fouriertle13}
with $p_0:=1$. By this and Lemma
\ref{fouriertll2}, we then
obtain Theorem \ref{Fouriert}.
\end{proof}

\chapter{Localized Generalized Herz--Hardy Spaces\label{sec7}}
\markboth{\scriptsize\rm\sc Localized Generalized Herz--Hardy Spaces}
{\scriptsize\rm\sc Localized Generalized Herz--Hardy Spaces}

Recall that, in 1979,
Goldberg \cite{Gold} introduced the local Hardy space $h^p(\rrn)$,
with $p\in(0,\infty)$, as a localization of the
classical Hardy spaces $H^p(\rrn)$
and showed some properties
which are different from the classical Hardy spaces.
Indeed, Goldberg proved that
the local Hardy space $h^p(\rrn)$ contains
$\mathcal{S}(\rrn)$ as a dense
subspace and established the boundedness of
pseudo-differential operators of order zero on $h^p(\rrn)$.
Moreover, due to \cite{Gold}, the localized Hardy space
is well defined on manifolds. After that, lots of
nice works have been done in the study of
localized Hardy spaces and their applications;
see, for instance, \cite{cky18,d05,m94,yyz21}.
Specially, some variants of the local
Hardy spaces $h^p(\rrn)$
have also been studied; see, for instance,
\cite{t12} for weighted
localized Hardy spaces,
\cite{wl15} for localized Herz--Hardy spaces,
\cite{yy11} for (weighted) localized
Orlicz--Hardy spaces, and \cite{yy12} for localized
Musielak--Orlicz--Hardy spaces.
Very recently, Sawano et al. \cite{SHYY} introduced
the local Hardy space $h_X(\rrn)$ associated with
the ball quasi-Banach function
space $X$ and gave various maximal
function characterizations of $h_X(\rrn)$.
Wang et al. \cite{WYY} then established
atomic, molecular, and various Littlewood--Paley
function characterizations of $h_{X}(\rrn)$
as well as showed the boundedness
of pseudo-differential operators
on these localized Hardy spaces.

The main target of this chapter
is devoted to introducing the
localized generalized Herz--Hardy spaces
and then establishing their complete
real-variable theory. To achieve
this, we begin with showing their various
maximal function characterizations
via the known maximal function characterizations
of the local Hardy space $h_X(\rrn)$
associated with the ball quasi-Banach function
space $X$. Then, to overcome the obstacle
caused by the deficiency of associate spaces
of global generalized Herz spaces, we
establish improved
atomic and molecular characterizations of $h_X(\rrn)$
as well as the boundedness of
pseudo-differential operators on $h_X(\rrn)$
without recourse to associate spaces (see Theorems
\ref{atomxxl}, \ref{molexl},
and \ref{pseudox} below). Via these improved conclusions,
we establish the atomic and molecular characterizations
of localized generalized Herz--Hardy spaces
and the boundedness of pseudo-differential
operators on localized generalized Herz--Hardy spaces.
In addition, to clarify the relation
between localized generalized Herz--Hardy spaces
and generalized Herz--Hardy spaces, we find the
relation between $h_X(\rrn)$ and
Hardy spaces $H_X(\rrn)$ associated with
the ball quasi-Banach function space $X$
(see Theorem \ref{relax} below).
This extends the results obtained
by Goldberg \cite[Lemma 4]{Gold} for
classical Hardy spaces and also Nakai and Sawano
\cite[Lemma 9.1]{NSa12} for variable
Hardy spaces. Applying this and some auxiliary lemmas about
generalized Herz spaces, we then
obtain the relation between
localized generalized Herz--Hardy spaces
and generalized Herz--Hardy spaces.
As applications, we also establish various
Littlewood--Paley function
characterizations of $h_X(\rrn)$ and hence, together
with the construction of the
quasi-norm $\|\cdot\|_{\HerzS}$,
we further obtain the Littlewood--Paley function
characterizations of localized generalized Herz--Hardy spaces.
In addition, we introduce localized generalized
Morrey--Hardy spaces. Using the
equivalence between generalized
Herz spaces and generalized Morrey spaces, we
also conclude the
corresponding real-variable
characterizations and applications
of localized generalized Morrey--Hardy spaces
in this chapter.

For any given $N\in\mathbb{N}$ and
for any $f\in\mathcal{S}'(\rrn)$,
the \emph{local
grand maximal
function}\index{local grand maximal function}
$m_{N}(f)$ is defined by setting,
for any $x\in\rrn$,\index{$m_N$}
\begin{equation}\label{sec7e1}
m_{N}(f)(x):=\sup\left\{|f\ast\phi_{t}(y)|:\ t\in(0,1),
\ |x-y|<t,\ \phi\in\mathcal{F}_{N}(\rrn)\right\},
\end{equation}
where $\mathcal{F}_N(\rrn)$ is as in \eqref{defn}.
Then we introduce localized generalized Herz--Hardy
spaces\index{localized generalized Herz--Hardy space} as follows.

\begin{definition}\label{lghh}
Let $p$, $q\in(0,\infty)$, $\omega\in M(\rp)$,
and $N\in\mathbb{N}$.
Then
\begin{enumerate}
  \item[{\rm(i)}] the \emph{local generalized Herz--Hardy space}
  $\HaSaHol$\index{$\HaSaHol$},
  associated with the local generalized
  Herz space $\HerzSo$,
  is defined to be the set of all the
  $f\in\mathcal{S}'(\rrn)$ such that $$
\|f\|_{\HaSaHol}:=\|m_{N}(f)\|_{\HerzSo}<\infty;$$
  \item[{\rm(ii)}] the \emph{local generalized Herz--Hardy
  space} $\HaSaHl$\index{$\HaSaHl$},
  associated with the global generalized Herz space $\HerzS$,
  is defined to be the set of all the
  $f\in\mathcal{S}'(\rrn)$ such that $$
\|f\|_{\HaSaHl}:=\|m_{N}(f)\|_{\HerzS}<\infty.$$
\end{enumerate}
\end{definition}

\begin{remark}
In Definition \ref{lghh}, for any given $\alpha\in\rr$
and for any $t\in(0,\infty)$, let
$\omega(t):=t^\alpha$. Then, in this case,
the local generalized Herz--Hardy
space $\HaSaHol$ goes back to the classical
\emph{homogeneous local
Herz-type Hardy space}\index{homogeneous local
Herz-type Hardy \\space}
$h\dot{K}_p^{\alpha,q}(\rrn)$\index{$h\dot{K}_p^{\alpha,q}(\rrn)$}
which was originally introduced in \cite[Definition 1.2]{fy97}
(see also \cite[Section 2.6]{LYH}).
However, to the best of our knowledge, even
in this case, the local generalized
Herz--Hardy space $\HaSaHl$ is also new.
\end{remark}

Recall that local and global generalized Morrey spaces
are given in Remark \ref{remhs}(iv). We now introduce the
following concepts of localized generalized Morrey--Hardy
spaces\index{localized generalized Morrey--Hardy \\space}.

\begin{definition}\label{lgm}
Let $p$, $q\in(0,\infty)$, $\omega\in M(\rp)$, and $N\in\mathbb{N}$.
\begin{enumerate}
  \item[{\rm(i)}] The \emph{local
  generalized Morrey--Hardy space}
  $h\MorrSo$\index{$h\MorrSo$},
  associated with the local generalized Morrey space $\MorrSo$,
  is defined to be the set of all
  the $f\in\mathcal{S}'(\rrn)$ such that $$
\|f\|_{h\MorrSo}:=\|m_{N}(f)\|_{\MorrSo}<\infty;$$
  \item[{\rm(ii)}] The \emph{local
  generalized Morrey--Hardy
  space} $h\MorrS$\index{$h\MorrS$},
  associated with the global generalized Morrey space $\MorrS$,
  is defined to be the set of all the
  $f\in\mathcal{S}'(\rrn)$ such that $$
\|f\|_{h\MorrS}:=\|m_{N}(f)\|_{\MorrS}<\infty.$$
\end{enumerate}
\end{definition}

\begin{remark}\label{remark5.0.9}
\begin{enumerate}
  \item[{\rm(i)}] We should point out that,
  to the best of our knowledge,
  in Definition \ref{lgm}, even when
  $\omega(t):=t^\alpha$ for any $t\in(0,\infty)$ and
  for any given $\alpha\in\rr$,
  the local generalized Morrey--Hardy
  spaces $h\MorrSo$ and $h\MorrS$
  are also new.
  \item[{\rm(ii)}] In Definition \ref{lgm},
  let $p$, $q\in[1,\infty)$
  and $\omega$ satisfy $$\Mw\in(-\infty,0).$$
  Then, by Remark \ref{remhs}(iv),
  we find that, in this case, the
  local generalized Morrey--Hardy spaces
  $h\MorrSo$ and $h\MorrS$ coincide, respectively, with
  the local generalized Herz--Hardy
  spaces $\HaSaHol$ and $\HaSaHl$
  in the sense of equivalent norms.
\end{enumerate}
\end{remark}

\section{Maximal Function Characterizations}

In this section, we establish the
maximal function characterizations
of localized generalized Herz--Hardy spaces.
To begin with,
we recall the following definitions of
various localized maximal functions
(see, for instance, \cite[Definition 5.1]{SHYY}).

\begin{definition}\label{df511}
Let $N\in\mathbb{N},\
a,\ b\in(0,\infty)$,
$\phi\in\mathcal{S}(\rrn)$, and
$f\in\mathcal{S}'(\rrn)$.
\begin{enumerate}
\item[(i)] The \emph{local radial maximal
function}\index{local
radial maximal function}
$m(f,\phi)$\index{$m(f,\phi)$} is
defined by setting, for any $x\in\rrn$,
$$m(f,\phi)(x):=\sup_{t\in(0,1)}|f\ast\phi_{t}(x)|.$$
\item[(ii)] The \emph{local non-tangential
maximal function}\index{local
non-tangential maximal \\function}
$m^{*}_{a}(f,\phi)$\index{$m^{*}_{a}$}, with aperture
$a\in(0,\infty)$, is defined by setting,
for any $x\in\rrn$,
$$
m_{a}^{*}(f,\phi)(x):=\sup_{t\in(0,1)}\sup_{
\{y\in\rrn:\ |y-x|<at\}}|f\ast\phi_{t}(y)|.
$$
\item[(iii)] The \emph{local maximal
function $m_b^{**}(f,\phi)$ of
Peetre type}\index{local maximal
function of Peetre type}\index{$m_b^{**}$}
is defined by setting, for any $x\in\rrn$,
$$
m_b^{**}(f,\phi)(x):=\sup_{(y,t)\in\rrn
\times(0,1)}\frac{|f\ast\phi_{t}(x-y)|}
{(1+t^{-1}|y|)^{b}}.
$$
\item[(iv)] The \emph{local grand maximal
function $m^{**}_{b,N}(f)$
of Peetre type}\index{local grand maximal function
of Peetre type}\index{$m^{**}_{b,N}$} is defined
by setting, for any $x\in\rrn$,
$$m^{**}_{b,N}(f)(x):=\sup_{\phi\in
\mathcal{F}_{N}(\rrn)}
\sup_{(y,t)\in\rrn\times(0,1)}\frac{
|f\ast\phi_{t}(x-y)|}
{(1+t^{-1}|y|)^{b}},$$
where, for any $N\in\mathbb{N}$,
$\mathcal{F}_N(\rrn)$ is defined as in
\eqref{defn}.
\end{enumerate}
\end{definition}

Now, we establish the following maximal function characterizations
of the local generalized Herz--Hardy space
$\HaSaHol$.\index{maximal function characterization}

\begin{theorem}\label{Th6.1}
Let $p,\ q,\ a,\ b\in(0,\infty),\ \omega\in M(\rp),\
N\in\mathbb{N}$,
and $\phi\in\mathcal{S}(\rrn)$
satisfy $\int_{\rrn}
\phi(x)\,dx\neq0.$
\begin{enumerate}
\item[\rm{(i)}] Let $N\in\mathbb{N}\cap[\lfloor b+1\rfloor,\infty)$
and $\omega$ satisfy $\m0(\omega)\in(-\frac{n}{p},\infty)$.
Then, for any $f\in\mathcal{S}'(\rrn)$,
$$
\|m(f,\phi)\|_{\HerzSo}\lesssim\|m^*_{a}(f,\phi)\|_{\HerzSo}
\lesssim\|m^{**}_{b}(f,\phi)\|_{\HerzSo},$$
\begin{align*}
\|m(f,\phi)\|_{\HerzSo}&\lesssim\|m_{N}(f)\|_{\HerzSo}
\lesssim\|m_{\lfloor b+1\rfloor}(f)\|_{\HerzSo}\\
&\lesssim\|m^{**}_{b}(f,\phi)\|_{\HerzSo},
\end{align*}
and
$$
\|m^{**}_{b}(f,\phi)\|_{\HerzSo}
\sim\|m^{**}_{b,N}(f)\|_{\HerzSo},
$$
where the implicit positive
constants are independent of $f$.
\item[\rm{(ii)}] Let $\omega$ satisfy
$\m0(\omega)\in(-\frac{n}{p},\infty)$
and $\mi(\omega)\in(-\frac{n}{p},\infty)$,
and $$b\in\left(2\max\left\{\frac{n}{p},
\frac{n}{q},\Mw+\frac{n}{p}\right\},\infty\right).$$
Then, for any $f\in\mathcal{S}'(\rrn),$
$$\|m^{**}_{b}(f,\phi)\|_{\HerzSo}
\lesssim\|m(f,\phi)\|_{\HerzSo},$$
where the implicit positive constant
is independent of $f$. In particular,
when $N\in\mathbb{N}\cap[\lfloor b+1\rfloor,\infty)$,
if one of the quantities
$$\|m(f,\phi)\|_{\HerzSo},\ \|m^*_{a}(f,\phi)\|
_{\HerzSo},\ \|m_{N}
(f)\|_{\HerzSo},$$
$$\|m^{**}_{b}(f,\phi)\|_{\HerzSo},
\ \text{and }\|m^{**}_{b,N}(f)\|_{\HerzSo}$$
is finite, then the others are also finite
and mutually equivalent with the
positive equivalence constants independent of $f$.
\end{enumerate}
\end{theorem}

To prove the above maximal function
characterizations, we require the following
auxiliary inequality which was obtained in
\cite[Remark 4.4]{WYY}.

\begin{lemma}\label{Th6.1l2}
Let $X$ be a ball quasi-Banach function
space and $r\in(0,\infty)$ such that
$\mc^{(r)}$ is bounded on $X$.
Then there exists a positive constant
$C$ such that, for any
$f\in X$ and $z\in\rrn$,
$$
\left\|\left\{
\int_{z+[0,1]^{n}}
|f(\cdot-y)|^{r}\,dy\right\}^{\frac{1}{r}}
\right\|_{X}\leq C
(1+|z|)^{\frac{n}{r}}\|f\|_{X}.
$$
\end{lemma}

We also need the maximal function
characterizations of localized Hardy spaces
associated with ball quasi-Banach function
spaces as follows, which is just
\cite[Theorem 5.3]{SHYY} and
plays a key role in the proof of
Theorem \ref{Th6.1}.

\begin{lemma}\label{Th6.1l1}
Let $X$ be a ball quasi-Banach function
space, $a,\ b\in(0,\infty)$,
$N\in\mathbb{N}$,
and $\phi\in\mathcal{S}(\rrn)$
satisfy $\int_{\rrn}\phi(x)\,dx\neq0.$
\begin{enumerate}
\item[\rm{(i)}] Let $N\in\mathbb{N}
\cap[\lfloor b+1\rfloor,\infty)$.
Then, for any $f\in\mathcal{S}'(\rrn)$,
$$
\|m(f,\phi)\|_{X}\lesssim\|m^*_{a}(f,\phi)\|_{X}
\lesssim\|m^{**}_{b}(f,\phi)\|_{X},$$
\begin{align*}
\|m(f,\phi)\|_{X}\lesssim\|m_{N}(f)\|_{X}
\lesssim\|m_{\lfloor b+1\rfloor}(f)\|_{X}
\lesssim\|m^{**}_{b}(f,\phi)\|_{X},
\end{align*}
and
$$
\|m^{**}_{b}(f,\phi)\|_{X}
\sim\|m^{**}_{b,N}(f)\|_{X},
$$
where the implicit positive
constants are independent of $f$.
\item[\rm{(ii)}] Let $r$, $A\in
(0,\infty)$ satisfy $(b-A)r>n$.
Assume that $X$ is
strictly $r$-convex and
there exists a positive constant $C$
such that, for any $g\in X$
and $z\in\rrn$,
\begin{equation}\label{Th6.1e3}
\left\|\left\{
\int_{z+[0,1]^{n}}
|g(\cdot-y)|^{r}\,dy\right\}^{\frac{1}{r}}
\right\|_{X}\leq C
(1+|z|)^{A}\|g\|_{X}.
\end{equation}
Then, for any $f\in\mathcal{S}'(\rrn),$
$$\|m^{**}_{b}(f,\phi)\|_{X}
\lesssim\|m(f,\phi)\|_{X},$$
where the implicit positive constant
is independent of $f$. In particular,
when $N\in\mathbb{N}\cap[\lfloor
b+1\rfloor,\infty)$,
if one of the quantities
$$\|m(f,\phi)\|_{X},\ \|m^*_{a}(f,\phi)\|
_{X},\ \|m_{N}
(f)\|_{X},$$
$$\|m^{**}_{b}(f,\phi)\|_{X},
\ \text{and }\|m^{**}_{b,N}(f)\|_{X}$$
is finite, then the others are also finite
and mutually equivalent with the
positive equivalence constants independent of $f$.
\end{enumerate}
\end{lemma}

\begin{remark}\label{Th6.1r}
Let $X$ be a ball quasi-Banach
function space satisfying
Assumption \ref{assfs}
for some $0<\theta<s\leq1$ and
Assumption \ref{assas} for the same $s$.
Then we claim that $X$ satisfies all
the assumptions of Lemma \ref{Th6.1l1}(ii).
Indeed, by Assumption \ref{assas}
and \cite[Lemma 2.6]{CWYZ} with
$X$ and $p$ therein replaced, respectively,
by $X^{1/s}$ and $\frac{s}{\theta}$, we find
that $X^{1/\theta}$ is a ball Banach function space
and hence $X$ is strictly $\theta$-convex.
On the other hand, from both Assumption \ref{assfs}
and Lemma \ref{Th6.1l2} with $r:=\theta$,
it follows that \eqref{Th6.1e3}
holds true with $r:=\theta$
and $A\in(\frac{n}{\theta},\infty)$.
This then implies that all the assumptions
of Lemma \ref{Th6.1l1}(ii) hold true for
the above $X$ and hence finishes the proof of
the above claim.
\end{remark}

Via the above two lemmas, we now show
Theorem \ref{Th6.1}.

\begin{proof}[Proof of Theorem \ref{Th6.1}]
Let all the symbols be as in the present theorem.
Then, from the assumption
$\m0(\omega)\in(-\frac{n}{p},\infty)$ and
Theorem \ref{Th3}, it follows that
the local generalized
Herz space $\HerzSo$ under consideration
is a BQBF space. This,
together with Lemma \ref{Th6.1l1}(i),
then finishes the proof of (i).

Next, we show (ii).
Indeed, by Lemma \ref{Th6.1l1}(ii),
we conclude that it suffices
to prove that all the assumptions
of Lemma \ref{Th6.1l1}(ii) hold true for
the Herz space $\HerzSo$ under consideration.
Namely, there exist $r$, $A\in
(0,\infty)$ with $(b-A)r>n$ satisfying
that $\HerzSo$ is
strictly $r$-convex and,
for any $f\in\HerzSo$
and $z\in\rrn$,
\begin{equation}\label{Th6.1e2}
\left\|\left\{
\int_{z+[0,1]^{n}}
|f(\cdot-y)|^{r}\,dy\right\}^{\frac{1}{r}}
\right\|_{\HerzSo}\lesssim
(1+|z|)^{A}\|f\|_{\HerzSo}.
\end{equation}
To achieve this, let
$$
r\in\left(
\frac{2n}{b},\min\left\{p,q,
\frac{n}{\Mw+n/p}\right\}\right).
$$
Then, from Theorem \ref{strc} with
$s:=r$, we deduce that
the local generalized Herz space
$\HerzSo$
is strictly $r$-convex.
On the other hand, let
$A\in(\frac{n}{r},b-\frac{n}{r})$.
We then prove that \eqref{Th6.1e2}
holds true for this $A$.
Indeed, by Lemma \ref{mbhl}
and Remark \ref{power}(i), we find that,
for any $f\in L^1_{\mathrm{loc}}(\rrn)$,
\begin{equation*}
\left\|\mc^{(r)}(f)\right\|_{\HerzSo}
\lesssim\|f\|_{\HerzSo}.
\end{equation*}
Applying this, Lemma \ref{Th6.1l2},
and the assumption $A>\frac{n}{r}$,
we further conclude that, for any
$f\in\HerzSo$ and $z\in\rrn$,
\begin{align*}
&\left\|\left\{
\int_{z+[0,1]^{n}}
|f(\cdot-y)|^{r}\,dy\right\}^
{\frac{1}{r}}\right\|_{\HerzSo}\\
&\quad\lesssim
(1+|z|)^{\frac{n}{r}}\|f\|_{\HerzSo}\lesssim
(1+|z|)^{A}\|f\|_{\HerzSo}.
\end{align*}
This finishes the proof
of \eqref{Th6.1e2}.
Thus, under the assumptions
of the present theorem, $\HerzSo$
satisfies all the assumptions
of Lemma \ref{Th6.1l1}(ii),
which completes the proof of (ii),
and hence of Theorem \ref{Th6.1}.
\end{proof}

\begin{remark}
We point out that, in Theorem \ref{Th6.1},
if $\omega(t):=t^{\alpha}$ for any $t\in(0,\infty)$
and for any given $\alpha\in[n(1-\frac{1}{p}),
\infty)$, then Theorem \ref{Th6.1} goes back to
\cite[Theorem 2.6.1]{LYH}.
\end{remark}

Via Theorem \ref{Th6.1} and Remark \ref{remark5.0.9}(ii),
we immediately obtain the maximal function characterizations
\index{maximal function characterization}of
the local generalized Morrey--Hardy space $h\MorrSo$ as follows;
we omit the details.

\begin{corollary}
Let $a$, $b\in(0,\infty)$, $p$, $q\in[1,\infty)$,
$\omega\in M(\rp)$, $N\in\mathbb{N}$,
and $\phi\in\mathcal{S}(\rrn)$
satisfy $\int_{\rrn}
\phi(x)\,dx\neq0.$
\begin{enumerate}
\item[\rm{(i)}] Let $N\in\mathbb{N}
\cap[\lfloor b+1\rfloor,\infty)$
and $\omega$ satisfy $\MI(\omega)\in(-\infty,0)$
and
$$
-\frac{n}{p}<\m0(\omega)\leq\M0(\omega)<0.
$$
Then, for any $f\in\mathcal{S}'(\rrn)$,
$$
\|m(f,\phi)\|_{\MorrSo}\lesssim
\|m^*_{a}(f,\phi)\|_{\MorrSo}
\lesssim\|m^{**}_{b}(f,\phi)
\|_{\MorrSo},$$
\begin{align*}
\|m(f,\phi)\|_{\MorrSo}&\lesssim
\|m_{N}(f)\|_{\MorrSo}
\lesssim\|m_{\lfloor b+1\rfloor}
(f)\|_{\MorrSo}\\
&\lesssim\|m^{**}_{b}(f,\phi)\|_{\MorrSo},
\end{align*}
and
$$
\|m^{**}_{b}(f,\phi)\|_{\MorrSo}
\sim\|m^{**}_{b,N}(f)\|_{\MorrSo},
$$
where the implicit positive
constants are independent of $f$.
\item[\rm{(ii)}] Let $\omega$ satisfy
$$
-\frac{n}{p}<\m0(\omega)\leq\M0(\omega)<0
$$
and
$$
-\frac{n}{p}<\mi(\omega)\leq\MI(\omega)<0,
$$
and $b\in(2n\max\{\frac{1}{p},\frac{1}{q}\},\infty)$.
Then, for any $f\in\mathcal{S}'(\rrn),$
$$\|m^{**}_{b}(f,\phi)\|_{\MorrSo}
\lesssim\|m(f,\phi)\|_{\MorrSo},$$
where the implicit positive constant
is independent of $f$. In particular,
when $N\in\mathbb{N}\cap[\lfloor b+1
\rfloor,\infty)$, if one of the quantities
$$\|m(f,\phi)\|_{\MorrSo},\ \|m^*_{a}
(f,\phi)\|_{\MorrSo},\ \|m_{N}
(f)\|_{\MorrSo},$$
$$\|m^{**}_{b}(f,\phi)\|_{\MorrSo},
\ \text{and }\|m^{**}_{b,N}(f)\|_{\MorrSo}$$
is finite, then the others are also finite
and mutually equivalent with the
positive equivalence constants independent of $f$.
\end{enumerate}
\end{corollary}

Next, we establish the following maximal function
characterizations\index{maximal function characterization} of
the local generalized Herz--Hardy space $\HaSaHl$
with the help of the known maximal function
characterizations of the local Hardy
space $h_X(\rrn)$ presented
in Lemma \ref{Th6.1l1}.

\begin{theorem}\label{Th6.3}
Let $p,\ q,\ a,\ b\in(0,\infty),\ \omega\in M(\rp),\
N\in\mathbb{N}$,
and $\phi\in\mathcal{S}(\rrn)$
satisfy $\int_{\rrn}
\phi(x)\,dx\neq0.$
\begin{enumerate}
\item[\rm{(i)}] Let $N\in\mathbb{N}
\cap[\lfloor b+1\rfloor,\infty)$
and $\omega$ satisfy $\m0(\omega)
\in(-\frac{n}{p},\infty)$ and
$\MI(\omega)\in(-\infty,0)$.
Then, for any $f\in\mathcal{S}'(\rrn)$,
$$
\|m(f,\phi)\|_{\HerzS}\lesssim
\|m^*_{a}(f,\phi)\|_{\HerzS}
\lesssim\|m^{**}_{b}(f,\phi)\|_{\HerzS},$$
\begin{align*}
\|m(f,\phi)\|_{\HerzS}&\lesssim
\|m_{N}(f)\|_{\HerzS}
\lesssim\|m_{\lfloor b+1\rfloor}
(f)\|_{\HerzS}\\
&\lesssim\|m^{**}_{b}(f,\phi)\|_{\HerzS},
\end{align*}
and
$$
\|m^{**}_{b}(f,\phi)\|_{\HerzS}\sim
\|m^{**}_{b,N}(f)\|_{\HerzS},
$$
where the implicit positive
constants are independent of $f$.
\item[\rm{(ii)}] Let $\omega$ satisfy
$\m0(\omega)\in(-\frac{n}{p},\infty)$
and
$$
-\frac{n}{p}<\mi(\omega)\leq\MI(\omega)<0,
$$
and $$b\in\left(2\max\left\{\frac{n}{p},
\frac{n}{q},\Mw+\frac{n}{p}\right\},
\infty\right).$$
Then, for any $f\in\mathcal{S}'(\rrn),$
$$\|m^{**}_{b}(f,\phi)\|_{\HerzS}
\lesssim\|m(f,\phi)\|_{\HerzS},$$
where the implicit positive constant
is independent of $f$. In particular,
when $N\in\mathbb{N}\cap[\lfloor
b+1\rfloor,\infty)$,
if one of the quantities
$$\|m(f,\phi)\|_{\HerzS},\
\|m^*_{a}(f,\phi)\|_{\HerzS},\ \|m_{N}
(f)\|_{\HerzS},$$
$$\|m^{**}_{b}(f,\phi)\|_{\HerzS},
\ \text{and }\|m^{**}_{b,N}(f)\|_{\HerzS}$$
is finite, then the others are also finite
and mutually equivalent with the
positive equivalence constants
independent of $f$.
\end{enumerate}
\end{theorem}

\begin{proof}
Let all the symbols be as in the present theorem.
Note that $\omega$ satisfies
$\m0(\omega)\in(-\frac{n}{p},\infty)$
and $\MI(\omega)\in(-\infty,0)$. From this and
Theorem \ref{Th2}, it follows that
the global generalized
Herz space $\HerzS$ under consideration
is a BQBF space. Thus,
applying Lemma \ref{Th6.1l1}(i),
we then complete
the proof of (i).

Now, we prove (ii).
Indeed, we only need to show
that the Herz space $\HerzS$
under consideration satisfies all
the assumptions of Lemma \ref{Th6.1l1}(ii).
For this purpose, let
$r,\ A\in(0,\infty)$ be such that
$$
r\in\left(
\frac{2n}{b},\min\left
\{p,q,\frac{n}{\Mw+n/p}\right\}\right)
$$
and $A\in(\frac{n}{r},b-\frac{n}{r})$.
Then, applying Theorem \ref{Th3.1} with
$s:=r$, we conclude that
the global generalized Herz space
$\HerzS$ is strictly $r$-convex.
On the other hand, from Lemma \ref{mbhg}
and Remark \ref{power}(i), we
deduce that, for any
$f\in L^1_{\mathrm{loc}}(\rrn)$,
\begin{equation*}
\left\|\mc^{(r)}(f)\right\|_
{\HerzS}\lesssim\|f\|_{\HerzS}.
\end{equation*}
By this, Lemma \ref{Th6.1l2},
and the assumption $A>\frac{n}{r}$,
we further find that,
for any $f\in\HerzS$ and $z\in\rrn$,
\begin{align*}
&\left\|\left\{
\int_{z+[0,1]^{n}}
|f(\cdot-y)|^{r}\,dy\right\}^
{\frac{1}{r}}\right\|_{\HerzS}\\
&\quad\lesssim
(1+|z|)^{\frac{n}{r}}\|f\|_{\HerzS}\lesssim
(1+|z|)^{A}\|f\|_{\HerzS}.
\end{align*}
Combining this and
the fact that $\HerzS$ is
strictly $r$-convex, we conclude that
all the assumptions of Lemma \ref{Th6.1l2}(ii)
hold true for $\HerzS$.
This finishes the proof of (ii), and hence
of Theorem \ref{Th6.3}.
\end{proof}

Using Theorem \ref{Th6.3} and Remark \ref{remark5.0.9}(ii),
we immediately obtain the following maximal function
characterizations\index{maximal function characterization}
of the local generalized Morrey--Hardy space $h\MorrS$;
we omit the details.

\begin{corollary}
Let $a$, $b\in(0,\infty)$, $p$, $q
\in[1,\infty)$,
$\omega\in M(\rp)$, $N\in\mathbb{N}$,
and $\phi\in\mathcal{S}(\rrn)$
satisfy $\int_{\rrn}
\phi(x)\,dx\neq0.$
\begin{enumerate}
\item[\rm{(i)}] Let $N\in\mathbb{N}
\cap[\lfloor b+1\rfloor,\infty)$
and $\omega$ satisfy $\MI(\omega)\in(-\infty,0)$
and
$$
-\frac{n}{p}<\m0(\omega)\leq\M0(\omega)<0.
$$
Then, for any $f\in\mathcal{S}'(\rrn)$,
$$
\|m(f,\phi)\|_{\MorrS}\lesssim\|
m^*_{a}(f,\phi)\|_{\MorrS}
\lesssim\|m^{**}_{b}(f,\phi)\|_{\MorrS},$$
\begin{align*}
\|m(f,\phi)\|_{\MorrS}&\lesssim\|m_{N}(f)\|_{\MorrS}
\lesssim\|m_{\lfloor b+1\rfloor}
(f)\|_{\MorrS}\\
&\lesssim\|m^{**}_{b}(f,\phi)\|_{\MorrS},
\end{align*}
and
$$
\|m^{**}_{b}(f,\phi)\|_{\MorrS}\sim
\|m^{**}_{b,N}(f)\|_{\MorrS},
$$
where the implicit positive
constants are independent of $f$.
\item[\rm{(ii)}] Let $\omega$ satisfy
$$
-\frac{n}{p}<\m0(\omega)\leq\M0(\omega)<0
$$
and
$$
-\frac{n}{p}<\mi(\omega)\leq\MI(\omega)<0,
$$
and $b\in(2n\max\{\frac{1}{p},\frac{1}{q}\},\infty)$.
Then, for any $f\in\mathcal{S}'(\rrn),$
$$\|m^{**}_{b}(f,\phi)\|_{\MorrS}\lesssim
\|m(f,\phi)\|_{\MorrS},$$
where the implicit positive constant
is independent of $f$. In particular,
when $N\in\mathbb{N}\cap[\lfloor b+1
\rfloor,\infty)$, if one of the quantities
$$\|m(f,\phi)\|_{\MorrS},\ \|m^*_{a}
(f,\phi)\|_{\MorrS},\ \|m_{N}
(f)\|_{\MorrS},$$
$$\|m^{**}_{b}(f,\phi)\|_{\MorrS},
\ \text{and }\|m^{**}_{b,N}(f)\|_{\MorrS}$$
is finite, then the others are also finite
and mutually equivalent with the
positive equivalence constants
independent of $f$.
\end{enumerate}
\end{corollary}

\section{Relations with Generalized Herz--Hardy \\ Spaces}

The target of this section is to investigate the
relation between
localized generalized Herz--Hardy spaces and
generalized Herz--Hardy spaces.
For this purpose, we first
establish a relation between
localized Hardy spaces and
Hardy spaces associated with
ball quasi-Banach function spaces,
which extends the results obtained
by Goldberg \cite[Lemma 4]{Gold} for
classical Hardy spaces and Nakai and Sawano
\cite[Lemma 9.1]{NSa12} for variable
Hardy spaces. Using this general conclusion,
we then obtain the relation between
localized generalized Herz--Hardy spaces
and generalized Herz--Hardy spaces.

First, we give the relation between
localized Hardy spaces and Hardy spaces
associated with ball quasi-Banach function spaces.
To this end, we first recall the definition
of localized Hardy spaces associate with
ball quasi-Banach function space
as follows, which was introduced in
\cite[Definition 5.2]{SHYY}.

\begin{definition}\label{hardyxl}
Let $X$ be a ball quasi-Banach
function space and $N\in
\mathbb{N}$. Then the \emph{local Hardy
space} $h_X(\rrn)$\index{$h_X(\rrn)$} is defined to be
the set of all the $f\in\mathcal{S}'(\rrn)$
such that
$$
\left\|f\right\|_{h_X(\rrn)}
:=\left\|m_N(f)\right\|_{X}<\infty.
$$
\end{definition}

Then we have the following conclusion.

\begin{theorem}\label{relax}
Let $X$ be a ball quasi-Banach
function space and let $r\in(0,\infty)$
be such that $X$ is strictly $r$-convex
and $\mc$ is bounded on $X^{1/r}$.
Assume $\varphi\in\mathcal{S}(\rrn)$
satisfying that
$$
\1bf_{B(\0bf,1)}\leq
\widehat{\varphi}\leq\1bf_{B(\0bf,2)}.
$$
Then there exists a
constant $C\in[1,\infty)$ such that,
for any $f\in\mathcal{S}'(\rrn)$,
$$
C^{-1}\|f\|_{h_X(\rrn)}
\leq\left\|f\ast\varphi\right\|_{X}
+\left\|f-f\ast\varphi\right\|_{H_X(\rrn)}
\leq C\|f\|_{h_X(\rrn)}.
$$
\end{theorem}

To prove this theorem,
we need the following Plancherel--P\'olya--Nikol'skij
inequality\index{Plancherel--P\'olya--Nikol'skij \\
inequality} which is a consequence
of \cite[p.\,16, Theorem]{t83}
(see also \cite[Lemma 2.6]{NSa12}).

\begin{lemma}\label{relaxl1}
Let $\varphi\in\mathcal{S}(\rrn)$ be such that
$\widehat{\varphi}$ has compact support,
and let $r\in(0,\infty)$.
Then there exists a positive
constant $C$, independent
of $\varphi$, such that, for any
$f\in\mathcal{S}'(\rrn)$ and $x$, $y\in\rrn$,
$$
\left|f\ast\varphi(y)\right|
\leq C\left(1+|x-y|\right)^{
\frac{n}{r}}\mc^{(r)}\left(
f\ast\varphi\right)(x).
$$
\end{lemma}

Via Lemma \ref{relaxl1}, we next show
Theorem \ref{relax}.

\begin{proof}[Proof of Theorem \ref{relax}]
Let all the symbols be as in the
present theorem and $f\in\mathcal{S}'(\rrn)$.
We first prove that
\begin{equation}\label{relaxe0}
\left\|f\right\|_{h_X(\rrn)}
\lesssim\left\|f\ast\varphi\right\|
_{X}+\left\|f-f\ast\varphi\right\|_{H_X(\rrn)}.
\end{equation}
Indeed, from Definitions
\ref{hardyx} and \ref{hardyxl},
\eqref{sec6e1}, and
\eqref{sec7e1}, we deduce that
\begin{equation}\label{relaxe1}
\left\|f-f\ast\varphi\right\|_{h_X(\rrn)}\leq
\left\|f-f\ast\varphi\right\|_{H_X(\rrn)}.
\end{equation}
On the other hand, by
the assumption that $\mc$ is
bounded on $X^{1/r}$ and Remark
\ref{power}(i), we conclude
that $\mc^{(r)}$ is bounded
on $X$. This, together with
Lemma \ref{Th6.1l2},
further implies that,
for any
$g\in X$ and $z\in\rrn$,
$$
\left\|\left\{
\int_{z+[0,1]^{n}}
|g(\cdot-y)|^{r}\,dy\right\}^{\frac{1}{r}}
\right\|_{X}\lesssim
(1+|z|)^{\frac{n}{r}}\|g\|_{X}.
$$
Applying this, the assumption
that $X$ is strictly $r$-convex, and
Lemma \ref{Th6.1l1}(ii),
we find that
\begin{equation*}
\left\|f\ast\varphi\right\|_{h_X(\rrn)}
\sim\left\|m(f\ast\varphi,\varphi)
\right\|_{X},
\end{equation*}
which, combined with \eqref{relaxe1},
implies that
\begin{align}\label{relaxe2}
\left\|f\right\|_{h_X(\rrn)}
&\lesssim\left\|f\ast\varphi\right\|_{h_X(\rrn)}
+\left\|f-f\ast\varphi\right\|_{h_X(\rrn)}\notag\\
&\lesssim\left\|m(f\ast\varphi,\varphi)
\right\|_{X}+\left\|f-f\ast\varphi\right\|_{H_X(\rrn)}.
\end{align}
Now, we estimate $m(f\ast\varphi,\varphi)$.
Indeed, using Lemma
\ref{relaxl1} with $y$ replaced by $x-y$, we
conclude that, for any $t\in(0,1)$ and $x\in\rrn$,
\begin{align*}
&\left|\left(f\ast\varphi
\right)\ast\varphi_t(x)\right|\notag\\
&\quad=\left|\int_{\rrn}f\ast\varphi(x-y)
\varphi_t(y)\,dy\right|\notag\\
&\quad\lesssim
\int_{\rrn}\left(1+|y|\right)^{\frac{n}{r}}
\mc^{(r)}\left(
f\ast\varphi\right)(x)\left|\varphi_t(y)
\right|\,dy\notag\\
&\quad\lesssim\frac{1}{t^n}
\mc^{(r)}\left(
f\ast\varphi\right)(x)\int_{\rrn}
\left(1+\left|\frac{y}{
t}\right|\right)^{\frac{n}{r}}
\left|\varphi\left(\frac{y}{t}
\right)\right|\,dy\notag\\
&\quad\sim\mc^{(r)}\left(
f\ast\varphi\right)(x)\int_{\rrn}
\left(1+|y|\right)^{\frac{n}{r}}
\left|\varphi(y)\right|\,dy\notag\\
&\quad\lesssim\mc^{(r)}\left(
f\ast\varphi\right)(x)\int_{\rrn}
\frac{1}{(1+|y|)^{n+1}}\,dy
\sim\mc^{(r)}\left(
f\ast\varphi\right)(x).
\end{align*}
Therefore, from \eqref{sec7e1},
it follows that
\begin{equation*}
m\left(f\ast\varphi,\varphi
\right)\lesssim\mc^{(r)}\left(
f\ast\varphi\right).
\end{equation*}
Then, combining this, Definition \ref{Df1}(ii),
the assumption that $\mc$ is bounded
on $X^{1/r}$, and Remark \ref{power}(i),
we conclude that
\begin{align*}
\left\|m\left(f\ast\varphi,
\varphi\right)\right\|_{X}
\lesssim\left\|\mc^{(r)}\left(
f\ast\varphi\right)\right\|_{X}
\lesssim\left\|f\ast\varphi\right\|_X.
\end{align*}
By this and \eqref{relaxe2},
we further obtain
$$
\left\|f\right\|_{h_X(\rrn)}
\lesssim\left\|f\ast\varphi\right\|
_{X}+\left\|f-f\ast\varphi\right\|_{H_X(\rrn)},
$$
which completes the proof of
\eqref{relaxe0}.

Conversely, we show
\begin{equation}\label{relaxe4}
\left\|f\ast\varphi\right\|
_{X}+\left\|f-f\ast\varphi\right\|_{H_X(\rrn)}
\lesssim\left\|f\right\|_{h_X(\rrn)}.
\end{equation}
To achieve this, we
first estimate the
term $\|f-f\ast
\varphi\|_{H_X(\rrn)}$.
Indeed, applying Definition \ref{smax}(i),
we find that
\begin{align}\label{relaxe5}
M\left(f-f\ast\varphi,\varphi\right)
&=\sup_{t\in(0,\infty)}
\left|\left(f-f\ast\varphi\right)
\ast\varphi_t\right|\notag\\
&=\sup_{t\in(0,\infty)}
\left|f\ast\varphi_t-f\ast\left(
\varphi\ast\varphi_t\right)\right|.
\end{align}
We next claim that, for any
$t\in[2,\infty)$, $f\ast\varphi_
t-f\ast(\varphi\ast\varphi_t)=0$.
Indeed, notice that, for any
$t\in(0,\infty)$,
$$
f\ast\varphi_t-f\ast\left(
\varphi\ast\varphi_t\right)\in\mathcal{S}'(\rrn).
$$
Thus, for any $t\in(0,\infty)$, we have
\begin{equation}\label{relaxe6}
\mathcal{F}\left(f\ast\varphi_t-f\ast\left(
\varphi\ast\varphi_t\right)\right)
=\widehat{f}\widehat{\varphi}(t\cdot)
\left(1-\widehat{\varphi}\right).
\end{equation}
Using the assumption
$\widehat{\varphi}\1bf_{B(\0bf,1)}=1$,
we conclude that, for any $t\in(0,\infty)$
and $x\in B(\0bf,1)$,
$\widehat{\varphi}(tx)[1-\widehat{
\varphi}(x)]=0$.
On the other hand, for any
$t\in[2,\infty)$ and $x\in [B(\0bf,1)]^{
\complement}$, we have
$|tx|\geq t\geq2$. This,
together with the assumption
that $\widehat{\varphi}
\1bf_{[B(\0bf,2)]^{\complement}}=0$,
further implies that, for any
$t\in[2,\infty)$ and $x\in[B(\0bf,1)]^{
\complement}$,
$$\widehat{\varphi}(tx)[1-\widehat{
\varphi}(x)]=0.$$
Therefore,
for any $t\in[2,\infty)$, we have
$$\widehat{\varphi}(t\cdot)\left(1-\widehat{
\varphi}\right)=0.$$
By this and \eqref{relaxe6},
we find that, for any
$t\in[2,\infty)$,
$$
\mathcal{F}\left(f\ast\varphi_t-f\ast\left(
\varphi\ast\varphi_t\right)\right)=0
$$
in $\mathcal{S}'(\rrn)$, and hence
$f\ast\varphi_t-f\ast(
\varphi\ast\varphi_t)=0$, which completes the
proof of the above claim.
Now, we define $\psi\in\mathcal{S}(\rrn)$
by setting, for any $x\in\rrn$,
\begin{equation}\label{relaxe12}
\psi(x):=\frac{1}{2^n}\varphi\left(\frac{x}
{2}\right).
\end{equation}
Then, from \eqref{relaxe5} and
the above claim, we deduce that
\begin{align*}
M\left(f-f\ast\varphi,\varphi\right)
&=\sup_{t\in(0,2)}
\left|f\ast\varphi_t-f\ast\left(
\varphi\ast\varphi_t\right)\right|\\
&\leq\sup_{t\in(0,1)}
\left|f\ast\psi_t\right|
+\sup_{t\in(0,2)}\left|f\ast\left(
\varphi\ast\varphi_t\right)\right|\\
&=m(f,\psi)+\sup_{t\in(0,2)}\left|f\ast\left(
\varphi\ast\varphi_t\right)\right|.
\end{align*}
Combining this and both Lemmas \ref{Th5.4l1}(ii)
and \ref{Th6.1l1}(ii), we find that
\begin{align}\label{relaxe8}
\left\|f-f\ast\varphi\right\|_{H_X(\rrn)}
&\sim\left\|M(f-f\ast\varphi,\varphi)\right\|_{X}
\notag\\&\lesssim
\left\|m(f,\psi)\right\|_X
+\left\|\sup_{t\in(0,2)}
\left|f\ast\left(\varphi\ast\varphi_t
\right)\right|\right\|_X\notag\\
&\sim\left\|f\right\|_{h_X(\rrn)}
+\left\|\sup_{t\in(0,2)}
\left|f\ast\left(\varphi\ast\varphi_t
\right)\right|\right\|_X.
\end{align}

Next, we prove that
\begin{equation}\label{relaxe9}
\left\|\sup_{t\in(0,2)}
\left|f\ast\left(\varphi\ast\varphi_t
\right)\right|\right\|_X\lesssim
\left\|f\right\|_{h_X(\rrn)}.
\end{equation}
Indeed, for any
$t\in(0,\infty)$, we have
$$
\varphi\ast\varphi_t=
\mathcal{F}\left(\widehat{\varphi}(-\cdot)
\widehat{\varphi}(-t\cdot)\right).
$$
Thus, for any $\alpha$,
$\beta\in\zp^n$, $t\in(0,2)$,
and $x\in\rrn$, we obtain
\begin{align}\label{relaxe10}
&\left|x^{\beta}\partial^{\alpha}
\left(\varphi\ast\varphi_t\right)(x)\right|
\notag\\&\quad\sim
\left|\mathcal{F}\left(
\partial^{\beta}\left(
\left(\cdot\right)^
{\alpha}\widehat{\varphi}
(-\cdot)\widehat{\varphi}
(-t\cdot)\right)\right)(x)\right|\notag\\
&\quad\lesssim\int_{\rrn}\left|
\partial^{\beta}\left(
\left(\cdot\right)^
{\alpha}\widehat{\varphi}
(-\cdot)\widehat{\varphi}
(-t\cdot)\right)(x)\right|
\,dx\notag\\
&\quad\lesssim\sum_{\gfz{\gamma\in\zp^n}{\gamma\leq\beta}}\int_{\rrn}
\left[\left|x^{\alpha}\widehat{\varphi}
(-x)\partial^{\gamma}\left(
\widehat{\varphi}(-t\cdot)\right)(x)\right|\r.\notag\\
&\qquad\lf.+\left|\widehat{\varphi}(-tx)
\partial^{\beta-\gamma}
\left(\left(\cdot\right)^{\alpha}
\widehat{\varphi}(-\cdot)\right)(x)\right|\right]\,dx
\notag\\&\quad\lesssim\sum_{\gfz{\gamma\in\zp^n}
{\gamma\leq\beta}}\int_{\rrn}\lf[t^{|\gamma|}
\left|x^{\alpha}\widehat{\varphi}(-x)\right|
+\left|\partial^{\beta-\gamma}
\left(\left(\cdot\right)^{\alpha}
\widehat{\varphi}(-\cdot)\right)(x)\right|\r]
\,dx\notag\\
&\quad\lesssim\int_{\rrn}
\frac{1}{(1+|x|)^{n+1}}\,dx\sim1.
\end{align}
Fix an $N\in\mathbb{N}\cap
(\frac{n}{r}+1,\infty)$.
Then, applying \eqref{relaxe10}, we conclude that,
for any $\alpha\in\zp^n$ with $|\alpha|\leq N$,
any $t\in(0,2)$, and $x\in\rrn$,
\begin{align*}
&\left(1+|x|\right)^N\left|\partial^{\alpha}
\left(\frac{1}{2^n}\left(\varphi\ast\varphi_t\right)
\left(\frac{\cdot}{2}\right)\right)(x)\right|\notag\\
&\quad\lesssim\sum_{\beta\in\zp^n,\,|\beta|\leq N}
\left|x^{\beta}\partial^{\alpha}
\left(\frac{1}{2^n}\left(\varphi\ast\varphi_t\right)
\left(\frac{\cdot}{2}\right)\right)(x)\right|\notag\\
&\quad\lesssim\sum_{\beta\in\zp^n,\,|\beta|\leq
N}\left|x^{\beta}\partial^{\alpha}\left(
\varphi\ast\varphi_t\right)\left(\frac{x}
{2}\right)\right|\lesssim1.
\end{align*}
This, together
with \eqref{defp},
further implies that, for any $t\in(0,2)$,
$$
p_N\left(\frac{1}{2^n}\left(\varphi\ast\varphi_t\right)
\left(\frac{\cdot}{2}\right)\right)\lesssim1,
$$
where the implicit positive
constant is independent of $t$.
From this, \eqref{defn}, and \eqref{sec7e1},
it follows that
$$
\sup_{t\in(0,2)}
\left|f\ast\left(\varphi\ast\varphi_t
\right)\right|\lesssim m_N(f).
$$
Combining this and Definition \ref{Df1}(ii),
we further conclude that
$$
\left\|\sup_{t\in(0,2)}
\left|f\ast\left(\varphi\ast\varphi_t
\right)\right|\right\|_X\lesssim
\left\|m_N(f)\right\|_{X}\sim
\left\|f\right\|_{h_X(\rrn)}.
$$
This finishes the estimation of
\eqref{relaxe9}. From this and
\eqref{relaxe8}, we further infer that
\begin{equation}\label{relaxe11}
\left\|f-f\ast\varphi\right\|_{H_X(\rrn)}
\lesssim\left\|f\right\|_{h_X(\rrn)}.
\end{equation}

In addition, applying \eqref{relaxe12}
and Definition \ref{df511}(i), we
find that
$$
\left|f\ast\varphi\right|
=\left|f\ast\psi_{\frac{1}{2}}\right|
\leq m\left(f,\psi\right).
$$
From this, Definition \ref{Df1}(ii),
and Lemma \ref{Th6.1l1}(ii), it follows that
$$
\left\|f\ast\varphi\right\|_X
\leq\left\|m\left(f,\psi\right)\right\|_X
\sim\left\|f\right\|_{h_X(\rrn)}.
$$
This, together with \eqref{relaxe11},
further implies that
$$
\left\|f\ast\varphi\right\|
_{X}+\left\|f-f\ast\varphi\right\|_{H_X(\rrn)}
\lesssim\left\|f\right\|_{h_X(\rrn)},
$$
which completes the proof of \eqref{relaxe4}.
Combining both \eqref{relaxe0} and \eqref{relaxe4},
we obtain
$$
\left\|f\right\|_{h_X(\rrn)}\sim
\left\|f\ast\varphi\right\|
_{X}+\left\|f-f\ast\varphi\right\|_{H_X(\rrn)}.
$$
This finishes the proof of Theorem \ref{relax}.
\end{proof}

Based on the above relation
between $h_X(\rrn)$ and $H_X(\rrn)$,
we now show the following relation
between the local generalized Herz--Hardy space
$\HaSaHol$ and the generalized Herz--Hardy space
$\HaSaHo$.

\begin{theorem}\label{relal}
Let $p,\ q\in(0,\infty)$, $\omega\in M(\rp)$ with
$\m0(\omega)\in(-\frac{n}{p},\infty)$ and $\mi(\omega)
\in(-\frac{n}{p},\infty)$, and $\varphi
\in\mathcal{S}(\rrn)$ with
$$\1bf_{B(\0bf,1)}
\leq\widehat{\varphi}
\leq\1bf_{B(\0bf,2)}.$$
Then there exists a
constant $C\in[1,\infty)$ such that,
for any $f\in\mathcal{S}'(\rrn)$,
\begin{align*}
C^{-1}\|f\|_{\HaSaHol}
&\leq\left\|f\ast\varphi\right\|_{\HerzSo}
+\left\|f-f\ast\varphi\right\|_{\HaSaHo}\\
&\leq C\|f\|_{\HaSaHol}.
\end{align*}
\end{theorem}

\begin{proof}
Let all the symbols be as in the present theorem.
Then, combining the assumption
$\m0(\omega)\in(-\frac{n}{p},\infty)$
and Theorem \ref{Th3}, we
find that, under the assumptions
of the present theorem, the local
generalized Herz space $\HerzSo$
is a BQBF space. This implies
that, to complete the proof of
the present theorem, we only need to show that
the Herz space $\HerzSo$
under consideration satisfies all
the assumptions of Theorem \ref{relax}.

Indeed, let
$$
r\in\left(0,\min\left\{
p,q,\frac{n}{\Mw+n/p}\right\}\right).
$$
Then, applying Theorem \ref{strc}
with $s:=r$, we conclude
that $\HerzSo$ is strictly $r$-convex.
On the other hand, from
Lemma \ref{mbhl}, it follows
that, for any $f\in L^1_{\mathrm{loc}}(\rrn)$,
$$
\left\|\mc(f)\right\|_{[\HerzSo]^{1/r}}
\lesssim\left\|f\right\|_{[\HerzSo]^{1/r}}.
$$
Therefore, all the
assumptions of Theorem \ref{relax} hold
true for $\HerzSo$ and hence
the proof of Theorem \ref{relal} is completed.
\end{proof}

Theorem \ref{relal}, together with Remark \ref{remark5.0.9}(ii),
further implies that the following
conclusion holds true for the local
generalized Morrey--Hardy space $h\MorrSo$; we omit the details.

\begin{corollary}\label{relalm}
Let $p,\ q\in[1,\infty)$, $\omega\in M(\rp)$ with
$$
-\frac{n}{p}<\m0(\omega)\leq\M0(\omega)<0
$$
and
$$
-\frac{n}{p}<\mi(\omega)\leq\MI(\omega)<0,
$$
and $\varphi\in\mathcal{S}(\rrn)$ with
$\1bf_{B(\0bf,1)}\leq\widehat{\varphi
}\leq\1bf_{B(\0bf,2)}$.
Then there exists a constant $C\in[1,\infty)$
such that, for any $f\in\mathcal{S}'(\rrn)$,
\begin{align*}
C^{-1}\|f\|_{h\MorrSo}&\leq\left\|
f\ast\varphi\right\|_{\MorrSo}+
\left\|f-f\ast\varphi\right\|_{H\MorrSo}\\
&\leq C\|f\|_{h\MorrSo}.
\end{align*}
\end{corollary}

We next show the relation between
local generalized Herz--Hardy space
$\HaSaHl$ and generalized Herz--Hardy
space $\HaSaH$ as follows.

\begin{theorem}\label{relag}
Let $p,\ q\in(0,\infty)$, $\omega\in M(\rp)$ with
$\m0(\omega)\in(-\frac{n}{p},\infty)$ and
$$-\frac{n}{p}<\mi(\omega)
\leq\MI(\omega)<0,$$
and $\varphi\in\mathcal{S}(\rrn)$ with
$$\1bf_{B(\0bf,1)}
\leq\widehat{\varphi}
\leq\1bf_{B(\0bf,2)}.$$
Then there exists a
constant $C\in[1,\infty)$ such that,
for any $f\in\mathcal{S}'(\rrn)$,
\begin{align*}
C^{-1}\|f\|_{\HaSaHl}
&\leq\left\|f\ast\varphi\right\|_{\HerzS}
+\left\|f-f\ast\varphi\right\|_{\HaSaH}\\
&\leq C\|f\|_{\HaSaHl}.
\end{align*}
\end{theorem}

\begin{proof}
Let all the symbols be
as in the present theorem.
Then, since $\omega\in M(\rp)$
satisfies that
$\m0(\omega)\in(-\frac{n}{p},\infty)$
and $\MI(\omega)\in(-\infty,0)$,
from Theorem \ref{Th2},
it follows that the global
generalized Herz space $\HerzS$
under consideration is
a BQBF space. Thus,
to finish the proof of the
present theorem, it suffices to prove that
all the assumptions of Theorem \ref{relax} hold
true for $\HerzS$.
Namely, there exists an $r\in(0,\infty)$
such that $\HerzS$ is strictly
$r$-convex and the Hardy--Littlewood
maximal operator $\mc$ is bounded on
$[\HerzS]^{1/r}$.

To this end, let
$$
r\in\left(0,\min\left\{
p,q,\frac{n}{\Mw+n/p}\right\}\right).
$$
Then, by Theorem \ref{Th3.1}
with $s$ therein replaced by $r$, we find
that $\HerzS$ is strictly $r$-convex.
On the other hand, applying
Lemma \ref{mbhg}, we conclude
that, for any $f\in L^1_{\mathrm{loc}}(\rrn)$,
$$
\left\|\mc(f)\right\|_{[\HerzS]^{1/r}}
\lesssim\left\|f\right\|_{[\HerzS]^{1/r}}.
$$
This further implies that all the
assumptions of Theorem \ref{relax} hold
true for $\HerzS$ and
then finishes the proof of Theorem \ref{relag}.
\end{proof}

Via Theorem \ref{relag} and Remark \ref{remark5.0.9}(ii),
we immediately obtain the following relation between
the local Hardy space $h\MorrS$ and
the Hardy space $H\MorrS$ associated with
the global generalized Morrey space $\MorrS$;
we omit the details.

\begin{corollary}
Let $p$, $q$, $\omega$, and $\varphi$ be as
in Corollary \ref{relalm}.
Then there exists a constant
$C\in[1,\infty)$
such that, for any $f\in\mathcal{S}'(\rrn)$,
\begin{align*}
C^{-1}\|f\|_{h\MorrS}&\leq\left\|f\ast\varphi
\right\|_{\MorrS}+
\left\|f-f\ast\varphi\right\|_{H\MorrS}\\
&\leq C\|f\|_{h\MorrS}.
\end{align*}
\end{corollary}

\section{Atomic Characterizations}

The main target of this section is
to establish the atomic characterization
of the local generalized Herz--Hardy spaces
$\HaSaHol$ and $\HaSaHl$. Indeed,
by the known atomic characterization
of localized Hardy spaces associated with
ball quasi-Banach function spaces,
we obtain the atomic characterization
of $\HaSaHol$. On the other hand,
in order to show the atomic characterization
of $\HaSaHl$, recall that the associate
spaces of global generalized Herz spaces
are still unknown. To overcome this difficulty,
we first prove an improved atomic characterization
of the local Hardy space $h_X(\rrn)$,
with $X$ being a ball quasi-Banach function
space, without recourse to the associate space
$X'$. From this improved conclusion
of $h_X(\rrn)$, we deduce the desired atomic
characterization of $\HaSaHl$.

To begin with, we establish the atomic characterization
of the local generalized Herz--Hardy space
$\HaSaHol$. For this purpose,
we first introduce the definitions
of local atoms and the local atomic Hardy space
$\aHaSaHol$ associated
with the local generalized Herz space $\HerzSo$
as follows.

\begin{definition}\label{loatom}
Let $p$, $q\in(0,\infty)$, $\omega\in M(\rp)$
with $\m0(\omega)\in(-\frac{n}{p},\infty)$,
$r\in[1,\infty]$, and $d\in\zp$. Then a measurable
function $a$ is called
a \emph{local-$(\HerzSo,
\,r,\,d)$-atom}\index{local-$(\HerzSo,\,r,\,d)$-atom}
if
\begin{enumerate}
  \item[{\rm(i)}] there exists a ball
  $B(x_0,r_0)\in\mathbb{B}$, with $x_0\in\rrn$
  and $r_0\in(0,\infty)$, such that
  $\supp(a):=\{x\in\rrn:\ a(x)\neq0\}
  \subset B(x_0,r_0)$;
  \item[{\rm(ii)}] $\|a\|_{L^r(\rrn)}
  \leq\frac{|B(x_0,r_0)|^{1/r}}{\|\1bf_{B
  (x_0,r_0)}\|_{\HerzSo}}$;
  \item[{\rm(iii)}] when $r_0\in(0,1)$, then,
  for any $\alpha\in\zp^n$
  satisfying $|\alpha|\leq d$, $$
  \int_{\rrn}a(x)x^\alpha\,dx=0.$$
\end{enumerate}
\end{definition}

\begin{definition}\label{latsl}
Let $p$, $q\in(0,\infty)$, $\omega\in M(\rp)$ with
$\m0(\omega)\in(-\frac{n}{p},\infty)$
and $\mi(\omega)\in(-\frac{n}{p},\infty)$,
$$
s\in\left(0,\min\left\{1,p,q,
\frac{n}{\Mw+n/p}\right\}\right),
$$
$d\geq\lfloor n(1/s-1)\rfloor$ be a fixed integer,
and
$$
r\in\left(\max\left\{1,p,\frac{n}
{\mw+n/p}\right\},\infty\right].
$$
Then the \emph{local generalized
atomic Herz--Hardy space}
\index{local generalized atomic Herz--Hardy \\ space}
$\aHaSaHol$\index{$\aHaSaHol$},
associated with the local
generalized Herz space $\HerzSo$,
is defined to be the set of
all the $f\in\mathcal{S}'(\rrn)$ such that there exist
a sequence $\{a_{j}\}_{j\in\mathbb{N}}$ of
local-$(\HerzSo,\,r,\,d)$-atoms supported,
respectively, in the balls
$\{B_{j}\}_{j\in\mathbb{N}}\subset\mathbb{B}$
and a sequence $\{\lambda_{j}
\}_{j\in\mathbb{N}}\subset[0,\infty)$ satisfying that
$$
f=\sum_{j\in\mathbb{N}}\lambda_{j}a_{j}
$$
in $\mathcal{S}'(\rrn)$, and
$$
\left\|\left\{\sum_{j\in\mathbb{N}}
\left[\frac{\lambda_{j}}
{\|\1bf_{B_{j}}\|_{\HerzSo}}
\right]^{s}\1bf_{B_{j}}\right\}^{\frac{1}{s}}
\right\|_{\HerzSo}<\infty.
$$
Moreover, for any $f\in\aHaSaHol$,
$$
\|f\|_{\aHaSaHol}:=\inf\left\{
\left\|\left\{\sum_{j\in\mathbb{N}}
\left[\frac{\lambda_{j}}
{\|\1bf_{B_{j}}\|_{\HerzSo}}
\right]^{s}\1bf_{B_{j}}\right\}^{\frac{1}{s}}
\right\|_{\HerzSo}\right\},
$$
where the infimum is taken over
all the decompositions of $f$ as above.
\end{definition}

Then we have the following atomic
characterization\index{atomic characterization} of the local
generalized Herz--Hardy space $\HaSaHol$.

\begin{theorem}\label{lMolel}
Let $p$, $q$, $\omega$, $d$, $s$, and $r$ be as
in Definition \ref{latsl}.
Then $$\HaSaHol=\aHaSaHol$$ with equivalent quasi-norms.
\end{theorem}

To show the above atomic characterization,
we require the known atomic characterization
of the local Hardy space $h_X(\rrn)$
associated with ball quasi-Banach function
space $X$. First, we recall the following definition
of local-$(X,\,r,\,d)$-atoms,
which is just \cite[Definition 4.6]{WYY}.

\begin{definition}\label{loatomx}
Let $X$ be a ball quasi-Banach
function space,
$r\in[1,\infty]$, and $d\in\zp$.
Then a measurable
function $a$ is called
a \emph{local-$(X,\,r,
\,d)$-atom}\index{local-$(X,\,r,\,d)$-atom} if
\begin{enumerate}
  \item[{\rm(i)}] there exists a ball
  $B(x_0,r_0)\in\mathbb{B}$, with $x_0\in\rrn$
  and $r_0\in(0,\infty)$, such that
  $\supp(a):=\{x\in\rrn:\ a(x)\neq0\}\subset B(x_0,r_0)$;
  \item[{\rm(ii)}] $\|a\|_{L^r(\rrn)}
  \leq\frac{|B(x_0,r_0)|^{1/r}}{\|\1bf_{B
  (x_0,r_0)}\|_{X}}$;
  \item[{\rm(iii)}] when $r_0\in(0,1)$, then,
  for any $\alpha\in\zp^n$
  satisfying $|\alpha|\leq d$, $$
  \int_{\rrn}a(x)x^\alpha\,dx=0.$$
\end{enumerate}
\end{definition}

\begin{remark}\label{loatomxr}
Let $X$ be a ball quasi-Banach function space,
$r,\,t\in[1,\infty]$, and $d\in\zp$.
\begin{enumerate}
  \item[{\rm(i)}] Assume that
  $a$ is a local-$(X,\,t,\,d)$-atom
  supported in a ball $B\in\mathbb{B}$ and
  $r\leq t$. Then it is obviously
  that $a$ is also
  a local-$(X,\,r,\,d)$-atom supported in $B$.
  \item[{\rm(ii)}] Obviously, for
  any $(X,\,r,\,d)$-atom $a$ supported
  in the ball $B\in\mathbb{B}$, $a$
  is also a local-$(X,\,r,\,d)$-atom
  supported in $B$.
\end{enumerate}
\end{remark}

Via local-$(X,\,r,\,d)$-atoms,
we now present the following
definition of local atomic Hardy spaces
associated with ball quasi-Banach function
spaces (see, for instance, \cite{WYY}).

\begin{definition}\label{atom-hardyl}
Let $X$ be a ball quasi-Banach
function space, $r\in(1,\infty]$,
$d\in\zp$, and $s\in(0,1]$.
Then the \emph{local atomic Hardy space}
$h^{X,r,d,s}(\rrn)$\index{$h^{X,r,d,s}(\rrn)$},
associated with $X$, is defined to be the set of
all the $f\in\mathcal{S}'(\rrn)$ such that
there exist a sequence
$\{a_{j}\}_{j\in\mathbb{N}}$ of
local-$(X,\,r,\,d)$-atoms supported,
respectively, in the balls
$\{B_{j}\}_{j\in\mathbb{N}}\subset\mathbb{B}$
and a sequence $\{\lambda_{j}
\}_{j\in\mathbb{N}}\subset[0,\infty)$
satisfying that
$$
f=\sum_{j\in\mathbb{N}}\lambda_{j}a_{j}
$$
in $\mathcal{S}'(\rrn)$, and
$$
\left\|\left[\sum_{j\in\mathbb{N}}
\left(\frac{\lambda_{j}}
{\|\1bf_{B_{j}}\|_{X}}
\right)^{s}\1bf_{B_{j}}\right]^{\frac{1}{s}}
\right\|_{X}<\infty.
$$
Moreover, for any $f\in
h^{X,r,d,s}(\rrn)$,
$$
\|f\|_{h^{X,r,d,s}(\rrn)}:=\inf\left\{
\left\|\left[\sum_{j\in\mathbb{N}}
\left(\frac{\lambda_{j}}
{\|\1bf_{B_{j}}\|_{X}}
\right)^{s}\1bf_{B_{j}}\right]^{\frac{1}{s}}
\right\|_{X}\right\},
$$
where the infimum is taken over
all the decompositions of $f$ as above.
\end{definition}

Then we state the following
atomic characterization of
the local Hardy space $h_X(\rrn)$,
which was obtained in
\cite[Theorem 4.8]{WYY} and plays
an essential role in the proof of
Theorem \ref{lMolel}.

\begin{lemma}\label{atomxl}
Let $X$ be a ball quasi-Banach
function space satisfying both
Assumption \ref{assfs} with
$0<\theta<s\leq1$ and
Assumption \ref{assas} with the
same $s$ and $r\in(1,\infty]$, and
let $d\geq\lfloor n(1/\theta-1)\rfloor$
be a fixed integer.
Then
$$
h_X(\rrn)=
h^{X,r,d,s}(\rrn)
$$
with equivalent quasi-norms.
\end{lemma}

Via Lemma \ref{atomxl},
we now show Theorem \ref{lMolel}.

\begin{proof}[Proof of Theorem \ref{lMolel}]
Let $p$, $q$, $\omega$,
$d$, $s$, and $r$ be as in the present theorem.
Then, from the assumption $\m0(\omega)
\in(-\frac{n}{p},\infty)$
and Theorem \ref{Th3}, we deduce
that the local generalized Herz space
$\HerzSo$ under consideration is a
BQBF space. This implies that,
to complete the proof of the present theorem,
we only need to prove that the
Herz space $\HerzSo$ under
consideration satisfies all the assumptions of
Lemma \ref{atomxl}.

First, let $\theta\in(0,s)$
satisfy
\begin{equation}\label{lMolele1}
\left\lfloor n\left(\frac{1}{s}-1\right)
\right\rfloor\leq n\left(\frac{1}{\theta}
-1\right)<\left\lfloor n\left(
\frac{1}{s}-1\right)
\right\rfloor+1.
\end{equation}
We now show that $\HerzSo$
satisfies Assumption \ref{assfs}
with the above $\theta$ and $s$.
Indeed, by Lemma \ref{Atoll3},
we conclude that, for any
$\{f_{j}\}_{j\in\mathbb{N}}
\subset L^1_{{\rm loc}}(\rrn)$,
\begin{equation}\label{lMolele2}
\left\|\left\{\sum_{j\in\mathbb{N}}
\left[\mc^{(\theta)}
(f_{j})\right]^s\right\}^{1/s}
\right\|_{\HerzSo}\lesssim
\left\|\left(\sum_{j\in\mathbb{N}}
|f_{j}|^{s}\right)^{1/s}\right\|_{\HerzSo},
\end{equation}
which implies that
Assumption \ref{assfs} holds
true for $\HerzSo$ with the
above $\theta$ and $s$.

Next, we prove that $\HerzSo$
satisfies Assumption \ref{assas}
with the above $s$ and $r$.
Indeed, from Lemma \ref{mbhal},
we deduce that $[\HerzSo]^{1/s}$
is a BBF space and, for any $f\in
L_{{\rm loc}}^{1}(\rrn)$,
\begin{equation}\label{lMolele3}
\left\|\mc^{((r/s)')}(f)\right\|_
{([\HerzSo]^{1/s})'}\lesssim
\left\|f\right\|_{([\HerzSo]^{1/s})'}.
\end{equation}
This implies that, for the above
$s$ and $r$, Assumption
\ref{assas} holds true with $X:=\HerzSo$.
Finally, using \eqref{lMolele1},
we find that $d\geq\lfloor n(1/\theta-1)\rfloor$.
This, together with
\eqref{lMolele2}, the fact that
$[\HerzSo]^{1/s}$ is a BBF space,
and \eqref{lMolele3}, further implies
that, under the assumptions
of the present theorem,
$\HerzSo$ satisfies all the
assumptions of Lemma \ref{atomxl}.
Therefore, we have
$$
\HaSaHol=\aHaSaHol
$$
with equivalent quasi-norms,
which completes the proof of
Theorem \ref{lMolel}.
\end{proof}

As an application, we now
establish the atomic characterization of the local
generalized Morrey--Hardy space $h\MorrSo$. To achieve this,
we first introduce the definition of the local atoms
associated with the local generalized Morrey space $\MorrSo$ as follows.

\begin{definition}\label{loatomm}
Let $p$, $q\in[1,\infty)$, $\omega\in M(\rp)$
with $\MI(\omega)\in(-\infty,0)$ and
$$
-\frac{n}{p}<\m0(\omega)\leq\M0(\omega)<0,
$$
$r\in[1,\infty]$, and $d\in\zp$. Then a measurable
function $a$ on $\rrn$ is called
a \emph{local-$(\MorrSo,\,r,
\,d)$-atom}\index{local-$(\MorrSo,\,r,\,d)$-atom}
if
\begin{enumerate}
  \item[{\rm(i)}] there exists a ball
  $B(x_0,r_0)\in\mathbb{B}$, with $x_0\in\rrn$
  and $r_0\in(0,\infty)$, such that
  $\supp(a):=\{x\in\rrn:\ a(x)\neq0\}\subset B(x_0,r_0)$;
  \item[{\rm(ii)}] $\|a\|_{L^r(\rrn)}
  \leq\frac{|B(x_0,r_0)|^{1/r}}{\|\1bf_{B
  (x_0,r_0)}\|_{\MorrSo}}$;
  \item[{\rm(iii)}] when $r_0\in(0,1)$, then,
  for any $\alpha\in\zp^n$
  satisfying that $|\alpha|\leq d$, $$
  \int_{\rrn}a(x)x^\alpha\,dx=0.$$
\end{enumerate}
\end{definition}

Now, we give the following atomic characterization
\index{atomic characterization}of $h\MorrSo$, which
is a corollary of both Theorem \ref{lMolel} and
Remark \ref{remhs}(iv); we omit the details.

\begin{corollary}\label{lMolelm}
Let $p$, $q\in[1,\infty)$, $\omega\in M(\rp)$ with
$$
-\frac{n}{p}<\m0(\omega)\leq\M0(\omega)<0
$$
and
$$
-\frac{n}{p}<\mi(\omega)\leq\MI(\omega)<0,
$$
$s\in(0,1)$, $d\geq\lfloor
n(1/s-1)\rfloor$ be a fixed integer,
and
$$
r\in\left(\frac{n}{\mw+n/p},\infty\right].
$$
Then the \emph{local generalized
atomic Morrey--Hardy space}
\index{local generalized
atomic Morrey--Hardy space}$h\aMorrSo$\index{$h\aMorrSo$},
associated with the local generalized Morrey space $\MorrSo$,
is defined to be the set of all the
$f\in\mathcal{S}'(\rrn)$ such that there exist
a sequence $\{a_{j}\}_{j\in\mathbb{N}}$ of
local-$(\MorrSo,\,r,\,d)$-atoms supported,
respectively, in the balls
$\{B_{j}\}_{j\in\mathbb{N}}\subset\mathbb{B}$
and a sequence $\{\lambda_{j}\}
_{j\in\mathbb{N}}\subset[0,\infty)$
satisfying that
$$
f=\sum_{j\in\mathbb{N}}\lambda_{j}a_{j}
$$
in $\mathcal{S}'(\rrn)$,
and
$$
\left\|\left\{\sum_{j\in\mathbb{N}}
\left[\frac{\lambda_{j}}
{\|\1bf_{B_{j}}\|_{\MorrSo}}
\right]^{s}\1bf_{B_{j}}\right\}^{\frac{1}{s}}
\right\|_{\MorrSo}<\infty
$$
Moreover,
for any $f\in h\aMorrSo$,
$$
\|f\|_{h\aMorrSo}:=\inf\left\{
\left\|\left\{\sum_{j\in\mathbb{N}}\left[\frac{\lambda_{j}}
{\|\1bf_{B_{j}}\|_{\MorrSo}}\right]
^{s}\1bf_{B_{j}}\right\}^{\frac{1}{s}}
\right\|_{\MorrSo}\right\},
$$
where the infimum is
taken over all the decompositions of $f$ as above.
Then
$$
h\MorrSo=h\aMorrSo
$$
with equivalent quasi-norms.
\end{corollary}

The remainder of this section
is devoted to establishing the
atomic characterization of
the local generalized Herz--Hardy
space $\HaSaHl$. To this end, we first
introduce the following local-$(\HerzS,\,r,\,d)$-atoms
and the local atomic Hardy space associated with
$\HerzS$.

\begin{definition}\label{goatom}
Let $p$, $q\in(0,\infty)$, $\omega\in M(\rp)$
with $\m0(\omega)\in(-\frac{n}{p},\infty)$
and $\MI(\omega)\in(-\infty,0)$,
$r\in[1,\infty]$, and $d\in\zp$. Then a measurable
function $a$ is called
a \emph{local-$(\HerzS,
\,r,\,d)$-atom}\index{local-$(\HerzS,\,r,\,d)$-atom}
if
\begin{enumerate}
  \item[{\rm(i)}] there exists a ball
  $B(x_0,r_0)\in\mathbb{B}$, with $x_0\in\rrn$
  and $r_0\in(0,\infty)$, such that
  $\supp(a):=\{x\in\rrn:\ a(x)\neq0\}\subset B(x_0,r_0)$;
  \item[{\rm(ii)}] $\|a\|_{L^r(\rrn)}
  \leq\frac{|B(x_0,r_0)|^{1/r}}{\|\1bf_{B
  (x_0,r_0)}\|_{\HerzS}}$;
  \item[{\rm(iii)}] when $r_0\in(0,1)$, then,
  for any $\alpha\in\zp^n$
  satisfying $|\alpha|\leq d$, $$
  \int_{\rrn}a(x)x^\alpha\,dx=0.$$
\end{enumerate}
\end{definition}

\begin{definition}\label{latsg}
Let $p$, $q\in(0,\infty)$, $\omega\in M(\rp)$ with
$\m0(\omega)\in(-\frac{n}{p},\infty)$
and $$-\frac{n}{p}<
\mi(\omega)\leq\MI(\omega)<0,$$
$$
s\in\left(0,\min\left\{1,p,q,
\frac{n}{\Mw+n/p}\right\}\right),
$$
$d\geq\lfloor n(1/s-1)\rfloor$ be a fixed integer,
and
$$
r\in\left(\max\left\{1,p,\frac{n}
{\mw+n/p}\right\},\infty\right].
$$
Then the \emph{local generalized
atomic Herz--Hardy space}
\index{local generalized atomic Herz--Hardy \\ space}
$\aHaSaHl$\index{$\aHaSaHl$},
associated with the global
generalized Herz space $\HerzS$,
is defined to be the set of
all the $f\in\mathcal{S}'(\rrn)$ such that there exist
a sequence $\{a_{j}\}_{j\in\mathbb{N}}$ of
local-$(\HerzS,\,r,\,d)$-atoms supported,
respectively, in the balls
$\{B_{j}\}_{j\in\mathbb{N}}\subset\mathbb{B}$
and a sequence $\{\lambda_{j}
\}_{j\in\mathbb{N}}\subset[0,\infty)$
satisfying that
$$
f=\sum_{j\in\mathbb{N}}\lambda_{j}a_{j}
$$
in $\mathcal{S}'(\rrn)$, and
$$
\left\|\left\{\sum_{j\in\mathbb{N}}
\left[\frac{\lambda_{j}}
{\|\1bf_{B_{j}}\|_{\HerzS}}
\right]^{s}\1bf_{B_{j}}\right\}^{\frac{1}{s}}
\right\|_{\HerzS}<\infty.
$$
Moreover, for any $f\in\aHaSaHl$,
$$
\|f\|_{\aHaSaHl}:=\inf\left\{
\left\|\left\{\sum_{j\in\mathbb{N}}
\left[\frac{\lambda_{j}}
{\|\1bf_{B_{j}}\|_{\HerzS}}
\right]^{s}\1bf_{B_{j}}\right\}^{\frac{1}{s}}
\right\|_{\HerzS}\right\},
$$
where the infimum is taken over
all the decompositions of $f$ as above.
\end{definition}

Then we have the following atomic
characterization\index{atomic characterization} of the local
generalized Herz--Hardy space $\HaSaHl$.

\begin{theorem}\label{lMoleg}
Let $p$, $q$, $\omega$, $d$, $s$, and $r$ be
as in Definition \ref{latsg}.
Then $$\HaSaHl=\aHaSaHl$$ with equivalent quasi-norms.
\end{theorem}

To prove this atomic characterization,
we first establish an atomic characterization
of the local Hardy space $h_X(\rrn)$
via borrowing some ideas from \cite[Theorem 4.8]{WYY}
and get rid of the usage
of associate spaces. Namely, we have the following
conclusion.

\begin{theorem}\label{atomxxl}
Let $X$ be a ball quasi-Banach function space
satisfy:
\begin{enumerate}
  \item[{\rm(i)}] Assumption \ref{assfs}
  holds true with $0<\theta<s\leq1$;
  \item[{\rm(ii)}] for the above $s$,
  $X^{1/s}$ is a ball Banach function space and
  there exists a linear space
  $Y\subset\Msc(\rrn)$
  equipped with a seminorm $\|\cdot\|_{Y}$
  such that, for any $f\in\Msc(\rrn)$,
  \begin{equation}\label{atomxxle0}
  \|f\|_{X^{1/s}}
  \sim\sup\left\{\|fg\|_{L^1(\rrn)}:\
  \|g\|_{Y}=1\right\},
  \end{equation}
  where the positive equivalence
  constants are independent of $f$;
  \item[{\rm(iii)}] for the above
  $s$ and $Y$, there exist an $r\in(1,\infty]$
  and a positive constant $C$ such that,
  for any $f\in L^1_{\mathrm{loc}}(\rrn)$,
  $$
  \left\|\mc^{((r/s)')}
  (f)\right\|_{Y}\leq C\|f\|_{Y}.
  $$
\end{enumerate}
Then
$$
h_X(\rrn)=h^{X,r,d,s}(\rrn)
$$
with equivalent quasi-norms.
\end{theorem}

\begin{remark}
We should point out that
Theorem \ref{atomxxl} is an
improved version of the known
atomic characterization of $h_X(\rrn)$
obtained in \cite[Theorem 4.8]{WYY}.
Indeed, if $Y\equiv(X^{1/s})'$ in Theorem
\ref{atomxxl}, then this theorem goes back to
\cite[Theorem 4.8]{WYY}.
\end{remark}

To show Theorem \ref{atomxxl}, we
require the following
atomic decomposition of the local
Hardy space $h_X(\rrn)$, which
was obtained in \cite[p.\,37-39]{WYY}.

\begin{lemma}\label{atomxxll1}
Let $X$ be a
ball quasi-Banach function space
satisfying Assumption \ref{assfs}
with $0<\theta<s\leq1$,
$d\geq\lfloor n(1/\theta-1)\rfloor$
be a fixed integer, and $f\in h_X(\rrn)$.
Then there exist $\{\lambda_j\}_{j\in\mathbb{N}}
\subset[0,\infty)$ and
$\{a_j\}_{j\in\mathbb{N}}$
of local-$(X,\,\infty,\,d)$-atoms
supported, respectively, in the balls
$\{B_j\}_{j\in\mathbb{N}}\subset\mathbb{B}$
satisfying that
$$
f=\sum_{j\in\mathbb{N}}\lambda_{j}a_{j}
$$
in $\mathcal{S}'(\rrn)$, and
$$
\left\|\left[\sum_{j\in\mathbb{N}}
\left(\frac{\lambda_{j}}
{\|\1bf_{B_{j}}\|_{X}}
\right)^{s}\1bf_{B_{j}}\right]^{\frac{1}{s}}
\right\|_{X}\lesssim\left\|f\right\|_{h_X(\rrn)},
$$
where the implicit positive constant
is independent of $f$.
\end{lemma}

Now, we turn to prove Theorem \ref{atomxxl}.

\begin{proof}[Proof of Theorem \ref{atomxxl}]
Let $X$, $r$, $d$, and $s$ be as in
the present theorem. We first prove
that $h_X(\rrn)\subset h^{X,r,d,s}(\rrn)$.
For this purpose, let $f\in h_X(\rrn)$.
Then, applying the assumption (i) of the present
theorem and
Lemma \ref{atomxxll1}, we find that
there exist $\{\lambda_j\}_{j\in\mathbb{N}}
\subset[0,\infty)$ and
$\{a_j\}_{j\in\mathbb{N}}$
of local-$(X,\,\infty,\,d)$-atoms
supported, respectively, in the balls
$\{B_j\}_{j\in\mathbb{N}}\subset\mathbb{B}$
satisfying that
\begin{equation}\label{atomxxle1}
f=\sum_{j\in\mathbb{N}}\lambda_{j}a_{j}
\end{equation}
in $\mathcal{S}'(\rrn)$, and
\begin{equation}\label{atomxxle2}
\left\|\left[\sum_{j\in\mathbb{N}}
\left(\frac{\lambda_{j}}
{\|\1bf_{B_{j}}\|_{X}}
\right)^{s}\1bf_{B_{j}}\right]^{\frac{1}{s}}
\right\|_{X}\lesssim\left\|f\right\|_{h_X(\rrn)}.
\end{equation}
In addition, for any $j\in\mathbb{N}$,
from Remark \ref{loatomxr}(i)
with $a:=a_j$ and $t:=\infty$, it follows that
$a_j$ is a local-$(X,\,r,\,d)$-atom supported in
the ball $B_j$. This, combined with
\eqref{atomxxle1}, \eqref{atomxxle2},
and Definition \ref{atom-hardyl},
further implies that $f\in h^{X,r,d,s}(\rrn)$
and
\begin{equation}\label{atomxxle3}
\left\|f\right\|_{h^{X,r,d,s}(\rrn)}
\leq\left\|\left[\sum_{j\in\mathbb{N}}
\left(\frac{\lambda_{j}}
{\|\1bf_{B_{j}}\|_{X}}
\right)^{s}\1bf_{B_{j}}\right]^{\frac{1}{s}}
\right\|_{X}\lesssim\left\|f\right\|_{h_X(\rrn)},
\end{equation}
which completes the proof that
$h_X(\rrn)\subset h^{X,r,d,s}(\rrn)$.

Conversely, we next show that $h^{X,r,d,s}(\rrn)
\subset h_X(\rrn)$. Indeed, let
$f\in h^{X,r,d,s}(\rrn)$,
$\{\lambda_{l,j}\}_{l\in\{1,2\},\,j\in\mathbb{N}}
\subset[0,\infty)$, and $\{a_{l,j}\}
_{l\in\{1,2\},\,j\in\mathbb{N}}$
be a sequence of
local-$(X,\,r,\,d)$-atoms supported,
respectively, in the balls
$\{B_{l,j}\}_{l\in\{1,2\},\,j\in\mathbb{N}}
\subset\mathbb{B}$
satisfying that, for any
$j\in\mathbb{N}$, $r(B_{1,j})\in(0,1)$
and $r(B_{2,j})\in[1,\infty)$,
\begin{equation}\label{atomxxle4}
f=\sum_{l=1}^{2}\sum_{j\in\mathbb{N}}
\lambda_{l,j}a_{l,j}
\end{equation}
in $\mathcal{S}'(\rrn)$, and
\begin{equation}\label{atomxxle5}
\left\|\left[\sum_{l=1}^{2}\sum_{j\in\mathbb{N}}
\left(\frac{\lambda_{l,j}}
{\|\1bf_{B_{l,j}}\|_{X}}
\right)^{s}\1bf_{B_{l,j}}\right]^{\frac{1}{s}}
\right\|_{X}<\infty.
\end{equation}
In what follows, let $\phi\in\mathcal{S}(\rrn)$
be such that $\supp(\phi)\subset B(\0bf,1)$
and $\int_{\rrn}\phi(x)\,dx\neq0$.
Then, using \eqref{atomxxle4}
and repeating an argument similar to
that used in the estimation of \eqref{Atogxe10}
with $\{\lambda_j\}_{j\in\mathbb{N}}$
and $\{a_j\}_{j\in\mathbb{N}}$ therein
replaced, respectively,
by $\{\lambda_{l,j}\}_{l\in\{1,2\},
\,j\in\mathbb{N}}$, and $\{a_{l,j}\}
_{l\in\{1,2\},\,j\in\mathbb{N}}$,
we conclude that, for any $t\in(0,\infty)$,
\begin{align*}
\left|f\ast\phi_t\right|
\leq\sum_{l=1}^{2}\sum_{j\in\mathbb{N}}
\lambda_{l,j}\left|a_{l,j}\ast\phi_t\right|.
\end{align*}
This, together with
Definition \ref{df511}(i), further implies that
\begin{equation*}
m(f,\phi)\leq\sum_{l=1}^{2}\sum_{j\in\mathbb{N}}
\lambda_{l,j}m(a_{l,j},\phi).
\end{equation*}
Using this and Definition \ref{Df1}(ii),
we find that
\begin{align}\label{atomxxle6}
\left\|m(f,\phi)\right\|_X
&\lesssim\left\|\sum_{j\in\mathbb{N}}
\lambda_{1,j}m(a_{1,j},\phi)\right\|_X
+\left\|\sum_{j\in\mathbb{N}}
\lambda_{2,j}m(a_{2,j},\phi)\right\|_X\notag\\
&=:\mathrm{IV_1}+\mathrm{IV_2}.
\end{align}

We first estimate $\mathrm{IV}_1$. Indeed,
for any given $j\in\mathbb{N}$,
by the assumption $r(B_{1,j})\in(0,1)$
and both Definitions \ref{loatomx} and \ref{atom},
we conclude that $a_j$ is an
$(X,\,r,\,d)$-atom supported in $B_j$.
Thus, using both Definitions \ref{df511}(i)
and \ref{smax}(i) and
repeating the arguments similar to that used
in the estimations of \eqref{Atogxe2},
\eqref{Atogxe5}, and \eqref{Atogxe6}
with $\{\lambda_j\}_{j\in\mathbb{N}}$,
$\{a_j\}_{j\in\mathbb{N}}$, and
$\{B_j\}_{j\in\mathbb{N}}$ therein
replaced, respectively,
by $\{\lambda_{1,j}\}_{j\in
\mathbb{N}}$, $\{a_{1,j}\}
_{j\in\mathbb{N}}$,
and $\{B_{1,j}\}_{j\in\mathbb{N}}$,
we find that
\begin{align}\label{atomxxle7}
\mathrm{IV}_1\lesssim
\left\|\sum_{j\in\mathbb{N}}
\lambda_{1,j}M(a_{1,j},\phi)\right\|_X\lesssim
\left\|\left[\sum_{j\in\mathbb{N}}
\left(\frac{\lambda_{1,j}}
{\|\1bf_{B_{1,j}}\|_{X}}
\right)^{s}\1bf_{B_{1,j}}\right]^{\frac{1}{s}}
\right\|_{X}.
\end{align}
This is the desired estimate of $\mathrm{IV}_1$.

On the other hand, we deal with $\mathrm{IV}_2$.
To this end, we first claim that,
for any given $j\in\mathbb{N}$ and for any $x
\in(2B_{2,j})^{\complement}$,
$m(a_{2,j},\phi)(x)=0$.
Indeed, fix a $j\in\mathbb{N}$.
Then, applying the assumption
$r(B_{2,j})\in[1,\infty)$, we find
that, for any $t\in(0,1)$,
$x\in(2B_{2,j})^{\complement}$,
and $y\in B_{2,j}$,
$$
|x-y|\geq r(B_{2,j})\geq1>t,
$$
which implies that
$y\in[B(x,t)]^{\complement}$.
This, together with the
assumption $\supp(\phi)
\subset B(\0bf,1)$, further implies that,
for any $t\in(0,1)$ and $x\in(2B_{2,j})^{\complement}$,
\begin{align*}
a_{2,j}\ast\phi_t(x)
&=\int_{\rrn}\phi_t(x-y)a_{2,j}(y)\,dy\\
&=\int_{B_{2,j}\cap B(x,t)}\phi_t(x-y)a_{2,j}(y)\,dy=0.
\end{align*}
Therefore, by Definition \ref{df511}(i),
we further conclude that,
for any given $j\in\mathbb{N}$ and for any $x
\in(2B_{2,j})^{\complement}$,
$m(a_{2,j},\phi)(x)=0$,
which completes the proof
of the above claim.
Combining this claim and both
Definitions \ref{df511}(i)
and \ref{smax}(i), we obtain
\begin{align*}
\mathrm{IV}_2&\sim
\left\|\sum_{j\in\mathbb{N}}
\lambda_{2,j}m(a_{2,j},\phi)\1bf_{
2B_{2,j}}\right\|_X
\lesssim\left\|\sum_{j\in\mathbb{N}}
\lambda_{2,j}M(a_{2,j},\phi)\1bf_{
2B_{2,j}}\right\|_X.
\end{align*}
From this and an argument similar to that
used in the estimation of
$\mathrm{II}_1$ in the proof of
Theorem \ref{Atogx}
with $\{\lambda_j\}_{j\in\mathbb{N}}$,
$\{a_j\}_{j\in\mathbb{N}}$, and
$\{B_j\}_{j\in\mathbb{N}}$ therein
replaced, respectively,
by $\{\lambda_{2,j}\}_{j\in
\mathbb{N}}$, $\{a_{2,j}\}
_{j\in\mathbb{N}}$,
and $\{B_{2,j}\}_{j\in\mathbb{N}}$,
it follows that
\begin{align}\label{atomxxle8}
\mathrm{IV}_2\lesssim
\left\|\left[\sum_{j\in\mathbb{N}}
\left(\frac{\lambda_{2,j}}
{\|\1bf_{B_{2,j}}\|_{X}}
\right)^{s}\1bf_{B_{2,j}}\right]^{\frac{1}{s}}
\right\|_{X},
\end{align}
which is the desired estimate of
$\mathrm{IV}_2$. Thus,
by both the assumptions (i) and (ii)
of the present theorem,
Remark \ref{Th6.1r}, Lemma
\ref{Th6.1l1}(ii), \eqref{atomxxle6},
\eqref{atomxxle7},
\eqref{atomxxle8}, and \eqref{atomxxle5},
we find that
\begin{equation}\label{atomxxle9}
\left\|f\right\|_{h_X(\rrn)}
\sim\left\|m(f,\phi)\right\|_{X}
\lesssim\left\|\left[\sum_{l=1}^{2}\sum_{j\in\mathbb{N}}
\left(\frac{\lambda_{l,j}}
{\|\1bf_{B_{l,j}}\|_{X}}
\right)^{s}\1bf_{B_{l,j}}\right]^{\frac{1}{s}}
\right\|_{X}<\infty.
\end{equation}
This further implies that $f\in h_X(\rrn)$
and hence $h_X(\rrn)\subset h^{X,r,d,s}(\rrn)$.
Moreover, from \eqref{atomxxle9},
the choice of $\{\lambda_{l,j}\}_{
l\in\{1,2\},j\in\mathbb{N}}$, and
Definition \ref{atom-hardyl},
we deduce that
$$
\left\|f\right\|_{h_X(\rrn)}
\lesssim\left\|f\right\|_{h^{X,r,d,s}(\rrn)},
$$
which, combined with \eqref{atomxxle3},
implies that
$$
\left\|f\right\|_{h_X(\rrn)}
\sim\left\|f\right\|_{h^{X,r,d,s}(\rrn)}.
$$
Thus, we have
$$
h_X(\rrn)=h^{X,r,d,s}(\rrn)
$$
with equivalent quasi-norms, which
completes the proof of Theorem \ref{atomxxl}.
\end{proof}

Via the above atomic
characterization of local
Hardy spaces associated with ball
quasi-Banach function spaces,
we now prove Theorem \ref{lMoleg}.

\begin{proof}[Proof of Theorem \ref{lMoleg}]
Let $p$, $q$, $\omega$, $r$, $s$, and $d$
be as in the present theorem. Then, combining
the assumptions $\m0(\omega)\in(-\frac{n}{p},
\infty)$ and $\MI(\omega)\in(-\infty,0)$,
and Theorem \ref{Th2}, we conclude that
the global generalized Herz space $\HerzS$
under consideration is a BQBF space.
From this and Theorem \ref{atomxxl},
it follows that, to finish
the proof of the present theorem, we only need
to show that
the assumptions (i) through (iii)
of Theorem \ref{atomxxl} hold true for
$\HerzS$.

First, we show that Theorem
\ref{atomxxl}(i) holds true for $\HerzS$.
To this end, let $\theta\in(0,s)$ be such that
\begin{equation}\label{lMolege1}
\left\lfloor
n\left(\frac{1}{s}-1\right)
\right\rfloor\leq n\left(
\frac{1}{\theta}-1
\right)<\left\lfloor
n\left(\frac{1}{s}-1\right)
\right\rfloor+1.
\end{equation}
Then, applying Lemma \ref{Atogl5},
we find that, for any $\{f_j\}_{
j\in\mathbb{N}}\subset L^1_{\mathrm{loc}}(
\rrn)$,
\begin{equation*}
\left\|\left\{
\sum_{j\in\mathbb{N}}
\left[
\mc^{(\theta)}(f_j)
\right]^{s}\right\}^{1/s}\right\|_{\HerzS}
\lesssim\left\|\left(\sum_{j\in\mathbb{N}}
|f_j|^s\right)^{1/s}\right\|_{\HerzS}.
\end{equation*}
This implies that,
for the above $\theta$ and
$s$, $\HerzS$ satisfies Assumption
\ref{assfs} and hence
Theorem \ref{atomxxl}(i) holds
true.

Next, we prove that $\HerzS$ satisfies Theorem
\ref{atomxxl}(ii). Indeed, from
the assumptions $\m0(\omega)
\in(-\frac{n}{p},\infty)$
and $\MI(\omega)\in(-\infty,0)$, and
Lemma \ref{rela}, it follows that
\begin{equation*}
\m0\left(\omega^s\right)
=s\m0(\omega)>-\frac{n}{p/s}
\end{equation*}
and
\begin{equation*}
\MI\left(\omega^s\right)=s
\MI(\omega)<0.
\end{equation*}
Applying these,
the assumptions $p/s,\ q/s\in(1,\infty)$,
and Theorem \ref{ball}
with $p$, $q$, and $\omega$
replaced, respectively, by $p/s$,
$q/s$, and $\omega^s$, we find that
$\Kmp^{p/s,q/s}_{\omega^s}(\rrn)$
is a BBF space. Moreover, by
Lemma \ref{convexl}, we
conclude that
$$\lf[\HerzS\r]^{1/s}=
\Kmp^{p/s,q/s}_{\omega^s}(\rrn).$$
Thus, $[\HerzS]^{1/s}$ is a BBF space.
On the
other hand, from \eqref{atoge3},
it follows that,
for any $f\in\Msc(\rrn)$,
\begin{equation*}
\|f\|_{[\HerzS]^{1/s}}\sim
\sup\left\{\|fg\|_{L^1(\rrn)}:\
\|g\|_{\dot{\mathcal{B}}_{1/\omega^s}
^{(p/s)',(q/s)'}(\rrn)}=1\right\}.
\end{equation*}
Therefore, $[\HerzS]^{1/s}$ is
a BBF space and \eqref{atomxxle0}
holds true with $$Y:=
\dot{\mathcal{B}}_{1/\omega^s}
^{(p/s)',(q/s)'}(\rrn).$$
These further imply that
the Herz space $\HerzS$
under consideration satisfies Theorem
\ref{atomxxl}(ii).

Finally, we show that Theorem
\ref{atomxxl}(iii) holds true for $\HerzS$.
Indeed, using Lemma \ref{mbhag},
we conclude that, for any
$f\in L_{{\rm loc}}^{1}(\rrn)$,
\begin{equation*}
\left\|\mc^{((r/s)')}(f)\right\|_{\dot{\mathcal{B}}
_{1/\omega^{s}}^{(p/s)',(q/s)'}(\rrn)}\lesssim
\|f\|_{\dot{\mathcal{B}}_
{1/\omega^{s}}^{(p/s)',(q/s)'}(\rrn)},
\end{equation*}
which implies that Theorem \ref{atomxxl}(iii)
holds true for $\HerzS$ with
$$Y:=\dot{\mathcal{B}}_{1/\omega^s}
^{(p/s)',(q/s)'}(\rrn).$$
Moreover, by \eqref{lMolege1},
we find that
$d\geq\lfloor
n(1/\theta-1)\rfloor$.
This, together with the fact that
the assumptions (i)
through (iii) of Theorem \ref{atomxxl}
hold true for the Herz
space $\HerzS$ under consideration,
further implies that
$$
\HaSaHl=\aHaSaHl
$$
with equivalent quasi-norms,
which completes the proof of
Theorem \ref{lMoleg}.
\end{proof}

As an application of Theorem
\ref{lMoleg}, we now establish
the atomic characterization of the local generalized
Morrey--Hardy space $h\MorrS$. To this end,
we first introduce the definition
of local-$(\MorrS,\,r,\,d)$-atoms
as follows.

\begin{definition}
Let $p$, $q$, $\omega$, $r$, and $d$ be as in
Definition \ref{loatomm}.
Then a measurable
function $a$ on $\rrn$ is called
a \emph{local-$(\MorrS,\,r,\,d)$-atom}\index{local-$(\MorrS,\,r,\,d)$-atom}
if
\begin{enumerate}
  \item[{\rm(i)}] there exists a ball $B(x_0,r_0)
  \in\mathbb{B}$, with $x_0\in\rrn$ and $r_0\in(0,\infty)$,
  such that $\supp(a):=
  \{x\in\rrn:\ a(x)\neq0\}\subset B(x_0,r_0)$;
  \item[{\rm(ii)}] $\|a\|_{L^r(\rrn)}\leq
  \frac{|B(x_0,r_0)|^{1/r}}{\|\1bf_{B(x_0,r_0)}\|_{\MorrS}}$;
  \item[{\rm(iii)}] when $r_0\in(0,1)$, then, for any $\alpha\in\zp^n$
  satisfying $|\alpha|\leq d$, $$\int_{\rrn}a(x)x^\alpha\,dx=0.$$
\end{enumerate}
\end{definition}

Then we have the following atomic
characterization\index{atomic characterization} of $h\MorrS$,
which can be deduced from Theorem \ref{lMoleg}
and Remark \ref{remhs}(iv) immediately; we omit the details.

\begin{corollary}
Let $p$, $q$, $\omega$,
$r$, $d$, and $s$ be as in
Corollary \ref{lMolelm}.
Then the \emph{local generalized
atomic Morrey--Hardy space}
\index{local generalized
atomic Morrey--Hardy space}$h\aMorrS$\index{$h\aMorrS$},
associated with the global generalized
Morrey space $\MorrS$,
is defined to be the set of all the
$f\in\mathcal{S}'(\rrn)$ such that there exist
a sequence $\{a_{j}\}_{j\in\mathbb{N}}$ of
local-$(\MorrS,\,r,\,d)$-atoms supported,
respectively, in the balls
$\{B_{j}\}_{j\in\mathbb{N}}\subset\mathbb{B}$
and a sequence $\{\lambda_{j}\}
_{j\in\mathbb{N}}\subset[0,\infty)$
satisfying that
$
f=\sum_{j\in\mathbb{N}}\lambda_{j}a_{j}
$
in $\mathcal{S}'(\rrn)$,
and
$$
\left\|\left\{\sum_{j\in\mathbb{N}}
\left[\frac{\lambda_{j}}
{\|\1bf_{B_{j}}\|_{\MorrS}}
\right]^{s}\1bf_{B_{j}}\right\}^{\frac{1}{s}}
\right\|_{\MorrS}<\infty.
$$
Moreover,
for any $f\in h\aMorrS$,
$$
\|f\|_{h\aMorrS}:=\inf\left\{
\left\|\left\{\sum_{j\in\mathbb{N}}
\left[\frac{\lambda_{j}}
{\|\1bf_{B_{j}}\|_{\MorrS}}\right]
^{s}\1bf_{B_{j}}\right\}^{\frac{1}{s}}
\right\|_{\MorrS}\right\},
$$
where the infimum is
taken over all the decompositions of $f$ as above.
Then
$$
h\MorrS=h\aMorrS
$$
with equivalent quasi-norms.
\end{corollary}

\section{Molecular Characterizations}

In this section, we investigate
the molecular characterization
of localized generalized Herz--Hardy spaces
via viewing generalized Herz spaces as
special cases of ball quasi-Banach function spaces.
Precisely, we establish the molecular
characterization of the local generalized
Herz--Hardy space $\HaSaHol$ via using the
known molecular characterization
of localized Hardy spaces associated with
ball quasi-Banach function spaces obtained
in \cite[Theorem 5.2]{WYY} (see also Lemma
\ref{llmolell1} below).
On the other hand, to prove the molecular
characterization of the local
generalized Herz--Hardy space $\HaSaHl$, recall
that the associate spaces
of the global generalized Herz spaces
are still unknown. To overcome
this obstacle, we
first establish an improved molecular
characterization of localized Hardy spaces
associated with ball quasi-Banach function
spaces (see Theorem \ref{molexl} below)
without recourse to associate spaces. Combining
this molecular characterization and
the fact that the global generalized
Herz space $\HerzS$ is a special
ball quasi-Banach function space,
we then obtain the
desired molecular characterization of $\HaSaHl$.

We first investigate the molecular characterization
of $\HaSaHol$. To begin with, we introduce the local molecules
associated with the local generalized Herz space $\HerzSo$
as follows.

\begin{definition}\label{lomole}
Let $p$, $q\in(0,\infty)$, $\omega\in M(\rp)$
with $\m0(\omega)\in(-\frac{n}{p},\infty)$,
$r\in[1,\infty]$,
$d\in\mathbb{Z}_{+}$, and $\tau\in
(0,\infty)$. Then a measurable
function $m$ on $\rrn$ is called a
\emph{local-$(\HerzSo,\,r,\,d,\,\tau)$-molecule}
\index{local-$(\HerzSo,\,r,\,d,
\,\tau)$-molecule}centered at a ball
$B(x_0,r_0)\in\mathbb{B}$, with $x_0\in
\rrn$ and $r_0\in(0,\infty)$, if
\begin{enumerate}
  \item[{\rm(i)}] for any $i\in\mathbb{Z}_{+}$,
$$\left\|m\1bf_{S_{i}(B(x_0,r_0))}
\right\|_{L^{r}(\rrn)}\leq2^{-\tau i}
\frac{|B(x_0,r_0)|^{1/r}}{\|\1bf_{B(x_0,r_0)}
\|_{\HerzSo}};$$
  \item[{\rm(ii)}] when $r_0\in(0,1)$, then, for any
  $\alpha\in\zp^n$ with $|\alpha|\leq d$,
  $$\int_{\rrn}m(x)x^\alpha\,dx=0.$$
\end{enumerate}
\end{definition}

Then we establish the following molecular characterization
\index{molecular \\characterization}of
the local generalized Herz--Hardy space $\HaSaHol$.

\begin{theorem}\label{llmolel}
Let $p,\ q\in(0,\infty)$,
$\omega\in M(\rp)$ with $\m0(\omega)
\in(-\frac{n}{p},\infty)$
and $\mi(\omega)\in(-\frac{n}{p},\infty)$,
$$
s\in\left(0,\min\left\{1,p,q,\frac{n}
{\Mw+n/p}\right\}\right),
$$
$d\geq\lfloor n(1/s-1)\rfloor$ be a fixed integer,
$$r\in\left(\max\left\{1,p,\frac{n}{
\mw+n/p}\right\},\infty\right],$$
and $\tau\in(0,\infty)$ with $\tau>
n(1/s-1/r)$.
Then $f\in\HaSaHol$ if and only if $f\in\mathcal{S}'(\rrn)$
and there exist a sequence $\{m_{j}\}_{j\in\mathbb{N}}$
of local-$(\HerzSo,\,r,\,d,
\,\tau)$-molecules centered, respectively,
at the balls
$\{B_{j}\}_{j\in\mathbb{N}}\subset\mathbb{B}$ and
a sequence $\{\lambda_{j}\}_{j\in\mathbb{N}}\subset
[0,\infty)$ such that
\begin{equation*}
f=\sum_{j\in\mathbb{N}}\lambda_{j}a_{j}
\end{equation*}
in $\mathcal{S}'(\rrn)$, and
$$
\left\|\left\{\sum_{j\in\mathbb{N}}\left[\frac{\lambda_{j}}
{\|\1bf_{B_{j}}\|_{\HerzSo}}\right]^s\1bf_
{B_{j}}\right\}^{\frac{1}{s}}\right
\|_{\HerzSo}<\infty.
$$
Moreover, there exists a
constant $C\in[1,\infty)$ such that, for
any $f\in\HaSaHol$,
\begin{align*}
C^{-1}\|f\|_{\HaSaHol}&\leq\inf\left\{
\left\|\left\{\sum_{j\in\mathbb{N}}
\left[\frac{\lambda_{i}}
{\|\1bf_{B_{j}}\|_{\HerzSo}}\right]^s\1bf_
{B_{j}}\right\}^{\frac{1}{s}}
\right\|_{\HerzSo}\right\}\\
&\leq C\|f\|_{\HaSaHol},
\end{align*}
where the infimum is taken
over all the decompositions of $f$ as above.
\end{theorem}

To obtain this molecular characterization,
we first recall the following
definition of local-$(X,\,r,\,d,\,\tau)$-molecules
with $X$ being a ball quasi-Banach function space,
which is just \cite[Definition 5.1]{WYY}.

\begin{definition}\label{lomolex}
Let $X$ be a ball quasi-Banach
function space, $r\in[1,\infty]$,
$d\in\mathbb{Z}_{+}$, and $\tau\in
(0,\infty)$. Then a measurable
function $m$ on $\rrn$ is called a
\emph{local-$(X,\,r,\,d,\,\tau)$-molecule}
\index{local-$(X,\,r,\,d,
\,\tau)$-molecule}centered at a ball
$B(x_0,r_0)\in\mathbb{B}$, with $x_0\in
\rrn$ and $r_0\in(0,\infty)$, if
\begin{enumerate}
  \item[{\rm(i)}] for any $i\in\mathbb{Z}_{+}$,
$$\left\|m\1bf_{S_{i}(B(x_0,r_0))}
\right\|_{L^{r}(\rrn)}\leq2^{-\tau i}
\frac{|B(x_0,r_0)|^{1/r}}{\|\1bf_{B(x_0,r_0)}
\|_{X}};$$
  \item[{\rm(ii)}] when $r_0\in(0,1)$,
  then, for any
  $\alpha\in\zp^n$ with $|\alpha|\leq d$,
  $$\int_{\rrn}m(x)x^\alpha\,dx=0.$$
\end{enumerate}
\end{definition}

\begin{remark}\label{lomolexr}
Let $X$ be a ball quasi-Banach function
space, $r\in[1,\infty]$, $d\in\zp$,
and $\tau\in(0,\infty)$. Then
it is obvious that
\begin{enumerate}
  \item[{\rm(i)}] for any local-$(X,\,r,\,d)$-atom
$a$ supported in $B\in\mathbb{B}$,
$a$ is a local-$(X,\,r,\,d,\,\tau)$-molecule centered
at $B$;
  \item[{\rm(ii)}] for any $(X,\,r,\,d,\,\tau)$-molecule
$m$ centered at a ball $B\in\mathbb{B}$,
$m$ is a local-$(X,\,r,\,d,\,\tau)$-molecule centered
at $B$.
\end{enumerate}
\end{remark}

The following atomic characterization
of $h_X(\rrn)$ was established in
\cite[Theorem 5.2]{WYY}, which is
an essential tool in the proof of
Theorem \ref{llmolel}.

\begin{lemma}\label{llmolell1}
Let $X$ be a ball quasi-Banach
function space satisfying both
Assumption \ref{assfs} with
$0<\theta<s\leq1$ and Assumption
\ref{assas} with the same $s$
and $r\in(1,\infty]$, and let
$d\geq\lfloor n(1/\theta-1)\rfloor$
be a fixed integer.
Then $f\in h_X(\rrn)$ if and only
if $f\in\mathcal{S}'(\rrn)$
and there exist a sequence
$\{m_{j}\}_{j\in\mathbb{N}}$
of local-$(X,\,r,\,d,\,\tau)$-molecules
centered, respectively,
at the balls $\{B_{j}\}_{j\in\mathbb{N}}
\subset\mathbb{B}$ and
a sequence $\{\lambda_{j}\}
_{j\in\mathbb{N}}\subset
[0,\infty)$ such that
$
f=\sum_{j\in\mathbb{N}}\lambda_{j}a_{j}
$
in $\mathcal{S}'(\rrn)$, and
$$
\left\|\left[\sum_{j\in\mathbb{N}}
\left(\frac{\lambda_{j}}
{\|\1bf_{B_{j}}\|_{X}}
\right)^{s}\1bf_{B_{j}}
\right]^{\frac{1}{s}}
\right\|_{X}<\infty.
$$
Moreover, for any $f\in h_X(\rrn)$,
\begin{align*}
\|f\|_{h_X(\rrn)}\sim
\inf\left\{
\left\|\left[\sum_{j\in\mathbb{N}}
\left(\frac{\lambda_{i}}
{\|\1bf_{B_{j}}\|_{X}}\right)^s\1bf_
{B_{j}}\right]^{\frac{1}{s}}\right\|_{X}
\right\}
\end{align*}
with the positive equivalence
constants independent of $f$,
where the infimum is taken
over all the decompositions of $f$ as above.
\end{lemma}

Based on the above lemma, we next
show Theorem \ref{llmolel}.

\begin{proof}[Proof of Theorem \ref{llmolel}]
Let $p$, $q$, $\omega$,
$r$, $d$, and $s$
be as in the present theorem. Since
$\omega$ satisfies
that $\m0(\omega)\in(-\frac{n}{p},
\infty)$, from Theorem \ref{Th3},
it follows that the local
generalized Herz space $\HerzSo$
under consideration is a BQBF space.
Thus, to complete the present theorem,
it suffices to prove that
all the assumptions of Lemma \ref{llmolell1}
hold true for the Herz space $\HerzSo$
under consideration. Namely,
$\HerzSo$ satisfies both Assumption \ref{assfs}
with the above $s$ and some $\theta\in(0,s)$,
and Assumption \ref{assas} with the above
$s$ and $r$.

Indeed, let $\theta\in(0,s)$ satisfy
that \begin{equation}\label{llmolele1}
\left\lfloor n\left(\frac{1}{s}-1\right)
\right\rfloor\leq n\left(\frac{1}{\theta}
-1\right)<\left\lfloor n\left(
\frac{1}{s}-1\right)
\right\rfloor+1
\end{equation}
and
\begin{equation}\label{llmolele2}
\tau>n\left(\frac{1}{\theta}-
\frac{1}{r}\right).
\end{equation}
Then, by Lemma \ref{Atogl5},
we find that,
for any $\{f_{j}\}_{
j\in\mathbb{N}}\subset L^{1}_{\text{loc}}(\rrn)$,
\begin{equation*}
\left\|\left\{\sum_{j\in\mathbb{N}}\left[\mc^{(\theta)}
(f_{j})\right]^s\right\}^{1/s}
\right\|_{\HerzSo}\lesssim
\left\|\left(\sum_{j\in\mathbb{N}}
|f_{j}|^{s}\right)^{1/s}\right\|_{\HerzSo},
\end{equation*}
which implies that, for the above $\theta$
and $s$, Assumption \ref{assfs} holds
true for $\HerzSo$.
On the other hand, applying Lemma \ref{mbhal},
we conclude that $[\HerzSo]^{1/s}$ is a BBF space
and, for any $f\in L_{{\rm loc}}^{1}(\rrn)$,
\begin{equation*}
\left\|\mc^{((r/s)')}(f)\right\|
_{([\HerzSo]^{1/s})'}\lesssim
\left\|f\right\|_{([\HerzSo]^{1/s})'}.
\end{equation*}
This implies that $\HerzSo$ satisfies
Assumption \ref{assas} with the above
$s$ and $r$.
Finally, from \eqref{llmolele1},
we deduce that $d\geq\lfloor n(1/\theta-1)
\rfloor$. Combining this,
\eqref{llmolele2}, and the facts
that $\HerzSo$ satisfies both Assumptions
\ref{assfs} and \ref{assas},
we find that the Herz space
$\HerzSo$ under consideration
satisfies all the assumptions of
Lemma \ref{llmolell1}
and then complete the
proof of Theorem \ref{llmolel}.
\end{proof}

As an application of Theorem \ref{llmolel}, we now establish
the molecular characterization of
the local generalized Morrey--Hardy space
$h\MorrSo$ via introducing the definition
of the local-$(\MorrSo,\,r,\,d,\,\tau)$-molecules as follows.

\begin{definition}\label{lomolem}
Let $p$, $q\in[1,\infty)$, $\omega\in M(\rp)$
with $\MI(\omega)\in(-\infty,0)$ and
$$
-\frac{n}{p}<\m0(\omega)\leq\M0(\omega)<0,
$$
$r\in[1,\infty]$, $d\in\mathbb{Z}_{+}$, and $\tau\in
(0,\infty)$. Then a measurable
function $m$ on $\rrn$ is called a
\emph{local-$(\MorrSo,\,r,\,d,\,\tau)$-molecule}
\index{local-$(\MorrSo,\,r,\,d,
\,\tau)$-molecule}centered at a ball
$B(x_0,r_0)\in\mathbb{B}$, with $x_0\in
\rrn$ and $r_0\in(0,\infty)$, if
\begin{enumerate}
  \item[{\rm(i)}] for any $i\in\mathbb{Z}_{+}$,
$$\left\|m\1bf_{S_{i}(B(x_0,r_0))}
\right\|_{L^{r}(\rrn)}\leq2^{-\tau i}
\frac{|B(x_0,r_0)|^{1/r}}{\|\1bf_{B(x_0,r_0)}
\|_{\MorrSo}};$$
  \item[{\rm(ii)}] when $r_0\in(0,1)$, then, for any
  $\alpha\in\zp^n$ with $|\alpha|\leq d$,
  $$\int_{\rrn}m(x)x^\alpha\,dx=0.$$
\end{enumerate}
\end{definition}

Then, combining Remark \ref{remark5.0.9}(ii)
and Theorem \ref{llmolel}, we immediately
conclude the following molecular
characterization of $h\MorrSo$;
we omit the details.

\begin{corollary}\label{llmolelm}
Let $p,\ q\in[1,\infty)$,
$\omega\in M(\rp)$ with
$$
-\frac{n}{p}<\m0(\omega)\leq\M0(\omega)<0
$$
and
$$
-\frac{n}{p}<\mi(\omega)\leq\MI(\omega)<0,
$$
$s\in(0,1)$, $d\geq\lfloor n(1/s-1)\rfloor$
be a fixed integer,
$$r\in\left(\frac{n}{\mw+n/p},\infty\right],$$
and $\tau\in(0,\infty)$ with $\tau>
n(1/s-1/r)$.
Then $f\in h\MorrSo$ if and only if $f\in\mathcal{S}'(\rrn)$
and there exist a sequence $\{m_{j}\}_{j\in\mathbb{N}}$
of local-$(\MorrSo,\,r,\,d,\,\tau)$-molecules
centered, respectively,
at cubes $\{B_{j}\}_{j\in\mathbb{N}}\subset\mathbb{B}$ and
a sequence $\{\lambda_{j}\}_{j\in\mathbb{N}}\subset
[0,\infty)$ such that
\begin{equation*}
f=\sum_{j\in\mathbb{N}}\lambda_{j}a_{j}
\end{equation*}
in $\mathcal{S}'(\rrn)$, and
$$
\left\|\left\{\sum_{j\in\mathbb{N}}\left[\frac{\lambda_{j}}
{\|\1bf_{B_{j}}\|_{\MorrSo}}\right]^s\1bf_
{B_{j}}\right\}^{\frac{1}{s}}\right
\|_{\MorrSo}<\infty.
$$
Moreover, there exists a
constant $C\in[1,\infty)$ such that, for
any $f\in h\MorrSo$,
\begin{align*}
C^{-1}\|f\|_{h\MorrSo}&\leq\inf\left\{
\left\|\left\{\sum_{j\in\mathbb{N}}\left[\frac{\lambda_{i}}
{\|\1bf_{B_{j}}\|_{\MorrSo}}\right]^s\1bf_
{B_{j}}\right\}^{\frac{1}{s}}\right\|_{\MorrSo}\right\}\\
&\leq C\|f\|_{h\MorrSo},
\end{align*}
where the infimum is taken
over all the decompositions of $f$ as above.
\end{corollary}

Next, we are devoted to establishing
the molecular characterization
of the local generalized
Herz--Hardy space $\HaSaHl$.
For this purpose, we
first introduce the definition
of local-$(\HerzS,\,r,\,d,\,\tau)$-molecules
as follows.

\begin{definition}
Let $p$, $q\in(0,\infty)$, $\omega\in M(\rp)$
with $\m0(\omega)\in(-\frac{n}{p},\infty)$ and
$$
-\frac{n}{p}<\mi(\omega)\leq
\MI(\omega)<0,
$$
$r\in[1,\infty]$,
$d\in\mathbb{Z}_{+}$, and $\tau\in
(0,\infty)$. Then a measurable
function $m$ on $\rrn$ is called a
\emph{local-$(\HerzS,\,r,\,d,\,\tau)$-molecule}
\index{local-$(\HerzS,\,r,\,d,
\,\tau)$-molecule}centered at a ball
$B(x_0,r_0)\in\mathbb{B}$, with $x_0\in
\rrn$ and $r_0\in(0,\infty)$, if
\begin{enumerate}
  \item[{\rm(i)}] for any $i\in\mathbb{Z}_{+}$,
$$\left\|m\1bf_{S_{i}(B(x_0,r_0))}
\right\|_{L^{r}(\rrn)}\leq2^{-\tau i}
\frac{|B(x_0,r_0)|^{1/r}}{\|\1bf_{B(x_0,r_0)}
\|_{\HerzS}};$$
  \item[{\rm(ii)}] when $r_0\in(0,1)$, then, for any
  $\alpha\in\zp^n$ with $|\alpha|\leq d$,
  $$\int_{\rrn}m(x)x^\alpha\,dx=0.$$
\end{enumerate}
\end{definition}

Via local-$(\HerzS,\,r,\,d,\,\tau)$-molecules,
we have the following molecular
characterization of the local generalized
Herz--Hardy space $\HaSaHl$.

\begin{theorem}\label{llmoleg}
Let $p,\ q\in(0,\infty)$,
$\omega\in M(\rp)$ with $\m0(\omega)
\in(-\frac{n}{p},\infty)$
and
$$-\frac{n}{p}
<\mi(\omega)\leq\MI(\omega)<0,$$
$$
s\in\left(0,\min\left\{1,p,q,\frac{n}
{\Mw+n/p}\right\}\right),
$$
$d\geq\lfloor n(1/s-1)\rfloor$ be a fixed integer,
$$r\in\left(\max\left\{1,p,\frac{n}{
\mw+n/p}\right\},\infty\right],$$
and $\tau\in(0,\infty)$ with $\tau>
n(1/s-1/r)$.
Then $f\in\HaSaHl$ if and only if $f\in\mathcal{S}'(\rrn)$
and there exist a sequence $\{m_{j}\}_{j\in\mathbb{N}}$
of local-$(\HerzS,\,r,\,d,
\,\tau)$-molecules centered, respectively,
at the balls
$\{B_{j}\}_{j\in\mathbb{N}}\subset\mathbb{B}$ and
a sequence $\{\lambda_{j}\}_{j\in\mathbb{N}}\subset
[0,\infty)$ such that
\begin{equation*}
f=\sum_{j\in\mathbb{N}}\lambda_{j}a_{j}
\end{equation*}
in $\mathcal{S}'(\rrn)$, and
$$
\left\|\left\{\sum_{j\in\mathbb{N}}\left[\frac{\lambda_{j}}
{\|\1bf_{B_{j}}\|_{\HerzS}}\right]^s\1bf_
{B_{j}}\right\}^{\frac{1}{s}}\right
\|_{\HerzS}<\infty.
$$
Moreover, there exists a
constant $C\in[1,\infty)$ such that, for
any $f\in\HaSaHl$,
\begin{align*}
C^{-1}\|f\|_{\HaSaHl}&\leq\inf\left\{
\left\|\left\{\sum_{j\in\mathbb{N}}\left[\frac{\lambda_{i}}
{\|\1bf_{B_{j}}\|_{\HerzS}}\right]^s\1bf_
{B_{j}}\right\}^{\frac{1}{s}}\right\|_{\HerzS}\right\}\\
&\leq C\|f\|_{\HaSaHl},
\end{align*}
where the infimum is taken
over all the decompositions of $f$ as above.
\end{theorem}

To establish the above molecular characterization,
we first show the following improved
molecular characterization
of localized Hardy spaces associated with
ball quasi-Banach function spaces.

\begin{theorem}\label{molexl}
Let $X$ be a ball quasi-Banach function space
satisfy:
\begin{enumerate}
  \item[{\rm(i)}] Assumption \ref{assfs}
  holds true with $0<\theta<s\leq1$;
  \item[{\rm(ii)}] for the above $s$,
  $X^{1/s}$ is a ball Banach function space and
  there exists a linear space
  $Y\subset\Msc(\rrn)$
  equipped with a seminorm $\|\cdot\|_{Y}$
  such that, for any $f\in\Msc(\rrn)$,
  \begin{equation*}
  \|f\|_{X^{1/s}}
  \sim\sup\left\{\|fg\|_{L^1(\rrn)}:\
  \|g\|_{Y}=1\right\},
  \end{equation*}
  where the positive equivalence
  constants are independent of $f$;
  \item[{\rm(iii)}] for the above
  $s$ and $Y$, there exist an $r\in(1,\infty]$
  and a positive constant $C$ such that,
  for any $f\in L^1_{\mathrm{loc}}(\rrn)$,
  $$
  \left\|\mc^{((r/s)')}
  (f)\right\|_{Y}\leq C\|f\|_{Y}.
  $$
\end{enumerate}
Then $f\in h_X(\rrn)$ if and only
if $f\in\mathcal{S}'(\rrn)$
and there exist a sequence
$\{m_{j}\}_{j\in\mathbb{N}}$
of local-$(X,\,r,\,d,\,\tau)$-molecules
centered, respectively,
at the balls $\{B_{j}\}_{j\in\mathbb{N}}
\subset\mathbb{B}$ and
a sequence $\{\lambda_{j}\}
_{j\in\mathbb{N}}\subset
[0,\infty)$ such that
\begin{equation}\label{molexle0}
f=\sum_{j\in\mathbb{N}}\lambda_{j}a_{j}
\end{equation}
in $\mathcal{S}'(\rrn)$, and
\begin{equation}\label{molexle00}
\left\|\left[\sum_{j\in\mathbb{N}}
\left(\frac{\lambda_{j}}
{\|\1bf_{B_{j}}\|_{X}}
\right)^{s}\1bf_{B_{j}}
\right]^{\frac{1}{s}}
\right\|_{X}<\infty.
\end{equation}
Moreover, there exists a
positive constant $C$ such that,
for any $f\in h_X(\rrn)$,
\begin{align}\label{molexle01}
C^{-1}\|f\|_{h_X(\rrn)}
&\leq\inf\left\{
\left\|\left[\sum_{j\in\mathbb{N}}
\left(\frac{\lambda_{i}}
{\|\1bf_{B_{j}}\|_{X}}\right)^s\1bf_
{B_{j}}\right]^{\frac{1}{s}}\right\|_{X}
\right\}\notag\\
&\leq C\left\|f\right\|_{h_X(\rrn)},
\end{align}
where the infimum is taken
over all the decompositions of $f$ as above.
\end{theorem}

\begin{proof}
Let all the symbols
be as in the present theorem.
We first prove the necessity.
To achieve this, let $f\in h_X(\rrn)$.
Then, from the assumptions
(i) through (iii) of the present theorem,
and Theorem \ref{atomxxl}, we deduce that
$f\in h^{X,r,d,s}(\rrn)$.
Thus, by Definition \ref{atom-hardyl},
we conclude that there exist
a sequence $\{\lambda_j\}_{j\in\mathbb{N}}
\subset[0,\infty)$ and a sequence
$\{a_j\}_{j\in\mathbb{N}}$ of
local-$(X,\,r,\,d)$-atoms supported,
respectively, in the balls
$\{B_j\}_{j\in\mathbb{N}}$ such that
\begin{equation}\label{molexle1}
f=\sum_{j\in\mathbb{N}}\lambda_{j}a_{j}
\end{equation}
in $\mathcal{S}'(\rrn)$, and
\begin{equation}\label{molexle2}
\left\|\left[\sum_{j\in\mathbb{N}}
\left(\frac{\lambda_{j}}
{\|\1bf_{B_{j}}\|_{X}}
\right)^{s}\1bf_{B_{j}}\right]^{\frac{1}{s}}
\right\|_{X}<\infty.
\end{equation}
Using Remark \ref{lomolexr}(i),
we find that, for any $j\in\mathbb{N}$,
$a_j$ is a local-$(X,\,r,\,r,\,\tau)$-molecule
centered at $B_j$.
This, combined with both \eqref{molexle1}
and \eqref{molexle2}, implies
that the necessity holds true.
Furthermore, by the choice of $\{\lambda_j\}_{
j\in\mathbb{N}}$, Definition \ref{atom-hardyl},
and Theorem \ref{atomxxl} again,
we conclude that
\begin{align}\label{molexle3}
\inf\left\{
\left\|\left[\sum_{j\in\mathbb{N}}
\left(\frac{\lambda_{i}}
{\|\1bf_{B_{j}}\|_{X}}\right)^s\1bf_
{B_{j}}\right]^{\frac{1}{s}}\right\|_{X}
\right\}\leq\left\|f\right\|_{h^{X,r,d,s}(\rrn)}
\sim\left\|f\right\|_{h_X(\rrn)},
\end{align}
where the infimum is taken
over all the sequences $\{\lambda_j\}_{
j\in\mathbb{N}}\subset[0,\infty)$
and $\{m_j\}_{j\in\mathbb{N}}$
of local-$(X,\,r,\,d,\,\tau)$-molecules
centered, respectively, at
the balls $\{B_j\}_{j\in\mathbb{N}}
\subset\mathbb{B}$
such that both \eqref{molexle0}
and \eqref{molexle00} hold true.

Conversely, we show the sufficiency.
Let $f\in\mathcal{S}'(\rrn)$,
$\{\lambda_{l,j}\}_{l\in\{1,2\},\,j
\in\mathbb{N}}\subset[0,\infty)$,
and $\{m_{l,j}\}_{l\in\{1,2\},\,j
\in\mathbb{N}}$ be a sequence of
local-$(X,\,r,\,d,\,\tau)$-molecules
centered, respectively, at the
balls $\{B_{l,j}\}_{l\in\{1,2\},\,j
\in\mathbb{N}}\in\mathbb{B}$ such that,
for any $j\in\mathbb{N}$,
$r(B_{1,j})\in(0,1)$
and $r(B_{2,j})\in[1,\infty)$,
\begin{equation}\label{molexle4}
f=\sum_{l=1}^{2}\sum_{j\in\mathbb{N}}
\lambda_{l,j}m_{l,j}
\end{equation}
in $\mathrm{S}'(\rrn)$, and
\begin{equation}\label{molexle5}
\left\|\left[\sum_{l=1}^{2}\sum_{j\in\mathbb{N}}
\left(\frac{\lambda_{l,j}}
{\|\1bf_{B_{l,j}}\|_{X}}
\right)^{s}\1bf_{B_{l,j}}\right]^{\frac{1}{s}}
\right\|_{X}<\infty.
\end{equation}
In what follows,
fix a $\phi\in\mathcal{S}(\rrn)$
satisfying that $\supp(\phi)\subset
B(\0bf,1)$ and $$\int_{\rrn}\phi(x)\,dx
\neq0.$$ Then, applying \eqref{molexle4}
and an argument
similar to that used in the
estimation of \eqref{Atogxe10}
with $\{\lambda_j\}_{j\in\mathbb{N}}$
and $\{a_j\}_{j\in\mathbb{N}}$
therein replaced, respectively, by
$\{\lambda_{l,j}\}_{l\in\{1,2\},\,j\in\mathbb{N}}$
and $\{m_{l,j}\}_{l\in\{1,2\},\,j\in\mathbb{N}}$,
we find that, for any $t\in(0,\infty)$,
\begin{align*}
\left|f\ast\phi_t\right|
\leq\sum_{l=1}^{2}\sum_{j\in\mathbb{N}}
\lambda_{l,j}\left|m_{l,j}\ast\phi_t\right|.
\end{align*}
Combining this and
Definition \ref{df511}(i), we further obtain
\begin{equation*}
m(f,\phi)\leq\sum_{l=1}^{2}\sum_{j\in\mathbb{N}}
\lambda_{l,j}m(m_{l,j},\phi).
\end{equation*}
This, together with Definition
\ref{Df1}(ii), implies that
\begin{align}\label{molexle6}
\left\|m(f,\phi)\right\|_X
&\lesssim\left\|\sum_{j\in\mathbb{N}}
\lambda_{1,j}m(m_{1,j},\phi)\right\|_X
+\left\|\sum_{j\in\mathbb{N}}
\lambda_{2,j}m(m_{2,j},\phi)\right\|_X\notag\\
&=:\mathrm{V}_1+\mathrm{V}_2.
\end{align}
We next estimate
$\mathrm{V}_1$ and $\mathrm{V}_2$,
respectively.

First, we deal with $\mathrm{V}_1$.
Indeed, for any $j\in\mathbb{N}$,
by Remark \ref{lomolexr}(ii) with
$m:=m_{1,j}$, we conclude that
$m_{1,j}$ is an $(X,\,r,\,d,\,\tau)$-molecule
centered at $B_{1,j}$.
Therefore, from Definitions \ref{df511}(i),
\ref{smax}(i), and \eqref{Df1}(ii),
and some arguments similar to those used in
the estimations of \eqref{Molegxe5},
\eqref{Molegxe6}, \eqref{Molegxe9},
and \eqref{Molegxe10}, it follows that
\begin{align}\label{molexle7}
\mathrm{V}_1\lesssim\left\|
\sum_{j\in\mathbb{N}}\lambda_j
M(m_{1,j},\phi)\right\|_X\lesssim\left\|\left[
\sum_{j\in\mathbb{N}}\left(\frac{\lambda_{1,j}}{
\|\1bf_{B_{1,j}}\|_X}
\right)^s\1bf_{B_{1,j}}\right]^{1/s}\right\|_X,
\end{align}
which is the desired estimate of $\mathrm{V}_1$.

Conversely, we next estimate
$\mathrm{V}_2$. To this end,
for any $j\in\mathbb{N}$ and $k
\in\zp$, let
$$
\mu_{j,k}:=2^{-k(\tau+\frac{n}{r})}
\frac{\|\1bf_{\tk B_{2,j}}\|_X}{\|
\1bf_{B_{2,j}}\|_X}
$$
and
$$
a_{j,k}:=2^{k(\tau+\frac{n}{r})}
\frac{\|\1bf_{B_{2,j}}\|_X}{
\|\1bf_{\tk B_{2,j}\|_X}}m_{2,j}
\1bf_{S_k(B_{2,j})}.
$$
Then, for any $j\in\mathbb{N}$,
we have
\begin{equation}\label{molexle8}
m_{2,j}=\sum_{k\in\zp}m_{2,j}
\1bf_{S_k(B_{2,j})}
=\sum_{k\in\zp}\mu_{j,k}
a_{j,k}
\end{equation}
almost everywhere in $\rrn$.
In addition, by the Tonelli theorem,
the H\"{o}lder inequality, and
Definition \ref{lomolex}(i)
with $m$ therein replaced by
$m_{2,j}$ with $j\in\mathbb{N}$, we
conclude that, for any $t\in(0,\infty)$
and $x\in\rrn$,
\begin{align*}
&\sum_{k\in\zp}
\int_{\rrn}\mu_{j,k}\left|
a_{j,k}(y)\right|\left|
\phi_t(x-y)\right|\,dy\\
&\quad=\sum_{k\in\zp}\int_{S_k(B_{2,j})}
\left|m_{2,j}(y)\right|\left|
\phi_t(x-y)\right|\,dy\\
&\quad\leq\sum_{k\in\zp}
\left\|m_{2,j}\1bf_{
S_k(B_{2,j})}\right\|_{L^r(\rrn)}
\left\|\phi_t\right\|_{L^{r'}(\rrn)}\\
&\quad\leq\frac{|B_{2,j}|^{1/r}}{\|\1bf_{B_{2,j}}\|
_{X}}\left\|\phi_t\right\|_{L^{r'}(\rrn)}
\sum_{k\in\zp}2^{-k\tau}\sim1.
\end{align*}
From this, \eqref{molexle8},
and the Fubini theorem, we deduce that,
for any $j\in\mathbb{N}$, $t\in(0,\infty)$,
and $x\in\rrn$,
\begin{align*}
m_{2,j}\ast\phi_t(x)
&=\int_{\rrn}\sum_{k\in\zp}
\mu_{j,k}a_{j,k}(y)\phi_{
t}(x-y)\,dy\\
&=\sum_{k\in\zp}\mu_{j,k}
\int_{\rrn}a_{j,k}\phi_{j,k}(x-y)\,dy
\\&=\sum_{k\in\zp}\mu_{j,k}a_{j,k}\ast
\phi_t(x).
\end{align*}
This, together with Definition \ref{df511},
further implies that, for any
$j\in\mathbb{N}$,
$$
m(m_{2,j},\phi)
\leq\sum_{k\in\zp}
\mu_{j,k}m(a_{j,k},\phi).
$$
By this and Definition \ref{Df1}(ii),
we conclude that
\begin{equation}\label{molexle10}
\mathrm{V}_2\lesssim
\left\|\sum_{j\in\mathbb{N}}
\sum_{k\in\zp}\lambda_{2,j}\mu_{j,k}
m(a_{j,k},\phi)\right\|_{X}.
\end{equation}

We now claim that, for any $j\in\mathbb{N}$
and $k\in\zp$, $a_{j,k}$ is
a local-$(X,\,r,\,d)$-atom
supported in $\tk B_{2,j}$.
Indeed, applying
Definition \ref{lomolex}(i),
we find that, for any $j\in\mathbb{N}$
and $k\in\zp$,
\begin{equation}\label{molexle9}
\left\|a_{j,k}\right\|_{L^r(\rrn)}
\leq2^{\frac{nk}{r}}\frac{|B_{2,j}|^{1/r}}
{\|\1bf_{\tk B_{2,j}}\|_X}=\frac{
|\tk B_{2,j}|^{1/r}}{\|\1bf_{\tk B_{2,j}}}.
\end{equation}
On the other hand, for any $j\in\mathbb{N}$
and $k\in\zp$, observe that
$\supp(a_{j,k})
\subset 2^k B_{2,j}$ and
$$r(\tk B_{2,j})=\tk r(B_{2,j})\geq\tk\geq1,$$
which, combined with \eqref{molexle9}
and Definition \ref{loatomx}, imply
that $a_{j,k}$ is a
local-$(X,\,r,\,d)$-atom supported
in $\tk B_{2,j}$. This
finishes the proof of the
above claim.
Combining this claim,
the fact that, for
any $j\in\mathbb{N}$ and
$k\in\zp$, $r(\tk B_{2,j})\geq1$,
the claim obtained in the proof
of Theorem \ref{atomxxl}, \eqref{molexle10},
and Definitions \ref{df511}(i) and \ref{smax}(i),
we conclude that
\begin{align*}
\mathrm{V}_2&\lesssim
\left\|\sum_{j\in\mathbb{N}}
\sum_{k\in\zp}\lambda_{2,j}\mu_{j,k}
m(a_{j,k},\phi)\1bf_{\tka B_{2,j}}\right\|_{X}\\
&\lesssim\left\|\sum_{j\in\mathbb{N}}
\sum_{k\in\zp}\lambda_{2,j}\mu_{j,k}
M(a_{j,k},\phi)\1bf_{\tka B_{2,j}}\right\|_{X}.
\end{align*}
Using this and an argument similar
to that used in the estimation
of $\mathrm{II}_1$ in the proof of Theorem
\ref{Atogx} with $\{\lambda_j\}_{j\in\mathbb{N}}$
and $\{a_{j}\}_{j\in\mathbb{N}}$
therein replaced, respectively, by
$\{\lambda_{2,j}
\mu_{j,k}\}_{j\in\mathbb{N},\,k\in\zp}$
and $\{a_{j,k}\}_{j\in\mathbb{N},\,k\in\zp}$,
we find that
\begin{align*}
\mathrm{V}_2&\lesssim
\left\|\left[
\sum_{j\in\mathbb{N}}
\sum_{k\in\zp}
\left(\frac{\lambda_{2,j}\mu_{j,k}}{
\|\1bf_{\tk B_{2,j}}\|_X}\right)^s
\1bf_{\tk B_{2,j}}
\right]^{\frac{1}{s}}\right\|_X\\
&\sim\left\|\left[\sum_{j\in\mathbb{N}}
\left(\frac{\lambda_{2,j}}{
\|\1bf_{B_{2,j}}\|_X}\right)^s
\sum_{k\in\zp}2^{-ks(\tau+
\frac{n}{r})}\1bf_{\tk B_{2,j}}
\right]^{\frac{1}{s}}\right\|_X.
\end{align*}
This, together with an argument similar
to that used in the estimation
of \eqref{Molegxe9} with $\{\lambda_j\}_{j\in\mathbb{N}}$
and $\{B_{j}\}_{j\in\mathbb{N}}$
therein replaced, respectively, by
$\{\lambda_{2,j}\}_{j\in\mathbb{N}}$
and $\{B_{2,j}\}_{j\in\mathbb{N}}$, further implies
that
\begin{equation}\label{molexle11}
\mathrm{V}_2\lesssim\left\|\left[
\sum_{j\in\mathbb{N}}\left(\frac{\lambda_{2,j}}{
\|\1bf_{B_{2,j}}\|_X}
\right)^s\1bf_{B_{1,j}}\right]^{1/s}\right\|_X,
\end{equation}
which is the desired estimate of $\mathrm{V}_2$.
Then, from both the assumptions (i) and (ii)
of the present theorem,
Remark \ref{Th6.1r}, Lemma
\ref{Th6.1l1}(ii), \eqref{molexle6}, \eqref{molexle7},
\eqref{molexle11}, and \eqref{molexle5},
it follows that
\begin{equation}\label{molexle12}
\left\|f\right\|_{h_X(\rrn)}
\sim\left\|m(f,\phi)\right\|_{X}
\lesssim\left\|\left[
\sum_{l=1}^{2}\sum_{j\in\mathbb{N}}
\left(\frac{\lambda_{l,j}}
{\|\1bf_{B_{l,j}}\|_{X}}
\right)^{s}\1bf_{B_{l,j}}\right]^{\frac{1}{s}}
\right\|_{X}<\infty.
\end{equation}
This further implies that $f\in h_X(\rrn)$,
and hence finishes the proof of the sufficiency.
Moreover, by \eqref{molexle12} and
the choice of $\{\lambda_{l,j}\}_{
l\in\{1,2\},\,j\in\mathbb{N}}$, we conclude that
$$
\left\|f\right\|_{h_X(\rrn)}\lesssim
\inf\left\{
\left\|\left[\sum_{j\in\mathbb{N}}
\left(\frac{\lambda_{i}}
{\|\1bf_{B_{j}}\|_{X}}\right)^s\1bf_
{B_{j}}\right]^{\frac{1}{s}}\right\|_{X}
\right\},
$$
where the infimum is taken
over all the sequences $\{\lambda_j\}_{
j\in\mathbb{N}}\subset[0,\infty)$
and $\{m_j\}_{j\in\mathbb{N}}$
of local-$(X,\,r,\,d,\,\tau)$-molecules
centered, respectively, at
the balls $\{B_j\}_{j\in\mathbb{N}}
\subset\mathbb{B}$
such that both
\eqref{molexle0} and \eqref{molexle00}
hold true.
From this and \eqref{molexle3}, we deduce that
\eqref{molexle01} holds true,
which completes the proof of Theorem
\ref{molexl}.
\end{proof}

\begin{remark}
We should point out that Theorem
\ref{molexl} is an improved version
of the known molecular
characterization of $h_X(\rrn)$ established by
Wang et al.\ in \cite[Theorem 5.2]{WYY}.
Indeed, if $Y\equiv(X^{1/s})'$ in
Theorem \ref{molexl}, then this theorem
goes back to \cite[Theorem 5.2]{WYY}.
\end{remark}

With the help of the above improved
molecular characterization of $h_X(\rrn)$,
we now prove Theorem \ref{llmoleg}.

\begin{proof}[Proof of Theorem \ref{llmoleg}]
Let $p$, $q$, $\omega$, $r$,
$d$, and $s$ be as in the present
theorem. Then, using the
assumptions $\m0(\omega)\in
(-\frac{n}{p},\infty)$ and
$\MI(\omega)\in(-\infty,0)$, and Theorem \ref{Th2},
we find that the global generalized
Herz space $\HerzS$ under consideration
is a BQBF space.
Therefore, to complete
the proof of the present theorem, we only need to
show that $\HerzS$ satisfies
(i) through (iii) of Theorem \ref{molexl}.
Indeed, let $\theta\in(0,s)$ satisfy that
\begin{equation}\label{llmolege1}
\left\lfloor n
\left(\frac{1}{s}-1\right)\right\rfloor
\leq n\left(\frac{1}{\theta}-1\right)
<\left\lfloor n
\left(\frac{1}{s}-1\right)\right\rfloor+1
\end{equation}
and
\begin{equation}\label{llmolege2}
\tau>n\left(\frac{1}{\theta}-\frac{1}{r}
\right).
\end{equation}
For the above $\theta$, $s$, and
$r$, by the proof of
Theorem \ref{lMoleg}, we conclude
that the following
three statements hold true:
\begin{enumerate}
  \item[{\rm(i)}] for any $\{f_{j}\}_{
j\in\mathbb{N}}\subset L^{1}_{\text{loc}}(\rrn)$,
\begin{equation*}
\left\|\left\{\sum_{j\in\mathbb{N}}\left[\mc^{(\theta)}
(f_{j})\right]^s\right\}^{1/s}
\right\|_{\HerzS}\lesssim
\left\|\left(\sum_{j\in\mathbb{N}}
|f_{j}|^{s}\right)^{1/s}\right\|_{\HerzS};
\end{equation*}
  \item[{\rm(ii)}] $[\HerzS]^{1/s}$ is a BBF space
and, for any $f\in\Msc(\rrn)$,
\begin{equation*}
\|f\|_{[\HerzS]^{1/s}}\sim
\sup\left\{\|fg\|_{L^1(\rrn)}:\
\|g\|_{\dot{\mathcal{B}}_{1/\omega^s}
^{(p/s)',(q/s)'}(\rrn)}=1\right\}
\end{equation*}
with the positive equivalence constants
independent of $f$;
  \item[{\rm(iii)}] for any
$f\in L_{{\rm loc}}^{1}(\rrn)$,
\begin{equation*}
\left\|\mc^{((r/s)')}(f)\right\|_{\dot{\mathcal{B}}
_{1/\omega^{s}}^{(p/s)',(q/s)'}(\rrn)}\lesssim
\|f\|_{\dot{\mathcal{B}}_
{1/\omega^{s}}^{(p/s)',(q/s)'}(\rrn)}.
\end{equation*}
\end{enumerate}
These further imply that
the Herz space $\HerzS$ under consideration
satisfies the assumptions
(i) through (iii) of Theorem
\ref{molexl}. In addition, from
\eqref{llmolege1} and
the assumption $d\geq\lfloor
n(1/s-1)\rfloor$, it follows that
$d\geq\lfloor n(1/s-1)\rfloor=
\lfloor n(1/\theta-1)\rfloor$.
This, together with \eqref{llmolege2} and
the fact that (i) through (iii)
of Theorem \ref{molexl} hold true, then
finishes the proof of Theorem \ref{llmoleg}.
\end{proof}

As an application of Theorem
\ref{llmoleg}, we next establish
the molecular characterization of
the local generalized Morrey--Hardy
space $h\MorrS$. To begin with,
we introduce the following concept of
local molecules associated with
the global generalized Morrey space $\MorrS$.

\begin{definition}
Let $p$, $q$, $\omega$,
$r$, $d$, and $\tau$ be as in
Definition \ref{lomolem}.
Then a measurable
function $m$ on $\rrn$ is called a
\emph{local-$(\MorrS,\,r,\,d,\,\tau)$-molecule}
\index{local-$(\MorrS,\,r,\,d,
\,\tau)$-molecule}centered at a ball
$B(x_0,r_0)\in\mathbb{B}$, with $x_0\in
\rrn$ and $r_0\in(0,\infty)$, if
\begin{enumerate}
  \item[{\rm(i)}] for any $i\in\mathbb{Z}_{+}$,
$$\left\|m\1bf_{S_{i}(B(x_0,r_0))}
\right\|_{L^{r}(\rrn)}\leq2^{-\tau i}
\frac{|B(x_0,r_0)|^{1/r}}{\|\1bf_{B(x_0,r_0)}
\|_{\MorrS}};$$
  \item[{\rm(ii)}] when $r_0\in(0,1)$, then, for any
  $\alpha\in\zp^n$ with $|\alpha|\leq d$,
  $$\int_{\rrn}m(x)x^\alpha\,dx=0.$$
\end{enumerate}
\end{definition}

Then, by Remark \ref{remark5.0.9}
and Theorem \ref{llmoleg}, we conclude
the molecular characterization
of the local generalized
Morrey--Hardy space $h\MorrS$ as follows;
we omit the details.

\begin{corollary}
Let $p$, $q$, $\omega$,
$r$, $d$, and $s$ be
as in Corollary \ref{llmolelm}.
Then $f\in h\MorrS$ if and only if $f\in\mathcal{S}'(\rrn)$
and there exist a sequence $\{m_{j}\}_{j\in\mathbb{N}}$
of local-$(\MorrS,\,r,\,d,\,\tau)$-molecules
centered, respectively,
at cubes $\{B_{j}\}_{j\in\mathbb{N}}\subset\mathbb{B}$ and
a sequence $\{\lambda_{j}\}_{j\in\mathbb{N}}\subset
[0,\infty)$ such that
$$
f=\sum_{j\in\mathbb{N}}\lambda_{j}a_{j}
$$
in $\mathcal{S}'(\rrn)$, and
$$
\left\|\left\{\sum_{j\in\mathbb{N}}\left[\frac{\lambda_{j}}
{\|\1bf_{B_{j}}\|_{\MorrS}}\right]^s\1bf_
{B_{j}}\right\}^{\frac{1}{s}}\right
\|_{\MorrS}<\infty.
$$
Moreover, there exists a
constant $C\in[1,\infty)$ such that, for
any $f\in h\MorrS$,
\begin{align*}
C^{-1}\|f\|_{h\MorrS}&\leq\inf\left\{
\left\|\left\{\sum_{j\in\mathbb{N}}\left[\frac{\lambda_{i}}
{\|\1bf_{B_{j}}\|_{\MorrS}}\right]^s\1bf_
{B_{j}}\right\}^{\frac{1}{s}}\right\|_{\MorrS}\right\}\\
&\leq C\|f\|_{h\MorrS},
\end{align*}
where the infimum is taken
over all the decompositions of $f$ as above.
\end{corollary}

\section{Littlewood--Paley Function Characterizations}

The main target of this section
is to prove the Littlewood--Paley
function characterizations of the local
generalized Herz--Hardy spaces $\HaSaHol$
and $\HaSaHl$. To achieve this, we first establish
various Littlewood--Paley function
characterizations of the
local Hardy space $h_X(\rrn)$
with $X$ being a ball quasi-Banach
function space.
To begin with, we present the following
concepts of the local Littlewood--Paley
functions.

\begin{definition}\label{df551}
Let $\varphi_0\in\mathcal{S}(\rrn)$
satisfy
$$
\1bf_{B(\0bf,1)}\leq\widehat{\varphi}_0
\leq\1bf_{B(\0bf,2)},
$$
and $\varphi\in\mathcal{S}(\rrn)$ satisfy
$$
\1bf_{B(\0bf,4)\setminus
B(\0bf,2)}\leq\widehat{\varphi}
\leq\1bf_{B(\0bf,8)\setminus B(\0bf,1)}.
$$
Then, for any $f\in\mathcal{S}'(\rrn)$,
the \emph{local Lusin area
function}\index{local Lusin area function}
$S_{\mathrm{loc}}(f)$\index{$S_{\mathrm{loc}}(f)$},
the \emph{local Littlewood--Paley
$g$-function}\index{local Littlewood--Paley $g$-function}
$g_{\mathrm{loc}}(f)$\index{$g_{\mathrm{loc}}(f)$}, and the
\emph{local Littlewood--Paley
$g^{*}_{\lambda}$-function}\index{local Littlewood--Paley $g^{*}_{\lambda}$-function}
$(g^{*}_{\lambda})_{\mathrm{loc}}
(f)$\index{$(g^{*}_{\lambda})_{\mathrm{loc}}$} with
$\lambda\in(0,\infty)$ are defined,
respectively, by setting, for any $x\in\rrn$,
$$S_{\mathrm{loc}}(f)(x):=\left|f\ast\varphi_0(x)\right|
+\left[\int_{0}^{1}\int_{B(x,t)}|f\ast\varphi_t(y)|^{2}
\,\frac{dy\,dt}{t^{n+1}}\right]^{\frac{1}{2}},$$
$$
g_{\mathrm{loc}}(f)(x):=\left|f\ast\varphi_0(x)\right|+
\left[\int_{0}^{1}|f\ast\varphi_{t}(
x)|^{2}\,\frac{dt}{t}\right]^{\frac{1}{2}},
$$
and
\begin{align*}
&\left(g_{\lambda}^{*}\right)
_{\mathrm{loc}}(f)(x)\\
&\quad:=\left|f\ast\varphi_0(x)\right|
+\left[\int_{0}
^{1}\int_{\rrn}
\left(\frac{t}{t+|x-y|}\right)^{\lambda n}|
f\ast\varphi_{t}(y)|^{2}
\,\frac{dy\,dt}{t^{n+1}}\right]^{\frac{1}{2}}.
\end{align*}
\end{definition}

We now establish the following
Lusin area function characterization
of the local Hardy space $h_X(\rrn)$.

\begin{theorem}\label{lplx}
Let $X$ be a ball quasi-Banach
function space satisfying both
Assumptions \ref{assfs} and \ref{assas}
with the same $s$. Then
$f\in h_X(\rrn)$ if and only if
$f\in\mathcal{S}'(\rrn)$ and
$S_{\mathrm{loc}}(f)\in X$.
Moreover, there exists a constant $C\in[1,\infty)$
such that, for any $f\in h_X(\rrn)$,
$$
C^{-1}\left\|f\right\|_{h_X(\rrn)}
\leq\left\|S_{\mathrm{loc}}(f)\right\|_X
\leq C\left\|f\right\|_{h_X(\rrn)}.
$$
\end{theorem}

In order to show this theorem, we
require some auxiliary lemmas and concepts.
First, we need the following Calder\'{o}n reproducing
formula\index{Calder\'{o}n reproducing formula}
which was given in \cite[(5.3)]{WYY} (see also
\cite[Proposition 1.1.6]{LGMF}).

\begin{lemma}\label{icrff}
Let $\varphi_0$ and $\varphi$
be as in Theorem \ref{lplx}.
Then there exist $a,\ b,\ c\in(0,\infty)$ and $\psi_0$,
$\psi\in\mathcal{S}(\rrn)$ such that
$$
\supp(\widehat{\psi}_0)
\subset B(\0bf,a),$$
$$
\supp(\widehat{\psi})
\subset B(\0bf,c)\setminus B(\0bf,b),$$
and, for any $f\in\mathcal{S}'(\rrn)$,
$$
f=f\ast\varphi_0\ast\psi_0+
\int_{0}^{1}f\ast\varphi_t\ast
\psi_t\,\frac{dt}{t}
$$
in $\mathcal{S}'(\rrn)$,
namely,
$$
f=f\ast\varphi_0\ast\psi_0+\lim\limits_{
\eps\to0^+}
\int_{\eps}^{1}f\ast\varphi_t\ast
\psi_t\,\frac{dt}{t}
$$
in $\mathcal{S}'(\rrn)$.
\end{lemma}

We also need the following
auxiliary estimate about
convolutions, which is just
\cite[(5.9)]{WYY}.

\begin{lemma}\label{esti-con}
Let $f\in\mathcal{S}'(\rrn)$ and $
\varphi\in\mathcal{S}(\rrn)$. Then
there exist a positive integer $m$,
depending only on $f$, and a positive constant $C$,
independent of $f$, such that,
for any $t\in(0,1]$ and $x\in\rrn$,
$$
\left|f\ast\varphi_t(x)\right|\leq Ct^{-n-m}
\left(1+|x|\right)^m.
$$
\end{lemma}

Let $S$ be the Lusin area operator
as in Definition \ref{df451}. Then
the following $L^p$ boundedness of $S$
plays an important role in the proof of
Theorem \ref{lplx}, which was obtained in
\cite[Theorem 7.8]{fs82}.

\begin{lemma}\label{lusinlp}
Let $p\in(1,\infty)$ and
$S$ be as in Definition \ref{df451}.
Then there exists a positive constant
$C$ such that, for any
$f\in L^p(\rrn)$,
$$
\left\|S(f)\right\|_{L^p(\rrn)}\leq
C\left\|f\right\|_{L^p(\rrn)}.
$$
\end{lemma}

Moreover, the following conclusion
is useful in the proof of Theorem
\ref{lplx}, which
characterizes the Hardy space $H_X(\rrn)$
via the Lusin area function $S$ and
was showed in \cite[Theorem 3.21]{SHYY}.

\begin{lemma}\label{hardy-lusin}
Let $X$ be a ball quasi-Banach function
space satisfying both Assumptions \ref{assfs}
and \ref{assas} with the same $s\in(0,1]$,
$\varphi$ be as in Definition \ref{df551}, and
$S$ be as in Definition \ref{df451} with
the above $\varphi$.
Then $f\in H_X(\rrn)$ if and only if
$f\in\mathcal{S}'(\rrn)$, $f$ vanishes weakly
at infinity, and $S(f)\in X$.
Moreover, for any $f\in H_X(\rrn)$,
$$
\left\|f\right\|_{H_X(\rrn)}
\sim\left\|S(f)\right\|_{X}
$$
with the positive equivalence constants
independent of $f$.
\end{lemma}

In what follows, for any
$f\in\mathcal{S}'(\rrn)$ and $x\in\rrn$, let
\begin{equation}\label{lusinl}
\widetilde{S}_{\mathrm{loc}}(f)
(x):=\left[\int_{0}^{1}\int_{B(x,t)}
\left|f\ast\varphi_t(y)\right|^2\,\frac{dy\,dt}{
t^{n+1}}\right]^{\frac{1}{2}}.
\end{equation}
Then, in order to prove Theorem \ref{lplx},
we also need a technique
estimate about $\widetilde{S}_{\mathrm{loc}}$
as follows, which was obtained in
\cite[p.\,53]{WYY}.

\begin{lemma}\label{lusinl-atom}
Let $X$ be a ball quasi-Banach function
space, $d\in\zp$, and $\theta\in(0,\infty)$.
Then there exists a positive constant
$C$ such that, for any local-$(X,\,\infty,\,d)$-atom
$a$ supported in the ball $B\in\mathbb{B}$
with $r(B)\in[1,\infty)$,
$$
\widetilde{S}_{\mathrm{loc}}(a)\1bf_{
(4B)^{\complement}}\leq C\frac{1}{\|\1bf_{B}
\|_X}\mc^{(\theta)}\left(\1bf_B\right).
$$
\end{lemma}

Furthermore, to show Theorem \ref{lplx},
we require some conclusions about $X$-tent spaces.
We first recall some basic concepts.
Let $\alpha\in(0,\infty)$ and
$f:\ \mathbb{R}^{n+1}_+\to \mathbb{C}$ be
a measurable function.
Then the \emph{Lusin area function}
$\mathcal{A}_{\alpha}(f)$\index{$\mathcal{A}_{\alpha}$},
with aperture $\alpha$,
is defined by setting, for any $x\in\rrn$,
$$
\mathcal{A}^{(\alpha)}(f)(x):=
\left\{\int_{\Gamma_{\alpha}(x)}
\left|f\ast\varphi_t(y)
\right|^2\,\frac{dy\,dt}{t^{n+1}}
\right\}^{\frac{1}{2}},
$$
where, for any $x\in\rrn$ and
$\alpha\in(0,\infty)$, $\Gamma_{\alpha}(x)$
is defined as in \eqref{gammaa}.
For any given ball quasi-Banach
function space $X$, the \emph{$X$-tent space}
\index{$X$-tent space}$T_X^{\alpha}(\mathbb{R}^{n+1}_+
)$\index{$T_X^{\alpha}(\mathbb{R}^{n+1}_+)$},
with aperture $\alpha$,  is defined
to be the set of all the measurable
functions $f:\ \mathbb{R}^{n+1}_+\to\mathbb{C}$ such that
$$\left\|f\right\|_{
T_X^{\alpha}(\mathbb{R}^{n+1}_+)}:=
\left\|\mathcal{A}_{\alpha}(f)\right\|_X
<\infty.$$
In addition, let $\alpha\in(0,\infty)$.
Then, for any ball $B(x,r)\in\mathbb{B}$
with $x\in\rrn$ and $r\in(0,\infty)$, let
$$T_{\alpha}(B)\index{$T_{\alpha}(B)$}:=\left\{
(y,t)\in\mathbb{R}^{n+1}_+:\ 0<t<r/\alpha,\
|y-x|<r-\alpha t\right\}.$$
When $\alpha:=1$, we denote
$T_{\alpha}(B)$ simply by $T(B)$\index{$T(B)$}.
Then the following definition
of atoms is just \cite[Definition 3.17]{SHYY}.

\begin{definition}
Let $X$ be a ball quasi-Banach
function space, $p\in(1,\infty)$, and
$\alpha\in(0,\infty)$. A measurable
function $a:\ \mathbb{R}^{n+1}_+\to\mathbb{C}$
is called a \emph{$(T_X^{\alpha},\,p)$-atom}
\index{$(T_X^{\alpha},\,p)$-atom}if
there exists a ball $B\in\mathbb{B}$ such that
\begin{enumerate}
\item[{\rm{(i)}}] $\supp(a):=
\{(y,t)\in\rrnp:\ a(y,t)\neq0\}\subset T_{\alpha}(B)$;
\item[{\rm{(ii)}}] $\|\mathcal{A}_{\alpha}
(a)\|_{L^p(\rrn)}\leq\frac{|B|^{1/p}}{\|\mathbf{1}_B\|_X}$.
\end{enumerate}
Moreover, if $a$ is a $(T_X^\alpha,p)$-atom
for any $p\in(1,\infty)$,
then $a$ is called a
\emph{$(T_X^\alpha,\infty
)$-atom}\index{$(T_X^{\alpha},\,\infty)$-atom}.
\end{definition}

Repeating an argument similar to
that used in the proof of
\cite[Lemma 4.8]{hyy13} with
$L^{\varphi}(\rrn)$ therein replaced
by $X$, we obtain the following
auxiliary conclusion about
$(T_{X}^1,\,\infty)$-atoms, which
plays a key role in the proof of
the Lusin area function characterization
of $h_X(\rrn)$; we omit the details.

\begin{lemma}\label{tent-con}
Let $d\in\zp$,
$\psi\in\mathcal{S}(\rrn)$
satisfy that, for any $\gamma\in\zp^n$
with $|\gamma|\leq d$,
$
\int_{\rrn}\psi(x)x^{\gamma}\,dx=0,
$
and $X$ be a ball quasi-Banach function space.
Assume that $a$ is a $(T_{X}^1,\,\infty)$-atom
supported in $T(B)$ with $B\in\mathbb{B}$.
Then, for any $\tau\in(0,\infty)$,
$\int_{0}^{\infty} a(\cdot,t)\ast\psi_t\,\frac{dt}{t}$ is
a harmless constant multiple of
an $(X,\,\infty,\,d,\,\tau)$-molecule
centered at the ball $B$.
\end{lemma}

Moreover, the following atomic characterization
of $X$-tent spaces obtained in
\cite[Theorem 3.19]{SHYY} is also an essential
tool in the proof of Theorem \ref{lplx}.

\begin{lemma}\label{tent-atom}
Let $f:\ \mathbb{R}^{n+1}_+\to\mathbb{C}$
be a measurable function
and $X$ a ball quasi-Banach function space.
Assume that $X$ satisfies both
Assumptions \ref{assfs} and \ref{assas}
with the same $s\in(0,1]$. Then
$f\in T^1_X(\mathbb{R}^{n+1}_+)$ if and only if there exist
a sequence $\{\lambda_j\}_{j\in\mathbb{N}}\subset
[0,\infty)$ and a sequence $\{a_j\}_{j\in\mathbb{N}}$
of $(T_X^1,\infty)$-atoms supported, respectively,
in $\{T(B_j)\}_{j\in\mathbb{N}}$ with $\{B_j\}_{j
\in\mathbb{N}}\subset\mathbb{B}$
such that, for almost every
$(x,t)\in\mathbb{R}^{n+1}_+$,
$$f(x,t)=\sum_{j\in\mathbb{N}}\lambda_j a_j(x,t)
$$
and
$$
|f(x,t)|=\sum_{j\in\mathbb{N}}\lambda_j |a_j(x,t)|.$$
Moreover,
$$\|f\|_{T^1_X(\mathbb{R}^{n+1}_+)}
\sim\left\|\left[\sum_{j\in\mathbb{N}}
\left(\frac{\lambda_j}{\|\mathbf{1}_{B_j}\|_X}\right)^s
\mathbf{1}_{B_j}\right]^{1/s}\right\|_X,$$
where the positive equivalence
constants are independent of $f$.
\end{lemma}

Via these preparations, we now show Theorem
\ref{lplx}.

\begin{proof}[Proof of Theorem \ref{lplx}]
Let all the symbols be as in the
present theorem. We first prove the
sufficiency. To this end, let $f\in\mathcal{S}'
(\rrn)$ satisfy $S_{\mathrm{loc}}(f)\in X$.
Then, applying Lemma \ref{icrff},
we find that there exist
$a,\ b,\ c\in(0,\infty)$ and $\psi_0$,
$\psi\in\mathcal{S}(\rrn)$ such that
$
\supp(\widehat{\psi}_0)
\subset B(\0bf,a)
$, $\supp(\widehat{\psi})
\subset B(\0bf,c)\setminus B(\0bf,b)$,
and
\begin{equation}\label{lplxe1}
f=f\ast\varphi_0\ast\psi_0+
\int_{0}^{1}f\ast\varphi_t\ast
\psi_t\,\frac{dt}{t}
\end{equation}
in $\mathcal{S}'(\rrn)$.
From both Definitions \ref{df551} and
\ref{Df1}(ii), it follows that
\begin{align*}
\left\|f\ast\varphi_t\1bf_{\{\tau\in(0,1)\}}
\right\|_{T^1_X(\rrnp)}&=\left\|\left[\int_{\Gamma(\cdot)}
\left|f\ast\varphi_t(y)\1bf_{\{\tau:\,\tau\in(0,1)\}}(t)
\right|^2\,\frac{dy\,dt}{t^{n+1}}\right]
^{\frac{1}{2}}\right\|_X\\
&=\left\|\left[\int_0^1\int_{B(\cdot,t)}
\left|f\ast\varphi_t(y)\right|^2\,
\frac{dy\,dt}{t^{n+1}}\right]^{\frac{1}{2}}\right\|_X\\
&\leq\left\|S_{\mathrm{loc}}(f)\right\|_X<\infty.
\end{align*}
This, combined with Lemma \ref{tent-atom},
further implies that there exist
a sequence $\{\lambda_{1,j}\}_{j\in\mathbb{N}}
\subset[0,\infty)$ and a sequence
$\{a_{1,j}\}_{j\in\mathbb{N}}$ of $(T_X^1,\,\infty)$-atoms
supported, respectively, in $\{T(B_{1,j})\}_{
j\in\mathbb{N}}$ with $\{B_{1,j}\}_{j\in\mathbb{N}}
\subset\mathbb{B}$ such that, for almost
every $(x,t)\in\rrnp$
\begin{equation}\label{lplxe2}
f\ast\varphi_t(x)\1bf_{
\{\tau:\,\tau\in(0,1)\}}(t)=\sum_{j\in\mathbb{N}}
\lambda_{1,j}a_{1,j}(x,t),
\end{equation}
\begin{equation}\label{lplxe3}
\left|f\ast\varphi_t(x)\right|\1bf_{
\{\tau:\,\tau\in(0,1)\}}(t)=\sum_{j\in\mathbb{N}}
\lambda_{1,j}a_{1,j}\left|(x,t)\right|,
\end{equation}
and
\begin{align}\label{lplxe4}
\left\|\left[\sum_{j=1}^\infty
\left(\frac{\lambda_{1,j}}{
\|\mathbf{1}_{B_{1,j}}\|_X}\right)^s
\mathbf{1}_{B_{1,j}}\right]^{1/s}\right\|_X
&\sim\left\|f\ast\varphi_{t}
\1bf_{\{\tau\in(0,1)\}}
\right\|_{T^1_{X}(\rrnp)}\notag\\
&\lesssim\left\|
S_{\mathrm{loc}}(f)\right\|_{X}.
\end{align}

Next, we show that
\begin{equation}\label{lplxe5}
\int_0^1f\ast\varphi_t\ast
\psi_t\,\frac{dt}{t}=\sum_{j\in\mathbb{N}}
\lambda_{1,j}\int_0^1a_j(\cdot,t)\ast\psi_t
\,\frac{dt}{t}
\end{equation}
in $\mathcal{S}'(\rrn)$.
Indeed, by Lemma \ref{esti-con},
we find that there exists an $m\in\mathbb{N}$,
depending only on $f$, such that, for any
$t\in(0,1]$ and $x\in\rrn$,
\begin{equation}\label{lplxe6}
\left|f\ast\varphi_t(x)\right|
\lesssim t^{-n-m}\left(1+|x|\right)^m.
\end{equation}
This further implies that, for any
$t\in(0,1]$ and $x\in\rrn$,
\begin{align*}
&\left|f\ast\varphi_t\right|\ast
\left|\psi_t\right|(x)\\&\quad=
\int_{\rrn}\left|f\ast\varphi_t(y)\right|
\left|\psi_t(x-y)\right|\,dy\\
&\quad\lesssim t^{-n-m}\int_{\rrn}
\left(1+|y|\right)^m\frac{1}{t^n}
\left|\psi\left(\frac{x-y}{t}\right)\right|\,dy\\
&\quad\sim t^{-n-m}\int_{\rrn}\left(
1+\left|x-ty\right|\right)^m\left|
\psi(y)\right|\,dy\\
&\quad\lesssim t^{-n-m}\left(
1+|x|\right)^m\left[
\int_{\rrn}\left|\psi(y)\right|\,dy
+\int_{\rrn}\left(1+|y|\right)^m
\left|\psi(y)\right|\,dy\right]\\
&\quad\lesssim t^{-n-m}\left(
1+|x|\right)^m\int_{\rrn}\frac{1}{(1+|y|)^{n+1}}\,dy
\sim t^{-n-m}\left(
1+|x|\right)^m.
\end{align*}
From this and \eqref{lplxe3}, we
deduce that
\begin{align*}
&\int_{0}^{1}\int_{\rrn}\sum_{j\in\mathbb{N}}
\lambda_{1,j}\left|a_{1,j}\right|\ast\left|
\psi_t\right|(x)\left|\eta(x)\right|\,\frac{dx\,dt}{t}\\
&\quad=\int_{0}^{1}\int_{\rrn}\sum_{j\in\mathbb{N}}
\lambda_{1,j}\left|a_{1,j}(x,t)\right|
\left|\psi_t(-\cdot)\right|\ast
\left|\eta\right|(x)\,\frac{dx\,dt}{t}\\
&\quad=\int_{0}^{1}\int_{\rrn}\left|f\ast\varphi_t(x)
\right|\left|\psi_t(-\cdot)\right|\ast
\left|\eta\right|(x)\,\frac{dx\,dt}{t}\\
&\quad=\int_{0}^{1}\int_{\rrn}\left|f\ast\varphi_t
\right|\ast\left|\psi_t\right|(x)
\left|\eta(x)\right|\,\frac{dx\,dt}{t}\\
&\quad\lesssim\int_{0}^{1}\int_{\rrn}
t^{-n-m}\left(1+|x|\right)^m\left(1+|x|\right)^{-n-m-1}
\,\frac{dx\,dt}{t}\\
&\quad\sim\int_{0}^{1}\frac{1}{t^{n+m+1}}\,dt
\int_{\rrn}\frac{1}{(1+|x|)^{n+1}}\,dx<\infty.
\end{align*}
Combining this, \eqref{lplxe1},
\eqref{lplxe2}, and the Fubini theorem, we further
conclude that
\begin{align*}
&\int_{\rrn}\left[
\int_{0}^{1}f\ast\varphi_t\ast\psi_t(x)
\,\frac{dt}{t}\right]\eta(x)\,dx\\
&\quad=\int_{0}^{1}\int_{\rrn}
f\ast\varphi_t\ast\psi_t(x)\eta(x)\,\frac{dx\,dt}{t}\\
&\quad=\int_{0}^{1}\int_{\rrn}
f\ast\varphi_t(x)\psi_t(-\cdot)\ast\eta(x)
\,\frac{dx\,dt}{t}\\
&\quad=\int_{0}^{1}\int_{\rrn}
\sum_{j\in\mathbb{N}}\lambda_{1,j}
a_{1,j}(x,t)\psi_t(-\cdot)\ast\eta(x)
\,\frac{dx\,dt}{t}\\
&\quad=\int_{\rrn}\left[\sum_{j\in\mathbb{N}}
\lambda_{1,j}\int_{0}^{1}a_{1,j}\ast
\psi_t(x)\,\frac{dt}{t}\right]\eta(x)\,dx,
\end{align*}
which implies that \eqref{lplxe5} holds true.

In addition, observe that $\supp(
\widehat{\psi})\subset
B(\0bf,c)\setminus B(\0bf,b)$.
This implies that, for any $\gamma\in\zp^n$,
\begin{align*}
\int_{\rrn}\psi(x)x^{\gamma}\,dx
&=\left(-2\pi i\right)^{-|\gamma|}
\int_{\rrn}\psi(x)\left(2\pi ix\right)^{\gamma}\,dx\\
&=\left(-2\pi i\right)^{-|\gamma|}\mathcal{F}
\left(\psi\left[2\pi i\left(\cdot\right)
\right]^{\gamma}\right)(\0bf)\\
&=\left(-2\pi i\right)^{-|\gamma|}\partial^{\gamma}
\widehat{\psi}(\0bf)=0.
\end{align*}
In what follows, let
$d\geq\lfloor n(1/\theta-1)\rfloor$ be
a fixed integer and $\tau\in(n(\frac{1}{\theta}-
\frac{1}{r}),\infty)$ a fixed positive constant.
Then, from Lemma \ref{tent-con},
we deduce that, for any $j\in\mathbb{N}$,
$\int_{0}^{\infty}a_{1,j}(\cdot,t)\ast\psi_t
\,\frac{dt}{t}$ is a
harmless constant multiple of an
$(X,\,\infty,\,d,\,\tau)$-molecule centered
at $B_{1,j}$.
On the other hand, by \eqref{lplxe3},
we find that, for any $j\in\mathbb{N}$ and
$t\in[1,\infty)$ and for
almost every $x\in\rrn$,
$a_{1,j}(x,t)=0$. This further implies that,
for any $j\in\mathbb{N}$,
$$\int_0^1a_{1,j}(\cdot,t)\ast\psi_t(y)\,\frac{dt}{t}
=\int_0^{\infty}a_{1,j}(\cdot,t)\ast\psi_t(y)\,\frac{dt}{t}.$$
Therefore, for any
$j\in\mathbb{N}$,
$\int_{0}^{1}a_{1,j}(\cdot,t)\ast\psi_t
\,\frac{dt}{t}$ is a
harmless constant multiple of an
$(X,\,\infty,\,d,\,\tau)$-molecule centered
at $B_{1,j}$. Combining this,
\eqref{lplxe5}, Lemma \ref{llmolell1},
and \eqref{lplxe4}, we further conclude that
\begin{align}\label{lplxe7}
&\left\|\int_{0}^{1}f\ast\varphi_t\ast
\psi_t(\cdot)\,\frac{dt}{t}\right\|_{h_X(\rrn)}
\notag\\&\quad\lesssim\left\|\left[\sum_{j=1}^\infty
\left(\frac{\lambda_{1,j}}{
\|\mathbf{1}_{B_{1,j}}\|_X}\right)^s
\mathbf{1}_{B_{1,j}}\right]^{1/s}\right\|_X
\lesssim\left\|S_{\mathrm{loc}}(f)\right\|_{X}.
\end{align}

Next, we prove that
\begin{equation}\label{lplxe8}
\left\|f\ast\varphi_0\ast\psi_0
\right\|_{h_X(\rrn)}\lesssim\left\|
S_{\mathrm{loc}}(f)\right\|_{X}.
\end{equation}
To achieve this, for any $
\alpha\in\mathbb{Z}^n$, let $Q_{\alpha}:=
2\alpha+[0,2)^n$, $\mathbf{e}:=
(1,\ldots,1)\in\rrn$ denote the
\emph{unit} of $\rrn$, and
$B_{\alpha}:=B(\mathbf{e}+\alpha,\sqrt{n})$.
Take arrangements of $\{Q_{\alpha}\}_{
\alpha\in\mathbb{Z}^n}$ and $\{B_{\alpha}\}_{
\alpha\in\mathbb{Z}^n}$, which are denoted,
respectively, by $\{Q_{2,j}\}_{j\in\mathbb{N}}$
and $\{B_{2,j}\}_{j\in\mathbb{N}}$.
Then it holds true that, for any $j\in\mathbb{N}$,
$Q_{2,j}\subset B_{2,j}$, $\{Q_j\}_{j\in\mathbb{N}}$
are pairwise disjoint cubes,
$$\bigcup\limits_{j\in\mathbb{N}}Q_{2,j}=\rrn,
\text{ and }
\sum_{j\in\mathbb{N}}\1bf_{B_{2,j}}\leq 3^n.$$
For any $j\in\mathbb{N}$, if
$\|f\ast\varphi_0\ast\psi_0\1bf_{
Q_{2,j}}\|_{L^{\infty}(\rrn)}=0$,
define $\lambda_{2,j}:=0$ and $a_{2,j}:=0$;
if $\|f\ast\varphi_0\ast\psi_0
\1bf_{Q_{2,j}}\|_{L^{\infty}(\rrn)}\neq0$,
define $$
\lambda_{2,j}:=
\left\|f\ast\varphi_0\ast\psi_0
\1bf_{Q_{2,j}}\right\|_{L^{\infty}(\rrn)}
\left\|\1bf_{B_{2,j}}\right\|_X
$$
and
$$
a_{2,j}:=\frac{f\ast\varphi_0\ast\psi_0\1bf_{Q_{2,j}}}{
\|f\ast\varphi_0\ast\psi_0\1bf_{Q_{2,j}}\|_{L^{\infty}(\rrn)}
\|\1bf_{B_{2,j}}\|_X}.
$$
Then, for any $x\in\rrn$,
\begin{equation}\label{lplxe10}
f\ast\varphi_0\ast\psi_0(x)=
\sum_{j\in\mathbb{N}}\lambda_{2,j}
a_{2,j}(x)
\end{equation}
and
\begin{equation}\label{lplxe11}
\left|f\ast\varphi_0\ast\psi_0(x)\right|=
\sum_{j\in\mathbb{N}}\lambda_{2,j}
\left|a_{2,j}(x)\right|.
\end{equation}
We now prove that $f\ast\varphi_0\ast\psi_0=
\sum_{j\in\mathbb{N}}\lambda_{2,j}
a_{2,j}$
in $\mathcal{S}'(\rrn)$.
Indeed, applying \eqref{lplxe6}
with $\varphi$ and $t$ therein
replaced, respectively, by
$\varphi_0\ast\psi_0$ and $1$, we find that,
for any $x\in\rrn$,
$$
\left|f\ast\varphi_0\ast\psi_0(x)\right|
\lesssim\left(1+|x|\right)^m,
$$
which, together with \eqref{lplxe11}, further
implies that, for any given $\eta\in
\mathcal{S}(\rrn)$,
\begin{align*}
&\int_{\rrn}\sum_{j\in\mathbb{N}}
\lambda_{2,j}\left|a_{2,j}(x)\right|
\left|\eta(x)\right|\,dx\\
&\quad=\int_{\rrn}\left|f\ast\varphi_0\ast
\psi_0(x)\right|\left|\eta(x)\right|\,dx\\
&\quad\lesssim\int_{\rrn}\left(1+|x|\right)^m
\left|\eta(x)\right|\,dx\sim\int_{\rrn}
\frac{1}{(1+|x|)^{n+1}}\,dx<\infty.
\end{align*}
From this, \eqref{lplxe10}, and the Fubini
theorem, we deduce that, for any given
$\eta\in\mathcal{S}(\rrn)$,
\begin{align*}
&\int_{\rrn}f\ast\varphi_0\ast\psi_0
(x)\eta(x)\,dx\\
&\quad=\int_{\rrn}\sum_{j\in\mathbb{N}}
\lambda_{2,j}a_{2,j}(x)\eta(x)\,dx
=\sum_{j\in\mathbb{N}}\lambda_{2,j}
\int_{\rrn}a_{2,j}(x)\eta(x)\,dx.
\end{align*}
This further implies that
$f\ast\varphi_0\ast\psi_0
=\sum_{j\in\mathbb{N}}\lambda_{2,j}a_{2,j}$
holds true in $\mathcal{S}'(\rrn)$.

In addition, notice that, for any
$j\in\mathbb{N}$, $\supp(a_{2,j})
\subset B_{2,j}$ and
$$\left\|a_{2,j}\right\|_{
L^{\infty}(\rrn)}\leq\frac{1}{\|\1bf_{B_{2,j}}\|_X}.$$
These further imply that, for any $j\in\mathbb{N}$,
$$
\left\|a_{2,j}\right\|_{L^r(\rrn)}
\leq\frac{|B_{2,j}|^{1/r}}{\|\1bf_{B_{2,j}}\|_X}.
$$
Moreover, for any $j\in\mathbb{N}$,
we have $r(B_{2,j})=\sqrt{n}\geq1$.
Therefore, applying Definition \ref{loatomx},
we find that, for any $j\in\mathbb{N}$,
$a_{2,j}$ is a local-$(X,\,r,\,d)$-atom
supported in the ball $B_{2,j}$.
Combining this, the fact that
$f\ast\varphi_0\ast\psi_0
=\sum_{j\in\mathbb{N}}\lambda_{2,j}a_{2,j}$
holds true in $\mathcal{S}'(\rrn)$,
Lemma \ref{atomxl}, and Definition \ref{atom-hardyl},
we conclude that
\begin{align}\label{lplxe12}
\left\|f\ast\varphi_0\ast
\psi_0\right\|_{h_X(\rrn)}
&\sim\left\|f\ast\varphi_0\ast
\psi_0\right\|_{h^{X,r,d,s}(\rrn)}\notag\\
&\lesssim\left\|\left[\sum_{j\in\mathbb{N}}
\left(\frac{\lambda_{2,j}}
{\|\1bf_{B_{2,j}}\|_{X}}
\right)^{s}\1bf_{B_{2,j}}\right]^{\frac{1}{s}}
\right\|_{X}\notag\\
&\sim\left\|\left[\sum_{j\in\mathbb{N}}
\left\|f\ast\varphi_0\ast\psi_0
\1bf_{Q_{2,j}}
\right\|_{L^{\infty}(\rrn)}^s\1bf_{B_{
2,j}}\right]^{\frac{1}{s}}\right\|_X.
\end{align}

Next, we estimate $\|f\ast\varphi_0\ast\psi_0
\1bf_{Q_{2,j}}\|_{L^{\infty}(\rrn)}\1bf_{B_{2,j}}$
for any $j\in\mathbb{N}$. Indeed, for any
given $j\in\mathbb{N}$ and $x\in B_{2,j}$,
and for any $y\in Q_{2,j}$, we have
$|y-x|<2r(B_{2,j})=2\sqrt{n}$. This,
combined with Lemma \ref{relaxl1} with
$\varphi$, $r$, and $y$ therein replaced,
respectively, by $\varphi_0$, $\theta$,
and $z$, implies that, for any
$j\in\mathbb{N}$ and $x\in B_{2,j}$,
\begin{align}\label{lplxe13}
&\left\|f\ast\varphi_0\ast\psi_0
\1bf_{Q_{2,j}}\right\|_{L^{\infty}(\rrn)}\notag\\
&\quad\leq\sup_{y\in\rrn,\,
|y-x|<2\sqrt{n}}\left|f\ast\varphi_0
\ast\psi_0(y)\right|\notag\\
&\quad\leq\sup_{y\in\rrn,\,
|y-x|<2\sqrt{n}}\int_{\rrn}
\left|f\ast\varphi_0(z)\right|
\left|\psi_0(y-z)\right|\,dz\notag\\
&\quad\lesssim\mc^{(\theta)}\left(
f\ast\varphi_0\right)(x)
\sup_{y\in\rrn,\,|y-x|<2\sqrt{n}}\int_{\rrn}
\left(1+|x-z|\right)^{\frac{n}{\theta}}
\left|\psi_0(y-z)\right|\,dz.
\end{align}
Observe that, for any $x$, $y$, $z\in\rrn$
with $|y-x|<2\sqrt{n}$,
$$|x-z|\leq|y-x|+|y-z|<2\sqrt{n}+|y-z|.$$
From this and \eqref{lplxe13}, we deduce that,
for any $j\in\mathbb{N}$ and $x\in B_{2,j}$,
\begin{align*}
&\left\|f\ast\varphi_0\ast\psi_0
\1bf_{Q_{2,j}}\right\|_{L^{\infty}(\rrn)}\\
&\quad\lesssim\mc^{(\theta)}\left(
f\ast\varphi_0\right)(x)
\sup_{y\in\rrn,\,|y-x|<2\sqrt{n}}\int_{\rrn}
\left(1+|y-z|\right)^{\frac{n}{\theta}}
\left|\psi_0(y-z)\right|\,dz\\
&\quad\lesssim\mc^{(\theta)}\left(
f\ast\varphi_0\right)(x)\int_{\rrn}
\frac{1}{(1+|z|)^{n+1}}\,dz\sim
\mc^{(\theta)}\left(f\ast\varphi_0\right)(x).
\end{align*}
This further implies that, for any $j\in\mathbb{N}$,
$$
\left\|f\ast\varphi_0\ast\psi_0
\1bf_{Q_{2,j}}\right\|_{L^{\infty}(\rrn)}\1bf_{B_{2,j}}
\lesssim\mc^{(\theta)}\left(f\ast\varphi_0\right)
\1bf_{B_{2,j}}.
$$
Applying this, Definition \ref{Df1}(ii),
\eqref{lplxe12}, the fact that
$\sum_{j\in\mathbb{N}}\1bf_{B_{2,j}}\leq 3^n$,
Assumption \ref{assfs} with
$f_1:=f\ast\varphi_0$ and $f_j:=
0$ for any $j\in\mathbb{N}\cap[2,\infty)$, and
Definition \ref{df511}, we find that
\begin{align*}
\left\|f\ast\varphi_0\ast
\psi_0\right\|_{h_X(\rrn)}&\lesssim\left\|
\mc^{(\theta)}\left(f\ast\varphi_0\right)
\left(\sum_{j\in\mathbb{N}}
\1bf_{B_{2,j}}\right)^{\frac{1}{s}}\right\|_X\\
&\lesssim\left\|
\mc^{(\theta)}\left(f\ast\varphi_0\right)\right\|_X
\lesssim\left\|f\ast\varphi_0\right\|_X\lesssim\left\|
S_{\mathrm{loc}}(f)\right\|_X,
\end{align*}
which completes the proof of \eqref{lplxe8}.
Combining \eqref{lplxe1}, \eqref{lplxe7},
and \eqref{lplxe8}, we have
\begin{align}\label{lplxe14}
\left\|f\right\|_{h_X(\rrn)}
&\lesssim\left\|f\ast\varphi_0\ast
\psi_0\right\|_{h_X(\rrn)}+
\left\|\int_{0}^{1}f\ast\varphi_t\ast
\psi_t(\cdot)\,\frac{dt}{t}\right\|_{h_X(\rrn)}
\notag\\&\lesssim\left\|S_{\mathrm{loc}}(f)\right\|_X
<\infty.
\end{align}
This further implies that $f\in h_X(\rrn)$ and
then finishes the proof of the sufficiency.

Next, we prove the necessity. For this purpose,
let $f\in h_X(\rrn)$. Then, by Remark \ref{Th6.1r}
and Theorem \ref{relax}, we find that
\begin{equation}\label{lplxe15}
\left\|f\ast\varphi_0\right\|_{X}\lesssim
\left\|f\right\|_{h_X(\rrn)}.
\end{equation}
On the other hand, we estimate
$\widetilde{S}_{\mathrm{loc}}(f)$
with $\widetilde{S}_{\mathrm{loc}}$ as
in \eqref{lusinl}.
Indeed, from both \eqref{lusinl} and
Definition \ref{df451}, we infer that
\begin{align}\label{lplxe16}
\left\|\widetilde{S}_{\mathrm{loc}}(f)\right\|_X
&\lesssim\left\|\widetilde{S}_{\mathrm{loc}}
\left(f\ast\varphi_0\right)\right\|_{X}
+\left\|\widetilde{S}_{\mathrm{loc}}
\left(f-f\ast\varphi_0\right)\right\|_{X}\notag\\
&\lesssim\left\|\widetilde{S}_{\mathrm{loc}}
\left(f\ast\varphi_0\right)\right\|_{X}
+\left\|S\left(f-f\ast\varphi_0\right)\right\|_{X},
\end{align}
where $S(f-f\ast\varphi_0)$ is
defined as in Definition \ref{df451} with
$f$ therein replaced by $f-f\ast\varphi_0$.
By Remark \ref{Th6.1r} and Theorem \ref{relax},
we conclude that
$$f-f\ast\varphi_0\in H_X(\rrn).$$
From this, Lemma \ref{hardy-lusin},
and Theorem \ref{relax} again, we further
deduce that
\begin{equation}\label{lplxe17}
\left\|S\left(f-f\ast\varphi_0\right)\right\|_{X}
\sim\left\|f-f\ast\varphi_0\right\|_{H_X(\rrn)}
\lesssim\left\|f\right\|_{h_X(\rrn)}.
\end{equation}

Now, let $\{Q_{2,j}\}_{j\in\mathbb{N}}$
and $\{B_{2,j}\}_{j\in\mathbb{N}}$ be the
same as in the proof of the sufficiency.
For any $j\in\mathbb{N}$,
if
$\|f\ast\varphi_0\1bf_{Q_{2,j}}
\|_{L^{\infty}(\rrn)}=0$,
define $$\mu_{2,j}:=0\text{ and }b_{2,j}:=0;$$
if $\|f\ast\varphi_0\1bf_{Q_{2,j}}
\|_{L^{\infty}(\rrn)}\neq0$,
define $$
\mu_{2,j}:=
\left\|f\ast\varphi_0\1bf_{Q_{2,j}}
\right\|_{L^{\infty}(\rrn)}
\left\|\1bf_{B_{2,j}}\right\|_X
$$
and
$$
b_{2,j}:=\frac{f\ast\varphi_0\1bf_{Q_{2,j}}}{
\|f\ast\varphi_0\1bf_{Q_{2,j}}\|_{L^{\infty}(\rrn)}
\|\1bf_{B_{2,j}}\|_X}.
$$
Then, repeating an argument similar to that
used in the proof of the sufficiency of
the present theorem with $\varphi_0\ast
\psi_0$ therein replaced by $\varphi_0$, we
find that, for any $j\in\mathbb{N}$,
$b_{2,j}$ is a local-$(X,\,r,\,d)$-atom
supported in the ball $B_{2,j}$, and
$$
f\ast\varphi_0=\sum_{j\in\mathbb{N}}\mu_{2,j}
b_{2,j}
$$
in $\mathcal{S}'(\rrn)$.
This further implies that, for any $t\in(0,\infty)$
and $x\in\rrn$,
\begin{align*}
f\ast\varphi_0\ast\varphi_0(x)
&=\int_{\rrn}f\ast\varphi_0(y)
\varphi_t(x-y)\,dy\\
&=\sum_{j\in\mathbb{N}}\mu_{2,j}
\int_{\rrn}b_{2,j}(y)\varphi_t(x-y)\,dy\\
&=\sum_{j\in\mathbb{N}}\mu_{2,j}
b_{2,j}\ast\varphi_t(x).
\end{align*}
By this, \eqref{lusinl}, the
Minkowski inequality, and Definition \ref{Df1}(ii),
we obtain
\begin{align}\label{lplxe18}
&\left\|\widetilde{S}_{\mathrm{loc}}\left(f
\ast\varphi_0\right)\right\|_X\notag\\
&\quad\leq\left\|\sum_{j\in\mathbb{N}}
\mu_{2,j}\widetilde{S}_{\mathrm{loc}}
\left(b_{2,j}\right)\right\|_X\notag\\
&\quad\lesssim\left\|\sum_{j\in\mathbb{N}}
\mu_{2,j}\widetilde{S}_{\mathrm{loc}}
\left(b_{2,j}\right)\1bf_{2B_{2,j}}\right\|_X
+\left\|\sum_{j\in\mathbb{N}}
\mu_{2,j}\widetilde{S}_{\mathrm{loc}}
\left(b_{2,j}\right)\1bf_{(2B_{2,j})^{\complement}}
\right\|_X\notag\\
&\quad=:\mathrm{VI}_1+\mathrm{VI}_2.
\end{align}
For any $j\in\mathbb{N}$,
from \eqref{lusinl},
Lemma \ref{lusinlp} with $p:=r$
and $f:=b_{2,j}$, and
the fact that $b_{2,j}$
is a local-$(X,\,\infty,\,r)$-atom
supported in the ball $B_{2,j}$, it follows that
\begin{align*}
\left\|\widetilde{S}_{\mathrm{loc}}
\left(b_{2,j}\right)\right\|_{L^r(\rrn)}
\leq\left\|S\left(b_{2,j}\right)\right\|_{L^r(\rrn)}
\lesssim\left\|b_{2,j}\right\|_{L^r(\rrn)}
\lesssim\frac{|B_{2,j}|^{1/r}}{\|\1bf_{B_{2,j}}\|_{X}}.
\end{align*}
This, together with an argument similar
to that used in the estimation of $\mathrm{II}_1$,
further implies that
\begin{align}\label{lplxe19}
\mathrm{VI}_1&\lesssim
\left\|\left[\sum_{j\in\mathbb{N}}
\left(\frac{\mu_{2,j}}
{\|\1bf_{B_{2,j}}\|_{X}}
\right)^{s}\1bf_{B_{2,j}}\right]^{\frac{1}{s}}
\right\|_{X}\notag\\
&\sim\left\|\left[\sum_{j\in\mathbb{N}}
\left\|f\ast\varphi_0\1bf_{Q_{2,j}}
\right\|_{L^{\infty}(\rrn)}^s\1bf_{B_{
2,j}}\right]^{\frac{1}{s}}\right\|_X.
\end{align}
On the other hand, applying Lemma \ref{lusinl-atom}
with $a:=b_{2,j}$ for any $j\in\mathbb{N}$
and \eqref{Atogxe6} with
$\{\lambda_{j}\}_{j\in\mathbb{N}}$ and
$\{B_{j}\}_{j\in\mathbb{N}}$ therein replaced,
respectively, by $\{\mu_{2,j}\}_{j\in\mathbb{N}}$
and $\{B_{2,j}\}_{j\in\mathbb{N}}$, we find that
\begin{align*}
\mathrm{VI}_2&\lesssim
\left\|\left\{
\sum_{j\in\mathbb{N}}
\left[\frac{\mu_{2,j}}{\|\1bf_
{B_{2,j}}\|_{X}}\mc^{(\theta)}(\1bf_
{B_{2,j}})\right]^s
\right\}^{\frac{1}{s}}
\right\|_{X}\\
&\lesssim\left\|\left[
\sum_{j\in\mathbb{N}}\left(\frac{\mu_{2,j}}
{\|\1bf_{B_{2,j}}\|_{X}}\right)^s\1bf_
{B_{2,j}}\right]^{\frac{1}{s}}\right\|_{X}\\
&\sim\left\|\left[\sum_{j\in\mathbb{N}}
\left\|f\ast\varphi_0\1bf_{Q_{2,j}}
\right\|_{L^{\infty}(\rrn)}^s\1bf_{B_{
2,j}}\right]^{\frac{1}{s}}\right\|_X.
\end{align*}
Combining this, \eqref{lplxe18}, and
\eqref{lplxe19}, we further obtain
\begin{equation}\label{lplxe20}
\left\|\widetilde{S}_{\mathrm{loc}}\left(f
\ast\varphi_0\right)\right\|_X
\lesssim\left\|\left[\sum_{j\in\mathbb{N}}
\left\|f\ast\varphi_0\1bf_{Q_{2,j}}
\right\|_{L^{\infty}(\rrn)}^s\1bf_{B_{
2,j}}\right]^{\frac{1}{s}}\right\|_X.
\end{equation}

Next, in order to complete the estimation of
$\|\widetilde{S}_{\mathrm{loc}}(f
\ast\varphi_0)\|_X$,
we first deal with
$\|f\ast\varphi_0\1bf_{Q_{2,j}}
\|_{L^{\infty}(\rrn)}\1bf_{B_{2,j}}$ for any
$j\in\mathbb{N}$. To this end,
define a function $\phi_0$ by setting,
for any $x\in\rrn$,
$$
\phi_0(x):=\frac{1}{2^n}\varphi_0\left(\frac{x}{
2}\right).
$$
Then, by Definition \ref{df511}(ii) with
$\phi:=\phi_0$ and $a:=4\sqrt{n}$, we conclude that,
for any $j\in\mathbb{N}$ and $x\in B_{2,j}$,
\begin{align*}
\left\|f\ast\varphi_0\1bf_{Q_{2,j}}
\right\|_{L^{\infty}(\rrn)}&=
\sup_{y\in Q_{2,j}}\left|f\ast\varphi_0(y)\right|
\leq\sup_{y\in\rrn,\,|y-x|<2\sqrt{n}}
\left|f\ast\left(\phi_0\right)_{\frac{1}{2}}
(y)\right|\\&\leq m^*_{4\sqrt{n}}
\left(f,\phi_0\right)(x),
\end{align*}
where $\left(\phi_0\right)_{\frac{1}{2}}$
is defined as in \eqref{smaxe1} with
$\phi:=\phi_0$ and $t:=\frac{1}{2}$.
This implies that, for any $j\in\mathbb{N}$,
\begin{equation}\label{lplxe21}
\left\|f\ast\varphi_0\1bf_{Q_{2,j}}
\right\|_{L^{\infty}(\rrn)}\1bf_{B_{2,j}}
\leq m^*_{4\sqrt{n}}
\left(f,\phi_0\right)\1bf_{B_{2,j}}.
\end{equation}
On the other hand, notice that
$$
\int_{\rrn}\phi_0(x)\,dx
=\int_{\rrn}\varphi_0(x)\,dx
=\widehat{\varphi}_0(\0bf)\neq0.
$$
From this, \eqref{lplxe20}, Definition \ref{Df1}(ii),
\eqref{lplxe21}, the fact that
$\sum_{j\in\mathbb{N}}\1bf_{B_{2,j}}\leq 3^n$,
Remark \ref{Th6.1r}, and Lemma \ref{Th6.1l1},
it further follows that
\begin{align*}
\left\|\widetilde{S}_{\mathrm{loc}}
\left(f\ast\varphi_0\right)\right\|_{X}
&\lesssim\left\|m^*_{4\sqrt{n}}
\left(f,\phi_0\right)\left(
\sum_{j\in\mathbb{N}}\1bf_{
B_{2,j}}\right)^{\frac{1}{s}}\right\|_X\\
&\lesssim\left\|m^*_{4\sqrt{n}}
\left(f,\phi_0\right)\right\|_X\sim
\left\|f\right\|_{h_X(\rrn)},
\end{align*}
which then completes the estimation of
$\|\widetilde{S}_{\mathrm{loc}}(f
\ast\varphi_0)\|_{X}$. Combining this,
\eqref{lplxe16}, \eqref{lplxe17}, and
\eqref{lplxe15}, we further conclude that
\begin{equation*}
\left\|S_{\mathrm{loc}}(f)\right\|_{X}
\lesssim\left\|f\ast\varphi_0\right\|_X
+\left\|\widetilde{S}_{\mathrm{loc}}(f)\right\|_{X}
\lesssim\left\|f\right\|_{h_X(\rrn)}<\infty.
\end{equation*}
This implies that $S_{\mathrm{loc}}(f)\in X$,
and hence finishes the proof
of the necessity.
Moreover, by both \eqref{lplxe14} and
\eqref{lplxe21}, we find that
$$
\left\|f\right\|_{h_X(\rrn)}
\sim\left\|S_{\mathrm{loc}}(f)\right\|_{X},
$$
where the positive equivalence
constants are independent of $f$,
which then completes the proof of Theorem \ref{lplx}.
\end{proof}

Next, we are devoted to establishing
the Littlewood--Paley $g$-function
characterization of the local Hardy
space $h_X(\rrn)$ as follows.

\begin{theorem}\label{lplgx}
Let $X$ be a ball quasi-Banach function
space satisfying Assumption \ref{assas}
with some $s\in(0,1]$.
Assume that, for the above $s$, Assumption
\ref{assfs} holds true for both $X$ and $X^{s/2}$.
Then $f\in h_X(\rrn)$ if and only if
$f\in\mathcal{S}'(\rrn)$ and
$g_{\mathrm{loc}}(f)\in X$.
Moreover, there exists a constant $C\in[1,\infty)$
such that, for any $f\in h_X(\rrn)$,
$$
C^{-1}\left\|f\right\|_{h_X(\rrn)}
\leq\left\|g_{\mathrm{loc}}(f)\right\|_X
\leq C\left\|f\right\|_{h_X(\rrn)}.
$$
\end{theorem}

To show this theorem,
we need a technique estimate about the
local Littlewood--Paley $g$-function.
For this purpose,
let $\varphi\in\mathcal{S}(\rrn)$
be as in Definition \ref{df551},
$b,\ t\in(0,\infty)$, and, for
any $f\in\mathcal{S}'(\rrn)$
and $x\in\rrn$,
\begin{equation*}
\left(\varphi^{*}_{t}f\right
)_{b}(x):=\sup_{y\in\rrn}
\frac{|f\ast\varphi_t(y)|}{(1+t^{-1}
|x-y|)^{b}}
\end{equation*}
and
\begin{equation*}
g_{b,*}^{(l)}(f)(x):=\left\{
\int_{0}^{1}\left[\left(\varphi_t^*f
\right)_b(x)\right]^2\,\frac{dt}{t}
\right\}^{\frac{1}{2}}.
\end{equation*}
Then we have the following estimate
about $g_{\mathrm{loc}}$ and $g_{b,*}^{(l)}$,
which plays an essential role in
the proof of the Littlewood--Paley
$g$-function characterization of $h_X(\rrn)$.

\begin{proposition}\label{peetreg}
Let $0<\theta<s\leq1$, $b\in(\frac{ns}{2\theta},
\infty)$, and
$X$ be a ball quasi-Banach function space satisfying
that there exists a positive constant $C$ such that,
for any $\{f_{j}\}_{j\in\mathbb{N}}
\subset L^1_{{\rm loc}}(\rrn)$,
\begin{equation*}
\left\|\left\{\sum_{j\in\mathbb{N}}
\left[\mc^{(\theta)}
(f_{j})\right]^s\right\}^{1/s}
\right\|_{X^{s/2}}\leq C
\left\|\left(\sum_{j\in\mathbb{N}}
|f_{j}|^{s}\right)^{1/s}\right\|_{X
^{s/2}}.
\end{equation*}
Then there exists a positive constant
$C$ such that, for any $f\in\mathcal{S}'(\rrn)$,
$$
\left\|g_{b,*}^{(l)}(f)\right\|_{X}
\leq C\left\|g_{\mathrm{loc}}(f)\right\|_X.
$$
\end{proposition}

To prove Proposition \ref{peetreg},
we first present the following
auxiliary lemma which was obtained in
\cite[Lemma 3.21]{WYYZ}.

\begin{lemma}\label{lplgxl1}
Let $\varphi\in\mathcal{S}(\rrn)$ satisfy
that, for any $\xi\in\rrn\setminus\{\0bf\}$,
there exists a $j\in\mathbb{Z}$ such that
$\widehat{\varphi}(2^j\xi)\neq0$.
Then, for any given $N_0\in\mathbb{N}$
and $r\in(0,\infty)$, there exists
a positive constant $C_{(n,N_0,r,\varphi)}$,
depending only on $n$, $N_0$, $r$,
and $\varphi$, such that, for any
$t\in[1,2]$, $b\in(0,N_0]$, $j\in\mathbb{Z}$,
$f\in\mathcal{S}'(\rrn)$, and $x\in\rrn$,
\begin{align*}
&\left[\left(\varphi^*
_{2^{-j}t}f\right)_b(x)\right]^{r}\\
&\quad\leq C_{(n,N_0,r,\varphi)}
\sum_{k=0}^\infty2^{-kN_0r}2^{(k+j)n}\int_{\rrn}
\frac{|\varphi_{2^{-k-j}t}\ast f(y)|^r}{(1+2^jt^{-1}|x-y|)
^{br}}\,dy.
\end{align*}
\end{lemma}

Furthermore, we need the following
well-known Aoki--Rolewicz theorem
(see \cite[Exercise 1.46]{LGCF} and
also \cite{a42,r57}).

\begin{lemma}\label{glamsl1}
Let $X$ be a ball quasi-Banach function
space. Then there exists a $\nu\in(0,1]$
such that, for any $\{f_{j}\}_{j\in\mathbb{N}}
\subset\Msc(\rrn)$,
$$
\left\|\sum_{j\in\mathbb{N}}
f_j\right\|_X^{\nu}\leq
4\sum_{j\in\mathbb{N}}\left\|
f_j\right\|_X^{\nu}.
$$
\end{lemma}

Based on the above two
conclusions, we now show
Proposition \ref{peetreg}.

\begin{proof}[Proof of Proposition
\ref{peetreg}]
Let all the symbols be as
in the present proposition and
$f\in\mathcal{S}'(\rrn)$. Fix
an $N_0\in\mathbb{N}\cap[b,\infty)$.
Then, from Lemma \ref{lplgxl1} with
$r:=\frac{2\theta}{s}$ and the
Minkowski inequality, it follows that,
for any $x\in\rrn$,
\begin{align}\label{peetrege1}
&\left[g_{b,*}(f)(x)\right]^{2}\notag\\&\quad=
\sum_{j\in\mathbb{N}}\int_{2^{-j}}^{
2^{-j+1}}\left|\left(\varphi_t^{*}f
\right)(x)\right|^2\,\frac{dt}{t}\notag\\
&\quad=\sum_{j\in\mathbb{N}}\int_{1}^{2}
\left|\left(\varphi_{2^{-j}t}^{*}f
\right)(x)\right|^2\,\frac{dt}{t}\notag\\
&\quad\lesssim\sum_{j\in\mathbb{N}}
\left[
\sum_{k=0}^{\infty}2^{-\frac{2kN_0\theta}{s}}
2^{(k+j)n}\int_{\rrn}\frac{[
\int_{1}^{2}|\varphi_{2^{-k-j}t}\ast
f(y)|^2\,\frac{dt}{t}
]^{\frac{\theta}{s}}}{(1+2^j
|x-y|)^{\frac{2b\theta}{s}}}\,dy
\right]^{\frac{s}{\theta}}.
\end{align}
Moreover, by the assumption
$b>\frac{ns}{2\theta}$, we find that,
for any $j\in\mathbb{N}$, $k\in\zp$, and $x\in\rrn$,
\begin{align*}
&\int_{\rrn}\frac{[
\int_{1}^{2}|\varphi_{2^{-k-j}t}\ast
f(y)|^2\,\frac{dt}{t}
]^{\frac{\theta}{s}}}{(1+2^j
|x-y|)^{\frac{2b\theta}{s}}}\,dy\\
&\quad=\int_{|y-x|<2^{-j}}\frac{[
\int_{1}^{2}|\varphi_{2^{-k-j}t}\ast
f(y)|^2\,\frac{dt}{t}
]^{\frac{\theta}{s}}}{(1+2^j
|x-y|)^{\frac{2b\theta}{s}}}\,dy
+\sum_{j\in\mathbb{N}}\int_{
2^{i-j-1}\leq|y-x|<2^{i-j}}\cdots\\
&\quad\lesssim\sum_{i=0}^{\infty}
2^{-\frac{2ib\theta}{s}}\int_{|y-x|<2^{i-j}}
\left[\int_{1}^{2}|\varphi_{2^{-k-j}t}\ast
f(y)|^2\,\frac{dt}{t}
\right]^{\frac{\theta}{s}}\,dy\\
&\quad\lesssim2^{-jn}\mc\left(
\left[\int_{1}^{2}\left|\varphi_{
2^{-k-j}t}\ast f(\cdot)\right|^2\,\frac{dt}{t}
\right]^{\frac{\theta}{s}}\right)(x)
\sum_{i=0}^{\infty}2^{-\frac{2i\theta}{s}(b-
\frac{ns}{2\theta})}\\
&\quad\sim2^{-jn}\mc\left(
\left[\int_{1}^{2}\left|\varphi_{
2^{-k-j}t}\ast f(\cdot)\right|^2\,\frac{dt}{t}
\right]^{\frac{\theta}{s}}\right)(x).
\end{align*}
This, combined with \eqref{peetrege1}
and the Minkowski inequality again, further
implies that, for any $x\in\rrn$,
\begin{align*}
&\left[g^{(l)}_{b,*}(f)(x)
\right]^{\frac{2\theta}{s}}\\
&\quad\lesssim\left\{
\sum_{j\in\mathbb{N}}\left[
\sum_{k=0}^{\infty}2^{-\frac{2k\theta}{s}
(N_0-\frac{ns}{2\theta})}\mc\left(
\left[\int_{1}^{2}\left|\varphi_{
2^{-k-j}t}\ast f(\cdot)\right|^2\,\frac{dt}{t}
\right]^{\frac{\theta}{s}}\right)(x)
\right]^{\frac{s}{\theta}}
\right\}^{\frac{\theta}{s}}\\
&\quad\lesssim\sum_{k=0}^{\infty}
2^{-\frac{2k\theta}{s}(N_0-\frac{ns}{2\theta})}
\left\{\sum_{j\in\mathbb{N}}\left[
\mc\left(\left[\int_{1}^{2}\left|\varphi_{
2^{-k-j}t}\ast f(\cdot)\right|^2\,\frac{dt}{t}
\right]^{\frac{\theta}{s}}\right)(x)
\right]^{\frac{s}{\theta}}
\right\}^{\frac{\theta}{s}}.
\end{align*}
Using this, Definition \ref{Df1}(ii),
Lemma \ref{glamsl1} with
$X$ replaced by $X^{\frac{s}{2\theta}}$,
and Remark \ref{power}(ii) with
$X$ replaced by $X^{\frac{s}{2}}$,
we find that there exists
a constant $\nu\in(0,1]$ such that, for any $x\in\rrn$,
\begin{align}\label{peetrege2}
&\left\|g_{b,*}^{(l)}(f)
\right\|_{X}^{\frac{2\theta\nu}{s}}\notag\\
&\quad=\left\|\left|g_{b,*}^{(l)}
(f)\right|^{\frac{2\theta}{s}}\right\|^{\nu}_{X^{
\frac{s}{2\theta}}}\notag\\
&\quad\lesssim\sum_{k=0}^{\infty}
2^{-\frac{2k\theta\nu}{s}(N_0-\frac{ns}{2\theta})}
\left\|\left\{\sum_{j\in\mathbb{N}}
\left[\mc\left(\left[\int_{1}^{2}\left|\varphi_{
2^{-k-j}t}\ast f(\cdot)\right|^2\,\frac{dt}{t}
\right]^{\frac{\theta}{s}}\right)
\right]^{\frac{s}{\theta}}
\right\}^{\frac{\theta}{s}}
\right\|_{X^{\frac{s}{2\theta}}}^{\nu}\notag\\
&\quad\lesssim\sum_{k=0}^{\infty}
2^{-\frac{2k\theta\nu}{s}(N_0-\frac{ns}{2\theta})}
\left\|\left[\sum_{j\in\mathbb{N}}
\int_{1}^{2}\left|\varphi_{
2^{-k-j}t}\ast f(\cdot)\right|^2\,\frac{dt}{t}
\right]^{\frac{\theta}{s}}
\right\|_{X^{\frac{s}{2\theta}}}^{\nu}\notag\\
&\quad\sim\sum_{k=0}^{\infty}
2^{-\frac{2k\theta\nu}{s}(N_0-\frac{ns}{2\theta})}
\left\|\left[\sum_{j\in\mathbb{N}}
\int_{1}^{2}\left|\varphi_{
2^{-k-j}t}\ast f(\cdot)\right|^2\,\frac{dt}{t}
\right]^{\frac{1}{2}}
\right\|_{X}^{\frac{2\theta\nu}{s}}.
\end{align}
Observe that, for any $k\in\zp$ and $x\in\rrn$,
\begin{align*}
&\sum_{j\in\mathbb{N}}
\int_{1}^{2}\left|\varphi_{
2^{-k-j}t}\ast f(x)\right|^2\,\frac{dt}{t}\\
&\quad=\sum_{j\in\mathbb{N}}\int_{2^{-k-j}}^{
2^{-k-j+1}}\left|\varphi_t\ast f(x)
\right|^2\,\frac{dt}{t}\\
&\quad=\int_{0}^{2^{-k}}\left|\varphi_t\ast f(x)
\right|^2\,\frac{dt}{t}
\leq g_{\mathrm{loc}}(f)(x).
\end{align*}
Thus, combining this, \eqref{peetrege2},
Definition \ref{Df1}(ii),
and the assumption $N_0\geq b>
\frac{ns}{2\theta}$, we further
obtain
$$
\left\|g^{(l)}_{
b,*}(f)\right\|_{X}^{\frac{2\theta\nu}{s}}
\lesssim\left\|g_{\mathrm{loc}}\right\|_
X^{\frac{2\theta\nu}{s}}\sum_{k=0}^{\infty}
2^{-\frac{2k\theta\nu}{s}(N_0-\frac{ns}{2\theta})}
\sim\left\|g_{\mathrm{loc}}\right\|_
X^{\frac{2\theta\nu}{s}},
$$
which then completes the proof of
Proposition \ref{peetreg}.
\end{proof}

We then turn to show Theorem \ref{lplgx}.
To this end, we need the
$L^p$ boundedness of the Littlewood--Paley
$g$-function as follows, which was obtained
in \cite[Theorem 7.7]{fs82}.

\begin{lemma}\label{glp}
Let $p\in(1,\infty)$ and
$g$ be as in Definition \ref{df452}.
Then there exists a positive constant
$C$ such that, for any
$f\in L^p(\rrn)$,
$$
\left\|g(f)\right\|_{L^p(\rrn)}\leq
C\left\|f\right\|_{L^p(\rrn)}.
$$
\end{lemma}

Moreover, let $\varphi$ be as in
Definition \ref{df551}. For any
$f\in\mathcal{S}'(\rrn)$ and $x\in\rrn$,
let
$$
\widetilde{g}_{\mathrm{loc}}(f)(x):=
\left[\int_{0}^{1}\left|f\ast
\varphi_t(x)\right|^2\,\frac{dt}{t}\right]^{
\frac{1}{2}}.
$$
Then, repeating an argument similar
to that used in the proof of
\cite[p.\,53]{WYY} with
$\widetilde{S}_l$ therein replaced
by $\widetilde{g}_{\mathrm{loc}}$,
we then obtain the following
technique lemma about $\widetilde{g}_{\mathrm{loc}}$,
which is important in the
proof of Theorem \ref{lplgx}; we omit
the details.

\begin{lemma}\label{gl-atom}
Let $X$ be a ball quasi-Banach function
space, $d\in\zp$ and $\theta\in(0,\infty)$.
Then there exists a positive constant
$C$ such that, for any local-$(X,\,\infty,\,d)$-atom
$a$ supported in the ball $B\in\mathbb{B}$
with $r(B)\in[1,\infty)$,
$$
\widetilde{g}_{\mathrm{loc}}(a)\1bf_{
(2B)^{\complement}}\leq C\frac{1}{\|\1bf_{B}
\|_X}\mc^{(\theta)}\left(\1bf_B\right).
$$
\end{lemma}

Via the above preparations, we now
prove Theorem \ref{lplgx}.

\begin{proof}[Proof of Theorem \ref{lplgx}]
Let all the symbols be
as in the present theorem. We first
prove the necessity. For this purpose,
let $f\in h_X(\rrn)$.
Then, repeating an argument similar to
that used in the proof of the necessity of
Theorem \ref{lplx} with Lemmas
\ref{hardy-lusin}, \ref{lusinlp}, and
\ref{lusinl-atom}, therein replaced, respectively,
by Lemmas \ref{lusinl}, \ref{glp}, and
\ref{gl-atom}, we find that
$g_{\mathrm{loc}}(f)\in X$ and
\begin{equation}\label{lplgxe1}
\left\|g_{\mathrm{loc}}(f)\right\|
\lesssim\left\|f\right\|_{h_X(\rrn)},
\end{equation}
which then completes the proof of the necessity.

Now, we show the sufficiency.
To this end, let $f\in\mathcal{S}'(\rrn)$
with $g_{\mathrm{loc}}(f)\in X$ and
$b\in(\frac{ns}{2\theta},\infty)$.
Notice that, for any given $t\in(0,\infty)$
and $x\in\rrn$, and for any $y\in B(x,t)$,
we have $\frac{1}{1+t^{-1}|x-y|}\sim1$.
This further implies that, for any
$x\in\rrn$,
\begin{align*}
\widetilde{S}_{\mathrm{loc}}(f)(x)
&\sim\left\{\int_{0}^{1}\int_{B(x,t)}
\left[\frac{|f\ast\varphi_t(y)|}{
(1+t^{-1}|x-y|)^b}\right]^2
\,\frac{dy\,dt}{t^{n+1}}
\right\}^{\frac{1}{2}}\\
&\lesssim\left\{\int_{0}^{1}
\left[\left(\varphi^*_tf\right)_b(x)
\right]^2\left|B(x,t)\right|
\,\frac{dt}{t^{n+1}}
\right\}^{\frac{1}{2}}\sim g^{(l)}_{b,*}(f)(x).
\end{align*}
From this, Definition \ref{df551}, \eqref{lusinl},
Definition \ref{Df1}(ii), and Proposition \ref{peetreg},
it follows that
\begin{align*}
\left\|S_{\mathrm{loc}}(f)\right\|_X
&\lesssim\left\|f\ast\varphi_0\right\|_X
+\left\|\widetilde{S}_{\mathrm{loc}}(f)\right\|_X\\
&\lesssim\left\|f\ast\varphi_0\right\|_X
+\left\|g_{b,*}^{(l)}(f)\right\|_X\\&\lesssim
\left\|g_{\mathrm{loc}}(f)\right\|_X<\infty,
\end{align*}
which, combined with Theorem \ref{lplx},
further implies that $f\in h_X(\rrn)$ and
\begin{equation}\label{lplgxe2}
\left\|f\right\|_{h_X(\rrn)}
\sim\left\|S_{\mathrm{loc}}(f)\right\|_X
\lesssim\left\|g_{\mathrm{loc}}(f)\right\|_X.
\end{equation}
This implies that the sufficiency holds true.
Therefore, by both \eqref{lplgxe1} and \eqref{lplgxe2},
we conclude that
$$
\left\|f\right\|_{h_X(\rrn)}\sim
\left\|g_{\mathrm{loc}}(f)\right\|_X
$$
with the positive equivalence constants
independent of $f$, which then completes the
proof of Theorem \ref{lplgx}.
\end{proof}

We next investigate the
Littlewood--Paley $g_{\lambda}^*$-function
characterization of the Hardy space $h_X(\rrn)$.
Namely, we have the following conclusion.

\begin{theorem}\label{lplglx}
Let $X$ be a ball quasi-Banach function
space satisfying both Assumptions
\ref{assfs} and \ref{assas} with the same
$s\in(0,1]$, and $\lambda\in(\max\{1,
\frac{1}{s}\},\infty)$.
Then $f\in h_X(\rrn)$ if and only if
$f\in\mathcal{S}'(\rrn)$ and
$(g_{\lambda}^*)_{\mathrm{loc}}(f)\in X$.
Moreover, there exists a constant $C\in[1,\infty)$
such that, for any $f\in h_X(\rrn)$,
$$
C^{-1}\left\|f\right\|_{h_X(\rrn)}
\leq\left\|\left(g_{\lambda}^*
\right)_{\mathrm{loc}}(f)\right\|_X
\leq C\left\|f\right\|_{h_X(\rrn)}.
$$
\end{theorem}

In order to prove this theorem, we
first show the following auxiliary estimate
about $(g_{\lambda}^*)_{\mathrm{loc}}$ and
$S_{\mathrm{loc}}$.

\begin{proposition}\label{glams}
Let $X$ be a ball quasi-Banach function,
$s\in(0,1]$, and $\lambda\in(\max\{
1,\frac{2}{s}\},\infty)$. Assume that
$X^{1/s}$ is a ball Banach function space
and the Hardy--Littlewood maximal operator
$\mc$ is bounded on $(X^{1/s})'$.
Then there exists a positive constant $C$
such that, for any $f\in\mathcal{S}'(\rrn)$,
$$
\left\|\left(g_{\lambda}^*\right)_
{\mathrm{loc}}(f)\right\|_X\leq
C\left\|S_{\mathrm{loc}}(f)\right\|_X.
$$
\end{proposition}

To obtain this proposition,
we need the following
estimate on the change of angles
in $X$-tent spaces, which is just
\cite[Theorem 3.3]{CWYZ}.

\begin{lemma}\label{glamsl2}
Let $X$ be a ball quasi-Banach function and
$s\in(0,1]$. Assume that
$X^{1/s}$ is a ball Banach function space
and the Hardy--Littlewood maximal operator
$\mc$ is bounded on $(X^{1/s})'$.
Then there exists a positive constant
$C$ such that, for any $\alpha\in[1,\infty)$
and any measurable function $F$ on $\rrnp$,
$$
\left\|\mathcal{A}^{(\alpha)}(F)\right\|_X
\leq C\alpha^{n\max\{\frac{1}{2},\frac{1}{s}\}}
\left\|\mathcal{A}^{(1)}(F)\right\|_X.
$$
\end{lemma}

Now, we show Proposition \ref{glams}.

\begin{proof}[Proof of Proposition \ref{glams}]
Let all the symbols be as in the
present proposition, $f\in\mathcal{S}'(\rrn)$,
and, for any $\alpha$ and $x\in\rrn$,
\begin{equation*}
\widetilde{S}^{(\alpha)}_{\mathrm{loc}}(f)
(x):=\left[\int_{0}^{1}\int_{B(x,\alpha t)}
\left|f\ast\varphi_t(y)\right|^2\,\frac{dy\,dt}{
t^{n+1}}\right]^{\frac{1}{2}}.
\end{equation*}
Then, applying Lemma \ref{glamsl2}
with $F:=f\ast\varphi_t\1bf_{\{\tau\in(0,1)\}}$,
we find that, for any $\alpha\in[1,\infty)$,
\begin{align}\label{glamse1}
\left\|\widetilde{S}^{(\alpha)}
_{\mathrm{loc}}(f)\right\|_{X}&=
\left\|\left\{\int_{\Gamma(\cdot)}
\left[\left|f\ast
\varphi_t(y)\right|\1bf_{\{\tau:\,
\tau\in(0,1)\}}
(t)\right]^2\,\frac{dy\,dt}{t^{n+1}}
\right\}^{\frac{1}{2}}\right\|_X\notag\\
&=\left\|\mathcal{A}^{(\alpha)}
\left(f\ast\varphi_t\1bf_{
\{\tau\in(0,1)\}}\right)\right\|_X\notag\\
&\lesssim\alpha^{n\max\{\frac{1}{2},
\frac{1}{s}\}}\left\|\mathcal{A}^{(1)}
\left(f\ast\varphi_t\1bf_{
\{\tau\in(0,1)\}}\right)\right\|_X\notag\\
&\sim\alpha^{n\max\{\frac{1}{2},
\frac{1}{s}\}}\left\|\widetilde{S}_
{\mathrm{loc}}(f)\right\|_X.
\end{align}
In addition, for any $x\in\rrn$, we have
\begin{align*}
&\left(g_{\lambda}^*
\right)_{\mathrm{loc}}(f)(x)\\
&\quad\leq
\left|f\ast\varphi_0(x)\right|
+\left\{\int_{0}^{1}
\int_{|y-x|<t}\left(
\frac{t}{t+|x-y|}\right)^{\lambda n}
\left|f\ast\varphi_t(x)\right|^2
\,\frac{dy\,dt}{t^{n+1}}\right.\\
&\qquad\left.
+\sum_{j\in\mathbb{N}}\int_{0}^{1}
\int_{2^{j-1}t\leq|y-x|
<2^jt}\cdots\right\}^{
\frac{1}{2}}\\
&\quad\lesssim\left|f\ast\varphi_0(x)\right|
+\left\{\int_{0}^{1}
\int_{|y-x|<t}\left(
\frac{t}{t+|x-y|}\right)^{\lambda n}
\left|f\ast\varphi_t(x)\right|^2
\,\frac{dy\,dt}{t^{n+1}}\right\}^{\frac{1}{2}}\\
&\qquad+\sum_{j\in\mathbb{N}}
2^{-\frac{j\lambda n}{2}}
\left\{\int_{0}^{1}\int_{|y-x|<2^{j}t}\cdots
\right\}^{\frac{1}{2}}\\
&\quad\sim\left|f\ast\varphi_0(x)\right|
+\sum_{j=0}^{\infty}2^{-\frac{j\lambda n}{2}}
\widetilde{S}^{(2^j)}_{\mathrm{loc}}(f)(x).
\end{align*}
This, combined with Definition \ref{Df1}(ii),
Lemma \ref{glamsl1} with $f_1$ therein
replaced by $|f\ast\varphi_0|$ and
$f_{j+2}$ therein replaced by
$2^{-\frac{j\lambda n}{2}}\widetilde{S}^{
(2^j)}_{\mathrm{loc}}(f)$
for any $j\in\zp$, \eqref{glamse1}
with $\alpha$ therein replaced by $2^j$ for
any $j\in\zp$, and the assumption
$\lambda>\max\{1,\frac{2}{s}\}$,
further implies that
there exists a $\nu\in(0,1]$ such that
\begin{align*}
\left\|\left(g^*_{\lambda}
\right)_{\mathrm{loc}}(f)\right\|_X^{\nu}
&\lesssim\left\|f\ast\varphi_0\right\|_X^{\nu}
+\sum_{j=0}^{\infty}2^{-\frac{j\lambda n\nu}{2}}
\left\|\widetilde{S}^{(2^j)}
_{\mathrm{loc}}(f)\right\|_X^{\nu}\\
&\lesssim\left\|f\ast\varphi_0\right\|_X^{\nu}
+\left\|\widetilde{S}_{\mathrm{loc}}(f)
\right\|_X^{\nu}\sum_{j=0}^{\infty}
2^{-\frac{jn\nu}{2}\max\{\lambda-1,\lambda-
\frac{2}{s}\}}\\
&\sim\left\|f\ast\varphi_0\right\|_X^{\nu}
+\left\|\widetilde{S}
_{\mathrm{loc}}(f)\right\|_X^{\nu}
\lesssim\left\|S_{\mathrm{loc}}
(f)\right\|_X^{\nu},
\end{align*}
which then completes the proof of
Proposition \ref{glams}.
\end{proof}

Via Proposition \ref{glams}, we
next prove Theorem \ref{lplglx}.

\begin{proof}[Proof of Theorem \ref{lplglx}]
Let all the symbols be as in the
present theorem. We first show the
necessity. For this purpose, let $f\in
\mathcal{S}'(\rrn)$.
Then, by Proposition \ref{glams} and Theorem
\ref{lplx}, we conclude that
\begin{equation}\label{lplglxe1}
\left\|\left(g_{\lambda}
^{*}\right)_{\mathrm{loc}}(f)\right\|_X
\lesssim\left\|S_{\mathrm{loc}}(f)\right\|_X
\sim\left\|f\right\|_{h_X(\rrn)}<\infty,
\end{equation}
which implies that $(g_{\lambda}^{*})_{
\mathrm{loc}}(f)\in X$
and then completes the proof of the
necessity.

Next, we prove the sufficiency. To
this end, let $f\in\mathcal{S}'(\rrn)$
satisfy that $(g_{\lambda}^{*})_{
\mathrm{loc}}(f)\in X$.
Observe that, for any given $t\in(0,\infty)$
and $x\in\rrn$, and for any $y\in B(x,t)$,
$$\frac{t}{t+|x-y|}\sim1.$$ This, together
with Definition \ref{df551},
further implies that, for any $x\in\rrn$,
\begin{align*}
&S_{\mathrm{loc}}(f)(x)\\
&\quad\sim\left|f\ast\varphi_0(x)\right|
+\left[\int_{0}^{1}
\int_{B(x,t)}
\left(\frac{t}{
t+|x-y|}\right)^{\lambda n}\left|
f\ast\varphi_t(y)\right|^2
\,\frac{dy\,dt}{t^{n+1}}\right]^{\frac{1}{2}}\\
&\quad\lesssim\left(g_{\lambda}
^*\right)_{\mathrm{loc}}(f)(x).
\end{align*}
From this and Definition \ref{Df1}(ii),
it follows that
\begin{equation*}
\left\|S_{\mathrm{loc}}(f)\right\|_X
\lesssim\left\|\left(
g_{\lambda}^*\right)_{\mathrm{loc}}(f)\right\|_X<\infty.
\end{equation*}
This, combined with Theorem \ref{lplx},
further implies that $f\in h_X(\rrn)$
and
$$
\left\|f\right\|_{h_X(\rrn)}\sim
\left\|S_{\mathrm{loc}}(f)\right\|_X
\lesssim\left\|\left(
g_{\lambda}^*\right)_{\mathrm{loc}}(f)\right\|_X,
$$
which completes the proof of the sufficiency.
By this and \eqref{lplglxe1},
we further conclude that
$$
\left\|f\right\|_{h_X(\rrn)}\sim
\left\|\left(
g_{\lambda}^*\right)_{\mathrm{loc}}(f)\right\|_X,
$$
where the positive equivalence
constants are independent of $f$, which then
completes the proof of Theorem \ref{lplglx}.
\end{proof}

Via the various Littlewood--Paley function
characterizations of $h_X(\rrn)$
obtained above, in
the remainder of this section,
we turn to establish
the Lusin area function, the
Littlewood--Paley $g$-function,
and the Littlewood--Paley $g_{\lambda}^*$-function
characterization of local generalized
Herz--Hardy spaces $\HaSaHol$ and
$\HaSaHl$. Indeed, we
have the following Littlewood--Paley
function characterizations of the
local generalized Herz--Hardy space
$\HaSaHol$.

\begin{theorem}\label{lusl}
Let $p,\ q\in(0,\infty)$, $\omega\in
M(\rp)$ satisfy
$\m0(\omega)\in(-\frac{n}{p},\infty)$ and
$\mi(\omega)\in(-\frac{n}{p},\infty)$,
$$
s_{0}:=\min\left\{1,p,q,\frac{n}{\Mw+n/p}\right\},
$$
and $\lambda\in(\max\{1,2/s_{0}\},\infty)$.
Then the following four statements are
mutually equivalent:
\begin{enumerate}
  \item[{\rm(i)}] $f\in\HaSaHol$;
  \item[{\rm(ii)}] $f\in\mathcal{S}'(\rrn)$
  and $S_{\mathrm{loc}}(f)\in\HerzSo$;
  \item[{\rm(iii)}] $f\in\mathcal{S}'(\rrn)$
  and $g_{\mathrm{loc}}(f)\in\HerzSo$;
  \item[{\rm(iv)}] $f\in\mathcal{S}'(\rrn)$
  and $(g_{\lambda}^*
  )_{\mathrm{loc}}(f)\in\HerzSo$.
\end{enumerate}
Moreover, for any $f\in\HaSaHol$,
\begin{align*}
\left\|f\right\|_{\HaSaHol}&\sim\left\|
S_{\mathrm{loc}}(f)\right\|_{\HerzSo}
\sim\left\|g_{\mathrm{loc}}(f)\right\|_{\HerzSo}
\\&\sim\left\|\left(g_{\lambda}^*\right)_{\mathrm{loc}}
(f)\right\|_{\HerzSo},
\end{align*}
where the positive equivalence
constants are independent of $f$.
\end{theorem}

\begin{proof}
Let all the symbols be as in the present
theorem, $s\in(\frac{2}{\lambda},s_0)$,
and $\theta\in(0,\min\{s,\frac{s^2}{2}\})$.
Then, by the proof of Theorem \ref{lusin},
we find that the following three
statements hold true:
\begin{enumerate}
  \item[{\rm(i)}] $\HerzSo$ is the BQBF space
  and $[\HerzSo]^{1/s}$ is a BBF space;
  \item[{\rm(ii)}] for any $\{f_j\}_{j\in\mathbb{N}}
  \subset L^1_{\mathrm{loc}}(\rrn)$,
  \begin{equation*}
  \left\|
\left\{
\sum_{j\in\mathbb{N}}
\left[\mc^{(\theta)}(f_j)\right]^{s}
\right\}^{1/s}
\right\|_{\HerzSo}
\lesssim\left\|
\left(\sum_{j\in\mathbb{N}}
|f_j|^{s}
\right)^{1/s}
\right\|_{\HerzSo}
  \end{equation*}
  and
  \begin{equation*}
\left\|
\left\{
\sum_{j\in\mathbb{N}}
\left[\mc^{(\theta)}(f_j)\right]^{s}
\right\}^{1/s}
\right\|_{[\HerzSo]^{s/2}}
\lesssim\left\|
\left(\sum_{j\in\mathbb{N}}
|f_j|^{s}
\right)^{1/s}
\right\|_{[\HerzSo]^{s/2}};
  \end{equation*}
  \item[{\rm(iii)}] for any $f\in L^1_{\mathrm
  {loc}}(\rrn)$,
\begin{equation*}
\left\|\mc^{((r/s)')}(f)\right\|_
{([\HerzSo]^{1/s})'}\lesssim
\left\|f\right\|_{([\HerzSo]^{1/s})'}.
\end{equation*}
\end{enumerate}
Combining these and Theorems \ref{lplx},
\ref{lplgx}, and \ref{lplglx}, we further
conclude that (i), (ii), (iii), and (iv)
are mutually equivalent and, for any
$f\in\HaSaHol$,
\begin{align*}
\left\|f\right\|_{\HaSaHol}&\sim\left\|
S_{\mathrm{loc}}(f)\right\|_{\HerzSo}
\sim\left\|g_{\mathrm{loc}}(f)\right\|_{\HerzSo}
\\&\sim\left\|\left(g_{\lambda}^*\right)_{\mathrm{loc}}
(f)\right\|_{\HerzSo}
\end{align*}
with the positive constants independent
of $f$.
This finishes the proof of Theorem \ref{lusl}.
\end{proof}

From Theorem \ref{lusl} and Remark \ref{remark5.0.9}(ii),
we immediately deduce the following
Littlewood--Paley function characterizations of
the local generalized Morrey--Hardy space $h\MorrSo$;
we omit the details.

\begin{corollary}\label{luslm}
Let $p,\ q\in[1,\infty)$, $\omega\in
M(\rp)$ satisfy
$$-\frac{n}{p}<
\m0(\omega)\leq\M0(\omega)<0$$ and
$$-\frac{n}{p}<
\mi(\omega)\leq\MI(\omega)<0,$$
and $\lambda\in(2,\infty)$.
Then the following four statements are
mutually equivalent:
\begin{enumerate}
  \item[{\rm(i)}] $f\in h\MorrSo$;
  \item[{\rm(ii)}] $f\in\mathcal{S}'(\rrn)$
  and $S_{\mathrm{loc}}(f)\in\MorrSo$;
  \item[{\rm(iii)}] $f\in\mathcal{S}'(\rrn)$
  and $g_{\mathrm{loc}}(f)\in\MorrSo$;
  \item[{\rm(iv)}] $f\in\mathcal{S}'(\rrn)$
  and $(g_{\lambda}^*)_{\mathrm{loc}}
  (f)\in\MorrSo$.
\end{enumerate}
Moreover, for any $f\in h\MorrSo$,
$$
\left\|f\right\|_{h\MorrSo}\sim\left\|
S_{\mathrm{loc}}(f)\right\|_{\MorrSo}
\sim\left\|g_{\mathrm{loc}}(f)\right\|_{\MorrSo}\sim
\left\|\left(g_{\lambda}^*\right)_{\mathrm{loc}}
(f)\right\|_{\MorrSo},
$$
where the positive equivalence
constants are independent of $f$.
\end{corollary}

On the other hand, the
following theorem gives
the Littlewood--Paley function characterizations
of the local generalized Herz--Hardy
space $\HaSaHl$.

\begin{theorem}\label{lusg}
Let $p,\ q\in(0,\infty)$, $\omega\in
M(\rp)$ satisfy
$\m0(\omega)\in(-\frac{n}{p},\infty)$ and
$$
-\frac{n}{p}<\mi(\omega)\leq\MI(\omega)<0,
$$
$$
s_{0}:=\min\left\{1,p,q,\frac{n}{\Mw+n/p}\right\},
$$
and $\lambda\in(\max\{1,2/s_{0}\},\infty)$.
Then the following four statements are
mutually equivalent:
\begin{enumerate}
  \item[{\rm(i)}] $f\in\HaSaHl$;
  \item[{\rm(ii)}] $f\in\mathcal{S}'(\rrn)$
  and $S_{\mathrm{loc}}(f)\in\HerzS$;
  \item[{\rm(iii)}] $f\in\mathcal{S}'(\rrn)$
  and $g_{\mathrm{loc}}(f)\in\HerzS$;
  \item[{\rm(iv)}] $f\in\mathcal{S}'(\rrn)$
  and $(g_{\lambda}^*
  )_{\mathrm{loc}}(f)\in\HerzS$.
\end{enumerate}
Moreover, for any $f\in\HaSaHl$,
\begin{align*}
\left\|f\right\|_{\HaSaHl}&\sim\left\|
S_{\mathrm{loc}}(f)\right\|_{\HerzS}
\sim\left\|g_{\mathrm{loc}}(f)\right\|_{\HerzS}\\&\sim
\left\|\left(g_{\lambda}^*\right)_{\mathrm{loc}}
(f)\right\|_{\HerzS},
\end{align*}
where the positive equivalence
constants are independent of $f$.
\end{theorem}

To show this theorem, we need the following
auxiliary lemma about translations,
whose proof is just to repeat the proof of Lemma
\ref{lusingl1}(iii) with
$\{S,g,g_{\lambda}^*\}$ therein replaced by
$\{S_{\mathrm{loc}},
g_{\mathrm{loc}},(g_{\lambda}^*
)_{\mathrm{loc}}\}$;
we omit the details.

\begin{lemma}\label{lusgl1}
Let $f\in\mathcal{S}'(\rrn)$, $\xi\in\rrn$,
and $\lambda\in(0,\infty)$.
Then, for any $A\in\{S_{\mathrm{loc}},
g_{\mathrm{loc}},(g_{\lambda}^*
)_{\mathrm{loc}}\}$,
$$
A\left(\tau_{\xi}(f)\right)
=\tau_{\xi}\left(A(f)\right).
$$
\end{lemma}

With the help of Lemma \ref{lusgl1},
we now show Theorem \ref{lusg}.

\begin{proof}[Proof of Theorem \ref{lusg}]
Let all the symbols be as in the present theorem.
Then, repeating an argument similar
to that used in the proof of Theorem \ref{lusing}
with Lemma \ref{lusingl1}(iii)
replaced by Lemma \ref{lusgl1}, we
find that (i), (ii), (iii), and (iv) are
mutually equivalent and, for any
$f\in\HaSaHl$,
\begin{align*}
\left\|f\right\|_{\HaSaHl}&\sim\left\|
S_{\mathrm{loc}}(f)\right\|_{\HerzS}
\sim\left\|g_{\mathrm{loc}}(f)\right\|_{\HerzS}\\&\sim
\left\|\left(g_{\lambda}^*\right)_{\mathrm{loc}}
(f)\right\|_{\HerzS}
\end{align*}
with the positive equivalence
constants independent of $f$.
This then finishes the proof of Theorem \ref{lusg}.
\end{proof}

As an application of Theorem \ref{lusg},
we have the Littlewood--Paley function
characterizations of the local generalized
Morrey--Hardy space
$h\MorrS$ as follows, which
can be deduced from
Theorem \ref{lusg}
and Remark \ref{remark5.0.9}(ii)
directly; we omit the details.

\begin{corollary}
Let $p$, $q$, $\omega$, and $\lambda$
be as in Corollary \ref{luslm}.
Then the following four statements are
mutually equivalent:
\begin{enumerate}
  \item[{\rm(i)}] $f\in h\MorrS$;
  \item[{\rm(ii)}] $f\in\mathcal{S}'(\rrn)$
  and $S_{\mathrm{loc}}(f)\in\MorrS$;
  \item[{\rm(iii)}] $f\in\mathcal{S}'(\rrn)$
  and $g_{\mathrm{loc}}(f)\in\MorrS$;
  \item[{\rm(iv)}] $f\in\mathcal{S}'(\rrn)$
  and $(g_{\lambda}^*)_{\mathrm{loc}}
  (f)\in\MorrS$.
\end{enumerate}
Moreover, for any $f\in h\MorrS$,
\begin{align*}
\left\|f\right\|_{h\MorrS}&\sim\left\|
S_{\mathrm{loc}}(f)\right\|_{\MorrS}
\sim\left\|g_{\mathrm{loc}}(f)\right\|_{\MorrS}\\&\sim
\left\|\left(g_{\lambda}^*\right)_{\mathrm{loc}}
(f)\right\|_{\MorrS},
\end{align*}
where the positive equivalence
constants are independent of $f$.
\end{corollary}

\section{Boundedness of Pseudo-Differential \\ Operators}

The target of this section is to show the boundedness
of pseudo-differential operators on the local generalized
Herz--Hardy spaces $\HaSaHol$ and $\HaSaHl$.
Recall that the \emph{H\"{o}rmander class}\index{H\"{o}rmander class}
$S_{1,0}^{0}(\rrn)$\index{$S_{1,0}^{0}(\rrn)$} is
defined to be the set of all the
infinitely differentiable functions
$\sigma$ on $\rrn\times\rrn$ such that,
for any $\alpha,\ \beta\in\mathbb{Z}_{+}^{n}$,
there exists a positive constant $C_{(\alpha,\beta)}$,
depending on both $\alpha$ and $\beta$,
such that, for any $x,\ \xi\in\rrn$,
\begin{equation}\label{hormander}
\left|\partial_{(1)}^{\alpha}
\partial_{(2)}^{\beta}\sigma(x,\xi)
\right|\leq C_{(\alpha,\beta)}(1+|\xi|)^{-|\beta|}.
\end{equation}
Via the H\"{o}rmander class $S_{1,0}^{0}(\rrn)$,
we present the definition
of pseudo-differential
operators as follows (see, for instance,
\cite[Subsection 4.5.4]{LGMF}).

\begin{definition}\label{df561}
Let $\sigma\in S_{1,0}^0(\rrn)$.
Then the
\emph{pseudo-differential operator}\index{pseudo-differential operator}
$T_\sigma$ is defined by setting,
for any $f\in\mathcal{S}(\rrn)$ and $x\in\rrn$,
$$
T_\sigma(f)(x):=\int_{\rrn}\sigma(x,\xi)e^{2\pi
ix\cdot\xi}\widehat{f}(\xi)\,d\xi,
$$
where $\sigma$ is called the
\emph{symbol} of $T_{\sigma}$.
\end{definition}

\begin{remark}
Let $T_{\sigma}$ be a
pseudo-differential operator
with symbol $\sigma\in S^0_{1,0}(\rrn)$.
Then, applying \cite[p.\,250,\ (43)]{s93},
we find that there exists a function
$K$ on $\rrn\times\rrn$ such that,
for any $f\in L^2(\rrn)$
with compact support, and for almost
every $x\notin\supp(f)$,
\begin{equation}\label{pseudore}
T_{\sigma}(f)(x)
=\int_{\rrn}K(x,x-y)f(y)\,dy.
\end{equation}
Moreover, from \cite[p.\,235,\ (9)]{s93},
it follows that the above
$K$ satisfies that, for any given
$M\in(0,\infty)$, there exists a positive
constant $C$ such that, for any $z\in\rrn$
with $|z|\geq1$,
\begin{equation}\label{pseudo9}
\left|K(x,z)\right|\leq C|z|^{-M}.
\end{equation}
The function $K$ satisfying both \eqref{pseudore}
and \eqref{pseudo9} is called the
\emph{kernel} of $T_{\sigma}$.
\end{remark}

Then the following theorem shows that
the pseudo-differential operator $T_\sigma$,
with $\sigma\in S_{1,0}^{0}(\rrn)$, is bounded
on the local generalized Herz--Hardy space $\HaSaHol$.

\begin{theorem}\label{posedo}
Let $p,\ q\in(0,\infty)$, $\omega\in M(\rp)$ satisfy
$\m0(\omega)\in(-\frac{n}{p},\infty)$ and
$\mi(\omega)\in(-\frac{n}{p},\infty)$,
and $T_\sigma$ be a pseudo-differential
operator with $\sigma\in S_{1,0}^{0}(\rrn)$.
Then $T_{\sigma}$ is well
defined on $\HaSaHol$ and
there exists a positive constant
$C$ such that, for any $f\in
\HaSaHol$,
$$
\left\|T_\sigma(f)\right\|_{\HaSaHol}
\leq C\|f\|_{\HaSaHol}.
$$
\end{theorem}

To prove Theorem \ref{posedo}, we
need the following boundedness
of pseudo-differential operators
on the local Hardy space
$h_X(\rrn)$, which is just
\cite[Theorem 4.5]{WYY}.

\begin{lemma}\label{posedol1}
Let $X$ be a ball quasi-Banach
function space satisfying both
Assumptions \ref{assfs} and
\ref{assas} with the same
$s\in(0,1]$, and
having an absolutely
continuous quasi-norm.
Assume that $T_\sigma$
is a pseudo-differential
operator with $\sigma\in S_{1,0}^{0}(\rrn)$.
Then $T_{\sigma}$
is well defined on $H_X(\rrn)$ and
there exists a positive constant
$C$ such that, for any $f\in H_X(\rrn)$,
$$
\left\|T_\sigma(f)\right\|_{
H_X(\rrn)}\leq C\|f\|_{H_X(\rrn)}.
$$
\end{lemma}

Via the above lemma, we now show
Theorem \ref{posedo}.

\begin{proof}[Proof of
Theorem \ref{posedo}]
Let all the symbols
be as in the present theorem.
Then, combining the assumption
$\m0(\omega)\in(-\frac{n}{p},
\infty)$ and both Theorems \ref{Th3}
and \ref{abso},
we conclude that the local generalized
Herz space $\HerzSo$ under consideration
is a BQBF space having
an absolutely continuous
quasi-norm. This, combined with
Lemma \ref{posedol1}, implies that,
to finish the proof of the present theorem,
we only need to show that
the Herz space $\HerzSo$ satisfies both
Assumptions \ref{assfs} and
\ref{assas} with the same $s\in(0,1]$.

To achieve this, let
$$
s\in\left(0,\min\left\{1,p,q,
\frac{n}{\Mw+n/p}\right\}\right)
$$
and $\theta\in(0,s)$.
Then, applying Lemma \ref{Atogl5},
we find that, for any $\{f_{j}\}_{
j\in\mathbb{N}}\subset L^{1}_{\text{loc}}(\rrn)$,
\begin{equation*}
\left\|\left\{\sum_{j\in\mathbb{N}}\left[\mc^{(\theta)}
(f_{j})\right]^s\right\}^{1/s}
\right\|_{\HerzSo}\lesssim
\left\|\left(\sum_{j\in\mathbb{N}}
|f_{j}|^{s}\right)^{1/s}\right\|_{\HerzSo},
\end{equation*}
which further implies that
Assumption \ref{assfs} holds true
for $\HerzSo$ with the above $\theta$ and $s$.

On the other hand, from Lemma \ref{mbhal}
with $r:=\infty$, it follows that
$[\HerzSo]^{1/s}$ is a BBF space and,
for any $f\in L_{{\rm loc}}^{1}(\rrn)$,
\begin{equation*}
\left\|\mc(f)\right\|_{([\HerzSo]^{1/s})'}\lesssim
\left\|f\right\|_{([\HerzSo]^{1/s})'}.
\end{equation*}
This implies that
the Herz space $\HerzSo$ under consideration
satisfies Assumption \ref{assas}
with the above $s$ and $r:=\infty$.
Therefore, both Assumptions \ref{assfs} and
\ref{assas} hold true for $\HerzSo$
with the same $s$.
Combining this, the assumptions
that $\HerzSo$ under consideration
is a BQBF space having
an absolutely continuous
quasi-norm, and Lemma \ref{posedo}
with $X:=\HerzSo$,
we further conclude that $T_{\sigma}$ is well
defined on $\HaSaHol$ and, for any $f\in
\HaSaHol$,
$$
\left\|T_\sigma(f)\right\|_{\HaSaHol}
\lesssim\|f\|_{\HaSaHol},
$$
which completes the proof of Theorem
\ref{posedo}.
\end{proof}

\begin{remark}
We should point out that,
in Theorem \ref{posedo}, when $p\in(1,\infty)$ and
$\omega(t):=t^\alpha$ for any $t\in(0,\infty)$
and for any given $\alpha\in[n(1-\frac{1}{p}),
\infty)$, then
Theorem \ref{posedo} goes
back to \cite[Theorem 2.6]{wl15}.
\end{remark}

Via Theorem \ref{posedo}
and Remark \ref{remark5.0.9}(ii),
we immediately obtain the boundedness of
the pseudo-differential operators $T_\sigma$
on the local generalized Morrey--Hardy space $h\MorrSo$ as follows;
we omit the details.

\begin{corollary}\label{posedom}
Let $p,\ q\in[1,\infty)$, $\omega\in M(\rp)$ with
$$
-\frac{n}{p}<\m0(\omega)\leq\M0(\omega)<0
$$
and
$$
-\frac{n}{p}<\mi(\omega)\leq\MI(\omega)<0,
$$
and $T_\sigma$ be a pseudo-differential
operator with $\sigma\in S_{1,0}^{0}(\rrn)$.
Then $T_{\sigma}$ is well defined on
$h\MorrSo$ and
there exists a positive
constant $C$ such that, for any $f\in
h\MorrSo$,
$$
\left\|T_\sigma(f)\right\|_{h\MorrSo}
\leq C\|f\|_{h\MorrSo}.
$$
\end{corollary}

We next turn to establish
the boundedness of pseudo-differential
operators on the local generalized
Herz--Hardy space $\HaSaHl$.
Namely, we have the following conclusion.

\begin{theorem}\label{pseudog}
Let $p,\ q\in(0,\infty)$, $\omega\in M(\rp)$ satisfy
$\m0(\omega)\in(-\frac{n}{p},\infty)$ and
$$-\frac{n}{p}
<\mi(\omega)\leq\MI(\omega)<0,$$
and $T_\sigma$ be a pseudo-differential
operator with symbol $\sigma\in S_{1,0}^{0}(\rrn)$.
Then $T_{\sigma}$ is well
defined on $\HaSaHl$ and
there exists a positive constant
$C$ such that, for any $f\in
\HaSaHl$,
$$
\left\|T_\sigma(f)\right\|_
{\HaSaHl}\leq C\|f\|_{\HaSaHl}.
$$
\end{theorem}

Due to the deficiency of associate spaces
and the absolutely continuity of
quasi-norms of global generalized Herz spaces,
we can not show Theorem \ref{pseudog}
using Lemma \ref{posedol1}
directly. To overcome these difficulties,
we first establish a new
boundedness criterion of pseudo-differential
operators on the local
Hardy space $h_X(\rrn)$ associated with
the ball quasi-Banach function space
$X$ as follows.

\begin{theorem}\label{pseudox}
Let $X$ be a ball quasi-Banach function
space, $Y$ a linear space equipped with a
quasi-seminorm $\|\cdot\|_Y$, and
$Y_0$ a linear space equipped with
a quasi-seminorm $\|\cdot\|_{Y_0}$, and let
$\eta\in(1,\infty)$
and $0<\theta<s<s_0\leq1$ be such that
\begin{enumerate}
  \item[{\rm(i)}] for the above $\theta$ and
  $s$, Assumption \ref{assfs} holds true;
  \item[{\rm(ii)}] both $\|\cdot\|_{Y}$
  and $\|\cdot\|_{Y_0}$
  satisfy Definition \ref{Df1}(ii);
  \item[{\rm(iii)}] $\1bf_{B(\0bf,1)}
  \in Y_0$;
  \item[{\rm(iv)}] $X^{1/s}$ is a ball Banach
  function space and, for any $f\in\Msc(\rrn)$,
  $$
  \|f\|_{X^{1/s}}
  \sim\sup\left\{\|fg\|_{L^1(\rrn)}:\
  \|g\|_{Y}=1\right\}
  $$
  and
  $$
  \|f\|_{X^{1/s_{0}}}
  \sim\sup\left\{\|fg\|_{L^1(\rrn)}:\
  \|g\|_{Y_0}=1\right\}
  $$
  with the positive equivalence constants
  independent of $f$;
  \item[{\rm(v)}] $\mc^{(\eta)}$ is bounded
  on $Y$ and $Y_0$.
\end{enumerate}
Assume that $T_{\sigma}$ is a
pseudo-differential operator
with symbol $\sigma\in S^0_{1,0}(\rrn)$.
Then $T_{\sigma}$ is well defined
on $h_X(\rrn)$ and there exists a positive
constant $C$ such that,
for any $f\in h_X(\rrn)$,
$$
\left\|T_{\sigma}(f)\right\|_{h_X(\rrn)}
\leq C\left\|f\right\|_{h_X(\rrn)}.
$$
\end{theorem}

To prove this theorem, we require
some auxiliary conclusions. First,
the following lemma shows
that the local Hardy space $h_X(\rrn)$
can be continuously embedded in
$\mathcal{S}'(\rrn)$, which
plays a vital role in the proof of
Theorem \ref{pseudox}.

\begin{lemma}\label{pseudoxl1}
Let $X$ be a ball quasi-Banach
function space and $N\in\mathbb{N}$.
Then the local Hardy space
$h_X(\rrn)$ embeds continuously
into $\mathcal{S}'(\rrn)$. Namely,
there exists a positive constant
$C$ such that, for any $f\in h_X(\rrn)$
and $\phi\in\mathcal{S}(\rrn)$,
$$
|\langle f,\phi\rangle|\leq
Cp_N(\phi)\|f\|_{h_X(\rrn)},
$$
where $p_N$ is defined as in
\eqref{defp}.
\end{lemma}

\begin{proof}
Let all the symbols be as in
the present lemma, $f\in h_X(\rrn)$,
and $\phi\in\mathcal{S}(\rrn)$.
We define $\psi$ by setting, for any
$x\in(\rrn)$,
$$
\psi(x):=\frac{1}{2^n}\phi\left(
\frac{x}{2}\right).
$$
Now, we show the present lemma
by considering the following two
cases on $\psi$.

\emph{Case 1)} $p_N(\psi)=0$.
In this case, applying \ref{defp} with
$\phi$ therein replaced by $\psi$, we find
that $\psi=0$. This further implies that
$\phi=0$. Thus, we have
\begin{equation}\label{pseudoxl1e1}
\left|\left\langle
f,\phi\right\rangle\right|=0=p_N(\phi)\left\|
f\right\|_{h_X(\rrn)},
\end{equation}
which completes the proof of the present
lemma in this case.

\emph{Case 2)} $p_N(\psi)\neq0$. In this case,
we have $\frac{\psi(-\cdot)}{p_N(\psi)}\in\mathcal{S}
(\rrn)$,
\begin{equation}\label{pseudoxl1e2}
\left(\frac{\psi(-\cdot)}{p_N(\psi)}\right)_{
\frac{1}{2}}=\frac{\phi(-\cdot)}{p_N(\psi)},
\end{equation}
and
\begin{equation}\label{pseudoxl1e3}
p_N\left(\frac{\psi(-\cdot)}
{p_N(\psi)}\right)=\frac{
p_N(\psi)}{p_N(\psi)}=1,
\end{equation}
where $(\frac{\psi(-\cdot)}{p_N(\psi)})_{\frac{1}{2}}$
is defined as in \eqref{smaxe1} with
$\phi$ and $t$ therein replaced,
respectively, by $\frac{\psi(-\cdot)}{p_N(\psi)}$
and $\frac{1}{2}$. From
\eqref{pseudoxl1e3}, it follows that
$\frac{\psi(-\cdot)}{p_N(\psi)}\in\mathcal{F}_N(\rrn)$.
This, combined with \eqref{pseudoxl1e2} and
\eqref{sec7e1}, further implies that,
for any $x\in B(\0bf,1)$,
\begin{align}\label{pseudoxl1e4}
\left|\left\langle
f,\phi\right\rangle\right|&=
\left|f\ast\left[\phi(-\cdot)\right](\0bf)\right|
\notag\\
&=p_N(\psi)\left|f\ast\left[
\left(\frac{\psi(-\cdot)}{p_N(\psi)}\right)_{
\frac{1}{2}}\right](\0bf)\right|\notag\\
&\leq p_{N}(\psi)m_N(f)(x).
\end{align}
In addition, notice that
\begin{align*}
p_N(\psi)&=\sum_{\alpha\in\zp^n,\,|\alpha
|\leq N}\sup_{x\in\rrn}\left(1+|x|\right)^N
\left|\partial^{\alpha}\left(
\frac{1}{2^n}\phi\left(\frac{\cdot}{2}
\right)\right)(x)\right|\\
&\sim\sum_{\alpha\in\zp^n,\,|\alpha
|\leq N}\sup_{x\in\rrn}\left(1+|x|\right)^N
\left|\partial^{\alpha}\phi\left(\frac{x}{2}
\right)\right|\\
&\sim\sum_{\alpha\in\zp^n,\,|\alpha
|\leq N}\sup_{x\in\rrn}\left(1+|x|\right)^N
\left|\partial^{\alpha}\phi(x)\right|\sim
p_N(\phi).
\end{align*}
By this and \eqref{pseudoxl1e4}, we conclude
that
\begin{equation*}
\left|\left\langle
f,\phi\right\rangle\right|\1bf_{B(\0bf,1)}
\lesssim p_N(\phi)m_N(f).
\end{equation*}
Combining this, an argument similar to that
used in the estimation of \eqref{bddol2e3}
with $\mc_N(f)$ therein replaced by
$m_N(f)$, and Definition \ref{hardyxl},
we further obtain
\begin{equation}\label{pseudoxl1e6}
\left|\left\langle
f,\phi\right\rangle\right|
\lesssim p_N(\phi)\left\|f\right\|_{h_X(\rrn)},
\end{equation}
which completes the proof of
the present lemma in this case.
Thus, by \eqref{pseudoxl1e1} and
\eqref{pseudoxl1e6}, we then
complete the proof of Lemma \ref{pseudoxl1}.
\end{proof}

To prove Theorem \ref{pseudox},
we also need the following
technique estimate about
local atoms.

\begin{proposition}\label{pseudoae}
Let $T_{\sigma}$ be a pseudo-differential
operator with $\sigma\in S^0_{1,0}(\rrn)$,
$X$ a ball quasi-Banach function space,
$\phi\in\mathcal{S}(\rrn)$, $\theta\in(0,1]$,
$d\geq n(\frac{1}{\theta}-1)-1$
be a fixed integer, and $r\in[1,\infty]$.
Then there exists a positive constant $C$
such that, for any
local-$(X,\,r,\,d)$-atom $a$ supported
in the ball $B\in\mathbb{B}$,
$$
m\left(T_{\sigma}(a),\phi\right)\1bf_{(2B)^{
\complement}}\leq C\frac{1}{\|\1bf_B\|_X}
\mc^{(\theta)}(\1bf_B).
$$
\end{proposition}

In order to show this estimate,
we need the following auxiliary lemma
about symbols and kernels of pseudo-differential
operators, which is just
\cite[Lemma 6]{Gold}.

\begin{lemma}\label{pseudot}
Let $T_{\sigma}$ be a pseudo-differential
operator with symbol $\sigma\in
S^0_{1,0}(\rrn)$,
$\phi\in\mathcal{S}(\rrn)$, and $t\in(0,1)$.
The operator $T_{\sigma}^{(t)}$ is defined
by setting, for any $f\in\mathcal{S}(\rrn)$,
$$T_{\sigma}^{(t)}(f):=
\left[T(f)\right]\ast\phi_t.$$
Then $T_{\sigma}^{(t)}$
is a pseudo-differential operator
with symbol $\sigma_t$ and kernel $K_t$
satisfying that, for any
$\alpha$, $\beta\in\zp^n$,
there exists a positive constant $C_{(\alpha,\beta)}$,
depending on $\alpha$ and $\beta$ but
independent of $t$, such that, for any $x$,
$\xi\in\rrn$,
\begin{equation*}
\left|\partial^{\alpha}_{(1)}\partial
^{\beta}_{(2)}\sigma_t\left(x,\xi\right)\right|
\leq C_{(\alpha,\beta)}\left(
1+|\xi|\right)^{-|\beta|}
\end{equation*}
and, for any $z\in\rrn\setminus\{\0bf\}$,
\begin{equation}\label{pseudoe2}
\left|\partial^{\alpha}_{(1)}\partial
^{\beta}_{(2)}K_t\left(x,z\right)\right|
\leq C_{(\alpha,\beta)}|z|^{-n-|\beta|}.
\end{equation}
\end{lemma}

Applying the above lemma, we next
prove Proposition \ref{pseudoae}.

\begin{proof}[Proof of Proposition
\ref{pseudoae}]
Let all the symbols be as in the
present proposition and $a$
a local-$(X,\,r,\,d)$-atom supported
in the ball $B(x_0,r_0)$ with $x_0\in\rrn$
and $r_0\in(0,\infty)$.
Then we claim that, for any
$t\in(0,1)$,
\begin{equation}\label{pseudoaee1}
\left|T_{\sigma}(a)\ast\phi_t\right|
\1bf_{[B(x_0,2r_0)]^{\complement}}
\lesssim\frac{r_0^{d+1}}{|x-x_0|^{n+d+1}}
\left\|a\right\|_{L^1(\rrn)},
\end{equation}
where the implicit positive
constant is independent of $a$.
Assume that this claim holds true for the
moment. Then, by Definition \ref{df511}(i),
\eqref{pseudoaee1}, and an argument
similar to that used in the estimation
of \eqref{atogxl1e1}, we conclude that
$$
m\left(T_{\sigma}(a),\phi\right)\1bf_{[B(x_0,
2r_0)]^{\complement}}
\lesssim\frac{1}{\|\1bf_(B(x_0,r_0))\|_X}
\mc^{(\theta)}\left(\1bf_{B(x_0,r_0)}\right)
$$
with the implicit positive constant independent
of $a$, which completes the proof of the
present proposition. Therefore,
to complete the whole proof, it suffices
to show the above claim. To achieve this, we
now prove the above claim by considering the
following two cases on $r_0$.

\emph{Case 1)} $r_0\in(0,1)$.
In this case, for any $t\in(0,\infty)$
and $x\in[B(x_0,2r_0)]^{\complement}$,
from Lemma \ref{pseudot}, \eqref{pseudore}
with $T_{\sigma}$, $K$, and
$f$ therein replaced, respectively,
by $T_{\sigma}^{(t)}$, $K_t$,
and $a$, Definition
\ref{loatomx}(iii), and
the Taylor remainder theorem,
it follows that, for any
$y\in B(x_0,r_0)$, there exists a
$t_y\in(0,1)$ such that
\begin{align}\label{pseudoaee2}
&\left[T_{\sigma}(a)\right]\ast\phi_t(x)
\notag\\&\quad=T_{\sigma}^{(t)}(a)(x)
=\int_{B(x_0,r_0)}
K_t(x,x-y)a(y)\,dy\notag\\
&\quad=\int_{B(x_0,r_0)}\left[
K_t(x,x-y)-\sum_{\gfz{\gamma\in\zp^n}{|\gamma|\leq d}}\frac{\partial^{\gamma}_{
(2)}K(x,x-x_0)}{\gamma!}(y-x_0)^{\gamma}
\right]a(y)\,dy\notag\\
&\quad=\int_{B(x_0,r_0)}\sum_{\gfz{\gamma\in\zp^n}
{|\gamma|=d+1}}\frac{\partial^{\gamma}_{
(2)}K_t(x,x-t_yy-(1-t_y)x_0)}{\gamma!}a(y)\,dy.
\end{align}
In addition, notice that,
for any $x\in[B(x_0,2r_0)]^{
\complement}$ and $y\in B(x_0,r_0)$,
$$
|x-x_0|\geq2r_0>2|y-x_0|,
$$
which further implies that
\begin{align*}
\left|x-t_yy-\left(1-t_y\right)x_0\right|
&=\left|x-x_0-t_y\left(y-x_0\right)\right|\\
&\geq\left|x-x_0\right|-\left|y-x_0\right|
>\frac{1}{2}\left|x-x_0\right|.
\end{align*}
Combining this, \eqref{pseudoaee2},
and \eqref{pseudoe2} with
$\alpha$, $\beta$, and $z$ therein
replaced, respectively, by $\0bf$,
$\gamma$ for any $\gamma\in\zp^n$
satisfying that $|\gamma|=d+1$, and
$x-t_yy-(1-t_y)x_0$ for any
$x\in[B(x_0,2r_0)]^{\complement}$
and $y\in B(x_0,r_0)$,
we find that, for any $t\in(0,1)$
and $x\in[B(x_0,2r_0)]^{\complement}$,
\begin{align*}
\left|T_{\sigma}(a)\ast\phi_t(x)\right|
&\lesssim\int_{B(x_0,r_0)}
\frac{|y-x_0|^{d+1}}{|x-t_yy-(1-t_y)
x_0|^{n+d+1}}|a(y)|\,dy\\
&\lesssim\frac{r_0^{d+1}}{|x-x_0|^{n+d+1}}
\left\|a\right\|_{L^1(\rrn)},
\end{align*}
where the implicit positive constant
is independent of $a$. This finishes
the proof of the above claim in this case.

\emph{Case 2)} $r_0\in[1,\infty)$. In this case,
for any $y\in B(x_0,r_0)$, we have
$$
\left|x-y\right|
\geq\left|x-x_0\right|-\left|y-y_0\right|
>\frac{1}{2}\left|x-x_0\right|
\geq r_0\geq1.
$$
Using this, \eqref{pseudore} with
$T_{\sigma}$, $K$, and $f$ therein
replaced, respectively, by
$T_{\sigma}^{(t)}$ and
$K_t$ for any $t\in(0,1)$, and $a$,
Lemma \ref{pseudot}, and \eqref{pseudo9} with
$K$, $z$, and $M$ therein replaced, respectively,
by $K_t$ for any $t\in(0,1)$, $x-y$, and $n+d+1$,
we conclude that, for any $t\in(0,1)$ and
for almost every $x\in[B(x_0,2r_0)]^{\complement}$,
\begin{align}\label{pseudoaee4}
&\left|T_{\sigma}(a)\ast\phi_t(x)\right|\notag\\
&\quad=\left|T_{\sigma}^{(t)}(a)(x)\right|
\leq\int_{B(x_0,r_0)}\left|K_t(x,x-y)\right|
\left|a(y)\right|\,dy\notag\\
&\quad\lesssim\int_{B(x_0,r_0)}\frac{|a(y)|}{
|x-y|^{n+d+1}}\,dy
\lesssim\frac{r_0^{d+1}}{|x-x_0|^{n+d+1}}
\left\|a\right\|_{L^1(\rrn)},
\end{align}
where the implicit positive constants
are independent of $t$ and $a$.
This then finishes the proof of the above
claim in this case.
Therefore, combining \eqref{pseudoaee2} and
\eqref{pseudoaee4}, we further
find that \eqref{pseudoaee1} holds
true and hence complete the proof
of the above claim.
This further finishes the proof of Proposition
\ref{pseudoae}.
\end{proof}

Furthermore, the following
proposition gives the atomic decomposition
of the local Hardy space $h_X(\rrn)$
with convergence in $h_X(\rrn)$,
which is an essential tool
in the proof of Theorem \ref{pseudox}.

\begin{proposition}\label{atoma}
Let $X$ be a ball quasi-Banach
function space satisfying
both Assumption \ref{assfs}
with $0<\theta<s\leq1$ and
Assumption \ref{assas} with
the same $s$, $d\geq
\lfloor n(1/\theta-1)\rfloor$
be a fixed integer,
and $r\in(1,\infty]$.
Assume that $X$ has an absolutely
continuous quasi-norm.
Let $\{a_j\}_{j\in\mathbb{N}}$
be a sequence of local-$(X,
\,r,\,d)$-atoms supported, respectively,
in the balls $\{B_j\}_{j\in\mathbb{N}}
\subset\mathbb{B}$ and $\{\lambda_j\}
_{j\in\mathbb{N}}\subset[0,\infty)$
such that
$$f:=\sum_{j\in\mathbb{N}}
\lambda_ja_j$$
in $\mathcal{S}'(\rrn)$,
and
\begin{equation}\label{atomae1}
\left\|\left[
\sum_{j\in\mathbb{N}}
\left(\frac{\lambda_j}{\|\1bf_{B_j}\|_{X}}
\right)^s\1bf_{B_j}
\right]^{\frac{1}{s}}
\right\|_{X}<\infty.
\end{equation}
Then $f\in h_X(\rrn)$
and $f=\sum_{j\in\mathbb{N}}\lambda_j
a_j$ holds true in $h_X(\rrn)$.
\end{proposition}

In order to prove this proposition,
we first present the following
dominated convergence theorem of
ball quasi-Banach function spaces,
which was obtained in \cite[Lemma 6.3]{yhyy22}
(see also \cite[Chapter 1,\ Proposition 3.6]{BSIO}).

\begin{lemma}\label{absox}
Let $X$ be a ball quasi-Banach function
space having an absolutely
continuous quasi-norm.
Assume that $g\in X$ and
$\{f_k\}_{k\in\mathbb{N}}$
is a sequence of measurable functions
satisfying that
$|f_k|\leq|g|$ for any $k\in\mathbb{N}$
and $\lim_{k\to\infty}f_k=f$
almost everywhere in
$\rrn$. Then
\begin{equation*}
\lim\limits_{
k\to\infty}\left\|f_k-f\right\|_{X}=0.
\end{equation*}
\end{lemma}

We now prove Proposition \ref{atoma}.

\begin{proof}[Proof of Proposition
\ref{atoma}]
Let all the symbols be as in
the present theorem. Then,
using both Definition \ref{atom-hardyl}
and Lemma \ref{atomxl},
we find that $f\in h_X(\rrn)$ and
$f=\sum_{j\in\mathbb{N}}\lambda_ja_j$
holds true in $\mathcal{S}'(\rrn)$.
Next, we show that $f=\sum_{
j\in\mathbb{N}}\lambda_ja_j$
also holds true in $h_X(\rrn)$.
To this end, for any $k\in\mathbb{N}$,
let $$f_k:=\sum_{j=1}^{k}
\lambda_ja_j.$$ Then, in order to
show that $f=\sum_{
j\in\mathbb{N}}\lambda_ja_j$
holds true in $h_X(\rrn)$,
it suffices to prove
that
\begin{equation}\label{atomae2}
\lim\limits_{k\to\infty}
\left\|f_k-f\right\|_{h_X(\rrn)}=0.
\end{equation}
Indeed, for any $k\in\mathbb{N}$,
we have
$$
f_k-f=\sum_{j=k+1}^{\infty}
\lambda_ja_j
$$
in $\mathcal{S}'(\rrn)$.
Applying this, Lemma \ref{atomxl},
Definition \ref{atom-hardyl},
\eqref{atomae1},
and Lemma \ref{absox}
with $\{f_k\}_{k\in\mathbb{N}}$,
$f$, and $g$ therein replaced,
respectively, by
$$\left\{\left[
\sum_{j=1}^k\left(
\frac{\lambda_j}{\|\1bf_{B_j}\|_X}
\right)^s\1bf_{B_j}\right]^{\frac{1}{s}}
\right\}_{k\in\mathbb{N}},$$
$$
\left[
\sum_{j\in\mathbb{N}}\left(
\frac{\lambda_j}{\|\1bf_{B_j}\|_X}
\right)^s\1bf_{B_j}\right]^{\frac{1}{s}},
$$
and
$$
\left[
\sum_{j\in\mathbb{N}}\left(
\frac{\lambda_j}{\|\1bf_{B_j}\|_X}
\right)^s\1bf_{B_j}\right]^{\frac{1}{s}},
$$
we conclude that
\begin{align*}
\left\|f_k-f\right\|_{h_X(\rrn)}
&\sim\left\|f_k-f\right\|_{h
^{X,r,d,s}(\rrn)}\\
&\lesssim\left\|
\left[\sum_{j=k+1}^{\infty}\left(
\frac{\lambda_j}{\|\1bf_{B_j}\|_X}
\right)^s\1bf_{B_j}\right]^{\frac{1}{s}}
\right\|_{X}\to0
\end{align*}
as $k\to\infty$. This finishes the proof
of \eqref{atomae2}, and further implies that
$f=\sum_{j\in\mathbb{N}}
\lambda_ja_j$ holds true in
$h_X(\rrn)$, which then completes the
proof of Proposition \ref{atoma}.
\end{proof}

In addition, to show Theorem \ref{pseudox},
we also need the boundedness
of pseudo-differential operators
on Lebesgue spaces and the
localized weighted Hardy spaces.
Recall that the
\emph{local weighted Hardy space}
$h_{\upsilon}^p(\rrn)$\index{$h_{\upsilon}^p(\rrn)$},
with $p\in(0,\infty)$ and $\upsilon\in A_{\infty}(\rrn)$,
is defined as in Definition \ref{hardyxl}
with $X:=L^p_{\upsilon}(\rrn)$.
Then we have the following
two lemmas which can be found,
respectively, in
\cite[p.\,250,\ Proposition 4]{s93}
and \cite[Corollary 4.14(b)]{WYY}.

\begin{lemma}\label{bdd-pseudo}
Let $p\in(1,\infty)$ and
$T_{\sigma}$ be a pseudo-differential
operator with $\sigma\in S^0_{1,0}(\rrn)$.
Then $T_{\sigma}$ is well defined on
$L^p(\rrn)$ and there exists a positive
constant $C$ such that, for any $f\in
L^p(\rrn)$,
$$
\left\|T_\sigma(f)\right\|_{L^p(\rrn)}
\leq C\|f\|_{L^p(\rrn)}.
$$
\end{lemma}

\begin{lemma}\label{pseudow}
Let $p\in(0,1]$, $\upsilon\in
A_{\infty}(\rrn)$,
and $T_\sigma$ be a pseudo-differential
operator with symbol $\sigma\in S_{1,0}^{0}(\rrn)$.
Then $T_{\sigma}$ is well
defined on $h^p_{\upsilon}(\rrn)$ and
there exists a positive constant
$C$ such that, for any $f\in
h^p_{\upsilon}(\rrn)$,
$$
\left\|T_\sigma(f)\right\|_
{h^p_{\upsilon}(\rrn)}\leq C\|f\|_{
h^p_{\upsilon}(\rrn)}.
$$
\end{lemma}

Via above preparations, we next
prove Theorem \ref{pseudox}.

\begin{proof}[Proof of Theorem \ref{pseudox}]
Let all the symbols be as in
the present theorem, $f\in
h_X(\rrn)$, and $d\geq\lfloor n(
1/\theta-1)\rfloor$ be a fixed integer.
Then, by the assumption (i) of the present
theorem and Lemma \ref{atomxxll1},
we find that there exist
$\{\lambda_j\}_{j\in\mathbb{N}}
\subset[0,\infty)$ and $\{a_j\}_{j\in
\mathbb{N}}$ of local-$(X,\,\infty,\,d)$-atoms
supported, respectively, in
the balls $\{B_j\}_{j\in\mathbb{N}}\subset
\mathbb{B}$ such that
\begin{equation}\label{pseudoxe1}
f=\sum_{j\in\mathbb{N}}\lambda_ja_j
\end{equation}
in $\mathcal{S}'(\rrn)$, and
\begin{equation}\label{pseudoxe2}
\left\|\left[\sum_{j\in\mathbb
N}\left(\frac{\lambda_j}{\|\1bf_{B_j}\|_{X}}
\right)^s\1bf_{B_j}
\right]^{\frac{1}{s}}\right\|_{X}\lesssim
\|f\|_{h_X(\rrn)}.
\end{equation}

In addition, from the assumptions
(ii) through (v) of the present theorem
and Lemma \ref{bddol1} with
$Y:=Y_0$, $\theta:=\eta$, and $s:=s_0$,
we deduce that there exists an $\eps\in(0,1)$
such that, for any $g\in\Msc(\rrn)$,
\begin{equation*}
\|g\|_{L^{s_0}_{\upsilon}(\rrn)}
\lesssim\|g\|_X,
\end{equation*}
where $\upsilon:=[\mc(\1bf_{B(\0bf,1)})]
^{\eps}$. This, combined with
both \eqref{pseudoxe1}
and \eqref{pseudoxe2}, further implies that
\begin{equation}\label{pseudoxe3}
f=\sum_{j\in\mathbb{N}}\lambda_ja_j
=\sum_{j\in\mathbb{N}}\left[\lambda_j
\frac{\|\1bf_{B_j}\|_{L^{s_0}_{\upsilon}(\rrn)}}
{\|\1bf_{B_j}\|_X}\right]\left[\frac
{\|\1bf_{B_j}\|_X}{\|\1bf_{B_j}\|_{L^{s_0}
_{\upsilon}(\rrn)}}a_j\right]
\end{equation}
in $\mathcal{S}'(\rrn)$, and
\begin{align}\label{pseudoxe4}
&\left\|\left\{
\sum_{j\in\mathbb{N}}\left[
\frac{\lambda_j\frac{\|\1bf_{B_j}
\|_{L^{s_0}_{\upsilon}(\rrn)}}
{\|\1bf_{B_j}\|_X}}{\|\1bf_{B_j}\|
_{L^{s_0}_{\upsilon}(\rrn)}}\right]^s
\1bf_{B_j}\right\}^{\frac{1}{s}}
\right\|_{L^{s_0}_{\upsilon}(\rrn)}\notag\\
&\quad\lesssim\left\|\left[\sum_{j\in\mathbb
N}\left(\frac{\lambda_j}{\|\1bf_{B_j}\|_{X}}
\right)^s\1bf_{B_j}
\right]^{\frac{1}{s}}\right\|_{X}\lesssim
\|f\|_{h_X(\rrn)}<\infty.
\end{align}
Notice that, for any
$j\in\mathbb{N}$, applying
Definition \ref{loatomx} with
$X:=L^{s_0}_{\upsilon}(\rrn)$,
we conclude that
$\frac
{\|\1bf_{B_j}\|_X}{\|\1bf_{B_j}\|_{L^{s_0}
_{\upsilon}(\rrn)}}a_j$ is a
local-$(L^{s_0}_{\upsilon}(\rrn),\,\infty,\,d)$-atom
supported in $B_j$.

Now, we prove that $f\in h_{\upsilon}^{s_0}(\rrn)$
and $f=\sum_{j\in\mathbb{N}}\lambda_j
a_j$ holds true in $h_{\upsilon}^{s_0}(\rrn)$.
Indeed, using \cite[Remarks 2.4(b),
2.7(b), and 3.4(i)]{WYY}, we find that
the following four statements hold true:
\begin{enumerate}
  \item[{\rm(i)}] $L^{s_0}_{
  \upsilon}(\rrn)$ is a BQBF space;
  \item[{\rm(ii)}] for any
$\{f_{j}\}_{j\in\mathbb{N}}
\subset L^1_{{\rm loc}}(\rrn)$,
\begin{equation*}
\left\|\left\{\sum_{j\in\mathbb{N}}
\left[\mc^{(\theta)}
(f_{j})\right]^s\right\}^{1/s}
\right\|_{L^{s_0}_{
\upsilon}(\rrn)}\lesssim
\left\|\left(\sum_{j\in\mathbb{N}}
|f_{j}|^{s}\right)^{1/s}\right\|_{
L^{s_0}_{\upsilon}(\rrn)};
\end{equation*}
  \item[{\rm(iii)}] $[L^{s_0}_{\upsilon}(\rrn)]
  ^{1/s}$ is a BBF space and,
for any $f\in L^1_{\mathrm{loc}}(\rrn)$,
\begin{equation*}
\left\|\mc(f)
\right\|_{([L^{s_0}_{\upsilon}(\rrn)]
^{1/s})'}\lesssim
\left\|f\right\|_{([L^{s_0}_{\upsilon}(\rrn)]
^{1/s})'};
\end{equation*}
  \item[{\rm(iv)}] $L^{s_0}_{\upsilon}(\rrn)$
  has an absolutely continuous quasi-norm.
\end{enumerate}
These, combined with Proposition \ref{atoma}
with $X:=L^{s_0}_{\upsilon}(\rrn)$
and $r:=\infty$ , \eqref{pseudoxe3}, and
\eqref{pseudoxe4}, further imply that
$f\in h^{s_0}_{\upsilon}
(\rrn)$ and $f=\sum_{j\in\mathbb{N}}
\lambda_ja_j$ holds true in
$h^{s_0}_{\upsilon}(\rrn)$.
From this and Lemma \ref{pseudow}
with $p:=s_0$, it follows that
$T_{\sigma}(f)=\sum_{j\in\mathbb{N}}
\lambda_jT_{\sigma}(a_j)$ holds true in $h_{\upsilon}
^{s_0}(\rrn)$.
This, together with Lemma \ref{pseudoxl1}
with $X:=L^{s_0}_{
\upsilon}(\rrn)$, further implies that
\begin{equation}\label{pseudoxe5}
T_{\sigma}(f)=\sum_{j\in\mathbb{N}}
\lambda_jT_{\sigma}(a_j)
\end{equation}
in $\mathcal{S}'(\rrn)$.

Now, we show that $T_{\sigma}(f)\in h_X(\rrn)$
and
$$
\left\|T_{\sigma}(f)\right\|_{h_X(\rrn)}
\lesssim\left\|f\right\|_{h_X(\rrn)}.
$$
To achieve this, fix a $\phi\in\mathcal{S}(\rrn)$
satisfying that $\supp(\phi)\subset
B(\0bf,1)$ and $$\int_{\rrn}\phi(x)\,dx\neq0.$$
Then, by \eqref{pseudoxe5} and
repeating an argument similar to that
used in the estimation of \ref{Atogxe10}
with $f$ and $\{a_j\}_{j\in\mathbb{N}}$
therein replaced, respectively, by
$T_{\sigma}(f)$ and $\{T_{\sigma}
(a_j)\}_{j\in\mathbb{N}}$, we find that, for
any $t\in(0,\infty)$,
$$
\left|T_{\sigma}(f)\ast\phi_t
\right|\leq\sum_{j\in\mathbb{N}}
\lambda_j\left|T_{\sigma}(a_j)\ast\phi_t
\right|,
$$
which, combined with Definition \ref{df511},
further implies that
$$
m\left(T_{\sigma}(f),\phi\right)
\leq\sum_{j\in\mathbb{N}}
\lambda_jm\left(T_{\sigma}(a_j),\phi\right).
$$
Applying this and Definition \ref{Df1}(ii),
we conclude that
\begin{align}\label{pseudoxe6}
&\left\|m\left(T_{\sigma}(f),\phi\right)
\right\|_X\notag\\&\quad\leq\left\|\sum_{j\in\mathbb{N}}
\lambda_jm\left(T_{\sigma}(a_j),\phi\right)
\right\|_X\notag\\
&\quad\lesssim\left\|\sum_{j\in\mathbb{N}}
\lambda_jm\left(T_{\sigma}(a_j),\phi\right)\1bf_{2B_j}
\right\|_X+\left\|\sum_{j\in\mathbb{N}}
\lambda_jm\left(T_{\sigma}(a_j),\phi\right)\1bf_{(2B_j)
^{\complement}}\right\|_X\notag\\
&\quad=:\mathrm{VII}_1+\mathrm{VII}_2.
\end{align}
We then estimate $\mathrm{VII}_1$ and
$\mathrm{VII}_2$ respectively.

First, we deal with $\mathrm{VII}_1$.
To this end, we first estimate
$m(T_{\sigma}(a_j),\phi)$ for any $j\in\mathbb{N}$.
Indeed, for any $j\in\mathbb{N}$,
from Definition \ref{df511}(i) and
Lemma \ref{Atogxl3} with $f:=T_{\sigma}(a_j)$
and $\Phi:=\|\phi\|_{L^{\infty}(\rrn)}\1bf_{B(
\0bf,1)}$, we infer that
\begin{equation}\label{pseudoxe7}
m\left(T_{\sigma}(a_j),\phi\right)
\leq\sup_{t\in(0,\infty)}\left|
T_{\sigma}(a_j)\ast\phi_t\right|
\lesssim\mc\left(T_{\sigma}(a_j)\right).
\end{equation}
Let $r\in(\max\{1,s\eta'\},\infty)$ with
$\frac{1}{\eta}+\frac{1}{\eta'}=1$.
Then, using \eqref{pseudoxe7},
the $L^r(\rrn)$ boundedness
of the Hardy--littlewood maximal operator
$\mc$, Lemma \ref{bdd-pseudo}
with $p:=r$, and Definition \ref{loatomx}(ii),
we find that, for any $j\in\mathbb{N}$,
\begin{align}\label{pseudoxe8}
\left\|m\left(T_{\sigma}
(a_j,\phi)\right)\right\|_{L^r(\rrn)}
&\lesssim\left\|\mc\left(T_{\sigma}(a_j)\right)
\right\|_{L^r(\rrn)}
\lesssim\left\|T_{\sigma}(a_j)\right\|_{L^r(\rrn)}
\notag\\&\lesssim\left\|a_j\right\|_{L^r(\rrn)}
\lesssim\frac{|B_j|^{1/r}}{\|\1bf_{B_j}\|_{X}}.
\end{align}
In addition, from the assumption (ii) of
the present theorem and
an argument similar to that
used in the proof of \eqref{bddoxl4e1}
with $(r_0/s)'$ and $(X^{1/s})'$
therein replaced, respectively, by
$\eta$ and $Y$, it follows that,
for any $f\in L^1_{\mathrm{loc}}(\rrn)$,
\begin{equation*}
\left\|\mc^{((r/s)')}(f)\right\|_{Y}
\lesssim\|f\|_Y.
\end{equation*}
This, together with \eqref{pseudoxe8},
an argument similar to that used
in the estimation of $\mathrm{II}_1$
in the proof of Theorem \ref{Atogx} with
$\{\mc(a_j)\}_{j\in\mathbb{N}}$
therein replaced by $$\left\{m\left(T_{\sigma}(a_j),
\phi\right)\right\}_{j\in\mathbb{N}},$$
and \eqref{pseudoxe2}, implies that
\begin{equation}\label{pseudoxe10}
\mathrm{VII}_1\lesssim
\left\|\left[
\sum_{j\in\mathbb{N}}\left(
\frac{\lambda_j}{\|\1bf_{B_j}\|_X}\right)^s\1bf_{
B_j}\right]^{\frac{1}{s}}\right\|_X\lesssim
\left\|f\right\|_{h_X(\rrn)},
\end{equation}
which is the desired estimate of $\mathrm{VII}_1$.

On the other hand, we estimate $\mathrm{VII}_2$.
Indeed, applying Proposition \ref{pseudoae} with
$a:=a_j$ for any $j\in\mathbb{N}$
and Definition \ref{Df1}(ii), we
find that
$$
\mathrm{VII}_2\lesssim\left\|
\sum_{j\in\mathbb{N}}
\frac{\lambda_j}{\|\1bf_{B_j}\|_X}
\mc^{(\theta)}(\1bf_{B_j})\right\|_X.
$$
From this, \eqref{Atogxe6}, and \eqref{pseudoxe2},
we further deduce that
\begin{equation}\label{pseudoxe11}
\mathrm{VII}_2\lesssim
\left\|\left[
\sum_{j\in\mathbb{N}}\left(
\frac{\lambda_j}{\|\1bf_{B_j}\|_X}\right)^s\1bf_{
B_j}\right]^{\frac{1}{s}}\right\|_X\lesssim
\left\|f\right\|_{h_X(\rrn)},
\end{equation}
which completes the estimation of
$\mathrm{VII}_2$. Moreover,
combining both the assumptions (i) and (iv)
of the present theorem,
Remark \ref{Th6.1r}, Lemma
\ref{Th6.1l1}(ii), \eqref{pseudoxe6},
\eqref{pseudoxe10}, and \eqref{pseudoxe11},
we obtain
$$
\|T_{\sigma}(f)\|_{h_X(\rrn)}
\sim\left\|m\left(T_{\sigma}(f),\phi\right)\right\|_{X}
\lesssim\|f\|_{h_X(\rrn)}.
$$
This finishes the proof of Theorem \ref{pseudox}
\end{proof}

Next, we show Theorem \ref{pseudog}
with the help of Theorem \ref{pseudox}.

\begin{proof}[Proof of Theorem \ref{pseudog}]
Let all the symbols
be as in the present theorem. Then,
applying the assumptions $\m0(\omega)
\in(-\frac{n}{p},\infty)$ and $\MI(\omega)
\in(-\infty,0)$, and Theorem \ref{Th2},
we find that the Herz space $\HerzS$
under consideration is a BQBF space.
In addition, let
$$
s\in\left(0,\min\left\{1,p,q,\frac{n}{\Mw+n/p}
\right\}\right).
$$
Then, from the assumptions $\m0(\omega)
\in(-\frac{n}{p},\infty)$
and $\MI(\omega)\in(-\infty,0)$, and
Lemma \ref{rela}, it follows that
\begin{equation*}
\m0\left(\omega^s\right)
=s\m0(\omega)>-\frac{n}{p/s}
\end{equation*}
and
\begin{equation*}
\MI\left(\omega^s\right)=s
\MI(\omega)<0.
\end{equation*}
Combining these,
the assumptions $p/s,\ q/s\in(1,\infty)$,
and Theorem \ref{ball}
with $p$, $q$, and $\omega$
replaced, respectively, by $p/s$,
$q/s$, and $\omega^s$, we conclude that
the Herz space
$\Kmp^{p/s,q/s}_{\omega^s}(\rrn)$ is a BBF space.
From this and Lemma \ref{convexl}, it follows that
$[\HerzS]^{1/s}$ is a BBF space.
Moreover, let
$\theta\in(0,s)$,
$$
s_0\in\left(s,\min\left\{1,p,q,\frac{n}{
\Mw+n/p}\right\}\right),
$$
and $\eta\in(1,\infty)$ satisfy
$$
\eta<\min\left\{\frac{n}{n(1-s/p)-s\mw},\left(
\frac{p}{s}\right)'\right\}.
$$ Then,
for the above $\theta$, $s$, $s_0$, and $\eta$,
repeating an argument similar to that
used in the proof of Theorem \ref{bddog},
we find that the
following five statements hold true:
\begin{enumerate}
  \item[{\rm(i)}] for any $\{f_j\}_{j
  \in\mathbb{N}}\subset L^1_{\mathrm{loc}}(\rrn)$,
$$
\left\|\left\{
\sum_{j\in\mathbb{N}}
\left[\mc^{(\theta)}(f_j)\right]^{s}
\right\}^{1/s}
\right\|_{\HerzS}
\lesssim\left\|
\left(\sum_{j\in\mathbb{N}}
|f_j|^{s}
\right)^{1/s}
\right\|_{\HerzS};
$$
  \item[{\rm(ii)}] for any $f\in\Msc(\rrn)$,
  $$
  \|f\|_{[\HerzS]^{1/s}}
  \sim\sup\left\{\|fg\|_{L^1(\rrn)}:\
  \|g\|_{\HerzScsd}=1\right\}
  $$
  and
  $$
  \|f\|_{[\HerzS]^{1/s_0}}
  \sim\sup\left\{\|fg\|_{L^1(\rrn)}:\
  \|g\|_{\HerzScsod}=1\right\}
  $$
  with the positive equivalence constants
  independent of $f$;
  \item[{\rm(iii)}] both $\|\cdot\|_{\HerzScsd}$
  and $\|\cdot\|_{\HerzScsod}$ satisfy
  Definition \ref{Df1}(ii);
  \item[{\rm(iv)}] $\1bf_{B(\0bf,1)}\in\HerzScsod$;
  \item[{\rm(v)}] $\mc^{(\eta)}$ is bounded
  on $\HerzScsd$ and $\HerzScsod$.
\end{enumerate}
These, together with the
facts that $\HerzS$ is a BQBF space
and $[\HerzS]^{1/s}$ is a BBF space, and
Theorem \ref{pseudox} with $X:=\HerzS$,
$Y:=\HerzScsd$, and $Y_0:=\HerzScsod$,
further imply that $T_{\sigma}$ is well defined
on $\HaSaHl$ and, for any $f\in\HaSaHl$,
$$
\left\|T_{\sigma}(f)\right\|_{\HaSaHl}
\lesssim\left\|f\right\|_{\HaSaHl},
$$
which then completes the proof of
Theorem \ref{pseudog}.
\end{proof}

By both Theorem \ref{pseudog} and
Remark \ref{remark5.0.9}(ii),
we immediately conclude that the pseudo-differential
operator $T_\sigma$, with $\sigma\in S_{1,0}^{0}(\rrn)$,
is bounded on the local generalized Morrey--Hardy space
$h\MorrS$ as follows;
we omit the details.

\begin{corollary}
Let $p$, $q$, and $\omega$
be as in Corollary \ref{posedom}
and $T_\sigma$ be a pseudo-differential
operator with symbol
$\sigma\in S_{1,0}^{0}(\rrn)$.
Then $T_{\sigma}$ is well defined on $h\MorrS$ and
there exists a positive
constant $C$ such that, for any $f\in
h\MorrS$,
$$
\left\|T_\sigma(f)\right\|_{h\MorrS}
\leq C\|f\|_{h\MorrS}.
$$
\end{corollary}

\chapter{Weak Generalized
Herz--Hardy Spaces\label{sec8}}
\markboth{\scriptsize\rm\sc Weak Generalized Herz--Hardy Spaces}
{\scriptsize\rm\sc Weak Generalized Herz--Hardy Spaces}

The main target of this chapter
is to establish a complete real-variable theory of
weak generalized Herz--Hardy spaces.
Note that the classical weak Hardy space $WH^p(\rrn)$, with
$p\in(0,1]$, plays an important role
in the study of the boundedness of operators in harmonic
analysis. Indeed, in 1986, to find the biggest function space $X$
such that Calder\'{o}n--Zygmund operators
are bounded from $X$ to the weak Lebesgue space $WL^1(\rrn)$,
Fefferman and Soria \cite{FSo} introduced the weak Hardy space $WH^1(\rrn)$.
Via establishing the $\infty$-atomic characterization
of $WH^1(\rrn)$, Fefferman and Soria showed that
the convolutional Calder\'{o}n--Zygmund
operators with kernels satisfying the Dini condition were bounded from
$WH^1(\rrn)$ to $WL^1(\rrn)$.
In addition, let $\delta\in(0,1]$ and
$T$ be a convolutional $\delta$-type Calder\'{o}n--Zygmund operator.
It is well known that, for any given $p\in(n/(n+\delta),1]$, $T$ is bounded
on the classical Hardy space $H^p(\rrn)$ (see \cite{AM}).
However, $T$ is not bounded on $H^{n/(n+\delta)}(\rrn)$,
which is called the \emph{critical case}\index{critical case}
or the \emph{endpoint case}\index{endpoint case}.
In 1988, to deal with this critical case,
Liu \cite{Li} first introduced the weak Hardy space $WH^p(\rrn)$ with
$p\in(0,1]$, and showed that the
convolutional $\delta$-type Calder\'{o}n--Zygmund operators are
bounded from $H^{n/(n+\delta)}(\rrn)$
to $WH^{n/(n+\delta)}(\rrn)$ via establishing the $\infty$-atomic
characterization of $WH^p(\rrn)$.

Furthermore, the classical weak Hardy space
is the real interpolation space between the classical Hardy space
and the Lebesgue space $L^\infty(\rrn)$
(see \cite{FRS}), which is another
motivation to develop the real-variable theory of weak Hardy spaces.
Later on, a lot of works have
been done in the study of weak Hardy spaces and their variants;
see, for instance, \cite{qy00} for weighted weak Hardy spaces,
\cite{yyyz16} for variable weak Hardy spaces, \cite{ho17}
for weak Hardy--Morrey spaces, \cite{lyj16,YLK}
for weak Musielak--Orlicz Hardy spaces, and
\cite{gh15,h14} for vector-valued weak Hardy spaces.
Particularly, Zhang et al. \cite{ZWYY} introduced
the weak Hardy space $WH_X(\rrn)$ associated with
the ball quasi-Banach function space $X$ and characterized
$WH_{X}(\rrn)$ via various maximal functions, atoms,
and molecules. Moreover, they also established the boundedness
of Calder\'{o}n--Zygmund operators
from the Hardy space $H_X(\rrn)$ to
the weak Hardy space $WH_{X}(\rrn)$ in \cite{ZWYY}.
After that, various Littlewood--Paley function characterizations
of $WH_X(\rrn)$ and the real interpolation theorems have
been studied by Wang et al. in \cite{WYYZ}.

In this chapter, we first introduce weak
generalized Herz--Hardy spaces and then
establish their complete
real-variable theory. For this purpose,
recall that Zhang et al.\ \cite{ZWYY}
and Wang et al.\ \cite{WYYZ}
investigated the real-variable
theory of the weak Hardy space $WH_X(\rrn)$
associated with the
ball quasi-Banach function space $X$.
Combining these known results and
establishing some improved atomic and molecular
characterizations of $WH_X(\rrn)$
(see Propositions \ref{watomgd},
\ref{watomgr}, and \ref{wmolegr} below)
as well as new real interpolation
between the Hardy space $H_X(\rrn)$ and the Lebesgue
space $L^{\infty}(\rrn)$ (see Proposition
\ref{intergx} below),
we obtain various maximal function,
atomic, and molecular characterizations
of weak generalized Herz--Hardy spaces
and also show that the real interpolation spaces
between generalized Herz--Hardy spaces and
the Lebesgue space $L^{\infty}(\rrn)$
are just the new introduced weak
generalized Herz--Hardy spaces.
In addition, by establishing a
technique lemma about the
quasi-norm $\|\cdot\|_{W\HerzS}$ and
the Littlewood--Paley
function characterizations of $WH_X(\rrn)$
obtained in \cite{WYYZ}, we show
various Littlewood--Paley
function characterizations of weak generalized
Herz--Hardy spaces.
Furthermore, we establish
two boundedness criteria of
Calder\'{o}n--Zygmund operators
from the Hardy space $H_X(\rrn)$
to the weak Hardy space $WH_X(\rrn)$ and,
as a consequence, we finally deduce the boundedness
of Clader\'{o}n--Zygmund operators
from generalized Herz--Hardy spaces
to weak generalized Herz--Hardy spaces
even in the critical case.
In addition, we introduce weak generalized Morrey--Hardy spaces.
From the fact that, under some reasonable and sharp assumptions
on exponents, generalized Morrey spaces are included
into the scale of generalized Herz spaces, we also obtain
the corresponding real-variable characterizations
and applications of weak generalized Morrey--Hardy spaces
in this chapter.

We now introduce weak generalized Herz spaces as
follows.

\begin{definition}\label{df6.0.7}
Let $p,\ q\in(0,\infty)$ and $\omega\in M(\rp)$.
\begin{enumerate}
  \item[{\rm(i)}] The \emph{weak local generalized
  Herz space}\index{weak local generalized Herz space}
  $\WHerzSo$\index{$\WHerzSo$}
  is defined to be the set of all the measurable
  functions $f$ on $\rrn$ such that
  $$
  \|f\|_{\WHerzSo}:=\sup_{\alpha\in(0,\infty)}\left\{
  \alpha\left\|\1bf_{\{x\in\rrn:\ |f(x)|>\alpha\}}\right\|_{\HerzSo}
  \right\}<\infty.
  $$
  \item[{\rm(ii)}] The \emph{weak global generalized
  Herz space}\index{weak global generalized Herz space}
  $\WHerzS$\index{$\WHerzS$}
  is defined to be the set of all the measurable functions
  $f$ on $\rrn$ such that
  $$
  \|f\|_{\WHerzS}:=\sup_{\alpha\in(0,\infty)}\left\{
  \alpha\left\|\1bf_{\{x\in\rrn:\ |f(x)|>\alpha\}}\right\|_{\HerzS}
  \right\}<\infty.
  $$
\end{enumerate}
\end{definition}

\begin{remark}
In Definition \ref{df6.0.7}, when $\omega(t):=t^\alpha$
for any $t\in(0,\infty)$ and
for any given $\alpha\in\rr$, then
the weak local generalized Herz space $\WHerzSo$
coincides with the classical \emph{homogeneous weak Herz space
$W\dot{K}_p^{\alpha,q}(\rrn)$}\index{homogeneous weak Herz space}\index{$W\dot{K}_p^{\alpha,q}(\rrn)$}
which was originally introduced in \cite{hly97} (see also \cite[Section 1.4]{LYH}).
However, we should point out that,
even in this case, the weak global generalized Herz space $\WHerzS$
is also new.
\end{remark}

For any $N\in\mathbb{N}$ and $f\in\mathcal{S}'(\rrn)$, via
the non-tangential grand maximal function $\mathcal{M}_N(f)$
of $f$ as in \eqref{sec6e1}, we then introduce the definitions
of weak generalized Herz--Hardy spaces as
follows.\index{weak generalized Herz--Hardy space}

\begin{definition}\label{df6.11}
Let $p$, $q\in(0,\infty)$, $\omega\in M(\rp)$,
and $N\in\mathbb{N}$.
\begin{enumerate}
  \item[{\rm(i)}] The \emph{weak generalized
  Herz--Hardy space
  $\WHaSaHo$\index{$\WHaSaHo$}}, associated with
  the weak local generalized Herz space $W\HerzSo$,
  is defined to be the set of all the
  $f\in\mathcal{S}'(\rrn)$ such that
  $$
  \|f\|_{\WHaSaHo}:=\left\|\mathcal{M
  }_{N}(f)\right\|_{\WHerzSo}<\infty.
  $$
  \item[{\rm(ii)}] The \emph{weak generalized
  Herz--Hardy space
  $\WHaSaH$\index{$\WHaSaH$}}, associated with
  the weak global generalized Herz space $W\HerzS$,
  is defined to be the set of all the
  $f\in\mathcal{S}'(\rrn)$ such that
  $$
  \|f\|_{\WHaSaH}:=\left\|\mathcal{M}_{
  N}(f)\right\|_{\WHerzS}<\infty.
  $$
\end{enumerate}
\end{definition}

\begin{remark}
We point out that, in Definition \ref{df6.11},
if $\omega(t):=t^\alpha$ for any $t\in(0,\infty)$
and for any given $\alpha\in\rr$,
then the weak generalized Herz--Hardy space $\WHaSaHo$
coincides with the classical \emph{homogeneous weak Herz-type Hardy space}
\index{homogeneous weak Herz-type Hardy \\ space}$WH
\dot{K}_p^{\alpha,q}(\rrn)$\index{$WH\dot{K}_p^{\alpha,q}(\rrn)$}
which was originally introduced by Hu et al.\ in
\cite[Definition 3]{hly971}
(see also \cite[Definition 2.1.2]{LYH}).
However, even in this case,
the weak generalized Herz--Hardy space $W\HaSaH$ is also new.
\end{remark}

We next introduce the following concepts of
both weak generalized Morrey spaces
and associated Hardy spaces.

\begin{definition}\label{df6.13}
Let $p$, $q\in(0,\infty)$ and $\omega\in M(\rp)$.
\begin{enumerate}
  \item[{\rm(i)}] The \emph{weak local generalized
  Morrey space}\index{weak local generalized Morrey space}
  $W\MorrSo$\index{$W\MorrSo$} is defined to
  be the set of all the measurable
  functions $f$ on $\rrn$ such that
  $$
  \|f\|_{W\MorrSo}:=\sup_{\alpha\in(0,\infty)}\left\{
  \alpha\left\|\1bf_{\{x\in\rrn:\ |f(x)|>\alpha\}}\right\|_{\MorrSo}
  \right\}<\infty.
  $$
  \item[{\rm(ii)}] The \emph{weak global generalized
  Morrey space}\index{weak global generalized Morrey \\space}
  $W\MorrS$\index{$W\MorrS$} is defined to
  be the set of all the measurable functions
  $f$ on $\rrn$ such that
  $$
  \|f\|_{W\MorrS}:=\sup_{\alpha\in(0,\infty)}\left\{
  \alpha\left\|\1bf_{\{x\in\rrn:\ |f(x)|>\alpha\}}\right\|_{\MorrS}
  \right\}<\infty.
  $$
\end{enumerate}
\end{definition}

\begin{remark}\label{remark6.14}
\begin{enumerate}
  \item[{\rm(i)}] We should point out that,
  in Definition \ref{df6.13}, even when $\omega(t):=t^\alpha$
  for any $t\in(0,\infty)$ and for any given $\alpha\in\rr$,
  the weak generalized Morrey spaces $W\MorrSo$ and $W\MorrS$
  are also new.
  \item[{\rm(ii)}] In Definition \ref{df6.13}, let $p$, $q\in[1,\infty)$
  and $\omega$ satisfy $$\Mw\in(-\infty,0).$$
  Then, by Remark \ref{remhs}(iv), we conclude that,
  in this case, the weak generalized Morrey spaces
  $W\MorrSo$ and $W\MorrS$ coincide, respectively, with the
  weak generalized Herz spaces $\WHerzSo$ and $\WHerzS$ in the sense
  of equivalent quasi-norms.
\end{enumerate}
\end{remark}

\begin{definition}\label{df6.15}
Let $p$, $q\in(0,\infty)$, $\omega\in M(\rp)$,
and $N\in\mathbb{N}$.
\begin{enumerate}
  \item[{\rm(i)}] The \emph{weak generalized
  Morrey--Hardy space\index{weak generalized Morrey--Hardy \\space}
  $WH\MorrSo$\index{$WH\MorrSo$}}, associated with
  the weak local generalized Morrey space $W\MorrSo$,
  is defined to be the set of all
  the $f\in\mathcal{S}'(\rrn)$ such that
  $$
  \|f\|_{WH\MorrSo}:=\left\|\mathcal{M}_{
  N}(f)\right\|_{W\MorrSo}<\infty.
  $$
  \item[{\rm(ii)}] The \emph{weak
  generalized Morrey--Hardy space $WH\MorrS$},
  associated with the weak
  global generalized Morrey space $W\MorrSo$\index{$\MorrSo$},
  is defined to be the set of all the
  $f\in\mathcal{S}'(\rrn)$ such that
  $$
  \|f\|_{WH\MorrS}:=\left\|\mathcal{M}_{N}
  (f)\right\|_{W\MorrS}<\infty.
  $$
\end{enumerate}
\end{definition}

\begin{remark}\label{remark6.16}
In Definition \ref{df6.15}, let $p$, $q\in[1,\infty)$
and $\omega\in M(\rp)$ satisfy $\Mw\in(-\infty,0)$.
Then, from Remark \ref{remark6.14}(ii), it follows that,
in this case, the weak generalized Morrey--Hardy spaces
$WH\MorrSo$ and $WH\MorrS$ coincide, respectively,
with the weak generalized Herz--Hardy spaces
$\WHaSaHo$ and $\WHaSaH$ with equivalent quasi-norms.
\end{remark}

\section{Maximal Function Characterizations}

In this section, we establish
the maximal function characterizations
of the weak generalized
Herz--Hardy spaces $W\HaSaHo$ and $W\HaSaH$.
Recall that various radial
and non-tangential maximal functions
are defined as in Definition \ref{smax}.
Via these maximal functions, we then
show the following maximal
function characterizations
\index{maximal function characterization}of
the weak generalized Herz--Hardy space $W\HaSaHo$.

\begin{theorem}\label{macl}
Let $p$, $q$, $a,\ b\in(0,\infty)$, $\omega\in M(\rp)$,
$N\in\mathbb{N}$, and $\phi\in\mathcal{S}(\rrn)$ satisfy
$\int_{\rrn}\phi(x)\,dx\neq0.$
\begin{enumerate}
\item[\rm{(i)}] Let $N\in\mathbb{N}
\cap[\lfloor b+1\rfloor,\infty)$
and $\omega$ satisfy $\m0(\omega)
\in(-\frac{n}{p},\infty)$.
Then, for any $f\in\mathcal{S}'(\rrn)$,
$$
\|M(f,\phi)\|_{W\HerzSo}\lesssim
\|M^*_{a}(f,\phi)\|_{W\HerzSo}
\lesssim\|M^{**}_{b}(f,\phi)\|_{W\HerzSo},
$$
\begin{align*}
\|M(f,\phi)\|_{W\HerzSo}&\lesssim
\|\mathcal{M}_{N}(f)\|_{W\HerzSo}
\lesssim\|\mathcal{M}_{\lfloor b+1\rfloor}(f)
\|_{W\HerzSo}\\&\lesssim
\|M^{**}_{b}(f,\phi)\|_{W\HerzSo},
\end{align*}
and
$$
\|M^{**}_{b}(f,\phi)\|_{W\HerzSo}\sim
\|\mathcal{M}^{**}_{b,N}(f)\|_{W\HerzSo},
$$
where the implicit positive
constants are independent of $f$.
\item[\rm{(ii)}] Let
$\omega\in M(\rp)$ satisfy
$\m0(\omega)\in(-\frac{n}{p},\infty)$
and $\mi(\omega)\in(-\frac{n}{p},\infty)$.
Assume $b\in(\max\{\frac{n}{p},\Mw+\frac{n}{p}\},\infty)$.
Then, for any $f\in\mathcal{S}'(\rrn),$
$$\|M^{**}_{b}(f,\phi)\|_{W\HerzSo}
\lesssim\|M(f,\phi)\|_{W\HerzSo},$$
where the implicit positive
constant is independent of $f$.
In particular,
when $N\in\mathbb{N}\cap[\lfloor b+1\rfloor,\infty)$,
if one of the quantities
$$\|M(f,\phi)\|_{W\HerzSo},\
\|M^*_{a}(f,\phi)\|_{W\HerzSo},\ \|\mathcal{M}_{N}
(f)\|_{W\HerzSo},$$
$$
\|M^{**}_{b}(f,\phi)\|_{W\HerzSo},
\ \text{and}\ \|\mathcal{M}^{**}_{b,N}(f)\|_{W\HerzSo}$$
is finite, then the others are also finite and mutually
equivalent with the
positive equivalence
constants independent of $f$.
\end{enumerate}
\end{theorem}

To prove Theorem \ref{macl},
we first recall the following
concept of weak ball quasi-Banach
function spaces introduced in
\cite[Definition 2.8]{ZWYY}.

\begin{definition}\label{de-bqbfwx}
Let $X$ be a ball quasi-Banach
function space. The
\emph{weak ball quasi-Banach function
space}\index{weak ball quasi-Banach function space}
$WX$\index{$WX$} is defined to be the set
of all the measurable functions $f$ on $\rrn$
such that
$$
\left\|f\right\|_{WX}
:=\sup_{\lambda\in(0,\infty)}
\left\{\lambda\left\|\1bf_{\{x\in\rrn:\,
|f(x)|>\lambda\}}
\right\|_X\right\}<\infty.
$$
\end{definition}

\begin{remark}\label{wball}
Let $X$ be a ball quasi-Banach function
space and $E\subset\rrn$ a measurable
set. Then it is easy to show that
$$
\left\|\1bf_{E}\right\|_{WX}
=\left\|\1bf_E\right\|_X
$$
\end{remark}

The following conclusion
obtained in \cite[Lemma 2.13]{ZWYY}
shows that the weak ball
quasi-Banach function space
is also a ball quasi-Banach function
space. This plays a key
role in the proof of Theorem
\ref{macl}.

\begin{lemma}\label{bqbfwx}
Let $X$ be a ball quasi-Banach
function space. Then the weak ball
quasi-Banach function space $WX$
is also a ball quasi-Banach function
space.
\end{lemma}

To show the maximal function
characterizations of $W\HaSaHo$,
we also require a lemma about
the boundedness of the Hardy--Littlewood
maximal operator on weak local generalized
Herz spaces as follows.

\begin{lemma}\label{wmbhl}
Let $p,\ q\in(0,\infty)$ and $\omega\in M(\rp)$ satisfy
$\m0(\omega)\in(-\frac{n}{p},\infty)$ and $\mi(\omega)\in(-\frac{n}{p},
\infty)$.
Then, for any given $$r\in\left(0,\min
\left\{p,\frac{n}{\Mw+n/p}\right\}\right),$$
there exists a positive constant $C$ such that,
for any $f\in L^1_{\mathrm{loc}}(\rrn)$,
$$
\left\|\mc(f)\right\|_{[W\HerzSo]^{1/r}}
\leq C\|f\|_{[W\HerzSo]^{1/r}}.
$$
\end{lemma}

To prove Lemma \ref{wmbhl},
we first state the following
auxiliary conclusion
in ball quasi-Banach function spaces,
which was obtained in
\cite[Theorem 4.4]{syy22}.

\begin{lemma}\label{vwmbhll1}
Let $X$ be a ball quasi-Banach
function space and there exists a
$p_-\in(0,\infty)$ such that,
for any given $r\in(0,p_-)$ and
$u\in(1,\infty)$,
there exists a positive constant
$C$ such that, for any
$\{f_j\}_{j\in\mathbb{N}}\subset
L^1_{\mathrm{loc}}(\rrn)$,
$$
\left\|\left\{\sum_{j\in\mathbb{N}}
\left[\mc\left(f_j\right)\right]^u
\right\}^{\frac{1}{u}}\right\|_{X^{1/r}}
\leq C\left\|\left(\sum_{j\in\mathbb{N}}
\left|f_j\right|^u
\right)^{\frac{1}{u}}\right\|_{X^{1/r}}.
$$
Then, for any given $r\in(0,p_-)$
and $u\in(1,\infty)$, there
exists a positive constant $C$ such that,
for any $\{f_j\}_{j\in\mathbb{N}}
\subset L^1_{\mathrm{loc}}(\rrn)$,
$$
\left\|\left\{\sum_{j\in\mathbb{N}}
\left[\mc\left(f_j\right)\right]^u
\right\}^{\frac{1}{u}}\right\|_{(WX)^{1/r}}
\leq C\left\|\left(\sum_{j\in\mathbb{N}}
\left|f_j\right|^u
\right)^{\frac{1}{u}}\right\|_{(WX)^{1/r}}.
$$
\end{lemma}

Via this conclusion, we now prove
Lemma \ref{wmbhl}.

\begin{proof}[Proof of Lemma \ref{wmbhl}]
Let all the symbols be as in the present lemma.
Then, by the assumption
$\m0(\omega)\in(-\frac{n}{p},\infty)$
and Theorem \ref{Th3}, we
conclude that the local generalized
Herz space $\HerzSo$ under consideration
is a BQBF space.
Moreover, let $$
p_-:=\min\left\{p,\frac{n}{\Mw+n/p}\right\}.
$$
Then, for any given $r\in(0,p_-)$ and
$u\in(1,\infty)$,
applying Lemma \ref{vmbhl},
we find that, for any
$\{f_j\}_{j\in\mathbb{N}}\subset
L^1_{\mathrm{loc}}(\rrn)$,
$$
\left\|\left\{\sum_{j\in\mathbb{N}}
\left[\mc\left(f_j\right)\right]^u
\right\}^{\frac{1}{u}}\right\|_{[\HerzSo]^{1/r}}
\lesssim\left\|\left(\sum_{j\in\mathbb{N}}
\left|f_j\right|^u
\right)^{\frac{1}{u}}\right\|_{[\HerzSo]^{1/r}}.
$$
This, combined with the fact
that $\HerzSo$ is a BQBF space and
Lemma \ref{vwmbhll1} with
$X:=\HerzSo$,
$f_1:=f$, and $f_j:=0$ for any
$j\in\mathbb{N}\cap[2,\infty)$,
further implies that,
for any given $r\in(0,p_-)$ and
for any $f\in L^1_{\mathrm{loc}}(\rrn)$,
$$
\left\|\mc(f)\right\|_{[W\HerzSo]^{1/r}}
\lesssim\left\|f\right\|_{[W\HerzSo]^{1/r}},
$$
which then completes the
proof of Lemma \ref{wmbhl}.
\end{proof}

Via the above lemmas and the known maximal
function characterizations of
Hardy spaces associated with
ball quasi-Banach function spaces
proved in \cite[Theorem 3.1]{SHYY}
(see also Lemma \ref{Th5.4l1} above),
we next show Theorem \ref{macl}.

\begin{proof}[Proof of Theorem \ref{macl}]
Let all the symbols be as in the
present theorem.
Then, from the assumption
$\m0(\omega)\in(-\frac{n}{p},\infty)$
and Theorem \ref{Th3}, it follows that
the local generalized Herz space $\HerzSo$
under consideration is a BQBF space.
By this and Lemma \ref{bqbfwx} with
$X:=\HerzSo$,
we conclude that the weak space
$W\HerzSo$ is also a BQBF space. This,
combined with Lemma \ref{Th5.4l1}(i),
finishes the proof of (i).

Next, we show (ii).
Indeed, let
$$
r\in\left(\frac{n}{b},\min\left\{
p,\frac{n}{\Mw+n/p}\right\}\right).
$$
Then, using
Lemma \ref{wmbhl}, we find that,
for any
$f\in L^1_{\mathrm{loc}}(\rrn)$,
$$
\left\|\mc(f)\right\|_{[W\HerzSo]^{1/r}}
\lesssim\left\|f\right\|_{[W\HerzSo]^{1/r}}.
$$
From this, the fact that $b\in(\frac{n}{r},\infty)$,
and Lemma \ref{Th5.4l1}(ii),
we deduce that (ii) holds true and hence
complete the proof of Theorem \ref{macl}.
\end{proof}

Via Theorem \ref{macl} and Remark \ref{remark6.16},
we immediately obtain the following maximal function
characterizations of the weak generalized Morrey--Hardy spaces
$WH\MorrSo$; we omit the
details.\index{maximal function characterization}

\begin{corollary}
Let $a,\ b\in(0,\infty)$, $p$, $q\in[1,\infty)$,
$\omega\in M(\rp)$,
$N\in\mathbb{N}$, and $\phi\in\mathcal{S}(\rrn)$ satisfy
$\int_{\rrn}\phi(x)\,dx\neq0.$
\begin{enumerate}
\item[\rm{(i)}] Let
$N\in\mathbb{N}\cap[\lfloor b+1\rfloor,\infty)$
and $\omega$ satisfy $\MI(\omega)\in(-\infty,0)$ and
$$
-\frac{n}{p}<\m0(\omega)\leq\M0(\omega)<0.
$$
Then, for any $f\in\mathcal{S}'(\rrn)$,
$$
\|M(f,\phi)\|_{W\MorrSo}\lesssim
\|M^*_{a}(f,\phi)\|_{W\MorrSo}
\lesssim\|M^{**}_{b}(f,\phi)\|_{W\MorrSo},
$$
\begin{align*}
\|M(f,\phi)\|_{W\MorrSo}&\lesssim
\|\mathcal{M}_{N}(f)\|_{W\MorrSo}
\lesssim\|\mathcal{M}_{\lfloor
b+1\rfloor}(f)\|_{W\MorrSo}\\&\lesssim
\|M^{**}_{b}(f,\phi)\|_{W\MorrSo},
\end{align*}
and
$$
\|M^{**}_{b}(f,\phi)\|_{W\MorrSo}\sim\|
\mathcal{M}^{**}_{b,N}(f)\|_{W\MorrSo},
$$
where the implicit positive
constants are independent of $f$.
\item[\rm{(ii)}] Let $\omega\in M(\rp)$ satisfy
$$
-\frac{n}{p}<\m0(\omega)\leq\M0(\omega)<0
$$
and
$$
-\frac{n}{p}<\m0(\omega)\leq\M0(\omega)<0.
$$
Assume $b\in(\frac{n}{p},\infty)$.
Then, for any $f\in\mathcal{S}'(\rrn),$
$$\|M^{**}_{b}(f,\phi)\|_{W\MorrSo}
\lesssim\|M(f,\phi)\|_{W\MorrSo},$$
where the implicit positive
constant is independent of $f$.
In particular, when $N\in\mathbb{N}
\cap[\lfloor b+1\rfloor,\infty)$,
if one of the quantities
$$\|M(f,\phi)\|_{W\MorrSo},\
\|M^*_{a}(f,\phi)\|_{W\MorrSo},\ \|\mathcal{M}_{N}
(f)\|_{W\MorrSo},$$
$$
\|M^{**}_{b}(f,\phi)\|_{W\MorrSo},
\ \text{and}\ \|\mathcal{M}^{**}
_{b,N}(f)\|_{W\MorrSo}$$
is finite, then the others are also
finite and mutually
equivalent with the
positive equivalence
constants independent of $f$.
\end{enumerate}
\end{corollary}

We now establish the maximal function characterizations
of the weak generalized Herz--Hardy space $W\HaSaH$ as
follows.\index{maximal function characterization}

\begin{theorem}\label{Th7.11}
Let $p$, $q$, $a,\ b\in(0,\infty)$, $\omega\in M(\rp)$,
$N\in\mathbb{N}$, and $\phi\in\mathcal{S}(\rrn)$ satisfy
$\int_{\rrn}\phi(x)\,dx\neq0.$
\begin{enumerate}
\item[\rm{(i)}] Let $N\in\mathbb{N}
\cap[\lfloor b+1\rfloor,\infty)$
and $\omega$ satisfy $\m0(\omega)
\in(-\frac{n}{p},\infty)$ and
$\MI(\omega)\in(-\infty,0)$.
Then, for any $f\in\mathcal{S}'(\rrn)$,
$$
\|M(f,\phi)\|_{W\HerzS}\lesssim
\|M^*_{a}(f,\phi)\|_{W\HerzS}
\lesssim\|M^{**}_{b}(f,\phi)\|_{W\HerzS},
$$
\begin{align*}
\|M(f,\phi)\|_{W\HerzS}&\lesssim
\|\mathcal{M}_{N}(f)\|_{W\HerzS}
\lesssim\|\mathcal{M}_{\lfloor
b+1\rfloor}(f)\|_{W\HerzS}\\&\lesssim
\|M^{**}_{b}(f,\phi)\|_{W\HerzS},
\end{align*}
and
$$
\|M^{**}_{b}(f,\phi)\|_{W\HerzS}\sim
\|\mathcal{M}^{**}_{b,N}(f)\|_{W\HerzS},
$$
where the implicit positive
constants are independent of $f$.
\item[\rm{(ii)}] Let $\omega\in M(\rp)$ satisfy
$\m0(\omega)\in(-\frac{n}{p},\infty)$ and
$$-\frac{n}{p}<\mi(\omega)\leq\MI(\omega)<0.$$
Assume $b\in(\max\{\frac{n}{p},\Mw+\frac{n}{p}\},\infty)$.
Then, for any $f\in\mathcal{S}'(\rrn),$
$$\|M^{**}_{b}(f,\phi)\|_{W\HerzS}
\lesssim\|M(f,\phi)\|_{W\HerzS},$$
where the implicit positive
constant is independent of $f$.
In particular, when $N\in\mathbb{N}
\cap[\lfloor b+1\rfloor,\infty)$,
if one of the quantities
$$\|M(f,\phi)\|_{W\HerzS},\
\|M^*_{a}(f,\phi)\|_{W\HerzS},\ \|\mathcal{M}_{N}
(f)\|_{W\HerzS},$$
$$
\|M^{**}_{b}(f,\phi)\|_{W\HerzS},
\ \text{and}\ \|\mathcal{M}^{**}_{b,N}(f)\|_{W\HerzS}$$
is finite, then the others are also finite and mutually
equivalent with the
positive equivalence constants independent of $f$.
\end{enumerate}
\end{theorem}

To prove this maximal function characterizations,
we first show the boundedness of the
Hardy--Littlewood maximal operator on
weak global generalized Herz spaces as follows.

\begin{lemma}\label{wmbhg}
Let $p,\ q\in(0,\infty)$ and $\omega\in M(\rp)$ satisfy
$\m0(\omega)\in(-\frac{n}{p},\infty)$ and $\mi(\omega)\in(-\frac{n}{p},
\infty)$.
Then, for any given $$r\in\left(0,\min\left\{p,\frac{n}{\Mw+n/p}\right\}\right),$$
there exists a positive constant $C$ such that,
for any $f\in L^1_{\mathrm{loc}}(\rrn)$,
$$
\left\|\mc(f)\right\|_{[W\HerzS]^{1/r}}\leq C\|f\|_{[W\HerzS]^{1/r}}.
$$
\end{lemma}

\begin{proof}
Let all the symbols be as in the present lemma.
Then, from the assumptions
$\m0(\omega)\in(-\frac{n}{p},\infty)$
and $\MI(\omega)\in(-\infty,0)$,
and Theorem \ref{Th2}, it follows
that the global generalized
Herz space $\HerzS$ under consideration
is a BQBF space.
Let $$
p_-:=\min\left\{p,\frac{n}{\Mw+n/p}\right\}.
$$
Then, for any given $r\in(0,p_-)$ and
$u\in(1,\infty)$,
applying Lemma \ref{vmbhg},
we find that, for any
$\{f_j\}_{j\in\mathbb{N}}\subset
L^1_{\mathrm{loc}}(\rrn)$,
$$
\left\|\left\{\sum_{j\in\mathbb{N}}
\left[\mc\left(f_j\right)\right]^u
\right\}^{\frac{1}{u}}\right\|_{[\HerzS]^{1/r}}
\lesssim\left\|\left(\sum_{j\in\mathbb{N}}
\left|f_j\right|^u
\right)^{\frac{1}{u}}\right\|_{[\HerzS]^{1/r}}.
$$
This, combined with the fact
that $\HerzS$ is a BQBF space and
Lemma \ref{vwmbhll1} with
$X:=\HerzS$,
$f_1:=f$, and $f_j:=0$ for any
$j\in\mathbb{N}\cap[2,\infty)$,
further implies that,
for any given $r\in(0,p_-)$ and for
any $f\in L^1_{\mathrm{loc}}(\rrn)$,
$$
\left\|\mc(f)\right\|_{[W\HerzS]^{1/r}}
\lesssim\left\|f\right\|_{[W\HerzS]^{1/r}},
$$
which then completes the
proof of Lemma \ref{wmbhg}.
\end{proof}

Using Lemma \ref{wmbhg}, we next
show Theorem \ref{Th7.11}.

\begin{proof}[Proof of Theorem \ref{Th7.11}]
Let all the symbols be as in the present theorem.
Then, applying the assumptions
$\m0(\omega)\in(-\frac{n}{p},\infty)$
and $\MI(\omega)\in(-\infty,0)$, and
Theorem \ref{Th2}, we conclude that
the global generalized Herz space $\HerzS$
under consideration is a BQBF space.
This, together with
Lemma \ref{bqbfwx}, further implies
that the weak Herz space
$W\HerzS$ is a BQBF space.
Thus, by Lemma \ref{Th5.4l1}(i),
we then complete the proof of (i).

Next, we prove (ii).
Indeed, let
\begin{equation}\label{th7.11e2}
r\in\left(\frac{n}{b},\min
\left\{p,\frac{n}{\Mw+n/p}\right\}\right).
\end{equation}
Then, from Lemma \ref{wmbhg},
it follows that, for any
$f\in L^1_{\mathrm{loc}}(\rrn)$,
\begin{equation}\label{th7.11e1}
\left\|\mc(f)\right\|_{[W\HerzS]^{1/r}}
\lesssim\left\|f\right\|_{[W\HerzS]^{1/r}}.
\end{equation}
In addition, by \eqref{th7.11e2},
we find that $b\in(\frac{n}{r},\infty)$.
Combining this, \eqref{th7.11e1},
and Lemma \ref{Th5.4l1}(ii),
we then complete the proof of (ii),
and hence of Theorem \ref{Th7.11}.
\end{proof}

As an application of Theorem \ref{Th7.11}, we next
present the following result about the maximal
function characterizations of
the weak generalized Morrey--Hardy space $WH\MorrS$,
which is just a immediately corollary of both
Theorem \ref{Th7.11}
and Remark \ref{remark6.16}; we omit the
details.\index{maximal function characterization}

\begin{corollary}
Let $a,\ b\in(0,\infty)$, $p$,
$q\in[1,\infty)$, $\omega\in M(\rp)$,
$N\in\mathbb{N}$, and
$\phi\in\mathcal{S}(\rrn)$ satisfy
$\int_{\rrn}\phi(x)\,dx\neq0.$
\begin{enumerate}
\item[\rm{(i)}] Let $N\in\mathbb{N}
\cap[\lfloor b+1\rfloor,\infty)$
and $\omega$ satisfy
$\MI(\omega)\in(-\infty,0)$ and
$$
-\frac{n}{p}<\m0(\omega)\leq\M0(\omega)<0.
$$
Then, for any $f\in\mathcal{S}'(\rrn)$,
$$
\|M(f,\phi)\|_{W\MorrS}\lesssim
\|M^*_{a}(f,\phi)\|_{W\MorrS}
\lesssim\|M^{**}_{b}(f,\phi)\|_{W\MorrS},
$$
\begin{align*}
\|M(f,\phi)\|_{W\MorrS}&\lesssim
\|\mathcal{M}_{N}(f)\|_{W\MorrS}
\lesssim\|\mathcal{M}_{\lfloor
b+1\rfloor}(f)\|_{W\MorrS}\\&\lesssim
\|M^{**}_{b}(f,\phi)\|_{W\MorrS},
\end{align*}
and
$$
\|M^{**}_{b}(f,\phi)\|_{W\MorrS}\sim\|
\mathcal{M}^{**}_{b,N}(f)\|_{W\MorrS},
$$
where the implicit positive
constants are independent of $f$.
\item[\rm{(ii)}] Let $\omega\in M(\rp)$ satisfy
$$
-\frac{n}{p}<\m0(\omega)\leq\M0(\omega)<0
$$
and
$$
-\frac{n}{p}<\m0(\omega)\leq\M0(\omega)<0.
$$
Assume $b\in(\frac{n}{p},\infty)$.
Then, for any $f\in\mathcal{S}'(\rrn),$
$$\|M^{**}_{b}(f,\phi)\|_{W\MorrS}
\lesssim\|M(f,\phi)\|_{W\MorrS},$$
where the implicit positive
constant is independent of $f$.
In particular, when
$N\in\mathbb{N}\cap[\lfloor b+1\rfloor,\infty)$,
if one of the quantities
$$\|M(f,\phi)\|_{W\MorrS},\
\|M^*_{a}(f,\phi)\|_{W\MorrS},\ \|\mathcal{M}_{N}
(f)\|_{W\MorrS},$$
$$
\|M^{**}_{b}(f,\phi)\|_{W\MorrS},
\ \text{and}\
\|\mathcal{M}^{**}_{b,N}(f)\|_{W\MorrS}$$
is finite, then the others are
also finite and mutually
equivalent with the
positive equivalence constants
independent of $f$.
\end{enumerate}
\end{corollary}

\section{Relations with Weak Generalized Herz \\Spaces}

In this section,
via the relation between
ball quasi-Banach function
spaces and associated Hardy spaces obtained
in \cite[Theorem 3.4]{SHYY}
(see also Lemma \ref{Th5.7l1} above),
we investigate the relations between
weak generalized Herz spaces and
associated Hardy spaces.
We first establish the relation between
$W\HerzSo$ and the associated Hardy
space $\WHaSaHo$ as follows. Indeed, the following
theorem shows that, under some
reasonable and sharp assumptions,
$W\HerzSo=W\HaSaHo$ with equivalent quasi-norms.

\begin{theorem}\label{wrelal}
Let $p\in(1,\infty)$,
$q\in(0,\infty)$, and
$\omega\in M(\rp)$ satisfy
$$-\frac{n}{p}<\m0(\omega)\leq\M0(\omega)<
\frac{n}{p'}$$
and
$$-\frac{n}{p}<\mi(\omega)\leq\MI(\omega)<\frac{n}{p'},$$
where $\frac{1}{p}+\frac{1}{p'}=1$. Then
\begin{enumerate}
\item[\rm{(i)}] $W\HerzSo\hookrightarrow\mathcal{S}'(\rrn)$.
\item[\rm{(ii)}] If $f\in\WHerzSo$, then $f\in\WHaSaHo$
and there exists a positive constant $C$, independent
of $f$, such that $$\left\|f\right\|_{\WHaSaHo}
\leq C\left\|f\right\|_{\WHerzSo}.$$
\item[\rm{(iii)}] If $f\in\WHaSaHo$, then there exists a
locally integrable function $g$ belonging to $\WHerzSo$
such that $g$ represents $f$, which means that $f=g$ in
$\mathcal{S}'(\rrn)$, $$\left\|f\right\|_{
\WHaSaHo}=\left\|g\right\|_{\WHaSaHo},$$
and there exists a positive constant $C$, independent of $f$,
such that $$\left\|g\right\|_{\WHerzSo}
\leq C\left\|f\right\|_{\WHaSaHo}.$$
\end{enumerate}
\end{theorem}

\begin{proof}
Let all the symbols be as in the present theorem.
Then, from the assumption
$\m0(\omega)\in(-\frac{n}{p},\infty)$
and Theorem \ref{Th3}, it follows that the
local generalized Herz space $\HerzSo$
under consideration is a BQBF space.
This, together with Lemma \ref{bqbfwx},
further implies that the weak local
generalized Herz space $W\HerzSo$ is a
BQBF space. Therefore,
to finish the proof of the present theorem,
we only need to show that $W\HerzSo$
satisfies all the assumptions of
Lemma \ref{Th5.7l1}.
Namely, there exists an $r\in(1,\infty)$
such that the Hardy--Littlewood
maximal operator $\mc$ is bounded on
$[W\HerzSo]^{1/r}$.

Indeed, applying the assumptions
$$\mw\in\left(-\frac{n}{p},\infty\right)$$ and
$$\Mw\in\left(-\infty,\frac{n}{p'}\right),$$
and Remark \ref{remark 1.1.3}(iii),
we conclude that
$$
\frac{n}{\Mw+n/p}\in(1,\infty).
$$
Combining this and the assumption
$p\in(1,\infty)$, we further find that
$$
\min\left\{p,\frac{n}
{\Mw+n/p}\right\}\in(1,\infty).
$$
Thus, we can choose an
$$
r\in\left(1,\min\left\{
p,\frac{n}{\Mw+n/p}\right\}\right).
$$
For this $r$, from Lemma \ref{wmbhl}, we infer that,
for any $f\in L^1_{\mathrm{loc}}(\rrn)$,
$$
\left\|\mc(f)\right\|_{[W\HerzSo]^{1/r}}
\lesssim\left\|f\right\|_{[W\HerzSo]^{1/r}}.
$$
This further implies that
there exists an $r\in(1,\infty)$
such that $\mc$ is bounded on
$[W\HerzSo]^{1/r}$.
Thus, all the assumptions of
Lemma \ref{Th5.7l1} hold
true for $W\HerzSo$, and hence
the proof of Theorem \ref{wrelal} is
then completed.
\end{proof}

Using Theorem \ref{wrelal} and
Remark \ref{remark6.16}, we immediately
obtain the following conclusion, which shows that
$W\MorrSo=WH\MorrSo$ with
equivalent quasi-norms
under some reasonable and sharp
assumptions; we omit the details.

\begin{corollary}\label{coro6.3.2}
Let $p\in(1,\infty)$, $q\in[1,\infty)$, and
$\omega\in M(\rp)$ satisfy
$$-\frac{n}{p}<\m0(\omega)\leq\M0(\omega)<0$$
and
$$-\frac{n}{p}<\mi(\omega)\leq\MI(\omega)<0.$$
Then
\begin{enumerate}
\item[\rm{(i)}] $W\MorrSo\hookrightarrow\mathcal{S}'(\rrn)$.
\item[\rm{(ii)}] If $f\in W\MorrSo$, then $f\in WH\MorrSo$
and there exists a positive constant $C$, independent
of $f$, such that $$\left\|f
\right\|_{WH\MorrSo}\leq C\left\|f\right\|_{W\MorrSo}.$$
\item[\rm{(iii)}] If $f\in WH\MorrSo$, then there exists a
locally integrable function $g$ belonging to $W\MorrSo$
such that $g$ represents $f$, which means that $f=g$ in
$\mathcal{S}'(\rrn)$, $$\left\|f\right\|_{WH\MorrSo}=
\left\|g\right\|_{WH\MorrSo},$$
and there exists a positive constant $C$, independent of $f$,
such that $$\left\|g\right\|_{W\MorrSo}
\leq C\left\|f\right\|_{WH\MorrSo}.$$
\end{enumerate}
\end{corollary}

Similarly, we next prove that, under some reasonable and
sharp assumptions,
the weak generalized Herz--Hardy space
$\WHaSaH$ coincides with the weak global generalized
Herz space $\WHerzS$ with equivalent quasi-norms as follows.

\begin{theorem}\label{wrelag}
Let $p\in(1,\infty)$,
$q\in(0,\infty)$, and
$\omega\in M(\rp)$ with
$$-\frac{n}{p}<\m0(\omega)\leq\M0(\omega)<\frac{n}{p'}$$
and
$$-\frac{n}{p}<\mi(\omega)\leq\MI(\omega)<0.$$
Then
\begin{enumerate}
\item[\rm{(i)}] $W\HerzS\hookrightarrow\mathcal{S}'(\rrn)$.
\item[\rm{(ii)}] If $f\in\WHerzS$, then $f\in\WHaSaH$
and there exists a positive constant $C$, independent
of $f$, such that $$\left\|f\right\|_{\WHaSaH}\leq
C\left\|f\right\|_{\WHerzS}.$$
\item[\rm{(iii)}] If $f\in\WHaSaH$, then there exists a
locally integrable function $g$ belonging to $\WHerzS$
such that $g$ represents $f$, which means that $f=g$ in
$\mathcal{S}'(\rrn)$, $$\left\|f\right\|_{\WHaSaH}
=\left\|g\right\|_{\WHaSaH},$$
and there exists a positive constant $C$, independent of $f$,
such that $$\left\|g\right\|_{\WHerzS}\leq
C\left\|f\right\|_{\WHaSaH}.$$
\end{enumerate}
\end{theorem}

\begin{proof}
Let all the symbols be as in the present theorem.
Notice that $\omega$ satisfies
both $\m0(\omega)\in(-\frac{n}{p},\infty)$ and
$\MI(\omega)\in(-\infty,0)$.
From these and Theorem \ref{Th2},
it follows that the global generalized
Herz space $\HerzS$ under consideration
is a BQBF space. By this and
Lemma \ref{bqbfwx}, we conclude that
the weak global generalized Herz space
$W\HerzS$ is also a BQBF space.
On the other hand, let
$$
r\in\left(1,\min\left\{
p,\frac{n}{\Mw+n/p}\right\}\right).
$$
Then, applying Lemma \ref{wmbhg},
we find that,
for any $f\in L^1_{\mathrm{loc}}(\rrn)$,
$$
\left\|\mc(f)\right\|_{[W\HerzS]^{1/r}}
\lesssim\left\|f\right\|_{[W\HerzS]^{1/r}}.
$$
This, together with the facts that
$r\in(1,\infty)$ and $W\HerzS$ is a
BQBF space, and Lemma \ref{Th5.7l1}
with $X:=W\HerzS$,
finishes the proof of Theorem \ref{wrelal}.
\end{proof}

Via both Theorem \ref{wrelag} and Remark \ref{remark6.16},
we immediately obtain the following relation
between the weak generalized Morrey--Hardy space $WH\MorrS$
and the weak global generalized Morrey space $W\MorrS$; we omit the
details.

\begin{corollary}
Let $p$, $q$, and $\omega$ be as in Corollary \ref{coro6.3.2}.
Then
\begin{enumerate}
\item[\rm{(i)}] $W\MorrS\hookrightarrow\mathcal{S}'(\rrn)$.
\item[\rm{(ii)}] If $f\in W\MorrS$, then $f\in WH\MorrS$
and there exists a positive constant $C$, independent
of $f$, such that $$\left\|f\right\|_{WH\MorrS}\leq
C\left\|f\right\|_{W\MorrS}.$$
\item[\rm{(iii)}] If $f\in WH\MorrS$, then there exists a
locally integrable function $g$ belonging to $W\MorrS$
such that $g$ represents $f$, which means that $f=g$ in
$\mathcal{S}'(\rrn)$, $$\left\|f\right\|_{WH\MorrS}=
\left\|g\right\|_{WH\MorrS},$$
and there exists a positive constant $C$, independent of $f$,
such that $$\left\|g\right\|_{W\MorrS}
\leq C\left\|f\right\|_{WH\MorrS}.$$
\end{enumerate}
\end{corollary}

\section{Atomic Characterizations}

The main target of this section
is to characterize weak generalized Herz--Hardy spaces
via atoms. For this purpose, we first establish the
atomic characterization of the weak generalized
Herz--Hardy space $W\HaSaHo$ via
the known atomic characterization
of weak Hardy spaces associated with
ball quasi-Banach function spaces directly.
However, due to the deficiency
of associate spaces of global generalized
Herz spaces, the atomic characterization
of Hardy spaces associated with
ball quasi-Banach function spaces mentioned
above is not applicable to establish the atomic
characterization of the weak
generalized Herz--Hardy space $W\HaSaH$.
To overcome this obstacle,
we first establish an improved
atomic characterization
of weak Hardy spaces associated with ball
quasi-Banach function spaces (see both Propositions
\ref{watomgd} and \ref{watomgr} below)
without recourse to associate spaces.
Then, using this improved conclusion,
we obtain the desired atomic characterization
of $W\HaSaH$.

Recall that the concept
of $(\HerzSo,\,r,\,d)$-atoms
is given in Definition \ref{atom}.
We now introduce the
weak generalized atomic Herz--Hardy spaces associated
with the local generalized Herz spaces as follows.

\begin{definition}\label{dfwatom}
Let $p,\ q\in(0,\infty)$, $\omega\in M(\rp)$ with
$\m0(\omega)\in(-\frac{n}{p},\infty)$ and $\mi(\omega)\in(
-\frac{n}{p},\infty)$,
$$
p_-\in\left(0,\frac{\min\{p,q,
\frac{n}{\Mw+n/p}\}}{\max\{1,p,q\}}\right),
$$
$d\geq\lfloor n(1/p_--1)\rfloor$ be a fixed
integer, and
$$
r\in\left(\max\left\{1,p,\frac{n}
{\mw+n/p}\right\},\infty\right].
$$
Then the \emph{weak generalized atomic Herz--Hardy space}
\index{weak generalized atomic Herz--Hardy \\ space}$WH
\Kmp_{\omega,\0bf}^{p,q,r,d}
(\rrn)$\index{$WH\Kmp_{\omega,\0bf}^{p,q,r,d}(\rrn)$},
associated with
the weak local generalized
Herz space $W\HerzSo$, is defined to
be the set of all the
$f\in\mathcal{S}'(\rrn)$ such that
there exist a sequence $\{a_{i,j}\}_{i\in
\mathbb{Z},\,j\in\mathbb{N}}$ of
$(\HerzSo,\,r,\,d)$-atoms supported,
respectively, in the balls $\{B_{i,j}
\}_{i\in\mathbb{Z},\,j\in\mathbb{N}}
\subset\mathbb{B}$ and three
positive constants $c\in(0,1]$,
$A$, and $\widetilde{A}$,
independent of $f$, satisfying that,
for any $i\in\mathbb{Z}$,
$$
\sum_{j\in\mathbb{N}}\1bf_{cB_{i,j}}\leq A,
$$
$$
f=\sum_{i\in\mathbb{Z}}\sum_{j\in\mathbb{N}}
\widetilde{A}2^i\left\|\1bf_{B_{i,j}}
\right\|_{\HerzSo}a_{i,j}
$$
in $\mathcal{S}'(\rrn)$, and
$$
\sup_{i\in\mathbb{Z}}\left\{2^i
\left\|\sum_{j\in\mathbb{N}}
\1bf_{B_{i,j}}\right\|_{\HerzSo}\right\}<\infty.
$$
Moreover, for any $f\in
WH\Kmp_{\omega,\0bf}^{p,q,r,d}(\rrn)$,
$$
\|f\|_{WH\Kmp_{\omega,\0bf}^{p,q,r,d}(\rrn)}:=
\inf\sup_{i\in\mathbb{Z}}2^i\left\{
\left\|\sum_{j\in\mathbb{N}}
\1bf_{B_{i,j}}\right\|_{\HerzSo}\right\},
$$
where the infimum is taken over
all the decompositions of
$f$ as above.
\end{definition}

Then we have the following atomic characterization
of the weak generalized Herz--Hardy space
$\WHaSaHo$.\index{atomic characterization}

\begin{theorem}\label{watom}
Let $p$, $q$, $\omega$,
$r$, and $d$ be as in Definition \ref{dfwatom}.
Then
$$
WH\HerzSo=WH\Kmp_{\omega,\0bf}^{p,q,r,d}(\rrn)
$$
with equivalent quasi-norms.
\end{theorem}

To show this theorem,
we first recall the following
definition of weak Hardy spaces
associated with ball quasi-Banach
function spaces, which is just
\cite[Definition 2.21]{ZWYY}.

\begin{definition}\label{hardyxw}
Let $X$ be a ball quasi-Banach
function space and $N\in
\mathbb{N}$. Then the \emph{weak Hardy
space} $WH_X(\rrn)$\index{$WH_X(\rrn)$}
is defined to be the set of
all the $f\in\mathcal{S}'(\rrn)$ such that
$$
\left\|f\right\|_{WH_X(\rrn)}
:=\left\|\mc_N(f)\right\|_{WX}
<\infty.
$$
\end{definition}

The following conclusion shows that,
under some assumptions of the ball
quasi-Banach function space $X$,
distributions in the weak Hardy space
$WH_X(\rrn)$ can be decomposed
into linear combinations of atoms,
which was proved in
\cite[Theorem 4.2]{ZWYY} and
plays an essential role
in the proof of
Theorem \ref{watom} above.

\begin{lemma}\label{watomd}
Let $X$ be a ball quasi-Banach function
space and let $r\in(0,1)$,
$r_0,\ s,\ p_-\in(0,\infty)$, and
$\theta\in(1,\infty)$ be such that
the following four statements hold true:
\begin{enumerate}
  \item[{\rm(i)}] there exists a positive constant
  $C$ such that, for any $\{f_j\}_{j\in\mathbb{N}}\subset
  L^1_{\mathrm{loc}}(\rrn)$,
  $$
  \left\|\left\{\sum_{j\in\mathbb{N}}
  \left[\mc\left(f_j\right)\right]^{1/r}
  \right\}^r\right\|_{X^{1/r}}
  \leq C\left\|\left(\sum_{j\in\mathbb{N}}
  \left|f_j\right|^{1/r}\right)^{r}\right\|_{X^{1/r}};
  $$
  \item[{\rm(ii)}] $X$ is $\theta$-concave and
  $\mc$ bounded on $X^{\frac{1}{\theta p_-}}$;
  \item[{\rm(iii)}] $X^{1/s}$ is a
  ball Banach function space and $\mc$ bounded on
  $(X^{1/s})'$;
  \item[{\rm(iv)}] $\mc$ is bounded on $(WX)^{1/r_0}$.
\end{enumerate}
Assume that $d\geq\lfloor n(1/p_--1)\rfloor$
is a fixed nonnegative integer and $f\in WH_X(\rrn)$.
Then there exist
$\{a_{i,j}\}_{i\in\mathbb{Z},\,j\in\mathbb{N}}$
of $(X,\,\infty,\,d)$-atoms supported,
respectively, in the balls $\{B_{i,j}\}_{
i\in\mathbb{Z},\,j\in\mathbb{N}}\subset\mathbb{B}$
and three positive constants $c\in(0,1]$,
$A$, and $\widetilde{A}$, independent of $f$,
such that, for any
$i\in\mathbb{Z}$, $$\sum_{j\in\mathbb{N}}
\1bf_{cB_{i,j}}\leq A,$$
$$f=\sum_{i\in\mathbb{Z}}\sum_{j\in\mathbb{N}}
\widetilde{A}2^i\left\|\1bf_{B_{i,j}}\right\|_Xa_{i,j}$$
in $\mathcal{S}'(\rrn)$, and
$$
\sup_{i\in\mathbb{Z}}\left\{
2^i\left\|\sum_{j\in\mathbb{N}}
\1bf_{B_{i,j}}\right\|_X\right\}\lesssim
\left\|f\right\|_{WH_X(\rrn)},
$$
where the implicit positive constant is independent
of $f$.
\end{lemma}

The following atomic reconstruction
theorem of the weak Hardy space $WH_X(\rrn)$
is just \cite[Theorem 4.7]{ZWYY},
which also plays a vital role in
the proof of Theorem \ref{watom}.

\begin{lemma}\label{watomr}
Let $X$ be a ball quasi-Banach function space
and let $p_-\in(0,1)$ be such that the following
three statements hold true:
\begin{enumerate}
  \item[{\rm(i)}] for any given
  $\theta\in(0,p_-)$ and $u\in(1,\infty)$,
  there exists a positive constant $C$ such that,
for any $\{f
_{j}\}_{j\in\mathbb{N}}\subset L_{{\rm loc}}^{1}(\rrn)$,
\begin{equation*}
\left\|\left\{\sum_{j\in\mathbb{N}}\left[\mc
(f_{j})\right]^{u}\right\}^
{\frac{1}{u}}\right\|_{X^{1/\theta}}
\leq C\left\|\left\{
\sum_{j\in\mathbb{N}}|f_{j}|^{u}\right\}
^{\frac{1}{u}}\right\|_{X^{1/\theta}};
\end{equation*}
  \item[{\rm(ii)}] for any $s\in(0,p_-)$,
  $X^{1/s}$ is a ball Banach function space;
  \item[{\rm(iii)}] there exist an $s_0\in(0,p_-)$,
  an $r_0\in(s_0,\infty)$, and a $C\in(0,\infty)$
  such that, for any $f\in L^1_{\mathrm{loc}}(\rrn)$,
  $$
  \left\|\mc^{((r_0/s_0)')}(f)\right\|_{(X^{1/s_0})'}
  \leq C\left\|f\right\|_{(X^{1/s_0})'}.
  $$
\end{enumerate}
Let $d\geq\lfloor n(1/p_--1)\rfloor$
be a fixed integer, $c\in(0,1]$,
$r\in(\max\{1,r_0\},\infty]$, and
$A,\ \widetilde{A}\in(0,\infty)$.
Assume that
$\{a_{i,j}\}_{i\in\mathbb{Z},\,j\in\mathbb{N}}$
is a sequence of $(X,\,r,\,d)$-atoms supported,
respectively, in the balls $\{B_{i,j}\}_{
i\in\mathbb{Z},\,j\in\mathbb{N}}\subset\mathbb{B}$
such that, for any $i\in\mathbb{Z}$,
$$\sum_{j\in\mathbb{N}}\1bf_{cB_{i,j}}\leq A,$$
$f:=\sum_{i\in\mathbb{Z}}\sum_{j\in\mathbb{N}}
\widetilde{A}2^i\|\1bf_{B_{i,j}}\|_Xa_{i,j}$
converges in $\mathcal{S}'(\rrn)$, and
$$
\sup_{i\in\mathbb{Z}}\left\{
2^i\left\|\sum_{j\in\mathbb{N}}
\1bf_{B_{i,j}}\right\|_X\right\}<\infty.
$$
Then $f\in WH_X(\rrn)$ and
$$
\left\|f\right\|_{WH_X(\rrn)}\lesssim
\sup_{i\in\mathbb{Z}}2^i\left\|\sum_{j\in\mathbb{N}}
\1bf_{B_{i,j}}\right\|_X,
$$
where the implicit positive constant is
independent of $f$.
\end{lemma}

Based on the above two lemmas, we
now show the atomic characterization
of the weak generalized
Herz--Hardy space $W\HaSaHo$.

\begin{proof}[Proof of Theorem \ref{watom}]
Let $p$, $q$, $\omega$, $r$,
$p_-$, and $d$ be as in the present theorem. We first prove
$W\HaSaHo\subset WH\Kmp_{\omega,\0bf}^{p,q,r,d}(\rrn)$.
To this end, let $f\in W\HaSaHo$.
Now, we claim that, for the
above $r$ and $p_-$, all the assumptions of
Lemma \ref{watomd} hold true
with $X:=\HerzSo$. Assume that this claim holds
true for the moment. Then, from
Lemma \ref{watomd} with $X$ therein replaced
by $\HerzSo$, it follows that
there exist
$\{a_{i,j}\}_{i\in\mathbb{Z},\,j\in\mathbb{N}}$
of $(\HerzSo,\,\infty,\,d)$-atoms supported,
respectively, in the balls $\{B_{i,j}\}_{
i\in\mathbb{Z},\,j\in\mathbb{N}}\subset\mathbb{B}$
and three positive constants $c\in(0,1]$,
$A$, and $\widetilde{A}$, independent of $f$,
such that, for any
$i\in\mathbb{Z}$, $$\sum_{j\in\mathbb{N}}
\1bf_{cB_{i,j}}\leq A,$$
\begin{equation}\label{watome1}
f=\sum_{i\in\mathbb{Z}}\sum_{j\in\mathbb{N}}
\widetilde{A}2^i\left\|\1bf_{B_{i,j}}
\right\|_{\HerzSo}a_{i,j}
\end{equation}
in $\mathcal{S}'(\rrn)$, and
\begin{equation}\label{watome2}
\sup_{i\in\mathbb{Z}}\left\{2^i\left\|
\sum_{j\in\mathbb{N}}\1bf_{B_{i,j}}
\right\|_{\HerzSo}\right\}
\lesssim\left\|f\right\|_{W\HaSaHo}
<\infty.
\end{equation}
In addition, by Lemma \ref{Atoll2}
with $t:=\infty$, we conclude that, for any
$i\in\mathbb{Z}$ and $j\in\mathbb{N}$,
$a_{i,j}$ is a $(\HerzSo,\,r,\,d)$-atom
supported in the ball $B_{i,j}$.
This, combined with \eqref{watome1},
\eqref{watome2}, and Definition \ref{dfwatom},
further implies that $f\in WH\Kmp_{\omega,\0bf}
^{p,q,r,d}(\rrn)$ and then
finishes the proof that
$W\HaSaHo\subset
WH\Kmp_{\omega,\0bf}^{p,q,r,d}(\rrn)$.

Thus, to complete the proof that
$W\HaSaHo\subset WH\Kmp_{\omega,\0bf}^{p,q,r,d}(\rrn)$,
it remains to show the above claim.
Indeed, applying the assumption $\m0
\in(-\frac{n}{p},\infty)$ and
Theorem \ref{Th3}, we find that the local
generalized Herz space $\HerzSo$
under consideration is a BQBF space.
Using Lemma \ref{vmbhl} with $u:=\frac{1}{r}$,
we conclude that,
for any $\{f_j\}_{j\in\mathbb{N}}\subset
L^1_{\mathrm{loc}}(\rrn)$,
\begin{equation*}
\left\|\left\{\sum_{j\in\mathbb{N}}
\left[\mc\left(f_j\right)\right]^{1/r}
\right\}^r\right\|_{[\HerzSo]^{1/r}}
\lesssim\left\|\left(\sum_{j\in\mathbb{N}}
\left|f_j\right|^{1/r}\right)^{r}
\right\|_{[\HerzSo]^{1/r}},
\end{equation*}
which implies that the assumption (i)
of Lemma \ref{watomd} holds true for $\HerzSo$.
We next prove that there exist
a $\theta$, an $s$, and an $r_0$
such that (ii), (iii),
and (iv) of Lemma \ref{watomd} hold
true. Indeed, let
$$
\theta\in\left(\max\{1,p,q\},
\frac{1}{p_-}\min
\left\{p,q,\frac{n}{\Mw+n/p}\right\}\right).
$$
Then, from the assumptions
$p/\theta,\ q/\theta\in(0,1)$,
the reverse Minkowski inequality,
and Lemma \ref{convexll},
we deduce that, for any $\{f_{j}\}_{j\in
\mathbb{N}}\subset\Msc(\rrn)$,
\begin{equation}\label{watome4}
\sum_{j\in\mathbb{N}}\|f_{j}\|_{[\HerzSo]^{1/\theta}}\leq
\left\|\sum_{j\in\mathbb{N}}|f_{j}|
\right\|_{[\HerzSo]^{1/\theta}},
\end{equation}
which implies that $\HerzSo$
is $\theta$-concave. On the other hand,
using Lemma \ref{mbhl} with
$r:=\theta p_-$, we conclude that,
for any $f\in L_{{\rm loc}}^{1}(\rrn)$,
\begin{equation*}
\left\|\mc(f)\right\|_{[\HerzSo]^{\frac{1}{\theta p_-}}}
\lesssim\left\|f\right\|_{[\HerzSo]^{\frac{1}{\theta p_-}}}.
\end{equation*}
From this and \eqref{watome4},
we further infer that,
for the above $\theta$ and $p_-$,
Lemma \ref{watomd}(ii) holds true.
Let
$$
s\in\left(0,\min\left\{1,
p,q,\frac{n}{\Mw+n/p}\right\}\right).
$$
Then, by Lemma \ref{mbhal} with
$r:=\infty$, we find that $[\HerzSo]^{1/s}$
is a BBF space and, for any
$f\in L^1_{\mathrm{loc}}(\rrn)$,
$$
\left\|\mc(f)\right\|_{([\HerzSo]^{1/s})'}
\lesssim\left\|f\right\|_{([\HerzSo]^{1/s})'},
$$
which imply that Lemma \ref{watomd}(iii) holds true
with the above $s$.
Finally, let
$$r_0\in
\left(0,\min\left\{p,\frac{n}{\Mw+n/p}\right\}\right).
$$
Then, from Lemma \ref{wmbhl}
with $r:=r_0$, we deduce that, for any
$f\in L^1_{\mathrm{loc}}(\rrn)$,
\begin{equation*}
\left\|\mc(f)\right\|_{[W\HerzSo]^{1/r_0}}\lesssim
\left\|f\right\|_{[W\HerzSo]^{1/r_0}}.
\end{equation*}
This implies that Lemma \ref{watomd}(iv)
holds true with the above $r_0$.
Therefore, the Herz space
$\HerzSo$ under consideration
satisfies all the assumptions of
Lemma \ref{watomd} with
$r$ and $p_-$ in the present theorem.
This then finishes the proof of the above
claim and further implies that
$$
W\HaSaHo\subset WH\Kmp_{\omega,\0bf}^{
p,q,r,d}(\rrn).
$$
Moreover, combining \eqref{watome1},
\eqref{watome2}, and Definition \ref{dfwatom}
again, we conclude that, for any $f\in W\HaSaHo$,
\begin{equation}\label{watome5}
\left\|f\right\|_{WH\Kmp_{\omega,\0bf}^{
p,q,r,d}(\rrn)}\lesssim\left\|f\right\|_{W\HaSaHo}.
\end{equation}

Conversely, we now prove that
$WH\Kmp_{\omega,\0bf}^{
p,q,r,d}(\rrn)\subset W\HaSaHo$.
Indeed, from the proof that $W\HaSaHo
\subset WH\Kmp_{\omega,\0bf}^{
p,q,r,d}(\rrn)$,
it follows that the local generalized
Herz space $\HerzSo$ under consideration
is a BQBF space. Thus, in order to
finish the proof that
$WH\Kmp_{\omega,\0bf}^{
p,q,r,d}(\rrn)\subset W\HaSaHo$,
we only need to show that the
assumptions (i) through (iii) of
Lemma \ref{watomr} hold true
for $\HerzSo$. Indeed, for any
given $\theta\in(0,p_-)$ and
$u\in(1,\infty)$, applying
Lemma \ref{vmbhl} with $r:=\theta$, we find
that, for any $\{f_{j}\}_{j\in\mathbb{N}}
\subset L_{{\rm loc}}^{1}(\rrn)$,
\begin{equation*}
\left\|\left\{\sum_{j\in\mathbb{N}}\left[\mc
(f_{j})\right]^{u}\right\}^{\frac{1}{u}}
\right\|_{[\HerzSo]^{1/\theta}}
\lesssim\left\|\left(
\sum_{j\in\mathbb{N}}|f_{j}|
^{u}\right)^{\frac{1}{u}}\right\|_
{[\HerzSo]^{1/\theta}}.
\end{equation*}
This implies that Lemma \ref{watomr}(i)
holds true for $\HerzSo$.
In addition, from Lemma \ref{mbhal},
we deduce that, for any $s\in(0,p_-)$,
$[\HerzSo]^{1/s}$ is a BBF space.
This implies that Lemma \ref{watomr}(ii)
holds true.
Finally, we show that there exist
an $s_0$ and an $r_0$ such that Lemma
\ref{watomr}(iii) holds true with
these $s_0$ and $r_0$.
To do this, let $s_0\in(0,p_-)$ and
$$
r_0\in\left(\max\left\{1,p,
\frac{n}{\Mw+n/p}\right\},r\right).
$$
Then, using Lemma \ref{mbhal}
again with $s$ and $r$ therein
replaced, respectively, by
$s_0$ and $r_0$, we conclude that,
for any $f\in L_{{\rm loc}}^{1}(\rrn)$,
\begin{equation*}
\left\|\mc^{((r_0/s_0)')}(f)\right\|_
{([\HerzSo]^{1/s_0})'}\lesssim
\left\|f\right\|_{([\HerzSo]^{1/s_0})'},
\end{equation*}
which implies that, for the above
$s_0$ and $r_0$, Lemma \ref{watomr}(iii)
holds true for $\HerzSo$. Therefore,
all the assumptions of Lemma \ref{watomr}
hold true for the Herz space $\HerzSo$
under consideration. This then implies
that $WH\Kmp_{\omega,\0bf}^{
p,q,r,d}(\rrn)\subset W\HaSaHo$.
Moreover, using Definition \ref{dfwatom}
and Lemma \ref{watomr} again,
we further find that, for any $f\in
WH\Kmp_{\omega,\0bf}^{p,q,r,d}(\rrn)$,
$$
\left\|f\right\|_{W\HaSaHo}
\lesssim\left\|f\right\|_
{WH\Kmp_{\omega,\0bf}^{p,q,r,d}(\rrn)},
$$
which, together with \eqref{watome5},
further implies that
$$
W\HaSaHo=WH\Kmp_{\omega,\0bf}^{p,q,r,d}(\rrn)
$$
with equivalent quasi-norms.
This finishes the proof of Theorem
\ref{watom}.
\end{proof}

As an application,
we now establish the atomic
characterization of the weak
generalized Morrey--Hardy space
$WH\MorrSo$. Recall that the definition
of the $(\MorrSo,\,r,\,d)$-atoms
is as in Definition \ref{atomm}.
Then, applying Theorem \ref{watom}
and Remark \ref{remark6.16},
we obtain the following atomic
characterization of the weak
generalized Morrey--Hardy space
$WH\MorrSo$ immediately; we omit the
details.\index{atomic characterization}

\begin{corollary}\label{watomm}
Let $p,\ q\in[1,\infty)$,
$p_-\in(0,\min\{p,q\}/\max\{p,q\})$,
$\omega\in M(\rp)$ with
$$
-\frac{n}{p}<\m0(\omega)\leq\M0(\omega)<0
$$
and
$$
-\frac{n}{p}<\mi(\omega)\leq\MI(\omega)<0,
$$
$d\geq\lfloor n(1/p_--1)\rfloor$ be a fixed
integer, and
$$
r\in\left(\frac{n}{\mw+n/p},\infty\right].
$$
Then the \emph{weak generalized atomic Morrey--Hardy space}
\index{weak generalized atomic Morrey--Hardy space}$WH
\textsl{\textbf{M}}_{\omega,\0bf}^{p,q,r,d}(\rrn)$\index{$WH
\textsl{\textbf{M}}_{\omega,\0bf}^{p,q,r,d}(\rrn)$}, associated with
the weak local generalized
Morrey space $W\MorrSo$, is defined to
be the set of all the
$f\in\mathcal{S}'(\rrn)$ such that
there exist a sequence $\{a_{i,j}\}_{i\in
\mathbb{Z},\,j\in\mathbb{N}}$ of
$(\MorrSo,\,r,\,d)$-atoms supported,
respectively, in the balls
$\{B_{i,j}\}_{i\in\mathbb{Z},\,j\in\mathbb{N}}
\subset\mathbb{B}$ and three
positive constants $c\in(0,1]$,
$A$, and $\widetilde{A}$,
independent of $f$, satisfying that,
for any $i\in\mathbb{Z}$,
$$
\sum_{j\in\mathbb{N}}\1bf_{cB_{i,j}}\leq A,
$$
$$
f=\sum_{i\in\mathbb{Z}}\sum_{j\in\mathbb{N}}
\widetilde{A}2^i\left\|\1bf_{B_{i,j}}\right\|_{
\MorrSo}a_{i,j}
$$
in $\mathcal{S}'(\rrn)$, and
$$
\sup_{i\in\mathbb{Z}}
\left\{2^i\left\|\sum_{j\in\mathbb{N}}
\1bf_{B_{i,j}}\right\|_{\MorrSo}\right\}<\infty.
$$
Moreover, for any
$f\in WH\textsl{\textbf{M}}_{
\omega,\0bf}^{p,q,r,d}(\rrn)$,
$$
\|f\|_{WH\textsl{\textbf{M}}_{
\omega,\0bf}^{p,q,r,d}(\rrn)}:=
\inf\sup_{i\in\mathbb{Z}}\left\{
2^i\left\|\sum_{j\in\mathbb{N}}
\1bf_{B_{i,j}}\right\|_{\MorrSo}\right\},
$$
where the infimum is taken over
all the decompositions of
$f$ as above.
Then $$
WH\MorrSo=WH\textsl{\textbf{M}
}_{\omega,\0bf}^{p,q,r,d}(\rrn)
$$
with equivalent quasi-norms.
\end{corollary}

The remainder of this section is
devoted to establishing the atomic
characterization of the weak
generalized Herz--Hardy space $W\HaSaH$.
To this end, we first
introduce the following definition of weak
generalized atomic Herz--Hardy spaces
via $(\HerzS,\,r,\,d)$-atoms introduced
in Definition \ref{deatomg}.

\begin{definition}\label{dfwatomg}
Let $p,\ q\in(0,\infty)$, $\omega\in M(\rp)$ with
$\m0(\omega)\in(-\frac{n}{p},\infty)$ and
$$-\frac{n}{p}<\mi(\omega)
\leq\MI(\omega)<0,
$$
$$
p_-\in\left(0,\frac{\min\{p,q,\frac{n}
{\Mw+n/p}\}}{\max\{1,p,q\}}\right),
$$
$d\geq\lfloor n(1/p_--1)\rfloor$ be a fixed
integer, and
$$
r\in\left(\max\left\{1,p,\frac{n}{\mw+n/p}\right\},\infty\right].
$$
Then the \emph{weak generalized atomic Herz--Hardy space}
\index{weak generalized atomic Herz--Hardy \\ space}$WH
\Kmp_{\omega}^{p,q,r,d}(\rrn)$\index{$WH\Kmp_{\omega}^{p,q,r,d}(\rrn)$},
associated with
the weak global generalized Herz space $W\HerzS$, is defined to
be the set of all the
$f\in\mathcal{S}'(\rrn)$ such that
there exist a sequence $\{a_{i,j}\}_{i\in
\mathbb{Z},\,j\in\mathbb{N}}$ of
$(\HerzS,\,r,\,d)$-atoms supported,
respectively, in the balls
$\{B_{i,j}\}_{i\in\mathbb{Z},\,j
\in\mathbb{N}}\subset\mathbb{B}$
and three positive constants $c\in(0,1]$,
$A$, and $\widetilde{A}$,
independent of $f$, satisfying that,
for any $i\in\mathbb{Z}$,
$$
\sum_{j\in\mathbb{N}}\1bf_{cB_{i,j}}\leq A,
$$
$$
f=\sum_{i\in\mathbb{Z}}\sum_{j\in\mathbb{N}}
\widetilde{A}2^i\left\|\1bf_{B_{i,j}}
\right\|_{\HerzS}a_{i,j}
$$
in $\mathcal{S}'(\rrn)$, and
$$
\sup_{i\in\mathbb{Z}}\left\{2^i
\left\|\sum_{j\in\mathbb{N}}
\1bf_{B_{i,j}}\right\|_{\HerzS}\right\}<\infty,
$$
Moreover, for any $f\in
WH\Kmp_{\omega}^{p,q,r,d}(\rrn)$,
$$
\|f\|_{WH\Kmp_{\omega}^{p,q,r,d}(\rrn)}:=
\inf\sup_{i\in\mathbb{Z}}\left\{2^i
\left\|\sum_{j\in\mathbb{N}}
\1bf_{B_{i,j}}\right\|_{\HerzS}\right\},
$$
where the infimum is taken over all the
decompositions of $f$ as above.
\end{definition}

Then we have the following
atomic characterization of the
weak generalized Herz--Hardy space $W\HaSaH$.

\begin{theorem}\label{watomg}
Let $p$, $q$, $\omega$, $r$, and $d$
be as in Definition \ref{dfwatomg}.
Then
$$
W\HaSaH=WH\Kmp_{\omega}^{p,q,r,d}(\rrn)
$$
with equivalent quasi-norms.
\end{theorem}

Due to the deficiency of associate spaces
of global generalized Herz spaces, we can
not prove this theorem via
the atomic characterization of
weak Hardy spaces associated with
ball quasi-Banach function spaces
obtained in \cite[Theorems 4.2 and 4.7]{ZWYY}
(see also both Lemmas \ref{watomd} and
\ref{watomr} above) directly.
To overcome this difficulty, we
establish an improved
atomic characterization
of $WH_X(\rrn)$ associated with the
ball quasi-Banach function space $X$
without recourse to the associate
space $X'$. Indeed, we have
the following atomic decomposition theorem.

\begin{proposition}\label{watomgd}
Let $X$ be a ball quasi-Banach function
space and $Y\subset\Msc(\rrn)$ a linear
space equipped with a quasi-seminorm
$\|\cdot\|_Y$, and let $r\in(0,1)$,
$r_0,\ s,\ p_-\in(0,\infty)$, and
$\theta,\ \theta_0\in(1,\infty)$ be such that
the following seven statements hold true:
\begin{enumerate}
  \item[{\rm(i)}] there exists a positive constant
  $C$ such that, for any $\{f_j\}_{j\in\mathbb{N}}\subset
  L^1_{\mathrm{loc}}(\rrn)$,
  $$
  \left\|\left\{\sum_{j\in\mathbb{N}}
  \left[\mc\left(f_j\right)\right]^{1/r}
  \right\}^r\right\|_{X^{1/r}}
  \leq C\left\|\left(\sum_{j\in\mathbb{N}}
  \left|f_j\right|^{1/r}\right)^{r}\right\|_{X^{1/r}};
  $$
  \item[{\rm(ii)}] $X$ is $\theta$-concave and
  $\mc$ bounded on $X^{\frac{1}{\theta p_-}}$;
  \item[{\rm(iii)}] $\mc$ is bounded on $(WX)^{1/r_0}$;
  \item[{\rm(iv)}] $\|\cdot\|_Y$ satisfies
  Definition \ref{Df1}(ii);
  \item[{\rm(v)}] $\1bf_{B(\0bf,1)}\in Y$;
  \item[{\rm(vi)}] for any $f\in\Msc(\rrn)$,
  $$
  \|f\|_{X^{1/s}}
  \sim\sup\left\{\|fg\|_{L^1(\rrn)}:\
  \|g\|_{Y}=1\right\},
  $$
  where the positive equivalence
  constants are independent of $f$;
  \item[{\rm(vii)}] $\mc^{(\theta_0)}$ is
  bounded on $Y$.
\end{enumerate}
Assume that $d\geq\lfloor n(1/p_--1)\rfloor$
is a fixed nonnegative integer and $f\in WH_X(\rrn)$.
Then there exist
$\{a_{i,j}\}_{i\in\mathbb{Z},\,j\in\mathbb{N}}$
of $(X,\,\infty,\,d)$-atoms supported,
respectively, in the balls $\{B_{i,j}\}_{
i\in\mathbb{Z},\,j\in\mathbb{N}}\subset\mathbb{B}$
and three positive constants $c\in(0,1]$,
$A$, and $\widetilde{A}$, independent of $f$,
such that, for any
$i\in\mathbb{Z}$, $$\sum_{j\in\mathbb{N}}
\1bf_{cB_{i,j}}\leq A,$$
$$f=\sum_{i\in\mathbb{Z}}\sum_{j\in\mathbb{N}}
\widetilde{A}2^i\left\|\1bf_{B_{i,j}}
\right\|_Xa_{i,j}$$
in $\mathcal{S}'(\rrn)$, and
$$
\sup_{i\in\mathbb{Z}}\left\{
2^i\left\|\sum_{j\in\mathbb{N}}
\1bf_{B_{i,j}}\right\|_X\right\}\lesssim
\left\|f\right\|_{WH_X(\rrn)},
$$
where the implicit positive constant is independent
of $f$.
\end{proposition}

\begin{remark}
We should point out that
Proposition \ref{watomgd}
is an improved version of the
known atomic decomposition
of $WH_X(\rrn)$ obtained
in \cite[Theorem 4.2]{ZWYY}. Indeed,
if $Y\equiv(X^{1/s})'$ in Proposition
\ref{watomgd}, then this proposition goes
back to \cite[Theorem 4.2]{ZWYY}.
\end{remark}

To show this proposition, we first
prove the following technique
lemma by borrowing some
ideas from the proof of
\cite[Lemma 4.3]{ZWYY}, which shows
a crucial fact that
distributions in weak Hardy spaces vanish
weakly at infinity.

\begin{lemma}\label{WX-weak vanish}
Let $X$ be a ball quasi-Banach
function space, $Y\subset\Msc(\rrn)$ a
linear space equipped with a
quasi-seminorm $\|\cdot\|_Y$,
$\theta\in(1,\infty)$, and $s\in(0,\infty)$
satisfy the following four statements:
\begin{enumerate}
  \item[{\rm(i)}] $\|\cdot\|_Y$ satisfies
  Definition \ref{Df1}(ii);
  \item[{\rm(ii)}] $\1bf_{B(\0bf,1)}\in Y$;
  \item[{\rm(iii)}] for any $f\in\Msc(\rrn)$,
  $$
  \|f\|_{X^{1/s}}
  \sim\sup\left\{\|fg\|_{L^1(\rrn)}:\
  \|g\|_{Y}=1\right\},
  $$
  where the positive equivalence
  constants are independent of $f$;
  \item[{\rm(iv)}] $\mc^{(\theta)}$ is
  bounded on $Y$.
\end{enumerate}
Assume $f\in WH_X(\rrn)$. Then
$f$ vanishes weakly at infinity.
\end{lemma}

In order to prove Lemma
\ref{WX-weak vanish}, we require
the following auxiliary estimate about
$A_p$-weights
(see, for instance,
\cite[Corollary 7.6(3)]{d01}).

\begin{lemma}\label{reverse Holder}
Let $p\in[1,\infty)$ and
$\upsilon\in A_p(\rrn)$. Then
there exist two positive
constants $\delta$ and $C$
such that, for any ball $B_1\in\mathbb{B}$
and any ball $B_2\subset B_1$,
$$
\frac{\upsilon(B_2)}{\upsilon(B _1)}
\leq C\left(\frac{|B_2|}{|B_1|}\right)^{
\delta}.
$$
\end{lemma}

We next show Lemma \ref{WX-weak vanish}.

\begin{proof}[Proof of Lemma
\ref{WX-weak vanish}]
Let all the symbols be as in the present
lemma. Then, from the assumptions
(i) through (iv) of the present proposition
and Lemma \ref{bddol1},
we deduce that there exists
an $\eps\in(0,1)$ such that, for any
$g\in\Msc(\rrn)$,
\begin{equation}\label{wxe1}
\left\|g\right\|_{L^{s}_{\upsilon}(\rrn)}
\lesssim\left\|g\right\|_X,
\end{equation}
where $\upsilon:=
[\mc(\1bf_{B(\0bf,1)})]^{\eps}$.
Let $\varphi\in\mathcal{S}(\rrn)$ and $N\in\mathbb{N}$.
Then, by \eqref{sec6e1},
we conclude that,
for any $t\in(0,\infty)$ and
$x,\ y\in\rrn$ with $|y-x|<t$,
$$
\left|f\ast\varphi_t(x)\right|
\lesssim\mc_N(f)(y),
$$
which further implies that
$$
\left|f\ast\varphi_t(x)\right|
\1bf_{B(x,t)}\lesssim\mc_N(f)\1bf_{B(x,t)}.
$$
Applying this, Definitions \ref{df237}(ii)
and \ref{de-bqbfwx}, \eqref{wxe1},
Remark \ref{wball} with
$X:=L^{s}_{\upsilon}(\rrn)$, and Definition
\ref{hardyxw}, we find that,
for any $t\in(0,\infty)$ and $x\in\rrn$,
\begin{align*}
\left|f\ast\varphi_t(x)\right|
&\lesssim\frac{\|\mc_N(f)
\1bf_{B(x,t)}\|_{WL^s_{\upsilon}(\rrn)}}
{\|\1bf_{B(x,t)}\|_{WL^s_{\upsilon}(\rrn)}}\\
&\lesssim\frac{\|\mc_N(f)\|_{WX}}{\|\1bf_{B(x,t)}
\|_{L^s_{\upsilon}(\rrn)}}
\sim\frac{\|f\|_{WH_X(\rrn)}}{\|\1bf_{B(x,t)}
\|_{L^s_{\upsilon}(\rrn)}}.
\end{align*}
Therefore, for any given $\phi\in\mathcal{S}(\rrn)$
and for any $t\in(0,\infty)$,
we have
\begin{align}\label{wxe2}
&\left|\int_{\rrn}f\ast\varphi_t(x)
\phi(x)\,dx\right|\notag\\&\quad
\lesssim\int_{\rrn}\frac{1}
{\|\1bf_{B(x,t)}\|_{L^s_{\upsilon}(\rrn)}}\left|
\phi(x)\right|\,dx\notag\\
&\quad\sim\int_{\rrn}\frac{1}{\|\1bf_{B(x,1)}
\|_{L^s_{\upsilon}(\rrn)}}
\frac{\|\1bf_{B(x,1)}\|_{L^s_{\upsilon}(\rrn)}}
{\|\1bf_{B(x,t)}\|_{L^s_{\upsilon}(\rrn)}}
\left|\phi(x)\right|\,dx.
\end{align}
Moreover, using Lemma \ref{2.1.6}
with $r:=1$, we find that, for any $x\in\rrn$,
\begin{align}\label{wxe3}
\left\|\1bf_{B(x,1)}
\right\|_{L^s_{\upsilon}(\rrn)}
&=\left\{\int_{B(x,1)}\left[\mc\left(
\1bf_{B(\0bf,1)}\right)(y)
\right]^{\eps}\,dy\right\}^{\frac{1}{s}}\notag\\
&\sim\left[\int_{B(x,1)}\left(1+|y|
\right)^{-n\eps}\,dy
\right]^{\frac{1}{s}}\notag\\
&\gtrsim\left(1+|x|
\right)^{-\frac{n\eps}{s}}
\left|B(x,1)\right|^{\frac{1}{s}}
\sim\left(1+|x|
\right)^{-\frac{n\eps}{s}}.
\end{align}
On the other hand, from
Lemma \ref{czool2} with $f:=\1bf_{B(\0bf,1)}$
and $\delta:=\eps$, it follows that
$\upsilon\in A_1(\rrn)$. Then,
for any $t\in[1,\infty)$
and $x\in\rrn$,
by Lemma \ref{reverse Holder}
with $p:=1$, $B_1:=B(x,t)$, and
$B_2:=B(x,1)$, we find that there
exists a $\delta\in(0,\infty)$,
independent of $t$ and $x$, such that
\begin{align*}
\frac{\|\1bf_{B(x,1)}\|_{L^s_{\upsilon}(\rrn)}}
{\|\1bf_{B(x,t)}\|_{L^s_{\upsilon}(\rrn)}}
\sim\left[\frac{\upsilon(B(x,1))}{\upsilon(
B(x,t))}\right]^{\frac{1}{s}}
\lesssim\left[\frac{|B(x,1)|}{|B(x,t)|}
\right]^{\frac{\delta}{s}}
\sim t^{-\frac{n\delta}{s}}.
\end{align*}
Combining this, \eqref{wxe2}, and
\eqref{wxe3}, we further conclude that,
for any given $\phi\in\mathcal{S}(\rrn)$,
\begin{align*}
&\left|\int_{\rrn}f\ast\varphi_t(x)
\phi(x)\,dx\right|\\
&\quad\lesssim t^{-\frac{n\delta}{s}}
\int_{\rrn}\left(1+|x|\right)^{\frac{n\eps}{s}}
\left|\phi(x)\right|\,dx
\lesssim t^{-\frac{n\delta}{s}}\to0
\end{align*}
as $t\to\infty$. This implies that
$f\ast\varphi_t\to0$ in $\mathcal{S}'(\rrn)$
as $t\to\infty$, and hence $f$ vanishes weakly
at infinity, which then completes the proof
of Lemma \ref{WX-weak vanish}.
\end{proof}

Applying Lemma \ref{WX-weak vanish},
we now prove the atomic
decomposition of $WH_X(\rrn)$ as follows.

\begin{proof}[Proof of Proposition \ref{watomgd}]
Let all the symbols be as in the
present proposition and $f\in WH_X(\rrn)$.
Then, by the assumptions (iv) through (vii)
of the present proposition and
Lemma \ref{WX-weak vanish},
we find that $f$ vanishes weakly
at infinity. Using this,
(i) through (iii) of the present
proposition, and repeating
the proof of \cite[Theorem 4.2]{ZWYY},
we conclude that there exist a sequence
$\{a_{i,j}\}_{i\in\mathbb{Z},\,j\in\mathbb{N}}$
of $(X,\,\infty,\,d)$-atoms supported,
respectively, in the balls $\{B_{i,j}\}_{
i\in\mathbb{Z},\,j\in\mathbb{N}}\subset\mathbb{B}$
and three positive constants $c\in(0,1]$,
$A$, and $\widetilde{A}$, independent of $f$,
such that, for any
$i\in\mathbb{Z}$, $$\sum_{j\in\mathbb{N}}
\1bf_{cB_{i,j}}\leq A,$$
$$f=\sum_{i\in\mathbb{Z}}\sum_{j\in\mathbb{N}}
\widetilde{A}2^i\left\|\1bf_{B_{i,j}}
\right\|_Xa_{i,j}$$
in $\mathcal{S}'(\rrn)$, and
$$
\sup_{i\in\mathbb{Z}}\left\{
2^i\left\|\sum_{j\in\mathbb{N}}
\1bf_{B_{i,j}}\right\|_X\right\}\lesssim
\left\|f\right\|_{WH_X(\rrn)}.
$$
This finishes the proof of Proposition
\ref{watomgd}.
\end{proof}

On the other hand, we establish
the following atomic reconstruction
theorem of $WH_X(\rrn)$ without
recourse to associate spaces.

\begin{proposition}\label{watomgr}
Let $X$ be a ball quasi-Banach function space
and let $p_-\in(0,1)$ be such that the following
three statements hold true:
\begin{enumerate}
  \item[{\rm(i)}] for any given
  $\theta\in(0,p_-)$ and $u\in(1,\infty)$,
  there exists a positive constant $C$ such that,
for any $\{f
_{j}\}_{j\in\mathbb{N}}\subset L_{{\rm loc}}^{1}(\rrn)$,
\begin{equation*}
\left\|\left\{\sum_{j\in\mathbb{N}}\left[\mc
(f_{j})\right]^{u}\right\}^
{\frac{1}{u}}\right\|_{X^{1/\theta}}
\leq C\left\|\left\{
\sum_{j\in\mathbb{N}}|f_{j}|^{u}\right\}
^{\frac{1}{u}}\right\|_{X^{1/\theta}};
\end{equation*}
  \item[{\rm(ii)}] for any $s\in(0,p_-)$,
  $X^{1/s}$ is a ball Banach function space;
  \item[{\rm(iii)}] there exist an $s_0\in(0,p_-)$,
  an $r_0\in(s_0,\infty)$, a $C\in(0,\infty)$,
  and a linear space $Y_0\subset\Msc(\rrn)$
  equipped with a quasi-seminorm $\|\cdot\|_{Y_0}$
  such that, for any
  $f\in\Msc(\rrn)$,
  $$
  \left\|f\right\|_{X^{1/s_0}}
  \sim\sup\left\{\left\|fg\right\|_{L^1(\rrn)}:\
  \left\|g\right\|_{Y_0}=1\right\}
  $$
  with the positive equivalence constants
  independent of $f$
  and, for any $f\in L^1_{\mathrm{loc}}(\rrn)$,
  $$
  \left\|\mc^{((r_0/s_0)')}(f)\right\|_{Y_0}
  \leq C\left\|f\right\|_{Y_0}.
  $$
\end{enumerate}
Let $d\geq\lfloor n(1/p_--1)\rfloor$
be a fixed integer, $c\in(0,1]$,
$r\in(\max\{1,r_0\},\infty]$, and
$A,\ \widetilde{A}\in(0,\infty)$.
Assume that
$\{a_{i,j}\}_{i\in\mathbb{Z},\,j\in\mathbb{N}}$
is a sequence of $(X,\,r,\,d)$-atoms supported,
respectively, in the balls $\{B_{i,j}\}_{
i\in\mathbb{Z},\,j\in\mathbb{N}}\subset\mathbb{B}$
such that, for any $i\in\mathbb{Z}$,
$$\sum_{j\in\mathbb{N}}
\1bf_{cB_{i,j}}\leq A,$$
$$f:=\sum_{i\in\mathbb{Z}}\sum_{j\in\mathbb{N}}
\widetilde{A}2^i\left\|\1bf_{B_{i,j}}
\right\|_Xa_{i,j}$$
converges in $\mathcal{S}'(\rrn)$, and
$$
\sup_{i\in\mathbb{Z}}\left\{
2^i\left\|\sum_{j\in\mathbb{N}}
\1bf_{B_{i,j}}\right\|_X\right\}<\infty.
$$
Then $f\in WH_X(\rrn)$ and
$$
\left\|f\right\|_{WH_X(\rrn)}\lesssim
\sup_{i\in\mathbb{Z}}2^i\left\|\sum_{j\in\mathbb{N}}
\1bf_{B_{i,j}}\right\|_X,
$$
where the implicit positive constant is
independent of $f$.
\end{proposition}

\begin{remark}
We should point out that Proposition
\ref{watomgr} is an improved version
of the known atomic reconstruction
of $WH_X(\rrn)$ established in
\cite[Theorem 4.8]{ZWYY}. Indeed,
if $Y_0\equiv(X^{1/s_0})'$ in Proposition
\ref{watomgr}, then this proposition
goes back to \cite[Theorem 4.8]{ZWYY}.
\end{remark}

To prove Proposition \ref{watomgr},
we first establish the following
technique lemma via
borrowing some ideas from \cite[Lemma 4.8]{ZWYY}.

\begin{lemma}\label{at-es}
Let $s\in(0,\infty)$, $r\in(s,\infty]$,
$X$ be a ball quasi-Banach function
space, and $Y\subset\Msc(\rrn)$ a linear
space equipped with a quasi-seminorm $\|\cdot\|_Y$.
Assume that the following two statements hold true:
\begin{enumerate}
  \item[{\rm(i)}] for any $f\in\Msc(\rrn)$,
  $$
  \|f\|_{X^{1/s}}
  \sim\sup\left\{\|fg\|_{L^1(\rrn)}:\
  \|g\|_{Y}=1\right\},
  $$
  where the positive equivalence
  constants are independent of $f$;
  \item[{\rm(ii)}] there exists a
  positive constant $C$ such that,
  for any $f\in L^1_{
  \mathrm{loc}}(\rrn)$,
  $$
  \left\|\mc^{((r/s)')}(f)\right\|_{Y}
  \leq C\left\|f\right\|_Y.
  $$
\end{enumerate}
Let $\{\lambda_j\}
_{j\in\mathbb{N}}\subset[0,\infty)$,
$\{B_j\}_{j\in\mathbb{N}}\subset\mathbb{B}$,
and $\{a_j\}_{j\in\mathbb{N}}\subset\Msc(\rrn)$
satisfy that, for any $j\in\mathbb{N}$,
$\supp(a_j)\subset B_j$ and $\|a_j\|_{L^r(\rrn)}
\leq|B_j|^{1/r}$. Then there
exists a positive constant $C$, independent of
$\{\lambda_j\}_{j\in\mathbb{N}}$,
$\{B_j\}_{j\in\mathbb{N}}$, and
$\{a_j\}_{j\in\mathbb{N}}$, such that
$$
\left\|\left(\sum_{j\in\mathbb{N}}
\left|\lambda_ja_j\right|^s
\right)^{\frac{1}{s}}\right\|_X
\leq C\left\|\left(\sum_{j\in\mathbb{N}}
\left|\lambda_j\1bf_{B_j}\right|^s
\right)^{\frac{1}{s}}\right\|_X
$$
\end{lemma}

\begin{proof}
Let all the symbols be as in the
present lemma and $g\in\Msc(\rrn)$
with $\|g\|_Y=1$. Then, by
the Tonelli theorem, the assumption
that, for any $j\in\mathbb{N}$,
$\supp(a_j)\subset B_j$, the H\"{o}lder
inequality, and the assumption that,
for any $j\in\mathbb{N}$, $\|a_j\|_{L^r(\rrn)}
\leq|B_j|^{1/r}$,
we find that
\begin{align}\label{at-ese1}
&\left|\int_{\rrn}\sum_{j\in\mathbb{N}}
\left|\lambda_ja_j(x)\right|^sg(x)
\,dx\right|\notag\\
&\quad\leq\sum_{j\in\mathbb{N}}
\lambda_j^s\int_{B_j}\left|a_j(x)\right|^s
\left|g(x)\right|\,dx\notag\\
&\quad\leq\sum_{\rrn}\lambda_j^s
\left\|\left|a_j\right|^s\right\|_{L^{r/s}(\rrn)}
\left\|g\1bf_{B_j}\right\|_{L^{(r/s)'}(\rrn)}\notag\\
&\quad=\sum_{j\in\mathbb{N}}\lambda_j^s
\left\|a_j\right\|^s_{L^r(\rrn)}
\left\|g\1bf_{B_j}\right\|_{L^{(r/s)'}(\rrn)}\notag\\
&\quad\leq\sum_{j\in\mathbb{N}}
\lambda_j^s\left|B_j\right|^{\frac{s}{r}}
\left\|g\1bf_{B_j}\right\|_{L^{(r/s)'}(\rrn)}.
\end{align}
Observe that, from \eqref{hlmax} and
\eqref{hlmaxp} with $\theta:=(r/s)'$ and
$f:=g$, it follows that, for any
$j\in\mathbb{N}$ and $x\in B_j$,
\begin{align}\label{at-ese2}
\mc^{((r/s)')}(g)(x)&\geq
\left[\frac{1}{|B_j|}\int_{B_j}
\left|g(y)\right|^{(r/s)'}
\,dy\right]^{1/(r/s)'}\notag\\
&\sim\left|B_j\right|^{-1/(r/s)'}
\left\|g\1bf_{B_j}
\right\|_{L^{(r/s)'}(\rrn)}.
\end{align}

On the other hand, by the assumption
(ii), we find that, for any $f$,
$g\in\Msc(\rrn)$,
$$
\left\|fg\right\|_{L^1(\rrn)}
\lesssim\left\|f\right\|_{X^{1/s}}
\left\|g\right\|_Y.
$$
This, combined with \eqref{at-ese1},
\eqref{at-ese2}, the Tonelli theorem,
the assumption (ii) of
the present lemma, and the assumption that
$\|g\|_Y=1$, further implies that
\begin{align*}
&\left|\int_{\rrn}\sum_{j\in\mathbb{N}}
\left|\lambda_ja_j(x)\right|^sg(x)
\,dx\right|\\
&\quad\lesssim\sum_{j\in\mathbb{N}}
\lambda_j^s\int_{B_j}\left|B_{j}
\right|^{-1/(r/s)'}\left\|g\1bf_{B_j}
\right\|_{L^{(r/s)'}(\rrn)}\,dx\\
&\quad\lesssim\sum_{j\in\mathbb{N}}
\lambda_j^s\int_{B_j}\mc^{((r/s)')}(g)(x)\,dx\\
&\quad\sim\int_{\rrn}\sum_{j\in\mathbb{N}}
\left|\lambda_j\1bf_{B_j}(x)\right|^s
\mc^{((r/s)')}(g)(x)\,dx\\
&\quad\lesssim\left\|\sum_{j\in\mathbb{N}}
\left|\lambda_j\1bf_{B_j}\right|^s
\right\|_{X^{1/s}}\left\|\mc^{((r/s)')}
(g)\right\|_Y
\lesssim\left\|\left(\sum_{j\in\mathbb{N}}
\left|\lambda_j\1bf_{B_j}\right|^s
\right)^{\frac{1}{s}}\right\|_{X}^s.
\end{align*}
Applying this, Definition \ref{convex}(i)
with $p:=\frac{1}{s}$,
the arbitrariness of
$g$, and the assumption (i) of the present lemma,
we conclude that
\begin{align*}
\left\|\left(\sum_{j\in\mathbb{N}}
\left|\lambda_ja_{j}\right|^s
\right)^{\frac{1}{s}}\right\|_{X}
=\left\|\sum_{j\in\mathbb{N}}
\left|\lambda_ja_{j}\right|^s
\right\|_{X^{1/s}}^{\frac{1}{s}}
\lesssim\left\|\left(\sum_{j\in\mathbb{N}}
\left|\lambda_j\1bf_{B_j}\right|^s
\right)^{\frac{1}{s}}\right\|_{X}.
\end{align*}
This then finishes the proof of
Lemma \ref{at-es}.
\end{proof}

Via Lemma \ref{at-es}, we
now prove the atomic reconstruction
theorem.

\begin{proof}[Proof of Proposition
\ref{watomgr}]
Let $X$, $r$, $d$, $c$,
$A$, and $\widetilde{A}$
be as in the present proposition and
$\{a_{i,j}\}_{i\in\mathbb{Z},\,
j\in\mathbb{N}}$ a sequence of $(X,\,r,\,d)$-atoms
supported, respectively, in the balls
$\{B_{i,j}\}_{i\in\mathbb{Z},\,j\in\mathbb{N}}
\subset\mathbb{B}$ such that,
for any $i\in\mathbb{Z}$, $$\sum_{j\in\mathbb{N}}
\1bf_{cB_{i,j}}\leq A,$$
$$\sum_{i\in\mathbb{Z}}\sum_{j\in\mathbb{N}}
\widetilde{A}2^i\left\|\1bf_{B_{i,j}}
\right\|_Xa_{i,j}$$
converges in $\mathcal{S}'(\rrn)$, and
$$
\sup_{i\in\mathbb{Z}}\left\{
2^i\left\|\sum_{j\in\mathbb{N}}\1bf_{B_{i,j}}
\right\|_X\right\}<\infty.
$$
Then, repeating the proof
of \cite[Theorem 4.7]{ZWYY} via
replacing \cite[Lemma 4.8]{ZWYY} therein
by Lemma \ref{at-es} here,
we obtain $f\in WH_X(\rrn)$
and
$$
\left\|f\right\|_{WH_X(\rrn)}
\lesssim\sup_{i\in\mathbb{Z}}
2^i\left\|\sum_{j\in\mathbb{N}}\1bf_{B_{i,j}}
\right\|_X,
$$
which then completes the proof
of Proposition \ref{watomgr}.
\end{proof}

Via the improved atomic characterization
of $WH_X(\rrn)$ established in both
Propositions \ref{watomgd} and
\ref{watomgr} above, we next
show the atomic characterization of
the weak generalized Herz--Hardy space
$W\HaSaH$.

\begin{proof}[Proof of Theorem
\ref{watomg}]
Let $p$, $q$, $\omega$, $r$,
$p_-$, and $d$ be as in the
present theorem. We first
prove that $W\HaSaH\subset
WH\Kmp_{\omega}^{p,q,r,d}(\rrn)$.
To this end, we show that,
the global generalized
Herz space $\HerzS$ under consideration
satisfies all the assumptions of Proposition
\ref{watomgd}. Indeed, from the assumptions
$\m0(\omega)\in(-\frac{n}{p},\infty)$
and $\MI(\omega)\in(-\infty,0)$,
and Theorem \ref{Th2}, it
follows that $\HerzS$ is a BQBF space.
Applying Lemma \ref{vmbhg} with
$u:=\frac{1}{r}$, we conclude that,
for any $\{f_j\}_{j\in\mathbb{N}}\subset
L^1_{\mathrm{loc}}(\rrn)$,
$$
\left\|\left\{\sum_{j\in\mathbb{N}}
\left[\mc\left(f_j\right)\right]^{1/r}
\right\}^r\right\|_{[\HerzS]^{1/r}}
\lesssim\left\|\left(\sum_{j\in\mathbb{N}}
\left|f_j\right|^{1/r}\right)^{r}\right\|_{
[\HerzS]^{1/r}},
$$
which implies that Proposition \ref{watomgd}(i)
holds true.
Now, we choose a
$$
\theta\in\left(\max\{1,p,q\},
\frac{1}{p_-}\min
\left\{p,q,\frac{n}{\Mw+n/p}\right\}\right).
$$
Then, by the assumptions
$p/\theta,\ q/\theta\in(0,1)$,
the reverse Minkowski inequality,
and Lemma \ref{convexl},
we find that, for any $\{f_{j}\}_{j\in
\mathbb{N}}\subset\Msc(\rrn)$,
\begin{equation*}
\sum_{j\in\mathbb{N}}\|f_{j}\|_{[\HerzS]^{1/\theta}}\leq
\left\|\sum_{j\in\mathbb{N}}|f_{j}|
\right\|_{[\HerzS]^{1/\theta}}.
\end{equation*}
This implies that $\HerzS$
is $\theta$-concave.
Moreover, from Lemma \ref{mbhg} with
$r:=\theta p_-$, we deduce that,
for any $f\in L_{{\rm loc}}^{1}(\rrn)$,
\begin{equation*}
\left\|\mc(f)\right\|_{[\HerzS]^{\frac{1}{\theta p_-}}}
\lesssim\left\|f\right\|_{[\HerzS]^{\frac{1}{\theta p_-}}},
\end{equation*}
which, together with the fact
that $\HerzS$ is $\theta$-concave,
further implies that Proposition
\ref{watomgd}(ii) holds true.
In addition, let
$$
r_0\in\left(0,\min\left\{p,\frac{n}
{\Mw+n/p}\right\}\right).
$$
Then, using Lemma \ref{wmbhg} with
$r:=r_0$, we conclude that,
for any $f\in L^1_{\mathrm{loc}}(\rrn)$,
$$
\left\|\mc(f)\right\|_{[W\HerzS]^{1/r_0}}
\lesssim\left\|f\right\|_{[W\HerzS]^{1/r_0}},
$$
which implies that Proposition \ref{watomgd}(iii)
holds true with the above $r_0$.
Finally, let
$$
s\in\left(0,\min\left\{1,p,q,
\frac{n}{\Mw+n/p}\right\}\right)
$$
and $\theta_0\in(1,\infty)$ satisfy
$$
\theta_0<\min\left\{\frac{n}{n(1-s/p)-s\mw},
\left(\frac{p}{s}\right)'\right\}.
$$
Then, repeating an argument similar to that
used in the proof of Theorem \ref{bddog}
with $\eta$ therein replaced by $\theta_0$,
we find that
\begin{enumerate}
  \item[{\rm(i)}] for any $f\in\Msc(\rrn)$,
  $$
  \|f\|_{[\HerzS]^{1/s}}
  \sim\sup\left\{\|fg\|_{L^1(\rrn)}:\
  \|g\|_{\HerzScsd}=1\right\}
  $$
  with positive equivalence
  constants independent of $f$;
  \item[{\rm(ii)}] $\|\cdot\|_{\HerzScsd}$ satisfies
  Definition \ref{Df1}(ii);
  \item[{\rm(iii)}] $\1bf_{B(\0bf,1)}\in\HerzScsd$;
  \item[{\rm(iv)}] $\mc^{(\theta_0)}$ is bounded
  on $\HerzScsd$.
\end{enumerate}
These imply that the assumptions
(iv) through (vii) of Proposition \ref{watomgd}
hold true with the above $s$ and $\theta_0$.
Thus, all the assumptions of
Proposition \ref{watomgd} hold
true for the Herz space $\HerzS$ under consideration.

Let $f\in W\HaSaH$. Then,
from Proposition \ref{watomgd}
with $X:=\HerzS$, we infer that there exists
$\{a_{i,j}\}_{i\in\mathbb{Z},\,j\in\mathbb{N}}$
of $(\HerzS,\,\infty,\,d)$-atoms supported,
respectively, in the balls $\{B_{i,j}\}_{
i\in\mathbb{Z},\,j\in\mathbb{N}}\subset\mathbb{B}$
and three positive constants $c\in(0,1]$,
$A$, and $\widetilde{A}$, independent of $f$,
such that, for any
$i\in\mathbb{Z}$, $$\sum_{j\in\mathbb{N}}
\1bf_{cB_{i,j}}\leq A,$$
\begin{equation}\label{watomge1}
f=\sum_{i\in\mathbb{Z}}\sum_{j\in\mathbb{N}}
\widetilde{A}2^i\left\|\1bf_{B_{i,j}}
\right\|_{\HerzS}a_{i,j}
\end{equation}
in $\mathcal{S}'(\rrn)$, and
\begin{equation}\label{watomge2}
\sup_{i\in\mathbb{Z}}\left\{
2^i\left\|\sum_{j\in\mathbb{N}}
\1bf_{B_{i,j}}\right\|_{\HerzS}\right\}\lesssim
\left\|f\right\|_{W\HaSaH}.
\end{equation}
In addition, using Lemma \ref{Atoll2} with
$X:=\HerzS$ and $t:=\infty$, we find that,
for any $i\in\mathbb{Z}$
and $j\in\mathbb{N}$, $a_{i,j}$ is
a $(\HerzS,\,r,\,d)$-atom
supported in the ball $B_{i,j}$.
Combining this, \eqref{watomge1},
\eqref{watomge2}, and Definition \ref{dfwatomg},
we conclude that $f\in WH\Kmp_{\omega}^{p,q,r,d}(\rrn)$
and
\begin{equation}\label{watomge3}
\left\|f\right\|_{WH\Kmp_{\omega}^{p,q,r,d}(\rrn)}
\lesssim\left\|f\right\|_{W\HaSaH}.
\end{equation}
This then finishes the proof that $W\HaSaH\subset
WH\Kmp_{\omega}^{p,q,r,d}(\rrn)$.

Conversely, we next show that
$WH\Kmp_{\omega}^{p,q,r,d}(\rrn)\subset W\HaSaH$.
For this purpose, we first prove that
the Herz space $\HerzS$ under consideration
satisfies the assumptions (i), (ii),
and (iii) of Proposition \ref{watomgr} with
$p_-$ as in the present theorem.
Indeed, for any given $\theta\in(0,p_-)$
and $u\in(1,\infty)$, from Lemma \ref{vmbhg} with
$r:=\theta$, it follows that,
for any $\{f_j\}_{j\in\mathbb{N}}\subset
\Msc(\rrn)$,
$$
\left\|\left\{
\sum_{j\in\mathbb{N}}\left[\mc\left(f_j\right)\right]^u
\right\}^{\frac{1}{u}}\right\|_{[\HerzS]^{1/\theta}}
\lesssim\left\|\left(\sum_{j\in\mathbb{N}}\left|f_j\right|^u
\right)^{\frac{1}{u}}\right\|_{[\HerzS]^{1/\theta}},
$$
which implies that Proposition \ref{watomgr}(i)
holds true. Next, we show that,
for any $s\in(0,p_-)$, $[\HerzS]^{1/s}$
is a BBF. Indeed, fix an $s\in(0,p_-)$. Then,
from the assumptions $\m0(\omega)
\in(-\frac{n}{p},\infty)$
and $\MI(\omega)\in(-\infty,0)$, and
Lemma \ref{rela}, it follows that
\begin{equation*}
\m0\left(\omega^s\right)
=s\m0(\omega)>-\frac{n}{p/s}
\end{equation*}
and
\begin{equation*}
\MI\left(\omega^s\right)=s
\MI(\omega)<0.
\end{equation*}
Applying these,
the assumptions $p/s,\ q/s\in(1,\infty)$,
and Theorem \ref{ball}
with $p$, $q$, and $\omega$
replaced, respectively, by $p/s$,
$q/s$, and $\omega^s$, we find that
$\Kmp^{p/s,q/s}_{\omega^s}(\rrn)$
is a BBF space. From this and Lemma \ref{convexl},
we further infer that $[\HerzS]^{1/s}$
is a BBF. This then implies that
Proposition \ref{watomgr}(ii)
holds true.

Finally, let $s_0\in(0,p_-)$ and
$$
r_0\in\left(\max\left\{1,p,
\frac{n}{\Mw+n/p}\right\},\infty\right).
$$
Then, repeating an argument similar to that
used in the proof of Theorem \ref{Atog} with
$s$ and $r$ therein replaced, respectively,
by $s_0$ and $r_0$, we conclude that,
for any $f\in\Msc(\rrn)$,
\begin{equation*}
\|f\|_{[\HerzS]^{1/s_0}}\sim
\sup\left\{\|fg\|_{L^1(\rrn)}:\
\|g\|_{\dot{\mathcal{B}}_{1/\omega^{s_0}}
^{(p/s_0)',(q/s_0)'}(\rrn)}=1\right\}
\end{equation*}
and, for any
$f\in L_{{\rm loc}}^{1}(\rrn)$,
\begin{equation*}
\left\|\mc^{((r_0/s_0)')}(f)
\right\|_{\dot{\mathcal{B}}
_{1/\omega^{s_0}}^{(p/s_0)',(q/s_0)'}(\rrn)}\lesssim
\|f\|_{\dot{\mathcal{B}}_
{1/\omega^{s_0}}^{(p/s_0)',(q/s_0)'}(\rrn)}.
\end{equation*}
These imply that the global generalized Herz
space $\HerzS$ under consideration
satisfies Proposition \ref{watomgr}(iii)
with the above $s_0$, $r_0$, and
$
Y_0:=\dot{\mathcal{B}}_
{1/\omega^{s_0}}^{(p/s_0)',(q/s_0)'}(\rrn).
$
Therefore, the assumptions (i) through
(iii) of Proposition \ref{watomgr} hold
true for $\HerzS$.

In addition, by the
proof that $W\HaSaH\subset
WH\Kmp_{\omega}^{p,q,r,d}(\rrn)$,
we find that $\HerzS$ is a BQBF space.
Combining this, the fact that
$\HerzS$ satisfies (i), (ii), and (iii)
of Proposition \ref{watomgr},
Proposition \ref{watomgr} with $X:=\HerzS$,
and Definition \ref{dfwatomg},
we conclude that
$WH\Kmp_{\omega}^{p,q,r,d}(\rrn)\subset W\HaSaH$
and, for any $f\in
WH\Kmp_{\omega}^{p,q,r,d}(\rrn)$,
\begin{equation*}
\left\|f\right\|_{W\HaSaH}\lesssim
\left\|f\right\|_{WH\Kmp_{\omega}^{p,q,r,d}(\rrn)}.
\end{equation*}
This, combined with \eqref{watomge3},
further implies that
$$
W\HaSaH=WH\Kmp_{\omega}^{p,q,r,d}(\rrn)
$$
with equivalent quasi-norms,
which then completes the proof of Theorem \ref{watomg}.
\end{proof}

Finally, combining Theorem
\ref{watomg} and Remark \ref{remark6.16},
we have the following atomic characterization
of $WH\MorrS$; we omit the
details.\index{atomic characterization}

\begin{corollary}
Let $p$, $q$, $\omega$,
$r$, and $d$ be as in Corollary \ref{watomm}.
Then the \emph{weak generalized atomic Morrey--Hardy space}
\index{weak generalized atomic Morrey--Hardy space}$WH
\textsl{\textbf{M}}_{\omega}^{p,q,r,d}(\rrn)$\index{$WH
\textsl{\textbf{M}}_{\omega}^{p,q,r,d}(\rrn)$}, associated with
the weak local generalized Morrey space $W\MorrS$, is defined to
be the set of all the $f\in\mathcal{S}'(\rrn)$ such that
there exist a sequence $\{a_{i,j}\}_{i\in
\mathbb{Z},\,j\in\mathbb{N}}$ of $(\MorrS,\,r,\,d)$-atoms supported,
respectively, in the balls $\{B_{i,j}\}_{i\in\mathbb{Z},\,j\in\mathbb{N}}
\subset\mathbb{B}$ and three positive constants $c\in(0,1]$,
$A$, and $\widetilde{A}$, independent of $f$, satisfying that,
for any $i\in\mathbb{Z}$,
$$
\sum_{j\in\mathbb{N}}\1bf_{cB_{i,j}}\leq A,
$$
$$
f=\sum_{i\in\mathbb{Z}}\sum_{j\in\mathbb{N}}
\widetilde{A}2^i\left\|\1bf_{B_{i,j}}\right\|_{
\MorrS}a_{i,j}
$$
in $\mathcal{S}'(\rrn)$, and
$$
\sup_{i\in\mathbb{Z}}
\left\{2^i\left\|\sum_{j\in\mathbb{N}}
\1bf_{B_{i,j}}\right\|_{\MorrS}\right\}<\infty,
$$
Moreover, for any
$f\in WH\textsl{\textbf{M}}_{\omega}^{p,q,r,d}(\rrn)$,
$$
\|f\|_{WH\textsl{\textbf{M}}_{\omega}^{p,q,r,d}(\rrn)}:=
\inf\sup_{i\in\mathbb{Z}}\left\{
2^i\left\|\sum_{j\in\mathbb{N}}
\1bf_{B_{i,j}}\right\|_{\MorrS}\right\},
$$
where the infimum is taken over all the decompositions of
$f$ as above.
Then $$
WH\MorrS=WH\textsl{\textbf{M}}_{\omega}^{p,q,r,d}(\rrn)
$$
with equivalent quasi-norms.
\end{corollary}

\section{Molecular Characterizations}

In this section, we establish
the molecular characterization
of weak generalized Herz--Hardy spaces.
Indeed, applying the
molecular characterization of
weak Hardy spaces $WH_X(\rrn)$ associated with
ball quasi-Banach function spaces $X$
obtained in \cite{ZWYY},
we immediately obtain the
molecular characterization of
the weak generalized Herz--Hardy
space $W\HaSaHo$. However, due to the
deficiency of the associate space
of the global generalized Herz
space $\HerzS$, we can not
prove the molecular characterization
of the weak generalized Herz--Hardy
space $W\HaSaH$ via
the known molecular characterization
of $WH_X(\rrn)$ directly.
To overcome this difficulty,
we first prove an improved
molecular reconstruction theorem of
weak Hardy spaces associated with
ball quasi-Banach function spaces
(see Proposition \ref{wmolegr} below)
without recourse to associate spaces.
Then, combining this improved conclusion
and the atomic
characterization of $W\HaSaH$ established
in the last section, we obtain the desired
molecular characterization of $W\HaSaH$.

Recall that the definition of
$(\HerzSo,\,r,\,d,\,\tau)$-molecules
are given in Definition \ref{mole}. We
now characterize
the weak generalized Herz--Hardy space $WH\HerzSo$
via these molecules
as follows.\index{molecular \\characterization}

\begin{theorem}\label{wmole}
Let $p,\ q\in(0,\infty)$, $\omega\in M(\rp)$ satisfy
$\m0(\omega)\in(-\frac{n}{p},\infty)$ and $\mi(\omega)\in(
-\frac{n}{p},\infty)$,
$$
p_-\in\left(0,\frac{\min\{p,q,\frac{n}{\Mw+n/p}\}}
{\max\{1,p,q\}}\right),
$$
$d\geq\lfloor n(1/p_--1)\rfloor$ be a fixed
integer,
$$
r\in\left(\max\left\{1,p,\frac{n}{\mw+n/p}
\right\},\infty\right],
$$
and $\tau\in(n(\frac{1}{p_-}-
\frac{1}{r}),\infty)$.
Then $f$ belongs to the weak generalized Herz--Hardy space
$WH\HerzSo$ if and only if $f\in\mathcal{S}'(\rrn)$
and there exist a sequence $\{m_{i,j}\}_{i\in
\mathbb{Z},\,j\in\mathbb{N}}$ of
$(\HerzSo,\,r,\,d,\,\tau)$-molecules centered,
respectively, at the balls $\{B_{i,j}
\}_{i\in\mathbb{Z},\,j\in\mathbb{N}}\subset\mathbb{B}$
and three positive constants $c\in(0,1]$,
$A$, and $\widetilde{A}$, independent of $f$, such that,
for any $i\in\mathbb{Z}$,
$$
\sum_{j\in\mathbb{N}}\1bf_{cB_{i,j}}\leq A,
$$
$$
f=\sum_{i\in\mathbb{Z}}\sum_{j\in\mathbb{N}}
\widetilde{A}2^i\left\|\1bf_{
B_{i,j}}\right\|_{\HerzSo}m_{i,j}
$$
in $\mathcal{S}'(\rrn)$, and
$$
\sup_{i\in\mathbb{Z}}\left\{2^i
\left\|\sum_{j\in\mathbb{N}}
\1bf_{B_{i,j}}\right\|_{\HerzSo}\right\}<\infty.
$$
Moreover, there exists a constant
$C\in[1,\infty)$ such that,
for any $f\in\WHaSaHo$,
\begin{align*}
C^{-1}\|f\|_{\WHaSaHo}
&\leq\inf\left\{\sup_{i\in\mathbb{Z}}
2^i\left\|\sum_{j\in\mathbb{N}}
\1bf_{B_{i,j}}\right\|_{\HerzSo}\right\}\\
&\leq C\|f\|_{\WHaSaHo},
\end{align*}
where the infimum is taken over
all the decompositions of
$f$ as above.
\end{theorem}

To show this molecular characterization,
we first recall the following
molecular characterization of weak
Hardy spaces $WH_X(\rrn)$ associated
with ball quasi-Banach
function spaces $X$, which was obtained by
Zhang et al.\ in
\cite[Theorem 5.3]{ZWYY}.

\begin{lemma}\label{wmoler}
Let $X$ be a ball quasi-Banach function space
and let $p_-\in(0,1)$ and $p_+\in[p_-,\infty)$
be such that the following
three statements hold true:
\begin{enumerate}
  \item[{\rm(i)}] for any given
  $\theta\in(0,p_-)$ and $u\in(1,\infty)$,
  there exists a positive constant $C$ such that,
for any $\{f
_{j}\}_{j\in\mathbb{N}}\subset L_{{\rm loc}}^{1}(\rrn)$,
\begin{equation*}
\left\|\left\{\sum_{j\in\mathbb{N}}\left[\mc
(f_{j})\right]^{u}\right\}^
{\frac{1}{u}}\right\|_{X^{1/\theta}}
\leq C\left\|\left\{
\sum_{j\in\mathbb{N}}|f_{j}|^{u}\right\}
^{\frac{1}{u}}\right\|_{X^{1/\theta}};
\end{equation*}
  \item[{\rm(ii)}] for any $s\in(0,p_-)$,
  $X^{1/s}$ is a ball Banach function space;
  \item[{\rm(iii)}] for any given
  $s\in(0,p_-)$ and $r\in(p_+,\infty)$,
  there exists a positive constant
  $C$ such that,
  for any $f\in L^1_{\mathrm{loc}}(\rrn)$,
  $$
  \left\|\mc^{((r/s)')}(f)\right\|_{(X^{1/s})'}
  \leq C\left\|f\right\|_{(X^{1/s})'}.
  $$
\end{enumerate}
Let $d\geq\lfloor n(1/p_--1)\rfloor$
be a fixed integer, $r\in(\max\{1,p_+\},\infty]$,
$\tau\in(n(\frac{1}{p_-}-\frac{n}{r}),\infty)$,
$c\in(0,1]$, and $A,\ \widetilde{A}\in(0,\infty)$.
Assume that
$\{a_{i,j}\}_{i\in\mathbb{Z},\,j\in\mathbb{N}}$
is a sequence of $(X,\,r,\,d,\,\tau)$-molecules
centered, respectively, at the balls $\{B_{i,j}\}_{
i\in\mathbb{Z},\,j\in\mathbb{N}}\subset\mathbb{B}$
such that, for any $i\in\mathbb{Z}$, $$\sum_{j\in\mathbb{N}}
\1bf_{cB_{i,j}}\leq A,$$
$$f:=\sum_{i\in\mathbb{Z}}\sum_{j\in\mathbb{N}}
\widetilde{A}2^i\left\|\1bf_{B_{i,j}}
\right\|_Xa_{i,j}$$
converges in $\mathcal{S}'(\rrn)$, and
$$
\sup_{i\in\mathbb{Z}}\left\{
2^i\left\|\sum_{j\in\mathbb{N}}
\1bf_{B_{i,j}}\right\|_X\right\}<\infty.
$$
Then $f\in WH_X(\rrn)$ and
$$
\left\|f\right\|_{WH_X(\rrn)}\lesssim
\sup_{i\in\mathbb{Z}}\left\{2^i\left\|\sum_{j\in\mathbb{N}}
\1bf_{B_{i,j}}\right\|_X\right\},
$$
where the implicit positive constant is
independent of $f$.
\end{lemma}

With the help of the above lemma,
we next show Theorem \ref{wmole}.

\begin{proof}[Proof of Theorem \ref{wmole}]
Let all the symbols be as in the
present theorem. We first prove the necessity.
To achieve this, let $f\in W\HaSaHo$.
Then, applying Theorem \ref{watom},
we find that $f\in
WH\Kmp_{\omega,\0bf}^{p,q,r,d}(\rrn)$.
From this and Definition \ref{dfwatom},
it follows that there exist a sequence
$\{a_{i,j}\}_{i\in\mathbb{Z},\,j\in\mathbb{N}}$
of $(\HerzSo,\,r,\,d)$-atoms supported,
respectively, in the balls $\{B_{i,j}\}_{
i\in\mathbb{Z},\,j\in\mathbb{N}}\subset\mathbb{B}$
and three positive constants $c\in(0,1]$,
$A$, and $\widetilde{A}$, independent of $f$,
such that, for any $i\in\mathbb{Z}$,
$$\sum_{j\in\mathbb{N}}\1bf_{cB_{i,j}}\leq A,$$
\begin{equation}\label{wmolee1}
f=\sum_{i\in\mathbb{Z}}\sum_{j\in\mathbb{N}}
\widetilde{A}2^i\left\|
\1bf_{B_{i,j}}\right\|_{\HerzSo}a_{i,j}
\end{equation}
in $\mathcal{S}'(\rrn)$, and
\begin{equation}\label{wmolee2}
\sup_{i\in\mathbb{Z}}
\left\{2^i\left\|\sum_{j\in\mathbb{N}}
\1bf_{B_{i,j}}\right\|_{\HerzSo}\right\}<\infty.
\end{equation}
On the other hand, by Remark \ref{molexr}
with $X:=\HerzSo$, we conclude that,
for any $i\in\mathbb{Z}$ and $j\in\mathbb{N}$,
$a_{i,j}$ is a $(\HerzSo,\,r,\,d,\,\tau)$-molecule
centered at $B_{i,j}$.
Combining this, \eqref{wmolee1}, and
\eqref{wmolee2}, we then complete the
proof the necessity. Moreover,
using both Definition \ref{dfwatom} and Theorem
\ref{watom} again, we find that
\begin{align}\label{wmolee3}
\inf\left\{\sup_{i\in\mathbb{Z}}
2^i\left\|\sum_{j\in\mathbb{N}}
\1bf_{B_{i,j}}\right\|_{\HerzSo}\right\}
&\lesssim\left\|f\right\|_{WH\Kmp_{\omega,\0bf}^{
p,q,r,d,}(\rrn)}\notag\\&\sim\left\|f\right\|_{W\HaSaHo},
\end{align}
where the infimum is taken over all the
decompositions of $f$ as in the present theorem.

Conversely, we next prove the sufficiency.
For this purpose, assume that $f\in\mathcal{S}'(\rrn)$,
$\{m_{i,j}\}_{i\in
\mathbb{Z},\,j\in\mathbb{N}}$ is a sequence of
$(\HerzSo,\,r,\,d,\,\tau)$-molecules centered,
respectively, at the balls $\{B_{i,j}
\}_{i\in\mathbb{Z},\,j\in\mathbb{N}}\subset\mathbb{B}$,
$c\in(0,1]$, and $A,\ \widetilde{A}\in(0,\infty)$
satisfying that, for any $i\in\mathbb{Z}$,
$$
\sum_{j\in\mathbb{N}}\1bf_{cB_{i,j}}\leq A,
$$
$$
f=\sum_{i\in\mathbb{Z}}\sum_{j\in\mathbb{N}}
\widetilde{A}2^i\left\|\1bf_{
B_{i,j}}\right\|_{\HerzSo}m_{i,j}
$$
in $\mathcal{S}'(\rrn)$, and
$$
\sup_{i\in\mathbb{Z}}\left\{2^i
\left\|\sum_{j\in\mathbb{N}}
\1bf_{B_{i,j}}\right\|_{\HerzSo}\right\}<\infty.
$$
Let
$$
p_+:=\max\left\{1,p,\frac{n}{\mw+n/p}\right\}.
$$
Then, repeating an argument similar
to that used in the proof of Theorem \ref{watom}
with $s_0$ and $r_0$ therein replaced,
respectively, by $s$ and $r$, we conclude
that the local generalized
Herz space $\HerzSo$ under consideration
satisfies the following four statements:
\begin{enumerate}
  \item[{\rm(i)}] $\HerzSo$ is a BQBF space;
  \item[{\rm(i)}] for any given
  $\theta\in(0,p_-)$ and $u\in(1,\infty)$, and
for any $\{f
_{j}\}_{j\in\mathbb{N}}\subset L_{{\rm loc}}^{1}(\rrn)$,
\begin{equation*}
\left\|\left\{\sum_{j\in\mathbb{N}}\left[\mc
(f_{j})\right]^{u}\right\}^
{\frac{1}{u}}\right\|_{[\HerzSo]^{1/\theta}}
\lesssim\left\|\left\{
\sum_{j\in\mathbb{N}}|f_{j}|^{u}\right\}
^{\frac{1}{u}}\right\|_{[\HerzSo]^{1/\theta}};
\end{equation*}
  \item[{\rm(ii)}] for any $s\in(0,p_-)$,
  $[\HerzSo]^{1/s}$ is a BBF space;
  \item[{\rm(iii)}] for any given
  $s\in(0,p_-)$ and $r\in(p_+,\infty)$,
and for any $f\in L^1_{\mathrm{loc}}(\rrn)$,
  $$
  \left\|\mc^{((r/s)')}(f)\right\|_{([\HerzSo]^{1/s})'}
  \lesssim\left\|f\right\|_{([\HerzSo]^{1/s})'}.
  $$
\end{enumerate}
These, together with Lemma \ref{wmoler}
with $X:=\HerzSo$, further imply
that $f\in W\HaSaHo$. This then finishes
the proof of the sufficiency.
Moreover, applying Lemma \ref{wmoler} again
with $X:=\HerzSo$ and the choice of
$\{m_{i,j}\}_{i\in\mathbb{Z},\,
j\in\mathbb{N}}$,
we find that
$$
\left\|f\right\|_{W\HaSaHo}
\lesssim\inf\sup_{i\in\mathbb{Z}}
\left\{2^i\left\|\sum_{j\in\mathbb{N}}
\1bf_{B_{i,j}}\right\|_{\HerzSo}\right\},
$$
where the infimum is taken over
all the decompositions of $f$ as in
the present theorem. From this and \eqref{wmolee3},
we further infer that, for any $f\in W\HaSaHo$,
$$
\left\|f\right\|_{W\HaSaHo}
\sim\inf\sup_{i\in\mathbb{Z}}\left\{
2^i\left\|\sum_{j\in\mathbb{N}}
\1bf_{B_{i,j}}\right\|_{\HerzSo}\right\},
$$
where the infimum is taken over
all the decompositions of $f$ as in
the present theorem. This finishes
the proof of Theorem \ref{wmole}.
\end{proof}

From
Theorem \ref{wmole} and Remark \ref{remark6.16} above, we deduce
the following molecular characterization of
the weak generalized Morrey--Hardy space $WH\MorrSo$; we
omit the details.

\begin{corollary}\label{wmolem}
Let $p,\ q\in[1,\infty)$, $p_-\in
(0,\min\{p,q\}/\max\{p,q\})$,
$\omega\in M(\rp)$ with
$$
-\frac{n}{p}<\m0(\omega)\leq\M0(\omega)<0
$$
and
$$
-\frac{n}{p}<\mi(\omega)\leq\MI(\omega)<0,
$$
$d\geq\lfloor n(1/p_--1)\rfloor$ be a fixed
integer,
$$
r\in\left(\frac{n}{\mw+n/p},\infty\right],
$$
and $\tau\in(n(\frac{1}{p_-}-
\frac{1}{r}),\infty)$.
Then $f$ belongs to the weak generalized Morrey--Hardy space
$WH\MorrSo$ if and only if $f\in\mathcal{S}'(\rrn)$
and there exist a sequence $\{m_{i,j}\}_{i\in
\mathbb{Z},\,j\in\mathbb{N}}$ of
$(\MorrSo,\,r,\,d,\,\tau)$-molecules centered,
respectively, at the balls $\{B_{i,j}
\}_{i\in\mathbb{Z},\,j\in\mathbb{N}}\subset\mathbb{B}$
and three positive constants $c\in(0,1]$,
$A$, and $\widetilde{A}$, independent of $f$, such that,
for any $i\in\mathbb{Z}$,
$$
\sum_{j\in\mathbb{N}}\1bf_{cB_{i,j}}\leq A,
$$
$$
f=\sum_{i\in\mathbb{Z}}\sum_{j\in\mathbb{N}}
\widetilde{A}2^i\left\|\1bf_{
B_{i,j}}\right\|_{\MorrSo}m_{i,j}
$$
in $\mathcal{S}'(\rrn)$, and
$$
\sup_{i\in\mathbb{Z}}\left\{2^i
\left\|\sum_{j\in\mathbb{N}}
\1bf_{B_{i,j}}\right\|_{\MorrSo}\right\}<\infty.
$$
Moreover, there exists a constant
$C\in[1,\infty)$ such that,
for any $f\in WH\MorrSo$,
\begin{align*}
C^{-1}\|f\|_{WH\MorrSo}
&\leq\inf\sup_{i\in\mathbb{Z}}\left\{
2^i\left\|\sum_{j\in\mathbb{N}}
\1bf_{B_{i,j}}\right\|_{\MorrSo}\right\}\\
&\leq C\|f\|_{WH\MorrSo},
\end{align*}
where the infimum is taken over all
the decompositions of
$f$ as above.
\end{corollary}

Now, we turn to establish the molecular
characterization of the weak generalized
Herz--Hardy space $W\HaSaH$. Recall that
the concept of $(\HerzS,\,r,
\,d,\,\tau)$-molecules is
introduced in Definition \ref{moleg} above.
Then we show the following conclusion.

\begin{theorem}\label{wmoleg}
Let $p,\ q\in(0,\infty)$, $\omega\in M(\rp)$ satisfy
$\m0(\omega)\in(-\frac{n}{p},\infty)$ and
$$-\frac{n}{p}<
\mi(\omega)\leq\MI(\omega)<0,$$
$
r\in(\max\{1,p,\frac{n}{\mw+n/p}\},\infty],
$
$$
p_-\in\left(0,\frac{\min\{p,q,\frac{n}{\Mw+n/p}\}}
{\max\{1,p,q\}}\right),
$$
$d\geq\lfloor n(1/p_--1)\rfloor$ be a fixed
integer,
and $\tau\in(n(\frac{1}{p_-}-
\frac{1}{r}),\infty)$.
Then $f$ belongs to the weak generalized Herz--Hardy space
$WH\HerzS$ if and only if $f\in\mathcal{S}'(\rrn)$
and there exist a sequence $\{m_{i,j}\}_{i\in
\mathbb{Z},\,j\in\mathbb{N}}$ of
$(\HerzS,\,r,\,d,\,\tau)$-molecules centered,
respectively, at the balls $\{B_{i,j}
\}_{i\in\mathbb{Z},\,j\in\mathbb{N}}\subset\mathbb{B}$
and three positive constants $c\in(0,1]$,
$A$, and $\widetilde{A}$, independent of $f$, such that,
for any $i\in\mathbb{Z}$,
$$
\sum_{j\in\mathbb{N}}\1bf_{cB_{i,j}}\leq A,
$$
$$
f=\sum_{i\in\mathbb{Z}}\sum_{j\in\mathbb{N}}
\widetilde{A}2^i\left\|\1bf_{
B_{i,j}}\right\|_{\HerzS}m_{i,j}
$$
in $\mathcal{S}'(\rrn)$, and
$$
\sup_{i\in\mathbb{Z}}\left\{2^i
\left\|\sum_{j\in\mathbb{N}}
\1bf_{B_{i,j}}\right\|_{\HerzS}\right\}<\infty.
$$
Moreover, there exists a constant
$C\in[1,\infty)$ such that,
for any $f\in\WHaSaH$,
\begin{align*}
C^{-1}\|f\|_{\WHaSaH}
&\leq\inf\sup_{i\in\mathbb{Z}}
\left\{2^i\left\|\sum_{j\in\mathbb{N}}
\1bf_{B_{i,j}}\right\|_{\HerzS}\right\}\\
&\leq C\|f\|_{\WHaSaH},
\end{align*}
where the infimum is taken over all
the decompositions of
$f$ as above.
\end{theorem}

In order to show this theorem, we first
establish the following molecular
reconstruction theorem of the weak Hardy space
$WH_X(\rrn)$, which improves
\cite[Theorem 5.3]{ZWYY} via removing the assumption
about associate spaces.

\begin{proposition}\label{wmolegr}
Let $X$ be a ball quasi-Banach function space and
let $p_-\in(0,1)$ and $p_+\in[p_-,\infty)$
be such that the following
three statements hold true:
\begin{enumerate}
  \item[{\rm(i)}] for any given
  $\theta\in(0,p_-)$ and $u\in(1,\infty)$,
  there exists a positive constant $C$ such that,
for any $\{f
_{j}\}_{j\in\mathbb{N}}\subset L_{{\rm loc}}^{1}(\rrn)$,
\begin{equation*}
\left\|\left\{\sum_{j\in\mathbb{N}}\left[\mc
(f_{j})\right]^{u}\right\}^
{\frac{1}{u}}\right\|_{X^{1/\theta}}
\leq C\left\|\left\{
\sum_{j\in\mathbb{N}}|f_{j}|^{u}\right\}
^{\frac{1}{u}}\right\|_{X^{1/\theta}};
\end{equation*}
  \item[{\rm(ii)}] for any $s\in(0,p_-)$,
  $X^{1/s}$ is a ball Banach function space
  and there exists a linear space $Y_s\subset
  \Msc(\rrn)$
  equipped with a quasi-seminorm $\|\cdot\|_{Y_s}$
  such that, for any $f\in\Msc(\rrn)$,
  $$
  \left\|f\right\|_{X^{1/s}}\sim
  \sup\left\{\left\|fg\right\|_{L^1(\rrn)}:\
  \left\|g\right\|_{Y_s}=1\right\}
  $$
  with the positive equivalence constants
  independent of $f$;
  \item[{\rm(iii)}] for any given
  $s\in(0,p_-)$ and $r\in(p_+,\infty)$,
  there exists a positive constant
  $C$ such that,
  for any $f\in L^1_{\mathrm{loc}}(\rrn)$,
  $$
  \left\|\mc^{((r/s)')}(f)\right\|_{Y_s}
  \leq C\left\|f\right\|_{Y_s}.
  $$
\end{enumerate}
Let $d\geq\lfloor n(1/p_--1)\rfloor$
be a fixed integer, $r\in(\max\{1,p_+\},\infty]$,
$\tau\in(n(\frac{1}{p_-}-\frac{n}{r}),\infty)$,
$c\in(0,1]$, and $A,\ \widetilde{A}\in(0,\infty)$.
Assume that
$\{a_{i,j}\}_{i\in\mathbb{Z},\,j\in\mathbb{N}}$
is a sequence of $(X,\,r,\,d,\,\tau)$-molecules
centered, respectively, at the balls $\{B_{i,j}\}_{
i\in\mathbb{Z},\,j\in\mathbb{N}}\subset\mathbb{B}$
such that, for any $i\in\mathbb{Z}$, $$\sum_{j\in\mathbb{N}}
\1bf_{cB_{i,j}}\leq A,$$
$$f:=\sum_{i\in\mathbb{Z}}\sum_{j\in\mathbb{N}}
\widetilde{A}2^i\left\|\1bf_{B_{i,j}}
\right\|_Xa_{i,j}$$
converges in $\mathcal{S}'(\rrn)$, and
$$
\sup_{i\in\mathbb{Z}}\left\{
2^i\left\|\sum_{j\in\mathbb{N}}
\1bf_{B_{i,j}}\right\|_X\right\}<\infty.
$$
Then $f\in WH_X(\rrn)$ and
$$
\left\|f\right\|_{WH_X(\rrn)}\lesssim
\sup_{i\in\mathbb{Z}}2^i\left\|\sum_{j\in\mathbb{N}}
\1bf_{B_{i,j}}\right\|_X,
$$
where the implicit positive constant is
independent of $f$.
\end{proposition}

\begin{proof}
Let $X$, $r$, $d$, $\tau$,
$c$, $A$, and $\widetilde{A}$ be as in the present
proposition and $\{a_{i,j}\}_{i\in\mathbb{Z},\,j\in\mathbb{N}}$
a sequence of $(X,\,r,\,d,\,\tau)$-molecules
centered, respectively, at the balls $\{B_{i,j}\}_{
i\in\mathbb{Z},\,j\in\mathbb{N}}\subset\mathbb{B}$
such that, for any $i\in\mathbb{Z}$, $$\sum_{j\in\mathbb{N}}
\1bf_{cB_{i,j}}\leq A,$$
the summation
$$\sum_{i\in\mathbb{Z}}\sum_{j\in\mathbb{N}}
\widetilde{A}2^i\left\|\1bf_{B_{i,j}}\right\|_Xa_{i,j}$$
converges in $\mathcal{S}'(\rrn)$, and
$$
\sup_{i\in\mathbb{Z}}\left\{
2^i\left\|\sum_{j\in\mathbb{N}}
\1bf_{B_{i,j}}\right\|_X\right\}<\infty.
$$
Then, repeating the proof of
\cite[Theorem 5.3]{ZWYY} via replacing
\cite[Lemma 4.8]{ZWYY} therein by Lemma
\ref{at-es} here, we conclude that
$$f:=\sum_{i\in\mathbb{Z}}\sum_{j\in\mathbb{N}}
\widetilde{A}2^i
\left\|\1bf_{B_{i,j}}\right\|_{X}a_{i,j}\in
WH_X(\rrn)$$ and
$$
\left\|f\right\|_{WH_X(\rrn)}
\lesssim\sup_{i\in\mathbb{Z}}\left\{
2^i\left\|\sum_{j\in\mathbb{N}}\1bf_{B_{i,j}}
\right\|_X\right\},
$$
which then completes the proof of
Proposition \ref{wmolegr}.
\end{proof}

\begin{remark}
We should point out that
Proposition \ref{wmolegr} is an
improved version of the known molecular
reconstruction theorem obtained
by Zhang et al.\ in \cite[Theorem 5.3]{ZWYY}.
Indeed, if $Y_s\equiv(X^{1/s})'$ in
Proposition \ref{wmolegr}, then
this proposition is just \cite[Theoreem 5.3]{ZWYY}.
\end{remark}

Using this proposition and
the atomic characterization of
the weak generalized Herz--Hardy
space $W\HaSaH$ obtained in the last
section, we now prove the molecular
characterization of $W\HaSaH$.

\begin{proof}[Proof of Theorem \ref{wmoleg}]
Let all the symbols be as in the present theorem.
We first show the necessity. Indeed,
let $f\in W\HaSaH$.
Then, applying Theorem \ref{watomg},
we find that $f\in
WH\Kmp_{\omega}^{p,q,r,d}(\rrn)$.
This, together with Definition \ref{dfwatomg},
further implies that there exist a sequence
$\{a_{i,j}\}_{i\in\mathbb{Z},\,j\in\mathbb{N}}$
of $(\HerzS,\,r,\,d)$-atoms supported,
respectively, in the balls $\{B_{i,j}\}_{
i\in\mathbb{Z},\,j\in\mathbb{N}}\subset\mathbb{B}$
and three positive constants $c\in(0,1]$,
$A$, and $\widetilde{A}$, independent of $f$,
such that, for any $i\in\mathbb{Z}$,
$$\sum_{j\in\mathbb{N}}\1bf_{cB_{i,j}}\leq A,$$
\begin{equation}\label{wmolege1}
f=\sum_{i\in\mathbb{Z}}\sum_{j\in\mathbb{N}}
\widetilde{A}2^i\left\|
\1bf_{B_{i,j}}\right\|_{\HerzS}a_{i,j}
\end{equation}
in $\mathcal{S}'(\rrn)$, and
\begin{equation}\label{wmolege2}
\sup_{i\in\mathbb{Z}}
\left\{2^i\left\|\sum_{j\in\mathbb{N}}
\1bf_{B_{i,j}}\right\|_{\HerzS}\right\}<\infty.
\end{equation}
In addition, from Remark \ref{molexr}
with $X:=\HerzS$, it follows that,
for any $i\in\mathbb{Z}$ and $j\in\mathbb{N}$,
$a_{i,j}$ is a $(\HerzS,\,r,\,d,\,\tau)$-molecule
centered at $B_{i,j}$.
This, combined with \eqref{wmolege1} and
\eqref{wmolege2}, finishes the proof of the necessity.
Moreover, using both
Definition \ref{dfwatomg} and Theorem
\ref{watomg} again, we conclude that
\begin{align}\label{wmolege3}
\inf\left\{\sup_{i\in\mathbb{Z}}
2^i\left\|\sum_{j\in\mathbb{N}}
\1bf_{B_{i,j}}\right\|_{\HerzS}\right\}
&\lesssim\left\|f\right\|_{WH\Kmp_{\omega}^{
p,q,r,d,}(\rrn)}\notag\\&\sim\left\|f\right\|_{W\HaSaH},
\end{align}
where the infimum is taken over all the
decompositions of $f$ as in the present theorem.

Conversely, we prove the sufficiency.
To this end, let $f\in\mathcal{S}'(\rrn)$,
$\{m_{i,j}\}_{i\in\mathbb{Z},\,j\in\mathbb{N}}$
be a sequence of $(\HerzS,\,r,\,d,\,\tau)$-molecules
centered, respectively, at the balls $\{B_{i,j}
\}_{i\in\mathbb{Z},\,j\in\mathbb{N}}\subset\mathbb{B}$,
$c\in(0,1]$, and $A,\ \widetilde{A}\in(0,\infty)$
such that, for any $i\in\mathbb{Z}$,
$$
\sum_{j\in\mathbb{N}}\1bf_{cB_{i,j}}\leq A,
$$
$$
f=\sum_{i\in\mathbb{Z}}\sum_{j\in\mathbb{N}}
\widetilde{A}2^i\left\|\1bf_{
B_{i,j}}\right\|_{\HerzS}m_{i,j}
$$
in $\mathcal{S}'(\rrn)$, and
$$
\sup_{i\in\mathbb{Z}}\left\{2^i
\left\|\sum_{j\in\mathbb{N}}
\1bf_{B_{i,j}}\right\|_{\HerzS}\right\}<\infty.
$$
We next claim that the global generalized
Herz space $\HerzS$ under consideration
satisfies all the assumptions of Proposition
\ref{wmolegr}. Assume that this
claim holds true for the moment. Then,
applying Proposition \ref{wmolegr} with
$X:=\HerzS$, we conclude that $f\in W\HaSaH$,
which completes the proof of the sufficiency.

Thus, to finish the proof of the sufficiency,
we only need to show the above claim.
For this purpose, let
$$
p_+:=\max\left\{1,p,\frac{n}{\mw+n/p}\right\}.
$$
Then, repeating an argument similar
to that used in the proof of Theorem \ref{watomg}
with $s_0$ and $r_0$ therein replaced,
respectively, by $s$ and $r$, we find
that the global generalized
Herz space $\HerzS$ under consideration
satisfies the following four statements:
\begin{enumerate}
  \item[{\rm(i)}] $\HerzS$ is a BQBF space;
  \item[{\rm(i)}] for any given
  $\theta\in(0,p_-)$ and $u\in(1,\infty)$, and
for any $\{f
_{j}\}_{j\in\mathbb{N}}\subset L_{{\rm loc}}^{1}(\rrn)$,
\begin{equation*}
\left\|\left\{\sum_{j\in\mathbb{N}}\left[\mc
(f_{j})\right]^{u}\right\}^
{\frac{1}{u}}\right\|_{[\HerzS]^{1/\theta}}
\lesssim\left\|\left\{
\sum_{j\in\mathbb{N}}|f_{j}|^{u}\right\}
^{\frac{1}{u}}\right\|_{[\HerzS]^{1/\theta}};
\end{equation*}
  \item[{\rm(ii)}] for any $s\in(0,p_-)$,
  $[\HerzS]^{1/s}$ is a BBF space and, for any
  $f\in\Msc(\rrn)$,
  $$
  \left\|f\right\|_{[\HerzS]^{1/s}}
  \sim\sup\left\{\left\|fg\right\|_{L^1(\rrn)}:\
  \left\|g\right\|_{\dot{\mathcal{B}}_{1/\omega^s}
  ^{(p/s)',(q/s)'}(\rrn)}\right\}
  $$
  with the positive equivalence constants
  independent of $f$;
  \item[{\rm(iii)}] for any given
  $s\in(0,p_-)$ and $r\in(p_+,\infty)$,
and for any $f\in L^1_{\mathrm{loc}}(\rrn)$,
  $$
  \left\|\mc^{((r/s)')}(f)\right\|_{
  \dot{\mathcal{B}}_{1/\omega^{s}}
  ^{(p/s)',(q/s)'}(\rrn)}
  \lesssim\left\|f\right\|_{\dot{\mathcal{B}}_{1/\omega^{s}}
  ^{(p/s)',(q/s)'}(\rrn)}.
  $$
\end{enumerate}
These further imply that
the global generalized Herz space $\HerzS$
under consideration satisfies
all the assumptions of Proposition \ref{wmoleg}
with the above $p_-$, $p_+$,
and $
Y_s:=\dot{\mathcal{B}}_{1/\omega^{s}}
  ^{(p/s)',(q/s)'}(\rrn)
$
for any $s\in(0,p_-)$. This then
finishes the above claim and further implies
that the sufficiency holds true.

In addition, from Proposition \ref{wmolegr} again
with $X:=\HerzS$ and the choice of
$\{m_{i,j}\}_{i\in\mathbb{Z},\,
j\in\mathbb{N}}$,
we deduce that
$$
\left\|f\right\|_{W\HaSaH}
\lesssim\inf\sup_{i\in\mathbb{Z}}
\left\{2^i\left\|\sum_{j\in\mathbb{N}}
\1bf_{B_{i,j}}\right\|_{\HerzS}\right\},
$$
where the infimum is taken over
all the decompositions of $f$ as in
the present theorem. This, together
with \eqref{wmolege3},
further implies that, for any $f\in
W\HaSaH$,
$$
\left\|f\right\|_{W\HaSaH}
\sim\inf\sup_{i\in\mathbb{Z}}\left\{
2^i\left\|\sum_{j\in\mathbb{N}}
\1bf_{B_{i,j}}\right\|_{\HerzS}\right\},
$$
where the infimum is taken over
all the decompositions of $f$ as in
the present theorem. This then finishes
the proof of Theorem \ref{wmoleg}.
\end{proof}

As an application, we now establish
the following molecular
characterization of weak Hardy spaces
associated with global generalized Morrey spaces,
which is just a simple corollary of Theorem \ref{wmoleg}
and Remark \ref{remark6.16}; we omit the
details.

\begin{corollary}
Let $p$, $q$, $\omega$, $r$, $d$, and $\tau$ be
as in Corollary \ref{wmolem}.
Then $f$ belongs to the weak generalized Morrey--Hardy space
$WH\MorrS$ if and only if $f\in\mathcal{S}'(\rrn)$
and there exist a sequence $\{m_{i,j}\}_{i\in
\mathbb{Z},\,j\in\mathbb{N}}$ of
$(\MorrS,\,r,\,d,\,\tau)$-molecules centered,
respectively, at the balls $\{B_{i,j}
\}_{i\in\mathbb{Z},\,j\in\mathbb{N}}\subset\mathbb{B}$
and three positive constants $c\in(0,1]$,
$A$, and $\widetilde{A}$, independent of $f$, such that,
for any $i\in\mathbb{Z}$,
$$
\sum_{j\in\mathbb{N}}\1bf_{cB_{i,j}}\leq A,
$$
$$
f=\sum_{i\in\mathbb{Z}}\sum_{j\in\mathbb{N}}
\widetilde{A}2^i\left\|\1bf_{
B_{i,j}}\right\|_{\MorrS}m_{i,j}
$$
in $\mathcal{S}'(\rrn)$, and
$$
\sup_{i\in\mathbb{Z}}\left\{2^i
\left\|\sum_{j\in\mathbb{N}}
\1bf_{B_{i,j}}\right\|_{\MorrS}\right\}<\infty.
$$
Moreover, there exists a constant
$C\in[1,\infty)$ such that,
for any $f\in WH\MorrS$,
\begin{align*}
C^{-1}\|f\|_{WH\MorrS}
&\leq\inf\left\{\sup_{i\in\mathbb{Z}}
2^i\left\|\sum_{j\in\mathbb{N}}
\1bf_{B_{i,j}}\right\|_{\MorrS}\right\}\\
&\leq C\|f\|_{WH\MorrS},
\end{align*}
where the infimum is taken over all
the decompositions of
$f$ as above.
\end{corollary}

\section{Littlewood--Paley Function Characterizations}

The main target of this section is to
characterize weak generalized Herz--Hardy spaces
via various Littlewood--Paley functions.
Precisely, using the Lusin area function
$S$, the $g$-function $g$, and
the $g_{\lambda}^*$-function $g_{\lambda}^*$
presented in Definitions \ref{df451} and \ref{df452},
we establish several equivalent characterizations
of the weak generalized Herz--Hardy spaces
$\WHaSaHo$ and $\WHaSaH$.
To begin with, we state the following
Littlewood--Paley function characterizations
of the Hardy space $\WHaSaHo
$.\index{Littlewood--Paley function characterization}

\begin{theorem}\label{wlusin}
Let $p,\ q\in(0,\infty)$, $\omega\in M(\rp)$
with $\m0(\omega)\in(-\frac{n}{p},\infty)$
and $\mi(\omega)\in(-\frac{n}{p},\infty)$,
$$
p_-:=\min\left\{p,q,\frac{n}{\Mw+n/p}\right\},
$$
$$
p_+:=\max\left\{p,\frac{n}{\mw+n/p}\right\},
$$
and
$$
\lambda\in\left(\max\left\{
\frac{2}{\min\{1,p_-\}},\,1-\frac{2}{\max\{1,p_+\}}+
\frac{2}{\min\{1,p_-\}}\right\},\infty\right).
$$
Assume that, for any
$f\in\mathcal{S}'(\rrn)$, $S(f)$ and
$g_{\lambda}^*(f)$
are as in Definition \ref{df451}, and
$g(f)$ is as in Definition \ref{df452}.
Then the following four statements are mutually
equivalent:
\begin{enumerate}
  \item[{\rm(i)}] $f\in W\HaSaHo$;
  \item[{\rm(ii)}] $f\in\mathcal{S}'(\rrn)$,
  $f$ vanishes weakly at infinity,
  and $S(f)\in W\HerzSo$;
  \item[{\rm(iii)}] $f\in\mathcal{S}'(\rrn)$,
  $f$ vanishes weakly at infinity,
  and $g(f)\in W\HerzSo$;
  \item[{\rm(iv)}] $f\in\mathcal{S}'(\rrn)$,
  $f$ vanishes weakly at infinity,
  and $g_{\lambda}^*(f)\in W\HerzSo$.
\end{enumerate}
Moreover, for any $f\in W\HaSaHo$,
\begin{align*}
\left\|f\right\|_{W\HaSaHo}
&\sim\left\|S(f)\right\|_{W\HerzSo}
\sim\left\|g(f)\right\|_{W\HerzSo}\\
&\sim\left\|g_{\lambda}^*(f)\right\|_{W\HerzSo},
\end{align*}
where the positive equivalence
constants are independent of $f$.
\end{theorem}

In order to prove this theorem,
we need the following
Littlewood--Paley function characterizations
of weak Hardy spaces associated with ball
quasi-Banach function spaces, which are
just
\cite[Theorems 3.12, 3.16, and 3.21]{WYYZ}.

\begin{lemma}\label{wlpx}
Let $X$ be a ball quasi-Banach function space
and let $p_-\in(0,\infty)$, $p_+\in[p_-,\infty)$, and
$\theta_0\in(1,\infty)$ be
such that the following four statements hold true:
\begin{enumerate}
  \item[{\rm(i)}] for any given
  $\theta\in(0,p_-)$ and $u\in(1,\infty)$,
  there exists a positive constant $C$ such that,
for any $\{f
_{j}\}_{j\in\mathbb{N}}\subset L_{{\rm loc}}^{1}(\rrn)$,
\begin{equation*}
\left\|\left\{\sum_{j\in\mathbb{N}}\left[\mc
(f_{j})\right]^{u}\right\}^
{\frac{1}{u}}\right\|_{X^{1/\theta}}
\leq C\left\|\left\{
\sum_{j\in\mathbb{N}}|f_{j}|^{u}\right\}
^{\frac{1}{u}}\right\|_{X^{1/\theta}}
\end{equation*}
and
\begin{equation*}
\left\|\left\{\sum_{j\in\mathbb{N}}\left[\mc
(f_{j})\right]^{u}\right\}^
{\frac{1}{u}}\right\|_{(WX)^{1/\theta}}
\leq C\left\|\left\{
\sum_{j\in\mathbb{N}}|f_{j}|^{u}\right\}
^{\frac{1}{u}}\right\|_{(WX)^{1/\theta}};
\end{equation*}
  \item[{\rm(ii)}] $X$ is $\theta_0$-concave;
  \item[{\rm(iii)}] for any $s\in(0,\min\{1,p_-\})$,
  $X^{1/s}$ is a ball Banach function space;
  \item[{\rm(iv)}] for any given
  $s\in(0,\min\{1,p_-\})$ and $r\in(p_+,\infty)$,
  there exists a positive constant
  $C$ such that,
  for any $f\in L^1_{\mathrm{loc}}(\rrn)$,
  $$
  \left\|\mc^{((r/s)')}(f)\right\|_{(X^{1/s})'}
  \leq C\left\|f\right\|_{(X^{1/s})'}.
  $$
\end{enumerate}
Let $$
\lambda\in\left(\max\left\{
\frac{2}{\min\{1,p_-\}},\,1-\frac{2}{\max\{1,p_+\}}+
\frac{2}{\min\{1,p_-\}}\right\},\infty\right).
$$
Then the following four statements are mutually
equivalent:
\begin{enumerate}
  \item[{\rm(i)}] $f\in WH_X(\rrn)$;
  \item[{\rm(ii)}] $f\in\mathcal{S}'(\rrn)$,
  $f$ vanishes weakly at infinity,
  and $S(f)\in WX$;
  \item[{\rm(iii)}] $f\in\mathcal{S}'(\rrn)$,
  $f$ vanishes weakly at infinity,
  and $g(f)\in WX$;
  \item[{\rm(iv)}] $f\in\mathcal{S}'(\rrn)$,
  $f$ vanishes weakly at infinity,
  and $g_{\lambda}^*(f)\in WX$.
\end{enumerate}
Moreover, for any $f\in WH_X(\rrn)$,
$$
\left\|f\right\|_{WH_X(\rrn)}
\sim\left\|S(f)\right\|_{WX}
\sim\left\|g(f)\right\|_{WX}
\sim\left\|g_{\lambda}^*(f)\right\|_{WX}
$$
with the positive equivalence
constants independent of $f$.
\end{lemma}

Moreover, to prove Theorem \ref{wlusin},
we still require
the following auxiliary lemma about
the Fefferman--Stein vector-valued
inequality on weak local generalized
Herz spaces.

\begin{lemma}\label{vwmbhl}
Let $p,\ q\in(0,\infty)$ and $\omega\in M(\rp)$ satisfy
$\m0(\omega)\in(-\frac{n}{p},\infty)$
and $\mi(\omega)\in(-\frac{n}{p},
\infty)$.
Then, for any given
$$r\in\left(0,\min\left\{p,\frac{n}{\Mw+n/p}\right\}\right)$$
and $u\in(1,\infty)$,
there exists a positive constant $C$ such that,
for any $\{f_{j}\}_{
j\in\mathbb{N}}\subset L^1_{\mathrm{loc}}(\rrn)$,
$$
\left\|\left\{\sum_{j\in\mathbb{N}}\left[\mc
(f_{j})\right]^{u}\right\}^
{\frac{1}{u}}\right\|_{[W\HerzSo]^{1/r}}
\leq C\left\|\left\{
\sum_{j\in\mathbb{N}}|f_{j}|^{u}\right\}
^{\frac{1}{u}}\right\|_{[W\HerzSo]^{1/r}}.
$$
\end{lemma}

\begin{proof}
Let all the symbols be as in the present lemma.
Then, by the assumption
$\m0(\omega)\in(-\frac{n}{p},\infty)$
and Theorem \ref{Th3}, we
conclude that the local generalized
Herz space $\HerzSo$ under consideration
is a BQBF space.
Moreover, let $$
p_-:=\min\left\{p,\frac{n}{\Mw+n/p}\right\}.
$$
Then, for any given $r\in(0,p_-)$ and
$u\in(1,\infty)$,
applying Lemma \ref{vmbhl},
we find that, for any
$\{f_j\}_{j\in\mathbb{N}}\subset
L^1_{\mathrm{loc}}(\rrn)$,
$$
\left\|\left\{\sum_{j\in\mathbb{N}}
\left[\mc\left(f_j\right)\right]^u
\right\}^{\frac{1}{u}}\right\|_{[\HerzSo]^{1/r}}
\lesssim\left\|\left(\sum_{j\in\mathbb{N}}
\left|f_j\right|^u
\right)^{\frac{1}{u}}\right\|_{[\HerzSo]^{1/r}}.
$$
This, combined with the fact
that $\HerzSo$ is a BQBF space and
Lemma \ref{vwmbhll1} with
$X:=\HerzSo$,
further implies that,
for any given $r\in(0,p_-)$
and $u\in(1,\infty)$, and
for any $\{f_j\}_{j\in\mathbb{N}}
\subset L^1_{\mathrm{loc}}(\rrn)$,
$$
\left\|\left\{\sum_{j\in\mathbb{N}}
\left[\mc\left(f_j\right)\right]^u
\right\}^{\frac{1}{u}}\right\|_{[W\HerzSo]^{1/r}}
\lesssim\left\|\left(\sum_{j\in\mathbb{N}}
\left|f_j\right|^u
\right)^{\frac{1}{u}}\right\|_{[W\HerzSo]^{1/r}},
$$
which then completes the
proof of Lemma \ref{vwmbhl}.
\end{proof}

Via the above two lemmas, we next
prove Theorem \ref{wlusin}.

\begin{proof}[Proof of Theorem \ref{wlusin}]
Let all the symbols as in the present theorem.
Then, by the assumption $\m0(\omega)
\in(-\frac{n}{p},\infty)$ and
Theorem \ref{Th3}, we conclude that the local
generalized Herz space $\HerzSo$
under consideration is a BQBF space.
Therefore, from Lemma \ref{wlpx}, it follows
that, to finish the proof of the
present theorem, we only
need to show that $\HerzSo$ satisfies all
the assumptions of Lemma \ref{wlpx}.
Indeed, for any given
$\theta\in(0,p_-)$ and $u\in(1,\infty)$,
using Lemma \ref{vwmbhl} with
$r:=\theta$, we find that,
for any $\{f_{j}\}_{
j\in\mathbb{N}}\subset L^1_{\mathrm{loc}}(\rrn)$,
\begin{equation}\label{wlusine1}
\left\|\left\{\sum_{j\in\mathbb{N}}\left[\mc
(f_{j})\right]^{u}\right\}^
{\frac{1}{u}}\right\|_{[W\HerzSo]^{1/\theta}}
\lesssim\left\|\left\{
\sum_{j\in\mathbb{N}}|f_{j}|^{u}\right\}
^{\frac{1}{u}}\right\|_{[W\HerzSo]^{1/\theta}}.
\end{equation}
In addition, let $\theta_0\in
(\max\{1,p,q\},\infty)$. Then,
repeating an argument similar
to that used in the proof of Theorem \ref{watom}
with $s_0$ and $r_0$ therein replaced,
respectively, by $s$ and $r$, we
conclude that the following four statements
hold true:
\begin{enumerate}
  \item[{\rm(i)}] for any given
  $\theta\in(0,p_-)$ and $u\in(1,\infty)$, and
for any $\{f
_{j}\}_{j\in\mathbb{N}}\subset L_{{\rm loc}}^{1}(\rrn)$,
\begin{equation*}
\left\|\left\{\sum_{j\in\mathbb{N}}\left[\mc
(f_{j})\right]^{u}\right\}^
{\frac{1}{u}}\right\|_{[\HerzSo]^{1/\theta}}
\lesssim\left\|\left\{
\sum_{j\in\mathbb{N}}|f_{j}|^{u}\right\}
^{\frac{1}{u}}\right\|_{[\HerzSo]^{1/\theta}};
\end{equation*}
  \item[{\rm(ii)}] $\HerzSo$ is $\theta_0$-concave;
  \item[{\rm(iii)}] for any $s\in(0,\min\{1,p_-\})$,
  $[\HerzSo]^{1/s}$ is a BBF space;
  \item[{\rm(iv)}] for any given
  $s\in(0,\min\{1,p_-\})$ and $r\in(p_+,\infty)$,
  and for any $f\in L^1_{\mathrm{loc}}(\rrn)$,
  $$
  \left\|\mc^{((r/s)')}(f)\right\|_{([\HerzSo]^{1/s})'}
  \lesssim\left\|f\right\|_{([\HerzSo]^{1/s})'}.
  $$
\end{enumerate}
These, combined with \eqref{wlusine1}, further
imply that the Herz space $\HerzSo$
under consideration satisfies all the assumptions
of Lemma \ref{wlpx} with $p_-$ and $p_+$
as in the present theorem. This then implies that
(i), (ii), (iii), and (iv) of the present
theorem are mutually equivalent and,
for any $f\in W\HaSaHo$,
\begin{align*}
\left\|f\right\|_{W\HaSaHo}
&\sim\left\|S(f)\right\|_{W\HerzSo}
\sim\left\|g(f)\right\|_{W\HerzSo}\\
&\sim\left\|g_{\lambda}^*(f)\right\|_{W\HerzSo}
\end{align*}
with the positive equivalence constants
independent of $f$, which completes
the proof of Theorem \ref{wlusin}.
\end{proof}

Combining Theorem \ref{wlusin}
and Remark \ref{remark6.16},
we immediately conclude
the following littlewood--Paley
function characterizations of
the weak generalized Morrey--Hardy
space $WH\MorrSo$;
we omit the details.

\begin{corollary}\label{wlusinm}
Let $p,\ q\in[1,\infty)$, $\omega\in M(\rp)$
with
$$
-\frac{n}{p}<\m0(\omega)\leq\M0(\omega)<0
$$
and
$$
-\frac{n}{p}<\mi(\omega)\leq\MI(\omega)<0,
$$
and
$$
\lambda\in\left(\max\left\{
2,3-\frac{2(\mw+n/p)}{n}\right\},\infty\right).
$$
Then the following four statements are mutually
equivalent:
\begin{enumerate}
  \item[{\rm(i)}] $f\in W\MorrSo$;
  \item[{\rm(ii)}] $f\in\mathcal{S}'(\rrn)$,
  $f$ vanishes weakly at infinity,
  and $S(f)\in W\MorrSo$;
  \item[{\rm(iii)}] $f\in\mathcal{S}'(\rrn)$,
  $f$ vanishes weakly at infinity,
  and $g(f)\in W\MorrSo$;
  \item[{\rm(iv)}] $f\in\mathcal{S}'(\rrn)$,
  $f$ vanishes weakly at infinity,
  and $g_{\lambda}^*(f)\in W\MorrSo$.
\end{enumerate}
Moreover, for any $f\in W\MorrSo$,
\begin{align*}
\left\|f\right\|_{W\MorrSo}
&\sim\left\|S(f)\right\|_{W\MorrSo}
\sim\left\|g(f)\right\|_{W\MorrSo}\\&
\sim\left\|g_{\lambda}^*(f)\right\|_{W\MorrSo},
\end{align*}
where the positive equivalence
constants are independent of $f$.
\end{corollary}

Now, we turn to investigate
the Littlewood--Paley function
characterizations of the weak generalized
Herz--Hardy space $W\HaSaH$. Indeed,
we have the following conclusion.

\begin{theorem}\label{wlusing}
Let $p,\ q\in(0,\infty)$, $\omega\in M(\rp)$
with $\m0(\omega)\in(-\frac{n}{p},\infty)$
and
$$
-\frac{n}{p}<\mi(\omega)\leq\MI(\omega)<0,
$$
$$
p_-:=\min\left\{p,q,\frac{n}{\Mw+n/p}\right\},
$$
$$
p_+:=\max\left\{p,\frac{n}{\mw+n/p}\right\},
$$
and
$$
\lambda\in\left(\max\left\{
\frac{2}{\min\{1,p_-\}},\,1-\frac{2}{\max\{1,p_+\}}+
\frac{2}{\min\{1,p_-\}}\right\},\infty\right).
$$
Assume that, for any
$f\in\mathcal{S}'(\rrn)$, $S(f)$
and $g_{\lambda}^*(f)$ are as in Definition
\ref{df451}, and $g(f)$ is as in Definition
\ref{df452}.
Then the following four statements are mutually
equivalent:
\begin{enumerate}
  \item[{\rm(i)}] $f\in W\HaSaH$;
  \item[{\rm(ii)}] $f\in\mathcal{S}'(\rrn)$,
  $f$ vanishes weakly at infinity,
  and $S(f)\in W\HerzS$;
  \item[{\rm(iii)}] $f\in\mathcal{S}'(\rrn)$,
  $f$ vanishes weakly at infinity,
  and $g(f)\in W\HerzS$;
  \item[{\rm(iv)}] $f\in\mathcal{S}'(\rrn)$,
  $f$ vanishes weakly at infinity,
  and $g_{\lambda}^*(f)\in W\HerzS$.
\end{enumerate}
Moreover, for any $f\in W\HaSaH$,
\begin{align*}
\left\|f\right\|_{W\HaSaH}
&\sim\left\|S(f)\right\|_{W\HerzS}
\sim\left\|g(f)\right\|_{W\HerzS}\\
&\sim\left\|g_{\lambda}^*(f)\right\|_{W\HerzS},
\end{align*}
where the positive equivalence
constants are independent of $f$.
\end{theorem}

To prove these Littlewood--Paley
function characterizations, we first
establish the following representation
formula of the quasi-norm $\|\cdot\|_{W\HerzS}$.

\begin{lemma}\label{wlusingl1}
Let $p$, $q\in(0,\infty)$ and
$\omega\in M(\rp)$. Then, for any
$f\in\Msc(\rrn)$,
$$
\left\|f\right\|_{W\HerzS}
=\sup_{\xi\in\rrn}\left\|\tau_{\xi}\left(f\right)
\right\|_{W\HerzSo},
$$
where, for any
$\xi\in\rrn$, the operator $\tau_{\xi}$ is
defined as in \eqref{translation}.
\end{lemma}

\begin{proof}
Let all the symbols be as in the
present lemma. Then, using
Definition \ref{df6.0.7}(ii),
Remark \ref{remhs}(ii), and Definition
\ref{df6.0.7}(i), we find that,
for any $f\in\Msc(\rrn)$,
\begin{align*}
\left\|f\right\|_{W\HerzS}
&=\sup_{\lambda\in(0,\infty)}
\left[\lambda\left\|\1bf_{\{y\in\rrn:\
|f(y)|>\lambda\}}\right\|_{\HerzS}\right]\\
&=\sup_{\lambda\in(0,\infty)}
\left\{\lambda\sup_{\xi\in\rrn}\left[
\left\|\1bf_{\{y\in\rrn:\
|f(y)|>\lambda\}}\left(\cdot-\xi\right)
\right\|_{\HerzSo}\right]\right\}\\
&=\sup_{\lambda\in(0,\infty),\,\xi\in\rrn}
\left[\lambda\left\|\1bf_{\{y\in\rrn:\
|\tau_{\xi}(f)(y)|>\lambda\}}\right\|_{\HerzSo}
\right]\\&=\sup_{\xi\in\rrn}\left\{\sup_{\lambda\in(0,\infty)}
\left[\lambda\left\|\1bf_{\{y\in\rrn:\
|\tau_{\xi}(f)(y)|>\lambda\}}\right\|_{\HerzSo}
\right]\right\}\\
&=\sup_{\xi\in\rrn}\left\|\tau_{\xi}\left(f\right)
\right\|_{\HerzSo}.
\end{align*}
This finishes the proof of Lemma \ref{wlusingl1}.
\end{proof}

Based on this formula and the Littlewood--Paley
function characterizations of $W\HaSaHo$
obtained in Theorem \ref{wlusin} above, we next prove
Theorem \ref{wlusing}.

\begin{proof}[Proof of Theorem \ref{wlusing}]
Let all the symbols be as in the present
theorem, $f\in\mathcal{S}'(\rrn)$, and
$\phi\in\mathcal{S}(\rrn)$ satisfy
$\int_{\rrn}\phi(x)\,dx\neq0$.
We first show that (i) implies (ii). Indeed, assume
$f\in W\HaSaH$. Let $M$ denote the radial
maximal function defined as in Definition
\ref{smax}(i). Then, for any $\xi\in\rrn$,
using Theorem \ref{macl}(ii)
with $f$ therein replaced by
$\tau_{\xi}(f)$, Lemma \ref{lusingl1}(ii),
Lemma \ref{wlusingl1} with $f$
therein replaced by $M(f,\phi)$, and
Theorem \ref{Th7.11}(ii), we conclude that
\begin{align}\label{wlusinge1}
&\left\|\tau_{\xi}\left(
f\right)\right\|_{W\HaSaHo}\notag\\
&\quad\sim\left\|M\left(\tau_{\xi}\left(f\right),
\phi\right)\right\|_{W\HerzSo}\sim\left\|\tau_{\xi}\left(M\left(
f,\phi\right)\right)\right\|_{W\HerzSo}\notag\\
&\quad \lesssim\left\|M\left(f,\phi\right)\right\|_{W\HerzS}
\sim\left\|f\right\|_{W\HaSaH}<\infty,
\end{align}
which further implies that
$\tau_{\xi}(f)\in W\HaSaHo$.
Applying this with $\xi:=\0bf$, we further find
that $f\in W\HaSaHo$. This, together with
Theorem \ref{wlusin}, implies that
$f$ vanishes weakly at infinity.
On the other hand, from Lemma \ref{lusingl1}(iii)
with $A:=S$, the fact
that $\tau_{\xi}(f)\in W\HaSaHo$ for any
$\xi\in\rrn$, Theorem \ref{wlusin}, and \eqref{wlusinge1},
it follows that
\begin{align*}
\left\|\tau_{\xi}\left(S\left(f\right)
\right)\right\|_{W\HerzSo}&
\sim\left\|S\left(\tau_{\xi}\left(
f\right)\right)\right\|_{W\HerzSo}\\
&\sim\left\|\tau_{\xi}\left(
f\right)\right\|_{W\HaSaHo}\lesssim
\left\|f\right\|_{W\HaSaH}.
\end{align*}
Using this and Lemma \ref{wlusingl1} again
with $f$ replaced by $S(f)$, we find that
\begin{equation}\label{wlusinge2}
\left\|S\left(f\right)\right\|_{W\HerzS}
\lesssim\left\|f\right\|_{W\HaSaH}<\infty,
\end{equation}
which completes the proof that (i) implies (ii).

Conversely, we next prove that (ii) implies
(i). Indeed, assume that $f$ vanishes weakly at infinity
and $S(f)\in W\HerzS$.
Then, by Lemma \ref{lusingl1}(iv), we conclude that,
for any $\xi\in\rrn$, $\tau_{\xi}(f)$ vanishes
weakly at infinity. In addition,
from Lemma \ref{lusingl1}(iii) with $A:=S$
and Lemma \ref{wlusingl1} with
$f$ therein replaced by $S(f)$, we deduce that,
for any $\xi\in\rrn$,
\begin{align}\label{wlusinge3}
\left\|S\left(\tau_{\xi}\left(f\right)\right)
\right\|_{W\HerzSo}&=\left\|\tau_{\xi}
\left(S\left(f\right)\right)\right\|_{W\HerzSo}
\notag\\&\leq\left\|S\left(f\right)
\right\|_{W\HerzS}<\infty,
\end{align}
which, combined with the fact that
$\tau_{\xi}(f)$ vanishes weakly at infinity
and Theorem \ref{wlusin}, further implies that
$\tau_{\xi}(f)\in W\HaSaHo$.
Moreover, for any $\xi\in\rrn$,
using Lemma \ref{lusingl1}(ii),
Theorem \ref{macl}(ii) with $f$ replaced
by $\tau_{\xi}(f)$, Theorem \ref{wlusin} again,
and \eqref{wlusinge3}, we conclude that
\begin{align*}
\left\|\tau_{\xi}\left(
M\left(f,\phi\right)\right)\right\|_{W\HerzSo}
&=\left\|M\left(\tau_{\xi}
\left(f\right),\phi\right)\right\|_{W\HerzSo}\\
&\sim\left\|\tau_{\xi}\left(
f\right)\right\|_{W\HaSaHo}
\sim\left\|S\left(\tau_{\xi}\left(
f\right)\right)\right\|_{W\HerzSo}\\
&\lesssim\left\|S\left(f\right)\right\|_{W\HerzS},
\end{align*}
which, together with Theorem \ref{Th7.11}(ii)
and Lemma \ref{wlusingl1}
with $f$ replaced by $M(f,\phi)$, further implies that
\begin{align}\label{wlusinge4}
\left\|f\right\|_{W\HaSaH}
&\sim\left\|M\left(f,\phi\right)\right\|
_{W\HerzS}\notag\\
&\sim\sup_{\xi\in\rrn}
\left\{\left\|\tau_{\xi}\left(M\left(
f,\phi\right)\right)\right\|_{W\HerzSo}\right\}\notag\\
&\lesssim\left\|S\left(f\right)\right\|_{W\HerzS}<\infty,
\end{align}
which completes the proof that (ii) implies (i).
Therefore, (i) is equivalent to (ii). Moreover,
combining \eqref{wlusinge2} and \eqref{wlusinge4},
we find that
$$
\left\|f\right\|_{W\HerzS}
\sim\left\|S\left(f\right)\right\|_{W\HaSaH}.
$$

Now, repeating an argument similar
to that used in the estimations
of both \eqref{wlusinge2} and
\eqref{wlusinge4} above with $S$ therein
replaced by $g$ or $g_{\lambda}^*$, we conclude that
(i) is equivalent to (iii) and (i)
is also equivalent to (iv). Moreover, for
any $A\in\{g,g_{\lambda}^*\}$ and
$f\in W\HaSaH$, it holds true that
$$
\left\|f\right\|_{W\HaSaH}
\sim\left\|A\left(f\right)\right\|_{W\HerzS}.
$$
This then finishes the proof of Theorem \ref{wlusing}.
\end{proof}

Using the above theorem and Remark \ref{remark6.16},
we immediately obtain the following Littlewood--Paley
function characterizations of
the weak generalized Morrey--Hardy
space $WH\MorrS$; we omit the
details.

\begin{corollary}
Let $p$, $q$, $\omega$, and $\lambda$
be as in Corollary \ref{wlusinm}.
Then the following four statements are mutually
equivalent:
\begin{enumerate}
  \item[{\rm(i)}] $f\in W\MorrS$;
  \item[{\rm(ii)}] $f\in\mathcal{S}'(\rrn)$,
  $f$ vanishes weakly at infinity,
  and $S(f)\in W\MorrS$;
  \item[{\rm(iii)}] $f\in\mathcal{S}'(\rrn)$,
  $f$ vanishes weakly at infinity,
  and $g(f)\in W\MorrS$;
  \item[{\rm(iv)}] $f\in\mathcal{S}'(\rrn)$,
  $f$ vanishes weakly at infinity,
  and $g_{\lambda}^*(f)\in W\MorrS$.
\end{enumerate}
Moreover, for any $f\in W\MorrS$,
\begin{align*}
\left\|f\right\|_{W\MorrS}
&\sim\left\|S(f)\right\|_{W\MorrS}
\sim\left\|g(f)\right\|_{W\MorrS}\\
&\sim\left\|g_{\lambda}^*(f)\right\|_{W\MorrS},
\end{align*}
where the positive equivalence
constants are independent of $f$.
\end{corollary}

\section{Boundedness of Calder\'{o}n--Zygmund Operators}

Let $\delta\in(0,1)$,
$d\in\zp$, $K$ be a
standard kernel as in Definition \ref{def-s-k}
with $\delta\in(0,1)$, and
$T$ a $d$-order Calder\'{o}n--Zygmund
operator as in Definition \ref{defin-C-Z-s}
with kernel $K$.
The main target of this
section is to investigate the
boundedness of operator $T$ from
generalized Herz--Hardy spaces to weak generalized
Herz--Hardy spaces.
First, we establish two boundedness
criteria of $T$ from
Hardy spaces $H_X(\rrn)$
to weak Hardy spaces $WH_X(\rrn)$,
which improve the
boundedness of convolutional type
Calder\'{o}n--Zygmund operators from
$H_X(\rrn)$ to $WH_X(\rrn)$ obtained in
\cite[Theorem 6.3]{ZWYY}.
As applications, we obtain the
boundedness of $T$ from generalized
Herz--Hardy spaces to weak generalized
Herz--Hardy spaces
even in the critical case $p=\frac{n}{n+d+\delta}$.

Next, we show the following
boundedness of Calder\'{o}n--Zygmund
operators about Hardy-type spaces
associated with ball quasi-Banach function
spaces.

\begin{theorem}\label{wczx}
Let $d\in\zp$, $K$ be a $d$-order standard kernel
with $\delta\in(0,1)$, and
$T$ a $d$-order Calder\'{o}n--Zygmund
operator with kernel $K$
having the vanishing moments up to $d$.
Let $X$ be a ball quasi-Banach function
space, $Y$ a linear space equipped with a
quasi-seminorm $\|\cdot\|_Y$,
and $Y_0$ a linear space equipped with
a quasi-seminorm $\|\cdot\|_{Y_0}$,
and let $r_1\in(0,\infty)$, $\eta\in(1,\infty)$,
and $0<\theta<s<s_0\leq1$ be such that
\begin{enumerate}
  \item[{\rm(i)}] for the above $\theta$ and
  $s$, Assumption \ref{assfs} holds true;
  \item[{\rm(ii)}] both $\|\cdot\|_{Y}$
  and $\|\cdot\|_{Y_0}$
  satisfy Definition \ref{Df1}(ii);
  \item[{\rm(iii)}] $\1bf_{B(\0bf,1)}
  \in Y_0$;
  \item[{\rm(iv)}] for any $f\in\Msc(\rrn)$,
  $$
  \|f\|_{X^{1/s}}
  \sim\sup\left\{\|fg\|_{L^1(\rrn)}:\
  \|g\|_{Y}=1\right\}
  $$
  and
  $$
  \|f\|_{X^{1/s_{0}}}
  \sim\sup\left\{\|fg\|_{L^1(\rrn)}:\
  \|g\|_{Y_0}=1\right\}
  $$
  with the positive equivalence constants
  independent of $f$;
  \item[{\rm(v)}] $\mc^{(\eta)}$ is bounded
  on both $Y$ and $Y_0$;
  \item[{\rm(vi)}] $\mc$ is bounded
  on $(WX)^{1/r_1}$;
  \item[{\rm(vii)}] there exists a positive
  constant $C$ such that, for any $\{f_j\}_{j\in\mathbb{N}}
  \subset L^1_{\mathrm{loc}}(\rrn)$,
  $$
  \left\|\left\{\sum_{j\in\mathbb{N}}
  \left[\mc\left(f_j
  \right)\right]^{\frac{n+d+\delta}{n}}
  \right\}^{\frac{n}{n+d+\delta}}
  \right\|_{(WX)^{\frac{n+d+\delta}{n}}}
  \leq C\left\|\left(\sum_{j\in\mathbb{N}}
  \left|f_j\right|^{\frac{n+d+\delta}{n}}
  \right)^{\frac{n}{n+d+\delta}}
  \right\|_{X^{\frac{n+d+\delta}{n}}}.
  $$
\end{enumerate}
If $\theta\in(0,\frac{n}{n+d}]$, then
$T$ is well defined on $H_X(\rrn)$ and
there exists a positive
constant $C$ such that,
for any $f\in H_X(\rrn)$,
$$
\|T(f)\|_{WH_X(\rrn)}\leq C\|f\|_{H_X(\rrn)}.
$$
\end{theorem}

\begin{remark}
We should point out that
Theorem \ref{wczx} is an improved
version of the known boundedness
of convolutional type Calder\'{o}n--Zygmund
operators from $H_X(\rrn)$ to
$WH_X(\rrn)$ established by Zhang et al.\
in \cite[Theorem 6.3]{ZWYY}. Indeed, in Theorem
\ref{wczx}, if $d=0$, $K(x,y)\equiv K_1(x-y)$
for some $K_1\in L^1_{\mathrm{loc}}(\rrn\setminus
\{\0bf\})$, $Y\equiv(X^{1/s})'$,
and $Y_0\equiv(X^{1/s_0})'$, then
$T$ coincides with the \emph{convolution
$\delta$-type Calder\'{o}n--Zygmund operator}
as in \cite[Theorem 6.3]{ZWYY}, and
this proposition goes back to
\cite[Theorem 6.3]{ZWYY}.
\end{remark}

To establish this boundedness
criterion, we first show the
following boundedness of Calder\'{o}n--Zygmund
operators.

\begin{proposition}\label{wczxa}
Let $d\in\zp$, $K$ be a $d$-order
standard kernel with $\delta\in(0,1)$,
and $T$ a $d$-order
Calder\'{o}n--Zygmund operator with kernel
$K$ having the vanishing moments up to order $d$.
Let $X$ be a ball quasi-Banach
function space satisfying Assumption
\ref{assfs} with some
$0<\theta<s\leq1$ and Assumption \ref{assas}
with the same $s$ and some $r_0\in(1,\infty)$.
Assume that there exists an $r_1\in(0,\infty)$
such that $\mc$ is bounded on $(WX)^{1/r_1}$, and
there exists a positive constant
$C$ such that,
for any $\{f_j\}_{j\in\mathbb{N}}
\subset L^1_{\mathrm{loc}}(\rrn)$,
\begin{align}\label{wczxae0}
&\left\|\left\{\sum_{j\in\mathbb{N}}
\left[\mc\left(f_j
\right)\right]^{\frac{n+d+\delta}{n}}
\right\}^{\frac{n}{n+d+\delta}}
\right\|_{(WX)^{\frac{n+d+\delta}{n}}}\notag\\
&\quad\leq C\left\|\left(\sum_{j\in\mathbb{N}}
\left|f_j\right|^{\frac{n+d+\delta}{n}}
\right)^{\frac{n}{n+d+\delta}}
\right\|_{X^{\frac{n+d+\delta}{n}}}.
\end{align}
If $\theta\in(0,\frac{n}{n+d}]$ and $X$ has
an absolutely continuous quasi-norm,
then $T$ has a unique extension on $H_X(\rrn)$
and there exists a positive constant $C$
such that, for any $f\in H_X(\rrn)$,
$$
\|T(f)\|_{WH_X(\rrn)}\leq C\|f\|_{H_X(\rrn)}.
$$
\end{proposition}

To show this proposition, we
need the following technique
estimate about radial maximal
functions.

\begin{lemma}\label{wczx es-at}
Let $d\in\zp$, $K$ be a
$d$-order standard kernel with
$\delta\in(0,1)$, and $T$ a $d$-order
Calder\'{o}n--Zygmund with kernel $K$
having the vanishing moments up to order $d$.
Assume that $X$ is a ball quasi-Banach function space,
$\phi\in\mathcal{S}(\rrn)$, and
$r\in[2,\infty)$. Then there exists
a positive constant $C$ such that,
for any $(X,\,r,\,d)$-atom $a$ supported
in the ball $B\in\mathbb{B}$,
$$
M\left(T\left(a\right),\phi\right)\1bf_{
(4B)^\complement}\leq C\frac{1}{\|\1bf_B\|_X}
\left[\mc\left(\1bf_B\right)
\right]^{\frac{n+d+\delta}{n}},
$$
where the radial maximal function
$M$ is defined as in Definition \ref{smax}(i).
\end{lemma}

\begin{proof}
Let all the symbols
be as in the present lemma and
$a$ an $(X,\,r,\,d)$-atom supported
in the ball $B:=B(x_0,r_0)$ with
$x_0\in\rrn$ and $r_0\in(0,\infty)$.
Then, combining Definition \ref{atomx}(iii),
the fact that $T$ has the vanishing
moments up to order $d$, and Definition
\ref{Def-T-s-v}, we conclude that,
for any $\gamma\in\zp^n$ with $|\gamma|\leq d$,
\begin{equation}\label{wczx es-ate1}
\int_{\rrn}T(a)(x)x^{\gamma}\,dx=0.
\end{equation}
Fix a $t\in(0,\infty)$ and
an $x\in[B(x_0,4r_0)]^{\complement}$.
Then, from \eqref{wczx es-ate1}, it follows that
\begin{align}\label{wczx es-ate2}
&\left|T(a)\ast\phi_t(x)\right|\notag\\
&\quad=\left|\int_{\rrn}\phi_t(x-y)T(a)(y)\,dy\right|
\notag\\&\quad=\left|
\frac{1}{t^n}\int_{\rrn}\left[\phi\left(
\frac{x-y}{t}\right)-\sum_{\gfz{\gamma\in\zp^n}{|\gamma|
\leq d}}\frac{\partial^{\gamma}(\phi(\frac{x-\cdot}
{t}))(x_0)}{\gamma!}(y-x_0)^{\gamma}\right]
T(a)(y)\,dy\right|\notag\\
&\quad\leq\frac{1}{t^n}\left\{
\int_{|y-x_0|<2r_0}\left|\phi\left(
\frac{x-y}{t}\right)-\sum_{\gfz{\gamma\in\zp^n}{|\gamma|
\leq d}}\frac{\partial^{\gamma}(\phi(\frac{x-\cdot}
{t}))(x_0)}{\gamma!}(y-x_0)^{\gamma}\right|
\right.\notag\\&\qquad\left.\times\left|
T(a)(y)\right|\,dy+\int_{2r_0\leq|y-x_0|<
\frac{|x-x_0|}{2}}\cdots+\int_{|y-x_0|\geq
\frac{|x-x_0|}{2}}\cdots\right\}\notag\\
&\quad=:\mathrm{VIII}_1+\mathrm{VIII}_2+\mathrm{VIII}_3.
\end{align}
We next estimate $\mathrm{VIII}_1$, $\mathrm{VIII}_2$,
and $\mathrm{VIII}_3$, respectively.

First, we deal with $\mathrm{VIII}_1$. Indeed,
applying the Taylor remainder theorem, we find
that, for any $y\in\rrn$ with
$|y-x_0|<2r_0$, there exists a $t_y\in(0,1)$
such that
\begin{align}\label{wczx es-ate3}
\mathrm{VIII}_1&=\frac{1}{t^n}
\int_{|y-x_0|<2r_0}\left|\sum_{\gfz{\gamma\in\zp^n}{|\gamma|=d+1}}\frac{\partial^{
\gamma}(\phi(\frac{x-\cdot}{t}))(t_yy+(1-t_y)x_0)}
{\gamma!}\right.\notag\\
&\quad\times(y-x_0)^{\gamma}\Bigg|\left|
T(a)(y)\right|\,dy\notag\\&\lesssim
\frac{1}{t^{n+d+1}}\int_{|y-x_0|<2r_0}
\sum_{\gfz{\gamma\in\zp^n}{|\gamma|=d+1}}\left|
\partial^{\gamma}\phi\left(
\frac{x-t_yy-(1-t_y)x_0}
{t}\right)\right|\notag\\&\quad
\times\left|y-x_0\right|^{d+1}\left|
T(a)(y)\right|\,dy\notag\\
&\lesssim\int_{|y-x_0|<2r_0}\frac{|y-x_0|^{d+1}}
{|x-t_yy-(1-t_y)x_0|^{n+d+1}}\left|T(a)(y)\right|\,dy.
\end{align}
Notice that, for any $y\in\rrn$ with
$|y-x_0|<2r_0$, we have
$$|y-x_0|<\frac{1}{2}(4r_0)\leq
\frac{1}{2}|x-x_0|,$$ which,
together with the fact that $t_y\in(0,1)$,
further implies that
\begin{equation}\label{wczx es-ate8}
\left|x-t_yy-(1-t_y)x_0\right|
\geq\left|x-x_0\right|-t_y\left|y-x_0\right|
>\frac{1}{2}\left|x-x_0\right|.
\end{equation}
From this, \eqref{wczx es-ate3},
the H\"{o}lder inequality,
Lemma \ref{czop} with $p:=r$, and Definition
\ref{atomx}(ii), we deduce that
\begin{align}\label{wczx es-ate4}
\mathrm{VIII}_1&\lesssim\frac{r_0^{d+1}}{|x-x_0|^{n+d+1}}
\left\|T(a)\1bf_{B(x_0,2r_0)}
\right\|_{L^1(\rrn)}\notag\\
&\lesssim\frac{r_0^{d+1}}{|x-x_0|^{n+d+1}}
\left|B(x_0,2r_0)\right|^{1-\frac{1}{r}}
\left\|T(a)\1bf_{B(x_0,2r_0)}
\right\|_{L^r(\rrn)}\notag\\
&\lesssim\frac{r_0^{d+1}}{|x-x_0|^{n+d+1}}
\left|B(x_0,2r_0)\right|^{1-\frac{1}{r}}
\left\|a\right\|_{L^r(\rrn)}\notag\\
&\lesssim\frac{1}{\|\1bf_{B}\|_{X}}
\frac{r_0^{n+d+1}}{|x-x_0|^{n+d+1}}.
\end{align}
Moreover, by the assumption
$|x-x_0|\geq4r_0$, we find that
$\frac{r_0}{|x-x_0|}\in(0,\frac{1}{4}]$.
This, combined with \eqref{wczx es-ate4} and
the assumption $\delta\in(0,1)$,
further implies that
\begin{align}\label{wczx es-ate5}
\mathrm{VIII}_1\lesssim\frac{1}{\|\1bf_{B}\|_{X}}
\frac{r_0^{n+d+\delta}}{|x-x_0|^{n+d+\delta}},
\end{align}
which is the desired estimate of
$\mathrm{VIII}_1$.

Now, we estimate $\mathrm{VIII}_2$.
For this purpose, we first estimate
$T(a)(y)$ with $y\in[B(x_0,2r_0)]^{\complement}$.
Indeed, fix a $y\in[B(x_0,2r_0)]^{\complement}$.
By \eqref{czp} with $f:=a$, Definition
\ref{atomx}(iii), and the Taylor remainder
theorem, we conclude that, for any
$z\in B(x_0,r_0)$, there exists a
$\tau_z\in[0,1]$ such that
\begin{align}\label{wczx es-ate6}
&\left|T(a)(y)\right|\notag\\
&\quad=\left|\int_{\rrn}K(y,z)a(z)\,dz\right|
\notag\\&\quad=\left|\int_{\rrn}
\left[K(y,z)-\sum_{\gfz{\beta\in\zp^n}
{|\beta|\leq d}}\frac{\partial_{(2)}^{
\beta}K(y,x_0)}{\beta!}(z-x_0)^{\beta}
\right]a(z)\,dz\right|\notag\\
&\quad=\left|\int_{\rrn}
\sum_{\gfz{\beta\in\zp^n}{|\beta|=d}}
\frac{\partial_{(2)}^{\beta}K(y,\tau_zz+(1-\tau_z)
x_0)-\partial_{(2)}^{\beta}
K(y,x_0)}{\beta!}\right.\notag\\
&\qquad\times(z-x_0)^{\beta}a(z)\,dz\Bigg|.
\end{align}
Observe that, for any $z\in B(x_0,r_0)$,
$$
\left|x_0-\tau_zz
-(1-\tau_z)x_0\right|\leq|z-x_0|<r_0
\leq\frac{1}{2}|y-x_0|.
$$
Combining this, \eqref{wczx es-ate6}, \eqref{regular2-s}
with $\gamma$, $x$, $y$, and $z$ therein
replaced, respectively, by
$\beta$, $y$, $x_0$, and $\tau_zz+(1-\tau_z)x_0$
with $z\in B(x_0,r_0)$, the H\"{o}lder inequality,
and Definition \ref{atomx}(ii), we further find
that, for any $y\in[B(x_0,2r_0)]^{\complement}$,
\begin{align}\label{wczx es-ate7}
\left|T(a)(y)\right|
&\lesssim\int_{|z-x_0|<r_0}\frac{|
x_0-\tau_zz-(1-\tau_z)x_0|^{\delta}}
{|y-x_0|^{n+d+\delta}}\left|z-x_0\right|^{d}
\left|a(z)\right|\,dz\notag\\
&\lesssim\frac{r_0^{d+\delta}}{|y-x_0|^{n+d+\delta}}
\left\|a\1bf_{B(x_0,r_0)}\right\|_{L^1(\rrn)}\notag\\
&\lesssim\frac{r_0^{d+\delta}}{|y-x_0|^{n+d+\delta}}
\left|B(x_0,r_0)\right|^{1-\frac{1}{r}}
\left\|a\1bf_{B(x_0,r_0)}\right\|_{L^r(\rrn)}\notag\\
&\lesssim\frac{1}{\|\1bf_B\|_X}
\frac{r_0^{n+d+\delta}}{|y-x_0|^{n+d+\delta}},
\end{align}
which is the desired estimate of $T(a)(y)$
with $y\in[B(x_0,2r_0)]^{\complement}$.
This, together with the Taylor remainder
theorem, \eqref{wczx es-ate8}, and
the assumption that $\delta\in(0,1)$,
further implies that, for any $y\in\rrn$
satisfying $2r_0\leq|y-x_0|<\frac{|x-x_0|}{2}$,
there exists a $t_y\in(0,1)$ such that
\begin{align}\label{wczx es-ate9}
\mathrm{VIII}_2&=
\frac{1}{t^n}\int_{2r_0\leq|y-x_0|
<\frac{|x-x_0|}{2}}
\left|\sum_{\gfz{\gamma\in\zp^n}
{|\gamma|=d+1}}\frac{\partial^{
\gamma}(\phi(\frac{x-\cdot}{t}))(t_yy+(1-t_y)x_0)}
{\gamma!}(y-x_0)^{\gamma}\right|\notag\\
&\quad\times\left|T(a)(y)\right|\,dy\notag\\
&\lesssim\frac{1}{t^{n+d+1}}\int_{2r_0\leq|y-x_0|
<\frac{|x-x_0|}{2}}\sum_{\gfz{\gamma\in\zp^n}
{|\gamma|=d+1}}\left|\partial^{\gamma}
\phi\left(\frac{x-t_yy-(1-t_y)x_0}{t}\right)
\right|\notag\\&\quad\times
|y-x_0|^{d+1}\frac{1}{\|\1bf_B\|_X}
\frac{r_0^{n+d+\delta}}{|y-
x_0|^{n+d+\delta}}\,dy\notag\\
&\lesssim\frac{r_0^{n+d+\delta}}{
\|\1bf_B\|_X}\int_{2r_0\leq|y-x_0|<\frac{|x-x_0|}{2}}
\frac{1}{|x-t_yy-(1-t_y)x_0|^{n+d+1}}\notag\\
&\quad\times\frac{1}{|y-x_0|^{n+\delta-1}}\,dy\notag\\
&\lesssim\frac{1}{\|\1bf_B\|_X}
\frac{r_0^{n+d+\delta}}{|x-x_0|^{n+d+1}}
\int_{|y-x_0|<\frac{|x-x_0|}{2}}
\frac{1}{|y-x_0|^{n+\delta-1}}\,dy\notag\\
&\sim\frac{1}{\|\1bf_B\|_X}\frac{r_0^{n+d+\delta}}
{|x-x_0|^{n+d+\delta}}.
\end{align}
This is the desired estimate of $\mathrm{VIII}_2$.

Finally, we deal with $\mathrm{VIII}_3$.
Indeed, for any $y\in\rrn$ satisfying
$|y-x_0|\geq\frac{|x-x_0|}{2}$, we have
$
|y-x_0|\geq2r_0.
$
By this, \eqref{wczx es-ate7},
and the assumption $\delta\in(0,1)$,
we conclude that
\begin{align*}
\mathrm{VIII}_3&\lesssim
\int_{|y-x_0|\geq\frac{|x-x_0|}{2}}
\left|\phi_t(x-y)\right|\left|
T(a)(y)\right|\,dy\notag\\
&\quad+\int_{|y-x_0|\geq\frac{|x-x_0|}{2}}
\sum_{\gfz{\gamma\in\zp^n}
{|\gamma|\leq d}}\frac{1}{t^{n+|\gamma|}}
\left|\partial^{\gamma}\phi\left(
\frac{x-x_0}{t}\right)\right||y-x_0|
^{|\gamma|}\left|T(a)(y)\right|\,dy\notag\\
&\lesssim\int_{|y-x_0|\geq\frac{|x-x_0|}{2}}
\left|\phi_t(x-y)\right|\frac{1}{
\|\1bf_B\|_X}\,dy\frac{r_0^{n+d+\delta}}
{|y-x_0|^{n+d+\delta}}\,dy\notag\\
&\quad+\sum_{\gfz{\gamma\in\zp}
{|\gamma|\leq d}}\int_{|y-x_0|\geq\frac{|x-x_0|}{2}}
\frac{1}{|x-x_0|^{n+|\gamma|}}
|y-x_0|^{|\gamma|}\notag\\
&\quad\times\frac{1}{\|\1bf_B\|_X}
\frac{r_0^{n+d+\delta}}{|y-x_0|
^{n+d+\delta}}\,dy\notag\\
&\lesssim\frac{r_0^{n+d+\delta}}{\|\1bf_B\|_X}
\left[\frac{1}{|x-x_0|^{n+d+\delta}}
\left\|\phi\right\|_{L^1(\rrn)}\right.\notag\\
&\quad\left.+\sum_{\gfz{\gamma\in\zp^n}
{|\gamma|\leq d}}\frac{1}{|x-x_0|^{n+|\gamma|}}
\int_{|y-x_0|\geq\frac{|x-x_0|}{2}}
\frac{1}{|y-x_0|^{n+d+\delta-|\gamma|}}
\,dy\right]\notag\\
&\sim\frac{1}{\|\1bf_B\|_X}\frac{r_0^{n+d+\delta}}
{|x-x_0|^{n+d+\delta}},
\end{align*}
which is the desired estimate of $\mathrm{VIII}_3$.
Combining this, \eqref{wczx es-ate2},
\eqref{wczx es-ate5}, and \eqref{wczx es-ate9},
we further find that,
for any $t\in(0,\infty)$ and
$x\in[B(x_0,4r_0)]^{\complement}$,
$$
\lf|T(a)\ast\phi_t(x)\r|\lesssim
\frac{1}{\|\1bf_B\|_X}\frac{r_0^{n+d+\delta}}
{|x-x_0|^{n+d+\delta}},
$$
which, together with the
arbitrariness of $t$ and an argument
similar to that used in the estimation
of \eqref{atogxl1e1} with $n+d+1$
therein replaced by $n+d+\delta$, implies that
$$
M\left(T\left(a\right),\phi\right)(x)
\lesssim\frac{1}{\|\1bf_B\|_{X}}
\left[\mc\left(\1bf_B
\right)(x)\right]^{\frac{n+d+\delta}{n}}.
$$
This then finishes the proof
of Lemma \ref{wczx es-at}.
\end{proof}

Applying Lemma \ref{wczx es-at}, we next
show Proposition \ref{wczxa}.

\begin{proof}[Proof of Proposition
\ref{wczxa}]
Let all the symbols be as in
the present proposition and
$r:=\max\{2,r_0\}$.
Then, by the
H\"{o}lder inequality, we find that,
for any $f\in L^1_{\mathrm{loc}}(\rrn)$,
$x\in\rrn$, and $B\in\mathbb{B}$ satisfying
that $x\in B$,
\begin{align*}
\left[\frac{1}{|B|}
\int_{B}|f(y)|^{(r/s)'}
\,dy\right]^{\frac{1}{(r/s)'}}
&\leq\left[\frac{1}{|B|}
\int_{B}|f(y)|^{(r_0/s)'}
\,dy\right]^{\frac{1}{(r_0/s)'}}\\
&\leq\mc^{((r_0/s)')}(f)(x),
\end{align*}
which further implies that
$$
\mc^{((r/s)')}(f)(x)\leq
\mc^{((r_0/s)')}(f)(x).
$$
From this, the fact
that Assumption \ref{assas}
holds true with $s$ and $r_0$
as in the present proposition,
Remark \ref{remark127},
and Definition \ref{Df1}(ii),
it follows that, for any
$f\in L^1_{\mathrm{loc}}(\rrn)$,
\begin{equation}\label{wczxae1}
\left\|\mc^{((r/s)')}(f)
\right\|_{(X^{1/s})'}\leq
\left\|\mc^{((r_0/s)')}(f)
\right\|_{(X^{1/s})'}\lesssim
\left\|f\right\|_{(X^{1/s})'}.
\end{equation}
On the other hand, using the assumption
$\theta\in(0,\frac{n}{n+d}]$,
we find that
$
d\leq n(\frac{1}{\theta}-1).
$
This further implies that
\begin{equation}\label{wczxae2}
d\leq\left\lfloor n\left(
\frac{1}{\theta}-1\right)\right\rfloor
=:d_X.
\end{equation}
Then, combining
\eqref{wczxae1}, the
assumption that $X$ satisfies Assumption
\ref{assfs} for $\theta$ and $s$ as
in the present proposition,
and Lemma \ref{finiatomll1},
we conclude that, for any
$g\in H^{X,r,d_X,s}_{\mathrm{fin}}(\rrn)$,
\begin{equation}\label{wczxae3}
\|g\|_{H^{X,r,d_X,s}_{\mathrm{fin}}(\rrn)}
\sim\|g\|_{H_X(\rrn)}
\end{equation}
with the positive equivalence constants
independent of $g$.

Let $\phi\in\mathcal{S}(\rrn)$
satisfy $\supp(\phi)\subset B(\0bf,1)$
and $\int_{\rrn}\phi(x)\,dx\neq0$,
$f\in H^{X,r,d_X,s}_{\mathrm{fin}}(\rrn)$,
$m\in\mathbb{N}$,
$\{\lambda_{j}\}_{j=1}^m\subset
[0,\infty)$, and $\{a_j\}_{j=1}^m$
of $(X,\,r,\,d_X)$-atoms supported,
respectively, in the balls $\{
B_j\}_{j=1}^m\subset\mathbb{B}$
such that
$f=\sum_{j=1}^{m}\lambda_ja_j$.
Using this and the linearity
of $T$, we find that
$
T(f)=\sum_{j=1}^{m}\lambda_jT(a_j).
$
This implies that
\begin{align*}
M\left(T\left(f\right),\phi\right)
&=\sup_{t\in(0,\infty)}
\left|T\left(f\right)\ast\phi_t\right|\\
&=\sup_{t\in(0,\infty)}\left|
\left[\sum_{j=1}^{m}
\lambda_jT(a_j)\right]\ast\phi_t\right|\\
&\leq\sum_{j=1}^{m}\lambda_j
\sup_{t\in(0,\infty)}\left|T(a_j)
\ast\phi_t\right|=
\sum_{j=1}^{m}\lambda_jM\left(
T\left(a_j\right),\phi\right).
\end{align*}
From this, Lemma \ref{bqbfwx}, and
Definition \ref{Df1}(ii),
we deduce that
\begin{align}\label{wczxae4}
&\left\|M\left(T\left(f\right),
\phi\right)\right\|_{WX}\notag\\&\quad\leq
\left\|\sum_{j=1}^{m}\lambda_j
M\left(T\left(a_j\right),
\phi\right)\right\|_{WX}\notag\\
&\quad\lesssim\left\|\sum_{j=1}^{m}\lambda_j
M\left(T\left(a_j\right),
\phi\right)\1bf_{4B_j}\right\|_{WX}
+\left\|\sum_{j=1}^{m}\lambda_j
M\left(T\left(a_j\right),
\phi\right)\1bf_{(4B_j)^{\complement}}
\right\|_{WX}\notag\\
&\quad=:\mathrm{IX}_1+\mathrm{IX}_2.
\end{align}

Next, we first estimate $\mathrm{IX}_1$.
Indeed, for any $j\in\{1,\ldots,m\}$,
applying Lemma \ref{Atogxl3}
with $\Phi:=\|\phi\|_{L^{\infty}(\rrn)}
\1bf_{B(\0bf,1)}$, $f:=T(a_j)$, the $L^r$
boundedness of $\mc$, Lemma \ref{czop}
with $p:=r$, and Definition \ref{atomx}(ii),
we conclude that
\begin{align*}
\left\|M\left(T\left(a_j\right),\phi
\right)\right\|_{L^r(\rrn)}
\lesssim\left\|\mc\left(T\left(
a_j\right)\right)\right\|_{L^r(\rrn)}
\lesssim\left\|a_j\right\|_{L^r(\rrn)}
\lesssim\frac{|B_j|^{1/r}}{\|\1bf_{B_j}\|_X}.
\end{align*}
By this, Definition \ref{hardyxw},
and an argument similar to that used
in the estimation of $\mathrm{II}_1$
in the proof of Theorem \ref{Atogx} with
$\{\mc(a_j)\1bf_{2B_j}\}_{j\in\mathbb{N}}$
therein replaced by $\{M(T(a_j),\phi)\1bf_{
4B_j}\}_{j=1}^m$,
we obtain
\begin{align}\label{wczxae5}
\mathrm{IX}_1&\lesssim
\left\|\sum_{j=1}^{m}\lambda_j
M\left(T\left(a_j\right),
\phi\right)\1bf_{4B_j}\right\|_{X}
\lesssim\left\|\left[\sum_{j=1}^{m}
\left(\frac{\lambda_j}{\|\1bf_{B_j}
\|_{X}}\right)^s\1bf_{B_j}
\right]^{\frac{1}{s}}\right\|_X,
\end{align}
which is the desired estimate of $\mathrm{IX}_1$.

Now, we deal with $\mathrm{IX}_2$.
Indeed, from \eqref{wczxae2} and Definition
\ref{atomx}, it follows that,
for any $j\in\{1,\ldots,m\}$, $a_j$
is an $(X,\,r,\,d)$-atom supported
in the ball $B_j$. Combining
this, Lemma \ref{wczx es-at}
with $a:=a_j$ for any $j\in\{1,\ldots,m\}$,
Definition \ref{convex}(i), \eqref{wczxae0},
and Lemma \ref{l4320}
with $r$ and $\{a_j\}_{j\in\mathbb{N}}$
therein replaced, respectively,
by $s$ and $\{\frac{\lambda_j}{\|
\1bf_{B_j}\|_X}\1bf_{B_j}\}_{j=1}^m$,
we further conclude that
\begin{align}\label{wczxae7}
\mathrm{IX}_2&\lesssim\left\|
\sum_{j=1}^{m}\frac{\lambda_j}
{\|\1bf_{B_j}\|_X}\left[\mc\left(\1bf_{B_j}\right)
\right]^{\frac{n+d+\delta}{n}}\right\|_{WX}\notag\\
&\sim\left\|\left\{\sum_{j=1}^{m}\frac{\lambda_j}
{\|\1bf_{B_j}\|_X}\left[\mc\left(\1bf_{B_j}\right)
\right]^{\frac{n+d+\delta}{n}}\right\}^{
\frac{n}{n+d+\delta}}\right\|_{(WX)^{\frac{n+d+
\delta}{n}}}^{\frac{n+d+\delta}{n}}\notag\\
&\lesssim\left\|\left(\sum_{j=1}^{m}
\frac{\lambda_j}{\|\1bf_{B_j}\|_X}
\1bf_{B_j}\right)^{\frac{n}
{n+d+\delta}}\right\|_{X^{\frac{n+d+\delta}{n}}}^{
\frac{n+d+\delta}{n}}\notag\\
&\sim\left\|\sum_{j=1}^{m}
\frac{\lambda_j}{\|\1bf_{B_j}\|_X}
\1bf_{B_j}\right\|_{X}
\lesssim\left\|\left[\sum_{j=1}^{m}
\left(\frac{\lambda_j}{\|\1bf_{B_j}\|_X}
\right)^s\1bf_{B_j}
\right]^{\frac{1}{s}}\right\|_{X},
\end{align}
which is the desired estimate of $\mathrm{IX}_2$.
From this, \eqref{wczxae5}, and \eqref{wczxae4},
it follows that
$$
\left\|M\left(T\left(f\right),
\phi\right)\right\|_{WX}
\lesssim\left\|\left[\sum_{j=1}^{m}
\left(\frac{\lambda_j}{\|\1bf_{B_j}\|_X}
\right)^s\1bf_{B_j}
\right]^{\frac{1}{s}}\right\|_{X}.
$$
This, combined with the fact that there exists
an $r_1\in(0,\infty)$ such that
$\mc$ is bounded on $(WX)^{1/r_1}$, Lemmas \ref{bqbfwx}
and \ref{Th5.4l1}(ii), Definition \ref{hardyxw},
the choice of $\{\lambda_j\}_{j=1}^m$,
\eqref{fi}, and \eqref{wczxae3} with $g:=f$,
further implies that
\begin{align}\label{wczxae6}
\left\|T(f)\right\|_{WH_X(\rrn)}
&\sim\left\|M\left(T\left(f\right),
\phi\right)\right\|_{WX}\notag\\
&\lesssim\left\|f\right\|_{H_{\mathrm{fin}}^{
X,r,d_X,s}(\rrn)}\sim\left\|f\right\|_{H_X(\rrn)}.
\end{align}

Furthermore, using the assumption that
$X$ has an absolutely continuous
quasi-norm and \cite[Remark 3.12]{SHYY},
we conclude that
the finite atomic Hardy space $H^{X,r,d_X,s}_
{\mathrm{fin}}(\rrn)$ is dense in the
Hardy space $H_X(\rrn)$. By this,
the fact that \eqref{wczxae6} holds
true for any $f\in H_{\mathrm{fin}}^{
X,r,d_X,s}(\rrn)$,
and a standard density argument, we find that
$T$ has a unique extension on
$H_X(\rrn)$ and, for any $f\in H_X(\rrn)$,
$$
\|T(f)\|_{WH_X(\rrn)}\lesssim\|f\|_{H_X(\rrn)},
$$
which then completes the proof of Proposition
\ref{wczxa}.
\end{proof}

Now, to show Theorem \ref{wczx},
we still need the boundedness
of Calder\'{o}n--Zygmund operators
from weighted Hardy spaces
to weighted weak Hardy spaces.
Let $p\in(0,\infty)$ and
$\upsilon\in A_{\infty}(\rrn)$.
Recall that the \emph{weighted
weak Hardy space} $WH^p_{\upsilon}
(\rrn)$\index{$WH^p_{\upsilon}(\rrn)$} is defined
as in Definition \ref{hardyxw} with
$X:=L^p_{\upsilon}(\rrn)$.
Then we have the following boundedness
of Calder\'{o}n--Zygmund operators from
weighted Hardy spaces to weighted weak
Hardy spaces.

\begin{proposition}\label{wcz-weight}
Let $s_0\in(0,1]$, $\upsilon\in A_1(\rrn)$,
$d\in\zp$, $K$ be a $d$-order
standard kernel with $\delta\in(0,1)$,
and $T$ a $d$-order
Calder\'{o}n--Zygmund operator with kernel $K$
having the vanishing moments up to order $d$.
If $s_0\in[\frac{n}{n+d+\delta},1]$,
then $T$ has a unique extension on
$H^{s_0}_{\upsilon}(\rrn)$
and there exists a positive constant $C$
such that, for any $f\in
H^{s_0}_{\upsilon}(\rrn)$,
$$
\|T(f)\|_{WH^{s_0}_{\upsilon}(\rrn)}
\leq C\|f\|_{
H^{s_0}_{\upsilon}(\rrn)}.
$$
\end{proposition}

To prove Proposition \ref{wcz-weight},
we need some auxiliary lemmas.
The following one is the weighted weak
vector-valued inequality which was
obtained in \cite[Theorem 3.1(a)]{aj80}.

\begin{lemma}\label{wcz-weightl1}
Let $p\in[1,\infty)$, $r\in(1,\infty)$,
and $\upsilon\in A_p(\rrn)$. Then
there exists a positive constant $C$
such that, for any $\{f_j
\}_{j\in\mathbb{N}}\subset L^1_{\mathrm{
loc}}(\rrn)$,
$$
\left\|\left\{\sum_{j\in\mathbb{N}}
\left[\mc\left(f_j
\right)\right]^{r}
\right\}^{\frac{1}{r}}
\right\|_{WL^p_{\upsilon}(\rrn)}
\leq C\left\|\left(\sum_{j\in\mathbb{N}}
\left|f_j\right|^{r}
\right)^{\frac{1}{r}}
\right\|_{L^p_{\upsilon}(\rrn)}.
$$
\end{lemma}

In addition, the following
lemma gives the boundedness of
the Hardy--Littlewood maximal operator
on weighted weak Lebesgue
spaces (see, for instance,
\cite[Corollary 2.6]{lyj16}),
which plays a key role in the
proof of proof of
Proposition \ref{wcz-weight}.

\begin{lemma}\label{wcz-weightl2}
Let $p\in(1,\infty)$ and
$\upsilon\in A_p(\rrn)$. Then
there exists a positive constant
$C$ such that, for any $f\in L^1_{
\mathrm{loc}}(\rrn)$,
$$
\left\|\mc(f)\right\|_{WL^p_{\upsilon}(\rrn)}
\leq C\left\|f\right\|_{WL^p_{\upsilon}(\rrn)}.
$$
\end{lemma}

Via the above two lemmas, we next
prove Proposition \ref{wcz-weight}.

\begin{proof}[Proof of Proposition
\ref{wcz-weight}]
Let all the symbols be as in
the present proposition.
We first claim that all the assumptions
of Proposition \ref{wczxa} hold true
for the weighted Lebesgue space
$L_{\upsilon}^{s_0}(\rrn)$ under consideration.
To this end, let $s\in(0,s_0)$ and
$$
\theta\in\left(0,\min\left\{
s,\frac{n}{n+d}\right\}\right).
$$
Then, by \cite[Remarks 2.4(b),
2.7(b), and 3.4(i)]{WYY}, we conclude that
the following four statements hold true:
\begin{enumerate}
  \item[{\rm(i)}] $L^{s_0}_{
  \upsilon}(\rrn)$ is a BQBF space;
  \item[{\rm(ii)}] for any
$\{f_{j}\}_{j\in\mathbb{N}}
\subset L^1_{{\rm loc}}(\rrn)$,
\begin{equation*}
\left\|\left\{\sum_{j\in\mathbb{N}}
\left[\mc^{(\theta)}
(f_{j})\right]^s\right\}^{1/s}
\right\|_{L^{s_0}_{
\upsilon}(\rrn)}\lesssim
\left\|\left(\sum_{j\in\mathbb{N}}
|f_{j}|^{s}\right)^{1/s}\right\|_{
L^{s_0}_{\upsilon}(\rrn)};
\end{equation*}
  \item[{\rm(iii)}] $[L^{s_0}_{\upsilon}(\rrn)]
  ^{1/s}$ is a BBF space and
there exists an $r_0\in(1,\infty)$ such that,
for any $f\in L^1_{\mathrm{loc}}(\rrn)$,
\begin{equation*}
\left\|\mc^{((r_0/s)')}(f)
\right\|_{([L^{s_0}_{\upsilon}(\rrn)]
^{1/s})'}\lesssim
\left\|f\right\|_{([L^{s_0}_{\upsilon}(\rrn)]
^{1/s})'};
\end{equation*}
  \item[{\rm(iv)}] $L^{s_0}_{\upsilon}(\rrn)$
  has an absolutely continuous quasi-norm.
\end{enumerate}
In addition, let $r_1\in(0,s_0)$.
Then, from the fact that
$\upsilon\in A_1(\rrn)$ and Lemma
\ref{wcz-weightl2} with $p:=\frac{s_0}{r_1}$,
it follows that, for any
$f\in L^1_{\mathrm{loc}}(\rrn)$,
\begin{equation}\label{wcz-weighte1}
\left\|\mc(f)\right\|_{WL^{s_0/r_1}_{\upsilon}(\rrn)}
\lesssim\left\|f\right\|_{WL^{s_0/
r_1}_{\upsilon}(\rrn)}.
\end{equation}
Combining Definitions \ref{convex}(ii) and
\ref{df237}, we find that $$WL^{s_0/
r_1}_{\upsilon}(\rrn)=\left[WL^{s_0}_{\upsilon}
(\rrn)\right]^{1/r_1}.$$
This, together with \eqref{wcz-weighte1},
further implies that
the Hardy--Littlewood maximal operator
$\mc$ is bounded on $[WL^{s_0}_{\upsilon}
(\rrn)]^{1/r_1}$.
Finally, we prove that
\eqref{wczxae0} holds true
with $X:=L^{s_0}_{\upsilon}(\rrn)$.
Indeed, applying Definitions \ref{convex}(i)
and \ref{df237}, and Lemma \ref{wcz-weightl1}
with $r:=\frac{n+d+\delta}{n}$ and
$p:=\frac{n+d+\delta}{n}s_0$, we conclude
that, for any
$\{f_j\}_{j\in\mathbb{N}}
\subset L^1_{\mathrm{loc}}(\rrn)$,
\begin{align*}
&\left\|\left\{\sum_{j\in\mathbb{N}}
\left[\mc\left(f_j
\right)\right]^{\frac{n+d+\delta}{n}}
\right\}^{\frac{n}
{n+d+\delta}}\right\|_{[WL_{\upsilon}^{s_0}
(\rrn)]^{\frac{n+d+\delta}{n}}}\\
&\quad=\left\|\left\{\sum_{j\in\mathbb{N}}
\left[\mc\left(f_j
\right)\right]^{\frac{n+d+\delta}{n}}
\right\}^{\frac{n}
{n+d+\delta}}\right\|_{WL_{\upsilon}^
{\frac{n+d+\delta}{n}s_0}(\rrn)}\\
&\quad\lesssim\left\|\left\{\sum_{j\in\mathbb{N}}
\left[\mc\left(f_j
\right)\right]^{\frac{n+d+\delta}{n}}
\right\}^{\frac{n}
{n+d+\delta}}\right\|_{L_{\upsilon}^
{\frac{n+d+\delta}{n}s_0}(\rrn)}\\&\quad\sim
\left\|\left(\sum_{j\in\mathbb{N}}\left|
f_j\right|^{
\frac{n}{n+d+\delta}}
\right)^{\frac{n}{n+d+\delta}}
\right\|_{[L^{s_0
}_{\upsilon}(\rrn)]^{\frac{n+d+\delta}{n}}},
\end{align*}
which implies that \eqref{wczxae0} holds true
with $X:=L^{s_0}_{\upsilon}(\rrn)$.
This, together with (i) through (iv)
in the proof of the present proposition
and the fact that
$\mc$ is bounded on $[WL_{\upsilon}^{s_0}(\rrn)]^{1/r_1}$,
further implies that
the weighted Lebesgue space $L^{s_0}_{\upsilon}(\rrn)$
under consideration satisfies all the assumptions
of Proposition \ref{wczxa}.
This then finishes the proof of
the above claim. Thus, by Proposition
\ref{wczxa} with $X:=L^{s_0}_{\upsilon}(\rrn)$,
we conclude that $T$ has a unique extension
on $H_{\upsilon}^{s_0}(\rrn)$ and,
for any $f\in H_{\upsilon}^{s_0}(\rrn)$,
$$
\left\|T\left(f\right)\right\|_{WH_{
\upsilon}^{s_0}(\rrn)}
\lesssim\left\|f\right\|_{H_{\upsilon}^{s_0}(\rrn)},
$$
which completes the proof of Proposition
\ref{wcz-weight}.
\end{proof}

Via above preparations, we now
show Theorem \ref{wczx}.

\begin{proof}[Proof of Theorem \ref{wczx}]
Let all the symbols be as in
the present theorem and $f\in H_X(\rrn)$.
Then, by the assumption
$\theta\in(0,
\frac{n}{n+d}]$,
we conclude that
$
d\leq n\left(\frac{1}{\theta}-1\right).
$
This further implies that
\begin{equation}\label{wczxe3}
d\leq\left\lfloor n\left(
\frac{1}{\theta}-1\right)\right\rfloor
=:d_X.
\end{equation}
From this, the assumption (i) of the present
theorem, and Lemma \ref{Atogl4},
it follows that there exist
$\{\lambda_j\}_{j\in\mathbb{N}}
\subset[0,\infty)$ and $\{a_j\}_{j\in
\mathbb{N}}$ of $(X,\,\infty,\,d_X)$-atoms
supported, respectively, in
the balls $\{B_j\}_{j\in\mathbb{N}}\subset
\mathbb{B}$ such that
$f=\sum_{j\in\mathbb{N}}\lambda_ja_j$
in $\mathcal{S}'(\rrn)$, and
\begin{equation}\label{wczxe2}
\left\|\left[\sum_{j\in\mathbb
N}\left(\frac{\lambda_j}{\|\1bf_{B_j}\|_{X}}
\right)^s\1bf_{B_j}
\right]^{\frac{1}{s}}\right\|_{X}\lesssim
\|f\|_{H_X(\rrn)}.
\end{equation}
Combining these and an
argument similar to that used
in the proof of Theorem \ref{bddox}
with Proposition \ref{bddoxl5}
therein replaced by Proposition
\ref{wcz-weight} here,
we find that
\begin{equation}\label{wczxe4}
T(f)=\sum_{j\in\mathbb{N}}
\lambda_jT(a_j)
\end{equation}
in $\mathcal{S}'(\rrn)$.

Next, we prove that
$T(f)\in WH_X(\rrn)$ and
$$\left\|T(f)\right\|_{WH_X(\rrn)}\lesssim
\left\|f\right\|_{H_X(\rrn)}.$$
For this purpose, let
$r\in(\max\{2,s\eta'\},\infty)$
with $\frac{1}{\eta}+\frac{1}{\eta'}=1$.
Then, from \eqref{wczxe3}
and Definition \ref{atomx},
we infer that, for any $j\in\mathbb{N}$,
$a_j$ is an $(X,\,r,\,d)$-atom
supported in the ball $B_j$.
On the other hand, using the assumption (ii) of
the present theorem and
an argument similar to that
used in the estimation of \eqref{bddoxl4e1}
with $(r_0/s)'$ and $(X^{1/s})'$
therein replaced, respectively, by
$\eta$ and $Y$, we conclude that,
for any $f\in L^1_{\mathrm{loc}}(\rrn)$,
\begin{equation}\label{wczxe5}
\left\|\mc^{((r/s)')}(f)\right\|_{Y}
\lesssim\|f\|_Y.
\end{equation}
In addition, let $\phi\in\mathcal{S}(\rrn)$
satisfy $\supp(\phi)\subset B(\0bf,1)$
and $\int_{\rrn}\phi(x)\,dx\neq0$.
Then, from \eqref{wczxe4}, \eqref{wczxe5},
the arguments similar to those used
in the estimations of \eqref{wczxae4},
\eqref{wczxae5}, and \eqref{wczxae7}
with $\{\lambda_j\}_{j=1}^m$
and $\{a_j\}_{j=1}^m$ therein
replaced, respectively,
by $\{\lambda_j\}_{j\in\mathbb{N}}$
and $\{a_j\}_{j\in\mathbb{N}}$,
and \eqref{wczxe2},
we deduce that
\begin{align*}
\left\|M\left(T\left(f\right),
\phi\right)\right\|_{WX}
&\lesssim\left\|\sum_{j\in\mathbb{N}}
\lambda_jM(T(a_j),\phi)\right\|_{WX}\\
&\lesssim\left\|\left[\sum_{j\in\mathbb{N}}
\left(\frac{\lambda_j}{\|\1bf_{B_{j}}\|_X}\right)^s
\1bf_{B_j}\right]^{\frac{1}{s}}\right\|_{X}
\lesssim\left\|f\right\|_{H_X(\rrn)},
\end{align*}
which, combined with the assumption
(vi) of the present theorem,
Lemmas \ref{bqbfwx} and \ref{Th5.4l1}(ii),
and Definition \ref{hardyxw},
further implies that
$$
\left\|T\left(f\right)\right\|_{WH_X(\rrn)}
\sim\left\|M\left(T\left(f\right),
\phi\right)\right\|_{WX}
\lesssim\left\|f\right\|_{H_X(\rrn)}.
$$
This then finishes the proof of Theorem \ref{wczx}.
\end{proof}

Next, we turn to establish
the boundedness of Calder\'{o}n--Zygmund
operators from generalized
Herz--Hardy spaces to weak generalized
Herz--Hardy spaces. Indeed,
the following conclusion gives
the boundedness of Calder\'{o}n--Zygmund
operators from $\HaSaHo$ to $W\HaSaHo$.

\begin{theorem}\label{wczoo}
Let $d\in\zp$, $\delta\in(0,1)$,
$p\in[\frac{n}{n+d+\delta},\infty)$,
$q\in(0,\infty)$,
$K$ be a $d$-order standard kernel as
in Definition \ref{def-s-k}, $T$ a
$d$-order Calder\'{o}n--Zygmund operator
with kernel $K$ having the vanishing moments
up to order $d$, and $\omega\in
M(\rp)$ with
$$
-\frac{n}{p}<\m0(\omega)\leq
\M0(\omega)<n-\frac{n}{p}+d+\delta
$$
and
$$
-\frac{n}{p}<\mi(\omega)\leq
\MI(\omega)<n-\frac{n}{p}+d+\delta.
$$
Then $T$ has a unique extension
on $\HaSaHo$ and there exists
a positive constant $C$ such that,
for any $f\in\HaSaHo$,
$$
\|T(f)\|_{W\HaSaHo}\leq C\|f\|_{
\HaSaHo}.
$$
\end{theorem}

To prove this theorem, we need the following weak type
inequality about the Hardy--Littlewood maximal operator,
which is a part of \cite[Theorem 5.6.6]{LGCF}.

\begin{lemma}\label{wczool10}
Let $r\in(1,\infty)$. Then there exists a positive constant
$C$ such that, for any $\{f_{j}\}_{j\in\mathbb{N}}
\subset L^1_{\mathrm{loc}}(\rrn)$ and $\lambda\in(0,\infty)$,
$$
\lambda\left\|
\1bf_{\{y\in\rrn:\ \mvrys>\lambda\}}\right\|_{L^1(\rrn)}
\leq C\left\|\vr\right\|_{L^1(\rrn)}.
$$
\end{lemma}

Via this estimate, we obtain the weak type
vector-valued inequalities on both local and
global generalized
Herz spaces as follows, which
play an important role in the
proof of the boundedness of
Calder\'{o}n--Zygmund operators
from generalized Herz--Hardy spaces
to weak generalized Herz--Hardy spaces.

\begin{proposition}\label{wczool1}
Let $p\in[1,\infty)$, $q\in(0,\infty)$,
$r\in(1,\infty)$, and $\omega\in M(\rp)$
with
$$
-\frac{n}{p}<\m0(\omega)\leq\M0(\omega)<\frac{n}{p'}
$$
and
$$
-\frac{n}{p}<\mi(\omega)\leq\MI(\omega)<\frac{n}{p'},
$$
where $\frac{1}{p}+\frac{1}{p'}=1$.
Then there exists a positive constant $C$ such that, for any
$\{f_{j}\}_{j\in\mathbb{N}}\subset L^{1}_{{\rm loc}}(\rrn)$,
\begin{equation}\label{wczool1e0}
\left\|\mvr\right\|_{\WHerzSo}\leq C\left\|\vr\right\|_{\HerzSo}
\end{equation}
and
\begin{equation}\label{wczol1e00}
\left\|\mvr\right\|_{\WHerzS}\leq C\left\|\vr\right\|_{\HerzS}.
\end{equation}
\end{proposition}

\begin{proof}
Let all the symbols be as in the present
proposition and $\{f_j\}_{j\in\mathbb{N}}$ be any
given sequence of local integrable functions on $\rrn$.
Then we show the present proposition by considering
the following two cases on $p$.

\emph{Case 1)} $p\in(1,\infty)$.
In this case, by Definition
\ref{df6.0.7} and Theorems \ref{Th3.4}
and \ref{Th3.3},
we find that
\begin{align*}
\left\|\mvr\right\|_{W\HerzSo}&\leq
\left\|\mvr\right\|_{\HerzSo}\\
&\lesssim\left\|\vr\right\|_{\HerzSo}
\end{align*}
and
\begin{align*}
\left\|\mvr\right\|_{W\HerzS}&\leq
\left\|\mvr\right\|_{\HerzS}\\
&\lesssim\left\|\vr\right\|_{\HerzS}.
\end{align*}
These finish the proof of Proposition
\ref{wczool1} in this case.

\emph{Case 2)} $p=1$. In this case,
fix a point $\xi\in\rrn$.
For any $j\in\mathbb{N}$ and
$k\in\mathbb{Z}$,
let
$$
f_{j,k,1}:=f_j\1bf_{B(\xi,\tkmd)},
$$
$$
f_{j,k,2}:=f_j\1bf_{B(\xi,\tka)\setminus B(\xi,\tkmd)},
$$
and
$$
f_{j,k,3}:=f_j\1bf_{[B(\xi,\tka)]^\complement}.
$$
Then, obviously, for any $j\in
\mathbb{N}$ and $k\in\mathbb{Z}$,
$f=f_{j,k,1}+f_{j,k,2}+f_{j,k,3}$.
From this, we deduce that,
for any $\lambda\in(0,\infty)$,
\begin{align}\label{wczool1e13}
&\left\|\1bf_{\{y\in\rrn:\ \mvrys>
\lambda\}}(\cdot+\xi)\right\|_{
\Kmp_{\omega,\0bf}^{1,q}(\rrn)}\notag\\
&\quad\lesssim\left\{\sum_{k\in\mathbb{Z}}
\left[\omega(\tk)\right]^q\right.\notag\\
&\qquad\times\left.\left\|
\1bf_{\{y\in\rrn:\ \mvryo>\lambda/3\}}
\1bf_{B(\xi,\tk)\setminus B(\xi,\tkm)}
\right\|_{L^1(\rrn)}^q\right\}^{\frac{1}{q}}\notag\\
&\qquad+\left\{\sum_{k\in\mathbb{Z}}
\left[\omega(\tk)\right]^q\right.\notag\\
&\left.\qquad\times\left\|
\1bf_{\{y\in\rrn:\ \mvrytw>\lambda/3\}}
\1bf_{B(\xi,\tk)\setminus B(\xi,\tkm)}
\right\|_{L^1(\rrn)}^q\right\}^{\frac{1}{q}}\notag\\
&\qquad+\left\{\sum_{k\in\mathbb{Z}}
\left[\omega(\tk)\right]^q\right.\notag\\
&\left.\qquad\times\left\|
\1bf_{\{y\in\rrn:\ \mvryth>\lambda/3\}}
\1bf_{B(\xi,\tk)\setminus B(\xi
,\tkm)}\right\|_{L^1(\rrn)}^q\right\}^{
\frac{1}{q}}\notag\\
&\quad=:\mathrm{II}_{\lambda,1}+\mathrm{II}_{
\lambda,2}+\mathrm{II}_{\lambda,3}.
\end{align}
We next estimate $\mathrm{II}_{\lambda,1}$,
$\mathrm{II}_{\lambda,2}$, and $\mathrm{II}_{
\lambda,3}$, respectively.
To achieve this, let
$$
\eps\in\left(0,\left\{\mw+n,-\Mw\right\}\right)
$$
be a fixed positive constant.
Then, applying Lemma \ref{La3.5}, we conclude that,
for any $0<t<\tau<\infty$,
\begin{equation}\label{wczoott}
\frac{\omega(t)}{\omega(\tau)}\lesssim\left(
\frac{t}{\tau}\right)
^{\min\{\m0(\omega),\mi(\omega)\}-\eps}
\end{equation}
and, for any $0<\tau<t<\infty$,
\begin{equation}\label{wczootd}
\frac{\omega(t)}{\omega(\tau)}\lesssim\left(
\frac{t}{\tau}\right)
^{\max\{\M0(\omega),\MI(\omega)\}+\eps}.
\end{equation}
For the simplicity of the presentation,
let $m:=\mw-\eps$ and $M:=\Mw+\eps$. Thus,
we have $m\in(-n,\infty)$ and $M\in(-\infty,0)$.

Now, we deal with $\mathrm{II}_{\lambda,1}$. Indeed,
using the fact that the
Hardy--Littlewood maximal operator
$\mc$ satisfies the size
condition \eqref{size} (see, for instance,
\cite[Remark 4.4]{HSamko}),
together with the Minkowski integral inequality,
we find that, for any $k\in\mathbb{Z}$ and $x\in\rrn$
satisfying $\tkm\leq|x-\xi|<\tk$,
\begin{align}\label{wczool1e1}
\left\{\sum_{j\in\mathbb{N}}
\left[\mc(f_{j,k,1})(x)\right]^r\right\}^{\frac{1}{r}}
&\lesssim\left\{\sum_{j\in\mathbb{N}}
\left[\int_{\rrn}\frac{|f_{j,k,1}(y)|}{|x-y|^n}
\,dy\right]^r\right\}^{\frac{1}{r}}\notag\\
&\lesssim\int_{\rrn}\left\{\sum_{j\in\mathbb{N}}
\left[\frac{|f_{j,k,1}(y)|}{|x-y|^n}
\right]^r\right\}^{\frac{1}{r}}\,dy\notag\\
&\sim\int_{B(\xi,\tkmd)}\frac{\vrys}{|x-y|^n}\,dy\notag\\
&\sim\sum_{i=-\infty}^{k-2}\int_{B(\xi,2^i)
\setminus B(\xi,2^{i-1})}\frac{\vrys}{|x-y|^n}\,dy.
\end{align}
On the other hand, for any $k\in\mathbb{Z}$,
$i\in\mathbb{Z}\cap(-\infty,k-2]$, and $x,\ y\in\rrn$
satisfying $\tkm\leq|x-\xi|<\tk$ and $2^{i-1}\leq|y-\xi|<2^i$,
we have
\begin{align*}
|x-y|\geq|x-\xi|-|y-\xi|>\tkm-2^i>2^{k-2},
\end{align*}
which, together with \eqref{wczool1e1}, implies that,
for any $k\in\mathbb{Z}$ and $x\in\rrn$
satisfying $\tkm\leq|x-\xi|<\tk$,
\begin{align}\label{wczool1e2}
&\left\{\sum_{j\in\mathbb{N}}
\left[\mc(f_{j,k,1})(x)\right]^r\right\}^
{\frac{1}{r}}\notag\\
&\quad\lesssim\sum_{i=-\infty
}^{k-2}2^{-nk}\left\|
\vr\1bf_{B(\xi,2^i)
\setminus B(\xi,2^{i-1})}\right\|_{L^1(\rrn)}.
\end{align}
In addition, by Lemma \ref{mono}
and \eqref{wczootd}, we conclude that, for any
$k\in\mathbb{Z}$ and $i\in\mathbb{Z}\cap(-\infty,k-2]$,
\begin{align*}
\frac{\omega(\tk)}{\omega(2^i)}\sim
\frac{\omega(\tkmd)}{\omega(2^i)}
\lesssim\left(\frac{\tkmd}{2^i}\right)^M\sim2^{(k-i)M}.
\end{align*}
From this and \eqref{wczool1e2}, it follows that,
for any $k\in\mathbb{Z}$,
\begin{align}\label{wczool1e3}
&\omega(\tk)\left\|\mvro\1bf_{B(\xi,\tk)
\setminus B(\xi,\tkm)}\right\|_{L^1(\rrn)}\notag\\
&\quad\lesssim\omega(\tk)\sum_{i=-\infty}^{k-2}2^{-nk}
\left\|\1bf_{B(\xi,\tk)\setminus
B(\xi,\tkm)}\right\|_{L^1(\rrn)}\notag\\
&\qquad\times\left\|
\vr\1bf_{B(\xi,2^i)
\setminus B(\xi,2^{i-1})}\right\|_{L^1(\rrn)}\notag\\
&\quad\sim\sum_{i=-\infty}^{k-2}
\frac{\omega(\tk)}{\omega(2^i)}\omega(2^i)\left\|
\vr\1bf_{B(\xi,2^i)
\setminus B(\xi,2^{i-1})}\right\|_{L^1(\rrn)}\notag\\
&\quad\lesssim\sum_{i=-\infty}^{k-2}
2^{(k-i)M}\omega(2^i)\left\|
\vr\1bf_{B(\xi,2^i)
\setminus B(\xi,2^{i-1})}\right\|_{L^1(\rrn)}.
\end{align}
Applying this, the fact that, for any $\lambda\in(0,\infty)$,
\begin{equation}\label{wczool1e31}
\1bf_{\{y\in\rrn:\ \mvryo>\lambda/3\}}\leq\frac{3}{\lambda}
\mvro,
\end{equation}
Lemma \ref{l4320},
and the assumption $M\in(-\infty,0)$,
we further conclude that, for any $q\in(0,1]$ and
$\lambda\in(0,\infty)$,
\begin{align}\label{wczool1e4}
\mathrm{II}_{\lambda,1}&\lesssim
\lambda^{-1}\left(\sum_{k\in\mathbb{Z}}
\left[\omega(\tk)\right]^q\left\|
\mvro\1bf_{B(\xi,\tk)\setminus B(\xi,\tkm)}
\right\|_{L^1(\rrn)}^q\right)^{\frac{1}{q}}
\notag\\
&\lesssim\lambda^{-1}\left\{
\sum_{k\in\mathbb{Z}}\left[
\sum_{i=-\infty}^{k-2}
2^{(k-i)M}\omega(2^i)\right.\right.\notag\\
&\left.\left.\quad\times\left\|
\left(\sum_{j\in\mathbb{N}}
|f_j|^r\right)^{\frac{1}{r}}\1bf_{B(\xi,2^i)
\setminus B(\xi,2^{i-1})}\right\|_{L^1(\rrn)}
\right]^q\right\}^{\frac{1}{q}}\notag\\
&\lesssim\lambda^{-1}\left\{
\sum_{k\in\mathbb{Z}}
\sum_{i=-\infty}^{k-2}
2^{(k-i)Mq}\left[\omega(2^i)\right]^q\right.\notag\\
&\lf.\quad\times\left\|
\left(\sum_{j\in\mathbb{N}}
|f_j|^r\right)^{\frac{1}{r}}\1bf_{B(\xi,2^i)
\setminus B(\xi,2^{i-1})}\right\|_{
L^1(\rrn)}^q\right\}^{\frac{1}{q}}\notag\\
&\sim\lambda^{-1}\left\{
\sum_{i\in\mathbb{Z}}\left[\omega(\tk)\right]^q
\left\|\left(\sum_{j\in\mathbb{N}}
|f_j|^r\right)^{\frac{1}{r}}\1bf_{B(\xi,2^i)
\setminus B(\xi,2^{i-1})}\right\|_{L^1(\rrn)}^q\right.\notag\\
&\left.\quad\times\sum_{k=i+2}^{\infty}2^{(k-i)Mq}\right\}^{\frac{1}{q}}\notag\\
&\sim\lambda^{-1}\left\|\left(\sum_{j\in\mathbb{N}}
|f_j|^r\right)^{\frac{1}{r}}(\cdot+\xi)
\right\|_{\Kmp_{\omega,\0bf}^{1,q}(\rrn)}.
\end{align}
Moreover, using \eqref{wczool1e31}, \eqref{wczool1e3},
the H\"{o}lder inequality, and the assumption $M\in(-\infty,0)$,
we find that,
for any $q\in(1,\infty)$ and $\lambda\in(0,\infty)$,
\begin{align}\label{wczool1e5}
\mathrm{II}_{\lambda,1}&\lesssim
\lambda^{-1}\left(\sum_{k\in\mathbb{Z}}
\left[\omega(\tk)\right]^q\left\|
\mvro\1bf_{B(\xi,\tk)\setminus B(\xi,\tkm)}
\right\|_{L^1(\rrn)}^q\right)^{\frac{1}{q}}
\notag\\
&\lesssim\lambda^{-1}\left\{
\sum_{k\in\mathbb{Z}}\left[
\sum_{i=-\infty}^{k-2}
2^{(k-i)M}\omega(2^i)\right.\right.\notag\\
&\left.\left.\quad\times\left\|
\left(\sum_{j\in\mathbb{N}}
|f_j|^r\right)^{\frac{1}{r}}\1bf_{B(\xi,2^i)
\setminus B(\xi,2^{i-1})}\right\|_{L^1(\rrn)}
\right]^q\right\}^{\frac{1}{q}}\notag\\
&\lesssim\lambda^{-1}\left\{
\sum_{k\in\mathbb{Z}}\left[\sum_{i=-\infty}^{k-2}
2^{\frac{(k-i)Mq'}{2}}\right]^{\frac{q}{q'}}\right.\notag\\
&\quad\left.\times\sum_{i=-\infty}^{k-2}
2^{\frac{(k-i)Mq}{2}}\left[\omega(2^i)\right]^q\left\|
\left(\sum_{j\in\mathbb{N}}
|f_j|^r\right)^{\frac{1}{r}}\1bf_{B(\xi,2^i)
\setminus B(\xi,2^{i-1})}\right\|_{L^1(\rrn)}
^q\right\}^{\frac{1}{q}}\notag\\
&\sim\lambda^{-1}\left\{
\sum_{i\in\mathbb{Z}}\left[\omega(2^i)\right]^q
\left\|\left(\sum_{j\in\mathbb{N}}
|f_j|^r\right)^{\frac{1}{r}}\1bf_{B(\xi,2^i)
\setminus B(\xi,2^{i-1})}\right\|_{L^1(\rrn)}^q\right.\notag\\
&\left.\quad\times\sum_{k=i+2}^{\infty}2^{\frac{(k-i)Mq}{2}}
\right\}^{\frac{1}{q}}\notag\\
&\sim\lambda^{-1}\left\|\left(\sum_{j\in\mathbb{N}}
|f_j|^r\right)^{\frac{1}{r}}(\cdot+\xi)
\right\|_{\Kmp_{\omega,\0bf}^{1,q}(\rrn)},
\end{align}
which, together with \eqref{wczool1e4}, then
completes the estimate of $\mathrm{II}_{\lambda,1}$.

We next deal with $\mathrm{II}_{\lambda,2}$.
Indeed, from Lemma \ref{wczool10}, it follows that,
for any $k\in\mathbb{Z}$ and $\lambda\in(0,\infty)$,
\begin{align}\label{wczool1e11}
&\left\|\1bf_{\{y\in\rrn:\
\mvrytw>\lambda/3\}}\1bf_{B(\xi,\tk)\setminus
B(\xi,\tkm)}\right\|_{L^1(\rrn)}\notag\\
&\quad\leq\left\|\1bf_{\{y\in\rrn:\
\mvrytw>\lambda/3\}}\right\|_{L^1(\rrn)}\notag\\
&\quad\lesssim\lambda^{-1}\left\|
\vrtw\right\|_{L^1(\rrn)}\notag\\
&\quad\lesssim\lambda^{-1}\left[\left\|\vr\1bf_{B(\xi,\tka)
\setminus B(\xi,\tk)}\right\|_{L^1(\rrn)}\right.\notag\\
&\qquad+\left\|\vr\1bf_{B(\xi,\tk)
\setminus B(\xi,\tkm)}\right\|_{L^1(\rrn)}\notag
\\&\qquad\left.+
\left\|\vr\1bf_{B(\xi,\tkm)
\setminus B(\xi,\tkmd)}\right\|_{L^1(\rrn)}\right].
\end{align}
In addition, by Lemma \ref{mono}, we conclude that,
for any $k\in\mathbb{Z}$, $$\omega(\tka)\sim
\omega(\tk)\sim\omega(\tkm).$$
Therefore, using \eqref{wczool1e11}, we find that,
for any $\lambda\in(0,\infty)$,
\begin{align}\label{wczool1e12}
\mathrm{II}_{\lambda,2}&\lesssim\lambda^{-1}\left(
\left\{\sum_{k\in\mathbb{Z}}\left[\omega(\tka)\right]^q
\left\|\vr\1bf_{B(\xi,2^k)
\setminus B(\xi,2^{k-1})}\right\|_
{L^1(\rrn)}^q\right\}^{\frac{1}{q}}\right.\notag\\
&\quad+\left\{\sum_{k
\in\mathbb{Z}}\left[\omega(\tk)\right]^q
\left\|\vr\1bf_{B(\xi,2^k)
\setminus B(\xi,2^{k-1})}\right\|
_{L^1(\rrn)}^q\right\}^{\frac{1}{q}}\notag\\
&\quad\left.+
\left\{\sum_{k\in\mathbb{Z}}\left[
\omega(\tkm)\right]^q
\left\|\vr\1bf_{B(\xi,2^k)
\setminus B(\xi,2^{k-1})}\right\|_
{L^1(\rrn)}^q\right\}^{\frac{1}{q}}
\right)\notag\\
&\sim\lambda^{-1}\left\|\vr(\cdot+\xi)
\right\|_{\Kmp_{\omega,\0bf}^{1,q}(\rrn)},
\end{align}
which completes the estimate of $\mathrm{II}_{\lambda,2}$.

Finally, we estimate $\mathrm{II}_{\lambda,3}$. Indeed,
applying the fact that the
Hardy--Littlewood maximal operator
$\mc$ satisfies the size
condition \eqref{size} (see, for instance,
\cite[Remark 4.4]{HSamko}),
together with the Minkowski integral inequality,
we conclude that, for any $k\in\mathbb{Z}$ and $x\in\rrn$
satisfying $\tkm\leq|x-\xi|<\tk$,
\begin{align}\label{wczool1e6}
\left\{\sum_{j\in\mathbb{N}}
\left[\mc(f_{j,k,3})(x)\right]^r\right\}^{\frac{1}{r}}
&\lesssim\left\{\sum_{j\in\mathbb{N}}
\left[\int_{\rrn}\frac{|f_{j,k,3}(y)|}{|x-y|^n}
\,dy\right]^r\right\}^{\frac{1}{r}}\notag\\
&\lesssim\int_{\rrn}\left\{\sum_{j\in\mathbb{N}}
\left[\frac{|f_{j,k,3}(y)|}{|x-y|^n}
\right]^r\right\}^{\frac{1}{r}}\,dy\notag\\
&\sim\int_{[B(\xi,\tka)]^\complement}
\frac{\vrys}{|x-y|^n}\,dy\notag\\
&\sim\sum_{i=k+2}^{\infty}\int_{B(\xi,2^i)
\setminus B(\xi,2^{i-1})}\frac{\vrys}{|x-y|^n}\,dy.
\end{align}
On the other hand, for any $k\in\mathbb{Z}$,
$i\in\mathbb{Z}\cap[k+2,\infty)$, and $x,\ y\in\rrn$
satisfying $\tkm\leq|x-\xi|<\tk$ and $2^{i-1}\leq|y-\xi|<2^i$,
we have
\begin{align*}
|x-y|\geq|y-\xi|-|x-\xi|
>2^{i-1}-\tk\geq2^{i-2},
\end{align*}
which, together with \eqref{wczool1e6}, implies that,
for any $k\in\mathbb{Z}$ and $x\in\rrn$
satisfying $\tkm\leq|x-\xi|<\tk$,
\begin{align}\label{wczool1e7}
&\left\{\sum_{j\in\mathbb{N}}
\left[\mc(f_{j,k,3})(x)\right]^r\right\}^{\frac{1}{r}}
\notag\\&\quad
\lesssim\sum_{i=k+2}^{\infty}2^{-ni}\left\|
\vr\1bf_{B(\xi,2^i)
\setminus B(\xi,2^{i-1})}\right\|_{L^1(\rrn)}.
\end{align}
In addition, from Lemma \ref{mono}
and \eqref{wczoott}, we deduce that, for any
$k\in\mathbb{Z}$ and $i\in\mathbb{Z}\cap[k+2,\infty)$,
\begin{align*}
\frac{\omega(\tk)}{\omega(2^i)}\sim
\frac{\omega(2^{k+2})}{\omega(2^i)}
\lesssim\left(\frac{2^{k+2}}{2^i}\right)^m
\sim2^{(k-i)m}.
\end{align*}
Applying this and \eqref{wczool1e7},
we conclude that,
for any $k\in\mathbb{Z}$,
\begin{align*}
&\omega(\tk)\left\|\mvrth\1bf_{B(\xi,\tk)
\setminus B(\xi,\tkm)}\right\|_{L^1(\rrn)}\\
&\quad\lesssim\omega(\tk)\sum_{i=k+2}^{
\infty}2^{-ni}\left\|
\1bf_{B(\xi,\tk)\setminus
B(\xi,\tkm)}\right\|_{L^1(\rrn)}\\
&\qquad\times\left\|\vr\1bf_{B(\xi,2^i)
\setminus B(\xi,2^{i-1})}\right\|_{L^1(\rrn)}\\
&\quad\sim\sum_{i=k+2}^{\infty}2^{(k-i)n}
\frac{\omega(\tk)}{\omega(2^i)}\omega(2^i)\left\|
\vr\1bf_{B(\xi,2^i)
\setminus B(\xi,2^{i-1})}\right\|_{L^1(\rrn)}\\
&\quad\lesssim\sum_{i=k+2}^{\infty}
2^{(k-i)(m+n)}\omega(2^i)\left\|
\vr\1bf_{B(\xi,2^i)
\setminus B(\xi,2^{i-1})}\right\|_{L^1(\rrn)}.
\end{align*}
From this, the fact that, for any $\lambda\in(0,\infty)$,
\begin{equation*}
\1bf_{\{y\in\rrn:\ \mvryth>\lambda/3\}}
\leq\frac{3}{\lambda}
\mvrth
\end{equation*}
and an argument similar to that
used in the estimations of both \eqref{wczool1e4}
and \eqref{wczool1e5}, it follows that,
for any $\lambda\in(0,\infty)$,
\begin{equation}\label{wczool1e10}
\mathrm{II}_{\lambda,3}\lesssim\lambda^{-1}
\left\|\vr(\cdot+\xi)\right\|_{\Kmp_{\omega,\0bf}^{1,q}(\rrn)}.
\end{equation}
This finishes the estimate of $\mathrm{II}_{\lambda,3}$.

Combining \eqref{wczool1e13}, \eqref{wczool1e4},
\eqref{wczool1e5}, \eqref{wczool1e12}, and \eqref{wczool1e10},
we further conclude that, for any $\lambda\in(0,\infty)$,
\begin{equation}\label{wczool1e14}
\lambda\left\|\mvr(\cdot+\xi)\right\|_{\Kmp_{\omega,\0bf}^{1,q}(\rrn)}
\lesssim\left\|\vr(\cdot+\xi)\right\|_{\Kmp_{\omega,\0bf}^{1,q}(\rrn)}.
\end{equation}
Thus, from the arbitrariness of $\xi$, it follows
that
$$
\left\|\mvr\right\|_{W\Kmp_{\omega}^{1,q}(\rrn)}\lesssim
\left\|\vr\right\|_{\Kmp_{\omega}^{1,q}(\rrn)}.
$$
In particular, letting $\xi:=\0bf$ in \eqref{wczool1e14}, we have
$$
\left\|\mvr\right\|_{W\Kmp_{\omega,\0bf}^{1,q}(\rrn)}\lesssim
\left\|\vr\right\|_{\Kmp_{\omega,\0bf}^{1,q}(\rrn)}.
$$
Therefore, both \eqref{wczool1e0} and \eqref{wczol1e00}
hold true in this case, and hence
we complete the proof of
Proposition \ref{wczool1}.
\end{proof}

In addition, to show Theorem
\ref{wczoo}, we also need
the following lemma about the
convexification of weak local
generalized Herz spaces, which is just a
simple corollary of both \cite[Remark 2.14]{ZWYY}
and Lemma \ref{convexll}; we omit the details.

\begin{lemma}\label{wconvexl}
Let $p,\ q,\ s\in(0,\infty)$ and $\omega\in M(\rp)$.
Then
$$\left[\WHerzSo\right]^{1/s}=W\HerzSocs$$
with the same quasi-norms.
\end{lemma}

Then we show Theorem \ref{wczoo}.

\begin{proof}[Proof of Theorem \ref{wczoo}]
Let all the symbols be as in the
present theorem. Then, combining
the assumption $\m0(\omega)
\in(-\frac{n}{p},\infty)$ and
Theorems \ref{Th3} and \ref{abso},
we find that the local
generalized Herz space $\HerzSo$ under
consideration is a BQBF space having
an absolutely continuous quasi-norm.
Thus, to complete the proof
of the present theorem, we only need to show
that all the assumptions of Proposition
\ref{wczxa} hold true for $\HerzSo$.

Indeed, let $$
s\in\left(0,\min\left\{
1,p,q,\frac{n}{\Mw+n/p}\right\}\right)
$$
and
\begin{equation}\label{wczooe1}
\theta\in\left(0,\min\left\{
s,\frac{n}{n+d}\right\}\right).
\end{equation}
Then, applying Lemma \ref{Atogl5},
we conclude that, for any $\{f_j\}
_{j\in\mathbb{N}}\subset L^1_{\mathrm{loc}}
(\rrn)$,
$$
\left\|\left\{\sum_{j\in\mathbb{N}}
\left[\mc^{(\theta)}\left(f_j\right)\right]^{s}
\right\}^{\frac{1}{s}}\right\|_{\HerzSo}
\lesssim\left\|\left(\sum_{j\in\mathbb{N}}
\left|f_j\right|^s
\right)^{\frac{1}{s}}\right\|_{\HerzSo}.
$$
This implies that Assumption \ref{assfs}
holds true with the above $\theta$ and $s$.

Next, we prove that Assumption \ref{assas}
holds true with the above $s$ and some $r_0\in
(1,\infty)$. Indeed,
let $$
r_0\in\left(\max\left\{
1,p,\frac{n}{\mw+n/p}\right\},\infty\right).
$$
Then, from Lemma \ref{mbhal}
with $r:=r_0$, we deduce that
$[\HerzSo]^{1/s}$ is a BBF space
and, for any $f\in L^1_{\mathrm{loc}}(\rrn)$,
$$
\left\|\mc^{((r_0/s)')}(f)
\right\|_{([\HerzSo]^{1/s})'}
\lesssim\left\|f\right\|_{([\HerzSo]^{1/s})'},
$$
which implies that Assumption
\ref{assas} holds true.

Let $$
r_1\in\left(0,\min\left\{p,
\frac{n}{\Mw+n/p}\right\}\right).
$$
Then, using Lemma \ref{wmbhl} with
$r:=r_1$, we find that,
for any $f\in L^1_{\mathrm{loc}}(\rrn)$,
$$
\left\|\mc(f)
\right\|_{[W\HerzSo]^{1/r_1}}
\lesssim\left\|f\right\|_{
[W\HerzSo]^{1/r_1}},
$$
which implies that, for the above
$r_1$, the Hardy--Littlewood
maximal operator $\mc$ is bounded
on $[W\HerzSo]^{1/r_1}$.

Finally, we prove that \eqref{wczxae0}
holds true with $X:=\HerzSo$.
Indeed, from the assumptions
$$\mw\in\left(-\frac{n}{p},\infty\right)$$
and $$\Mw\in\left(-\infty,n-\frac{n}{p}
+d+\delta\right),$$
and Lemma \ref{rela}, it follows that
\begin{align}\label{wczooe2}
&\min\left\{\m0\left(\omega^{\frac{n}{
n+d+\delta}}\right),\,\mi\left(\omega^{\frac{n}{
n+d+\delta}}\right)\right\}\notag\\
&\quad=\frac{n}{n+d+\delta}\mw>-\frac{n}{
(n+d+\delta)p/n}
\end{align}
and
\begin{align}\label{wczooe3}
&\max\left\{\M0\left(\omega^{\frac{n}{
n+d+\delta}}\right),\,\MI\left(\omega^{\frac{n}{
n+d+\delta}}\right)\right\}\notag\\
&\quad=\frac{n}{n+d+\delta}\Mw\notag\\
&\quad<n-\frac{n}{(n+d+\delta)p/n}=
\frac{n}{((n+d+\delta)/p)'}.
\end{align}
These, together with Lemma \ref{wconvexl}
with $s:=\frac{n}{n+d+\delta}$,
the assumption $p\in[\frac{n}{n+d+\delta},\infty)$,
\eqref{wczool1e0} with
$p$, $q$, and $\omega$ therein
replaced, respectively, by $\frac{n+d+\delta}{n}p$,
$\frac{n+d+\delta}{n}q$, and $\omega^{\frac{n}
{n+d+\delta}}$, and
Lemma \ref{convexll} with $s:=\frac{n}{n+d+\delta}$,
further imply that,
for any $\{f_j\}_{j\in\mathbb{N}}\subset
L^1_{\mathrm{loc}}(\rrn)$,
\begin{align*}
&\left\|\left\{\sum_{j\in\mathbb{N}}
\left[\mc\left(f_j\right)
\right]^{\frac{n+d+\delta}{n}}
\right\}^{\frac{n}{n+d+\delta}}
\right\|_{[W\HerzSo]^{\frac{n+d+\delta}{n}}}\\
&\quad=\left\|\left\{\sum_{j\in\mathbb{N}}
\left[\mc\left(f_j\right)
\right]^{\frac{n+d+\delta}{n}}
\right\}^{\frac{n}{n+d+\delta}}
\right\|_{W\Kmp_{\omega^{\frac{n}{n+d+\delta}},
\0bf}^{\frac{n+d+\delta}{n}p,
\frac{n+d+\delta}{n}q}(\rrn)}\\
&\quad\lesssim\left\|\left(\sum_{j\in\mathbb{N}}
\left|f_j\right|^{\frac{n+d+\delta}{n}}
\right)^{\frac{n}{n+d+\delta}}
\right\|_{\Kmp_{\omega^{\frac{n}{n+d+\delta}},
\0bf}^{\frac{n+d+\delta}{n}p,
\frac{n+d+\delta}{n}q}(\rrn)}\\
&\quad\sim\left\|\left(\sum_{j\in\mathbb{N}}
\left|f_j\right|^{\frac{n+d+\delta}{n}}
\right)^{\frac{n}{n+d+\delta}}
\right\|_{[\HerzSo]^{\frac{n+d+\delta}{n}}}.
\end{align*}
This then implies that \eqref{wczxae0}
holds true with $X:=\HerzSo$.
Combining this, the facts that
$\HerzSo$ satisfies both Assumptions
\ref{assfs} and \ref{assas} with
same $s$ and that $\mc$ is bounded on
$[W\HerzSo]^{1/r_1}$, we
find that the local generalized
Herz space $\HerzSo$ under consideration
satisfies all the assumptions of
Proposition \ref{wczxa}.
By this, \eqref{wczooe1}, and
Proposition \ref{wczxa} with $X:=\HerzSo$,
we conclude that $T$ has a unique
extension on $\HaSaHo$ and, for any $f\in
\HaSaHo$,
$$
\left\|T(f)\right\|_{W\HaSaHo}
\lesssim\left\|f\right\|_{\HaSaHo},
$$
which completes the proof of
Theorem \ref{wczoo}.
\end{proof}

Via Theorem \ref{wczoo}, Remark \ref{remark4.10}(ii),
and Remark \ref{remark6.16}, we immediately obtain the
following boundedness of
Calder\'{o}n--Zygmund operators from $H\MorrSo$
to $WH\MorrSo$; we omit the details.

\begin{corollary}\label{wczoom}
Let $d\in\zp$, $p,\ q\in[1,\infty)$,
$K$ be a $d$-order standard kernel as
in Definition \ref{def-s-k} with
some $\delta\in(0,1)$, $T$ a
$d$-order Calder\'{o}n--Zygmund operator
with kernel $K$ having the vanishing moments
up to order $d$, and $\omega\in
M(\rp)$ with
$$
-\frac{n}{p}<\m0(\omega)\leq
\M0(\omega)<0
$$
and
$$
-\frac{n}{p}<\mi(\omega)\leq
\MI(\omega)<0.
$$
Then $T$ has a unique extension
on $H\MorrSo$ and there exists
a positive constant $C$ such that,
for any $f\in H\MorrSo$,
$$
\|T(f)\|_{WH\MorrSo}\leq C\|f\|_{
H\MorrSo}.
$$
\end{corollary}

Finally, we turn to show the boundedness
of Calde\'{o}n--Zygmund operators
from the generalized Herz--Hardy space
$\HaSaH$ to the weak generalized
Herz--Hardy space $W\HaSaH$ as follows.

\begin{theorem}\label{wczoog}
Let $d\in\zp$, $\delta\in(0,1)$,
$p\in[\frac{n}{n+d+\delta},\infty)$,
$q\in(0,\infty)$,
$K$ be a $d$-order standard kernel as in
Definition \ref{def-s-k}, $T$ a
$d$-order Calder\'{o}n--Zygmund operator
with kernel $K$ having the vanishing moments
up to order $d$, and $\omega\in
M(\rp)$ with
$$
-\frac{n}{p}<\m0(\omega)\leq
\M0(\omega)<n-\frac{n}{p}+d+\delta
$$
and
$$
-\frac{n}{p}<\mi(\omega)\leq
\MI(\omega)<0.
$$
Then $T$ is well defined
on $\NHaSaH$ and there exists
a positive constant $C$ such that,
for any $f\in\HaSaH$,
$$
\|T(f)\|_{W\HaSaH}\leq C\|f\|_{
\HaSaH}.
$$
\end{theorem}

To show this theorem, we
need the following
auxiliary conclusion about
the convexification of
weak global generalized Herz spaces,
which is just a simple consequence of both
\cite[Remark 2.14]{ZWYY} and Lemma \ref{convexl};
we omit the details.

\begin{lemma}\label{wconvex}
Let $p,\ q,\ s\in(0,\infty)$ and $\omega\in M(\rp)$.
Then
$$\left[\WHerzS\right]^{1/s}=W\HerzScs$$
with the same quasi-norms.
\end{lemma}

Via the above lemma, we now
prove Theorem \ref{wczoog}.

\begin{proof}[Proof of Theorem \ref{wczoog}]
Let all the symbols be as in
the present theorem. Then, from
the assumptions $\m0(\omega)\in
(-\frac{n}{p},\infty)$ and
$\MI(\omega)\in(-\infty,0)$,
and Theorem \ref{Th2}, we infer
that the global generalized
Herz space $\HerzS$ under consideration
is a BQBF space. Therefore,
in order to finish the
proof of the present theorem,
it suffices to prove that
the assumptions (i) through (vii)
hold true for $\HerzS$.
For this purpose, let
$$
s_0\in\left(0,\min\left\{1,p,q,\frac{n}{
\Mw+n/p}\right\}\right),
$$
$s\in(0,s_0)$,
\begin{equation}\label{wczooge1}
\theta\in\left(0,\min\left\{s,\frac{n}{n+d}
\right\}\right),
\end{equation}
and $\eta\in(1,\infty)$ satisfy
$$
\eta<\min\left\{\frac{n}{n(1-s/p)-s\mw},\left(
\frac{p}{s}\right)'\right\}.
$$
Then, repeating the proof of Theorem \ref{bddog},
we find that the
following five statements hold true:
\begin{enumerate}
  \item[{\rm(i)}] for any $\{f_j\}_{j
  \in\mathbb{N}}\subset L^1_{\mathrm{loc}}(\rrn)$,
$$
\left\|\left\{
\sum_{j\in\mathbb{N}}
\left[\mc^{(\theta)}(f_j)\right]^{s}
\right\}^{1/s}
\right\|_{\HerzS}
\lesssim\left\|
\left(\sum_{j\in\mathbb{N}}
|f_j|^{s}
\right)^{1/s}
\right\|_{\HerzS};
$$
  \item[{\rm(ii)}] for any $f\in\Msc(\rrn)$,
  $$
  \|f\|_{[\HerzS]^{1/s}}
  \sim\sup\left\{\|fg\|_{L^1(\rrn)}:\
  \|g\|_{\HerzScsd}=1\right\}
  $$
  and
  $$
  \|f\|_{[\HerzS]^{1/s_0}}
  \sim\sup\left\{\|fg\|_{L^1(\rrn)}:\
  \|g\|_{\HerzScsod}=1\right\}
  $$
  with the positive equivalence constants
  independent of $f$;
  \item[{\rm(iii)}] both $\|\cdot\|_{\HerzScsd}$
  and $\|\cdot\|_{\HerzScsod}$ satisfy
  Definition \ref{Df1}(ii);
  \item[{\rm(iv)}] $\1bf_{B(\0bf,1)}\in\HerzScsod$;
  \item[{\rm(v)}] $\mc^{(\eta)}$ is bounded
  on both $\HerzScsd$ and $\HerzScsod$.
\end{enumerate}
These imply that the assumptions (i)
through (v) of Theorem \ref{wczx} hold
true for $\HerzS$.

Next, we show that Theorem
\ref{wczx}(vi) holds true.
Namely, there exists an $r_1\in(0,\infty)$
such that $\mc$ is bounded on $[
W\HerzS]^{1/r_1}$. Indeed,
let $$
r_1\in\left(0,\min\left\{p,
\frac{n}{\Mw+n/p}\right\}\right).
$$
From this and Lemma \ref{wmbhg} with
$r:=r_1$, it follows that,
for any $f\in L^1_{\mathrm{loc}}(\rrn)$,
$$
\left\|\mc\left(f
\right)\right\|_{[W\HerzS]^{1/r_1}}
\lesssim\left\|f\right\|_{[W\HerzS]^{1/r_1}},
$$
which implies that Theorem \ref{wczx}(iv)
holds true for $\HerzS$ with the above $r_1$.

Finally, we show that the assumption
(vii) of Theorem \ref{wczx} holds true.
Indeed, applying Lemma \ref{wconvex}
with $s:=\frac{n}{n+d+\delta}$,
\eqref{wczooe2}, \eqref{wczooe3},
the assumption $p\in[\frac{n}{n+d+\delta},\infty)$,
\eqref{wczol1e00} with
$p$, $q$, and $\omega$ therein
replaced, respectively, by $\frac{n+d+\delta}{n}p$,
$\frac{n+d+\delta}{n}q$, and $\omega^{\frac{n}
{n+d+\delta}}$, and
Lemma \ref{convexl} with $s:=\frac{n}{n+d+\delta}$,
we find that,
for any $\{f_j\}_{j\in\mathbb{N}}\subset
L^1_{\mathrm{loc}}(\rrn)$,
\begin{align*}
&\left\|\left\{\sum_{j\in\mathbb{N}}
\left[\mc\left(f_j\right)
\right]^{\frac{n+d+\delta}{n}}
\right\}^{\frac{n}{n+d+\delta}}
\right\|_{[W\HerzS]^{\frac{n+d+\delta}{n}}}\\
&\quad=\left\|\left\{\sum_{j\in\mathbb{N}}
\left[\mc\left(f_j\right)
\right]^{\frac{n+d+\delta}{n}}
\right\}^{\frac{n}{n+d+\delta}}
\right\|_{W\Kmp_{\omega^{\frac{n}
{n+d+\delta}}}^{\frac{n+d+\delta}{n}p,
\frac{n+d+\delta}{n}q}(\rrn)}\\
&\quad\lesssim\left\|\left(\sum_{j\in\mathbb{N}}
\left|f_j\right|^{\frac{n+d+\delta}{n}}
\right)^{\frac{n}{n+d+\delta}}
\right\|_{\Kmp_{\omega^{\frac{n}
{n+d+\delta}}}^{\frac{n+d+\delta}{n}p,
\frac{n+d+\delta}{n}q}(\rrn)}\\
&\quad\sim\left\|\left(\sum_{j\in\mathbb{N}}
\left|f_j\right|^{\frac{n+d+\delta}{n}}
\right)^{\frac{n}{n+d+\delta}}
\right\|_{[\HerzS]^{\frac{n+d+\delta}{n}}}.
\end{align*}
This further implies that Theorem
\ref{wczx}(vii) holds true with
$X:=\HerzS$.
Thus, the assumptions (i) through
(vii) of Theorem \ref{wczx} hold
true for the global generalized
Herz space $\HerzS$ under consideration.
Combining this, the fact that
$\HerzS$ is a BQBF space, \eqref{wczooge1},
and Theorem \ref{wczx} with $X:=\HerzS$,
we conclude that $T$ is well defined
on $\HaSaH$ and, for any $f\in\HaSaH$,
$$
\left\|T\left(
f\right)\right\|_{W\HaSaH}
\lesssim\left\|f\right\|_{\HaSaH}.
$$
This then finishes the proof of Theorem
\ref{wczoog}.
\end{proof}

Via Theorem \ref{wczoog},
Remark \ref{remark4.10}(ii), and Remark \ref{remark6.16},
we immediately obtain the following
boundedness of Calder\'{o}n--Zygmund operators
from the generalized Morrey--Hardy space $H\MorrS$
to the weak generalized Morrey--Hardy space $WH\MorrS$;
we omit the details.

\begin{corollary}
Let $d$, $p$, $q$, $\omega$, $K$, and $T$
be as in Corollary \ref{wczoom}.
Then $T$ is well defined
on $H\MorrS$ and there exists
a positive constant $C$ such that,
for any $f\in H\MorrS$,
$$
\|T(f)\|_{WH\MorrS}\leq C\|f\|_{
H\MorrS}.
$$
\end{corollary}

\section{Real Interpolations}

As applications,
we investigate the real interpolation
theorems about generalized Herz--Hardy
spaces in this section.
Precisely,
we show that the real interpolation
spaces between the generalized Herz--Hardy spaces
and the Lebesgue space $L^\infty(\rrn)$ are just the
weak generalized Herz--Hardy spaces.

Let $X_0$ and $X_1$ be two quasi-Banach spaces.
Recall that the pair $(X_{0},\,X_{1})$ is said to be
\emph{compatible}\index{compatible}
(see, for instance, \cite[p.\,20]{BL})
if there exists a Hausdorff space $\mathbb{X}$ such that
$X_{0}\subset\mathbb{X}$ and
$X_{1}\subset\mathbb{X}$.
For any compatible pair $(X_{0},\,X_{1})$
of quasi-Banach spaces, let
$$\index{$X_0+X_1$}
X_{0}+X_{1}:=\{a\in\mathbb{X}:\ \exists\,a_{0}\in
X_{0},\ \exists\,
a_{1}\in X_{1}\text{ such that }a=a_{0}+a_{1}\}.
$$
Moreover, for any $t\in(0,\infty)$, the
\emph{Peetre K-functional}\index{Peetre K-functional}
$K(t,f;X_{0},X_{1})$\index{$K(t,f;X_{0},X_{1})$}
(see, for instance, \cite[p.\,38]{BL})
is defined by setting, for any $f\in X_{0}+X_{1}$,
$$
K(t,f;X_{0},X_{1}):=\inf\left\{
\|f_{0}\|_{X_{0}}+t\|f_{1}\|_{X_{1}}:\
f=f_{0}+f_{1},\ f_{i}\in X_{i},\ i
\in\{0,1\}\right\}.
$$
We next recall the definition of the real interpolation
space between two quasi-Banach spaces
$X_0$ and $X_1$ as follows
(see, for instance, \cite[p.\,40]{BL}).

\begin{definition}\label{df681}
Let $X_0$, $X_1$ be two quasi-Banach
spaces satisfy that
the pair $(X_0,\,X_1)$ is compatible,
$\theta\in(0,1)$, and $q\in(0,\infty]$.
Then the\emph{
real interpolation space}\index{real interpolation space}
$(X_{0},\,X_{1})_{\theta,q}$\index{$(X_{0},\,X_{1})_{\theta,q}$}
between $X_{0}$ and $X_{1}$ is defined to be
the set of all the $f\in X_0+X_1$ such that
$$
\|f\|_{\theta,q}<\infty,
$$
where, for any $f\in X_{0}+X _{1}$,
$$
\|f\|_{\theta,q}:=\left\{
\begin{aligned}
&\left\{\int_{0}^{\infty}\left[
t^{-\theta}K(t,f;X_{0},X_{1})\right]^{q}
\,\frac{dt}{t}\right\}^{\frac{1}{q}},
\hspace{.3cm}\text{when }q\in(0,\infty),\\
&\sup_{t\in(0,\infty)}t^{-\theta}K(t,f;X_{0},X_{1}),
\hspace{2cm}\text{when }q=\infty.
\end{aligned}
\right.
$$
\end{definition}

Then the following conclusion shows that
the real interpolation space between
$\HaSaHo$ and $L^\infty(\rrn)$ is just
$WH\Kmp_{\omega^{1-\theta},\0bf}^
{p/(1-\theta),q/(1-\theta)}(\rrn)$.

\begin{theorem}\label{interl}
Let $p,\ q\in(0,\infty)$ and $\omega
\in M(\rp)$ satisfy
$\m0(\omega)\in(-\frac{n}{p},\infty)$
and $\mi(\omega)
\in(-\frac{n}{p},\infty)$.
Then, for any $\theta\in(0,1)$,
$$
\left(\HaSaHo,\,L^{\infty}(\rrn)
\right)_{\theta,\infty}=
WH\Kmp_{\omega^{1-\theta},\0bf}
^{p/(1-\theta),q/(1-\theta)}(\rrn).
$$
\end{theorem}

To show this real interpolation
theorem, we need the following
real interpolation result between
Hardy spaces associated with
ball quasi-Banach function spaces
and the Lebesgue space $L^{\infty}(\rrn)$,
which is just
\cite[Theorem 4.5]{WYYZ}.

\begin{lemma}\label{interx}
Let $X$ be a ball quasi-Banach function space and
let $p_-\in(0,1)$, $\theta_0\in(1,\infty)$,
and $r\in(0,\infty)$ be such that the following
four statements hold true:
\begin{enumerate}
  \item[{\rm(i)}] for any given
  $\theta\in(0,p_-)$ and $u\in(1,\infty)$,
  there exists a positive constant $C$ such that,
for any $\{f
_{j}\}_{j\in\mathbb{N}}\subset L_{{\rm loc}}^{1}(\rrn)$,
\begin{equation*}
\left\|\left\{\sum_{j\in\mathbb{N}}\left[\mc
(f_{j})\right]^{u}\right\}^
{\frac{1}{u}}\right\|_{X^{1/\theta}}
\leq C\left\|\left\{
\sum_{j\in\mathbb{N}}|f_{j}|^{u}\right\}
^{\frac{1}{u}}\right\|_{X^{1/\theta}};
\end{equation*}
  \item[{\rm(ii)}] $X$ is $\theta_0$-concave;
  \item[{\rm(iii)}] there exist an $s_0\in(0,p_-)$,
  an $r_0\in(s_0,\infty)$, and a $C\in(0,\infty)$
  such that $X^{1/s_0}$ is a ball Banach function
  space and, for any $f\in L^1_{\mathrm{loc}}(\rrn)$,
  $$
  \left\|\mc^{((r_0/s_0)')}(f)\right\|_{(X^{1/s_0})'}
  \leq C\left\|f\right\|_{(X^{1/s_0})'};
  $$
  \item[{\rm(iv)}] $\mc$ is bounded on $(WX)^{1/r}$.
\end{enumerate}
Assume that $\theta\in(0,1)$.
Then
$$
\left(H_X(\rrn),\,L^{\infty}(\rrn)
\right)_{\theta,\infty}
=WH_{X^{1/(1-\theta)}}(\rrn).
$$
\end{lemma}

We now show Theorem \ref{interl}.

\begin{proof}[Proof of Theorem \ref{interl}]
Let all the symbols be as in the present
theorem. Then, by the assumption
$\m0(\omega)\in(-\frac{n}{p},\infty)$
and Theorem \ref{Th3},
we conclude that the local generalized
Herz space $\HerzSo$ under consideration
is a BQBF space.
Thus, in order to finish the proof
of the present theorem,
we only need to show that $\HerzSo$ satisfies
the assumptions (i) through (iv)
of Lemma \ref{interx}.

Indeed, let
$$
p_-:=\min\left\{1,p,q,\frac{n}{\Mw+n/p}\right\}.
$$
Then, for any given $\theta\in(0,p_-)$
and $u\in(1,\infty)$,
from Lemma \ref{vmbhl} with $r:=\theta$,
it follows that, for any
$\{f_{j}\}_{j\in\mathbb{N}}\subset
L_{{\rm loc}}^{1}(\rrn)$,
\begin{equation}\label{interle1}
\left\|\left\{\sum_{j\in\mathbb{N}}\left[\mc
(f_{j})\right]^{u}\right\}^{\frac{1}{u}}
\right\|_{[\HerzSo]^{1/\theta}}
\lesssim\left\|\left(
\sum_{j\in\mathbb{N}}|f_{j}|^{u}\right)
^{\frac{1}{u}}\right\|_{[\HerzSo]^{1/\theta}}.
\end{equation}
This implies that Lemma \ref{interx}(i)
holds true with the above $p_-$.

Next, we prove that Lemma \ref{interx}(ii)
holds true. To this end, let
$$\theta_0\in\left(\max\left\{1,p,q
\right\},\infty\right).$$
Then, applying the reverse Minkowski
inequality and Lemma \ref{convexll} with
$s:=\theta_0$,
we conclude that, for any $\{f_{j}\}_{j\in
\mathbb{N}}\subset\Msc(\rrn)$,
\begin{equation}\label{interle2}
\sum_{j\in\mathbb{N}}\|f_{j}\|_{[\HerzSo]
^{1/\theta_0}}\leq
\left\|\sum_{j\in\mathbb{N}}|f_{j}|
\right\|_{[\HerzSo]^{1/\theta_0}},
\end{equation}
which then implies that
the Herz space $\HerzSo$ under consideration
is $\theta_0$-concave and hence
Lemma \ref{interx}(ii) holds
true.

In addition, let $s_0\in(0,p_-)$ and
$$
r_0\in\left(\min\left\{1,p,
\frac{n}{\mw+n/p}\right\},\infty\right).
$$
We now prove that, for this $s_0$ and
$r_0$, the assumption (iii) of Lemma \ref{interx}
holds true.
Indeed, by Lemma \ref{mbhal} with
$s:=s_0$ and $r:=r_0$, we find that
$[\HerzSo]^{1/s_0}$ is a BBF space and,
for any $f\in L_{{\rm loc}}^{1}(\rrn)$,
\begin{equation}\label{interle3}
\left\|\mc^{((r_0/s_0)')}(f)
\right\|_{([\HerzSo]^{1/s_0})'}\lesssim
\left\|f\right\|_{([\HerzSo]^{1/s_0})'}.
\end{equation}
This finishes the proof of Lemma \ref{interx}(iii).

Finally, we show that Lemma \ref{interx}(iv)
holds true. Namely, there exists
an $r\in(0,\infty)$ such that
$\mc$ is bounded on $[W\HerzSo]^{1/r}$.
Indeed, let
$$
r\in\left(
0,\min\left\{p,\frac{n}
{\Mw+n/p}\right\}\right).
$$
Using this and Lemma \ref{wmbhl}, we conclude
that, for any $f\in L^1_\mathrm{loc}(\rrn)$,
\begin{equation*}
\left\|\mc(f)\right\|_{[W\HerzSo]^{1/r}}
\lesssim\|f\|_{[W\HerzSo]^{1/r}},
\end{equation*}
which implies that Lemma \ref{interx}(iv)
holds true for $\HerzSo$.
Combining this, \eqref{interle1},
\eqref{interle2}, and \eqref{interle3},
we further find that
the local generalized Herz space
$\HerzSo$ under consideration satisfies
the assumptions (i) through (iv)
of Lemma \ref{interx}. This further implies
that, for any $\theta\in(0,1)$,
$$
\left(\HaSaHo,\,L^{\infty}(\rrn)
\right)_{\theta,\infty}=
WH\Kmp_{\omega^{1-\theta},\0bf}
^{p/(1-\theta),q/(1-\theta)}(\rrn),
$$
which completes the proof of Theorem
\ref{interl}.
\end{proof}

Via Theorem \ref{interl}, Remark \ref{remark4.10}(ii),
and Remark \ref{remark6.16}, we obtain the following
real interpolation theorem
which shows that the real interpolation space between
$H\MorrSo$ and $L^\infty(\rrn)$ is just the weak space
$WH\textsl{\textbf{M}}_{\omega^{1-\theta},
\0bf}^{p/(1-\theta),q/(1-\theta)}(\rrn)$; we omit the details.

\begin{corollary}\label{interlm}
Let $p,\ q\in[1,\infty)$ and $\omega\in M(\rp)$ with
$$
-\frac{n}{p}<\m0(\omega)\leq\M0(\omega)<0
$$
and
$$
-\frac{n}{p}<\mi(\omega)\leq\MI(\omega)<0.
$$
Then, for any $\theta\in(0,1)$,
$$
\left(H\MorrSo,\,L^{\infty}(\rrn)\right)_{\theta,\infty}=
WH\textsl{\textbf{M}}_{\omega^
{1-\theta},\0bf}^{p/(1-\theta),q/(1-\theta)}(\rrn).
$$
\end{corollary}

The remainder of this section is
devoted to establishing
a real interpolation theorem
between the generalized Herz--Hardy
space $\HaSaH$ and the Lebesgue space
$L^{\infty}(\rrn)$. Indeed,
we have the following conclusion which
shows that the real interpolation space
between $\HaSaH$ and $L^{\infty}(\rrn)$
is just the weak generalized Herz--Hardy space
$WH\Kmp_{\omega^{1-\theta}}
^{p/(1-\theta),q/(1-\theta)}(\rrn)$.

\begin{theorem}\label{interg}
Let $p,\ q\in(0,\infty)$ and $\omega\in M(\rp)$ satisfy
$\m0(\omega)\in(-\frac{n}{p},\infty)$ and
$$
-\frac{n}{p}<\mi(\omega)\leq\MI(\omega)<0.
$$
Then, for any $\theta\in(0,1)$,
$$
\left(\HaSaH,\,L^{\infty}(\rrn)\right)_{\theta,\infty}=
WH\Kmp_{\omega^{1-\theta}}^{p/(1-\theta),q/(1-\theta)}(\rrn).
$$
\end{theorem}

To show this theorem, we first establish
the following real interpolation
theorem about the Hardy spaces
associated with ball quasi-Banach function
spaces without recourse to
associate spaces, which
improves the results obtained in \cite{WYYZ}
(see also Lemma \ref{interx} above).

\begin{proposition}\label{intergx}
Let $X$ be a ball quasi-Banach function space and
$Y\subset\Msc(\rrn)$ a linear space equipped
with a quasi-seminorm $\|\cdot\|_Y$,
and let $p_-\in(0,1)$, $\theta_0\in(1,\infty)$,
and $r\in(0,\infty)$ be such that the following
eight statements hold true:
\begin{enumerate}
  \item[{\rm(i)}] for any given
  $\theta\in(0,p_-)$ and $u\in(1,\infty)$,
  there exists a positive constant $C$ such that,
for any $\{f
_{j}\}_{j\in\mathbb{N}}\subset L_{{\rm loc}}^{1}(\rrn)$,
\begin{equation*}
\left\|\left\{\sum_{j\in\mathbb{N}}\left[\mc
(f_{j})\right]^{u}\right\}^
{\frac{1}{u}}\right\|_{X^{1/\theta}}
\leq C\left\|\left\{
\sum_{j\in\mathbb{N}}|f_{j}|^{u}\right\}
^{\frac{1}{u}}\right\|_{X^{1/\theta}};
\end{equation*}
  \item[{\rm(ii)}] $X$ is $\theta_0$-concave;
  \item[{\rm(iii)}] $\mc$ is bounded on $(WX)^{1/r}$;
  \item[{\rm(iv)}] for any $s\in(0,p_-)$, $X^{1/s}$
  is a ball Banach function space;
  \item[{\rm(v)}] $\|\cdot\|_{Y}$ satisfies
  Definition \ref{Df1}(ii);
  \item[{\rm(vi)}] $\1bf_{B(\0bf,1)}\in Y$;
  \item[{\rm(vii)}] there exists an $s_0\in(0,p_-)$
  such that, for any $f\in\Msc(\rrn)$,
  $$
  \left\|f\right\|_{X^{1/s_0}}
  \sim\sup\left\{\left\|fg\right\|_{L^1(\rrn)}:\
  \left\|g\right\|_Y=1\right\},
  $$
  where the positive equivalence
  constants are independent of $f$;
  \item[{\rm(viii)}] there exist
  an $r_0\in(s_0,\infty)$ and a $C\in(0,\infty)$
  such that, for any $f\in L^1_{\mathrm{loc}}(\rrn)$,
  $$
  \left\|\mc^{((r_0/s_0)')}(f)\right\|_{Y}
  \leq C\left\|f\right\|_{Y}.
  $$
\end{enumerate}
Then, for any $\theta\in(0,1)$,
$$
\left(H_X(\rrn),\,L^{\infty}(\rrn)
\right)_{\theta,\infty}
=WH_{X^{1/(1-\theta)}}(\rrn).
$$
\end{proposition}

\begin{remark}
We should point out that
Proposition \ref{intergx} is an
improved version of the known real
interpolation theorem between
the Hardy space $H_X(\rrn)$
associated with the ball quasi-Banach
function space $X$ and the
Lebesgue space $L^{\infty}(\rrn)$
established in \cite[Theorem 4.5]{WYYZ}.
Indeed, if $Y\equiv(X^{1/s_0})'$ in
Proposition \ref{intergx}, then
this proposition goes back to
\cite[Theorem 4.5]{WYYZ}.
\end{remark}

In order to prove this proposition,
we need the following Calder\'{o}n reproducing
formula\index{Calder\'{o}n reproducing formula}
which was obtained by
Calder\'{o}n
\cite[p.\,219,\ Section 3]{c77}
(see also \cite[Lemma 4.1]{ct75} or
\cite[Lemma 4.6]{yyyz16}).

\begin{lemma}\label{crff}
Let $\psi\in\mathcal{S}(\rrn)$
satisfy that $\supp(\psi)\subset B(\0bf,1)$
and $\int_{\rrn}\psi(x)\,dx=0$.
Then there exist a $\phi\in\mathcal{S}(\rrn)$
and $0<a<b<\infty$ such that
$$
\supp(\widehat{\phi})
\subset B(\0bf,b)\setminus B(\0bf,a)
$$
and, for any $x\in\rrn\setminus\{\0bf\}$,
$$
\int_{0}^{\infty}\widehat{\psi}(tx)
\widehat{\phi}(tx)\,\frac{dt}{t}=1.
$$
Moreover, let $\eta$ be defined by setting
$\widehat{\eta}(\0bf):=1$
and, for any $x\in\rrn\setminus\{\0bf\}$,
$$
\widehat{\eta}(x):=\int_{1}^{\infty}
\widehat{\psi}(tx)\widehat{\phi}(tx)\,
\frac{dt}{t}.
$$
Then
\begin{enumerate}
  \item[{\rm(i)}] $\eta$ is well defined;
  \item[{\rm(ii)}] $\eta\in C_{\mathrm{c}}^{\infty}
  (\rrn)$;
  \item[{\rm(iii)}] there exists a $\delta\in(0,\infty)$
  such that, for any $x\in B(\0bf,\delta)$,
  $\widehat{\eta}(x)=1$.
\end{enumerate}
\end{lemma}

We also require the following
technique lemma which is just
\cite[Theorem 1.64]{fs82}.

\begin{lemma}\label{crf}
Let $\phi\in\mathcal{S}(\rrn)$
satisfy $\int_{\rrn}\phi(x)\,dx=0$
and, for any $x\in\rrn\setminus\{\0bf\}$,
$$
\int_{\rrn}\widehat{\phi}(tx)\,\frac{dt}{t}=1.
$$
Assume that $f\in\mathcal{S}'(\rrn)$ vanishes
weakly at infinity. Then
$$
f=\int_{0}^{\infty}f\ast\phi_t\,\frac{dt}{t}
$$
in $\mathcal{S}'(\rrn)$, namely,
$$
f=\lim\limits_{\eps\to0^+,\,A\to\infty}\int_{\eps}^{A}
f\ast\phi_t
\,\frac{dt}{t}
$$
in $\mathcal{S}'(\rrn)$.
\end{lemma}

Via both Lemmas \ref{crff} and \ref{crf},
we next prove Proposition \ref{intergx}.

\begin{proof}[Proof of Proposition
\ref{intergx}]
Let all the symbols be as in the present
proposition and $\theta\in(0,1)$.
Then, by the proof of
\cite[Theorem 4.5]{WYYZ}, we find that
\begin{equation}\label{intergxe3}
\left(H_X(\rrn),\,L^{\infty}(\rrn)
\right)_{\theta,\infty}\subset
WH_{X^{1/(1-\theta)}}(\rrn).
\end{equation}
Conversely, we next show that
\begin{equation}\label{intergxe2}
WH_{X^{1/(1-\theta)}}(\rrn)\subset
\left(H_X(\rrn),\,L^{\infty}(\rrn)
\right)_{\theta,\infty}.
\end{equation}
To this end, let
$f\in WH_{X^{1/(1-\theta)}}(\rrn)$.
Then, from the assumptions
(v) through (viii) of the present proposition
and Lemma \ref{WX-weak vanish}
with $X$, $\theta$, and $s$ therein
replaced, respectively, by $X^{1/(1-\theta)}$,
$(r_0/s_0)'$, and $s_0/(1-\theta)$,
it follows that $f$ vanishes weakly
at infinity. In addition, let
$d\geq\lfloor n(1/\min\{\frac{p_-}{\theta_0},
r_0\}-1)\rfloor$ be a
fixed nonnegative integer and
$\psi\in\mathcal{S}(\rrn)$ satisfy that
$\supp(\widehat{\psi})\subset B(\0bf,1)$
and, for any $\gamma\in\zp^n$ with
$|\gamma|\leq d$,
$$\int_{\rrn}\psi(x)x^{\gamma}\,dx=0.$$
Then, using Lemma \ref{crff}, we
find that there exist a $\phi\in\mathcal{S}(\rrn)$
and $0<a<b<\infty$ such that
$$\supp(\widehat{\phi})\subset
B(\0bf,b)\setminus B(\0bf,a)$$
and, for any $x\in\rrn\setminus\{\0bf\}$,
$$
\int_{0}^{\infty}\widehat{\psi}(tx)
\widehat{\phi}(tx)\,\frac{dt}{t}=1.
$$
These, combined with Lemma \ref{crf}
with $\phi$ therein replaced
by $\psi\ast\phi$ and the fact
that $f$ vanishes weakly at infinity,
further imply that
\begin{equation}\label{intergxe1}
f=\int_{0}^{\infty}f\ast\psi_t\ast\phi_t\,
\frac{dt}{t}
\end{equation}
in $\mathcal{S}'(\rrn)$.

Let $x_0:=(2,\ldots,2)\in\rrn$ and $\eta$
be the same as in Lemma \ref{crff}, namely,
$\widehat{\eta}(\0bf):=1$ and,
for any $x\in\rrn\setminus\{\0bf\}$,
$$
\widehat{\eta}(x):=\int_{1}^{\infty}
\widehat{\psi}(tx)\widehat{\phi}(tx)\,
\frac{dt}{t}.
$$
By (i), (ii), and (iii) of
Lemma \ref{crff}, we
find that $\eta\in\mathcal{S}(\rrn)$.
For any $x\in\rrn$ and $t\in(0,\infty)$,
let
$$\widetilde{\phi}(x):=\phi(x-x_0),\
\widetilde{\psi}(x):=\psi(x+x_0),$$
$$F(x,t):=f\ast\widetilde{\phi}_t(x),\
G(x,t):=f\ast\eta_t(x),$$ and
$$
M_{\bigtriangledown}(f)(x):=
\sup_{(y,t)\in\Gamma_{3(|x_0|+1)}(x)}
\left[\left|F(y,t)\right|+
\left|G(y,t)\right|\right],
$$
where, for any $x\in\rrn$,
$\Gamma_{3(|x_0|+1)}(x)$ is defined as in
\eqref{gammaa} with $a:=3(|x_0|+1)$.
Then, applying \eqref{intergxe1},
we conclude that
$$
f(\cdot)=\int_{0}^{\infty}f\ast\widetilde{\phi}_t
\ast\widetilde{\psi}_t(\cdot)\,\frac{dt}{t}
=\int_{0}^{\infty}\int_{\rrn}
F(y,t)\widetilde{\psi}_t(\cdot-y)\,\frac{dy\,dt}{t}
$$
in $\mathcal{S}'(\rrn)$.
From this and the proof of \cite[Lemma 4.4]{WYYZ},
we deduce that, for any
$\alpha\in(0,\infty)$, there exist a $g_{\alpha}
\in L^{\infty}(\rrn)$ and
a $b_{\alpha}\in\mathcal{S}'(\rrn)$
such that $f=g_{\alpha}+b_{\alpha}$
in $\mathcal{S}'(\rrn)$, $\|g_{\alpha}\|_{L^{\infty}}
\lesssim\alpha$, and
$$
\left\|b_{\alpha}\right\|_{H_X(\rrn)}
\lesssim\left\|M_{\bigtriangledown}(f)
\1bf_{\{x\in\rrn:\
M_{\bigtriangledown}(f)(x)>\alpha\}}\right\|_{X},
$$
where the implicit positive constants
are independent of both $f$ and $\alpha$.
This, together with the proof of
\cite[Theorem 4.5]{WYYZ},
further implies
that, for any $t\in(0,\infty)$,
\begin{align*}
&t^{-\theta}K\left(t,f;H_X(\rrn),L^{\infty}
(\rrn)\right)\\&\quad\lesssim\left\|
M_\bigtriangledown
(f)\right\|_{X^{1/(1-\theta)}}
\lesssim\left\|M_{\bigtriangledown}
(f)\right\|_{WX^{1/(1-\theta)}}
\sim\left\|f
\right\|_{WH_{X^{1/(1-\theta)}}(\rrn)}.
\end{align*}
Using this and Definition \ref{df681},
we find that
$$
\left\|f\right\|_{\theta,\infty}
\lesssim\left\|f
\right\|_{WH_{X^{1/(1-\theta)}}(\rrn)}<\infty.
$$
This further implies that
$f\in (H_X(\rrn),\,L^{\infty}(\rrn))_{\theta,\infty}$
and then finishes the proof of
\eqref{intergxe2}.
Combining this and \eqref{intergxe3},
we conclude that
$$
\left(H_X(\rrn),\,L^{\infty}(\rrn)\right)_{\theta,\infty}
=WH_{X^{1/(1-\theta)}}(\rrn)
$$
and hence complete the proof of Proposition
\ref{intergx}.
\end{proof}

Via the improved real interpolation
result obtained above, we now prove Theorem
\ref{interg}.

\begin{proof}[Theorem \ref{interg}]
Let all the symbols be as in the present
theorem. Then, by the assumptions
$\m0(\omega)\in(-\frac{n}{p},\infty)$ and
$\MI(\omega)\in(-\infty,0)$, and
Theorem \ref{Th2}, we find that
the global generalized Herz space
$\HerzS$ under consideration is a BBF space.
Therefore, in order to finish
the proof of the present theorem,
we only need to show
that the assumptions (i) through
(viii) of Proposition \ref{intergx} hold
true for $\HerzS$.

Indeed, let
$$
p_-:=\left\{1,p,q,\frac{n}{\Mw+n/p}\right\}.
$$
Then, for any given $\theta\in(0,p_-)$ and
$u\in(1,\infty)$, from Lemma \ref{vmbhl}
with $r:=\theta$, it follows that,
for any $\{f_{j}\}_{j\in\mathbb{N}}\subset
L^1_{\mathrm{loc}}(\rrn)$,
$$
\left\|\left\{\sum_{j\in\mathbb{N}}
\left[\mc\left(f_j\right)\right]^u
\right\}^{\frac{1}{u}}\right\|_{[\HerzS]^{1/\theta}}
\lesssim\left\|\left(\sum_{j\in\mathbb{N}}
\left|f_j\right|^u
\right)^{\frac{1}{u}}\right\|_{[\HerzS]^{1/\theta}}.
$$
This implies that Proposition \ref{intergx}(i)
holds true with the above $p_-$.

Next, we prove that Proposition \ref{intergx}(ii)
holds true with some $\theta_0\in(1,\infty)$.
Indeed, let $\theta_0\in(\max\{1,p,q\},\infty)$.
Combining this, the reverse Minkowski inequality,
and Lemma \ref{convexl} with $s:=\theta_0$,
we conclude that, for any
$\{f_{j}\}_{j\in\mathbb{N}}\subset
\Msc(\rrn)$,
$$
\sum_{j\in\mathbb{N}}\left\|f_j\right\|_{
[\HerzS]^{1/\theta_0}}
\leq\left\|\sum_{j\in\mathbb{N}}
\left|f_j\right|\right\|_{[\HerzS]^{1/\theta_0}},
$$
which further implies that the Herz space
$\HerzS$ under consideration is $\theta_0$-concave
and hence Proposition \ref{intergx}(ii)
holds true.

Let
$$
r\in\left(0,\min\left\{p,\frac{n}{\Mw+n/p}
\right\}\right).
$$
Then, from Lemma \ref{wmbhg}, it follows
that, for any $f\in L^1_{\mathrm{loc}}(\rrn)$,
$$
\left\|\mc\left(f\right)\right\|_{[W\HerzS]^{1/r}}
\lesssim\left\|f\right\|_{[W\HerzS]^{1/r}}.
$$
This implies that the Hardy--Littlewood
maximal operator $\mc$ is bounded on
the weak Herz space $[W\HerzS]^{1/r}$
and hence further implies that,
for the above $r$,
Proposition \ref{intergx}(iii) holds true.

Next, we prove that, for any $s\in(0,p_-)$,
$[\HerzS]^{1/s}$ is a BBF space.
Indeed, for any $s\in(0,p_-)$, by
the assumptions $\m0(\omega)\in(-\frac{n}{p},
\infty)$ and $\MI(\omega)\in(-\infty,0)$,
and Lemma \ref{rela}, we find that
$$
\m0\left(\omega^s\right)=
s\m0(\omega)>-\frac{n}{p/s}
$$
and
$$
\MI\left(\omega^s\right)
=s\MI(\omega)<0,
$$
which, combined with the assumptions
$p/s,\ q/s\in(1,\infty)$, Theorem \ref{ball}
with $p$, $q$, and $\omega$ therein replaced,
respectively, by $p/s$, $q/s$, and $\omega^s$,
and Lemma \ref{convexl},
further imply that
$[\HerzS]^{1/s}$ is a BBF space.
This finishes the proof that
$\HerzS$ satisfies Proposition \ref{intergx}(iv).

Finally, we show that there exist
a linear space $Y\subset\Msc(\rrn)$,
an $s_0\in(0,p_-)$, and
an $r_0\in(s_0,\infty)$ such that
(v) through (viii) of Proposition
\ref{intergx} hold true.
To this end, let $s_0\in(0,p_-)$ and
$$
r_0\in\left(\max\left\{1,
p,\frac{n}{\mw+n/p}\right\},\infty\right).
$$
Then, repeating an argument similar to that
used in the proof of Theorem \ref{bddog}
with $\eta$ therein replaced by
$(r_0/s_0)'$, we conclude that the following
four statements hold true:
\begin{enumerate}
  \item[{\rm(i)}] $\|\cdot\|_{\HerzScsod}$ satisfies
  Definition \ref{Df1}(ii);
  \item[{\rm(ii)}] $\1bf_{B(\0bf,1)}\in\HerzScsod$;
  \item[{\rm(iii)}] for any $f\in\Msc(\rrn)$,
  $$
  \|f\|_{[\HerzS]^{1/s_0}}
  \sim\sup\left\{\|fg\|_{L^1(\rrn)}:\
  \|g\|_{\HerzScsod}=1\right\}
  $$
  with the positive equivalence
  constants independent of $f$;
  \item[{\rm(iv)}] for any $f\in L^1_{\mathrm{loc}}
  (\rrn)$,
  $$
  \left\|\mc^{((r_0/s_0)')}
  \left(f\right)\right\|_{\HerzScsod}
  \lesssim\left\|f\right\|_{\HerzScsod}.
  $$
\end{enumerate}
These imply that (v) through (viii) of Proposition
\ref{intergx} hold true with the above
$s_0$ and $r_0$, and
$Y:=\HerzScsod$. Therefore,
the assumptions (i) through (viii) of
Proposition \ref{intergx} hold true for
the global generalized Herz space
$\HerzS$ under consideration. This then
implies that, for any $\theta\in(0,1)$,
$$
\left(\HaSaH,\,L^{\infty}(\rrn)\right)_{\theta,\infty}=
WH\Kmp_{\omega^{1-\theta}}^{p/(1-\theta),q/(1-\theta)}(\rrn),
$$
which completes the proof of Theorem
\ref{interg}.
\end{proof}

As an application, by Theorem \ref{interg},
Remark \ref{remark4.10}(ii),
and Remark \ref{remark6.16}, we conclude the
following real interpolation theorem
which shows that the real interpolation
space between $H\MorrS$ and $L^\infty(\rrn)$ is
just the weak generalized Morrey--Hardy space
$WH\textsl{\textbf{M}}_{\omega^{1-\theta}}
^{p/(1-\theta),q/(1-\theta)}(\rrn);$
we omit the details.

\begin{corollary}
Let $p$, $q$, and $\omega$ be as in
Corollary \ref{interlm}.
Then, for any $\theta\in(0,1)$,
$$
\left(H\MorrS,\,L^{\infty}(\rrn)\right)_{\theta,\infty}=
WH\textsl{\textbf{M}}_{\omega^{1-\theta}}^{p/(1-\theta),q/(1-\theta)}(\rrn).
$$
\end{corollary}

\chapter{Inhomogeneous Generalized Herz
Spaces and Inhomogeneous \\ Block
Spaces\label{sec9}}
\markboth{\scriptsize\rm\sc Inhomogeneous Generalized Herz Spaces and
Inhomogeneous Block Spaces}
{\scriptsize\rm\sc Inhomogeneous Generalized Herz Spaces and
Inhomogeneous Block Spaces}

The targets of this chapter are threefold.
The first one is to introduce the
inhomogeneous counterparts
of generalized Herz spaces as in
Definition \ref{gh} and then find
both dual spaces and associate spaces of these
inhomogeneous local generalized Herz spaces.
The second one is to introduce
inhomogeneous block spaces and
prove the duality between
inhomogeneous global generalized Herz spaces
and these inhomogeneous block spaces.
Moreover, we also establish the boundedness
of some important operators on these block spaces.
The last one is, as applications,
to investigate the boundedness and the
compactness characterizations
of commutators on inhomogeneous
generalized Herz spaces.

\section{Inhomogeneous Generalized
Herz Spaces}

In this section, we first introduce
the inhomogeneous counterparts of local and global
generalized Herz spaces
and investigate their basic properties.
Indeed, we obtain the relation between
these Herz spaces and the corresponding
homogeneous generalized Herz spaces.
Then, under some reasonable
and sharp assumptions, we show that
inhomogeneous generalized Herz spaces
are ball (quasi-)Banach function spaces.
Furthermore, we show the convexity and
the absolutely continuity of
quasi-norms of these Herz spaces as well as
establish a boundedness criterion of sublinear operators
and the Fefferman--Stein vector-valued inequalities
on these Herz spaces. Finally, we find both the
dual and the associate space
of the inhomogeneous local generalized
Herz space $\NHerzSo$
when $p$, $q\in(1,\infty)$,
and establish the extrapolation theorems
of both the inhomogeneous local
and the inhomogeneous global generalized Herz
spaces.

We begin this section with the following
definitions of inhomogeneous generalized Herz spaces.

\begin{definition}\label{igh}
Let $p,\ q\in(0,\infty)$ and $\omega\in M(\rp)$.
\begin{enumerate}
  \item[(i)] The \emph{inhomogeneous local generalized
  Herz space}\index{inhomogeneous local generalized \\ Herz space}
  $\NHerzSo$\index{$\NHerzSo$} is defined to be
  the set of all the measurable
  functions $f$ on $\rrn$ such that
  \begin{align*}
\|f\|_{\NHerzSo}:=\left\{
\left\|f\1bf_{B(\0bf,1)}\right\|_{
L^{p}(\rrn)}^{q}
+\sum_{k\in\mathbb{N}}\left[\omega(\tk)
\right]^{q}\left\|f\1bf_{B(\0bf,\tk)\setminus
B(\0bf,\tkm)}\right\|_{L^{p}(\rrn)}
^{q}\right\}^{\frac{1}{q}}
\end{align*}
is finite.
  \item[(ii)] The \emph{inhomogeneous global generalized
  Herz space}\index{inhomogeneous global generalized \\ Herz space}
  $\NHerzS$\index{$\NHerzS$} is defined to be the
  set of all the measurable functions $f$ on $\rrn$ such that
  \begin{align*}
\|f\|_{\NHerzS}:&=\sup_{\xi\in\rrn}
\left\{\left\|f\1bf_{B(\xi,1)}
\right\|_{L^{p}(\rrn)}^{q}\right.
\\&\quad\left.+\sum_{k\in\mathbb{N}}
\left[\omega(\tk)\right]^{q}\left\|f
\1bf_{B(\xi,\tk)\setminus
B(\xi,\tkm)}\right\|_{L^{p}(\rrn)}
^{q}\right\}^{\frac{1}{q}}
\end{align*}
is finite.
\end{enumerate}
\end{definition}

\begin{remark}\label{ighr}
\begin{enumerate}
  \item[(i)] Obviously, by Definition
  \ref{igh}, we conclude that,
  for any $f\in\Msc(\rrn)$,
\begin{equation*}
\|f\|_{\NHerzS}=\sup_{\xi\in\rrn}\left\|
f(\cdot+\xi)\right\|_{\NHerzSo}.
\end{equation*}
  \item[(ii)] In Definition
  \ref{igh}(i), for any given $\alpha\in\rr$
  and for any $t\in(0,\infty)$, let $\omega(t):=t^\alpha$.
  Then, in this case, the inhomogeneous
  local generalized Herz space $\NHerzSo$ coincides with the
classical \emph{inhomogeneous Herz space}\index{inhomogeneous Herz space}
  $K_{p}^{\alpha,q}(\rrn)$\index{$K_{p}^{\alpha,q}(\rrn)$},
  which was originally introduced in \cite[Definition 1.1(b)]{LuY95}
  (see also \cite[Chapter 1]{LYH}),
  with the same quasi-norms. In particular,
  when $p\in(1,\infty)$, $q=1$,
  and $\alpha=n(1-\frac{1}{p})$,
  the inhomogeneous local generalized Herz space
  $\NHerzSo$ is just the \emph{Beurling
  algebra}\index{Beurling algebra} $A^p(\rrn)$\index{$A^p(\rrn)$}
  which was originally introduced in \cite{Beu}
  (see also \cite[Definition 1.1]{GCJ}).
  \item[(iii)] We should point out that, even
  when $\omega(t):=t^\alpha$ for any $t\in(0,\infty)$
  and for any given $\alpha\in\rr$,
  the inhomogeneous global generalized
  Herz space $\NHerzS$ is also new.
\end{enumerate}
\end{remark}

The following conclusion gives the relation
between the local generalized Herz space and its inhomogeneous
counterpart.

\begin{theorem}\label{rehi}
Let $p,\ q\in(0,\infty)$ and $\omega\in M(\rp)$ satisfy
$
\m0(\omega)\in(0,\infty)$ and $\mi(\omega)\in(0,\infty)$.
Then
$$
\NHerzSo=\HerzSo\cap L^{p}(\rrn).
$$
Moreover, there exists a constant $C\in[1,\infty)$ such that,
for any $f\in\NHerzSo$,
$$
C^{-1}\|f\|_{\NHerzSo}\leq\|f\|_{\HerzSo}
+\|f\|_{L^p(\rrn)}\leq C\|f\|_{\NHerzSo}.
$$
\end{theorem}

\begin{proof}
Let all the symbols be as in the present theorem.
We first show that
\begin{equation}\label{rehie1}
\NHerzSo\cap L^p(\rrn)\subset
\HerzSo.
\end{equation}
To this end, let $f\in\NHerzS\cap L^p(\rrn)$.
Then, applying both Definitions \ref{gh}
and \ref{igh}, we find that
\begin{align}\label{er1}
&\|f\|_{\NHerzSo}\notag\\&\quad\lesssim\left\|f
\1bf_{B(\0bf,1)}\right\|_{L^p(\rrn)}+
\left[\sum_{k\in\mathbb{N}}
\left[\omega(\tk)\right]^q\left\|
f\1bf_{B(\0bf,\tk)\setminus
B(\0bf,\tkm)}\right\|_{L^p(\rrn)}
^q\right]^{\frac{1}{q}}\notag\\
&\quad\lesssim\|f
\|_{L^p(\rrn)}+\|f\|_{\HerzSo}<\infty.
\end{align}
This finishes the proof of \eqref{rehie1}.

Conversely, we prove
$\NHerzSo\subset\HerzSo\cap L^p(\rrn)$.
Indeed, from Lemma \ref{Th1},
it follows that, for any $k\in\mathbb{N}$,
$$\omega(\tk)\gtrsim2^{k[\mi(\omega)-\eps]}$$
and, for any $k\in\zmn$, $$\omega(2^{k})
\lesssim2^{k[\m0(\omega)-\eps]},$$
where $\eps\in(0,\mw)$ is a fixed positive constant.

We now claim that, for any $f\in\NHerzSo$,
\begin{equation}\label{er5}
\left\|f\1bf_{[B(\0bf,1)]^{\complement}}
\right\|_{L^p(\rrn)}\lesssim
\|f\|_{\NHerzSo}.
\end{equation}
Indeed, by the assumption that,
for any $k\in\mathbb{N}$,
$\omega(\tk)\gtrsim2^{k[\mi(\omega)-\eps]}$,
we conclude that, for any $f\in\NHerzSo$,
\begin{align}\label{er2}
&\left\|f\1bf_{[B(\0bf,1)]^\complement}
\right\|_{L^p(\rrn)}^p\notag\\&\quad=
\sum_{k\in\mathbb{N}}\left[\omega(\tk)\right]^{-p}
\left[\omega(\tk)\right]^p
\left\|f\1bf_{B(\0bf,\tk)\setminus
B(\0bf,\tkm)}\right\|_{L^p(\rrn)}^p\notag\\
&\quad\lesssim\sum_{k\in\mathbb{N}}
2^{kp[-\mi(\omega)+\eps]}
\left[\omega(\tk)\right]^p
\left\|f\1bf_{B(\0bf,\tk)\setminus
B(\0bf,\tkm)}\right\|_{L^p(\rrn)}^p,
\end{align}
which, together with Lemma \ref{l4320}
and $-\mi(\omega)+\eps\in(-\infty,0)$,
further implies that, for any $q\in(0,p]$,
\begin{align}\label{er6}
&\left\|f\1bf_{[B(\0bf,1)]^\complement}
\right\|_{L^p(\rrn)}^p\notag\\
&\quad\lesssim\left\{\sum_{k\in\mathbb{N}}
2^{kq[-\mi(\omega)+
\eps]}\left[\omega(\tk)\right]^q
\left\|f\1bf_{B(\0bf,\tk)\setminus B(\0bf,\tkm)}
\right\|_{L^p(\rrn)}^q\right\}^{\frac{p}{q}}\notag\\
&\quad\lesssim\left[\sum_{k\in\mathbb{N}}
\left[\omega(\tk)\right]^q\left\|f\1bf_{B(\0bf,\tk)
\setminus B(\0bf,\tkm)}\right\|_{
L^p(\rrn)}^q\right]^{\frac{p}{q}}\notag\\
&\quad\lesssim\|f\|_{\NHerzSo}^p.
\end{align}
This finishes the proof of \eqref{er5}
when $q\in(0,p]$.
On the other hand, applying \eqref{er2},
the H\"{o}lder inequality,
and the assumption $-\mi(\omega)+
\eps\in(-\infty,0)$, we find that,
for any $q\in(p,\infty)$ and $f\in\NHerzSo$,
\begin{align*}
&\left\|f\1bf_{[B(\0bf,1)]^\complement}
\right\|_{L^p(\rrn)}^p\\&\quad
\lesssim\left\{\sum_{k\in\mathbb{N}}
2^{kp[-\mi(\omega)+\eps](q/p)'}\right\}^{[(q/p)']^{-1}}\\
&\qquad\times\left\{\sum_{k\in\mathbb{N}}
\left[\omega(\tk)\right]^q\left\|f\1bf_{B(\0bf,\tk)\setminus B(\0bf,\tkm)}
\right\|_{L^p(\rrn)}^q\right\}^{\frac{p}{q}}\\
&\quad\sim\left\{\sum_{k\in\mathbb{N}}
\left[\omega(\tk)\right]^q\left\|
f\1bf_{B(\0bf,\tk)\setminus
B(\0bf,\tkm)}\right\|_{L^p(\rrn)}^q\right\}^{
\frac{p}{q}}\lesssim\|f\|_{\NHerzSo}^p.
\end{align*}
Combining this and \eqref{er6}, we further
conclude that \eqref{er5}
holds true for any $q\in(0,\infty)$
and hence complete the prove of the above claim.
Thus, from Definition \ref{igh}(i), we deduce that, for any
$f\in\NHerzSo$,
\begin{align}\label{er3}
\|f\|_{L^p(\rrn)}^p=\left\|f\1bf_{B(\0bf,1)}
\right\|_{L^p(\rrn)}^p+
\left\|f\1bf_{[B(\0bf,1)]^\complement}
\right\|_{L^p(\rrn)}^p
\lesssim\|f\|_{\NHerzSo}^p<\infty,
\end{align}
which further implies that
$\NHerzSo\subset L^p(\rrn)$.

We next show that $\NHerzSo\subset\HerzSo$. To achieve this,
let $f\in\NHerzSo$.
Then, by the assumption that,
for any $k\in\zmn$,
$\omega(2^{k})\lesssim2^{k[\m0(\omega)-\eps]}$, Lemma
\ref{l4320}, and
the assumption $\m0(\omega)-\eps\in(0,\infty)$,
we conclude that,
for any $q\in[p,\infty)$,
\begin{align}\label{er9}
&\sum_{k\in\zmn}
\left[\omega(\tk)\right]^q
\left\|f\1bf_{B(\0bf,\tk)\setminus
B(\0bf,\tkm)}\right\|_{L^p(\rrn)}^q\notag\\
&\quad\lesssim
\sum_{k\in\zmn}2^{kq[\m0(\omega)-\eps]}
\left\|f\1bf_{B(\0bf,\tk)\setminus
B(\0bf,\tkm)}\right\|_{L^p(\rrn)}^q\notag\\
&\quad\lesssim\left\{\sum_{k\in\zmn}
2^{kp[\m0(\omega)-\eps]}\left\|
f\1bf_{B(\0bf,\tk)\setminus
B(\0bf,\tkm)}\right\|_{L^p(\rrn)}
^p\right\}^{\frac{q}{p}}\notag\\
&\quad\lesssim\left[\sum_{k\in\zmn}
\left\|f\1bf_{B(\0bf,\tk)\setminus
B(\0bf,\tkm)}\right\|_{L^p(\rrn)}^p
\right]^{\frac{q}{p}}
\sim\left\|f\1bf_{B(\0bf,1)}
\right\|_{L^p(\rrn)}^q.
\end{align}
On the other hand, using the
assumption that, for any $k\in\zmn$,
$\omega(2^{k})\lesssim2^{
k[\m0(\omega)-\eps]}$, the H\"{o}lder
inequality, and the
assumption $\m0(\omega)-\eps
\in(0,\infty)$ again, we
find that, for any $q\in(0,p)$,
\begin{align*}
&\sum_{k\in\zmn}
\left[\omega(\tk)\right]^q
\left\|f\1bf_{B(\0bf,\tk)\setminus
B(\0bf,\tkm)}\right\|_{L^p(\rrn)}^q\\
&\quad\lesssim
\sum_{k\in\zmn}2^{kq[\m0(\omega)-\eps]}
\left\|f\1bf_{B(\0bf,\tk)\setminus
B(\0bf,\tkm)}\right\|_{L^p(\rrn)}^q\notag\\
&\quad\lesssim\left\{\sum_{k\in\zmn}
2^{kq[\m0(\omega)-\eps](p/q)'}\right\}^{
(p/q)'^{-1}}
\left[\sum_{k\in\zmn}\left\|
f\1bf_{B(\0bf,\tk)\setminus
B(\0bf,\tkm)}\right\|_{L^p(\rrn)}^p
\right]^{\frac{q}{p}}\notag\\
&\quad\sim\left[\sum_{k\in\zmn}
\left\|f\1bf_{B(\0bf,\tk)\setminus
B(\0bf,\tkm)}\right\|_{
L^p(\rrn)}^p\right]^{\frac{q}{p}}
\sim\left\|f\1bf_{B(\0bf,1)}
\right\|_{L^p(\rrn)}^q,
\end{align*}
which, combined with \eqref{er9}, further implies that
\begin{align}\label{er4}
&\|f\|_{\HerzSo}\notag\\&\quad=\left\{\sum_{k\in\zmn}
\left[\omega(\tk)\right]^q\left\|f\1bf_{B(\0bf,\tk)
\setminus B(\0bf,\tkm)}\right\|_{L^p(\rrn)}^q+
\sum_{k\in\mathbb{N}}\cdots\right\}^{\frac{1}{q}}\notag\\
&\quad\lesssim\left\{\left\|f\1bf_{
B(\0bf,1)}\right\|_{L^p(\rrn)}^q+\sum_{k\in\mathbb{N}}
\left[\omega(\tk)\right]^q\left\|
f\1bf_{B(\0bf,\tk)\setminus
B(\0bf,\tkm)}\right\|_{L^p(\rrn)}^q
\right\}^\frac{1}{q}\notag\\
&\quad\sim\|f\|_{\NHerzSo}<\infty.
\end{align}
This implies that $\NHerzSo\subset\HerzSo$.
Therefore, we complete the proof that
$$\NHerzSo=\HerzSo\cap L^p(\rrn).$$

Moreover, from \eqref{er1}, \eqref{er3},
and \eqref{er4},
it follows that, for any $f\in\NHerzSo$,
$$
\|f\|_{\NHerzSo}\sim\|f\|_{\HerzSo}+\|f\|_{L^p(\rrn)}
$$
with positive equivalence constants independent
of $f$. This finishes
the proof of Theorem \ref{rehi}.
\end{proof}

\begin{remark}
We should point out that,
in Theorem \ref{rehi}, when $\omega(t):=t^\alpha$
for any $t\in(0,\infty)$
and for any given
$\alpha\in\rr$, Theorem \ref{rehi}
goes back to \cite[Proposition 1.1.2]{LYH}.
\end{remark}

For the relation between the global generalized
Herz space $\HerzS$ and its
inhomogeneous counterpart $\NHerzS$,
we have the following conclusion.

\begin{theorem}\label{rehig}
Let $p,\ q\in(0,\infty)$ and $\omega\in M(\rp)$ satisfy
$\mi(\omega)\in(0,\infty)$. Then
$$
\NHerzS=\HerzS=\left\{f\in\Msc(\rrn):\
f(x)=0\text{ for almost
every $x\in\rrn$}\right\}.
$$
\end{theorem}

\begin{proof}
Let all the symbols be as in the present theorem
and $f$ be a measurable function on $\rrn$ such that
\begin{equation}\label{rehige1}
\left|\left\{x\in\rrn:\ |f(x)|>0\right\}\right|\neq0.
\end{equation}
We first show that there exists a $k_{0}\in
\mathbb{N}$ such that
$$
|\{x\in B(\0bf,2^{k_{0}}):\ |f(x)|>0\}|\neq0.
$$
Indeed, if, for any
$k\in\mathbb{N}$, $|\{x\in B(\0bf,\tk):\ |f(x)|>0\}|=0$,
from this and the
fact that
$$\left\{x\in\rrn:\
|f(x)|>0\right\}=\bigcup_{k\in\mathbb{N}}
\left\{x\in B(\0bf,\tk):\
|f(x)|>0\right\},$$
we deduce that
\begin{align*}
\left|\left\{x\in\rrn:\ |f(x)|>0\right\}\right|
\leq\sum_{k\in\mathbb{N}}
\left|\left\{x\in B(\0bf,\tk):\ |f(x)|>0\right\}\right|=0,
\end{align*}
which contradicts \eqref{rehige1}. Similarly
to the argument above, we conclude that there
exists a $k_1\in\mathbb{N}$ such that
$|E_{f,k_0,k_1}|\in(0,\infty)$, where
$$
E_{f,k_0,k_1}:=\left\{x\in B(\0bf,2^{k_{0}}):\
|f(x)|>\frac{1}{k_{1}}\right\}.
$$

Next, for any $i\in\mathbb{N}\cap[k_0+1,\infty)$,
let $\xi_{i}\in\rrn$ satisfy $|\xi_{i}|=2^{i}+2^{k_{0}}$.
Then, for any $x\in B(\0bf,2^{k_{0}})$ and
$i\in\mathbb{N}\cap[k_{0}+1,\infty)$, we have
\begin{align*}
|x-\xi_{i}|\leq|x|+|\xi_{i}|
<2^{k_{0}}+2^i+2^{k_{0}}\leq2^{i+1}\end{align*} and
\begin{align*}
|x-\xi_{i}|\geq|\xi_{i}|-|x|>2^i+2^{k_{0}}-2^{k_{0}}
=2^i.\end{align*}
This further implies that,
for any $i\in\mathbb{N}\cap[k_0+1,\infty)$,
$$B(\0bf,2^{k_{0}})\subset B(\xi_{i},2^{i+1})\setminus
B(\xi_{i},2^i).$$
Applying this and both Definitions \ref{gh} and \ref{igh},
we find that,
for any $i\in\mathbb{N}\cap[k_0+1,\infty),$
\begin{align}\label{ee1}
\|f\|_{\NHerzS}&\geq\left\|f\1bf_{B(\0bf,2^{k_{0}})}(
\cdot+\xi_{i})\right\|_{\NHerzSo}
=\omega(2^{i+1})\left\|f\1bf_{B(\0bf,2^{k_{0}})}
\right\|_{L^p(\rrn)}\notag\\&\geq\omega(2^{i+1})
\frac{|E_{f,k_0,k_1}|^{1/p}}{k_{1}}
\end{align}
and
\begin{align}\label{ee2}
\|f\|_{\HerzS}&\geq\left\|f\1bf_{B(\0bf,2^{k_{0}})}(
\cdot+\xi_{i})\right\|_{\HerzSo}
=\omega(2^{i+1})\left\|f\1bf_{B(\0bf,2^{k_{0}})}\right\|_{
L^p(\rrn)}\notag\\&\geq\omega
(2^{i+1})\frac{|E_{f,k_0,k_1}|^{1/p}}{k_{1}}.
\end{align}
On the other hand, from Lemma \ref{Th1},
it follows that, for any $k\in\mathbb{N}$,
$$\omega(2^{k})\gtrsim2^{k[\mi(\omega)-\eps]},$$
where $\eps\in(0,\mi(\omega))$ is a fixed positive constant.
This, together with
\eqref{ee1} and \eqref{ee2},
further implies that,
for any $i\in\mathbb{N}\cap[k_{0}+1,\infty)$,
$$
\|f\|_{\NHerzS}\gtrsim2^{i[\mi(\omega)-\eps]}
$$
and
$$
\|f\|_{\HerzS}\gtrsim2^{i[\mi(\omega)-\eps]}.
$$
Using this and the assumption
$\mi(\omega)\in(0,\infty)$, we further find that
\begin{align*}
\|f\|_{\NHerzS}=\infty=\|f\|_{\HerzS},
\end{align*}
which completes the proof of Theorem \ref{rehig}.
\end{proof}

Now, we show that the
inhomogeneous local generalized
Herz space $\NHerzSo$ is
a ball quasi-Banach function space.

\begin{theorem}\label{Th3n}
Let $p,\ q\in(0,\infty)$ and
$\omega\in M(\rp)$. Then
the inhomogeneous local generalized
Herz space $\NHerzSo$
is a ball quasi-Banach function space.
\end{theorem}

\begin{proof}
Let $p$, $q\in(0,\infty)$ and $\omega\in M(\rp)$.
Obviously, the inhomogeneous
local generalized Herz space $\NHerzSo$ is a
quasi-normed linear space
satisfying both (i) and (ii) of Definition \ref{Df1}.
Moreover, Definition \ref{Df1}(iii)
is a simple corollary of
the monotone convergence theorem.
Therefore, using Proposition \ref{xcom},
we further find that
$\NHerzSo$ is complete and hence $\NHerzSo$ is a
quasi-Banach space.

Thus, to finish the proof of the present theorem,
it remains to show that the inhomogeneous
local generalized Herz space $\NHerzSo$
satisfies Definition \ref{Df1}(iv). To this end,
let $B(x_0,r)\in\mathbb{B}$ with $x_0\in\rrn$
and $r\in(0,\infty)$.
Then, for any $k\in\mathbb{N}\cap(\ln(r+|x_{0}|)/\ln2+1,\infty)$,
we have $2^{k-1}>r+|x_{0}|$.
This implies that,
for any $x\in B(x_{0},r)$,
\begin{align*}
|x|\leq|x-x_{0}|+|x_{0}|<r+|x_{0}|<\tkm.
\end{align*}
Therefore, $x\in B(\0bf,\tkm)$ and hence
$B(x_{0},r)\subset B(\0bf,\tkm)$.
By this, we further conclude that
$$B(x_{0},r)\cap\left[B(\0bf,\tk)\setminus
B(\0bf,\tkm)\right]=\emptyset.$$
Combining this and Definition \ref{igh}(i), we find that
\begin{align*}
&\left\|\1bf_{B(x_0,r)}\right\|_{\NHerzSo}^q\\
&\quad=\left\|\1bf_{B(x_0,r)}
\1bf_{B(\0bf,1)}\right\|_{L^p(\rrn)}^q\\
&\qquad+\sum_{k\in\mathbb{N}
\cap[1,\frac{\ln(r+|x_0|)}{\ln2}]}
\left[\omega(\tk)\right]^q\left\|\1bf_{B(x_0,r)}
\1bf_{B(\0bf,\tk)\setminus B(
\0bf,\tkm)}\right\|_{L^p(\rrn)}^q\\&\quad<\infty,
\end{align*}
which implies that $\1bf_{B(x_0,r)}\in\NHerzSo$ and hence
finishes the proof of Theorem \ref{Th3n}.
\end{proof}

However, the following example shows
that inhomogeneous global generalized
Herz spaces may not be ball quasi-Banach
function spaces.

\begin{example}\label{Ex2n}
Let $p$, $q\in(0,\infty)$, $\alpha_1\in\rr$,
$\alpha_2\in[0,\infty)$, and
$$
\omega(t):=\left\{\begin{aligned}
&t^{\alpha_{1}}(1-\ln t)
\hspace{0.3cm}\text{when }t\in(0,1],\\
&t^{\alpha_{2}}(1+\ln t)
\hspace{0.3cm}\text{when }t\in(1,\infty).
\end{aligned}\right.
$$
Then $\1bf_{B(\0bf,1)}\notin\NHerzS$
and hence, in this case, the
inhomogeneous global generalized Herz space
$\NHerzS$ is not a
ball quasi-Banach function space.
\end{example}

\begin{proof}
Let all the symbols be as in the present example
and, for any $k\in\mathbb{N}$, $\xi_k\in\rrn$
satisfy $|\xi_k|=2^{k}+1$. Then,
by the proof of Example \ref{Ex2},
we find that,
for any $k\in\mathbb{N}$,
$$
B(\0bf,1)\subset
B(\xi_k,2^{k+1})\setminus B(\xi_{k},\tk).
$$
From this, Remark \ref{ighr}, and Definition
\ref{igh}(i), we deduce that
\begin{align*}
\left\|\1bf_{B(\0bf,1)}\right\|_{\NHerzS}
&\geq\left\|\1bf_{B(\0bf,1)}
(\cdot+\xi_k)\right\|_{\NHerzSo}
=\omega(2^{k+1})\left\|
\1bf_{B(\0bf,1)}\right\|_{L^p(\rrn)}\\
&\sim2^{(k+1)\alpha_2}
\left[1+(k+1)\ln2\right]\to\infty
\end{align*}
as $k\to\infty$. This further implies that
$\1bf_{B(\0bf,1)}\notin\NHerzS$
and hence the inhomogeneous global
generalized Herz space $\NHerzS$
under consideration
is not a BQBF space.
This then finishes
the proof of Example \ref{Ex2n}.
\end{proof}

Under a reasonable and
sharp assumption, we next
prove that the inhomogeneous
global generalized
Herz space $\NHerzS$ is also a
ball quasi-Banach function space.
Namely, the following conclusion holds true.

\begin{theorem}\label{Th2n}
Let $p,\ q\in(0,\infty)$ and $\omega\in M(\rp)$ satisfy
$\MI(\omega)\in(-\infty,0)$. Then the inhomogeneous
global generalized Herz space $\NHerzS$
is a ball quasi-Banach function space.
\end{theorem}

\begin{proof}
Let $p$, $q\in(0,\infty)$
and $\omega\in M(\rp)$ with
$\MI(\omega)\in(-\infty,0)$.
Obviously, the inhomogeneous global
generalized Herz space $\NHerzS$
is a quasi-normed linear space
satisfying both (i) and (ii)
of Definition \ref{Df1}.
We now show that Definition
\ref{Df1}(iii) holds true for $\NHerzS$.
Indeed, for any given
$\{f_{m}\}_{m\in\mathbb{N}}\subset
\Msc(\rrn)$ and $f\in\Msc(\rrn)$
satisfying $0\leq f_{m}\uparrow f$
almost everywhere as $m\to\infty$,
and for any given
$\alpha\in(0,\|f\|_{\NHerzS})$,
from Remark \ref{ighr}(i),
we deduce that there exists a
$\xi\in\rrn$ such that
$$\left\|f(\cdot+\xi)
\right\|_{\NHerzSo}>\alpha.$$
Using this and the monotone convergence theorem,
we conclude that there exists an $N\in\mathbb{N}$ such that,
for any $m\in\mathbb{N}\cap(N,\infty)$,
$$\left\|f_{m}(\cdot+\xi)\right\|_{\NHerzSo}>\alpha.$$
By this and Remark \ref{ighr}(i) again, we find that,
for any $m\in\mathbb{N}\cap(N,\infty)$,
\begin{align*}
\alpha<\sup_{\xi\in\rrn}\left\|f_{m}(\cdot+\xi)\right\|_{\NHerzSo}
=\|f_{m}\|_{\NHerzS}.
\end{align*}
Letting $\alpha\to\|f\|_{\NHerzS}$ and $m\to\infty$,
we have
$$
\|f\|_{\NHerzS}\leq\lim\limits_{m\to\infty}
\|f_m\|_{\NHerzS}.
$$
On the other hand, it is easy to
show that
$$
\lim\limits_{m\to\infty}
\left\|f_m\right\|_{\NHerzS}\leq
\left\|f\right\|_{\NHerzS}.
$$
Therefore, we conclude that $$\left\|f_{m}
\right\|_{\NHerzS}
\uparrow\left\|f\right\|_{\NHerzS}$$
as $m\to\infty$
and hence Definition \ref{Df1}(iii) holds true
for $\NHerzS$. From this and
Proposition \ref{xcom}, it follows that
the inhomogeneous global generalized
Herz space $\NHerzS$
is quasi-Banach space.

Next, we prove that Definition \ref{Df1}(iv) holds true
for $\NHerzS$. Indeed, by
Lemma \ref{Th1}, we find that, for any $k\in\mathbb{N}$,
$$\omega(\tk)\lesssim2^{k[\MI(\omega)+\eps]},$$ where
$\eps\in(0,-\MI(\omega))$ is a fixed positive constant.
Combining this, Definition \ref{igh}(i),
and the assumption $\MI(\omega)+\eps\in(-\infty,0)$,
we further conclude
that, for any $B(x_0,r)\in\mathbb{B}$
with $x_0\in\rrn$ and $r\in(0,\infty)$,
and $\xi\in\rrn$,
\begin{align*}
&\left\|\1bf_{B(x_0,r)}(\cdot+\xi)
\right\|_{\NHerzSo}^q\\
&\quad=\left\|\1bf_{B(x_0,r)}
\1bf_{B(\xi,1)}\right\|_{L^p(\rrn)}^q\\
&\qquad+\sum_{k\in\mathbb{N}}
\left[\omega(\tk)\right]^q
\left\|\1bf_{B(x_0,r)}\1bf_{B(\xi,\tk)
\setminus B(\xi,\tkm)}\right\|_{L^p(\rrn)}^q\\
&\quad\lesssim\left\|\1bf_{B(x_0,r)}
\right\|_{L^p(\rrn)}^q
\left\{1+\sum_{k\in\mathbb{N}}
2^{kq[\MI(\omega)+\eps]}\right\}
<\infty.
\end{align*}
This, together with Remark \ref{ighr}(i), further implies
that, for any $x_0\in\rrn$ and $r\in(0,\infty)$,
$$
\left\|\1bf_{B(x_0,r)}\right\|_{\NHerzS}\lesssim1.
$$
Therefore, $\1bf_{B(x_0,r)}\in\NHerzS$ for
any $x_0\in\rrn$ and $r\in(0,\infty)$.
This implies that Definition \ref{Df1}(iv) holds true
for $\NHerzS$, and hence finishes
the proof of Theorem \ref{Th2n}.
\end{proof}

\begin{remark}
By Examples \ref{ex117} and \ref{Ex2n},
 we find that the assumption
$\MI(\omega)\in(-\infty,0)$ in Theorem
\ref{Th2n} is sharp. Indeed, if $\MI
(\omega)\in[0,\infty)$,
then both Examples \ref{ex117} and \ref{Ex2n}
show that there exists an $\omega\in M(\rp)$
such that $\NHerzS$ is not a ball
quasi-Banach function space.
\end{remark}

Moreover, the following
theorem indicates that the inhomogeneous
local generalized Herz space $\NHerzSo$
is a ball Banach function space
when $p,$ $q\in[1,\infty)$.

\begin{theorem}\label{ballln}
Let $p,\ q\in[1,\infty)$ and $\omega\in M(\rp)$.
Then the inhomogeneous local
generalized Herz space $\NHerzSo$ is a ball Banach function space.
\end{theorem}

\begin{proof}
Let $p$, $q\in[1,\infty)$ and $\omega\in M(\rp)$.
Then, using Theorem \ref{Th3n}, we find that
the inhomogeneous local generalized Herz space $\NHerzSo$
is a BQBF space. Moreover,
$\NHerzSo$ obviously
satisfies the triangle inequality
due to $p$, $q\in[1,\infty)$.

Thus, to complete
the proof of the present theorem,
it remains to show that \eqref{bal}
holds true with $X:=\NHerzSo$.
To achieve this, we first claim that,
for any given $k_0\in\mathbb{Z}$,
\eqref{bal} holds true with $X:=\NHerzSo$
and $B:=B(\0bf,2^{k_0})$.
We show this claim by
considering the following two cases
on $k_0$.

\emph{Case 1)} $k_0\in\zmn$. In this case, we have
$B(\0bf,2^{k_0})\subset B(\0bf,1)$. By this, the
H\"{o}lder inequality, and
Definition \ref{igh}(i), we conclude
that, for any $f\in\NHerzSo$,
\begin{align}\label{th2ne1}
\int_{B(\0bf,2^{k_0})}|f(y)|\,dy&\lesssim\left\|
f\1bf_{B(\0bf,2^{k_0})}\right\|_{L^p(\rrn)}
\lesssim\left\|
f\1bf_{B(\0bf,1)}\right\|_{L^p(\rrn)}\notag\\
&\lesssim\|f\|_{\NHerzS},
\end{align}
which implies that the above claim holds true in this case.

\emph{Case 2)} $k_0\in\mathbb{N}$. In this case, from
the H\"{o}lder inequality, the fact that,
for any given $\alpha\in(0,\infty)$ and
$m\in\mathbb{N}$, and for any $\{a_j\}_{j\in\mathbb{N}}
\subset\mathbb{C}$,
$$
\left(\sum_{j=1}^m|a_j|\right)^\alpha
\leq\max\left\{1,m^{\alpha-1}\right\}\sum_{j=1}^{m}
|a_j|^{\alpha},
$$
and Definition \ref{igh}(i),
we deduce that, for any $f\in\NHerzSo$,
\begin{align*}
&\int_{B(\0bf,2^{k_0})}|f(y)|\,dy\\&\quad\lesssim\left\|
f\1bf_{B(\0bf,2^{k_0})}\right\|_{L^p(\rrn)}\\
&\quad\sim\left[\left\|f\1bf_{B
(\0bf,1)}\right\|_{L^p(\rrn)}^p+\sum_{k=1}^{k_0}
\left\|f\1bf_{B(\0bf,\tk)\setminus
B(\0bf,\tkm)}\right\|_{L^p(\rrn)}^p\right]^{\frac{1}{p}}\\
&\quad\lesssim\left\{\left\|
f\1bf_{B(\0bf,1)}\right\|_{L^p(\rrn)}^q
+\sum_{k=1}^{k_0}\left[\omega(\tk)\right]^q
\left\|f\1bf_{B(\0bf,\tk)
\setminus B(\0bf,\tkm)}
\right\|_{L^p(\rrn)}^q\right\}^{\frac{1}{q}}\\
&\quad\lesssim\|f\|_{\NHerzSo}.
\end{align*}
This, combined with \eqref{th2ne1}, then finishes
the proof of the above claim.

Applying the above claim, we
then show that \eqref{bal} holds
true for any $B\in\mathbb{B}$.
Notice that, for any $B\in\mathbb{B}$,
there exists a $k\in\mathbb{Z}$ such that
$B\subset B(\0bf,2^{k})$. By this and the above claim,
we further conclude that,
for any $f\in\NHerzSo$,
\begin{align*}
\int_{B}|f(y)|\,dy\leq\int_{B(\0bf,\tk)}|f(y)|\,dy
\lesssim\|f\|_{\NHerzSo},
\end{align*}
which implies that \eqref{bal} holds true with
$X:=\NHerzSo$, and hence finishes the
proof of Theorem \ref{Th2n}.
\end{proof}

Similarly, we now show that the inhomogeneous
global generalized Herz space $\NHerzS$ is also a ball Banach
function space when $p$, $q\in[1,\infty)$.

\begin{theorem}\label{balln}
Let $p,\ q\in[1,\infty)$ and $\omega\in M(\rp)$ satisfy
$\MI(\omega)\in(-\infty,0)$.
Then the inhomogeneous global generalized Herz
space $\NHerzS$ is a ball Banach function space.
\end{theorem}

\begin{proof}
Let all the symbols be as in the present theorem.
Then, by the assumption $\MI(\omega)\in(-\infty,0)$
and Theorem \ref{Th2n}, we find that the inhomogeneous
global generalized Herz space $\NHerzS$ is a BQBF
space. Moreover, notice that
$\NHerzS$ satisfies the triangle inequality
when $p$, $q\in[1,\infty)$.
Thus, to finish the proof of the present theorem,
we only need to show that,
for any given $B\in\mathbb{B}$,
\eqref{bal} holds true with $X:=\NHerzS$.
Indeed, from Definition \ref{igh}(ii), we deduce that,
for any $f\in\NHerzS$,
$$\|f\|_{\NHerzSo}\leq\|f\|_{\NHerzS}$$
and hence $f\in\NHerzSo$.
Using this, Theorem \ref{ballln},
and Definition \ref{igh}(ii) again,
we further conclude that, for any $f\in\NHerzS$,
\begin{align*}
\int_{B}|f(y)|\,dy\lesssim\|f\|_{\NHerzSo}\lesssim\|f\|_{\NHerzS}.
\end{align*}
This implies that \eqref{bal} holds true with
$X:=\NHerzS$, and then finishes the
proof of Theorem \ref{balln}.
\end{proof}

\subsection{Convexities}

The main target of this subsection is to show
the convexity of the inhomogeneous
generalized Herz spaces.
For this purpose, we first investigate
the relations between these Herz spaces
and their convexifications as follows,
which are useful in the study of the inhomogeneous
generalized Herz--Hardy spaces in next chapters.

\begin{lemma}\label{convexlln}
Let $p,\ q,\ s\in(0,\infty)$ and $\omega\in M(\rp)$.
Then
$$\left[\NHerzSo\right]^{1/s}=\NHerzSocs$$
with the same quasi-norms.
\end{lemma}

\begin{proof}
Let all the symbols be as in the present lemma.
Then we find that, for any $f\in\Msc(\rrn)$,
\begin{align}\label{convexllne1}
&\|f\|_{[\NHerzSo]^{1/s}}\notag\\&\quad=\left\|
|f|^{\frac{1}{s}}\right\|_{\NHerzSo}^s\notag\\
&\quad=\left\{\left[\int_{|y|<1}|f(y)|^{
\frac{p}{s}}\,dy\right]
^{\frac{q}{p}}+\sum_{k\in\mathbb{N}}
\left[\omega(\tk)\right]^{q}
\left[\int_{\tkm\leq|y|<\tk}|f(y)|^{
\frac{p}{s}}\,dy\right]
^{\frac{q}{p}}\right\}^{\frac{s}{q}}\notag\\
&\quad=\left\{\left[\int_{|y|<1}|f(y)|^{
\frac{p}{s}}\,dy\right]
^{\frac{q/s}{p/s}}\right.\notag\\
&\left.\qquad+\sum_{k\in\mathbb{N}}
\left[\omega^s(\tk)\right]
^{\frac{q}{s}}\left[\int_{\tkm\leq|y|<\tk}
|f(y)|^{\frac{p}{s}}\,dy\right]
^{\frac{q/s}{p/s}}\right\}^{\frac{s}{q}}\notag\\
&\quad=\|f\|_{\NHerzSocs},
\end{align}
which completes the proof of Lemma \ref{convexlln}.
\end{proof}

\begin{lemma}\label{convexln}
Let $p,\ q,\ s\in(0,\infty)$ and $\omega\in M(\rp)$.
Then
$$\left[\NHerzS\right]^{1/s}=\NHerzScs$$
with the same quasi-norms.
\end{lemma}

\begin{proof}
Let $p$, $q$, $s\in(0,\infty)$ and $\omega\in M(\rp)$.
Then, applying an argument similar to
that used in the proof of
\eqref{convexllne1}, we conclude that,
for any $f\in\Msc(\rrn)$ and $\xi\in\rrn$,
$$
\left\|f(\cdot+\xi)\right\|_{[\NHerzSo]^{1/s}}
=\left\|f(\cdot+\xi)\right\|_{\NHerzSocs},
$$
which, together with Remark \ref{ighr}(i),
further implies that
$$\|f\|_{[\NHerzS]^{1/s}}=\|f\|_{\NHerzScs}.$$
This finishes the proof of Lemma \ref{convexln}.
\end{proof}

Via both Lemmas \ref{convexlln} and
\ref{convexln}, we next prove
the convexity of inhomogeneous
generalized Herz spaces. Namely,
we have the following two conclusions.

\begin{theorem}\label{strcn}
Let $p$, $q\in(0,\infty)$, $s\in(0,
\min\{p,q\}]$, and $\omega\in M(\rp)$.
Then the inhomogeneous
local generalized Herz space $\NHerzSo$ is
strictly $s$-convex.\index{strictly $s$-convex}
\end{theorem}

\begin{proof}
Let all the symbols be as in the present theorem.
Then, from the assumption $s\in(0,\min\{p,q\}]$,
it follows that $p/s,\ q/s\in[1,\infty)$.
This, combined with Lemma \ref{convexlln} and
the Minkowski inequality, implies that,
for any sequence
$\{f_{j}\}_{j\in\mathbb{N}}$ of measurable functions in
$[\NHerzSo]^{1/s}$, and for any $N\in\mathbb{N}$,
\begin{align*}
\left\|\sum_{j=1}^{N}
|f_{j}|\right\|_{[\NHerzSo]^{1/s}}\leq\sum_{j=1}^{N}
\|f_{j}\|_{[\NHerzSo]^{1/s}}\leq\sum_{j\in\mathbb{N}}
\|f_{j}\|_{[\NHerzSo]^{1/s}}.
\end{align*}
By this and the monotone convergence theorem, we further find that
\begin{equation}\label{strcne1}
\left\|\sum_{j\in\mathbb{N}}|f_{j}|\right\|_{[\NHerzSo]^{1/s}}\leq
\sum_{j\in\mathbb{N}}\|f_{j}\|_{[\NHerzSo]^{1/s}},
\end{equation}
which completes the proof of Theorem \ref{strcn}.
\end{proof}

\begin{theorem}\label{Th3.1n}
Let $p,\ q\in(0,\infty)$, $s\in(0,\min\{p,q\}]$,
and $\omega\in M(\rp)$.
Then the inhomogeneous
global generalized Herz space $\NHerzS$ is strictly
$s$-convex.\index{strictly $s$-convex}
\end{theorem}

\begin{proof}
Let all the symbols be as in the present theorem.
Then, similarly to the estimation
of \eqref{strcne1}, using Lemma \ref{convexln},
we conclude that, for any given sequence
$\{f_{j}\}_{j\in\mathbb{N}}$ of
measurable functions in $[\NHerzS]^{1/s}$,
$$
\left\|\sum_{j\in\mathbb{N}}|f_{j}|\right\|_{[\NHerzS]^{1/s}}\leq
\sum_{j\in\mathbb{N}}\|f_{j}\|_{[\NHerzS]^{1/s}}.
$$
This finishes the proof of Theorem \ref{Th3.1n}.
\end{proof}

\subsection{Absolutely Continuous Quasi-Norms}

In this subsection, we investigate the absolutely
continuity of the quasi-norms of
inhomogeneous generalized
Herz spaces. To be precise, we show that the
inhomogeneous local generalized Herz space
$\NHerzSo$ has an
absolutely continuous quasi-norm but find that
the inhomogeneous global generalized
Herz space may
not have an absolutely continuous
quasi-norm via a counterexample.

We first prove that $\NHerzSo$ has an absolutely continuous
quasi-norm as follows.

\begin{theorem}\label{abson}
Let $p,\ q\in(0,\infty)$ and $\omega\in M(\rp)$.
Then the inhomogeneous
local generalized Herz space
$\NHerzSo$ has an absolutely continuous quasi-norm.
\end{theorem}

\begin{proof}
Let $p,\ q\in(0,\infty)$, $\omega\in M(\rp)$,
and $f$ be a given measurable function
in $\NHerzSo$. We now show that $f$ has
an absolutely continuous quasi-norm in $\NHerzSo$.
To this end, let $\{E_{i}\}_{i\in\mathbb{N}}$
be a sequence of measurable sets satisfying $\1bf_{E_{i}}\to0$
almost everywhere as $i\to\infty$.
Then, for any $k$, $i\in\mathbb{N}$, let
$$
a_0:=\left\|f\1bf_{B(\0bf,1)}\right\|_{L^p(\rrn)},
$$
$$
a_{0,i}:=\left\|f\1bf_{B(\0bf,1)}\1bf_{E_i}\right\|_{L^p(\rrn)},
$$
$$
a_{k}:=\omega(\tk)\left\|f\1bf_{B(\0bf,
\tk)\setminus B(\0bf,\tkm)}\right\|_{L^{p}(\rrn)},
$$
and
$$
a_{k,i}:=\omega(\tk)\left\|f
\1bf_{B(\0bf,
\tk)\setminus
B(\0bf,\tkm)}\1bf_{E_{i}}\right\|_{L^{p}(\rrn)}.
$$
From the assumption
$f\in\NHerzSo$, it follows that
\begin{equation}\label{len}
\left(\sum_{k\in\zp}|a_{k}|^{q}
\right)^{\frac{1}{q}}<\infty.
\end{equation}
This implies that, for any $k\in\zp$,
$a_k\in[0,\infty)$. Therefore,
$f\1bf_{B(\0bf,1)}\in L^p(\rrn)$
and, for any $k\in\mathbb{N}$,
$$f\1bf_{B(\0bf,\tk)\setminus
B(\0bf,\tkm)}\in L^{p}(\rrn).$$
By this, the facts that, for any $k,\ i\in\mathbb{N}$,
$$
\left|f\1bf_{B(\0bf,1)}\1bf_{E_i}\right|
\leq\left|f\1bf_{B(\0bf,1)}\right|
$$
and
$$\left|f\1bf_{B(\0bf,\tk)\setminus
B(\0bf,\tkm)}\1bf_{E_{i}}\right|\leq
\left|f\1bf_{B(\0bf,\tk)\setminus
B(\0bf,\tkm)}\right|,$$
and the dominated convergence theorem,
we conclude that, for any $k\in\mathbb{Z}_+$,
$a_{k,i}\to0$ as $i\to\infty$.
Applying this, the fact that,
for any $k,\ i\in\mathbb{N}$,
$|a_{k,i}|\leq|a_{k}|$, \eqref{len}, and
the dominated convergence theorem again, we find that
\begin{align*}
\lim\limits_{i\to\infty}\|f\1bf_{E_{i}}\|_{\NHerzSo}
=\lim\limits_{i\to\infty}
\left(\sum_{k\in\mathbb{Z}_+}|a_{k,i}|^{q}\right)^{\frac{1}{q}}
=\left(\sum_{k\in\mathbb{Z}_+}
\lim\limits_{i\to\infty}|a_{k,i}|^{q}\right)^{\frac{1}{q}}=0,
\end{align*}
which implies that $f$ has an absolutely continuous
quasi-norm in $\NHerzSo$.
Thus, by the arbitrariness of $f$,
we further conclude that the inhomogeneous
local generalized Herz space $\NHerzSo$ has an absolutely
continuous quasi-norm, and hence
complete the proof of Theorem \ref{abson}.
\end{proof}

On the other hand, via the following counterexample,
we find that inhomogeneous global generalized
Herz spaces may not have absolutely continuous quasi-norms.

\begin{example}\label{countern}
Let $p,\ q\in(0,\infty)$, $\alpha\in(-\infty,0)$,
and $$E:=\bigcup\limits_{k\in\mathbb{N}}\left(
k-1+k^{-\frac{2}{\alpha p}},k+k^{-\frac{2}{\alpha p}}\right).$$
For any $t\in(0,\infty)$, let $\omega(t):=t^{\alpha}$.
Then the characteristic function
$\1bf_{E}\in\mathcal{K}_{\omega}^{p,q}(\rr)$,
but $\1bf_{E}$
does not have an absolutely continuous quasi-norm in
$\mathcal{K}_{\omega}^{p,q}(\rr)$.
This implies that the inhomogeneous
global generalized Herz space
$\mathcal{K}_{\omega}^{p,q}(\rr)$ does not
have an absolutely continuous quasi-norm.
\end{example}

\begin{proof}
Let all the symbols be as in the present example.
From Definition \ref{igh}(i), we deduce that,
for any $\xi\in\rr$,
\begin{align}\label{counterne1}
&\left\|\1bf_{E}(\cdot+\xi)
\right\|_{\mathcal{K}_{\omega,
\0bf}^{p,q}(\rr)}\notag\\
&\quad=\left[\left\|\1bf_{E}
\1bf_{B(\xi,1)}\right\|_{L^p(\rr)}^q+
\sum_{k\in\mathbb{N}}2^{k\alpha q}\left\|
\1bf_{E}\1bf_{B(\xi,\tk)\setminus
B(\xi,\tkm)}\right\|_{L^p(\rr)}^q
\right]^{\frac{1}{q}}\notag\\
&\quad\leq\left[\left\|\1bf_{B(\xi,1)}
\right\|_{L^p(\rr)}^q+
\sum_{k\in\mathbb{N}}2^{k\alpha q}\left\|
\1bf_{E}\1bf_{B(\xi,\tk)}\right\|_{
L^p(\rr)}^q\right]^{\frac{1}{q}}\notag\\
&\quad\lesssim\left\|\1bf_{B(\xi,1)}
\right\|_{L^p(\rr)}+
\left[\sum_{k\in\mathbb{N}}2^{k\alpha q}\left\|
\1bf_{E}\1bf_{B(\xi,\tk)}
\right\|_{L^p(\rr)}^q\right]^{\frac{1}{q}}.
\end{align}
Applying this, we find that,
for any given $\alpha\in(-\infty,-\frac{1}{p})$
and for any $\xi\in\rr$,
\begin{align}\label{counterne2}
\left\|\1bf_{E}(\cdot+\xi)\right\|_{\mathcal{K}_{\omega,
\0bf}^{p,q}(\rr)}\lesssim1+\left(\sum_{k\in\mathbb{N}}
2^{k\alpha q}2^{\frac{k\alpha q}{p}}\right)^{\frac{1}{q}}\sim1.
\end{align}

On the other hand, by the proof
of Example \ref{conter}, we conclude that
$$
\mathrm{J}_{\xi,1}:\sim\left[
\sum_{k\in\mathbb{N}}2^{k\alpha q}\left\|
\1bf_{E}\1bf_{(\xi,\xi+\tk)}\right\|_{
L^p(\rr)}^q\right]^{\frac{1}{q}}
\lesssim1
$$
and
$$
\mathrm{J}_{\xi,2}:\sim\left[\sum_{
k\in\mathbb{N}}2^{k\alpha q}\left\|
\1bf_{E}\1bf_{(\xi-\tk,\xi)}\right\|_{
L^p(\rr)}^q\right]^{\frac{1}{q}}
\lesssim1
$$
in Example \ref{conter}
still hold true for $\alpha=-\frac{1}{p}$.
Thus, from Example \ref{conter} again, it
follows that, for any given
$\alpha\in[-\frac{1}{p},0)$ and
for any $\xi\in\rr$,
$\mathrm{J}_{\xi,1}\lesssim1$ and
$\mathrm{J}_{\xi,2}
\lesssim1$.
This, together with \eqref{counterne1},
further implies that, for any
given $\alpha\in[-\frac{1}{p},0)$ and
for any $\xi\in\rr$,
\begin{align*}
&\left\|\1bf_{E}(\cdot+\xi)
\right\|_{\mathcal{K}_{\omega,
\0bf}^{p,q}(\rr)}\\&\quad
\lesssim\left\|\1bf_{B(\xi,1)}\right\|_{L^p(\rr)}
+\left[\sum_{k\in\mathbb{N}}2^{k\alpha q}\left\|
\1bf_{E}\1bf_{(\xi,\xi+\tk)}
\right\|_{L^p(\rr)}^q\right]^{\frac{1}{q}}\\
&\qquad+\left[\sum_{k\in\mathbb{N}}
2^{k\alpha q}\left\|
\1bf_{E}\1bf_{(\xi-\tk,\xi)}\right\|_{
L^p(\rr)}^q\right]^{\frac{1}{q}}\\
&\quad\sim\left\|\1bf_{B(\xi,1)}\right\|_{L^p(\rr)}
+\mathrm{J_{\xi,1}}+\mathrm{J_{\xi,2}}\lesssim1.
\end{align*}
Combining this, \eqref{counterne2},
and Remark \ref{ighr}(i),
we further find that $\|\1bf_E\|_{
\mathcal{K}_{\omega}^{p,q}(\rr)}<\infty$,
and hence $\1bf_{E}\in\mathcal{K}_{\omega}^{p,q}(\rr)$.

We now show that the characteristic
function $\1bf_E$ does
not have an absolutely continuous
quasi-norm in the inhomogeneous
global generalized Herz space
$\mathcal{K}_{\omega}^{p,q}(\rr)$
under consideration.
Indeed, for any $\widetilde{k}\in\mathbb{N}$,
let $$F_{\widetilde{k}}:=\left(\widetilde{k},
\infty\right)\text{ and }
\xi_{\widetilde{k}}:=\widetilde{k}
-2+\widetilde{k}^{-\frac{2}{\alpha p}}.$$
Then, obviously, for any $x\in\rr$,
we have $\1bf_{F_{\widetilde{k}}}(x)\to0$
as $\widetilde{k}\to\infty$.
In addition,
from the assumption that, for any
$\widetilde{k}\in\mathbb{N}$,
$\widetilde{k}-1+\widetilde{k}^{-\frac{2}{
\alpha p}}\geq\widetilde{k}$,
it follows that, for any
$\widetilde{k}\in\mathbb{N}$,
$$\left(\widetilde{k}-1+\widetilde{k}^{
-\frac{2}{\alpha p}},
\widetilde{k}+\widetilde{k}^{
-\frac{2}{\alpha p}}\right)
\subset F_{\widetilde{k}}.$$
Using this, Remark \ref{ighr}(i),
and Definition \ref{igh}(i),
we further conclude that,
for any $\widetilde{k}\in\mathbb{N}$,
\begin{align*}
\left\|\1bf_{E}\1bf_{F_{\widetilde{k}}}
\right\|_{\mathcal{K}_{\omega}^{p,q}(\rr)}
&\geq\left\|\1bf_{(\widetilde{k}-1+\widetilde{k}
^{-\frac{2}{\alpha p}},
\widetilde{k}+\widetilde{k}^{-\frac{2}{\alpha p}})}
\right\|_{\mathcal{K}_{\omega}^{p,q}(\rr)}\\
&\geq\left\|\1bf_{(\widetilde{k}-1+
\widetilde{k}^{-\frac{2}{\alpha p}},
\widetilde{k}+\widetilde{k}^{-\frac{2}{\alpha p}})}
(\cdot+\xi_{\widetilde{k}})\right\|_{
\mathcal{K}_{\omega,\0bf}^{p,q}(\rr)}\\
&\geq2^{\alpha}\left\|\1bf_{(\widetilde{k}-1
+\widetilde{k}^{-\frac{2}{\alpha p}},
\widetilde{k}+\widetilde{k}^{-\frac{2}{\alpha p}})}
\1bf_{(\xi_{\widetilde{k}}+1,
\xi_{\widetilde{k}}+2)}\right\|_{L^{p}(\rr)}\\
&=2^{\alpha}\left|\left(\widetilde{k}-1+
\widetilde{k}^{-\frac{2}{\alpha p}},\widetilde{k}+
\widetilde{k}^{-\frac{2}{\alpha p}}\right)
\right|^{\frac{1}{p}}=2^{\alpha},
\end{align*}
which implies that $\1bf_{E}$ does not
have an absolutely continuous
quasi-norm in $\mathcal{K}_{
\omega}^{p,q}(\rr)$, and hence
the inhomogeneous global generalized Herz space
$\mathcal{K}_{\omega}^{p,q}(\rr)$ under consideration
does not have an absolutely continuous quasi-norm.
This finishes the proof of Example \ref{countern}.
\end{proof}

\begin{remark}
We should point out
that the assumption $\alpha\in(-\infty,0)$
in Example \ref{countern} is reasonable,
which means that, under this assumption,
the space $\NHerzS$ with $\omega(t):=t^{\alpha}$
for any $t\in[0,\infty)$ is a ball
quasi-Banach function space. Indeed,
combining Example \ref{ex113}
and Theorem \ref{Th2n}, we find that,
when $\omega(t):=t^\alpha$ for any $t\in(0,\infty)$
and for any given $\alpha\in(-\infty,0)$,
the inhomogeneous global generalized
Herz space $\mathcal{K}_{\omega}^{p,q}(\rr)$
is a ball quasi-Banach function
space.
\end{remark}

\subsection{Boundedness of Sublinear Operators
and \\ Fefferman--Stein Vector-Valued Inequalities}

In this subsection, we first
establish a boundedness criterion
of sublinear operators on inhomogeneous
generalized Herz spaces. As applications, we obtain
the boundedness of both the Hardy--Littlewood maximal operator
and Calder\'{o}n--Zygmund operators on these spaces immediately.
Then we establish the Fefferman--Stein
vector-valued inequalities
on inhomogeneous generalized Herz spaces.
The conclusions obtained in
this subsection play important
roles in the study of
inhomogeneous generalized Herz--Hardy spaces
in next chapter.

To begin with, we
show the following boundedness criterion
of sublinear operators
on both inhomogeneous local
and inhomogeneous global
generalized Herz spaces.

\begin{theorem}\label{vmaxhn}
Let $p\in(1,\infty),\
q\in(0,\infty)$, and $\omega\in M(\rp)$ satisfy
\begin{equation*}
-\frac{n}{p}<\mi(\omega)\leq\MI(\omega)<\frac{n}{p'},
\end{equation*}
where $\frac{1}{p}+\frac{1}{p'}=1$.
Assume that $T$ is a
sublinear operator bounded on $L^{p}(\rrn)$ and there
exists a positive constant $\widetilde{C}$ such that,
for any $f\in\NHerzSo$ and $x\notin\supp(f)$,
\begin{equation}\label{sizen}
|T(f)(x)|\leq\widetilde{C}\int_{\rrn}\frac{|f(y)|}{|x-y|^{n}}\,dy.
\end{equation}
Then there exists a
positive constant $C$, independent
of $f$, such that
\begin{enumerate}
  \item[{\rm(i)}] for any $f\in\NHerzSo$,
  \begin{equation}\label{vmaxhne0}
  \left\|T(f)\right\|_{\NHerzSo}\leq
  C\left[\widetilde{C}+
  \left\|T\right\|_{L^p(\rrn)
  \to L^p(\rrn)}\right]\|f\|_{\NHerzSo};
  \end{equation}
  \item[{\rm(ii)}] for any $f\in\NHerzS$,
  \begin{equation}\label{vmaxhne00}
  \left\|T(f)\right\|_{\NHerzS}\leq
  C\left[\widetilde{C}+
  \left\|T\right\|_{L^p(\rrn)
  \to L^p(\rrn)}\right]\|f\|_{\NHerzS}.
  \end{equation}
\end{enumerate}
\end{theorem}

\begin{proof}
Let all the symbols be as in the present theorem.
We first prove (i). To this end, let $f\in\NHerzSo$.
Then, by Definition \ref{igh}(i), we find that
\begin{align}\label{vmaxhne1}
&\left\|T(f)\right\|_{\NHerzSo}\notag\\&\quad\lesssim
\left\|T(f)\1bf_{B(\0bf,1)}\right\|_{L^p(\rrn)}\notag\\
&\qquad+\left\{\sum_{k\in\mathbb{N}}
\left[\omega(\tk)\right]^{q}\left\|
T(f)\1bf_{B(\0bf,\tk)\setminus
B(\0bf,\tkm)}\right\|_{L^p(\rrn)}^q
\right\}^{\frac{1}{q}}\notag\\
&\quad=:\mathrm{A}+\mathrm{B}.
\end{align}
We next deal with $\mathrm{A}$
and $\mathrm{B}$, respectively.
Indeed, applying
both the sublinearity and
the $L^p(\rrn)$ boundedness of the
operator $T$ and Definition \ref{igh}(i),
we conclude that
\begin{align}\label{vmaxhne2}
\mathrm{A}&
\lesssim\left\|T\left(f\1bf_{B(\0bf,2)}
\right)\1bf_{B(\0bf,1)}\right\|_{L^p(\rrn)}
+\left\|T\left(f\1bf_{[B(\0bf,2)]
^\complement}\right)\1bf_{B(\0bf,1)}\right\|_{L^p(\rrn)}
\notag\\&\lesssim\left\|T\right\|_{L^p(\rrn)\to
L^p(\rrn)}
\left\|f\1bf_{B(\0bf,2)}\right\|_{L^p(\rrn)}
+\left\|T\left(f\1bf_{[B(\0bf,2)]^\complement}\right)
\1bf_{B(\0bf,1)}\right\|_{L^p(\rrn)}\notag\\
&\lesssim\left\|T\right\|_{L^p(\rrn)\to
L^p(\rrn)}\left[
\left\|f\1bf_{B(\0bf,1)}\right\|_{L^p(\rrn)}
+\left\|f\1bf_{B(\0bf,2)\setminus
B(\0bf,1)}\right\|_{L^p(\rrn)}\right]\notag\\
&\quad+\left\|T\left(f\1bf_{[B(\0bf,2)]
^\complement}\right)
\1bf_{B(\0bf,1)}\right\|_{L^p(\rrn)}\notag\\
&\lesssim\left\|T\right\|_{L^p(\rrn)\to
L^p(\rrn)}\left\|f\right\|_{\NHerzSo}+
\left\|T\left(f\1bf_{[B(\0bf,2)]^\complement}\right)
\1bf_{B(\0bf,1)}\right\|_{L^p(\rrn)}.
\end{align}
Notice that, for any $k\in\mathbb{N}
\cap[3,\infty)$ and $x$, $y\in\rrn$
satisfying $|x|<1$ and $\tkm\leq|y|<\tk$, we have
\begin{align*}
|x-y|\geq|y|-|x|>2^{k-1}-1=\tkmd,
\end{align*}
which, combined with \eqref{sizen} and
the H\"{o}lder inequality,
implies that, for any $x\in\rrn$ satisfying $|x|<1$,
\begin{align}\label{vmaxhne3}
\left|T\left(f\1bf_{[B(\0bf,2)]
^\complement}\right)(x)\right|
&\leq\widetilde{C}
\int_{|y|\geq2}\frac{|f(y)|}{|x-y|^n}\,dy\notag\\
&\lesssim\widetilde{C}\sum_{k=3}^{\infty}2^{-nk}
\int_{\tkm\leq|y|<\tk}|f(y)|\,dy\notag\\
&\lesssim\widetilde{C}\sum_{k=3}^{\infty}
2^{-\frac{nk}{p}}\left\|f
\1bf_{B(\0bf,\tk)\setminus
B(\0bf,\tkm)}\right\|_{L^p(\rrn)}.
\end{align}
In addition, from Lemma \ref{Th1}, it follows that,
for any $k\in\mathbb{N}$,
$$\omega(\tk)\gtrsim2^{k[\mi(\omega)-\eps_1]},$$
where $\eps_1\in(0,\mi(\omega)
+\frac{n}{p})$ is a fixed positive constant.
Thus, using \eqref{vmaxhne3}, Lemma
\ref{l4320},
and the assumption $-\mi(\omega)-\frac{n}{p
}+\eps_1\in(-\infty,0)$, we find that,
for any $q\in(0,1]$,
\begin{align}\label{vmaxhne4}
&\left\|T\left(f\1bf_{[B(\0bf,2)]^\complement}
\right)\1bf_{B(\0bf,1)}\right\|_{L^p(\rrn)}
\notag\\&\quad\lesssim\widetilde{C}
\sum_{k=3}^{\infty}2^{
k[-\mi(\omega)-\frac{n}{p}+\eps_1]}
\omega(\tk)\left\|f\1bf_{B(\0bf,\tk)\setminus
B(\0bf,\tkm)}\right\|_{L^p(\rrn)}\notag\\
&\quad\lesssim\widetilde{C}
\left\{\sum_{k=3}^{\infty}
2^{k[-\mi(\omega)-\frac{n}{p}+\eps_1]
q}\left[\omega(\tk)\right]^q
\left\|f\1bf_{B(\0bf,\tk)\setminus
B(\0bf,\tkm)}\right\|_{L^p(\rrn)
}^q\right\}^{\frac{1}{q}}\notag\\
&\quad\lesssim\widetilde{C}\left\{\sum_{k\in
\mathbb{N}}\left[\omega(\tk)\right]^q
\left\|f\1bf_{B(\0bf,\tk)\setminus
B(\0bf,\tkm)}\right\|_{L^p(\rrn)
}^q\right\}^{\frac{1}{q}}
\lesssim\widetilde{C}\|f\|_{\NHerzSo}.
\end{align}
Moreover, by the assumption that, for any $k\in\mathbb{N}$,
$\omega(\tk)\gtrsim2^{k[\mi(\omega)-\eps_1]}$,
the H\"{o}lder inequality, and the assumption
$-\mi(\omega)-\frac{n}{p}+\eps_1\in(-\infty,0)$ again,
we conclude that, for any $q\in(1,\infty)$,
\begin{align*}
&\left\|T\left(f\1bf_{[B(\0bf,2)]^\complement}
\right)\1bf_{B(\0bf,1)}\right\|_{L^p(\rrn)}\\
&\quad\lesssim\widetilde{C}
\sum_{k=3}^{\infty}2^{k[-\mi
(\omega)-\frac{n}{p}+\eps_1]}
\omega(\tk)\left\|f\1bf_{B(\0bf,\tk)\setminus
B(\0bf,\tkm)}\right\|_{L^p(\rrn)}\\
&\quad\lesssim\widetilde{C}\left\{\sum_{k=3}^{\infty}
2^{k[-\mi(\omega)-\frac{n}{p}+\eps_1]q'}
\right\}^{\frac{1}{q'}}\\
&\qquad\times\left\{\sum_{k\in\mathbb{N}}
\left[\omega(\tk)\right]^q
\left\|f\1bf_{B(\0bf,\tk)\setminus
B(\0bf,\tkm)}\right\|_{
L^p(\rrn)}^q\right\}^{\frac{1}{q}}\\
&\quad\sim\widetilde{C}
\left\{\sum_{k\in\mathbb{N}}
\left[\omega(\tk)\right]^q
\left\|f\1bf_{B(\0bf,\tk)\setminus
B(\0bf,\tkm)}\right\|_{L^p
(\rrn)}^q\right\}^{\frac{1}{q}}
\lesssim\widetilde{C}\|f\|_{\NHerzSo}.
\end{align*}
Combining this, \eqref{vmaxhne2}, and \eqref{vmaxhne4},
we further obtain
\begin{align}\label{vmaxhne5}
\mathrm{A}&\lesssim\left[\widetilde{C}+
\left\|T\right\|_{L^p(\rrn)\to
L^p(\rrn)}\right]\|f\|_{\NHerzSo},
\end{align}
where the implicit positive
constant is independent of both $T$ and $f$,
which completes the estimate of $\mathrm{A}$.

We now deal with $\mathrm{B}$. To achieve this,
for any $k\in\mathbb{N}$,
let
$$
f_{k,0}:=f\1bf_{B(\0bf,2^{-2})},
$$
$$
f_{k,1}:=f\1bf_{B(\0bf,\tkmd)\setminus B(\0bf,2^{-2})},
$$
$$
f_{k,2}:=f\1bf_{B(\0bf,\tka)\setminus B(\0bf,\tkmd)},
$$
and
$$
f_{k,3}:=f\1bf_{[B(\0bf,\tka)]^\complement}.
$$
Then, from the sublinearity of $T$,
we deduce that
\begin{align}\label{vmaxhne6}
\mathrm{B}&\lesssim\left\{
\sum_{k\in\mathbb{N}}\left[\omega(\tk)\right]^q
\left\|T\left(f_{k,0}\right)\1bf_{B(\0bf,\tk)\setminus
B(\0bf,\tkm)}\right\|_{L^p(\rrn)}^q\right\}^{\frac{1}{q}}\notag\\
&\quad+\left\{\sum_{k\in\mathbb{N}}\left[\omega(\tk)\right]^q
\left\|T\left(f_{k,1}\right)\1bf_{B(\0bf,\tk)\setminus
B(\0bf,\tkm)}\right\|_{L^p(\rrn)}^q\right\}^{\frac{1}{q}}\notag\\
&\quad+\left\{\sum_{k\in\mathbb{N}}\left[\omega(\tk)\right]^q
\left\|T\left(f_{k,2}\right)\1bf_{B(\0bf,\tk)\setminus
B(\0bf,\tkm)}\right\|_{L^p(\rrn)}^q\right\}^{\frac{1}{q}}\notag\\
&\quad+\left\{\sum_{k\in\mathbb{N}}\left[\omega(\tk)\right]^q
\left\|T\left(f_{k,3}\right)\1bf_{B(\0bf,\tk)\setminus
B(\0bf,\tkm)}\right\|_{L^p(\rrn)}^q\right\}^{\frac{1}{q}}\notag\\
&=:\widetilde{\mathrm{J}_0}+\widetilde{\mathrm{J}_1}+
\widetilde{\mathrm{J}_2}+\widetilde{\mathrm{J}_3}.
\end{align}
Then we estimate $\widetilde{\mathrm{J}_0}$,
$\widetilde{\mathrm{J}_1}$, $\widetilde{\mathrm{J}_2}$,
and $\widetilde{\mathrm{J}_3}$, respectively.
Indeed, using Lemma \ref{Th1}, we find that,
for any $k\in\mathbb{N}$,
$$\omega(\tk)\lesssim2^{k[\MI(\omega)+\eps_2]},$$
where $\eps_2\in(0,-\MI(\omega)+\frac{n}{p'})$
is a fixed positive constant.
Moreover, for any $k\in\mathbb{N}$
and $x$, $y\in\rrn$ satisfying
$\tkm\leq|x|<\tk$ and $|y|<2^{-2}$,
we have
\begin{align*}
|x-y|\geq|x|-|y|>\tkm-2^{-2}>\tkmd,
\end{align*}
which, together with \eqref{sizen} and the H\"{o}lder
inequality, further implies that,
for any $k\in\mathbb{N}$ and
$x\in\rrn$ satisfying $\tkm\leq|x|<\tk$,
\begin{align*}
\left|T\left(f_{k,0}\right)(x)\right|
\leq\widetilde{C}
\int_{|y|<2^{-2}}\frac{|f(y)|}{|x-y|^n}\,dy
\lesssim\widetilde{C}2^{-nk}\left\|f
\1bf_{B(\0bf,1)}\right\|_{L^p(\rrn)}.
\end{align*}
From this, the assumptions that,
for any $k\in\mathbb{N}$,
$\omega(\tk)\lesssim2^{k[\MI(\omega)+\eps_2]}$ and
$\MI(\omega)-\frac{n}{p'}+\eps_
2\in(-\infty,0)$, and Definition \ref{igh}(i),
it follows that
\begin{align}\label{vmaxhne7}
\widetilde{\mathrm{J}_0}&\lesssim
\widetilde{C}\left\{\sum_{k\in\mathbb{N}}
2^{k[\MI(\omega)+\eps_2]q}2^{-nkq}\left\|
f\1bf_{B(\0bf,1)}\right\|_{L^p(\rrn)
}^q2^{\frac{nkq}{p}}\right\}^{\frac{1}{q}}\notag\\
&\sim\widetilde{C}\left\|f\1bf_{B(\0bf,1)}
\right\|_{L^p(\rrn)}
\left\{\sum_{k\in\mathbb{N}}2^{k[\MI(\omega)
-\frac{n}{p'}+\eps_2]q}
\right\}^{\frac{1}{q}}\notag\\
&\sim\widetilde{C}\left\|f\1bf_{
B(\0bf,1)}\right\|_{L^p(\rrn)}
\lesssim\widetilde{C}\|f\|_{\NHerzSo}.
\end{align}
This finishes the estimate of $\widetilde{\mathrm{J}_0}$.

In addition, by Lemma \ref{Th1},
we conclude that, for any $0<t<\tau<\infty$,
\begin{align*}
\left(\frac{\tau}{t}\right)^{\mi(\omega)-\eps_3}
\lesssim\frac{\omega(\tau)}{\omega(t)}\lesssim\left(
\frac{\tau}{t}\right)^{\MI(\omega)+\eps_3},
\end{align*}
where
$$\eps_3\in\left(0,\min\left\{\mi(\omega)+\frac{n}{p},
-\MI(\omega)+\frac{n}{p'}\right\}\right)$$
is a fixed positive constant.
From this, repeating an argument similar
to that used in the
estimations of $\mathrm{J}_1$, $\mathrm{J}_2$,
and $\mathrm{J}_3$ in the proof of Theorem \ref{vmaxh}
with $\mathbb{Z}$ therein replaced
by $\mathbb{N}$, and
Definition \ref{igh}(i), we deduce that
\begin{equation*}
\widetilde{\mathrm{J}_1}\lesssim
\widetilde{C}\|f\|_{\NHerzSo},
\end{equation*}
\begin{equation*}
\widetilde{\mathrm{J}_2}\lesssim
\left\|T\right\|_{L^p(\rrn)\to
L^p(\rrn)}\|f\|_{\NHerzSo},
\end{equation*}
and
\begin{equation*}
\widetilde{\mathrm{J}_3}
\lesssim\widetilde{C}\|f\|_{\NHerzSo}.
\end{equation*}
Combining these, \eqref{vmaxhne6}, \eqref{vmaxhne7},
we further find that
\begin{equation}\label{vmaxhne10}
\mathrm{B}\lesssim\left[\widetilde{C}+
\left\|T\right\|_{L^p(\rrn)\to
L^p(\rrn)}\right]\|f\|_{\NHerzSo},
\end{equation}
where the implicit positive constant
is independent of both $T$ and $f$.
This then finishes the estimate of $\mathrm{B}$.
Therefore, applying \eqref{vmaxhne1}, \eqref{vmaxhne5},
and \eqref{vmaxhne10}, we obtain
$$
\|T(f)\|_{\NHerzSo}\lesssim
\left[\widetilde{C}+
\left\|T\right\|_{L^p(\rrn)\to
L^p(\rrn)}\right]\|f\|_{\NHerzSo},
$$
where the implicit positive constant
is independent of both $T$ and $f$.
This implies that \eqref{vmaxhne0} holds true.

We next show \eqref{vmaxhne00}. Indeed,
for any given $f\in\NHerzS$ and for any $\xi\in\rrn$,
from Remark \ref{ighr}(i), it follows that
$\|f(\cdot+\xi)\|_{\NHerzSo}<\infty$ and hence
$f(\cdot+\xi)\in\NHerzSo$. Thus, by
\eqref{vmaxhne0}, we conclude that
$$
\left\|T(f)(\cdot+\xi)\right\|_{\NHerzSo}
\lesssim\left[\widetilde{C}+
\left\|T\right\|_{L^p(\rrn)\to
L^p(\rrn)}\right]
\left\|f(\cdot+\xi)\right\|_{\NHerzSo}.
$$
By this, the arbitrariness of
$\xi$, and Remark \ref{ighr}(i) again,
we find that
$$
\|T(f)\|_{\NHerzS}\lesssim
\left[\widetilde{C}+
\left\|T\right\|_{L^p(\rrn)\to
L^p(\rrn)}\right]\|f\|_{\NHerzS}
$$
with the implicit positive constant
independent of $f$.
This finishes the proof of
\eqref{vmaxhne00} and hence
of Theorem \ref{vmaxhn}.
\end{proof}

\begin{remark}
We should point out that,
in Theorem \ref{vmaxhn}, if $\omega(t):=t^\alpha$
for any $t\in(0,\infty)$ and for
any given $\alpha\in\rr$, then Theorem \ref{vmaxhn}(i)
coincides with the result for
classical inhomogeneous Herz spaces
obtained in \cite[Corollary 2.1]{LiY}.
\end{remark}

Applying Theorem \ref{vmaxhn} and the facts that
the Hardy--Littlewood maximal operator satisfies
the condition \eqref{sizen}
(see, for instance, \cite[Remark 4.4]{HSamko})
and is bounded on the
Lebesgue space $L^p(\rrn)$ when $p\in(1,\infty)$
(see, for instance, \cite[Theorem 2.1.6]{LGCF}),
we immediately obtain the following two conclusions
which give the boundedness
of the Hardy--Littlewood maximal operator,
respectively, on inhomogeneous local and
inhomogeneous global generalized
Herz spaces; we omit the details.

\begin{corollary}\label{bhlon}
Let $p$, $q$, and $\omega$ be as
in Theorem \ref{vmaxhn} and $\mc$
the Hardy--Littlewood maximal operator
as in \eqref{hlmax}. Then there exists
a positive constant $C$ such that,
for any $f\in L^1_{\mathrm{loc}}(\rrn)$,
$$
\left\|\mc(f)\right\|_{\NHerzSo}\leq
C\|f\|_{\NHerzSo}.
$$
\end{corollary}

\begin{remark}
We should point out that, in Corollary \ref{bhlon},
if $\omega(t):=t^\alpha$ for any $t\in(0,\infty)$
and for any given $\alpha\in\rr$,
the conclusion obtained in this corollary goes back to
\cite[Theorem 4.1]{IM}.
\end{remark}

\begin{corollary}\label{Th4n}
Let $p$, $q$, and $\omega$ be as
in Theorem \ref{vmaxhn} and $\mc$
the Hardy--Littlewood maximal operator
as in \eqref{hlmax}. Then there exists
a positive constant $C$ such that,
for any $f\in L^1_{\mathrm{loc}}(\rrn)$,
$$
\left\|\mc(f)\right\|_{\NHerzS}\leq
C\|f\|_{\NHerzS}.
$$
\end{corollary}

Furthermore, let $d\in\zp$. Recall that
the $d$-order Calder\'{o}n--Zygmund operator
is defined as in Definition \ref{defin-C-Z-s}.
Then, repeating an argument
similar to that used in the proof of
Corollary \ref{czoh} with Theorem \ref{vmaxh}
therein replaced by Theorem \ref{vmaxhn},
we obtain the following two
boundedness criteria of $d$-order
Calder\'{o}n--Zygmund operators,
respectively, on inhomogeneous local and
inhomogeneous global
generalized Herz spaces; we omit the details.

\begin{corollary}\label{czohn}
Let $p$, $q$, and $\omega$ be as
in Theorem \ref{vmaxhn}, $d\in\zp$,
$\delta\in(0,1]$, $K$ be a $d$-order standard
kernel as in Definition
\ref{def-s-k} with the above $\delta$,
and $T$ a $d$-order Calder\'{o}n--Zygmund
operator as in Definition \ref{defin-C-Z-s}
with kernel $K$.
Then $T$ is well defined on $\NHerzSo$ and
there exists a positive
constant $C$ such that,
for any $f\in\NHerzSo$,
$$
\left\|T(f)\right\|_{\NHerzSo}\leq
C\|f\|_{\NHerzSo}.
$$
\end{corollary}

\begin{remark}
In Corollary \ref{czohn},
for any given $\alpha\in\rr$
and for any $t\in(0,\infty)$,
let $\omega(t):=t^\alpha$.
Then we point out that, in this case,
Corollary \ref{czohn} goes back to \cite[Remark 5.1.1]{LYH}.
\end{remark}

\begin{corollary}\label{czohgn}
Let $p$, $q$, and $\omega$ be as
in Theorem \ref{vmaxhn}, $d\in\zp$,
$\delta\in(0,1]$, $K$ be a $d$-order standard
kernel as in Definition
\ref{def-s-k} with the above $\delta$,
and $T$ a $d$-order Calder\'{o}n--Zygmund
operator as in Definition \ref{defin-C-Z-s}
with kernel $K$.
Then $T$ is well defined on $\NHerzS$
and there exists a positive constant $C$ such that,
for any $f\in\NHerzS$,
$$
\left\|T(f)\right\|_{\NHerzS}\leq
C\|f\|_{\NHerzS}.
$$
\end{corollary}

Now, we turn to investigate
the Fefferman--Stein vector-valued inequality
on inhomogeneous generalized
Herz spaces. Indeed, repeating an
argument similar to that used in the
proof of Theorem \ref{Th3.4} via replacing
Theorem \ref{vmaxh}(i)
therein by Theorem \ref{vmaxhn}(i),
we conclude the following Fefferman--Stein
vector-valued inequality
on inhomogeneous local generalized Herz spaces;
we omit the details.\index{Fefferman--Stein vector-valued \\inequality}

\begin{theorem}\label{Th3.4n}
Let $p,\ r\in(1,\infty)$, $q\in(0,\infty)$, and
$\omega\in M(\rp)$ satisfy
$$-\frac{n}{p}<\mi(\omega)\leq\MI(\omega)<\frac{n}{p'},$$
where $\frac{1}{p}+\frac{1}{p'}=1$.
Then there exists a positive constant $C$ such that,
for any $\{f_{j}\}_{j\in\mathbb{N}}\subset L^1_{\mathrm{loc}}(\rrn)$,
$$
\left\|\left\{\sum_{j\in\mathbb{N}}\left[\mc
(f_{j})\right]^{r}\right\}^{\frac{1}{r}}\right\|_{\NHerzSo}
\leq C\left\|\left(
\sum_{j\in\mathbb{N}}|f_{j}|^{r}\right)^{\frac{1}{r}}\right\|_{\NHerzSo}.
$$
\end{theorem}

\begin{remark}
We should point out that,
in Theorem \ref{Th3.4n}, when $\omega(t):=t^\alpha$
for any $t\in(0,\infty)$ and for any given $\alpha\in\rr$,
then Theorem \ref{Th3.4n} coincides with the known
Fefferman--Stein vector-valued inequality
on the classical inhomogeneous Herz
spaces obtained in \cite[Corollary 4.5]{i10}.
\end{remark}

Similarly, we also have the following
Fefferman--Stein vector-valued inequality on the
inhomogeneous global generalized Herz space $\NHerzS$.
The proof of this theorem is just to
repeat the proof of Theorem \ref{Th3.4}
with Theorem \ref{vmaxh}(i) therein replaced
by Theorem \ref{vmaxhn}(ii); we omit the
details.\index{Fefferman--Stein vector-valued \\inequality}

\begin{theorem}\label{Th3.3n}
Let $p,\ r\in(1,\infty)$, $q\in(0,\infty)$, and
$\omega\in M(\rp)$ satisfy
$$-\frac{n}{p}<\mi(\omega)\leq\MI(\omega)<\frac{n}{p'},$$
where $\frac{1}{p}+\frac{1}{p'}=1$.
Then there exists a positive constant $C$
such that, for any $\{f_{j}\}_{j\in\mathbb{N}}
\subset L_{{\rm loc}}^{1}(\rrn)$,
$$
\left\|\left\{\sum_{j\in\mathbb{N}}\left[\mc
(f_{j})\right]^{r}\right\}^{\frac{1}{r}}\right\|_{\NHerzS}
\leq C\left\|\left(
\sum_{j\in\mathbb{N}}|f_{j}|^{r}
\right)^{\frac{1}{r}}\right\|_{\NHerzS}.
$$
\end{theorem}

\subsection{Dual and Associate
Spaces of Inhomogeneous
Local Generalized Herz Spaces}

The target of this subsection is to
investigate dual spaces and
associate spaces
of inhomogeneous local generalized
Herz spaces. To this end, we first
introduce the concept of the
inhomogeneous local generalized Herz
space $\HerzSx$ as follows.

\begin{definition}\label{ighx}
Let $p$, $q\in(0,\infty)$, $\omega\in M(\rp)$,
and $\xi\in\rrn$.
Then the \emph{inhomogeneous
local generalized Herz space} $\NHerzSx$
is defined to be the set of all the measurable
functions $f$ on $\rrn$ such that
\begin{align*}
\|f\|_{\NHerzSx}:=
\left\{\left\|
f\1bf_{B(\xi,1)}\right\|_{L^{p}(\rrn)}^{q}
+\sum_{k\in\mathbb{N}}\left[
\omega(\tk)\right]^{q}\left\|f\1bf_{B(\xi,\tk)\setminus
B(\xi,\tkm)}\right\|_{L^{p}(\rrn)}
^{q}\right\}^{\frac{1}{q}}
\end{align*}
is finite.
\end{definition}

Next, we establish the following
dual theorem of inhomogeneous local
generalized Herz spaces $\HerzSx$ when
$p$, $q\in(1,\infty)$.

\begin{theorem}\label{dualn}
Let $p,\ q\in(1,\infty)$, $\omega\in M(\rp)$,
and $\xi\in\rrn$.
Then the dual space of $\NHerzSx$,
denoted by $(\NHerzSx)^*$,
is $\NHerzSxd$ in the following sense:
\begin{enumerate}
 \item[{\rm(i)}]
 Let $g\in\NHerzSxd$. Then the
 linear functional
 \begin{equation}\label{funldn}
 \phi_{g}:\ f\mapsto\phi_{g}(f):=
 \int_{\rrn}f(y)g(y)\,dy,
 \end{equation}
 defined for any $f\in\NHerzSx$,
 is bounded on $\NHerzSx$.
 \item[{\rm(ii)}]
 Conversely, any continuous linear functional
 on $\NHerzSx$ arises as in \eqref{funldn}
 with a unique $g\in\NHerzSxd$.
\end{enumerate}
Moreover, for any $g\in\NHerzSxd$,
$$\left\|g\right\|_{
\NHerzSxd}=\left\|\phi_{g}
\right\|_{(\NHerzSx)^*}.$$
\end{theorem}

\begin{proof}
Let all the symbols be as in
the present theorem.
We first show (i). Indeed, from
the Tonelli theorem
and the H\"{o}lder inequality,
we deduce that, for any
$f\in\NHerzSx$ and $g\in\NHerzSxd$,
\begin{align}\label{dualne0}
&\left\|fg\right\|_{L^1(\rrn)}\notag\\&\quad=
\int_{B(\xi,1)}|f(y)g(y)|\,dy+\sum_{k\in\mathbb{N}}
\int_{B(\xi,\tk)\setminus
B(\xi,\tkm)}|f(y)g(y)|\,dy\notag\\
&\quad\leq\left\|f\1bf_{B(\xi,1)}\right\|_{L^{p}(\rrn)}
\left\|g\1bf_{B(\xi,1)}\right\|_{L^{p'}
(\rrn)}\notag\\&\qquad+\sum_{k\in\mathbb{N}}
\left\|f\1bf_{B(\xi,\tk)\setminus
B(\xi,\tkm)}\right\|_{L^{p}(\rrn)}
\left\|g\1bf_{B(\xi,\tk)\setminus
B(\xi,\tkm)}\right\|_{L^{p'}(\rrn)}\notag\\
&\quad\leq\left\{\left\|f\1bf_{B(\xi,1)}
\right\|_{L^{p}(\rrn)}^{q}
+\sum_{k\in\mathbb{N}}\left[\omega(\tk)\right]^{q}
\left\|f\1bf_{B(\xi,\tk)\setminus
B(\xi,\tkm)}\right\|_{L^{p}(\rrn)}^{q}
\right\}^{\frac{1}{q}}\notag\\
&\qquad\times\left\{\left\|g
\1bf_{B(\xi,1)}\right\|_{L^{p'}(\rrn)}^{q'}
+\sum_{k\in\mathbb{N}}\left[\omega(\tk)\right]^{-q'}
\left\|g\1bf_{B(\xi,\tk)\setminus
B(\xi,\tkm)}\right\|_{L^{p'}(\rrn)}^{q'}
\right\}^{\frac{1}{q'}}\notag\\
&\quad=\|f\|_{\NHerzSx}\|g\|_{\NHerzSxd}.
\end{align}
This then implies that the linear functional $\phi_g$
as in \eqref{funldn} is bounded on $\NHerzSx$ and
\begin{equation}\label{dualne1}
\left\|\phi_g\right\|_{(\NHerzSx)^*}\leq
\left\|g\right\|_{\NHerzSxd},
\end{equation}
which completes the proof of (i).

Conversely, we show (ii).
To this end, let $\phi\in(\NHerzSx)^*$,
$A_{\xi,0}:=B(\xi,1)$, and, for any
$k\in\mathbb{N}$,
$$
A_{\xi,k}:=B(\xi,\tk)\setminus
B(\xi,\tkm).
$$
Then, by Definition \ref{ighx},
we find that, for any given
$k\in\zp$ and for any $f\in
L^p(A_{\xi,k})$,
\begin{align*}
\|f\|_{\NHerzSx}&\sim
\left\|f\1bf_{A_{\xi,k}}
\right\|_{L^p(\rrn)}
\sim\left\|f\right\|_{L^p(A_{\xi,k})}
<\infty,
\end{align*}
which further implies that
\begin{align*}
\left|\phi(f)\right|&\leq
\left\|\phi\right\|_{(\NHerzSx)^*}
\|f\|_{\NHerzSx}\sim
\left\|\phi\right\|_{(\NHerzSx)^*}
\left\|f\right\|_{
L^p(A_{\xi,k})}.
\end{align*}
Thus, $\phi\in(L^p(A_{\xi,k}))^*$
for any $k\in\zp$.
Combining this and the Riesz
representation theorem
(see \cite[Theorem 4.11]{b11}),
we conclude that,
for any $k\in\mathbb{Z}_+$,
there exists a unique $g_k\in L^{p'}(A_{\xi,k})$
such that, for any
$f\in L^p(A_{\xi,k})$,
\begin{equation}\label{dualne2}
\phi(f)=\int_{A_{\xi,k}}
f(y)g_k(y)\,dy.
\end{equation}
Let $$
g:=\sum_{k\in\zp}g_k.
$$
Then, from \eqref{dualne2} and an
argument similar to that
used in the proof of
Theorem \ref{dual} with
$\{g_k\}_{k\in\mathbb{Z}}$
and $\{B(\xi,\tk)\setminus
B(\xi,\tkm)\}_{k\in\mathbb{Z}}$
therein replaced, respectively,
by $\{g_k\}_{k\in\zp}$ and
$\{A_{\xi,k}\}_{k\in\zp}$ here,
it follows that
$g\in\NHerzSxd$,
$$
\left\|g\right\|_{\NHerzSxd}
\leq\left\|\phi\right\|_{(\NHerzSx)^*},
$$
and, for any $f\in\NHerzSx$,
$$
\phi(f)=\int_{\rrn}f(y)g(y)\,dy.
$$
These, together with (i) and \eqref{dualne1},
further imply that the
linear functional $\phi_g$ as in
\eqref{funldn}
is bounded on $\NHerzSx$, $\phi=\phi_g$,
and \begin{align*}
\|\phi\|_{(\NHerzSx)^*}=\|\phi_g\|_{(\NHerzSx)^*}
\leq\|g\|_{\NHerzSxd}\leq\|\phi\|_{(\NHerzSx)^*}.
\end{align*}
Therefore, we have
$$\left\|\phi_g\right\|_{(
\NHerzSx)^*}=\left\|g\right\|_{\NHerzSx}.$$
Using this and the linearity of $\phi_g$ about $g$,
we further find that $g$ is unique.
This finishes the proof of (ii),
and hence of Theorem \ref{dualn}.
\end{proof}

Using the above dual theorem,
we now prove the following
conclusion which shows that,
when $p$, $q\in(1,\infty)$,
the associate space of the inhomogeneous
local generalized Herz space $\HerzSo$
is just $\Kmc_{1/\omega,
\0bf}^{p',q'}(\rrn)$.

\begin{theorem}\label{associaten}
Let $p,\ q\in(1,\infty)$ and
$\omega\in M(\rp)$.
Then
$$
\left(\NHerzSo\right)'=\Kmc_{1/\omega,
\0bf}^{p',q'}(\rrn)
$$
with the same norms,
where $(\NHerzSo)'$ denotes
the associate space
of the inhomogeneous
local generalized Herz space
$\NHerzSo$.
\end{theorem}

\begin{proof}
Let all the symbols be as in the
present theorem. Then, applying an argument
similar to that used in the
proof of Theorem \ref{associate}
with Theorems \ref{Th3} and \ref{abso}
therein replaced, respectively,
by Theorems \ref{Th3n} and \ref{abson}, we
find that
$$\left(\NHerzSo\right)'=\Kmc_{1/\omega,
\0bf}^{p',q'}(\rrn)$$
and, for any $f\in(\NHerzSo)'$,
$$
\|f\|_{(\NHerzSo)'}=\|f\|_{\Kmc_{1/\omega,
\0bf}^{p',q'}(\rrn)}.
$$
This finishes the proof of
Theorem \ref{associaten}.
\end{proof}

\subsection{Extrapolation Theorems}

In this subsection, we establish the
extrapolation theorems of inhomogeneous
generalized Herz spaces, which
plays an essential role
in the study of the boundedness
and the compactness characterizations
of commutators on these Herz spaces.

We first give the following extrapolation
theorem of inhomogeneous local generalized
Herz spaces.

\begin{theorem}\label{th175n}
Let $p$, $q\in(1,\infty)$,
$r_0\in[1,\infty)$, and $\omega\in M(\rp)$
satisfy
$$
-\frac{n}{p}<\mi(\omega)
\leq\MI(\omega)<\frac{n}{p'}.
$$
Assume that $\mathcal{F}$ is a set of
all pairs of nonnegative
measurable functions $(F,\,G)$ such that,
for any given
$\upsilon\in A_{r_0}(\rrn)$,
$$
\int_{\rrn}\left[F(x)\right]^{r_0}\upsilon(x)\,dx
\leq C_{(r,[\upsilon]_{A_{r_0}(\rrn)})}
\int_{\rrn}\left[G(x)\right]^{r_0}\upsilon(x)\,dx,
$$
where the positive constant
$C_{(r_0,[\upsilon]_{A_{r_0}(\rrn)})}$
is independent of $(F,\,G)$, but depending on
both $r_0$ and $[\upsilon]_{A_{r_0}(\rrn)}$.
Then there exists a
positive constant $C$ such that, for any
$(F,\,G)\in\mathcal{F}$ with
$\|F\|_{\NHerzSo}<\infty$,
$$
\left\|F\right\|_{\NHerzSo}\leq
C\left\|G\right\|_{\NHerzSo}.
$$
\end{theorem}

To show this theorem, we
require two technique
lemmas. The following one is related
to the boundedness of the Hardy--Littlewood maximal
operator on inhomogeneous local generalized
Herz spaces.

\begin{lemma}\label{mbhln}
Let $p,\ q\in(0,\infty)$ and $\omega\in M(\rp)$ satisfy
$\mi(\omega)\in(-\frac{n}{p},\infty)$.
Then, for any given
$$r\in\left(0,\min\left\{p,\frac{n}
{\MI(\omega)+n/p}\right\}\right),$$
there exists a positive constant $C$ such that,
for any $f\in L^1_{\mathrm{loc}}(\rrn)$,
$$
\left\|\mc(f)\right\|_{[\NHerzSo]^{1/r}}\leq
C\|f\|_{[\NHerzSo]^{1/r}}.
$$
\end{lemma}

\begin{proof}
Let all the symbols be as in the
present lemma and
$$r\in\left(0,\min\left\{p,\frac{n}
{\MI(\omega)+n/p}\right\}\right).$$
Then, from an argument similar to
that used in the proof of Lemma \ref{mbhl} with Lemma
\ref{convexll} and Corollary \ref{bhlo}
therein replaced, respectively,
by Lemma \ref{convexlln} and Corollary \ref{bhlon},
it follows that, for any $f\in L^1_{\mathrm{loc}}(\rrn)$,
$$
\left\|\mc(f)\right\|_{[\NHerzSo]^{1/r}}\lesssim
\|f\|_{[\NHerzSo]^{1/r}}.
$$
This then finishes the proof of Lemma \ref{mbhln}.
\end{proof}

On the other hand,
we need the following
conclusion about the boundedness of
the Hardy--Littlewood maximal operator
on associate spaces of inhomogeneous
local generalized Herz spaces.

\begin{lemma}\label{mbhaln}
Let $p,\ q\in(0,\infty)$ and $\omega\in M(\rp)$ satisfy
$\mi(\omega)\in(-\frac{n}{p},\infty)$.
Then, for any given $$s\in\left(0,\min\left\{p,q,
\frac{n}{\MI(\omega)+n/p}\right\}\right)$$
and $$r\in\left(\max\left\{p,\frac{n}{
\mi(\omega)+n/p}\right\},\infty\right],$$
the Herz space $[\NHerzSo]^{1/s}$
is a ball Banach function space and there exists a positive
constant $C$ such that, for any $f\in
L_{{\rm loc}}^{1}(\rrn)$,
\begin{equation*}
\left\|\mc^{((r/s)')}(f)
\right\|_{([\NHerzSo]^{1/s})'}\leq C
\left\|f\right\|_{([\NHerzSo]^{1/s})'}.
\end{equation*}
\end{lemma}

\begin{proof}
Let all the symbols be as in
the present lemma, $$s\in\left(0,\min\left\{p,q,
\frac{n}{\MI(\omega)+n/p}\right\}\right),$$
and $$r\in\left(\max\left\{p,\frac{n}{
\mi(\omega)+n/p}\right\},\infty\right].$$
Then, by the proof of
Lemma \ref{mbhal} via replacing
Lemma \ref{convexll}, Corollary
\ref{bhlo}, and Theorem \ref{associate}
therein, respectively, by Lemma
\ref{convexlln}, Corollary
\ref{bhlon}, and Theorem \ref{associaten},
we conclude that the Herz space $[\NHerzSo]^{1/s}$
is a BBF space and,
for any $f\in L_{{\rm loc}}^{1}(\rrn)$,
\begin{equation*}
\left\|\mc^{((r/s)')}(f)
\right\|_{([\NHerzSo]^{1/s})'}\lesssim
\left\|f\right\|_{([\NHerzSo]^{1/s})'},
\end{equation*}
which completes the proof of Lemma \ref{mbhaln}.
\end{proof}

Via the above two lemmas and
the known extrapolation theorem of
ball Banach function spaces obtained
by Tao et al.\ \cite[Lemma 2.13]{TYYZ}
(see also Lemma \ref{th175l1} above),
we next prove Theorem \ref{th175n}.

\begin{proof}[Proof of Theorem \ref{th175n}]
Let all the symbols be as in
the present theorem. Then, repeating the
proof of Theorem \ref{th175} with
Theorem \ref{balll} and Lemmas
\ref{mbhl} and \ref{mbhal} used
therein replaced, respectively,
by Theorem \ref{ballln} and Lemmas
\ref{mbhln} and \ref{mbhaln}, we
find that the inhomogeneous
local generalized Herz space $\NHerzSo$
under consideration is a BBF space and
the Hardy--Littlewood maximal operator
$\mc$ as in \eqref{hlmax} is bounded
on both $\NHerzSo$ and $(\NHerzSo)'$.
From these and Lemma \ref{th175l1} with
$X:=\NHerzSo$, we deduce that, for any
$(F,\,G)\in\mathcal{F}$ and
$\|F\|_{\NHerzSo}<\infty$,
$$
\left\|F\right\|_{\NHerzSo}
\lesssim\left\|G\right\|_{\NHerzSo}.
$$
Thus, the proof
of Theorem \ref{th175n} is completed.
\end{proof}

Finally, we establish the extrapolation
theorem of inhomogeneous global
generalized Herz spaces as follows.

\begin{theorem}\label{extragn}
Let $p$, $q\in(1,\infty)$,
$r_0\in[1,\infty)$, and $\omega\in M(\rp)$
satisfy
$$
-\frac{n}{p}<\mi(\omega)
\leq\MI(\omega)<\frac{n}{p'}.
$$
Assume that $\mathcal{F}$ is a set of
all pairs of nonnegative
measurable functions $(F,\,G)$ such that,
for any given
$\upsilon\in A_{r_0}(\rrn)$,
\begin{equation*}
\int_{\rrn}\left[F(x)\right]^{r_0}\upsilon(x)\,dx
\leq C_{(r,[\upsilon]_{A_{r_0}(\rrn)})}
\int_{\rrn}\left[G(x)\right]^{r_0}\upsilon(x)\,dx,
\end{equation*}
where the positive constant
$C_{(r_0,[\upsilon]_{A_{r_0}(\rrn)})}$
is independent of $(F,\,G)$, but depending on
both $r_0$ and $[\upsilon]_{A_{r_0}(\rrn)}$.
Then there exists a
positive constant $C$ such that, for any
$(F,\,G)\in\mathcal{F}$ with
$\|F\|_{\NHerzS}<\infty$,
$$
\left\|F\right\|_{\NHerzS}\leq
C\left\|G\right\|_{\NHerzS}.
$$
\end{theorem}

\begin{proof}
Let all the symbols be as in the
present theorem and $(F,\,G)\in\mathcal{F}$
with $\|F\|_{\NHerzS}<\infty$.
Then, repeating the proof of Theorem
\ref{extrag} via replacing Remark \ref{remhs}(ii)
and Theorem \ref{th175} therein,
respectively, by Remark \ref{ighr}(i)
and Theorem \ref{th175n}, we obtain
$$
\left\|F\right\|_{\NHerzS}
\lesssim\left\|G\right\|_{\NHerzS}.
$$
This finishes the proof of Theorem
\ref{extragn}.
\end{proof}

\section{Inhomogeneous Block Spaces and Their Applications}

Let $p$, $q\in(0,\infty)$ and $\omega\in M(\rp)$.
In this section, we first introduce the
inhomogeneous block spaces
$\Nbspace$ and establish the duality between these
spaces and inhomogeneous
global generalized Herz spaces. As applications,
we obtain the boundedness of both sublinear operators
and Calder\'{o}n--Zygmund
operators on inhomogeneous block spaces.
In particular, the boundedness
of powered Hardy--Littlewood
maximal operators on inhomogeneous block spaces can
be deduced directly
from the boundedness of sublinear operators,
which plays an important role
in the next chapter.

\subsection{Inhomogeneous Block Spaces}

In this subsection, we introduce the
inhomogeneous counterparts of block spaces
as in Chapter \ref{sec4}.
To this end, we first present some necessary
concepts. For any $k\in\mathbb{Z}$,
we use $\mathcal{D}_k$ as in \eqref{dyadic-k}
to denote the set of all the standard dyadic cubes
on $\rrn$ of level $-k$. Moreover,
for any $k\in\mathbb{N}$ and $\xi\in\rrn$,
let
$$
\mathscr{Q}_k^{(\0bf)}:=
\mathcal{D}_{k-1}\cap\left[
Q(\0bf,2^{k+1})\setminus
Q(\0bf,2^k)\right],
$$
$$
\mathscr{Q}_{k}^{(\xi)}
:=\left\{Q\in\mathcal{Q}:\ Q-\{\xi\}\in
\mathscr{Q}_k^{(\0bf)}\right\},
$$
and $\mathscr{Q}_{0}^{(\xi)}:=\{Q(\xi,2)\}$.

\begin{remark}\label{remark215n}
Obviously, for any $\xi\in\rrn$
and $k\in\mathbb{N}$, $$\mathscr{Q}_{k}
^{(k)}=\mathcal{Q}_{k}^{(\xi)},$$
where $\mathcal{Q}_{k}^{(\xi)}$ is defined
as in \eqref{qmc}.
Combining this and Remark \ref{remark215},
we further find that, for any $\xi\in\rrn$
and $k\in\mathbb{N}$,
$$\sharp\mathscr{Q}_k^{(\xi)}=2^{2n}-2^n.$$
\end{remark}

Recall that the concept of $(\omega,\,p)$-blocks
is introduced in Definition \ref{blo}.
Then we give the definition of inhomogeneous
block spaces via $(\omega,\,p)$-blocks as follows.

\begin{definition}\label{blosn}
Let $p,\ q\in(0,\infty)$ and $\omega\in M(\rp)$.
Then the \emph{inhomogeneous block
space}\index{inhomogeneous block space}
$\Nbspace$\index{$\Nbspace$} is
defined to be the set of all the
measurable functions $f$ on $\rrn$ such that
\begin{equation*}
f=\sum_{l\in\mathbb{N}}\sum_{k\in\mathbb{Z}_+}
\sum_{Q\in\mathscr{Q}_{k}^{(\xil)}}\lambda_{
\xil,k,Q}b_{\xil,k,Q}
\end{equation*}
almost everywhere in $\rrn$, and
$$
\left\{\sum_{l\in\mathbb{N}}
\left[\sum_{k\in\mathbb{Z}_+}
\sum_{Q\in\mathscr{Q}_k^{(\xil)}}
\lambda_{\xil,k,Q}^{q}
\right]^{\frac{1}{q}}\right\}<\infty,
$$
where $\{\xil\}_{l\in\mathbb{N}}
\subset\rrn$,
$\{\lambda_{\xil,k,Q}\}_{l\in\mathbb{N},
\,k\in\mathbb{Z}_+,\,Q\in\mathscr{Q}
_{k}^{(\xil)}}\subset[0,\infty)$, and, for any
$l\in\mathbb{N}$, $k\in\mathbb{Z}_+$, and
$Q\in\mathscr{Q}_k^{(\xil)}$, $b_{\xil,k,Q}$
is an $(\omega,\,p)$-block supported in the cube
$Q$. Moreover, for any $f\in\Msc(\rrn)$,
\begin{displaymath}
\|f\|_{\Nbspace}:=\inf\left\{\sum_{l\in\mathbb{N}}
\left[\sum_{k\in\mathbb{Z}_+}
\sum_{Q\in\mathscr{Q}_k^{(\xil)}}
\lambda_{\xil,k,Q}^{q}
\right]^{\frac{1}{q}}\right\},
\end{displaymath}
where the infimum is taken over all
the decompositions of $f$ as above.
\end{definition}

\begin{remark}
We point out that, for any $p$, $q\in(0,\infty)$
and $\omega\in M(\rp)$,
the inhomogeneous block space $\Nbspace$ is a
linear space equipped with the quasi-seminorm
$\|\cdot\|_{\bspace}$ and, in particular,
when $p$, $q\in[1,\infty)$,
$\|\cdot\|_{\bspace}$ is a seminorm.
\end{remark}

\subsection{Duality Between Inhomogeneous Block Spaces
and \\ Global Generalized Herz Spaces}

The main target of this subsection is to establish
the duality between
inhomogeneous block spaces and inhomogeneous global
generalized Herz spaces. Namely, we have the
following dual result.

\begin{theorem}\label{pren}
Let $p,\ q\in(1,\infty)$ and $\omega\in M(\rp)$.
Then the dual space of $\Nbspace$ is
$\NHerzSea$ in the following sense:
\begin{enumerate}
 \item[{\rm(i)}] Let $g\in\NHerzSea$.
 Then the linear functional
\begin{equation}\label{funln}
\phi_{g}:\ f\mapsto\phi_{g}(f):=
\int_{\rrn}f(y)g(y)\,dy,
\end{equation}
defined for any $f\in\Nbspace$,
is bounded on $\Nbspace$.
 \item[{\rm(ii)}]
 Conversely, any continuous
 linear functional on $\Nbspace$ arises as
 in \eqref{funln} with a unique $g\in\NHerzSea$.
\end{enumerate}
Moreover, there exists a constant
$C\in[1,\infty)$ such that,
for any $g\in\NHerzSea$,
$$C^{-1}\left\|g\right\|_{\NHerzSea}\leq
\left\|\phi_{g}\right\|_{
(\mathcal{B}_{\omega}^{p,q}(\rrn))^*
}\leq C\left\|g\right\|_{\NHerzSea},$$
where $(\mathcal{B}_{\omega}^{p,q}(\rrn))^*$
denotes the dual space of $\Nbspace$.
\end{theorem}

To show this dual theorem, we require
the following equivalent characterizations
of both inhomogeneous local generalized Herz
spaces $\NHerzSx$ and inhomogeneous block
spaces $\Nbspace$.

\begin{lemma}\label{prenl1}
Let $p$, $q\in(0,\infty)$ and $\omega\in M(\rp)$.
Then
\begin{enumerate}
  \item[{\rm(i)}] for any given
  $\xi\in\rrn$, a measurable function
  $f$ belongs to $\NHerzSx$ if and only if
  $$
  f=\sum_{k\in\mathbb{Z}_+}\sum_{Q\in
  \mathscr{Q}_k^{(\xi)}}
  \lambda_{k,Q}b_{k,Q}
  $$
  almost everywhere in $\rrn$, and
  $$\left[
  \sum_{k\in\mathbb{Z}_+}\sum_{Q\in\mathscr{Q}_k
  ^{(\xi)}}\lambda_{k,Q}^q\right]
  ^{\frac{1}{q}}<\infty,$$
where $\{\lambda_{k,Q}\}_{
  k\in\mathbb{Z}_+,\,Q\in\mathscr{Q}_k^{(\xi)}}
  \subset[0,\infty)$ and,
  for any $k\in\mathbb{Z}$ and
  $Q\in\mathscr{Q}_k^{(\xi)}$, $b_{k,Q}$ is an
  $(\omega,\,p)$-block
  supported in the cube $Q$. Moreover,
  for any $f\in\NHerzSx$,
  $$
  \left\|f\right\|_{\NHerzSx}\sim\left[
  \sum_{k\in\mathbb{Z}_+}\sum_{Q\in\mathscr{Q}_k
  ^{(\xi)}}\lambda_{k,Q}^q\right]^{\frac{1}{q}},
  $$
  where the positive equivalence constants are independent
  of both $\xi$ and $f$;
  \item[{\rm(ii)}] a measurable function
  $f$ belongs to $\Nbspace$
  if and only if
  $$
  \normmm{f}_{\Nbspace}:=\inf\left\{
  \sum_{l\in\mathbb{N}}\left\|f_{\xil}
  \right\|_{\Kmc_{\omega,\xil}^{p,q}
  (\rrn)}\right\}<\infty,
  $$
  where the infimum is taken over all the sequences
  $\{\xil\}_{l\in\mathbb{N}}\subset\rrn$ and
  $\{f_{\xil}\}_{l\in\mathbb{N}}\subset\Msc(\rrn)$
  such that, for any $l\in\mathbb{N}$,
  $f_{\xil}\in\Kmp_{\omega,\xil}^{p,q}(\rrn)$ and
  $
  f=\sum_{l\in\mathbb{N}}f_{\xil}
  $
  almost everywhere in $\rrn$. Moreover, for any
  $f\in\Nbspace$,
  $$
  \|f\|_{\Nbspace}\sim\normmm{f}_{\Nbspace},
  $$
  where the positive equivalence constants are independent
  of $f$.
\end{enumerate}
\end{lemma}

\begin{proof}
Let all the symbols be as in the present
lemma. We first prove the necessity of (i).
For this purpose, fix a $\xi\in\rrn$ and
an $f\in\NHerzSx$. For any
$k\in\mathbb{Z}_+$ and $Q\in\mathscr{
Q}_k^{(\xi)}$,
let $$\lambda_{k,Q}:=\omega\left(|Q|^{1/n}\right)
\left\|f\1bf_Q\right\|_{L^p(\rrn)}.$$
Moreover, for any
$k\in\mathbb{Z}_+$ and
$Q\in\mathscr{Q}_k^{(\xi)}$,
let $b_{k,Q}:=0$ when $\|f\1bf_Q\|_{L^p(\rrn)}=0$;
otherwise, let
$$
b_{k,Q}:=\left[\omega\left(|Q|^{\frac{1}{n}}\right)
\right]^{-1}
\left[\left\|f\1bf_{Q}
\right\|_{L^p(\rrn)}\right]^{-1}f\1bf_Q.
$$
Then, for any $k\in\zp$ and
$Q\in\mathscr{Q}_k^{(\xi)}$,
$b_{k,Q}$ is an $(\omega,\,p)$-block
supported in $Q$ and
\begin{align*}
f&=f\1bf_{Q(\xi,2)}+\sum_{k\in\mathbb{N}}
f\1bf_{Q(\xi,\tka)\setminus Q(\xi,\tk)}\\
&=\sum_{k\in\zp}\sum_{Q\in\mathscr{Q}_k^{(\xi)}}
f\1bf_Q=\sum_{k\in\zp}\sum_{Q\in\mathscr{Q}_k^{(\xi)}}
\lambda_{k,Q}b_{k,Q}.
\end{align*}
Since $Q(\xi,2)\subset B(\xi,2)$,
from this and Definition \ref{ighx},
we deduce that
\begin{align}\label{prenl1e1}
\sum_{Q\in\mathscr{Q}_0^{(\xi)}}
\lambda_{k,Q}^q&=\left[\omega(2)\right]^q
\left\|f\1bf_{Q(\xi,2)}
\right\|_{L^p(\rrn)}^q\notag\\
&\lesssim\left\|f\1bf_{B(\xi,1)}
\right\|_{L^p(\rrn)}^q+\left[
\omega(2)\right]^q
\left\|f\1bf_{B(\xi,2)\setminus
B(\xi,1)}\right\|_{L^p(\rrn)}^q\notag\\
&\lesssim\left\|f\right\|_{\NHerzSx}^q<\infty.
\end{align}
On the other hand, by the fact that,
for any $k\in\mathbb{N}$, $$\mathscr{Q}_{k}^{(\xi)}
=\mathcal{Q}_{k}^{(\xi)},$$
the estimation
of \eqref{prel1e1}, and Definition \ref{ighx}
again, we obtain
\begin{align*}
\sum_{k\in\mathbb{N}}\sum_{Q\in\mathscr{Q}_k^{(\xi)}}
\lambda_{k,Q}^q&\lesssim
\sum_{k\in\mathbb{N}}\left[\omega(\tk)\right]^q
\left\|f\1bf_{B(\xi,\tk)\setminus
B(\xi,\tkm)}\right\|_{L^p(\rrn)}^q
\lesssim\left\|f\right\|_{\NHerzSx}^q.
\end{align*}
This, combined with \eqref{prenl1e1},
further implies that
\begin{align}\label{prenl1e2}
\left[\sum_{k\in\zp}
\sum_{Q\in\mathscr{Q}_k^{(\xi)}}
\lambda_{k,Q}^q\right]^{\frac{1}{q}}
\lesssim\left\|f\right\|_{\NHerzSx}<\infty,
\end{align}
where the implicit positive constant
is independent of both $\xi$ and $f$, which
completes the proof of the
necessity of (i).

Conversely, we now show the sufficiency
of (i).  To this end, let
$f\in\Msc(\rrn)$ satisfy that
$$
  f=\sum_{k\in\mathbb{Z}_+}\sum_{
  Q\in\mathscr{Q}_k^{(\xi)}}
  \lambda_{k,Q}b_{k,Q}
  $$
  almost everywhere in $\rrn$, and
  $$\left[
  \sum_{k\in\mathbb{Z}_+}\sum_{Q\in
  \mathscr{Q}_k
  ^{(\xi)}}\lambda_{k,Q}^q\right]
  ^{\frac{1}{q}}<\infty,$$
where $\{\lambda_{k,Q}\}_{
  k\in\mathbb{Z}_+,\,Q\in\mathscr{Q}_k^{(\xi)}}
  \subset[0,\infty)$ and,
  for any $k\in\mathbb{Z}_+$ and
  $Q\in\mathscr{Q}_k^{(\xi)}$, $b_{k,Q}$ is an
  $(\omega,\,p)$-block supported in the cube $Q$.
Then, using the fact that
$B(\xi,2)\subset Q(\xi,4)$, Remark \ref{remark215n},
and Definition \ref{blo}, we conclude that
\begin{align}\label{prenl1e3}
\left\|f\1bf_{B(\xi,2)}
\right\|_{L^p(\rrn)}^q&\leq
\left\|\sum_{k=0}^{1}\sum_{Q\in\mathscr{Q}_k^{(\xi)}}
\lambda_{k,Q}b_{k,Q}\right\|_{L^p(\rrn)}^q\notag\\
&\lesssim\sum_{k=0}^{1}\sum_{Q\in\mathscr{Q}_{k}
^{(\xi)}}\lambda_{k,Q}^q\left\|b_{k,Q}
\right\|_{L^p(\rrn)}^q\lesssim\sum_{k=0}^{1}
\sum_{Q\in\mathscr{Q}_k^{(\xi)}}\lambda_{k,Q}^q.
\end{align}
On the other hand, by the fact that,
for any $k\in\mathbb{N}$, $$\mathscr{Q}_{k}^{(\xi)}
=\mathcal{Q}_{k}^{(\xi)}$$
and the estimation
of \eqref{prel1e6}, we find that,
for any $j\in\mathbb{N}\cap[2,\infty)$,
$$
\left\|f\1bf_{B(\xi,2^j)
\setminus B(\xi,2^{j-1})}\right\|_{L^p(\rrn)}^q
\lesssim\left[\omega(2^j)\right]^{-q}
\sum_{k=j}^{j+1}\sum_{Q\in\mathscr{Q}_k^{(\xi)}}
\lambda_{k,Q}^q.
$$
Combining this, \eqref{prenl1e3}, and
Definition \ref{ighx}, we obtain
\begin{align*}
\left\|f\right\|_{\NHerzSx}&\lesssim
\left\{\left\|f\1bf_{
B(\xi,2)}\right\|_{L^p(\rrn)}^q+
\sum_{j=2}^{\infty}
\left[\omega(2^j)\right]^q
\left\|f\1bf_{B(\xi,2^j)
\setminus B(\xi,2^{j-1})}
\right\|_{L^p(\rrn)}^q\right\}^{\frac{1}{q}}\\
&\lesssim\left[\sum_{k\in\zp}
\sum_{Q\in\mathscr{Q}_{k}^{(\xi)}}
\lambda_{k,Q}^q\right]^{\frac{1}{q}}<\infty,
\end{align*}
where the implicit positive constants
are independent of both $\xi$ and $f$,
which further implies that
$f\in\NHerzSx$ and hence completes the
proof of the sufficiency of (i).
Moreover, from this
and \eqref{prenl1e2}, it follows that,
for any $f\in\NHerzSx$,
$$
\left\|f\right\|_{\NHerzSx}\sim\left[
\sum_{k\in\mathbb{Z}_+}
\sum_{Q\in\mathscr{Q}_k
^{(\xi)}}\lambda_{k,Q}^q
\right]^{\frac{1}{q}},
$$
where the positive equivalence constants are independent
of both $\xi$ and $f$. This finishes the proof of (i).

Next, repeating the proof of
Theorem \ref{pre}(ii) via replacing
Theorem \ref{pre}(i) therein
by (i), we find that
(ii) holds true, which then completes the proof
of Theorem \ref{pren}.
\end{proof}

From Lemma \ref{prenl1}, we immediately
deduce the following relation between
inhomogeneous local generalized
Herz spaces and inhomogeneous block spaces,
which is useful in the below roof of
Theorem \ref{pren}; we omit the details.

\begin{lemma}\label{een}
Let $p$, $q\in(0,\infty)$
and $\omega\in M(\rp)$.
Then, for any given $\xi\in\rrn$,
$\NHerzSx\subset\Nbspace$.
Moreover, there exists a
positive constant $C$, independent
of $\xi$, such that,
for any $f\in\NHerzSx$,
$$
\|f\|_{\Nbspace}\leq
C\|f\|_{\NHerzSx}.
$$
\end{lemma}

Via the above preparations, we
now show Theorem \ref{pren}.

\begin{proof}[Proof of Theorem \ref{pren}]
Let all the symbols be as in
the present theorem. Then,
repeating the proof of Theorem \ref{pre}
with Lemmas \ref{prel1} and \ref{ee}
therein replaced, respectively, by
Lemmas \ref{prenl1} and \ref{een},
we find that both (i) and (ii) of the present
theorem hold true. This
then finishes the proof of Theorem
\ref{pren}
\end{proof}

\subsection{Boundedness of Sublinear Operators}

The main target of this subsection is
to establish the boundedness of
both sublinear operators and
Calder\'{o}n--Zygmund operators on
inhomogeneous block spaces
(see Theorems \ref{subn} and
\ref{czoon} below). In particular,
the boundedness of powered Hardy--Littlewood
maximal operators on inhomogeneous block
spaces (see Corollary
\ref{maxbln} below) can be concluded directly
by the boundedness of sublinear operators
on inhomogeneous
block spaces, which plays a key role in
the study of inhomogeneous generalized
Herz--Hardy spaces in the next chapter.

First, we establish the following
boundedness criterion
of sublinear operators on inhomogeneous
block spaces.

\begin{theorem}\label{subn}
Let $p\in(1,\infty)$, $q\in(0,\infty)$,
and $\omega\in M(\rp)$ satisfy
$$-\frac{n}{p}<\mi(\omega)
\leq\MI(\omega)<\frac{n}{p'},$$
where $\frac{1}{p}+\frac{1}{p'}=1$.
Let $T$ be a bounded sublinear operator
on $L^p(\rrn)$ satisfying that
there exists a positive constant
$\widetilde{C}$ such that, for any
$f\in\NHerzSo$ and $x\notin\supp(f)$,
\begin{equation*}
|T(f)(x)|\leq\widetilde{C}
\int_{\rrn}\frac{|f(y)|}{|x-y|^n}\,dy
\end{equation*}
and, for any $\{f_{j}\}_{j\in\mathbb{N}}
\subset\Msc(\rrn)$ and
almost every $x\in\rrn$,
\begin{equation}\label{contin}
\left|T\left(\sum_{j\in\mathbb{N}}
f_{j}\right)(x)\right|\leq\sum_{j
\in\mathbb{N}}|T(f_{j})(x)|.
\end{equation}
If $T$ is well defined on $\Nbspace$, then
there exists a positive constant
$C$ such that, for any
$f\in\Nbspace$,
$$
\left\|T(f)\right\|_{\Nbspace}
\leq C\|f\|_{\Nbspace}.
$$
\end{theorem}

In order to prove this theorem,
we first establish the following
lattice property of inhomogeneous
block spaces.

\begin{lemma}\label{equan}
Let $p,\ q\in(0,\infty)$ and
$\omega\in M(\rp)$.
Then a measurable function $f$
on $\rrn$ belongs to the block space
$\Nbspace$ if and only if
there exists a measurable function $g\in
\Nbspace$ such that $|f|\leq g$
almost everywhere in $\rrn$. Moreover,
for these $f$ and $g$,
$$\|f\|_{\Nbspace}\leq\|g\|_{\Nbspace}.$$
\end{lemma}

\begin{proof}
Let all the symbols be as in
the present lemma and $f\in\Msc(\rrn)$.
Then, applying an argument similar to that used
in the proof of Lemma \ref{ee}
with Definition \ref{blos}
therein replaced by Definition \ref{blosn},
we conclude that $f\in\Nbspace$ if and only if
there exists a measurable function $g\in
\Nbspace$ such that $|f|\leq g$
almost everywhere in $\rrn$, and
$$\left\|f\right\|_{\Nbspace}
\leq\left\|g\right\|_{\Nbspace},$$
which complete the proof of Lemma \ref{equan}.
\end{proof}

\begin{remark}
Let $p$, $q\in(0,\infty)$,
$\omega\in M(\rp)$, and $f\in\Msc(\rrn)$.
Then, by the argument similar to that
used in Remark \ref{remark233}, we find that
$f\in\Nbspace$ if and only if
$|f|\in\Nbspace$ and, for any
$g\in\Nbspace$,
$$
\left\|g\right\|_{\Nbspace}
=\left\|\,|g|\,\right\|_{\Nbspace}.
$$
\end{remark}

Then, with the help of the above lemma
and the equivalent
characterization of inhomogeneous
block spaces obtained in Lemma \ref{prenl1}(ii),
we obtain the following boundedness
criterion of sublinear operators
on inhomogeneous block spaces, which
plays a key role in the proof of Theorem
\ref{subn}.

\begin{proposition}\label{subrn}
Let $p$, $q\in(0,\infty)$, $\omega\in M(\rp)$,
and $T$ be a bounded sublinear
operator on $\NHerzSo$ such that
\eqref{contin} holds true.
If $T$ is well defined on $\Nbspace$, then
there exists a positive constant $C$ such that,
for any $f\in\Nbspace$,
$$
\left\|T(f)\right\|_{\Nbspace}\leq C\|f\|_{\Nbspace}.
$$
\end{proposition}

\begin{proof}
Let all the symbols be as in
the present proposition and
$f\in\Nbspace$. Then,
repeating the proof of Proposition \ref{subr}
via replacing Lemmas \ref{prel1}(ii)
and \ref{equa} therein,
respectively, by Lemmas \ref{prenl1}(ii)
and \ref{equan}, we obtain
$$
\left\|T(f)\right\|_{\Nbspace}
\lesssim\left\|f\right\|_{\Nbspace},
$$
which completes the proof of Proposition
\ref{subrn}.
\end{proof}

We next show Theorem \ref{subn}.

\begin{proof}[Proof of Theorem \ref{subn}]
Let all the symbols be as in the present theorem.
Then, by both Theorem \ref{vmaxhn}
and Proposition \ref{subrn},
we find that,
for any $f\in\Nbspace$,
$$
\left\|T(f)\right\|_{\Nbspace}\lesssim
\|f\|_{\Nbspace}.
$$
This then finishes the proof of
Theorem \ref{subn}.
\end{proof}

As an application,
we now establish the following
boundedness of powered Hardy--Littlewood
maximal operators on inhomogeneous block spaces.

\begin{corollary}\label{maxbln}
Let $p\in(1,\infty)$, $q\in(0,\infty)$, $r\in[1,p)$, and
$\omega\in M(\rp)$ satisfy
$$
-\frac{n}{p}<\mi(\omega)\leq\MI(\omega)<n
\left(\frac{1}{r}-\frac{1}{p}\right).
$$
Then there exists a positive constant $C$ such that,
for any $f\in L^1_{\mathrm{loc}}(\rrn)$,
$$
\left\|\mc^{(r)}(f)\right\|_{\Nbspace}\leq
C\|f\|_{\Nbspace}.
$$
\end{corollary}

\begin{proof}
Let all the symbols be as in
the present corollary. Then,
repeating an argument similar
to that used in the proof of Corollary
\ref{maxbl} with Lemma \ref{convexll},
Corollary \ref{bhlo}, and Proposition
\ref{subr} therein replaced,
respectively, by Lemma \ref{convexlln},
Corollary \ref{bhlon}, and Proposition
\ref{subrn}, we find that, for any
$f\in L^1_{\mathrm{loc}}(\rrn)$,
$$
\left\|\mc^{(r)}(f)\right\|_{\Nbspace}\lesssim
\|f\|_{\Nbspace}.
$$
This implies that the powered Hardy--Littlewood
maximal operator $\mc^{(r)}$ is bounded
on $\Nbspace$, and hence finishes the proof of
Corollary \ref{maxbln}.
\end{proof}

Finally, we turn to establish
the following boundedness of
Calder\'{o}n--Zygmund operators
on inhomogeneous block spaces.

\begin{theorem}\label{czoon}
Let $p,\ q\in(1,\infty)$
and $\omega\in M(\rp)$ with
$$
-\frac{n}{p}<\mi(\omega)
\leq\MI(\omega)<\frac{n}{p'},
$$
where $\frac{1}{p}+\frac{1}{p'}=1.$
Assume $d\in\zp$
and $T$ is a $d$-order Calder\'{o}n--Zygmund
operator as in Definition \ref{defin-C-Z-s}.
Then $T$ is well defined on $\Nbspace$
and there exists a positive constant
$C$ such that, for any $f\in\Nbspace$,
$$
\left\|T(f)\right\|_{\Nbspace}
\leq C\|f\|_{\Nbspace}.
$$
\end{theorem}

\begin{proof}
Let all the symbols be as
in the present theorem. Then,
from an argument similar to
that used in the proof of
Theorem \ref{czoo} with
Lemma \ref{prel1}(ii) and Corollary
\ref{czoh} therein
replaced, respectively, by
Lemma \ref{prenl1}(ii) and Corollary
\ref{czohn}, we deduce that
$T$ is well defined on $\Nbspace$ and,
for any $f\in\Nbspace$,
$$
\left\|T(f)\right\|_{\Nbspace}
\lesssim\|f\|_{\Nbspace},
$$
which completes the proof of Theorem
\ref{czoon}.
\end{proof}

\section{Boundedness and Compactness Characterizations
of Commutators}

Let $b\in L^1_{\mathrm{loc}}(\rrn)$
and $T_{\Omega}$ be a singular integral
operator defined as in \eqref{czcc}
with homogeneous kernel $\Omega$ satisfying
\eqref{degrz}.
Recall that the commutator $[b,T_{\Omega}]$
is defined as in \eqref{comu}.
In this section, we establish the
boundedness and the compactness
characterizations of the commutator
$[b,T_\Omega]$ on inhomogeneous
generalized Herz spaces via
the boundedness and the compactness
characterizations of commutators
on ball Banach function spaces
obtained in \cite{TYYZ}
and Chapter \ref{sec5}.

\subsection{Boundedness Characterizations}

The target of this subsection
is to establish the boundedness
characterization of commutators
on both inhomogeneous
local and inhomogeneous
global generalized Herz spaces.
We first consider commutators
on inhomogeneous local generalized Herz spaces.
Namely, we have the following conclusion.

\begin{theorem}\label{iffbddln}
Let $p,\ q\in(1,\infty)$ and
$\omega\in M(\rp)$ satisfy
$$
-\frac{n}{p}<\mi(\omega)
\leq\MI(\omega)<\frac{n}{p'},
$$
where $\frac{1}{p}+\frac{1}{p'}=1.$
Assume that $\Omega$ is a homogeneous Lipschitz
function of degree zero on $\mathbb{S}^{n-1}$
satisfying \eqref{mz},
$T_{\Omega}$ a singular integral
operator with homogeneous kernel $\Omega$,
and $b\in L^1_{{\rm loc}}(\rrn)$.
Then the commutator $[b,T_{\Omega}]$
is bounded on the inhomogeneous
local generalized Herz space
$\NHerzSo$ if and only if $b\in\BMO$.
Moreover, there exists a constant
$C\in[1,\infty)$, independent of $b$, such that
$$
C^{-1}\|b\|_{\BMO}\leq\|[b,T_{\Omega}]
\|_{\NHerzSo
\to\NHerzSo}\leq C\|b\|_{\BMO}.
$$
\end{theorem}

To show this theorem, we first
establish the following boundedness
of commutators on inhomogeneous
local generalized Herz spaces.

\begin{proposition}\label{suffbddln}
Let $p,\ q\in(1,\infty)$, $\omega\in M(\rp)$
be as in Theorem \ref{iffbddln},
and $r\in(1,\infty]$.
Assume that $b\in\BMO$, $\Omega\in
L^r(\mathbb{S}^{n-1})$ satisfies
both \eqref{degrz} and \eqref{mz},
and $T_{\Omega}$ is a singular integral
operator with homogeneous kernel $\Omega$.
Then there exists a positive constant $C$
such that, for any $f\in\NHerzSo$,
$$
\left\|[b,T_{\Omega}](f)\right\|_{\NHerzSo}
\leq C\|b\|_{\BMO}\|f\|_{\NHerzSo}.
$$
\end{proposition}

\begin{proof}
Let all the symbols be as in
the present proposition.
Then, from Theorem \ref{ballln},
we deduce that the inhomogeneous
local generalized
Herz space $\NHerzSo$ under consideration
is a BBF space.

Therefore, in order to complete
the proof of the present proposition,
it suffices to show that $\NHerzSo$ satisfies
all the assumptions of Lemma \ref{suffbddll1}.
Indeed, applying Corollary \ref{bhlon},
we conclude that the Hardy--Littlewood maximal
operator $\mc$ is bounded on $\NHerzSo$, namely,
for any $f\in L^1_{\mathrm{loc}}(\rrn)$,
\begin{equation}\label{suffbddlne1}
\left\|\mc(f)\right\|_{\NHerzSo}
\lesssim\left\|f\right\|_{\NHerzSo}.
\end{equation}
On the other hand, from the assumption
$\MI(\omega)\in(-\frac{n}{p},\frac{n}{p'}),$
it follows that
\begin{equation*}
\frac{n}{\MI(\omega)+n/p}>\frac{n}{n(1/p'+1/p)}=1.
\end{equation*}
Combining this and the
assumptions $p$, $q\in(1,\infty)$,
we conclude that
$$
\min\left\{p,q,\frac{n}{
\MI(\omega)+n/p}\right\}\in(1,\infty).
$$
This, together with Lemma \ref{mbhaln},
further implies that
the Hardy--Littlewood maximal
operator $\mc$ is bounded on $(\NHerzSo)'$,
namely, for any
$f\in L^1_{\mathrm{loc}}(\rrn)$,
\begin{equation}\label{suffbddlne2}
\left\|\mc(f)\right\|_{(\NHerzSo)'}
\lesssim\left\|f\right\|_{(\NHerzSo)'}.
\end{equation}
Using this and \eqref{suffbddlne1},
we find that the
inhomogeneous local
generalized Herz space $\NHerzSo$ under
consideration satisfies all
the assumptions of Lemma \ref{suffbddll1},
which completes the proof of
Proposition \ref{suffbddln}.
\end{proof}

To show Theorem \ref{iffbddln},
we also require the following
necessity of the boundedness
of commutators on inhomogeneous local
generalized Herz spaces.

\begin{proposition}\label{necebddln}
Let $p\in(1,\infty)$, $q\in[1,\infty)$,
and $\omega\in M(\rp)$
be as in Theorem \ref{iffbddln}.
Assume that $b\in L^1_{{\rm loc}}(\rrn)$
and $\Omega\in
L^\infty(\mathbb{S}^{n-1})$ satisfies that
there exists an open set $\Lambda\subset
\mathbb{S}^{n-1}$ such that $\Omega$
never vanishes and
never changes sign on $\Lambda$.
If the commutator $[b,T_{\Omega}]$
is bounded on the inhomogeneous
local generalized Herz space $\NHerzSo$
and, for any bounded
measurable set $F\subset\rrn$ and
almost every $x\in\rrn\setminus\overline{F}$,
\begin{equation}\label{necebddlne1}
\left[b,T_\Omega\right]\left(\1bf_{F}\right)(x)
=\int_{F}\left[b(x)-b(y)\right]\frac{
\Omega(x-y)}{|x-y|^n}\,dy,
\end{equation}
then $b\in
\BMO$ and there exists a positive constant $C$,
independent of $b$, such that
$$
\|b\|_{\BMO}\leq C\|[b,T_{\Omega}]
\|_{\NHerzSo\to\NHerzSo}.
$$
\end{proposition}

\begin{proof}
Let all the symbols be as in the
present proposition.
Then, combining the assumptions
$p$, $q\in[1,\infty)$,
and Theorem \ref{ballln},
we conclude that the inhomogeneous
local generalized Herz space
$\NHerzSo$ under consideration
is a BBF space.
On the other hand,
from Corollary \ref{bhlon}, it follows
that the Hardy--Littlewood maximal operator
$\mc$ is bounded on the
Herz space $\NHerzSo$. This implies
that $\NHerzSo$ satisfies
all the assumptions of
Lemma \ref{necebddll1} and
hence finishes the proof
of Proposition \ref{necebddln}.
\end{proof}

We next prove Theorem \ref{iffbddln}.

\begin{proof}[Proof of Theorem
\ref{iffbddln}]
Let all the symbols be as in
the present theorem. Then, repeating
the proof of Theorem \ref{iffbddl}
via replacing Propositions
\ref{suffbddl} and \ref{necebddl}
therein, respectively,
by Propositions \ref{suffbddln} and
\ref{necebddln}, we complete
the proof of Theorem \ref{iffbddln}.
\end{proof}

In the remainder of this subsection,
we are devoted to establishing the following
boundedness characterization of
commutators on inhomogeneous global
generalized Herz spaces.

\begin{theorem}\label{iffbddgn}
Let $p,\ q\in(1,\infty)$ and
$\omega\in M(\rp)$ satisfy
$$
-\frac{n}{p}<\mi(\omega)
\leq\MI(\omega)<0,
$$
where $\frac{1}{p}+\frac{1}{p'}=1.$
Assume that $\Omega$ is a homogeneous
function of degree zero on $\mathbb{S}^{n-1}$
satisfying both \eqref{lips}
and \eqref{mz}, and $T_{\Omega}$ a
singular integral
operator with homogeneous kernel
$\Omega$. Then, for any
$b\in L^1_{{\rm loc}}(\rrn)$, the commutator
$[b,T_{\Omega}]$ is bounded on
the inhomogeneous global generalized Herz space
$\NHerzS$ if and only if $b\in\BMO$. Moreover,
there exists a constant $C\in[1,\infty)$,
independent of $b$, such that
$$
C^{-1}\|b\|_{\BMO}\leq\left\|[b,T_{\Omega}]
\right\|_{\NHerzS\to\NHerzS}
\leq C\|b\|_{\BMO}.
$$
\end{theorem}

To show this boundedness characterization,
we first give the following
boundedness of commutators on
inhomogeneous global generalized
Herz spaces.

\begin{proposition}\label{suffbddgn}
Let $p$, $q$,
and $\omega$ be as in Theorem
\ref{iffbddgn}
and $r\in(1,\infty]$.
Assume that $\Omega\in
L^r(\mathbb{S}^{n-1})$ satisfies
both \eqref{degrz} and \eqref{mz},
and $T_{\Omega}$ is a singular integral
operator with homogeneous kernel $\Omega$.
Then there exists a positive constant $C$
such that, for any $b\in\BMO$ and $f\in\NHerzS$,
$$
\|[b,T_{\Omega}](f)\|_{\NHerzS}\leq
C\|b\|_{\BMO}\|f\|_{\NHerzS}.
$$
\end{proposition}

\begin{proof}
Let all the symbols be as in
the present proposition. Then,
repeating the proof of Proposition
\ref{suffbddg} with both Definition
\ref{gh} and Theorem \ref{extrag}
therein replaced, respectively,
by both Definition \ref{igh} and Theorem
\ref{extragn}, we conclude that,
for any $b\in\BMO$ and $f\in\NHerzS$,
$$
\left\|\left[b,T_{\Omega}
\right](f)\right\|_{\NHerzS}
\lesssim\|b\|_{\BMO}\|f\|_{\NHerzS}
$$
with the implicit
positive constant independent of both
$b$ and $f$. This then finishes the proof of
Proposition \ref{suffbddln}.
\end{proof}

On the other hand, we now
establish the following necessity
of the boundedness of commutators
on inhomogeneous global generalized
Herz spaces.

\begin{proposition}\label{necebddgn}
Let $p\in(1,\infty)$,
$q\in[1,\infty)$, and $\omega\in M(\rp)$
be as in Theorem \ref{iffbddgn}.
Assume that $b\in L^1_{{\rm loc}}
(\rrn)$ and $\Omega\in L^\infty(\mathbb{S}^{n-1})$
satisfies that there exists an open set
$\Lambda\subset\mathbb{S}^{n-1}$
such that $\Omega$
never vanishes and never changes sign
on $\Lambda$.
If the commutator $[b,T_{\Omega}]$
is bounded on the inhomogeneous
global generalized Herz space $\NHerzS$
and, for any bounded
measurable set $F\subset\rrn$
and almost every
$x\in\rrn\setminus\overline{F}$,
$$
\left[b,T_\Omega\right]
\left(\1bf_{F}\right)(x)
=\int_{F}\left[b(x)-b(y)\right]\frac{
\Omega(x-y)}{|x-y|^n}\,dy,
$$
then $b\in
\BMO$ and there exists
a positive constant $C$,
independent of $b$, such that
$$
\|b\|_{\BMO}\leq C\|[b,T_{\Omega}]
\|_{\NHerzS\to\NHerzS}.
$$
\end{proposition}

\begin{proof}
Let all the symbols
be as in the present proposition.
Then, combining the assumption
$\MI(\omega)\in(-\infty,0)$,
Theorem \ref{balln}, and
Corollary \ref{Th4n},
we find that the inhomogeneous
global generalized Herz space
$\NHerzS$ under consideration
is a BBF space and the
Hardy--Littlewood maximal operator
$\mc$ is bounded on $\NHerzS$.
From this and Lemma \ref{necebddll1},
it follows that $b\in\BMO$ and
$$
\|b\|_{\BMO}\lesssim
\left\|
\left[b,T_{\Omega}\right]
\right\|_{\NHerzS\to\NHerzS},
$$
which completes the proof of
Proposition \ref{necebddgn}.
\end{proof}

Via Propositions \ref{suffbddgn} and
\ref{necebddgn}, we
next show Theorem \ref{iffbddgn}.

\begin{proof}[Proof of Theorem
\ref{iffbddgn}]
Let all the symbols be as in
the present theorem. Then, repeating
the proof of Theorem \ref{iffbddg}
via replacing Propositions
\ref{suffbddg} and \ref{necebddg}
therein, respectively,
by Propositions \ref{suffbddgn} and
\ref{necebddgn}, we complete
the proof of Theorem \ref{iffbddgn}.
\end{proof}

\subsection{Compactness Characterizations}

In this subsection, we establish the
compactness characterizations
of commutators on inhomogeneous
generalized Herz spaces.
To begin with, we give
the compactness characterization
of commutators
on the inhomogeneous local generalized
Herz space $\NHerzSo$ as follows.

\begin{theorem}\label{iffcomln}
Let $p,\ q\in(1,\infty)$ and
$\omega\in M(\rp)$ satisfy
$$
-\frac{n}{p}<\mi(\omega)\leq
\MI(\omega)<\frac{n}{p'},
$$
where $\frac{1}{p}+\frac{1}{p'}=1.$
Assume that $\Omega$ is a homogeneous
function satisfying \eqref{lips}, \eqref{degrz},
and \eqref{mz}, and $T_{\Omega}$ a singular integral
operator with homogeneous kernel $\Omega$. Then,
for any $b\in L^1_{{\rm loc}}(\rrn)$,
the commutator $[b,T_{\Omega}]$
is compact on the inhomogeneous
local generalized Herz space
$\NHerzSo$ if and only if $b\in\CMO$.
\end{theorem}

In order to prove this theorem, we
first show the following sufficiency of
the compactness of commutators
on $\NHerzSo$.

\begin{proposition}\label{compsuln}
Let $p$, $q$, and $\omega$ be
as in Theorem \ref{iffcomln}.
Assume that $b\in L^1_{\mathrm{loc}}(\rrn)$,
$\Omega\in L^\infty(\mathbb{S}^{n-1})$ satisfies
\eqref{degrz}, \eqref{mz},
and the $L^{\infty}$-Dini condition,
and $T_{\Omega}$ is a singular
integral operator with homogeneous kernel
$\Omega$. If $b\in\CMO$,
then the commutator $[b,T_{\Omega}]$ is compact on the
inhomogeneous local generalized Herz space $\NHerzSo$.
\end{proposition}

\begin{proof}
Let all the symbols be as
in the present proposition and $b\in\CMO$.
Then, repeating an argument similar to
that used in the proof of Proposition
\ref{compsul} with Theorem \ref{balll},
\eqref{th175e1}, and \eqref{th175e2}
therein replaced, respectively,
by Theorem \ref{ballln},
\eqref{suffbddlne1}, and \eqref{suffbddlne2},
we find that the commutator $[b,T_{\Omega}]$
is compact on the inhomogeneous
local generalized Herz space $\NHerzSo$
under consideration. This then finishes the
proof of Proposition \ref{compsuln}.
\end{proof}

Conversely, we establish the
following necessity of the compactness
of commutators on inhomogeneous local
generalized Herz spaces.

\begin{proposition}\label{compneceln}
Let $p$, $q$, and $\omega$ be as
in Theorem \ref{iffcomln}.
Assume that $b\in L^1_{{\rm loc}}(\rrn)$ and $\Omega\in
L^\infty(\mathbb{S}^{n-1})$ satisfies that
there exists an open set $\Lambda\subset
\mathbb{S}^{n-1}$ such that
$\Omega$ never vanishes and never
changes sign on $\Lambda$. If
the commutator $[b,T_{\Omega}]$ is compact
on the inhomogeneous local
generalized Herz space $\NHerzSo$ and
\eqref{necebddlne1} holds true,
then $b\in\CMO$.
\end{proposition}

\begin{proof}
Let all the symbols be as
in the present proposition, and let $b\in
L^1_{\mathrm{loc}}(\rrn)$ be such
that the commutator $[b,T_{\Omega}]$
is compact on the inhomogeneous
local generalized Herz space $\NHerzSo$
under consideration.
Then, repeating the proof of Proposition
\ref{compnecel} via replacing Theorem \ref{balll},
\eqref{th175e1}, and \eqref{th175e2}
therein, respectively, by Theorem \ref{ballln},
\eqref{suffbddlne1}, and \eqref{suffbddlne2},
we conclude that $b\in\CMO$, and hence
complete the proof of Proposition \ref{compneceln}.
\end{proof}

Applying Propositions
\ref{compsuln} and \ref{compneceln}, we now
give the proof of Theorem \ref{iffcomln}.

\begin{proof}[Proof of Theorem
\ref{iffcomln}]
Let all the symbols be as in
the present theorem. Then, from
an argument similar to that
used in the proof of Theorem \ref{iffcoml}
with both Propositions
\ref{compsul} and \ref{compnecel}
therein replaced, respectively,
by both Propositions \ref{compsuln}
and \ref{compneceln}, we deduce
that Theorem \ref{iffcomln} holds
true, and hence complete the
proof of the present theorem.
\end{proof}

Next, we are devoted to establishing
the compactness characterization of
commutators on inhomogeneous global
generalized Herz spaces. Indeed,
we have the following conclusion.

\begin{theorem}\label{iffcomgn}
Let $p,\ q\in(1,\infty)$ and
$\omega\in M(\rp)$ satisfy
$$
-\frac{n}{p}<\mi(\omega)\leq\MI(\omega)<0,
$$
where $\frac{1}{p}+\frac{1}{p'}=1.$
Assume that $\Omega$ is a homogeneous
function of degree zero on $\mathbb{S}^{n-1}$
satisfying both \eqref{lips} and \eqref{mz},
and $T_{\Omega}$ a singular integral
operator with homogeneous kernel $\Omega$. Then, for
any $b\in L^1_{{\rm loc}}(\rrn)$,
the commutator $[b,T_{\Omega}]$
is compact on the inhomogeneous global generalized
Herz space $\NHerzS$ if and only if $b\in\CMO$.
\end{theorem}

To show this theorem, we first give
the following sufficiency of
the compactness of commutators
on the inhomogeneous global generalized
Herz space $\NHerzS$.

\begin{proposition}\label{compsugn}
Let $p$, $q$, and $\omega$ be
as in Theorem \ref{iffcomgn}.
Assume that
$b\in L^1_{\mathrm{loc}}(\rrn)$,
$\Omega\in L^\infty(\mathbb{S}^{n-1})$
satisfies \eqref{degrz}, \eqref{mz},
and the $L^{\infty}$-Dini condition,
and $T_{\Omega}$ is a singular
integral operator with homogeneous
kernel $\Omega$. If $b\in\CMO$,
then the commutator $[b,T_{\Omega}]$
is compact on $\NHerzS$.
\end{proposition}

\begin{proof}
Let all the symbols
be as in the present proposition
and $b\in\CMO$. Then, from
the assumptions $p$,
$q\in(1,\infty)$ and $\MI(\omega)\in(-\infty,0)$,
and Theorem \ref{balln}, it follows
that the inhomogeneous global
generalized Herz space $\NHerzS$
under consideration is a BBF space.
Thus, to finish the proof of the present
proposition, we only need to show that the assumptions
(i) through (iii) of Proposition \ref{compsugl1}
hold true for $\NHerzS$.

First, we show that Proposition \ref{compsugl1}(i)
holds true for $\NHerzS$, namely,
both $[b,T_{\Omega}]$ and $T_{\Omega}^*$
are well defined on $\NHerzS$.
Indeed, applying Definition \ref{igh},
we find that, for any $f\in\NHerzS$,
$$
\left\|f\right\|_{\NHerzSo}
\leq\left\|f\right\|_{\NHerzS}<\infty,
$$
which implies that $f\in\NHerzSo$
and hence $$\NHerzS\subset\NHerzSo.$$
Combining this and Proposition
\ref{suffbddln}, we conclude that
$[b,T_{\Omega}](f)\in\NHerzSo$ and hence
$[b,T_{\Omega}]$ is well defined on $\NHerzS$.
Now, we prove that $T_{\Omega}^*$
is also well defined on $\NHerzS$. Notice that,
by the fact that $\NHerzS\subset\NHerzSo$,
we only need to show that
$T_{\Omega}^*$ is well defined on $\NHerzSo$.
To this end, repeating an argument
similar to that used in the proof of
Lemma \ref{compsugl2} with
Theorem \ref{balll}, \eqref{th175e1},
and \eqref{th175e2} therein replaced,
respectively, by Theorem
\ref{ballln}, \eqref{suffbddlne1},
and \eqref{suffbddlne2}, we find that, for any
$f\in\NHerzSo$,
$$
\left\|T_{\Omega}^*(f)\right\|_{\NHerzSo}
\lesssim\left\|f\right\|_{\NHerzSo},
$$
which then implies that $T_{\Omega}^*$
is well defined on $\NHerzSo$ and hence
on $\NHerzS$. Therefore, Proposition
\ref{compsugl1}(i) holds
true for $\NHerzS$.

Next, we prove that $\NHerzS$ satisfies
Proposition \ref{compsugl1}(ii). Indeed,
from Corollary \ref{Th4n}, it follows
that, for any $f\in L^1_{\mathrm{loc}}(\rrn)$,
$$
\left\|\mc(f)\right\|_{\NHerzS}
\lesssim\left\|f\right\|_{\NHerzS},
$$
which implies that the Hardy--Littlewood
maximal operator $\mc$ is bounded on
$\NHerzS$ and hence Proposition \ref{compsugl1}(ii)
holds true.

Finally, applying Theorem \ref{extragn},
we find that the inhomogeneous
global generalized Herz space $\NHerzS$
under consideration satisfies
Proposition \ref{compsugl1}(iii).
Thus, by the fact that $\NHerzS$
is a BBF space and Proposition \ref{compsugl1} with
$X:=\NHerzS$, we find that
the commutator $[b,T_{\Omega}]$ is compact
on $\NHerzS$, which completes the
proof of Proposition \ref{compsugn}.
\end{proof}

To establish the compactness
characterization Theorem \ref{iffcomgn},
we also need the following necessity
of the compactness of commutators on $\NHerzS$.

\begin{proposition}\label{compnecegn}
Let $p$, $q$, and $\omega$ be as
in Theorem \ref{iffcomgn}.
Assume that $b\in L^1_{{\rm loc}}(\rrn)$
and $\Omega\in
L^\infty(\mathbb{S}^{n-1})$ satisfies that
there exists an open set $\Lambda\subset
\mathbb{S}^{n-1}$ such that
$\Omega$ never vanishes and never
changes sign on $\Lambda$. If
the commutator $[b,T_{\Omega}]$ is compact
on the inhomogeneous global
generalized Herz space $\NHerzS$ and
satisfies \eqref{necebddlne1},
then $b\in\CMO$.
\end{proposition}

\begin{proof}
Let all the symbols be as
in the present proposition, and let $b\in
L^1_{\mathrm{loc}}(\rrn)$ be such
that the commutator $[b,T_{\Omega}]$
is compact on the inhomogeneous
global generalized Herz space $\NHerzS$
under consideration.
Then, repeating the proof of Proposition
\ref{compneceg} via Theorems \ref{ball}
and \ref{extrag} and Corollary \ref{Th4}
therein replaced,
respectively, by Theorems \ref{balln}
and \ref{extragn} and Corollary \ref{Th4n},
we find that $b\in\CMO$, and hence
complete the proof of
Proposition \ref{compnecegn}.
\end{proof}

Via the above two propositions, we
now show Theorem \ref{iffcomgn}.

\begin{proof}[Proof of Theorem
\ref{iffcomgn}]
Let all the symbols be as in
the present theorem. Then, repeating
an argument similar to that
used in the proof of Theorem \ref{iffcomg}
with both Propositions
\ref{compsug} and \ref{compneceg}
therein replaced, respectively,
by both Propositions \ref{compsugn}
and \ref{compnecegn}, we then complete the proof
of Theorem \ref{iffcomgn}.
\end{proof}

\chapter{Hardy Spaces Associated
with Inhomogeneous Generalized Herz Spaces\label{sec10}}
\markboth{\scriptsize\rm\sc Hardy Spaces
Associated with Inhomogeneous
Generalized Herz Spaces}{\scriptsize\rm\sc Hardy
Spaces Associated with Inhomogeneous Generalized Herz Spaces}

The target of this chapter is
devoted to establishing a
complete real-variable theory of
Hardy spaces associated with
inhomogeneous generalized
Herz spaces.
Precisely, based on the inhomogeneous
generalized Herz spaces studied
in Chapter \ref{sec9},
we first introduce inhomogeneous generalized
Herz--Hardy spaces, inhomogeneous
localized generalized Herz--Hardy spaces,
and inhomogeneous weak generalized Herz--Hardy
spaces in this chapter. Then, via using
the known results about Hardy spaces $H_X(\rrn)$
associated with ball quasi-Banach
function spaces $X$ as well as some
improved characterizations
of $H_X(\rrn)$ established in Chapters \ref{sec6}
through \ref{sec8}, we investigate
both the real-variable characterizations
and their applications of these Hardy spaces associated with
inhomogeneous generalized
Herz spaces, which include various
maximal function, the (finite) atomic,
the molecular, and
the Littlewood--Paley function
characterizations as well as
the boundedness of Calder\'{o}n--Zygmund
operators or pseudo-differential operators, duality,
the properties of the Fourier transform,
and the real interpolation theorems
about these Herz--Hardy spaces.

\section{Inhomogeneous Generalized Herz--Hardy \\
Spaces}

In this section, we first introduce the
inhomogeneous generalized Herz--Hardy spaces
$\NHaSaHo$ and $\NHaSaH$, where $p$,
$q\in(0,\infty)$ and $\omega\in M(\rp)$.
Then we establish some real-variable characterizations
and also give some applications of them.
To be precise, we show the maximal
function, the atomic,
the finite atomic, the molecular,
and the Littlewood--Paley function
characterizations of $\NHaSaHo$ and
$\NHaSaH$. As applications, we
find the dual space of $\NHaSaHo$
under some reasonable and sharp assumptions
and establish the boundedness of the
Calder\'{o}n--Zygmund operators on both
$\NHaSaHo$ and $\NHaSaH$ as well
as investigate the Fourier
transform properties of these
generalized Herz--Hardy spaces.

To begin with, recall that,
for any given $N\in\mathbb{N}$ and
for any $f\in\mathcal{S}'(\rrn)$,
the non-tangential grand maximal function
$\mathcal{M}_N(f)$ of $f$ is defined
as in \eqref{sec6e1}. Then we introduce the
inhomogeneous generalized Herz--Hardy spaces via
the operator $\mathcal{M}_N$ as
follows.\index{inhomogeneous generalized Herz--Hardy space}

\begin{definition}\label{herzhn}
Let $p$, $q\in(0,\infty)$, $\omega\in M(\rp)$, and
$N\in\mathbb{N}$.
Then
\begin{enumerate}
  \item[{\rm(i)}] the \emph{inhomogeneous
  generalized Herz--Hardy space}
  $\NHaSaHo$\index{$\NHaSaHo$}, associated with
  the inhomogeneous local generalized Herz space $\NHerzSo$,
  is defined to be the set of all
  the $f\in\mathcal{S}'(\rrn)$ such that
  $$\left\|f\right\|_{\NHaSaHo}:=\left\|
  \mathcal{M}_{N}(f)\right\|_{\NHerzSo}<\infty;$$
  \item[{\rm(ii)}] the \emph{inhomogeneous
  generalized Herz--Hardy space}
  $\NHaSaH$\index{$\NHaSaH$}, associated with
  the inhomogeneous global generalized
  Herz space $\NHerzS$,
  is defined to be the set of all the
  $f\in\mathcal{S}'(\rrn)$ such that $$
\|f\|_{\NHaSaH}:=\left\|\mathcal{M}_{N}
(f)\right\|_{\NHerzS}<\infty.$$
\end{enumerate}
\end{definition}

\begin{remark}\label{sec101r}
We should point out that, in Definition \ref{herzhn},
when $\omega(t):=t^\alpha$ for any
$t\in(0,\infty)$ and for any given $\alpha\in\rr$,
the inhomogeneous generalized
Herz--Hardy space $\NHaSaHo$ coincides with the classical
\emph{inhomogeneous Herz-type
Hardy space}\index{inhomogeneous Herz-type Hardy space}
$HK_p^{\alpha,q}(\rrn)$\index{$HK_p^{\alpha,q}(\rrn)$}
which was originally introduced in \cite[Definition 2.1]{LuY95}
(see also \cite[Definition 2.1.1]{LYH}).
\end{remark}

\subsection{Maximal Function
Characterizations}

The main target of this subsection
is to characterize
inhomogeneous generalized
Herz--Hardy spaces via various
maximal functions. Via using the radial
and the non-tangential
maximal functions presented in
Definition \ref{smax},
we first show the following maximal
function characterizations
of the inhomogeneous generalized
Herz--Hardy space
$\NHaSaHo$.\index{maximal function characterization}

\begin{theorem}\label{Th5.4n}
Let $p$, $q$, $a,\ b\in(0,\infty)$,
$\omega\in M(\rp)$,
$N\in\mathbb{N}$, and
$\phi\in\mathcal{S}(\rrn)$ satisfy
$\int_{\rrn}\phi(x)\,dx\neq0.$
\begin{enumerate}
\item[\rm{(i)}] Let $N\in\mathbb{N}
\cap[\lfloor b+1\rfloor,\infty)$.
Then, for any $f\in\mathcal{S}'(\rrn)$,
$$
\|M(f,\phi)\|_{\NHerzSo}\lesssim\|
M^*_{a}(f,\phi)\|_{\NHerzSo}
\lesssim\|M^{**}_{b}(f,\phi)\|_{\NHerzSo},
$$
\begin{align*}
\|M(f,\phi)\|_{\NHerzSo}&\lesssim\|
\mathcal{M}_{N}(f)\|_{\NHerzSo}
\lesssim\|\mathcal{M}_{\lfloor b+1
\rfloor}(f)\|_{\NHerzSo}\\&\lesssim
\|M^{**}_{b}(f,\phi)\|_{\NHerzSo},
\end{align*}
and
$$
\|M^{**}_{b}(f,\phi)\|_{\NHerzSo}\sim\|
\mathcal{M}^{**}_{b,N}(f)\|_{\NHerzSo},
$$
where the implicit positive
constants are independent of $f$.
\item[\rm{(ii)}] Let $\omega\in M(\rp)$
satisfy $\mi(\omega)\in(-\frac{n}{p},\infty)$.
Assume $$b\in\left(\max\left\{\frac{n}{p},
\MI(\omega)+\frac{n}{p}\right\},\infty\right).$$
Then, for any $f\in\mathcal{S}'(\rrn),$
$$\|M^{**}_{b}(f,\phi)\|_{\NHerzSo}
\lesssim\|M(f,\phi)\|_{\NHerzSo},$$
where the implicit positive
constant is independent of $f$.
In particular, when $N\in\mathbb{N}
\cap[\lfloor b+1\rfloor,\infty)$,
if one of the quantities
$$\|M(f,\phi)\|_{\NHerzSo},\
\|M^*_{a}(f,\phi)\|_{\NHerzSo},\
\|\mathcal{M}_{N}
(f)\|_{\NHerzSo},$$
$$
\|M^{**}_{b}(f,\phi)\|_{\NHerzSo},
\ \text{and}\ \|\mathcal{M}^{**}_{b,N
}(f)\|_{\NHerzSo}$$
is finite, then the others are
also finite and mutually
equivalent with the
positive equivalence constants independent of $f$.
\end{enumerate}
\end{theorem}

\begin{proof}
Let all the symbols be as in
the present theorem. Then, repeating
the proof of Theorem \ref{Th5.4}
with $\Mw$, Theorem \ref{Th3}, and
Lemma \ref{mbhl} therein replaced,
respectively, by
$\MI(\omega)$, Theorem \ref{Th3n}, and
Lemma \ref{mbhln}, we find that
both (i) and (ii) of the present theorem
hold true. This finishes
the proof of Theorem \ref{Th5.4n}.
\end{proof}

\begin{remark}
We should point out that, in Theorem \ref{Th5.4n},
when $p\in(1,\infty)$ and $\omega(t):=t^\alpha$
for any $t\in(0,\infty)$ and for any given $\alpha\in\rr$,
then Theorem \ref{Th5.4n} goes back to the maximal
function characterizations of $HK_p^{\alpha,q}(\rrn)$
obtained in \cite[p.\,27]{LYH}.
\end{remark}

Next, we establish the following maximal
function characterizations
of the inhomogeneous generalized
Herz--Hardy space
$\NHaSaH$.\index{maximal function characterization}

\begin{theorem}\label{Th5.6n}
Let $p$, $q$, $a,\ b\in(0,\infty)$,
$\omega\in M(\rp)$,
$N\in\mathbb{N}$, and $\phi\in
\mathcal{S}(\rrn)$ satisfy
$\int_{\rrn}\phi(x)\,dx\neq0.$
\begin{enumerate}
\item[\rm{(i)}] Let $N\in\mathbb{N}\cap
[\lfloor b+1\rfloor,\infty)$
and $\omega$ satisfy $\MI(\omega)\in(-\infty,0)$.
Then, for any $f\in\mathcal{S}'(\rrn)$,
$$
\|M(f,\phi)\|_{\NHerzS}\lesssim
\|M^*_{a}(f,\phi)\|_{\NHerzS}
\lesssim\|M^{**}_{b}(f,\phi)\|_{\NHerzS},
$$
\begin{align*}
\|M(f,\phi)\|_{\NHerzS}&\lesssim\|
\mathcal{M}_{N}(f)\|_{\NHerzS}
\lesssim\|\mathcal{M}_{\lfloor
b+1\rfloor}(f)\|_{\NHerzS}\\&\lesssim
\|M^{**}_{b}(f,\phi)\|_{\NHerzS},
\end{align*}
and
$$
\|M^{**}_{b}(f,\phi)\|_{\NHerzS}\sim\|
\mathcal{M}^{**}_{b,N}(f)\|_{\NHerzS},
$$
where the implicit positive
constants are independent of $f$.
\item[\rm{(ii)}] Let $\omega\in M(\rp)$ satisfy
$$-\frac{n}{p}<\mi(\omega)\leq\MI(\omega)<0.$$
Assume $b\in(\frac{n}{p},\infty)$.
Then, for any $f\in\mathcal{S}'(\rrn),$
$$\|M^{**}_{b}(f,\phi)\|_{\NHerzS}
\lesssim\|M(f,\phi)\|_{\NHerzS},$$
where the implicit positive
constant is independent of $f$.
In particular, when $N\in\mathbb{N}
\cap[\lfloor b+1\rfloor,\infty)$,
if one of the quantities
$$\|M(f,\phi)\|_{\NHerzS},\
\|M^*_{a}(f,\phi)\|_{\NHerzS},\ \|\mathcal{M}_{N}
(f)\|_{\NHerzS},$$
$$
\|M^{**}_{b}(f,\phi)\|_{\NHerzS},
\ \text{and}\ \|\mathcal{M}^{**}_{b,N}(f)\|_{\NHerzS}$$
is finite, then the others are also finite and mutually
equivalent with the
positive equivalence constants independent of $f$.
\end{enumerate}
\end{theorem}

To prove this theorem, we first
show the following lemma about the boundedness of the
Hardy--Littlewood maximal operator on
inhomogeneous generalized Herz spaces.

\begin{lemma}\label{mbhgn}
Let $p,\ q\in(0,\infty)$ and $\omega\in M(\rp)$ satisfy
$\mi(\omega)\in(-\frac{n}{p},\infty)$.
Then, for any given
$r\in(0,\min\{p,\frac{n}{\MI(\omega)+n/p}\})$,
there exists a positive constant $C$ such that,
for any $f\in L^1_{\mathrm{loc}}(\rrn)$,
$$
\left\|\mc(f)\right\|_{[\NHerzS]^{1/r}}
\leq C\|f\|_{[\NHerzS]^{1/r}}.
$$
\end{lemma}

\begin{proof}
Let all the symbols be as
in the present lemma. Then, for any
given $r\in(0,\min\{p,\frac{n}{\MI(\omega)
+n/p}\})$, using
Lemma \ref{rela} and
repeating the proof of Lemma \ref{mbhg} with
$\mw$, $\Mw$, Lemma \ref{convexl}, and
Corollary \ref{Th4} therein replaced,
respectively, by
$\mi(\omega)$, $\MI(\omega)$, Lemma \ref{convexln}, and
Corollary \ref{Th4n}, we conclude that,
for any $f\in L^1_{\mathrm{loc}}(\rrn)$,
$$
\left\|\mc(f)\right\|_{[\NHerzS]^{1/r}}
\lesssim\|f\|_{[\NHerzS]^{1/r}}.
$$
This finishes the proof of Lemma \ref{mbhgn}.
\end{proof}

Via this lemma, we now show
Theorem \ref{Th5.6n}.

\begin{proof}[Proof of Theorem \ref{Th5.6n}]
Let all the symbols be as
in the present theorem. Then, repeating
the proof of Theorem \ref{Th5.6n}
via replacing Theorem \ref{Th2} and
Lemma \ref{mbhg} therein, respectively,
by Theorem \ref{Th2n} and
Lemma \ref{mbhgn}, we find that
both (i) and (ii) of the present
theorem hold true, and hence
complete the proof of Theorem \ref{Th5.6n}.
\end{proof}

\subsection{Relations with Inhomogeneous Generalized
Herz \\ Spaces}

In this subsection, we investigate the relations
between inhomogeneous
generalized Herz spaces and
associated Hardy spaces. To begin with,
the following conclusion shows that,
under some reasonable and sharp assumptions,
$\NHaSaHo=\NHerzSo$
in the sense of equivalent quasi-norms.

\begin{theorem}\label{Th5.7n}
Let $p\in(1,\infty)$,
$q\in(0,\infty)$, and
$\omega\in M(\rp)$ satisfy
$$-\frac{n}{p}<\mi(\omega)\leq\MI(\omega)<\frac{n}{p'},$$
where $\frac{1}{p}+\frac{1}{p'}=1$.
Then
\begin{enumerate}
\item[\rm{(i)}] $\NHerzSo\hookrightarrow\mathcal{S}'(\rrn)$.
\item[\rm{(ii)}] If $f\in\NHerzSo$, then $f\in\NHaSaHo$
and there exists a positive constant $C$, independent
of $f$, such that $$\left\|f\right\|_{\NHaSaHo}
\leq C\left\|f\right\|_{\NHerzSo}.$$
\item[\rm{(iii)}] If $f\in\NHaSaHo$, then there exists a
locally integrable function $g\in\NHerzSo$
such that $g$ represents $f$, which means that $f=g$ in
$\mathcal{S}'(\rrn)$, $$\left\|f\right\|_{\NHaSaHo}
=\left\|g\right\|_{\NHaSaHo},$$
and there exists a positive constant $C$, independent of $f$,
such that $$\left\|g\right\|_{\NHerzSo}
\leq C\left\|f\right\|_{\NHaSaHo}.$$
\end{enumerate}
\end{theorem}

\begin{proof}
Let all the symbols be as
in the present theorem. Then,
repeating an argument similar to
that used in the proof of Theorem \ref{Th5.7}
with both
Theorem \ref{Th3} and Lemma \ref{mbhl} therein
replaced, respectively, by
both Theorem \ref{Th3n} and Lemma \ref{mbhln},
we find that (i) through (iii)
of the present theorem hold true.
This finishes the proof of Theorem \ref{Th5.7n}.
\end{proof}

\begin{remark}
We should point out that, in Theorem \ref{Th5.7n},
if $\omega(t):=t^\alpha$ for any
$t\in(0,\infty)$ and for any given $\alpha\in\rr$,
then the conclusion obtained in this theorem
goes back to \cite[Proposition 2.1.1(1)]{LYH}.
\end{remark}

Similarly, the following conclusion shows that,
under some reasonable and sharp assumptions,
$\NHaSaH=\NHerzS$
with equivalent quasi-norms.

\begin{theorem}\label{Th5.8n}
Let $p\in(1,\infty)$, $q\in(0,\infty)$, and $\omega\in M(\rp)$
satisfy
$$-\frac{n}{p}<\mi(\omega)\leq\MI(\omega)<0.$$
Then
\begin{enumerate}
\item[\rm{(i)}] $\NHerzS\hookrightarrow\mathcal{S}'(\rrn)$.
\item[\rm{(ii)}] If $f\in\NHerzS$, then $f\in\NHaSaH$
and there exists a positive constant $C$, independent
of $f$, such that $$\left\|f\right\|_{\NHaSaH}\leq
C\left\|f\right\|_{\NHerzS}.$$
\item[\rm{(iii)}] If $f\in\NHaSaH$, then there exists a
locally integrable function $g\in\NHerzS$
such that $g$ represents $f$, which means that $f=g$ in
$\mathcal{S}'(\rrn)$, $$\left\|f\right\|_{\NHaSaH}=
\left\|g\right\|_{\NHaSaH},$$
and there exists a positive constant $C$, independent of $f$,
such that $$\left\|g\right\|_{\NHerzS}
\leq C\left\|f\right\|_{\NHaSaH}.$$
\end{enumerate}
\end{theorem}

\begin{proof}
Let all the symbols be as in
the present theorem. Then, repeating
the proof of Theorem \ref{Th5.8} via
replacing both Theorem \ref{Th2} and Lemma
\ref{mbhg} therein, respectively,
by both Theorem \ref{Th2n} and Lemma \ref{mbhgn},
we conclude that (i), (ii), and (iii)
of the present theorem hold true. This then
finishes the proof of Theorem \ref{Th5.8n}.
\end{proof}

\subsection{Atomic Characterizations}

The target of this subsection is to
establish the atomic characterization
of inhomogeneous generalized
Herz--Hardy spaces.
To begin with,
we show the atomic
characterization of the
inhomogeneous generalized
Herz--Hardy space $\NHaSaHo$. For this
purpose, we first introduce
the concepts of both $(\NHerzSo,\,r,\,d)$-atoms
and the atomic Hardy space $\NaHaSaHo$
associated with the
inhomogeneous local generalized
Herz space $\NHerzSo$ as follows.

\begin{definition}\label{atomn}
Let $p$, $q\in(0,\infty)$, $\omega\in M(\rp)$,
$r\in[1,\infty]$, and $d\in\zp$.
Then a measurable function
$a$ on $\rrn$ is called a
\emph{$(\NHerzSo,\,r,
\,d)$-atom}\index{$(\NHerzSo,\,r,\,d)$-atom}
if there exists a ball $B\in\mathbb{B}$
such that
\begin{enumerate}
  \item[{\rm(i)}] $\supp(a):=
  \{x\in\rrn:\ a(x)\neq0\}\subset B$;
  \item[{\rm(ii)}] $\|a\|_{L^r(\rrn)}
  \leq\frac{|B|^{1/r}}{\|\1bf_{B}\|_{\NHerzSo}}$;
  \item[{\rm(iii)}] for any $\alpha\in\zp^n$
  with $|\alpha|\leq d$,
  $$\int_{\rrn}a(x)x^\alpha\,dx=0.$$
\end{enumerate}
\end{definition}

\begin{definition}\label{atsln}
Let $p$, $q\in(0,\infty)$, $\omega\in M(\rp)$ with
$\mi(\omega)\in(-\frac{n}{p},\infty)$,
$r\in(\max\{1,p,\frac{n}{\mi(\omega)+n/p}\},\infty]$,
$$s\in\left(0,\min\left\{
1,p,q,\frac{n}{\MI(\omega)+n/p}\right\}\right),$$
and $d\geq\lfloor n(
1/s-1)\rfloor$ be a fixed integer.
Then the \emph{inhomogeneous generalized atomic
Herz--Hardy space}\index{inhomogeneous generalized
atomic Herz--Hardy \\space} $\NaHaSaHo$\index{$\NaHaSaHo$}
is defined to be the set of all the
$f\in\mathcal{S}'(\rrn)$ such that there exist
$\{\lambda_j\}_{j\in\mathbb{N}}
\subset[0,\infty)$ and
a sequence $\{a_{j}\}_{j\in\mathbb{N}}$ of
$(\NHerzSo,\,r,\,d)$-atoms supported,
respectively, in the balls $\{B_{j}\}_{j\in\mathbb{N}}
\subset\mathbb{B}$ such that
$$
f=\sum_{j\in\mathbb{N}}\lambda_{j}a_{j}
$$
in $\mathcal{S}'(\rrn)$,
and
$$
\left\|\left\{
\sum_{j\in\mathbb{N}}
\left[\frac{\lambda_j}{\|
\1bf_{B_j}\|_{\NHerzSo}}\right]^s
\1bf_{B_j}\right\}^
{\frac{1}{s}}\right\|_{\NHerzSo}<\infty.
$$
Moreover, for any $f\in\NaHaSaHo$,
$$
\|f\|_{\NaHaSaHo}:=\inf\left\{
\left\|\left\{\sum_{j\in\mathbb{N}}
\left[\frac{\lambda_{j}}
{\|\1bf_{B_{j}}\|_{\NHerzSo}}\right]^{s}
\1bf_{B_{j}}\right\}^{\frac{1}{s}}
\right\|_{\NHerzSo}\right\},
$$
where the infimum is taken over all
the decompositions of $f$ as above.
\end{definition}

Next, we give the atomic characterization
of the Hardy space $\NHaSaHo$ as follows.

\begin{theorem}\label{Atoln}
Let $p$, $q$, $\omega$, $d$, $s$, and $r$
be as in Definition \ref{atsln}.
Then $$\NHaSaHo=\NaHaSaHo$$ with equivalent quasi-norms.
\end{theorem}

To establish this atomic characterization,
we need the following two auxiliary
lemmas about the
Fefferman--Stein vector-valued
inequality on inhomogeneous
local generalized Herz spaces.

\begin{lemma}\label{vmbhln}
Let $p,\ q\in(0,\infty)$ and $\omega\in M(\rp)$ satisfy
$\mi(\omega)\in(-\frac{n}{p},\infty)$.
Then, for any given $u\in(1,\infty)$ and
$$r\in\left(0,\min\left\{p,\frac{n}{
\MI(\omega)+n/p}\right\}\right),$$ there exists
a positive constant $C$ such that, for any $\{f
_{j}\}_{j\in\mathbb{N}}\subset
L_{{\rm loc}}^{1}(\rrn)$,
\begin{equation*}
\left\|\left\{\sum_{j\in\mathbb{N}}
\left[\mc
(f_{j})\right]^{u}\right\}^{\frac{1}{u}}
\right\|_{[\NHerzSo]^{1/r}}
\leq C\left\|\left(
\sum_{j\in\mathbb{N}}|f_{j}|^{u}\right)^{
\frac{1}{u}}\right\|_{[\NHerzSo
]^{1/r}}.
\end{equation*}
\end{lemma}

\begin{proof}
Let all the symbols be as in the
present lemma, $u\in(1,\infty)$,
and $$r\in\left(0,\min\left\{p,\frac{n}{
\MI(\omega)+n/p}\right\}\right).$$
Then, repeating the proof of Lemma
\ref{vmbhl} with both
Theorem \ref{Th3.4} and Lemma \ref{convexl}
therein replaced, respectively,
by both Theorem \ref{Th3.4n}
and Lemma \ref{convexln},
we find that, for any $\{f
_{j}\}_{j\in\mathbb{N}}\subset
L_{{\rm loc}}^{1}(\rrn)$,
\begin{equation*}
\left\|\left\{\sum_{j\in\mathbb{N}}
\left[\mc
(f_{j})\right]^{u}\right\}^{\frac{1}{u}}
\right\|_{[\NHerzSo]^{1/r}}
\lesssim\left\|\left(
\sum_{j\in\mathbb{N}}|f_{j}|^{u}\right)^{
\frac{1}{u}}\right\|_{[\NHerzSo
]^{1/r}}.
\end{equation*}
This finishes the proof of Lemma \ref{vmbhln}.
\end{proof}

\begin{lemma}\label{Atolnl3}
Let $p$, $q\in(0,\infty)$,
$\omega\in M(\rp)$ satisfy
$\mi(\omega)\in(-\frac{n}{p},
\infty)$, $s\in(0,\infty)$,
and $$
\theta\in\left(0,\min\left\{
s,p,\frac{n}{\MI(\omega)+n/p}\right\}\right).
$$
Then there exists a positive
constant $C$ such that, for any
$\{f_j\}_{j\in\mathbb{N}}\subset
L^1_{\mathrm{loc}}(\rrn)$,
$$
\left\|
\left\{
\sum_{j\in\mathbb{N}}
\left[\mc^{(\theta)}(f_j)\right]^{s}
\right\}^{1/s}
\right\|_{\NHerzSo}
\leq C\left\|
\left(\sum_{j\in\mathbb{N}}
|f_j|^{s}
\right)^{1/s}
\right\|_{\NHerzSo}.
$$
\end{lemma}

\begin{proof}
Let all the symbols be as
in the present lemma. Then,
from Remark \ref{power}(ii) and
repeating an argument similar
to that used in the proof of
Lemma \ref{Atoll3} with
Lemma \ref{vmbhl} therein replaced
by Lemma \ref{vmbhln},
we infer that, for any
$\{f_j\}_{j\in\mathbb{N}}\subset
L^1_{\mathrm{loc}}(\rrn)$,
$$
\left\|
\left\{
\sum_{j\in\mathbb{N}}
\left[\mc^{(\theta)}(f_j)\right]^{s}
\right\}^{1/s}
\right\|_{\NHerzSo}
\lesssim\left\|
\left(\sum_{j\in\mathbb{N}}
|f_j|^{s}
\right)^{1/s}
\right\|_{\NHerzSo}.
$$
This finishes the
proof of Lemma \ref{Atolnl3}.
\end{proof}

Via both Lemmas \ref{vmbhln}
and \ref{Atolnl3} above, we
now show Theorem \ref{Atoln}.

\begin{proof}[Proof of Theorem
\ref{Atoln}]
Let all the symbols be as in the
present theorem. Then, applying both
Lemmas \ref{Atoll1} and \ref{Atogl4}
and repeating
the proof of Theorem \ref{Atol}
with Theorem \ref{Th3} and
Lemmas \ref{mbhal} and \ref{Atoll3}
therein replaced, respectively,
by Theorem \ref{Th3n} and
Lemmas \ref{mbhaln} and \ref{Atolnl3},
we conclude that
$$\NHaSaHo=\NaHaSaHo$$ with equivalent
quasi-norms. This finishes the proof of
Theorem \ref{Atoln}.
\end{proof}

Next, we turn to establish the atomic
characterization of the inhomogeneous
generalized Herz--Hardy space $\NHaSaH$.
To do this, we first introduce the following
definition of $(\NHerzS,\,r,\,d)$-atom.

\begin{definition}\label{deatomgn}
Let $p$, $q\in(0,\infty)$, $\omega\in
M(\rp)$ with $\MI(\omega)\in(-\infty,0)$,
$r\in[1,\infty]$, and $d\in\zp$.
Then a measurable function
$a$ on $\rrn$ is called a
\emph{$(\NHerzS,
\,r,\,d)$-atom}\index{$(\NHerzS,\,r,\,d)$-atom}
if there exists a ball $B\in\mathbb{B}$ such that
\begin{enumerate}
  \item[{\rm(i)}] $\supp(a)
  :=\{x\in\rrn:\ a(x)\neq0\}\subset B$;
  \item[{\rm(ii)}] $\|a\|_{L^r(\rrn)}
  \leq\frac{|B|^{1/r}}{\|\1bf_{B}\|_{\NHerzS}}$;
  \item[{\rm(iii)}] for any $\alpha\in\zp^n$
  with $|\alpha|\leq d$,
  $$\int_{\rrn}a(x)x^\alpha\,dx=0.$$
\end{enumerate}
\end{definition}

Via these atoms, we now give the following
concept of the inhomogeneous generalized atomic
Herz--Hardy space $\NaHaSaH$.

\begin{definition}\label{atsgn}
Let $p$, $q\in(0,\infty)$,
$\omega\in M(\rp)$ satisfy
$$-\frac{n}{p}<\mi(\omega)\leq\MI(\omega)<0,$$
$s\in(0,\min\{1,p,q\})$,
$d\geq\lfloor n(
1/s-1)\rfloor$ be a fixed integer, and
$$r\in\left(\max\left\{1,\frac{n}{
\mi(\omega)+n/p}\right\},\infty\right].$$
Then the \emph{inhomogeneous generalized atomic
Herz--Hardy space}\index{inhomogeneous generalized
atomic Herz--Hardy \\space} $\NaHaSaH$\index{$\NaHaSaH$}
is defined to be the set of all the
$f\in\mathcal{S}'(\rrn)$ such that there exist
$\{\lambda_j\}_{j\in\mathbb{N}}
\subset[0,\infty)$ and
a sequence $\{a_{j}\}_{j\in\mathbb{N}}$ of
$(\NHerzS,\,r,\,d)$-atoms supported,
respectively, in the balls $\{B_{j}\}_{j\in\mathbb{N}}
\subset\mathbb{B}$ such that
$$
f=\sum_{j\in\mathbb{N}}\lambda_{j}a_{j}
$$
in $\mathcal{S}'(\rrn)$,
and
$$
\left\|\left\{
\sum_{j\in\mathbb{N}}
\left[\frac{\lambda_j}{\|
\1bf_{B_j}\|_{\NHerzS}}\right]^s
\1bf_{B_j}\right\}^
{\frac{1}{s}}\right\|_{\NHerzS}<\infty.
$$
Moreover, for any $f\in\NaHaSaH$,
$$
\|f\|_{\NaHaSaH}:=\inf\left\{
\left\|\left\{\sum_{j\in\mathbb{N}}
\left[\frac{\lambda_{j}}
{\|\1bf_{B_{j}}\|_{\NHerzS}}\right]^{s}
\1bf_{B_{j}}\right\}^{\frac{1}{s}}
\right\|_{\NHerzS}\right\},
$$
where the infimum is taken over all
the decompositions of $f$ as above.
\end{definition}

Then we have the following atomic
characterization of the inhomogeneous
generalized Herz--Hardy space $\NHaSaH$.

\begin{theorem}\label{Atogn}
Let $p$, $q$, $\omega$, $r$, $s$, and $d$
be as in Definition \ref{atsgn}.
Then $$\NHaSaH=\NaHaSaH$$ with
equivalent quasi-norms.
\end{theorem}

To prove this theorem, we require some
preliminary lemmas about both
inhomogeneous global generalized Herz
spaces and inhomogeneous block spaces.
First, the following one is related to
the boundedness of powered Hardy--Littlewood
maximal operators on inhomogeneous block spaces.

\begin{lemma}\label{mbhagn}
Let $p,\ q\in(0,\infty)$ and
$\omega\in M(\rp)$ satisfy
$\mi(\omega)\in(-\frac{n}{p},\infty)$.
Then, for any given $s\in
(0,\min\{p,q,\frac{n}{\MI(\omega)+n/p}\})$
and $$r\in\left(\max\left\{p,
\frac{n}{\mi(\omega)+n/p}
\right\},\infty\right],$$
there exists a positive constant
$C$ such that,
for any $f\in L_{{\rm loc}}^{1}(\rrn)$,
\begin{equation*}
\left\|\mc^{((r/s)')}(f)\right\|_{\mathcal{B}
_{1/\omega^{s}}^{(p/s)',(q/s)'}(\rrn)}\leq C
\|f\|_{\mathcal{B}_
{1/\omega^{s}}^{(p/s)',(q/s)'}(\rrn)}.
\end{equation*}
\end{lemma}

\begin{proof}
Let all the symbols be as in
the present lemma,
$s\in(0,\min\{p,q,\frac{n}
{\MI(\omega)+n/p}\})$,
and $$r\in\left(\max\left\{p,
\frac{n}{\mi(\omega)+n/p}
\right\},\infty\right].$$
Then, using Lemma \ref{rela} and
repeating an argument similar
to that used in the proof
Lemma \ref{mbhag} with
Corollary \ref{maxbl} therein
replaced by Corollary \ref{maxbln},
we find that, for any $f\in L_{{\rm loc}}^{1}(\rrn)$,
\begin{equation*}
\left\|\mc^{((r/s)')}(f)\right\|_{\mathcal{B}
_{1/\omega^{s}}^{(p/s)',(q/s)'}(\rrn)}\lesssim
\|f\|_{\mathcal{B}_
{1/\omega^{s}}^{(p/s)',(q/s)'}(\rrn)}.
\end{equation*}
This finishes the proof of Lemma \ref{mbhagn}.
\end{proof}

We also need the following two
technical lemmas about
the Fefferman--Stein vector-valued
inequality on inhomogeneous global generalized
Herz spaces.

\begin{lemma}\label{vmbhgn}
Let $p,\ q\in(0,\infty)$ and $\omega\in M(\rp)$ satisfy
$\mi(\omega)\in(-\frac{n}{p},\infty)$.
Then, for any given $u\in(1,\infty)$ and
$$r\in\left(0,\min\left\{p,\frac{n}
{\MI(\omega)+n/p}\right\}\right),$$ there exists
a positive constant $C$ such that, for any $\{f
_{j}\}_{j\in\mathbb{N}}\subset L_{{\rm loc}}^{1}(\rrn)$,
\begin{equation*}
\left\|\left\{\sum_{j\in\mathbb{N}}\left[\mc
(f_{j})\right]^{u}\right\}^{\frac{1}{u}}
\right\|_{[\NHerzS]^{1/r}}
\leq C\left\|\left(
\sum_{j\in\mathbb{N}}|f_{j}|^{u}\right)^
{\frac{1}{u}}\right\|_{[\NHerzS]^{1/r}}.
\end{equation*}
\end{lemma}

\begin{proof}
Let all the symbols be as in the present lemma,
$u\in(1,\infty)$, and
$$r\in\left(0,\min\left\{p,\frac{n}
{\MI(\omega)+n/p}\right\}\right).$$
Then, repeating the proof of
Lemma \ref{vmbhg} with both
Theorem \ref{Th3.3} and
Lemma \ref{convexl} therein replaced,
respectively, by both Theorem \ref{Th3.3n} and
Lemma \ref{convexln},
we conclude that, for any $\{f
_{j}\}_{j\in\mathbb{N}}\subset
L_{{\rm loc}}^{1}(\rrn)$,
\begin{equation*}
\left\|\left\{\sum_{j\in\mathbb{N}}\left[\mc
(f_{j})\right]^{u}\right\}^{\frac{1}{u}}
\right\|_{[\NHerzS]^{1/r}}
\lesssim\left\|\left(
\sum_{j\in\mathbb{N}}|f_{j}|^{u}\right)^
{\frac{1}{u}}\right\|_{[\NHerzS]^{1/r}}.
\end{equation*}
This then finishes the proof of Lemma
\ref{vmbhgn}.
\end{proof}

\begin{lemma}\label{Atognl5}
Let $p$, $q\in(0,\infty)$,
$\omega\in M(\rp)$ satisfy
$\mi(\omega)\in(-\frac{n}{p},
\infty)$, $s\in(0,\infty)$,
and $$
\theta\in\left(0,\min\left\{
s,p,\frac{n}{\MI(\omega)+n/p}\right\}\right).
$$
Then there exists a positive
constant $C$ such that, for any
$\{f_j\}_{j\in\mathbb{N}}\subset
L^1_{\mathrm{loc}}(\rrn)$,
$$
\left\|
\left\{
\sum_{j\in\mathbb{N}}
\left[\mc^{(\theta)}(f_j)\right]^{s}
\right\}^{1/s}
\right\|_{\NHerzS}
\leq C\left\|
\left(\sum_{j\in\mathbb{N}}
|f_j|^{s}
\right)^{1/s}
\right\|_{\NHerzS}.
$$
\end{lemma}

\begin{proof}
Let all the symbols be
as in the present lemma.
Then, applying Remark \ref{power}(ii)
and repeating the proof of
Lemma \ref{Atogl5} via replacing
Lemma \ref{vmbhg} therein by
Lemma \ref{vmbhgn},
we find that,
for any $\{f_{j}\}_{
j\subset\mathbb{N}}
\subset L^1_{\mathrm{loc}}(\rrn)$,
$$
\left\|
\left\{
\sum_{j\in\mathbb{N}}
\left[\mc^{(\theta)}(f_j)\right]^{s}
\right\}^{1/s}
\right\|_{\NHerzS}
\lesssim\left\|
\left(\sum_{j\in\mathbb{N}}
|f_j|^{s}
\right)^{1/s}
\right\|_{\NHerzS}.
$$
This finishes the proof of
Lemma \ref{Atognl5}.
\end{proof}

Moreover, to
show Theorem \ref{Atogn}, we require the
following variant of
the H\"{o}lder inequality
of inhomogeneous global
generalized Herz spaces.

\begin{lemma}\label{Atognl1}
Let $p$, $q\in(1,\infty)$ and $\omega\in M(\rp)$.
Then there exists a positive constant $C$ such that,
for any $f$, $g\in\Msc(\rrn)$,
$$
\|fg\|_{L^1(\rrn)}\leq C\|f\|_{\NHerzS}
\|g\|_{\NBspace}.
$$
\end{lemma}

\begin{proof}
Let all the symbols be as in the
present lemma and $f$, $g\in\Msc(\rrn)$.
Then, repeating
an argument similar to that used
in the estimation of \eqref{ope}
with both \eqref{duale0} and
Lemma \ref{prel1}(ii) therein replaced,
respectively, by \eqref{dualne0} and
Lemma \ref{prenl1}(ii), we obtain
$$
\left\|fg\right\|_{L^1(\rrn)}
\lesssim\left\|f\right\|_{\NHerzS}
\left\|f\right\|_{\NBspace},
$$
which completes the proof of Lemma \ref{Atognl1}.
\end{proof}

In addition, the following
equivalent characterization
of inhomogeneous
global generalized Herz spaces
is also an essential tool
in the proof of Theorem \ref{Atogn}.

\begin{lemma}\label{Atognl2}
Let $p$, $q\in(1,\infty)$ and $\omega
\in M(\rp)$ satisfy
$\MI(\omega)\in(-\infty,0)$. Then
a measurable function $f$ belongs to
the inhomogeneous
global generalized Herz space $\NHerzS$ if and
only if
$$
\|f\|_{\NHerzS}^{\star}:=\sup
\left\{\|fg\|_{L^1(\rrn)}:\
\|g\|_{\NBspace}=1\right\}<\infty.
$$
Moreover, there exists a constant
$C\in[1,\infty)$ such that,
for any $f\in\NHerzS$,
$$
C^{-1}\|f\|_{\NHerzS}\leq
\|f\|_{\NHerzS}^{\star}\leq C\|f\|_{\NHerzS}.
$$
\end{lemma}

\begin{proof}
Let all the symbols be as in the present lemma
and $f\in\Msc(\rrn)$.
We first show the necessity. For
this purpose, assume $f\in\NHerzS$. Then,
repeating the estimation of \eqref{Atogl2e1} with
Lemma \ref{Atogl1} therein replaced
by Lemma \ref{Atognl1}, we obtain
\begin{equation}\label{Atognl2e1}
\|f\|_{\NHerzS}^{\star}
\lesssim\|f\|_{\NHerzS}<\infty.
\end{equation}
This then finishes the proof of
the necessity.

Conversely, we show the sufficiency.
To this end, assume $$\|f\|_
{\HerzS}^{\star}<\infty.$$
Then, repeating an argument similar
to that used in the estimation of
\eqref{Atogl2e3} via replacing both
Theorems \ref{Th2} and \ref{pre} therein,
respectively, by both Theorems \ref{Th2n} and
\ref{pren}, we find that
\begin{equation}\label{Atognl2e2}
\|f\|_{\NHerzS}\lesssim\|f\|_{\NHerzS}
^{\star}<\infty.
\end{equation}
This implies that $f\in\HerzS$ and
hence finishes the proof of the sufficiency.
Moreover, combining both \eqref{Atognl2e1} and
\eqref{Atognl2e2}, we conclude that
$$
\|f\|_{\NHerzS}\sim
\|f\|_{\NHerzS}^{\star}
$$
with the positive equivalence constants
independent of $f$.
This finishes the proof of Lemma \ref{Atognl2}.
\end{proof}

Via the above lemmas, we next
prove Theorem \ref{Atogn}.

\begin{proof}[Proof of Theorem
\ref{Atogn}]
Let all the symbols be as in
the present theorem. Then, repeating
the proof of Theorem \ref{Atog}
via replacing Theorem \ref{Th2}
and Lemmas \ref{mbhag}, \ref{Atogl5},
and \ref{Atogl2} therein, respectively,
by Theorem \ref{Th2n}
and Lemmas \ref{mbhagn}, \ref{Atognl5},
and \ref{Atognl2}, we obtain
$$\NHaSaH=\NaHaSaH$$
with equivalent quasi-norms. This then
finishes the proof of
Theorem \ref{Atogn}.
\end{proof}

\subsection{Finite Atomic Characterizations}

In this subsection, we establish the finite atomic
characterizations of both the inhomogeneous
generalized Herz--Hardy spaces
$\NHaSaHo$ and $\NHaSaH$. To this end,
we first introduce the following inhomogeneous
generalized finite atomic Herz--Hardy
space $\NiaHaSaHo$.

\begin{definition}\label{finiatn}
Let $p$, $q\in(0,\infty)$, $\omega\in M(\rp)$ with
$\mi(\omega)\in(-\frac{n}{p},\infty)$,
$$s\in\left(0,\min\left\{
1,p,q,\frac{n}{\MI(\omega)+n/p}\right\}\right),$$
$d\geq\lfloor n(
1/s-1)\rfloor$ be a fixed integer, and
$$r\in\left(\max\left\{1,p,\frac{n}
{\mi(\omega)+n/p}\right\},\infty\right].$$
Then the \emph{inhomogeneous generalized finite atomic
Herz--Hardy space\index{inhomogeneous generalized
finite atomic Herz--Hardy \\ space}
$\NiaHaSaHo$\index{$\NiaHaSaHo$}},
associated with $\NHerzSo$, is defined
to be the set of all
finite linear combinations of
$(\NHerzSo,\,r,\,d)$-atoms.
Moreover,
for any $f\in\NiaHaSaHo$,
$$
\|f\|_{\NiaHaSaHo}:=\left\{
\inf\left\|\left\{\sum_{j=1}
^{N}\left[\frac{\lambda_{j}}{
\|\1bf_{B_{j}}\|_{\NHerzSo}}
\right]^{s}\1bf_{
B_{j}}\right\}^{\frac{1}{s}}
\right\|_{\NHerzSo}\right\},
$$
where the infimum is taken over all
finite linear combinations of $f$,
namely, $N\in\mathbb{N}$, $f=\sum_{j=1}^{N}
\lambda_ja_j$,
$\{\lambda_j\}_{j=1}^N\subset[0,\infty)$,
and $\{a_j\}_{j=1}^N$
being $(\NHerzSo,\,r,\,\,d)$-atoms supported,
respectively, in the balls
$\{B_j\}_{j=1}^N\subset\mathbb{B}$.
\end{definition}

Then we establish the following finite
atomic characterization
of $\NHaSaHo$. Recall that $C(\rrn)$ denotes the set of all
continuous functions on $\rrn$.

\begin{theorem}\label{finiatomln}
Let $p$, $q$, $\omega$, $d$, $s$, and $r$ be as in Definition
\ref{finiatn}.
\begin{enumerate}
  \item[{\rm(i)}] If $$r\in\left(\max\left\{1,p,\frac{n}{
      \mi(\omega)+n/p}\right\},\infty\right),$$ then
  $\|\cdot\|_{\NiaHaSaHo}$ and
  $\|\cdot\|_{\NHaSaHo}$ are equivalent quasi-norms
on the inhomogeneous generalized finite
atomic Herz--Hardy space $\NiaHaSaHo$.
  \item[{\rm(ii)}] If $r=\infty$, then
  $\|\cdot\|_{H\Kmc_{\omega,\0bf,\mathrm{fin}
  }^{p,q,\infty,d,s}(\rrn)}$
  and $\|\cdot\|_{\NHaSaHo}$ are
  equivalent quasi-norms
  on $H\Kmc_{\omega,\0bf,\mathrm{fin}
  }^{p,q,\infty,d,s}(\rrn)
  \cap C(\rrn)$.
\end{enumerate}
\end{theorem}

\begin{proof}
Let all the symbols be as in
the present theorem. Then,
from Lemma \ref{finiatomll1} and
an argument similar to
that used in the proof of
Theorem \ref{finiatoml} with
both Lemmas \ref{mbhal} and \ref{Atoll3} therein
replaced, respectively, by
both Lemmas \ref{mbhaln} and \ref{Atolnl3},
we infer that
both (i) and (ii) of the present theorem
hold true, which completes
the proof of Theorem \ref{finiatomln}.
\end{proof}

Now, we turn to establish the finite
atomic characterization of the inhomogeneous
generalized Herz--Hardy space $\NHaSaH$.
For this purpose, we first
introduce the following finite atomic
Hardy spaces.

\begin{definition}\label{finiatgn}
Let $p$, $q\in(0,\infty)$, $\omega\in M(\rp)$ with
$$
-\frac{n}{p}<\mi(\omega)\leq
\MI(\omega)<0,
$$
$s\in(0,\min\{1,p,q\}),$
$d\geq\lfloor n(
1/s-1)\rfloor$ be a fixed integer,
and
$$r\in\left(\max\left\{1,\frac{n}{\mi(\omega)
+n/p}\right\},\infty\right].$$
Then the \emph{inhomogeneous generalized finite atomic
Herz--Hardy space\index{inhomogeneous generalized
finite atomic Herz--Hardy \\ space}
$\NiaHaSaH$\index{$\NiaHaSaH$}},
associated with $\NHerzS$, is defined
to be the set of all
finite linear combinations of
$(\NHerzS,\,r,\,d)$-atoms.
Moreover,
for any $f\in\NiaHaSaH$,
$$
\|f\|_{\NiaHaSaH}:=\left\{
\inf\left\|\left\{\sum_{j=1}
^{N}\left[\frac{\lambda_{j}}{
\|\1bf_{B_{j}}\|_{\NHerzS}}
\right]^{s}\1bf_{
B_{j}}\right\}^{\frac{1}{s}}
\right\|_{\NHerzS}\right\},
$$
where the infimum is taken over all
finite linear combinations of $f$,
namely, $N\in\mathbb{N}$, $f=\sum_{j=1}^{N}
\lambda_ja_j$,
$\{\lambda_j\}_{j=1}^N\subset[0,\infty)$,
and $\{a_j\}_{j=1}^N$
being $(\NHerzS,\,r,\,\,d)$-atoms supported,
respectively, in the balls
$\{B_j\}_{j=1}^N\subset\mathbb{B}$.
\end{definition}

Via this concept, we
show the finite atomic characterization
of $\NHaSaH$ as follows.

\begin{theorem}\label{finiatomgn}
Let $p$, $q$, $\omega$, $d$, $s$,
and $r$ be as in Definition
\ref{finiatgn}.
\begin{enumerate}
  \item[{\rm(i)}] If $$r\in\left(\max
  \left\{1,\frac{n}{\mi(\omega)+n/p}
  \right\},\infty\right),$$ then
  $\|\cdot\|_{\NiaHaSaH}$ and $\|\cdot\|_{\NHaSaH}$
  are equivalent quasi-norms
on the inhomogeneous generalized
finite atomic Herz--Hardy space $\NiaHaSaH$.
  \item[{\rm(ii)}] If $r=\infty$, then
  $\|\cdot\|_{H\Kmc_{\omega,\mathrm{fin}}
  ^{p,q,\infty,d,s}(\rrn)}$
  and $\|\cdot\|_{\NHaSaH}$ are equivalent quasi-norms
  on $H\Kmc_{\omega,\mathrm{fin}}
  ^{p,q,\infty,d,s}(\rrn)
  \cap C(\rrn)$.
\end{enumerate}
\end{theorem}

\begin{proof}
Let all the symbols be as in the
present lemma and $f\in\NiaHaSaH$.
Then, using Lemma \ref{finiatomgl1}
and repeating the proof of
Theorem \ref{finiatomg} via replacing both
Theorem \ref{Atog} and Lemma \ref{Atogl5}
therein, respectively, by both Theorem
\ref{Atogn} and Lemma \ref{Atognl5},
we conclude that, if
$r\in(1,\infty)$, then
$$
\left\|f\right\|_{\NiaHaSaH}
\sim\left\|f\right\|_{\NHaSaH}
$$
and, if $r=\infty$ and $f\in C(\rrn)$, then
$$
\left\|f\right\|_{H\Kmc_{\omega,\mathrm{fin}}
^{p,q,\infty,d,s}(\rrn)}
\sim\left\|f\right\|_{\NHaSaH}.
$$
These further imply that both
(i) and (ii) of the present theorem hold true,
and then finish the proof of Theorem
\ref{finiatomgn}.
\end{proof}

\subsection{Molecular Characterizations}

This subsection aims to establish the
molecular characterizations of inhomogeneous
generalized Herz--Hardy spaces.
To do this,
recall that, for any $j\in\mathbb{N}$,
$S_{j}(B):=(2^{j}B)\setminus(2^{j-1}B)$
and $S_{0}(B):=B$. Then we introduce the
$(\NHerzSo,\,r,\,d,\,\tau)$-molecules as follows.

\begin{definition}\label{molen}
Let $p$, $q\in(0,\infty)$, $\omega\in M(\rp)$,
$\tau\in(0,\infty)$,
$r\in[1,\infty]$, and $d\in\zp$.
Then a measurable function
$m$ on $\rrn$ is called a
\emph{$(\NHerzSo,\,r,
\,d,\,\tau)$-molecule}\index{$(\NHerzSo,\,r,\,d,\,\tau)$-molecule}
centered at a ball $B\in\mathbb{B}$ if
\begin{enumerate}
  \item[{\rm(i)}] for any $j\in\zp$,
  $$\left\|m\1bf_{S_{j}(B)}
  \right\|_{L^{r}(\rrn)}\leq2^{-\tau j}
\frac{|B|^{1/r}}{\|\1bf_{B}\|_{\NHerzSo}};$$
 \item[{\rm(ii)}] for any
 $\alpha\in\zp^n$ with $|\alpha|\leq d$,
  $$\int_{\rrn}m(x)x^\alpha\,dx=0.$$
\end{enumerate}
\end{definition}

We next show the following molecular
characterization of $\NHaSaHo$.

\begin{theorem}\label{Moleln}
Let $p,\ q\in(0,\infty)$,
$\omega\in M(\rp)$ with
$\mi(\omega)\in(-\frac{n}{p},\infty)$,
$$s\in\left(0,\min\left\{1,p,q,
\frac{n}{\MI(\omega)+n/p}\right\}\right),$$
$d\geq\lfloor n(1/s-1)\rfloor$ be a fixed integer,
$$r\in\left(\max\left\{1,p,\frac{n}
{\mi(\omega)+n/p}\right\},\infty\right],$$ and
$\tau\in(0,\infty)$ with $\tau>n(1/s-1/r)$.
Then $f\in\NHaSaHo$ if and only
if $f\in\mathcal{S}'(\rrn)$
and there exist a sequence
$\{m_{j}\}_{j\in\mathbb{N}}$
of $(\NHerzSo,\,r,\,d,\,\tau)$-molecules
centered, respectively,
at the balls $\{B_{j}\}_{j\in\mathbb{N}}
\subset\mathbb{B}$ and
a sequence $\{\lambda_{j}\}
_{j\in\mathbb{N}}\subset
[0,\infty)$ such that
$
f=\sum_{j\in\mathbb{N}}\lambda_{j}a_{j}
$
in $\mathcal{S}'(\rrn)$, and
$$
\left\|\left\{\sum_{j\in\mathbb{N}}
\left[\frac{\lambda_{j}}
{\|\1bf_{B_{j}}\|_{\NHerzSo}}
\right]^{s}\1bf_{B_{j}}
\right\}^{\frac{1}{s}}
\right\|_{\NHerzSo}<\infty.
$$
Moreover, there exists
a constant $C\in[1,\infty)$
such that, for
any $f\in\NHaSaHo$,
\begin{align*}
C^{-1}\|f\|_{\NHaSaHo}&\leq\inf\left\{
\left\|\left\{\sum_{j\in\mathbb{N}}
\left[\frac{\lambda_{i}}
{\|\1bf_{B_{j}}\|_{\NHerzSo}}\right]^s\1bf_
{B_{j}}\right\}^{\frac{1}{s}}\right\|_{\NHerzSo}
\right\}\\&\leq C\|f\|_{\NHaSaHo},
\end{align*}
where the infimum is taken
over all the decompositions of $f$ as above.
\end{theorem}

\begin{proof}
Let all the symbols be as in
the present theorem and $f\in\mathcal{S}'(\rrn)$.
Then, applying Lemma \ref{Molell1} and
repeating an argument similar to that used
in the proof of Theorem \ref{Molel} with
Theorem \ref{Th3} and Lemmas \ref{mbhal}
and \ref{Atoll3} therein replaced,
respectively, by
Theorem \ref{Th3n} and Lemmas \ref{mbhaln}
and \ref{Atolnl3}, we
find that $f\in\NHaSaHo$ if and only if
there exist a sequence
$\{m_{j}\}_{j\in\mathbb{N}}$
of $(\NHerzSo,\,r,\,d,\,\tau)$-molecules
centered, respectively,
at the balls $\{B_{j}\}_{j\in\mathbb{N}}
\subset\mathbb{B}$ and
a sequence $\{\lambda_{j}\}
_{j\in\mathbb{N}}\subset
[0,\infty)$ such that
$$
f=\sum_{j\in\mathbb{N}}\lambda_{j}a_{j}
$$
in $\mathcal{S}'(\rrn)$, and
$$
\left\|\left\{\sum_{j\in\mathbb{N}}
\left[\frac{\lambda_{j}}
{\|\1bf_{B_{j}}\|_{\NHerzSo}}
\right]^{s}\1bf_{B_{j}}
\right\}^{\frac{1}{s}}
\right\|_{\NHerzSo}<\infty,
$$
and, moreover,
$$
\left\|f\right\|_{\NHaSaHo}
\sim\inf\left\{
\left\|\left\{\sum_{j\in\mathbb{N}}
\left[\frac{\lambda_{j}}
{\|\1bf_{B_{j}}\|_{\NHerzSo}}
\right]^{s}\1bf_{B_{j}}
\right\}^{\frac{1}{s}}
\right\|_{\NHerzSo}\right\},
$$
where the infimum is taken over
all the decompositions of $f$ as
in the present theorem. This finishes
the proof of Theorem \ref{Moleln}.
\end{proof}

Now, we establish the molecular characterization
of the inhomogeneous
generalized Herz--Hardy space
$\NHaSaH$. For this purpose,
we first introduce the following
concept of $(\NHerzS,\,r,\,d,\,\tau)$-molecules.

\begin{definition}\label{molegn}
Let $p$, $q\in(0,\infty)$, $\omega\in M(\rp)$
with $\MI(\omega)\in(-\infty,0)$,
$\tau\in(0,\infty)$, $r\in[1,\infty]$,
and $d\in\zp$.
Then a measurable function $m$ on
$\rrn$ is called a
\emph{$(\NHerzS,\,r,\,d,\,\tau)
$-molecule}\index{$(\NHerzS,\,r,\,d,\,\tau)$-molecule}
centered at a ball $B\in\mathbb{B}$ if
\begin{enumerate}
  \item[{\rm(i)}] for any
  $j\in\zp$,
  $$
  \left\|m\1bf_{S_j(B)}\right\|_{L^r(\rrn)}
  \leq2^{-j\tau}\frac{|B|^{1/r}}
  {\|\1bf_B\|_{\NHerzS}};
  $$
 \item[{\rm(ii)}] for any $\alpha\in
 \zp^n$ with $|\alpha|\leq d$,
  $$\int_{\rrn}m(x)x^\alpha\,dx=0.$$
\end{enumerate}
\end{definition}

Via these molecules, we characterize the
Hardy space $\NHaSaH$ as follows.

\begin{theorem}\label{Molegn}
Let $p,\ q\in(0,\infty)$,
$\omega\in M(\rp)$ with
$$-\frac{n}{p}<\mi
(\omega)\leq\MI(\omega)<0,$$
$s\in(0,\min\{1,p,q\})$,
$d\geq\lfloor n(1/s-1)\rfloor$ be a
fixed integer,
$$r\in\left(\max\left\{1,
\frac{n}{\mi(\omega)+n/p}\right\},\infty\right],$$
and
$
\tau\in(n(
\frac{1}{s}-\frac{1}{r}),\infty).
$
Then $f\in\NHaSaH$ if and only if
$f\in\mathcal{S}'(\rrn)$
and there exist $\{\lambda_{j}\}
_{j\in\mathbb{N}}\subset[0,\infty)$ and
a sequence $\{m_{j}\}_{j\in\mathbb{N}}$
of $(\NHerzS,\,r,\,d,\,\tau)$-molecules
centered, respectively,
at the balls $\{B_{j}\}_{j\in\mathbb{N}}
\subset\mathbb{B}$
such that
\begin{equation*}
f=\sum_{j\in\mathbb{N}}\lambda_{j}m_{j}
\end{equation*}
in $\mathcal{S}'(\rrn)$, and
\begin{equation*}
\left\|\left\{\sum_{j\in\mathbb{N}}
\left[\frac{\lambda_{j}}
{\|\1bf_{B_{j}}\|_{\NHerzS}}
\right]^{s}\1bf_{B_{j}}\right\}^{\frac{1}{s}}
\right\|_{\NHerzS}<\infty.
\end{equation*}
Moreover, there exists a
constant $C\in[1,\infty)$ such that, for
any $f\in\NHaSaH$,
\begin{align*}
C^{-1}\|f\|_{\NHaSaH}&\leq\inf\left\{
\left\|\left\{\sum_{j\in\mathbb{N}}
\left[\frac{\lambda_{i}}
{\|\1bf_{B_{j}}\|_{\NHerzS}}\right]^s\1bf_
{B_{j}}\right\}^{\frac{1}{s}}\right\|_{\NHerzS}
\right\}\\&\leq C\|f\|_{\NHaSaH},
\end{align*}
where the infimum is taken
over all the decompositions of $f$ as above.
\end{theorem}

\begin{proof}
Let all the symbols be as in
the present theorem and $f\in\mathcal{S}'(\rrn)$.
Then, using Theorem \ref{Molegx} and
repeating the proof of Theorem
\ref{Atog} via replacing
Theorem \ref{Th2} and Lemmas \ref{mbhag},
\ref{Atogl5}, and \ref{Atogl2} therein,
respectively, by Theorem \ref{Th2n} and
Lemmas \ref{mbhagn},
\ref{Atognl5}, and \ref{Atognl2}, we
conclude that $f\in\NHaSaH$ if and only if
there exist a sequence
$\{m_{j}\}_{j\in\mathbb{N}}$
of $(\NHerzS,\,r,\,d,\,\tau)$-molecules
centered, respectively,
at the balls $\{B_{j}\}_{j\in\mathbb{N}}
\subset\mathbb{B}$ and
a sequence $\{\lambda_{j}\}
_{j\in\mathbb{N}}\subset
[0,\infty)$ such that
$
f=\sum_{j\in\mathbb{N}}\lambda_{j}a_{j}
$
in $\mathcal{S}'(\rrn)$, and
$$
\left\|\left\{\sum_{j\in\mathbb{N}}
\left[\frac{\lambda_{j}}
{\|\1bf_{B_{j}}\|_{\NHerzS}}
\right]^{s}\1bf_{B_{j}}
\right\}^{\frac{1}{s}}
\right\|_{\NHerzS}<\infty,
$$
and, moreover,
$$
\left\|f\right\|_{\NHaSaH}
\sim\inf\left\{
\left\|\left\{\sum_{j\in\mathbb{N}}
\left[\frac{\lambda_{j}}
{\|\1bf_{B_{j}}\|_{\NHerzS}}
\right]^{s}\1bf_{B_{j}}
\right\}^{\frac{1}{s}}
\right\|_{\NHerzS}\right\},
$$
where the infimum is taken over
all the decompositions of $f$ as
in the present theorem.
This then finishes the
proof of Theorem \ref{Molegn}.
\end{proof}

\subsection{Littlewood--Paley Function
Characterizations}

In this subsection, we establish the
equivalent characterizations of both the
inhomogeneous generalized
Herz--Hardy spaces $\NHaSaHo$ and $\NHaSaH$
via the Lusin area function,
the Littlewood--Paley $g$-function, and the Littlewood--Paley
$g_\lambda^*$-function defined as in Definitions \ref{df451}
and \ref{df452}.
First, we show these Littlewood--Paley
function characterizations of $\NHaSaHo$ as follows.

\begin{theorem}\label{lusinn}
Let $p,\ q\in(0,\infty)$, $\omega\in
M(\rp)$ satisfy
$\mi(\omega)\in(-\frac{n}{p},\infty)$,
$$
s_{0}:=\min\left\{1,p,q,\frac{n}{\MI(\omega)+n/p}\right\},
$$
and $\lambda\in(\max\{1,2/s_{0}\},\infty)$.
Then the following four statements are
mutually equivalent:
\begin{enumerate}
  \item[{\rm(i)}] $f\in\NHaSaHo$;
  \item[{\rm(ii)}] $f\in\mathcal{S}'(\rrn)$,
  $f$ vanishes weakly
  at infinity, and $S(f)\in\NHerzSo$;
  \item[{\rm(iii)}] $f\in\mathcal{S}'(\rrn)$,
  $f$ vanishes weakly
  at infinity, and $g(f)\in\NHerzSo$;
  \item[{\rm(iv)}] $f\in\mathcal{S}'(\rrn)$,
  $f$ vanishes weakly
  at infinity, and $g_{\lambda}^*(f)\in\NHerzSo$.
\end{enumerate}
Moreover, for any $f\in\NHaSaHo$,
$$
\|f\|_{\NHaSaHo}\sim\|S(f)\|_{\NHerzSo}
\sim\|g(f)\|_{\NHerzSo}\sim
\left\|g_{\lambda}^*(f)\right\|_{\NHerzSo},
$$
where the positive equivalence
constants are independent of $f$.
\end{theorem}

\begin{proof}
Let all the symbols be as in the
present theorem. Then, applying
Lemma \ref{lusinl1} and repeating
an argument similar to that used in
the proof of Theorem \ref{lusin}
with $\mw$, $\Mw$, Theorem \ref{Th3},
and Lemmas \ref{mbhal} and \ref{Atoll3}
therein replaced, respectively, by
$\mi(\omega)$, $\MI(\omega)$, Theorem \ref{Th3n},
and Lemmas \ref{mbhaln} and \ref{Atolnl3},
we find that (i), (ii), (iii), and (iv)
of the present theorem are mutually equivalent
and, for any $f\in\NHaSaHo$,
$$
\|f\|_{\NHaSaHo}\sim\|S(f)\|_{\NHerzSo}
\sim\|g(f)\|_{\NHerzSo}\sim
\left\|g_{\lambda}^*(f)\right\|_{\NHerzSo}
$$
with positive equivalence constants
independent of $f$. This finishes the proof of
Theorem \ref{lusinn}.
\end{proof}

\begin{remark}
Let $p\in(1,\infty)$, $q\in[1,\infty)$,
and $\omega(t):=t^\alpha$ for any $t\in(0,\infty)$
and for any given $\alpha\in
(-\frac{n}{p},\frac{n}{p'})$
in Theorem \ref{lusinn}. Then, by Remark \ref{sec101r}
and Theorem \ref{Th5.7n}, we find that,
in this case, the inhomogeneous
generalized Herz--Hardy space $\NHaSaHo$
coincides with the classical
inhomogeneous Herz space $K_p^{\alpha,q}(\rrn)$,
and hence Theorem \ref{lusinn}
goes back to \cite[Theorem 1.1.1]{LYH}.
\end{remark}

In addition, the following conclusion
gives various Littlewood--Paley function
characterizations of the inhomogeneous
generalized Herz--Hardy space $\NHaSaH$.

\begin{theorem}\label{lusingn}
Let $p,\ q\in(0,\infty)$, $\omega\in
M(\rp)$ satisfy
$$
-\frac{n}{p}<\mi(\omega)\leq\MI(\omega)<0,
$$
and $\lambda\in(2\max\{1,\frac{1}{p},
\frac{1}{q}\},\infty)$.
Then the following four statements are
mutually equivalent:
\begin{enumerate}
  \item[{\rm(i)}] $f\in\NHaSaH$;
  \item[{\rm(ii)}] $f\in\mathcal{S}'(\rrn)$,
  $f$ vanishes weakly
  at infinity, and $S(f)\in\NHerzS$;
  \item[{\rm(iii)}] $f\in\mathcal{S}'(\rrn)$,
  $f$ vanishes weakly
  at infinity, and $g(f)\in\NHerzS$;
  \item[{\rm(iv)}] $f\in\mathcal{S}'(\rrn)$,
  $f$ vanishes weakly
  at infinity, and $g_{\lambda}^*(f)\in\NHerzS$.
\end{enumerate}
Moreover, for any $f\in\NHaSaH$,
$$
\|f\|_{\NHaSaH}\sim\|S(f)\|_{\NHerzS}
\sim\|g(f)\|_{\NHerzS}\sim
\left\|g_{\lambda}^*(f)\right\|_{\NHerzS},
$$
where the positive equivalence
constants are independent of $f$.
\end{theorem}

\begin{proof}
Let all the symbols be as in
the present theorem. Then,
using Lemma \ref{lusingl1} and
repeating the proof of
Theorem \ref{lusing} via
replacing Theorems \ref{Th5.4}
and \ref{lusin} therein, respectively,
by Theorems \ref{Th5.4n}
and \ref{lusinn}, we conclude that
(i), (ii), (iii), and (iv) of the
present theorem are mutually
equivalent and, for any $f\in\NHaSaH$,
$$
\|f\|_{\NHaSaH}\sim\|S(f)\|_{\NHerzS}
\sim\|g(f)\|_{\NHerzS}\sim
\left\|g_{\lambda}^*(f)\right\|_{\NHerzS},
$$
where the positive equivalence
constants are independent of $f$. This
then finishes the proof of Theorem \ref{lusing}.
\end{proof}

\subsection{Dual Space of $\NHaSaHo$}

The main target of this subsection
is to establish the dual theorem of
the inhomogeneous generalized Herz--Hardy space
$\NHaSaHo$. To this end, we first introduce
the following Campanato-type function space
$\NHaSaHoc$.

\begin{definition}\label{camsn}
Let $p,\ q,\ s\in(0,\infty)$, $r\in[1,\infty)$, $d\in\mathbb{Z}_{+}$,
and $\omega\in M(\rp)$. Then the \emph{Campanato-type
function space}\index{Campanato-type function \\space}
$\NHaSaHoc$\index{$\NHaSaHoc$},
associated with the inhomogeneous local generalized Herz space $\NHerzSo$,
is defined to be the set of all the $f\in L^{r}_{\rm loc}(\rrn)$
such that
\begin{align*}
&\|f\|_{\NHaSaHoc}\\&\quad:=\sup\left\|
\left\{\sum_{i=1}^{m}\left[\frac{\lambda_{i}}{
\|\1bf_{B_{i}}\|_{\NHerzSo}}\right]^{s}\1bf_{B_{i}}\right\}^{\frac{1}{s}}
\right\|_{\NHerzSo}^{-1}\\
&\qquad\times\sum_{j=1}^{m}\left\{\frac{\lambda_{j}|B_{j}|}{
\|\1bf_{B_{j}}\|_{\NHerzSo}}\left[\frac{1}{|B_{j}|}\int_{
B_{j}}\left|f(x)-P_{B_{j}}^{d}f(x)\right|^{r}\,dx\right]
^{\frac{1}{r}}\right\}
\end{align*}
is finite,
where the supremum is taken over
all $m\in\mathbb{N}$, $\{B_{j}\}_{j=1}^{m}\subset\mathbb{B}$,
and $\{\lambda_{j}\}_{j=1}^{m}
\subset[0,\infty)$ with $\sum_{j=1}^{m}\lambda_{
j}\neq0$.
\end{definition}

\begin{remark}
Let all the symbols be as in Definition \ref{camsn}.
Then, obviously, $\mathcal{P}_d(\rrn)\subset
\NHaSaHoc$ and, for any $f\in\NHaSaHoc$,
$\|f\|_{\NHaSaHoc}=0$ if and only if $f\in\mathcal{P}_d(\rrn)$.
Thus, in what follows,
we always identify $f\in\NHaSaHoc$ with $\{f+P:\ P\in\mathcal{P}_d(\rrn)\}$.
\end{remark}

Then, using \cite[Remark 3.3(iii)
and Proposition 3.4]{HYYZ} with
$X$ therein replaced by $\NHerzSo$,
we immediately obtain the following
equivalent characterizations of the Campanato-type
function space $\NHaSaHoc$; we
omit the details.

\begin{proposition}
Let $p$, $q$, $\omega$, $r$, $d$, and
$s$ be as in Definition \ref{camsn}.
Then the following three statements are equivalent:
\begin{enumerate}
  \item[{\rm(i)}] $f\in\NHaSaHoc$;
  \item[{\rm(ii)}] $f\in L^r_\mathrm{loc}(\rrn)$ and
  \begin{align*}
  &\normmm{f}_{\NHaSaHoc}\\&\quad:=\sup\inf\left\|
\left\{\sum_{i=1}^{m}\left[\frac{\lambda_{i}}{
\|\1bf_{B_{i}}\|_{\NHerzSo}}\right]^{s}\1bf_{B_{i}}\right\}^{\frac{1}{s}}
\right\|_{\NHerzSo}^{-1}\\
&\qquad\times\sum_{j=1}^{m}\left\{\frac{\lambda_{j}|B_{j}|}{
\|\1bf_{B_{j}}\|_{\NHerzSo}}\left[\frac{1}{|B_{j}|}\int_{
B_{j}}\left|f(x)-P(x)\right|^{r}\,dx\right]
^{\frac{1}{r}}\right\}
\end{align*}
is finite, where the supremum is the same as in Definition \ref{camsn}
and the infimum is taken over all $P\in\mathcal{P}_d(\rrn)$;
  \item[{\rm(iii)}] $f\in L^r_\mathrm{loc}(\rrn)$ and
  \begin{align*}
&\widetilde{\|f\|}_{\NHaSaHoc}\\&\quad:=\sup\left\|
\left\{\sum_{i\in\mathbb{N}}\left[\frac{\lambda_{i}}{
\|\1bf_{B_{i}}\|_{\NHerzSo}}\right]^{s}\1bf_{B_{i}}\right\}^{\frac{1}{s}}
\right\|_{\NHerzSo}^{-1}\\
&\qquad\times\sum_{j\in\mathbb{N}}\left\{\frac{\lambda_{j}|B_{j}|}{
\|\1bf_{B_{j}}\|_{\NHerzSo}}\left[\frac{1}{|B_{j}|}\int_{
B_{j}}\left|f(x)-P_{B_{j}}^{d}f(x)\right|^{r}\,dx\right]
^{\frac{1}{r}}\right\}
\end{align*}
is finite, where the supremum is taken over all $\{B_j\}_{
j\in\mathbb{N}}\subset\mathbb{B}$ and $\{\lambda_j\}_{j\in\mathbb{N}}
\subset[0,\infty)$ satisfying
$$
\left\|
\left\{\sum_{i\in\mathbb{N}}\left[\frac{\lambda_{i}}{
\|\1bf_{B_{i}}\|_{\NHerzSo}}\right]^{s}\1bf_{B_{i}}\right\}^{\frac{1}{s}}
\right\|_{\NHerzSo}\in(0,\infty).
$$
\end{enumerate}
Moreover, there exists a constant $C\in[1,\infty)$ such that,
for any $f\in L^r_\mathrm{loc}(\rrn)$,
$$
C^{-1}\|f\|_{\NHaSaHoc}\leq\normmm{f}_{\NHaSaHoc}\leq
C\|f\|_{\NHaSaHoc}
$$
and
$$
\widetilde{\|f\|}_{\NHaSaHoc}=\|f\|_{\NHaSaHoc}.
$$
\end{proposition}

Next, we give the following dual theorem
which shows that
the dual space of the inhomogeneous
generalized Herz--Hardy space
$\NHaSaHo$ is just the Campanato-type
function space $\NHaSaHod$.

\begin{theorem}\label{dualhn}
Let $p$, $q\in(0,\infty)$, $\omega\in M(\rp)$ with
$\mi(\omega)\in(-\frac{n}{p},\infty)$,
$$p_{-}:=\min\left\{1,p,\frac{n}
{\MI(\omega)+n/p}\right\},$$
$d\geq\lfloor n(1/p_{-}-1)\rfloor$ be a
fixed integer, $s\in(0,\min\{p_{-},q\})$, and $$r\in\left(\max\left\{1,p,\frac{n}{\mi(\omega)+n/p}
\right\},\infty\right].$$
Then $\NHaSaHod$ is the dual space of
$\NHaSaHo$ in the following sense:
\begin{enumerate}
  \item[{\rm(i)}]
  Let $g\in\NHaSaHod$. Then the linear
  functional
  \begin{equation}\label{dualh2n}
  L_{g}:\ f\mapsto L_{g}(f):=\int_{\rrn}f(x)g(x)\,dx,
  \end{equation}
  initially defined for any $f\in\NiaHaSaHo$,
  has a bounded extension to
  the inhomogeneous generalized Herz--Hardy space $\NHaSaHo$.
  \item[{\rm(ii)}] Conversely, any continuous linear functional
  on $\NHaSaHo$ arises as in \eqref{dualh2n} with a unique
  $g\in\NHaSaHod$.
\end{enumerate}
Moreover, there exists a constant $C\in[1,\infty)$ such that,
for any $g\in\NHaSaHod$,
$$
C^{-1}\|g\|_{\NHaSaHod}\leq\|L_{g}\|_{(\NHaSaHo)^*}
\leq C\|g\|_{\NHaSaHod},
$$
where $(\NHaSaHo)^*$ denotes the dual space of $\NHaSaHo$.
\end{theorem}

\begin{proof}
Let all the symbols be as in
the present theorem. Then, from Lemma
\ref{dualhl1} and an argument similar
to that used in the proof of Theorem \ref{dualh} with
Theorem \ref{Th3} and
Lemmas \ref{mbhal} and \ref{vmbhl} therein replaced,
respectively, by Theorem \ref{Th3n} and
Lemmas \ref{mbhaln} and \ref{vmbhln}, we infer
that both (i) and (ii) of the present theorem
hold true. This finishes
the proof of Theorem \ref{dualhn}.
\end{proof}

Finally, applying this dual theorem,
we immediately obtain
the following equivalence of
Campanato-type
function spaces $\NHaSaHoc$; we omit the details.

\begin{corollary}
Let $p$, $q$, $\omega$, $p_-$, $d$,
and $s$ be as in Theorem \ref{dualhn},
$$
p_+:=\max\left\{1,p,\frac{n}{\mi(\omega)+n/p}\right\},
$$
$r\in[1,p_+')$, $d_0:=\lfloor n(1/p_--1)\rfloor$, and $s_0
\in(0,\min\{p_-,q\})$. Then
$$
\NHaSaHoc=\mathcal{L}_{\omega,\0bf}^{p,q,1,d_0,s_0}(\rrn)
$$
with equivalent quasi-norms.
\end{corollary}

\subsection{Boundedness of Calder\'{o}n--Zygmund Operators}

Let $d\in\zp$. Recall that
the standard kernel $K$ and
the $d$-order Calder\'{o}n--Zygmund
operator with kernel $K$ are defined,
respectively, in Definitions \ref{def-s-k}
and \ref{defin-C-Z-s}.
In this subsection, we establish
the boundedness of
$d$-order Calder\'{o}n--Zygmund
operators on inhomogeneous
generalized Herz--Hardy spaces.
First, the following conclusion
gives the boundedness of these operators
on the inhomogeneous generalized Herz--Hardy
space $\NHaSaHo$.

\begin{theorem}\label{bddoln}
Let $d\in\zp$, $\delta\in(0,1]$,
$p$, $q\in(\frac{n}{n+d+\delta},\infty)$,
$K$ be a $d$-order standard kernel defined
as in Definition \ref{def-s-k}, $T$ a
$d$-order Calder\'{o}n--Zygmund operator
with kernel $K$ having the vanishing moments
up to order $d$, and $\omega\in
M(\rp)$ with
$$
-\frac{n}{p}<\mi(\omega)\leq
\MI(\omega)<n-\frac{n}{p}+d+\delta.
$$
Then $T$ has a unique extension
on $\NHaSaHo$ and there exists
a positive constant $C$ such that,
for any $f\in\NHaSaHo$,
$$
\|T(f)\|_{\NHaSaHo}\leq C\|f\|_{
\NHaSaHo}.
$$
\end{theorem}

\begin{proof}
Let all the symbols be
as in the present theorem. Then,
applying Proposition \ref{bddoxl4}
and repeating the proof of Theorem
\ref{bddol} with $\mw$, $\Mw$, Theorems \ref{Th3}
and \ref{abso}, and Lemmas \ref{mbhal} and \ref{Atoll3}
therein replaced, respectively, by
$\mi(\omega)$, $\MI(\omega)$, Theorems \ref{Th3n} and
\ref{abson}, and Lemmas \ref{mbhaln} and \ref{Atolnl3},
we find that $T$ has a unique extension on
$\NHaSaHo$ and, for any $f\in\NHaSaHo$,
$$
\left\|T(f)\right\|_{\NHaSaHo}
\lesssim\left\|f\right\|_{\NHaSaHo}.
$$
This then finishes the proof of Theorem
\ref{bddoln}.
\end{proof}

\begin{remark}
We should point out that, in Theorem
\ref{bddoln}, when $d=0$, $p\in(1,\infty)$,
and $\omega(t):=t^{\alpha}$ for any $t\in(0,\infty)$
and for any given $\alpha\in[n(1-\frac{1}{p}),
n(1-\frac{1}{p})+\delta)$,
then this theorem goes back to \cite[Theorem 1]{ll97}.
\end{remark}

On the other hand, we prove the following
boundedness of Calder\'{o}n--Zygmund
operators on the Hardy space $\NHaSaH$.

\begin{theorem}\label{bddogn}
Let $d\in\zp$, $\delta\in(0,1]$,
$p$, $q\in(\frac{n}
{n+d+\delta},\infty)$,
$K$ be a
$d$-order standard kernel defined as in
Definition \ref{def-s-k}, $T$ a
$d$-order Calder\'{o}n--Zygmund operator
with kernel $K$ having the vanishing moments
up to order $d$, , and $\omega\in
M(\rp)$ satisfy
$$
-\frac{n}{p}<\mi(\omega)\leq
\MI(\omega)<0.
$$
Then $T$ can be extended into a bounded
linear operator on $\NHaSaH$, namely,
there exists
a positive constant $C$ such that,
for any $f\in\NHaSaH$,
$$
\|T(f)\|_{\NHaSaH}\leq C\|f\|_{
\NHaSaH}.
$$
\end{theorem}

\begin{proof}
Let all the symbols be as
in the present theorem. Then, from
Theorem \ref{bddox} and an argument similar to
that used in the proof of Theorem
\ref{bddog} with $\mw$, $\Mw$,
Theorem \ref{Th2},
Lemmas \ref{ee}, \ref{equa}, \ref{Atogl5},
and \ref{Atogl2}, and Corollary \ref{maxbl}
therein replaced, respectively, by
$\mi(\omega)$, $\MI(\omega)$, Theorem \ref{Th2n},
Lemmas \ref{een}, \ref{equan}, \ref{Atognl5},
and \ref{Atognl2}, and Corollary \ref{maxbln},
it follows that $T$ is well defined on
$\NHaSaH$ and, for any $f\in\NHaSaH$,
$$
\left\|T(f)\right\|_{\NHaSaH}
\lesssim\left\|f\right\|_{\NHaSaH}.
$$
This finishes the proof of Theorem
\ref{bddogn}.
\end{proof}

\subsection{Fourier Transform}

The target of this subsection is to
investigate the Fourier transform properties
of distributions in the
inhomogeneous generalized Herz--Hardy spaces
$\NHaSaHo$ and $\NHaSaH$. We first study
the Fourier transform in $\NHaSaHo$ as follows.

\begin{theorem}\label{Fouriertln}
Let $p,\ q\in(0,1]$, $\omega\in M(\rp)$ with
$\mi(\omega)\in(0,\infty)$,
and $p_{-}\in(0,\frac{n}{\MI(\omega)+n/p})$.
Then, for any $f\in\NHaSaHo$, there
exists a continuous function
$g$ on $\rrn$ such that $\widehat{f}=g$
in $\mathcal{S}'(\rrn)$
and
$$
\lim\limits_{|x|\to0^{+}}\frac{
|g(x)|}{|x|^{n(\frac{1}{p_-}-1)}}=0.
$$
Moreover, there exists a positive
constant $C$, independent of both $f$ and $g$,
such that, for any $x\in\rrn$,
$$
|g(x)|\leq C\|f\|_{\NHaSaHo}\max
\left\{1,|x|^{n(\frac{1}{p_-}-1)}\right\}
$$
and
$$
\int_{\rrn}|g(x)|\min\left\{|x|^{-\frac{n}{p_{-}}},
|x|^{-n}\right\}\,dx\leq C\|f\|_{\NHaSaHo}.
$$
\end{theorem}

In order to prove this theorem, we first
establish the following technique lemma
about the quasi-norms
of the characteristic function
of balls on $\NHerzSo$.

\begin{lemma}\label{Fouriertll1n}
Let $p,\ q\in(0,1]$, $\omega\in M(\rp)$ with
$\mi(\omega)\in(0,\infty)$,
and $p_{-}\in(0,\frac{n}{\MI(\omega)+n/p})$.
Then there exists a positive constant $C$
such that, for any $B\in\mathbb{B}$,
\begin{equation*}
\|\1bf_{B}\|_{\NHerzSo}\geq C\min
\left\{|B|,|B|^{\frac{1}{p_{-}}}\right\}.
\end{equation*}
\end{lemma}

\begin{proof}
Let all the symbols be as in the present lemma.
We first claim that, for any $B(x_0,2^{k_0})\in\mathbb{B}$ with
$x_0\in\rrn$ and $k_0\in\mathbb{Z}\cap(-\infty,0)$,
$$
\left|B(x_0,2^{k_0})\right|^{\frac{1}{p_-}}\lesssim
\left\|\1bf_{B(x_0,2^{k_0})}\right\|_{\NHerzSo}.
$$
To show this inequality,
we consider the following five cases on $x_0$.

\emph{Case 1)} $x_0=\0bf$. In this case,
notice that $B(\0bf,2^{k_0})\subset B(\0bf,1)$.
By this and Definition \ref{igh}(i),
we find that
\begin{align}\label{Fouriertll1ne2}
\left|B(\0bf,2^{k_0})\right|^{\frac{1}{p}}=
\left\|\1bf_{B(\0bf,2^{k_0})}\right\|_{L^p(\rrn)}
\leq\left\|\1bf_{B(\0bf,2^{k_0})}\right\|_{\NHerzSo}.
\end{align}
On the other hand,
from the assumption $\mi(\omega)\in(0,\infty)$
and Remark \ref{remark 1.1.3}(iii), it follows that
\begin{align}\label{Fouriertll1ne3}
p_-<\frac{n}{\MI(\omega)+n/p}
\leq\frac{n}{\mi(\omega)+n/p}<p,
\end{align}
which, combined with \eqref{Fouriertll1ne2},
further implies that
\begin{align}\label{Fouriertll1ne4}
\left|B(\0bf,2^{k_0})\right|^{\frac{1}{p_-}}<
\left|B(\0bf,2^{k_0})\right|^{\frac{1}{p}}
\leq\left\|\1bf_{B(\0bf,2^{k_0})}\right\|_{\NHerzSo}.
\end{align}
This finishes the proof of the above claim in this case.

\emph{Case 2)} $|x_0|\in(0,3\cdot2
^{k_0-1})$. In this case,
applying both \eqref{fouriertle11}
and the assumption $k_0\in\mathbb{Z}
\cap(-\infty,0)$, we find that
\begin{align*}
B\left(\frac{3\cdot2^{k_0-1}}{|x_0|
}x_0,2^{k_0-2}\right)
&\subset\left[B(x_0,2^{k_0})\cap
\left(B(\0bf,2^{k_0})\setminus
B(\0bf,2^{k_0-1})\right)\right]\\
&\subset\left[B(x_0,2^{k_0})\cap B(\0bf,1)\right].
\end{align*}
By this, \eqref{Fouriertll1ne3},
and Definition \ref{igh}(i),
we conclude that
\begin{align}\label{Fouriertll1ne5}
\left|B(x_0,2^{k_0})\right|^{\frac{1}{p_-}}
&<\left|B(x_0,2^{k_0})\right|^{\frac{1}{p}}
\sim\left|B\left(\frac{3\cdot2^{k_0-1}}
{|x_0|}x_0,2^{k_0-2}\right)\right|^
{\frac{1}{p}}\notag\\
&\lesssim\left|B(x_0,2^{k_0})\cap
B(\0bf,1)\right|^{\frac{1}{p}}
\lesssim\left\|\1bf_{B(x_0,2^{k_0})
}\right\|_{\NHerzSo},
\end{align}
which implies that the above claim holds
true in this case.

\emph{Case 3)} $|x_0|\in[3\cdot2^{k_0-1},
2^{k_0+1})$. In this case,
for any $y\in B(x_0,2^{k_0})$, we have
\begin{align*}
|y|\leq|y-x_0|+|x_0|<2^{k_0}+2^{k_0+1}<2.
\end{align*}
This implies that $B(x_0,2^{k_0})\subset
B(\0bf,2)$. Therefore, from \eqref{Fouriertll1ne3}
and Definition \ref{igh}(i), we deduce that
\begin{align}\label{Fouriertll1ne6}
\left|B(x_0,2^{k_0})\right|^{\frac{1}{p_-}}
&<\left|B(x_0,2^{k_0})\right|^{\frac{1}{p}}\notag\\
&\lesssim\left|B(x_0,2^{k_0})\cap
B(\0bf,1)\right|^{\frac{1}{p}}\notag\\
&\quad+\omega(2)\left|B(x_0,2^{k_0})\cap
\left(B(\0bf,2)\setminus B(\0bf,1)\right)\right|^{\frac{1}{p}}\notag\\
&\lesssim\left\|\1bf_{B(x_0,2^{k_0})}\right\|_{\NHerzSo},
\end{align}
which completes the proof of the above claim in this case.

\emph{Case 4)} $|x_0|\in[2^{k},2^{k+1})$ with $k\in\mathbb{Z}
\cap[k_0+1,0]$. In this case,
for any $y\in B(x_0,2^{k_0})$, we have
\begin{align*}
|y|\leq|y-x_0|+|x_0|<2^{k_0}+2^{k+1}<2^2,
\end{align*}
which implies that $B(x_0,2^{k_0})\subset
B(\0bf,2^2)$. Thus, using \eqref{Fouriertll1ne3} and Definition \ref{igh}(i),
we conclude that
\begin{align}\label{Fouriertll1ne7}
\left|B(x_0,2^{k_0})\right|^{\frac{1}{p_-}}
&<\left|B(x_0,2^{k_0})\right|^{\frac{1}{p}}\notag\\
&\lesssim\left|B(x_0,2^{k_0})\cap
B(\0bf,1)\right|^{\frac{1}{p}}\notag\\
&\quad+\omega(2)\left|B(x_0,2^{k_0})\cap
\left(B(\0bf,2)\setminus B(\0bf,1)\right)\right|^{\frac{1}{p}}\notag\\
&\quad+\omega(2^2)\left|B(x_0,2^{k_0})\cap
\left(B(\0bf,2^2)\setminus B(\0bf,2)\right)\right|^{\frac{1}{p}}\notag\\
&\lesssim\left\|\1bf_{B(x_0,2^{k_0})}\right\|_{\NHerzSo}.
\end{align}
This implies that, in this case, the above claim holds true.

\emph{Case 5)} $|x_0|\in[2^{k},2^{k+1})$ with $k\in\mathbb{N}$.
In this case, by
the assumption $k_0\in\mathbb{Z}
\cap(-\infty,0)$ and \eqref{fouriertle15}, we obtain
$$
B(x_0,2^{k_0})\subset\left[B(\0bf,2^{k+2})
\setminus B(\0bf,2^{k-1})\right].
$$
This, together with \eqref{Fouriertll1ne3}
and Definition \ref{igh}(i), further implies that
\begin{align*}
\left|B(x_0,2^{k_0})\right|^{\frac{1}{p_-}}
&<\left|B(x_0,2^{k_0})\right|^{\frac{1}{p}}\\
&\lesssim\omega(\tk)\left|B(x_0,2^{k_0})\cap
\left(B(\0bf,2^k)\setminus B(\0bf,2^{k-1})\right)\right|^{\frac{1}{p}}\\
&\quad+\omega(2^{k+1})\left|B(x_0,2^{k_0})\cap
\left(B(\0bf,2^{k+1})\setminus B(\0bf,2^k)\right)\right|^{\frac{1}{p}}\\
&\quad+\omega(2^{k+2})\left|B(x_0,2^{k+2})\cap
\left(B(\0bf,2^{k+1})\setminus B(\0bf,2)\right)\right|^{\frac{1}{p}}\\
&\lesssim\left\|\1bf_{B(x_0,2^{k_0})}\right\|_{\NHerzSo}.
\end{align*}
From this, \eqref{Fouriertll1ne4},
\eqref{Fouriertll1ne5}, \eqref{Fouriertll1ne6},
and \eqref{Fouriertll1ne7}, it follows that the
above claim holds true.

In addition, for any $k_0\in\zp$, repeating
an argument similar
to that used in
the estimation of \eqref{fouriertle16}
with $r$, $\mw$, and Lemma \ref{La3.5}
therein replaced, respectively, by
$2^{k_0}$, $\mi(\omega)$, and Lemma \ref{Th1},
we conclude that
\begin{equation}\label{Fouriertll1ne1}
\left|B(x_0,2^{k_0})\right|\lesssim\left\|
\1bf_{B(x_0,2^{k_0})}\right\|_{\NHerzS}.
\end{equation}
Assume $B(x_0,r)\in\mathbb{B}$ with $x_0\in\rrn$
and $r\in(0,\infty)$. Then there exists a $k\in\mathbb{Z}$
such that $r\in[2^{k},2^{k+1})$.
Thus, by the above claim and
\eqref{Fouriertll1ne1},
we find that
\begin{align*}
&\min\left\{\left|B(x_0,r)\right|,\left|
B(x_0,r)\right|^{\frac{1}{p_-}}\right\}\\
&\quad\sim\min\left\{\left|B(x_0,2^{k_0})\right|,\left|
B(x_0,2^{k_0})\right|^{\frac{1}{p_-}}\right\}\\
&\quad\lesssim\left\|\1bf_{B(x_0,2^{k_0})}\right\|_{\NHerzSo}
\lesssim\left\|\1bf_{B(x_0,r)}\right\|_{\NHerzSo}.
\end{align*}
This finishes the proof of Lemma \ref{Fouriertll1n}.
\end{proof}

Via the above lemma, we now prove
Theorem \ref{Fouriertln}.

\begin{proof}[Proof of Theorem \ref{Fouriertln}]
Let all the symbols be as in
the present theorem and $f\in\NHaSaHo$.
Then, using Lemma \ref{fouriertll2}
and repeating the proof of Theorem
\ref{Fouriertl} with
$\Mw$, Theorem \ref{Th3}, and Lemmas
\ref{vmbhl} and \ref{Fouriertll1}
therein replaced, respectively, by
$\MI(\omega)$, Theorem \ref{Th3n}, and Lemmas
\ref{vmbhln} and \ref{Fouriertll1n},
we conclude that there
exists a continuous function
$g$ on $\rrn$ such that $\widehat{f}=g$
in $\mathcal{S}'(\rrn)$,
$$
\lim\limits_{|x|\to0^{+}}\frac{
|g(x)|}{|x|^{n(\frac{1}{p_-}-1)}}=0,
$$
and, for any $x\in\rrn$,
$$
|g(x)|\lesssim\|f\|_{\NHaSaHo}\max
\left\{1,|x|^{n(\frac{1}{p_-}-1)}\right\}
$$
and
$$
\int_{\rrn}|g(x)|\min\left\{|x|^{-\frac{n}{p_{-}}},
|x|^{-n}\right\}\,dx\lesssim\|f\|_{\NHaSaHo},
$$
where the implicit positive constants
are independent of both $f$ and $g$.
This finishes the proof of Theorem \ref{Fouriertln}.
\end{proof}

Next, we investigate the Fourier transform of the Hardy space
$\NHaSaH$. Indeed, we have the following
conclusion.

\begin{theorem}\label{Fouriertn}
Let $p\in(0,1),\ q\in(0,1]$,
$\omega\in M(\rp)$ with
$$
n\left(1-\frac{1}{p}\right)<
\mi(\omega)\leq\MI(\omega)<0,
$$
and $p_{-}\in(0,p)$.
Then, for any $f\in\NHaSaH$, there
exists a continuous function
$g$ on $\rrn$ such that $\widehat{f}=g$
in $\mathcal{S}'(\rrn)$,
and
$$
\lim\limits_{|x|\to0^{+}}\frac{|g(x)|}
{|x|^{n(\frac{1}{p_-}-1)}}=0.
$$
Moreover, there exists a positive
constant $C$, independent of both $f$ and $g$,
such that, for any $x\in\rrn$,
$$
|g(x)|\leq C\|f\|_{\NHaSaH}\max\left\{1,
|x|^{n(\frac{1}{p_-}-1)}\right\}
$$
and
$$
\int_{\rrn}|g(x)|\min\left\{|x|^{-\frac{n}{p_{-}}},
|x|^{-n}\right\}\,dx\leq C\|f\|_{\NHaSaH}.
$$
\end{theorem}

To show this theorem,
we require the following auxiliary
estimate for the quasi-norm
$\|\cdot\|_{\NHaSaH}$ of the characteristic
function of balls.

\begin{lemma}\label{estafn}
Let $p\in(0,1),\ q\in(0,1]$, $\omega\in M(\rp)$ with
$$
n\left(1-\frac{1}{p}\right)<\mi(\omega)\leq\MI(\omega)<0,
$$
and $p_{-}\in(0,p)$.
Then there exists a positive constant $C$ such that,
for any $B\in\mathbb{B}$,
\begin{equation*}
\|\1bf_{B}\|_{\NHerzS}
\geq C\min\left\{|B|,|B|^{\frac{1}{p_{-}}}\right\}.
\end{equation*}
\end{lemma}

\begin{proof}
Let all the symbols be as in the present lemma and
$B(x_0,r)\in\mathbb{B}$ with $x_0\in\rrn$ and $r\in(0,\infty)$.
We show Lemma \ref{estafn} by considering the following
three cases on $r$.

\emph{Case 1)} $r\in(0,1]$. In this case, we have
$B(x_0,r)\subset B(x_0,1)$. From this, the assumption
$p_-\in(0,p)$, Definition \ref{igh}(i), and
Remark \ref{ighr}(i), we deduce that
\begin{align}\label{estafne1}
\left|B(x_0,r)\right|^{\frac{1}{p_-}}
&\leq\left|B(x_0,r)\right|^{\frac{1}{p}}
\sim\left\|\1bf_{B(x_0,r)}\right\|_{L^p(\rrn)}\notag\\
&\lesssim\left\|\1bf_{B(x_0,r)}(\cdot+x_0)\right\|_{\NHerzSo}
\lesssim\left\|\1bf_{B(x_0,r)}\right\|_{\NHerzS}.
\end{align}
This implies that Lemma \ref{estafn} holds true in this case.

\emph{Case 2)} $r\in(1,2]$. In this case,
we have
$B(x_0,1)\subset B(x_0,r)$.
Applying this, the H\"{o}lder inequality,
Definition \ref{igh}(i), and Remark \ref{ighr}(i), we
conclude that
\begin{align}\label{estafne2}
\left|B(x_0,r)\right|&\sim\left|B(x_0,1)\right|
\lesssim\left\|\1bf_{B(x_0,1)}\right\|_{L^p(\rrn)}
\lesssim\left\|\1bf_{B(x_0,1)}(\cdot+x_0)\right\|_{\NHerzSo}\notag\\
&\lesssim\left\|\1bf_{B(x_0,1)}\right\|_{\NHerzS}
\lesssim\left\|\1bf_{B(x_0,r)}\right\|_{\NHerzS},
\end{align}
which completes the proof of Lemma \ref{estafn} in this case.

\emph{Case 3)} $r\in(2,\infty)$. In this case,
there exists a $k\in\mathbb{N}$
such that $r\in(\tk,\tka]$.
This implies that
\begin{equation}\label{estafne3}
\left[B(x_0,\tk)\setminus B(x_0,\tkm)\right]\subset B(x_0,r).
\end{equation}
In addition, let $\eps\in(0,-n+\mi(\omega)+\frac{n}{p})$
be a fixed positive constant. Then, using \eqref{estafne3},
Lemmas \ref{Th1} and \ref{mono}, Definition \ref{igh}(i),
and Remark \ref{ighr}(i), we find that
\begin{align*}
\left|B(x_0,r)\right|&\sim r^n
\lesssim r^{\mi(\omega)+\frac{n}{p}-\eps}
\lesssim r^{\frac{n}{p}}\omega(r)
\sim2^{\frac{nk}{p}}\omega(\tk)\\
&\sim\omega(\tk)\left\|\1bf_{B(x_0,\tk)
\setminus B(x_0,\tkm)}\right\|_{L^p(\rrn)}\\
&\lesssim\left\|\1bf_{B(x_0,r)}(\cdot+x_0)
\right\|_{\NHerzSo}\notag\\
&\lesssim\left\|\1bf_{B(x_0,r)}\right\|_{\NHerzS}.
\end{align*}
This, together with \eqref{estafne1} and \eqref{estafne2},
further implies that
$$
\min\left\{\left|B(x_0,r)\right|,
\left|B(x_0,r)\right|^{\frac{1}{p_{-}}}\right\}\lesssim
\left\|\1bf_{B(x_0,r)}\right\|_{\NHerzS},
$$
which completes the proof of the present lemma
in this case, and hence the whole proof
of Lemma \ref{estafn}.
\end{proof}

We now prove Theorem \ref{Fouriertn}
via Lemma \ref{estafn}.

\begin{proof}[Proof of Theorem \ref{Fouriertn}]
Let all the symbols be as
in the present theorem and $f\in\NHaSaH$.
Then, from Lemma \ref{fouriertll2} and
an argument similar to
that used in the proof of
Theorem \ref{Fouriert} via replacing Theorems
\ref{Th2} and Lemmas \ref{vmbhg} and
\ref{estaf} therein, respectively, by Theorem
\ref{Th2n} and Lemmas \ref{vmbhgn} and
\ref{estafn}, we infer that
there exists a continuous function
$g$ on $\rrn$ such that $\widehat{f}=g$
in $\mathcal{S}'(\rrn)$,
$$
\lim\limits_{|x|\to0^{+}}\frac{|g(x)|}
{|x|^{n(\frac{1}{p_-}-1)}}=0,
$$
and, for any $x\in\rrn$,
$$
|g(x)|\lesssim\|f\|_{\NHaSaH}\max\left\{1,
|x|^{n(\frac{1}{p_-}-1)}\right\}
$$
and
$$
\int_{\rrn}|g(x)|\min\left\{|x|^{-\frac{n}{p_{-}}},
|x|^{-n}\right\}\,dx\lesssim\|f\|_{\NHaSaH}
$$
with the implicit positive constants
independent of both $f$ and $g$.
This then finishes the proof of
Theorem \ref{Fouriertn}.
\end{proof}

\section{Inhomogeneous Localized Generalized \\
Herz--Hardy Spaces}

In this section, we introduce inhomogeneous
localized generalized Herz--Hardy spaces and then
establish their maximal function,
atomic, molecular as well
as various Littlewood--Paley function
characterizations.
As applications, we also obtain the boundedness
of the pseudo-differential operators on these
localized Hardy spaces as well as
investigate the relation between
inhomogeneous generalized Herz--Hardy spaces and
inhomogeneous localized generalized Herz--Hardy
spaces.

Recall that, for any given
$N\in\mathbb{N}$ and for any $f\in\mathcal{S}'(\rrn)$,
the local grand maximal function $m_N(f)$
of $f$ is defined as
in \eqref{sec7e1}. Then we introduce the
inhomogeneous
localized generalized Herz--Hardy spaces
via this maximal function
as follows.\index{localized generalized Herz--Hardy space}

\begin{definition}\label{lghhn}
Let $p$, $q\in(0,\infty)$, $\omega\in
M(\rp)$, and $N\in\mathbb{N}$.
Then
\begin{enumerate}
  \item[{\rm(i)}] the \emph{inhomogeneous
  local generalized Herz--Hardy space}
  $\NHaSaHol$\index{$\NHaSaHol$}, associated with
  the inhomogeneous local generalized Herz space $\NHerzSo$,
  is defined to be the set of all
  the $f\in\mathcal{S}'(\rrn)$ such that $$
\|f\|_{\NHaSaHol}:=\|m_{N}(f)\|_{\NHerzSo}<\infty;$$
  \item[{\rm(ii)}] the \emph{inhomogeneous
  local generalized Herz--Hardy
  space} $\NHaSaHl$\index{$\NHaSaHl$}, associated with
  the inhomogeneous global
  generalized Herz space $\NHerzS$,
  is defined to be the set of all the
  $f\in\mathcal{S}'(\rrn)$ such that $$
\|f\|_{\NHaSaHl}:=\|m_{N}(f)\|_{\NHerzS}<\infty.$$
\end{enumerate}
\end{definition}

\begin{remark}
In Definition \ref{lghhn},
for any given $\alpha\in\rr$
and for any $t\in(0,\infty)$, let
$\omega(t):=t^\alpha$. Then, in this case,
the inhomogeneous local generalized
Herz--Hardy space $\NHaSaHol$
goes back to the classical \emph{inhomogeneous local
Herz-type Hardy space}\index{inhomogeneous local Herz-type Hardy \\ space}
$hK_p^{\alpha,q}(\rrn)$\index{$hK_p^{\alpha,q}(\rrn)$}
which was originally introduced in
\cite[Definition 1.2]{fy97} (see also
\cite[Section 2.6]{LYH}).
\end{remark}

\subsection{Maximal Function Characterizations}

In this subsection, we characterize the inhomogeneous
localized generalized Herz--Hardy
spaces via various maximal functions.
Recall that several localized radial and
localized non-tangential maximal functions
are given in Definition \ref{df511}.
Using these maximal functions, we first
establish the following
maximal function characterizations of
the local Hardy space $\NHaSaHol$.

\begin{theorem}\label{Th6.1n}
Let $p,\ q,\ a,\ b\in(0,\infty),\ \omega\in M(\rp),\
N\in\mathbb{N}$,
and $\phi\in\mathcal{S}(\rrn)$
satisfy $\int_{\rrn}
\phi(x)\,dx\neq0.$
\begin{enumerate}
\item[\rm{(i)}] Let $N\in\mathbb{N}
\cap[\lfloor b+1\rfloor,\infty)$.
Then, for any $f\in\mathcal{S}'(\rrn)$,
$$
\|m(f,\phi)\|_{\NHerzSo}\lesssim
\|m^*_{a}(f,\phi)\|_{\NHerzSo}
\lesssim\|m^{**}_{b}
(f,\phi)\|_{\NHerzSo},$$
\begin{align*}
\|m(f,\phi)\|_{\NHerzSo}&
\lesssim\|m_{N}(f)\|_{\NHerzSo}
\lesssim\|m_{\lfloor b+1\rfloor}(f)\|_{\NHerzSo}\\
&\lesssim\|m^{**}_{b}(f,\phi)\|_{\NHerzSo},
\end{align*}
and
$$
\|m^{**}_{b}(f,\phi)\|_{\NHerzSo}
\sim\|m^{**}_{b,N}(f)\|_{\NHerzSo},
$$
where the implicit positive
constants are independent of $f$.
\item[\rm{(ii)}] Let $\omega$ satisfy
$\mi(\omega)\in(-\frac{n}{p},\infty)$,
and $$b\in\left(2\max\left\{\frac{n}{p},
\frac{n}{q},\MI(\omega)+\frac{n}
{p}\right\},\infty\right).$$
Then, for any $f\in\mathcal{S}'(\rrn),$
$$\|m^{**}_{b}(f,\phi)\|_{\NHerzSo}
\lesssim\|m(f,\phi)\|_{\NHerzSo},$$
where the implicit positive constant
is independent of $f$. In particular,
when $N\in\mathbb{N}\cap[\lfloor
b+1\rfloor,\infty)$, if one of the quantities
$$\|m(f,\phi)\|_{\NHerzSo},\
\|m^*_{a}(f,\phi)\|_{\NHerzSo},\ \|m_{N}
(f)\|_{\NHerzSo},$$
$$\|m^{**}_{b}(f,\phi)\|_{\NHerzSo},
\ \text{and }\|m^{**}_{b,N}(f)\|_{\NHerzSo}$$
is finite, then the others are also finite
and mutually equivalent with the
positive equivalence constants independent of $f$.
\end{enumerate}
\end{theorem}

\begin{proof}
Let all the symbols be as in the
present theorem. Then, applying
both Lemmas \ref{Th6.1l2} and \ref{Th6.1l1}
and repeating the proof of Theorem \ref{Th6.1}
via replacing $\Mw$, Theorems \ref{Th3} and
\ref{strc}, and Lemma \ref{mbhl} therein,
respectively, by
$\MI(\omega)$, Theorems \ref{Th3n} and
\ref{strcn}, and Lemma \ref{mbhln},
we find that both (i) and (ii) of the present
theorem hold true. This then finishes the
proof of Theorem \ref{Th6.1n}.
\end{proof}

\begin{remark}
In Theorem \ref{Th6.1n}, let $p\in(1,\infty)$
and, for any given
$\alpha\in[n(1-\frac{1}{p}),\infty)$ and
for any $t\in(0,\infty)$,
$\omega(t):=t^\alpha$. Then, in this case,
the conclusion obtained in Theorem \ref{Th6.1n} goes back
to \cite[Theorem 2.2]{wl15}.
\end{remark}

Next, we turn to prove
the maximal function characterizations
of the inhomogeneous local generalized Herz--Hardy space
$\NHaSaHl$ as follows.

\begin{theorem}\label{Th6.3n}
Let $p,\ q,\ a,\ b\in(0,\infty),\ \omega\in M(\rp),\
N\in\mathbb{N}$,
and $\phi\in\mathcal{S}(\rrn)$
satisfy $\int_{\rrn}
\phi(x)\,dx\neq0.$
\begin{enumerate}
\item[\rm{(i)}] Let $N\in\mathbb{N}\cap
[\lfloor b+1\rfloor,\infty)$
and $\omega$ satisfy $\MI(\omega)\in(-\infty,0)$.
Then, for any $f\in\mathcal{S}'(\rrn)$,
$$
\|m(f,\phi)\|_{\NHerzS}\lesssim\|
m^*_{a}(f,\phi)\|_{\NHerzS}
\lesssim\|m^{**}_{b}(f,\phi)\|_{\NHerzS},$$
\begin{align*}
\|m(f,\phi)\|_{\NHerzS}&
\lesssim\|m_{N}(f)\|_{\NHerzS}
\lesssim\|m_{\lfloor b+1\rfloor}(f)\|_{\NHerzS}\\
&\lesssim\|m^{**}_{b}(f,\phi)\|_{\NHerzS},
\end{align*}
and
$$
\|m^{**}_{b}(f,\phi)\|_{\NHerzS}
\sim\|m^{**}_{b,N}(f)\|_{\NHerzS},
$$
where the implicit positive
constants are independent of $f$.
\item[\rm{(ii)}] Let $\omega$ satisfy
$$
-\frac{n}{p}<\mi(\omega)\leq\MI(\omega)<0,
$$
and $b\in(2n\max\{\frac{1}{p},\frac{1}{q}\},\infty)$.
Then, for any $f\in\mathcal{S}'(\rrn),$
$$\|m^{**}_{b}(f,\phi)\|_{\NHerzS}
\lesssim\|m(f,\phi)\|_{\NHerzS},$$
where the implicit positive constant
is independent of $f$. In particular,
when $N\in\mathbb{N}\cap[\lfloor
b+1\rfloor,\infty)$, if one of the quantities
$$\|m(f,\phi)\|_{\NHerzS},\ \|m^*_{a}
(f,\phi)\|_{\NHerzS},\ \|m_{N}
(f)\|_{\NHerzS},$$
$$\|m^{**}_{b}(f,\phi)\|_{\NHerzS},
\ \text{and }\|m^{**}_{b,N}(f)\|_{\NHerzS}$$
is finite, then the others are also finite
and mutually equivalent with the
positive equivalence constants independent of $f$.
\end{enumerate}
\end{theorem}

\begin{proof}
Let all the symbols be as in the
present theorem. Then, from Lemmas
\ref{Th6.1l2} and \ref{Th6.1l1}
and an argument similar to that
used in the proof of Theorem \ref{Th6.3} with
$\Mw$, Theorems \ref{Th2} and \ref{Th3.1},
and Lemma \ref{mbhg} therein replaced,
respectively, by
$\MI(\omega)$, Theorems \ref{Th2n} and \ref{Th3.1n},
and Lemma \ref{mbhgn}, we deduce that
both (i) and (ii) of the present theorem hold true.
This then finishes the proof of
Theorem \ref{Th6.3n}.
\end{proof}

\subsection{Relations with
Inhomogeneous Generalized \\ Herz--Hardy Spaces}

The main target of
this subsection is to investigate the
relation between
inhomogeneous generalized Herz--Hardy spaces and inhomogeneous
localized generalized Herz--Hardy spaces.
First, we give the following relation
between $\NHaSaHo$ and $\NHaSaHol$ associated
with the inhomogeneous local generalized
Herz space $\HerzSo$.

\begin{theorem}\label{relaln}
Let $p,\ q\in(0,\infty)$, $\omega\in M(\rp)$ satisfy
$\mi(\omega)\in(-\frac{n}{p},\infty)$, and $\varphi
\in\mathcal{S}(\rrn)$ satisfy
$$\1bf_{B(\0bf,1)}
\leq\widehat{\varphi}
\leq\1bf_{B(\0bf,2)}.$$
Then there exists a
constant $C\in[1,\infty)$ such that,
for any $f\in\mathcal{S}'(\rrn)$,
\begin{align*}
C^{-1}\|f\|_{\NHaSaHol}
&\leq\left\|f\ast\varphi\right\|_{\NHerzSo}
+\left\|f-f\ast\varphi\right\|_{\NHaSaHo}\\
&\leq C\|f\|_{\NHaSaHol}.
\end{align*}
\end{theorem}

\begin{proof}
Let $p$, $q$, $\omega$, and $\varphi$
be as in the present theorem and
$f\in\mathcal{S}'(\rrn)$. Then, using
Theorem \ref{relax}
and repeating the proof of Theorem \ref{relal}
with $\Mw$, Theorems \ref{Th3} and
\ref{strc}, and Lemma \ref{mbhl} therein
replaced, respectively, by
$\MI(\omega)$, Theorems \ref{Th3n} and
\ref{strcn}, and Lemma \ref{mbhln},
we obtain
$$
\left\|f\right\|_{\NHaSaHol}
\sim\left\|f\ast\varphi\right\|_{\NHerzSo}
+\left\|f-f\ast\varphi\right\|_{\NHaSaHo}
$$
with positive equivalence constants independent of
$f$. This then finishes the
proof of Theorem \ref{relaln}.
\end{proof}

On the other hand,
we obtain the relation
between the inhomogeneous generalized
Herz--Hardy space $\NHaSaH$ and the inhomogeneous
local generalized Herz--Hardy space $\NHaSaHl$
as follows.

\begin{theorem}\label{relagn}
Let $p,\ q\in(0,\infty)$, $\omega\in M(\rp)$ satisfy
$$-\frac{n}{p}<\mi(\omega)
\leq\MI(\omega)<0,$$
and $\varphi\in\mathcal{S}(\rrn)$ satisfy
$$\1bf_{B(\0bf,1)}
\leq\widehat{\varphi}
\leq\1bf_{B(\0bf,2)}.$$
Then there exists a
constant $C\in[1,\infty)$ such that,
for any $f\in\mathcal{S}'(\rrn)$,
\begin{align*}
C^{-1}\|f\|_{\NHaSaHl}
&\leq\left\|f\ast\varphi\right\|_{\NHerzS}
+\left\|f-f\ast\varphi\right\|_{\NHaSaH}\\
&\leq C\|f\|_{\NHaSaHl}.
\end{align*}
\end{theorem}

\begin{proof}
Let all the symbols be
as in the present theorem.
Then, applying Theorem \ref{relax}
and repeating an argument
similar to that used in
the proof of Theorem \ref{relag}
with $\Mw$, Theorems \ref{Th2} and
\ref{Th3.1}, and Lemma \ref{mbhg} therein
replaced, respectively, by
$\MI(\omega)$, Theorems \ref{Th2n} and
\ref{Th3.1n}, and Lemmas \ref{mbhgn},
we find that, for any $f\in\mathcal{S}'(\rrn)$,
$$
\left\|f\right\|_{\NHaSaHl}
\sim\left\|f\ast\varphi\right\|_{\NHerzS}
+\left\|f-f\ast\varphi\right\|_{\NHaSaH}.
$$
This then finishes the proof of
Theorem \ref{relagn}.
\end{proof}

\subsection{Atomic Characterizations}

In this subsection, we establish the atomic
characterizations of the inhomogeneous local
generalized Herz--Hardy spaces
$\NHaSaHol$ and $\NHaSaHl$. For this purpose, we first
introduce both local-$(\NHerzSo,\,r,\,d)$-atoms
and local atomic Hardy spaces associated
with the inhomogeneous local generalized
Herz space $\NHerzSo$
as follows.

\begin{definition}\label{loatomn}
Let $p$, $q\in(0,\infty)$, $\omega\in M(\rp)$,
$r\in[1,\infty]$, and $d\in\zp$. Then a measurable
function $a$ is called
a \emph{local-$(\NHerzSo,
\,r,\,d)$-atom}\index{local-$(\NHerzSo,\,r,\,d)$-atom}
if
\begin{enumerate}
  \item[{\rm(i)}] there exists a ball
  $B(x_0,r_0)\in\mathbb{B}$, with $x_0\in\rrn$
  and $r_0\in(0,\infty)$, such that
  $$\supp(a):=\left\{x\in\rrn:\
  a(x)\neq0\right\}\subset B(x_0,r_0);$$
  \item[{\rm(ii)}] $\|a\|_{L^r(\rrn)}
  \leq\frac{|B(x_0,r_0)|^{1/r}}{\|\1bf_{B
  (x_0,r_0)}\|_{\NHerzSo}}$;
  \item[{\rm(iii)}] when $r_0\in(0,1)$, then,
  for any $\alpha\in\zp^n$
  such that $|\alpha|\leq d$,
  $$\int_{\rrn}a(x)x^\alpha\,dx=0.$$
\end{enumerate}
\end{definition}

\begin{definition}\label{latsln}
Let $p$, $q\in(0,\infty)$, $\omega\in M(\rp)$ with
$\mi(\omega)\in(-\frac{n}{p},\infty)$,
$$
s\in\left(0,\min\left\{1,p,q,
\frac{n}{\MI(\omega)+n/p}\right\}\right),
$$
$d\geq\lfloor n(1/s-1)\rfloor$ be a fixed integer,
and
$$
r\in\left(\max\left\{1,p,\frac{n}
{\mi(\omega)+n/p}\right\},\infty\right].
$$
Then the \emph{inhomogeneous local generalized
atomic Herz--Hardy space}
\index{inhomogeneous local generalized atomic Herz--Hardy \\ space}
$\NaHaSaHol$\index{$\NaHaSaHol$},
associated with the inhomogeneous local
generalized Herz space $\NHerzSo$,
is defined to be the set of
all the $f\in\mathcal{S}'(\rrn)$ such that there exist
a sequence $\{a_{j}\}_{j\in\mathbb{N}}$ of
local-$(\NHerzSo,\,r,\,d)$-atoms supported,
respectively, in the balls
$\{B_{j}\}_{j\in\mathbb{N}}\subset\mathbb{B}$
and a sequence $\{\lambda_{j}
\}_{j\in\mathbb{N}}\subset[0,\infty)$ such that
$$
f=\sum_{j\in\mathbb{N}}\lambda_{j}a_{j}
$$
in $\mathcal{S}'(\rrn)$, and
$$
\left\|\left\{\sum_{j\in\mathbb{N}}
\left[\frac{\lambda_{j}}
{\|\1bf_{B_{j}}\|_{\NHerzSo}}
\right]^{s}\1bf_{B_{j}}\right\}^{\frac{1}{s}}
\right\|_{\NHerzSo}<\infty.
$$
Moreover, for any $f\in\NaHaSaHol$,
$$
\|f\|_{\NaHaSaHol}:=\inf\left\{
\left\|\left\{\sum_{j\in\mathbb{N}}
\left[\frac{\lambda_{j}}
{\|\1bf_{B_{j}}\|_{\NHerzSo}}
\right]^{s}\1bf_{B_{j}}\right\}^{\frac{1}{s}}
\right\|_{\NHerzSo}\right\},
$$
where the infimum is taken over
all the decompositions of $f$ as above.
\end{definition}

Then we have the following atomic characterization of
$\NHaSaHol$.

\begin{theorem}\label{lMoleln}
Let $p$, $q$, $\omega$, $d$, $s$, and $r$ be as
in Definition \ref{latsln}.
Then $$\NHaSaHol=\NaHaSaHol$$ with equivalent quasi-norms.
\end{theorem}

\begin{proof}
Let $p$, $q$, $\omega$, $d$, $s$, and $r$
be as in Definition \ref{latsln}. Then,
from Lemma \ref{atomxl} and
an argument similar to
that used in the proof of Theorem
\ref{lMolel} via replacing Theorem
\ref{Th3} and Lemmas \ref{mbhal} and \ref{Atoll3}
therein, respectively, by
Theorem \ref{Th3n} and Lemmas
\ref{mbhaln} and \ref{Atolnl3}, we infer
that
$$\NHaSaHol=\NaHaSaHol$$ with equivalent quasi-norms.
This then finishes the proof of
Theorem \ref{lMoleln}.
\end{proof}

In the remainder of this subsection,
we turn to prove the atomic characterization
of the inhomogeneous local generalized
Herz--Hardy space $\NHaSaHl$. To do this,
we first introduce the following
definition of local-$(\NHerzS,\,r,\,d)$-atoms.

\begin{definition}\label{goatomn}
Let $p$, $q\in(0,\infty)$, $\omega\in M(\rp)$
with $\MI(\omega)\in(-\infty,0)$,
$r\in[1,\infty]$, and $d\in\zp$. Then a measurable
function $a$ is called
a \emph{local-$(\NHerzS,
\,r,\,d)$-atom}\index{local-$(\NHerzS,\,r,\,d)$-atom}
if
\begin{enumerate}
  \item[{\rm(i)}] there exists a ball
  $B(x_0,r_0)\in\mathbb{B}$, with $x_0\in\rrn$
  and $r_0\in(0,\infty)$, such that
  $$\supp(a):=\left\{x\in\rrn:\
  a(x)\neq0\right\}\subset B(x_0,r_0);$$
  \item[{\rm(ii)}] $\|a\|_{L^r(\rrn)}
  \leq\frac{|B(x_0,r_0)|^{1/r}}{\|\1bf_{B
  (x_0,r_0)}\|_{\NHerzS}}$;
  \item[{\rm(iii)}] when $r_0\in(0,1)$, then,
  for any $\alpha\in\zp^n$
  such that $|\alpha|\leq d$,
  $$\int_{\rrn}a(x)x^\alpha\,dx=0.$$
\end{enumerate}
\end{definition}

Via these atoms, we now introduce
the following atomic Hardy spaces.

\begin{definition}\label{latsgn}
Let $p$, $q\in(0,\infty)$, $\omega\in M(\rp)$ with
$$-\frac{n}{p}<
\mi(\omega)\leq\MI(\omega)<0,$$
$
s\in(0,\min\{1,p,q\}),
$
$d\geq\lfloor n(1/s-1)\rfloor$ be a fixed integer,
and
$$
r\in\left(\max\left\{1,\frac{n}
{\mi(\omega)+n/p}\right\},\infty\right].
$$
Then the \emph{inhomogeneous local generalized
atomic Herz--Hardy space}
\index{inhomogeneous local generalized atomic Herz--Hardy \\ space}
$\NaHaSaHl$\index{$\NaHaSaHl$},
associated with the inhomogeneous global
generalized Herz space $\NHerzS$,
is defined to be the set of
all the $f\in\mathcal{S}'(\rrn)$ such that there exist
a sequence $\{a_{j}\}_{j\in\mathbb{N}}$ of
local-$(\NHerzS,\,r,\,d)$-atoms supported,
respectively, in the balls
$\{B_{j}\}_{j\in\mathbb{N}}\subset\mathbb{B}$
and a sequence $\{\lambda_{j}
\}_{j\in\mathbb{N}}\subset[0,\infty)$ such that
$
f=\sum_{j\in\mathbb{N}}\lambda_{j}a_{j}
$
in $\mathcal{S}'(\rrn)$, and
$$
\left\|\left\{\sum_{j\in\mathbb{N}}
\left[\frac{\lambda_{j}}
{\|\1bf_{B_{j}}\|_{\NHerzS}}
\right]^{s}\1bf_{B_{j}}\right\}^{\frac{1}{s}}
\right\|_{\NHerzS}<\infty.
$$
Moreover, for any $f\in\NaHaSaHl$,
$$
\|f\|_{\NaHaSaHl}:=\inf\left\{
\left\|\left\{\sum_{j\in\mathbb{N}}
\left[\frac{\lambda_{j}}
{\|\1bf_{B_{j}}\|_{\NHerzS}}
\right]^{s}\1bf_{B_{j}}\right\}^{\frac{1}{s}}
\right\|_{\NHerzS}\right\},
$$
where the infimum is taken over
all the decompositions of $f$ as above.
\end{definition}

Then we establish the following atomic
characterization\index{atomic characterization} of the
inhomogeneous local
generalized Herz--Hardy space $\HaSaHl$.

\begin{theorem}\label{lMolegn}
Let $p$, $q$, $\omega$, $d$, $s$, and $r$ be
as in Definition \ref{latsgn}.
Then $$\NHaSaHl=\NaHaSaHl$$ with equivalent quasi-norms.
\end{theorem}

\begin{proof}
Let all the symbols be as in
the present theorem. Then, applying
Theorems \ref{atomxxl}
and \ref{balln} and repeating
the proof of Theorem \ref{Atog}
with Theorem \ref{Th2}
and Lemmas \ref{mbhag}, \ref{Atogl5},
and \ref{Atogl2} therein replaced, respectively,
by Theorem \ref{Th2n}
and Lemmas \ref{mbhagn}, \ref{Atognl5},
and \ref{Atognl2}, we conclude that
$$\NHaSaHl=\NaHaSaHl$$
with equivalent quasi-norms.
This then
finishes the proof of
Theorem \ref{lMolegn}.
\end{proof}

\subsection{Molecular Characterizations}

This subsection is devoted to
establishing the molecular
characterizations of both the inhomogeneous local
generalized Herz--Hardy spaces
$\NHaSaHol$ and $\NHaSaHl$. To achieve this, we
first show the molecular
characterization of $\NHaSaHol$ via
introducing the following
local-$(\NHerzSo,\,r,\,d,\,\tau)$-molecules.

\begin{definition}\label{lomolen}
Let $p$, $q\in(0,\infty)$, $\omega\in M(\rp)$,
$r\in[1,\infty]$,
$d\in\mathbb{Z}_{+}$, and $\tau\in
(0,\infty)$. Then a measurable
function $m$ on $\rrn$ is called a
\emph{local-$(\NHerzSo,\,r,\,d,\,\tau)$-molecule}
\index{local-$(\NHerzSo,\,r,\,d,
\,\tau)$-molecule}centered at a ball
$B(x_0,r_0)\in\mathbb{B}$, with $x_0\in
\rrn$ and $r_0\in(0,\infty)$, if
\begin{enumerate}
  \item[{\rm(i)}] for any $i\in\mathbb{Z}_{+}$,
$$\left\|m\1bf_{S_{i}(B(x_0,r_0))}
\right\|_{L^{r}(\rrn)}\leq2^{-\tau i}
\frac{|B(x_0,r_0)|^{1/r}}{\|\1bf_{B(x_0,r_0)}
\|_{\NHerzSo}};$$
  \item[{\rm(ii)}] when $r_0\in(0,1)$, then, for any
  $\alpha\in\zp^n$ with $|\alpha|\leq d$,
  $$\int_{\rrn}m(x)x^\alpha\,dx=0.$$
\end{enumerate}
\end{definition}

Then we establish the following
molecular characterization
\index{molecular \\characterization}of
the inhomogeneous
local generalized Herz--Hardy space $\NHaSaHol$.

\begin{theorem}\label{llmoleln}
Let $p,\ q\in(0,\infty)$,
$\omega\in M(\rp)$ with
$\mi(\omega)\in(-\frac{n}{p},\infty)$,
$$
s\in\left(0,\min\left\{1,p,q,\frac{n}
{\MI(\omega)+n/p}\right\}\right),
$$
$d\geq\lfloor n(1/s-1)\rfloor$ be a fixed integer,
$$r\in\left(\max\left\{1,p,\frac{n}{
\mi(\omega)+n/p}\right\},\infty\right],$$
and $\tau\in(n(1/s-1/r),\infty)$.
Then $f\in\NHaSaHol$ if and only
if $f\in\mathcal{S}'(\rrn)$
and there exist a sequence $\{m_{j}\}_{j\in\mathbb{N}}$
of local-$(\NHerzSo,\,r,\,d,
\,\tau)$-molecules centered, respectively,
at the balls
$\{B_{j}\}_{j\in\mathbb{N}}
\subset\mathbb{B}$ and
a sequence $\{\lambda_{j}\}_{j\in\mathbb{N}}\subset
[0,\infty)$ such that
$
f=\sum_{j\in\mathbb{N}}\lambda_{j}a_{j}$
in $\mathcal{S}'(\rrn)$, and
$$
\left\|\left\{\sum_{j\in\mathbb{N}}
\left[\frac{\lambda_{j}}
{\|\1bf_{B_{j}}\|_{\NHerzSo}}\right]^s\1bf_
{B_{j}}\right\}^{\frac{1}{s}}\right
\|_{\NHerzSo}<\infty.
$$
Moreover, there exists a
constant $C\in[1,\infty)$ such that, for
any $f\in\NHaSaHol$,
\begin{align*}
C^{-1}\|f\|_{\NHaSaHol}&\leq\inf\left\{
\left\|\left\{\sum_{j\in\mathbb{N}}
\left[\frac{\lambda_{i}}
{\|\1bf_{B_{j}}\|_{\NHerzSo}}\right]^s\1bf_
{B_{j}}\right\}^{\frac{1}{s}}
\right\|_{\NHerzSo}\right\}\\
&\leq C\|f\|_{\NHaSaHol},
\end{align*}
where the infimum is taken
over all the decompositions of $f$ as above.
\end{theorem}

\begin{proof}
Let all the symbols be as
in the present theorem and $f\in\mathcal{S}'(\rrn)$.
Then, from Lemma \ref{llmolell1}
and an argument similar to
that used in the proof of Theorem \ref{llmolel}
with Theorem \ref{Th3} and Lemmas \ref{mbhal}
and \ref{Atoll3} therein replaced,
respectively, by Theorem \ref{Th3n}
and Lemmas \ref{mbhaln} and \ref{Atolnl3},
we infer that $f\in\NHaSaHol$ if and only
if there exist a sequence
$\{m_{j}\}_{j\in\mathbb{N}}$
of local-$(\NHerzSo,\,r,\,d,
\,\tau)$-molecules centered, respectively,
at the balls $\{B_{j}\}_{j\in\mathbb{N}}
\subset\mathbb{B}$ and
a sequence $\{\lambda_{j}\}_{j\in\mathbb{N}}\subset
[0,\infty)$ such that
$
f=\sum_{j\in\mathbb{N}}\lambda_{j}a_{j}$
in $\mathcal{S}'(\rrn)$, and
$$
\left\|\left\{\sum_{j\in\mathbb{N}}
\left[\frac{\lambda_{j}}
{\|\1bf_{B_{j}}\|_{\NHerzSo}}\right]^s\1bf_
{B_{j}}\right\}^{\frac{1}{s}}\right
\|_{\NHerzSo}<\infty,
$$
and, moreover,
\begin{align*}
\|f\|_{\NHaSaHol}\sim\inf\left\{
\left\|\left\{\sum_{j\in\mathbb{N}}
\left[\frac{\lambda_{i}}
{\|\1bf_{B_{j}}\|_{\NHerzSo}}\right]^s\1bf_
{B_{j}}\right\}^{\frac{1}{s}}
\right\|_{\NHerzSo}\right\}
\end{align*}
with the positive equivalence
constants independent of $f$,
where the infimum is taken
over all the decompositions of $f$ as in
the present theorem. This then
finishes the proof of Theorem \ref{llmoleln}.
\end{proof}

We next characterize the
inhomogeneous local generalized
Herz--Hardy space $\NHaSaHl$ via molecules.
To this end, we
first introduce the concept of
local-$(\NHerzS,\,r,\,d,\,\tau)$-molecules
as follows.

\begin{definition}
Let $p$, $q\in(0,\infty)$,
$\omega\in M(\rp)$ with
$$
-\frac{n}{p}<\mi(\omega)\leq
\MI(\omega)<0,
$$
$r\in[1,\infty]$,
$d\in\mathbb{Z}_{+}$, and $\tau\in
(0,\infty)$. Then a measurable
function $m$ on $\rrn$ is called a
\emph{local-$(\NHerzS,\,r,\,d,\,\tau)$-molecule}
\index{local-$(\NHerzS,\,r,\,d,
\,\tau)$-molecule}centered at a ball
$B(x_0,r_0)\in\mathbb{B}$, with $x_0\in
\rrn$ and $r_0\in(0,\infty)$, if
\begin{enumerate}
  \item[{\rm(i)}] for any $i\in\mathbb{Z}_{+}$,
$$\left\|m\1bf_{S_{i}(B(x_0,r_0))}
\right\|_{L^{r}(\rrn)}\leq2^{-\tau i}
\frac{|B(x_0,r_0)|^{1/r}}{\|\1bf_{B(x_0,r_0)}
\|_{\NHerzS}};$$
  \item[{\rm(ii)}] when $r_0\in(0,1)$, then, for any
  $\alpha\in\zp^n$ with $|\alpha|\leq d$,
  $$\int_{\rrn}m(x)x^\alpha\,dx=0.$$
\end{enumerate}
\end{definition}

Now, we prove the following molecular
characterization of $\NHaSaHl$.

\begin{theorem}\label{llmolegn}
Let $p,\ q\in(0,\infty)$,
$\omega\in M(\rp)$ with
$$-\frac{n}{p}
<\mi(\omega)\leq\MI(\omega)<0,$$
$
s\in(0,\min\{1,p,q\}),
$
$d\geq\lfloor n(1/s-1)\rfloor$ be a fixed integer,
$$r\in\left(\max\left\{1,\frac{n}{
\mi(\omega)+n/p}\right\},\infty\right],$$
and $\tau\in(0,\infty)$ with $\tau>
n(1/s-1/r)$.
Then $f\in\NHaSaHl$ if and only if $f\in\mathcal{S}'(\rrn)$
and there exist a sequence $\{m_{j}\}_{j\in\mathbb{N}}$
of local-$(\NHerzS,\,r,\,d,
\,\tau)$-molecules centered, respectively,
at the balls
$\{B_{j}\}_{j\in\mathbb{N}}\subset\mathbb{B}$ and
a sequence $\{\lambda_{j}\}_{j\in\mathbb{N}}\subset
[0,\infty)$ such that
$
f=\sum_{j\in\mathbb{N}}\lambda_{j}a_{j}
$
in $\mathcal{S}'(\rrn)$, and
$$
\left\|\left\{\sum_{j\in\mathbb{N}
}\left[\frac{\lambda_{j}}
{\|\1bf_{B_{j}}\|_{\NHerzS}}\right]^s\1bf_
{B_{j}}\right\}^{\frac{1}{s}}\right
\|_{\NHerzS}<\infty.
$$
Moreover, there exists a
constant $C\in[1,\infty)$ such that, for
any $f\in\NHaSaHl$,
\begin{align*}
C^{-1}\|f\|_{\NHaSaHl}&\leq\inf\left\{
\left\|\left\{\sum_{j\in\mathbb{N}}\left[\frac{\lambda_{i}}
{\|\1bf_{B_{j}}\|_{\NHerzS}}\right]^s\1bf_
{B_{j}}\right\}^{\frac{1}{s}}\right\|_{\HerzS}\right\}\\
&\leq C\|f\|_{\NHaSaHl},
\end{align*}
where the infimum is taken
over all the decompositions of $f$ as above.
\end{theorem}

\begin{proof}
Let all the symbols be as in
the present theorem and $f\in\mathcal{S}'(\rrn)$.
Then, using Theorems \ref{molexl}
and \ref{balln} and repeating
the proof of Theorem \ref{Atog}
with Theorem \ref{Th2}
and Lemmas \ref{mbhag}, \ref{Atogl5},
and \ref{Atogl2} therein replaced, respectively,
by Theorem \ref{Th2n}
and Lemmas \ref{mbhagn}, \ref{Atognl5},
and \ref{Atognl2}, we find that
$f\in\NHaSaHl$ if and only
if there exist a sequence
$\{m_{j}\}_{j\in\mathbb{N}}$
of local-$(\NHerzS,\,r,\,d,\,\tau)$-molecules
centered, respectively,
at the balls $\{B_{j}\}_{j\in\mathbb{N}}
\subset\mathbb{B}$ and
a sequence $\{\lambda_{j}\}_{j\in\mathbb{N}}\subset
[0,\infty)$ such that
$
f=\sum_{j\in\mathbb{N}}\lambda_{j}a_{j}$
in $\mathcal{S}'(\rrn)$, and
$$
\left\|\left\{\sum_{j\in\mathbb{N}}
\left[\frac{\lambda_{j}}
{\|\1bf_{B_{j}}\|_{\NHerzS}}\right]^s\1bf_
{B_{j}}\right\}^{\frac{1}{s}}\right
\|_{\NHerzS}<\infty,
$$
and, moreover,
\begin{align*}
\|f\|_{\NHaSaHl}\sim\inf\left\{
\left\|\left\{\sum_{j\in\mathbb{N}}
\left[\frac{\lambda_{i}}
{\|\1bf_{B_{j}}\|_{\NHerzS}}\right]^s\1bf_
{B_{j}}\right\}^{\frac{1}{s}}
\right\|_{\NHerzS}\right\}
\end{align*}
with the positive equivalence
constants independent of $f$,
where the infimum is taken
over all the decompositions of $f$ as in
the present theorem. This then finishes
the proof of Theorem \ref{llmolegn}.
\end{proof}

\subsection{Littlewood--Paley Function Characterizations}

The main target of this subsection is to
characterize both the
inhomogeneous local generalized
Herz--Hardy spaces $\NHaSaHol$
and $\NHaSaHl$ via the local
Lusin area function, the local
$g$-function, and the local
$g_\lambda^*$-function defined
as in Definition \ref{df511}.
To begin with, we establish the following
Littlewood--Paley characterizations
of $\NHaSaHol$.

\begin{theorem}\label{lusln}
Let $p,\ q\in(0,\infty)$, $\omega\in
M(\rp)$ satisfy
$\mi(\omega)\in(-\frac{n}{p},\infty)$,
$$
s_{0}:=\min\left\{1,p,q,\frac{n}{
\MI(\omega)+n/p}\right\},
$$
and $\lambda\in(\max\{1,2/s_{0}\},\infty)$.
Then the following four statements are
mutually equivalent:
\begin{enumerate}
  \item[{\rm(i)}] $f\in\NHaSaHol$;
  \item[{\rm(ii)}] $f\in\mathcal{S}'(\rrn)$
  and $S_{\mathrm{loc}}(f)\in\NHerzSo$;
  \item[{\rm(iii)}] $f\in\mathcal{S}'(\rrn)$
  and $g_{\mathrm{loc}}(f)\in\NHerzSo$;
  \item[{\rm(iv)}] $f\in\mathcal{S}'(\rrn)$
  and $(g_{\lambda}^*
  )_{\mathrm{loc}}(f)\in\NHerzSo$.
\end{enumerate}
Moreover, for any $f\in\NHaSaHol$,
\begin{align*}
\left\|f\right\|_{\NHaSaHol}&\sim\left\|
S_{\mathrm{loc}}(f)\right\|_{\NHerzSo}
\sim\left\|g_{\mathrm{loc}}(f)\right\|_{\NHerzSo}\\&\sim
\left\|\left(g_{\lambda}^*\right)_{\mathrm{loc}}
(f)\right\|_{\NHerzSo},
\end{align*}
where the positive equivalence
constants are independent of $f$.
\end{theorem}

\begin{proof}
Let all the symbols be as in
the present theorem. Then, from
Theorems \ref{lplx}, \ref{lplgx}, and
\ref{lplglx} and an argument similar
to that used in the proof of
Theorem \ref{lusin} via replacing
Theorem \ref{Th3} and Lemmas \ref{mbhal}
and \ref{Atoll3} therein, respectively,
by Theorem \ref{Th3n} and Lemmas \ref{mbhaln}
and \ref{Atolnl3}, we deduce that
(i) through (iv) of the present theorem
are mutually equivalent and, for any
$f\in\NHaSaHol$,
\begin{align*}
\left\|f\right\|_{\NHaSaHol}&\sim\left\|
S_{\mathrm{loc}}(f)\right\|_{\NHerzSo}
\sim\left\|g_{\mathrm{loc}}(f)\right\|_{\NHerzSo}\\&\sim
\left\|\left(g_{\lambda}^*\right)_{\mathrm{loc}}
(f)\right\|_{\NHerzSo}
\end{align*}
with the positive equivalence constants
independent of $f$. This then finishes
the proof of Theorem \ref{lusln}.
\end{proof}

Next, we turn to characterize
the inhomogeneous local generalized
Herz--Hardy space $\NHaSaHl$ via various
Littlewood--Paley functions as follows.

\begin{theorem}\label{lusgn}
Let $p,\ q\in(0,\infty)$, $\omega\in
M(\rp)$ satisfy
$$
-\frac{n}{p}<\mi(\omega)\leq\MI(\omega)<0,
$$
and $\lambda\in(2\max\{1,\frac{1}{p},
\frac{1}{q}\},\infty)$.
Then the following four statements are
mutually equivalent:
\begin{enumerate}
  \item[{\rm(i)}] $f\in\NHaSaHl$;
  \item[{\rm(ii)}] $f\in\mathcal{S}'(\rrn)$
  and $S_{\mathrm{loc}}(f)\in\NHerzS$;
  \item[{\rm(iii)}] $f\in\mathcal{S}'(\rrn)$
  and $g_{\mathrm{loc}}(f)\in\NHerzS$;
  \item[{\rm(iv)}] $f\in\mathcal{S}'(\rrn)$
  and $(g_{\lambda}^*
  )_{\mathrm{loc}}(f)\in\NHerzS$.
\end{enumerate}
Moreover, for any $f\in\NHaSaHl$,
\begin{align*}
\left\|f\right\|_{\NHaSaHl}&\sim\left\|
S_{\mathrm{loc}}(f)\right\|_{\NHerzS}
\sim\left\|g_{\mathrm{loc}}(f)\right\|_{\NHerzS}\\&\sim
\left\|\left(g_{\lambda}^*\right)_{\mathrm{loc}}
(f)\right\|_{\NHerzS},
\end{align*}
where the positive equivalence
constants are independent of $f$.
\end{theorem}

\begin{proof}
Let all the symbols be as in the present theorem.
Then, repeating the proof of Theorem
\ref{lusing} with Theorems \ref{Th5.4},
\ref{Th5.6}, and \ref{lusin}, and
Lemma \ref{lusingl1}(iii) therein
replaced, respectively, by Theorems \ref{Th5.4n},
\ref{Th5.6n}, and \ref{lusinn}, and
Lemma \ref{lusgl1}, we
find that (i) through (iv) of the
present theorem are
mutually equivalent and, for any
$f\in\NHaSaHl$,
\begin{align*}
\left\|f\right\|_{\NHaSaHl}&\sim\left\|
S_{\mathrm{loc}}(f)\right\|_{\NHerzS}
\sim\left\|g_{\mathrm{loc}}(f)\right\|_{\NHerzS}\\&\sim
\left\|\left(g_{\lambda}^*\right)_{\mathrm{loc}}
(f)\right\|_{\NHerzS}
\end{align*}
with the positive equivalence
constants independent of $f$.
This then finishes the
proof of Theorem \ref{lusgn}.
\end{proof}

\subsection{Boundedness of Pseudo-Differential Operators}

Recall that the
H\"{o}rmander class $S^0_{1,0}(\rrn)$ is defined
to be the set of all the infinitely differentiable
functions $\sigma$ on $\rrn\times\rrn$
satisfying \eqref{hormander}. Let
$\sigma\in S^0_{1,0}(\rrn)$ and $T_\sigma$ be
a pseudo-differential operator
defined as in Definition \ref{df561}.
In this subsection, we investigate the boundedness of
$T_\sigma$ on inhomogeneous
localized generalized Herz--Hardy spaces.
First, the following
conclusion gives the boundedness of
operator $T_\sigma$ on
the inhomogeneous local generalized
Herz--Hardy space $\NHaSaHol$.

\begin{theorem}\label{posedon}
Let $p,\ q\in(0,\infty)$, $\omega\in M(\rp)$ satisfy
$\mi(\omega)\in(-\frac{n}{p},\infty)$,
and $T_\sigma$ be a pseudo-differential
operator with the symbol $\sigma\in S_{1,0}^{0}(\rrn)$.
Then $T_{\sigma}$ is well
defined on $\NHaSaHol$ and
there exists a positive constant
$C$ such that, for any $f\in
\NHaSaHol$,
$$
\left\|T_\sigma(f)\right\|_{\NHaSaHol}
\leq C\|f\|_{\NHaSaHol}.
$$
\end{theorem}

\begin{proof}
Let all the symbols be as
in the present theorem. Then, applying
Lemma \ref{posedol1} and repeating
an argument similar to that
used in the proof of Theorem
\ref{posedo} with $\Mw$, Theorems
\ref{Th3} and \ref{abso}, and Lemmas
\ref{mbhal} and \ref{Atoll3} therein
replaced, respectively, by $\MI(\omega)$, Theorems
\ref{Th3n} and \ref{abson}, and Lemmas
\ref{mbhaln} and \ref{Atolnl3}, we conclude
that $T_{\sigma}$ is well defined
on $\NHaSaHol$ and, for any $f\in
\NHaSaHol$,
$$
\left\|T_\sigma(f)\right\|_{\NHaSaHol}
\lesssim\|f\|_{\NHaSaHol}.
$$
This then finishes the proof of Theorem
\ref{posedon}.
\end{proof}

\begin{remark}
We should point out that,
in Theorem \ref{posedon},
when $\omega(t):=t^\alpha$
for any $t\in(0,\infty)$
and for any given $\alpha\in\rr$,
the result obtained in Theorem \ref{posedon}
goes back to \cite[Corollary 4.14(c)]{WYY}.
\end{remark}

In addition, we establish the following
boundedness of
pseudo-differential operators on the
local Hardy space $\NHaSaHl$.

\begin{theorem}\label{pseudogn}
Let $p,\ q\in(0,\infty)$, $\omega\in M(\rp)$ satisfy
$$-\frac{n}{p}
<\mi(\omega)\leq\MI(\omega)<0,$$
and $T_\sigma$ be a pseudo-differential
operator with symbol $\sigma\in S_{1,0}^{0}(\rrn)$.
Then $T_{\sigma}$ is well
defined on $\NHaSaHl$ and
there exists a positive constant
$C$ such that, for any $f\in
\NHaSaHl$,
$$
\left\|T_\sigma(f)\right\|_
{\NHaSaHl}\leq C\|f\|_{\NHaSaHl}.
$$
\end{theorem}

\begin{proof}
Let all the symbols be as
in the present theorem. Then,
using Theorems \ref{bddox} and \ref{balln}
and repeating the proof of Theorem
\ref{bddog} with $$\mw,\ \Mw,$$
Theorem \ref{Th2},
Lemmas \ref{ee}, \ref{equa}, \ref{Atogl5},
and \ref{Atogl2}, and Corollary \ref{maxbl}
therein replaced, respectively, by
$\mi(\omega)$, $\MI(\omega)$, Theorem \ref{Th2n},
Lemmas \ref{een}, \ref{equan}, \ref{Atognl5},
and \ref{Atognl2}, and Corollary \ref{maxbln},
we conclude that $T_{\sigma}$ is well defined on
$\NHaSaHl$ and, for any $f\in\NHaSaHl$,
$$
\left\|T(f)\right\|_{\NHaSaHl}
\lesssim\left\|f\right\|_{\NHaSaHl},
$$
which then completes the proof of Theorem
\ref{pseudogn}.
\end{proof}

\section{Inhomogeneous Weak Generalized \\
Herz--Hardy Spaces}

In this section, we introduce both inhomogeneous weak
generalized Herz spaces and
inhomogeneous weak generalized Herz--Hardy spaces.
Then we establish the maximal function, the atomic,
the molecular, and the Littlewood--Paley function
characterizations of these Hardy spaces. As applications,
we also obtain the boundedness of
Calder\'{o}n--Zygmund operators
from inhomogeneous generalized
Herz--Hardy spaces to inhomogeneous
weak generalized Herz--Hardy spaces
even in the critical case as well as
show that the real interpolation spaces
between inhomogeneous generalized
Herz--Hardy spaces and the Lebesgue
space $L^\infty(\rrn)$ are just the new
introduced inhomogeneous weak generalized
Herz--Hardy spaces.

To begin with, we introduce inhomogeneous weak
generalized Herz spaces as follows.

\begin{definition}\label{df6.0.7n}
Let $p,\ q\in(0,\infty)$ and $\omega\in M(\rp)$.
\begin{enumerate}
  \item[{\rm(i)}] The \emph{inhomogeneous weak local generalized
  Herz space} $W\NHerzSo$\index{$W\NHerzSo$} is defined to be
  the set of all the measurable
  functions $f$ on $\rrn$ such that
  $$
  \|f\|_{W\NHerzSo}:=\sup_{\alpha\in(0,\infty)}\left\{
  \alpha\left\|\1bf_{\{x\in\rrn:\ |f(x)|>\alpha\}}\right\|_{\NHerzSo}
  \right\}<\infty.
  $$
  \item[{\rm(ii)}] The \emph{inhomogeneous weak global generalized
  Herz space}\index{inhomogeneous weak global generalized
  Herz space} $W\NHerzS$\index{$W\NHerzS$} is defined to
  be the set of all the measurable functions
  $f$ on $\rrn$ such that
  $$
  \|f\|_{W\NHerzS}:=\sup_{\alpha\in(0,\infty)}\left\{
  \alpha\left\|\1bf_{\{x\in\rrn:\ |f(x)|>\alpha\}}\right\|_{\NHerzS}
  \right\}<\infty.
  $$
\end{enumerate}
\end{definition}

\begin{remark}
In Definition \ref{df6.0.7n}, when $\omega(t):=t^\alpha$
for any $t\in(0,\infty)$ and for any given $\alpha\in\rr$,
the inhomogeneous weak local generalized Herz space $W\NHerzSo$
coincides with the classical \emph{inhomogeneous weak
Herz space\index{inhomogeneous weak Herz space}
$WK_p^{\alpha,q}(\rrn)$\index{$WK_p^{\alpha,q}(\rrn)$}}
which was originally introduced in \cite{hly97}
(see also \cite[Section 1.4]{LYH}).
\end{remark}

Let $N\in\mathbb{N}$ and the non-tangential grand maximal
function $\mathcal{M}_N(f)$ of $f\in\mathcal{S}'(\rrn)$ be
as in \eqref{sec6e1}. We now introduce
the following inhomogeneous weak generalized
Herz--Hardy spaces
via the operator $\mc_N$.\index{inhomogeneous weak generalized Herz--Hardy space}

\begin{definition}\label{df6.11n}
Let $p$, $q\in(0,\infty)$, $\omega\in M(\rp)$,
and $N\in\mathbb{N}$.
\begin{enumerate}
  \item[{\rm(i)}] The \emph{inhomogeneous
  weak generalized Herz--Hardy space
  $\WNHaSaHo$\index{$\WNHaSaHo$}}, associated with
  the inhomogeneous weak local generalized
  Herz space $W\NHerzSo$,
  is defined to be the set of all
  the $f\in\mathcal{S}'(\rrn)$ such that
  $$
  \|f\|_{\WNHaSaHo}:=\left\|
  \mathcal{M}_{N}(f)\right\|_{W\NHerzSo}<\infty.
  $$
  \item[{\rm(ii)}] The \emph{inhomogeneous
  weak generalized Herz--Hardy space
  $\WNHaSaH$\index{$\WNHaSaH$}}, associated with
  the inhomogeneous weak global
  generalized Herz space $W\NHerzS$,
  is defined to be the set of all the
  $f\in\mathcal{S}'(\rrn)$ such that
  $$
  \|f\|_{\WNHaSaH}:=\left\|
  \mathcal{M}_{N}(f)\right\|_{W\NHerzS}<\infty.
  $$
\end{enumerate}
\end{definition}

\begin{remark}
We point out that, in Definition \ref{df6.11n},
if $\omega(t):=t^\alpha$ for any
$t\in(0,\infty)$ and for any given $\alpha\in\rr$,
then the inhomogeneous weak
generalized Herz--Hardy space $\WNHaSaHo$
coincides with the classical
\emph{inhomogeneous weak Herz-type Hardy space}
\index{inhomogeneous weak Herz-type Hardy \\ space}$WHK_p^{\alpha,q}
(\rrn)$\index{$WHK_p^{\alpha,q}(\rrn)$}
(see, for instance, \cite[Definition 2.1.2]{LYH}).
\end{remark}

\subsection{Maximal Function Characterizations}

The main target of this subsection
is to establish the equivalent characterizations
of inhomogeneous weak generalized Herz--Hardy spaces
via various radial and non-tangential maximal functions
given in Definition \ref{smax}.
To begin with, we prove the following
maximal function characterizations of
the inhomogeneous weak generalized Herz--Hardy space
$W\NHaSaHo$.

\begin{theorem}\label{macln}
Let $p$, $q$, $a,\ b\in(0,\infty)$, $\omega\in M(\rp)$,
$N\in\mathbb{N}$, and $\phi\in\mathcal{S}(\rrn)$ satisfy
$\int_{\rrn}\phi(x)\,dx\neq0.$
\begin{enumerate}
\item[\rm{(i)}] Let $N\in\mathbb{N}
\cap[\lfloor b+1\rfloor,\infty)$.
Then, for any $f\in\mathcal{S}'(\rrn)$,
$$
\|M(f,\phi)\|_{W\NHerzSo}\lesssim
\|M^*_{a}(f,\phi)\|_{W\NHerzSo}
\lesssim\|M^{**}_{b}(f,\phi)\|_{W\NHerzSo},
$$
\begin{align*}
\|M(f,\phi)\|_{W\NHerzSo}&\lesssim
\|\mathcal{M}_{N}(f)\|_{W\NHerzSo}
\lesssim\|\mathcal{M}_{\lfloor b+1\rfloor}(f)
\|_{W\NHerzSo}\\&\lesssim
\|M^{**}_{b}(f,\phi)\|_{W\NHerzSo},
\end{align*}
and
$$
\|M^{**}_{b}(f,\phi)\|_{W\NHerzSo}\sim
\|\mathcal{M}^{**}_{b,N}(f)\|_{W\NHerzSo},
$$
where the implicit positive
constants are independent of $f$.
\item[\rm{(ii)}] Let
$\omega\in M(\rp)$ satisfy
$\mi(\omega)\in(-\frac{n}{p},\infty)$ and
$$b\in\left(\max\left\{\frac{n}{p},\MI
(\omega)+\frac{n}{p}\right\},\infty\right).$$
Then, for any $f\in\mathcal{S}'(\rrn),$
$$\|M^{**}_{b}(f,\phi)\|_{W\NHerzSo}
\lesssim\|M(f,\phi)\|_{W\NHerzSo},$$
where the implicit positive
constant is independent of $f$.
In particular,
when $N\in\mathbb{N}\cap[\lfloor b+1\rfloor,\infty)$,
if one of the quantities
$$\|M(f,\phi)\|_{W\NHerzSo},\
\|M^*_{a}(f,\phi)\|_{W\NHerzSo},\ \|\mathcal{M}_{N}
(f)\|_{W\NHerzSo},$$
$$
\|M^{**}_{b}(f,\phi)\|_{W\NHerzSo},
\ \text{and}\ \|\mathcal{M}^{**}_{b,N}(f)\|_{W\NHerzSo}$$
is finite, then the others are also finite and mutually
equivalent with the
positive equivalence
constants independent of $f$.
\end{enumerate}
\end{theorem}

To prove this theorem,
we first show the following auxiliary
lemma about the boundedness of
the Hardy--Littlewood maximal operator on inhomogeneous
weak local generalized Herz spaces.

\begin{lemma}\label{wmbhln}
Let $p,\ q\in(0,\infty)$ and
$\omega\in M(\rp)$ satisfy
$\mi(\omega)\in(-\frac{n}{p},\infty)$.
Then, for any given $$r\in\left(0,\min
\left\{p,\frac{n}{\MI(\omega)+n/p}\right\}\right),$$
there exists a positive constant $C$ such that,
for any $f\in L^1_{\mathrm{loc}}(\rrn)$,
$$
\left\|\mc(f)\right\|_{[W\NHerzSo]^{1/r}}
\leq C\|f\|_{[W\NHerzSo]^{1/r}}.
$$
\end{lemma}

\begin{proof}
Let all the symbols be as in
the present lemma and $$r\in\left(0,\min
\left\{p,\frac{n}{\MI(\omega)+n/p}\right\}\right).$$
Then, repeating an argument similar
to that used in the proof of
Lemma \ref{wmbhl} with both Theorem \ref{Th3}
and Lemma \ref{vmbhl} therein replaced,
respectively, by both
Theorem \ref{Th3n} and Lemma \ref{vmbhln},
we find that, for any
$f\in L^1_{\mathrm{loc}}(\rrn)$,
$$
\left\|\mc(f)\right\|_{[W\NHerzSo]^{1/r}}
\lesssim\|f\|_{[W\NHerzSo]^{1/r}}.
$$
This then finishes the proof of
Lemma \ref{wmbhln}.
\end{proof}

With the help of this lemma, we next
show Theorem \ref{macln}.

\begin{proof}[Proof of Theorem \ref{macln}]
Let all the symbols be as in the present
theorem. Then, from Lemma \ref{Th5.4l1} and
an argument similar to
that used in the proof of Theorem \ref{macl} with
$\Mw$, Theorem \ref{Th3}, and Lemma \ref{wmbhl}
therein replaced, respectively, by
$\MI(\omega)$, Theorem \ref{Th3n},
and Lemma \ref{wmbhln}, we infer that
both (i) and (ii) of the
present theorem hold true,
and hence complete the proof of
Theorem \ref{macln}.
\end{proof}

Now, we consider the maximal function characterizations
of the inhomogeneous weak generalized Herz--Hardy space
$W\NHaSaH$. Indeed, we have the following
conclusion.

\begin{theorem}\label{Th7.11n}
Let $p$, $q$, $a,\ b\in(0,\infty)$,
$\omega\in M(\rp)$,
$N\in\mathbb{N}$, and
$\phi\in\mathcal{S}(\rrn)$ satisfy
$\int_{\rrn}\phi(x)\,dx\neq0.$
\begin{enumerate}
\item[\rm{(i)}] Let $N\in\mathbb{N}
\cap[\lfloor b+1\rfloor,\infty)$
and $\omega$ satisfy
$\MI(\omega)\in(-\infty,0)$.
Then, for any $f\in\mathcal{S}'(\rrn)$,
$$
\|M(f,\phi)\|_{W\NHerzS}\lesssim
\|M^*_{a}(f,\phi)\|_{W\NHerzS}
\lesssim\|M^{**}_{b}(f,\phi)\|_{W\NHerzS},
$$
\begin{align*}
\|M(f,\phi)\|_{W\NHerzS}&\lesssim
\|\mathcal{M}_{N}(f)\|_{W\NHerzS}
\lesssim\|\mathcal{M}_{\lfloor
b+1\rfloor}(f)\|_{W\NHerzS}\\&\lesssim
\|M^{**}_{b}(f,\phi)\|_{W\NHerzS},
\end{align*}
and
$$
\|M^{**}_{b}(f,\phi)\|_{W\NHerzS}\sim
\|\mathcal{M}^{**}_{b,N}(f)\|_{W\NHerzS},
$$
where the implicit positive
constants are independent of $f$.
\item[\rm{(ii)}] Let $\omega\in M(\rp)$ satisfy
$-\frac{n}{p}<\mi(\omega)\leq\MI(\omega)<0$
and $b\in(\frac{n}{p},\infty)$.
Then, for any $f\in\mathcal{S}'(\rrn),$
$$\|M^{**}_{b}(f,\phi)\|_{W\NHerzS}
\lesssim\|M(f,\phi)\|_{W\NHerzS},$$
where the implicit positive
constant is independent of $f$.
In particular, when $N\in\mathbb{N}
\cap[\lfloor b+1\rfloor,\infty)$,
if one of the quantities
$$\|M(f,\phi)\|_{W\NHerzS},\
\|M^*_{a}(f,\phi)\|_{W\NHerzS},\ \|\mathcal{M}_{N}
(f)\|_{W\NHerzS},$$
$$
\|M^{**}_{b}(f,\phi)\|_{W\NHerzS},
\ \text{and}\ \|\mathcal{M}^{**}_{b,N}(f)\|_{W\NHerzS}$$
is finite, then the others are also finite and mutually
equivalent with the
positive equivalence constants independent of $f$.
\end{enumerate}
\end{theorem}

To prove this theorem, we first establish the
boundedness of the Hardy--Littlewood maximal operator
on the convexification of $W\NHerzS$ as follows

\begin{lemma}\label{wmbhgn}
Let $p,\ q\in(0,\infty)$ and
$\omega\in M(\rp)$ satisfy
$\mi(\omega)\in(-\frac{n}{p},\infty)$.
Then, for any given $$r\in\left(0,\min
\left\{p,\frac{n}{\MI(\omega)+n/p}\right\}\right),$$
there exists a positive constant $C$ such that,
for any $f\in L^1_{\mathrm{loc}}(\rrn)$,
$$
\left\|\mc(f)\right\|_{[W\NHerzS]^{1/r}}
\leq C\|f\|_{[W\NHerzS]^{1/r}}.
$$
\end{lemma}

\begin{proof}
Let all the symbols be as in
the present lemma.
Then, using Lemma \ref{vwmbhll1} and
repeating the proof of
Lemma \ref{wmbhg} with both Theorem \ref{Th2}
and Lemma \ref{vmbhg} therein replaced,
respectively, by both
Theorem \ref{Th2n} and Lemma \ref{vmbhgn},
we conclude that, for any given
$$r\in\left(0,\min
\left\{p,\frac{n}{\MI(\omega)+n/p}
\right\}\right)$$
and for any
$f\in L^1_{\mathrm{loc}}(\rrn)$,
$$
\left\|\mc(f)\right\|_{[W\NHerzS]^{1/r}}
\lesssim\|f\|_{[W\NHerzS]^{1/r}}.
$$
This then finishes the proof of
Lemma \ref{wmbhgn}.
\end{proof}

We next turn to show the maximal
function characterizations Theorem \ref{Th7.11n}

\begin{proof}[Proof of Theorem \ref{Th7.11n}]
Let all the symbols be as in the present
theorem. Then, applying Lemma \ref{Th5.4l1} and
repeating an argument similar to
that used in the proof of Theorem \ref{Th7.11} via
replacing $\Mw$, Theorem \ref{Th2}, and
Lemma \ref{wmbhg} therein, respectively, by
$\MI(\omega)$, Theorem \ref{Th2n},
and Lemma \ref{wmbhgn}, we find that
both (i) and (ii) of the present theorem hold true.
This finishes the proof of
Theorem \ref{Th7.11n}.
\end{proof}

\subsection{Relations with
Inhomogeneous Weak Generalized \\ Herz Spaces}

In this subsection, we investigate the relations
between inhomogeneous weak generalized Herz spaces
and associated Hardy spaces.
We first study the relation between
$W\NHerzSo$ and $\WNHaSaHo$.
Indeed, the following conclusion
shows that, under some reasonable
and sharp assumptions on the exponents, then
$$\WNHaSaHo=W\NHerzSo$$
in the sense of equivalent quasi-norms.

\begin{theorem}\label{wrelaln}
Let $p\in(1,\infty)$,
$q\in(0,\infty)$, and
$\omega\in M(\rp)$ satisfy
$$-\frac{n}{p}<\mi(\omega)
\leq\MI(\omega)<\frac{n}{p'},$$
where $\frac{1}{p}+\frac{1}{p'}=1$. Then
\begin{enumerate}
\item[\rm{(i)}] $W\NHerzSo
\hookrightarrow\mathcal{S}'(\rrn)$.
\item[\rm{(ii)}] If
$f\in W\NHerzSo$, then $f\in W\NHaSaHo$
and there exists a positive
constant $C$, independent
of $f$, such that $$\left\|f\right\|_{
W\NHaSaHo}\leq C\left\|f\right\|_{W\NHerzSo}.$$
\item[\rm{(iii)}] If $f
\in W\NHaSaHo$, then there exists a
locally integrable function
$g$ belonging to $W\NHerzSo$
such that $g$ represents $f$,
which means that $f=g$ in
$\mathcal{S}'(\rrn)$, $$
\left\|f\right\|_{W\NHaSaHo}=
\left\|g\right\|_{W\NHaSaHo},$$
and there exists a positive
constant $C$, independent of $f$,
such that $$\left\|g\right\|_{W\NHerzSo
}\leq C\left\|f\right\|_{W\NHaSaHo}.$$
\end{enumerate}
\end{theorem}

\begin{proof}
Let all the symbols be as in
the present theorem. Then, from
Lemma \ref{Th5.7l1} and an
argument similar to that
used in the proof of Theorem \ref{wrelal}
with $$\mw,\ \Mw,$$ Theorem \ref{Th3},
and Lemma \ref{wmbhl} therein replaced,
respectively, by $\mi(\omega)$,
$\MI(\omega)$, Theorem \ref{Th3n},
and Lemma \ref{wmbhln}, we deduce that
(i) through (iii) of the present theorem
hold true. This then finishes the proof of
Theorem \ref{wrelaln}.
\end{proof}

Similarly, we can show
that $\WNHaSaH=W\NHerzS$
with equivalent quasi-norms as follows.

\begin{theorem}\label{wrelagn}
Let $p\in(1,\infty)$, $q\in(0,\infty)$, and
$\omega\in M(\rp)$ with
$$-\frac{n}{p}<\mi(\omega)\leq\MI(\omega)<0.$$
Then
\begin{enumerate}
\item[\rm{(i)}] $W\NHerzS\hookrightarrow\mathcal{S}'(\rrn)$.
\item[\rm{(ii)}] If $f\in W\NHerzS$, then $f\in W\NHaSaH$
and there exists a positive constant $C$, independent
of $f$, such that $$\left\|f\right\|_{W\NHaSaH}\leq
C\left\|f\right\|_{W\NHerzS}.$$
\item[\rm{(iii)}] If $f\in W\NHaSaH$, then there exists a
locally integrable function $g$ belonging to $W\NHerzS$
such that $g$ represents $f$, which means that $f=g$ in
$\mathcal{S}'(\rrn)$, $$\left\|f\right\|_{W\NHaSaH}=
\left\|g\right\|_{W\NHaSaH},$$
and there exists a positive constant $C$, independent of $f$,
such that $$\left\|g\right\|_{W\NHerzS}\leq
C\left\|f\right\|_{W\NHaSaH}.$$
\end{enumerate}
\end{theorem}

\begin{proof}
Let all the symbols be as in
the present theorem. Then, applying
Lemma \ref{Th5.7l1} and repeating
the proof of Theorem \ref{wrelag} via
replacing $\Mw$, Theorem \ref{Th2},
and Lemma \ref{wmbhg} therein,
respectively, by $\MI(\omega)$,
Theorem \ref{Th2n}, and Lemma \ref{wmbhln},
we conclude that
(i), (ii), and (iii) of the present theorem
hold true. This further finishes the proof of
Theorem \ref{wrelagn}.
\end{proof}

\subsection{Atomic Characterizations}

The main target of this
subsection is to characterize both the
inhomogeneous weak generalized
Herz--Hardy spaces $\WNHaSaHo$ and $\WNHaSaH$
by means of atoms.
For this purpose, we first introduce the following
inhomogeneous weak generalized atomic Herz--Hardy
space $WH\Kmc_{\omega,\0bf}^{p,q,r,d}(\rrn)$
via $(\NHerzSo,\,r,\,d)$-atoms
introduced in Definition \ref{atomn}.

\begin{definition}\label{dfwatomn}
Let $p,\ q\in(0,\infty)$, $\omega\in M(\rp)$ with
$\mi(\omega)\in(-\frac{n}{p},\infty)$,
$$
p_-\in\left(0,\frac{\min\{p,q,
\frac{n}{\MI(\omega)+n/p}\}}{\max\{1,p,q\}}\right),
$$
$d\geq\lfloor n(1/p_--1)\rfloor$ be a fixed
integer, and
$$
r\in\left(\max\left\{1,p,\frac{n}
{\mi(\omega)+n/p}\right\},\infty\right].
$$
The \emph{inhomogeneous
weak generalized atomic Herz--Hardy space}
\index{inhomogeneous weak generalized atomic Herz--Hardy \\ space}$WH
\Kmc_{\omega,\0bf}^{p,q,r,d}
(\rrn)$\index{$WH\Kmc_{\omega,\0bf}^{p,q,r,d}(\rrn)$},
associated with $W\NHerzSo$, is defined to
be the set of all the
$f\in\mathcal{S}'(\rrn)$ such that
there exist a sequence $\{a_{i,j}\}_{i\in
\mathbb{Z},\,j\in\mathbb{N}}$ of
$(\NHerzSo,\,r,\,d)$-atoms supported,
respectively, in the balls $\{B_{i,j}
\}_{i\in\mathbb{Z},\,j\in\mathbb{N}}
\subset\mathbb{B}$ and three
positive constants $c\in(0,1]$,
$A$, and $\widetilde{A}$,
independent of $f$, satisfying that,
for any $i\in\mathbb{Z}$,
$$
\sum_{j\in\mathbb{N}}\1bf_{cB_{i,j}}\leq A,
$$
$$
f=\sum_{i\in\mathbb{Z}}\sum_{j\in\mathbb{N}}
\widetilde{A}2^i\left\|\1bf_{B_{i,j}}
\right\|_{\NHerzSo}a_{i,j}
$$
in $\mathcal{S}'(\rrn)$, and
$$
\sup_{i\in\mathbb{Z}}\left\{2^i
\left\|\sum_{j\in\mathbb{N}}
\1bf_{B_{i,j}}\right\|_{\NHerzSo}\right\}<\infty.
$$
Moreover, for any $f\in
WH\Kmc_{\omega,\0bf}^{p,q,r,d}(\rrn)$,
$$
\|f\|_{WH\Kmc_{\omega,\0bf}^{p,q,r,d}(\rrn)}:=
\inf\left\{\sup_{i\in\mathbb{Z}}2^i
\left\|\sum_{j\in\mathbb{N}}
\1bf_{B_{i,j}}\right\|_{\NHerzSo}\right\},
$$
where the infimum is taken over
all the decompositions of
$f$ as above.
\end{definition}

Then we have the following atomic characterization
of the inhomogeneous weak generalized Herz--Hardy space
$W\NHaSaHo$.\index{atomic characterization}

\begin{theorem}\label{watomn}
Let $p$, $q$, $\omega$,
$r$, and $d$ be as in Definition \ref{dfwatomn}.
Then
$$
WH\NHerzSo=WH\Kmc_{\omega,\0bf}^{p,q,r,d}(\rrn)
$$
with equivalent quasi-norms.
\end{theorem}

\begin{proof}
Let all the symbols be as in the
present theorem. Then, from
Lemmas \ref{watomd} and \ref{watomr}
and an argument similar
to that used in the proof of Theorem
\ref{watom} with $\Mw$, Theorem \ref{Th3},
and Lemmas \ref{convexll}, \ref{mbhl},
\ref{mbhal}, \ref{vmbhl}, and \ref{wmbhl} therein
replaced, respectively, by $\MI(\omega)$,
Theorem \ref{Th3n}, and Lemmas \ref{convexlln},
\ref{mbhln}, \ref{mbhaln},
\ref{vmbhln}, and \ref{wmbhln}, it follows that
$$
WH\NHerzSo=WH\Kmc_{\omega,\0bf}^{p,q,r,d}(\rrn)
$$
with equivalent quasi-norms. This then
finishes the proof of Theorem \ref{watomn}.
\end{proof}

Now, we turn to show the atomic
characterization of the inhomogeneous
weak generalized Herz--Hardy space $\WNHaSaH$.
For this purpose, we
first introduce the following
weak atomic Hardy space
$WH\Kmc_{\omega}^{p,q,r,d}(\rrn)$ via
$(\NHerzS,\,r,\,d)$-atoms
given in Definition \ref{deatomgn}.

\begin{definition}\label{dfwatomgn}
Let $p,\ q\in(0,\infty)$, $\omega\in M(\rp)$ with
$$-\frac{n}{p}<\mi(\omega)
\leq\MI(\omega)<0,$$
$
p_-\in(0,\min\{p,q\}/
\max\{1,p,q\}),
$
$d\geq\lfloor n(1/p_--1)\rfloor$ be a fixed
integer, and
$$
r\in\left(\max\left\{1,
\frac{n}{\mi(\omega)+n/p}\right\},\infty\right].
$$ The \emph{inhomogeneous
weak generalized atomic Herz--Hardy space}
\index{weak generalized atomic Herz--Hardy \\ space}$WH
\Kmc_{\omega}^{p,q,r,d}(
\rrn)$\index{$WH\Kmc_{\omega}^{p,q,r,d}(\rrn)$},
associated with $W\NHerzS$, is defined to
be the set of all the
$f\in\mathcal{S}'(\rrn)$ such that
there exist a sequence $\{a_{i,j}\}_{i\in
\mathbb{Z},\,j\in\mathbb{N}}$ of
$(\NHerzS,\,r,\,d)$-atoms supported,
respectively, in the balls
$\{B_{i,j}\}_{i\in\mathbb{Z},\,j
\in\mathbb{N}}\subset\mathbb{B}$
and three positive constants $c\in(0,1]$,
$A$, and $\widetilde{A}$,
independent of $f$, satisfying that,
for any $i\in\mathbb{Z}$,
$$
\sum_{j\in\mathbb{N}}\1bf_{cB_{i,j}}\leq A,
$$
$$
f=\sum_{i\in\mathbb{Z}}\sum_{j\in\mathbb{N}}
\widetilde{A}2^i\left\|\1bf_{B_{i,j}}
\right\|_{\NHerzS}a_{i,j}
$$
in $\mathcal{S}'(\rrn)$, and
$$
\sup_{i\in\mathbb{Z}}\left\{2^i
\left\|\sum_{j\in\mathbb{N}}
\1bf_{B_{i,j}}\right\|_{\NHerzS}\right\}<\infty,
$$
Moreover, for any $f\in
WH\Kmc_{\omega}^{p,q,r,d}(\rrn)$,
$$
\|f\|_{WH\Kmc_{\omega}^{p,q,r,d}(\rrn)}:=
\inf\left\{\sup_{i\in\mathbb{Z}}2^i
\left\|\sum_{j\in\mathbb{N}}
\1bf_{B_{i,j}}\right\|_{\NHerzS}\right\},
$$
where the infimum is taken over all the
decompositions of $f$ as above.
\end{definition}

Then we have the following
atomic characterization of $W\NHaSaH$.

\begin{theorem}\label{watomgn}
Let $p$, $q$, $\omega$, $r$, and $d$
be as in Definition \ref{dfwatomgn}.
Then
$$
W\NHaSaH=WH\Kmc_{\omega}^{p,q,r,d}(\rrn)
$$
with equivalent quasi-norms.
\end{theorem}

\begin{proof}
Let all the symbols be as in
the present theorem. Then, using
Propositions \ref{watomgd} and \ref{watomgr}
and repeating an argument similar to
that used in the proof of Theorem
\ref{watomg} via replacing $\mw$, $\Mw$,
Theorem \ref{Th2}, Lemmas \ref{convexl},
\ref{ee}, \ref{equa}, \ref{mbhg}, \ref{mbhag},
\ref{Atogl2}, and \ref{wmbhg} therein, respectively,
by $\mi(\omega)$, $\MI(\omega)$,
Theorem \ref{Th2n}, Lemmas \ref{convexln},
\ref{een}, \ref{equan}, \ref{mbhgn},
\ref{mbhagn}, \ref{Atognl2}, and \ref{wmbhgn},
we find that $$
W\NHaSaH=WH\Kmc_{\omega}^{p,q,r,d}(\rrn)
$$
with equivalent quasi-norms. This finishes
the proof of Theorem \ref{watomgn}.
\end{proof}

\subsection{Molecular Characterizations}

In this subsection,  we establish the molecular
characterizations of inhomogeneous weak generalized
Herz--Hardy spaces. To this end, via
$(\NHerzSo,\,r,\,d,\,\tau)$-molecules introduced in Definition
\ref{molen}, we first show the
following molecular characterization
of the inhomogeneous
weak generalized Herz--Hardy space $\WNHaSaHo$.

\begin{theorem}\label{wmolen}
Let $p,\ q\in(0,\infty)$, $\omega\in M(\rp)$ satisfy
$\mi(\omega)\in(-\frac{n}{p},\infty)$,
$$
p_-\in\left(0,\frac{\min\{p,q,\frac{n}{\MI(\omega)+n/p}\}}
{\max\{1,p,q\}}\right),
$$
$d\geq\lfloor n(1/p_--1)\rfloor$ be a fixed
integer,
$$
r\in\left(\max\left\{1,p,\frac{n}{\mi(\omega)+n/p}
\right\},\infty\right],
$$
and $\tau\in(n(\frac{1}{p_-}-
\frac{1}{r}),\infty)$.
Then $f$ belongs to the inhomogeneous
weak generalized Herz--Hardy space
$WH\NHerzSo$ if and only if $f\in\mathcal{S}'(\rrn)$
and there exist a sequence $\{m_{i,j}\}_{i\in
\mathbb{Z},\,j\in\mathbb{N}}$ of
$(\NHerzSo,\,r,\,d,\,\tau)$-molecules centered,
respectively, at the balls $\{B_{i,j}
\}_{i\in\mathbb{Z},\,j\in\mathbb{N}}\subset\mathbb{B}$
and three positive constants $c\in(0,1]$,
$A$, and $\widetilde{A}$, independent of $f$, such that,
for any $i\in\mathbb{Z}$,
$$
\sum_{j\in\mathbb{N}}\1bf_{cB_{i,j}}\leq A,
$$
$$
f=\sum_{i\in\mathbb{Z}}\sum_{j\in\mathbb{N}}
\widetilde{A}2^i\left\|\1bf_{
B_{i,j}}\right\|_{\NHerzSo}m_{i,j}
$$
in $\mathcal{S}'(\rrn)$, and
$$
\sup_{i\in\mathbb{Z}}\left\{2^i
\left\|\sum_{j\in\mathbb{N}}
\1bf_{B_{i,j}}\right\|_{\NHerzSo}\right\}<\infty.
$$
Moreover, there exists a constant
$C\in[1,\infty)$ such that,
for any $f\in W\NHaSaHo$,
\begin{align*}
C^{-1}\|f\|_{W\NHaSaHo}
&\leq\inf\left\{\sup_{i\in\mathbb{Z}}
2^i\left\|\sum_{j\in\mathbb{N}}
\1bf_{B_{i,j}}\right\|_{\NHerzSo}\right\}\\
&\leq C\|f\|_{W\NHaSaHo},
\end{align*}
where the infimum is taken over
all the decompositions of
$f$ as above.
\end{theorem}

\begin{proof}
Let all the symbols be as in the
present theorem and $f\in\mathcal{S}'(\rrn)$.
Then, from Lemma \ref{wmoler} and
an argument similar to that
used in the proof of Theorem \ref{wmole}
via replacing $\mw$, Definition
\ref{dfwatom}, Theorems \ref{Th3} and
\ref{watom}, and Lemmas \ref{mbhal} and
\ref{vmbhl} therein, respectively,
by $\mi(\omega)$, Definition
\ref{dfwatomn}, Theorems \ref{Th3n} and
\ref{watomn}, and Lemmas \ref{mbhaln} and
\ref{vmbhln}, it follows that
$f\in W\NHaSaHo$ if and only if
there exist a sequence $\{m_{i,j}\}_{i\in
\mathbb{Z},\,j\in\mathbb{N}}$ of
$(\NHerzSo,\,r,\,d,\,\tau)$-molecules centered,
respectively, at the balls $\{B_{i,j}
\}_{i\in\mathbb{Z},\,j\in\mathbb{N}}\subset\mathbb{B}$
and three positive constants $c\in(0,1]$,
$A$, and $\widetilde{A}$, independent of $f$, such that,
for any $i\in\mathbb{Z}$,
$$
\sum_{j\in\mathbb{N}}\1bf_{cB_{i,j}}\leq A,
$$
$$
f=\sum_{i\in\mathbb{Z}}\sum_{j\in\mathbb{N}}
\widetilde{A}2^i\left\|\1bf_{
B_{i,j}}\right\|_{\NHerzSo}m_{i,j}
$$
in $\mathcal{S}'(\rrn)$, and
$$
\sup_{i\in\mathbb{Z}}\left\{2^i
\left\|\sum_{j\in\mathbb{N}}
\1bf_{B_{i,j}}\right\|_{\NHerzSo}\right\}<\infty.
$$
Moreover, it also holds true that
\begin{align*}
\|f\|_{W\NHaSaHo}
\sim\inf\left\{\sup_{i\in\mathbb{Z}}
2^i\left\|\sum_{j\in\mathbb{N}}
\1bf_{B_{i,j}}\right\|_{\NHerzSo}\right\}
\end{align*}
with the positive equivalence constants
independent of $f$,
where the infimum is taken over
all the decompositions of
$f$ as in the present theorem. This
finishes the proof of Theorem \ref{molen}.
\end{proof}

In addition, recall that the concept
of $(\NHerzS,\,r,\,d,\,\tau)$-molecules
is given in Definition \ref{molegn}. Then,
via these molecules,
the following conclusion shows the
molecular characterization of the inhomogeneous
weak generalized Herz--Hardy space $W\NHaSaH$.

\begin{theorem}\label{wmolegn}
Let $p,\ q\in(0,\infty)$, $\omega\in M(\rp)$ satisfy
$$-\frac{n}{p}<
\mi(\omega)\leq\MI(\omega)<0,$$
$
p_-\in(0,\min\{p,q\}/\max\{1,p,q\}),
$
$d\geq\lfloor n(1/p_--1)\rfloor$ be a fixed
integer,
$$
r\in\left(\max\left\{1,
\frac{n}{\mi(\omega)+n/p}\right\},\infty\right],
$$
and $\tau\in(n(\frac{1}{p_-}-
\frac{1}{r}),\infty)$.
Then $f\in WH\NHerzS$ if and
only if $f\in\mathcal{S}'(\rrn)$
and there exist a sequence $\{m_{i,j}\}_{i\in
\mathbb{Z},\,j\in\mathbb{N}}$ of
$(\NHerzS,\,r,\,d,\,\tau)$-molecules centered,
respectively, at the balls $\{B_{i,j}
\}_{i\in\mathbb{Z},\,j\in\mathbb{N}}\subset\mathbb{B}$
and three positive constants $c\in(0,1]$,
$A$, and $\widetilde{A}$, independent of $f$, such that,
for any $i\in\mathbb{Z}$,
$$
\sum_{j\in\mathbb{N}}\1bf_{cB_{i,j}}\leq A,
$$
$$
f=\sum_{i\in\mathbb{Z}}\sum_{j\in\mathbb{N}}
\widetilde{A}2^i\left\|\1bf_{
B_{i,j}}\right\|_{\NHerzS}m_{i,j}
$$
in $\mathcal{S}'(\rrn)$, and
$$
\sup_{i\in\mathbb{Z}}\left\{2^i
\left\|\sum_{j\in\mathbb{N}}
\1bf_{B_{i,j}}\right\|_{\NHerzS}\right\}<\infty.
$$
Moreover, there exists a constant
$C\in[1,\infty)$ such that,
for any $f\in W\NHaSaH$,
\begin{align*}
C^{-1}\|f\|_{W\NHaSaH}
&\leq\inf\left\{\sup_{i\in\mathbb{Z}}
2^i\left\|\sum_{j\in\mathbb{N}}
\1bf_{B_{i,j}}\right\|_{\NHerzS}\right\}\\
&\leq C\|f\|_{W\NHaSaH},
\end{align*}
where the infimum is taken over all
the decompositions of
$f$ as above.
\end{theorem}

\begin{proof}
Let all the symbols be as in
the present theorem and $f\in\mathcal{S}'(\rrn)$.
Then, using Proposition \ref{wmolegr}
and repeating an argument similar to that
used in the proof of Theorem \ref{wmoleg} via
replacing $\mw$, Definition \ref{watomg},
Theorems \ref{Th2}, \ref{ball}, and
\ref{watomg}, and Lemmas \ref{mbhag}, \ref{vmbhg},
and \ref{Atogl2} therein, respectively, by
$\mi(\omega)$, Definition \ref{watomgn},
Theorems \ref{Th2n}, \ref{balln}, and
\ref{watomgn}, and Lemmas \ref{mbhagn}, \ref{vmbhgn},
and \ref{Atognl2}, we conclude that
$f\in WH\NHerzS$ if and
only if there exist a sequence $\{m_{i,j}\}_{i\in
\mathbb{Z},\,j\in\mathbb{N}}$ of
$(\NHerzS,\,r,\,d,\,\tau)$-molecules centered,
respectively, at the balls $\{B_{i,j}
\}_{i\in\mathbb{Z},\,j\in\mathbb{N}}\subset\mathbb{B}$
and three positive constants $c\in(0,1]$,
$A$, and $\widetilde{A}$, independent of $f$, such that,
for any $i\in\mathbb{Z}$,
$$
\sum_{j\in\mathbb{N}}\1bf_{cB_{i,j}}\leq A,
$$
$$
f=\sum_{i\in\mathbb{Z}}\sum_{j\in\mathbb{N}}
\widetilde{A}2^i\left\|\1bf_{
B_{i,j}}\right\|_{\NHerzS}m_{i,j}
$$
in $\mathcal{S}'(\rrn)$, and
$$
\sup_{i\in\mathbb{Z}}\left\{2^i
\left\|\sum_{j\in\mathbb{N}}
\1bf_{B_{i,j}}\right\|_{\NHerzS}\right\}<\infty.
$$
Moreover, we also obtain
\begin{align*}
\|f\|_{W\NHaSaH}
\sim\inf\left\{\sup_{i\in\mathbb{Z}}
2^i\left\|\sum_{j\in\mathbb{N}}
\1bf_{B_{i,j}}\right\|_{\NHerzS}\right\}
\end{align*}
with the positive equivalence constants
independent of $f$,
where the infimum is taken over all
the decompositions of
$f$ as in the present theorem. This
then finishes the proof of Theorem \ref{wmolegn}.
\end{proof}

\subsection{Littlewood--Paley Function Characterizations}

This subsection is devoted to
characterizing both the
inhomogeneous weak generalized Herz--Hardy spaces
$\WNHaSaHo$ and $\WNHaSaH$ via the Lusin area function,
the $g$-function, and the $g_\lambda^*$-function
defined as in Definitions \ref{df451}
and \ref{df452}. First, we give
the following Littlewood--Paley function
characterization of $\WNHaSaHo$.

\begin{theorem}\label{wlusinn}
Let $p,\ q\in(0,\infty)$, $\omega\in M(\rp)$
with $\mi(\omega)\in(-\frac{n}{p},\infty)$,
$$
p_-:=\min\left\{p,q,\frac{n}{\MI(\omega)+n/p}\right\},
$$
$$
p_+:=\max\left\{p,\frac{n}{\mi(\omega)+n/p}\right\},
$$
and
$$
\lambda\in\left(\max\left\{
\frac{2}{\min\{1,p_-\}},\,1-\frac{2}{\max\{1,p_+\}}+
\frac{2}{\min\{1,p_-\}}\right\},\infty\right).
$$
Assume that, for any
$f\in\mathcal{S}'(\rrn)$, $S(f)$ and
$g_{\lambda}^*(f)$
are as in Definition \ref{df451}, and
$g(f)$ is as in Definition \ref{df452}.
Then the following four statements are mutually
equivalent:
\begin{enumerate}
  \item[{\rm(i)}] $f\in W\NHaSaHo$;
  \item[{\rm(ii)}] $f\in\mathcal{S}'(\rrn)$,
  $f$ vanishes weakly at infinity,
  and $S(f)\in W\NHerzSo$;
  \item[{\rm(iii)}] $f\in\mathcal{S}'(\rrn)$,
  $f$ vanishes weakly at infinity,
  and $g(f)\in W\NHerzSo$;
  \item[{\rm(iv)}] $f\in\mathcal{S}'(\rrn)$,
  $f$ vanishes weakly at infinity,
  and $g_{\lambda}^*(f)\in W\NHerzSo$.
\end{enumerate}
Moreover, for any $f\in W\NHaSaHo$,
\begin{align*}
\left\|f\right\|_{W\NHaSaHo}
&\sim\left\|S(f)\right\|_{W\NHerzSo}
\sim\left\|g(f)\right\|_{W\NHerzSo}\\
&\sim\left\|g_{\lambda}^*(f)\right\|_{W\NHerzSo},
\end{align*}
where the positive equivalence
constants are independent of $f$.
\end{theorem}

To prove this theorem, we need the following
auxiliary lemma related to the Fefferman--Stein
vector-valued inequality on
inhomogeneous weak local generalized Herz spaces.

\begin{lemma}\label{vwmbhln}
Let $p,\ q\in(0,\infty)$ and $\omega\in M(\rp)$ satisfy
$\mi(\omega)\in(-\frac{n}{p},\infty)$.
Then, for any given
$$r\in\left(0,\min\left\{p,
\frac{n}{\MI(\omega)+n/p}\right\}\right)$$
and $u\in(1,\infty)$,
there exists a positive constant $C$ such that,
for any $\{f_{j}\}_{
j\in\mathbb{N}}\subset L^1_{\mathrm{loc}}(\rrn)$,
$$
\left\|\left\{\sum_{j\in\mathbb{N}}\left[\mc
(f_{j})\right]^{u}\right\}^
{\frac{1}{u}}\right\|_{[W\NHerzSo]^{1/r}}
\leq C\left\|\left\{
\sum_{j\in\mathbb{N}}|f_{j}|^{u}\right\}
^{\frac{1}{u}}\right\|_{[W\NHerzSo]^{1/r}}.
$$
\end{lemma}

\begin{proof}
Let all the symbols be as
in the present lemma, $r\in(0,\min\{p,
\frac{n}{\MI(\omega)+n/p}\}),$
and $u\in(1,\infty)$. Then,
using Lemma \ref{vwmbhll1} and repeating
the proof of Lemma \ref{vwmbhl} with
$\Mw$, Theorem \ref{Th3}, and Lemma
\ref{vmbhl} therein replaced, respectively,
by $\MI(\omega)$, Theorem \ref{Th3n}, and Lemma
\ref{vmbhln}, we find that, for any $\{f_{j}\}_{
j\in\mathbb{N}}\subset L^1_{\mathrm{loc}}(\rrn)$,
$$
\left\|\left\{\sum_{j\in\mathbb{N}}\left[\mc
(f_{j})\right]^{u}\right\}^
{\frac{1}{u}}\right\|_{[W\NHerzSo]^{1/r}}
\lesssim\left\|\left\{
\sum_{j\in\mathbb{N}}|f_{j}|^{u}\right\}
^{\frac{1}{u}}\right\|_{[W\NHerzSo]^{1/r}}.
$$
This then finishes the proof of Lemma \ref{vwmbhln}.
\end{proof}

Via Lemma \ref{vwmbhln}, we next
prove Theorem \ref{wlusinn}.

\begin{proof}[Proof of Theorem
\ref{wlusinn}]
Let all the symbols be as in
the present theorem. Then, repeating
an argument similar to that used
in the proof of Theorem \ref{wlusin}
with Theorem \ref{Th3} and
Lemmas \ref{convexll}, \ref{mbhal},
\ref{vmbhl}, and \ref{vwmbhl} therein replaced,
respectively, by Theorem \ref{Th3n} and
Lemmas \ref{convexlln}, \ref{mbhaln},
\ref{vmbhln}, and \ref{vwmbhln}, we infer
that (i) through (iv) of the present
theorem are mutually equivalent. Moreover,
we also find that, for any $f\in W\NHaSaHo$,
\begin{align*}
\left\|f\right\|_{W\NHaSaHo}
&\sim\left\|S(f)\right\|_{W\NHerzSo}
\sim\left\|g(f)\right\|_{W\NHerzSo}\\
&\sim\left\|g_{\lambda}^*(f)\right\|_{W\NHerzSo},
\end{align*}
which then completes the proof of Theorem
\ref{wlusinn}.
\end{proof}

Now, we turn to establish the following
Littlewood--Paley function characterization
of the inhomogeneous weak generalized Herz--Hardy
space $W\NHaSaH$.

\begin{theorem}\label{wlusingn}
Let $p,\ q\in(0,\infty)$, $\omega\in M(\rp)$
satisfy
$$
-\frac{n}{p}<\mi(\omega)\leq\MI(\omega)<0,
$$
$
p_-:=\min\{p,q\},
$
$
p_+:=\max\{p,\frac{n}{\mi(\omega)+n/p}\},
$
and
$$
\lambda\in\left(\max\left\{
\frac{2}{\min\{1,p_-\}},\,1-\frac{2}{\max\{1,p_+\}}+
\frac{2}{\min\{1,p_-\}}\right\},\infty\right).
$$
Assume that, for any $f\in\mathcal{S}'(\rrn)$, $S(f)$
and $g_{\lambda}^*(f)$ are as in Definition
\ref{df451}, and $g(f)$ is as in Definition
\ref{df452}.
Then the following four statements are mutually
equivalent:
\begin{enumerate}
  \item[{\rm(i)}] $f\in W\NHaSaH$;
  \item[{\rm(ii)}] $f\in\mathcal{S}'(\rrn)$,
  $f$ vanishes weakly at infinity,
  and $S(f)\in W\NHerzS$;
  \item[{\rm(iii)}] $f\in\mathcal{S}'(\rrn)$,
  $f$ vanishes weakly at infinity,
  and $g(f)\in W\NHerzS$;
  \item[{\rm(iv)}] $f\in\mathcal{S}'(\rrn)$,
  $f$ vanishes weakly at infinity,
  and $g_{\lambda}^*(f)\in W\NHerzS$.
\end{enumerate}
Moreover, for any $f\in W\NHaSaH$,
\begin{align*}
\left\|f\right\|_{W\NHaSaH}
&\sim\left\|S(f)\right\|_{W\NHerzS}
\sim\left\|g(f)\right\|_{W\NHerzS}\\
&\sim\left\|g_{\lambda}^*(f)\right\|_{W\NHerzS},
\end{align*}
where the positive equivalence
constants are independent of $f$.
\end{theorem}

To show the above theorem, we first
establish the following representation
formula of the quasi-norm $\|\cdot\|_{W\NHerzS}$.

\begin{lemma}\label{wlusingnl1}
Let $p$, $q\in(0,\infty)$ and
$\omega\in M(\rp)$. Then, for any
$f\in\Msc(\rrn)$,
$$
\left\|f\right\|_{W\NHerzS}
=\sup_{\xi\in\rrn}\left\|\tau_{\xi}\left(f\right)
\right\|_{W\NHerzSo},
$$
where, for any $\xi\in\rrn$,
the operator $\tau_{\xi}$ is defined
as in \eqref{translation}.
\end{lemma}

\begin{proof}
Let all the symbols be as in the
present lemma and $f\in\Msc(\rrn)$.
Then, from Definition
\ref{df6.0.7n}(ii), Remark \ref{ighr}(i),
and Definition \ref{df6.0.7n}(i), it follows that
\begin{align*}
\left\|f\right\|_{W\NHerzS}
&=\sup_{\lambda\in(0,\infty)}
\left[\lambda\left\|\1bf_{\{y\in\rrn:\
|f(y)|>\lambda\}}\right\|_{\NHerzS}\right]\\
&=\sup_{\lambda\in(0,\infty)}
\left\{\lambda\sup_{\xi\in\rrn}\left[
\left\|\1bf_{\{y\in\rrn:\
|f(y)|>\lambda\}}\left(\cdot-\xi\right)
\right\|_{\NHerzSo}\right]\right\}\\
&=\sup_{\lambda\in(0,\infty),\,\xi\in\rrn}
\left[\lambda\left\|\1bf_{\{y\in\rrn:\
|\tau_{\xi}(f)(y)|>\lambda\}}\right\|_{\NHerzSo}
\right]\\&=\sup_{\xi\in\rrn}\left\{\sup_{\lambda\in(0,\infty)}
\left[\lambda\left\|\1bf_{\{y\in\rrn:\
|\tau_{\xi}(f)(y)|>\lambda\}}\right\|_{\NHerzSo}
\right]\right\}\\
&=\sup_{\xi\in\rrn}\left\|\tau_{\xi}\left(f\right)
\right\|_{\NHerzSo}.
\end{align*}
This then finishes the proof of Lemma \ref{wlusingnl1}.
\end{proof}

We next give the proof of Theorem \ref{wlusingn}
as follows.

\begin{proof}[Proof of Theorem \ref{wlusingn}]
Let all the symbols be as
in the present theorem. Then, repeating
the proof of Theorem \ref{wlusing} via
replacing Theorems \ref{macl}(ii), \ref{Th7.11}(ii),
and \ref{wlusin}, and Lemma \ref{wlusingl1} therein,
respectively, by Theorems \ref{macln}(ii), \ref{Th7.11n}(ii),
and \ref{wlusinn}, and Lemma \ref{wlusingnl1},
we infer that (i), (ii), (iii), and (iv)
are mutually equivalent and, for any $f\in W\NHaSaH$,
\begin{align*}
\left\|f\right\|_{W\NHaSaH}
&\sim\left\|S(f)\right\|_{W\NHerzS}
\sim\left\|g(f)\right\|_{W\NHerzS}\\
&\sim\left\|g_{\lambda}^*(f)\right\|_{W\NHerzS}.
\end{align*}
This further finishes the proof of
Theorem \ref{wlusingn}.
\end{proof}

\subsection{Boundedness of Calder\'{o}n--Zygmund Operators}

In this subsection, we establish the boundedness
of Calder\'{o}n--Zygmund operators from
inhomogeneous generalized Herz--Hardy spaces
to inhomogeneous weak generalized
Herz--Hardy spaces even in the critical case.
Let $d\in\zp$. Recall that
the standard kernel $K$ and
the $d$-order Calder\'{o}n--Zygmund
operator $T$ with kernel $K$ are defined,
respectively, in Definitions \ref{def-s-k}
and \ref{defin-C-Z-s}.
To begin with, we show the following boundedness
of the operator $T$ from
$\NHaSaHo$ to $W\NHaSaHo$.

\begin{theorem}\label{wczoon}
Let $d\in\zp$, $\delta\in(0,1)$,
$p\in[\frac{n}{n+d+\delta},\infty)$,
$q\in(0,\infty)$,
$K$ be a $d$-order standard kernel as
in Definition \ref{def-s-k}, $T$ a
$d$-order Calder\'{o}n--Zygmund operator
with kernel $K$ having the vanishing moments
up to order $d$, and $\omega\in
M(\rp)$ satisfy
$$
-\frac{n}{p}<\mi(\omega)\leq
\MI(\omega)<n-\frac{n}{p}+d+\delta.
$$
Then $T$ has a unique extension
on $\NHaSaHo$ and there exists
a positive constant $C$ such that,
for any $f\in\NHaSaHo$,
$$
\|T(f)\|_{W\NHaSaHo}\leq C\|f\|_{
\NHaSaHo}.
$$
\end{theorem}

To prove this theorem, we first show
the following technique estimate
about the Hardy--Littlewood
maximal operator.

\begin{lemma}\label{sizev}
Let $r\in[1,\infty)$ and
$\mc$ be the Hardy--Littlewood
maximal operator defined
as in \eqref{hlmax}. Then
there exists a positive constant $C$
such that, for any $\{f_{j}\}
_{j\in\mathbb{N}}\subset\Msc(\rrn)$ and
$x\notin\bigcup_{j\in\mathbb{N}}\supp(f_j)$,
$$
\left\{\sum_{j\in\mathbb{N}}\left[
\mc\left(f_j\right)(x)\right]^r
\right\}^{\frac{1}{r}}\leq C
\int_{\rrn}\frac{[\sum_{j\in\mathbb{N}}
|f_j(y)|^r]^{\frac{1}{r}}}{|x-y|^n}\,dy.
$$
\end{lemma}

\begin{proof}
Let all the symbols be as in the present lemma,
and let $\{f_j\}_{j\in\mathbb{N}}
\subset\Msc(\rrn)$.
Then, by the fact that the
Hardy--Littlewood maximal operator
$\mc$ satisfies the size
condition \eqref{size} (see, for instance,
\cite[Remark 4.4]{HSamko}) and
the Minkowski integral inequality,
we conclude that, for any
$x\notin\bigcup_{j\in\mathbb{N}}\supp(f_j)$,
\begin{align*}
\left\{\sum_{j\in\mathbb{N}}
\left[\mc(f_j)(x)\right]^{r}\right\}^{\frac{1}{r}}
&\lesssim\left\{\sum_{j\in\mathbb{N}}
\left[\int_{\rrn}\frac{|f_j(y)|}{|x-y|^{n}}
\,dy\right]^r\right\}^{\frac{1}{r}}\\
&\lesssim\int_{\rrn}\left[
\sum_{j\in\mathbb{N}}
\frac{|f_j(y)|^r}{|x-y|^{nr}}
\right]^{\frac{1}{r}}\,dy\\
&\sim\int_{\rrn}\frac{[\sum_{j\in\mathbb{N}}
|f_j(y)|^r]^{\frac{1}{r}}}{|x-y|^n}\,dy.
\end{align*}
This then finishes the proof of Lemma \ref{sizen}.
\end{proof}

Via this lemma, we now establish
the following weak type vector-valued inequalities
on inhomogeneous generalized Herz spaces,
which paly essential roles in the proof of
the boundedness of Calder\'{o}n--Zygmund
operators from inhomogeneous generalized
Herz--Hardy spaces to inhomogeneous
weak generalized Herz--Hardy spaces.

\begin{proposition}\label{wczool1n}
Let $p\in[1,\infty)$, $q\in(0,\infty)$,
$r\in(1,\infty)$, and $\omega\in M(\rp)$
with
$$
-\frac{n}{p}<\mi(\omega)\leq\MI
(\omega)<\frac{n}{p'},
$$
where $\frac{1}{p}+\frac{1}{p'}=1$.
Then there exists a positive
constant $C$ such that, for any
$\{f_{j}\}_{j\in\mathbb{N}}\subset
L^{1}_{{\rm loc}}(\rrn)$,
\begin{equation}\label{wczoonl1e0}
\left\|\mvr\right\|_{W\NHerzSo}
\leq C\left\|\vr\right\|_{\NHerzSo}
\end{equation}
and
\begin{equation}\label{wczoonl1e1}
\left\|\mvr\right\|_{W\NHerzS}\leq
C\left\|\vr\right\|_{\NHerzS}.
\end{equation}
\end{proposition}

\begin{proof}
Let all the symbols be as in
the present proposition,
and let $\{f_j\}_{j\in\mathbb{N}}\subset
\Msc(\rrn)$.
We next show the present proposition
by considering the following
two cases on $p$.

\emph{Case 1)} $p\in(1,\infty)$.
In this case, from Definition
\ref{df6.0.7n} and Theorems \ref{Th3.4n}
and \ref{Th3.3n}, it follows that
\begin{align*}
&\left\|\mvr\right\|_{W\NHerzSo}\\
&\quad\leq\left\|\mvr\right\|_{\NHerzSo}
\lesssim\left\|\vr\right\|_{\NHerzSo}
\end{align*}
and
\begin{align*}
&\left\|\mvr\right\|_{W\NHerzS}\\&\quad\leq
\left\|\mvr\right\|_{\NHerzS}
\lesssim\left\|\vr\right\|_{\NHerzS}.
\end{align*}
These finish the proof of
the present proposition in this case.

\emph{Case 2)} $p=1$. In this case,
fix a $\xi\in\rrn$. Then, using Definition
\ref{igh}(i), we find that, for any
$\lambda\in(0,\infty)$,
\begin{align}\label{wczoonl1e2}
&\left\|\1bf_{\{y\in\rrn:\ \{\sum_{j\in\mathbb{N}}
[\mc(f_j)(y)]^r\}^{\frac{1}{r}}>\lambda\}
}\left(\cdot+\xi\right)
\right\|_{\Kmc_{\omega,\0bf}^{1,q}(\rrn)}\notag\\
&\quad\lesssim\left\|\1bf_{\{y\in\rrn:\
\{\sum_{j\in\mathbb{N}}
[\mc(f_j)(y)]^r\}^{\frac{1}{r}}>\lambda\}
}\1bf_{B(\xi,1)}\right\|_{L^1(\rrn)}\notag\\
&\qquad+\left\{\sum_{k\in\mathbb{N}}
\left[\omega(\tk)\right]^q\right.\notag\\
&\qquad\left.\times
\left\|\1bf_{\{y\in\rrn:\ \{\sum_{j\in\mathbb{N}}
[\mc(f_j)(y)]^r\}^{\frac{1}{r}}>\lambda\}
}\1bf_{B(\xi,\tk)\setminus B(\xi,\tkm)}
\right\|_{L^1(\rrn)}^q
\right\}^{\frac{1}{q}}\notag\\
&\quad=:\mathrm{A}_{\lambda}+
\mathrm{B}_{\lambda}.
\end{align}
Now, we deal with $\mathrm{A}_{\lambda}$
and $\mathrm{B}_{\lambda}$, respectively.
Indeed, for any $\lambda\in(0,\infty)$,
applying the fact that
$$
\1bf_{\{y\in\rrn:\ \{\sum_{j\in\mathbb{N}}
[\mc(f_j\1bf_{[B(\xi,2)]
^{\complement}})(y)]^r\}^{\frac{1}{r}}>\lambda\}}
\leq\lambda^{-1}\left\{\sum_{j\in\mathbb{N}}
\left[\mc\left(f_j\1bf_{[B(\xi,2)]
^{\complement}}\right)\right]^r
\right\}^{\frac{1}{r}}
$$
and Lemma \ref{wczool10}
with $\{f_j\}_{j\in\mathbb{N}}$
therein replaced by $\{f_j
\1bf_{B(\xi,2)}\}_{j\in\mathbb{N}}$,
we conclude that
\begin{align}\label{wczoonl1e3}
\mathrm{A}_{\lambda}&\lesssim
\left\|\1bf_{\{y\in\rrn:\
\{\sum_{j\in\mathbb{N}}
[\mc(f_j\1bf_{B(\xi,2)})
(y)]^r\}^{\frac{1}{r}}>\lambda\}
}\1bf_{B(\xi,1)}\right\|_{L^1(\rrn)}\notag\\
&\quad+\left\|\1bf_{\{y\in\rrn:\
\{\sum_{j\in\mathbb{N}}
[\mc(f_j\1bf_{[B(\xi,2)]^{\complement}})
(y)]^r\}^{\frac{1}{r}}>\lambda\}
}\1bf_{B(\xi,1)}\right\|_{L^1(\rrn)}\notag\\
&\lesssim\left\|\1bf_{\{y\in\rrn:\
\{\sum_{j\in\mathbb{N}}
[\mc(f_j\1bf_{B(\xi,2)})
(y)]^r\}^{\frac{1}{r}}>\lambda\}}
\right\|_{L^1(\rrn)}\notag\\
&\quad+\lambda^{-1}\left\|
\left\{\sum_{j\in\mathbb{N}}\left[
\mc\left(f_j\1bf_{[B(\xi,2)]^{\complement}}
\right)\right]^{r}
\right\}^{\frac{1}{r}}\right\|_{L^1(\rrn)}\notag\\
&\lesssim\lambda^{-1}\left[
\left\|\left(\sum_{j\in\mathbb{N}}|f_j|^r
\right)^{\frac{1}{r}}\1bf_{B(\xi,2)}
\right\|_{L^1(\rrn)}\right.\notag\\&\quad\left.+
\left\|\left\{\sum_{j\in\mathbb{N}}\left[
\mc\left(f_j\1bf_{[B(\xi,2)]^{\complement}}
\right)\right]^{r}
\right\}^{\frac{1}{r}}\right\|_{L^1(\rrn)}\right].
\end{align}
In addition,
from Lemma \ref{sizev} with $\{f_j\}_{j\in\mathbb{N}}$
therein replaced by $\{f_j\1bf_{
[B(\xi,2)]^{\complement}}\}_{j\in\mathbb{N}}$,
we infer that, for any
$x\in B(\xi,1)$,
\begin{equation}\label{wczoonl1e4}
\left\{\sum_{j\in\mathbb{N}}
\left[\mc\left(f_j
\1bf_{[B(\xi,2)]^{\complement}}
\right)(x)\right]^r\right\}^{\frac{1}{r}}
\lesssim\int_{|y-\xi|\geq2}
\frac{[\sum_{j\in\mathbb{N}}|f_j(y)|^r]^{
\frac{1}{r}}}{|x-y|^n}\,dy.
\end{equation}
Therefore, using \eqref{wczoonl1e3}
and repeating an argument similar to
that used in the estimations of
\eqref{vmaxhne3}, \eqref{vmaxhne4},
and \eqref{vmaxhne5} with \eqref{sizen},
$\0bf$, $|f|$, and $|T(f)|$ therein replaced,
respectively, by \eqref{wczoonl1e4}, $\xi$,
$(\sum_{j\in\mathbb{N}}|f_j|^r)^{\frac{1}{r}}$,
and $\{\sum_{j\in\mathbb{N}}
[\mc(f_j)]^r\}^{\frac{1}{r}}$, we find
that, for any $\lambda\in(0,\infty)$,
\begin{align}\label{wczoonl1e5}
\mathrm{A}_{\lambda}\lesssim\lambda^{-1}\left\|
\vr\left(\cdot+\xi\right)
\right\|_{\Kmc_{\omega,\0bf}^{1,q}(\rrn)},
\end{align}
which is the desired estimate of
$\mathrm{A}_{\lambda}$.

Next, we deal with $\mathrm{B}_{\lambda}$.
For this purpose, for any $k\in\mathbb{N}$,
let
$$f_{j,k,0}:=f_j\1bf_{B(\xi,2^{-2})},$$
$$f_{j,k,1}:=f_j\1bf_{B(\xi,2^{k-2})\setminus
B(\xi,2^{-2})},$$
$$f_{j,k,2}:=f_j\1bf_{B(\xi,2^{k+1})\setminus
B(\xi,2^{k-2})},$$
and
$$f_{j,k,3}:=f_j\1bf_{[B(\xi,2^{k+1})]^{\complement}}.$$
Then, for any $\lambda\in(0,\infty)$, we have
\begin{align}\label{wczoonl1e6}
\mathrm{B}_{\lambda}&\lesssim
\left\{\sum_{k\in\mathbb{N}}
\left[\omega(\tk)\right]^q\right.\notag\\
&\quad\left.\times
\left\|\1bf_{\{y\in\rrn:\ \{\sum_{j\in\mathbb{N}}
[\mc(f_{j,k,0})(y)]^r\}^{\frac{1}{r}}>\lambda\}
}\1bf_{B(\xi,\tk)\setminus B(\xi,\tkm)}
\right\|_{L^1(\rrn)}^q
\right\}^{\frac{1}{q}}\notag\\
&\quad+\left\{\sum_{k\in\mathbb{N}}
\left[\omega(\tk)\right]^q\right.\notag\\
&\quad\left.\times
\left\|\1bf_{\{y\in\rrn:\ \{\sum_{j\in\mathbb{N}}
[\mc(f_{j,k,1})(y)]^r\}^{\frac{1}{r}}>\lambda\}
}\1bf_{B(\xi,\tk)\setminus B(\xi,\tkm)}
\right\|_{L^1(\rrn)}^q
\right\}^{\frac{1}{q}}\notag\\
&\quad+\left\{\sum_{k\in\mathbb{N}}
\left[\omega(\tk)\right]^q\right.\notag\\
&\quad\left.\times
\left\|\1bf_{\{y\in\rrn:\ \{\sum_{j\in\mathbb{N}}
[\mc(f_{j,k,2})(y)]^r\}^{\frac{1}{r}}>\lambda\}
}\1bf_{B(\xi,\tk)\setminus B(\xi,\tkm)}
\right\|_{L^1(\rrn)}^q
\right\}^{\frac{1}{q}}\notag\\
&\quad+\left\{\sum_{k\in\mathbb{N}}
\left[\omega(\tk)\right]^q\right.\notag\\
&\quad\left.\times
\left\|\1bf_{\{y\in\rrn:\ \{\sum_{j\in\mathbb{N}}
[\mc(f_{j,k,3})(y)]^r\}^{\frac{1}{r}}>\lambda\}
}\1bf_{B(\xi,\tk)\setminus B(\xi,\tkm)}
\right\|_{L^1(\rrn)}^q
\right\}^{\frac{1}{q}}\notag\\
&=:\widetilde{\mathrm{II}}_{\lambda,0}+
\widetilde{\mathrm{II}}_{\lambda,1}
+\widetilde{\mathrm{II}}_{\lambda,2}+
\widetilde{\mathrm{II}}_{\lambda,3}.
\end{align}
Using Lemma \ref{sizev} with $\{f_j\}_{j\in\mathbb{N}}$
therein replaced by $\{f_{j,k,0}\}_{j\in\mathbb{N}}$,
we conclude that, for any $x\in\rrn$ with $|x|\geq1$,
\begin{align}\label{wczoonl1e7}
\left\{\sum_{j\in\mathbb{N}}
\left[\mc\left(f_{j,k,0}
\right)(x)\right]^r\right\}^{\frac{1}{r}}
\lesssim\int_{|y-\xi|<2^{-2}}
\frac{[\sum_{j\in\mathbb{N}}|f_j(y)|^r]^{
\frac{1}{r}}}{|x-y|^n}\,dy.
\end{align}
Thus, for any $\lambda\in(0,\infty)$,
from the fact that
$$
\1bf_{\{y\in\rrn:\ \{\sum_{j\in\mathbb{N}}
[\mc(f_{j,k,0})(y)]^r\}^{\frac{1}{r}}>\lambda\}}
\leq\lambda^{-1}\left\{\sum_{j\in\mathbb{N}}
\left[\mc\left(f_{j,k,0}\right)\right]^r
\right\}^{\frac{1}{r}}
$$
and an argument similar to that
used in the estimation of $\widetilde{\mathrm{J}}_0$
in the proof of Theorem \ref{vmaxhn}
with \eqref{sizen},
$\0bf$, $|f|$, and $|T(f)|$ therein replaced,
respectively, by \eqref{wczoonl1e7}, $\xi$,
$(\sum_{j\in\mathbb{N}}|f_j|^r)^{\frac{1}{r}}$,
and $\{\sum_{j\in\mathbb{N}}
[\mc(f_j)]^r\}^{\frac{1}{r}}$,
we obtain
\begin{align}\label{wczoonl1e8}
\widetilde{\mathrm{II}}_{\lambda,0}
&\lesssim\lambda^{-1}\left\{\sum_{j\in\mathbb{N}}
\left[\omega(\tk)\right]^q\left\|
\mvr\1bf_{B(\xi,\tk)\setminus
B(\xi,\tkm)}\right\|_{L^1(\rrn)}^q
\right\}^{\frac{1}{q}}\notag\\
&\lesssim\lambda^{-1}\left\|\vr
\left(\cdot+\xi\right)\right\|_{\Kmc_{\omega,\0bf}^{
1,q}(\rrn)},
\end{align}
which completes the estimation of
$\widetilde{\mathrm{II}}_{\lambda,0}$.
On the other hand, repeating
the estimations of $\widetilde{\mathrm{II}}_{
\lambda,1}$, $\widetilde{\mathrm{II}}_{\lambda,2}$,
and $\widetilde{\mathrm{II}}_{\lambda,3}$
in the proof of Proposition
\ref{wczool1} via
replacing $\mw$, $\Mw$, and $\mathbb{Z}$
therein, respectively, by
$\mi(\omega)$, $\MI(\omega)$, and $\mathbb{N}$,
we find that, for any $\lambda
\in(0,\infty)$,
\begin{equation*}
\widetilde{\mathrm{II}}_{\lambda,1}
\lesssim\lambda^{-1}\left\|\vr\left(
\cdot+\xi\right)\right\|_{\Kmc_{\omega,\0bf}^{
1,q}(\rrn)},
\end{equation*}
\begin{equation*}
\widetilde{\mathrm{II}}_{\lambda,2}
\lesssim\lambda^{-1}\left\|\vr\left(
\cdot+\xi\right)\right\|_{\Kmc_{\omega,\0bf}^{
1,q}(\rrn)},
\end{equation*}
and
\begin{equation*}
\widetilde{\mathrm{II}}_{\lambda,3}
\lesssim\lambda^{-1}\left\|\vr\left(
\cdot+\xi\right)\right\|_{\Kmc_{\omega,\0bf}^{
1,q}(\rrn)}.
\end{equation*}
These, together with both \eqref{wczoonl1e8}
and \eqref{wczoonl1e6}, further imply that,
for any $\lambda\in(0,\infty)$,
$$
\mathrm{B}_{\lambda}\lesssim\lambda^{-1}
\left\|\vr\left(
\cdot+\xi\right)\right\|_{\Kmc_{\omega,\0bf}^{
1,q}(\rrn)},
$$
which is the desired estimate of $
\mathrm{B}_{\lambda}$.
From this, \eqref{wczoonl1e5}, and
\eqref{wczoonl1e2}, we deduce that,
for any $\lambda\in(0,\infty)$,
\begin{align*}
&\lambda\left\|\1bf_{\{y\in\rrn:\ \{\sum_{j\in\mathbb{N}}
[\mc(f_j)(y)]^r\}^{\frac{1}{r}}>\lambda\}
}\left(\cdot+\xi\right)
\right\|_{\Kmc_{\omega,\0bf}^{1,q}(\rrn)}\\
&\quad\lesssim\left\|\vr\left(
\cdot+\xi\right)\right\|_{\Kmc_{\omega,\0bf}^{
1,q}(\rrn)}.
\end{align*}
Combining this and Definition \ref{df6.0.7n}(i),
we conclude that
\begin{align}\label{wczoonl1e9}
&\left\|\mvr
\left(\cdot+\xi\right)\right\|_{W
\Kmc_{\omega,\0bf}^{1,q}(\rrn)}\notag\\
&\quad\lesssim\left\|\vr\left(\cdot+\xi\right)
\right\|_{\Kmc_{\omega,\0bf}^{1,q}(\rrn)}.
\end{align}
Therefore, by both Definition
\ref{df6.0.7n}(ii) and the arbitrariness of $\xi$,
we have
$$
\left\|\mvr\right\|_{W\Kmc_{\omega}^{1,q}(\rrn)}\lesssim
\left\|\vr\right\|_{\Kmc_{\omega}^{1,q}(\rrn)}.
$$
In particular, letting $\xi:=\0bf$ in
\eqref{wczoonl1e9}, we obtain
$$
\left\|\mvr\right\|_{W\Kmc_{\omega,\0bf}^{1,q}(\rrn)}\lesssim
\left\|\vr\right\|_{\Kmc_{\omega,\0bf}^{1,q}(\rrn)}.
$$
This further implies that
both \eqref{wczoonl1e0} and \eqref{wczoonl1e1}
hold true in this case, and hence finishes the proof of
Proposition \ref{wczool1n}.
\end{proof}

Moreover, in order to show Theorem
\ref{wczoon}, we also require the following
convexification lemma of inhomogeneous
weak local generalized Herz spaces, which
can be immediately concluded
by Lemma \ref{convexlln} and
\cite[Remark 2.14]{ZWYY}; we omit the details.

\begin{lemma}\label{wconvexln}
Let $p,\ q,\ s\in(0,\infty)$ and $\omega\in M(\rp)$.
Then
$$\left[W\NHerzSo\right]^{1/s}=W\NHerzSocs$$
with the same quasi-norms.
\end{lemma}

Using both Proposition \ref{wczool1n}
and Lemma \ref{wconvexln},
we now prove Theorem \ref{wczoon}.

\begin{proof}[Proof of Theorem \ref{wczoon}]
Let all the symbols be as in the present theorem.
Then, from Proposition \ref{wczxa}
and an argument similar to that
used in the proof of Theorem \ref{wczoo} via
replacing \eqref{wczool1e0}, Theorems \ref{Th3}
and \ref{abso}, and Lemmas \ref{convexll},
\ref{mbhal}, \ref{Atoll3}, \ref{wmbhl},
and \ref{wconvexl} therein, respectively,
by \eqref{wczoonl1e0}, Theorems \ref{Th3n}
and \ref{abson}, and Lemmas \ref{convexlln},
\ref{mbhaln}, \ref{Atolnl3}, \ref{wmbhln},
and \ref{wconvexln}, we infer that
$T$ has a unique extension on
$\NHaSaHo$ and, for any $f\in\NHaSaHo$,
$$
\|T(f)\|_{W\NHaSaHo}\lesssim\|f\|_{
\NHaSaHo}.
$$
This then finishes the proof of Theorem
\ref{wczoon}.
\end{proof}

Next, we are devoted to
establishing the boundedness of
Calder\'{o}n--Zygmund
operators from the inhomogeneous
generalized Herz--Hardy space $\NHaSaH$ to
the inhomogeneous weak generalized Herz--Hardy
space $\WNHaSaH$ as follows.

\begin{theorem}\label{wczoogn}
Let $d\in\zp$, $\delta\in(0,1)$,
$p\in[\frac{n}{n+d+\delta},\infty)$,
$q\in(0,\infty)$,
$K$ be a $d$-order standard kernel as in
Definition \ref{def-s-k}, $T$ a
$d$-order Calder\'{o}n--Zygmund operator
with kernel $K$ having the vanishing moments
up to order $d$, and $\omega\in M(\rp)$ satisfy
$$
-\frac{n}{p}<\mi(\omega)\leq
\MI(\omega)<0.
$$
Then $T$ is well defined
on $\NHaSaH$ and there exists
a positive constant $C$ such that,
for any $f\in\NHaSaH$,
$$
\|T(f)\|_{W\NHaSaH}\leq C\|f\|_{
\NHaSaH}.
$$
\end{theorem}

To show this theorem, we
need a preliminary lemma.
Indeed, applying both Lemma \ref{convexln}
and \cite[Remark 2.14]{ZWYY},
we immediately obtain the
following convexification lemma
of inhomogeneous weak global
generalized Herz spaces, which plays
a key role in the proof of Theorem
\ref{wczoogn}; we omit the details.

\begin{lemma}\label{wconvexn}
Let $p,\ q,\ s\in(0,\infty)$ and $\omega\in M(\rp)$.
Then
$$\left[W\NHerzS\right]^{1/s}=W\NHerzScs$$
with the same quasi-norms.
\end{lemma}

With the help of this lemma, we now
turn to prove Theorem \ref{wczoogn}.

\begin{proof}[Proof of Theorem
\ref{wczoogn}]
Let all the symbols be as in the present
theorem. Then, using Theorem \ref{wczx}
and repeating an argument similar to
that used in the proof of Theorem
\ref{wczoog} via replacing $\mw$, $\Mw$,
\eqref{wczol1e00}, Theorem \ref{Th2},
and Lemmas \ref{convexl}, \ref{ee},
\ref{equa}, \ref{mbhag}, \ref{Atogl5},
\ref{wmbhg}, and \ref{wconvex} therein, respectively,
by $\mi(\omega)$, $\MI(\omega)$,
\eqref{wczoonl1e1}, Theorem \ref{Th2n},
and Lemmas \ref{convexln}, \ref{een},
\ref{equan}, \ref{mbhagn}, \ref{Atognl5},
\ref{wmbhgn}, and \ref{wconvexn}, we conclude
that $T$ is well defined on $\NHaSaH$ and,
for any $f\in\NHaSaH$,
$$
\|T(f)\|_{W\NHaSaH}\lesssim\|f\|_{
\NHaSaH},
$$
which completes the proof of
Theorem \ref{wczoogn}.
\end{proof}

\subsection{Real Interpolations}

In this subsection, we establish a relation between
inhomogeneous generalized Herz--Hardy
spaces and inhomogeneous
weak generalized Herz--Hardy spaces
via the real interpolation.
To be precise, we show that
the real interpolation spaces between inhomogeneous
generalized Herz--Hardy spaces and the Lebesgue
space $L^\infty(\rrn)$ are just the inhomogeneous weak
generalized Herz--Hardy spaces new introduced in
the present section. To begin with, recall
that the real interpolation
space between two quasi-Banach spaces
is defined as in Definition \ref{df681}.

We first show the following
interpolation theorem
about the inhomogeneous generalized
Herz--Hardy space $\NHaSaHo$ and
the Lebesgue space $L^{\infty}(\rrn)$.

\begin{theorem}\label{interln}
Let $p,\ q\in(0,\infty)$ and
$\omega\in M(\rp)$ satisfy $\mi(\omega)
\in(-\frac{n}{p},\infty)$.
Then, for any $\theta\in(0,1)$,
$$
\left(\NHaSaHo,\,L^{\infty}(\rrn)
\right)_{\theta,\infty}=
WH\Kmc_{\omega^{1-\theta},
\0bf}^{p/(1-\theta),q/(1-\theta)}(\rrn).
$$
\end{theorem}

\begin{proof}
Let all the symbols be as in the
present theorem and $\theta\in(0,1)$.
Then, repeating the
proof of Theorem \ref{interl} with
$\mw$, $\Mw$, Theorem \ref{Th3},
and Lemmas \ref{convexll}, \ref{mbhal},
\ref{vmbhl}, and \ref{wmbhl} therein replaced,
respectively, by $\mi(\omega)$, $\MI(\omega)$,
Theorem \ref{Th3n}, and Lemmas
\ref{convexlln}, \ref{mbhaln},
\ref{vmbhln}, and \ref{wmbhln}, we find
that
$$
\left(\NHaSaHo,\,L^{\infty}(\rrn)
\right)_{\theta,\infty}=
WH\Kmc_{\omega^{1-\theta},
\0bf}^{p/(1-\theta),q/(1-\theta)}(\rrn).
$$
This finishes the proof of Theorem
\ref{interln}.
\end{proof}

Finally, we establish the following
interpolation theorem which
shows that the real interpolation space
between $\NHaSaH$ and $L^{\infty}(\rrn)$
is just the inhomogeneous
weak generalized Herz--Hardy space
$WH\Kmc_{\omega^{1-\theta}}
^{p/(1-\theta),q/(1-\theta)}(\rrn)$.

\begin{theorem}\label{intergn}
Let $p,\ q\in(0,\infty)$ and $\omega\in M(\rp)$ satisfy
$$
-\frac{n}{p}<\mi(\omega)\leq\MI(\omega)<0.
$$
Then, for any $\theta\in(0,1)$,
$$
\left(\NHaSaH,\,L^{\infty}(\rrn)
\right)_{\theta,\infty}=
WH\Kmc_{\omega^{1-\theta}}^{p/
(1-\theta),q/(1-\theta)}(\rrn).
$$
\end{theorem}

\begin{proof}
Let all the symbols be as in
the present theorem. Then, from
Proposition \ref{intergx} and
an argument similar to that
used in the proof Theorem \ref{interg} via
replacing $\mw$, $\Mw$,
Theorems \ref{Th2} and \ref{ball}, and
Lemmas \ref{convexl}, \ref{ee}, \ref{equa},
\ref{mbhag}, \ref{vmbhg}, \ref{Atogl2}, and
\ref{wmbhg} therein, respectively, by
$\mi(\omega)$, $\MI(\omega)$,
Theorems \ref{Th2n} and \ref{balln}, and
Lemmas \ref{convexln}, \ref{een}, \ref{equan},
\ref{mbhagn}, \ref{vmbhgn}, \ref{Atognl2}, and
\ref{wmbhgn}, we infer that,
for any $\theta\in(0,1)$,
$$
\left(\NHaSaH,\,L^{\infty}(\rrn)
\right)_{\theta,\infty}=
WH\Kmc_{\omega^{1-\theta}}^{p/
(1-\theta),q/(1-\theta)}(\rrn).
$$
This then finishes the proof of Theorem
\ref{intergn}.
\end{proof}

\newpage

\markboth{\scriptsize\rm\sc Index}
{\scriptsize\rm\sc Index}
\printindex
\addcontentsline{toc}{chapter}{Index}

\newpage

\markboth{\scriptsize\rm\sc Abstract}{\scriptsize\rm\sc Abstract}
\addcontentsline{toc}{chapter}{Abstract}

\chapter*{Abstract}

This book is devoted to exploring properties of
generalized Herz spaces and establishing a complete
real-variable theory of Hardy spaces associated with
local and global generalized Herz spaces
via a totally fresh perspective which
means that the authors view these generalized Herz spaces as
special cases of ball quasi-Banach function spaces.
To be precise, in this book, the authors first study
some basic properties of generalized Herz spaces
and obtain boundedness and compactness characterizations
of commutators on them. Then
the authors introduce the associated
Herz--Hardy spaces, local Herz--Hardy spaces,
and weak Herz--Hardy spaces,
and develop a complete real-variable
theory of these Herz--Hardy spaces,
including their various
maximal function, atomic, finite
atomic, molecular as well as
various Littlewood--Paley function
characterizations. As applications,
the authors establish the boundedness
of some important operators arising
from harmonic analysis on these Herz--Hardy spaces.
Finally, the inhomogeneous Herz--Hardy spaces and their
complete real-variable theory are also investigated.
Due to the deficiency of the associate space of
the global Herz space,
the known real-variable
characterizations about Hardy-type spaces
associated with ball quasi-Banach function spaces are
not applicable to Hardy spaces associated with
global generalized Herz spaces which
need an improved generalization of
the existing one, done by
the authors also in this book and having more
additional anticipating applications.
The authors should also point out
that, with the fresh perspective
and the improved conclusions on the
real-variable theory of Hardy spaces associated
with ball quasi-Banach function spaces,
the exponents in all the obtained results of this book
are sharp. Moreover, all of these results in this book
are new and have never been published before.

\end{document}